%% file: birarakelov.tex
\def\DateTime{23/March/2019}
\def\Version{Version $1.0$}
\def\DateChapZero{2019_03_23} 
\def\DateChapOne{2019_03_23} 
\def\DateChapTwo{2019_03_23} 
\def\DateChapThree{2019_03_23} 
\def\DateChapFour{2019_03_23} 
\def\DateChapFive{2019_03_23} 
\def\DateChapSix{2019_03_23} 
\def\DateChapSeven{2019_03_23} 
\def\DateChapAppendix{2019_03_23} 
\def\FrontCoverColor{natsumushiiro} 
\def\IfChapVersion{\if01} 
\def\NoColorNoquery{\if01} 

\documentclass[a4paper,leqno,10pt,english]{smfbook}
\usepackage{smfthm,appendix,CJK}
\usepackage[T2A,T1]{fontenc}
\usepackage[english,francais]{babel}
\usepackage[utf8]{inputenc}
\usepackage{latexsym,amscd,color}
\usepackage{amsmath,amsfonts,amssymb,mathrsfs,mathtools}
\usepackage{qsymbols,euscript,enumitem}
\usepackage{bm,comment}
\usepackage{url,xspace}
\usepackage{verbatim}
\usepackage{tikz}
\usepackage{tikz-cd}
\usetikzlibrary{positioning,arrows,shapes}
\usepackage[normalem]{ulem}
\usepackage{manfnt} 
\usepackage{etoolbox}

\usepackage{hyperref}
\usepackage{xstring}

\makeindex

\theoremstyle{plain}
\author{Huayi CHEN}
\address{Institut de Math\'ematiques de Jussieu - Paris Rive Gauche, Universit\'e Paris Diderot}
\email{huayi.chen@imj-prg.fr}
\author{Atsushi MORIWAKI}
\address{University of Kyoto}
\email{moriwaki@math.kyoto-u.ac.jp}
\title{Arakelov geometry \\ over adelic curves}
\date{\DateTime, \Version}
\input xy
\xyoption{all}
\numberwithin{equation}{chapter}

\def\lbr{[\![}
\def\rbr{]\!]}
\def\indic{1\hspace{-2.5pt}\mathrm{l}}
\def\esssup_#1{\underset{#1}{\mathrm{ess\,sup\, }}}
\newcommand{\ndot}{\mathord{\cdot}}
\newcommand{\sndot}{\hspace{0.05em}\mathclap{\cdot}\hspace{0.05em}}
\newcommand{\emptyinnprod}{\langle\kern.15em,\kern-.02em\rangle} 
\newcommand{\rest}[2]{\left.{#1}\right\vert_{{#2}}}
\def\colorsout#1{\bgroup\markoverwith{\textcolor{#1}{\rule[0.5ex]{2pt}{0.7pt}}}\ULon} 
\def\coloruline#1{\bgroup\markoverwith{\textcolor{#1}{\rule[-0.5ex]{2pt}{0.7pt}}}\ULon} 

\definecolor{rose}{rgb}{1.0, 0.0, 0.5}
\definecolor{rossocorsa}{rgb}{0.83, 0.0, 0.0}
\definecolor{rosewood}{rgb}{0.4, 0.0, 0.04}
\definecolor{ruby}{rgb}{0.88, 0.07, 0.37}
\definecolor{orange-red}{rgb}{1.0, 0.27, 0.0}
\definecolor{mred}{rgb}{0.83, 0.0, 0.0}
\definecolor{mgreen}{rgb}{0.33, 0.42, 0.18}
\definecolor{purpleheart}{rgb}{0.41, 0.21, 0.61}
\definecolor{gray}{rgb}{0.5, 0.5, 0.5}
\definecolor{pastelgray}{rgb}{0.81, 0.81, 0.77}
\definecolor{natsumushiiro}{RGB}{206, 228, 174}
\def\ifgray{\if01}
\ifgray
\definecolor{rossocorsa}{rgb}{0.5, 0.5, 0.5}
\definecolor{mred}{rgb}{0.5, 0.5, 0.5}
\definecolor{blue}{rgb}{0.5, 0.5, 0.5}
\fi
\NoColorNoquery
\definecolor{mred}{rgb}{.0, .0, .0}
\definecolor{blue}{rgb}{.0, .0, .0}
\fi

\def\ifquery{\if00}
\ifquery
\def\lquery#1#2{\setlength\marginparwidth{#1}
\marginpar{\raggedright\fontsize{7.81}{9} 
\selectfont\upshape\hrule\smallskip 
#2\par\smallskip\hrule}}
\def\query#1{\lquery{65pt}{#1}}

\else
\def\query#1{}
\fi
\NoColorNoquery
\def\query#1{}
\fi
\IfChapVersion
\def\ChapVersion#1{\setlength\marginparwidth{100pt}
\marginpar{\raggedright\fontsize{7.81}{9} 
\selectfont\upshape\smallskip 
#1\par\smallskip}}
\fi

\DeclarePairedDelimiter{\norm}{\lVert}{\rVert}
\DeclareSymbolFont{bbold}{U}{bbold}{m}{n}
\DeclareMathSymbol{\bbalpha}{\mathord}{bbold}{"0B}
\DeclareMathSymbol{\bbbeta}{\mathord}{bbold}{"0C}
\DeclareMathSymbol{\bbgamma}{\mathord}{bbold}{"0D}
\DeclareMathSymbol{\bbdelta}{\mathord}{bbold}{"0E}
\DeclareMathSymbol{\bbespilon}{\mathord}{bbold}{"0F}
\DeclareMathSymbol{\bbzeta}{\mathord}{bbold}{"10}
\DeclareMathSymbol{\bbeta}{\mathord}{bbold}{"11}
\DeclareMathSymbol{\bbtheta}{\mathord}{bbold}{"12}
\DeclareMathSymbol{\bbiota}{\mathord}{bbold}{"13}
\DeclareMathSymbol{\bbkappa}{\mathord}{bbold}{"14}
\DeclareMathSymbol{\bblambda}{\mathord}{bbold}{"15}
\DeclareMathSymbol{\bbmu}{\mathord}{bbold}{"16}
\DeclareMathSymbol{\bbnu}{\mathord}{bbold}{"17}
\DeclareMathSymbol{\bbxi}{\mathord}{bbold}{"18}
\DeclareMathSymbol{\bbpi}{\mathord}{bbold}{"19}
\DeclareMathSymbol{\bbrho}{\mathord}{bbold}{"1A}
\DeclareMathSymbol{\bbsigma}{\mathord}{bbold}{"1B}
\DeclareMathSymbol{\bbtau}{\mathord}{bbold}{"1C}
\DeclareMathSymbol{\bbupsilon}{\mathord}{bbold}{"1D}
\DeclareMathSymbol{\bbphi}{\mathord}{bbold}{"1E}
\DeclareMathSymbol{\bbchi}{\mathord}{bbold}{"1F}
\DeclareMathSymbol{\bbpsi}{\mathord}{bbold}{"20}
\DeclareMathSymbol{\bbomega}{\mathord}{bbold}{"7F}

\def\upint{\mathchoice%
    {\mkern13mu\overline{\vphantom{\intop}\mkern7mu}\mkern-20mu}%
    {\mkern7mu\overline{\vphantom{\intop}\mkern7mu}\mkern-14mu}%
    {\mkern7mu\overline{\vphantom{\intop}\mkern7mu}\mkern-14mu}%
    {\mkern7mu\overline{\vphantom{\intop}\mkern7mu}\mkern-14mu}%
  \int}
\def\lowint{\mkern3mu\underline{\vphantom{\intop}\mkern7mu}\mkern-10mu\int}

\DeclareMathOperator{\card}{card}
\DeclareMathOperator{\Div}{Div}
\DeclareMathOperator{\dist}{dist}
\DeclareMathOperator{\defp}{dp}
\DeclareMathOperator{\dpur}{dpur}

\DeclareMathOperator{\Hom}{Hom}
\DeclareMathOperator{\Ker}{Ker}
\DeclareMathOperator{\ord}{ord}
\DeclareMathOperator{\Pic}{Pic}
\DeclareMathOperator{\pr}{pr}
\DeclareMathOperator{\Id}{Id}
\DeclareMathOperator{\sgn}{sgn}

\DeclareMathOperator{\hdeg}{\widehat{\deg}}

\DeclareMathOperator{\sposdeg}{\widehat{\deg}^{\kern.1em\lower.7ex\hbox{$\scriptstyle\mathrm{a}$}}_{+}}

\DeclareMathOperator{\smaxslop}{\widehat{\mu}^{\kern.1em\lower.5ex\hbox{$\scriptstyle\mathrm{a}$}}_{\max}}

\DeclareMathOperator{\supdeg}{\widehat{\deg}^{\kern.1em\lower.7ex\hbox{$\scriptstyle\mathrm{s}$}}_{\mathrm{up}}}

\DeclareMathOperator{\Spec}{Spec}
\DeclareMathOperator{\rang}{rk}
\def\sbullet{{\scriptscriptstyle\bullet}}
\DeclareMathOperator{\vol}{vol}

\def\scirc{\mathop{\raise.3ex\hbox{$\scriptscriptstyle\circ$}}}

\newcommand{\intervalle}[4]{\mathopen{#1}#2,#3\mathclose{#4}}

\begin{document}
\frontmatter 
\maketitle

\tableofcontents
\mainmatter%
\include{ch0_\DateChapZero}

\include{ch1_\DateChapOne}

\include{ch2_\DateChapTwo}

\include{ch3_\DateChapThree}

\include{ch4_\DateChapFour}

\include{ch5_\DateChapFive}
\include{ch6_\DateChapSix}

\include{ch7_\DateChapSeven}

\include{cha_\DateChapAppendix}

\backmatter
\bibliography{birarakelov}
\bibliographystyle{smfplain_nobysame}

\printindex

\end{document}

%% file: ch0_2019_03_23.tex

\addtocounter{chapter}{-1}
\chapter{Introduction}
\addtocontents{Introduction}

\IfChapVersion
\ChapVersion{Version of Chapter 0 : \\ \StrSubstitute{\DateChapZero}{_}{\_}}
\fi

The purpose of this book is to build up the fundament of an Arakelov theory over adelic curves in order to provide a unified framework for the researches of arithmetic geometry in several directions.

Let us begin with a brief description of the main ideas of Arakelov geometry. In number theory, it is well known that number fields are similar to fields of rational functions over algebraic curves defined over a base field, which is often assumed to be finite. It is expected that the geometry of schemes of finite type over $\mathbb Z$ should be similar to the algebraic geometry of schemes of finite type over a regular projective curve. However, some nice properties, especially finiteness of cohomological groups (in the case where the base field is finite), fail to hold in the arithmetic setting, which prevents using geometrical methods to count the arithmetic objects. The core problem is that schemes over $\Spec\mathbb Z$, even projective, are not ``compact'', and in general it is not possible to ``compactify'' them in the category of schemes. The seminal works of Arakelov \cite{Arakelov74,Arakelov75} propose to ``compactify'' a   scheme of finite type over $\Spec\mathbb Z$ by transcendental objects, such as the associated complex analytic variety, Hermitian metrics, Green functions, and differential forms etc.  In the case of relative dimension zero, the idea of Arakelov corresponds to the classic approach in algebraic number theory to include the infinite places of a number field to obtain a product formula, and to introduce Hermitian norms on projective modules over an algebraic integer ring to study the geometry of numbers and counting problems. Most interestingly, the approach of Arakelov proposes an intersection theory for divisors on a projective arithmetic surface (relative dimension one case), which is similar to the intersection pairing of Cartier divisors on classic projective surfaces.

The works of Arakelov have opened a gate to a new geometric theory of arithmetic varieties (schemes of finite type over $\Spec\mathbb Z$). Inspired by the classic algebraic geometry, many results have been obtained and enriched Arakelov's geometry.  Among the wide literature, we can mention for example the arithmetic Hodge index theorem by Faltings \cite{Faltings84} and Hriljac \cite{Hriljac85}, arithmetic intersection theory of higher dimensional arithmetic varieties and arithmetic Riemann-Roch theorem by Gillet and Soul\'e \cite{Gillet_Soule84,Gillet-Soule}, see also \cite{GRS08}. Arakelov geometry also provides an alternative approach (compared to the classic Weil height theory) to the height theory in arithmetic geometry, see \cite{Arakelov74,Arakelov75,Szprio85,Faltings91,BGS94},  (see also the approach of Philippon \cite{Philippon86,Philippon94}, and \cite{Soule91} for the comparison between Philippon height and Arakelov height). The Arakelov height is often more precise than the Weil height machine since the choice of a Hermitian metric on a line bundle permits to construct an explicit height function associated with that line bundle (in the Weil heigh machine, the height function is defined only up to a bounded function).

These advancements have led to fruitful applications in number theory, such as the proof of Mordell's conjecture by Faltings \cite{Faltings83,Faltings83b} and the alternative proof by Vojta \cite{Vojta91} (see also the proof of Bombieri \cite{Bombieri90} and the generalisation of Vojta's approach to the study of subvarieties in an Abelian variety \cite{Faltings91}), equidistribution of algebraic points in an arithmetic variety and applications to Bogomolov's conjecture \cite{Szpiro_Ullmo_Zhang,Ullmo98,Zhang98}, algebraicity of formal leaves of algebraic foliation \cite{Bost2001} etc. 

Although the philosophy of Arakelov allows to inspire notions and results of algebraic geometry and has already led to a rich arithmetic theory, the realisation of Arakelov theory is rather different from that of the classic algebraic geometry and usually gets involved subtle tools in analysis. The transition of technics on the two sides is often obscure. For example, the $abc$ conjecture, which can be easily established in the function field setting by algebraic geometry tools (see \cite{Mason84}), turns out to be very deep in the number field setting. Conversely, the Bogomolov's conjecture has been resolved in the number field setting, before the adaptation of its proof in the function field setting by using Berkovich analytic spaces (see \cite{Gubler07,Yamaki16,Yamaki17}). It is therefore an interesting problem to provide a uniform fundament for Arakelov geometry, both in the function field and number field settings, and the adelic approach is a natural choice for this goal. We would like however to mention that  Durov \cite{Durov07} has proposed an approach of different nature to algebrify the Arakelov geometry over number fields.

The theory of ad\`eles in the study of global fields was firstly introduced by Chevalley \cite[Chapitre III]{Chevalley51} for function fields and by Weil \cite{Weil51} for number fields. This theory allows to consider all places of a global filed in a unified way. It also leads to a uniform approach in the geometry of numbers in global fields, either via the adelic version of Minkowski's theorems and Siegel's lemma developed by McFeat \cite{McFeat71}, Bombieri-Vaaler \cite{Bombieri-Vaaler83}, Thunder \cite{Thunder96}, Roy-Thunder \cite{Roy_Thunder96}, or via the study of adelic vector bundles developed by Gaudron \cite{Gaudron08}, generalising the slope theory introduced by Bost \cite{BostBour96,Bost2001}. 

Several works have been realised in the adelification of Arakelov theory. Besides the result of Gaudron on adelic vector bundles over global fields mentioned above, we can for example refer to \cite{Zhang95a} for adelic metrics on arithmetic line bundles and applications to the Bogomolov problem for cycles. Moreover, Moriwaki \cite{Moriwaki16} has studied the birational geometry of adelic line bundles over arithmetic varieties. The key point is to consider an arithmetic variety as a scheme of finite type over a global field, together with a family of analytic varieties (possibly equipped with metrised vector bundles) associated with the scheme, which is parametrised by the set of all places of the global field. Classic objects in Arakelov geometry can be naturally considered in this setting. For example, given a Hermitian line bundle over a classic arithmetic variety (scheme of finite type over $\Spec\mathbb Z$), the algebraic structure of the line bundle actually induces, for each finite place of $\mathbb Q$, a metric on the pull-back of the line bundle on the corresponding analytic space.   

In this book, we introduce the notion of \emph{adelic curves} and develop an Arakelov theory over them. By adelic curve we mean a field equipped with a family of absolute values parametrised by a measure space, such that the logarithmic absolute value of each non-zero element of the filed is an integrable function on the measure space, with $0$ as its integral. This property is called \emph{product formula}. Note that this notion has been studied by Gubler \cite{Gubler97} in the setting of height theory and is also considered by Ben Yaakov and Hrushovski \cite{Yaakov_Hrushovski,Hrushovski16}
in a recent work on model theory of global fields.  Clearly the notion of adelic curve generalises the classic one of global field, where the measure space is given by the set of all places of the global field equipped with the discrete measure of local degrees. However, this is certainly not the only motivation for the general notion of adelic curves. Our choice is rather inspired by several bunches of researches which are apparently  transversal to each other, which we will resume as follows (we will explain further the reason for the choice of terminology ``adelic curve''). 

\emph{1. Finitely generated extensions of a number field.} From a point of view of birational geometry, we expect that the field of rational functions of an algebraic variety determines the geometric properties of the variety. In Arakelov geometry, we consider integral schemes of finite type over $\Spec\mathbb Q$, whose function field is a finitely generated extension of $\mathbb Q$. Moriwaki \cite{Moriwaki00} has developed an Arakelov height theory for varieties over a finitely generated extension of a number field and applied it to the study of Bogomolov problem over such a field (see \cite{Moriwaki01}, see also  \cite{Moriwaki04} for a panoramic view). Burgos, Philippon and Sombra \cite{BPS14} have expressed the height of cycles in a projective variety over a finitely generated extension of $\mathbb Q$ as an integral of local heights over the set of places of the field.

\emph{2. Trivially valued field.} In number theory, we usually consider non-trivial absolute values on fields. Note that on any field there exists a trivial absolute value which takes value $1$ on each non-zero element of the field. Note that a trivial product formula is satisfied in this setting. Although the trivially valued fields are very simple, the corresponding geometry of numbers is rather rich, which has wide interactions with the classic geometry of lattices or Hermitian vector bundles. In fact, given a finite-dimensional vector space over a trivially valued field, the  ultrametric norms on it are canonically in bijection to the decreasing $\mathbb R$-filtrations on the vector space. The $\mathbb R$-filtration is a key method of the works \cite{Chen10b,Chen08,Chen10}, where the main idea consists in associating to each Hermitian vector bundle an $\mathbb R$-filtration on the generic fibre, which captures the arithmetic information such as successive minima or successive slopes.  

\emph{3. Harder-Narasimhan theory for vector bundles on higher dimensional varieties.} Harder and Narasimhan theory \cite{Harder-Nara} is an important tool in the study of vector bundles on a projective curve. In the geometry of Euclidean lattices, the counterpart of Harder-Narasimhan theory has been proposed by Stuhler \cite{Stuhler76} and Grayson \cite{Grayson84}. Later Bost \cite{BostBour96} has generalised their work in the setting of Hermitian vector bundles on the spectrum of the ring of algebraic integers in a number field. Moreover, he has developed the slope inequalities in this framework and applied them to the study of algebraicity of formal schemes \cite{Bost2001,Bost2004,BostICM}. Note that the slope function and the notion of semistability can be naturally defined for torsion-free coherent sheaves on a polarised projective variety \cite{Takemoto72}. This allows Shatz \cite{Shatz77} and Maruyama \cite{Maruyama81} to develop a Harder-Narasimhan theory for general torsion-free coherent sheaves. However, it seems that the analogue of their results in the arithmetic case is still missing.

\emph{4. Fields of algebraic numbers, Siegel fields.} The geometry of numbers for algebraic (not necessarily finite) extensions plays an important role both in Diophantine problems and in Arakelov geometry. Recall that the Minkowski's theorem and Siegel's lemma in geometry of numbers admit an adelic version for number fields, see \cite{Bombieri-Vaaler83,McFeat71}. They also have an absolute counterpart over $\overline{\mathbb Q}$, see \cite{Roy_Thunder96,Roy_Thunder99,Zhang95}. In Arakelov geometry, a notion of Hermitian vector bundle over $\overline{\mathbb Q}$ has been proposed in the work \cite{Bost_Chen} of Bost and Chen, on which the absolute Siegel's lemma applies and is useful in the study of tensorial semistability of classic Hermitian vector bundles. Similarly, in the approach of Gaudron and R\'emond \cite{Gaudron_Remond13} to the tensorial semistability, the absolute Siegel's lemma is also a key argument. In \cite{Gaudron_Remond14}, the notion of Siegel field has been proposed. A Siegel field is a subfield of $\overline{\mathbb Q}$ on which an analogue of Siegel's lemma is true. In order to formulate a geometry of numbers for a Siegel field, Gaudron and R\'emond have introduced a topology on the space of all places of such a field and a Borel measure on it.

\emph{5. Algebraic extensions of function fields.} 
In \cite{MR1491742}, Corvaja and Zannier have studied the arithmetic of algebraic extensions of function fields. They have characterised the infinite algebraic extensions of function fields of a curve which still satisfy a product formula. They have also discussed several examples of product formulas associated with algebraic surfaces in revealing the non-uniqueness of the extension of a product formula under finite field extensions.

The above results are obtained in various settings of arithmetic geometry. It turns out that these settings can be naturally included in the framework of adelic curves (see \S\ref{Sec: Examples of adelic curves} for details) in order to treat the geometry of various fields analogously to that of vector bundles on projective curves. For example, on the field $\mathbb Q(T)$ of rational functions with coefficients in $\mathbb Q$, three types of absolute values are defined (see \S\ref{subsec:fun:field:Q} for details): the valuation corresponding to closed points of $\mathbb P^1_{\mathbb Q}$, the natural extensions of $p$-adic absolute values, the Archimedean absolute value corresponding to divers embeddings of $\mathbb Q(T)$ in $\mathbb C$. Note that Jensen's formula for Mahler measure shows that these absolute values, once suitably parametrised by a measure space, satisfies a product formula. Thus we can consider it as an adelic curve. This is actually a particular case of polarised arithmetic projective varieties, where the polarisation provides a structure of adelic curve on the field of rational functions on the projective variety. Moreover, algebraic coverings of an adelic curve can be naturally constructed (see Section \ref{Sec:algebracexte}), which provides a framework for the study of the arithmetic of algebraic extensions.

Note that Gubler \cite{Gubler97} has introduced a similar notion of $M$-field and extended the Arakelov height theory to this setting. An $M$-field is a field equipped with a measure space and a family of functions parametrised by the measure space which are absolute values almost everywhere, and the height of an arithmetic variety is defined as the integration along the measure space of local heights. However, our main  concern is to build up a suitable geometry of numbers while the purpose of \cite{Gubler97} is to extend the Arakelov height theory in  a sufficiently general setting in order to include the theory of Nevanlinna. In Diophantine geometry the geometry of numbers is as important as a the height theory, particularly in the geometrisation of the method of ``auxiliary polynomials''. We propose the notion of \emph{adelic vector bundles} on adelic curves, which consist of a finite-dimensional vector space over the underlying field, equipped with a measurable family of norms parametrised by the measure space. Our  choice facilites the study of algebraic constructions of adelic vector bundles. The height of arithmetic varieties is described in a global way by the asymptotic behaviour of graded linear series equipped with structures of adelic vector bundles, rather than the integral of local heights.

In the framework of model theory, Ben Yaakov and Hrushovski \cite{Yaakov_Hrushovski,Hrushovski16} 
also consider the formalisme of a field equipped with a family of absolute value parametrised by a measure space, which satisfied a product formula (called \emph{globally valued field} in their terminologies). Their work permits to considered classic Diophantine geometry objects (in particular heights) in the model theory setting. 

In order to set up a theory of adelic vector bundles over adelic curves, we present in the first chapter various constructions and properties of seminormed vector spaces over a complete valued field. Although the constructions and results are basic, the subtleties in the interaction and the compatibility of divers algebraic constructions, such as restriction, quotient, dual, tensor product, exterior powers etc, have not been clarified in the literature in a systematic way. In particular, several classic results in the functional analysis over $\mathbb C$ are no longer true in the non-Archimedean setting. We choose carefully our approach of presentation to unify the treatment of non-Archimedean and Archimedean cases whenever possible, and specify the differences and highlight the subtleties in detail.   A particular attention is paid to the two constructions of tensor product seminorms: the $\pi$-tensor product and the $\varepsilon$-tensor product. These notions have been firstly introduced by Grothendieck \cite{Grothendieck53,Grothendieck55} in the setting of functional analysis over $\mathbb C$. It turns out that similar constructions can be defined more generally over an arbitrary complete valued field, and they are useful for example in the study of seminorms on exterior powers.

The orthogonality is another theme discussed in the first chapter. Classically the orthogonality is a natural notion in the study of inner product spaces. We consider an equivalent form of this notion, which can be defined in the setting of finitely generated seminormed vector spaces over an arbitrary complete valued field. This reformulation has been used in \cite{BMPS16} to study the arithmetic positivity on toric varieties. Here it will serve as a fundamental tool to study ultrametrically normed spaces, inner product spaces and the construction of orthogonal tensor products. In particular, an analogue of the Gram-Schmidt process holds for finite-dimensional ultrametrically seminormed spaces, which plays a key role in the compatibility of the determinant norm with respect to short exact sequences.

We also discuss extension of seminorms under a valued extension of scalars. We distinguish three extensions of seminorms, corresponding to the three types of tensor product. The compatibility of extension of scalars with respect to divers algebraic construction is also explained. These constructions are used in the pull-back of an adelic vector bundle by an algebraic covering of the adelic curve. 

Note that in the classic Arakelov theory, usually we consider a vector space over a global field equipped with a family of norms. However, from the point of view of birational geometry, it is natural to consider metrics which admits singularity, that is, degenerates on a closed subscheme (which is usually the base locus of a linear series) to a family of seminorms. Moreover, in the study of algebraicity of formal leaves of an arithmetic foliation, the canonical ``metrics'' on the tangent bundle are often seminorms. We refer the readers to \cite{Bost_Chambert09} for more details. Motivated by these observations, we choose to present a panoramic view on the tools about general seminormed vector spaces which could be useful in Arakelov geometry later.

The second chapter is devoted to a presentation of metrised line bundles on a projective scheme over a complete valued field. It could be considered as a higher dimensional version of the results presented in Chapter 1. We use Berkovich topology to define continuous metrics on a vector bundle. Note that in the case where the base field is $\mathbb C$, our definition coincides with the classic definition of continuous metric on a vector bundle over a complex analytic space (associated with a complex projective scheme).

The Fubini-Study metric is another important ingredient of Chapter 2. It is closely related to the positivity of metrics on line bundles. More precisely, a continuous metric on a line bundle over a projective scheme defined over a complete valued field is said to be \emph{semipositive} if it can be written as a uniform limit of Fubini-Study metrics. In the case where the absolute value is Archimedean, this definition is equivalent to the semipositivity of the curvature current of the metric. In the case where the absolute value is non-Archimedean and non-trivial, it is equivalent to the semipositivity condition proposed in \cite[Section 6.8]{Chamber_Ducros12} and \cite[Section 6]{Gubler_Kunnemann15}. However, in the trivial valuation case, it seems that our formulation is crucial to study the positivity of the metrics. 

In classic Hermitian geometry, the positivity is closely related to the extension of sections of an ample line bundle with a control on the supremum norms. We establish a non-Archimedean analogue of the extension property, generalising the main result of \cite{Chen_Moriwaki2017} to the non-necessarily reduced case.

The third chapter is devoted to the fundament of adelic curves. We first give the formal definition of this notion and illustrate by various examples. The algebraic coverings of adelic curves occupy an important part of the chapter. As mentioned above, an adelic curve is a field equipped with a family of absolute values parametrised by a measure space, which satisfies a product formula. Given an algebraic extension of the underlying field, there is a canonical family of absolute values parametrised by a measure space fibered on the initial measure space and equipped with a disintegration kernel.  This construction is important in the height theory for algebraic points and in the study of Siegel and Northcott properties.  Contrary to the approach of \cite{Gaudron_Remond14}, we do not assume the structural measurable space of an adelic curve to be a topological space and do not adopt the topological construction of algebraic coverings. Although it is possible to reduce the construction to the case of finite extensions by an argument of passage to projective limit, even for the simplest case of finite separable extension of the underlying field, the problem is highly non-trivial. The main subtleties come from the measurability of the fibre integral, which neither follows from the classic disintegration theory, nor from the property of extension of absolute values in algebraic number theory. The difficulty is resolved by using symmetric polynomials and Vandermonde matrix.

The analogue in the adelic curve setting of the geometry of numbers occupies the main part of the fourth chapter. Given an adelic curve, for any finite-dimensional vector space over the underlying field, we consider families of norms indexed by the structural measure space of the adelic curve. Natural measurability and dominancy conditions are defined for such norm families. 
An adelic vector bundle is a finite-dimensional vector space over the underlying field of the adelic curve, equipped with a measurable and dominated norm family. In the case where the adelic curve arises from a global field, this notion corresponds essentially the notion of adelic vector bundle in the work \cite{Gaudron08}. Note that in the classic global field case it is required that almost all norms in the structure of an adelic vector bundle come from a common integral model of the vector space. However, in our general setting of adelic curve, it is not adequate to discuss integral models since the integral ring in an adelic curves is not well defined. The condition of common integral model is replaced by the dominancy condition, which in the global field case can be considered as uniform limit of classic structure of adelic vector bundle.

The arithmetic invariants of adelic vector bundles are also discussed. For example, the Arakelov degree of an adelic vector bundle is defined as the integral of the logarithmic determinant norm of a non-zero maximal exterior power vector, similarly as in the classic case of Hermitian vector bundle over an arithmetic curve. Moreover, although the analogue of classic minima of lattices can not be reformulated in the adelic curve setting, due to the lack of integral models, the version of Roy and Thunder \cite{Roy_Thunder96}, which is based on the height function (or equivalently the Arakelov degree of the non-zero vectors), can be naturally generalised in our setting of adelic curves. However, it turns out that several fundamental results in geometry of numbers, such as Minkowski's theorems, are not true in the general setting, and the set of vectors in the adelic unit ball is not the good generalisation of lattice points of norm $\leqslant 1$.  This phenomenon suggests that the slope method of Bost \cite{BostBour96} might be more efficient in Diophantine geometry. In fact, inspired by the Harder-Narasimhan theory of vector bundles over curves, the notion of successive slopes has been proposed in \cite{Stuhler76,Grayson84} for Euclidean lattices and generalised in \cite{BostBour96,Bost2001} with applications to the period and isogenies of abelian varieties, and algebraicity of formal schemes. In the setting of adelic vector bundles on adelic curves, we build up an analogue of Harder-Narasimhan theory and the slope method. In this sense, adelic vector bundles on adelic curves have very similar properties as those of vector bundles on a regular projective curve, or Hermitian vector bundles over an arithmetic curve. It is for this reason that we have chosen the terminology of adelic curve. However, although the successive minima and the successive slopes are close in the number field case (see \cite{Borek05,chen17}), they can differ much in the general adelic curve setting, even for the simple case of a field equipped with several copies of the trivial absolute value. Note that the semistability of adelic vector bundles over such adelic curves plays an important role in Diophantine geometry of projective spaces, as for example in the work of Faltings and W\"ustholz \cite{Fal-Wusth} (although not written explicitly in the language of the slope theory). Our general setting of adelic vector bundles helps to understand the roles of different arithmetic invariants should play in a Diophantine argument. 

The adelic curve consisting of the trivial absolute value is also closely related to the geometric invariant theory. In the fifth chapter of the book, we explain this link and apply it to the estimation of the minimal slope of the tensor product of two adelic vector bundles.  In fact, an ultrametrically normed vector space over a trivially valued field can be considered as a decreasing $\mathbb R$-filtration of the vector space. In the geometric invariant theory, an action of the multiplicative group on a finite-dimensional vector space over a field corresponds to the decomposition of the vector space into the direct sum of eigensubspaces and thus determines an $\mathbb R$-filtration of the vector space by the eigenvalues. Therefore we can reformulate the Hilbert-Mumford criterion for general linear groups (or products of general linear groups) in terms of a slope inequality for adelic vector bundles on the adelic curve of one trivial absolute value.

Bogomolov (see \cite{Raynaud81}) has interpreted the semistability of a vector bundle over a projective curve as an inequality linking the $\mathbb R$-filtration and the Arakelov degree. This result can also be viewed as a link between the geometric invariant theory and the semistability in the theory of Harder-Narasimhan. Later Ramanan and Ramanathan \cite{Ramanan_Ramanathan} have given an algebraic proof of the semistability of the tensor product of two semistable vector bundles on a regular projective curve over a field of characteristic $0$. In the number field case, Bost \cite{Bost97} has conjectured that the arithmetic analogue of the tensorial semistability is also true. This conjecture is equivalent to the statement that the tensor product of two Hermitian vector bundles has a minimal slope which is bounded from below by the minimal slopes of the two Hermitian vector bundles.

In the setting of adelic vector bundles over adelic curves, we can consider the natural generalisation of Bost's conjecture stating that, if the underlying base field of the adelic curve is perfect, then the tensor product of two semistable \emph{Hermitian} adelic vector bundles is also semistable. Besides the function field case proved by Ramanan and Ramanathan, the generalised conjecture is also true in the case where the adelic curve is given by a perfect field equipped with a finite number of copies of the trivial absolute value (see \cite{Totaro96}). We prove here a weaker version of this conjecture, showing that the minimal slope of the tensor product of two (non-necessarily Hermitian)  adelic vector bundles is bounded from below by the sum of the minimal slopes of the two adelic vector bundles, minus three half of the logarithm of the rank of the tensor product bundle times the measure of Archimedean places. In particular, the conjecture is true if the base field is perfect and all absolute values in the adelic curve structure are non-Archimedean. This result is similar to the works \cite{Chen_pm,Andre11,Gaudron_Remond13,Bost_Chen} in the case where the adelic curve comes from a number field. However, the strategy of proof is different. In fact, the common point of the works cited above is a geometric version of Siegel's lemma proved by Zhang \cite{Zhang95}, which could be considered as an absolute version of Minkowski's second theorem, which is false for general adelic curves. Our method relies on the geometric invariant theory of grassmannian (with Pl\"ucker coordinates) and combines the technics of \cite{Chen_pm} and \cite{Bost_Chen}.

The sixth chapter is devoted to the study of metrised line bundles on arithmetic varieties over adelic curves. In the classic setting of adelic metrics such as \cite{Zhang95a,Moriwaki16}, it was required that an adelic metric should coincides with an integral model metric for all but finitely many places. Again the integral model metric is not adequate in our setting of adelic curves, the suitable notions of dominancy and measurability occupy thus an important part of the chapter. An adelic line bundle on a projective  variety is then defined to be an invertible sheaf equipped with a dominated and measurable family of metrics parametrised by the adelic curve. In the setting of global fields, our definition is slightly more general than the classic one, which includes the limits of classic adelic line bundles. The analogue of some classic geometric invariants, such as  height function, essential minimum, and arithmetic volume function is also discussed. In particular, in the definition of the arithmetic volume function, we use the positive degree instead of the logarithmic cardinal of the small sections since the latter is no longer adequate in the general setting. Note that the failure of Minkowski's first theorem brings several technical difficulties, notably the filtration by minima and the filtration by slopes do not lead to the same arithmetic invariants, on the contrary of the case of number fields as in \cite{chen17}. Our strategy consists in introducing a refinement of the method of arithmetic Newton-Okounkov bodies \cite{Boucksom_Chen},
which allows to treat the case of graded linear series equipped with filtrations which are not necessarily additive. 

In the seventh and the last chapter, we relate the asymptotic minimal slope to the absolute minimum of the height function of an adelic line bundle, which could be considered as a generalisation of Nakai-Moishezon's criterion in the setting of Arakelov geometry over an adelic curve. In the case where the analogue of a strong version of Minkowski's first theorem holds for the adelic curve, we deduce from the criterion an analogue of Siegel's lemma for adelic vector bundles on the adelic curve. Our work clarifies the arguments of geometric nature from several fundamental result in the classic geometry of numbers.

Limited by the volume of the monograph, many aspects are not included in the current text. First of all, an arithmetic intersection theory should be developed in the setting of Arakelov geometry over an adelic curve, which allows to interpret the height of arithmetic varieties as the arithmetic intersection number. Secondly, by using the adelic curve of several copies of the trivial absolute value, we expect to incorporate the conditions and results of geometric invariant theory into the arithmetic setting. Thirdly, the geometry of adelic vector bundles should lead to a Diophantine approximation theory of adelic curves. Finally, the fundamental works achieved in the monograph could be applied to the study of Nevanlinna theory of $M$-field proposed by Gubler.

{\bf Acknowledgements.} We are grateful to Jean-Beno{\^{i}}t Bost, Carlo Gasbarri  and Hugues Randriambololona for comments, and to Ita\"{i} 
Ben Yaacov and Ehud Hrushovski for having sent us their lecture notes and for letter communications.

%% file: ch1_2019_03_23.tex

\chapter{Metrized vector bundles: local theory}\label{chap:local:theory}

\IfChapVersion
\ChapVersion{Version of Chapter 1 : \\ \StrSubstitute{\DateChapOne}{_}{\_}}
\fi

The purpose of this chapter is to explain the constructions and properties of normed vector spaces over a complete valued field.   It will serve as the fundament for the global study of adelic vector bundles. Note that we need to consider both Archimedean and non-Archimedean cases. Hence we carefully choose the approach of presentation  to unify the statements  whenever possible, and to clarify the differences.

Throughout the chapter, let $k$ be a field equipped with an absolute value $|\ndot|$. We assume that $k$ is complete with respect to the topology induced by $|\ndot|$. We emphasise that $|\ndot|$ could be the trivial absolute value on $k$, namely $|a|=1$ for any $a\in k\setminus\{0\}$. If the absolute value $|\ndot|$ is Archimedean, then $k$ is either the field $\mathbb R$ of real numbers or the field $\mathbb C$ of complex numbers. For simplicity, we assume that $|\ndot|$ is the usual absolute value on $\mathbb R$ or $\mathbb C$ if it is Archimedean.

\section{Norms and seminorms}\label{Sec: Norms}

\begin{defi}\label{Def:norm}
Let $V$ be a vector space over $k$. 
A map $\norm{\ndot}:V\rightarrow\mathbb R_{\geqslant 0}$ is called a \emph{seminorm}\index{seminorm} on $V$ if the following conditions \ref{Item: positive homogeneous} and \ref{Item: triangle inequality} are satisfied:
\begin{enumerate}[label=(\alph*)]
\item\label{Item: positive homogeneous} for any $a\in k$ and any $x\in V$, one has $\|ax\|=|a|\cdot\|x\|$;
\item\label{Item: triangle inequality} the \emph{triangle inequality}\index{triangle inequality}: for {any  $(x,y)\in V\times V$}, one has
$\|x+y\|\leqslant\|x\|+\|y\|$.
\end{enumerate}
The couple $(V,\norm{\ndot})$ is called a \emph{seminormed vector space}\index{seminormed vector space} over $k$.
If in addition the following \emph{strong triangle inequality}\index{strong triangle inequality} is satisfied
\begin{equation*}
{\forall\,(x,y)\in V^2},\quad\|x+y\|\leqslant\max\{\|x\|,\|y\|\},
\end{equation*}
we say that the seminorm $\norm{\ndot}$ is \emph{ultrametric}\index{ultrametric}\index{seminorm!ultrametric ---}. 
Note that the existence of a non-identically vanishing ultrametric seminorm on $V$ implies that the absolute value $|\ndot|$ on $k$ is non-Archimedean. Furthermore, if the following additional condition \ref{Item: seminorm is a norm} is satisfied:
\begin{enumerate}[label=(\alph*)]
\setcounter{enumi}{2}
\item\label{Item: seminorm is a norm} for any $x\in V\setminus\{0\}$, one has $\|x\|>0$,
\end{enumerate}
the seminorm $\norm{\ndot}$ is called a \emph{norm}\index{norm} on $V$, and
the couple $(V,\norm{\ndot})$ is called a \emph{normed vector space}\index{normed vector space} over $k$.

If $(V,\norm{\ndot})$ is a seminormed vector space over $k$, then \[N_{\norm{\ndot}} := \{ x \in V \,:\, \| x \| = 0 \}\] is a vector subspace of $V$, called the \emph{null space}\index{null space} of $\norm{\ndot}$. Moreover, if we denote by $\pi : V \to V/N_{\norm{\ndot}}$ the linear map
of projection, then there is a unique norm  $\norm{\ndot}^{\sim}$
on $V/N_{\norm{\ndot}}$ such that $\norm{\ndot} = \norm{\ndot}^{\sim}\circ{\pi}$. The norm $\norm{\ndot}^{\sim}$ is called the \emph{norm associated with the seminorm $\norm{\ndot}$}\index{norm!--- associated with a seminorm}.
\end{defi}

\begin{defi}\label{Def:restriction}
Let $f : W \to V$ be a linear map 
of vector spaces over $k$ and $\norm{\ndot}$ be a seminorm on $V$.
We define $\norm{\ndot}_{f}:W\rightarrow\mathbb R_{\geqslant 0}$ to be
\[
\forall\,x\in W,\quad\|x\|_{f} := \|f(x)\|,
\]
which is a seminorm on $W$, called the seminorm \emph{induced} by $f$ and $\norm{\ndot}$\index{seminorm!--- induced by a map}.
Clearly, if $\norm{\ndot}$ is ultrametric, then also is $\norm{\ndot}_{f}$.
In the case where $f$ is injective, $\norm{\ndot}_{f}$ is often denoted by $\norm{\ndot}_{W \hookrightarrow V}$ and
is called the seminorm on $W$ \emph{induced} by $\norm{\ndot}$\index{seminorm!induced ---}\index{induced seminorm}, or the \emph{restriction} of $\norm{\ndot}$ {to} $W$\index{restriction of a seminorm}\index{seminorm!restriction}.
\end{defi}

\begin{enonce}[remark]{Notation}\label{Not:ball}
Let $(V,\norm{\ndot})$ be a seminormed vector space over $k$. If $\epsilon$ is a non-negative real number, we denote by $(V,\norm{\ndot})_{\leqslant\epsilon}$ or simply by $V_{\leqslant\varepsilon}$ the closed ball $\{x\in V\,:\,\|x\|\leqslant\epsilon\}$ of radius $\epsilon$ centered at the origin. Similarly, we denote by $(V,\norm{\ndot})_{<\epsilon}$ or by $V_{<\epsilon}$ the open ball $\{x\in V\,:\,\|x\|<\epsilon\}$.
\end{enonce}

\begin{prop}\label{Pro:dilatation} Assume that $|\ndot|$ is non-trivial. Let $\lambda\in \intervalle{]}{0}{1}{[}$ such that 
\[\lambda<\sup\{|a|\,:\,a\in k^{\times},\,|a|<1\}.\]
Let $(V,\norm{\ndot})$ be a seminormed vector space over $k$ and $x$ be a vector in $V$ such that $\|x\|>0$. There exists $b\in k^{\times}$ such that $\lambda\leqslant\|bx\|< 1$.
\end{prop}
\begin{proof}
Let $a$ be an element in $k^{\times}$ such that $\lambda<|a|<1$. We take $b=a^p$ with 
\[p=\left\lfloor\frac{\ln(\lambda)-\ln\|x\|}{\ln|a|}\right\rfloor.\]
By definition one has $p\leqslant (\ln(\lambda)-\ln\|x\|)/\ln|a|$. Hence
$|b|=|a|^p\geqslant \lambda/\|x\|$,
which leads to $\|bx\|=|b|\cdot\|x\|\geqslant\lambda$. Moreover, since $\lambda<|a|<1$ one has $\ln(\lambda)<\ln|a|<0$. Hence $\ln(\lambda)/\ln|a|>1$, which implies that $p>-\ln\|x\|/\ln|a|$. Hence
$|b|=|a|^p<\|x\|^{-1}$,
which leads to $\|bx\|<1$.
\end{proof}

\begin{prop}
\label{Pro:valeur}
Let $(V,\norm{\ndot})$ be an ultrametrically seminormed vector space over~$k$. 
\begin{enumerate}[label=\rm(\arabic*)]
\item\label{Item: sum of vectors of distinct norm} If $x_1,\ldots,x_n$ are vectors of $V$ such that the numbers $\|x_1\|,\ldots,\|x_n\|$ are distinct, then one has
$\|x_1+\cdots+x_n\|=\displaystyle\max_{i\in\{1,\ldots,n\}}\|x_i\|$.
\item\label{Item: quotient of norm values} The cardinal of the image of the composed map
\begin{equation}\label{equ:composedmap}\xymatrix{\relax V\setminus N_{\norm{\ndot}}\ar[r]^-{\norm{\ndot}}& \mathbb R_{>0}
\ar[r]& \mathbb R_{>0}
/|k^{\times}|}
\end{equation}
is not greater than the rank of $V/N_{\norm{\ndot}}$ over $k$, where $\mathbb R_{>0}$ 
denotes the multiplicative group of positive real numbers, and $|k^{\times}|$ is the image of $k^{\times}$ by $|\ndot|$.
\end{enumerate}
\end{prop}
\begin{proof}
\ref{Item: sum of vectors of distinct norm} The statement is trivial when $n=1$. Moreover, by induction it suffices to treat the case where $n=2$. Without loss of generality, we assume that $\|x_1\|<\|x_2\|$. Since $\norm{\ndot}$ is ultrametric, one has $\|x_1+x_2\|\leqslant\max\{\|x_1\|,\|x_2\|\}=\|x_2\|$. Moreover, 
\[\|x_2\|=\|x_1+x_2+(-x_1)\|\leqslant\max\{\|x_1+x_2\|,\|x_1\|\}.\]
Since $\|x_2\|>\|x_1\|$, one should have $\|x_2\|\leqslant\|x_1+x_2\|$. Therefore \[\|x_1+x_2\|=\|x_2\|=\max\{\|x_1\|,\|x_2\|\}.\]

\ref{Item: quotient of norm values} By replacing $V$ by $V/N_{\norm{\ndot}}$ and $\norm{\ndot}$ by the associated norm, we may assume that $\norm{\ndot}$ is actually a norm. Denote by $I$ the image of the composed map \eqref{equ:composedmap}. For each element $\alpha$ in $I$, we pick a vector $x_\alpha$ in $V\setminus\{0\}$ such that the image of $x_\alpha$ by the composed map is $\alpha$. We will show that the family $\{x_\alpha\}_{\alpha\in I}$ is linearly independent over $k$ and hence the cardinal of $I$ is not greater than the rank of $V$ over $k$. Assume that $\alpha_1,\ldots,\alpha_n$ are distinct elements of the set $I$ and $\lambda_1,\ldots,\lambda_n$ are non-zero elements of $k$. Then the values $\|\lambda_1x_{\alpha_1}\|,\ldots,\|\lambda_nx_{\alpha_n}\|$ are distinct. As the norm $\norm{\ndot}$ is ultrametric, by \ref{Item: sum of vectors of distinct norm} one has
\[\|\lambda_1x_{\alpha_1}+\cdots+\lambda_nx_{\alpha_n}\|=\max_{i\in\{1,\ldots,n\}}\|\lambda_ix_i\|>0.\]
Hence $\lambda_1x_{\alpha_1}+\cdots+\lambda_nx_{\alpha_n}$ is non-zero.
\end{proof}

\begin{coro}\label{Cor:finitevlaue} Let $(V,\norm{\ndot})$ be an ultrametrically seminormed vector space of finite rank over $k$.
Then we have the following:
\begin{enumerate}[label=\rm(\arabic*)]
\item\label{Item: value of discrete valuation}
If $|\ndot|$ is a discrete valuation (namely $|k^{\times}|$ is a discrete subgroup of $\mathbb R_{>0}$),  
then the image of $V\setminus N_{\norm{\ndot}}$ by $\norm{\ndot}$ is a discrete subset of $\mathbb R_{>0}$.  
\item\label{Item: norm value of trivial valuation}
If $|\ndot|$ is the trivial absolute value, then the image of $V$ by $\norm{\ndot}$ is a finite set, whose cardinal does not exceed $\dim_k(V/N_{\norm{\ndot}}) + 1$.
\end{enumerate}
\end{coro}

\subsection{Topology}\label{Subsec: topology}
Let $(V,\norm{\ndot})$ be a seminormed vector space over $k$. The seminorm $\norm{\ndot}$ induces a pseudometric $\dist(\ndot,\ndot)$ on $V$ such that $\mathrm{dist}(x,y):=\|x-y\|$ for any $(x,y)\in V^2$. We equip $V$ with the most coarse topology which makes continuous the functions $(y\in V)\mapsto \|x-y\|$ for any $x\in V$. 
In other words, a subset $U$ of $V$ is open if and only if, for any $x \in U$, there is a positive number
$\epsilon$ such that $\{ y \in V \,:\, \| y - x\| < \epsilon \} \subseteq U$.
This topology is said to be \emph{induced by} the seminorm $\norm{\ndot}$\index{topology induced by a seminorm}\index{seminorm!topology induced by ---}. The set $V$ equipped with this topology forms a topological vector space. For any vector subspace $W$ of $V$, the closure of $W$ is also a vector subspace of $V$. In particular, if $W$ is a \emph{hyperplane} in $V$ (namely the kernel of a linear form), then either $W$ is a closed vector subspace of $V$ or $W$ is dense in $V$. For any $x\in V$, the \emph{pseudodistance} between $W$ and $x$\index{pseudodistance} is defined as 
\[\mathrm{dist}(x,W):=\inf\{\|x-y\|\,:\,y\in W\}.\]
Then $\mathrm{dist}(x,W) = 0$ if and only if $x$ belongs to the closure of $W$. In particular, the null space of $(V,\norm{\ndot})$ is a closed subspace, which is the closure of the zero vector subspace $\{0\}$. Thus the topological vector space $V$ is separated if and only if $\norm{\ndot}$ is a norm.

\begin{prop}\label{Pro:applicationlinearcontinue}
Let $(V_1,\|\mathord{\cdot}\|_1)$ and $(V_2,\norm{\ndot}_2)$ be seminormed vector spaces over $k$, and $f:V_1\rightarrow V_2$ be a $k$-linear map.
Then we have the following:
\begin{enumerate}[label=\rm(\arabic*)]
\item\label{Item: continuous implies inclusion of null spaces}
If the map $f$ is continuous, then $f(N_{\|\mathord{\cdot}\|_1})\subseteq N_{\|\mathord{\cdot}\|_2}$.
\item\label{Item: continuous and bounded}
If there is a non-negative constant $C$ such that $\| f(x) \|_2 \leqslant C \| x \|_1$ for all $x \in V_1$,
then the map $f$ is continuous. The converse is true if  either (i) the absolute value $|\ndot|$ is non-trivial or (ii) $\dim_k (V_2/N_{\norm{\ndot}_2}) < \infty$.
\end{enumerate}
\end{prop}
\begin{proof}
\ref{Item: continuous implies inclusion of null spaces} Since $N_{\norm{\ndot}_2}$ is a closed subset of $V_2$, its inverse image by the continuous map $f$ is a closed subset of $V_1$, which clearly contains $0\in V_1$. Hence $f^{-1}(N_{\norm{\ndot}_2})$ contains $N_{\norm{\ndot}_1}$ since $N_{\norm{\ndot}_1}$ is the closure of $\{0\}$ in $V_1$.

\ref{Item: continuous and bounded}
Let $\{x_n\}_{n\in\mathbb N}$ be a sequence in $V_1$ which converges to a point $x\in V_1$. One has
\[\|f(x_n)-f(x)\|_{2}=\|f(x_n-x)\|_{2}\leqslant
C \|x_n-x\|_{1},\]
so that 
the sequence $\{f(x_n)\}_{n\in\mathbb N}$ converges to $f(x)$. Hence the map $f$ is continuous. 

Assume that $f$ is continuous.
First we consider the case where the absolute value $|\ndot|$ is not trivial.
The set $f^{-1}((V_2,\norm{\ndot}_2)_{<1})$ is an open subset of $V_1$ (see Notation \ref{Not:ball}). Hence there exists $\epsilon>0$ such that
$f^{-1}((V_2,\norm{\ndot}_2)_{<1})\supseteq (V_1,\norm{\ndot}_1)_{<\epsilon}$. As the absolute value $|\ndot|$ is not trivial, there exists $a\in k$ such that $0<|a|<1$. 
Let us see that $\| f(x) \|_2 \leqslant (\epsilon |a|)^{-1} \| x \|_1$ for all $x \in V_1$.
If $x \in N_{\norm{\ndot}_1}$, then the assertion is obvious by \ref{Item: continuous implies inclusion of null spaces}, so that we may assume that $x \not\in N_{\norm{\ndot}_1}$. Then 
there exists a unique integer $n$ such that \[\|a^nx\|_{1}<\epsilon\leqslant\|a^{n-1}x\|_{1}=|a|^{n-1}\cdot\|x\|_{1}.\]
Thus $\|f(a^nx)\|_{2}<1$ and hence
\[\|f(x)\|_{2}<|a|^{-n}\leqslant (\epsilon|a|)^{-1}\cdot\|x\|_{1},\]
as desired.

Next we assume that the absolute value $|\ndot|$ is trivial and $\dim_k (V_2/N_{\norm{\ndot}_2}) < \infty$.
By \ref{Item: norm value of trivial valuation} in Corollary~\ref{Cor:finitevlaue} there exist positive numbers $r$ and $\delta$ such that $\|y\|\leqslant r$ for any $y\in V_2$ and that $(V_2,\norm{\ndot}_2)_{<\delta}=N_{\norm{\ndot}_2}$.
If $f$ is continuous, then there exists $\epsilon>0$ such that \[f^{-1}(N_{\norm{\ndot}_2})=f^{-1}((V_2,\norm{\ndot}_2)_{<\delta})\supseteq(V_1,\norm{\ndot}_{1})_{<\epsilon}.\] Therefore one has
$\| f(x) \|_2 \leqslant (r/\epsilon) \| x\|_1$ for all $x \in V_1$.
\end{proof}

\begin{rema}\label{Rem:criterecontinue}
The hypothesis of non-triviality of the absolute value or $\dim_k(V_2/N_{\norm{\ndot}_2}) < \infty$ for the sufficiency part of the above proposition is essential. In fact, if $V$ is an infinite-dimensional vector space over a trivially valued field $k$, equipped with the norm $\norm{\ndot}$ such that $\|x\|=1$ for any $x\in V\setminus\{0\}$, then the topology on $V$ induced by the norm $\norm{\ndot}$ is discrete. In particular, any $k$-linear map from $V$ to a normed vector space over $k$ is continuous. However, one can take a basis $B$ of the vector space $V$ (which is an infinite set) and define a new norm $\norm{\ndot}'$ on $V$ such that
\[\Big\|\sum_{x\in B}n_xx\Big\|'=\displaystyle\max_{x\in B,\,n_x\neq 0}\varphi(x),\]
where $\varphi:B\rightarrow \intervalle{]}{0}{+\infty}{[}$ is a map which is not bounded. If $f$ is the identity map from $(V,\norm{\ndot})$ to $(V,\norm{\ndot}')$, then 
one can not find a non-negative constant $C$ such that $\|x \|' \leqslant C \| x \|$ for all $x \in V$.
\end{rema}

\subsection{Operator seminorm}\label{subsec:Operator norm}
Let $(V_1,\norm{\ndot}_1)$ and $(V_2,\norm{\ndot}_2)$ be seminormed vector spaces over $k$. 
Let $f:V_1\rightarrow V_2$ be a $k$-linear map.
We say that the linear map $f$ is \emph{bounded}\index{bounded linear map} if there is a non-negative constant $C$ such that
$\| f(x) \|_2 \leqslant C \| x \|_1$ for all $x \in V_1$.
Note that if $f$ is bounded, then $f$ is continuous and $f(N_{\norm{\ndot}_1})\subseteq N_{\norm{\ndot}_2}$ by Proposition~\ref{Pro:applicationlinearcontinue}.

If $f(N_{\norm{\ndot}_1})\subseteq N_{\norm{\ndot}_2}$, we denote by $\|f\|$ the element
\[\sup_{x\in V_1\setminus N_{\norm{\sndot}_1}}\frac{\|f(x)\|_{2}}{\|x\|_{1}}\in [0,+\infty].\]
If the relation $f(N_{\norm{\ndot}_1})\subseteq N_{\norm{\ndot}_2}$ does not hold, then by convention $\|f\|$ is defined to be $+\infty$. With this notation, the linear map $f$ is bounded if and only if $\norm{f}<+\infty$.

We denote by $\mathscr L(V_1,V_2)$ the set of all bounded $k$-linear maps from $V_1$ to $V_2$, which forms a vector space over $k$ since, for $(f, g) \in \mathscr L(V_1,V_2)^2$ and $x \in V_1 \setminus N_{\norm{\ndot}_1}$, 
\[
\frac{\| (f + g)(x) \|_2}{\|x\|_1} \leqslant \begin{cases}
{\displaystyle \max \left\{ \frac{\| f(x) \|_2}{\|x\|_1}, \frac{\| g(x) \|_2}{\|x\|_1} \right\}}  & \text{if $\norm{\ndot}_2$ is ultrametric}, \\[3ex]
{\displaystyle\frac{\| f(x) \|_2}{\|x\|_1} + \frac{\| g(x) \|_2}{\|x\|_1}}  & \text{otherwise}.
\end{cases}
\]
The map $\norm{\ndot}:\mathscr L(V_1,V_2)\longrightarrow \intervalle{[}{0}{+\infty}{[}$ defined above is a seminorm, called the \emph{operator seminorm}\index{operator seminorm}\index{seminorm!operator ---}.
Moreover, from the above formula, we observe that, if $\norm{\ndot}_2$ is ultrametric, then the operator seminorm is also ultrametric. If $\norm{\ndot}_2$ is a norm, then the operator seminorm is actually a norm, called \emph{operator norm}\index{operator norm}\index{norm!operator ---}.

In the case where either the absolute value $|\ndot|$ is non-trivial or $\dim_k(V_2/N_{\norm{\ndot}_2}) < \infty$, the space $\mathscr L(V_1,V_2)$ identifies with the vector space of all continuous $k$-linear maps from $V_1$ to $V_2$ (see Proposition \ref{Pro:applicationlinearcontinue}).

\subsection{Quotient seminorm}
\label{Subsec:Quotientnorm}

Let $g : V \to Q$ be a surjective linear map 
of vector spaces over $k$ and $\norm{\ndot}$ be a seminorm on $V$.
We define $\norm{\ndot}_{V \twoheadrightarrow Q}$ to be
\[
\forall\,y\in Q,\quad\|y\|_{V \twoheadrightarrow Q} := \inf \{ \| x \| \,:\, x \in V, \ g(x) = y \}.
\]
Then we have the following proposition:

\begin{prop}
\phantomsection\label{Pro:quotient}
\begin{enumerate}[label=\rm(\arabic*)]
\item\label{Item: quotient seminorm}
$\norm{\ndot}_{V \twoheadrightarrow Q}$ is a seminorm on $Q$. Moreover, if $\norm{\ndot}$ is ultrametric, then
also is $\norm{\ndot}_{V \twoheadrightarrow Q}$.

\item\label{Item: null space of quotient seminorm}
Let $N_{\norm{\ndot}_{V \twoheadrightarrow Q}}$ be the null space of $\norm{\ndot}_{V \twoheadrightarrow Q}$.
Then $g^{-1}(N_{\norm{\ndot}_{V \twoheadrightarrow Q}})$ coincides with the closure of $\Ker(g)$ with respect to the topology induced by the seminorm $\norm{\ndot}$.
In particlar, if $\Ker(g)$ is closed, then $\norm{\ndot}_{V \twoheadrightarrow Q}$ is a norm on $Q$.
\end{enumerate}
\end{prop}

\begin{proof}
\ref{Item: quotient seminorm} In order to see the condition \ref{Item: positive homogeneous} in Definition \ref{Def:norm}, we may assume that $a \neq 0$ since otherwise the assertion is obvious. Then
\begin{align*}
\|a y\|_{V \twoheadrightarrow Q} & = \inf \{ \| x' \| \,:\, x' \in V, \ g(x') = a y \} = \inf \{ \| a x \| \,:\, x \in V, \ g(x) = y \} \\
& = |a| \inf \{ \| x \| \,:\, x \in V, \ g(x) = y \} = |a|\cdot \|y\|_{V \twoheadrightarrow Q}.
\end{align*}
Fix $(y, y') \in Q^2$.
For any $\epsilon > 0$, we can find $(x, x') \in V^2$ such that $g(x) = y$, $g(x') = y'$,
$\|x \| \leqslant \|y\|_{V \twoheadrightarrow Q} + \epsilon$ and 
$\|x'\| \leqslant \|y'\|_{V \twoheadrightarrow Q} + \epsilon$. Then $g(x+x') = y+ y'$ and
\[
\|y + y'\|_{V \twoheadrightarrow Q} \leqslant \| x + x'\| \leqslant \|x \| + \|x'\| \leqslant \|y\|_{V \twoheadrightarrow Q} + \|y'\|_{V \twoheadrightarrow Q} + 2\epsilon,
\]
and hence \ref{Item: triangle inequality} holds. If $\norm{\ndot}$ is ultrametric, in {a} similar way we can see that $\norm{\ndot}_{V \twoheadrightarrow Q}$ is also ultrametric.

\medskip
\ref{Item: null space of quotient seminorm} {Let} $x \in V$ {and} $y = g(x)$.
It is easy to see 
$
\|y \|_{V \twoheadrightarrow Q} = \mathrm{dist}(x, \Ker(g))$.
Therefore
\[
x \in \overline{\Ker(g)} 
\quad\Longleftrightarrow\quad 
\mathrm{dist}(x, \Ker(g)) = 0 
\quad\Longleftrightarrow\quad 
\|y \|_{V \twoheadrightarrow Q} = 0,
\]
as required.
\end{proof}

Given a vector subspace $W$ of a seminormed vector space $(V,\norm{\ndot})$, the seminorm $\norm{\ndot}_{V \twoheadrightarrow V/W}$ on $V/W$ is called the \emph{quotient seminorm}\index{quotient seminorm}\index{seminorm!quotient ---} on $V/W$ of the seminorm $\norm{\ndot}$ on $V$. For simplicity, the seminorm $\norm{\ndot}_{V \twoheadrightarrow V/W}$ is often denoted by $\norm{\ndot}_{V/W}$. If the vector subspace $W$ is closed, then the seminorm $\norm{\ndot}_{V/W}$ is actually a norm, called the \emph{quotient norm}\index{quotient norm}\index{norm!quotient ---} of $\norm{\ndot}$ by the quotient map $V\twoheadrightarrow V/W$. Note that the norm $\norm{\ndot}^{\sim}$ identifies with the quotient norm of $\norm{\ndot}$ by the quotient map $V\twoheadrightarrow V/N_{\norm{\ndot}}$.

\begin{prop}\label{Pro: quotient norm and topology}
Let $(V,\norm{\ndot})$ be a seminormed vector space over $k$ and $W$ be a vector subspace of $V$. The topology on $V/W$ defined by the quotient seminorm coincides with the quotient topology. 
{In particular, the quotient map $V\rightarrow V/W$ is continuous if we equip $V/W$ with the quotient seminorm.}
\end{prop}
\begin{proof}
Recall that the quotient topology is the finest topology on $V/W$ which makes the quotient map $\pi:V\rightarrow V/W$ continuous. In other words, a subset $U$ of $V/W$ is open for the quotient topology if and only if $\pi^{-1}(U)$ is an open subset of $V$. If we equip $V/W$ with the topology induced by the quotient seminorm, then the quotient map is continuous since $\|\pi\|\leqslant 1$ (see Proposition \ref{Pro:applicationlinearcontinue}). Moreover, if $U$ is a subset of $V/W$ such that $\pi^{-1}(U)$ is an open subset of $V$, then, for any $u\in U$ and any $x_0\in V$ such that $\pi(x_0)=u$, there exists $\epsilon>0$ such that \[\{x\in V\,:\,\|x-x_0\|<\epsilon\}\subseteq\pi^{-1}(U).\] Hence for any $v\in V/W$ with $\|v-u\|<\epsilon$, there exists $x\in\pi^{-1}(U)$ such that $\pi(x)=v$. So $U$ is an open subset of $V/W$ for the topology defined by the quotient seminorm. The proposition is thus proved.
\end{proof}

\subsection{Topology of normed vector spaces of finite rank}\label{Subsec:Topologyofnormed}

If $V$ is a finite-dimensional $k$-vector space, then all norms on $V$ induce the same topology. More precisely, we have the following result (see \cite{Bourbaki81} Chapter I, \S2, no.3, Theorem 2 and the remark on the page I.15).

\begin{prop}\label{Pro:topologicalnormedspace}
Assume that the vector space $V$ is of finite rank over $k$. 
If $\norm{\ndot}$ and $\norm{\ndot}'$ are norms on $V$, then there are positive constants $C$ and $C'$ such that
$C \norm{\ndot}' \leqslant \norm{\ndot} \leqslant C' \norm{\ndot}'$ on $V$. In particular, $V$ is complete with respect to $\norm{\ndot}$.\footnote{That is, for any sequence $\{ x_n \}_{n\in\mathbb N}$ in $V$, if \[\lim_{N\rightarrow+\infty}\sup_{\begin{subarray}{c}(n,m)\in \mathbb N^2\\
n\geqslant N,\,m\geqslant N\end{subarray}}\norm{x_n-x_m}=0,\] then there exists $x \in V$ such that
$\lim_{n\to\infty} \| x_n - x \| = 0$.}
\end{prop}
\begin{proof}
Let $\{e_i\}_{i=1}^r$ be a basis of $V$ and $f : k^r \to V$ be the isomorphism  given by
$f(a_1, \ldots, a_r) = a_1 e_1 + \cdots + a_r e_r$. Here we consider the product topology on $k^r$ and the topology induced by any norm $\norm{\ndot}$ on $V$.
By Proposition~\ref{Pro:applicationlinearcontinue}, it is sufficient to show that $f$ is a homeomorphism. 
Since
\[
\| a_1 e_1 + \cdots + a_r e_r \| \leqslant \max \{ |a_1|, \ldots, |a_r| \} \sum_{i=1}^r \| e_i\|,
\]
$f$ is continuous. {It remains to show} that $f^{-1}$ is continuous.

We reason by induction on the rank $r$ of $V$. The case where $r=0$ is trivial. 
In the case where $r=1$, as $|a|/\| a e_1 \| = 1/\| e_1 \|$ for {any} $a \in k^{\times}$,
$f^{-1}$ is continuous by Proposition~\ref{Pro:applicationlinearcontinue}.

Assume that the proposition has been proved for vector spaces of rank $<r$. 
Let $W$ be the vector subspace of $V$ generated by $e_1,\ldots,e_{r-1}$. 
By the induction hypothesis, the map $g:k^{r-1}\rightarrow W$ sending $(a_1,\ldots,a_{r-1})\in k^{r-1}$ to $a_1e_1+\cdots+a_{r-1}e_{r-1}$ is a homeomorphism. In particular, the topological vector space $W$ is complete. As a consequence, $W$ is a closed vector subspace of $V$. By the rank $1$ case of the proposition proved above, the map $\overline{f}$ from $k$ to $V/W$ sending $a\in k$ to $a[e_r]$ is a homeomorphism.

In the following, we show that, if $U$ is an open neighbourhood of $(0,\ldots,0)\in k^r$, then there exists $\epsilon>0$ such that $f(U)$ contains all vectors $x\in V$ satisfying $\|x\|<\epsilon$. Without loss of generality, we may assume that $U$ is the open multidisc $B_\delta^r$, where $B_{\delta}=\{a\in k\,:\,|a|<\delta\}$ and $\delta>0$. Since the map $g:k^{r-1}\rightarrow W$ is a homeomorphism, there exists $\epsilon_1>0$ such that 
\[g(B_\delta^{r-1})\supseteq \{y\in W\,:\,\|y\|<\epsilon_1\}.\]
Let $\delta'=\min\{\epsilon_1/(2\|e_r\|),\delta\}$. Since the map $\overline{f}$ is a homeomorphism, there exists $\epsilon_2>0$ such that 
\[\overline{f}(B_{\delta'})\supseteq\{u\in V/W\,:\,\|u\|_{V/W}<\epsilon_2\},\]
where we consider the quotient norm on $V/W$. 
We claim that \[f(U)\supseteq\{x\in V\,:\,\|x\|<\epsilon\}\] with  $\epsilon=\frac 12\min \{\epsilon_1,\epsilon_2\}$. In fact, if $x$ is an element of $V$ such that $\|x\|<\epsilon$, then its class in $V/W$ has  norm $<\epsilon_2$. Hence there exists $a_r\in B_{\delta'}$ such that $[x]=a_r[e_r]$. Moreover, one has \[\|x-a_re_r\|\leqslant\|x\|+|a_r|\cdot\|e_r\|<\frac{1}{2}\epsilon_1+\delta'\|e_r\|\leqslant\epsilon_1.\]
Hence there exists $(a_1,\ldots,a_{r-1})\in B_\delta^{r-1}$ such that $g(a_1,\ldots,a_{r-1})=x-a_re_r$.
Thus $(a_1,\ldots,a_r)$ is an element in $B_\delta^r$ such that $f(a_1,\ldots,a_r)=x$. The proposition is proved.
\end{proof}

\begin{coro}\label{Coro:continuous linearform}
Let $f : V_1 \rightarrow V_2$ be a linear map 
of vector spaces over $k$, and
let $\norm{\ndot}_1$ and $\norm{\ndot}_2$ be seminorms on $V_1$ and $V_2$, respectively.
We assume that $f(N_{\norm{\ndot}_1})\subseteq N_{\norm{\ndot}_2}$ and $\dim_k(V_2/N_{\norm{\ndot}_2}) < \infty$. Then
the following conditions are equivalent:
\begin{enumerate}[label=\rm(\alph*)]
\item\label{Item: the map f is continuous} the map $f$ is continuous;
\item\label{Item: inverse image of null space is closed} $f^{-1}(N_{\norm{\ndot}_2})$ is a closed vector subspace of $V$;
\item\label{Item: norm of f is finite} $\|f\|$ is finite.
\end{enumerate}
\end{coro}
\begin{proof}
``\ref{Item: the map f is continuous}$\Longrightarrow$\ref{Item: inverse image of null space is closed}'': Since $f^{-1}(N_{\norm{\ndot}_2})$ is the inverse image by $f$ of the closed subset $N_{\norm{\ndot}_2}$ of $V_2$, if $f$ is continuous, then it is a closed subset of $V_1$.

``\ref{Item: inverse image of null space is closed}$\Longrightarrow$\ref{Item: norm of f is finite}'': The assertion is trivial when $f(V_1)\subseteq N_{\norm{\ndot}_2}$. In the following, we assume that $f(V_1)\not\subseteq N_{\norm{\ndot}_2}$. We set 
\[
Q := f(V_1)/(f(V_1)\cap N_{\norm{\ndot}_2}) \cong (f(V_1) + N_{\norm{\ndot}_2})/N_{\norm{\ndot}_2} \not= \{ 0 \}.
\]
Let $\norm{\ndot}_{Q}$ be the quotient seminorm on $Q$ induced by $V_1\to f(V_1) \to Q$ and $\norm{\ndot}_{1}$. By Proposition \ref{Pro:quotient} and the condition \ref{Item: inverse image of null space is closed}, the seminorm $\norm{\ndot}_Q$ is actually a norm.
Moreover, we can consider $Q$ as a vector subspace of $V_2/N_{\norm{\ndot}_2}$. Let $\norm{\ndot}'_Q$ be the restriction of $\norm{\ndot}_2^{\sim}$ {to} $Q$, where $\norm{\ndot}_2^\sim$ is the norm associated with the seminorm $\norm{\ndot}_2$.
By Proposition~\ref{Pro:topologicalnormedspace}, 
there is a constant $C$ with $\norm{\ndot}'_Q \leqslant C \norm{\ndot}_{Q}$.
Thus, for any $x\in V_1\setminus N_{\norm{\ndot}_1}$, one has
\[
\frac{\| f(x) \|_{2}}{\| x \|_{1}} = \frac{\| [f(x)] \|'_{Q}}{\| x \|_{1}} \leqslant \frac{C \|[f(x)]\|_{Q}}{\| x\|_{1}} \leqslant C,
\]
which implies $\| f \| \leqslant C$.

``\ref{Item: norm of f is finite}$\Longrightarrow$\ref{Item: the map f is continuous}'' follows from Proposition~\ref{Pro:applicationlinearcontinue}.  
\end{proof}

\begin{coro}\label{Coro:equivalenceofnrom}
Let $(V, \norm{\ndot})$ be a finite-dimensional seminormed vector space over $k$. Then we have the following:
\begin{enumerate}[label=\rm(\arabic*)]
\item\label{Item: vector space containing null space is closed}
Every vector subspace of $V$ containing $N_{\norm{\ndot}}$ is closed.

\item\label{Item: criterion of continuity}
Let $(V', \norm{\ndot}')$ be a seminormed vector space over $k$ and $f : V \to V'$ be a linear map 
of vector spaces over $k$ such that $f(N_{\norm{\ndot}})\subseteq N_{\norm{\ndot}'}$.
Then $f$ is continuous and $\| f \| < +\infty$.
\item\label{Item: criterion of continuous linear form}A linear form on $V$ is bounded if and only if its kernel contains $N_{\norm{\ndot}}$. 
\end{enumerate}
\end{coro}

\begin{proof}
\ref{Item: vector space containing null space is closed} Let $\pi:V\rightarrow V/N_{\norm{\ndot}}$ be the canonical projection map and $\norm{\ndot}^\sim$ be the norm on $V/N_{\norm{\ndot}}$ associated with $\norm{\ndot}$. 
By Proposition \ref{Pro: quotient norm and topology}, the linear map $\pi$ is continuous. If $W$ is a vector subspace of $V$ containing $N_{\norm{\ndot}}$, then one has $W=\pi^{-1}(\pi(W))$. By Proposition~\ref{Pro:topologicalnormedspace}, $\pi(W)$ is complete with respect to the induced norm of $\norm{\ndot}^\sim$ on $\pi(W)$, so that
$\pi(W)$ is closed. Hence $W$ is also closed since it is the inverse image of a closed subset of $V/N_{\norm{\ndot}}$ by a continuous linear map.

\ref{Item: criterion of continuity} By replacing $V'$ by $f(V)$, we may assume that $\dim_k(V') < \infty$.
Thus the assertion follows from \ref{Item: vector space containing null space is closed} and Corollary~\ref{Coro:continuous linearform}.

\ref{Item: criterion of continuous linear form} Let $f:V\rightarrow k$ be a linear form. If $f$ is bounded, by Corollary \ref{Coro:continuous linearform}, the kernel of $f$ is a closed vector subspace of $V$, hence it contains the closure of $\{0\}$, which is $N_{\norm{\ndot}}$. Conversely, if $\Ker(f)\supseteq N_{\norm{\ndot}}$, then by \ref{Item: criterion of continuity}, the linear form $f$ is bounded. 
\end{proof}

\begin{prop}\phantomsection\label{prop:quotient:norm:linear:map}
\begin{enumerate}[label=\rm(\arabic*)]
\item\label{Item: successive quotient seminorm}
Let $V \overset{\alpha}{\longrightarrow} W \overset{\beta}{\longrightarrow} Q$ be a sequence of surjective linear maps of finite-dimensional
 vector spaces over $k$. For any seminorm $\norm{\ndot}$ on $V$, one has
$\norm{\ndot}_{V \twoheadrightarrow W, W \twoheadrightarrow Q} = \norm{\ndot}_{V \twoheadrightarrow Q}$.

\item\label{Item: diagram of quotient seminorm}
Let 
\[\xymatrix{
V \ar[r]^-{f}\ar@{->>}[d]_-\alpha& W\ar@{->>}[d]^-\beta \\
V' \ar[r]_-{g}& W'}
\]
be a commutative diagram of linear maps of  finite-dimensional vector spaces over $k$ such that $\alpha$ and $\beta$ are surjective.
Then we have the following:
\begin{enumerate}[label=\rm(\arabic*)]
\renewcommand{\labelenumii}{\rm(\arabic{enumi}.\alph{enumii})}
\item\label{SubItem: diagram of quotient seminorm1}
Let $\norm{\ndot}_V$ and $\norm{\ndot}_W$ be seminorms on $V$ and $W$;
let $\norm{\ndot}_{V'}$ and $\norm{\ndot}_{W'}$ be the quotient seminorms of $\norm{\ndot}_V$ and $\norm{\ndot}_W$ on $V'$ and $W'$, respectively.
{If} $f(N_{\norm{\ndot}_V})\subseteq N_{\norm{\ndot}_W}$, then $g(N_{\norm{\ndot}_{V'}})\subseteq N_{\norm{\ndot}_{W'}}$ and $\| g \| \leqslant \| f\|$. 

\item\label{SubItem: diagram of quotient seminorm2}
We assume that $f$ and $g$ are injective. Let $\norm{\ndot}_W$ be a seminorm {on} $W$.
Then $\norm{\ndot}_{W, V \hookrightarrow W, V \twoheadrightarrow V'} \geqslant \norm{\ndot}_{W, W \twoheadrightarrow W', V' \hookrightarrow W'}$.
Moreover, if $\Ker(\beta) \subseteq f(V)$, then the equality $\norm{\ndot}_{W, V \hookrightarrow W, V \twoheadrightarrow V'} = \norm{\ndot}_{W, W \twoheadrightarrow W', V' \hookrightarrow W'}$ holds.
\end{enumerate}
\end{enumerate}
\end{prop}
\begin{proof}
\ref{Item: successive quotient seminorm} For $q\in Q$, one has
\[\begin{split}
\|q\|_{V\twoheadrightarrow W,W\twoheadrightarrow Q}&=\inf_{y\in W,\,\beta(y)=q}\|y\|_{V\twoheadrightarrow W}=\inf_{\begin{subarray}{c}y\in W\\\beta(y)=q\end{subarray}}\inf_{\begin{subarray}{c}x\in V\\\alpha(x)=y
\end{subarray}}\|x\|_V\\&=\inf_{x\in V,\,\beta(
\alpha(x))=q}\|x\|_V=\|q\|_{V\twoheadrightarrow Q},
\end{split}\]
as desired.

\medskip
\ref{SubItem: diagram of quotient seminorm1} By Proposition \ref{Pro:quotient}, $\alpha^{-1}(N_{\norm{\ndot}_{V'}})$ is the closure of $\Ker(\alpha)$ in $V$, hence is equal to $\Ker(\alpha)+N_{\norm{\ndot}_V}$. Let $y'$ be an element in $N_{\norm{\ndot}_{V'}}$. There then exists $y\in N_{\norm{\ndot}_V}$ such that $\alpha(y)=y'$. Therefore \[g(y')=g(\alpha(y))=\beta(f(y))\in N_{\norm{\ndot}_{W'}}\] since $f(y)\in N_{\norm{\ndot}_{W}}$.

It remains to prove that $\|g\|\leqslant\norm{f}$. Let $x'$ be an element of $V'$. For any $x\in V$ with $\alpha(x)=x'$, one has
\[\norm{g(x')}_{W'}=\norm{g(\alpha(x))}_{W'}=\norm{\beta(f(x))}_{W'}\leqslant\norm{f(x)}_{W}\leqslant\norm{f}\cdot\norm{x}_V,\]
which leads to
\[\norm{g(x')}_{W'}\leqslant\norm{f}\inf_{x\in V,\,\alpha(x)=x'}\norm{x}_V=\norm{f}\cdot\norm{x'}_{V'}.\]

\medskip
\ref{SubItem: diagram of quotient seminorm2} Note that $f(\Ker(\alpha)) = f(V) \cap \Ker(\beta)$. Therefore, for $v \in V$,
\begin{align*}
& \|\alpha(v) \|_{W, V \hookrightarrow W, V \twoheadrightarrow V'} = \inf \{ \| x \|_W \,:\, x \in f(v) + (f(V) \cap \Ker(\beta)) \} \\
\intertext{and}
& \norm{\alpha(v)}_{W, W \twoheadrightarrow W', V' \hookrightarrow W'} = \inf \{ \| x \|_W \,:\, x \in f(v) + \Ker(\beta) \},
\end{align*}
so that the first assertion follows. Moreover, if $\Ker(\beta) \subseteq f(V)$, then $f(V) \cap \Ker(\beta) = \Ker(\beta)$.
Thus the second assertion {holds}.
\end{proof}

\begin{prop}\phantomsection\label{lem:quotient:dim:1:comm}
\begin{enumerate}[label=\rm(\arabic*)]
\item\label{Item: norm of quotient operator}
Let $f : V \to W$ be a surjective linear map of vector spaces over $k$,
$\norm{\ndot}_V$ be a seminorm on $V$ and $\norm{\ndot}_W$ be the quotient seminorm of $\norm{\ndot}_V$ on $W$.
If the seminorm $\norm{\ndot}_W$ does not vanish,
then $\norm{f} = 1$.

\item\label{Item: continuit of diagram}
Let 
\[
\xymatrix{ V \ar[r]^{f_1} \ar[rd]_{f_2} & W_1 \ar[d]^{g} \\
& W_2}
\]
be a commutative digram of linear maps of finite-dimensional vector spaces over $k$
such that $g$ is an isomorphism and $\dim_k(W_1) = \dim_k(W_2) = 1$.
Let $\norm{\ndot}_V$, $\norm{\ndot}_{W_1}$ and $\norm{\ndot}_{W_2}$ be seminorms of $V$, $W_1$ and $W_2$, respectively.
Then $\|f_2\| = \|f_1\|\cdot\|g\|$ provided that $f_1$, $f_2$ and $g$ are continuous.
\end{enumerate}
\end{prop}

\begin{proof}
\ref{Item: norm of quotient operator} Since $\| v \|_V \geqslant \| f(v) \|_W$ for any $v\in V$, one has $\|f\| \leqslant 1$. Let $w$ be an element of $W$ such that $\|w\|_W>0$. Since 
\[\|w\|_W=\displaystyle\inf_{v\in V,\,f(v)=w}\|v\|_V,\] one has
\[1=\inf_{v\in V,\,f(v)=w}\frac{\|v\|_V}{\|f(v)\|_W}\geqslant\|f\|^{-1},\]
which leads to $\|f\|\geqslant 1$.

\ref{Item: continuit of diagram} As $g$ is an isomorphism and $\dim_k(W_1) = \dim_k(W_2) = 1$, for any $w_1\in W_1$,
$\| g \|\cdot\| w_1 \|_{W_1} = \| g(w_1) \|_{W_2}$, Therefore,
\[\begin{split}
\| f_2 \| &= \sup_{v \in V \setminus N_{\norm{\sndot}_V} } \frac{\|f_2(v)\|_{W_2}}{\| v \|_V} =
\sup_{v \in V \setminus N_{\norm{\sndot}_V}} \frac{\|g(f_1(v))\|_{W_2}}{\| v \|_V} \\& =
\sup_{v \in V \setminus N_{\norm{\sndot}_V}} \| g \| \frac{\|f_1(v)\|_{W_1}}{\| v \|_V}
= \| g \|\cdot\| f_1\|,
\end{split}
\]
as required.
\end{proof}

\begin{prop}\label{Pro: induce quotient zero ball}
Let $(V,\norm{\ndot}_V)$ be a finite-dimensional seminormed vector space over $k$, $W$ be a vector subspace of $V$ and $Q$ be the quotient vector space $V/W$. We denote by $i:W\rightarrow V$ and $\pi:V\rightarrow Q$ the inclusion map and the projection map, 
respectively. Let $\norm{\ndot}_W$ be the restriction of $\norm{\ndot}_V$ {to} $W$ and $\norm{\ndot}_Q$ be the quotient seminorm of $\norm{\ndot}_V$ on $Q$. Then one has $i(N_{\norm{\ndot}_W})\subseteq N_{\norm{\ndot}_V}$ and $\pi(N_{\norm{\ndot}_V})\subseteq N_{\norm{\ndot}_Q}$. Moreover, the linear maps $i$ and $\pi$ induce short exact sequences
\begin{equation}\label{Equ: exact sequence of null balls}\xymatrix{0\ar[r]&N_{\norm{\ndot}_W}\ar[r]&N_{\norm{\ndot}_V}\ar[r]&N_{\norm{\ndot}_Q}\ar[r]&0}\end{equation}
and
\begin{equation}\label{Equ: exact sequence of quotients}
\xymatrix{0\ar[r]&W/N_{\norm{\ndot}_W}\ar[r]&V/N_{\norm{\ndot}_V}\ar[r]&Q/N_{\norm{\ndot}_Q}\ar[r]&0},
\end{equation}
and the induced norm (resp. quotient norm) of $\norm{\ndot}_V^\sim$ on $W/N_{\norm{\ndot}_W}$ (resp. $Q/N_{\norm{\ndot}_Q}$) identifies with $\|\ndot\|^{\sim}_W$ (resp. $\norm{\ndot}_Q^\sim$).
\end{prop}
\begin{proof}
The relations $i(N_{\norm{\ndot}_W})\subseteq N_{\norm{\ndot}_V}$ and $\pi(N_{\norm{\ndot}_V})\subseteq N_{\norm{\ndot}_Q}$ follow directly from the definition of induced and quotient seminorms. Moreover, by definition one has $N_{\norm{\ndot}_W}=N_{\norm{\ndot}_V}\cap W$. 

For any element $x\in V$, $\pi(x)$ lies in $N_{\norm{\ndot}_Q}$ if and only if $x\in W+N_{\norm{\ndot}_V}$ since $W+N_{\norm{\ndot}_V}$ is the closure of $W$ in $V$. Therefore one has \[N_{\norm{\ndot}_Q}\cong (W+N_{\norm{\ndot}_V})/W\cong N_{\norm{\ndot}_V}/(W\cap N_{\norm{\ndot}_V})=N_{\norm{\ndot}_V}/N_{\norm{\ndot}_W},\]
which proves that \eqref{Equ: exact sequence of null balls} is an exact sequence. The exactness of \eqref{Equ: exact sequence of null balls} implies that of \eqref{Equ: exact sequence of quotients}. Moreover, if $x$ is an element of $W$, then 
\[\norm{[x]}_{W}^\sim=\norm{x}_W=\norm{x}_V=\norm{[x]}_V^\sim.\]
If $u$ is an element in $Q$, then
\[\norm{[u]}_Q^\sim=\norm{u}_Q=\inf_{y\in V,\,\pi(y)=u}\norm{y}_V=\inf_{v\in V/N_{\norm{\sndot}_V},\,\widetilde{\pi}(v)=[u]}\norm{v}_{V}^\sim,\]
where $\widetilde{\pi}:V/N_{\norm{\ndot}_V}\rightarrow Q/N_{\norm{\ndot}_Q}$ is the linear map induced by $\pi$. Hence $\norm{\ndot}_Q^\sim$ identifies with the quotient norm of $\norm{\ndot}_V^\sim$.
\end{proof}

\subsection{Dual norm}
\label{Subsec:dualnorm}
Let $(V,\norm{\ndot})$ be a seminormed vector space over $k$. We denote by $V^*$ the vector space $\mathscr L(V,k)$ (where we consider $|\ndot|$ as a norm on $k$) of bounded $k$-linear forms on $V$ (which necessarily vanish  on $N_{\norm{\ndot}}$), called the \emph{dual} normed vector space of $V$. The operator norm on $V^*$ is called the \emph{dual norm}\index{dual norm@dual norm}\index{norm!dual ---} of $\norm{\ndot}$, denoted by $\norm{\ndot}_*$. Note that in general $V^*$ is different from the (algebraic) dual vector space $V^\vee:=\Hom_{k}(V,k)$. {One has} \[V^{*} \subseteq (V/N_{\norm{\ndot}})^{\vee}
= \{ \varphi \in V^{\vee} \,:\, \rest{\varphi}{N_{\norm{\ndot}}} = 0 \}.\]

\begin{prop}
\label{Pro: equality Vstar}
Let $(V,\norm{\ndot})$ be a finite-dimensional seminormed vector space over $k$. Then the map $(V/N_{\norm{\ndot}})^\vee\rightarrow V^\vee$ sending $\varphi\in(V/N_{\norm{\ndot}})^\vee$ to its composition with the projection map $V\rightarrow V/N_{\norm{\ndot}}$ defines an isomorphism between $(V/N_{\norm{\ndot}})^\vee$ and $V^*$. In particular, the equality $V^*=V^\vee$ holds when $\|\ndot\|$ is a norm.
\end{prop}
\begin{proof}
By Corollary~\ref{Coro:equivalenceofnrom}, a linear form on $V$ is bounded if and only if its kernel contains $N_{\norm{\ndot}}$. Therefore  $V^*$ is canonically isomorphic to $(V/N_{\norm{\ndot}})^\vee$.
\end{proof}

If $x$ is an element of $V$, for any $\alpha\in V^*$ one has
\begin{equation}\label{Equ:dualnormformula}|\alpha(x)|\leqslant\|\alpha\|_*\cdot\|x\|.\end{equation}
Therefore the linear form on $V^*$ sending $\alpha\in V^*$ to $\alpha(x)\in k$ is bounded. Hence one obtains a $k$-linear map  from $V$ to the double dual space $V^{**}$ whose kernel contains $N_{\norm{\ndot}}$. It is called the \emph{canonical linear map}\index{canonical linear map} from $V$ to $V^{**}$. The \emph{double dual norm}\index{double dual norm}\index{norm!double dual ---} $\norm{\ndot}_{**}$ on $V^{**}$ induces by composition with the canonical $k$-linear map $V\rightarrow V^{**}$ a seminorm on $V$ which we still denote by $\norm{\ndot}_{**}$ by abuse of notation. 
Moreover, by \eqref{Equ:dualnormformula} we obtain
\begin{equation}\label{Equ:doubledual}\forall\,x\in V,\quad\|x\|_{**}\leqslant\|x\|.\end{equation}
We say that $(V,\norm{\ndot})$ is \emph{reflexive}\index{reflexive}\index{norm!reflexive ---} if the $k$-linear map $V\rightarrow V^{**}$ described above induces an isometric $k$-linear isomorphism between the normed vector spaces $(V/N_{\norm{\ndot}},\norm{\ndot}^\sim)$ and $(V^{**},\norm{\cdot}_{**})$.

The following proposition shows that, in the Archimedean case, the seminorm $\norm{\ndot}_{**}$ on $V$ identifies with $\norm{\ndot}$. In particular, a finite-dimensional seminormed vector space over an Archimedean complete field is always reflexive. We will see further in Corollary \ref{Cor:doubledual} that any finite-dimensional ultrametrically seminormed vector space over $k$ is also reflexive.

\begin{prop}\label{Pro:doubledualarch}
Assume that the absolute value $|\ndot|$ is Archimedean. Let $(V,\norm{\ndot})$ be a seminormed vector space over $k$. For any $x\in V$ one has $\|x\|=\|x\|_{**}$. 
\end{prop}
\begin{proof}
This is a direct consequence of Hahn-Banach theorem. In fact, if $x$ is a vector in $V\setminus N_{\|\ndot\|}$, then by Hahn-Banach theorem there exists a $k$-linear form $\tilde{f}:V/N_{\norm{\ndot}} \rightarrow k$ such that $\tilde{f}(\pi(x))=\|\pi(x)\|^{\sim}$ and that $|\tilde{f}(\pi(y))|\leqslant\|\pi(y)\|^{\sim}$ for any $y\in V$, where $\pi : V \to V/N_{\norm{\ndot}}$ is the canonical linear map. 
If we  set $f = \tilde{f} \circ \pi$, then
$f(x) = \|x\|$ and $|f(y)| \leqslant \| y \|$ for any $y \in V$. In particular, $\| f \|_* = 1$. Hence
\[\|x\|_{**}\geqslant \frac{|f(x)|}{\|f\|_*}= \|x\|.\]
\end{proof}

\begin{rema}\label{Rem:doubledual}
The above proposition is not true when the absolute value $|\ndot|$ is non-Archimedean. Let $(V,\norm{\ndot})$ be a normed vector space over $k$. If the absolute value $|\ndot|$ is non-Archimedean, then the dual norm $\norm{\ndot}_*$ is necessarily ultrametric
(cf. Subsection~\ref{subsec:Operator norm}). 
For the same reason, the double dual norm $\norm{\ndot}_{**}$ is ultrametric, and hence cannot identify with $\norm{\ndot}$ on $V$ once the norm $\norm{\ndot}$ is not ultrametric. In the next section, we will establish the analogue of the above proposition in the case where $V$ is of finite rank over $k$ and $\norm{\ndot}$ is ultrametric (see Corollary \ref{Cor:doubledual}). We refer to \cite{Cohen48} and \cite{Ingleton52} for more general results on non-Archimedean Hahn-Banach theorem.
\end{rema}

\begin{prop}\label{Pro:dualquotient}
Let $(V,\norm{\ndot}_V)$ be a seminormed vector space over $k$, $W$ be a 
vector subspace of $V$ and $Q=V/W$ be the quotient space.  Let $\norm{\ndot}_Q$ be the quotient seminorm on $Q$ induced by $\norm{\ndot}_V$. Then the map $Q^*\rightarrow V^*$ sending $\varphi\in Q^*$ to the composition of $\varphi$ with the projection map $V\rightarrow Q$ is an isometry from $Q^*$ to its image (equipped with the induced norm), where we consider the dual norms $\norm{\ndot}_{Q,*}$ and $\norm{\ndot}_{V,*}$ on $Q^*$ and $V^*$, respectively.
\end{prop}
\begin{proof}
Note that
$
{\|v\|_Q^{-1}} = \displaystyle\sup_{x\in V,\,[x]=v} {\| x \|_V^{-1}}
$
for $v \in Q \setminus N_{\norm{\ndot}_Q}$.
Thus, for $\varphi\in Q^*$,
\[
\| \varphi \|_{Q,*} = \sup_{v \in Q \setminus N_{\norm{\sndot}_Q}}
\frac{|\varphi(v)|}{\| v \|_Q} = \sup_{x \in V\setminus N_{\norm{\sndot}_V}} \frac{|\varphi([x])|}{\| x \|_V} ,
\]
as required.
\end{proof}

\begin{rema}\label{Rem:comparaisondualitye}
The dual statement of the above proposition for the dual of {a} restricted seminorm is much more subtle. Let $(V,\norm{\ndot}_V)$ be a seminormed vector space over $k$ and $W$ be a vector subspace of $V$. We denote by $\norm{\ndot}_W$ the restriction of the seminorm $\norm{\ndot}_V$ {to} $W$. Then the restriction {to} $W$ of bounded linear forms on $V$ defines a $k$-linear map $\pi$ from $V^*$ to $W^*$.  We are interested in the nature of the dual norm $\norm{\ndot}_{W,*}$. In the case where $k$ is Archimedean, the $k$-linear map $\pi:V^*\rightarrow W^*$ is surjective and the norm $\norm{\ndot}_{W,*}$ identifies with the quotient norm of $\norm{\ndot}_{V,*}$. This is a direct consequence of Hahn-Banach theorem which asserts that any bounded linear form on $W$ extends to $V$ with the same operator norm (see Lemma \ref{Lem:quotientequalsrestriction}
 for more details). However, the non-Archimedean analogue of this result is not true, even in the case where $V$ is finite-dimensional. In fact, assume that $(k,|\ndot|)$ is non-Archimedean and $V$ is a finite-dimensional vector space of rank $\geqslant 2$ over $k$, equipped with a norm $\norm{\ndot}$ which is not ultrametric. Then the double dual norm $\norm{\ndot}_{**}$ on $V$ is bounded from above by $\norm{\ndot}$ (see \eqref{Equ:doubledual}), and there exists at least an element $x\in V$ such that $\|x\|_{**}<\|x\|$ since $\norm{\ndot}_{**}$ is ultrametric but $\norm{\ndot}$ is not. However, both norms $\norm{\ndot}_{**}$ and $\norm{\ndot}$ induce the same dual norm on $V^\vee$ (see Proposition \ref{Pro:doubedualandquotient}). 
Therefore the quotient norm of $\norm{\ndot}_*$ on $(kx)^\vee$ can not identify with the dual norm of the restriction of $\norm{\ndot}$ {to} $kx$. We will show in Proposition \ref{Pro:quotientdualnonarch} that the non-Archimedean analogue of the above statement is true when the norm on $V$ is ultrametric.
\end{rema}

\subsection{Seminorm of the dual operator}\label{Subsec:Norm of the dual operator}

Let $(V_1,\norm{\ndot}_1)$ and $(V_2,\norm{\ndot}_2)$ be seminormed vector spaces over $k$, and $f:V_1\rightarrow V_2$ be a bounded linear map. 
Note that $f(N_{\norm{\ndot}_{1}}) \subseteq N_{\norm{\ndot}_2}$ by Proposition~\ref{Pro:applicationlinearcontinue}.
For any $\alpha\in V_2^*$, we let $f^*(\alpha)$ be the linear form on $V_1$ which sends $x\in V_1$ to $\alpha(f(x))$. Note that for $x\in V_1$ one has
\begin{equation}\label{Equ:dualoperatornorm} |f^*(\alpha)(x)|=|\alpha(f(x))|\leqslant\|\alpha\|_{2,*}\cdot\|f(x)\|_{2}\leqslant\|\alpha\|_{2,*}\cdot\|f\|\cdot\|x\|_1.\end{equation}
Therefore $f^*(\alpha)$ is a bounded linear form on $V_1$. Thus $f^*$ defines a  linear map from $V_2^*$ to $V_1^*$.

\begin{prop}\label{Pro:normdualopt}
Let $(V_1,\norm{\ndot}_1)$ and $(V_2,\norm{\ndot}_2)$ be seminormed vector spaces over $k$ and  $f:V_1\rightarrow V_2$ be a bounded linear map. Then one has $\|f^*\|\leqslant\|f\|$. The equality holds when {$\norm{\ndot}_2=\norm{\ndot}_{2,**}$ on $V_2$.}
\end{prop}
\begin{proof}
By \eqref{Equ:dualoperatornorm} we obtain that, if $\alpha$ is an element of $V_2^*$, then one has
\[\norm{f^*(\alpha)}_{1,*}\leqslant\|f\|\cdot\norm{\alpha}_{2,*}.\]
Hence $\|f^*\|\leqslant \|f\|$. If we apply this inequality to $f^*$, we obtain $\|f^{**}\|\leqslant\|f^{*}\|\leqslant\|f\|$. Let $\iota_1:V_1\rightarrow V_1^{**}$ and $\iota_2:V_2\rightarrow V_2^{**}$ be the canonical linear maps.
For any vector $x$ in $V_1$, one has $f^{**}(\iota_1(x))=\iota_2(f(x))$. Moreover, {if $\norm{\ndot}_2=\norm{\ndot}_{2,**}$ on $V_2$, then} \[\norm{f^{**}}\geqslant\sup_{x\in V_1,\,\|x\|_{1,**}>0}\frac{\|f^{**}(\iota_1(x))\|_{2,**}}{\|x\|_{1,**}}\geqslant \sup_{x\in V_1,\,\|x\|_{1}>0}\frac{\|f(x)\|_{2}}{\|x\|_{1}}=\|f\|, \]
as required.
\end{proof}

\subsection{Lattices and norms}
\label{Subsec:latticesandnorms} In this subsection, we assume that the absolute value $|\ndot|$ is non-Archimedean. Let $\mathfrak o_k:=\{a\in k\,:\,|a|\leqslant 1\}$ be the closed unit ball of $(k,|\ndot|)$. It is a valuation ring, namely for any $a\in k\setminus \mathfrak o_k$ one has $a^{-1}\in \mathfrak o_k$  (see \cite{Bourbaki64} Chapter IV, \S1, no.2). It is a discrete valuation ring (namely a Noetherian valuation ring) if and only if the absolute value $|\ndot|$ is discrete, namely the image of $k^{\times}$ by $|\ndot|$ is a discrete subgroup of $(\mathbb R_{>0},\times)$ (see \cite{Bourbaki64} Chapter IV, \S3, no.6). In this case, $\mathfrak o_k$ is a principal ideal domain. In particular, {its maximal} ideal $\{a\in k\,:\,|a|<1\}$ is generated by one element $\varpi$, called a \emph{uniformizing parameter}\index{uniformizing parameter} of $k$. Note that, if the absolute value $|\ndot|$ is not discrete, then $|k^{\times}|$ is a dense subgroup of $(\mathbb R_{>0},\times)$. This results from the facts that a subgroup of $(\mathbb R,+)$ is either discrete or dense (cf. \cite{Bourbaki06} Chapter V, \S1, no.1 and \S4, no.1) and that the exponential function defines an isomorphism between the topological groups $(\mathbb R,+)$ and $(\mathbb R_{>0},\times)$.

\begin{defi}\label{Def:lattice}
Let $V$ be a finite-dimensional vector space over $k$. 
A sub-$\mathfrak o_k$-module $\mathcal V$ of $V$ is called a \emph{lattice of $V$}\index{lattice} if $\mathcal V$ generates $V$ as a vector space over $k$ (i.e. the natural linear map 
$\mathcal V \otimes_{\mathfrak o_k} k \to V$ is surjective) and $\mathcal V$ is bounded in $V$ for {a} certain norm on $V$ (or equivalently for any norm on $V$, see Proposition~\ref{Pro:topologicalnormedspace}). 
In particular, if $\mathcal V$ is a sub-$\mathfrak o_k$-module of finite type of $V$, which generates $V$ as a vector space over $k$, then it is a lattice in $V$. If $\mathcal V$ is a lattice of $V$, we define a function $\norm{\ndot}_{\mathcal V}$ on $V$ as follows:\footnote{Note that, in the case where $|\ndot|$ is not the trivial absolute value, one has \[\inf\{|a|\,:\,a\in k^{\times},\;a^{-1}0\in\mathcal V\}=0.\] However, this equality does not hold when $|\ndot|$ is trivial.}
\[\forall\,x\in V\setminus\{0\},\quad\|x\|_{\mathcal V}:=\inf\{|a|\,:\,a\in k^{\times},\;a^{-1}x\in\mathcal V\},\quad\text{and $\norm{0}_{\mathcal V}:=0$}.\]
Clearly, {if $\mathcal V$ and $\mathcal V'$ are lattices of $V$ such that} $\mathcal V\subseteq\mathcal V'$, then one has $\norm{\ndot}_{\mathcal V}\geqslant\norm{\ndot}_{\mathcal V'}$. Note that, if the absolute value $|\ndot|$ is trivial, then $\mathfrak o_k=k$ and the only lattice of $V$ is $V$ itself.
\end{defi}

\begin{prop}\label{Pro:comparisonbyunitball}
Let $V$ be a finite-dimensional vector space over $k$, $\mathcal V$ be a lattice of $V$ and $\norm{\ndot}$ be a norm on $V$. Assume that $\mathcal V$ is contained in the unit ball of $(V,\norm{\ndot})$, then one has $\norm{\ndot}_{\mathcal V}\geqslant \norm{\ndot}$.
\end{prop}
\begin{proof}
Let $x\in V\setminus\{0\}$ and $a$ be an element of $k^{\times}$ such that $a^{-1}x\in\mathcal V$.  One has \[\|a^{-1}x\|=|a|^{-1}\cdot\|x\|\leqslant 1\] since $\mathcal V$ is contained in the unit ball of $(V,\norm{\ndot})$. Therefore $\|x\|\leqslant |a|$. Thus we deduce that 
\[\|x\|\leqslant\inf\{|a|\,:\,a\in k^{\times},\,a^{-1}x\in\mathcal V\}=\|x\|_{\mathcal V}.\]
\end{proof}

In the case where the absolute value $|\ndot|$ is non-trivial, the balls in an ultrametrically normed vector space are natural examples of lattices.

\begin{prop}\label{prop:ultrametric:ball:lattice}
Assume that the absolute value $|\ndot|$ is non-trivial. Let $V$ be a finite-dimensional vector space over $k$, equipped with an ultrametric norm $\norm{\ndot}$. For any $\epsilon>0$ the balls
$V_{\leqslant\epsilon}=\{x\in V\,:\,\|x\|\leqslant \epsilon\}$ and $V_{<\epsilon}=\{x\in V\,:\,\|x\|<\epsilon\}$
are both lattices of $V$. 
\end{prop}
\begin{proof}
Since the norm $\norm{\ndot}$ is ultrametric, both balls $V_{\leqslant\epsilon}$ and $V_{<\epsilon}$ are stable by addition. Clearly they are also stable by the multiplication by an element in $\mathfrak o_k$. Therefore they are sub-$\mathfrak o_k$-modules of $V$. Moreover, by definition they are bounded subsets of $V$. It remains to verify that they generate $V$ as a vector space over $k$. It suffices to treat the open ball case. Let $\{e_i\}_{i=1}^r$ be a basis of $V$ over $k$. Since the absolute value $|\ndot|$ is non-trivial, there exists a non-zero element $a\in k$ such that $|a|<1$. For sufficiently large integer $n\in\mathbb N_{>0}$, one has $\|a^ne_i\|<\epsilon$ for any $i\in\{1,\ldots,r\}$. Hence $V_{<\epsilon}$ contains a basis of the vector space $V$.
\end{proof}

The following proposition shows that each lattice defines a norm on the underlying vector space.

\begin{prop}\label{Pro:norminduiteparreseau}
Let $V$ be a finite-dimensional vector space over $k$ and $\mathcal V$ be a lattice of $V$. The map $\norm{\ndot}_{\mathcal V}$ is an ultrametric norm on $V$. Moreover, $\mathcal V$ is contained in the unit ball of $(V,\norm{\ndot}_{\mathcal V})$.
\end{prop}
\begin{proof} In the case where the absolute value $|\ndot|$ is trivial, one has $\mathcal V=V$ and the function $\norm{\ndot}_{\mathcal V}$ takes value $1$ on $V\setminus\{0\}$ and vanishes on $\{0\}$. The result is clearly true in this case. In the following, we assume that $|\ndot|$ is non-trivial.
For any $x\in V$, let $A_x$ be the set of all $a\in k^\times$ such that $a^{-1}x\in\mathcal V$. We claim that $A_x$ is non-empty and hence $\|x\|_{\mathcal V}$ is finite. Let $\{e_i\}_{i=1}^r$ be a subset of $\mathcal V$ which forms a basis of $V$ over $k$. We write $x$ in the form $x=a_1e_1+\cdots+a_re_r$ with $(a_1,\ldots,a_r)\in k^r$. Since $k$ is the fraction field of $\mathfrak o_k$, there exists $b\in k^{\times}$ such that $ba_1,\ldots,ba_r$ are all in $\mathfrak o_k$. Thus $bx\in \mathcal V$ and hence $b^{-1}\in A_x$. Therefore $\norm{\ndot}_{\mathcal V}$ is a map from $V$ to $\mathbb R_{\geqslant 0}$.

Let $x$ be an element of $V$ and $a\in k^{\times}$. The map $b\mapsto ab$ defines a bijection between $A_x$ and $A_{ax}$. Hence one has $\|ax\|_{\mathcal V}=|a|\cdot\|x\|_{\mathcal V}$.

Let $x$ and $y$ be elements of $V$, $a\in A_x$ and $b\in A_y$. One has $\{a^{-1}x,b^{-1}y\}\subseteq\mathcal V$. Note that
\begin{gather*}a^{-1}(x+y)=a^{-1}x+a^{-1}y=a^{-1}x+(a^{-1}b)(b^{-1}y),\\
b^{-1}(x+y)=b^{-1}x+b^{-1}y=(b^{-1}a)(a^{-1}x)+b^{-1}y.
\end{gather*} 
Since $\mathfrak o_k$ is a valuation ring, either $b^{-1}a\in \mathfrak o_k$, or $a^{-1}b\in \mathfrak o_k$. Hence, either $a\in A_{x+y}$ or $b\in A_{x+y}$. Therefore $\|x+y\|_{\mathcal V}\leqslant\max\{\|x\|_{\mathcal V},\|y\|_{\mathcal V}\}$.

It remains to verify that, if $\|x\|_{\mathcal V}=0$ then $x=0$. Assume that there exists a non-zero element $x\in V$ such that $\|x\|_{\mathcal V}=0$. Then there exists a sequence $\{a_n\}_{n\in\mathbb N}$ in $A_x$ such that $\lim_{n\rightarrow+\infty}|a_n|=0$. However, one has $a_n^{-1}x\in\mathcal V$ for any $n\in\mathbb N$. This contradicts the assumption that $\mathcal V$ is bounded. 

If $x$ is an element in $\mathcal V$, then $1$ belongs to $A_x$. Hence $\|x\|_{\mathcal V}\leqslant 1$.
\end{proof}

\begin{defi}\label{Def:inducedbylattice}
Let $V$ be a finite-dimensional vector space over $k$ and $\mathcal V$ be a lattice of $V$. We call $\norm{\ndot}_{\mathcal V}$ the norm on $V$ \emph{induced by the lattice $\mathcal V$}\index{norm!--- induced by a lattice}.
\end{defi}

\begin{prop}\label{Pro:finitetype}
Let $V$ be a finite-dimensional vector space over $k$ and $r$ be its rank over $k$. Let $\mathcal V$ be a lattice of $V$. Assume that $\mathcal V$ is an $\mathfrak o_k$-module of finite type. Then it is a free $\mathfrak o_k$-module of rank $r$.
\end{prop}
\begin{proof}
Since $\mathcal V$ is a sub-$\mathfrak o_k$-module of $V$, it is torsion-free. By \cite{Bourbaki64} Chapter VI, \S4, no.6, Lemma 1, any torsion-free module of finite type over a valuation ring is free. Hence $\mathcal V$ is a free $\mathfrak o_k$-module. Finally, since $\mathcal V$ generates $V$ as a vector space over $k$, any basis of $\mathcal V$ over $\mathfrak o_k$ is also a basis of $V$ over $k$. Hence the rank of $\mathcal V$ over $\mathfrak o_k$ is $r$.
\end{proof}

\begin{defi}\label{Def:pure}
Let $(V,\norm{\ndot})$ be a finite-dimensional seminormed vector space over $k$. We define the \emph{default of purity} of $\norm{\ndot}$ as 
\[\dpur(\norm{\ndot}):=\sup_{x\in V\setminus N_{\norm{\ndot}}}\dist(\ln\norm{x},\ln|k^{\times}|),\]
with \[\dist(\ln\norm{x},\ln|k^{\times}|):=\inf\{|\ln\norm{x}-\ln|a||\,:\,a\in k^{\times}\}.\] We say that the seminorm $\norm{\ndot}$ is \emph{pure}\index{pure}\index{seminorm!pure ---} if $\dpur(\norm{\ndot})=0$, or equivalently, the image of $V\setminus N_{\norm{\ndot}}$ by $\norm{\ndot}$ is contained in the closure of $|k^{\times}|$ in $\mathbb R_{>0}$. By definition, if the absolute value $|\ndot|$ is not discrete, then any seminorm on $V$ is pure; if $|\ndot|$ is discrete, then a seminorm $\norm{\ndot}$ on $V$ is pure if and only if its image is contained in that of $|\ndot|$. In the case where $|\ndot|$ is discrete, Moreover, for any lattice $\mathcal V$ of $V$, the norm $\norm{\ndot}_{\mathcal V}$ is pure.
\end{defi}

In the following, we study the correspondance between {ultrametric} norms and lattices of a finite-dimensional vector space over $k$. Note that the behaviour depends much on the discreteness of the absolute value $|\ndot|$.

\begin{prop}\label{Pro:normetreausauxdisc}
Assume that the absolute value $|\ndot|$ is discrete.
\begin{enumerate}[label=\rm(\arabic*)]
\item\label{Item: lattice norm unit ball} For any lattice $\mathcal V$ of $V$, one has $(V,\norm{\ndot}_{\mathcal V})_{\leqslant 1}=\mathcal V$ (see Notation \ref{Not:ball}).
\item\label{Item: lattice is free} Any lattice $\mathcal V$ of $V$ is a free $\mathfrak o_k$-module of rank $\dim_k(V)$.
\item\label{Item: criterion of lattice norm} Assume in addition that the absolute value $|\ndot|$ is non-trivial. Let $\norm{\ndot}$ be an ultrametric norm on $V$ and let $\mathcal V=(V,\norm{\ndot})_{\leqslant 1}$. Then one has $\norm{\ndot}\leqslant\norm{\ndot}_{\mathcal V}\leqslant |\varpi|^{-1}\norm{\ndot}$, where $\varpi$ is a uniformizing parameter of $k$. In particular, the default of purity of $\norm{\ndot}$ is bounded from above by {$-\ln|\varpi|$}. Moreover, if the norm $\norm{\ndot}$ is pure, then $\norm{\ndot}_{\mathcal V}=\norm{\ndot}$.
\end{enumerate} 
\end{prop}
\begin{proof}
\ref{Item: lattice norm unit ball} By Proposition \ref{Pro:norminduiteparreseau}, one has $\mathcal V \subseteq (V, \norm{\ndot}_{\mathcal V})_{\leqslant 1}$. 
Let $x$ be an element of $V$ such that $\| x \|_{\mathcal V} \leqslant 1$. In order to see that
$x \in \mathcal V$, we may assume that $x \not= 0$.
There is a sequence $\{\alpha_n\}_{n\in\mathbb N}$ in $k^{\times}$ such that $\alpha_n^{-1} x \in \mathcal V$ and
$\lim_{n\to\infty} |\alpha_n|= \| x \|_{\mathcal V}$. As $|\ndot|$ is discrete, there is $n\in\mathbb N$ such that
$|\alpha_n| = \| x \|_{\mathcal V}$, so that $\alpha_n \in \mathfrak o_k$ because $\| x \|_{\mathcal V} \leqslant 1$.
Therefore, $x \in \alpha_n \mathcal V \subseteq \mathcal V$, and hence $(V, \norm{\ndot}_{\mathcal V})_{\leqslant1} \subseteq \mathcal V$.

\ref{Item: lattice is free} Let $\{e_i\}_{i=1}^r$ be a basis of $V$ over $k$. We equip $V$ with the norm $\norm{\ndot}$ such that \[\|\lambda_1e_1+\cdots+\lambda_re_r\|=\max\{|\lambda_1|,\ldots,|\lambda_r|\}\] for any $(\lambda_1,\ldots,\lambda_r)\in k^r$. For any $\epsilon>0$, the ball \[(V,\norm{\ndot})_{\leqslant\epsilon}=\{a\in k\,:\,|a|\leqslant\epsilon\}^r\] is a free $\mathfrak o_k$-module of rank $r$ since $\mathfrak o_k$ is a principal ideal domain. Let $\mathcal V$ be a lattice. Since it is bounded, it is contained in certain ball $(V,\norm{\ndot})_{\leqslant\epsilon}$. Thus $\mathcal V$ is an $\mathfrak o_k$-module of finite type, and hence a free $\mathfrak o_k$-module of rank $r$ by Proposition \ref{Pro:finitetype}.

\ref{Item: criterion of lattice norm} By the definition of the uniformizing element, one has $|k^{\times}|=\{|\varpi|^{n}\,:\,n\in\mathbb Z\}$. If $x$ is a non-zero element in $V$ and if $A_x$ is the set of all $a\in k^{\times}$ such that \[\|a^{-1}x\|=|a|^{-1}\cdot\|x\|\leqslant 1,\] then one has
\[\{|a|\,:\,a\in A_x\}=\{|\varpi|^n\,:\,n\in\mathbb Z,\;|\varpi|^n\geqslant\|x\|\}.\]  
Since $\|x\|_{\mathcal V}=\inf\{|a|\,:\,a\in A_x\}$, one has $\|x\|_{\mathcal V}\geqslant\|x\|>|\varpi|\cdot\|x\|_{\mathcal V}$. Combined with the  fact that the norm $\norm{\ndot}_{\mathcal V}$ is pure, this implies the inequality {$\dpur(\norm{\ndot})\leqslant -\ln|\varpi|$}. If in addition the norm $\norm{\ndot}$ is pure, $\|x\|_{\mathcal V}$ belongs to $\{|\varpi|^n\,:\,n\in\mathbb Z\}$. Hence $\|x\|_{\mathcal V}=\|x\|$.
\end{proof}

\begin{rema}\label{Rem:reseaunor}
Let $V$ be a finite-dimensional vector space over $k$. We denote by $\mathrm{Lat}(V)$ the set of all lattices of $V$, and by $\mathrm{Nor}(V)$ that of all ultrametric norms on $V$. The correspondance $(\mathcal V\in\mathrm{Lat}(V))\mapsto\norm{\ndot}_{\mathcal V}$ defines a map from $\mathrm{Lat}(V)$ to $\mathrm{Nor}(V)$. Proposition \ref{Pro:normetreausauxdisc} shows that, if the absolute value $|\ndot|$ is discrete, then {this} map is injective, and its image is precisely the set of all pure ultrametric norms.
\end{rema}

\begin{prop}\label{Pro:notdiscret}
Assume that the absolute value $|\ndot|$ is not discrete.
\begin{enumerate}[label=\rm(\arabic*)]
\item\label{Item: not discrte lattice norm1} For any lattice $\mathcal V$ of $V$ one has 
$(V,\norm{\ndot}_{\mathcal V})_{<1}\subseteq\mathcal V\subseteq(V,\norm{\ndot}_{\mathcal V})_{\leqslant 1}$.
If in addition there exists {an} ultrametric norm $\norm{\ndot}$ on $V$ such that $\mathcal V=(V,\norm{\ndot})_{\leqslant 1}$, then one has $\mathcal V=(V,\norm{\ndot}_{\mathcal V})_{\leqslant 1}$.
\item\label{Item: not discrte lattice norm2} Let $\norm{\ndot}$ be an ultrametric norm on $V$ and $\mathcal V=(V,\norm{\ndot})_{\leqslant 1}$. Then $\norm{\ndot}=\norm{\ndot}_{\mathcal V}$.
\end{enumerate}
\end{prop}
\begin{proof}
\ref{Item: not discrte lattice norm1} If $x$ is an element of $\mathcal V$, by the relation $1x=x\in\mathcal V$ we obtain that $\|x\|_{\mathcal V}\leqslant 1$. Hence $\mathcal V\subseteq (V,\norm{\ndot}_{\mathcal V})_{\leqslant 1}$. In the following, we prove the inclusion  relation $(V,\norm{\ndot}_{\mathcal V})_{<1}\subseteq\mathcal V$. Let $x$ be an element in $V$ such that $\|x\|_{\mathcal V}<1$. By definition there exists $a\in k^{\times}$, $|a|<1$, such that $a^{-1}x\in\mathcal V$. Since $|a|<1$ one has $a\in \mathfrak o_k$. Therefore $x=a(a^{-1}x)\in\mathcal V$.

The second assertion of \ref{Item: not discrte lattice norm1} is a direct consequence of \ref{Item: not discrte lattice norm2}. In the following, we prove the statement \ref{Item: not discrte lattice norm2}. 
Let $x$ be an element of $V$ and 
\[
A_x=\{a\in k^{\times}\,:\,a^{-1}x\in\mathcal V\}=\{a\in k^{\times}\,:\,\|x\|\leqslant|a|\}.
\]
Since the image of $|\ndot|$ is dense in $\mathbb R_+$, one has 
$\|x\|_{\mathcal V}=\inf\{|a|\,:\,a\in A_x\}=\|x\|$.
Hence $\norm{\ndot}_{\mathcal V}=\norm{\ndot}$.
\end{proof}

\begin{rema}
Let $V$ be a finite-dimensional vector space over $k$. Proposition \ref{Pro:notdiscret} shows that, if the absolute value $|\ndot|$ is not discrete, the map $\mathrm{Lat}(V)\rightarrow\mathrm{Nor}(V)$, sending any lattice $\mathcal V$ of $V$ to the norm $\norm{\ndot}_{\mathcal V}$, is surjective (compare with Remark \ref{Rem:reseaunor}). 
\end{rema}

\begin{prop}\label{coro:dual:lattice:norm}
Let $V$ be a finite-dimensional vector space over $k$ and $\mathcal V$ be a lattice of $V$. Let $\mathcal V^\vee=\Hom_{\mathfrak o_k}(\mathcal V,\mathfrak o_k)$ be the dual $\mathfrak o_k$-module of $\mathcal V$. Then one has $\norm{\ndot}_{\mathcal V,*}=\norm{\ndot}_{\mathcal V^\vee}$ on $V^\vee$.
\end{prop}
\begin{proof}
Let $f$ be a non-zero element of $V^\vee$. Assume that $a$ is an element of $k^{\times}$ such that $a^{-1}f\in\mathcal V^\vee$. Then for any $x\in V$ and any $b\in k^{\times}$ such that $b^{-1}x\in\mathcal V$ one has $a^{-1}f(b^{-1}x)=(ab)^{-1}f(x)\in\mathfrak o_k$ and hence $|b|\geqslant |f(x)|/|a|$. Since $b$ is arbitrary one has $\norm{x}_{\mathcal V}\geqslant|f(x)|/|a|$ for any $x\in V$ and hence $|a|\geqslant\norm{f}_{\mathcal V,*}$. Since $a$ is arbitrary we obtain $\norm{f}_{\mathcal V^\vee}\geqslant\norm{f}_{\mathcal V,*}$.

Conversely, suppose that the operator norm of a non-zero linear form $f:V\rightarrow k$ is bounded from above by $1$, where we consider the norm $\norm{\ndot}_{\mathcal V}$ on $V$. Then for any $x\in\mathcal V$ one has $|f(x)|\leqslant\|x\|_{\mathcal V}\leqslant 1$ and hence $f(x)\in\mathfrak o_k$. This shows that $f\in\mathcal V^\vee$ and hence $\|f\|_{\mathcal V^\vee}\leqslant 1$. Therefore the unit ball of $\norm{\ndot}_{\mathcal V^\vee}$ contains that of $\norm{\ndot}_{\mathcal V,*}$. Moreover, since the norm $\norm{\ndot}_{\mathcal V}$ is pure, also is its dual norm $\norm{\ndot}_{\mathcal V,*}$. Therefore, the norm $\norm{\ndot}_{\mathcal V,*}$ coincides with the norm induced by its unit ball (see Propositions \ref{Pro:normetreausauxdisc} and \ref{Pro:notdiscret}). Therefore, $\norm{\ndot}_{\mathcal V^\vee}\leqslant\norm{\ndot}_{\mathcal V,*}$. The proposition is thus proved.
\end{proof}

\subsection{Trivial valuation case}\label{Subsection: Trivial valuation}

In this subsection, we study ultrametrically normed vector spaces over a trivially valued field. We fix a field $k$ equipped with the trivial absolute value $|\ndot|$. If $V$ is a vector space over $k$, we denote by $\Theta(V)$ the set of all non-zero vector subspaces of $V$. The set $\Theta(V)$ is equipped with the partial order of inclusion. If $\norm{\ndot}$ is an ultrametric norm on $V$, we denote by $\Psi(V,\norm{\ndot})$ the set of closed balls of $V$ (centered at the origin) which do not reduce to one point, namely (see Notation \ref{Not:ball})
\[\Psi(V,\norm{\ndot})=\Big\{(V,\norm{\ndot})_{\leqslant r}\,:\,r>0,\;(V,\norm{\ndot})_{\leqslant r}\neq\{0\}\Big\}.\]

\begin{prop}\label{Pro:characterisation of filtration and nomr}
Let $V$ be a finite-dimensional vector space equipped with an ultrametric norm $\norm{\ndot}$.
The set $\Psi(V,\norm{\ndot})$ is a totally ordered subset of $\Theta(V)$, whose cardinal does not exceed the rank of $V$ over $k$.
\end{prop}
\begin{proof}
By definition the set $\Psi(V,\norm{\ndot})$ is totally ordered with respect to the partial order of inclusion. In the following, we show that any element $W\in\Psi(V,\norm{\ndot})$ is a vector subspace of $V$ and hence belongs to $\Theta(V)$. Assume that $W=(V,\norm{\ndot})_{\leqslant r}$ with $r>0$. Since the absolute value on $k$ is trivial, for any $x\in W$ and any $a\in k$ one has $\|ax\|\leqslant\|x\|\leqslant r $. Moreover, since the norm $\norm{\ndot}$ is ultrametric, $W$ is stable by addition. Hence $\Psi(V,\norm{\ndot})$ is a totally ordered subset of $\Theta(V)$. In particular, the function $\rang_k(\ndot):\Psi(V,\norm{\ndot})\rightarrow\mathbb N_{\geqslant 1}$ is injective, which is bounded from above by $\rang_k(V)$. Therefore the cardinal of $\Psi(V,\norm{\ndot})$ does not exceed $\rang_k(V)$.
\end{proof}

The above proposition shows that the set $\Psi(V,\norm{\ndot})$ actually forms an increasing flag of non-zero vector subspaces of $V$. For any $W\in\Psi(V,\norm{\ndot})$, let \[\varphi_{\norm{\ndot}}(W):=\sup\{t\in\mathbb R\,:\,W\subseteq(V,\norm{\ndot})_{\leqslant \mathrm{e}^{-t}}\}.\]
Then $\varphi_{\norm{\ndot}}$ is a strictly decreasing function on $\Psi(V,\norm{\ndot})$ in the sense that, if $W_1$ and $W_2$ are two elements of $\Psi(V,\norm{\ndot})$ such that $W_1\subsetneq W_2$, then one has $\varphi_{\norm{\ndot}}(W_1)>\varphi_{\norm{\ndot}}(W_2)$. The following proposition shows that the norm $\norm{\ndot}$ is completely determined by the increasing flag $\Psi(V,\norm{\ndot})$ and the function $\varphi_{\norm{\ndot}}$. 

\begin{prop}\label{Pro: charageruasiont of filtration by flage}
Let $\Psi$ be a totally ordered subset of $\Theta(V)$ and $\varphi:(\Psi,\supseteq)\rightarrow (\mathbb R,\leqslant)$ be a function which preserves strictly the orders, that is, for any $(W_1, W_2) \in \Psi^2$ with $W_1 \subsetneq W_2$, one has $\varphi(W_1) > \varphi(W_2)$. Then there exists a unique ultrametric norm $\norm{\ndot}$ on $V$ such that $\Psi(V,\norm{\ndot})=\Psi$ and $\varphi_{\norm{\ndot}}=\varphi$.
\end{prop}
\begin{proof}
We write $\Psi$ in the form of an increasing flag $V_1\subsetneq \ldots\subsetneq V_n$. For $i\in\{1,\ldots, n\}$, let $a_i=\varphi(V_i)$. Since $\varphi$ preserves strictly the orders, one has $a_1>\ldots>a_n$. Let $\boldsymbol{e}=\{e_j\}_{j=1}^m$ be a basis of $V$ which is compatible with the flag $\Psi$ (namely $\card(\mathbf{e}\cap V_i)=\mathrm{rk}_k(V_i)$
for any $i\in\{1,\ldots,n\}$). For any $j\in\{1,\ldots,m\}$, there exists a unique $i\in\{1,\ldots,n\}$ such that $e_j\in V_i\setminus V_{i-1}$ (where $V_0=\{0\}$ by convention) and we let $r_j=\mathrm{e}^{-a_i}$. Let $\norm{\ndot}$ be the ultrametric norm on $V$ defined as  
\[\forall\,(\lambda_1,\ldots,\lambda_m)\in k^m,\quad\|\lambda_1e_1+\cdots+\lambda_me_m\|=\max_{\begin{subarray}{c}
j\in\{1,\ldots,m\}\\
\lambda_j\neq 0
\end{subarray}}r_j.\]
Note that for $r\geqslant 0$ the ball $(V,\norm{\ndot})_{\leqslant r}$ identifies with the vector subspace generated by those $e_j$ with $r_j\leqslant r$. Hence one has $\Psi(V,\norm{\ndot})=\Psi$. Moreover, for any $i\in\{1,\ldots,n\}$ one has
\[\varphi_{\norm{\ndot}}(V_i)=\sup\{t\in\mathbb R\,:\,V_i\subseteq(V,\norm{\ndot})_{\leqslant \mathrm{e}^{-t}}\}=a_i=\varphi(V_i).\]

Let $\norm{\ndot}'$ be another ultrametric norm on $V$ verifying the relations $\Psi(V,\norm{\ndot}')=\Psi$ and $\varphi_{\norm{\ndot}'}=\varphi$. For any $r\geqslant 0$, $(V,\norm{\ndot}')_{\leqslant r}=V_i$ if and only if $r\in \intervalle{[}{\mathrm{e}^{-a_{i}}}{\mathrm{e}^{-a_{i+1}}}{[}$, with the convention $a_0=+\infty$ and $a_{n+1}=-\infty$. Therefore one has $(V,\norm{\ndot})_{\leqslant r}=(V,\norm{\ndot}')_{\leqslant r}$ for any $r\geqslant 0$, which leads to $\norm{\ndot}=\norm{\ndot}'$. 
\end{proof}

\begin{defi}
Let $V$ be a finite-dimensional vector space over $k$. 
A family $\mathcal F=\{\mathcal F^t(V)\}_{t\in\mathbb R}$ of vector subspaces of $V$ parametrised by $\mathbb R$ is called
an \emph{$\mathbb R$-filtration}\index{R-filtration@$\mathbb R$-filtration} of $V$ if
it is separated ($\mathcal F^t(V)=\{0\}$ for sufficiently positive $t$), exhaustive ($\mathcal F^t(V)=V$ for sufficiently negative $t$) and left-continuous (the function $(t\in\mathbb R)\rightarrow\rang_k(\mathcal F^t(V))$ is left-continuous).
\end{defi}

\begin{defi}
Let $V$ be a finite-dimensional vector space over $k$ and $\mathcal F$ be an $\mathbb R$-filtration on $V$. Let $r$ be the rank of $V$ over $k$. We define a map $Z_{\mathcal F}:\{1,\ldots,r\}\rightarrow\mathbb R$ as follows:
\[\forall\,i\in\{1,\ldots,r\},\quad Z_{\mathcal F}(i):=\sup\{t\in\mathbb R\,:\,\rang_k(\mathcal F^t(V))\geqslant i\}.\]
By definition, for any $t\in\mathbb R$ and any $i\in\{1,\ldots,r\}$ one has
\begin{equation}
\label{Equ: critere de fonction Z}
Z_{\mathcal F}(i)\geqslant t\Longleftrightarrow\rang_k(\mathcal F^t(V))\geqslant i.
\end{equation}
\end{defi}

\begin{prop}
\label{Pro: comparaison des minima}Let $V$ be a finite-dimensional non-zero vector space over $k$ and $\mathcal F$ and $\mathcal G$ be $\mathbb R$-filtrations on $V$. Let $a\in\mathbb R$ such that, for any $t\in\mathbb R$ one has $\mathcal F^t(V)\subseteq\mathcal G^{t-a}(V)$. Then one has $Z_{\mathcal F}(i)\leqslant Z_{\mathcal G}(i)+a$ for any $i\in\{1,\ldots,\rang_k(V)\}$. 
\end{prop}
\begin{proof}
By the relation \eqref{Equ: critere de fonction Z}, for any $i\in\{1,\ldots,\rang_k(V)\}$, if $Z_{\mathcal F}(i)\geqslant t$, then $\rang_k(\mathcal F^t(V))\geqslant i$, which implies that $\rang_k(\mathcal G^{t-a}(V))\geqslant i$ and hence (still by the relation \eqref{Equ: critere de fonction Z}) $Z_{\mathcal G}(i)\geqslant t-a$. Therefore we obtain $Z_{\mathcal F}(i)-a\leqslant Z_{\mathcal G}(i)$.
\end{proof}

\begin{rema}\label{Rem: R-filtration as flag plus slopes} Let $V$ be a finite-dimensional vector space over $k$.
There are canonical bijections between the following three sets:
\begin{enumerate}[label=(\Alph*)]
\item\label{Item: flag+jumps}
the set of all pairs
\[
\big( 0=V_0\subsetneq V_1\subsetneq\ldots\subsetneq V_n=V,\ \mu_1>\ldots>\mu_n \big)
\]
such that $0=V_0\subsetneq V_1\subsetneq\ldots\subsetneq V_n=V$ is an increasing sequence of {vector} subspaces of $V$ and $\mu_1>\ldots>\mu_n$ is a decreasing sequence of real numbers.

\item\label{Item: R filtrations}
the set of all $\mathbb R$-filtrations $\mathcal F$ of $V$.

\item\label{Item: ultrametric norms}
the set of all ultrametric norms $\|\ndot\|$ of $V$ over $k$.
\end{enumerate}

In the following, we explain the construction of these canonical maps.

\medskip\noindent
$\bullet$ \ref{Item: flag+jumps}$\to $\ref{Item: R filtrations}:
The associated $\mathbb R$-filtration $\mathcal F$ on $V$ with the data
$\big( V_0\subsetneq\ldots\subsetneq V_n,\ \mu_1>\ldots>\mu_n \big)$
is defined by
$\mathcal F^t(V) := V_i$ if $t \in\intervalle{]}{\mu_{i+1}}{\mu_i}{]} \cap \mathbb R$, where $\mu_0 =+ \infty$ and $\mu_{n+1} = -\infty$ by convention.

\begin{center}
\begin{tikzpicture}[domain=0:6, samples=100, very thick]
\draw[thin] (0,0)--(6.92,0);
\draw[thin, ->] (7.08,0)--(9,0) node[right] {$t$};
\draw[very thick] (0,3.5)--(1,3.5) ; \node at(1.0,-0.2){{\small $\mu_{n}$}};
\draw[very thick] (1.08,3.0)--(2,3.0) ; \node at(2.0,-0.2){{\small $\mu_{n-1}$}};
\draw[very thick] (2.08,2.5)--(3,2.5) ; \node at(3.0,-0.2){{\small $\mu_{n-3}$}};
\draw[very thick] (5.08,1.0)--(6,1.0) ; \node at(5.0,-0.2){{\small $\mu_{3}$}};
\draw[very thick] (6.08,0.5)--(7,0.5) ; \node at(6.0,-0.2){{\small $\mu_{2}$}};
\draw[very thick] (7.08,0)--(9,0) ; \node at(7.0,-0.2){{\small $\mu_{1}$}};
\draw[very thick] (4.5,1.5)--(5.0,1.5) ;
\draw[very thick] (3.08,2.0)--(3.5,2.0) ;
\node at(4.0,-0.3){{\small $\cdots$}};
\node at(0.5,3.7){{\small $V_n$}};
\node at(1.5,3.2){{\small $V_{n-1}$}};
\node at(2.5,2.7){{\small $V_{n-2}$}};
\node at(5.5,1.2){{\small $V_{2}$}};
\node at(6.5,0.7){{\small $V_{1}$}};
\node at(7.5,0.2){{\small $V_{0}$}};
\node at(4.0,1.77){{\huge $\ddots$}};
\fill (1.0,3.5) circle (2.3pt);
\fill (2.0,3.0) circle (2.3pt);
\fill (3.0,2.5) circle (2.3pt);
\fill (5.0,1.5) circle (2.3pt);
\fill (6.0,1.0) circle (2.3pt);
\fill (7.0,0.5) circle (2.3pt);
\draw[thick] (1.0,3.0) circle (2pt);
\draw[thick] (2.0,2.5) circle (2pt);
\draw[thick] (3.0,2.0) circle (2pt);
\draw[thick] (5.0,1.0) circle (2pt);
\draw[thick] (6.0,0.5) circle (2pt);
\draw[thick] (7.0,0.0) circle (2pt);
\draw [very thin, dashed](1,3.5)--(1,0);
\draw [very thin, dashed](2,3.0)--(2,0);
\draw [very thin, dashed](3,2.5)--(3,0);
\draw [very thin, dashed](5,1.5)--(5,0);
\draw [very thin, dashed](6,1.0)--(6,0);
\draw [very thin, dashed](7,0.5)--(7,0);
\end{tikzpicture}
\end{center}

\medskip\noindent
$\bullet$ \ref{Item: R filtrations}$\to$\ref{Item: flag+jumps}:
One has a sequence $0=V_0\subsetneq V_1\subsetneq\ldots\subsetneq V_n=V$ such that
$\left\{ \mathcal F^t(V) \,:\, t \in \mathbb R \right\} = \{ V_0, V_1, \ldots, V_n \}$.
A sequence $\mu_1>\ldots>\mu_n$ in $\mathbb R$ is given by
$\mu_i = \sup \{ t \,:\, \mathcal F^t(V) = V_i \}$ for $i\in\{1, \ldots, n\}$.

\medskip\noindent
$\bullet$ \ref{Item: flag+jumps}$\to$\ref{Item: ultrametric norms}:
The corresponding norm $\|\ndot\|$ to the data
$\big( V_0\subsetneq\ldots\subsetneq V_n,\ \mu_1>\ldots>\mu_n \big)$ is given by 
\[
\| x \| = \begin{cases}
\mathrm{e}^{-\mu_i} & \text{if $x \in V_i \setminus V_{i-1}$}, \\
0 & \text{if $x = 0$}.
\end{cases}
\]

\medskip\noindent
$\bullet$ \ref{Item: ultrametric norms}$\to$\ref{Item: flag+jumps}:
By Proposition~\ref{Pro:characterisation of filtration and nomr},
there is an increasing sequence 
\[
0=V_0\subsetneq V_1\subsetneq\ldots\subsetneq V_n=V
\]
of subspaces of $V$
such that $\Psi(V,\norm{\ndot}) = \{ V_1, \ldots, V_n \}$.
A decreasing sequence of real numbers is given by
$\mu_i = \varphi_{\norm{\ndot}}(V_i)$ for $i\in\{1, \ldots, n\}$.

\medskip\noindent
$\bullet$ \ref{Item: R filtrations}$\to$\ref{Item: ultrametric norms}:
We define a function $\lambda_{\mathcal F}:V\rightarrow\mathbb R\cup\{+\infty\}$ such that
\[\forall\,x\in V,\quad \lambda_{\mathcal F}(x):=\sup\{t\in\mathbb R\,:\,x\in\mathcal F^t(V)\}.\]
Then the ultrametric norm $\norm{\ndot}$ on $V$ corresponding to $\mathcal F$ is given by 
\[\forall\,x\in V,\quad \|x\|=\mathrm{e}^{-\lambda_{\mathcal F}(x)}.\]

\medskip\noindent
$\bullet$ \ref{Item: ultrametric norms}$\to $\ref{Item: R filtrations}:
The corresponding filtration $\mathcal F$ to the norm $\|\ndot\|$ is given by $\mathcal F^t(V) = (V,\norm{\ndot})_{\leqslant\mathrm{e}^{-t}}$.

\medskip
Let  
$\mathcal F$ be an $\mathbb R$-filtration on $V$, which corresponds to an increasing flag
$0=V_0\subsetneq V_1\subsetneq\ldots\subsetneq V_n=V$ together with
a decreasing sequence
$\mu_1>\ldots>\mu_n$ of real numbers.
Note that the sets $\{\mu_1,\ldots,\mu_n\}$ and $\{Z_{\mathcal F}(1),\ldots,Z_{\mathcal F}(r)\}$ are actually equal, where $r$ denotes the rank of $V$ {over $k$}. Moreover, the value $\mu_i$ appears exactly $\rang_k(V_i/V_{i-1})$ times in the sequence $Z_{\mathcal F}(1),\ldots,Z_{\mathcal F}(r)$.
\end{rema}

\subsection{Metric on the space of norms}\label{Subsec:Distance}
Let $V$ be a finite-dimensional vector space over $k$. We denote by $\mathcal N_V$ the set of all norms on $V$. If $\norm{\ndot}_1$ and $\norm{\ndot}_2$ are norms on $V$, by Proposition~\ref{Pro:topologicalnormedspace}  
we obtain that 
\[\sup_{s\in V\setminus\{0\}}\Big|\ln\|s\|_1-\ln\|s\|_2\Big|\]
is finite. We denote by $d(\norm{\ndot}_1,\norm{\ndot}_2)$ this number, called the \emph{distance}\index{distance}\index{norm!distance between two norms} between $\norm{\ndot}_1$ and $\norm{\ndot}_2$. 
It is easy to see that the function $d : \mathcal N_V \times \mathcal N_V \to \mathbb R_{\geqslant 0}$ satisfies the axioms of metric.

\begin{rema}\label{Rem:localcompleteness}
Let $V$ be a finite-dimensional vector space over $k$ and $\norm{\ndot}_0$ be a norm on $V$, which is the trivial norm if the absolute value $|\ndot|$ is trivial (namely $\|x\|_0=1$ for any $x\in V\setminus\{0\}$). Let $\lambda$ be a real number in $\intervalle{]}01[$. If the absolute value $|\ndot|$ is non-trivial, we require in addition that $\lambda<\sup\{|a|\,:\,a\in k^{\times},\;|a|<1\}$.
We denote by $C$ the  annulus $\{x\in V\,:\,\lambda\leqslant\|x\|_0\leqslant 1\}$. Note that one has $C=V\setminus\{0\}$ when $|\ndot|$ is trivial.
For any norm $\norm{\ndot}$ on $V$, the restriction of the function $\ln\norm{\ndot}$ {to} $C$ is bounded, and the norm $\norm{\ndot}$ is uniquely determined by its restriction {to} $C$ (this is a consequence of Proposition \ref{Pro:dilatation} when $|\ndot|$ is non-trivial). Thus we can identify $\mathcal N_V$ with a closed subset of $\mathcal C_b(C)$, the space of bounded and continuous functions on $C$ equipped with the sup norm. In particular, $\mathcal N_V$ is a complete metric space.
\end{rema}

\begin{prop}\label{Pro:distancequotientandsub}
Let $V$ be a finite-dimensional vector space over $k$, and $\norm{\ndot}_1$ and $\norm{\ndot}_2$ be norms on $V$.
\begin{enumerate}[label=\rm(\arabic*)]
\item\label{Item: induced norm metric} Let $U$ be a vector subspace of $V$, $\norm{\ndot}_{U,1}$ and $\norm{\ndot}_{U,2}$ be the restrictions of $\norm{\ndot}_1$ and $\norm{\ndot}_2$ {to} $U$, respectively. Then one has $d(\norm{\ndot}_{U,1},\norm{\ndot}_{U,2})\leqslant d(\norm{\ndot}_{1},\norm{\ndot}_{2})$.
\item\label{Item: quotient norm metric} Let $W$ be a quotient vector space of $V$, $\norm{\ndot}_{W,1}$ and $\norm{\ndot}_{W,2}$ be quotient norms of $\norm{\ndot}_1$ and $\norm{\ndot}_2$ on $W$, respectively. Then one has $d(\norm{\ndot}_{W,1},\norm{\ndot}_{W,2})\leqslant d(\norm{\ndot}_{1},\norm{\ndot}_{2})$.
\end{enumerate}
\end{prop}
\begin{proof}
\ref{Item: induced norm metric} follows directly from the definition of the distance function.

\ref{Item: quotient norm metric} It is sufficient to show that
\[\Big|\ln\|x\|_{W,1}-\ln\|x\|_{W,2}\Big| \leqslant d(\norm{\ndot}_1, \norm{\ndot}_2).\]
for $x \in W \setminus \{ 0 \}$. Clearly we may assume that $\|x\|_{W,1} > \|x\|_{W,2}$.
For $\epsilon > 0$, one can choose $s \in V$ such that $[s] = x$ and $\|s \|_2 \leqslant \mathrm{e}^{\epsilon} \|x\|_{2, W}$. Then
\begin{align*}
0 & < \ln\|x\|_{W,1}-\ln\|x\|_{W,2} \leqslant \ln \|s\|_1 - \ln \big(\mathrm{e}^{-\epsilon} \|s \|_2\big) \\
& = (\ln \|s\|_1 - \ln \|s \|_2) + \epsilon \leqslant d(\norm{\ndot}_1, \norm{\ndot}_2) + \epsilon,
\end{align*}
as desired.
\end{proof}

\begin{prop}\label{Pro:distanceofoperatornorms}
Let $V$ and $W$ be finite-dimensional vector spaces over $k$, $\norm{\ndot}_{V,1}$ and $\norm{\ndot}_{V,2}$ be norms on $V$, and $\norm{\ndot}$ be a norm on $W$. Let $\norm{\ndot}_1$ and $\norm{\ndot}_2$ be the operator norms on $\mathscr L(V,W)$, where we consider the norm $\norm{\ndot}_W$ on $W$, and the norms $\norm{\ndot}_{V,1}$ and $\norm{\ndot}_{V,2}$ on $V$, respectively. Then one has
\begin{equation}\label{Equ:distanceopemaj}d(\norm{\ndot}_1,\norm{\ndot}_2)\leqslant d(\norm{\ndot}_{V,1},\norm{\ndot}_{V,2}).\end{equation}
In particular, one has
\begin{equation}\label{Equ:distancddualmaj}d(\norm{\ndot}_{V,1,*},\norm{\ndot}_{V,2,*})\leqslant d(\norm{\ndot}_{V,1}.\norm{\ndot}_{V,2}),\end{equation}
Moreover, the equality in \eqref{Equ:distancddualmaj} holds when both norms $\norm{\ndot}_{V,1}$ and $\norm{\ndot}_{V,2}$ are reflexive.
\end{prop}
\begin{proof}
For \eqref{Equ:distanceopemaj}, it is sufficient to show
\[\Big|\ln\|f\|_1-\ln\|f\|_2\Big|\leqslant  d(\norm{\ndot}_{V,1},\norm{\ndot}_{V,2})\]
for $f \in {\mathscr L}(V,W) \setminus \{ 0 \}$. Clearly we may assume that $\|f\|_1 > \|f\|_2$.
By definition one has
\[\|f\|_1=\sup_{x\in V\setminus\{0\}}\frac{\|f(x)\|_{W}}{\|x\|_{V,1}}\quad\text{ and }\quad\|f\|_2=\sup_{x\in V\setminus\{0\}}\frac{\|f(x)\|_{W}}{\|x\|_{V,2}},\]
so that, for $\epsilon > 0$, one can find $x \in V \setminus \{ 0 \}$ such that $\mathrm{e}^{-\epsilon} \|f\|_1 \leqslant \|f(x)\|_{W}/\|x\|_{V,1}$.
Therefore,
\begin{align*}
0 & < \ln\|f\|_{1}-\ln\|f\|_{2} \leqslant \ln \left( \mathrm{e}^{\epsilon} \frac{ \|f(x)\|_{W}}{\|x\|_{V,1}}\right)  - \ln \left( \frac{\|f(x)\|_{W}}{\|x\|_{V,2}}\right) \\
& = (\ln \|x\|_{V,2} - \ln \|x \|_{V,1}) + \epsilon \leqslant d(\norm{\ndot}_{V,1}, \norm{\ndot}_{V,2}) + \epsilon,
\end{align*}
as desired.

In order to obtain \eqref{Equ:distancddualmaj}, it suffices to apply \eqref{Equ:distanceopemaj} to the case where $(W,\norm{\ndot})=(k,|\ndot|)$. If in addition both norms $\norm{\ndot}_{V,1}$ and $\norm{\ndot}_{V,2}$ are reflexive, then one has
\[d(\norm{\ndot}_{V,1,*},\norm{\ndot}_{V,2,*})\geqslant d(\norm{\ndot}_{V,1,**},\norm{\ndot}_{V,2,**})=d(\norm{\ndot}_{V,1},\norm{\ndot}_{V,2}).\]
Hence the equality holds.
\end{proof}

\subsection{Direct sums}\label{Subsec:directsums}
Let $\mathscr S$ be the set of all convex and continuous functions $\psi:[0,1]\rightarrow [0,1]$ such that 
$\max\{t,1-t\}\leqslant \psi(t)$ for any $t\in[0,1]$.

\medskip

\begin{center}
\begin{tikzpicture}[domain=0:3, samples=100, very thick]
\draw[thick, ->] (-0.2,0)--(3.2,0) node[right] {$t$};
\draw[thick, ->] (0,-0.2)--(0,3.2);
\draw [very thin, dashed](0,3) node [left]{{\tiny $1$}}--(3,3);
\draw [very thin, dashed](3,0) node [below]{{\tiny $1$}}--(3,3);
\draw [very thin, dashed](0,1.5) node [left]{{\tiny ${\frac 12}$}}--(1.5,1.5) ;
\draw [very thin, dashed](1.5,0) node [below]{{\tiny ${\frac 12}$}}--(1.5,1.2);
\node at(-0.2,-0.2){{\tiny ${0}$}};
\draw [thick] (0,3)--(1.5,1.5);
\draw [thick] (1.5,1.5)--(3,3);
\node at(1.45,1.35){{\tiny $\max\{t,1-t\}$}};
\draw [thick] plot(\x, {3 * pow(pow(\x/3, 3/2) + pow(1 - \x/3, 3/2), 2/3)});
\node at(1.47,2.58) {{\tiny $\psi(t)$}};
\end{tikzpicture}
\end{center}

\medskip
\noindent Let $\norm{\ndot}$ be a norm on $\mathbb R^2$, where we consider the usual absolute value $|\ndot|_\infty$ on $\mathbb R$. We say that $\norm{\ndot}$ is an \emph{absolute normalised norm}\index{absolute normalised norm}\index{norm!absolute normalised ---} if $\|(1,0)\|=\|(0,1)\|=1$ and if
\[\forall\,(x,y)\in\mathbb R^2,\quad\|(x,y)\|=\|(|x|_\infty,|y|_\infty)\|.\]
By \cite[\S21, Lemma 3]{Bonsall_Duncan73}, the set of all absolute normalised norms on $\mathbb R^2$ can be parametrised by the functional space $\mathscr S$ (see \cite{Saito_Kato_Takahashi00} for the higher dimensional generalisation of this result). If $\norm{\ndot}$ is the absolute normalised norm corresponding to $\psi\in\mathscr S$, one has
\[\|(x,y)\|=(|x|+|y|)\psi\Big(\frac{|x|}{|x|+|y|}\Big).\]
In particular, one always has
\begin{equation}\label{Equ: minoration norme (x,y)}\|(x,y)\|\geqslant\max(|x|_\infty,|y|_\infty)\end{equation}
Conversely, given an absolute normalised norm $\norm{\ndot}$ on $\mathbb R^2$, the corresponding function in $\mathscr S$ is  
\begin{equation}\label{Equ:fonctioncorrespondsant}
(t\in[0,1])\longmapsto\|(t,1-t)\|.
\end{equation}

For example, the function $\psi(t)=\max\{t,1-t\}$, $t\in[0,1]$ corresponds to the norm $(x,y)\mapsto \max\{|x|_\infty,|y|_\infty\}$ on $\mathbb R^2$. If $p\geqslant 1$ is a real number, the function $\psi_p(t)=(t^p+(1-t)^p)^{1/p}$, $t\in[0,1]$ belongs to $\mathscr S$; it corresponds to the $\ell^p$-norm $(x,y)\mapsto (|x|^p+|y|^p)^{1/p}$. 

Given a function $\psi$ in $\mathscr S$, or equivalently an absolute normalised norm on $\mathbb R^2$, for any couple of finite-dimensional seminormed vector spaces over $k$, one can naturally attach to the direct sum of the vector spaces a direct sum seminorm, which depends on the function $\psi$.

\begin{lemm}\label{Lem:croissancef}
Let $a$, $b$, $a'$ and $b'$ be  real numbers such that $0\leqslant a\leqslant a'$ and $0\leqslant b\leqslant b'$. We assume in addition that $a+b>0$. If $\psi$ is a function in $\mathscr S$, then
\begin{equation}\label{Equ:inequalitypsi}(a+b)\psi\Big(\frac{a}{a+b}\Big)\leqslant(a'+b')\psi\Big(\frac{a'}{a'+b'}\Big).\end{equation}
\end{lemm}
\begin{proof}
For any $t\in[0,1]$, the value $\psi(t)$ is bounded from below by $t$. Moreover, one has $\psi(1)=1$. The function $t\mapsto \psi(t)/t$ on $\intervalle{]}{0}{1}{]}$ is non-increasing. In fact, for $0<s\leqslant t$, by the convexity of the function $\psi$ one has
\[\begin{split}\psi(t)&=\psi\Big(\frac{t-s}{1-s}+\frac{1-t}{1-s}s\Big)\leqslant \frac{t-s}{1-s}\psi(1)+\frac{1-t}{1-s}\psi(s)\\
&\leqslant\frac{t-s}{1-s}\cdot\frac{\psi(s)}{s}+\frac{1-t}{1-s}\psi(s)=\frac{t}{s}\psi(s).
\end{split}\]
In particular, one has
\[(a+b)\psi\Big(\frac{a}{a+b}\Big)\leqslant(a+b')\psi\Big(\frac{a}{a+b'}\Big).\]
Moreover, the function from $[0,1]$ to itself sending $t\in[0,1]$ to $\psi(1-t)$ also belongs to $\mathscr S$. By the above argument, we obtain that the function $t\mapsto \psi(1-t)/t$ is also non-increasing. Therefore
\[(a+b')\psi\Big(\frac{a}{a+b'}\Big)=(a+b')\psi\Big(1-\frac{b'}{a+b'}\Big)\leqslant (a'+b')\psi\Big(1-\frac{b'}{a'+b'}\Big).\]
The inequality \eqref{Equ:inequalitypsi} is thus proved.
\end{proof}

\begin{prop}\label{prop:norm:direct:sum:psi}
Let $(V,\norm{\ndot}_V)$ and $(W,\norm{\ndot}_W)$ be finite-dimensional seminormed vector spaces over $k$. For any $\psi\in\mathscr S$, let $\norm{\ndot}_\psi:V\oplus W\rightarrow\mathbb R_{\geqslant 0}$ be the map such that $\|(v,w)\|_\psi=0$ for $(v,w)\in N_{\norm{\ndot}_V}\oplus N_{\norm{\ndot}_W}$ and that, for any $(x,y)\in (V\oplus W)\setminus (N_{\norm{\ndot}_V}\oplus N_{\norm{\ndot}_W})$,
\[\|(x,y)\|_\psi:=(\|x\|+\|y\|) \psi\Big(\frac{\|x\|}{\|x\|+\|y\|}\Big).\]
Then $\norm{\ndot}_\psi$ is a seminorm on $V\oplus W$ such that $N_{\norm{\ndot}_\psi}=N_{\norm{\ndot}_V}\oplus N_{\norm{\ndot}_W}$. Moreover, for any $(x,y)\in V\times W$ one has
\begin{equation}\label{Equ:encadrementsommedirec}\max\{\|x\|,\|y\|\}\leqslant\|(x,y)\|_\psi\leqslant\|x\|+\|y\|.\end{equation}
\end{prop}
\begin{proof}By definition, for any $(x,y)\in V\oplus W$ and any $a\in k$, one has
\[\|(ax,ay)\|_\psi=|a|\cdot\|(x,y)\|_\psi.\] Moreover, for $(x,y)\not\in N_{\norm{\ndot}_V}\oplus N_{\norm{\ndot}_W}$, one has $\|(x,y)\|_\psi>0$.
Thus it remains to verify the triangle inequality.

Let $(x_1,y_1)$ and $(x_2,y_2)$ be two elements in $V\oplus W$ such that $(x_1+x_2,y_1+y_2)$ does not belong to $N_{\norm{\ndot}_V}\oplus N_{\norm{\ndot}_W}$. One has 
\[\|(x_1+x_2,y_1+y_2)\|_\psi=(\|x_1+x_2\|+\|y_1+y_2\|)\psi(u),\]
where 
\[u=\frac{\|x_1+x_2\|}{\|x_1+x_2\|+\|y_1+y_2\|}.\]
Since 
$\|x_1+x_2\|\leqslant\|x_1\|+\|x_2\|$ and $\|y_1+y_2\|\leqslant\|y_1\|+\|y_2\|$, by Lemma \ref{Lem:croissancef} one obtains that $\|(x_1+x_2,y_1+y_2)\|_\psi$
is bounded from above by 
\[(\|x_1\|+\|x_2\|+\|y_1\|+\|y_2\|)\psi(v),\]
with
\[v=\frac{\|x_1\|+\|x_2\|}{\|x_1\|+\|x_2\|+\|y_1\|+\|y_2\|}=\Big({1+\frac{\|y_1\|+\|y_2\|}{\|x_1\|+\|x_2\|}}\Big)^{-1}\]
{if} $\norm{x_1}+\norm{x_2}>0$, and $v=0$ otherwise.
If $\norm{x_1}>0$, let
\[s=\frac{\|x_1\|}{\|x_1\|+\|y_1\|}=\Big(1+\frac{\|y_1\|}{\|x_1\|}\Big)^{-1},\]
otherwise let $s=0$. Similarly, if $\norm{x_2}>0$, let
\[
t=\frac{\|x_2\|}{\|x_2\|+\|y_2\|}=\Big(1+\frac{\|y_2\|}{\|x_2\|}\Big)^{-1},\]
otherwise let $t=0$.
In the case where $\norm{x_1}$ and $\norm{x_2}$ are both $>0$, one has
\[\displaystyle\min\Big\{\frac{\|y_1\|}{\|x_1\|},\frac{\|y_2\|}{\|x_2\|}\Big\}\leqslant\frac{\|y_1\|+\|y_2\|}{\|x_1\|+\|x_2\|}\leqslant\max\Big\{\frac{\|y_1\|}{\|x_1\|},\frac{\|y_2\|}{\|x_2\|}\Big\},\]
and therefore $\min\{s,t\}\leqslant v\leqslant\max\{s,t\}$. By the convexity of the function $\psi$ we obtain
\[\psi(v)\leqslant\frac{v-t}{s-t}\psi(s)+\frac{s-v}{s-t}\psi(t).\]
Note that
\[\frac{v-t}{s-t}=\frac{\|x_1\|+\|y_1\|}{\|x_1\|+\|x_2\|+\|y_1\|+\|y_2\|},\quad
\frac{s-v}{s-t}=\frac{\|x_2\|+\|y_2\|}{\|x_1\|+\|x_2\|+\|y_1\|+\|y_2\|}.\]
Thus we obtain the triangle inequality 
$\|(x_1+x_2,y_1+y_2)\|_\psi\leqslant\|(x_1,y_1)\|_\psi+\|(x_2,y_2)\|_\psi$.

We now proceed with the proof of the inequalities \eqref{Equ:encadrementsommedirec}. The second inequality comes from the fact that $\psi$ takes values $\leqslant 1$. The first inequality is a consequence of Lemma \ref{Lem:croissancef}. In fact, by \eqref{Equ:inequalitypsi}, when $\norm{x}>0$ one has
\[\|x\|=(\|x\|+0)\psi\Big(\frac{\|x\|}{\|x\|+0}\Big)\leqslant \|(x,y)\|_{\psi}.\]
Similarly, one has $\|y\|\leqslant\|(x,y)\|_{\psi}$. The proposition is thus proved.
\end{proof}

\begin{defi}
The seminorm $\norm{\ndot}_\psi$ constructed in the above proposition is called the \emph{$\psi$-direct sum}\index{psi-direct sum@$\psi$-direct sum} of the seminorms of $V$ and $W$. 
\end{defi}

\begin{prop}
Let $\norm{\ndot}$ be an absolute normalised norm on $\mathbb R^2$. Then the dual norm $\norm{\ndot}_*$ is also an absolute normalised norm,
where $\Hom_{\mathbb R}(\mathbb R^2, \mathbb R)$ is identified with $\mathbb R^2$ by using the isomorphism $\iota : \mathbb R^2 \to \Hom_{\mathbb R}(\mathbb R^2, \mathbb R)$ given by
$\iota(x, y)(a, b) = ax + by$.
\end{prop}
\begin{proof}
Let $(x,y)$ be an element of $\mathbb R^2$. One has {(recall that $|\ndot|_\infty$ denotes the usual absolute value on $\mathbb R$)}
\[\|(x,y)\|_*=\sup_{(0,0)\neq(a,b)\in\mathbb R^2}\frac{|ax+by|_\infty}{\|(a,b)\|}.\]
Since $\norm{\ndot}$ is an absolute normalised norm on $\mathbb R^2$, from the above formula we deduce that $\|(x,y)\|_*=\|(|x|_\infty,|y|_\infty)\|_*$ for any $(x,y)\in\mathbb R^2$. Moreover, by \eqref{Equ: minoration norme (x,y)} one has
\[\|(1,0)\|_*=\sup_{(0,0)\neq(a,b)\in\mathbb R}\frac{|a|_\infty}{\|(a,b)\|}=\sup_{0\neq a\in\mathbb R}\frac{|a|_\infty}{\|(a,0)\|}=1.\]
Similarly, $\|(0,1)\|_*=1$. Therefore, $\norm{\ndot}_*$ is an absolute normalised norm on $\mathbb R^2$.
\end{proof}

\begin{defi}\label{Def:dualdirecsum}
Let $\psi$ be an element of $\mathscr S$, which corresponds to an absolute normalised norm $\norm{\ndot}$ on $\mathbb R^2$. The above proposition 
shows that the dual norm $\norm{\ndot}_*$ is also an absolute normalised norm. We denote by $\psi_*$ the element of $\mathscr S$ corresponding to this dual norm. 
Note that $\psi_*$ is actually given by
\[
\psi_*(t) = \sup_{\lambda \in \intervalle{]}{0}{1}{[}} \left\{  \frac{\lambda t + (1-\lambda)(1-t)}{\psi(\lambda)} \right\}.
\]
\end{defi}

The following proposition studies the dual of a direct sum norm.

\begin{prop}\label{Pro:dualdirectsum}
Let $(V,\norm{\ndot}_V)$ and $(W,\norm{\ndot}_W)$ be finite-dimensional seminormed vector spaces over $k$, $\psi$ be an element in $\mathscr S$, and $\norm{\ndot}_{\psi}$ be the $\psi$-direct sum of $\norm{\ndot}_V$ and $\norm{\ndot}_W$. Let $\psi_0\in\mathscr S$  such that $\psi_0(t)=\max\{t,1-t\}$ for any $t\in[0,1]$.
\begin{enumerate}[label=\rm(\arabic*)]
\item\label{Item: dual direct sum} Assume that the absolute value $|\ndot|$ is non-Archimedean. Then the dual norm $\norm{\ndot}_{\psi,*}$ identifies with the $\psi_0$-direct sum of $\norm{\ndot}_{V,*}$ and $\norm{\ndot}_{W,*}$.
\item\label{Item: dual direct sum2} Assume that the absolute value $|\ndot|$ is Archimedean. Then the dual norm $\norm{\ndot}_{\psi,*}$ identifies with the $\psi_*$-direct sum of $\norm{\ndot}_{V,*}$ and $\norm{\ndot}_{W,*}$.
\end{enumerate}
\end{prop}
\begin{proof} Since the null space of the seminorm $\norm{\ndot}$ is $N_{\norm{\ndot}_V}\oplus N_{\norm{\ndot}_W}$, we obtain that a linear form $(f,g)\in V^\vee\oplus W^\vee$ vanishes on $N_{\norm{\ndot}_{\psi}}$ if and only if it belongs to $V^*\oplus W^*$. In other words, one has $(V\oplus W)^*=V^*\oplus W^*$.

\ref{Item: dual direct sum} Let $(f,g)$ be an element in $V^*\oplus W^*$, one has
\[\begin{split}\|(f,g)\|_{\psi,*}&=\sup_{\begin{subarray}{c}(s,t)\in V\oplus W\\
\max\{\norm{s}_V,\norm{t}_W\}>0
\end{subarray}}\frac{|f(s)+g(t)|}{\|(s,t)\|_\psi}\\&\quad\leqslant\sup_{\begin{subarray}{c}(s,t)\in V\oplus W\\
\max\{\norm{s}_V,\norm{t}_W\}>0
\end{subarray}}\frac{\max\{|f(s)|,|g(t)|\}}{\max\{\|s\|_V,\|t\|_W\}}\leqslant\max\{\|f\|_{V,*},\|g\|_{W,*}\},
\end{split}\]
where the first inequality comes from \eqref{Equ:encadrementsommedirec} and the fact that the absolute value $|\ndot|$ is non-Archimedean. Moreover, one has
\[\|(f,g)\|_{\psi,*}\geqslant\sup_{s\in V\setminus N_{\norm{\sndot}_V}}\frac{|f(s)+g(0)|}{\|(s,0)\|_{\psi}}=\|f\|_{V,*}.\]
Similarly, one has
\[\|(f,g)\|_{\psi,*}\geqslant\sup_{t\in W\setminus N_{\norm{\sndot}_W}}\frac{|f(0)+g(t)|}{\|(0,t)\|_{\psi}}=\|g\|_{W,*}.\]
Therefore $\|(f,g)\|_{\psi,*}=\max\{\|f\|_{V,*},\|g\|_{W,*}\}$.

\ref{Item: dual direct sum2} Let $\norm{\ndot}$ be the absolute normalised norm on $\mathbb R^2$ corresponding to $\psi$ and let $\norm{\ndot}_*$ be its dual norm. For any $(s,t)\in V\oplus W$, on has
$\|(s,t)\|_\psi=\|(\|s\|_V,\|t\|_W)\|$.
Let $(f,g)$ be an element in $V^*\oplus W^*$. One has
\[\begin{split}&\quad\;\|(f,g)\|_{\psi,*}=\sup_{\begin{subarray}{c}(s,t)\in V\oplus W\\
\max\{\norm{s}_V,\norm{t}_W\}>0
\end{subarray}}\frac{|f(s)+g(t)|}{\|(s,t)\|_\psi}
\\&\leqslant\sup_{\begin{subarray}{c}(s,t)\in V\oplus W\\
\max\{\norm{s}_V,\norm{t}_W\}>0
\end{subarray}}\frac{\|f\|_{V,*}\cdot\|s\|_V+\|g\|_{W,*}\cdot\|t\|_W}{\|(\|s\|_V,\|t\|_W)\|}=\|(\|f\|_{V,*},\|g\|_{W,*})\|_*.
\end{split}\]
Moreover, since $k=\mathbb R$ or $\mathbb C$, by Hahn-Banach theorem, for any $a>0$, there exists $s\in V$ such that $\|s\|_V=a$ and that $f(s)=\|f\|_{V,*}\cdot\|s\|_V$. Similarly, for any $b>0$, there exists $t\in W$ such that $\|t\|_W=b$ and $g(t)=\|g\|_{W,*}\cdot\|t\|_W$. Therefore the inequality in the above formula is actually an equality. 
\end{proof}

\begin{prop}\label{prop:norm:quot:norm:direct:sum}
Let $f : V \to V'$ and $g : W \to W'$ be surjective linear maps of finite-dimensional vector spaces over $k$.
Let $\norm{\ndot}_V$ and $\norm{\ndot}_W$ be seminorms on  $V$ and $W$, and let $\norm{\ndot}_{V'}$ and $\norm{\ndot}_{W'}$ be the quotient seminorms
of $\norm{\ndot}_V$ and $\norm{\ndot}_W$ on $V'$ and $W'$, respectively.
Then the quotient seminorm $\norm{\ndot}_{V \oplus W, \psi, V \oplus W \twoheadrightarrow V' \oplus W'}$ 
of $\norm{\ndot}_{V \oplus W, \psi}$ on $V' \oplus W'$ coincides with
$\norm{\ndot}_{V' \oplus W', \psi}$.
\end{prop}

\begin{proof}
It is sufficient to see that
\[
\|(x', y')\|_{V' \oplus W', \psi} = \|(x', y')\|_{V \oplus W, \psi, V \oplus W \twoheadrightarrow V' \oplus W'}
\]
for all $x' \in V'$ and $y' \in W'$ with $\|x'\|_{V'} + \| y' \|_{W'} > 0$.
Let $x \in V$ and $y \in W$ with $f(x) = x'$ and $g(y) = y'$. Then,
as $\|x\|_V \geqslant \|x'\|_{V'}$ and $\|y\|_W \geqslant \|y'\|_{W'}$, by Lemma~\ref{Lem:croissancef}, one has
$\|(x', y')\|_{V' \oplus W', \psi} \leqslant \|(x, y)\|_{V \oplus W, \psi}$, so that
\[
 \|(x', y')\|_{V' \oplus W', \psi} \leqslant \|(x', y')\|_{V \oplus W, \psi, V \oplus W \twoheadrightarrow V' \oplus W'}.
\]

Let us consider the converse inequality. We
choose sequences $\{ x_n \}_{n\in\mathbb N}$ and $\{ y_n \}_{n\in\mathbb N}$ in $V$ and $W$ such that
$f(x_n) = x'$, $g(y_n) = y'$, $\lim_{n\to\infty} \| x_n \|_V = \|x'\|_{V'}$ and
$\lim_{n\to\infty} \| y_n \|_{W} = \|y'\|_{W'}$. 

We assume that $\| x' \|_{V'} + \| y' \|_{W'} > 0$. Then as $\| x_n \|_{V} + \| y_n \|_{W} > 0$ for sufficiently large $n$ and $\psi$ is continuous, one has
\begin{align*}
\|(x', y')\|_{V' \oplus W', \psi} & = (\|x'\|_{V'} + \|y'\|_{W'}) \psi\left( \frac{\|x'\|_{V'}}{\|x'\|_{V'} + \|y'\|_{W'}}\right) \\
& = \lim_{n\to\infty} (\|x_n\|_{V} + \|y_n\|_{W}) \psi\left( \frac{\|x_n\|_{V}}{\|x_n\|_{V} + \|y_n\|_{W}}\right) \\
& = \lim_{n\to\infty} \| (x_n, y_n )\|_{V \oplus W, \psi}  \geqslant \|(x', y')\|_{V \oplus W, \psi, V \oplus W \twoheadrightarrow V' \oplus W'},
\end{align*}
as required. Otherwise, as \[
0 \leqslant \|(x', y')\|_{V \oplus W, \psi, V \oplus W \twoheadrightarrow V' \oplus W'} \leqslant
\|(x_n, y_n)\|_{V \oplus W, \psi} \leqslant \| x_n \|_V + \| y_n \|_W
\]
and $\lim_{n\to\infty} \| x_n \|_V + \| y_n \|_W = 0$, one has
\[ \|(x', y')\|_{V \oplus W, \psi, V \oplus W \twoheadrightarrow V' \oplus W'} = 0,\] as desired.
\end{proof}

\begin{rema}\label{rem:psi:assoc:comm}
Let $\psi$ be an element of $\mathscr S$. Let $\{ \psi_n \}_{n=1}^{\infty}$ be a sequence of functions given in the following ways:
\[
\begin{cases}
\forall a \in \mathbb R_{\geqslant 0},  & \psi_1(a) = a, \\[1ex]
\forall (a, b) \in \mathbb R_{\geqslant 0}^2, & \psi_2(a, b) = \begin{cases}
{\displaystyle (a+b) \psi\left(\frac{a}{a+ b} \right)} & \text{if $a+ b > 0$}, \\
0 & \text{if $a = b = 0$},
\end{cases}\\[1ex]
\forall (a_1, \ldots, a_n) \in \mathbb R_{\geqslant 0}^n, &
\psi_n(a_1, \ldots, a_n) = \psi_2(\psi_{n-1}(a_1, \ldots, a_{n-1}), a_n).
\end{cases}
\]

Let $(V_1, \norm{\ndot}_1), \ldots, (V_n, \norm{\ndot}_n)$ be finite-dimensional normed vector spaces over $k$.
If we define 
\[
\|(x_1, \ldots, x_n)\|_{\psi} := \psi_n(\|x_1\|_1, \ldots, \|x_n\|_n)
\]
for $(x_1, \ldots, x_n) \in V_1 \oplus \cdots \oplus V_{n}$, then, by Proposition~\ref{prop:norm:direct:sum:psi},
it yields a norm on $V_1 \oplus \cdots \oplus V_{n}$.

We assume that 
\begin{equation}\label{eqn:rem:psi:assoc:comm:01}
\forall\, a_1, a_2, a_3 \in \mathbb R_{\geqslant 0},\quad
\psi_2(a_1, \psi_2(a_2, a_3)) = \psi_2(\psi_2(a_1,a_2), a_3).
\end{equation}
Then it is easy to see that
$\psi_{n}(a_1, \ldots, a_n) =\psi_2(\psi_{i}(a_1, \ldots, a_i), \psi_{n-i}(a_{i+1}, \ldots, a_n))$
for $i\in\{1,\ldots,n-1\}$, so that
the construction of the norm $\norm{\ndot}_{\psi}$ is associative.
If we assume $\psi_2(a, b) = \psi_2(b, a)$ for all $(a, b) \in \mathbb R_{\geqslant 0}^2$ in addition to
\eqref{eqn:rem:psi:assoc:comm:01},
then $\psi_{n}(a_1, \ldots, a_n)$ is symmetric, that is,
for any permutation $\sigma$, $\psi_{n}(a_{\sigma(1)}, \ldots, a_{\sigma(n)}) = \psi_{n}(a_1, \ldots, a_n)$, which means that its construction is order independent.
\end{rema}

\subsection{Tensor product seminorms}\label{Subsec:tensorproduct}
Let $V$ and $W$ be seminormed vector spaces of finite rank over $k$. On the tensor product space $V\otimes_k W$ there are several natural ways to construct tensor product seminorms. We refer the readers to the original article \cite{Grothendieck53} of Grothendieck for different constructions. In this subsection, we recall the $\pi$-tensor product and the $\varepsilon$-tensor product. We refer to the book \cite{Ryan02} for a more detailed presentation in the Archimedean case. 

\begin{defi}\label{Def:tensorproducts}
Let $(V_1,\norm{\ndot}_1),\ldots,(V_n,\norm{\ndot}_n)$ be 
seminormed vector spaces over $k$. We define a map $\norm{\ndot}_\pi:V_1\otimes_k\cdots\otimes_k V_n\rightarrow \intervalle{[}{0}{+\infty}{[}$ such that, for any $\varphi\in V_1\otimes_k\cdots\otimes_k V_n$,
\begin{equation}\|\varphi\|_\pi:=\inf\bigg\{\sum_{i=1}^N\|x_1^{(i)}\|_1\cdots\|x_n^{(i)}\|_n\,:\,\varphi=\sum_{i=1}^Nx_1^{(i)}\otimes\cdots\otimes x_n^{(i)}\bigg\}.\end{equation}
Note that $\norm{\ndot}_\pi$ is a seminorm on $V_1\otimes_k\cdots\otimes_k V_n$. For example, the triangle inequality can be checked as follows: for $\varphi, \psi \in V_1\otimes_k\cdots\otimes_k V_n$ and
a positive number $\epsilon$, we choose
expressions \[\varphi =\sum_{i=1}^Nx_1^{(i)}\otimes\cdots\otimes x_n^{(i)}\quad\text{and} \quad
\psi =\sum_{j=1}^M y_1^{(j)}\otimes\cdots\otimes y_n^{(j)}\] such that
\[\sum_{i=1}^N\|x_1^{(i)}\|_1\cdots\|x_n^{(i)}\|_n \leqslant \| \varphi\|_{\pi} + \epsilon\quad\text{and}\quad\sum_{j=1}^M\|y_1^{(j)}\|_1\cdots\|y_n^{(j)}\|_n \leqslant \| \psi\|_{\pi} + \epsilon.\]
Then, as \[\varphi + \psi = \sum_{i=1}^Nx_1^{(i)}\otimes\cdots\otimes x_n^{(i)} + \sum_{j=1}^M y_1^{(j)}\otimes\cdots\otimes y_n^{(i)},\] one has
\[
\| \varphi + \psi\|_{\pi} \leqslant \sum_{i=1}^N\|x_1^{(i)}\|_1\cdots\|x_n^{(i)}\|_n  +
\sum_{j=1}^M\|y_1^{(j)}\|_1\cdots\|y_n^{(j)}\|_n \leqslant \| \varphi\|_{\pi} + \| \psi \|_{\pi} + 2\epsilon,
\]
as desired. 
We call $\norm{\ndot}_\pi$ the \emph{$\pi$-tensor product}\index{pi-tensor product@$\pi$-tensor product} of the seminorms $\norm{\ndot}_1,\ldots,\norm{\ndot}_n$.

Any element $\varphi$ in the tensor product space $V_1\otimes_k\cdots\otimes_k V_n$ can be considered as a multilinear form on $V_1^*\times\cdots\times V_n^{*}$. In particular, if $\varphi$ is of the form $x_1\otimes\cdots\otimes x_n$, the corresponding multilinear form sends $(f_1,\ldots,f_n)\in V_1^*\times\cdots\times V_n^{*}$ to $f_1(x_1)\cdots f_n(x_n)\in k$. For any  $\varphi\in V_1\otimes_k\cdots\otimes_k V_n$, viewed as a $k$-multilinear form on $V_1^*\times\cdots\times V_n^{*}$, let
\[\|\varphi\|_{\varepsilon} :=\sup_{\begin{subarray}{c}(f_1,\ldots,f_n)\in V_1^*\times\cdots\times V_n^{*}\\
\forall\,i\in\{1,\ldots,n\},\;f_i\neq 0\end{subarray}}\frac{|\varphi(f_1,\ldots,f_n)|}{\|f_1\|_{1,*}\cdots\|f_n\|_{n,*}}.\]
Then $\norm{\ndot}_\varepsilon$ is a seminorm on the tensor product space $V_1\otimes_k\cdots\otimes_k V_n$, called the  \emph{$\varepsilon$-tensor product}\index{epsilon-tensor product@$\varepsilon$-tensor product} of seminorms  $\norm{\ndot}_1,\ldots,\norm{\ndot}_n$. It is a norm once the seminorms $\norm{\ndot}_1,\ldots,\norm{\ndot}_n$ are norms. Similarly to the dual norm case, if the absolute value $|\ndot|$ is non-Archimedean, then the $\varepsilon$-tensor product $\norm{\ndot}_{\varepsilon}$ is ultrametric. By Proposition \ref{Pro:doubedualandquotient} in the next {section}, we obtain that, in the case where all $V_i$ are of finite type over $k$, the $\varepsilon$-tensor product of $\norm{\ndot}_1,\ldots,\norm{\ndot}_n$ identifies with that of $\norm{\ndot}_{1,**},\ldots,\norm{\ndot}_{n,**}$.
\end{defi}

\begin{rema}\label{Rem:operateureps}
Let $(V_1,\norm{\ndot}_1)$ and $(V_2,\norm{\ndot}_2)$ be seminormed vector spaces over $k$. Let $\norm{\ndot}$ be the operator seminorm on the vector space $\mathscr L(V_1^{*},V_2)$, where we consider the dual norm $\norm{\ndot}_{1,*}$ on $V_1^{*}$ and the double dual seminorm $\norm{\ndot}_{2,**}$ on $V_2$. One has a canonical $k$-linear map from $V_1\otimes_kV_2$ to $\mathscr L(V_1^{*},V_2)$ sending $x\otimes y\in V_1\otimes_kV_2$ to the bounded linear map $(\alpha\in V_1^*)\mapsto \alpha(x)y$.  We claim that the seminorm on $V_1\otimes_kV_2$ induced by $\norm{\ndot}$ and the above canonical map identifies with the $\varepsilon$-tensor product $\norm{\ndot}_{\varepsilon}$ of $\norm{\ndot}_1$ and $\norm{\ndot}_2$. In fact, for any $\varphi\in V_1\otimes_k V_2$ one has
\[\|\varphi\|=\sup_{f_1\in V_1^*\setminus\{0\}}\frac{\|\varphi(f_1)\|_{2,**}}{\|f_1\|_{1,*}}=\sup_{\begin{subarray}{c}f_1\in V_1^*\setminus\{0\}\\
f_2\in V_2^*\setminus\{0\}
\end{subarray}}\frac{|\varphi(f_1,f_2)|}{\|f_1\|_{1,*}\|f_2\|_{2,*}}=\|\varphi\|_\varepsilon.\]
In particular, if {$\norm{\ndot}_2=\norm{\ndot}_{2,**}$ on $V_2$}, then the $\varepsilon$-tensor product norm $\norm{\ndot}_{\varepsilon}$ identifies with the operator seminorm if we consider tensors in $V_1\otimes_kV_2$ as $k$-linear operators from $(V_1^*,\norm{\ndot}_{1,*})$ to $(V_2,\norm{\ndot}_{2})$.
\end{rema}

\begin{prop}
\label{Pro:normpiisanorm}We keep the notation of Definition \ref{Def:tensorproducts}. If $\norm{\ndot}$ is a seminorm on $V_1\otimes_k\cdots\otimes_kV_n$ such that $\norm{x_1\otimes\cdots\otimes x_n}\leqslant\norm{x_1}_1\cdots\norm{x_n}_n$ for any $(x_1,\ldots,x_n)\in V_1\times\cdots\times V_n$, then one has $\norm{\ndot}\leqslant\norm{\ndot}_{\pi}$. In particular, the seminorm $\norm{\ndot}_\varepsilon$ is bounded from above by $\norm{\ndot}_\pi$. Moreover, if $\|\ndot\|_1, \ldots, \|\ndot\|_n$ are norms, then $\|\ndot\|_{\pi}$ is also a norm.
\end{prop}
\begin{proof}
Let $\varphi$ be an element of $V_1\otimes_k\cdots\otimes_k V_n$. 
If $\varphi$ is written in the form 
\[\varphi=\sum_{i=1}^Nx_1^{(i)}\otimes\cdots\otimes x_n^{(i)},\]
where $x_j^{(i)}\in V_j$ for any $j\in\{1,\ldots,n\}$, then one has 
\[\norm{\varphi}\leqslant\sum_{i=1}^N\norm{x_1^{(i)}\otimes\cdots\otimes x_n^{(i)}}\leqslant\sum_{i=1}^N\norm{x_1^{(1)}}_1\cdots\norm{x_n^{(1)}}_n.\]
Therefore we obtain $\norm{\ndot}\leqslant\norm{\ndot}_\pi$. Note that, for any {$(x_1,\ldots,x_n)\in V_1\times\cdots\times V_n$} one has
\[
\begin{split}\|x_1\otimes\cdots\otimes x_n\|_{\varepsilon}&=\sup_{\begin{subarray}{c}
(f_1,\ldots,f_n)\in V_1^*\times\cdots\times V_n^*\\
\forall\,i\in\{1,\ldots,n\},\,f_i\neq 0
\end{subarray}}\frac{|f_1(x_1)|\cdots|f_n(x_n)|}{\norm{f_1}_{1,*}\cdots\norm{f_n}_{n,*}}\\&=\|x_1\|_{1,**}\cdots\|x_n\|_{n,**}\leqslant \|x_1\|_{1}\cdots\|x_n\|_{n}.\end{split}
\]
Therefore, one has $\norm{\ndot}_{\varepsilon}\leqslant\norm{\ndot}_\pi$. If the seminorms $\norm{\ndot}_i$ ($i\in\{1,\ldots,n\}$) are norms, then $\norm{\ndot}_{\varepsilon}$ is a norm and hence $\norm{\ndot}_{\pi}$ is also a norm.
\end{proof}

\begin{rema}
From the definition we observe that the $\varepsilon$-tensor product and $\pi$-tensor product are commutative. Namely, if $V_1$ and $V_2$ are finite-dimensional normed vector spaces over $k$, then the canonical isomorphism $V_1\otimes_k V_2\rightarrow V_2\otimes_k V_1$ is an isometry if we consider $\varepsilon$-tensor products or $\pi$-tensor product norms on both sides. The $\varepsilon$-tensor product and the $\pi$-tensor product are also associative. Namely, if $V_1$, $V_2$ and $V_3$ are finite-dimensional normed vector spaces over $k$, then the canonical isomorphisms $(V_1\otimes_k V_2)\otimes_k V_3\rightarrow V_1\otimes_k V_2\otimes_k V_3$ and $V_1\otimes_k(V_2\otimes_k V_3)\rightarrow V_1\otimes_k V_2\otimes_k V_3$ are both isometries.  
\end{rema}

\begin{rema}\label{Rem:produittenrk1}
Let $(V_1,\norm{\ndot}_1),\ldots,(V_n,\norm{\ndot}_n)$ be finite-dimensional seminormed vector spaces over $k$. From the definition, we observe that, if $(u_1,\ldots,u_n)$ is an element of $V_1\times\cdots\times V_n$, then one has
\begin{equation}
\label{equ:epsilontensorpodurc}
\|u_1\otimes\cdots\otimes u_n\|_{\varepsilon}=\|u_1\|_{1,**}\cdots\|u_n\|_{n,**}.
\end{equation} 

If the seminormed vector spaces $(V_1,\norm{\ndot}_1),\ldots,(V_n,\norm{\ndot}_n)$ are reflexive, by \eqref{equ:epsilontensorpodurc} and Proposition \ref{Pro:normpiisanorm}, we obtain that, for any $(u_1,\ldots,u_n)\in V_1\times\cdots\times V_n$, one has 
\[{\prod_{i=1}^n\|u_i\|_i =\|u_1\otimes\cdots\otimes u_n\|_{\varepsilon}\leqslant\|u_1\otimes\cdots\otimes u_n\|_\pi.}\]
Moreover, by definition one has $\|u_1\otimes\cdots\otimes u_n\|_\pi\leqslant\|u_1\|_1\cdots\|u_n\|_n$. Therefore 
\begin{equation}\label{Equ: split tensor}\|u_1\otimes\cdots\otimes u_n\|_{\varepsilon}=\|u_1\otimes\cdots\otimes u_n\|_{\pi}=\|u_1\|_1\cdots\|u_n\|_n.\end{equation}
In particular, if $V_1,\ldots,V_n$ are seminormed vector spaces of rank $1$ over $k$ (in this case they are necessarily reflexive), then their $\varepsilon$-tensor product and $\pi$-tensor product norms are the same. We simply call it the \emph{tensor product}\index{tensor product} of the seminorms $\norm{\ndot}_1,\ldots,\norm{\ndot}_n$.
\end{rema}

\begin{prop}
\label{Pro:dualitypiepsilon} Let $(V_1,\norm{\ndot}_1),\ldots,(V_n,\norm{\ndot}_n)$ be finite-dimensional seminormed vector spaces over $k$. Let $\norm{\ndot}_{*,\pi}$ and $\norm{\ndot}_{*,\varepsilon}$ be respectively the {$\pi$-tensor product} and the {$\varepsilon$-tensor product} of the dual norms $\norm{\ndot}_{1,*},\ldots,\norm{\ndot}_{n,*}$.
The $\varepsilon$-tensor product of $\norm{\ndot}_1,\ldots,\norm{\ndot}_n$ identifies with the seminorm induced by the dual norm $\norm{\ndot}_{*,\pi,*}$ on $(V_1^*\otimes_k\cdots\otimes_k V_n^*)^*$ by the natural linear map $V_1\otimes_k\cdots\otimes_k V_n\rightarrow (V_1^*\otimes_k\cdots\otimes_k V_n^*)^*$. If the absolute value $|\ndot|$ is Archimedean, then the $\pi$-tensor product of $\norm{\ndot}_1,\ldots,\norm{\ndot}_n$ identifies with seminorm induced by the dual norm $\norm{\ndot}_{*,\varepsilon,*}$ on $(V_1^*\otimes_k\cdots\otimes_k V_n^*)^*$ by the natural linear map $V_1\otimes_k\cdots\otimes_k V_n\rightarrow (V_1^*\otimes_k\cdots\otimes_k V_n^*)^*$.
\end{prop}
\begin{proof}
Let $\varphi$ be an element in $V_1\otimes_k\cdots\otimes_k V_n$, which can also be viewed as a $k$-multilinear form on $V_1^*\times\cdots\times V_n^*$ or a linear form on $V_1^*\otimes_k\cdots\otimes_k V_n^*$.  Let $\alpha$ be an element in $V_1^*\otimes_k\cdots\otimes_k V_n^*$. If $\alpha$ is written in the form
\[\alpha=\sum_{i=1}^Nf_1^{(i)}\otimes\cdots\otimes f_n^{(i)},\]
where $f_j^{(i)}\in V_j^*$, then one has
${\varphi}(\alpha)=\sum_{i=1}^N\varphi(f_1^{(i)},\ldots,f_n^{(i)})$
and hence
\[|{\varphi}(\alpha)|\leqslant\sum_{i=1}^N|\varphi(f_1^{(i)},\ldots,f_n^{(i)})|.\]
Thus we obtain
\[\frac{|{\varphi}(\alpha)|}{\sum_{i=1}^N\|f_1^{(i)}\|_{1,*}\cdots\|f_n^{(i)}\|_{n,*}}\leqslant\frac{\sum_{i=1}^N|\varphi(f_1^{(i)},\ldots,f_n^{(i)})|}{\sum_{i=1}^N\|f_1^{(i)}\|_{1,*}\cdots\|f_n^{(i)}\|_{n,*}}\leqslant\|\varphi\|_{\varepsilon}.\] 
Therefore $\varphi$ is a bounded linear form on $(V_1^*\otimes_k\cdots\otimes_k V_n^*,\norm{\ndot}_{*,\pi})$ and $\|\varphi\|_{*,\pi,*}\leqslant\|\varphi\|_{\varepsilon}$.

For any $(f_1,\ldots,f_n)\in (V_1^*\setminus\{0\})\times\cdots\times (V_n^*\setminus\{0\})$ one has
\[\frac{|{\varphi}(f_1\otimes\cdots\otimes f_n)|}{\|f_1\otimes\cdots\otimes f_n\|_{*,\pi}}=\frac{|\varphi(f_1,\ldots,f_n)|}{\|f_1\otimes\cdots\otimes f_n\|_{*,\pi}}\geqslant\frac{|\varphi(f_1,\ldots,f_n)|}{\|f_1\|_{1,*}\cdots\|f_n\|_{n,*}}.\]
Therefore one has $\|\varphi\|_\varepsilon\leqslant\|\varphi\|_{*,\pi,*}$. The first assertion is thus proved.

If $|\ndot|$ is Archimedean, any finite-dimensional normed {vector} space is reflexive. By the first assertion, the dual norm of the $\pi$-tensor product of $\norm{\ndot}_1,\ldots,\norm{\ndot}_n$ is the $\varepsilon$-tensor product of $\norm{\ndot}_{1,*},\ldots,\norm{\ndot}_{n,*}$. By taking the double dual seminorm we obtain that the $\pi$-tensor product of $\norm{\ndot}_1,\ldots,\norm{\ndot}_n$ identifies with the seminorm induced by $\norm{\ndot}_{*,\varepsilon,*}$.
\end{proof}

\begin{prop}\label{Pro:quotientavecpitensor}
Let $V$ and $W$ be seminormed vector spaces over $k$, and $Q$ be a quotient space of $V$, equipped with the quotient seminorm. Let $V_0$ be the kernel of the projection map $V\rightarrow Q$. Then the canonical isomorphism $(V\otimes_k W)/(V_0\otimes_k W)\rightarrow Q\otimes_k W$ is an isometry, where we consider the $\pi$-tensor product seminorms on $V\otimes_k W$ and $Q\otimes_k W$, and the quotient seminorm on $(V\otimes_k W)/(V_0\otimes_k W)$. 
\end{prop}
\begin{proof} 
Let $\psi$ be an element of $Q\otimes_kW$. One has
\[\begin{split}\|\psi\|_{\pi}&=\inf\Big\{\sum_{i=1}^N\|\alpha_i\|\cdot\| y_i\|\,:\, {N\in\mathbb N},\;\psi=\sum_{i=1}^N\alpha_i\otimes y_i\}\\
&=\inf_{N\in\mathbb N}\inf_{\begin{subarray}{c}(\alpha_i)_{i=1}^N\in Q^N\\
(y_i)_{i=1}^N\in W^N\\
\psi=\sum_{i=1}^N\alpha_i\otimes y_i\end{subarray}}\inf_{\begin{subarray}{c}(x_i)_{i=1}^N\in V^N\\
[x_i]=\alpha_i\end{subarray}}\sum_{i=1}^N\|x_i\|\cdot\|y_i\|=\inf_{\begin{subarray}{c}\varphi\in V\otimes W\\
[\varphi]=\psi\end{subarray}}\|\varphi\|_\pi.
\end{split}\]
\end{proof}

\begin{rema}\label{Rem:quotientepsilontensor}
We consider the $\varepsilon$-tensor product analogue of the above proposition. Let $f$ be an element of $V\otimes_kW$, viewed as a $k$-bilinear form on $V^*\times W^*$. Then its image $g$ in $Q\otimes_kW$ corresponds to the restriction of $f$ {to} $Q^*\times W^*$. By Proposition \ref{Pro:dualquotient}, the dual {norm} on $Q^*$ of the quotient {seminorm} identifies with the restriction {to $Q^*$} of the dual norm on $V^*$. Therefore, one has $\|g\|_\varepsilon\leqslant\|f\|_\varepsilon$. However, in the case where the absolute value $|\ndot|$ is Archimedean, in general the inequality 
\[\|g\|_{\varepsilon}\leqslant\inf_{\begin{subarray}{c}
f\in V\otimes_kW\\
f|_{Q^*\times W^*}=g
\end{subarray}}\|f\|_{\varepsilon}\]
is strict. In fact, this problem is closely related to the extension property of the normed vector space $V^*$, which consists of extending a linear operator defined on a vector subspace of $V^*$ and valued in another seminormed vector space while keeping the operator seminorm. In the case where the linear operator is a linear form (namely valued in $k$), it is just {a consequence of} Hahn-Banach theorem. However, in general the extension property does not hold, except in the cases where $\rang_k(V)\leqslant 2$ or the norm on $V$ comes from a symmetric semipositive bilinear form (see \S\ref{Subsec:innerproduct} for the notation). We refer the readers to \cite{Kakutani39,Saccoman01} for more details.

In the case where the absolute value $|\ndot|$ is non-Archimedean, any dual norm is ultrametric, we will give a proof for the $\varepsilon$-tensor product analogue of Proposition \ref{Pro:quotientavecpitensor}, by using the ultrametric Gram-Schmidt process (see Proposition \ref{Pro:quotientr1eps}).
\end{rema}

\begin{prop}\label{Pro: restriction and tensors}
Let $(V, \norm{\ndot}_V)$ and $(W, \norm{\ndot}_W)$ be seminormed
vector spaces over $k$, $V_0$ be a vector subspace of $V$ and $\norm{\ndot}_{V_0}$ be the restriction of $\norm{\ndot}_V$ {to} $V_0$. 
\begin{enumerate}[label=\rm(\arabic*)]
\item\label{Item: sub of pi tensor product} Let $\norm{\ndot}_\pi$ be the $\pi$-tensor product of $\norm{\ndot}_V$ and $\norm{\ndot}_W$, $\norm{\ndot}_{\pi,0}$ be the $\pi$-tensor product of $\norm{\ndot}_{V_0}$ and $\norm{\ndot}_{W}$. Then the seminorm $\norm{\ndot}_{\pi,0}$ is bounded from below by the restriction of $\norm{\ndot}_\pi$ {to} $V_0\otimes_kW$.
\item\label{Item: sub of epsilon tensor product} Let $\norm{\ndot}_{\varepsilon}$ be the $\varepsilon$-tensor product of $\norm{\ndot}_V$ and $\norm{\ndot}_W$, and $\norm{\ndot}_{\varepsilon,0}$ be the $\varepsilon$-tensor product of $\norm{\ndot}_{V_0}$ and $\norm{\ndot}_{W}$. Then the seminorm $\norm{\ndot}_{\varepsilon,0}$ is bounded from below by the restriction of $\norm{\ndot}_{\varepsilon}$ {to} $V_0\otimes_kW$.
\end{enumerate}
\end{prop}
\begin{proof}
\ref{Item: sub of pi tensor product} Let $\varphi$ be an element of $V_0\otimes_kW$. By definition, for any writing of $\varphi$ as $\sum_{i=1}^N x_i\otimes y_i$ with $\{x_1,\ldots,x_N\}\subseteq V_0$ and $\{y_1,\ldots,y_N\}\subseteq W$, one has
\[\|\varphi\|_\pi\leqslant\sum_{i=1}^N\|x_i\|_{V_0}\cdot\|y_i\|_W.\]
Therefore $\|\varphi\|_{\pi}\leqslant\|\varphi\|_{\pi,0}$.

\ref{Item: sub of epsilon tensor product} We consider the canonical linear map $V^*\rightarrow V_0^*$ sending a bounded linear form on $V$ to its restriction {to} $V_0$. Note that for any $f\in V^*$ one has $\norm{f_0}_{V_0,*}\leqslant\norm{f}_{V,*}$, where $f_0$ is the restriction of $f$ {to} $V_0$. Therefore, for any element $\varphi$ of $V_0\otimes_kW$, viewed as a bilinear form on $V_0^*\times W^*$ or as a bilinear form {on} $V^*\times W^*$ via the inclusion $V_0\otimes_kW\subseteq {V\otimes_kW}$, one has
\[\|\varphi\|_{\varepsilon,0}=\sup_{\begin{subarray}{c}(f_0,g)\in V_0^*\times W^*\\
f_0\neq 0,\,g\neq 0\end{subarray}}\frac{|\varphi(f_0,g)|}{\|f\|_{V_0,*}\|g\|_{W,*}}\geqslant\sup_{\begin{subarray}{c}(f,g)\in V^*\times W^*\\
f\neq 0,\,g\neq 0\end{subarray}}\frac{|\varphi(f,g)|}{\|f\|_{V,*}\|g\|_{W,*}}=\|\varphi\|_{\varepsilon}.\]
\end{proof}

\begin{prop}\label{Pro: operator norm of tensors}
Let $n$ be a positive integer 
and $\{(V_j,\norm{\ndot}_{V_j})\}_{j=1}^n$ and $\{(W_j,\norm{\ndot}_{W_j})\}_{j=1}^n$ be finite-dimensional seminormed vector spaces over $k$. For any $j\in\{1,\ldots,n\}$, let $f_j:V_j\rightarrow W_j$ be a bounded $k$-linear map. Let $f:V_1\otimes_k\cdots\otimes V_n\rightarrow W_1\otimes_k\cdots\otimes_kW_n$ be the $k$-linear map sending $x_1\otimes\cdots\otimes x_n$ to $f_1(x_1)\otimes\cdots\otimes f_n(x_n)$.
\begin{enumerate}[label=\rm(\arabic*)]
\item\label{Item: norm of tensor product operator} We equip the vector spaces {$V_1\otimes_k\cdots\otimes_k V_n$} and $W_1\otimes_k\cdots\otimes_kW_n$ with the $\pi$-tensor product seminorms of $\{\norm{\ndot}_{V_j}\}_{j=1}^n$ and of $\{\norm{\ndot}_{W_j}\}_{j=1}^n$, respectively. Then the operator seminorm of $f$ is bounded from above by $\|f_1\|\cdots\|f_n\|$.
\item\label{Item: tensor product of dual operator} We equip the vector spaces {$V_1\otimes_k\cdots\otimes_k V_n$} and $W_1\otimes_k\cdots\otimes_kW_n$ with the $\varepsilon$-tensor product seminorms of $\{\norm{\ndot}_{V_j}\}_{j=1}^n$ and of $\{\norm{\ndot}_{W_j}\}_{j=1}^n$, respectively. Then the operator seminorm of $f$ is bounded from above by $\|f_1^*\|\cdots\|f_n^*\|$.
\end{enumerate}
\end{prop}
\begin{proof}
\ref{Item: norm of tensor product operator} Let $\varphi$ be an element in $V_1\otimes_k\cdots\otimes_kV_n$, which is written as
\[\varphi=\sum_{i=1}^Nx_1^{(i)}\otimes\cdots\otimes x_n^{(i)}\]
where $x_j^{(i)}\in V_j$ for any $j\in\{1,\ldots,n\}$. By definition, one has 
\[f(\varphi)=\sum_{i=1}^Nf_1(x_1^{(i)})\otimes\cdots\otimes f_n(x_n^{(i)}).\]
Therefore
\[\|f(\varphi)\|_\pi\leqslant\sum_{i=1}^N\|f_1(x_1^{(i)})\|_{W_1}\cdots\|f_n(x_n^{(i)})\|_{W_n}\leqslant\bigg(\prod_{i=1}^n\|f_i\|\bigg)\sum_{i=1}^N\|x_1^{(i)}\|_{V_1}\cdots\|x_n^{(i)}\|_{V_n}.\]
Thus $\|f(\varphi)\|_\pi\leqslant \|f_1\|\cdots\|f_n\|\cdot\|\varphi\|_\pi$.

\ref{Item: tensor product of dual operator} Let $\varphi$ be an element in $V_1\otimes_k\cdots\otimes_nV_n$, which can be viewed as a multilinear form on $V_1^*\times\cdots\times V_n^*$. Then the element $f(\varphi)\in W_1\otimes_k\cdots\otimes_kW_n$, viewed as a multilinear form on $W_1^*\times\cdots\times W_n^*$, sends $(\beta_1,\ldots,\beta_n)\in W_1^*\times\cdots\times W_n^*$ to
$\varphi(f_1^*(\beta_1),\ldots,f_n^*(\beta_n))$. Thus for $(\beta_1,\ldots,\beta_n)\in (W_1^*\setminus\{0\})\times\cdots\times (W_n^*\setminus\{0\})$, one has
\[\frac{|f(\varphi)(\beta_1,\ldots,\beta_n)|}{\|\beta_1\|_{W_1,*}\cdots\|\beta_n\|_{W_n,*}}\leqslant\frac{\|\varphi\|_{\varepsilon}\|f_1^*(\beta_1)\|_{V_1,*}\cdots\|f_n^*(\beta_n)\|_{V_n,*}}{\|\beta_1\|_{W_1,*}\cdots\|\beta_n\|_{W_n,*}}\leqslant\|\varphi\|_{\varepsilon}\|f_1^*\|\cdots\|f_n^*\|,\]
so that $\|f(\varphi)\|_{\varepsilon} \leqslant \|\varphi\|_{\varepsilon}\| f_1^{*} \| \cdots \| f_n^{*} \|$, as required.
\end{proof}

\subsection{Exterior power seminorm}\label{Subsect: Exterior power norm}
Let $V$ be a vector space over $k$ and $r$ be the rank of $V$ over $k$. For any $i\in\mathbb N$, we let $\Lambda^iV$ be the $i^{\text{th}}$ exterior power of the vector space $V$. It is a quotient vector space of $V^{\otimes i}$. 
\begin{defi}\label{Def:exterior:power:norm}
Let $\norm{\ndot}$ be a seminorm on the vector space $V$ and $\norm{\ndot}_\pi$ be the $\pi$-tensor power of $\norm{\ndot}$ on $V^{\otimes i}$. The \emph{$i^{\text{th}}$ $\pi$-exterior power seminorm}\index{pi-exterior power seminorm@$\pi$-exterior power seminorm}\index{seminorm!pi-exterior power@$\pi$-exterior power} of $\norm{\ndot}$ on $\Lambda^iV$ is by definition the quotient seminorm on $\Lambda^iV$ of $\norm{\ndot}_\pi$ induced by the canonical projection map $V^{\otimes i}\rightarrow\Lambda^iV$ sending $x_1\otimes\cdots\otimes x_i$ to $x_1\wedge\cdots\wedge x_i$, denoted by $\norm{\ndot}_{\Lambda_\pi^i}$, or simply by $\norm{\ndot}_{\Lambda^i}$. Similarly, the $\varepsilon$-tensor product seminorm $\norm{\ndot}_{\varepsilon}$ on $V^{\otimes i}$ induces by quotient a seminorm on $\Lambda^iV$, called the \emph{$i^{\text{th}}$ $\varepsilon$-exterior power}\index{epsilon-exterior power seminorm@$\varepsilon$-exterior power seminorm}\index{seminorm!epsilon-exterior power@$\varepsilon$-exterior power} of $\norm{\ndot}$, denoted by $\norm{\ndot}_{\Lambda^i_{\varepsilon}}$.
\end{defi}

\begin{prop}
\label{Pro:Hadamard's inequality}
Let $(V,\norm{\ndot})$ be a seminormed vector space over $k$ and $i$ be a natural number.  For any $(x_1,\ldots,x_i)\in V^i$ one has
\[\|x_1\wedge\cdots\wedge x_i\|_{\Lambda_\varepsilon^i}\leqslant \|x_1\wedge\cdots\wedge x_i\|_{\Lambda^i_\pi}\leqslant\|x_1\|\cdots\|x_i\|.\]
\end{prop}
\begin{proof}The first inequality follows from Proposition \ref{Pro:normpiisanorm}.

Note that $x_1\wedge\cdots\wedge x_i$ is the image of $x_1\otimes\cdots\otimes x_i$ by the canonical projection map $V^{\otimes i}\rightarrow\Lambda^iV$. Therefore one has
\[\|x_1\wedge\cdots\wedge x_i\|_{\Lambda^i_\pi}\leqslant\|x_1\otimes\cdots\otimes x_i\|_\pi\leqslant\|x_1\|\cdots\|x_i\|,\]
where $\norm{\ndot}_{\pi}$ denotes the $\pi$-tensor power of $\norm{\ndot}$.
\end{proof}

\begin{prop}
\label{Pro:Hadamard}
Let $V$ and $W$ be seminormed vector spaces over $k$ and $f:V\rightarrow W$ be a bounded $k$-linear map. Let $i$ be a positive integer.  
The $k$-linear map $f$ induces by passing to the $i^{\text{th}}$ exterior power a $k$-linear map $\Lambda^if:\Lambda^iV\rightarrow\Lambda^iW$.
\begin{enumerate}[label=\rm(\arabic*)]
\item\label{Item: operator norm exterior power} If we equip $\Lambda^iV$ and $\Lambda^iW$ with the $i^{\mathrm{th}}$ $\pi$-exterior power seminorms, then the operator seminorm of $\Lambda^if$ is bounded from above by $\|f\|^i$.
\item\label{Item: operator norm exterior power2} If we equip $\Lambda^iV$ and $\Lambda^iW$ with the $i^{th}$ $\varepsilon$-exterior power seminorms, then the operator seminorm of $\Lambda^if$ is bounded from above by $\|f^*\|^i$.
\end{enumerate}
\end{prop}
\begin{proof}
Let us consider a commutative diagram:
\[
\xymatrix{\relax
V ^{\otimes i} \ar[r]^-{f^{\otimes i}}\ar[d]& W^{\otimes i}\ar[d] \\
\Lambda^i V \ar[r]_-{\Lambda^i f}& \Lambda^i W}
\]
By \ref{SubItem: diagram of quotient seminorm1} in Proposition~\ref{prop:quotient:norm:linear:map}, we obtain
$\| \Lambda^i f \| \leqslant \| f^{\otimes i} \|$.
Thus the assertions follow from Proposition~\ref{Pro: operator norm of tensors}.
\end{proof}

\subsection{Determinant seminorm}
\label{SubSec:determinant norm}
Let $V$ be a finite-dimensional vector space over $k$. Recall that the \emph{determinant}\index{determinant} of $V$ is defined as the maximal exterior power $\Lambda^rV$ of the vector space $V$, where $r$ is the rank of $V$ over $k$. It is a quotient space of rank $1$ of $V^{\otimes r}$. We denote by $\det(V)$ the determinant of $V$.

\begin{defi}\label{Def:determinant:norm}
Assume that the vector space $V$ is equipped with a seminorm $\norm{\ndot}$. We call the \emph{determinant seminorm of $\norm{\ndot}$}\index{determinant seminorm}\index{seminorm!determinant ---} on $\det(V)$ and we denote by $\norm{\ndot}_{\det}$ the $\pi$-exterior power seminorm of $\norm{\ndot}$, that is, quotient seminorm induced by the $\pi$-tensor power seminorm on $V^{\otimes r}$.
\end{defi}

\begin{prop}[Hadamard's inequality]\label{Pro:hadamard}
Let $(V,\norm{\ndot})$ be a finite-dimensional seminormed vector space of rank $r>0$ over $k$. For any $\eta\in\det(V)$,
\[\|\eta\|_{\det}=\inf\big\{\|x_1\|\cdots\|x_r\|\,:\,\eta=x_1\wedge\cdots\wedge x_r\big\}.\] 
In particular, the determinant seminorm is a norm if and only if $\norm{\ndot}$ is a norm.
\end{prop}
\begin{proof} If $\eta$ is written in the form $\eta=x_1\wedge\cdots\wedge x_r$, where $x_1,\ldots,x_r$ are elements in $V$, then it is the image of $x_1\otimes\cdots\otimes x_r$ by the canonical projection $V^{\otimes r}\rightarrow\det(V)$. Therefore one has
$\|\eta\|_{\det}\leqslant\|x_1\|\cdots\|x_r\|$. Thus we obtain 
\[\|\eta\|_{\det}\leqslant\inf\big\{\|x_1\|\cdots\|x_r\|\,:\,\eta=x_1\wedge\cdots\wedge x_r\big\}.\]

In the following, we prove the converse inequality. It suffices to treat the case where $\eta\neq 0$. By definition one has
\[\|\eta\|_{\det}=\inf\Big\{\sum_{i=1}^N\|x_1^{(i)}\|\cdots\|x_r^{(i)}\|\,:\,\eta=\sum_{i=1}^Nx_1^{(i)}\wedge\cdots\wedge x_r^{(i)}\Big\}.\]
Let $\{x_j^{(i)}\}_{i\in\{1,\ldots,N\},\,j\in\{1,\ldots,r\}}$ be elements in $V$ such that
$\eta=\sum_{i=1}^Nx_1^{(i)}\wedge\cdots\wedge x_r^{(i)}$.
Let $\{e_j\}_{j=1}^r$ be a basis of $V$ and $\eta_0=e_1\wedge\cdots\wedge e_r$. For any $i\in\{1,\ldots,N\}$, there exists $a_i\in k$ such that $x_1^{(i)}\wedge\cdots\wedge x_r^{(i)}=a_i\eta_0$. Without loss of generality, we may assume that all $a_i$ are non-zero and that
\[\frac{\|x_1^{(1)}\|\cdots\|x_r^{(1)}\|}{|a_1|}=\min_{i\in\{1,\ldots,N\}}\frac{\|x_1^{(i)}\|\cdots\|x_r^{(i)}\|}{|a_i|}.\]
Note that one has
\[\eta=(a_1+\cdots+a_N)\eta_0=\Big(1+\frac{a_2}{a_1}+\cdots+\frac{a_N}{a_1}\Big)x_1^{(1)}\wedge\cdots\wedge x_r^{(1)},\]
and 
\[\begin{split}&\quad\;\Big|1+\frac{a_2}{a_1}+\cdots+\frac{a_N}{a_1}\Big|\cdot\|x_1^{(1)}\|\cdots\|x_r^{(1)}\|\\ &\leqslant \Big(1+\Big|\frac{a_2}{a_1}\Big|+\cdots+\Big|\frac{a_N}{a_1}\Big|\Big)\|x_1^{(1)}\|\cdots\|x_r^{(1)}\|\leqslant\sum_{i=1}^N\|x_1^{(i)}\|\cdots\|x_r^{(i)}\|.
\end{split}\]
The proposition is thus proved.
\end{proof}

\begin{rema}\label{Rem:Hadamard basis}
Let $(V,\norm{\ndot})$ be a {non-zero} finite-dimensional normed vector space over $k$. Let $r$ be the rank of $V$ {over $k$}. Proposition \ref{Pro:hadamard} shows that 
\begin{equation}\label{Equ:DeltaV}\inf\bigg\{\frac{\|x_1\|\cdots\|x_r\|}{\|x_1\wedge\cdots\wedge x_r\|_{\det}}\,:\,(x_1,\ldots,x_r)\in V^r,\;x_1\wedge\cdots\wedge x_r\neq 0\bigg\}=1.\end{equation}
If the infimum is attained by some $(e_1,\ldots,e_r)\in V^r$, then $\{e_i\}_{i=1}^r$ is called an \emph{Hadamard basis}\index{Hadamard basis} of $(V,\norm{\ndot})$. {By convention, the empty subset of the zero normed vector space is considered as an Hadamard basis.} 
\end{rema}

\begin{coro}\label{Cor:exactsequencenorm}
Let $V$ be a finite-dimensional seminormed vector space over $k$ and $W$ be a vector subspace of $V$. Twhe canonical isomorphism (see \cite{Bourbaki70} Chapter III, \S 7, no.7)
\begin{equation}\label{Equ:isocanoniquesuiteexacte}\det(W)\otimes\det(V/W)\longrightarrow\det(V)\end{equation}
has seminorm $\leqslant 1$, where we consider the determinant seminorm of the induced seminorm on $\det(W)$ and that of the quotient seminorm on $\det(V/W)$, and the  tensor product seminorm on $\det(W)\otimes\det(V/W)$ (see Remark \ref{Rem:produittenrk1}).
\end{coro}
\begin{proof}
Let $\{x_1,\ldots,x_n\}$ be a basis of $W$ and $\{y_1,\ldots,y_m\}$ be elements in $V\setminus W$ whose image in $V/W$ forms a basis of $V/W$. By Proposition \ref{Pro:hadamard} one has
\[\|x_1\wedge\cdots\wedge x_n\wedge y_1\wedge\cdots\wedge y_m\|_{\det}\leqslant \|x_1\|\cdots\|x_n\|\cdot\|y_1\|\cdots\|y_m\|.\]
Note that if we replace each $y_i$ by an element $y_i'$ in the same equivalent class, one has
\[x_1\wedge\cdots\wedge x_n\wedge y_1\wedge\cdots\wedge y_m=x_1\wedge\cdots\wedge x_n\wedge y_1'\wedge\cdots\wedge y_m'.\]
Hence  we obtain
\[\|x_1\wedge\cdots\wedge x_n\wedge y_1\wedge\cdots\wedge y_m\|\leqslant \|x_1\|\cdots\|x_n\|\cdot\|[y_1]\|\cdots\|[ y_m]\|.\]
{Therefore, for any $\eta\in\det(W)$ and $\eta'\in\det(V/W)$ one has
\[\norm{\eta\wedge\eta'}_{\det}\leqslant \Big(\inf_{\begin{subarray}{c}(x_1,\ldots,x_n)\in W^n\\x_1\wedge\cdots\wedge x_n=\eta\end{subarray}}\norm{x_1}\cdots\norm{x_n}\Big)\Big(\inf_{\begin{subarray}{c}(y_1,\ldots,y_m)\in (V\setminus W)^m\\
{[y_1]}\wedge\cdots\wedge[y_m]=\eta'\end{subarray}}\norm{[y_1]}\cdots\norm{[y_m]}\Big),\] 
which leads to, by Proposition \ref{Pro:hadamard}, the inequality $\norm{\eta\wedge\eta'}_{\det}\leqslant\norm{\eta}_{\det}\cdot\norm{\eta'}_{\det}$.}
\end{proof}

\begin{prop}\label{Pro: tensor product and deteminant pi}
Let $(V,\norm{\ndot}_V)$ and $(W,\norm{\ndot}_W)$ be finite-dimensional seminormed vector spaces over $k$, and $n$ and $m$ be respectively the ranks of $V$ and $W$ over $k$. We equip $V\otimes_kW$ with the $\pi$-tensor product seminorm $\norm{\ndot}_{\pi}$.  Then the natural $k$-linear isomorphism $\det(V\otimes_kW)\cong\det(V)^{\otimes m}\otimes_k\det(W)^{\otimes n}$ is an isometry, where we consider the determinant seminorm of $\norm{\ndot}_{\pi}$ on $\det(V\otimes_kW)$ and the tensor product of determinant seminorms on $\det(V)^{\otimes m}\otimes_k\det(W)^{\otimes n}$.
\end{prop}
\begin{proof}
Let $\norm{\ndot}'$ be the seminorm on $\det(V)^{\otimes m}\otimes\det(W)^{\otimes n}$ given by tensor product of determinant seminorms. By Proposition \ref{Pro:quotientavecpitensor}, the seminorm $\norm{\ndot}'$ identifies with the quotient of the $\pi$-tensor {power on $(V\otimes_kW)^{\otimes nm}$ of the seminorm $\norm{\ndot}_{\pi}$ on $V\otimes_kW$.} In other words, $\norm{\ndot}'$ identifies with $\norm{\ndot}_{\pi,\det}$.
\end{proof}

\begin{prop}\label{Pro: determinant of exterior power}
Let $(V,\norm{\ndot})$ be a finite-dimensional seminormed vector space over $k$ and $r$ be the rank of $V$ over $k$. Let $i$ be a positive  
integer. Then the canonical $k$-linear isomorphism $\det(\Lambda^{i}V)\rightarrow\det(V)^{\otimes\binom{r-1}{i-1}}$ is an isometry, where we consider the $i^{\text{th}}$ $\pi$-exterior power seminorm on $\Lambda^iV$.  
\end{prop}
\begin{proof}
Consider the following commutative diagram
\[\xymatrix{V^{\otimes i\binom{r}{i}}\ar@{->>}[r]^-{p_1}\ar@{->>}[d]_-{p_2}&\det(V)^{\otimes\binom{r-1}{i-1}}
\\(\Lambda^iV)^{\otimes \binom{r}{i}}\ar@{->>}[r]_{p_3}&\det(\Lambda^iV)\ar[u]_-\simeq}\]
By definition, if we equip $V^{\otimes i\binom{r}{i}}$ with the $\pi$-tensor product seminorm, then its quotient seminorm on $(\Lambda^iV)^{\otimes\binom{r}{i}}$ identifies with the $\pi$-tensor product of the $\pi$-exterior power seminorm. Moreover, by Proposition \ref{Pro:quotientavecpitensor}, the quotient seminorm on $\det(\Lambda^iV)$ (induced by $p_3$) of the tensor product of the $\pi$-exterior power seminorm identifies with the determinant seminorm of the latter. Still by the same proposition, the quotient seminorm on $\det(V)^{\otimes\binom{r-1}{i-1}}$ induced by $p_1$ identifies with the tensor power of the determinant seminorm. Therefore the natural isomorphism $\det(\Lambda^{i}V)\rightarrow\det(V)^{\otimes\binom{r-1}{i-1}}$ preserves actually the seminorms by using \ref{Item: successive quotient seminorm} in Proposition~\ref{prop:quotient:norm:linear:map}. 
\end{proof}

\subsection{Seminormed graded algebra}\label{Subsec:graded:algebra:norms}

Let $R_\sbullet=\bigoplus_{n\in\mathbb N}R_n$ be a graded $k$-algebra such that,
for any $n\in\mathbb N$, $R_n$ is of finite rank over $k$. 
For any $n\in\mathbb N$, let $\norm{\ndot}_{n}$ be a seminorm on $R_n$. 
We say that $\overline R_{\sbullet}=\{(R_n,\norm{\ndot}_n)\}_{n\in\mathbb N}$ is a 
\emph{seminormed graded algebra over $k$}\index{seminormed graded algebra} if the following submultiplicativity condition is satisfied: for any $(n,m)\in\mathbb N^2$ and any $(a,b)\in R_{n}\times R_{m}$, one has
\[\|a\cdot b\|_{n+m}\leqslant\|a\|_{n}\cdot\|b\|_{m}.\]
Furthermore, we say that $\overline R_\sbullet$ is \emph{of finite type}\index{of finite type}\index{seminormed graded algebra!--- of finite type} if the underlying graded $k$-algebra $R_\sbullet$ is of finite type over $k$.

Let $M_\sbullet=\bigoplus_{m\in\mathbb Z}M_m$ be a $\mathbb Z$-graded $k$-linear space and $h$ be a positive integer. We say that $M_\sbullet$ is an \emph{$h$-graded $R_\sbullet$-module}\index{h-graded R-module@$h$-graded $R_\sbullet$-module} if $M_\sbullet$ is equipped with a structure of $R_\sbullet$-module such that
\[\forall\,(n,m)\in\mathbb N\times\mathbb Z,\quad\forall\, (a,x)\in R_n\times M_m,\quad ax\in M_{nh+m}.\]
Let $M_\sbullet$ be an $h$-graded $R_\sbullet$-module. Assume that each homogeneous component $M_m$ is of finite rank over $k$ and is equipped with a seminorm $\norm{\ndot}_{M_m}$. 
We say that $\overline{M}_\sbullet=\big\{(M_m,\norm{\ndot}_{M_m}) \big\}_{m\in\mathbb Z}$ is a \emph{seminormed $h$-graded $\overline R_\sbullet$-module}\index{seminormed h-graded R-module@seminormed $h$-graded $\overline R_\sbullet$-module} if the following condition is satisfied: for any $(n,m)\in\mathbb N\times\mathbb Z$ and any $(a,x)\in R_{n}\times M_{m}$, one has
\[\|a\cdot x\|_{M_{nh+m}}\leqslant\|a\|_{n}\cdot\|x\|_{M_m}.\] 
We say that an $h$-graded $\overline R_\sbullet$-module $\overline M_\sbullet$ is \emph{of finite type}\index{of finite type}
\index{seminormed h-graded R-module@seminormed $h$-graded $\overline R_\sbullet$-module!--- of finite type} if the underlying $h$-graded $R_\sbullet$-module $M_\sbullet$ is of finite type.

\begin{prop}\label{prop:norm:graded:module:quotient}
Let $\overline R_{\sbullet} = \{(R_n,\norm{\ndot}_n)\}_{n\in\mathbb N}$ 
be a seminormed graded algebra over $k$.
Let $I_{\sbullet}$ be a homogenous ideal of $R_{\sbullet}$ and
$R'_{\sbullet} := R_{\sbullet}/I_{\sbullet}$.
\begin{enumerate}[label=\rm(\arabic*)]
\item\label{Item: quotient algebra seminorms}
Let $\norm{\ndot}'_n$ be the quotient seminorm on $R'_n$ induced by $\norm{\ndot}_n$ and $R_n \to R'_n$.
Then ${\overline R}{}_{\sbullet}' = \{ (R'_n,\norm{\ndot}'_n)\}_{n\in\mathbb N}$ forms a seminormed graded
algebra over $k$.

\item\label{Item: quotient module seminorm}
Let $\overline M_{\sbullet} = \{ (M_m, \norm{\ndot}_{M_m}) \}_{m \in \mathbb N}$ be a normed
$h$-graded $\overline R_{\sbullet}$-module and
$f_{\sbullet} : M_{\sbullet} \to N_{\sbullet}$ be a homomorphism of $h$-graded modules over $R_{\sbullet}$.\footnote{That is, for each $m\in\mathbb Z$,  $f_m : M_m \to N_m$ is a $k$-linear map such that {$f_{nh + m}(a \cdot x) = [a] \cdot f_m(x)$} for all $a \in R_n$ and $x \in M_m$.}
We assume that $I_{\sbullet} \cdot N_{\sbullet} = 0$ and
$f_m : M_m \to N_m$ is surjective for all $m \in \mathbb Z$.
Let $\norm{\ndot}_{N_m}$ be the quotient seminorm on $N_m$ induced by $\norm{\ndot}_{M_m}$ and $f_m$.
Then
$\overline N_{\sbullet} = \{ (N_m, \norm{\ndot}_{N_m}) \}_{m \in \mathbb N}$ forms
a seminormed $h$-graded ${\overline R}{}'_\sbullet$-module.
\end{enumerate}
\end{prop}

\begin{proof}
First let us see the following:
\begin{equation}\label{eqn:prop:norm:graded:module:quotient:01}
\forall\, (n, m) \in \mathbb N \times \mathbb Z,\, (a', y) \in R'_n \times N_m,\quad
\| a' \cdot y \|_{N_{nh + m}} \leqslant \| a' \|'_n \cdot \| y \|_{N_m}.
\end{equation}
Indeed, for a fixed positive number $\epsilon$, one can find $a \in R_n$ and
$x \in M_m$ such that 
\[
\begin{cases}
{[a]=a'}, & \| a \|_{n} \leqslant \mathrm{e}^{\epsilon} \|a'\|'_{n},\\
f_m(x) = y, & \| x \|_{M_m} \leqslant \mathrm{e}^{\epsilon} \|y \|_{N_m}.
\end{cases}
\]
Then, as $f_m(a\cdot x) = a' \cdot y$,
\[
\| a' \cdot y \|_{N_{nh + m}} \leqslant \| a \cdot x \|_{M_{nh + m}} \leqslant
\| a \|_{n}\cdot \| x \|_{M_m} \leqslant \mathrm{e}^{2\epsilon} \| a' \|'_{n} \cdot\|  y \|_{N_m},
\]
which implies \eqref{eqn:prop:norm:graded:module:quotient:01} because $\epsilon$ is an arbitrary positive number. Applying \eqref{eqn:prop:norm:graded:module:quotient:01} to the case where
$\overline M_{\sbullet} = \overline R_{\sbullet}$ 
and $\overline N_{\sbullet} = \overline R{}'_{\sbullet}$, one has
\begin{equation}\label{eqn:prop:norm:graded:module:quotient:02}
\forall\, (n, n') \in \mathbb N^2,\quad\forall\,  (a', b')  \in R'_n \times R'_{n'},\quad
\| a' \cdot b' \|'_{n+n'} \leqslant \| a' \|'_n \cdot \| b' \|'_{n'}.
\end{equation}
Thus \ref{Item: quotient algebra seminorms} is proved, so that \ref{Item: quotient module seminorm} is also proved by \eqref{eqn:prop:norm:graded:module:quotient:01}.
\end{proof}

\subsection{Norm of polynomial}\label{Subsec:Norm:polynomial}
Let $k[X]$ be the polynomial ring of one variable over $k$.
For $f = a_n X^n + \cdots + a_1 X + a_0 \in k[X]$,
{We define $\| f \|$ to be}
\[
\| f \| := \max \{ |a_n|, \ldots, |a_1|, |a_0| \}.
\]
It is easy to see that $\|\ndot\|$ yields a norm of $k[X]$ over $k$.

\begin{prop}\label{prop:norm:polynomial}
For $f, g \in k[X]$, one has the following:
\begin{enumerate}[label=\rm(\arabic*)]
\item\label{Item: norm f g less than norm f times norm g}
If the absolute value of $k$ is Archimedean, then 
\[
\norm{f g} \leqslant \min \{ \deg(f) + 1, \deg(g) + 1 \} \norm{ f}\cdot \norm{ g},
\]
where the degree of the zero polynomial is defined to be $-1$ {by convention}.

\item\label{Item: norm fg equals norm f times norm g}
If the absolute value of $k$ is non-Archimedean, then $\norm{ f g} = \norm{ f}\cdot \norm{ g}$.
\end{enumerate}
\end{prop}

\begin{proof}
Clearly we may assume that $f \not= 0$, $g \not= 0$ and $\deg(f) \leqslant \deg(g)$. We set
\[
\begin{cases}
f = a_n X^n + \cdots + a_1 X + a_0, \\
g = b_m X^{m} + \cdots + b_1 X + b_0, \\
fg = c_{n+m} X^{n+m} + \cdots + c_1 X + c_0,
\end{cases}
\]
where $n = \deg(f)$ and $m = \deg(g)$. Then
\[
c_{l} = \sum_{(i, j) \in \Delta(l)} a_i b_j, 
\]
where \[\Delta(l) = \Big\{ (i, j) : i+j =l,\; i\in\{0,\ldots,n\},\; j\in\{0,\ldots, m\} \Big\},\]
so that, as $\card(\Delta(l)) \leqslant n+1$, one has
\[
| c_l | \leqslant \begin{cases}
{\displaystyle \sum_{(i,j) \in \Delta(l)} |a_i|\cdot |b_j| \leqslant (n+1) \| f \|\cdot \| g \|} & (\text{Archimedean case}), \\[1ex]
{\displaystyle \max_{(i,j) \in \Delta(l)} \{ |a_i|\cdot |b_j| \} \leqslant  \| f \|\cdot \| g \|} & (\text{non-Archimedean case}).
\end{cases}
\]
Thus \ref{Item: norm f g less than norm f times norm g} and the inequality $\| f g \| \leqslant \| f \|\cdot \| g \|$ in {the non-Archimedean case} are obtained.

Finally let us consider the converse inequality in {the non-Archimedean case}.
We set 
\[
\alpha = \min \{ i : |a_i| = \| f \| \}\quad\text{and}\quad
\beta = \min \{ j : |b_j| = \| g \|\}.
\]
Note that if $i+j = \alpha + \beta$ and $(i, j) \not= (\alpha, \beta)$, then 
$|a_i|\cdot |b_j| < \| f \|\cdot \|g \|$ because
either
$i < \alpha$ or $j < \beta$. Therefore,
$|c_{\alpha+\beta}| = \| f \|\cdot\| g \|$ by Proposition~\ref{Pro:valeur}
and hence $\| f g \| \geqslant \| f \|\cdot\| g \|$.
\end{proof}

\section{Orthogonality}

The orthogonality of bases plays an important role in the study of finite-dimensional normed vector spaces. In the classic functional analysis over $\mathbb R$ or $\mathbb C$, the orthogonality often refers to a property related to an inner product. This property actually has an equivalent form, which has an analogue in the non-Archimedean case. However, in a finite-dimensional normed vector space over a non-Archimedean valued field, there may not exist an orthogonal basis. One can remedy this problem by introducing an approximative variant of the orthogonality. This technic is useful in the study of determinant norms. 

\subsection{Inner product}\label{Subsec:innerproduct} In this subsection, we assume that the absolute value $|\ndot|$ is Archimedean. In this case the field $k$ is either $\mathbb R$ or $\mathbb C$ {and we assume that $|\ndot|$ is the usual absolute value}.

Let $V$ be a vector space over $k$. A map $\emptyinnprod
: V \times V \to k$ is called  a \emph{semidefinite inner product}\index{semidefinite inner product} on $V$
if the following conditions are satisfied:
\begin{enumerate}[label=(\roman*)]

\item\label{Item: linearlity in the second ariable}
$\langle x, ay + bz\rangle = a \langle x, y\rangle + b \langle x, z\rangle$ for all $(x, y, z) \in V^3$ and $(a, b) \in k^2$.

\item\label{Item: conj symmetry}
$\langle x, y \rangle = \overline{\langle y, x \rangle}$ for any $(x, y) \in V^2$,
where $\overline{\langle y, x \rangle}$ is the complex conjugation of $\langle y, x \rangle$.

\item
$\langle x, x \rangle \in \mathbb R_{\geqslant 0}$ for {any} $x \in V$.
\end{enumerate}
If $\langle x,x\rangle>0$ for any $x\in V\setminus\{0\}$, we just say that $\emptyinnprod$ is an \emph{inner product}\index{inner product}. Namely, an inner product means either a scalar product or a Hermitian product according to
$k = \mathbb R \text{ or } \mathbb C$. Note that the  semidefinite inner product $\emptyinnprod$
induces a seminorm $\norm{\ndot}$ on $V$ such that
$\|x\|=\langle x,x\rangle^{1/2}$ for any $x\in V$.

\begin{prop}
Let $V$ be a vector space over $k$, $\emptyinnprod$ be a  semidefinite inner product on $V$ and $\norm{\ndot}$ be the seminorm induced by $\emptyinnprod$. 
\begin{enumerate}[label=\rm(\arabic*)]
\item\label{Item: vanishing of N} For any $x\in N_{\norm{\ndot}}$ and any $y\in V$ one has $\langle x,y\rangle=\langle y,x\rangle =0$.
\item\label{Item: passage to quotient} The  semidefinite inner product $\emptyinnprod$ induces by passing to quotient an inner product $\emptyinnprod^\sim$ on $V/N_{\norm{\ndot}}$ such that \[\forall\,(x,y)\in V^2,\quad \langle [x],[y]\rangle^\sim=\langle x,y\rangle,\]
where $[x]$ and $[y]$ are the classes of $x$ and $y$ in $V/N_{\norm{\ndot}}$, respectively.
Moreover, one has $\left(\norm{\alpha}^\sim\right)^2=\langle\alpha,\alpha\rangle^\sim$ for any $\alpha\in V/N_{\norm{\ndot}}$.
\item\label{Item: Riesz representation} Assume that $V$ is of finite rank over $k$. For any bounded linear form {$f$ on $V$} there exists an element $y$ in $V$ such that $f(x)=\langle y,x\rangle$ for any $x\in V$. Moreover, the element $y$ is unique up to addition by an element in $N_{\norm{\ndot}}$. 
\end{enumerate}
\end{prop}
\begin{proof}

\ref{Item: vanishing of N} By Cauchy-Schwarz inequality, one has
$|\langle x,y\rangle|^2\leqslant \norm{x}^2\cdot\norm{y}^2=0$.
Hence $\langle x,y\rangle=0$. Similarly, $\langle y,x\rangle=0$.

\ref{Item: passage to quotient} By \ref{Item: vanishing of N} and the properties \ref{Item: linearlity in the second ariable} and \ref{Item: conj symmetry} of  semidefinite inner product, we obtain that, if $x$, $x'$, $y$ and $y'$ are vectors in $V$ such that $x-x'\in N_{\norm{\ndot}}$ and $y-y'\in N_{\norm{\ndot}}$, then $\langle x, y \rangle = \langle x', y' \rangle$. Therefore the  semidefinite inner product $\emptyinnprod$ induces by passing to quotient a function \[\emptyinnprod^\sim:(V/N_{\norm{\ndot}})\times(V/N_{\norm{\ndot}})\rightarrow k.\] From the definition it is straightforward to check that $\emptyinnprod^\sim$ is a  semidefinite inner product and $\langle\alpha,\alpha\rangle^\sim=(\norm{\alpha}^\sim)^{2}$ for any $\alpha\in V/N_{\norm{\ndot}}$. It remains to verify that $\emptyinnprod^\sim$ is definite. Let $x$ be an element in $V$ such that $\langle [x],[x]\rangle^\sim=0$. Then one has $\langle x,x\rangle=\norm{x}^2=0$. Hence $\norm{x}=0$, namely $x\in N_{\norm{\ndot}}$.
 
\ref{Item: Riesz representation} Since $f$ is a bounded linear form, it vanishes on $N_{\norm{\ndot}}$. Hence there exists a unique linear form $\widetilde f:V/N_{\norm{\ndot}}\rightarrow k$ such that $\widetilde f\scirc\pi=f$, where $\pi:V\rightarrow V/N_{\norm{\ndot}}$ is the projection map. Moreover, by Riesz's representation theorem for usual finite-dimensional inner product space, there exists a unique $\beta\in V/N_{\norm{\ndot}}$ such that $\widetilde f(\alpha)=\langle\beta,\alpha\rangle^{\sim}$ for any $\alpha\in V/N_{\norm{\ndot}}$. Hence we obtain that the equivalence class $\beta$ equals the set of $y\in V$ {such that} $f(x)=\langle y,x\rangle$ for any $x\in V$.
\end{proof}

Let $V$ be a finite-dimensional vector space over $k$ equipped with a seminorm $\norm{\ndot}$. We say that the seminorm $\norm{\ndot}$ is \emph{Euclidean}\index{Euclidean}\index{seminorm!Euclidean ---} (resp. \emph{Hermitian}\index{Hermitian}\index{seminorm!Hermitian ---}) if $k=\mathbb R$ (resp. $k=\mathbb C$) and if the seminorm $\norm{\ndot}$ is induced by a  semidefinite inner product.  Note that if a seminorm $\norm{\ndot}$ on $V$ is Euclidean (resp. Hermitian), then also is its {dual norm} on $V^*$. In fact, if 
$\emptyinnprod$ 
is a  semidefinite inner product on $V$ and $\norm{\ndot}$ is the corresponding seminorm, then it induces (by Riesz's representation theorem) an $\mathbb R$-linear isometry $\iota:(V/N_{\norm{\ndot}},\norm{\ndot}^\sim)\rightarrow (V^*,\norm{\ndot}_*)$ such that
\[\forall\,(x,y)\in V^2,\quad\iota([x])(y)=\langle x,y\rangle.\]
Moreover, for $a\in k$ and $x\in V$ one has
$\iota(ax)=\overline{a}\,\iota(x)$.
Then the dual norm on $V^*$ is induced by the following inner product $\emptyinnprod_*$: 
\[\forall\,(\alpha,\beta)\in (V^*)^2,\quad\langle\alpha,\beta\rangle_*=\overline{\langle \iota^{-1}(\alpha),\iota^{-1}(\beta)\rangle^{\sim}}.\]

\begin{rema}
Let $\psi:[0,1]\rightarrow[0,1]$ be the function $t\mapsto (t^2+(1-t)^2)^{1/2}$. If $V$ and $W$ are finite-dimensional vector spaces over $k$ equipped with  semidefinite inner products, then the direct sum seminorm $\norm{\ndot}_\psi$ on $V\oplus W$ as constructed in \S\ref{Subsec:directsums} is induced by the  semidefinite inner product on $V\oplus W$ defined as
$\big\langle (x,y),(x',y')\big\rangle:=\langle x,x'\rangle+\langle y,y'\rangle$.
The seminorm $\norm{\ndot}_\psi$ is called the \emph{orthogonal direct sum}\index{orthogonal direct sum}\index{seminorm!orthogonal direct sum} of the seminorms on $V$ and $W$ corresponding to their  semidefinite inner products.
\end{rema}

\subsection{Orthogonal basis of an inner product}\label{Subsec:orthogonalinnerproduct}In this subsection, we assume that the absolute value $|\ndot|$ is Archimedean. Let $V$ be a finite-dimensional vector space over $k$ equipped with a  semidefinite inner product $\emptyinnprod$. {Let $\norm{\ndot}$ be the seminorm induced by $\emptyinnprod$.} 
We say that a basis $\{e_1,\ldots,e_r\}$ of $V$ is \emph{orthogonal}\index{orthogonal basis} if $\langle e_i,e_j\rangle=0$ for distinct indices $i$ and $j$ in $\{1,\ldots,r\}$. If in addition $\langle e_i,e_i\rangle=1$ for any $i\in\{1,\ldots,r\}$ {such that $e_i\in V\setminus N_{\norm{\ndot}}$}, we say that $\{e_1,\ldots,e_r\}$ is an \emph{orthonormal basis}\index{orthonormal basis}. Note that, if $\{e_1,\ldots,e_r\}$ is an orthogonal basis, then\begin{equation}\label{Equ:orthogonalcond}\forall\,(\lambda_1,\ldots,\lambda_r)\in k^r,\quad\|\lambda_1e_1+\cdots+\lambda_re_r\|^{2}=\sum_{i=1}^r|\lambda_i|^{2}\cdot\|e_i\|^{2}.\end{equation}
Moreover, by the Gram-Schmidt process, there always exists an {orthonormal}
basis of $V$ (cf. the proof of Proposition~\ref{Pro:existenceepsorth}).

The following proposition provides an alternative form for the orthogonality condition of a basis in a finite-dimensional vector space equipped with a  semidefinite inner product.

\begin{prop}\label{Pro:orthogonality}
Let $V$ be a finite-dimensional vector space over $k$, equipped with a  semidefinite inner product $\emptyinnprod$. 
Let $\{e_i\}_{i=1}^r$ be a basis of $V$. Then it is an orthogonal basis if and only if the following condition is satisfied:
\begin{equation}\label{Equ:weaklyorhtgonal}\forall\,(\lambda_1,\ldots,\lambda_r)\in k^r,\quad \|\lambda_1e_1+\cdots+\lambda_re_r\|\geqslant\max_{i\in\{1,\ldots,r\}}\|\lambda_ie_i\|.\end{equation}
\end{prop}
\begin{proof}
If $\{e_i\}_{i=1}^r$ is an orthogonal basis of $V$, then by \eqref{Equ:orthogonalcond} we obtain that the inequality \eqref{Equ:weaklyorhtgonal} holds. Conversely, assume given a basis $\{e_i\}_{i=1}^r$ of $V$ which verifies the condition \eqref{Equ:weaklyorhtgonal}. Then for any $(\lambda_1,\ldots,\lambda_{r-1})\in k^{r-1}$, one has
\[\|\lambda_1e_1+\cdots+\lambda_{r-1}e_{r-1}+e_r\|\geqslant\|e_r\|,\]
which implies that $e_r$ is orthogonal to the vector subspace generated by $e_1,\ldots,e_{r-1}$. 
 Indeed, $\|(\pm\epsilon)e_i + e_r\| \geqslant \| e_r \|$
for $\epsilon > 0$, which implies that $\epsilon \| e_i \|^2 \pm 2 \langle e_i, e_r \rangle \geqslant 0$, and hence
$\pm  \langle e_i, e_r \rangle \geqslant 0$ by taking {the limit when} $\epsilon \to 0$, as required. Therefore by
induction we obtain that the basis $\{e_i\}_{i=1}^r$ is an orthogonal basis.
\end{proof}

\subsection{Orthogonality in general cases}
In this subsection, we consider a general valued field $(k,|\ndot|)$, which is not necessarily Archimedean. 
Let $V$ be a finite-dimensional vector space over $k$ and $\norm{\ndot}$ be a seminorm on $V$. We say that a basis $\{e_i\}_{i=1}^r$ of $V$ is \emph{orthogonal}\index{orthogonal basis} if for any $(a_1,\ldots,a_r)\in k^r$ one has
\[\|a_1e_1+\cdots+a_re_r\|\geqslant\max_{i\in\{1,\ldots,r\}}\|a_ie_i\|.\] 
If in addition {$\|e_i\|=1$ for any $i\in\{1,\ldots,r\}$ such that $e_i\in V\setminus N_{\norm{\ndot}}$}, we say that the basis $\{e_i\}_{i=1}^r$ is \emph{orthonormal}\index{orthonormal basis}. We have seen in Proposition \ref{Pro:orthogonality} that this definition is equivalent to the definition in \S\ref{Subsec:orthogonalinnerproduct} when the absolute value $|\ndot|$ is Archimedean and the seminorm $\norm{\ndot}$ is induced by a  semidefinite inner product.

The existence of an orthogonal basis in the non-Archimedean case is not always true. We refer the readers to \cite[Example 2.3.26]{Perez-Garcia_Schikhof10} for a counter-example. Thus we need a refinement of the notion of orthogonality.

\begin{defi} \label{Def:alphaorthogonal}
Let $(V,\norm{\ndot})$ be a finite-dimensional seminormed vector space over $k$, and $\alpha\in \intervalle]01]$. We say that a basis $\{e_1,\ldots,e_r\}$ of $V$ is $\alpha$-\emph{orthogonal}\index{alpha-orthogonal basis@$\alpha$-orthogonal basis} if for any $(\lambda_1,\ldots,\lambda_r)\in k^r$ one has
\[\|\lambda_1e_1+\cdots+\lambda_re_r\|\geqslant\alpha\max(|\lambda_1|\cdot\|e_1\|,\ldots,|\lambda_r|\cdot\|e_r\|).\]
Note that the $1$-orthogonality is just the orthogonality defined in the beginning of the subsection. We refer the readers to \cite[\S2.3]{Perez-Garcia_Schikhof10} for more details about this notion.
\end{defi}

{\begin{prop}
\label{Pro: orthogonal basis: vanishing}
Let $(V,\norm{\ndot})$ be a finite-dimensional seminormed vector space over $k$, $\alpha$ be an element in $\intervalle{]}{0}{1}{]}$, and $\boldsymbol{e}=\{e_i\}_{i=1}^r$ be an $\alpha$-orthogonal basis of $(V,\norm{\ndot})$. Then the intersection of $\boldsymbol{e}$ with $N_{\norm{\ndot}}$ forms a basis of $N_{\norm{\ndot}}$.
\end{prop}
\begin{proof}
Without loss of generality, we assume that $\boldsymbol{e}\cap N_{\norm{\ndot}}=\{e_1,\ldots,e_n\}$, where $n\in\mathbb N$, $n\leqslant r$. Suppose that $N_{\norm{\ndot}}$ is not generated by $\boldsymbol{e}\cap N_{\norm{\ndot}}$, then there exists an element $x=\lambda_1e_1+\cdots+\lambda_re_r$ in $N_{\norm{\ndot}}$ which does not belong to the vector subspace of $V$ generated by $\boldsymbol{e}\cap N_{\norm{\ndot}}$. Therefore there exists $i\in\{n+1,\ldots,r\}$ such that $\lambda_i\neq 0$. Since the basis $\boldsymbol{e}$ is $\alpha$-orthogonal, one has
\[0=\norm{x}\geqslant \alpha |\lambda_i|\cdot\norm{e_i}>0,\]
which leads to a contradiction. 
\end{proof}}

\begin{prop}\label{Pro:orthogonal basis in sub and quotient spaces}
Let $(V,\norm{\ndot})$ be a finite-dimensional seminormed vector space over $k$, $\alpha\in\intervalle{]}{0}{1}{]}$ and $\boldsymbol{e}$ be an $\alpha$-orthogonal basis of $V$. 
Let $\boldsymbol{e}'$ be a subset of $\boldsymbol{e}$ and $W$ be the vector subspace of $V$ generated by all vectors in $\boldsymbol{e}'$.
\begin{enumerate}[label=\rm(\arabic*)]
\item\label{Item: subspace generated by orthogonal basis} 
The set $\boldsymbol{e}'$  
is an $\alpha$-orthogonal basis of $W$ with respect to the restriction of $\norm{\ndot}$ {to} $W$.
\item\label{Item: quotient orhtogonal basis}  
The image of $\boldsymbol{e}\cap(V\setminus W)$ in $V/W$ forms an $\alpha$-orthogonal basis of $V/W$ with respect to the quotient seminorm of $\norm{\ndot}$. {Moreover, for any $x\in\boldsymbol{e}\cap (V\setminus W)$, the quotient seminorm of the class of $x$ is bounded from below by $\alpha\norm{x}$. In particular, if}   $\alpha=1$, namely $\boldsymbol{e}$ is an orthogonal basis, then for any element $x\in\boldsymbol{e}\cap (V\setminus W)$, the quotient seminorm of the class of $x$ in $V/W$ is equal to $\|x\|$.
\end{enumerate}
\end{prop}
\begin{proof}
\ref{Item: subspace generated by orthogonal basis} Assume that 
$\boldsymbol{e}' =\{e_1,\ldots,e_n\}$. 
Since $W$ is generated by the vectors in $\boldsymbol{e}'$, $\{e_1,\ldots,e_n\}$ is a basis of $W$. Since $\boldsymbol{e}$ is an $\alpha$-orthogonal basis of $V$, for any $(\lambda_1,\ldots,\lambda_n)\in k^n$ one has
\[\|\lambda_1e_1+\cdots+\lambda_ne_n\|\geqslant\alpha\max\{|\lambda_1|\cdot\|e_1\|,\ldots,|\lambda_n|\cdot\|e_n\|\}.\]
Therefore $\{e_1,\ldots,e_n\}$ is an $\alpha$-orthogonal basis of $W$.

\ref{Item: quotient orhtogonal basis} Assume that 
$\boldsymbol{e}' =\{e_1,\ldots,e_n\}$ and $\boldsymbol{e}\cap(V\setminus W)=\{e_{n+1},\ldots,e_r\}$. It is clear that the canonical  image of $\boldsymbol{e}\cap(V\setminus W)$ in $V/W$ forms a basis of $V/W$. It remains to show that it is an $\alpha$-orthogonal basis. Let $\norm{\ndot}'$ be the quotient seminorm of $\norm{\ndot}$ on $V/W$. For any $i\in\{n+1,\ldots,r\}$, let $y_i$ be the canonical image of $e_i$ in $V/W$. Let $(\lambda_1,\ldots,\lambda_r)$ be an element in $k^r$. Since $\boldsymbol{e}$ is an $\alpha$-orthogonal basis of $(V,\norm{\ndot})$, for any $(\lambda_1,\ldots,\lambda_r)\in k^r$ one has
\[\|\lambda_1e_1+\cdots+\lambda_re_r\|\geqslant\alpha\max_{i\in\{1,\ldots,r\}}|\lambda_i|\cdot\|e_i\|\geqslant\alpha\max_{i\in\{n+1,\ldots,r\}}|\lambda_i|\cdot\|e_i\|.\]
Therefore, for any $(\lambda_{n+1},\ldots,\lambda_r)\in k^{r-n}$ one has (note that $y_j=[e_j]$)
\[\|\lambda_{n+1}y_{n+1}+\cdots+\lambda_ry_r\|'\geqslant \alpha\max_{i\in\{n+1,\ldots,r\}}|\lambda_i|\cdot\|e_i\|\geqslant\alpha\max_{i\in\{n+1,\ldots,r\}}|\lambda_i|\cdot\|y_i\|'.                                                                                                                                                                                                                                                                                                                                                                                                                                                                                                                                                                                                                                                                                                                                                                                                                                                                                                                                                                                                                                                                                                                                                                                                                                                                                                                                                                                                                                                                                                                                                                                                                                                                                                                                                                                                                                                                                                                                                                                                                                                                                                            \]
Hence $\{y_i\}_{i=n+1}^r$ is an $\alpha$-orthogonal basis of $V/W$. {The first inequality also implies that $\norm{y_i}'\geqslant\alpha\norm{e_i}$ for any $i\in\{n+1,\ldots,r\}$. If $\alpha=1$, for any $i\in\{n+1,\ldots,r\}$ one has}
\[\|y_i\|'\geqslant \|e_i\|\geqslant \|y_i\|',\]
which leads to the equality $\|y_i\|'=\|e_i\|$.
\end{proof}

The following proposition shows that, in the Archimedean case, any finite-dimensional normed vector space admits an orthogonal basis. In general case, for any $\alpha\in\intervalle{]}{0}{1}{[}$, any finite-dimensional normed vector space admits an $\alpha$-orthogonal basis.

\begin{prop}\label{Pro:existenceoforthogonal}
Let $(V,\norm{\ndot})$ be a finite-dimensional normed vector space over $k$. Then we have the following:
\begin{enumerate}[label=\rm(\arabic*)]
\item\label{Item: existence of alpha orthogonal basis} 
For any $\alpha\in\intervalle{]}{0}{1}{[}$, there exists an $\alpha$-orthogonal basis of $V$.

\item\label{Item: Hadamard basis orthogonal}
Any Hadamard basis of $V$ is orthogonal (see Remark~\ref{Rem:Hadamard basis}).

\item\label{Item: locally compact field existence of Hadamard basis}
If the field $k$ is locally compact and the absolute value $|\ndot|$ is not trivial, then $V$ admits an Hadamard basis, which is also an orthogonal basis. 
\end{enumerate}
\end{prop}
\begin{proof}
\ref{Item: existence of alpha orthogonal basis}, \ref{Item: Hadamard basis orthogonal} By Proposition~\ref{Pro:hadamard}, we can choose a basis $\boldsymbol{e}=\{e_i\}_{i=1}^r$ such that
\[
{\frac{\|e_1\|\cdots\|e_r\|}{\|e_1\wedge\cdots\wedge e_r\|} \leqslant \alpha^{-1}.}
\]
We claim that $\boldsymbol{e}=\{e_i\}_{i=1}^r$ is $\alpha$-orthogonal.
Let $(\lambda_1,\ldots,\lambda_r)$ be an element in $k^r$ and $x=\lambda_1e_1+\cdots+\lambda_re_r$. For any $i\in\{1,\ldots,r\}$, by Proposition \ref{Pro:hadamard} one has
\[\frac{\|e_1\|\cdots\|e_r\|}{\|e_1\wedge\cdots\wedge e_r\|}\leqslant\alpha^{-1}\frac{\|e_1\|\cdots\|e_{i-1}\|\cdot\|x\|\cdot\|e_{i+1}\|\cdots\|e_r\|}{\|e_1\wedge\cdots e_{i-1}\wedge x\wedge e_{i+1}\wedge\cdots \wedge e_r\|}.\]
Therefore one has $\|x\|\geqslant\alpha|\lambda_i|\cdot\|e_i\|$. Since $i\in\{1,\ldots,r\}$ is arbitrary, we obtain that $\{e_i\}_{i=1}^r$ is an $\alpha$-orthogonal basis. A similar argument also shows that an Hadamard basis  
is necessarily orthogonal.

\ref{Item: locally compact field existence of Hadamard basis} We assume that $k$ is locally compact and $|\ndot|$ is not trivial. Then, as $V$ is locally compact,
there is $a_0 \in k^{\times}$ such that $(V,\norm{\ndot})_{\leqslant |a_0|}$ is compact. 
By Proposition~\ref{Pro:dilatation}, if we choose $\lambda$ with $\lambda<\sup\{|a|\,:\,a\in k^{\times},\,|a|<1\}$, then,
for any $x \in V \setminus \{ 0 \}$, there is $b \in k^{\times}$ such that $\lambda \leqslant \| bx \| < 1$.
Here we set 
\[
C = \{ x \in V \,:\, \lambda |a_0| \leqslant \| x \| \leqslant |a_0| \},
\]
which is a compact set in $V$. For $(x_1,\ldots,x_r)\in (V\setminus\{0\})^r$, there are $b_1, \ldots, b_r \in k^{\times}$ such that
$\lambda \leqslant \| b_i x_i \| < 1$ for all $i$, so that
$(a_0b_1x_1, \ldots, a_0 b_r x_r) \in C^{r}$ and
\[
\frac{\|(a_0 b_1 x_1)\wedge\cdots\wedge (a_0 b_rx_r)\|}{\|a_0 b_1 x_1\|\cdots\|a_0 b_r x_r\|} = \frac{\|x_1\wedge\cdots\wedge x_r\|}{\|x_1\|\cdots\|x_r\|}.
\]
Hence the function
\[(x_1,\ldots,x_r)\longmapsto \frac{\|x_1\wedge\cdots\wedge x_r\|}{\|x_1\|\cdots\|x_r\|}\]
attains its maximal value on $(V\setminus\{0\})^r$, which is equal to $1$. The proposition is thus proved. 
\end{proof}

\begin{rema}
(1) In the case where $|\ndot|$ is trivial and $\norm{\ndot}$ is ultrametric, 
there is an orthogonal basis $\boldsymbol{e}$ for $\norm{\ndot}$ by Proposition~\ref{Pro:existenceepsorth}.
Thus, by Proposition~\ref{Pro:orthogonalesthadamard}, $\boldsymbol{e}$ is an Hadamard basis of $(V,\norm{\ndot})$.

\medskip
(2) We assume that $k$ is an infinite field and the absolute value $|\ndot|$ is trivial.
Fix a map $\lambda : {\mathbb P}^1(k) \to [\frac 12, 1]$.
Let $\pi : k^2 \setminus \{ (0,0) \} \to {\mathbb P}^1(k)$ be the natural map.
We set
\[
\forall\,x\in k^2,\quad
\| x \|_{\lambda} := 
\begin{cases}
\lambda(\pi(x)) & \text{if $x \not= (0,0)$}, \\
0 & \text{if $x = (0,0)$}.
\end{cases}
\]
It is easy to see that $\norm{\ndot}_{\lambda}$ satisfies
the axioms of norm:
(1) $\|ax\|_{\lambda}=|a|\cdot\|x\|_{\lambda}$; (2) $\|x+y\|_{\lambda}\leqslant\|x\|_{\lambda}+\|y\|_{\lambda}$;
(3) $\|x\|_{\lambda} = 0 \Longleftrightarrow x = 0$.
Choosing an infinite subset $S = \{ \zeta_1, \zeta_2, \ldots, \zeta_n, \ldots \}$ of $\mathbb P^1(k)$, we consider
$\lambda$ given by
\[
\lambda(\zeta) := \begin{cases}
\frac 12 + (\frac 12)^n & \text{if $\zeta \in S$ and $\zeta = \zeta_n$}, \\
1 & \text{otherwise}.
\end{cases}
\]
Then $\lambda(\mathbb P^1(k)) \subseteq \intervalle{]}{\frac 12}{1}{]}$, and
for any $\epsilon > 1/2$ there is $\zeta \in {\mathbb P}^1(k)$ with
$\lambda(\zeta) < \epsilon$.
Obviously $\{ \| x \|_{\lambda} \,:\, x \in k^2 \}$ is an infinite set, which means that
(2) in Corollary~\ref{Cor:finitevlaue} does not holds without the assumption that the norm is ultrametric.
Moreover, let us see that there is no orthogonal basis for $\norm{\ndot}_{\lambda}$.
Indeed, we assume that $\{e_1, e_2\}$ is an orthogonal basis for $\norm{\ndot}_{\lambda}$. 
By the property of $\lambda$,  
there is $x \in k^2 \setminus \{ (0,0) \}$ such that
$\| x \|_{\lambda} <  \min \{ \|e_1\|_{\lambda}, \|e_2\|_{\lambda}\}$. If we set $x = ae_1 + be_2$ with $(a, b) \not= (0,0)$, then
\[
\| x \|_{\lambda} \geqslant  \max \{ |a|\cdot\|e_1\|_{\lambda}, |b|\cdot\| e_2 \|_{\lambda} \} \geqslant \min \{ \|e_1\|_{\lambda}, \|e_2\|_{\lambda}\},
\]
which is a contradiction.
\end{rema}

\begin{coro}\label{Cor: existence of alpha orthogonal}
Let $(V,\norm{\ndot})$ be a finite-dimensional seminormed vector space over $k$. For any $\alpha\in \intervalle{]}{0}{1}{[}$, there exists an $\alpha$-orthogonal basis of $(V,\norm{\ndot})$. If the absolute value $|\ndot|$ is non-trivial and $(k,|\ndot|)$ is locally compact, then $(V,\norm{\ndot})$ admits an orthogonal basis.
\end{coro}
\begin{proof} If the absolute value $|\ndot|$ is non-trivial and $(k,|\ndot|)$ is locally compact, let $\alpha$ be an element in $\intervalle]01]$, otherwise let $\alpha$ be an element in $\intervalle]01[$.  
Let $W$ be the quotient vector space $V/N_{\norm{\ndot}}$, equipped with the quotient norm $\norm{\ndot}^\sim$. By Proposition \ref{Pro:existenceoforthogonal}, the normed vector space $(W,\norm{\ndot}^\sim)$ admits an $\alpha$-orthogonal basis $\{x_i\}_{i=1}^n$. For any $i\in\{1,\ldots,n\}$, let $e_i$ be an element in the class $x_i$. We also choose a basis $\{e_{j}\}_{j=n+1}^r$ of $N_{\norm{\ndot}}$. Hence $\{e_i\}_{i=1}^r$ becomes a basis of $V$. For any $(\lambda_1,\ldots,\lambda_r)\in k^r$ one has
\begin{align*}\norm{\lambda_1e_1+\cdots+\lambda_re_r} & =\norm{\lambda_1x_1+\cdots+\lambda_nx_n}^\sim\geqslant\alpha\max_{i\in\{1,\ldots,n\}}|\lambda_i|\cdot\norm{x_i}^\sim\\
&=\alpha\max_{i\in\{1,\ldots,n\}}|\lambda_i|\cdot\norm{e_i}=\alpha\max_{i\in\{1,\ldots,r\}}|\lambda_i|\cdot\norm{e_i}.
\end{align*}
\end{proof}

\begin{lemm}\label{Lem:normofdualbasis}
Let $(V,\norm{\ndot})$ be a finite-dimensional seminormed vector space over $k$, and $\alpha\in\intervalle]01]$. If $\boldsymbol{e}=\{e_i\}_{i=1}^r$ is an $\alpha$-orthogonal basis of $V$ and if $\{e_i^\vee\}_{i=1}^r$ is its dual basis, then, for any $i\in\{1,\ldots,r\}$, $e_i\not\in N_{\norm{\ndot}}$ if and only if $e_i^\vee\in V^*$, and in this case one has
\begin{equation}
\label{Equ:normeduale}1\leqslant\|e_i^\vee\|_*\cdot\|e_i\|\leqslant\alpha^{-1}.\end{equation}
\end{lemm}
\begin{proof}
The hypothesis that $e_i\not\in N_{\norm{\ndot}}$ actually implies that $e_i^\vee$ vanishes on $N_{\norm{\ndot}}$ since $N_{\norm{\ndot}}$ is generated by $\boldsymbol{e}\cap N_{\norm{\ndot}}$ (see Proposition \ref{Pro: orthogonal basis: vanishing}). By Corollary \ref{Coro:equivalenceofnrom}, $e_i^\vee$ belongs to $V^*$. Conversely, if $e_i$ belongs to $N_{\norm{\ndot}}$ then $e_i^{\vee}$ is not a bounded linear form on $V$ since it takes non-zero value on $e_i\in N_{\norm{\ndot}}$.

The first inequality of \eqref{Equ:normeduale} comes from the formula \eqref{Equ:dualnormformula} in \S\ref{Subsec:dualnorm}. In the following, we prove the second inequality. For any $(\lambda_1,\ldots,\lambda_r)\in k^r$ one has
\[e_i^{\vee}(\lambda_1e_1+\cdots+\lambda_re_r)=\lambda_i.\]
Hence
\begin{equation*}\|e_i^\vee\|_*=\sup_{\begin{subarray}{c}(\lambda_1,\ldots,\lambda_r)\in k^r\\
\lambda_i\neq 0\end{subarray}}\frac{|\lambda_i|}{\|\lambda_1e_1+\cdots+\lambda_re_r\|}\leqslant \alpha^{-1}\|e_i\|^{-1},\end{equation*}
where the inequality comes from the hypothesis that the basis $\{e_i\}_{i=1}^r$ is $\alpha$-orthogonal (so that $\|\lambda_1e_1+\cdots+\lambda_re_r\|\geqslant\alpha |\lambda_i|\cdot\|e_i\|$).
\end{proof}

\begin{prop}\label{Pro:alphaorthogonale}
Let $V$ be a finite-dimensional seminormed vector space over $k$ and $\alpha\in\intervalle]01]$. If $\{e_i\}_{i=1}^r$ is an $\alpha$-orthogonal basis
of $V$ and if $\{e_i^\vee\}_{i=1}^r$ is the dual basis of $\{e_i\}_{i=1}^r$, then $\{e_i^\vee\}_{i=1}^r\cap V^*$ is an $\alpha$-orthogonal basis of $V^*$.
Moreover, $\{e_i\}_{i=1}^r$ is an $\alpha$-orthogonal basis of $(V,\norm{\ndot}_{**})$ and one has
\begin{equation}\label{Equ:majoration par doubledual}
\alpha\|e_i\|\leqslant\|e_i\|_{**}\leqslant\|e_i\|.
\end{equation}
\end{prop}
\begin{proof}
By Proposition \ref{Pro: orthogonal basis: vanishing} {and Lemma \ref{Lem:normofdualbasis}}, the cardinal of $\{e_i^\vee\}_{i=1}^r\cap V^*$, which is equal to that of $\{e_i\}_{i=1}^r\cap(V\setminus N_{\norm{\ndot}})$, is $\dim_k(V^*)=\dim_k(V/N_{\norm{\ndot}})$. Therefore $\{e_i^\vee\}_{i=1}^r\cap V^*$ is a basis of $V^*$.

Consider $\xi=a_1e_1^\vee+\cdots+a_re_r^\vee$ in $V^*$. As $\xi(e_i)=a_i$ we get that 
\[\|\xi\|_*\geqslant \frac{|a_i|}{\|e_i\|}\geqslant\alpha |a_i|\cdot\|e_i^{\vee}\|_*\]
for any $i\in\{1,\ldots,r\}$ such that $\norm{e_i}\neq 0$,
where the second inequality comes from Lemma \ref{Lem:normofdualbasis}.
This implies that $\{e_i^\vee\}_{i=1}^r\cap V^*$ is an $\alpha$-orthogonal basis of $V^*$.  

Let $x=\lambda_1e_1+\cdots+\lambda_re_r$ be an element in $V$. Without loss of generality, we assume that $\boldsymbol{e}\cap N_{\norm{\ndot}}=\{e_{n+1},\ldots,e_r\}$. By definition, for $i\in\{1,\ldots,n\}$ one has
\begin{equation}\label{Equ:encadrementdoubledual}\norm{x}_{**}\geqslant\frac{|e_i^\vee(x)|}{\norm{e_i^\vee}_*}=\frac{|\lambda_i|}{\norm{e_i^\vee}_*}.\end{equation}
By Lemma \ref{Lem:normofdualbasis}, one has $1\leqslant\|e_i^\vee\|_*\cdot\|e_i\|\leqslant\alpha^{-1}$ and {hence \begin{equation}\label{Equ: alpah orthoogona of double dual}\norm{x}_{**}\geqslant \alpha|\lambda_i|\cdot\norm{e_i}\geqslant\alpha|\lambda_i|\cdot\norm{e_i}_{**},\end{equation}
where the last inequality comes from \eqref{Equ:doubledual}. For $i\in\{n+1,\ldots,r\}$ one has $\norm{e_i}_{**}=\norm{e_i}=0$. Hence
\[\norm{x}_{**}\geqslant\alpha\max_{i\in\{1,\ldots,r\}}|\lambda_i|\cdot\norm{e_i}_{**},\]
which shows that $\{e_i\}_{i=1}^r$ is an $\alpha$-orthogonal basis of $(V,\norm{\ndot}_{**})$. Moreover, \eqref{Equ: alpah orthoogona of double dual} also implies that $\norm{e_i}_{**}\geqslant\alpha\norm{e_i}$ for $i\in\{1,\ldots,n\}$, which, joint with the relation
\[\forall\,i\in\{n+1,\ldots,r\}, \quad \norm{e_i}=\norm{e_i}_{**}=0,\] 
leads to the first inequality of 
 \eqref{Equ:majoration par doubledual}.} The second {inequality} of \eqref{Equ:majoration par doubledual} comes from \eqref{Equ:doubledual}.
\end{proof}

\begin{coro}\label{Cor:doubledual}We suppose that the absolute value $|\ndot|$ is non-Archimedean.
Let $(V,\norm{\ndot})$ be a finite-dimensional seminormed vector space over $k$. Then the double dual seminorm $\norm{\ndot}_{**}$ on $V$ is the largest ultrametric seminorm on $V$ which is bounded from above by $\norm{\ndot}$, and one has $\norm{\ndot}\leqslant\rang(V)\norm{\ndot}_{**}$. If the seminorm $\norm{\ndot}$ is ultrametric, then one has $\norm{\ndot}_{**}=\norm{\ndot}$.
\end{coro}
\begin{proof}
We have seen in Remark \ref{Rem:doubledual} that the double dual seminorm $\norm{\ndot}_{**}$ is ultrametric, and in the formula \eqref{Equ:doubledual} of \S\ref{Subsec:dualnorm} that it is bounded from above by $\norm{\ndot}$. Let $\norm{\ndot}'$ be an ultrametric seminorm on $V$ such that $\norm{\ndot}'\leqslant\norm{\ndot}$. We will show that $\norm{\ndot}'\leqslant\norm{\ndot}_{**}$ and $\norm{\ndot}\leqslant r\norm{\ndot}_{**}$, where $r$ is the rank of $V$ over $k$.

Let $\alpha\in\intervalle]01[$. By Proposition \ref{Pro:existenceoforthogonal}, there exists an $\alpha$-orthogonal basis $\{e_i\}_{i=1}^r$ of $(V,\norm{\ndot})$. 
For any vector $x=\lambda_1e_1+\cdots+\lambda_re_r$ in $V$ one has 
\[\alpha^2\|x\|'\leqslant\alpha^2\max_{i\in\{1,\ldots,r\}}|\lambda_i|\cdot\|e_i\|\leqslant\alpha\max_{i\in\{1,\ldots,r\}}|\lambda_i|\cdot\|e_i\|_{**}\leqslant\|x\|_{**},\]
where the second inequality comes from \eqref{Equ:majoration par doubledual} and the third inequality follows from the fact that $\{e_i\}_{i=1}^r$ is an $\alpha$-orthogonal basis for $\norm{\ndot}_{**}$ (see Proposition \ref{Pro:alphaorthogonale}). Moreover, by the triangle inequality one has
\[\alpha^2\|x\|\leqslant\alpha^2\sum_{i=1}^r|\lambda_i|\cdot\|e_i\|\leqslant r\alpha^2\max_{i\in\{1,\ldots,r\}}|\lambda_i|\cdot\|e_i\|\leqslant r\|x\|_{**}\] 
Since $\alpha\in\intervalle]01[$ is arbitrary, we obtain $\norm{\ndot}'\leqslant\norm{\ndot}_{**}$ and $\norm{\ndot}\leqslant r\norm{\ndot}_{**}$. The first assertion of the proposition is thus proved.

If $\norm{\ndot}$ is ultrametric, it is certainly the largest ultrametric norm bounded from above by $\norm{\ndot}$. Hence one has $\norm{\ndot}=\norm{\ndot}_{**}$.
\end{proof}

\begin{rema} In Corollary \ref{Cor:doubledual},
the constant $\rang(V)$ in the inequality $\norm{\ndot}\leqslant \rang(V)\norm{\ndot}_{**}$ is optimal. We can consider for example  the vector space $V=k^r$ equipped with the $\ell^1$-norm
\[\forall\,(a_1,\ldots,a_r)\in k^r,\quad\|(a_1,\ldots,a_r)\|_{\ell^1}=|a_1|+\cdots+|a_r|.\]
Then its double dual norm is given by
\[\forall\,(a_1,\ldots,a_r)\in k^r,\quad\|(a_1,\ldots,a_r)\|_{\ell^1,**}=\max\{|a_1|,\ldots|a_r|\}.\]
In particular, one has
$\|(1,\ldots,1)\|_{\ell^1}=r\|(1,\ldots,1)\|_{\ell^1,**}$.
\end{rema}

\begin{prop}\label{Pro:doubedualandquotient}
Let $(V,\norm{\ndot})$ be a finite-dimensional seminormed vector space over the field $k$.
\begin{enumerate}[label=\rm(\arabic*)]
\item\label{Item: seminorm and double dual induce the same dual norm} The seminorm $\norm{\ndot}$ and its double dual seminorm  $\norm{\ndot}_{**}$ induce the same dual norm on the vector space $V^*$ of bounded linear forms.
\item\label{Item: quotient seminorm of rank one} If $W$ is a quotient space of rank $1$ of $V$, then the seminorms $\norm{\ndot}$ and $\norm{\ndot}_{**}$ induce the same quotient seminorm on $W$. 
\end{enumerate}
\end{prop}
\begin{proof}
\ref{Item: seminorm and double dual induce the same dual norm} The Archimedean case follows from Proposition \ref{Pro:doubledualarch}. It suffices to treat the case where the absolute value $|\ndot|$ is non-Archimedean. By \eqref{Equ:doubledual} of \S\ref{Subsec:dualnorm}, one has $\norm{\ndot}_{**}\leqslant\norm{\ndot}$. Moreover, by Corollary \ref{Cor: existence of alpha orthogonal}, for any $\alpha\in\intervalle]01[$ there exists an $\alpha$-orthogonal basis $\boldsymbol{e}=\{e_i\}_{i=1}^r$ of $(V,\norm{\ndot})$. By Proposition \ref{Pro:alphaorthogonale}, this basis is also an $\alpha$-orthogonal basis of $(V,\norm{\ndot}_{**})$ and one has $\alpha\norm{e_i}\leqslant\norm{e_i}_{**}$ for any $i\in\{1,\ldots,r\}$. By Proposition \ref{Pro: orthogonal basis: vanishing}
, this leads to $N_{\norm{\ndot}}=N_{\norm{\ndot}_{**}}$ and hence the space of bounded linear forms on $(V,\norm{\ndot}_{**})$ identifies with $(V,\norm{\ndot})^*$, and the relation $\norm{\ndot}_{**}\leqslant\norm{\ndot}$ leads to $\norm{\ndot}_*\leqslant\norm{\ndot}_{**,*}$. Let $\{e_i^\vee\}_{i=1}^r$ be the dual basis of $\{e_i\}_{i=1}^r$ and assume that $\{e_i^\vee\}_{i=1}^r\cap V^*=\{e_i^\vee\}_{i=1}^n$. Then by Proposition \ref{Pro:alphaorthogonale} the familiy $\{e_i^\vee\}_{i=1}^n$ is an $\alpha$-orthogonal basis of $V^*$ both for $\norm{\ndot}_*$ and $\norm{\ndot}_{**,*}$. Moreover, by Lemma \ref{Lem:normofdualbasis}, for any $i\in\{1,\ldots,n\}$ one has
\[\norm{e_i^\vee}_*\geqslant\norm{e_i}^{-1}\geqslant\alpha\norm{e_i}_{**}^{-1}\geqslant \alpha^2\norm{e_i^\vee}_{**,*},\]
where the middle inequality comes from Proposition \ref{Pro:alphaorthogonale}. Since the absolute value $|\ndot|$ is non-Archimedean, the norm $\norm{\ndot}_{**,*}$ is ultrametric. Hence for any $\varphi=\lambda_1e_1^\vee+\cdots+\lambda_ne_n^\vee\in V^*$ one has
\[\norm{\varphi}_{**,*}\leqslant\max_{i\in\{1,\ldots,n\}}|\lambda_i|\cdot\norm{e_i^\vee}_{**,*}\leqslant\alpha^{-2}\max_{i\in\{1,\ldots,n\}}|\lambda_i|\cdot\norm{e_i^\vee}_{*}\leqslant\alpha^{-1}\norm{\varphi}_*,\] 
where the last inequality comes from the fact that $\{e_i^\vee\}_{i=1}^n$ is an $\alpha$-orthogonal basis of $V^*$. Since $\alpha\in\intervalle]01[$ is arbitrary, we obtain $\norm{\ndot}_{**,*}=\norm{\ndot}_*$.

\ref{Item: quotient seminorm of rank one} If the kernel of the quotient map $V\rightarrow W$ does not contain $N_{\norm{\ndot}}$, then the quotient seminorm of $\norm{\ndot}$ on $W$ vanishes
because $\dim_K W = 1$. The quotient seminorm of $\norm{\ndot}_{**}$ on $W$ also vanishes since we have observed in the proof of \ref{Item: seminorm and double dual induce the same dual norm} that $N_{\norm{\ndot}}=N_{\norm{\ndot}_{**}}$. In the following we treat the case where the kernel of the projection map $V\rightarrow W$ contains $N_{\norm{\ndot}}$, or equivalent, the quotient seminorms of $\norm{\ndot}$ and $\norm{\ndot}_{**}$ on $W$ are actually norms. Since $W$ is of rank $1$, any norm on $W$ is uniquely determined by its dual norm on $W^\vee$. Let $\norm{\ndot}_W$ be the quotient norm on $W$ induced by $\norm{\ndot}$. By Proposition \ref{Pro:dualquotient}, the dual norm $\norm{\ndot}_{W,*}$ identifies with the restriction of $\norm{\ndot}_*$ {to} $W^\vee$ (viewed as a vector subspace of $V^*$). By (1), the norm $\norm{\ndot}_*$ identifies with the dual norm of $\norm{\ndot}_{**}$. As a consequence,  $\norm{\ndot}_W$ coincides with the quotient norm of $\norm{\ndot}_{**}$. 
\end{proof}

\begin{prop}\label{Pro:doubledualdet} We assume that the absolute value $|\ndot|$ is non-Archimedean.
Let $(V,\norm{\ndot})$ be a finite-dimensional seminormed vector space over $k$ and let $r$ be the rank of $V$. Then the quotient norm of the $\varepsilon$-tensor product seminorm $\norm{\ndot}_{\varepsilon}$ on $V^{\otimes r}$ by the canonical quotient map $V^{\otimes r}\rightarrow\det(V)$ identifies with the determinant seminorm on $\det(V)$ induced by $\norm{\ndot}$. In particular, $\norm{\ndot}$ and $\norm{\ndot}_{**}$ induce the same determinant seminorm on $\det(V)$.
\end{prop}
\begin{proof} Denote by $\norm{\ndot}_{\det_\varepsilon}$ the quotient seminorm on $\det(V)$ of the $\varepsilon$-tensor product seminorm on $V^{\otimes r}$.
We have seen in  Proposition \ref{Pro:normpiisanorm} that the $\varepsilon$-tensor product seminorm is always bounded from above by the $\pi$-tensor product seminorm. Therefore, one has $\norm{\ndot}_{\det_\varepsilon}\leqslant\norm{\ndot}_{\det}$. Moreover, if $\norm{\ndot}$ is not a norm, then the seminorm $\norm{\ndot}_{\det}$ vanishes. Hence the seminorm $\norm{\ndot}_{\det_\varepsilon}$ also vanishes. To prove the first assertion of the proposition, it remains to verify the inequality $\norm{\ndot}_{\det}\leqslant\norm{\ndot}_{\det_{\varepsilon}}$ in the case where $\norm{\ndot}$ is a norm.

Consider a tensor vector $\varphi$ in $V^{\otimes r}$, which is also viewed as a $k$-multilinear form on $(V^\vee)^r$. By definition, one has
\[\|\varphi\|_{\varepsilon}=\sup_{\begin{subarray}{c}(f_1,\ldots,f_r)\in (V^\vee)^r\\
\forall\,i\in\{1,\ldots,r\},\;f_i\neq 0\end{subarray}}\frac{|{\varphi}(f_1,\ldots,f_r)|}{\|f_1\|_*\cdots\|f_r\|_*}.\]
Let $\alpha\in\intervalle{]}{0}{1}{[}$, $\{x_i\}_{i=1}^r$ be an $\alpha$-orthogonal basis of $V$, and $\{x^{\vee}_i\}_{i=1}^r$ be its dual basis of $V^\vee$. Assume that $\varphi$ is written in the form
\[\varphi=\sum_{I=(i_1,\ldots,i_r)\in\{1,\ldots,r\}^r}a_I(x_{i_1}\otimes\cdots\otimes x_{i_r}),\]
where $a_I\in k$. Then one has
\[\forall\,(i_1,\ldots,i_r)\in\{1,\ldots,r\}^r,\quad\varphi(x_{i_1}^\vee,\ldots,x_{i_r}^\vee)=a_{(i_1,\ldots,i_r)}.\]
In particular,
\[\forall\,(i_1,\ldots,i_r)\in\{1,\ldots,r\}^r,\quad\|\varphi\|_{\varepsilon}\geqslant\frac{|a_{(i_1,\ldots,i_r)}|}{\|x^{\vee}_{i_1}\|_*\cdots\|x^{\vee}_{i_r}\|_*}.\] 
Note that the canonical image $\eta$ of $\varphi$ in $\det(V)$ is
\[\Big(\sum_{\sigma\in\mathfrak S_r}\mathrm{sgn}(\sigma)a_{(\sigma(1),\ldots,\sigma(r))}\Big) x_{1}\wedge\cdots\wedge x_{r},\]
where $\mathfrak S_r$ is the symmetric group of order $r$, namely the group of all bijections from the set $\{1,\ldots,r\}$ to itself, and $\sgn(\ndot):\mathfrak S_r\rightarrow\{\pm 1\}$ denotes the character of signature.
Hence
\begin{equation}\label{Equ:detnormofeps}\begin{split}\|\eta\|_{\det}&=\Big|\sum_{\sigma\in\mathfrak S_r}\sgn(\sigma)a_{(\sigma(1),\ldots,\sigma(r))}\Big|\cdot\|x_1\wedge\cdots\wedge x_r\|_{\det},\\
&\leqslant\Big|\sum_{\sigma\in\mathfrak S_r}\sgn(\sigma)a_{(\sigma(1),\ldots,\sigma(r))}\Big|\cdot\|x_1\|\cdots\|x_r\|,
\end{split}\end{equation}
where the inequality follows from \eqref{Equ:DeltaV}.
 Since the absolute value $|\ndot|$ is non-Archimedean, one has
\[\|\eta\|_{\det}\leqslant\|\varphi\|_{\varepsilon}\cdot\|x_1\|\cdots\|x_r\|\cdot\|x^{\vee}_1\|_*\cdots\|x^{\vee}_r\|_*\leqslant\|\varphi\|_{\varepsilon}\alpha^{-r},\]
where the second inequality comes from Lemma \ref{Lem:normofdualbasis}. Since $\alpha\in\intervalle{]}{0}{1}{[}$ is arbitrary, the first assertion is proved.

We proceed with the proof of the second assertion. By Proposition \ref{Pro:doubedualandquotient}, the seminorms $\norm{\ndot}$ and $\norm{\ndot}_{**}$ induce the same dual norm on $V^*$, and hence induce the same $\varepsilon$-tensor product seminorm on $V^{\otimes r}$. Therefore, by the first assertion of the proposition, we obtain that they induce the same determinant seminorm on $\det(V)$.
\end{proof}

\begin{rema}
In the above proposition, the non-Archimedean assumption on the absolute value is essential. In the Archimedean case, the inequality \eqref{Equ:detnormofeps} only leads to {a} weaker estimate $\norm{\ndot}_{\det}\leqslant r!{\norm{\ndot}_{\det_{\varepsilon}}}$, {where $\norm{\ndot}_{\det_{\varepsilon}}$ is the quotient seminorm on $\det(V)$ induced by the $\varepsilon$-tensor product seminorm}.
\end{rema}


\begin{defi}
Let $(V,\norm{\ndot})$ be a finite-dimensional seminormed vector space over $k$ and $r$ be the dimension of $V$ over $k$. We denote by $\norm{\ndot}_{\det_\varepsilon}$ the quotient seminorm of the $\varepsilon$-tensor power of $\norm{\ndot}$ on $V^{\otimes r}$ by the canonical projection map $V^{\otimes r}\rightarrow\det(V)$, called the \emph{$\varepsilon$-determinant seminorm}\index{epsilon-determinant seminorm@$\varepsilon$-determinant seminorm}\index{seminorm!epsilon-determinant@$\varepsilon$-determinant} of $\norm{\ndot}$.
\end{defi}

\begin{prop}\label{Pro: varepsilon determinant}
We assume that the absolute value $|\ndot|$ is Archimedean. Let $(V,\norm{\ndot})$ be a finite-dimensional seminormed vector space over $k$ and let $r$ be the rank of $V/N_{\norm{\ndot}}$ over $k$. Then the $\varepsilon$-determinant norm of the dual norm $\norm{\ndot}_*$ on $\det(V^*)$  is bounded from below by $(r!)^{-1}\norm{\ndot}^\sim_{\det,*}$, where $\norm{\ndot}^\sim_{\det,*}$ is the dual norm of the determinant of the norm $\norm{\ndot}^\sim$. 
\end{prop}
\begin{proof}
{By Corollary~\ref{Coro:equivalenceofnrom}, one has $V^*=(V/N_{\norm{\ndot}})^\vee$. Moreover one has $\norm{\ndot}_*=\norm{\ndot}^\sim_*$. Hence, by replacing $(V,\norm{\ndot})$ by $(V/N_{\norm{\ndot}},\norm{\ndot}^\sim)$, we may assume without loss of generality that $\norm{\ndot}$ is a norm.} 

Let $\varphi$ be an element in $V^{\vee\otimes r}$. Viewed as a $k$-multilinear form on $V^r$, one has
\[\|\varphi\|_{*,\varepsilon}=\sup_{
(x_1,\ldots,x_r)\in (V\setminus\{0\})^r
}\frac{|\varphi(x_1,\ldots,x_r)|}{\|x_1\|\cdots\|x_r\|}.\]
Let $\{e_i\}_{i=1}^r$ be a basis of $V$ and $\{e_i^\vee\}_{i=1}^r$ be its dual basis of $V^\vee$. Assume that $\varphi$ is written in the form
\[\varphi=\sum_{I=(i_1,\ldots,i_r)\in\{1,\ldots,r\}^r}a_I(e^{\vee}_{i_1}\otimes\cdots\otimes e^{\vee}_{i_r}),\]
where $a_I\in k$. Then for any $(i_1,\ldots,i_r)$ one has
\[\varphi(e_{i_1},\ldots,e_{i_r})=a_{(i_1,\ldots,i_r)}.\]
In particular,
\[\|\varphi\|_{*,\varepsilon}\geqslant\frac{|a_{(i_1,\ldots,i_r)}|}{\|e_{i_1}\|\cdots\|e_{i_r}\|}.\]
Note that the canonical image $\eta$ of $\varphi$ in $\det(V)$ is
\[\Big(\sum_{\sigma\in\mathfrak S_r}\sgn(\sigma)a_{(\sigma(1),\ldots,\sigma(r))}\Big)e_1^\vee\wedge\cdots\wedge e_r^\vee.\]
Therefore,
\[\begin{split}\|\eta\|_{\det,*}&=\bigg|\sum_{\sigma\in\mathfrak S_r}\sgn(\sigma)a_{(\sigma(1),\ldots,\sigma(r))}\bigg|\cdot\|e_1^\vee\wedge\cdots\wedge e_r^\vee\|_{\det,*}\\
&=\bigg|\sum_{\sigma\in\mathfrak S_r}\sgn(\sigma)a_{(\sigma(1),\ldots,\sigma(r))}\bigg|\cdot\frac{1}{\|e_1\wedge\cdots\wedge e_r\|_{\det}}\leqslant r!\|\varphi\|_{*,\varepsilon}\frac{\|e_1\|\cdots\|e_r\|}{\|e_1\wedge\cdots\wedge e_r\|_{\det}}.
\end{split}\]
By \eqref{Equ:DeltaV}, we obtain $\|\eta\|_{\det,*}\leqslant r!\|\varphi\|_{*,\varepsilon}$. The proposition is thus proved.
\end{proof}

{
\begin{prop}\label{Pro:alphatenso}
Let $\{(V_j,\norm{\ndot}_j)\}_{j=1}^d$ be a finite family of finite-dimensional seminormed $k$-vector spaces and let $\alpha$ be a real number in $\intervalle]01]$. For any $j\in\{1,\ldots,d\}$, let $\{e_i^{(j)}\}_{i=1}^{n_j}$ be an $\alpha$-orthogonal basis of $V_j$. Then 
\[e_{i_1}^{(1)}\otimes\cdots\otimes e_{i_d}^{(d)},\quad (i_1,\ldots,i_d)\in\prod_{j=1}^d\{1,\ldots,n_j\}\] form an $\alpha^d$-orthogonal basis of  $V_1\otimes_k\cdots\otimes_kV_d$ for the $\varepsilon$-tensor product norm $\norm{\ndot}_{\varepsilon}$ of $\{\norm{\ndot}_j\}_{j=1}^d$. Moreover, if $\norm{\ndot}$ is an \emph{ultrametric} norm on $V_1\otimes_k\cdots\otimes_kV_d$ such that 
\[\norm{x_1\otimes\cdots\otimes x_d}\leqslant\norm{x_1}_{1,**}\cdots\norm{x_d}_{d,**}\] for any $(x_1,\ldots,x_d)\in V_1\times\cdots\times V_d$, then $\norm{\ndot}\leqslant\norm{\ndot}_{\varepsilon}$.
\end{prop}
\begin{proof}
Let 
\[T=\sum_{(i_1,\ldots,i_d)\in\prod_{j=1}^d\{1,\ldots,n_j\}} a_{i_1,\ldots,i_d}e_{i_1}^{(1)}\otimes\cdots e_{i_d}^{(d)}\] be a tensor in $V_1\otimes_k\cdots\otimes_k V_d$, where $a_{i_1,\ldots,i_d}\in k$. We consider it as a 
$k$-multilinear form on $V_1^{\vee}\times\cdots\times V_d^{\vee}$. For any $j\in\{1,\ldots,d\}$, let $\{\varphi^{(j)}_i\}_{i=1}^{n_j}$ be the dual basis of $\{e_i^{(j)}\}_{i=1}^{n_j}$ and assume that 
\[\{\varphi^{(j)}_i\}_{i=1}^{n_j}\cap V^*=\{\varphi^{(j)}_i\}_{i=1}^{n_j'}.\] By Proposition \ref{Pro:alphaorthogonale}, 
$\{\varphi_i^{(j)}\}_{i=1}^{n_j'}$ is an $\alpha$-orthogonal basis of $V^{(j),*}$. For any $(i_1,\ldots,i_d)\in\prod_{j=1}^d\{1,\ldots,n_j'\}$ we have  
$T(\varphi_{i_1}^{(1)},\ldots,\varphi_{i_d}^{(d)})=a_{i_1,\ldots,i_d}$,
which leads to
\[\begin{split}\|T\|_{\varepsilon}\geqslant \frac{|a_{i_1,\ldots,i_d}|}{\|\varphi_{i_1}^{(1)}\|_{1,*}\cdots\|\varphi_{i_d}^{(d)}\|_{d,*}}&\geqslant\alpha^d|a_{i_1,\ldots,i_d}|\cdot\|e_{i_1}^{(1)}\|_{1,**}\cdots\|e_{i_d}^{(d)}\|_{d,**}\\&=\alpha^d|a_{i_1,\ldots,i_d}|\cdot\|e_{i_1}^{(1)}\otimes\cdots\otimes e_{i_d}^{(d)}\|_{\varepsilon},
\end{split}
\]
 where the second inequality follows from \eqref{Equ:normeduale} and the equality comes from Remark~\ref{Rem:produittenrk1}.
This completes the proof of the proposition because 
\[\norm{e_{i_1}^{(1)}\otimes\cdots\otimes e_{i_d}^{(d)}}_\varepsilon=\norm{e_{i_1}^{(1)}}_{1,**}\cdots\norm{e_{i_d}^{(d)}}_{d,**}\] vanishes once $(i_1,\ldots,i_d)\not\in\prod_{j=1}^d\{1,\ldots,n_j'\}$ (see Lemma \ref{Lem:normofdualbasis}). Moreover, if $\norm{\ndot}$ is an ultrametric norm on $V_1\otimes_k\cdots\otimes_k V_d$ such that \[\norm{x_1\otimes\cdots\otimes x_d}\leqslant\norm{x_1}_{1,**}\cdots\norm{x_d}_{d,**}\] for any $(x_1,\ldots,x_d)\in V_1\times\cdots\times V_d$, then 
\[\begin{split}\norm{T}&\leqslant\max_{(i_1,\ldots,i_d)\in\prod_{j=1}^d\{1,\ldots,n_j\}}|a_{i_1,\ldots,i_d}|\cdot\norm{e_{i_1}^{(1)}\otimes\cdots\otimes e_{i_d}^{(d)}}\\
&\qquad \leqslant\max_{(i_1,\ldots,i_d)\in\prod_{j=1}^d\{1,\ldots,n_j\}}|a_{i_1,\ldots,i_d}|\cdot\norm{e_{i_1}^{(1)}}_{1,**}\cdot\norm{e_{i_d}}_{d,**}\leqslant\alpha^{-d}\norm{T}_{\varepsilon}\end{split}\]
By Proposition \ref{Pro:existenceoforthogonal}, for any $\alpha\in\intervalle]01[$, there exist $\alpha$-orthogonal bases of $V_1,\ldots,V_d$ respectively. Hence $\norm{\ndot}\leqslant\norm{\ndot}_{\varepsilon}$.
\end{proof}

\begin{coro}\label{Cor: dual tensor product}
Assume that the absolute value $|\ndot|$ is non-Archimedean. Let $\{(V_j,\norm{\ndot}_j)\}_{j=1}^d$ be a finite family of finite-dimensional seminormed vector spaces over $k$, and $\norm{\ndot}_{\varepsilon}$ be the $\varepsilon$-tensor product of the seminorms $\norm{\ndot}_1,\ldots,\norm{\ndot}_d$. Then the dual norm $\norm{\ndot}_{\varepsilon,*}$ coincides with the $\varepsilon$-tensor product of the dual norms $\norm{\ndot}_{1,*},\ldots,\norm{\ndot}_{d,*}$.
\end{coro}
\begin{proof}
Let $\alpha$ be an element of $\intervalle]01[$. For $j\in\{1,\ldots,d\}$, let $\{e_i^{(j)}\}_{i=1}^{n_j}$ be an $\alpha$-orthogonal basis of $V_j$ over $k$ (see Proposition~\ref{Pro:existenceoforthogonal}) and  $\{\varphi_i^{(j)}\}_{i=1}^{n_j}$ be the dual basis of $\{e_i^{(j)}\}_{i=1}^{n_j}$, and assume that $\{\varphi_i^{(j)}\}_{i=1}^{n_j}\cap V_j^*=\{\varphi_i^{(j)}\}_{i=1}^{n_j'}$. Note that for any $(i_1,\ldots,i_d)\in\prod_{j=1}^d\{1,\ldots,n_j\}$, 
\[\norm{e_{i_1}^{(1)}\otimes\cdots\otimes e_{i_d}^{(d)}}_{\varepsilon}=\norm{e_{i_1}}_{1,**}\cdot\norm{e_{i_d}}_{d,**}\neq 0\] 
if and only if $(i_1,\ldots,i_d)\prod_{j=1}^d\in\{1,\ldots,n_i'\}$. Therefore, \[\{\varphi_{i_1}^{(1)}\otimes\cdots\otimes \varphi_{i_d}^{(d)}\}_{(i_1,\ldots,i_d)\in\prod_{j=1}^d\{1,\ldots,n_j'\}}\] forms a basis of the vector space of bounded linear forms on $(V_1\otimes_k\cdots\otimes_kV_d,\norm{\ndot}_{\varepsilon})$, which shows that $(V_1\otimes_k\cdots\otimes_k V_d)^*\cong V_1^{*}\otimes_k\cdots\otimes_kV_d^{*}$.

Let $\norm{\ndot}'$ be the $\varepsilon$-tensor product of $\norm{\ndot}_{1,*},\ldots,\norm{\ndot}_{d,*}$. By definition, for any $T$ in $V_1^{*}\otimes_k\cdots\otimes_kV_d^{*}$, one has
\[\begin{split}\norm{T}'&=\sup_{\begin{subarray}{c}
(s_1,\ldots,s_d)\in V_1\times\cdots\times V_d\\
\min\{\norm{s_1}_1,\ldots,\norm{s_d}_d\}>0
\end{subarray}}\frac{|T(s_1,\ldots,s_d)|}{\norm{s_1}_{1,**}\cdots\norm{s_d}_{d,**}}\\
&=\sup_{\begin{subarray}{c}
(s_1,\ldots,s_d)\in V_1\times\cdots\times V_d\\
\min\{\norm{s_1}_1,\ldots,\norm{s_d}_d\}>0
\end{subarray}}\frac{|T(s_1,\ldots,s_d)|}{\norm{s_1\otimes\cdots\otimes s_d}_{\varepsilon}}\leqslant\norm{T}_{\varepsilon,*},
\end{split}\]
where the second equality comes from \eqref{equ:epsilontensorpodurc}. Conversely, if $T$ is of the form $\psi_1\otimes\cdots\otimes\psi_d$, where $(\psi_1,\ldots,\psi_d)\in V_1^{*}\times\cdots\times V_d^{*}$, then
\[\norm{T}_{\varepsilon,*}=\sup_{\begin{subarray}{c}f\in V_1\otimes_k\cdots\otimes_k V_d\\\norm{f}_\varepsilon\neq 0\end{subarray}}\frac{|f(\psi_1,\ldots,\psi_d)|}{\norm{f}_{\varepsilon}}\leqslant\norm{\psi_1}_{1,*}\cdots\norm{\psi_d}_{d,*}.\] By Proposition \ref{Pro:alphatenso}, we obtain $\norm{\ndot}_{\varepsilon,*}\leqslant \norm{\ndot}'$.
\end{proof}
}

\subsection{Orthogonality and lattice norms}

In this subsection, we assume that the absolute value $|\ndot|$ is non-Archimedean and we denote by $\mathfrak o_k$ the valuation ring of $(k,|\ndot|)$. Let $V$ be a finite-dimensional vector space over $k$ and let $r$ be the rank of $V$.

\begin{prop}\label{Pro:orthogonallattice}
Let $\mathcal V$ be a lattice of $V$ which is an $\mathfrak o_k$-module of finite type (and hence a free $\mathfrak o_k$-module of rank $r$). Then any basis of $\mathcal V$ over $\mathfrak o_k$ is an orthonormal basis of $(V,\norm{\ndot}_{\mathcal V})$. 
\end{prop}
\begin{proof}
Let $\{e_i\}_{i=1}^r$ be a basis of $\mathcal V$ over $\mathfrak o_k$. Note that an element $a_1e_1+\cdots+a_re_r$ of $V$ (with $(a_1,\ldots,a_r)\in k^r$) belongs to $\mathcal V$ if and only if all $a_i$ are in $\mathfrak o_k$. Let $x=\lambda_1e_1+\cdots+\lambda_re_r$ be an element of $V$. If $a$ is an element in $k^{\times}$ such that $a^{-1}x$ belongs to $\mathcal V$, then one has $a^{-1}\lambda_i\in \mathfrak o_k$ and hence $|\lambda_i|\leqslant|a|$ for any $i\in\{1,\ldots,r\}$. Therefore $\max\{|\lambda_1|,\ldots,|\lambda_r|\}\leqslant\|x\|_{\mathcal V}$. Conversely, if $j\in\{1,\ldots,r\}$ is such that \[|\lambda_j|=\max\{|\lambda_1|,\ldots,|\lambda_r|\}>0,\] then one has $\lambda_i\lambda_j^{-1}\in \mathfrak o_k$ for any $i\in\{1,\ldots,r\}$. Hence $\lambda_j^{-1}x\in\mathcal V$ and \[\|x\|_{\mathcal V}\leqslant|\lambda_j|=\max\{|\lambda_1|,\ldots,|\lambda_r|\}.\]
\end{proof}

\begin{prop}\label{Pro:approximation}
Assume that the absolute value $|\ndot|$ is non-trivial. Let $\lambda\in\intervalle]01[$ be a real number such that 
\begin{equation}\label{Equ:approximationcondition}\lambda<\sup\{|a|\,:\,a\in k^{\times},\,|a|<1\}.\end{equation}
Then, for any ultrametric norm $\norm{\ndot}$ on $V$ there exists a lattice of finite type $\mathcal V$ of $V$, such that $\norm{\ndot}\leqslant\norm{\ndot}_{\mathcal V}\leqslant\lambda^{-1}\norm{\ndot}$. 
\end{prop}
\begin{proof}
Let $\alpha\in \intervalle]01[$ such that $\lambda/\alpha<\sup\{|a|\,:\,a\in k^{\times},\,|a|<1\}$. Let $\{e_i\}_{i=1}^r$ be an $\alpha$-orthogonal basis of $V$ (the existence of which has been proved in Proposition \ref{Pro:existenceoforthogonal}). By Proposition \ref{Pro:dilatation}, by dilating the vectors $e_i$, $i\in\{1,\ldots,r\}$, we may assume that $\lambda/\alpha\leqslant\|e_i\|<1$ for any $i$. Let $\mathcal V$ be the free $\mathfrak o_k$-module generated by $\{e_i\}_{i=1}^r$. It is a lattice of $V$. Moreover, by Proposition \ref{Pro:orthogonallattice}, $\{e_i\}_{i=1}^r$ is an orthonormal basis of $(V,\norm{\ndot}_{\mathcal V})$. In particular, for any vector $x=a_1e_1+\cdots+a_re_r$ in $V$, one has
\[\|x\|_{\mathcal V}=\max_{i\in\{1,\ldots,r\}}|a_i|\geqslant\max_{i\in\{1,\ldots,r\}}|a_i|\cdot\|e_i\|\geqslant\|x\|.\]
Moreover, by the $\alpha$-orthogonality of $\{e_i\}_{i=1}^r$ one has
\[\|x\|\geqslant\alpha\max_{i\in\{1,\ldots,r\}}|a_i|\cdot\|e_i\|\geqslant\lambda\|x\|_{\mathcal V}.\]
The proposition is thus proved.
\end{proof}

\subsection{Orthogonality and Hadamard property} We have seen in Proposition \ref{Pro:existenceoforthogonal} that an Hadamard basis of a finite-dimensional normed vector space is necessarily orthogonal. The converse of this assertion is also true when the absolute value $|\ndot|$ is non-Archimedean.

\begin{prop}\label{Pro:orthogonalesthadamard}
We assume that the absolute value $|\ndot|$ is non-Archimedean. Let $(V,\norm{\ndot})$ be a finite-dimensional seminormed vector space over $k$, and let $r$ be the rank of $V$ over $k$. If $\alpha$ is an element in $\intervalle]01]$ and if $\{x_i\}_{i=1}^r$ is an $\alpha$-orthogonal basis of $V$, then one has
\begin{equation}\label{Equ: comparison Hadamard orthogonal}\|x_1\wedge\cdots\wedge x_r\|\geqslant\alpha^r\|x_1\|\cdots\|x_r\|.\end{equation}
In particular, if $\norm{\ndot}$ is a norm, any orthogonal basis of $V$ is an Hadamard basis.
\end{prop}
\begin{proof}
If $N_{\norm{\ndot}}$ is non-zero, then the interserction of $\{x_i\}_{i=1}^r$ with $N_{\norm{\ndot}}$ is not empty (see Proposition \ref{Pro: orthogonal basis: vanishing}) and hence the inequality \eqref{Equ: comparison Hadamard orthogonal} holds. In the following, we assume that $\norm{\ndot}$ is a norm. Note that the case where $V=\{0\}$ is trivial. Hence we may assume that $r>0$. Let $\{x_i\}_{i=1}^r$ be an $\alpha$-orthogonal basis of $V$. Let $\{y_i\}_{i=1}^r$ be an arbitrary basis of $V$ and $A=(a_{ij})_{i\in\{1,\ldots,r\},\,j\in\{1,\ldots,r\}}\in k^{r\times r}$ be the transition matrix from $\{x_i\}_{i=1}^r$ to $\{y_i\}_{i=1}^r$, namely $y_i=\sum_{j=1}^ra_{ij}x_j$ for any $i\in\{1,\ldots,r\}$. By the $\alpha$-orthogonality of the basis $\{x_i\}_{i=1}^r$ one has
\[|a_{ij}|\leqslant\alpha^{-1}\frac{\|y_i\|}{\|x_j\|}.\]
Thus, by the assumption that the absolute value $|\ndot|$ is non-Archimedean, one has
\[\begin{split}\|y_1\wedge\cdots\wedge y_r\|=|\det(A)|\cdot\|x_1\wedge\cdots\wedge x_r\|
\leqslant\alpha^{-r}\frac{\|y_1\|\cdots\|y_r\|}{\|x_1\|\cdots\|x_r\|}\cdot\|x_1\wedge\cdots\wedge x_r\|
\end{split}\]
and hence
\[\frac{\|x_1\wedge\cdots\wedge x_r\|}{\|x_1\|\cdots\|x_r\|}\geqslant\alpha^r\frac{\|y_1\wedge\cdots\wedge y_r\|}{\|y_1\|\cdots\|y_r\|}.\]
Since the basis $\{y_i\}_{i=1}^r$ is arbitrary, 
by Proposition \ref{Pro:hadamard} (see also Remark \ref{Rem:Hadamard basis})
we obtain that 
\[\frac{\|x_1\wedge\cdots\wedge x_r\|}{\|x_1\|\cdots\|x_r\|}\geqslant\alpha^r.\]
The proposition is thus proved.
\end{proof}

\begin{rema}
The Archimedean analogue of Proposition \ref{Pro:orthogonalesthadamard} is not true in general. One can consider for example the case where $V=\mathbb R^2$ equipped with the norm $\norm{\ndot}$ such that $\|(a,b)\|=\max\{|a|,|b|\}$, where $|\ndot|$ is the usual absolute value on $\mathbb R$. Let $e_1=(1,0)$ and $e_2=(0,1)$. The basis $\{e_1,e_2\}$ is orthonormal. However,
\[\|e_1\wedge e_2\|=\Big\|\frac 1{\sqrt{2}}(e_1+e_2)\wedge\frac 1{\sqrt{2}}(e_1-e_2)\Big\|=\frac {1}{2}.\]
Therefore $\{e_1,e_2\}$ is not an Hadamard basis.
\end{rema}

The following proposition shows that the Archimedean analogue of Proposition \ref{Pro:orthogonalesthadamard} is true provided that the norm is induced by an inner product.

\begin{prop}\label{Pro:ArchimedeanHadamard}
Let $(V,\norm{\ndot})$ be a finite-dimensional normed vector space over $k$. Assume that the absolute value $|\ndot|$ is Archimedean and that the norm $\norm{\ndot}$ is induced by an inner product. Then any orthogonal basis of $V$ is an Hadamard basis.
\end{prop}
\begin{proof}
The field $k$ is locally compact, therefore $V$ admits an Hadamard basis $\boldsymbol{e}=\{e_i\}_{i=1}^r$, which is necessarily an orthogonal basis (see Proposition \ref{Pro:existenceoforthogonal}). Without loss of generality, we may assume that $\boldsymbol{e}$ is an orthonormal basis. Let $\boldsymbol{e}'=\{e_i'\}_{i=1}^r$ be another orthonormal basis. There exists a unitary matrix $A$ such that $\boldsymbol{e}'=A\boldsymbol{e}$. One has $|\det(A)|=1$ and hence
\[\|e_1'\wedge\cdots\wedge e_r'\|_{\det}=\|e_1\wedge\cdots\wedge e_r\|_{\det}=1=\|e_1'\|\cdots\|e_r'\|.\]
Therefore the basis $\{e_i'\}_{i=1}^r$ is also an Hadamard basis. Thus we have proved that any orthonormal basis is an Hadamard basis. If $\{x_i\}_{i=1}^r$ is an orthogonal basis, and if {$e_i=\|x_i\|^{-1}x_i$} for any $i\in\{1,\ldots,r\}$, then $\{e_i\}_{i=1}^r$ is an orthonormal basis of $V$, which is an Hadamard basis. We then deduce that $\{x_i\}_{i=1}^r$ is also an Hadamard basis.
\end{proof}

{
\begin{prop}\label{Pro: alpha orth is orth}
We assume that the absolute value $|\ndot|$ is trivial. Let $(V,\norm{\ndot})$ be 
an $r$-dimensional $(r\in\mathbb N_{>0})$, ultrametrically normed vector space over $k$, which corresponds to an increasing sequence
\[0=V_0\subsetneq V_1\subsetneq \ldots\subsetneq V_n=V\]
of vector subspaces of $V$ and a decreasing sequence $\mu_1>\ldots>\mu_n$ of real numbers as described in Remark \ref{Rem: R-filtration as flag plus slopes}.
\begin{enumerate}[label=\rm(\arabic*)]
\item\label{Item: criterion of orthgonal basis} A basis $\{x_j\}_{j=1}^r$ of $V$ is orthogonal if and only if $\operatorname{\mathrm{card}}(\{x_j\}_{j=1}^r\cap V_i)=\dim_k(V_i)$ for any $i\in\{1,\ldots,n\}$.
\item\label{Item: alpha orth is oth}  Let $\alpha$ be an element of $\mathopen{]}0,1\mathclose{]}$ such that 
\[\forall\,i\in\{1,\ldots,n\},\quad \alpha>\mathrm{e}^{-(\mu_i-\mu_{i+1})/r},\]
where 
$\mu_{n+1}=-\infty$ by convention. Then any $\alpha$-orthogonal basis of $(V,\norm{\ndot})$ is orthogonal.
\end{enumerate}
\end{prop}
\begin{proof}
\ref{Item: criterion of orthgonal basis} Note that the restriction of $\norm{\ndot}$ on each $V_i\setminus V_{i-1}$ is constant and takes $\mathrm{e}^{-\mu_i}$  as its value for $i\in\{1,\ldots,n\}$. Let  
\[\lambda_1\leqslant\ldots\leqslant\lambda_r\]
be the increasing sequence of positive real numbers such that $\mathrm{e}^{-\mu_i}$ appears exactly $\dim_k(V_i)-\dim_{k}(V_{i-1})$ times. Let $\{x_j\}_{j=1}^r$ be a basis of $V$ such that $\norm{x_1}\leqslant\ldots\leqslant\norm{x_r}$. For each $i\in\{1,\ldots,n\}$, the cardinal of $\{x_j\}_{j=1}^r\cap V_i$  does not exceed $\dim_k(V_i)$, so that $\norm{x_j}\geqslant\lambda_j$ for any $j\in\{1,\ldots,r\}$ and hence 
\[\prod_{j=1}^r\norm{x_j}\geqslant\prod_{j=1}^r\lambda_j.\]
Moreover, if the equality $\operatorname{\mathrm{card}}(\{x_j\}_{j=1}^r\cap V_i)=\dim_k(V_i)$ holds for any $i\in\{1,\ldots,n\}$, then the basis $\{x_j\}_{j=1}^r$ is an Hadamard basis, and hence is an orthogonal basis  (by Proposition \ref{Pro:existenceoforthogonal} \ref{Item: Hadamard basis orthogonal}). If there exists an index $i\in\{1,\ldots,n\}$ such that $\operatorname{\mathrm{card}}(\{x_j\}_{j=1}^r\cap V_i)<\dim_k(V_i)$, then there exists an element 
\[x=\lambda_1 x_1+\cdots+\lambda_r x_r\] of $V_i$
and a $j\in\{1,\ldots,r\}$ such that $\lambda_j\neq 0$ and $x_j\in V\setminus V_i$. As $\norm{x}\leqslant\mu_i<\norm{x_j}$, the basis  $\{x_j\}_{j=1}^r$ is not orthogonal.

\ref{Item: alpha orth is oth}  Let $\{e_j\}_{j=1}^r$ be an $\alpha$-orthogonal basis of $(V,\norm{\ndot})$. Without loss of generality, we assume that $\norm{e_1}\leqslant\ldots\leqslant\norm{e_r}$.  One has
\begin{equation}\label{eqn:Pro: alpha orth is orth:01}
\{\norm{e_1},\ldots,\norm{e_r}\}\subseteq\{\lambda_1,\ldots,\lambda_r\}= \{ \mathrm{e}^{-\mu_1}, \ldots, \mathrm{e}^{-\mu_n} \}.
\end{equation}
Moreover, since 
\[\operatorname{card}(\{{e}_j\}_{j=1}^r\cap V_i)\leqslant \dim_k(V_i)\]
for any $i\in\{1,\ldots,n\}$, one has $\norm{{e}_j}\geqslant\lambda_j$ for any $j\in\{1,\ldots,r\}$. 
Therefore, if the strict inequality $\norm{e_j}>\lambda_j$ holds for some  $j\in\{1,\ldots,r\}$, then, by \eqref{eqn:Pro: alpha orth is orth:01}, one can find $m$ such that $m > j$, $\norm{e_j} \geqslant \lambda_m$ and $\lambda_m \lambda_j^{-1} = \mathrm{e}^{\mu_i - \mu_{i+1}}$ for some $i$, so that, by our assumption, \[
\norm{e_j} \geqslant \lambda_j (\lambda_m \lambda_j^{-1}) = \lambda_j \mathrm{e}^{\mu_i-\mu_{i+1}} > \lambda_j \alpha^{-r}. \]
Therefore,
\[\prod_{\ell=1}^r\norm{e_{\ell}}> \alpha^{-r}\prod_{\ell=1}^r\lambda_{\ell}=\alpha^{-r}\norm{e_1\wedge\cdots\wedge e_r}_{\det}.\]
This contradicts Proposition \ref{Pro:orthogonalesthadamard}. Therefore one has $\norm{e_j}=\lambda_j$ for any $j\in\{1,\ldots,r\}$ and hence, for each $i\in\{1,\ldots,n\}$, the cardinal of $\{e_j\}_{j=1}^r\cap V_i$  is equal to  $\dim_k(V_i)$, which implies that $\{e_j\}_{j=1}^r$ is an orthogonal basis. 
\end{proof}
}

\subsection{Ultrametric Gram-Schimdt process}
\label{subsec:Ultrametric Gram-Schimdt process}
In this subsection, we consider a refinement of Proposition \ref{Pro:existenceoforthogonal}.
First we recall the spherically completeness of a metric space.

\begin{defi}\label{def:spherically:complete}
We say that a metric space $(X,d)$
is \emph{spherically complete}\index{spherically complete} if, for any decreasing
sequence 
\[
B_1 \supseteq B_2 \supseteq \cdots \supseteq B_{n} \supseteq B_{n+1} \supseteq \cdots
\]
of non-empty closed balls  in $X$, one has
$\bigcap_{n=1}^{\infty} B_n \not= \varnothing$.
A normed vector space $(V,\norm{\ndot})$ over $k$ is said to be \emph{spherically complete}\index{spherically complete} 
if $(V,\norm{\ndot})$ is spherically complete as a metric space. If $(k,|\ndot|)$, viewed as a normed vector space over $k$, is spherically complete, we say that the valued field $(k,|\ndot|)$ is spherically complete.
\end{defi}

\begin{rema}
If $(k, |\ndot|)$ is a discrete valuation field, then
$(k, |\ndot|)$ is spherically complete by \cite[Proposition~20.2]{Schikhof06}.
In particular, any locally compact non-Archimedean valued field is  spherically complete.
\end{rema}

\begin{lemm}\label{lem:orthogonal:spherically:complete}
Let $(V,\norm{\ndot})$ be an ultrametrically normed vector space of finite rank over $k$.
Then we have the following:
\begin{enumerate}[label=\rm(\arabic*)]
\item\label{Item: spherically complete implies distance attained}
Let $W$ be a vector subspace of $V$ over $k$.
If $W$ equipped with the restriction $\norm{\ndot}_W$ of $\norm{\ndot}$ {to} $W$ is
spherically complete, then, for $x \in V$, there is $w \in W$ such that
$\mathrm{dist}(x, W) = \| x - w \|$.

\item\label{Item: field spherically complete implies vector space}
If $(V, \norm{\ndot})$ has an orthogonal basis $\{e_i\}_{i=1}^{r}$ and
$(k, |\ndot|)$ is spherically complete, then
$(V, \norm{\ndot})$ is also spherically complete.
\end{enumerate}
\end{lemm}

\begin{proof}
For $a \in V$ and $\delta \in \mathbb R_{\geqslant 0}$, we set
\[B(a; \delta) := \{ x \in V \,:\, \| x - a \| \leqslant \delta \}.\] As $\norm{\ndot}$ is
ultrametric, we can easily see that
\begin{equation}\label{eqn:lem:orthogonal:spherically:complete:01}
B(a; \delta) = B(a';\delta)
\end{equation}
for all $\delta \in \mathbb R_{\geqslant 0}$ and $a, a' \in V$ with $\| a - a' \| \leqslant \delta$.

\medskip
\ref{Item: spherically complete implies distance attained} We can choose a decreasing sequence $\{ \delta_n \}_{n=1}^{\infty}$ of positive numbers and
a sequence $\{ w_n \}_{n=1}^{\infty}$ in $W$ such that
$\| x - w_n \| \leqslant \delta_n$ and $\lim_{n\to\infty} \delta_n = \mathrm{dist}(x, W)$.
As $B(x;\delta_n) \cap W = B(w_n;\delta_n) \cap W$ by \eqref{eqn:lem:orthogonal:spherically:complete:01},
$\{ B(x;\delta_n) \cap W \}_{n=1}^{\infty}$ yields a decreasing sequence of non-empty closed balls in $W$.
Thus, by our assumption, there is $w \in \bigcap_{n=1}^{\infty} B(x;\delta_n) \cap W$, so that
$\| x - w\| \leqslant \delta_n$ for all $n$, that is, $\| x - w \| \leqslant \mathrm{dist}(x, W)$, as required.

\medskip
\ref{Item: field spherically complete implies vector space} Note that $\|a_1 e_1 + \cdots + a_{r} e_{r} \| = \max \{ |a_1|\cdot\|e_1\|, \ldots, |a_{r}|\cdot\|e_{r}\| \}$, so that \[
B(a; \delta) = B_k(a_1;\delta/\|e_1\|) e_1 + \cdots + B_k(a_{r};\delta/\|e_{r}\|) e_{r}
\]
for $a = a_1 e_1 + \cdots + a_{r} e_{r} \in V$ and $\delta \in \mathbb R_{\geqslant 0}$,
where \[B_k(\lambda; \delta') = \{ t \in k \,:\, |t - \lambda| \leqslant \delta' \}\] for $\lambda \in k$ and $\delta' \in \mathbb R_{\geqslant 0}$.
Therefore the assertion follows.
\end{proof}

In the case of an ultrametrically normed finite-dimensional vector space, Proposition \ref{Pro:existenceoforthogonal} has the following refined form. This could be considered as an ultrametric analogue of Gram-Schmidt orthogonalisation process.

\begin{prop}\label{Pro:existenceepsorth}
Let $(V,\norm{\ndot})$ be an ultrametrically seminormed $k$-vector space of rank $r\geqslant 1$. Let
\begin{equation}\label{Equ:flagcomplete}0=V_0\subsetneq V_1\subsetneq V_2\subsetneq\ldots\subsetneq V_r=V\end{equation}
be a complete flag of subspaces of $V$.
Fix a real number $\alpha$ such that
\[
\alpha \in \begin{cases}
\intervalle]01], & \text{if $(k, |\ndot|)$ is spherically complete}, \\
\intervalle]01[, & \text{otherwise}.
\end{cases}
\] 
Then there exists an $\alpha$-orthogonal basis $\boldsymbol{e}$ of $V$ such that, for any $i\in\{1,\ldots,r\}$, $\card(V_i\cap\boldsymbol{e})=i$.
\end{prop}
\begin{proof} If a basis $\boldsymbol{e}$ of $V$ is such that, for any $i\in\{1,\ldots,r\}$, $\card(V_i\cap\boldsymbol{e})=i$, we say that the basis $\boldsymbol{e}$ is \emph{compatible}\index{compatible with the flag} with the flag \eqref{Equ:flagcomplete}.

We begin with the proof of the particular case where $\norm{\ndot}$ is a norm by induction on $r$, the dimension of $V$ over $k$. The case where $r=1$ is trivial. Assume that the proposition holds for all
vector spaces of dimension $<r$, where $r\geqslant 2$. Applying the induction hypothesis to $V_{r-1}$ and the flag $0=V_0
\subsetneq V_1\subsetneq\ldots\subsetneq V_{r-1}$ we get a basis $\{e_1,\ldots,e_{r-1}\}$ of $V_{r-1}$ compatible with the flag such 
that, for any $(\lambda_1,\ldots,\lambda_{r-1})\in k^{r-1}$
\begin{equation}\label{Equ:hyprecurence}\|\lambda_1e_1+\cdots+\lambda_{r-1}e_{r-1}\|\geqslant \alpha^{1/2}\max_{i\in\{1,\ldots,r-1\}}|\lambda_i|\cdot\|e_i\|.\end{equation}
Let $x$ be an element of $V\setminus V_{r-1}$. The distance between $x$ 
and $V_{r-1}$ is strictly positive since $V_{r-1}$ is closed in $V$ (see Proposition \ref{Pro:topologicalnormedspace}). Hence there exists $y\in V_{r-1}$ such that \begin{equation}\label{Equ:distance}
\|x-y\|\leqslant \alpha^{-1/2}\mathrm{dist}(x,V_{r-1}).
\end{equation}
In the case where $\alpha=1$ and $(k,|\ndot|)$ is spherically complete, the existence of $y$ follows
from Lemma~\ref{lem:orthogonal:spherically:complete}.
We choose $e_r=x-y$. The basis $\{e_1,\ldots,e_r\}$ is compatible with the flag $0=V_0\subsetneq V_1\subsetneq \ldots
\subsetneq V_r=V$. 

Let $(\lambda_1,\ldots,\lambda_r)$ be an element of $k^r$. We wish to find a lower bound for the 
norm of $z=\lambda_1e_1+\cdots+\lambda_re_r$. By \eqref{Equ:distance} we have that 
\[\|z\|\geqslant |\lambda_r|\cdot\mathrm{dist}(x,V_{r-1})\geqslant \alpha^{1/2}|\lambda_r|\cdot\|e_r\|.\] 
This provides our lower bound when $\|\lambda_r e_r\|\geqslant \|\lambda_1e_1+\cdots+\lambda_{r-1}e_{r-1}\|$. If 
$\|\lambda_r e_r\|< \|\lambda_1e_1+\cdots+\lambda_{r-1}e_{r-1}\|$ then we have 
\[\|z\|=\|\lambda_1e_1+\cdots+\lambda_{r-1}e_{r-1}\|\] because the norm is ultrametric (see Proposition \ref{Pro:valeur}). By the induction hypothesis 
\eqref{Equ:hyprecurence} we have that $\|z\|\geqslant \alpha|\lambda_i|\cdot\|e_i\|$ for any 
$i\in\{1,\ldots,r-1\}$. This completes the proof of the proposition in the case where $\norm{\ndot}$ is a norm.

We now consider the general seminorm case. Let $W$ be the quotient vector space $V/N_{\norm{\ndot}}$. For each $i\in\{0,\ldots,r\}$, let $W_i$ be $(V_i+N_{\norm{\ndot}})/N_{\norm{\ndot}}$. Applying the particular case of the proposition to $(W,\norm{\ndot}^\sim)$, we obtain the existence of an $\alpha$-orthogonal basis $\widetilde{\boldsymbol{e}}$ of $W$ such that $\card(\widetilde{\boldsymbol{e}}\cap W_i)=\rang_k(W_i)$. 
We set \[I = \{ i \in \{ 1, \ldots, r \} \,:\, W_{i-1} \subsetneq W_{i} \}\quad\text{ and }\quad
J = \{ j \in \{ 1, \ldots, r \} \,:\, W_{j-1} = W_{j} \}.\]
If $i \in I$, then there is a unique element $u_i \in \widetilde{\boldsymbol{e}}\cap (W_{i}\setminus W_{i-1})$, so that we can choose $e_i \in V_i$ such that the class of $e_i$ in $V/N_{\|\ndot\|}$ is $u_i$.
If $j \in J$, then we can pick up $e_j \in (N_{\norm{\ndot}}\cap V_j)\setminus V_{j-1}$. Indeed,
as $V_{j} \setminus V_{j-1} \not= \varnothing$ and $V_j \subseteq V_{j-1} + N_{\|\ndot\|}$,
we can find $x \in V_{j} \setminus V_{j-1}$, $y \in V_{j-1}$ and $e_j \in N_{\|\ndot\|}$ with
$x = y + e_j$, and hence $e_j \in (N_{\norm{\ndot}}\cap V_j)\setminus V_{j-1}$.
By construction, $\boldsymbol{e}:=\{e_i\}_{i=1}^r$ satisfies $\card(V_i \cap \boldsymbol{e}) = i$ for $i\in\{0, \ldots, r\}$.
In particular, $\boldsymbol{e}$ forms a basis of $V$.

Let us see that $\boldsymbol{e}$ is $\alpha$-orthogonal.
For any $(\lambda_1,\ldots,\lambda_r)\in k^r$, if we let $x=\lambda_1e_1+\cdots+\lambda_re_r$ and 
$u=\sum_{i \in I}\lambda_iu_i$,
then one has
\[\norm{x}=\norm{u}^\sim\geqslant\alpha\max_{i \in I}
|\lambda_i|\cdot\norm{u_i}^\sim=\alpha\max_{i\in\{1,\ldots,r\}}|\lambda_i|\cdot\norm{e_i},\]
where the inequality comes from the $\alpha$-orthogonality of $\widetilde{\boldsymbol{e}}$, as required. 
\end{proof}

\begin{coro}\label{Cor: existence of orthogonal basis spherically complete}
Let $(V,\norm{\ndot})$ be an ultrametrically seminormed vector space of finite rank over $k$.
If $(k, |\ndot|)$ is spherically complete, then $(V, \norm{\ndot})$ has an orthogonal basis.
In particular, $(V, \norm{\ndot})$ is spherically complete.
\end{coro}

\begin{rema}\label{Rem:Gram-Schmidt} 
Assume that the absolute value $|\ndot|$ is Archimedean. Let $V$ be a finite-dimensional vector space over $k$ and $\norm{\ndot}$ be a seminorm on $V$ which is induced by an inner product. Given a complete flag
\[0=V_0\subsetneq V_1\subsetneq V_2\subsetneq\ldots\subsetneq V_r=V\]
of $V$, the Gram-Schmidt process permits to construct an orthogonal basis $\boldsymbol{e}$ of $V$ such that $\card(\boldsymbol{e}\cap V_i)=i$ for any $i\in\{1,\ldots,r\}$, along the same line as in the proof of Proposition \ref{Pro:existenceepsorth}. The main point is that the field $k$ is locally compact, and hence the distance in \eqref{Equ:distance} is actually attained by some point in $V_{r-1}$. For a general seminorm, even though an orthogonal basis always exists, it is not always possible to find an orthogonal basis which is compatible with a given flag.  
\end{rema}

Proposition \ref{Pro:existenceepsorth} and the usual Gram-Schmidt process lead to the following projection result.

\begin{coro}\label{Cor:projection}
Let $V$ be a finite-dimensional vector space over $k$ equipped with a seminorm $\norm{\ndot}$ which is either ultrametric or induced by a  semidefinite inner product. Let $V_0$ be a vector subspace of $V$. For any $\alpha\in\intervalle]01[$ there exists a $k$-linear projection $\pi:V\rightarrow V_0$ (namely $\pi$ is a $k$-linear map and its restriction {to} $V_0$ is the identity map) such that $\|\pi\|\leqslant\alpha^{-1}$. If {$(k,|\ndot|)$ is non-Archimedean and spherically complete, or if}   
$\norm{\ndot}$ is induced by a  semidefinite inner product, we can choose the $k$-linear projection $\pi$ such that $\|\pi\|\leqslant 1$.  
\end{coro}
\begin{proof}
We first consider the ultrametric case. By Proposition \ref{Pro:existenceepsorth}, there exists an $\alpha$-orthogonal basis $\boldsymbol{e}=\{e_i\}_{i=1}^n$ of $V$ such that $\card(\boldsymbol{e}\cap V_0)=\rang_k(V_0)$. Without loss of generality, we may assume that $\boldsymbol{e}\cap V_0=\{e_i\}_{i=1}^m$, where $m=\rang_k(V_0)$. Let $\pi:V\rightarrow V_0$ be the $k$-linear map sending $\lambda_1e_1+\cdots+\lambda_ne_n\in V$ to $\lambda_1e_1+\cdots+\lambda_me_m\in V_0$. Since the basis $\boldsymbol{e}$ is $\alpha$-orthogonal, one has 
\[\|\lambda_1e_1+\cdots+\lambda_ne_n\|\geqslant\alpha\max_{i\in\{1,\ldots,n\}}\|\lambda_ie_i\|\geqslant\alpha\|\lambda_1e_1+\cdots+\lambda_me_m\|,\]
which implies that $\|\pi\|\leqslant\alpha^{-1}$.

{If $(k,|\ndot|)$ is non-Archimedean and spherically complete, or if} $\norm{\ndot}$ is induced by an inner product, we use the existence of an orthogonal basis $\boldsymbol{e}$ such that $\card(\boldsymbol{e}\cap V_0)=\rang_k(V_0)$. {By the same agrument as above, we obtain the existence of a linear projection  $\pi:V\rightarrow V_0$ such that $\norm{\pi}\leqslant 1$.}
\end{proof}

\begin{coro}\label{Cor:comparaisondenormes}
Let $V$ be a finite-dimensional vector space over $k$ and $\norm{\ndot}_1$ and $\norm{\ndot}_2$ be two seminorms on $V$. We assume that $\norm{\ndot}_1\leqslant\norm{\ndot}_2$ and that the seminorm $\norm{\ndot}_2$ is either ultrametric or induced by a  semidefinite inner product. If there exists a vector $x\in V$ such that $\|x\|_1<\|x\|_2$, then one has $\norm{\ndot}_{1,\det}'<\norm{\ndot}^\sim_{2,\det}$ on $\det(V/N_{\norm{\ndot}_2})\setminus\{0\}$, where $\norm{\ndot}_1'$ denotes the quotient seminorm of $\norm{\ndot}_1$ on $V/N_{\norm{\ndot}_2}$.
\end{coro} 
\begin{proof}
{The condition $\norm{\ndot}_1\leqslant\norm{\ndot}_2$ implies that $N_{\norm{\ndot}_2}\subseteq N_{\norm{\ndot}_1}$. In particular, for any $x\in V$, if we denote by $[x]$ the class of $x$ in $V/N_{\norm{\ndot}_2}$, then one has
\[\norm{[x]}_1'=\norm{x}_1,\quad \norm{[x]}^\sim_2=\norm{x}_2.\]
Therefore, by replacing $V$ by $V/N_{\norm{\ndot}_2}$, $\norm{\ndot}_1$ by $\norm{\ndot}_1'$, and $\norm{\ndot}_2$ by $\norm{\ndot}_2^\sim$, we may assume without loss of generality that $\norm{\ndot}_2$ is a norm. Moreover, the case where $\norm{\ndot}_1$ is not a norm is trivial since the seminorm $\norm{\ndot}_{1,\det}$ vanishes. Hence it suffices to treat the case where both seminorms $\norm{\ndot}_1$ and $\norm{\ndot}_2$ are norms.}

Let $\lambda=\|x\|_2/\|x\|_1>1$.
We first consider the ultrametric case. By Proposition \ref{Pro:existenceepsorth}, for any $\alpha\in\intervalle]01[$ there exists an $\alpha$-orthogonal basis $\{e_i\}_{i=1}^r$ of $(V,\norm{\ndot}_2)$ such that $e_1=x$. Hence {See Proposition \ref{Pro:orthogonalesthadamard} for the first inequality}
\begin{multline*}\|e_1\wedge\cdots\wedge e_r\|_{2,\det}\geqslant\alpha^r\|e_1\|_2\cdots\|e_r\|_2 \\
\geqslant\lambda\alpha^r\|e_1\|_{1}\cdots\|e_r\|_1\geqslant\lambda\alpha^r\|e_1\wedge\cdots\wedge e_r\|_{1,\det}.
\end{multline*}
Since $\alpha\in\intervalle]01[$ is arbitrary, one has 
${\norm{\ndot}_{2,\det}}/{\norm{\ndot}_{1,\det}}\geqslant\lambda>1.$

The Archimedean case is very similar. There exists an orthogonal basis $\{e_i\}_{i=1}^r$ of $(V,\norm{\ndot})$ such that $e_1=x$. We then proceed as above in replacing $\alpha$ by $1$.
\end{proof}

\begin{prop}
\label{Pro:quotientdualnonarch}
Let $V$ be a finite-dimensional vector space over $k$, equipped with a seminorm $\norm{\ndot}_V$, which is ultrametric or induced by a {semidefinite} inner product. Let $(W,\norm{\ndot}_W)$ be a seminormed vector space over $k$. For any $k$-vector subspace $V_0$ of $V$, {the $k$-linear map $\mathscr L(V,W)\rightarrow\mathscr L(V_0,W)$, sending $f\in\mathscr L(V,W)$ to its restriction {to} $V_0$, is surjective, and} the operator seminorm on $\mathscr L(V_0,W)$ coincides with the quotient seminorm induced by the operator seminorm on $\mathscr L(V,W)$. In particular, the dual norm on $V_0^*$ identifies with the quotient of the dual norm on $V^*$ by the canonical quotient map $V^*\rightarrow V_0^*$.
\end{prop}
\begin{proof} For any $f\in\mathscr L(V,W)$, the operator seminorm of $f|_{V_0}$ does not exceed that of $f$. In the following, we show that, for any linear map $g\in\mathscr L(V_0,W)$ and any $\alpha\in\intervalle]01[$, there exists a $k$-linear map $f:V\rightarrow W$ extending $g$ such that $\alpha\|f\|\leqslant\|g\|$. By Corollary \ref{Cor:projection}, there exists a $k$-linear projection $\pi:V\rightarrow V_0$ such that $\|\pi\|\leqslant\alpha^{-1}$. Let $f=g\circ\pi$. Then  $\|f\|\leqslant\alpha^{-1}\|g\|$. The proposition is thus proved.
\end{proof}

\begin{prop}\label{Pro:quotientr1eps}
Let $V$ and $W$ be finite-dimensional seminormed vector spaces over $k$, $V_0$ be a $k$-vector subspace of $V$, and $Q$ be the quotient vector space $V/V_0$. We assume that, either the absolute value $|\ndot|$ is non-Archimedean, or the seminorm on $V$ is induced by a  semidefinite inner product. Then the canonique isomorphism $(V\otimes_kW)/(V_0\otimes_kW)\rightarrow Q\otimes_kW$ is an isometry, where we consider the $\varepsilon$-tensor product seminorms on $V\otimes_kW$ and $Q\otimes_kW$, and the quotient seminorm on $(V\otimes_kW)/(V_0\otimes_kW)$.
\end{prop}
\begin{proof}
We have seen in Remark \ref{Rem:quotientepsilontensor} that, for any $f\in V\otimes_kW$ viewed as a $k$-bilinear form on $V^*\times W^*$, its restriction {to} $Q^*\times W^*$ has an $\varepsilon$-tensor product norm not greater than that of $f$. In the following, we consider an element $g\in Q\otimes_kW$, viewed as a $k$-bilinear form on $Q^*\otimes_kW^*$. We will show that, for any $\alpha\in\intervalle]01[$, there exists a $k$-bilinear form $f$ on $V^*\times W^*$ extending $g$ such that $\|f\|_{\varepsilon}\leqslant\alpha^{-1}\|g\|_{\varepsilon}$. By Proposition \ref{Pro:dualquotient}, the dual norm on $Q^*$ identifies with the restriction of the dual norm on $V^*$. By Corollary \ref{Cor:projection}, there exists a $k$-linear projection $\pi:V^*\rightarrow Q^*$ such that $\|\pi\|\leqslant\alpha^{-1}$ (in the non-Archimedean case, we use the fact that any dual norm is ultrametric). We let $f$ be the $k$-bilinear form on $V^*\times W^*$ such that $f(\xi,\eta)=g(\pi(\xi),\eta)$. Then for $(\xi,\eta)\in (V^*\setminus\{0\})\times(W^*\setminus\{0\})$ one has
\[\frac{|f(\xi,\eta)|}{\|\xi\|_*\|\eta\|_*}=\frac{|g(\pi(\xi),\eta)|}{\|\xi\|_*\|\eta\|_*}\leqslant\alpha^{-1}\frac{|g(\pi(\xi),\eta)|}{\|\pi(\xi)\|_*\|\eta\|_*}\leqslant\alpha^{-1}\|g\|_{\varepsilon}.\]
The proposition is thus proved.   
\end{proof}

\begin{coro}\label{Cor: determinant of exterior power}
We assume that the absolute value $|\ndot|$ is non-Archimedean. Let $(V,\norm{\ndot})$ be a finite-dimensional seminormed vector space over $k$ and $r$ be the rank of $V$ over $k$. Let $i$ be a positive  
integer. Then the canonical $k$-linear isomorphism $\det(\Lambda^{i}V)\rightarrow\det(V)^{\otimes\binom{r-1}{i-1}}$ is an isometry, where we consider  the $i^{\text{th}}$ $\varepsilon$-exterior power seminorm on $\Lambda^iV$.  
\end{coro}
\begin{proof}
Consider the following commutative diagram
\[\xymatrix{V^{\otimes i\binom{r}{i}}\ar@{->>}[r]^-{p_1}\ar@{->>}[d]_-{p_2}&\det(V)^{\otimes\binom{r-1}{i-1}}
\\(\Lambda^iV)^{\otimes \binom{r}{i}}\ar@{->>}[r]_{p_3}&\det(\Lambda^iV)\ar[u]_-\simeq}\]
By {Proposition \ref{Pro:quotientr1eps}}, if we equip $V^{\otimes i\binom{r}{i}}$ with the $\varepsilon$-tensor product seminorm, then its quotient seminorm on $(\Lambda^iV)^{\otimes\binom{r}{i}}$ identifies with the $\varepsilon$-tensor product of the $\varepsilon$-exterior power seminorm. Moreover, by Proposition \ref{Pro:doubledualdet}, the quotient seminorm on $\det(\Lambda^iV)$ (induced by $p_3$) of the tensor product of the $\varepsilon$-exterior power seminorm identifies with the determinant seminorm of the latter. Still by the same proposition, the quotient seminorm on $\det(V)^{\otimes\binom{r-1}{i-1}}$ induced by $p_1$ identifies with the tensor power of the determinant seminorm. Therefore the natural isomorphism $\det(\Lambda^{i}V)\rightarrow\det(V)^{\otimes\binom{r-1}{i-1}}$ is actually an isometry
by using (1) in Proposition~\ref{prop:quotient:norm:linear:map}. 
\end{proof}

\begin{coro}
We assume that the absolute value $|\ndot|$ is non-Archimedean. Let $(V,\norm{\ndot})$ be a finite-dimensional normed vector space over $k$ and $i$ be a positive  
integer. Then the $i^{\text{th}}$ $\varepsilon$-exterior power norm on $\Lambda^iV$ is the double dual norm of the $i^{\text{th}}$ $\pi$-exterior power norm.
\end{coro}
\begin{proof}
By definition, the $i^{\text{th}}$ $\varepsilon$-exterior power norm $\norm{\ndot}_{\Lambda_\varepsilon^i}$ on $\Lambda^iV$ is ultrametric and is bounded from above by the $i^{\text{th}}$ $\pi$-exterior power norm. By Corollary \ref{Cor:doubledual}, $\norm{\ndot}_{\Lambda_\pi^i,**}$ is the largest ultrametric norm bounded from {above by} the $i^{\text{th}}$ $\pi$-exterior power norm $\norm{\ndot}_{\Lambda^i_\pi}$. In particular, one has $\norm{\ndot}_{\Lambda_\pi^i,**}\geqslant\norm{\ndot}_{\Lambda^i_\varepsilon}$ and hence
\[\norm{\ndot}_{\Lambda^i_\varepsilon ,\det}\leqslant\norm{\ndot}_{\Lambda_\pi^i,**,\det}\leqslant\norm{\ndot}_{\Lambda^i_\pi,\det}.\]  
By Corollary \ref{Cor: determinant of exterior power}, one has $\norm{\ndot}_{\Lambda^i_\varepsilon ,\det}=\norm{\ndot}_{\Lambda^i_\pi,\det}$  and hence the above inequalities are actually equalities. By Corollary \ref{Cor:comparaisondenormes}, we obtain $\norm{\ndot}_{\Lambda_\varepsilon^i}=\norm{\ndot}_{\Lambda_\pi^i,**}$. 
\end{proof}

\begin{prop}\label{Pro: tensor product and deteminant epsilon}
Assume that  $|\ndot|$ is non-Archimedean. Let $V$ and $W$ be finite-dimensional seminormed vector spaces over $k$ and $n$ and $m$ be respectively the ranks of $V$ and $W$ over $k$. We equip $V\otimes_kW$ with the $\varepsilon$-tensor product seminorm $\norm{\ndot}_{\varepsilon}$.  Then the natural $k$-linear isomorphism $\det(V\otimes_kW)\cong\det(V)^{\otimes m}\otimes_k\det(W)^{\otimes n}$ is an isometry, where we consider the determinant seminorm of $\norm{\ndot}_{\varepsilon}$ on $\det(V\otimes_kW)$ and the tensor product of determinant seminorms on $\det(V)^{\otimes m}\otimes_k\det(W)^{\otimes n}$.
\end{prop}
\begin{proof}
Let $\norm{\ndot}'$ be the seminorm on $\det(V)^{\otimes m}\otimes\det(W)^{\otimes n}$ induced by tensor product of determinant seminorms. {By Propositions \ref{Pro:quotientr1eps} and \ref{Pro:doubledualdet}, the seminorm $\norm{\ndot}'$ identifies with the quotient  of the $\varepsilon$-tensor power seminorm on $(V\otimes_kW)^{\otimes nm}$ of $\norm{\ndot}_{\varepsilon}$. Therefore, by Proposition \ref{Pro:doubledualdet} the seminorm $\norm{\ndot}'$ identifies with $\norm{\ndot}_{\varepsilon,\det}$.
}
\end{proof}

\begin{prop}\label{Pro:exactsequence}
Let $V$ be a finite-dimensional vector space over $k$ and $\norm{\ndot}$ be a seminorm on $V$. We assume that the seminorm $\norm{\ndot}$ is either ultrametric or induced by a  semidefinite inner product. For any vector subspace $W$ of $V$, the canonical isomorphism
\[\det(W)\otimes_k\det(V/W)\longrightarrow\det(V)\] is an isometry,  where we consider the determinant seminorm of the induced seminorm on $\det(W)$ and that of the quotient seminorm on $\det(V/W)$, and the  tensor product seminorm on the tensor product space $\det(W)\otimes_k\det(V/W)$ (see Remark \ref{Rem:produittenrk1}).
\end{prop}
\begin{proof}
{By Proposition \ref{Pro: induce quotient zero ball}, if the seminorm $\norm{\ndot}$ is not a norm, then either its restriction {to} $W$ is not a norm, or its quotient seminorm on $V/W$ is not a norm. In both cases, the seminorms on $\det(W)\otimes_k\det(V/W)$ and on $\det(V)$ vanish. Therefore we may assume without loss of generality that $\norm{\ndot}$ is a norm.}

Let $f:\det(W)\otimes_k\det(V/W)\rightarrow\det(V)$ be the canonical isomorphism. We have seen in Corollary \ref{Cor:exactsequencenorm} that the operator norm of $f$ is $\leqslant 1$. Since $f$ is an isomorphism between vector spaces of dimension $1$ over $k$, to prove that $f$ is an isometry, it suffices to verify that $\|f\|\geqslant 1$.

We first treat the case where the norm $\norm{\ndot}$ is ultrametric. By Proposition \ref{Pro:existenceepsorth}, for any $\alpha\in\intervalle]01[$, there exists an $\alpha$-orthogonal basis $\boldsymbol{e}=\{e_i\}_{i=1}^r$ of $V$ such that $\card(\boldsymbol{e}\cap W)=\rang(W)$. Without loss of generality, we assume that $\{e_1,\ldots,e_n\}$ forms a basis of $W$, and $e_{n+1},\ldots,e_r$ are vectors in $V\setminus W$. For any $i\in\{n+1,\ldots,r\}$, let $\overline{e}_i$ be the image of $e_i$ in $V/W$. By Proposition \ref{Pro:orthogonalesthadamard}, one has
\[\begin{split}\|e_1\wedge\cdots\wedge e_r\|_{\det}&\geqslant\alpha^r\cdot\|e_1\|\cdots\|e_r\|\geqslant\alpha^r\|e_1\|\cdots\|e_n\|\cdot\|\overline e_{n+1}\|\cdots\|\overline{e}_r\|\\
&\geqslant\alpha^r\|e_1\wedge\cdots\wedge e_n\|_{\det}\cdot\|\overline{e}_{n+1}\wedge\cdots\wedge\overline{e}_r\|_{\det},
\end{split}\]
where the last equality comes from Corollary \ref{Cor:exactsequencenorm}. Therefore the operator norm of $f$ is bounded from below by $\alpha^r$. Since $\alpha\in\intervalle]01[$ is arbitrary, one has $\|f\|\geqslant 1$.

For the Archimedean case where the norm $\norm{\ndot}$ is induced by an inner product, by the classic Gram-Schmidt process we can construct an orthonormal basis $\boldsymbol{e}$ of $V$ such that $\card(\boldsymbol{e}\cap W)=\rang(W)$. By Proposition \ref{Pro:orthogonalesthadamard}, $\boldsymbol{e}$ is an Hadamard basis. We then proceed as above in replacing $\alpha$ by $1$ to conclude.
\end{proof}

In Proposition \ref{Pro:exactsequence}, the assumption on the seminorm is crucial. In order to study the behaviour of the determinant seminorms of an exact sequence of general seminormed vector spaces, we introduce the following invariant.

\begin{defi}\label{Def:Delta:V:norm}
Let $(V,\norm{\ndot})$ be a finite-dimensional seminormed vector space over $k$. Let $\mathcal H(V,\norm{\ndot})$ be the set of all normes $\norm{\ndot}_h$ on $V/N_{\norm{\ndot}}$ which is either ultrametric or induced by an inner product, and such that $\norm{\ndot}_h\geqslant\norm{\ndot}^{\sim}$.  We define $\Delta(V,\norm{\ndot})$ to be the number {(if $\norm{\ndot}$ vanishe, by convention $\Delta(V,\norm{\ndot})$ is defined to be $1$)}
\[\inf\bigg\{\frac{\norm{\ndot}_{h,\det}}{\norm{\ndot}_{\det}^\sim}\,:\,\norm{\ndot}_h\in\mathcal H(V,\norm{\ndot})\bigg\}\in \intervalle[1{+\infty}[,\]
where $\norm{\ndot}_{h,\det}$ and $\norm{\ndot}_{\det}^\sim$ are respectively determinant norms on $\det(V/N_{\norm{\ndot}})$ induced by the norms $\norm{\ndot}_{h}$ and $\norm{\ndot}^\sim$. By definition, one has $\Delta(V,\norm{\ndot})=\Delta(V/N_{\norm{\ndot}},\norm{\ndot}^\sim)$. Moreover, if $\norm{\ndot}$ is ultrametric or induced by a  semidefinite inner product, then $\Delta(V,\norm{\ndot})=1$. \end{defi}

{\begin{prop}Assume that $|\ndot|$ is non-Archimedean. Let $(V,\norm{\ndot})$ be a finite-dimensional seminormed vector space over $k$. One has 
\begin{equation}\label{Equ: upper bound of log delta}\ln\Delta(V,\norm{\ndot})\leqslant\rang(V/N_{\norm{\ndot}})\sup_{x\in V\setminus N_{\norm{\ndot}}}\big(\ln\norm{x}-\ln\norm{x}_{**}\big).\end{equation}
In particular, $\ln\Delta(V,\norm{\ndot})\leqslant\rang(V/N_{\norm{\ndot}})\ln(\rang(V/N_{\norm{\ndot}}))$.
\end{prop}
\begin{proof}
By replacing $(V,\norm{\ndot})$ by $(V/N_{\norm{\ndot}},\norm{\ndot}^\sim)$, we may assume without loss of generality that $\norm{\ndot}$ is a norm. Let 
\[\lambda=\sup_{x\in V\setminus\{0\}}\big(\ln\norm{x}-\ln\norm{x}_{**}\big).\]
By definition one has $\norm{\ndot}\leqslant \mathrm{e}^\lambda\norm{\ndot}_{**}$. Note that the norm $\mathrm{e}^\lambda\norm{\ndot}_{**}$ is ultrametric.
Therefore
\[\Delta(V,\norm{\ndot})\leqslant\frac{\mathrm{e}^{r\lambda}\norm{\ndot}_{**,\det}}{\norm{\ndot}_{\det}}=\mathrm{e}^{r\lambda},\]
where $r$ is the rank of $V$ over $k$, and the equality comes from Proposition \ref{Pro:doubledualdet}. The inequality \eqref{Equ: upper bound of log delta} is thus proved. The last inequality results from \eqref{Equ: upper bound of log delta} and Corollary \ref{Cor:doubledual}.
\end{proof}

\begin{prop}\label{Pro:minormationdenorme}
Let $(V,\norm{\ndot})$ be a finite-dimensional normed vector space over $k$. For any vector subspace $W$ of $V$, the norm of the canonical isomorphism
\[f:\det(W)\otimes\det(V/W)\longrightarrow\det(V)\]
is bounded from below by
\[\frac{\Delta(W,\norm{\ndot}_W)\Delta({V/W},\norm{\ndot}_{V/W})}{\Delta(V,\norm{\ndot})}\geqslant\Delta(V,\norm{\ndot})^{-1},\]
where $\norm{\ndot}_{W}$ is the restriction of the norm $\norm{\ndot}$ {to} the vector subspace $W$ and $\norm{\ndot}_{V/W}$ is the quotient norm of $\norm{\ndot}$ on the quotient space $V/W$.
\end{prop}
\begin{proof}
Let $\norm{\ndot}_{h}$ be a norm in $\mathcal H(V,\norm{\ndot})$. Let $\norm{\ndot}_{h,W}$ and $\norm{\ndot}_{h,V/W}$ be respectively the restriction of $\norm{\ndot}_{h}$ {to} $W$ and the quotient norm of $\norm{\ndot}_{h}$ on $V/W$. By Proposition \ref{Pro:exactsequence}, the canonical isomorphism
\[\det(W,\norm{\ndot}_{h,W})\otimes\det(V/W,\norm{\ndot}_{h,V/W})\longrightarrow\det(V,\norm{\ndot}_{V,h})\]
is an isometry. Hence
\[\frac{\norm{\ndot}_{h,\det}}{\norm{\ndot}_{\det}}=\frac{\norm{\ndot}_{h,W,\det}\norm{\ndot}_{h,V/W,\det}}{\|f\|\cdot\norm{\ndot}_{W,\det}\norm{\ndot}_{V/W,\det}}\geqslant\frac{1}{\|f\|}\cdot{\Delta(W,\norm{\ndot}_W)}{\Delta({V/W},\norm{\ndot}_{V/W})}.\]
Since $\norm{\ndot}_h\in\mathcal H(V,\norm{\ndot})$ is arbitrary, we obtain the lower bound announced in the proposition.
\end{proof}

\begin{coro}\label{Cor: comparaison de Delta}
Let $(V,\norm{\ndot}_V)$ be a finite-dimensional seminormed vector space over $k$, $W$ be a vector subspace of $V$, $\norm{\ndot}_W$ be the restriction of $\norm{\ndot}_V$ {to} $W$, and $\norm{\ndot}_{V/W}$ be the quotient of $\norm{\ndot}_V$ on $V/W$. One has
\begin{equation}\label{Equ:produitdeDelta}\Delta(W,\norm{\ndot}_W)\Delta({V/W},\norm{\ndot}_{V/W})\leqslant\Delta(V,\norm{\ndot}).\end{equation}
In particular, $\Delta(W,\norm{\ndot}_W)\leqslant\Delta(V,\norm{\ndot})$ and $\Delta({V/W},\norm{\ndot}_{V/W})\leqslant\Delta(V,\norm{\ndot})$.
\end{coro}
\begin{proof}By Proposition \ref{Pro: induce quotient zero ball}, we can assume without loss of generality that $\norm{\ndot}_V$ is a norm. By Corollary \ref{Cor:exactsequencenorm}, if we denote by $f:\det(W)\otimes\det(V/W)\rightarrow\det(V)$ the canonical isomorphism, then $\|f\|\leqslant 1$. The inequality \eqref{Equ:produitdeDelta} thus follows from Proposition \ref{Pro:minormationdenorme}. Finally, by definition one has $\Delta(W,\norm{\ndot}_W)\geqslant 1$ and $\Delta({V/W},\norm{\ndot}_{V/W})\geqslant 1$, thus we deduce from \eqref{Equ:produitdeDelta} the last two inequalities stated in the corollary.  
\end{proof}

\begin{rema}\label{rema:Delta:upper:bound}
Let $(V,\norm{\ndot})$ be a finite-dimensional normed vector space over $k$. We assume that the rank $r$ of $V/N_{\norm{\ndot}}$ is positive.
In the case where the absolute value $|\ndot|$ is non-Archimedean, Corollary \ref{Cor:doubledual} provides the upper bound $\Delta(V,\norm{\ndot})\leqslant r^r$. This result is also true in the Archimedean case (which follows from the existence of an orthogonal basis, see the beginning of \S\ref{Sec:John-Lowner} for details). However, as we will see in the next subsection (cf. Theorem \ref{Thm:John}), in the Archimedean case one has a better upper bound $\Delta(V,\norm{\ndot})\leqslant r^{r/2}$.
\end{rema}

\subsection{Dual determinant norm}\label{Subsec:dualdet}

Let $(V,\norm{\ndot})$ be a finite-dimensional seminormed vector space over $k$. We denote by $\norm{\ndot}_{\det}^\sim$ the determinant norm on $\det(V/N_{\norm{\ndot}})$ induced by $\norm{\ndot}^\sim$, and denote by $\norm{\ndot}^\sim_{\det,*}$ the dual norm of $\norm{\ndot}^\sim_{\det}$. Let $\norm{\ndot}_{*,\det}$ be the determinant norm on $\det(V^*)\cong\det(V/N_{\norm{\ndot}})^*$ of the dual norm $\norm{\ndot}_*$ on $V^*$. The purpose of this subsection is to compare these two norms. We denote by $\delta(V,\norm{\ndot})$ the ratio
\[\delta(V,\norm{\ndot}):=\frac{\norm{\ndot}_{\det,*}^\sim}{\norm{\ndot}_{*,\det}}.\]
In the case where there is no ambiguity on the seminorm $\norm{\ndot}$ on $V$, we also use the abbreviate notation $\delta(V)$ to denote $\delta(V,\norm{\ndot})$.
By definition, if $\eta$ is a non-zero element in $\det(V/N_{\norm{\ndot}})$ and if $\eta^\vee$ is its dual element in $\det(V^*)$, then one has
\begin{equation}\label{Equ:definitiondedelta}\delta(V,\norm{\ndot})^{-1}=\|\eta\|_{\det}^\sim\cdot\|\eta^\vee\|_{*,\det}. \end{equation}
In particular, one has (see Proposition \ref{Pro:doubledualdet})
\begin{equation}
\label{Equ: delta equals dual}\delta(V,\norm{\ndot})=\delta(V/N_{\norm{\ndot}},\norm{\ndot}^\sim)=\delta(V^*,\norm{\ndot}_*).\end{equation}

\begin{prop}\label{Pro:minorationdelta}
Let $(V,\norm{\ndot})$ be a finite-dimensional seminormed vector space over $k$. One has $\delta(V,\norm{\ndot})\geqslant 1$.
\end{prop}
\begin{proof}By \eqref{Equ: delta equals dual} we may assume without loss of generality that $\norm{\ndot}$ is a norm.

Let $\{e_i\}_{i=1}^r$ be a basis of $V$, and $\{e_i^\vee\}_{i=1}^r$ be its dual basis. One has
\[\|e_1^\vee\wedge\cdots\wedge e_r^\vee\|_{\det,*}=\|e_1\wedge\cdots\wedge e_r\|^{-1}_{\det}.\]
Therefore
\[\begin{split}\delta(V,\norm{\ndot})^{-1}&=\|e_1\wedge\cdots\wedge e_r\|_{\det}\cdot\|e_1^\vee\wedge\cdots\wedge e_r^\vee\|_{*,\det}\\
&\leqslant\|e_1\|\cdots\|e_r\|\cdot\|e_1^\vee\|_*\cdots\|e_r^\vee\|_*,
\end{split}\]
where the inequality comes from Proposition \ref{Pro:hadamard}.
If the basis $\{e_i\}_{i=1}^r$ is $\alpha$-orthogonal, where $\alpha\in\intervalle]01[$, by Lemma \ref{Lem:normofdualbasis} one has $
\delta(V,\norm{\ndot})^{-1}\leqslant\alpha^{-r}$.
Since for any $\alpha\in\intervalle]01[$ there exists an $\alpha$-orthogonal basis (see Proposition \ref{Pro:existenceoforthogonal}), one has $\delta(V,\norm{\ndot})\geqslant 1$.
\end{proof}

\begin{prop}\label{Pro:delta=1} Let $(V,\norm{\ndot})$ be a finite-dimensional seminormed vector space over $k$. 
Assume that the absolute value $|\ndot|$ is non-Archimedean, or the seminorm $\norm{\ndot}$ is induced by a  semidefinite inner product. Then one has $\delta(V,\norm{\ndot})=1$.
\end{prop}
\begin{proof}
By \eqref{Equ: delta equals dual} we may assume without loss of generality that $\norm{\ndot}$ is a norm.

We first treat the case where the absolute value $|\ndot|$ is non-Archimedean. Let $\alpha\in \intervalle]01[$ and $\{e_i\}_{i=1}^r$ be an $\alpha$-orthogonal basis of $(V,\norm{\ndot})$ (see Proposition \ref{Pro:existenceoforthogonal}
 for the existence of an $\alpha$-orthogonal basis). Then the dual basis $\{e_i^\vee\}_{i=1}^r$ is $\alpha$-orthogonal with respect to the dual norm $\norm{\ndot}_*$ (see Proposition \ref{Pro:alphaorthogonale}). In particular, by Proposition \ref{Pro:orthogonalesthadamard} one has
\[\frac{\|e_1\wedge\cdots\wedge e_r\|_{\det}}{\|e_1\|\cdots\|e_r\|}\geqslant\alpha^r,\quad
\frac{\|e_1^\vee\wedge\cdots\wedge e_r^\vee\|_{*,\det}}{\|e_1^\vee\|_*\cdots\|e_r^\vee\|_*}\geqslant\alpha^r.\]
Therefore
\[\delta(V,\norm{\ndot})=\frac{\|e_1\wedge\cdots\wedge e_r\|_{\det}^{-1}}{\|e_1^\vee\wedge\cdots\wedge e_r^\vee\|_{*,\det}}\leqslant\alpha^{-2r}\frac{1}{\|e_1\|\cdots\|e_r\|\cdot\|e_1^\vee\|_*\cdots\|e_r^\vee\|_*}\leqslant\alpha^{-2r},\]
where the last inequality comes from Lemma \ref{Lem:normofdualbasis}.
Since $\alpha$ is arbitrary, one has $\delta(V,\norm{\ndot})\leqslant 1$.

The proof of the Archimedean case is quite similar, where we use the existence of an orthogonal basis, which is also an Hadamard basis (see Proposition \ref{Pro:ArchimedeanHadamard}). We omit the details.
\end{proof}

The following Lemma is the Archimedean counterpart of Proposition \ref{Pro:quotientdualnonarch} (see also the comparison in Remark \ref{Rem:comparaisondualitye}).

\begin{lemm}\label{Lem:quotientequalsrestriction}
Assume that the absolute value $|\ndot|$ is Archimedean.
Let $(V,\norm{\ndot}_V)$ be a finite-dimensional seminormed vector space over $k$, $W$ be a vector subspace of $V$, and $\norm{\ndot}_W$ be the restriction of the seminorm $\norm{\ndot}_V$ {to} $W$. Then the map $F:V^*\rightarrow W^*$, which sends $\varphi\in V^*$ to its restriction {to} $W$, is surjective. Moreover, the quotient norm on $W^*$ induced by the dual norm $\norm{\ndot}_{V,*}$ coincides with the norm $\norm{\ndot}_{W,*}$.
\end{lemm}
\begin{proof}
Let $\psi$ be an element in $W^*$. If $\varphi$ is an element in $V^*$ which extends $\psi$, then clearly one has $\|\varphi\|_{V,*}\geqslant\|\psi\|_{W,*}$. Moreover, by Hahn-Banach theorem, there exists $\varphi_0\in V^*$ which extends $\psi$ and such that $\|\varphi_0\|_{V,*}=\|\psi\|_{W,*}$. Therefore, the map $F$ is surjective and the quotient norm on $W^\vee$ induced by $\norm{\ndot}_{V,*}$ coincides with $\norm{\ndot}_{W,*}$.
\end{proof}

\begin{prop}\label{Pro: restriction of epsion tensors}
Let $(V,\norm{\ndot}_V)$ and $(W,\norm{\ndot}_W)$ be finite-dimensional seminormed vector spaces over $k$, $V_0$ be a $k$-vector subspace of $V$ and $\norm{\ndot}_{V_0}$ be the restriction of $\norm{\ndot}_V$ on $V_0$. Denote by $\norm{\ndot}_{\varepsilon}$ and $\norm{\ndot}_{\pi}$ the $\varepsilon$-tensor product and the $\pi$-tensor product of the seminorms $\norm{\ndot}_{V}$ and $\norm{\ndot}_W$, respectively.
{
\begin{enumerate}[label=\rm(\arabic*)]
\item\label{Item: restriction epsilon tensor}
Assume that, either the absolute value $|\ndot|$ is Archimedean, or the seminorm $\norm{\ndot}_V$ is ultrametric. Then the $\varepsilon$-tensor product $\norm{\ndot}_{\varepsilon,0}$ of $\norm{\ndot}_{V_0}$ and $\norm{\ndot}_W$ identifies with the restriction of $\norm{\ndot}_{\varepsilon}$ {to} $V_0\otimes_kW$.
\item\label{Item: restriction pi tensor} Assume that the {seminorm} $\norm{\ndot}_V$ is either ultrametric or induced by a semidefinite inner product. Then the $\pi$-tensor product $\norm{\ndot}_{\pi,0}$ of $\norm{\ndot}_{V_0}$ and $\norm{\ndot}_W$ coincides with the restriction of $\norm{\ndot}_{\pi}$ {to} $V_0\otimes_kW$.
\end{enumerate}
}
\end{prop}
\begin{proof} \ref{Item: restriction epsilon tensor} Let $\varphi$ be a tensor in $V_0\otimes_kW$, viewed as a bilinear form on $V_0^*\times W^*$. By definition, one has
\[\norm{\varphi}_{\varepsilon,0}=\sup_{\begin{subarray}{c}(f_0,g)\in V_0^*\times W^*\\
f_0\neq 0,\,g\neq 0
\end{subarray}}\frac{|\varphi(f_0,g)|}{\norm{f_0}_{V_0,*}\cdot\norm{g}_{W,*}}.\]
Since the absolute value $|\ndot|$ is Archimedean or the norm $\norm{\ndot}_V$ is ultrametric, by Proposition \ref{Pro:quotientdualnonarch} (for the ultrametric case) and Lemma \ref{Lem:quotientequalsrestriction} (for the Archimedean case), the norm $\norm{\ndot}_{V_0,*}$ identifies with the quotient of $\norm{\ndot}_{V,*}$ by the canonical surjective map $V^*\rightarrow V_0^*$. Therefore, one has 
\[\sup_{\begin{subarray}{c}(f_0,g)\in V_0^*\times W^*\\
f_0\neq 0,\,g\neq 0
\end{subarray}}\frac{|\varphi(f_0,g)|}{\norm{f_0}_{V_0,*}\cdot\norm{g}_{W,*}}=\sup_{\begin{subarray}{c}(f,g)\in V^*\times W^*\\
f\neq 0,\,g\neq 0
\end{subarray}}\frac{|\varphi(f,g)|}{\norm{f}_{V,*}\cdot\norm{g}_{W,*}},\]
which shows $\norm{\varphi}_{\varepsilon,0}=\norm{\varphi}_{\varepsilon}$. 

\ref{Item: restriction pi tensor} We have already seen in Proposition \ref{Pro: restriction and tensors} \ref{Item: sub of pi tensor product} that $\norm{\ndot}_{\pi,0}$ is bounded from below by the restriction of $\norm{\ndot}_{\pi}$ {to} $V_0\otimes_kW$. Let $T$ be an element of $V_0\otimes_kW$, which is written, as an element of $V\otimes_kW$, in the form
$T=\sum_{i=1}^Nx_i\otimes y_i$,
where $\{x_1,\ldots,x_N\}\subseteq V$ and $\{y_1,\ldots,{y_N}\}\subseteq W$. By Corollary \ref{Cor:projection}, for any $\alpha\in\intervalle{]}{0}{1}{[}$, there exists a linear projection $\pi_\alpha:V\rightarrow V_0$ such that $\norm{\pi_\alpha}\leqslant\alpha^{-1}$. Since $T$ belongs to $V_0\otimes_kW$ one has
$T=\sum_{i=1}^N\pi_\alpha(x_i)\otimes y_i$.
Moreover,
\[\norm{T}_{0,\pi}\leqslant\sum_{i=1}^N\norm{\pi_\alpha(x_i)}_{V_0}\cdot\norm{y_i}_W\leqslant\alpha^{-1}\sum_{i=1}^N\norm{x_i}_V\cdot\norm{y_i}_W.\]
Since $\alpha$ and the writing $T=\sum_{i=1}^Nx_i\otimes y_i$ are arbitrary, we obtain $\norm{T}_{0,\pi}\leqslant\norm{T}_\pi$.
\end{proof}

\begin{prop}\label{Pro: pi tensorial preserve split tensor} 
Let $(V,\norm{\ndot}_V)$ and $(W,\norm{\ndot}_W)$ be seminormed vector spaces over $k$, and $\norm{\ndot}_{\pi}$ be the $\pi$-tensor product norm of $\norm{\ndot}_V$. We assume that $\norm{\ndot}_W$ is ultrametric. For any $(x,y)\in V\times W$, one has $\norm{x\otimes y}_\pi=\norm{x}_V\cdot\norm{y}_W$.
\end{prop}
\begin{proof}
By definition on has $\norm{x\otimes y}_\pi\leqslant\norm{x}_V\cdot\norm{y}_W$. It then suffices to show that, for any writing of $x\otimes y$ as 
\[\sum_{i=1}^Nx_i\otimes y_i,\]
with $(x_1,\ldots,x_n)\in V^n$ and $(y_1,\ldots,y_n)\in W^n$, one has
\[\norm{x}_V\cdot\norm{y}_W\leqslant\sum_{i=1}^N\norm{x_i}_V\cdot\norm{y_i}_W.\]
Therefore we may assume without loss of generality that $V$ and $W$ are finite-dimensional vector spaces over $k$. Consider the $k$-linear map $\ell$ from $W^*$ to $V$ sending $\varphi\in W^*$ to 
\[\varphi(y)x=\sum_{i=1}^N\varphi(y_i)x_i.\] 
We equip $W^*$ with the dual norm $\norm{\ndot}_{W,*}$ and consider the operator norm of $\ell$. On one hand, one has
\[\norm{\ell}=\sup_{\varphi\in W^*\setminus\{0\}}\frac{\norm{\varphi(y)x}_V}{\norm{\varphi}_{W,*}}=\sup_{\varphi\in W^*\setminus\{0\}}\frac{|\varphi(y)|\cdot\norm{x}_V}{\norm{\varphi}_{W,*}}=\norm{y}_{W,**}\cdot\norm{x}_V=\norm{y}_{W}\cdot\norm{x}_V,\]
where the last equality comes from Corollary \ref{Cor:doubledual} and the hypothesis that $\norm{\ndot}_W$ is ultrametric. On the other hand, one has
\[\begin{split}\norm{\ell}&=\sup_{\varphi\in W^*\setminus\{0\}}\frac{\norm{\varphi(y_1)x_1+\cdots+\varphi(y_N)x_N}_V}{\norm{\varphi}_{W,*}}\\
&\leqslant\sup_{\varphi\in W^*\setminus\{0\}}\sum_{i=1}^N\frac{|\varphi(y_i)|\cdot\norm{x_i}_V}{\norm{\varphi}_{W,*}}\leqslant\sum_{i=1}^N\norm{x_i}_{V}\cdot\norm{y_i}_{W,**}=\sum_{i=1}^N\norm{x_i}_{V}\cdot\norm{y_i}_{W},
\end{split}\]
where the last equality follows from Corollary \ref{Cor:doubledual} and the hypothesis that $\norm{\ndot}_W$ is ultrametric again. The proposition is thus proved.
\end{proof}

The following result provides a variant of Proposition \ref{Pro:minormationdenorme}. Note that it generalises (by using Proposition \ref{Pro:delta=1}) Proposition \ref{Pro:exactsequence}.

\begin{prop}\label{Pro:exactesequenceanddelta}
Let $(V,\norm{\ndot})$ be a finite-dimensional normed vector space over $k$. Assume that the absolute value $|\ndot|$ is Archimedean or the norm $\norm{\ndot}$ is ultrametric. For any vector subspace $W$ of $V$, the norm of the canonical isomorphism
\[f:\det(W)\otimes\det(V/W)\longrightarrow\det(V)\]
is bounded from below by
\[\frac{\delta(W,\norm{\ndot}_W)\delta({V/W},\norm{\ndot}_{V/W})}{\delta(V,\norm{\ndot})}\geqslant\delta(V,\norm{\ndot})^{-1},\]
where we consider the restriction $\norm{\ndot}_{W}$ of the norm $\norm{\ndot}$ {to} the vector subspace $W$ and the quotient norm $\norm{\ndot}_{V/W}$ of $\norm{\ndot}$ on the quotient space $V/W$. In particular, one has
\[\max\big\{\delta(W,\norm{\ndot}_W),\delta({V/W},\norm{\ndot}_{V/W})\big\}\leqslant\delta(V,\norm{\ndot}).\]
\end{prop}
\begin{proof} Let $\norm{\ndot}_{V/W}$ be the quotient norm on $V/W$ induced by $\norm{\ndot}_V$.
By Proposition \ref{Pro:dualquotient}, the dual norm $\norm{\ndot}_{V/W,*}$ coincides with the restriction of the norm $\norm{\ndot}_*$ {to} $(V/W)^\vee$. Moreover, by Lemma \ref{Lem:quotientequalsrestriction} (for the Archimedean case) and Proposition \ref{Pro:quotientdualnonarch} (for the non-Archimedean case), the quotient norm on $W^\vee$ induced by $\norm{\ndot}_*$ identifies with the dual norm $\norm{\ndot}_{W,*}$. Let $\alpha$ and $\beta$ be respectively non-zero elements in $\det(W)$ and $\det(V/W)$. Let $\alpha^\vee\in\det(W^\vee)$ and $\beta^\vee\in\det((V/W)^\vee)$ be their dual elements, $\eta$ be the image of $\alpha\otimes\beta$ by the canonical isomorphism $\det(W)\otimes\det(V/W)\rightarrow\det(V)$, and $\eta^\vee$ be the image of $\alpha^\vee\otimes\beta^\vee$ by the canonical isomorphism $\det(W^\vee)\otimes\det((V/W)^\vee)\rightarrow\det(V^\vee)$. Then $\eta^\vee$ is the dual element of $\eta$.

By Proposition \ref{Cor:exactsequencenorm}, one has
\[\|\eta^\vee\|_{*,\det}\leqslant\|\alpha^\vee\|_{W,*,\det}\cdot\|{\beta}^\vee\|_{V/W,*,\det}.\]
Hence by \eqref{Equ:definitiondedelta} one has
\[\begin{split}\frac{\delta(W,\norm{\ndot}_W)\delta({V/W},\norm{\ndot}_{V/W})}{\delta(V,\norm{\ndot})}&=\frac{\|\eta^\vee\|_{*,\det}\cdot\|\eta\|_{\det}}{\|{\alpha}^\vee\|_{W,*,\det}\|\alpha\|_{W,\det}\cdot\|{\beta}\|_{V/W,\det}\|\beta^\vee\|_{V/W,*,\det}}
\\&\leqslant\frac{\|\eta\|_{\det}}{\|\alpha\|_{W,\det}\cdot\|{\beta}\|_{V/W,\det}}=\|f\|.
\end{split}\]
Finally, by Corollary \ref{Cor:exactsequencenorm}, we obtain
\[\delta(W,\norm{\ndot}_W)\delta({V/W},\norm{\ndot}_{V/W})\leqslant\delta(V,\norm{\ndot}).\]
Since $\delta({W},\norm{\ndot}_W)$ and $\delta({V/W},\norm{\ndot}_{V/W})$ are $\geqslant 1$ (see Proposition \ref{Pro:minorationdelta}), we obtain the last inequality.
\end{proof}

\begin{coro}\label{Cor:convexitenorme}
Let $V$ be a finite-dimensional vector space over $k$ and $\norm{\ndot}$ be a norm on $V$. We assume that, either the norm $\norm{\ndot}$ is ultrametric or the absolute value $|\ndot|$ is Archimedean. If $W_1$ and $W_2$ are two $k$-vector subspaces of $V$, then the canonical isomorphism
\begin{equation}\label{Equ: convex metric}\det(W_1)\otimes\det(W_2)\longrightarrow\det(W_1\cap W_2)\otimes\det(W_1+W_2)\end{equation}
induced by the short exact sequence
\[\xymatrix{0\ar[r]&W_1\cap W_2\ar[r]&W_1\oplus W_2\ar[r]&W_1+W_2\ar[r]&0}\]
has operator norm $\leqslant \min\{\delta(W_1),\delta(W_2)\}/\delta(W_1\cap W_2)$, where in the above formulae we consider the restricted norms on the vector subspaces of $V$. In particular, if $\norm{\ndot}$ is an ultrametric norm, then the linear map \eqref{Equ: convex metric} has norm $\leqslant 1$.
\end{coro}
\begin{proof}
Consider the short exact sequence
\[\xymatrix{0\ar[r]&W_1\cap W_2\ar[r]&W_1\ar[r]&W_1/(W_1\cap W_2)\ar[r]&0}.\]
By Proposition \ref{Pro:exactesequenceanddelta}, the canonical element $\eta$ in \[\det(W_1)^\vee\otimes\det(W_1\cap W_2)\otimes\det(G)\]
has norm $\leqslant\delta(W_1)/\delta(W_1\cap W_2)\delta(G)$, where $G$ denotes the vector space $W_1/(W_1\cap W_2)$ equipped with the quotient norm $\norm{\ndot}_G$.

Similarly, consider the short exact sequence
\[\xymatrix{0\ar[r]&W_2\ar[r]&W_1+W_2\ar[r]&(W_1+W_2)/W_2\ar[r]&0}.\]
By Corollary \ref{Cor:exactsequencenorm}, the canonical element $\eta'$ in
\[\det(W_2)^\vee\otimes\det(G')^\vee\otimes\det(W_1+W_2)\]
has norm $\leqslant 1$, where $G'$ denotes the vector space $(W_1+W_2)/W_2$ equipped with the quotient norm $\norm{\ndot}_{G'}$. Therefore we obtain
\[\|\eta\otimes \eta'\|\leqslant\delta(W_1)/\delta(W_1\cap W_2).\] 

Let $f:G\rightarrow G'$ be the canonical isomorphism. One has $\|f(x)\|_{G'}\leqslant \|x\|_G$ for any $x\in G$. In particular, the canonical element of $\det(G)\otimes\det(G')^\vee$ has norme $\geqslant 1$. We deduce that the canonical element of
\[\det(W_1)^\vee\otimes\det(W_2)^\vee\otimes\det(W_1\cap W_2)\otimes\det(W_1+W_2)\]
has norm $\leqslant\delta(W_1)/\delta(W_1\cap W_2)$. By the symmetry between $W_1$ and $W_2$, we  obtain the announced inequality. 
\end{proof}

\begin{rema}
We assume that the absolute value $|\ndot|$ is non-Archimedean. The result of Corollary \ref{Cor:convexitenorme} is not true in general if the norm $\norm{\ndot}$ is not ultrametic. However, we can combine the proof of Corollary \ref{Cor:convexitenorme} and Proposition \ref{Pro:minormationdenorme} to show that the canonical isomorphism \eqref{Equ: convex metric} in Corollary \ref{Cor:convexitenorme} has an operator norm bounded from above by
\[\frac{\min\{\Delta(W_1),\Delta(W_2)\}}{\Delta(W_1\cap W_2)}.\] The same argument also works in the Archimedean case.
\end{rema}

\subsection{Ellipsoid of John and L\"owner}\label{Sec:John-Lowner} We assume that the absolute value $|\ndot|$ is Archimedean.
Let $V$ be a finite-dimensional vector space over $k$, equipped with a norm $\norm{\ndot}$. In this subsection, we discuss the approximation of the norm $\norm{\ndot}$ by Euclidean or Hermitian norms. Note that Proposition \ref{Pro:existenceoforthogonal} provides a result in this direction. Let $\{e_i\}_{i=1}^r$ be an orthonormal basis of $V$. Let $\emptyinnprod$ 
be an inner product on $V$ such that $\{e_i\}_{i=1}^r$ is orthogonal with respect to the inner product, and that $\langle e_i,e_i\rangle=r$ for any $i\in\{1,\ldots,r\}$. If $\norm{\ndot}_h$ denotes the norm on $V$ induced by the inner product $\emptyinnprod$, then for any $x=\lambda_1e_1+\cdots+\lambda_re_r\in V$ one has\[\frac{1}{r}\|x\|_h=\Big(\frac{|\lambda_1|^{2}+\cdots+|\lambda_r|^{2}}{r}\Big)^{1/2}\leqslant\max\{|\lambda_1|,\ldots,|\lambda_r|\}\leqslant\|x\|\] 
and
$\|x\|\leqslant|\lambda_1|+\cdots+|\lambda_r|\leqslant r^{1/2}\big(|\lambda_1|^{2}+\cdots+|\lambda_r|^{2}\big)^{1/2}
=\|x\|_h$.

The works of John \cite{John48} and L\"owner provide a stronger result on the comparison of inner product norms and general norms. We refer to the expository article of Henk \cite{Henk12} for the history of this theory. For the convenience of the readers, we include the statement and the proof of this result.

\begin{theo}[John-L\"owner]\label{Thm:John}
Let $V$ be a non-zero finite-dimensional vector space over $k$, equipped with a norm $\norm{\ndot}$. There exists a unique Euclidean or Hermitian norm $\norm{\ndot}_J$ bounded from above by $\norm{\ndot}$ such that, for any Euclidean or Hermitian norm $\norm{\ndot}_h$ satisfying $\norm{\ndot}_h\leqslant\norm{\ndot}$, one has $\norm{\ndot}_{h,\det}\leqslant\norm{\ndot}_{J,\det}$. Moreover, for any $x\in V$, one has $\|x\|_h\leqslant\|x\|\leqslant {r^{1/2}}\|x\|_h$,
where $r$ is the rank of $V$ over $k$.
\end{theo}
\begin{proof}
We fix an arbitrary inner product $\emptyinnprod'$
on $V$ and denote by $\Theta$ the vector space (over $\mathbb R$) of all endomorphisms of $V$ which are self-adjoint with respect to the inner product $\emptyinnprod'$. 
Recall that a $k$-linear map $u:V\rightarrow V$ is said to be \emph{self-adjoint}\index{self-adjoint} with respect to $\emptyinnprod'$
if and only if
\[\forall\,x,y\in V,\quad \langle u(x),y\rangle'=\langle x,u(y)\rangle'.\]
Let $\Theta^+$ be the set of all positive definite self-adjoint operators. Since any pair of self-adjoint operator can be simultaneously diagonalised by a basis of $V$, we obtain that $\Theta^+$ is a convex open subset of $\Theta$ and that the function $\log\det(\ndot)$ is strictly concave on $\Theta^+$.

Let $B=\{x\in V\,:\,\|x\|\leqslant 1\}$ be the unit ball of the norm $\norm{\ndot}$. For any $u\in \Theta^+$, let $B_u=\{x\in V\,:\,\langle x,u(x)\rangle'\leqslant 1\}$, which is the unit ball of the Euclidean or Hermitian norm $\norm{\ndot}_u$ on $V$ defined as
\[\forall\,x\in V,\quad \|x\|_u^2=\langle x,u(x)\rangle'.\] Let $\Theta_0$ be the set of all $u\in\Theta^+$ such that $B_u\supseteq B$. Then for any $u_0\in \Theta_0$, the set \[\Theta(u_0):=\{u\in\Theta_0\,:\,\det(u)\geqslant\det(u_0)\}\]
is a convex and compact subset of $\Theta$. In fact, from the concavity and the continuity of the function $\log\det(\cdot)$ we obtain that the set $\Theta(u_0)$ is convex and closed. Moreover, the condition $B_u\supseteq B$ for $u\in\Theta_0$ implies that the set $\Theta(u_0)$ is bounded in $\Theta$. Therefore the restriction of the function $\det(\cdot)$ {to} $\Theta_0$ attains its maximal value on a unique point $u_1\in\Theta_0$.

Let $\emptyinnprod$  
be the inner product on $V$ such that
\[\forall\,(x,y)\in V\times V,\quad\langle x,y\rangle=\langle x,u_1(y)\rangle'.\]We call it the \emph{John inner product}\index{John inner product} associated with the norm $\norm{\ndot}$. The corresponding Euclidean or Hermitian norm $\norm{\ndot}_{J}$ is called the \emph{John norm}\label{Page:John norm}\index{John norm} associated with $\norm{\ndot}$. 

In the following, we prove the relation 
\[\forall\,x\in V,\quad \|x\|_J\leqslant\|x\|\leqslant r^{1/2}\|x\|_J\] 
under the supplementary assumption that the unit ball $B$ is the convex hull of finitely many orbits of the action of $\{a\in k\,:\,|a|=1\}$ on $V$.

Without loss of generality, we assume that $\emptyinnprod' = \emptyinnprod$.
For any $x\in V$ such that $\|x\|\leqslant 1$, let $\varphi_x:\Theta\rightarrow\mathbb R$ be the linear functional which sends $u\in\Theta$ to $\langle x,u(x)\rangle$. If $u:V\rightarrow V$ is a self-adjoint linear operator such that $\varphi_x(u)\leqslant 0$ for any $x\in B$ such that $\langle x,x\rangle=1$, then one has $\mathrm{Tr}(u)\leqslant 0$. In fact, the condition \[\forall\,x\in B,\quad\langle x,x\rangle=1\Longrightarrow \varphi_x(u)\leqslant 0\] implies that $\mathrm{Id}+\varepsilon u\in\Theta_0$ for sufficiently small $\varepsilon>0$ (here we use the supplementary assumption that the convex body $B$ is spanned by a finite number of orbits). Therefore one has $\det(\mathrm{Id}+\varepsilon u)\leqslant\det(\mathrm{Id})=1$, which leads to $\mathrm{Tr}(u)\leqslant 0$. Therefore, the linear form $\mathrm{Tr}(\cdot)$ lies in the closure of the positive cone of $\Theta^\vee$ generated by $\varphi_x(\cdot)$ ($x\in B$, $\langle x,x\rangle=1$), namely there exist a sequence of elements $\{x_n\}_{n\in\mathbb N}$ in \[B\cap\{x\in V\,:\,\langle x,x\rangle=1\}\] and a sequence $\{\lambda_n\}_{n\geqslant 0}$ of  real numbers such that 
\begin{equation}\label{Equ:convegencetrace}\mathrm{Tr}(u)=\sum_{n\in\mathbb N}\lambda_n\langle x_n,u(x_n)\rangle\end{equation}
for any $u\in\Theta$. If we apply the identity to $u=\mathrm{Id}$, we obtain
\begin{equation}\label{Equ:formular}r=\sum_{n\in\mathbb N}\lambda_n\langle x_n,x_n\rangle=\sum_{n\in\mathbb N}\lambda_n.\end{equation}
Let $y$ be an element in $V$ such that $\langle y,y\rangle=1$. We apply the identity \eqref{Equ:convegencetrace} to the linear map $u(x)=\langle y,x\rangle y $, and obtain
\[1=\sum_{n\in\mathbb N}\lambda_n|\langle x_n,y\rangle|_\infty^2.\]
Thus there should exist $n\in\mathbb N$ such that $|\langle x_n,y\rangle|_\infty\geqslant r^{-1/2}$ since otherwise we have
\[1<\sum_{n\in\mathbb N}\frac{1}{r}\lambda_n=\frac 1r\cdot r=1,\]
where the first equality comes from \eqref{Equ:formular}, 
which leads to a contradiction. Since the unit ball $B=\{x\in V\,:\,\|x\|=1\}$ is invariant by the multiplication by any $\lambda\in k$ with $|\lambda|=1$, we obtain that, for any $y\in V$ such that $\langle y,y\rangle= 1$, there exists $x\in B$ such that $\mathrm{Re}\langle y,x\rangle\geqslant r^{-1/2}$.

We claim that the unit ball $B=\{x\in V\,:\,\|x\|\leqslant 1\}$ contains the set of all $x\in V$ such that $\langle x,x\rangle \leqslant 1/r$. In fact, if $x_0\in V$ is a point such that $\langle x_0,x_0\rangle\leqslant 1/r$ and that $\|x\|>1$, we can choose an  $\mathbb R$-affine function $f:V\rightarrow \mathbb R$ such that $f(x_0)=0$ and that $f(x)<0$ for any $x\in B$. Note that $\mathrm{Re}\emptyinnprod$ 
defines an inner product on $V$, where $V$ is viewed as a vector space over $\mathbb R$ if $k=\mathbb C$. By Riesz's theorem there exists $y\in V$ such that 
\[\forall\,x\in V,\quad f(x)=\mathrm{Re}\langle y,x\rangle+f(0).\]
Without loss of generality, we may assume that $\langle y,y\rangle=1$. One has
\[0=f(x_0)=\mathrm{Re}\langle y,x_0\rangle+f(0)\leqslant\langle y,y\rangle^{1/2}\langle x_0,x_0\rangle^{1/2}+f(0)=\frac{1}{\sqrt{r}}+f(0).\]
Hence $f(0)\geqslant -r^{-1/2}$. However, the above argument shows that there exists $x\in B$ such that $\mathrm{Re}\langle y,x\rangle\geqslant r^{-1/2}$.
Hence one has
\[0>f(x)=\mathrm{Re}\langle y,x\rangle+f(0)\geqslant 0,\]
which leads a contradiction.

Since $B\subseteq\{x\in V\,:\,\langle x,x\rangle\leqslant 1\}$, one has $\|x\|_J\leqslant\|x\|$ for any $x\in V$. Moreover, the relation 
\[\{x\in V\,:\,\langle x,x\rangle\leqslant 1/r\}\subseteq B=\{x\in V\,:\,\|x\|\leqslant 1\}\]
implies that $\|x\|\leqslant r^{1/2}\|x\|_J$. The theorem is thus proved under the supplementary hypothesis.

For the general case, we can construct a decreasing sequence of norms $\{\norm{\ndot}_n\}_{n\in\mathbb N}$ such that each unit ball $\{x\in V\,:\,\|x\|_n\leqslant 1\}$ verifies the supplementary hypothesis mentioned above and that the sequence
\[\sup_{0\neq x\in V}\frac{\|x\|_n}{\|x\|}\]
converges to $1$ when $n\rightarrow+\infty$. For each $n\in\mathbb N$, let $\norm{\ndot}_{n,J}$ be the John norm associated to the norm $\norm{\ndot}_n$. If we identify the set of Euclidean or Hermitian norms on $V$ with $\Theta^+$, we obtain that these John norms actually lies in a bounded subset of $\Theta$. Therefore there exists a subsequence of $\{\norm{\ndot}_{n,J}\}_{n\in\mathbb N}$ which converges in $\Theta$, whose limite should be the John norm associated with $\norm{\ndot}$ by the uniqueness of the John norm. 
Without loss of generality we may assume that $\{\norm{\ndot}_{n,J}\}_{n\in\mathbb N}$ converges in $\Theta$. By what we have established above, for any $n\in\mathbb N$ one has
\[\forall\,x\in V,\quad\|x\|_{n,J}\leqslant\|x\|_{n}\leqslant r^{1/2}\|x\|_{n,J}.\]
By taking the limit when $n\rightarrow +\infty$, we obtain the result announced in the theorem.
\end{proof}

\begin{rema}\label{rema:Thm:John}
Let $(V,\norm{\ndot})$ be a finite-dimensional normed vector space over $\mathbb R$ or $\mathbb C$ (equipped with the usual absolute value). Since $(V,\norm{\ndot})$ is  reflexive (see Proposition \ref{Pro:doubledualarch}), we deduce that the dual norm $\norm{\ndot}_{J,*}$ is the unique norm on $V^\vee$ which is bounded from below by $\norm{\ndot}_*$ and such that the corresponding determinant norm $\norm{\ndot}_{J,*,\det}$ is minimal. In particular, one has
\[\Delta(V^\vee,\norm{\ndot}_*)=\frac{\norm{\ndot}_{J,*,\det}}{\norm{\ndot}_{*,\det}}.\]
Similarly, one has
\begin{equation}\label{Equ:deltacommequotient}\Delta(V,\norm{\ndot})=\frac{\norm{\ndot}_{L,\det}}{\norm{\ndot}_{\det}},\end{equation}
where $\norm{\ndot}_L$ is the unique Euclidean or Hermitian norm on $V$ which is bounded from below by $\norm{\ndot}$ and such that $\norm{\ndot}_{L,\det}$ is minimal (called the \emph{L\"{o}wner norm}\index{Lowner norm@L\"{o}wner norm} of $\norm{\ndot}$), which is also equal to $\norm{\ndot}_{*,J,*}$. Theorem \ref{Thm:John} then leads to
\begin{equation}\label{Equ:majorationdelta}
\max\{\Delta(V,\norm{\ndot}),\Delta(V^\vee,\norm{\ndot}_*)\}\leqslant \rang(V)^{\rang(V)/2}
\end{equation}
We denote by $\lambda(V,\norm{\ndot})$ the constant
$\|\eta\|_{L,\det}\cdot\|\eta^\vee\|_{J,*,\det}$,
where $\eta$ is an arbitrary non-zero element in $\det(V)$, and $\eta^\vee$ is its dual element in $\det(V^\vee)$. With this notation, by \eqref{Equ:definitiondedelta} in \S\ref{Subsec:dualdet} one has
\begin{equation}\label{Equ:lesdeuxdelta}\Delta(V^\vee,\norm{\ndot}_*)\Delta(V,\norm{\ndot})=\lambda(V,\norm{\ndot}){\delta(V,\norm{\ndot})}.\end{equation}
Note that one has 
$\norm{\ndot}_J\leqslant \norm{\ndot}$
by definition. Hence we obtain \begin{equation}\label{Equ:lambdaV}\lambda(V,\norm{\ndot})=\frac{\norm{\ndot}_{L,\det}}{\norm{\ndot}_{J,\det}}\geqslant\frac{\norm{\ndot}_{L,\det}}{\norm{\ndot}_{\det}}=\Delta(V,\norm{\ndot}),\end{equation}
where the first equality comes from Proposition \ref{Pro:delta=1}.
Therefore the relation \eqref{Equ:lesdeuxdelta} leads to $\Delta(V^\vee,\norm{\ndot}_*)\geqslant\delta(V,\norm{\ndot})$. Since $\delta(V,\norm{\ndot})$ and $\lambda(V,\norm{\ndot})$ are both invariant by duality, one obtains
\begin{equation}\label{Equ:majorationdelta2}\begin{split}\delta(V,\norm{\ndot})&\leqslant\min\{\Delta(V,\norm{\ndot}),\Delta(V^\vee,\norm{\ndot}_*)\}\\
&\leqslant\max\{\Delta(V,\norm{\ndot}),\Delta(V^\vee,\norm{\ndot}_*)\}\leqslant\lambda(V,\norm{\ndot}).
\end{split}
\end{equation}
\end{rema}

{
\begin{rema}
We can deduce from Theorem \ref{Thm:John} a similar result for seminorms. Let $(V,\norm{\ndot})$ be a finite-dimensional seminormed vector space over $k$. Let $\norm{\ndot}_J^\sim$ be the John norm associated with $\norm{\ndot}^\sim$. It is induced by an inner product on $V/N_{\norm{\ndot}}$. Let  $\norm{\ndot}_J$ be the seminorm on $V$ given by the composition of $\norm{\ndot}_J^\sim$ with the canonical projection $V\rightarrow V/\norm{\ndot}$. It is a seminorm induced by a  semidefinite inner product. Moreover, the following inequalities hold 
\[\norm{\ndot}_J\leqslant\norm{\ndot}\leqslant \rang_k(V/N_{\norm{\ndot}})^{1/2}{\norm{\ndot}_J}.\] 
\end{rema}}

\subsection{Hilbert-Schmidt tensor norm}\label{Subsec:Hilbert-Schmidt} In this subsection, we assume that the absolute value $|\ndot|$ is Archimedean.

Let $V$ and $W$ be finite-dimensional vector spaces over $k$, equipped with  semidefinite inner products. 
For $f \in \Hom_k(V^*,W)$, the adjoint operator $f^* : W \to V^*$ of $f$ is defined by 
$\langle f(\alpha), y \rangle = \langle \alpha, f^*(y) \rangle_{*}$
for all $\alpha \in V^{*}$ and $y \in W$. Note that the adjoint operator $f^*$ exists for any $f$ because
the product $\emptyinnprod_*$ on $V^*$ is positive definite.
We can equip $\Hom_k(V^*,W)$ with the following  semidefinite inner product $\emptyinnprod_{\mathrm{HS}}$:
\[\forall\,{f, g} 
\in\Hom_k (V^*,W),\quad \langle f,g\rangle_{\mathrm{HS}}:=\mathrm{Tr}(f^*\circ g).\]
This  semidefinite inner product defines a seminorm on $\Hom_k(V^*,W)$, which induces by the canonical linear map $V\otimes_kW\rightarrow\Hom_k(V^*,W)$ a seminorm $\norm{\ndot}_{\mathrm{HS}}$ on $V\otimes_KW$, called the \emph{orthogonal tensor product}\index{orthogonal tensor product}\index{seminorm!orthogonal tensor product} of the seminorms of $V$ and $W$, or \emph{Hilbert-Schmidt seminorm}\index{Hilbert-Schmidt seminorm}\index{seminorm!Hilbert-Schmidt ---}.
Note that if $\{x_i\}_{i=1}^n$ and $\{y_j\}_{j=1}^m$ are respectively orthogonal basis of $V$ and $W$, then $\{x_i\otimes y_j\}_{i\in\{1,\ldots,n\},\,j\in\{1,\ldots,m\}}$ is an orthogonal basis of $V\otimes_kW$ with respect to $\emptyinnprod_{\mathrm{HS}}$.  
Moreover, for $x\in V$ and $y\in W$ one has
\begin{equation}\label{Equ:HStensorprod}\|x\otimes y\|_{\mathrm{HS}}=\|x\|\cdot\|y\|.\end{equation} In particular, if $V$ and $W$ are both of rank $1$ over $k$, then the orthogonal tensor product seminorm on $V\otimes_kW$ coincides with the $\varepsilon$-tensor product and the $\pi$-tensor product seminorms. In this case we just call it the \emph{tensor product seminorm}\index{tensor product seminorm}\index{seminorm!tensor product}.

The dual norm on $V^*\otimes_k W^*$ of the Hilbert-Schmidt seminorm on $V\otimes_k W$ coincides with the orthogonal tensor product of the dual norms on $V^*$ and $W^*$. Moreover, the orthogonal tensor product is commutative, namely the isomorphism from $V\otimes_kW$ to $W\otimes_kV$ given by the transposition is actually an isometry under orthogonal tensor product seminorms. Similarly, the  orthogonal tensor product is associative. More precisely, given three finite-dimensional vector spaces $U$, $V$ and $W$ over $k$, equipped with  semidefinite inner products, the natural isomorphism from $(U\otimes_kV)\otimes_kW$ to $U\otimes_k(V\otimes_kW)$ is an isometry for orthogonal tensor product seminorms.

The following assertion, which is similar to Proposition \ref{Pro:quotientr1eps}, studies the quotient norm of the orthogonal tensor product.

\begin{prop}\label{Pro:quotientr1Hs}
Let $V$ and $W$ be finite-dimensional seminormed vector spaces over $k$, $V_0$ be a $k$-vector subspace of $V$, and $Q$ be the quotient vector space $V/V_0$ equipped with the quotient seminorm. We assume that the seminorms of $V$ and $W$ are induced by  semidefinite inner products. Then the canonical isomorphism $(V\otimes_kW)/(V_0\otimes_kW)\rightarrow Q\otimes_kW$ is an isometry, where we consider the orthogonal tensor product seminorms on $V\otimes_kW$ and $Q\otimes_kW$, and the quotient seminorm on $(V\otimes_kW)/(V_0\otimes_kW)$.
\end{prop}
\begin{proof} By the Gram-Schmidt process we can identify the quotient space $Q$ with the orthogonal {supplementary} of $V_0$ in $V$.
Let $\boldsymbol{e}=\{e_i\}_{i=1}^n$ be an orthogonal basis of $V$ such that $\card(\boldsymbol{e}\cap V_0)=\rang_k(V_0)$. Then the projection $V\rightarrow Q$ defines an isometry  between $Q$ and the vector subspace $V_1$ of $V$ generated by $\boldsymbol{e}\setminus V_0$. Let $\boldsymbol{f}=\{f_j\}_{j=1}^m$ be an orthogonal basis of $W$. Then the basis $\boldsymbol{e}\otimes\boldsymbol{f}=\{e_i\otimes f_j\}_{(i,j)\in\{1,\ldots,n\}\times\{1,\ldots,m\}}$  of $V\otimes_kW$ is orthogonal. Moreover, one has \[\card((\boldsymbol{e}\otimes\boldsymbol{f})\cap(V_0\otimes_kW))=\rang_k(V_0\otimes_kW).\]
Thus $(\boldsymbol{e}\setminus V_0)\otimes\boldsymbol{f}$ forms an orthogonal basis of $Q\otimes_kW$ equipped with the quotient seminorm (where we identify $Q$ with $V_0^\perp$). Hence the quotient seminorm on $Q\otimes_kW$ identifies with the orthogonal tensor product seminorm.
\end{proof}

\begin{prop}\label{Pro: pi tensor sub Archimedean} Let $(V,\norm{\ndot}_V)$ and $(W,\norm{\ndot}_W)$ be finite-dimensional seminormed vector space over $k$, $V_0$ be a $k$-vector subspace of $V$ and $\norm{\ndot}_{V_0}$ be the restriction of $\norm{\ndot}_V$ on $V_0$. We assume that the absolute value $|\ndot|$ is Archimedean and that the seminorms $\norm{\ndot}_V$ and $\norm{\ndot}_W$ are induced by semidefinite inner products. Let $\norm{\ndot}$ be the orthogonal tensor product of $\norm{\ndot}_V$ and $\norm{\ndot}_W$, and $\norm{\ndot}_0$ be the orthogonal tensor product of $\norm{\ndot}_{V_0}$ and $\norm{\ndot}_W$. Then $\norm{\ndot}_0$ identifies with the restriction of $\norm{\ndot}$ {to} $V_0\otimes_kW$.
\end{prop}
\begin{proof}
Note that $\norm{\ndot}_*$ identifies with the orthogonal tensor product of $\norm{\ndot}_{V,*}$ and $\norm{\ndot}_{W,*}$, and $\norm{\ndot}_{0,*}$ identifies with the orthogonal tensor product of $\norm{\ndot}_{V_0,*}$ and $\norm{\ndot}_{W,*}$. Moreover, by Lemma \ref{Lem:quotientequalsrestriction}, $\norm{\ndot}_{V_0,*}$ identifies with the quotient norm of $\norm{\ndot}_*$ by the canonical surjective map $V^*\rightarrow V_0^*$. By Proposition \ref{Pro:quotientr1Hs}, we obtain that $\norm{\ndot}_{0,*}$ identifies with the quotient norm of $\norm{\ndot}_*$ by the canonical surjective map $V^*\otimes_kW^*\rightarrow V_0^*\otimes_kW^*$. Therefore, by Proposition \ref{Pro:dualquotient}, $\norm{\ndot}_0$ is the restriction of $\norm{\ndot}$ {to} $V_0\otimes_kW$.
\end{proof}

The following proposition compares $\varepsilon$-tensor product to orthogonal tensor product.

\begin{prop}\label{Pro:comparisonofHSandepsnorms}
Let $V$ and $W$ be finite-dimensional vector spaces over $k$, equipped with  semidefinite inner products. Let $\norm{\ndot}_{\varepsilon}$ and $\norm{\ndot}_{\mathrm{HS}}$ be respectively the $\varepsilon$-tensor product seminorm and the orthogonal tensor product seminorm on $V\otimes_kW$. Then $\norm{\ndot}_{\varepsilon}\leqslant\norm{\ndot}_{\mathrm{HS}}\leqslant\min\{\rang_k(V^*),\rang_k(W^*)\}^{1/2}\norm{\ndot}_{\varepsilon}$. 
\end{prop}
\begin{proof} Without loss of generality, we may assume that $\rang_k(V^*)\leqslant\rang_k (W^*)$.
Let $\varphi$ be an element of $V\otimes_kW$, viewed as a $k$-linear map from $V^*$ to $W$. Let $\lambda_1\geqslant\ldots\geqslant\lambda_r$ be the eigenvalues of the positive semidefinite operator $\varphi^*\circ\varphi$. By definition, the Hilbert-Schmidt seminorm of $\varphi$ is 
$\|\varphi\|_{\mathrm{HS}}=(\lambda_1+\ldots+\lambda_r)^{1/2}$.
Moreover, the operator seminorm of $\varphi$ is $\lambda_1^{1/2}$. In fact, if $\alpha_1,\ldots,\alpha_r$ are eigenvectors of $\varphi^*\circ\varphi$ of eigenvalues $\lambda_1,\ldots,\lambda_r$, respectively, then for any $(a_1,\ldots,a_r)\in k^r$ one has
\[\begin{split}&\quad\;\langle \varphi(a_1\alpha_1+\ldots+a_r\alpha_r),\varphi(a_1\alpha_1+\ldots+a_r\alpha_r)\rangle\\&=\langle \varphi^*(\varphi(a_1\alpha_1+\ldots+a_r\alpha_r)),a_1\alpha_1+\ldots+a_r\alpha_r\rangle
=\sum_{i=1}^r|a_i|^2\lambda_i.\end{split} \]
Therefore one has
$\|\varphi\|_{\varepsilon}\leqslant\|\varphi\|_{\mathrm{HS}}\leqslant\sqrt{r}\|\varphi\|_{\varepsilon}$.
\end{proof}

By using the duality between the $\varepsilon$-tensor product and $\pi$-tensor product (see Proposition \ref{Pro:dualitypiepsilon}), we deduce from the previous proposition the following corollary.

\begin{coro}\label{Cor:comparisonofHSandpinorms}
Let $V$ and $W$ be finite-dimensional vector spaces over $k$, equipped with {semidefinite} inner products. Let $\norm{\ndot}_{\pi}$ and $\norm{\ndot}_{\mathrm{HS}}$ be respectively the $\pi$-tensor product and the orthogonal tensor product norms on $V\otimes_kW$. Then one has
\[\norm{\ndot}_{\pi}\geqslant\norm{\ndot}_{\mathrm{HS}}\geqslant\min\{\rang_k(V^*),\rang_k(W^*)\}^{-1/2}\norm{\ndot}_{\pi}.\]
\end{coro}

The following proposition expresses the Hilbert-Schmidt norm of endomorphisms in terms of the operator norm.

\begin{prop}\label{Pro:Hilbert-Schmidt}
Let $V$ be a vector space of finite rank $r$ over $k$, equipped with an inner product $\emptyinnprod$. 
Let $f:V\rightarrow V$ be an endomorphism of $V$. Then one has 
\[\langle f,f\rangle_{\mathrm{HS}}=\sum_{i=1}^r\inf_{\begin{subarray}{c}
g\in\mathrm{End}_k(V)\\
\rang(g)\leqslant i-1
\end{subarray}}\|f-g\|^{2},\]
where $\norm{\ndot}$ denotes the operator norm on $\mathrm{End}_k(V)$.
\end{prop} 
\begin{proof}
Let $\{e_i\}_{i=1}^r$ be an orthonormal basis of $V$ consisting of the eigenvectors of the self-adjoint operator $f^*\circ f$. For any $i\in\{1,\ldots,r\}$, let $\lambda_i$ be the eigenvalue of $f^*\circ f$ corresponding to the eigenvector $e_i$. Without loss of generality, we may assume that $\lambda_1\geqslant\ldots\geqslant\lambda_r$. Since the self-adjoint operator $f^*\circ f$ is positive semidefinite, one has $\lambda_r\geqslant 0$. By definition, one has
$\langle f,f\rangle_{\mathrm{HS}}=\sum_{i=1}^r\lambda_i$.
In the following, we prove that, for any $i\in\{1,\ldots, r\}$, one has
\[\inf_{\begin{subarray}{c}
g\in\mathrm{End}_k(V)\\
\rang(g)\leqslant i-1
\end{subarray}}\|f-g\|^{2}=\lambda_i.\]
Let $\pi$ be the orthogonal projection of $V$ to the vector subspace generated by $\{e_1,\ldots,e_{i-1}\}$. Then the endomorphism $f\circ\pi$ has rank $\leqslant i-1$. Moreover, since any orthogonal projection is self-adjoint, one has 
\[(f-f\pi)^*(f-f\pi)=(f^*-\pi f^*)(f-f\pi)=f^*f+\pi f^*f\pi-\pi f^*f-f^*f\pi.\]
In particular, the linear endomorphism $(f-f\pi)^*(f-f\pi)$ sends an element $a_1e_1+\cdots+a_re_r$ in $V$ to
$a_i\lambda_ie_i+\cdots+a_r\lambda_re_r$.
Hence the operator norm of $(f-f\pi)^*(f-f\pi)$, which is equal to the square of the operator norm of $f-f\pi$, is $\lambda_i$.

It remains to prove that, for any $k$-linear endomorphism $g\in\mathrm{End}_k(V)$ of rank $\leqslant i-1$, one has $\|f-g\|^{2}\geqslant\lambda_i$. Let $W$ be the vector subspace of $V$ generated by $\{e_1,\ldots,e_i\}$. Since $g$ has rank $\leqslant i-1$, one has $\mathrm{Ker}(g)\cap W\neq\{0\}$. 
Let $x$ be a non-zero vector in $\mathrm{Ker}(g)\cap W$. One has
\[\|(f-g)(x)\|^{2}=\|f(x)\|^{2}=\langle f(x),f(x)\rangle=\langle f^*(f(x)),x\rangle.\]
Since $x\in W$, one obtains
$\|(f-g)(x)\|^{2}\geqslant \lambda_i\|x\|^{2}$.
Therefore $\|f-g\|^{2}\geqslant\lambda_i$. The proposition is thus proved.
\end{proof}

\begin{prop}\label{Pro:determinant HS}
Let $V$ be a finite-dimensional vector space over $k$, equipped with a  semidefinite inner product $\emptyinnprod$, $r$ be the rank of $V$, and $\norm{\ndot}_{\det'}$ be the quotient seminorm of the orthogonal tensor product seminorm on $V^{\otimes r}$ by the canonical quotient map $V^{\otimes r}\rightarrow\det(V)$. Then one has $\norm{\ndot}_{\det}=(r!)^{1/2}\norm{\ndot}_{\det'}$.
\end{prop}
\begin{proof}If the seminorm associated with the  semidefinite inner product on $V$ is not a norm, then both seminorms $\norm{\ndot}_{\det}$ and $\norm{\ndot}_{\det'}$ vanish. It then suffices to treat the case where $\emptyinnprod$ is an inner product.

Let $\varphi$ be an element in $V^{\otimes r}$. Let $\{e_i\}_{i=1}^r$ be an orthonormal basis of $V$. We write $\varphi$ into the form
\[\varphi=\sum_{I=(i_1,\ldots,i_r)\in\{1,\ldots,r\}^r}a_I(e_{i_1}\otimes\cdots e_{i_r}).\]
Then the canonical image $\eta$ of $\varphi$ in $\det(V)$ is
\[\Big(\sum_{\sigma\in\mathfrak S_r}\mathrm{sgn}(\sigma)a_{(\sigma(1),\ldots,\sigma(r))}\Big) e_{1}\wedge\cdots\wedge e_{r},\]
where $\mathfrak S_r$ is the symmetric group of order $r$. Hence the Cauchy-Schwarz inequality leads to 
\[\|\eta\|_{\det}=\Big|\sum_{\sigma\in\mathfrak S_r}\mathrm{sgn}(\sigma)a_{(\sigma(1),\ldots,\sigma(r))}\Big|\leqslant (r!)^{1/2}\|\varphi\|_{\mathrm{HS}},\]
where $\norm{\ndot}_{\mathrm{HS}}$ denotes the orthogonal tensor product norm on $V^{\otimes r}$. The equality is attained when $\varphi$ is of the form
$\sum_{\sigma\in\mathfrak S_r}\sgn(\sigma)e_{\sigma(1)}\otimes\cdots\otimes e_{\sigma(r)}$.
The proposition is thus proved.
\end{proof}

\begin{prop}\label{Pro: tensor product and deteminant HS}
Let $V$ and $W$ be finite-dimensional seminormed vector spaces. We assume that the seminorms of $V$ and $W$ are induced by  semidefinite inner products. Let $n$ and $m$ be respectively the ranks of $V$ and $W$ over $k$. We equip the tensor product $V\otimes_kW$ with the orthogonal tensor product seminorm $\norm{\ndot}_{\mathrm{HS}}$. Then the canonical isomorphism $\det(V\otimes_kW)\rightarrow\det(V)^{\otimes m}\otimes\det(W)^{\otimes n}$ is an isometry, where we consider the determinant of the Hilbert-Schmidt seminorm  on $\det(V\otimes_kW)$, and the tensor product $\norm{\ndot}'$ of determinant seminorms on $\det(V)^{\otimes m}\otimes\det(W)^{\otimes n}$.
\end{prop}
\begin{proof}{The assertion is trivial when at least one of the seminorms of $V$ and $W$ is not a norm since in this case both seminorms $\norm{\ndot}_{\mathrm{HS},\det}$ and $\norm{\ndot}'$ vanish.}

In the following, we assume that $V$ and $W$ are equipped with inner products. Let $\{e_i\}_{i=1}^n$ and $\{f_j\}_{j=1}^m$ be respectively orthonormal bases of $V$ and $W$, which are also Hadamard bases (by Proposition \ref{Pro:ArchimedeanHadamard}). Then $\{e_i\otimes f_j\}_{(i,j)\in\{1,\ldots,n\}\times\{1,\ldots,m\}}$ is an orthonormal basis of $V\otimes_kW$. By Proposition \ref{Pro:ArchimedeanHadamard}, it is also an Hadamard basis. Hence one has
\[\bigg\|\bigwedge_{i=1}^n\bigwedge_{j=1}^m (e_i\otimes f_j)\bigg\|_{\mathrm{HS},\det}=1=\bigg\|(e_1\wedge\cdots\wedge e_n)^{\otimes m}\otimes(f_1\wedge\cdots\wedge f_m)^{\otimes n}\bigg\|'.\]
The proposition is thus proved.
\end{proof}

\section{Extension of scalars}\label{Subsec:extensionofscalars}
In this section, we suppose given a field extension $K$ of $k$ equipped with a complete absolute value which extends $|\ndot|$ on $k$. By abuse of notation, we still use the notation $|\ndot|$ to denote the extended absolute value on $K$. We can thus consider $K$ as a normed vector space over $k$, which is ultrametric if and only if the absolute value $|\ndot|$ on $k$ is non-Archimedean. 

Let $(V,\norm{\ndot})$ be a finite-dimensional seminormed vector space over $k$. We consider the natural $K$-linear map from $V_K=V\otimes_kK$ to $\mathscr L(V^*,K)$ which sends $x\otimes a\in V\otimes_kK$ (with $x\in V$ and $a\in K$) to the $k$-linear map $(f\in V^*)\mapsto a f(x)$. We equip $V^*$ with the dual norm and $\mathscr L(V^*,K)$ with the operator norm, which induces by this natural $K$-linear map a seminorm on $V\otimes_kK$ denoted by $\norm{\ndot}_{K,\varepsilon}$ and called the seminorm \emph{induced by $\norm{\ndot}$ by $\varepsilon$-extension of scalars}\index{epsilon-extension of scalars@$\varepsilon$-extension of scalars}\index{seminorm!epsilon-extension of scalars@$\varepsilon$-extension of scalars}. Note that the seminorm $\norm{\ndot}_{K,\varepsilon}$ is necessarily ultrametric if $k$ is non-Archimedean. 
Moreover, if $(K,|\ndot|)$ is reflexive as normed vector space over $k$ (this condition is satisfied notably when $K/k$ is a finite extension), then the seminorm $\norm{\ndot}_{K,\varepsilon}$ is the $\varepsilon$-tensor product of $\norm{\ndot}$ and the absolute value on $K$ (viewed as a norm on the $k$-vector space $K$), see Remark \ref{Rem:operateureps}.

We denote by $\norm{\ndot}_{K,\pi}$ the $\pi$-tensor product seminorm on $V\otimes_kK$ of the seminorm $\norm{\ndot}$ on $V$ and the absolute value $|\ndot|$ on $K$, called the seminorm \emph{induced by $\norm{\ndot}$ by $\pi$-extension of scalars}\index{pi-extension of scalars@$\pi$-extension of scalars}\index{seminorm!pi-extension of scalars@$\pi$-extension of scalars}. If $|\ndot|$ is Archimedean and if the seminorm $\norm{\ndot}$ is induced by a  semidefinite inner product, we denote by $\norm{\ndot}_{K,\mathrm{HS}}$ the orthogonal tensor product of the seminorm $\norm{\ndot}$ on $V$ and the absolute value $|\ndot|$ on $K$ (in the Archimedean case the extension $K/k$ is always finite), called the seminorm \emph{induced by $\norm{\ndot}$ by orthogonal extension of scalars}\index{orthogonal extension of scalars}\index{seminorm!orthogonal extension of scalars}.

In what follows, an element $x\in V$ is often considered as an element of $V_K=V\otimes_kK$ by the inclusion map $V\rightarrow V\otimes_kK$ sending $x$ to $x\otimes 1$.

\subsection{Basic properties} In this subsection, we discuss some basic behaviour of norms induced by extension of scalars.

\begin{prop}\label{Pro:extensiondecorps}
Let $(V,\norm{\ndot})$ be a finite-dimensional seminormed vector space  over $k$. 
\begin{enumerate}[label=\rm(\arabic*)]
\item\label{item: extension eps} For any $x\in V$ one has $\|x\|_{K,\varepsilon}=\|x\|_{**}$, where $\norm{\ndot}_{**}$ denotes the double dual seminorm of $\norm{\ndot}$. In particular, if either $(k,|\ndot|)$ is Archimedean or 
$(V,\norm{\ndot})$ is ultrametric, then one has $\|x\|_{K,\varepsilon}=\|x\|$ for any $x\in V$.
\item\label{item: extension pi} For any $x\in V$ one has $\|x\|_{K,\pi}=\|x\|$. If $|\ndot|$ is Archimedean and  $\norm{\ndot}$ is induced by a  semidefinite inner product, for any $x\in V$ one has $\|x\|_{K,\mathrm{HS}}=\|x\|$.
\item\label{item: comparison} For any $y\in V_K$ one has $\|y\|_{K,\varepsilon}\leqslant\|y\|_{K,\pi}$. If $(k,|\ndot|)$ is $\mathbb R$ equipped with the usual absolute value, $K=\mathbb C$, and  $\norm{\ndot}$ is induced by a  semidefinite inner product, then for any $y\in V_{\mathbb C}$ one has \begin{gather}\label{Equ:comparison HS}\|y\|_{\mathbb C,\varepsilon}\leqslant\|y\|_{\mathbb C,\mathrm{HS}}\leqslant\|y\|_{\mathbb C,\pi},\\ 
\label{Equ:comparison HS2}\min\{\rang_{\mathbb R}(V^*),2\}^{-1/2}\|y\|_{\mathbb C,\pi}\leqslant\|y\|_{\mathbb C,\mathrm{HS}}\leqslant \min\{\rang_{\mathbb R}(V^*),2\}^{1/2}\|y\|_{\mathbb C,\varepsilon}.\end{gather}
\end{enumerate}
\end{prop}
\begin{proof}
\ref{item: extension eps} Let $\ell_x:V^*\rightarrow k$ be the $k$-linear map sending any bounded linear form $f\in V^*$ to $f(x)$. Let $\widetilde{\ell}_x:V^*\rightarrow K$ be the composition $\ell_x$ with the inclusion map $k\rightarrow K$. The operator norms of $\ell_x$ and $\widetilde{\ell}_x$ are the same.  Therefore one has $\|x\|_K=\|x\|_{**}$. The last assertion comes from Proposition \ref{Pro:doubledualarch} and Corollary \ref{Cor:doubledual}.

The first assertion of \ref{item: extension pi} follow from Remark~\ref{Rem:produittenrk1} in the Archimedean case and from Proposition \ref{Pro: pi tensorial preserve split tensor} in the non-Archimedean case (note that the absolute value on $K$, viewed as a norm when we consider $K$ as a vector space over $k$, is ultrametric once $|\ndot|$ is non-Archimedean).
The second assertion follows from  \eqref{Equ:HStensorprod} in \S\ref{Subsec:Hilbert-Schmidt}.

\ref{item: comparison} The first assertion follows from (1) and Proposition \ref{Pro:normpiisanorm}.

In the case where $(k,|\ndot|)$ is $\mathbb R$ equipped with the usual absolute value, $K=\mathbb C$, and  $\norm{\ndot}$ is induced by an inner product,
the inequalities follow from Proposition \ref{Pro:comparisonofHSandepsnorms} and Corollary \ref{Cor:comparisonofHSandpinorms}. 
\end{proof}

\begin{rema}\phantomsection\label{Rem:extensiondoubledual}\quad
\begin{enumerate}[label=\rm(\arabic*)]
\item\label{Item: same epsilon extension of scalars} Note that $\norm{\ndot}$ and its double dual seminorm $\norm{\ndot}_{**}$ induce the same dual norm on $V^*$ (see Proposition \ref{Pro:doubedualandquotient}). Hence they induce the same seminorm on $V_K$ by $\varepsilon$-extension of scalars. Moreover, if $K=k$, then $\norm{\ndot}_{K,\varepsilon}$ identifies with the double dual seminorm of $\norm{\ndot}$ on $V$.
\item\label{Item:extensiondoubledual 2} 
Assume that $k=\mathbb R$, $K=\mathbb C$ and $|\ndot|$ is the usual absolute value on $\mathbb R$. Suppose that the norm $\norm{\ndot}$ is induced by a  semidefinite inner product $\emptyinnprod$. Note that  $\emptyinnprod$ 
induces a  semidefinite inner product $\emptyinnprod_{\mathbb C}$,   
given by
\[\qquad \forall\,x,y,x',y'\in V,\quad
\langle x+iy,x'+iy'\rangle_{\mathbb C}=\langle x,x'\rangle+\langle y,y'\rangle+i\big(\langle x,y'\rangle-\langle y,x'\rangle\big).\]
Note that the seminorm corresponding to $\emptyinnprod_{\mathbb C}$
identifies with the orthogonal tensor product $\norm{\ndot}_{\mathbb C,\mathrm{HS}}$ of $\norm{\ndot}$ and $|\ndot|$. 
Moreover, an orthogonal basis of $(V,\emptyinnprod)$ 
remains to be an orthogonal basis of $(V_{\mathbb C},\emptyinnprod)_{\mathbb C}$, which implies that $\norm{\ndot}_{\mathbb C,\mathrm{HS}}$ is the unique seminorm on $V_{\mathbb C}$ extending $\norm{\ndot}$ which is induced by a  semidefinite inner product. 
In particular, one has $\emptyinnprod_{*,\mathbb C} =  \emptyinnprod_{\mathbb C,*}$,  
where $\emptyinnprod_{*}$ 
denotes the dual inner product of $\emptyinnprod$ 
(see \S\ref{Subsec:innerproduct}), and hence $\norm{\ndot}_{\mathbb C,\mathrm{HS},*}=\norm{\ndot}_{*,\mathbb C,\mathrm{HS}}$.
\item\label{item: dimension 1 ext} Let $V$ be a seminormed vector space of rank $1$ on $k$. Then the norms $\norm{\ndot}_{K,\varepsilon}$ and $\norm{\ndot}_{K,\pi}$ are the same since they take the same value on a non-zero vector of $V$ (by Proposition \ref{Pro:extensiondecorps}). Similarly, if $|\ndot|$ is Archimedean then one has $\norm{\ndot}_{K,\varepsilon}=\norm{\ndot}_{K,\mathrm{HS}}=\norm{\ndot}_{K,\pi}$. We just call this seminorm the seminorm \emph{induced by $\norm{\ndot}$ by extension of scalars}\index{extension of scalars}\index{seminorm!extension of scalars} and denote it by $\norm{\ndot}_K$.
\end{enumerate}
\end{rema}

\begin{prop}\label{Pro: null space in the extension}
Let $(V,\norm{\ndot})$ be a finite-dimensional seminormed vector space over $k$ and $N=N_{\norm{\ndot}}$ be the null space of $\norm{\ndot}$. 
\begin{enumerate}[label=\rm(\arabic*)]
\item\label{Item: null space of extension} The null spaces of the seminorms $\norm{\ndot}_{K,\varepsilon}$ and $\norm{\ndot}_{K,\pi}$ are both equal to $N_K$.
\item\label{Item: criterion linear form bounded} A linear form on the $K$-vector space $V_K$ is bounded with respect to the seminorm $\norm{\ndot}_{K,\varepsilon}$ if and only if it is bounded with respect to $\norm{\ndot}_{K,\pi}$. Moreover, the {underlying} vector spaces of $(V_K,\norm{\ndot}_{K,\varepsilon})^*$ and $(V_K,\norm{\ndot}_{K,\pi})^*$ are both canonically isomorphic to $(V_K/N_K)^\vee$.
\item\label{Item: extension of scalars associated norm} The quotient norm on $V_K/N_K$ induced by $\norm{\ndot}_{K,\varepsilon}$ (resp. $\norm{\ndot}_{K,\pi}$) identifies with the $\varepsilon$-extension of scalars $\norm{\ndot}^\sim_{K,\varepsilon}$ (resp. the $\pi$-extension of scalars $\norm{\ndot}^\sim_{K,\pi}$) of the norm $\norm{\ndot}^\sim$.
\end{enumerate}
\end{prop}
\begin{proof}
\ref{Item: null space of extension} 
Note that the relation $\norm{\ndot}_{K,\varepsilon}\leqslant\norm{\ndot}_{K,\pi}$ {holds} (see Proposition \ref{Pro:extensiondecorps} \ref{item: comparison}), so that it is sufficient to see that
(i) $\norm{x}_{K,\pi} = 0$ for $x \in N_K$ and (ii) $\norm{x}_{K,\varepsilon}>0$ for $x\in V_K\setminus N_K$.
Let $\{e_i\}_{i=1}^n$ be a basis of $V$ such that $\{e_i\}_{i=1}^r$ forms a basis of $N$.

(i)  We write $x$ in the form $x=\lambda_1e_1+\cdots+\lambda_re_r$ with $(\lambda_1,\ldots,\lambda_r)\in K^r$. One has
\[0 \leqslant \norm{x}_{K,\pi}\leqslant\sum_{i=1}^r|\lambda_i|\cdot\norm{e_i}_{K,\pi}=\sum_{i=1}^r|\lambda_i|\cdot\norm{e_i}=0, \]
where the first equality comes from Proposition \ref{Pro:extensiondecorps} \ref{item: extension pi}.

(ii) We set $x=\lambda_1e_1+\cdots+\lambda_ne_n$ ($\lambda_1,\ldots,\lambda_n \in K$). If $x$ does not belong to $N_K$, then there exists $j\in\{r+1,\ldots,n\}$ such that $\lambda_j\neq 0$. Note that $e_j^\vee$ belongs to $V^*$. Hence
\[\norm{x}_{K,\varepsilon}\geqslant \frac{|\lambda_j|}{\norm{e_j^\vee}_*}>0.\]

\ref{Item: criterion linear form bounded} By Corollary \ref{Coro:equivalenceofnrom} \ref{Item: criterion of continuous linear form}, a linear form on a finite-dimensional seminormed vector space is bounded if and only if it vanishes on the null space of the seminorm. By \ref{Item: null space of extension} we obtain that both seminorms $\norm{\ndot}_{K,\varepsilon}$ and $\norm{\ndot}_{K,\pi}$ admit $N_K$ as the null space. Hence we obtain the required result.

\ref{Item: extension of scalars associated norm} We identify $V^*$ with $(V/N)^\vee$ and then the norm $\norm{\ndot}_*$ identifies with the dual norm of $\norm{\ndot}^\sim$. Therefore by definition for any $x\in V_K$ one has $\norm{x}_{K,\varepsilon}=\norm{[x]}_{K,\varepsilon}^\sim$, where $[x]$ denotes the class of $x$ in $V_K/N_K$. The case of $\pi$-extension of scalars comes from Proposition \ref{Pro:quotientavecpitensor}.  
\end{proof}

The following proposition proves a universal property of the $\pi$-extension of scalars.

\begin{prop}\label{Pro: maximality property pi}
Let $(V,\norm{\ndot})$ be a finite-dimensional seminormed vector space over $k$. If 
$\norm{\ndot}'$ is a seminorm on $V_K$ whose restriction {to} $V$ is bounded from above by $\norm{\ndot}$, then the seminorm $\norm{\ndot}'$ is bounded from above by $\norm{\ndot}_{K,\pi}$. In particular, $\norm{\ndot}_{K,\pi}$ is the largest seminorm on $V_K=V\otimes_kK$ extending $\norm{\ndot}$.
\end{prop}
\begin{proof}
For any $x\in V$ and $a\in K$ one has \[\norm{x\otimes a}'= |a|\cdot\norm{x\otimes 1}'\leqslant |a|\cdot\norm{x}.\] By Proposition \ref{Pro:normpiisanorm}, we obtain $\norm{\ndot}'\leqslant\norm{\ndot}_{K,\pi}$.
\end{proof}

\begin{prop}
\label{Pro: produit pi tensoriel extension}
Let $(V_1,\norm{\ndot}_1)$ and $(V_2,\norm{\ndot}_2)$ be finite-dimensional seminormed vector spaces over $k$, and $\norm{\ndot}$ be the $\pi$-tensor product {seminorm} of $\norm{\ndot}_1$ and $\norm{\ndot}_2$. Then the {norm} $\norm{\ndot}_{K,\pi}$ identifies with the $\pi$-tensor product of $\norm{\ndot}_{1,K,\pi}$ and $\norm{\ndot}_{2,K,\pi}$. 
\end{prop}
\begin{proof}
Let $\norm{\ndot}'$ be the $\pi$-tensor product of $\norm{\ndot}_{1,K,\pi}$ and $\norm{\ndot}_{2,K,\pi}$. If $s$ is an element of $V_1\otimes_k V_2$, which is written as $s=x_1\otimes y_1+\cdots+x_n\otimes y_n$, with $(x_1,\ldots,x_n)\in V_1^n$ and $(y_1,\ldots,y_2)\in V_2^n$. Then one has
\[\|s\|'\leqslant\sum_{i=1}^n\|x_i\|_{1,K,\pi}\cdot\|y_i\|_{2,K,\pi}=\sum_{i=1}^n\|x_i\|_{1}\cdot\|y_i\|_2.\]
Therefore one has $\|s\|'\leqslant\|s\|$. By Proposition \ref{Pro: maximality property pi}, the norm $\norm{\ndot}'$ is bounded from above by $\norm{\ndot}_{K,\pi}$. 

To prove the converse inequality, by Proposition \ref{Pro:normpiisanorm}, it suffices to show that, for any split tensor $u\otimes v$ in $V_{1,K}\otimes_KV_{2,K}$ one has $\|u\otimes v\|_{K,\pi}\leqslant\|u\|_{1,K,\pi}\cdot\| v\|_{2,K,\pi}$. Assume that $u$ and $v$ are written as $u=\lambda_1x_1+\cdots+\lambda_nx_n$ and $v=\mu_1 y_1+\cdots+\mu_my_m$ with $(\lambda_1,\ldots,\lambda_n)\in K^n$, $(x_1,\ldots,x_n)\in V_1^n$, $(\mu_1,\ldots,\mu_m)\in K^m$ and $(y_1,\ldots,y_m)\in V_2^m$. Then one has 
\[\begin{split}\|u\otimes v\|_{K,\pi}&\leqslant\sum_{i=1}^n\sum_{j=1}^m|\lambda_i\mu_j|\cdot\|x_i\otimes y_j\|_{K,\pi}=\sum_{i=1}^n\sum_{j=1}^m|\lambda_i\mu_j|\cdot\|x_i\otimes y_j\|\\
&=\sum_{i=1}^n\sum_{j=1}^m|\lambda_i\mu_j|\cdot\|x_i\|_{1}\cdot\| y_j\|_{2}=\Big(\sum_{i=1}^n|\lambda_i|\cdot\|x_i\|_1\Big)\Big(\sum_{j=1}^m|\mu_j|\cdot\|y_j\|_2\Big).
\end{split}\]
Since the decompositions $u=\lambda_1x_1+\cdots+\lambda_nx_n$ and $v=\mu_1 y_1+\cdots+\mu_my_m$ are arbitrary, we obtain
\[\|u\otimes v\|_{K,\pi}\leqslant\|u\|_{1,K,\pi}\cdot\| v\|_{2,K,\pi}.\]
\end{proof}

\begin{prop}
\label{Pro: produit HS tensoriel extension}
Assume that $(k,|\ndot|)$ is the field $\mathbb R$ equipped with the usual absolute value and $K=\mathbb C$.
Let $(V_1,\emptyinnprod_1)$ and $(V_2,\emptyinnprod_2)$ be finite-dimensional vector spaces over $\mathbb R$ equipped with  semidefinite inner products, $\norm{\ndot}_1$ and $\norm{\ndot}_2$ be seminorms corresponding to $\emptyinnprod_1$ and $\emptyinnprod_2$, respectively,  and $\norm{\ndot}$ be the orthogonal tensor product of $\norm{\ndot}_1$ and $\norm{\ndot}_2$. Then the seminorm $\norm{\ndot}_{\mathbb C,\mathrm{HS}}$ identifies with the orthogonal tensor product of $\norm{\ndot}_{1,\mathbb C,\mathrm{HS}}$ and $\norm{\ndot}_{2,\mathbb C,\mathrm{HS}}$. 
\end{prop}
\begin{proof}
Let $\{x_i\}_{i=1}^n$ and $\{y_j\}_{j=1}^m$ be orthogonal bases of $(V_1,\norm{\ndot}_1)$ and $(V_2,\norm{\ndot}_2)$, respectively. Then $\{x_i\otimes y_j\}_{(i,j)\in\{1,\ldots,n\}\times\{1,\ldots,m\}}$ is an orthogonal basis of $(V_1\otimes_{\mathbb R}V_2,\norm{\ndot})$ and hence is an orthogonal basis of $(V_{1,\mathbb C}\otimes_{\mathbb C}V_{2,\mathbb C},\norm{\ndot}_{\mathbb C,\mathrm{HS}})$. Moreover, $\{x_i\}_{i=1}^n$ and $\{y_j\}_{j=1}^m$ are also orthogonal bases of  $(V_{1,\mathbb C},\norm{\ndot}_{1,\mathbb C,\mathrm{HS}})$ and $(V_{2,\mathbb C},\norm{\ndot}_{2,\mathbb C,\mathrm{HS}})$, respectively. Hence $\{x_i\otimes y_j\}_{(i,j)\in\{1,\ldots,n\}\times\{1,\ldots,m\}}$ is an orthogonal basis of $V_{1,\mathbb C}\otimes_{\mathbb C}V_{2,\mathbb C}$ with respect to the orthogonal tensor product of $\norm{\ndot}_{1,\mathbb C,\mathrm{HS}}$ and $\norm{\ndot}_{2,\mathbb C,\mathrm{HS}}$. The proposition is thus proved.
\end{proof}

\begin{prop}\label{Pro:distanceextension}
Let $V$ be a finite-dimensional vector space over $k$, and $\norm{\ndot}_1$ and $\norm{\ndot}_2$ be two norms on $V$. 
\begin{enumerate}[label=\rm(\arabic*)]

\item\label{Item: distance K} One has 
\[d(\norm{\ndot}_{1,K,\varepsilon},\norm{\ndot}_{2,K,\varepsilon})=d(\norm{\ndot}_{1,**},\norm{\ndot}_{2,**})\leqslant d(\norm{\ndot}_1,\norm{\ndot}_{2}).\]
In particular, if both norms $\norm{\ndot}_1$ and $\norm{\ndot}_2$ are reflexive, then \[d(\norm{\ndot}_{1,K,\varepsilon},\norm{\ndot}_{2,K,\varepsilon})=d(\norm{\ndot}_1,\norm{\ndot}_{2}).\]
\item\label{Item: distance Kpi} One has $d(\norm{\ndot}_{1,K,\pi},\norm{\ndot}_{2,K,\pi})= d(\norm{\ndot}_1,\norm{\ndot}_{2})$.
\item\label{Item: distance KpiHS} Assume that $|\ndot|$ is Archimedean and that $\norm{\ndot}_1$ and $\norm{\ndot}_2$ are induced by inner products. Then $d(\norm{\ndot}_{1,K,\mathrm{HS}},\norm{\ndot}_{2,K,\mathrm{HS}})= d(\norm{\ndot}_1,\norm{\ndot}_{2})$.
\end{enumerate}
\end{prop}
\begin{proof}
\ref{Item: distance K} By Proposition \ref{Pro:distanceofoperatornorms}, one has $d(\norm{\ndot}_{1,*},\norm{\ndot}_{2,*})\leqslant d(\norm{\ndot}_1,\norm{\ndot}_2)$. {By the same argument as that of the proof of Proposition \ref{Pro:distanceofoperatornorms}, we can show that}  $d(\norm{\ndot}_{1,K,\varepsilon},\norm{\ndot}_{2,K,\varepsilon})\leqslant d(\norm{\ndot}_{1,*},\norm{\ndot}_{2,*})$. Hence we obtain the inequality $d(\norm{\ndot}_{1,K,\varepsilon},\norm{\ndot}_{2,K,\varepsilon})\leqslant 
d(\norm{\ndot}_{1},\norm{\ndot}_{2})$. {By Proposition \ref{Pro:doubedualandquotient}, for $i\in\{1,2\}$, $\norm{\ndot}_i$ and $\norm{\ndot}_{i,**}$ induce the same dual norm on $V^*$, and hence $\norm{\ndot}_{i,K,\varepsilon}=\norm{\ndot}_{i,**,K,\varepsilon}$. Therefore the above argument actually leads to $d(\norm{\ndot}_{1,K,\varepsilon},\norm{\ndot}_{2,K,\varepsilon})\leqslant d(\norm{\ndot}_{1,**},\norm{\ndot}_{2,**})$. Conversely, by Proposition \ref{Pro:extensiondecorps} \ref{item: extension eps} we obtain that $\norm{\ndot}_{1,K,\varepsilon}$ and $\norm{\ndot}_{2,K,\varepsilon}$ extend $\norm{\ndot}_{1,**}$ and $\norm{\ndot}_{2,**}$, respectively. Hence one has $d(\norm{\ndot}_{1,K,\varepsilon},\norm{\ndot}_{2,K,\varepsilon})\geqslant d(\norm{\ndot}_{1,**},\norm{\ndot}_{2,**})$. 

The inequality $d(\norm{\ndot}_{1,**},\norm{\ndot}_{2,**})\leqslant d(\norm{\ndot}_1,\norm{\ndot}_{2})$ comes from Proposition \ref{Pro:distanceofoperatornorms}. The equality holds when both norms $\norm{\ndot}_1$ and $\norm{\ndot}_2$ are reflexive.}

\ref{Item: distance Kpi} By Proposition \ref{Pro:extensiondecorps} \ref{item: extension pi}, $\norm{\ndot}_{1,K,\pi}$ and $\norm{\ndot}_{2,K,\pi}$ extend $\norm{\ndot}_1$ and $\norm{\ndot}_2$, respectively, and hence $d(\norm{\ndot}_{1,K,\pi},\norm{\ndot}_{2,K,\pi})\geqslant d(\norm{\ndot}_1,\norm{\ndot}_2)$. In the following, we prove the converse inequality. 
If we set $\delta =  d(\norm{\ndot}_1,\norm{\ndot}_2)$, then $\mathrm{e}^{-\delta} \leqslant \|s\|_1/\|s\|_2 \leqslant \mathrm{e}^{\delta}$ for
$s \in V \setminus \{ 0 \}$, that is, $\|\ndot\|_1 \leqslant \mathrm{e}^{\delta} \|\ndot\|_2$ and $\|\ndot\|_2 \leqslant \mathrm{e}^{\delta} \|\ndot\|_1$.
By Proposition \ref{Pro: maximality property pi}, one has $\norm{\ndot}_{1,K,\pi}\leqslant \mathrm{e}^{\delta} \norm{\ndot}_{2,K,\pi}$.
By the same reason, $\|\ndot\|_{2,K,\pi}\leqslant \mathrm{e}^{\delta}\|\ndot\|_{1,K,\pi}$.
Hence the inequality $d(\norm{\ndot}_{1,K,\pi},\norm{\ndot}_{2,K,\pi})\leqslant \delta = d(\norm{\ndot}_1,\norm{\ndot}_2)$ holds.

\ref{Item: distance KpiHS} It suffices to treat the case where $k=\mathbb R$ and $K=\mathbb C$.  By Proposition \ref{Pro:extensiondecorps} \ref{item: extension pi}, $\norm{\ndot}_{1,K,\mathrm{HS}}$ and $\norm{\ndot}_{2,K,\mathrm{HS}}$ extend $\norm{\ndot}_1$ and $\norm{\ndot}_2$, respectively, and hence $d(\norm{\ndot}_{1,K,\mathrm{HS}},\norm{\ndot}_{2,K,\mathrm{HS}})\geqslant d(\norm{\ndot}_1,\norm{\ndot}_2)$. 
As in \ref{Item: distance Kpi}, if we set $\delta =  d(\norm{\ndot}_1,\norm{\ndot}_2)$, then $\|\ndot\|_1 \leqslant \mathrm{e}^{\delta} \|\ndot\|_2$ and $\|\ndot\|_2 \leqslant \mathrm{e}^{\delta} \|\ndot\|_1$.
Let $z$ be an element of $V_{\mathbb C}$, which is written as $z=x+iy$, where $x$ and $y$ are vectors in $V$. Then one has \[\|z\|_{1,\mathbb C,\mathrm{HS}}^2=\|x\|_{1}^2+\|y\|_{1}^2\leqslant \mathrm{e}^{2\delta}(\|x\|_{2}^2+\|y\|_2^2)=\mathrm{e}^{2\delta}\|z\|_{2,\mathbb C,\mathrm{HS}}^2.\]
Therefore $\norm{\ndot}_{1,\mathrm{HS},\mathbb C}\leqslant \mathrm{e}^{\delta}\norm{\ndot}_{2,\mathrm{HS},\mathbb C}$. Similarly, $\norm{\ndot}_{2,\mathbb C,\mathrm{HS}}\leqslant \mathrm{e}^{\delta} \norm{\ndot}_{1,\mathbb C,\mathrm{HS}}$, so that 
the inequality $d(\norm{\ndot}_{1,\mathbb C,\mathrm{HS}},\norm{\ndot}_{2,\mathbb C,\mathrm{HS}})\leqslant \delta=d(\norm{\ndot}_1,\norm{\ndot}_2)$ holds.
\end{proof}

\begin{prop}
\label{Pro:extensioncorpsmor}Let $(V_1, \norm{\ndot}_1)$ and $(V_2,\norm{\ndot}_2)$ be finite-dimensional seminormed vector {spaces} over $k$, and $f:V_1\rightarrow V_2$ be a bounded $k$-linear map. Let $f_K:V_{1,K}\rightarrow V_{2,K}$ be the $K$-linear map induced by $f$.
\begin{enumerate}[label=\rm(\arabic*)]
\item\label{Item: compare operator norms} 
If we consider the seminorms $\norm{\ndot}_{1,K,\varepsilon}$ and $\norm{\ndot}_{2,K,\varepsilon}$ on $V_{1,K}$ and $V_{2,K}$, respectively, then the operator seminorm of $f_K$ is bounded from above by that of $f^*$ (which is bounded from above by $\norm{f}$, see Proposition \ref{Pro:normdualopt}). The equality 
$\norm{f_K}=\norm{f}$  holds when $(V_2,\norm{\ndot}_2)$ is reflexive.
\item\label{Item: compare operator norms pi} If we consider the seminorms $\norm{\ndot}_{1,K,\pi}$ and $\norm{\ndot}_{2,K,\pi}$  on $V_{1,K}$ and $V_{2,K}$,  respectively, then the operator seminorms of $f_K$ and $f$ are the same.
\item\label{Item: compare operator norms HN} Assume that $(k,|\ndot|)$ is $\mathbb R$ equipped with the usual absolute value, $K=\mathbb C$ and that $\norm{\ndot}_1$ and $\norm{\ndot}_2$ are induced by  semidefinite inner products. If we consider the norms $\norm{\ndot}_{1,K,\mathrm{HS}}$ and $\norm{\ndot}_{2,K,\mathrm{HS}}$  on $V_{1,K}$ and $V_{2,K}$, respectively, then the operator seminorms of $f_K$ and $f$ are the same.
\end{enumerate}
\end{prop}
\begin{proof}
\ref{Item: compare operator norms}   Let $\varphi$ be an element of $V_{1,K}$, viewed as a $k$-linear map from $V_1^*$ to $K$. Then the element $f_K(\varphi)\in V_{2,K}$, viewed as a $k$-linear form from $V_2^*$ to $K$, sends $\beta\in V_2^*$ to $\varphi(f^*(\beta))\in K$. One has
\[|\varphi(f^*(\beta))|\leqslant\|\varphi\|_{1,K,\varepsilon}\cdot\|f^*\|\cdot\|\beta\|_{2,*}.\]
Therefore $\|f_K(\varphi)\|_{2,K,\varepsilon}\leqslant\|f^*\|\cdot\|\varphi\|_{1,K,\varepsilon}$. Since $\varphi$ is arbitrary, one has $\|f_K\|\leqslant\|f^*\|$. The {first assertion} is thus proved.

Assume that $(V_2,\norm{\ndot}_2)$ is reflexive. For any element $x\in V_1$ one has 
\[
\begin{cases}
\|x\|_{1,K}=\|x\|_{1,**}\leqslant\|x\|_1, \\ 
\|f_K(x)\|_{2,K,\varepsilon}=\|f(x)\|_{2,K,\varepsilon}=\|f(x)\|_{2,**}=\|f(x)\|_2
\end{cases}
\]
since $(V_2,\norm{\ndot}_2)$ is reflexive. Therefore one has \[\|f_K\|\geqslant\sup_{x\in V_1\setminus N_{\norm{\sndot}_1}}\frac{\|f_K(x)\|_{2,K}}{\|x\|_{1,K}}\geqslant\sup_{x\in V_1\setminus N_{\norm{\sndot}_1}}\frac{\|f(x)\|_2}{\|x\|_1}=\|f\|.\]

\ref{Item: compare operator norms pi} Since the norms $\norm{\ndot}_{1,K,\pi}$ and $\norm{\ndot}_{2,K,\pi}$ extend $\norm{\ndot}_1$ and $\norm{\ndot}_2$, respectively (see Proposition \ref{Pro:extensiondecorps}), the operator seminorm $\|f\|$ is bounded from above by $\|f_K\|$. It suffices to prove the converse inequality. Let $y$ be an element in $V_{1,K}$, which is written as $y=x_1\otimes a_1+\cdots+x_n\otimes a_n$, where $(x_1,\ldots,x_n)\in V_1^n$ and $(a_1,\ldots,a_n)\in K^n$. Then one has $f_K(y)=f(x_1)\otimes a_1+\cdots+f(x_n)\otimes a_n$. Hence
\[\|f_K(y)\|_{2,K,\pi}\leqslant\sum_{i=1}^n|a_i|\cdot\|f(x_i)\|_{2}\leqslant \|f\|\sum_{i=1}^n|a_i|\cdot\|x_i\|_1.\]
As the decomposition $y=x_1\otimes a_1+\cdots+x_n\otimes a_n$ is arbitrary, we obtain
\[\|f(y)\|_{2,K,\pi}\leqslant\|f\|\cdot\|y\|_{1,K,\pi}.\] 

\ref{Item: compare operator norms HN} Since the seminorms $\norm{\ndot}_{1,\mathbb C,\mathrm{HS}}$ and $\norm{\ndot}_{2,\mathbb C,\mathrm{HS}}$ extend $\norm{\ndot}_1$ and $\norm{\ndot}_2$, respectively (see Proposition \ref{Pro:extensiondecorps}), the operator seminorm $\|f\|$ is bounded from above by $\|f_{\mathbb C}\|$. Let $z$ be an element of $V_{1,\mathbb C}$, written as $u+iv$, where $u$ and $v$ are vectors in $V_1$. Then one has $f_{\mathbb C}(z)=f(u)+if(v)$. Therefore
\[\|f_{\mathbb C}(z)\|^2=\|f(u)\|_2^2+\|f(v)\|_2^2\leqslant\|f\|^2(\|u\|_1^2+\|u\|_2^2)=\|f\|^2\cdot\|z\|_{1,\mathbb C,\mathrm{FS}}^2.\]
Hence $\|f_{\mathbb C}\|^2=\|f\|^2$.
\end{proof}

\subsection{Direct sums} In this subsection, we discuss the behaviour of direct sums under scalar extension. We fix a continuous and convex function $\psi:[0,1]\rightarrow[0,1]$  such that $\max\{t,1-t\}\leqslant\psi(t)$ for any $t\in[0,1]$ (cf. \S\ref{Subsec:directsums}).

\begin{prop}\label{Pro:sommedirectext}
Let $(V,\norm{\ndot}_V)$ and $(W,\norm{\ndot}_{W})$ be finite-dimensional seminormed vector spaces over $k$. Let $\norm{\ndot}_{\psi}$ be the $\psi$-direct sum of $\norm{\ndot}_V$ and $\norm{\ndot}_W$. Then for $(f,g)\in V_K\oplus W_K$ one has \begin{equation}\label{Equ:comparenormext}\max\{\|f\|_{V,K,\varepsilon},\|g\|_{W,K,\varepsilon}\}\leqslant\|(f,g)\|_{\psi,K,\varepsilon}.\end{equation} The equality holds if either $(k,|\ndot|)$ is non-Archimedean or 
$\psi(t)=\max\{t,1-t\}$ for any $t\in[0,1]$. Moreover, for any $(f,g)\in V_K\oplus W_K$ one has
\begin{equation}\label{Equ:marginales}\|(f,0)\|_{\psi,K,\varepsilon}=\|f\|_{V,K,\varepsilon},\quad\|(0,g)\|_{\psi,K,\varepsilon}=\|g\|_{W,K,\varepsilon}.\end{equation} 
\end{prop}
\begin{proof}
By Proposition \ref{Pro:dualdirectsum}, the dual norm $\norm{\ndot}_{\psi,*}$ is a certain direct sum of $\norm{\ndot}_{V,*}$ and $\norm{\ndot}_{W,*}$. Hence one has 
\begin{equation}\|\alpha\|_{V,*}+\|\beta\|_{W,*}\geqslant\|(\alpha,\beta)\|_{\psi,*}\geqslant\max\{\|\alpha\|_{V,*},\|\beta\|_{V,*}\}.\end{equation} Therefore, for any $(f,g)\in V_K\oplus W_K$ one has
\[\|(f,g)\|_{\psi,K,\varepsilon}\geqslant\sup_{\begin{subarray}{c}
(\alpha,\beta)\in V^*\oplus W^*\\
(\alpha,\beta)\neq (0,0)
\end{subarray}}\frac{|f(\alpha)+g(\beta)|}{\|\alpha\|_{V,*}+\|\beta\|_{W,*}}=\max\{\|f\|_{V,K,\varepsilon},\|g\|_{V,K,\varepsilon}\}\]
which proves \eqref{Equ:comparenormext}. Moreover, for any $f\in V_K$ one has
\[\|(f,0)\|_{\psi,K,\varepsilon}\leqslant\sup_{\begin{subarray}{c}
(\alpha,\beta)\in V^*\oplus W^*\\
(\alpha,\beta)\neq (0,0)
\end{subarray}}\frac{|f(\alpha)|}{\max\{\|\alpha\|_{V,*},\|\beta\|_{W,*}\}}=\|f\|_{V,K,\varepsilon}.\]
Therefore, by \eqref{Equ:comparenormext} we obtain the equality $\|(f,0)\|_{\psi,K,\varepsilon}=\|f\|_{V,K}$. Similarly, for any $g\in W_K$ one has $\|(0,g)\|_{\psi,K,\varepsilon}=\|g\|_{W,K,\varepsilon}$.

Finally, we proceed with the proof of the equality part of \eqref{Equ:comparenormext}. If $(k,|\ndot|)$ is non-Archimedean, then the seminorm $\norm{\ndot}_{\psi,K,\varepsilon}$ is ultrametric and hence by \eqref{Equ:marginales} one has  \[\forall\,(f,g)\in V_K\oplus W_K,\quad\|(f,g)\|_{\psi,K,\varepsilon}\leqslant\max\{\|f\|_{V,K,\varepsilon},\|g\|_{W,K,\varepsilon}\},\]
which leads to (by \eqref{Equ:comparenormext}) the equality
\[\forall\,(f,g)\in V_K\oplus W_K,\quad\|(f,g)\|_{\psi,K,\varepsilon}=\max\{\|f\|_{V,K,\varepsilon},\|g\|_{W,K,\varepsilon}\}.\] In the case where $k$ is Archimedean and $\psi(t)=\max\{t,1-t\}$ for any $t\in[0,1]$, one has $\|(\alpha,\beta)\|_{\psi,*}=\|\alpha\|_{V,*}+\|\beta\|_{W,*}$ for any $(\alpha,\beta)\in V^\vee\oplus W^\vee$. Therefore \[\|(f,g)\|_{\psi,K,\varepsilon}=\sup_{\begin{subarray}{c}
(\alpha,\beta)\in V^*\oplus W^*\\
(\alpha,\beta)\neq (0,0)
\end{subarray}}\frac{|f(\alpha)+g(\beta)|}{\|\alpha\|_{V,*}+\|\beta\|_{W,*}}=\max\{\|f\|_{V,K,\varepsilon},\|g\|_{V,K,\varepsilon}\}.\]
\end{proof}

\begin{rema}\label{Rem:comparisondirectsumext}
Let $\psi$ be an element in $\mathscr S$ (see \S\ref{Subsec:directsums}), which corresponds to an absolute normalised norm $\norm{\ndot}$ on $\mathbb R^2$. Let $\psi_*$ be the element in $\mathscr S$ corresponding to the dual norm $\norm{\ndot}_*$ (see Definition \ref{Def:dualdirecsum}). Suppose given finite-dimensional seminormed vector spaces $(V,\norm{\ndot}_V)$ and $(W,\norm{\ndot}_W)$ over $\mathbb R$ (equipped with the usual absolute value). By Proposition \ref{Pro:dualdirectsum} \ref{Item: dual direct sum2}, the dual norm of $\norm{\ndot}_\psi$ (the $\psi$-direct sum of $\norm{\ndot}_V$ and $\norm{\ndot}_W$) identifies with the $\psi_*$-direct sum of $\norm{\ndot}_{V,*}$ and $\norm{\ndot}_{W,*}$. Therefore, for any $(f,g)\in V_{\mathbb C}
\oplus W_{\mathbb C}$, one has
\begin{align*}
\|(f,g)\|_{\psi,\mathbb C,\varepsilon} & =\sup_{\begin{subarray}{c}(\alpha,\beta)\in V^*\oplus W^*\\
(\alpha,\beta)\neq(0,0)\end{subarray}}\frac{|f(\alpha)+g(\beta)|}{\|(\|\alpha\|_{V,*},\|\beta\|_{W,*})\|_*}\\
&\leqslant
\sup_{\begin{subarray}{c}(\alpha,\beta)\in V^*\oplus W^*\\
(\alpha,\beta)\neq(0,0)\end{subarray}}\frac{\|f\|_{V,\mathbb C,\varepsilon}\cdot\|\alpha\|_{V,*}+\|g\|_{W,\mathbb C,\varepsilon}\cdot\|\beta\|_{W,*}}{\|(\|\alpha\|_{V,*},\|\beta\|_{W,*})\|_*}\\
& =\|(\|f\|_{V,\mathbb C,\varepsilon},\|g\|_{W,\mathbb C,\varepsilon})\|.
\end{align*}
In other words, the {seminorm} $\norm{\ndot}_{\psi,\mathbb C,\varepsilon}$ is bounded from above by the $\psi$-direct sum of $\norm{\ndot}_{V,\mathbb C,\varepsilon}$ and $\norm{\ndot}_{W,\mathbb C,\varepsilon}$.
\end{rema}

\begin{prop}\label{Pro:sommedirectext pi}
Let $(V,\norm{\ndot}_V)$ and $(W,\norm{\ndot}_{W})$ be finite-dimensional seminormed vector spaces over $k$. Let $\norm{\ndot}_{\psi}$ be the $\psi$-direct sum of $\norm{\ndot}_V$ and $\norm{\ndot}_W$, and  $\norm{\ndot}_{K,\pi,\psi}$ be the $\psi$-direct sum of $\norm{\ndot}_{V,K,\pi}$ and $\norm{\ndot}_{W,K,\pi}$. Then  $\norm{\ndot}_{K,\pi,\psi}\leqslant\norm{\ndot}_{\psi,K,\pi}$. 
\end{prop}
\begin{proof}
Let $(x,y)$ be an element in  $V\oplus W$. One has
\[\begin{split}\|(x,y)\|_{K,\pi,\psi}&=(\|x\|_{V,K,\pi}+\|y\|_{W,K,\pi})\psi\Big(\frac{\|x\|_{V,K,\pi}}{\|x\|_{V,K,\pi}+\|y\|_{W,K,\pi}}\Big)\\
&=(\|x\|_V+\|y\|_W)\psi\Big(\frac{\|x\|_V}{\|x\|_V+\|y\|_W}\Big)=\|(x,y)\|_\psi,
\end{split}\]
where the second equality comes from Proposition \ref{Pro:extensiondecorps} \ref{item: extension pi}. Therefore the seminorm $\norm{\ndot}_{K,\pi,\psi}$ extends $\norm{\ndot}_{\psi}$. By Proposition \ref{Pro: maximality property pi}, it is bounded from above by $\norm{\ndot}_{\psi,K,\pi}$.
\end{proof}

\begin{prop}
Assume that $(k,|\ndot|)$ is the real field $\mathbb R$ equipped with the usual absolute value. Let $(V,\emptyinnprod_V)$ and $(W,\emptyinnprod_W)$ be finite-dimensional vector spaces over $\mathbb R$, equipped with  semidefinite inner products, $\norm{\ndot}_V$ and $\norm{\ndot}_W$ be seminorms associated with $\emptyinnprod_V$ and $\emptyinnprod_W$, respectively, and $\norm{\ndot}$ be the orthogonal direct sum of $\norm{\ndot}_{V}$ and $\norm{\ndot}_W$. Then $\norm{\ndot}_{\mathbb C,\mathrm{HS}}$ is the orthogonal direct sum of $\norm{\ndot}_{V,\mathbb C,\mathrm{HS}}$ and $\norm{\ndot}_{W,\mathbb C,\mathrm{HS}}$. 
\end{prop}
\begin{proof}
Let $\norm{\ndot}'$ be the orthogonal direct sum of $\norm{\ndot}_{V,\mathbb C,\mathrm{HS}}$ and $\norm{\ndot}_{W,\mathbb C,\mathrm{HS}}$. It is a seminorm on $V_{\mathbb C}\oplus W_{\mathbb C}$ which is induced by a  semidefinite inner product. Moreover, for any $(x,y)\in V\oplus W$ one has 
\[\|(x,y)\|'=(\|x\|_{V,\mathbb C,\mathrm{HS}}^2+\|y\|_{V,\mathbb C,\mathrm{HS}}^2)^{1/2}=(\|x\|_V^2+\|y\|_V^2)^{1/2}=\|(x,y)\|,\]
where the second equality comes from Proposition \ref{Pro:extensiondecorps} \ref{item: extension pi}. Therefore, $\norm{\ndot}'$ is a seminorm extending $\norm{\ndot}$ which is induced by a  semidefinite inner product and hence one has $\norm{\ndot}'=\norm{\ndot}_{\mathbb C,\mathrm{HS}}$ (see Remark \ref{Rem:extensiondoubledual}).
\end{proof}

\subsection{Orthogonality} In this subsection, we discuss the preservation of the orthogonality under extension of scalars, and its consequences. We have seen in Remark \ref{Rem:extensiondoubledual} \ref{Item:extensiondoubledual 2} that the orthonormality is preserved by the orthogonal extension of scalars.

\begin{prop}\label{Pro:alpha-orthgonalextension}
Let $(V,\norm{\ndot})$ be a finite-dimensional seminormed vector space over $k$, and $\alpha$ be a real number in $\intervalle]01]$. If $\boldsymbol{e}=\{e_i\}_{i=1}^r$ is an $\alpha$-orthogonal basis of $V$ with respect to the norm $\norm{\ndot}$, then it is also an $\alpha$-orthogonal basis of $V_K$ with respect to the norms $\norm{\ndot}_{K,\varepsilon}$ and $\norm{\ndot}_{K,\pi}$.
\end{prop}
\begin{proof}
Let $\boldsymbol{e}^\vee=\{e_i^\vee\}_{i=1}^r$ be the dual basis of $\boldsymbol{e}$. By Proposition \ref{Pro:alphaorthogonale}, the intersection $\boldsymbol{e}^\vee\cap V^*$ is an $\alpha$-orthogonal bases of $V^*$, and one has $\|e_i^\vee\|_*\leqslant\alpha^{-1}\|e_i\|^{-1}$ for any $e_i^\vee\in \boldsymbol{e}^\vee\cap V^*$ (see Lemma \ref{Lem:normofdualbasis}). If $x=a_1e_1+\cdots+a_re_r$ is an element in $V_K$, where $(a_1,\ldots,a_r)\in K^r$, and if $\ell_x:V^*\rightarrow K$ is the $k$-linear map sending $\varphi\in V^*$ to $a_1\varphi(e_1)+\cdots+a_r\varphi(e_r)$, then  for any $i\in\{1,\ldots,r\}$ such that $e_i^\vee\in V^*$ (or equivalently, $e_i\not\in N_{\norm{\ndot}}$) one has
\[\|\ell_x\|_{K,\varepsilon}\geqslant\frac{|\ell_x(e_i^\vee)|}{\|e_i^\vee\|_*}=\frac{|a_i|}{\|e_i^\vee\|_*}\geqslant\alpha |a_i|\cdot\|e_i\|\geqslant\alpha |a_i|\cdot\|e_i\|_{K,\varepsilon},\]
where the last inequality comes from Proposition \ref{Pro:extensiondecorps} and the relation \eqref{Equ:doubledual}. Therefore $\boldsymbol{e}$ is also an $\alpha$-orthogonal basis for $\norm{\ndot}_{K,\varepsilon}$. 

By Proposition \ref{Pro:extensiondecorps} \ref{item: comparison}, one has $\norm{\ndot}_{K,\varepsilon}\leqslant\norm{\ndot}_{K,\pi}$. Therefore \[\|x\|_{K,\pi}\geqslant\|x\|_{K,\varepsilon}\geqslant\alpha\max_{i\in\{1,\ldots,r\}}|a_i|\cdot\|e_i\|=\alpha\max_{i\in\{1,\ldots,r\}}|a_i|\cdot\|e_i\|_{K,\pi},\]
where the last equality comes from Proposition \ref{Pro:extensiondecorps} \ref{item: extension pi}.
\end{proof}

By using the preservation of orthogonality of bases, we prove an universal property of $\varepsilon$-extension of scalars, which is an  ultrametric analogue of Proposition \ref{Pro: maximality property pi}. 

\begin{prop}\label{Pro:extensioncorps}
Assume that the absolute value $|\ndot|$ is non-Archimedean. Let $V$ be a finite-dimensional vector space over $k$, equipped with a seminorm $\norm{\ndot}$. Let $\norm{\ndot}'_{K}$ be an ultrametric seminorm on $V_K$ whose restriction {to} $V$ is bounded from above by $\norm{\ndot}_{**}$. 
Then one has $\norm{\ndot}'_{K}\leqslant\norm{\ndot}_{K,\varepsilon}$. In particular, $\norm{\ndot}_{K,\varepsilon}$ is the largest ultrametric seminorm on $V_K$ which extends $\norm{\ndot}_{**}$.  
\end{prop}
\begin{proof}
By Proposition \ref{Pro:alpha-orthgonalextension}, if $\alpha$ is an element of $\intervalle]01[$ and if $\{e_i\}_{i=1}^r$ is an $\alpha$-orthogonal basis of $(V,\norm{\ndot})$, then $\{e_i\}_{i=1}^r$ is also an $\alpha$-orthogonal basis of $(V_K,\norm{\ndot}_{K,\varepsilon})$. 
In particular, for any $(\lambda_1,\ldots,\lambda_r)\in K^r$ one has
\[\|\lambda_1e_1+\cdots+\lambda_re_r\|'_{K}\leqslant\max_{i\in\{1,\ldots,r\}}|\lambda_i|\cdot\|e_i\|_{**}\leqslant \alpha^{-1}\|\lambda_1e_1+\cdots+\lambda_re_r\|_{K,\varepsilon}.\]
Since $(V,\norm{\ndot})$ admits an $\alpha$-orthogonal basis for any $\alpha\in\intervalle]01[$, we obtain $\norm{\ndot}'_{K}\leqslant\norm{\ndot}_{K,\varepsilon}$ for any ultrametric seminorm $\norm{\ndot}'_{K}$ with $\|\ndot\|'_K \leqslant \|\ndot\|_{**}$ on $V$.
\end{proof}

\begin{coro}\label{Cor:extsucc}
Let $K'$ be an extension of $K$ equipped with a complete absolute value extending that on $K$. Let $(V,\norm{\ndot})$ be a finite-dimensional seminormed vector space over $k$. One has $\norm{\ndot}_{K,\natural,K',\natural}=\norm{\ndot}_{K',\natural}$ on $V_{K'}$, where $\natural=\varepsilon$ or $\pi$.
\end{coro}
\begin{proof}The assertion is trivial when the absolute value $|\ndot|$ is Archimedean since in this case $k=\mathbb R$ or $\mathbb C$ and hence either $k=K$ or $K=K'$. In the following, we assume that the absolue value $|\ndot|$ is non-Archimedean.

By Proposition \ref{Pro:extensioncorps}, $\norm{\ndot}_{K',\varepsilon}$ is the largest ultrametric seminorm on $V_{K'}$ extending the seminorm $\norm{\ndot}_{**}$ on $V$. Moreover, by Proposition \ref{Pro:extensiondecorps}, $\norm{\ndot}_{K,\varepsilon}$ is an ultrametric seminorm on $V_K$ extending $\norm{\ndot}_{**}$, and the seminorm $\norm{\ndot}_{K,\varepsilon,K',\varepsilon}$ extends $\norm{\ndot}_{K,\varepsilon}$. Therefore one has $\norm{\ndot}_{K,\varepsilon,K',\varepsilon}\leqslant\norm{\ndot}_{K',\varepsilon}$. By the same reason, as the norm $\norm{\ndot}_{K',\varepsilon}$ extends $\norm{\ndot}_{**}$, its restriction {to} $V_{K}$ is bounded from above by $\norm{\ndot}_{K,\varepsilon}$ and hence the restriction of $\norm{\ndot}_{K',\varepsilon}$ {to} $V_K$ coincides with $\norm{\ndot}_{K,\varepsilon}$ (since we have already shown that $\norm{\ndot}_{K,\varepsilon,K',\varepsilon}\leqslant\norm{\ndot}_{K',\varepsilon}$). Therefore one has $\norm{\ndot}_{K,\varepsilon,K',\varepsilon}\geqslant\norm{\ndot}_{K',\varepsilon}$, still by the maximality property (for $\norm{\ndot}_{K,\varepsilon,K',\varepsilon}$) proved in Proposition \ref{Pro:extensioncorps}.

The case of $\pi$-extension of scalars is quite similar. By Proposition \ref{Pro:extensiondecorps} \ref{item: extension pi}, the seminorm $\norm{\ndot}_{K,\pi,K',\pi}$ extends $\norm{\ndot}_{K,\pi}$ on $V_K$ and hence extends $\norm{\ndot}$ on $V$. By the maximality property proved in Proposition \ref{Pro: maximality property pi}, we obtain that $\norm{\ndot}_{K,\pi,K',\pi}\leqslant\norm{\ndot}_{K',\pi}$. In particular, the restriction of $\norm{\ndot}_{K',\pi}$ {to} $V_{K}$ is bounded from below by $\norm{\ndot}_{K,\pi}$. Moreover, this restricted seminorm extends $\norm{\ndot}$. Still by the maximality property proved in Proposition \ref{Pro: maximality property pi}, we obtain that the restriction of $\norm{\ndot}_{K',\pi}$ {to} $V_K$ is bounded from above by $\norm{\ndot}_{K,\pi}$. Therefore the restriction of $\norm{\ndot}_{K',\pi}$ {to} $V_K$ coincides with $\norm{\ndot}_{K,\pi}$. By Proposition \ref{Pro: maximality property pi}, the norm $\norm{\ndot}_{K',\pi}$ is bounded from above by $\norm{\ndot}_{K,\pi,K',\pi}$. The proposition is thus proved. 
\end{proof}

{
\begin{prop}\label{Pro: Quotient extension}
Let $(V,\norm{\ndot}_V)$ be a finite-dimensional seminormed vector space over $k$, $Q$ be a quotient vector space of $V$, and $\norm{\ndot}_Q$ be the quotient seminorm of $\norm{\ndot}_V$ on $Q$.
\begin{enumerate}[label=\rm(\arabic*)]
\item \label{Item: extension quotient} The seminorm $\norm{\ndot}_{Q,K,\pi}$ identifies with the quotient of $\norm{\ndot}_{V,K,\pi}$ on $Q_K$.
\item \label{Item: epsilon extension quotient} 
The seminorm $\norm{\ndot}_{Q,K,\varepsilon}$ is bounded from above 
by the quotient seminorm of $\norm{\ndot}_{V,K,\varepsilon}$ on $Q_K$. The equality holds if one  
of the following conditions is satisfied: {\rm(i)} $|\ndot|$ is non-Archimedean; {\rm(ii)} $|\ndot|$ is Archimedean and $\norm{\ndot}_V$ is induced by a semidefinite inner product; {\rm(iii)} $Q$ is of dimension $1$ over $k$.
\item\label{Item: HS extension quotient} Assume that $|\ndot|$ is Archimedean and $\norm{\ndot}_V$ is induced by a semidefinite inner product. Then $\norm{\ndot}_{Q,K,\mathrm{HS}}$ identifies with the quotient of $\norm{\ndot}_{V,K,\mathrm{HS}}$ on $Q_K$. 
\end{enumerate}
\end{prop}
\begin{proof}
\ref{Item: extension quotient} follows directly from Proposition \ref{Pro:quotientavecpitensor}.

\ref{Item: epsilon extension quotient} Let $\norm{\ndot}_{Q,K}'$ be the quotient of the seminorm $\norm{\ndot}_{V,K,\varepsilon}$ on $Q_K$. 
Let $p : V \to Q$ be the canonical linear map. Note that $Q^{*} \subseteq V^{*}$ 
via $\psi \mapsto \psi \circ p$. 
Moreover, by Proposition~\ref{Pro:dualquotient}, {$\| \psi \circ p \|_{V,*} = \| \psi \|_{Q,*}$} for $\psi \in Q^*$.
Thus, for $s \in Q_K$, 
{%
\begin{align*}
\| s \|'_{Q,K} & = \inf_{\substack{x \in V_K\\ p_K(x) = s}} \sup_{\varphi \in V^{*} \setminus \{ 0 \}} \frac{|\varphi_K(x)|}{\| \varphi \|_{V,*}}
\geqslant \inf_{\substack{x \in V_K\\ p_K(x) = s }} \sup_{\psi \in Q^{*} \setminus \{ 0 \}} \frac{|\psi_K \circ p_K (x)|}{\| \psi \circ p \|_{V,*}} \\
& = \sup_{\psi \in Q^{*} \setminus \{ 0 \}} \frac{|\psi_K(s)|}{\| \psi \|_{Q,*}} = \| s \|_{Q,K,\varepsilon},
\end{align*}}
and hence the first assertion holds.

In the following, we prove the equality $\norm{\ndot}_{Q,K}'=\norm{\ndot}_{Q,K,\varepsilon}$ under each of the three conditions (i), (ii) and (iii). We first assume that the condition (i) or (ii) is satisfied. By Proposition \ref{Pro:dualquotient}, the dual norm $\norm{\ndot}_{Q,*}$ identifies with the restriction of $\norm{\ndot}_{V,*}$ {to} $Q^*$. By Proposition \ref{Pro:quotientdualnonarch}, we obtain that the seminorm $\norm{\ndot}_{Q,K,\varepsilon}$ identifies with the quotient seminorm of $\norm{\ndot}_{V,K,\varepsilon}$ on $Q_K$. 

Assume that the condition (iii) is satisfied and that the absolute value $|\ndot|$ is Archimedean (the non-Archimedean case has already been proved above). Let $f$ be a continuous $k$-linear operator from $Q^*$ to $K$. Since $Q$ is assumed to be of dimension $1$ over $k$, the image of $f$ is contained in a $k$-linear subspace of dimension $1$ in $K$. Therefore by Hahn-Banach theorem we obtain that there exists a continuous $k$-linear map $g:V^*\rightarrow K$ extending $f$ such that $f$ and $g$ have the same operator seminorm. Hence the seminorm $\norm{\ndot}_{Q,K,\varepsilon}$ identifies with the quotient seminorm of $\norm{\ndot}_{V,K,\varepsilon}$ on $Q_K$.

\ref{Item: HS extension quotient} follows directly from Proposition \ref{Pro:quotientr1Hs}.
\end{proof}}

\begin{prop}\label{Pro: restriction of epsion tensor}
Let $(V,\norm{\ndot}_V)$ be a finite-dimensional seminormed vector space over $k$ and $W$ be a vector subspaces of $V$. Let $\norm{\ndot}_W$ be the restriction of $\norm{\ndot}_V$ {to} $W$.
{
\begin{enumerate}[label=\rm(\arabic*)]
\item\label{Item: eps restriction extension}  The restriction of $\norm{\ndot}_{V,K,\varepsilon}$ {to} $W_K$ is bounded from above by $\norm{\ndot}_{W,K,\varepsilon}$. If {$|\ndot|$ is Archimedean or }  $\norm{\ndot}_V$ is ultrametric, then the restriction of $\norm{\ndot}_{V,K,\varepsilon}$ {to} $W_K$ coincides with $\norm{\ndot}_{W,K,\varepsilon}$.
\item\label{Item: pi restriction extension} The restriction of $\norm{\ndot}_{V,K,\pi}$ {to} $W_K$ is bounded from above by $\norm{\ndot}_{W,K,\pi}$. It identifies with $\norm{\ndot}_{W,K,\pi}$ if $\norm{\ndot}_V$ is ultrametric or induced by a  semidefinite inner product.
\item\label{Item: HS restriction extension} Assume that $|\ndot|$ is Archimedean and that $\norm{\ndot}_V$ is induced by {a semidefinite} inner product. Then the restriction of $\norm{\ndot}_{V,K,\mathrm{HS}}$ {to} $W_K$ identifies with $\norm{\ndot}_{W,K,\mathrm{HS}}$. 
\end{enumerate}}
\end{prop}
\begin{proof}
\ref{Item: eps restriction extension} Assume that $|\ndot|$ is non-Archimedean. By Proposition \ref{Pro:extensiondecorps} \ref{item: extension eps}, the seminorm $\norm{\ndot}_{V,K,\varepsilon}$ extends $\norm{\ndot}_{V,**}$. The restriction of $\norm{\ndot}_{V,K,\varepsilon}$ {to} $V$ is then bounded from above by $\norm{\ndot}_V$, which implies that the restriction of $\norm{\ndot}_{V,K,\varepsilon}$ {to} $W$ is bounded from above by $\norm{\ndot}_W$. Since $\norm{\ndot}_{W,**}$ is the largest ultrametric seminorm on $W$ which is bounded from above by $\norm{\ndot}_W$ (see Corollary \ref{Cor:doubledual}), we deduce from Proposition \ref{Pro:extensioncorps} that the restriction of $\norm{\ndot}_{V,K,\varepsilon}$ {to} $W$ is bounded from {above} by $\norm{\ndot}_{W,**}$. By Proposition \ref{Pro:extensioncorps}, we obtain that the restriction of $\norm{\ndot}_{V,K,\varepsilon}$ {to} $W_K$ is bounded from above by $\norm{\ndot}_{W,K,\varepsilon}$.

If $\norm{\ndot}_V$ is ultrametric or $|\ndot|$ is Archimedean, the dual norm $\norm{\ndot}_{W,*}$ coincides with the quotient norm of $\norm{\ndot}_{V,*}$ induced by the canonical quotient map $V^*\rightarrow W^*$ (see Proposition \ref{Pro:quotientdualnonarch} for the ultrametric case, and Remark \ref{Rem:comparaisondualitye} for the Archimedean case). Therefore, any $f\in W_K$, viewed as a $k$-linear operator from $W^*$ to $K$ or as a $k$-linear operator from $V^*$ to $K$, has the same operator norm. In other words, the restriction of $\norm{\ndot}_{V,K}$ {to} $W_K$ coincides with $\norm{\ndot}_{W,K}$.

{
\ref{Item: pi restriction extension} follows directly from Proposition \ref{Pro: restriction of epsion tensors} \ref{Item: restriction pi tensor}.

\ref{Item: HS restriction extension} follows directly from Proposition \ref{Pro: pi tensor sub Archimedean}.
}
\end{proof}

\begin{prop}\label{Pro: extension Hadamard}
Let $(V,\norm{\ndot})$ be a finite-dimensional normed vector space over $k$. We assume that either $|\ndot|$ is non-Archimedean or the norm $\norm{\ndot}$ is induced by an inner product. If $\{e_i\}_{i=1}^r$ is an Hadamard basis of $V$, then it is also an Hadamard basis of $V_K$ with respect to the norm $\norm{\ndot}_{K,\varepsilon}$.
\end{prop}
\begin{proof}
By Proposition \ref{Pro:existenceoforthogonal}, $\{e_i\}_{i=1}^r$ is an orthogonal basis with respect to $\norm{\ndot}$. By Proposition \ref{Pro:alpha-orthgonalextension}, it is also an orthogonal basis with respect to $\norm{\ndot}_{K,\varepsilon}$. Hence it is an Hadamard basis with respect to $\norm{\ndot}_{K,\varepsilon}$ (see Propositions \ref{Pro:orthogonalesthadamard} and \ref{Pro:ArchimedeanHadamard}). 
\end{proof}

\begin{prop}\label{Pro:extensionofdet}
Let $(V,\norm{\ndot})$ be a finite-dimensional seminormed vector space over $k$. 
Let $\norm{\ndot}_{\det,K}$ be the seminorm induced by the determinant seminorm $\norm{\ndot}_{\det}$ of $\norm{\ndot}$ by extension of scalars.
\begin{enumerate}[label=\rm(\arabic*)]
\item\label{Pro:extensionofdet:epsilon}
If either $|\ndot|$ is non-Archimedean or the seminorm $\norm{\ndot}$ is induced by a semidefinite inner product, then
the determinant seminorm $\norm{\ndot}_{K,\varepsilon,\det}$ of $\norm{\ndot}_{K,\varepsilon}$ on $\det(V_K)$ coincides with $\norm{\ndot}_{\det,K}$.

\item\label{Pro:extensionofdet:pi}
The determinant seminorm $\norm{\ndot}_{K,\pi,\det}$ of $\norm{\ndot}_{K,\pi}$ coincides with $\norm{\ndot}_{\det,K}$.

\item\label{Pro:extensionofdet:HS}
Assume that $(k,|\ndot|)$ is $\mathbb R$ equipped with the usual absolute value and
$\norm{\ndot}$ is a seminorm associated with a semidefinite inner product $\emptyinnprod$.
Then the determinant seminorm $\norm{\ndot}_{\mathbb C,\mathrm{HS},\det}$ of $\norm{\ndot}_{\mathbb C,\mathrm{HS}}$ coincides with $\norm{\ndot}_{\det,\mathbb C}$.
\end{enumerate}
\end{prop} 

\begin{proof}
\ref{Pro:extensionofdet:epsilon}
If $\norm{\ndot}$ is not a norm, then $\norm{\ndot}_{K,\varepsilon}$ is not a norm either. In this case both seminorms $\norm{\ndot}_{K,\varepsilon,\det}$ and $\norm{\ndot}_{\det,K}$ vanish. Hence we may assume without loss of generality that $\norm{\ndot}$ is a norm.

We first assume that $(V,\norm{\ndot})$ admits an Hadamard basis $\{e_i\}_{i=1}^r$. By Proposition \ref{Pro: extension Hadamard}, it is also an Hadamard basis of $(V_K,\norm{\ndot}_{K,\varepsilon})$. Moreover, by Propositions \ref{Pro:alphaorthogonale} and \ref{Pro:extensiondecorps}, for any $i\in\{1,\ldots,r\}$, one has $\|e_i\|=\|e_i\|_{K,\varepsilon}$. In particular, the vector $e_1\wedge\cdots\wedge e_r$ has the same length under the determinant norms induced by $\norm{\ndot}$ and $\norm{\ndot}_{K,\varepsilon}$. This establishes  the proposition in the particular case where $(V,\norm{\ndot})$ admits an Hadamard basis (and hence in the case where $\norm{\ndot}$ is induced by an inner product, see Proposition \ref{Pro:existenceoforthogonal}).

In the following, we assume that the absolute value $|\ndot|$ is non-Archimedean. Let $\alpha$ be an element in $\intervalle]01[$ and $\{e_i\}_{i=1}^r$ be an $\alpha$-orthogonal basis of $(V,\norm{\ndot})$. By Proposition \ref{Pro:alpha-orthgonalextension}, it is also an $\alpha$-orthogonal basis of $(V,\norm{\ndot}_{K,\varepsilon})$. By Proposition \ref{Pro:orthogonalesthadamard}, one has
\begin{align*}
\|e_1\wedge\cdots\wedge e_r\|_{K,\varepsilon,\det} & \geqslant \alpha^r \|e_1\|_{K,\varepsilon}\cdots\|e_r\|_{K,\varepsilon} \\
& \geqslant\alpha^{2r}\|e_1\|\cdots\|e_r\|\geqslant\alpha^{2r}\|e_1\wedge\cdots\wedge e_r\|_{\det},
\end{align*}
where the second inequality comes from Propositions \ref{Pro:extensiondecorps} and \ref{Pro:alphaorthogonale}.
Conversely, one has
\[\|e_1\wedge\cdots e_r\|_{K,\varepsilon\det}\leqslant\|e_1\|_{K,\varepsilon}\cdots\|e_r\|_{K,\varepsilon}\leqslant\|e_1\|\cdots\|e_r\|\leqslant\alpha^{-r}\|e_1\wedge\cdots\wedge e_r\|_{\det},\]
where the second inequality comes from Proposition \ref{Pro:extensiondecorps} and the formula \eqref{Equ:doubledual} in \S\ref{Subsec:dualnorm}, and the last inequality results from Proposition \ref{Pro:orthogonalesthadamard}
. Thus one has \[\alpha^{-r}\norm{\ndot}_{\det}\geqslant\norm{\ndot}_{K,\det}\geqslant\alpha^{2r}\norm{\ndot}_{\det}.\]
Since $\alpha\in\intervalle{]}{0}{1}{[}$ is arbitrary, we obtain $\norm{\ndot}_{\det,K}=\norm{\ndot}_{K,\varepsilon,\det}$.

\medskip
\ref{Pro:extensionofdet:pi}
Let $r$ be the rank of $V$ over $k$. Note that the $r$-th $\pi$-tensor power of the norm $\norm{\ndot}_{K,\pi}$ on $V_K^{\otimes_K r}\cong (V^{\otimes_k r})\otimes_kK$ coincides with the $\pi$-tensor product of $r$ copies of $\norm{\ndot}$ and the absolute value $|\ndot|$ on $K$ (see Proposition \ref{Pro: produit pi tensoriel extension}). Hence by Proposition \ref{Pro:quotientavecpitensor} its quotient norm on $\det(V_K)$ coincides with $\norm{\ndot}_{\det,K}$.  

\medskip
\ref{Pro:extensionofdet:HS}
Let $r$ be the rank of $V$ over $k$. Note that the $r$-th orthogonal tensor power of the norm $\norm{\ndot}_{\mathbb C,\pi}$ on $V_{\mathbb C}^{\otimes_{\mathbb C} r}\cong (V^{\otimes_{\mathbb R} r})\otimes_{\mathbb R}\mathbb C$ coincides with the orthogonal tensor product of $r$ copies of $\norm{\ndot}$ and the usual absolute value $|\ndot|$ on $\mathbb C$ (see Proposition \ref{Pro: produit HS tensoriel extension}). Hence by Proposition \ref{Pro:determinant HS} its quotient seminorm on $\det(V_{\mathbb C})$ coincides with $\norm{\ndot}_{\det,\mathbb C}$.  
\end{proof}

\if01
\begin{prop}\label{Pro:extensionofdet}
Let $(V,\norm{\ndot})$ be a finite-dimensional seminormed vector space over $k$. We assume that either $|\ndot|$ is non-Archimedean or the seminorm $\norm{\ndot}$ is induced by a  semidefinite inner product. Then the determinant seminorm of $\norm{\ndot}_{K,\varepsilon}$ on $\det(V\otimes_kK)$ coincides with the seminorm induced by the determinant seminorm of $\norm{\ndot}$ by extension of scalars.
\end{prop} 
\begin{proof}If $\norm{\ndot}$ is not a norm, then $\norm{\ndot}_{K,\varepsilon}$ is not a norm either. In this case both seminorms $\norm{\ndot}_{K,\varepsilon,\det}$ and $\norm{\ndot}_{\det,K}$ vanish. Hence we may assume without loss of generality that $\norm{\ndot}$ is a norm.

We first assume that $(V,\norm{\ndot})$ admits an Hadamard basis $\{e_i\}_{i=1}^r$. By Proposition \ref{Pro: extension Hadamard}, it is also an Hadamard basis of $(V_K,\norm{\ndot}_{K,\varepsilon})$. Moreover, by Propositions \ref{Pro:alphaorthogonale} and \ref{Pro:extensiondecorps}, for any $i\in\{1,\ldots,r\}$, one has $\|e_i\|=\|e_i\|_{K,\varepsilon}$. In particular, the vector $e_1\wedge\cdots\wedge e_r$ has the same length under the determinant norms induced by $\norm{\ndot}$ and $\norm{\ndot}_{K,\varepsilon}$. This establishes  the proposition in the particular case where $(V,\norm{\ndot})$ admits an Hadamard basis (and hence in the case where $\norm{\ndot}$ is induced by an inner product, see Proposition \ref{Pro:existenceoforthogonal}).

In the following, we assume that the absolute value $|\ndot|$ is non-Archimedean. Let $\alpha$ be an element in $\intervalle]01[$ and $\{e_i\}_{i=1}^r$ be an $\alpha$-orthogonal basis of $(V,\norm{\ndot})$. By Proposition \ref{Pro:alpha-orthgonalextension}, it is also an $\alpha$-orthogonal basis of $(V,\norm{\ndot}_{K,\varepsilon})$. By Proposition \ref{Pro:orthogonalesthadamard}, one has
\begin{align*}
\|e_1\wedge\cdots\wedge e_r\|_{K,\varepsilon,\det} & \geqslant \alpha^r \|e_1\|_{K,\varepsilon}\cdots\|e_r\|_{K,\varepsilon} \\
& \geqslant\alpha^{2r}\|e_1\|\cdots\|e_r\|\geqslant\alpha^{2r}\|e_1\wedge\cdots\wedge e_r\|_{\det},
\end{align*}
where the second inequality comes from Propositions \ref{Pro:extensiondecorps} and \ref{Pro:alphaorthogonale}, $\norm{\ndot}_{\det}$ and $\norm{\ndot}_{K,\varepsilon,\det}$ being the determinant norms induced by $\norm{\ndot}$ and $\norm{\ndot}_{K,\varepsilon}$, respectively. Conversely, one has
\[\|e_1\wedge\cdots e_r\|_{K,\varepsilon\det}\leqslant\|e_1\|_{K,\varepsilon}\cdots\|e_r\|_{K,\varepsilon}\leqslant\|e_1\|\cdots\|e_r\|\leqslant\alpha^{-r}\|e_1\wedge\cdots\wedge e_r\|_{\det},\]
where the second inequality comes from Proposition \ref{Pro:extensiondecorps} and the formula \eqref{Equ:doubledual} in \S\ref{Subsec:dualnorm}, and the last inequality results from Proposition \ref{Pro:orthogonalesthadamard}
. Thus one has \[\alpha^{-r}\norm{\ndot}_{\det}\geqslant\norm{\ndot}_{K,\det}\geqslant\alpha^{2r}\norm{\ndot}_{\det}.\]
Since $\alpha\in\intervalle{]}{0}{1}{[}$ is arbitrary, we obtain $\norm{\ndot}_{\det,K}=\norm{\ndot}_{K,\varepsilon,\det}$.
\end{proof}

\begin{prop}\label{Pro: det ext pi}
Let $(V,\norm{\ndot})$ be a finite-dimensional seminormed vector space over $k$. Then the determinant seminorm of $\norm{\ndot}_{K,\pi}$ coincides with $\norm{\ndot}_{\det,K}$.
\end{prop}
\begin{proof} 
Let $r$ be the rank of $V$ over $k$. Note that the $r$-th $\pi$-tensor power of the norm $\norm{\ndot}_{K,\pi}$ on $V_K^{\otimes_K r}\cong (V^{\otimes_k r})\otimes_kK$ coincides with the $\pi$-tensor product of $r$ copies of $\norm{\ndot}$ and the absolute value $|\ndot|$ on $K$ (see Proposition \ref{Pro: produit pi tensoriel extension}). Hence by Proposition \ref{Pro:quotientavecpitensor} its quotient norm on $\det(V_K)$ coincides with $\norm{\ndot}_{\det,K}$.  
\end{proof}

\begin{prop}
\label{Pro: det ext HS}Assume that $(k,|\ndot|)$ is $\mathbb R$ equipped with the usual absolute value. 
Let $(V,\emptyinnprod)$ be a finite-dimensional vector space over $\mathbb R$ equipped with a  semidefinite inner product and let $\norm{\ndot}$ be the seminorm associated with $\emptyinnprod$. Then the determinant seminorm of $\norm{\ndot}_{\mathbb C,\mathrm{HS}}$ coincides with $\norm{\ndot}_{\det,\mathbb C}$.
\end{prop}
\begin{proof} 
Let $r$ be the rank of $V$ over $k$. Note that the $r$-th orthogonal tensor power of the norm $\norm{\ndot}_{\mathbb C,\pi}$ on $V_{\mathbb C}^{\otimes_{\mathbb C} r}\cong (V^{\otimes_{\mathbb R} r})\otimes_{\mathbb R}\mathbb C$ coincides with the orthogonal tensor product of $r$ copies of $\norm{\ndot}$ and the usual absolute value $|\ndot|$ on $\mathbb C$ (see Proposition \ref{Pro: produit HS tensoriel extension}). Hence by Proposition \ref{Pro:determinant HS} its quotient seminorm on $\det(V_{\mathbb C})$ coincides with $\norm{\ndot}_{\det,\mathbb C}$.  
\end{proof}
\fi

\begin{prop}\label{Pro:comparisonofdualnormes:scalar:extension}
Let $(V,\norm{\ndot})$ be a finite-dimensional seminormed vector space over $k$. 
\begin{enumerate}[label=\rm(\arabic*)]
\item\label{Pro:comparisonofdualnormes:scalar:extension:epsilon}
Let $\norm{\ndot}_{K,\varepsilon,*}$ be the dual norm of $\norm{\ndot}_{K,\varepsilon}$ and $\norm{\ndot}_{*,K,\varepsilon}$ be the norm induced by $\norm{\ndot}_*$ by the $\varepsilon$-extension of scalars. Then we have $\norm{\ndot}_{K,\varepsilon,*}\geqslant\norm{\ndot}_{*,K,\varepsilon}$, and the restrictions {to} $V^*$ of these two norms are both equal to the dual norm $\norm{\ndot}_*$. Moreover, the equality $\norm{\ndot}_{K,\varepsilon,*}=\norm{\ndot}_{*,K,\varepsilon}$ holds if $|\ndot|$ is non-Archimedean or if $V$ is of rank $1$ over $k$.

\item\label{Pro:comparisonofdualnormes:scalar:extension:epsilon:pi}
The dual norm $\norm{\ndot}_{K,\pi,*}$ of $\norm{\ndot}_{K,\pi}$ is equal to $\norm{\ndot}_{*,K,\varepsilon}$ on $V_K$.
\end{enumerate}
\end{prop}
\begin{proof}
\ref{Pro:comparisonofdualnormes:scalar:extension:epsilon}
Let $\varphi$ be an element in $V_K^*$. By definition one has
\[\|\varphi\|_{K,\varepsilon,*}=\sup_{x\in V_K\setminus N_{\norm{\sndot}_{K,\varepsilon}}}\frac{|\varphi(x)|}{\|x\|_K}\geqslant\sup_{x\in V\setminus N_{\norm{\sndot}_{**}}}\frac{|\varphi(x)|}{\|x\|_{**}}=\|\varphi\|_{*,K,\varepsilon}.\]
Note that for any $x\in V_K$ one has
\[\|x\|_{K,\varepsilon}=\sup_{\alpha\in V^*\setminus\{0\}}\frac{|\alpha(x)|}{\|\alpha\|_*}.\]
Therefore, if $\varphi\in V^* \setminus\{0\}$ then one has 
$\|x\|_{K,\varepsilon}\geqslant{|\varphi(x)|}/{\|\varphi\|_*}$,
which leads to
\[\|\varphi\|_{K,\varepsilon,*}\leqslant\|\varphi\|_*=\|\varphi\|_{*,K,\varepsilon},\]
where the equality comes from Proposition \ref{Pro:extensiondecorps} (in the non-Archimedean case we use the fact that the norm $\norm{\ndot}_*$ is ultrametric).

In the following we prove the equality $\norm{\ndot}_{K,\varepsilon,*}=\norm{\ndot}_{*,K,\varepsilon}$ under the assumption that $|\ndot|$ is non-Archimedean or $\rang_k(V)=1$. We treat firstly the case where $\rang_k(V)=1$. In this case, either the seminorm $\norm{\ndot}$ vanishes and $V_K^*$ is the trivial vector space, which has only one norm, or the seminorm $\norm{\ndot}$ is a norm and for any non-zero element $\eta$ in $V$ one has
\[\|\eta^\vee\|_{K,\varepsilon,*}=\|\eta\|_{K,\varepsilon}^{-1}=\|\eta\|^{-1}=\|\eta^\vee\|_{*}=\|\eta^\vee\|_{*,K,\varepsilon},\]
where $\eta^\vee$ denotes the dual element of $\eta$ in $V^*=V^\vee$. Hence the equality $\norm{\ndot}_{K,\varepsilon,*}=\norm{\ndot}_{*,K,\varepsilon}$ always holds.

We now treat the case where the absolute value $|\ndot|$ is non-Archimedean. Note that $\norm{\ndot}_{K,\varepsilon,*}$ is an ultrametric norm on $V^*\otimes_kK$ extending $\norm{\ndot}_{*}$. Hence by Proposition \ref{Pro:extensioncorps} one has $\norm{\ndot}_{K,\varepsilon,*}\leqslant\norm{\ndot}_{*,K,\varepsilon}$. Therefore the equality $\norm{\ndot}_{*,K,\varepsilon}=\norm{\ndot}_{K,\varepsilon,*}$ holds.

\medskip
\ref{Pro:comparisonofdualnormes:scalar:extension:epsilon:pi}
If $(k,|\ndot|)$ is $\mathbb R$ equipped with the usual absolute value and if $K=\mathbb C$, then by Proposition \ref{Pro:dualitypiepsilon}, the norm $\norm{\ndot}_{*,\mathbb C,\varepsilon,*}$ identifies with the $\pi$-tensor product of $\norm{\ndot}_*$ and $|\ndot|$ (here we consider the absolute value $|\ndot|$ on $\mathbb C$ as a norm on a vector space over $\mathbb R$). Hence it is equal to the norm $\norm{\ndot}_{\mathbb C,\pi,**}$ on $V_{\mathbb C}^{**}$, which implies the equality $\norm{\ndot}_{*,\mathbb C,\varepsilon}=\norm{\ndot}_{\mathbb C,\pi,*}$ since any finite-dimensional normed vector space over $\mathbb R$ is reflexive.

Assume that $|\ndot|$ is non-Archimedean. By \ref{Pro:comparisonofdualnormes:scalar:extension:epsilon} and the fact that $\norm{\ndot}_{K,\varepsilon}\leqslant\norm{\ndot}_{K,\pi}$ (which results from Proposition \ref{Pro: maximality property pi} and Proposition \ref{Pro:extensiondecorps} \ref{item: extension eps}), one has \[\norm{\ndot}_{*,K,\varepsilon}=\norm{\ndot}_{K,\varepsilon,*}\geqslant\norm{\ndot}_{K,\pi,*},\] which leads to \[\norm{\ndot}_{K,\pi,**}\geqslant\norm{\ndot}_{K,\varepsilon,**}=\norm{\ndot}_{K,\varepsilon},\] where the 
equality comes from the fact that the norm $\norm{\ndot}_{K,\varepsilon}$ is ultrametric. Note the the restriction of $\norm{\ndot}_{K,\pi,**}$ {to} $V$ is bounded from above by $\norm{\ndot}$ since $\norm{\ndot}$ identifies with the restriction of $\norm{\ndot}_{K,\pi}$ {to} $V$ (see Proposition \ref{Pro:extensiondecorps} \ref{item: extension pi}). As $\norm{\ndot}_{K,\pi,**}$ is ultrametric, by Proposition \ref{Pro:extensioncorps} we obtain $\norm{\ndot}_{K,\pi,**}\leqslant\norm{\ndot}_{K,\varepsilon}$, which leads to the equality $\norm{\ndot}_{K,\pi,**}=\norm{\ndot}_{K,\varepsilon}$. By passing to the dual norms, using Proposition \ref{Pro:doubedualandquotient} \ref{Item: seminorm and double dual induce the same dual norm} we obtain $\norm{\ndot}_{K,\pi,*}=\norm{\ndot}_{K,\varepsilon,*}=\norm{\ndot}_{*,K,\varepsilon}$.
\end{proof}

\if01
\begin{prop}\label{Pro:comparisonofdualnormes}Let $(V,\norm{\ndot})$ be a finite-dimensional seminormed vector space over $k$. Let $\norm{\ndot}_{K,\varepsilon,*}$ be the dual norm of $\norm{\ndot}_{K,\varepsilon}$ and $\norm{\ndot}_{*,K,\varepsilon}$ be the norm induced by $\norm{\ndot}_*$ by the $\varepsilon$-extension of scalars. Then we have $\norm{\ndot}_{K,\varepsilon,*}\geqslant\norm{\ndot}_{*,K,\varepsilon}$, and the restrictions on $V^*$ of these two norms are both equal to the dual norm $\norm{\ndot}_*$. Moreover, the equality $\norm{\ndot}_{K,\varepsilon,*}=\norm{\ndot}_{*,K,\varepsilon}$ holds if $|\ndot|$ is non-Archimedean or if $V$ is of rank $1$ over $k$.
\end{prop}
\begin{proof}
Let $\varphi$ be an element in $V_K^*$. By definition one has
\[\|\varphi\|_{K,\varepsilon,*}=\sup_{x\in V_K\setminus N_{\norm{\sndot}_{K,\varepsilon}}}\frac{|\varphi(x)|}{\|x\|_K}\geqslant\sup_{x\in V\setminus N_{\norm{\sndot}_{**}}}\frac{|\varphi(x)|}{\|x\|_{**}}=\|\varphi\|_{*,K,\varepsilon}.\]
Note that for any $x\in V_K$ one has
\[\|x\|_{K,\varepsilon}=\sup_{\alpha\in V^*\setminus\{0\}}\frac{|\alpha(x)|}{\|\alpha\|_*}.\]
Therefore, if $\varphi\in V^* \setminus\{0\}$ then one has 
$\|x\|_{K,\varepsilon}\geqslant{|\varphi(x)|}/{\|\varphi\|_*}$,
which leads to
\[\|\varphi\|_{K,\varepsilon,*}\leqslant\|\varphi\|_*=\|\varphi\|_{*,K,\varepsilon},\]
where the equality comes from Proposition \ref{Pro:extensiondecorps} (in the non-Archimedean case we use the fact that the norm $\norm{\ndot}_*$ is ultrametric).

In the following we prove the equality $\norm{\ndot}_{K,\varepsilon,*}=\norm{\ndot}_{*,K,\varepsilon}$ under the assumption that $|\ndot|$ is non-Archimedean or $\rang_k(V)=1$. We treat firstly the case where $\rang_k(V)=1$. In this case, either the seminorm $\norm{\ndot}$ vanishes and $V_K^*$ is the trivial vector space, which has only one norm, or the seminorm $\norm{\ndot}$ is a norm and for any non-zero element $\eta$ in $V$ one has
$\|\eta^\vee\|_{K,\varepsilon,*}=\|\eta\|_{K,\varepsilon}^{-1}=\|\eta\|^{-1}=\|\eta^\vee\|_{*}=\|\eta^\vee\|_{*,K,\varepsilon}$,
where $\eta^\vee$ denotes the dual element of $\eta$ in $V^*=V^\vee$. Hence the equality $\norm{\ndot}_{K,\varepsilon,*}=\norm{\ndot}_{*,K,\varepsilon}$ always holds.

We now treat the case where the absolute value $|\ndot|$ is non-Archimedean. Note that $\norm{\ndot}_{K,\varepsilon,*}$ is an ultrametric norm on $V^*\otimes_kK$ extending $\norm{\ndot}_{*}$. Hence by Proposition \ref{Pro:extensioncorps} one has $\norm{\ndot}_{K,\varepsilon,*}\leqslant\norm{\ndot}_{*,K,\varepsilon}$. Therefore the equality $\norm{\ndot}_{*,K,\varepsilon}=\norm{\ndot}_{K,\varepsilon,*}$ holds.
\end{proof}

\begin{prop}
\label{Pro:dual of Kpi}
Let $(V,\norm{\ndot})$ be a finite-dimensional seminormed vector space over $k$. One has $\norm{\ndot}_{K,\pi,*}=\norm{\ndot}_{*,K,\varepsilon}$ on $V_K$.
\end{prop}
\begin{proof}
If $(k,|\ndot|)$ is $\mathbb R$ equipped with the usual absolute value and if $K=\mathbb C$, then by Proposition \ref{Pro:dualitypiepsilon}, the norm $\norm{\ndot}_{*,\mathbb C,\varepsilon,*}$ identifies with the $\pi$-tensor product of $\norm{\ndot}_*$ and $|\ndot|$ (here we consider the absolute value $|\ndot|$ on $\mathbb C$ as a norm on a vector space over $\mathbb R$). Hence it is equal to the norm $\norm{\ndot}_{\mathbb C,\pi,**}$ on $V_{\mathbb C}^{**}$, which implies the equality $\norm{\ndot}_{*,\mathbb C,\varepsilon}=\norm{\ndot}_{\mathbb C,\pi,*}$ since any finite-dimensional normed vector space over $\mathbb R$ is reflexive.

Assume that $|\ndot|$ is non-Archimedean. By Proposition \ref{Pro:comparisonofdualnormes} and the fact that $\norm{\ndot}_{K,\varepsilon}\leqslant\norm{\ndot}_{K,\pi}$ (which results from Proposition \ref{Pro: maximality property pi} and Proposition \ref{Pro:extensiondecorps} \ref{item: extension eps}), one has $\norm{\ndot}_{*,K,\varepsilon}=\norm{\ndot}_{K,\varepsilon,*}\geqslant\norm{\ndot}_{K,\pi,*}$, which leads to $\norm{\ndot}_{K,\pi,**}\geqslant\norm{\ndot}_{K,\varepsilon,**}=\norm{\ndot}_{K,\varepsilon}$, where the second equality comes from the fact that the norm $\norm{\ndot}_{K,\varepsilon}$ is ultrametric. Note the the restriction of $\norm{\ndot}_{K,\pi,**}$ on $V$ is bounded from above by $\norm{\ndot}$ since $\norm{\ndot}$ identifies with the restriction of $\norm{\ndot}_{K,\pi}$ on $V$ (see Proposition \ref{Pro:extensiondecorps} \ref{item: extension pi}). As $\norm{\ndot}_{K,\pi,**}$ is ultrametric, by Proposition \ref{Pro:extensioncorps} we obtain $\norm{\ndot}_{K,\pi,**}\leqslant\norm{\ndot}_{K,\varepsilon}$, which leads to the equality $\norm{\ndot}_{K,\pi,**}=\norm{\ndot}_{K,\varepsilon}$. By passing to the dual norms, using Proposition \ref{Pro:doubedualandquotient} \ref{Item: seminorm and double dual induce the same dual norm} we obtain $\norm{\ndot}_{K,\pi,*}=\norm{\ndot}_{K,\varepsilon,*}=\norm{\ndot}_{*,K,\varepsilon}$.\end{proof}
\fi

The following proposition is an $\varepsilon$-tensor analogue of Propositions \ref{Pro: produit pi tensoriel extension} and \ref{Pro: produit HS tensoriel extension}.

\begin{prop}\label{Pro:epstensoextcorps} We assume that the absolute value $|\ndot|$ on $k$ is non-Archimedean. 
Let $(V_1,\norm{\ndot}_1)$ and $(V_2,\norm{\ndot}_2)$ be finite-dimensional ultrametrically seminormed vector space over $k$, 
and $\norm{\ndot}$ be the $\varepsilon$-tensor product norm of $\norm{\ndot}_1$ and $\norm{\ndot}_2$. Then $\norm{\ndot}_{K,\varepsilon}$ identifies with the $\varepsilon$-tensor product of $\norm{\ndot}_{1,K,\varepsilon}$ and $\norm{\ndot}_{2,K,\varepsilon}$. 
\end{prop}
\begin{proof}
Let $\norm{\ndot}_{\varepsilon}'$ be the $\varepsilon$-tensor product of the norms $\norm{\ndot}_{1,K,\varepsilon}$ and $\norm{\ndot}_{2,K,\varepsilon}$. By Remark \ref{Rem:operateureps}, it identifies with the seminorm induced by the operator seminorm on the $K$-vector space $\Hom_{K}(V_{1,K}^*,V_{2,K})$ by the canonical $K$-linear map \[V_{1,K}\otimes_KV_{2,K}\longrightarrow \Hom_{K}(V_{1,K}^*,V_{2,K}),\] where we consider the dual norm of $\norm{\ndot}_{1,K,\varepsilon}$ on $V_{1,K}^*$, which identifies with the norm $\norm{\ndot}_{1,*,K,\varepsilon}$ induced by $\norm{\ndot}_{1,*}$ by $\varepsilon$-extension of scalars (see \ref{Pro:comparisonofdualnormes:scalar:extension:epsilon} in Proposition \ref{Pro:comparisonofdualnormes:scalar:extension}). By Proposition \ref{Pro:extensioncorpsmor}, for  any $f\in\Hom_{k}(V_1^*,V_2)$, the seminorm of $f_K$ identifies with that of $f$. Therefore $\norm{\ndot}_{\varepsilon}'$ is an ultrametric norm on $V_{1,K}\otimes_KV_{2,K}$ which extends the $\varepsilon$-tensor product {$\norm{\ndot}_{\varepsilon}$} of $\norm{\ndot}_1$ and $\norm{\ndot}_2$. By Proposition \ref{Pro:extensioncorps}, one has $\norm{\ndot}_{\varepsilon}'\leqslant\norm{\ndot}_{K,\varepsilon}$.

In the following, we prove the converse inequality $\norm{\ndot}_{K,\varepsilon}\leqslant\norm{\ndot}_{\varepsilon}'$. Let $\alpha\in\intervalle]01[$ and $\{e_i\}_{i=1}^n$ and $\{f_j\}_{j=1}^m$ be respectively $\alpha$-orthogonal bases of $(V_1,\norm{\ndot}_1)$ and $(V_2,\norm{\ndot}_2)$. By Proposition \ref{Pro:alphatenso}, they are also $\alpha$-orthogonal bases of $(V_{1,K},\norm{\ndot}_{1,K,\varepsilon})$ and $(V_{2,K},\norm{\ndot}_{2,K,\varepsilon})$, respectively. By Propsition \ref{Pro:alphatenso}, the basis $\{e_i\otimes f_j\}_{i\in\{1,\ldots,n\},\,j\in\{1,\ldots,m\}}$ of $V_{1,K}\otimes_KV_{2,K}$ is $\alpha^2$-orthogonal with respect to the seminorm $\norm{\ndot}_{\varepsilon}'$. Hence for $(a_{ij})_{i\in\{1,\ldots,n\},\,j\in\{1,\ldots,m\}}\in K^{n\times m}$ and $T=\sum_{i,j}a_{ij}e_i\otimes f_j\in V_{1,K}\otimes_KV_{2,K}$, one has
\[\|T\|_{\varepsilon}'\geqslant\alpha^2\max_{\begin{subarray}{c}i\in\{1,\ldots,n\}\\
j\in\{1,\ldots,m\}
\end{subarray}}|a_{ij}|\cdot\|e_i\otimes f_j\|_{\varepsilon}\geqslant\alpha^2\|T\|_{K,\varepsilon}.\]
Since $\alpha\in\intervalle]01[$ is arbitrary, we obtain the inequality $\norm{\ndot}_{\varepsilon}'\geqslant\norm{\ndot}_{K,\varepsilon}$. 
\end{proof}

\begin{prop}\label{prop:extension:trivial:Laurent}
We assume that the absolute value $|\ndot|$ of $k$ is trivial. 
Let $(V,\norm{\ndot})$ be an ultrametrically seminormed vector space of finite rank over $k$.
Let $(K,|\ndot|_K)$ be an extension of $(k,|\ndot|)$ such that $|\ndot|_K$ is non-trivial and complete.
Let $\mathfrak o_{K}$ be the valuation ring of $(K, |\ndot|_K)$ and $\mathfrak m_{K}$ be the maximal ideal of $\mathfrak o_{K}$. Suppose the following assumptions
\ref{Item: k into o_k induces an iso} and \ref{Item: set of quotient of norms}:
\begin{enumerate}[label=\rm(\arabic*)]
\item\label{Item: k into o_k induces an iso}
the natural map $k \to \mathfrak o_{K}$ induces
an isomorphism $k \overset{\sim}{\longrightarrow} \mathfrak o_{K}/\mathfrak m_{K}$,
\item\label{Item: set of quotient of norms} $\left\{ \|v'\|/\|v\| \,:\, v, v' \in V\setminus N_{\norm{\ndot}} \right\} \cap |{K}^{\times}|_K \subseteq \{ 1 \}$.
\end{enumerate}
Let $\norm{\ndot}_{K,\varepsilon}$ be the seminorm of $V_K$ induced by $\norm{\ndot}$ by $\varepsilon$-extension of scalars. Then $\norm{\ndot}_{K,\varepsilon}$ is the only ultrametric seminorm on $V_{K}$ extending $\norm{\ndot}$.
\end{prop}

\begin{proof}
{We prove the assertion by induction on the rank $n$ of $V$ over $k$. The case where $n=1$ is trivial. In the following, we suppose that the assertion has been proved for seminormed vector spaces of rank $<n$ over $k$. Since $\norm{\ndot}$ is ultrametric, one has $\norm{\ndot}=\norm{\ndot}_{**}$ (see Corollary \ref{Cor:doubledual}). 
Let $\norm{\ndot}'$ be another ultrametric seminorm on $V_{K}$ extending $\norm{\ndot}$. By Proposition \ref{Pro:extensioncorps}, one has $\|\ndot\|'\leqslant\norm{\ndot}_{K,\varepsilon}$.

Let $r$ be the rank of $V/N_{\norm{\ndot}}$ and $\{e_i\}_{i=1}^n$ be an orthogonal basis of $V$ such that $\{e_i\}_{i=r+1}^n$ forms a basis of $N_{\norm{\ndot}}$ (see Proposition \ref{Pro: orthogonal basis: vanishing}).
If the equality $\norm{\ndot}'=\norm{\ndot}_{K,\varepsilon}$ does not hold, then there exists a vector $x\in V_{K}$ such that $\|x\|'<\|x\|_{K,\varepsilon}$. 
We write $x$ in the form $x=a_1e_1+\cdots+a_ne_n$ with $(a_1,\ldots,a_n)\in K^n$. Note that \[\norm{a_{r+1}e_{r+1}+\cdots+a_{n}e_n}'\leqslant\max_{i\in\{r+1,\ldots,n\}}|a_i|\cdot\norm{e_i}=0.\]
For the same reason, $\norm{a_{r+1}e_{r+1}+\cdots+a_ne_n}_{K,\varepsilon}=0$. Therefore one has 
\[\norm{a_1e_1+\cdots+a_re_r}'=\norm{x}'<\norm{x}_{K,\varepsilon}=\norm{a_1e_1+\cdots+a_re_r}_{K,\varepsilon}.\]
By replacing $x$ by $a_1e_1+\cdots+a_re_r$ we many assume without loss of generality that $a_{r+1}=\cdots=a_n=0$.

We will prove that $|a_i|_K\cdot\|e_i\|$ are the same for $i\in\{1,\ldots,r\}$ by contradiction. Without loss of generality, we assume on the contrary that 
\[|a_1|_K\cdot\|e_1\|\leqslant\cdots\leqslant|a_j|_K\cdot\|e_j\|<|a_{j+1}|_K\cdot\|e_{j+1}\|=\cdots=|a_r|_K\cdot\|e_r\|\]
with $j\in\{1,\ldots,r-1\}$. 
Note that
\[\|x\|'<\|x\|_{K,\varepsilon}=\max_{i\in\{1,\ldots,r\}}|a_i|_K\cdot\|e_i\|=|a_r|_K\cdot\|e_r\|.\] 
Moreover, by the induction hypothesis, the norms $\norm{\ndot}'$ and $\norm{\ndot}_{K,\varepsilon}$ coincide on $Ke_{j+1}+\cdots+Ke_r$. In particular, one has
$\|a_{j+1}e_{j+1}+\cdots+a_re_r\|'=|a_r|_K\cdot\|e_r\|$.
Therefore, if we let $y=a_1e_1+\cdots+a_{j}e_{j}$, then one has 
\[\|y\|'=\|x-(a_{j+1}e_{j+1}+\cdots+a_re_r)\|'=|a_r|_K\cdot\|e_r\|>\max_{i\in\{1,\ldots,j\}}|a_i|_K\cdot\|e_i\|=\|y\|_{K,\varepsilon},\]
which leads to a contradiction since $\norm{\ndot}'\leqslant \norm{\ndot}_{K,\varepsilon}$. Hence we should have
\[|a_1|_K\cdot\|e_1\|=\cdots=|a_r|_K\cdot\|e_r\|.\]
By the condition \ref{Item: set of quotient of norms}, we have $\|e_1\|=\cdots=\|e_r\|$ (namely the function $\norm{\ndot}$ is constant on $V\setminus N_{\norm{\ndot}}$) and hence $|a_1|_K=\cdots=|a_r|_K>0$. 
As $|a_i/a_r|_K = 1$ for any $i\in\{1,\ldots,r\}$, by the assumption \ref{Item: k into o_k induces an iso}, there exists {a} $b_i\in k^{\times}$ such that $|a_i/a_r -b_i|_K<1$, that is, $|a_i -b_ia_r |_K<|a_r|_K$. Thus,
by Proposition~\ref{Pro:valeur},
\[{\|x\|'}=\bigg\|a_r\sum_{i=1}^rb_ie_i+\sum_{i=1}^r(a_i-b_ia_r)e_i\bigg\|'=|a_r|_K\cdot\|e_r\|=\|x\|_{K,\varepsilon}\]
because
\[\bigg\|a_r\sum_{i=1}^rb_ie_i\bigg\|'=|a_r|_K\bigg\|\sum_{i=1}^rb_ie_i\bigg\|'=|a_r|_K\|e_r\|\text{ and }\bigg\|\sum_{i=1}^r(a_i-b_ia_r)e_i\bigg\|'<|a_r|_K\|e_r\|.\]
This leads to a contradiction. The proposition is thus proved.}
\end{proof}

\begin{prop}\label{Pro:HStensorproducttex} We assume that $k=\mathbb R$ and that $|\ndot|$ is the usual absolute value. Let $\{(V_i,\norm{\ndot}_i)\}_{i=1}^n$ be  finite-dimensional seminormed vector spaces over $k$. We assume that the seminorms $\norm{\ndot}_i$  are induced by  semidefinite inner products $\langle\,,\,\rangle_i$ and we let 
$\norm{\ndot}_{\mathrm{HS}}$ be their orthogonal tensor product. For $i\in\{1,\ldots,n\}$, let $\pi_i:V_{i,\mathbb C}\rightarrow W_i$ be a quotient spaces of rank $1$ of $V_{i,\mathbb C}$, and $\norm{\ndot}_{W_i}$  be the quotient seminorm on $W_i$ induced by $\norm{\ndot}_{i,\mathbb C}$. Let $\norm{\ndot}_{W}$ be the quotient seminorm on $W=\bigotimes_{i=1}^nW_i$ induced by $\norm{\ndot}_{\mathrm{HS},\mathbb C}$ and let $\norm{\ndot}$ be the tensor product of $\norm{\ndot}_{W_i}$. Then one has \[\frac{1}{\sqrt{2}}\norm{\ndot}\leqslant\norm{\ndot}_{W}\leqslant (\sqrt 2)^n\norm{\ndot}.\]
\end{prop}
\begin{proof} For any $i\in\{1,\ldots,n\}$, let $\norm{\ndot}_i'$ be the seminorm on $V_{i,\mathbb C}$ induced by the  semidefinite inner product $\langle\,,\,\rangle_{i,\mathbb C}$. One has $\norm{\ndot}_{i,\mathbb C}\leqslant\norm{\ndot}_i'\leqslant \sqrt{2}\norm{\ndot}_{i,\mathbb C}$ (see Remark \ref{Rem:extensiondoubledual} \ref{Item:extensiondoubledual 2}). Let $\langle \,,\,\rangle_{\mathrm{HS}}$ be the  semidefinite inner product corresponding to the orthogonal tensor product seminorm $\norm{\ndot}_{\mathrm{HS}}$ and $\norm{\ndot}'$  be the seminorm on $\bigotimes_{i=1}^nV_{i,\mathbb C}$ induced by $\langle\,,\,\rangle_{\mathrm{HS},\mathbb C}$. Still by Remark \ref{Rem:extensiondoubledual}  \ref{Item:extensiondoubledual 2} one has $\norm{\ndot}_{\mathrm{HS},\mathbb C}\leqslant\norm{\ndot}'\leqslant\sqrt{2}\norm{\ndot}_{\mathrm{HS},\mathbb C}$. Moreover, $\norm{\ndot}'$ coincides with the orthogonal tensor product of the seminorms $\norm{\ndot}_i'$.

For $i\in\{1,\ldots, n\}$, let $\norm{\ndot}_{W_i}'$ be the quotient seminorms on $W_i$ induced by $\norm{\ndot}_i'$. Let $\norm{\ndot}_{W}'$ be the quotient seminorm on $W$ induced by $\norm{\ndot}'$. By Proposition \ref{Pro:quotientr1Hs}, $\norm{\ndot}_{W}'$ coincides with the tensor product of the seminorms $\norm{\ndot}_{W_i}'$. Moreover, by the relations  $\norm{\ndot}_{i,\mathbb C}\leqslant\norm{\ndot}_i'\leqslant \sqrt{2}\norm{\ndot}_{i,\mathbb C}$ we obtain $\norm{\ndot}_{W_i}\leqslant\norm{\ndot}_{W_i}'\leqslant\sqrt{2}\norm{\ndot}_{W_i}$, which implies
\begin{equation}\label{Equ:encadrementW12}\norm{\ndot}\leqslant\norm{\ndot}_W'\leqslant(\sqrt{2})^n\norm{\ndot};\end{equation}
by the relation $\norm{\ndot}_{\mathrm{HS},\mathbb C}\leqslant\norm{\ndot}'\leqslant\sqrt{2}\norm{\ndot}_{\mathrm{HS},\mathbb C}$ we obtain
\begin{equation}\label{Equ:encadrementW}
\norm{\ndot}_{W}\leqslant\norm{\ndot}_{W}'\leqslant\sqrt{2}\norm{\ndot}_{W}.
\end{equation} 
Combining \eqref{Equ:encadrementW12} and \eqref{Equ:encadrementW}, we obtain 
$\frac{1}{\sqrt{2}}\norm{\ndot}\leqslant\norm{\ndot}_{W}\leqslant (\sqrt{2})^n\norm{\ndot}$. The proposition is thus proved.
\end{proof}

\begin{prop}\label{Pro:quotientnormextensionscalar}
Let $(V,\norm{\ndot})$ be a finite-dimensional seminormed vector space over $k$ and $W$ be a quotient vector space of rank $1$ of $V$. Let $\norm{\ndot}_W$ be the quotient seminorm on $W$ induced by $\norm{\ndot}$. 
\begin{enumerate}[label=\rm(\arabic*)]
\item\label{item: quotient de rang 1}
The seminorm $\norm{\ndot}_{W,K}$ coincides with the quotient seminorm on $W_K$ induced by $\norm{\ndot}_{K,\natural}$, where $\natural=\varepsilon$ or $\pi$.
\item\label{item: quotient de rang 1 HS} Assume that $(k,|\ndot|)$ is $\mathbb R$ equipped with the usual absolute value, $K=\mathbb C$, and $\norm{\ndot}$ is induced by a  semidefinite inner product. Then the seminorm $\norm{\ndot}_{W,\mathbb C}$ coincides with the quotient seminorm on $W_{\mathbb C}$ induced by $\norm{\ndot}_{\mathbb C,\mathrm{HS}}$.
\end{enumerate}
\end{prop}
\begin{proof}
\ref{item: quotient de rang 1} The case where $\natural=\pi$ follows directly from Proposition \ref{Pro:quotientavecpitensor}. In the following, we consider the case where $\natural=\varepsilon$.

Let $\norm{\ndot}_{W_K}$ be the quotient seminorm on $W_K$ induced by the seminorm $\norm{\ndot}_{K,\varepsilon}$. If the kernel of the quotient map $V\rightarrow W$ does not contain $N_{\norm{\ndot}}$, then the quotient seminorm $\norm{\ndot}_W$ vanishes since $W$ is of dimension $1$ over $k$. In this case the quotient seminorm $\norm{\ndot}_{W_K}$ also vanishes since the kernel of the quotient map $V_K\rightarrow W_K$ does not contain $N_{\norm{\ndot}_{K,\varepsilon}}=N_{\norm{\ndot}}\otimes_kK$ (see Proposition \ref{Pro: null space in the extension}).

In the following, we assume that that seminorm $\norm{\ndot}_{W}$ is a norm. In this case $\norm{\ndot}_{W_K}$ is also a norm since the kernel of the quotient map $V_K\rightarrow W_K$ contains $N_{\norm{\ndot}_{K,\varepsilon}}=N_{\norm{\ndot}}\otimes_kK$. We will show that the dual norms $\norm{\ndot}_{W_K,*}$ and $\norm{\ndot}_{W,K,\varepsilon,*}$ on $W^\vee_K$ are equal. Since $W$ is a vector space of rank $1$, it suffices to show that the restrictions of these norms {to} $W^\vee$ are the same. We identify $W^\vee_K$ with a vector subspace of rank $1$ of $V^\vee_K$. By Proposition \ref{Pro:dualquotient}, the norm $\norm{\ndot}_{W_K,*}$ coincides with the restriction of $\norm{\ndot}_{K,\varepsilon,*}$ {to} $W_K^\vee$, where $\norm{\ndot}_{K,\varepsilon,*}$ denotes the dual seminorm of $\norm{\ndot}_{K,\varepsilon}$. By \ref{Pro:comparisonofdualnormes:scalar:extension:epsilon} in Proposition \ref{Pro:comparisonofdualnormes:scalar:extension}, the restriction of $\norm{\ndot}_{K,\varepsilon,*}$ {to} $V^\vee$ coincides with $\norm{\ndot}_*$. Therefore, the restriction of $\norm{\ndot}_{K,\varepsilon,*}$ {to} $W^\vee$ coincides with the restriction of $\norm{\ndot}_*$ {to} $W^\vee$, which identifies with the dual norm of $\norm{\ndot}_W$ (by Proposition \ref{Pro:dualquotient}). By \ref{Pro:comparisonofdualnormes:scalar:extension:epsilon} in Proposition \ref{Pro:comparisonofdualnormes:scalar:extension}, one has $\norm{\ndot}_{W,K,*}=\norm{\ndot}_{W,*,K}$. Finally, since $W$ is of rank $1$, if $k$ is non-Archimedean, then any norm on $W^\vee$ is ultrametric. Hence by Proposition \ref{Pro:extensiondecorps} (for both Archimedean and non-Archimedean cases), the restriction of $\norm{\ndot}_{W,K,*}=\norm{\ndot}_{W,*,K}$ {to} $W^\vee$ identifies with $\norm{\ndot}_{W,*}$. The assertion is thus proved.

\ref{item: quotient de rang 1 HS} follows directly from Proposition \ref{Pro:quotientr1Hs}.
\end{proof}

\begin{prop}\label{Pro: quotient metric pi ext2}
Let $(V,\norm{\ndot})$ be a finite dimension seminormed vector space over $k$ and $W$ be a quotient space of rank $1$ of $V_K=V\otimes_kK$. We equip $V_K$ with the seminorm $\norm{\ndot}_{K,\pi}$ induced by $\norm{\ndot}$ by $\pi$-extension of scalars and $W$ with the quotient seminorm $\norm{\ndot}_W$ of $\norm{\ndot}_{K,\pi}$. Then for any $\ell\in W$ one has 
\[\|\ell\|_{W}=\inf_{\begin{subarray}{c}
s\in V,\,\lambda\in K^{\times}\\
[s]=\lambda\ell
\end{subarray}}|\lambda|^{-1}\|s\|.\]
\end{prop}
\begin{proof}
By definition one has
\[\|\ell\|_{W}=\inf_{\begin{subarray}{c}s\in V_K,\,\lambda\in K^{\times}\\
[s]=\lambda\ell
\end{subarray}}|\lambda|^{-1}\cdot\|s\|_{K,\pi}\leqslant\inf_{\begin{subarray}{c}
s\in V,\,\lambda\in K^{\times}\\
[s]=\lambda\ell
\end{subarray}}|\lambda|^{-1}\cdot\|s\|_{K,\pi}=\inf_{\begin{subarray}{c}
s\in V,\,\lambda\in K^{\times}\\
[s]=\lambda\ell
\end{subarray}}|\lambda|^{-1}\cdot\|s\|,\]
where the last equality comes from Proposition \ref{Pro:extensiondecorps}. 

Without loss of generality, we may assume that $\ell\neq 0$. Let $s$ be an element in $V_K$, which is written as $s=a_1x_1+\cdots+a_nx_n$, where $(a_1,\ldots,a_n)\in K^n$ and $(x_1,\ldots,x_n)\in V$. For any $i\in\{1,\ldots,n\}$, let $\lambda_i$ be the element of $K$ such that $[x_i]=\lambda_i\ell$. Then $[s]=\lambda\ell$ with
$\lambda=a_1\lambda_1+\cdots+a_r\lambda_r$.
Let 
\[h=\inf_{\begin{subarray}{c}
t\in V,\,\lambda\in K^{\times}\\
[t]=\lambda\ell
\end{subarray}}|\lambda|^{-1}\cdot\|t\|.\]
For any $i\in\{1,\ldots,n\}$ one has $\|x_i\|\geqslant |\lambda_i|h$. Hence
\[|\lambda|^{-1}\sum_{i=1}^n|a_i|\cdot\|x_i\|\geqslant|\lambda|^{-1}\sum_{i=1}^n|a_i|\cdot|\lambda_i| h\geqslant h.\]
The proposition is thus proved.
\end{proof}

\begin{prop}\label{Pro:quotientderang1norm}
Let $(V,\norm{\ndot})$ be a finite-dimensional seminormed vector space over $k$. We assume one of the following conditions:
\begin{enumerate}[label=\rm(\roman*)]
\item
$(k,|\ndot|)$ is non-Archimedean; 
\item
$k=\mathbb C$ equipped with the usual absolute value.
\end{enumerate}  
Let $W$ be a quotient space of rank $1$ of $V\otimes_kK$. Let $\norm{\ndot}_{W}$ be the quotient seminorm on $W$ induced by $\norm{\ndot}_{K,\varepsilon}$ (the seminorm on $V\otimes_kK$ induced by $\norm{\ndot}$ by $\varepsilon$-extension of scalars). Then, for any $\ell\in W$ one has
\[\|\ell\|_{W}=\inf_{\begin{subarray}{c}
s\in V,\,\lambda\in K^{\times}\\
[s]=\lambda\ell
\end{subarray}}|\lambda|^{-1}\|s\|.\]
\end{prop}
\begin{proof}
The case where $k=\mathbb C$ equipped with the usual absolute value is trivial since $K=k$. In the following, we assume that $(k,|\ndot|)$ is non-Archimedean.

By definition one has
\[\begin{split}\|\ell\|_{W}&=\inf_{\begin{subarray}{c}s\in V_K,\,\lambda\in K^{\times}\\
[s]=\lambda\ell
\end{subarray}}|\lambda|^{-1}\cdot\|s\|_{K,\varepsilon}\leqslant\inf_{\begin{subarray}{c}
s\in V,\,\lambda\in K^{\times}\\
[s]=\lambda\ell
\end{subarray}}|\lambda|^{-1}\cdot\|s\|_{K,\varepsilon}\\&=\inf_{\begin{subarray}{c}
s\in V,\,\lambda\in K^{\times}\\
[s]=\lambda\ell
\end{subarray}}|\lambda|^{-1}\cdot\|s\|\leqslant\inf_{\begin{subarray}{c}
s\in V,\,\lambda\in K^{\times}\\
[s]=\lambda\ell
\end{subarray}}|\lambda|^{-1}\cdot\|s\|,\end{split}\]
where the last equality comes from Proposition \ref{Pro:extensiondecorps}. 

We then prove the converse inequality. Let $\alpha$ be a real number in $\intervalle]01[$. By Proposition \ref{Pro:existenceoforthogonal}
, there exists an $\alpha$-orthogonal basis $\{s_i\}_{i=1}^r$ of $(V,\norm{\ndot})$. By Proposition \ref{Pro:alpha-orthgonalextension}, $\{s_i\}_{i=1}^r$ is also an $\alpha$-orthogonal basis of $(V_K,\norm{\ndot}_K)$. For each $i\in\{1,\ldots,r\}$, let $\lambda_i\in K$ such that $[s_i]=\lambda_i\ell$. Let $s=a_1s_1+\cdots+a_rs_r$ be an element in $V\otimes_kK$, where $(a_1,\ldots,a_r)\in K^r$. Assume that $[s]$ is of the form $\lambda\ell$, where $\lambda\in K^{\times}$. Then one has $\lambda=a_1\lambda_1+\cdots+a_r\lambda_r$, which leads to
$|\lambda|\leqslant\max_{i\in\{1,\ldots,r\}}|a_i|\cdot|\lambda_i|$ since the absolute value is non-Archimedean.
By the $\alpha$-orthogonality of the basis $\{s_i\}_{i=1}^r$, we obtain
\[\begin{split}|\lambda|^{-1}&\cdot\|s\|_{K,\varepsilon}\geqslant\frac{\alpha}{|\lambda|}\max_{i\in\{1,\ldots,r\}}|a_i|\cdot\|s_i\|_{**}\\
&\geqslant\alpha\min_{i\in\{1,\ldots,r\}}|\lambda_i|^{-1}\|s_i\|_{**}\geqslant\alpha^2\min_{i\in\{1,\ldots,r\}}|\lambda_i|^{-1}\norm{s_i},\end{split}\]
where the last inequality comes from Proposition \ref{Pro:alphaorthogonale}.
The proposition is thus proved.  
\end{proof}

\begin{coro}\label{Cor:quotientnorminv}
We keep the notation and  hypotheses of Proposition \ref{Pro:quotientderang1norm}. Let $V'$ be a quotient $k$-vector space of $V$, equipped with the quotient seminorm $\norm{\ndot}'$ induced by $\norm{\ndot}$. We assume that the projection map $\pi:V_{K}\rightarrow W$ factorises through $V'_K$. Then the quotient seminorm on $W$ induced by $\norm{\ndot}'_{K,\varepsilon}$ coincides with $\norm{\ndot}_W$.
\end{coro}
\begin{proof}
Let $\norm{\ndot}_W'$ be the quotient seminorm on $W$ induced by $\norm{\ndot}'_{K,\varepsilon}$. We apply Proposition \ref{Pro:quotientderang1norm} to $(V',\norm{\ndot}')$ and $W$ to obtain that, for any $\ell\in W$, one has
\[\|\ell\|_W'=\inf_{\begin{subarray}{c}
t\in V',\,\lambda\in K^{\times}\\
[t]=\lambda\ell
\end{subarray}}|\lambda|^{-1}\|t\|'=\inf_{\begin{subarray}{c}
s\in V,\,\lambda\in K^{\times}\\
[s]=\lambda\ell
\end{subarray}}|\lambda|^{-1}\|s\|.\] Still by Proposition \ref{Pro:quotientderang1norm}, we obtain  $\|\ell\|_W'=\|\ell\|_W$.
\end{proof}

{

\subsection{Extension of scalars in the real case}
In this subsection, we assume that $(k,|\ndot|)$ is the field $\mathbb R$ of real numbers equipped with the usual absolute value.
\begin{defi}
Let $V$ be a vector space over $\mathbb R$. We say that a seminorm $\norm{\ndot}$ on $V_{\mathbb C}:=V\otimes_{\mathbb R}\mathbb C$ is \emph{invariant under the complex conjugation}\index{invariant under the complex conjugation}\index{seminorm!--- invariant under the complex conjugation} if the equality $\norm{x+iy}=\norm{x-iy}$ holds for any $(x,y)\in V^2$.
\end{defi}

\begin{prop}\label{Pro: invariant under complex conj bound}
Let $(V,\norm{\ndot})$ be a finite-dimensional seminormed vector space over $\mathbb R$. The seminorms $\norm{\ndot}_{\mathbb C,\varepsilon}$ and $\norm{\ndot}_{\mathbb C,\pi}$ are invariant under the complex conjugation. If $\norm{\ndot}$ is induced by a semidefinite inner product, then $\norm{\ndot}_{\mathbb C,\mathrm{HS}}$ is invariant under the complex conjugation.
\end{prop}
\begin{proof}
These statements follow directly from the definition of different tensor product 
seminorms and the fact that the absolute value on $\mathbb C$ is invariant under the complex conjugation (namely $|a+ib|=|a-ib|$ for any $(a,b)\in\mathbb R^2$).
\end{proof}

\begin{prop}\label{Pro: lower bound of norm extend invariant}
Let $(V,\norm{\ndot})$ be a finite-dimensional vector space over $\mathbb R$ (equipped with the usual absolute value) and $\norm{\ndot}'$ be a seminorm on $V_{\mathbb C}$ extending $\norm{\ndot}$. Assume that $\norm{\ndot}'$ is invariant under the complex conjugation. Then for any $(x,y)\in V^2$ one has $\max\{\norm{x},\norm{y}\} \leqslant \norm{x+iy}' \leqslant \norm{x}+\norm{y}$.
\end{prop}
\begin{proof} One has
\begin{gather*}2\norm{x}=\norm{2x}=\norm{2x}'\leqslant\norm{x+iy}'+\norm{x-iy}'=2\norm{x+iy}',\\
2\norm{y}=\norm{2y}=\norm{2iy}'\leqslant\norm{x+iy}'+\norm{iy-x}'=2\norm{x+iy}'.
\end{gather*}
Therefore $\norm{x+iy}'\geqslant\max\{\norm{x},\norm{y}\}$. The relation $\norm{x+iy}'\leqslant\norm{x}+\norm{y}$ comes from the triangle inequality.
\end{proof}
}

{
\begin{prop}
\label{Pro: lower bound for epsilon extension}
Let $(V,\norm{\ndot})$ be a seminormed vector space over $\mathbb R$. For any $(x,y)\in V^2$ one has \begin{equation}\label{Equ: bounds of epsilon product}\max\{\norm{x},\norm{y}\}\leqslant \norm{x+iy}_{\mathbb C,\varepsilon}\leqslant (\norm{x}^2+\norm{y}^2)^{1/2}.\end{equation}
Moreover, for any seminorm $\norm{\ndot}'$ on $V_{\mathbb C}$ extending $\norm{\ndot}$ which is invariant under the complex conjugation, one has 
\begin{gather}\label{Equ: encadrement norm ' epsilo}\frac 12\norm{\ndot}'\leqslant\norm{\ndot}_{\mathbb C,\varepsilon}\leqslant\sqrt{2}\norm{\ndot}',\\ \label{Equ: encadrement norm ' pi}\norm{\ndot}'\leqslant\norm{\ndot}_{\mathbb C,\pi}\leqslant2\norm{\ndot}'.\end{gather}
\end{prop}
\begin{proof} The first inequality of \eqref{Equ: bounds of epsilon product} comes from Propositions \ref{Pro: invariant under complex conj bound} and \ref{Pro: lower bound of norm extend invariant}. Moreover,
one has
\[\begin{split}\norm{x+iy}_{\mathbb C,\varepsilon}&=\sup_{\varphi\in V^*\setminus\{0\}}\frac{\sqrt{\varphi(x)^2+\varphi(y)^2}}{\norm{\varphi}_*}\\
&\leqslant\sup_{(\varphi_1,\varphi_2)\in (V^*\setminus\{0\})^2}\bigg({
\frac{\varphi_1(x)^2}{\norm{\varphi_1}_{*}^2}+\frac{\varphi_2(y)^2}{\norm{\varphi_2}_*^2}}\bigg)^{1/2}=(\norm{x}^2+\norm{y}^2)^{1/2},\end{split}\]which proves the second inequality of \eqref{Equ: bounds of epsilon product}.

By Proposition \ref{Pro: lower bound of norm extend invariant}
, for any $(x,y)\in V^2$, one has
\[\frac 12\norm{x+iy}'\leqslant
\frac 12(\norm{x}+\norm{y})\leqslant\max\{\norm{x},\norm{y}\}\leqslant\norm{x+iy}_{\mathbb C,\varepsilon},\]
where the last inequality comes from \eqref{Equ: bounds of epsilon product}. Moreover, still by \eqref{Equ: bounds of epsilon product} one has
\[\norm{x+iy}_{\mathbb C,\varepsilon}\leqslant(\norm{x}^2+\norm{y}^2)^{1/2}\leqslant \sqrt{2}\max\{\norm{x},\norm{y}\}\leqslant\norm{x+iy}',\]
where the last inequality comes from Proposition \ref{Pro: lower bound of norm extend invariant}. Hence \eqref{Equ: encadrement norm ' epsilo} is proved.

Since the seminorm $\norm{\ndot}'$ extends $\norm{\ndot}$, by Proposition \ref{Pro: maximality property pi} one has $\norm{\ndot}'\leqslant\norm{\ndot}_{\mathbb C,\pi}$. Moreover, for any $(x,y)\in V^2$ one has
\[\norm{x+iy}_{\mathbb C,\pi}\leqslant\norm{x}+\norm{y}\leqslant 2\max\{\norm{x},\norm{y}\}\leqslant {2\norm{x+iy}'},\]
where the last inequality comes from Proposition \ref{Pro: lower bound of norm extend invariant}. Hence \eqref{Equ: encadrement norm ' pi} is proved. 
\end{proof}

\begin{prop}\label{Pro: quotient of Galois invariant ext}
Let $(V,\norm{\ndot}_V)$ be a finite-dimensional seminormed vector space over $\mathbb R$, $Q$ be a quotient vector space of $V$ and $\norm{\ndot}_Q$ be the quotient seminorm of $\norm{\ndot}_V$ on $Q$. Let $\norm{\ndot}$ be a seminorm on $V_{\mathbb C}$ extending $\norm{\ndot}_V$, which is invariant under the complex conjugation. Then the quotient seminorm of $\norm{\ndot}$ on $Q_{\mathbb C}$ extends $\norm{\ndot}_{Q}$. It is moreover invariant under the complex conjugation.
\end{prop}
\begin{proof}
Denote by $\norm{\ndot}'$ the quotient {seminorm} of $\norm{\ndot}$ on $Q_{\mathbb C}$. For $q\in Q$ one has
\[\norm{q}'=\inf_{\substack{(x,y)\in V^2,\\ [x]=q,\,[y]=0}}\norm{x+iy} \leqslant \inf_{x \in V,\, [x]=q}\norm{x}.\]
Since $\norm{\ndot}$ is invariant under the complex conjugation, for any $(x,y)\in V^2$ one has $\norm{x+iy}\geqslant\norm{x}$. Hence
$\norm{q}' \geqslant \inf\limits_{x \in V,\, [x]=q}\norm{x}$, so that $\norm{q}' = \inf\limits_{x \in V,\, [x]=q}\norm{x}$. Therefore,
\[\norm{q}'=\inf_{x\in V,\,[x]=q}\norm{x}=\inf_{x\in V,\,[x]=q}\norm{x}_V=\norm{q}_Q.\]
Finally, for any $(p,q)\in Q^2$ one has
\[\norm{p+iq}'=\inf_{\begin{subarray}{c}
(x,y)\in V^2\\
([x],[y])=(p,q)
\end{subarray}}\norm{x+iy}=\inf_{\begin{subarray}{c}
(x,y)\in V^2\\
([x],[y])=(p,q)
\end{subarray}}\norm{x-iy}=\norm{p-iq}'.\]
\end{proof}
}

\begin{rema}\label{Rem:quotientrg1}
Let $(V,\norm{\ndot})$ be finite-dimensional seminormed vector space over $\mathbb R$, $W$ be a quotient vector space of rank one of $V_{\mathbb C}:=V\otimes_{\mathbb R}\mathbb C$. Let $\norm{\ndot}_W$ be the quotient seminorm on $W$ induced by $\norm{\ndot}_{\mathbb C,\varepsilon}$. If $\ell$ is a vector of $W$, then clearly one has
\[\|\ell\|_{W}\leqslant\inf_{\begin{subarray}{c}
s\in V,\,\lambda\in \mathbb C^{\times}\\
[s]=\lambda\ell
\end{subarray}}|\lambda|^{-1}\|s\|.
\]
The equality is in general not satisfied (see the counter-example in Remark \ref{Rem:couter-exampleextR}). However, we can show that
\begin{equation}\label{Equ:minorationnormeleW}2\|\ell\|_{W}\geqslant\inf_{\begin{subarray}{c}
s\in V,\,\lambda\in \mathbb C^{\times}\\
[s]=\lambda\ell
\end{subarray}}|\lambda|^{-1}\|s\|.
\end{equation}
 In fact, by definition one has
\[\|\ell\|_W=\inf_{\begin{subarray}{c}
s\in V_{\mathbb C},\,\lambda\in \mathbb C^{\times}\\
[s]=\lambda\ell
\end{subarray}}|\lambda|^{-1}\|s\|_{\mathbb C,\varepsilon}.\]
Let $s$ be an element in $V_{\mathbb C}$, which is written as $s=s_1 +is_2$, where $s_1$ and $s_2$ are vectors in $V$. Assume that $\lambda_1$ and $\lambda_2$ are complex numbers such that $[s_1]=\lambda_1\ell$ and $[s_2]=\lambda_2\ell$. Then one has $[s]=(\lambda_1+i\lambda_2)\ell$. By Proposition \ref{Pro: lower bound for epsilon extension},
\[\|s\|_{\mathbb C,\varepsilon}\geqslant\max\{\|s_1\|,\|s_2\|\}\geqslant\frac{1}{2}(\|s_1\|+\|s_2\|),\]
and $|\lambda_1+i\lambda_2|\leqslant|\lambda_1|+|\lambda_2|$. Hence
\[\frac{\|s\|_{\mathbb C,\varepsilon}}{|\lambda_1+i\lambda_2|}\geqslant\frac12\cdot\frac{\|s_1\|+\|s_2\|}{|\lambda_1|+|\lambda_2|}\geqslant\frac{1}{2}\inf_{\begin{subarray}{c}
s\in V,\,\lambda\in \mathbb C^{\times}\\
[s]=\lambda\ell
\end{subarray}}|\lambda|^{-1}\|s\|.\]
Thus we obtain \eqref{Equ:minorationnormeleW}.

In particular, if $V'$ is a quotient vector space of $V$ such that the projection map $V_{\mathbb C}\rightarrow W$ factorises through $V'_{\mathbb C}$, $\norm{\ndot}'$ is the quotient seminorm on $V'$ induced by $\norm{\ndot}$, and $\norm{\ndot}_W'$ is the quotient seminorm on $W$ induced by $\norm{\ndot}'_{\mathbb C,\varepsilon}$, then one has
\[\norm{\ndot}_W'\leqslant\norm{\ndot}_W\leqslant 2\norm{\ndot}_W'.\]
 In fact, by the above argument, for any non-zero element $\ell\in W$ one has
\[\|\ell\|_W'\leqslant\inf_{\begin{subarray}{c}
t\in V',\,\lambda\in\mathbb C^{\times}\\
[t]=\lambda\ell
\end{subarray}}|\lambda|^{-1}\|t\|=\inf_{\begin{subarray}{c}
s\in V,\,\lambda\in \mathbb C^{\times}\\
[s]=\lambda\ell
\end{subarray}}|\lambda|^{-1}\|s\|\leqslant 2\|\ell\|_W'.
\]
\end{rema}

The following proposition should be compared with \ref{Pro:extensionofdet:pi} and \ref{Pro:extensionofdet:HS} in Propositions \ref{Pro:extensionofdet}.

\begin{prop}\label{Pro:extdet}
Let $(V,\norm{\ndot})$ be a finite-dimensional seminormed vector space over $\mathbb R$. Denote by $\norm{\ndot}_{\det}$ and $\norm{\ndot}_{\mathbb C,\varepsilon,\det}$ the determinant seminorms induced by $\norm{\ndot}$ and $\norm{\ndot}_{\mathbb C,\varepsilon}$, respectively. Then one has
\begin{equation}\label{Equ:comparedetineq}\norm{\ndot}_{\mathbb C,\varepsilon,\det}\leqslant\norm{\ndot}_{\det,\mathbb C}\leqslant \frac{\delta(V_{\mathbb C},\norm{\ndot}_{\mathbb C,\varepsilon})}{\delta(V,\norm{\ndot})}\norm{\ndot}_{\mathbb C,\varepsilon,\det},\end{equation}
where $r$ is the rank of $V$ and $\norm{\ndot}_{\det,\mathbb C}$ is the seminorm on $\det(V)\otimes_{\mathbb R}\mathbb C$ induced by $\norm{\ndot}_{\det}$ by extension of scalars. 
\end{prop}
\begin{proof}
In the case where $\norm{\ndot}$ is not a norm, both seminorms $\norm{\ndot}_{\mathbb C,\varepsilon,\det}$ and $\norm{\ndot}_{\det,\mathbb C}$ vanish. In the following, we treat the case where $\norm{\ndot}$ is a norm.

Let $\{e_i\}_{i=1}^r$ be an Hadamard basis of $(V,\norm{\ndot})$. One has
\[\|e_1\wedge\cdots\wedge e_r\|_{\mathbb C,\varepsilon,\det}\leqslant\|e_1\|_{\mathbb C,\varepsilon}\cdots\|e_r\|_{\mathbb C,\varepsilon}=\|e_1\|\cdots\|e_r\|=\|e_1\wedge\cdots\wedge e_r\|_{\det},\]
where the first equality comes from Propositions \ref{Pro:extensiondecorps} and \ref{Pro:doubledualarch}. Hence we obtain \[\norm{\ndot}_{\mathbb C,\varepsilon,\det}\leqslant\norm{\ndot}_{\det,\mathbb C}.\] Similarly, if $\{\alpha_i\}_{i=1}^r$ is an Hadamard basis of $(V^\vee,\norm{\ndot}_*)$, one has
\[\|\alpha_1\wedge\cdots\wedge\alpha_r\|_{*,\det,\mathbb C}=\|\alpha_1\wedge\cdots\wedge\alpha_r\|_{*,\det}=\|\alpha_1\|_*\cdots\|\alpha_r\|_*,\]
where $\norm{\ndot}_{*,\det}$ denotes the determinant norm of $\norm{\ndot}_*$. Since $\alpha_1,\ldots,\alpha_r$ are elements in $V^\vee$, by \ref{Pro:comparisonofdualnormes:scalar:extension:epsilon} in Proposition \ref{Pro:comparisonofdualnormes:scalar:extension} one has $\|\alpha_i\|_*=\|\alpha_i\|_{\mathbb C,\varepsilon,*}$ for any $i\in\{1,\ldots,r\}$, where $\norm{\ndot}_{\mathbb C,\varepsilon,*}$ is the dual norm of $\norm{\ndot}_{\mathbb C,\varepsilon}$. Hence we obtain
\begin{equation}\label{Equ:minorationde*detc}\|\alpha_1\wedge\cdots\wedge\alpha_r\|_{*,\det,\mathbb C}=\|\alpha_1\|_{\mathbb C,*}\cdots\|\alpha_r\|_{\mathbb C,\varepsilon,*}\geqslant\|\alpha_1\wedge\cdots\wedge\alpha_r\|_{\mathbb C,\varepsilon,*,\det},\end{equation}
where $\norm{\ndot}_{\mathbb C,\varepsilon,*,\det}$ denotes the determinant norm of $\norm{\ndot}_{\mathbb C,\varepsilon,*}$. 

Let $\eta$ be a non-zero element of $\det(V)$, and $\eta^\vee$ be its dual element in $\det(V^\vee)$. By definition (see \S\ref{Subsec:dualdet}) one has
\begin{equation}\label{Equ:deltateta}\|\eta\|_{\det}=\delta(V,\norm{\ndot})^{-1}\|\eta^\vee\|_{*,\det}^{-1},\end{equation}
where $\norm{\ndot}_{*,\det}$ is the determinant norm of the dual norm $\norm{\ndot}_*$ on $V^{\vee}$. Since $\eta^\vee$ belongs to $V^\vee$, by \eqref{Equ:minorationde*detc} we obtain
$\|\eta^\vee\|_{*,\det}=\|\eta^\vee\|_{*,\det,\mathbb C}\geqslant\|\eta^\vee\|_{\mathbb C,\varepsilon,*,\det}$.
Hence we obtain
\[\|\eta\|_{\det}\leqslant\delta(V,\norm{\ndot})^{-1}\|\eta^\vee\|_{\mathbb C,\varepsilon,*,\det}^{-1}=\frac{\delta(V_{\mathbb C},\norm{\ndot}_{\mathbb C})}{\delta(V,\norm{\ndot})}\|\eta\|_{\mathbb C,\varepsilon,\det}.\]
The proposition is thus proved.
\end{proof}

The following proposition shows that, in the Archimedean case, the norm obtained by extension of scalars is ``almost the largest'' norm extending the initial one. 
\begin{prop}\label{Pro:normCalmostthelargest}
Let $(V,\norm{\ndot})$ be a finite-dimensional seminormed vector space over $\mathbb R$. Let $\norm{\ndot}'$ be a seminorm on $V_{\mathbb C}$ which extends $\norm{\ndot}$. Then one has {$\norm{\ndot}'\leqslant 2\norm{\ndot}_{\mathbb C,\varepsilon}$}.
\end{prop}
\begin{proof}
Let $s+it$ be an element of $V_{\mathbb C}$, where $(s,t)\in V^2$. One has $\|s+it\|'\leqslant\|s\|+\|t\|$.  By Proposition \ref{Pro: lower bound for epsilon extension}, $\max\{\|s\|,\|t\|\}\leqslant \|s+it\|_{\mathbb C,\varepsilon}$. Hence we obtain
$\|s+it\|'\leqslant 2\|s+it\|_{\mathbb C,\varepsilon}$.
\end{proof}

\begin{prop}\label{Pro:dualextensionscalar}
Let $(V,\norm{\ndot})$ be a finite-dimensional seminormed vector space over $\mathbb R$. Let $\norm{\ndot}_{\mathbb C,\varepsilon,*}$ be the dual norm of $\norm{\ndot}_{\mathbb C,\varepsilon}$ and $\norm{\ndot}_{*,\mathbb C,\varepsilon}$ be the norm on $E_{\mathbb C}^*$ induced by $\norm{\ndot}_*$ by $\varepsilon$-extension of scalars. One has
$\norm{\ndot}_{\mathbb C,\varepsilon,*}\leqslant 2\norm{\ndot}_{*,\mathbb C,\varepsilon}$.
\end{prop}
\begin{proof}
By \ref{Pro:comparisonofdualnormes:scalar:extension:epsilon} in Proposition \ref{Pro:comparisonofdualnormes:scalar:extension}, the restriction of $\norm{\ndot}_{\mathbb C,\varepsilon,*}$ {to} $V^*$ coincides with $\norm{\ndot}_*$. Hence Proposition \ref{Pro:normCalmostthelargest} leads to the inequality $\norm{\ndot}_{\mathbb C,\varepsilon,*}\leqslant 2\norm{\ndot}_{*,\mathbb C,\varepsilon}$.
\end{proof}

\begin{rema}\label{Rem:couter-exampleextR}
The results of Proposition \ref{Pro:extensionofdet} is not necessarily true for a general seminormed vector space over an Archimedean valued field. Consider the vector space $V=\mathbb R^2$ equipped with the norm $\norm{\ndot}$ such that
\[\forall\,(a,b)\in\mathbb R^2,\quad\|(a,b)\|=(\max\{a,b,0\}^2+\min\{a,b,0\}^2)^{1/2}.\]
In other words, if $a$ and $b$ have the same sign, one has $\|(a,b)\|=\max\{|a|,|b|\}$; otherwise $\|(a,b)\|=(a^2+b^2)^{1/2}$.
The unit disc of this norm is represented by Figure \ref{Fig: unit ball}.
\begin{figure}[ht]
\begin{center}
\caption{Unit ball of the norm $\|\!\cdot\!\|$}
\label{Fig: unit ball}

\vspace{2mm}
\begin{tikzpicture}[scale=0.8]
\filldraw[fill=gray!50,draw=black] (1,0)--(1,1)--(0,1) arc (90:180:1) -- (-1,-1) -- (0,-1) arc (270:360:1) --cycle;
\draw[->] (-2,0) to (2,0);
\node at (2.1,0.2) {$a$};
\draw[->] (0,-2) to (0,2);
\node at (0.3,2) {$b$};
\node at (-0.3,-0.3) {$O$};
\node at (1.2,-0.2) {$1$};
\node at (0.2,1.2) {$1$};
\end{tikzpicture}
\end{center}
\end{figure}
Let $\{e_1,e_2\}$ be the canonical basis of $\mathbb R^2$, where $e_1=(1,0)$ and $e_1=(0,1)$. One has
\[\|e_1\wedge e_2\|_{\det}=\inf_{ad-bc\neq 0}\frac{\|(a,b)\|\cdot\|(c,d)\|}{|ad-bc|},\]
where $\norm{\ndot}_{\det}$ is the determinant norm induced by $\norm{\ndot}$.
Note that if $a,b,c,d$ are four real numbers such that $\max\{|a|,|b|,|c|,|d|\}\leqslant 1$ and that $abcd\geqslant 0$, then one has $|ad-bc|\leqslant\max\{|ad|,|bc|\}\leqslant 1$ since $ad$ and $bc$ have the same sign. Hence
\[\|e_1\wedge e_2\|_{\det}=\inf_{\begin{subarray}{c}ad-bc\neq 0\\ abcd<0\end{subarray}}\frac{\|(a,b)\|\cdot\|(c,d)\|}{|ad-bc|}=\frac{1}{\sqrt{2}}.\]
Moreover, $(e_1+e_2,e_1-e_2)$ forms an Hadamard basis of $(V,\norm{\ndot})$.

The dual norm of $\norm{\ndot}$ is given by the following formula
\[\forall\,(\lambda,\mu)\in\mathbb R^2,\quad\|\lambda e_1^\vee+\mu e_2^\vee\|_*=\begin{cases}
|\lambda|+|\mu|,&\lambda\mu<0,\\
(\lambda^2+\mu^2)^{1/2},&\lambda\mu\geqslant 0.
\end{cases}\]
The unit disc of the dual norm is represented by Figure 
\ref{Fig: unit ball dual}.
\begin{figure}[ht]
\begin{center}
\caption{Unit ball of the norm $\|\!\cdot\!\|_*$}
\label{Fig: unit ball dual}

\vspace{2mm}
\begin{tikzpicture}[scale=0.8]
\filldraw[fill=gray!50,draw=black] (1,0) arc (0:90:1)--(-1,0) arc (180:270:1) -- cycle;
\draw[->] (-2,0) to (2,0);
\node at (2.1,0.3) {$\lambda$};
\draw[->] (0,-2) to (0,2);
\node at (0.3,2) {$\mu$};
\node at (-0.3,-0.3) {$O$};
\node at (1.2,-0.2) {$1$};
\node at (0.2,1.2) {$1$};
\end{tikzpicture}
\end{center}
\end{figure}

Consider now a vector $x+iy\in V\otimes_{\mathbb R}\mathbb C$, where $x$ and $y$ are vectors in $V$, and $i$ is the imaginary unit. One has
\[\|x+iy\|_{\mathbb C,\varepsilon}=\sup_{\varphi\in V^\vee\setminus\{0\}}\frac{|\varphi(x)+i\varphi(y)|}{\|\varphi\|_*}=\sup_{\varphi\in V^\vee\setminus\{0\}}\frac{(\varphi(x)^2+\varphi(y)^2)^{1/2}}{\|\varphi\|_*}.\]
In particular, one has
\[\|e_1+ie_2\|_{\mathbb C,\varepsilon}=\sup_{(\lambda,\mu)\neq(0,0)}\frac{(\lambda^2+\mu^2)^{1/2}}{f(\lambda,\mu)},\]
where
\[f(\lambda,\mu)=\begin{cases}
|\lambda|+|\mu|,&\lambda\mu<0,\\
(\lambda^2+\mu^2)^{1/2},&\lambda\mu\geqslant 0.
\end{cases}\]
Hence one has $\|e_1+ie_2\|_{\mathbb C,\varepsilon}=1$. Similarly, one has $\|ie_1+e_2\|_{\mathbb C,\varepsilon}=1$. Therefore
\[\|e_1\wedge e_2\|_{\mathbb C,\varepsilon,\det}=\frac 12\|(e_1+ie_2)\wedge (ie_1+e_2)\|_{\mathbb C,\varepsilon,\det}\leqslant \frac 12,\]
where $\norm{\ndot}_{\mathbb C,\varepsilon,\det}$ is the determinant norm associated with $\norm{\ndot}_{\mathbb C,\varepsilon}$. In particular, $(e_1+e_2,e_1-e_2)$ is no longer an Hadamard basis of $(V\otimes_{\mathbb R}\mathbb C,\norm{\ndot}_{\mathbb C,\varepsilon})$.

The above construction also provides a counter-example to the statement of Proposition \ref{Pro:quotientderang1norm} in the case where $(k,|\ndot|)$ is $\mathbb R$ equipped with the usual absolute value and $K=\mathbb C$. Consider the surjective $\mathbb C$-linear map $\pi$ from $\mathbb C^2$ to $\mathbb C$ which sends $(z_1,z_2)\in\mathbb C^2$ to $z_1-iz_2$. Let $\norm{\ndot}'$ be the quotient norm on $\mathbb C$ induced by $\norm{\ndot}_{\mathbb C,\varepsilon}$. Since $\pi(1,i)=2$ we obtain that 
\[\|1\|'\leqslant \frac{1}{2}\|(1,i)\|_{\mathbb C,\varepsilon}=\frac 12.\]
However, for any non-zero element $(\lambda,\mu)\in\mathbb R^2$ one has
\[\frac{\|(\lambda,\mu)\|}{|\pi(\lambda,\mu)|}=\frac{\|(\lambda,\mu)\|}{\sqrt{\lambda^2+\mu^2}}=\begin{cases}
\max(|\lambda|,|\mu|)/\sqrt{\lambda^2+\mu^2},&\lambda\mu\geqslant 0,\\
1,&\lambda\mu<0,
\end{cases}\]
which is bounded from below by $1/\sqrt{2}$.
\end{rema}

%% file: ch2_2019_03_23.tex

\chapter{Local metrics}

\IfChapVersion
\ChapVersion{Version of Chapter 2 : \\ \StrSubstitute{\DateChapTwo}{_}{\_}}
\fi

Throughout the chapter, let $k$ be a field equipped with an absolute value $|\ndot|$. We assume that $k$ is complete with respect to this absolute value. If $|\ndot|$ is Archimedean, we assume that it is the usual absolute value on $\mathbb R$ or $\mathbb C$.

\section{Metrised vector bundles}\label{Sec: Metrised vector bundles}
Let $\pi:X\rightarrow\Spec k$ be a $k$-scheme. Let $F_X$ be the functor from the category $\mathbf{A}_k$ of all $k$-algebras to that of sets, which sends any $k$-algebra $A$ to the set of all $k$-points of $X$ valued in $A$. Recall that a \emph{$k$-point of $X$ valued in $A$} is by definition a $k$-morphism from $\Spec A$ to $X$. If we denote by $\mathbf{E}_k$ the full subcategory of $\mathbf{A}_k$ of all field extensions of $k$, then the scheme $X$ identifies (as a set) with the colimit\footnote{One can fix two Grothendieck universes $\mathcal U$ and $\mathcal V$ such that $k\in\mathcal U\in\mathcal V$, and take $\mathbf{A}_k$ as the category of $k$-algebras whose underlying sets lie in $\mathcal U$. It is then a small category with respect to the univers $\mathcal V$. Thus we can consider the colimite of a functor from $\mathbf{A}_k$ to that of all sets in $\mathcal U$.} of the functor $F_X$ restricted to $\mathbf{E}_k$ (see \cite{Dem_Gab} page 18, th\'eor\`eme de comparaison). 

\subsection{Berkovich space associated to a scheme}
The Berkovich space associated with a $k$-scheme can also be constructed as a colimit. By \emph{valued extension}\index{valued extension} of $k$, we refer to a field extension $k'$ of $k$ equipped with an absolute value which extends $|\ndot|$ on $k$. If $(k_1,|\ndot|_1)$ and $(k_2,|\ndot|_2)$ are two valued extensions of $k$, we call \emph{morphism}\index{valued extension!morphism of@morphism of ---} from $(k_1,|\ndot|_1)$ to $(k_2,|\ndot|_2)$ any field homomorphism $k_1\rightarrow k_2$ which preserves the absolute values. The valued extensions of $k$ and morphisms between them form a category which we denote by $\mathbf{VE}_k$. One has a forgetful functor $w$  from $\mathbf{VE}_k$ to $\mathbf{A}_k$ which consists of forgetting the absolute values.

\begin{defi}\label{Def:specification} Let $X$ be a $k$-scheme.
The \emph{Berkovich space}\index{Berkovich space associated with a scheme} (as a set) associated with $X$ is defined as the colimit of the composed functor $F_X\circ w$ from $\mathbf{VE}_k$ to the category of sets, denoted by $X^{\mathrm{an}}$. The universal property of colimits defines a map $j$ from $X^{\mathrm{an}}$ to $X$, called the \emph{specification map}\index{specification map}.
\end{defi}

Let $x$ be a point of $X^{\mathrm{an}}$. We denote by $\kappa(x)$ the residue field of the point $j(x)$ of the scheme $X$, called the \emph{residue field}\index{Berkovich space associated with a scheme!residue field of a point} of $x$. We consider the point $j(x)$ as a $k$-morphism from $\Spec\kappa(x)$ to $X$. By definition, if $y:\Spec k'\rightarrow X$ is a $k$-point of $X$ taking values in some valued extension $(k',|\ndot|_y)$ of $k$, which represents the point $x\in X^{\mathrm{an}}$, then, as a $k$-morphism of schemes, it factorises through the point $j(x)$. Therefore the residue field $\kappa(x)$ is a subfield of $k'$. Note that the restriction of $|\ndot|_y$ {to} $\kappa(x)$ does not depend on the choice of the representative $y$. Hence we obtain an absolute value on $\kappa(x)$ extending $|\ndot|$ on $k$, denoted by $|\ndot|_x$, and called the absolute value \emph{induced} by $x$\index{Berkovich space associated with a scheme!absolute value induced by a point}. Note that two different points of $X^{\mathrm{an}}$ may have the same residue field. However in this case they induce different absolute values on the residue field. 

On the Berkovich space $X^{\mathrm{an}}$, one can naturally define the Zariski topology, which is the most coarse topology making the specification map $j:X^{\mathrm{an}}\rightarrow X$ continuous. Moreover, according to Berkovich \cite{Berkovich90}, the construction of $X^{\mathrm{an}}$ allows to define a finer topology, which we describe as follows. Let $U$ be an open subscheme of $X$. The ring $\mathcal O_X(U)$ 
of all regular functions on $U$ can be identified with the set of all $k$-morphisms from $U$ to the affine line $\mathbb A_k^1$. Let $f$ be a  regular function on $U$. If $k'$ is a valued extension of $k$ and $y:\Spec k'\rightarrow U$ is a $k$-point of $U$ valued in $k'$, then the evaluation of $f$ at $y$ determines an element in $k'$ which we denote by $f(y)$. Since $k'$ is equipped with an absolute value extending $|\ndot|$, we can evaluate the absolute value of the element $f(y)$, which we denote by $|f|(y)$. Note that the value of $|f|(y)$ only depends on the equivalence class of $y$ in $X^{\mathrm{an}}$. Thus we obtain a non-negative function $|f|$ defined on $U^{\mathrm{an}}=j^{-1}(U)$.

\begin{defi}\label{Def:mapassociatetoamor}
Let $X$ be a scheme over $\Spec k$. The \emph{Berkovich topology}\index{Berkovich topology} on $X^{\mathrm{an}}$ is defined as the most coarse topology on $X^{\mathrm{an}}$ which makes the specification map $j:X^{\mathrm{an}}\rightarrow X$ and all functions of the form $|f|$ continuous, where $f$ runs over the set of all regular functions on Zariski open subsets of the scheme $X$. We refer the readers to \cite[\S3.4]{Berkovich90} for more details.

The construction of Berkovich topological spaces associated with $k$-schemes is functorial. Let $X$ and $Y$ be  $k$-schemes and $\varphi:X\rightarrow Y$ be a $k$-morphism. It induces a morphism of functors from $F_X$ to $F_Y$, which determines, by passing to colimit, a map $\varphi^{\mathrm{an}}:X^{\mathrm{an}}\rightarrow Y^{\mathrm{an}}$, called the \emph{map associated with}\index{Berkovich space associated with a scheme!map associated with scheme morphism} the morphism of $k$-schemes $X\rightarrow Y$.
\end{defi}

\begin{prop}\label{Pro:fonctorialiteberko}
Let $\varphi:X\rightarrow Y$ be a morphism of $k$-schemes. Then the map $\varphi^{\mathrm{an}}$ between Berkovich spaces is continuous with respect to the Berkovich topologies.
\end{prop}
\begin{proof}
Clearly the map $\varphi^{\mathrm{an}}$ is continuous with respect to the Zariski topologies. It suffices to prove that, for any regular function $f$ on a Zariski open subset $U$ of $Y$, the function $|f|\circ\varphi^{\mathrm{an}}$ is continuous. Let $g$ be the image of $f$ by the morphism of sheaves $\mathcal O_Y\rightarrow\varphi_*(\mathcal O_X)$ in the structure of the morphism of schemes $\varphi$. It is a regular function on $\varphi^{-1}(U)$. For any $x\in\varphi^{-1}(U)^{\mathrm{an}}$, the residue field $\kappa(x)$ is a valued extension of $\kappa(y)$ with $y=\varphi^{\mathrm{an}}(x)$. Moreover, $g(x)$ is the canonical image of $f(y)$ in $\kappa(x)$. Therefore, one has $|f|\circ\varphi^{\mathrm{an}}=|g|$, which is a continuous function.
\end{proof}

\begin{rema}\label{Rem: Berkovich function}
Assume that $k$ is an Archimedean valued field, that is, $k=\mathbb R$ or $\mathbb C$. If $X$ is a $k$-scheme, then the Berkovich space $X^{\mathrm{an}}$ identifies (as a set) with the set $X(\mathbb C)$ of all complex points of $X$ modulo the action of the Galois group $\mathrm{Gal}(\mathbb C/k)$. In particular, if $X$ is the affine line $\mathbb A^1_k$, then the Berkovich space associated with $X$ is $\mathbb C$ when $k=\mathbb C$, and is $\mathbb C/\tau$ when $k=\mathbb R$, where $\tau$ denotes the complex conjugation. In this case the Berkovich topology on $X^{\mathrm{an}}$ is generated by functions of the form $|P(z)|$, where $P$ is a polynomial in $k[z]$. Therefore, it coincides with the usual topology on $\mathbb C$ or on $\mathbb C/\tau$. In fact, in the case where $k=\mathbb C$, the usual topology on $\mathbb C$ is generated by the functions $(z\in\mathbb C)\mapsto|z-a|$ (where $a\in\mathbb C$). In the case where $k=\mathbb R$, the usual topology on $\mathbb C/\tau$ is generated by the functions \[z\longmapsto |z-a|\cdot|z-\overline a|=|z^2-2\mathrm{Re}(a)z+a\bar{a}|, \text{ where $a\in\mathbb C$}.\]

For a general $k$-scheme $X$, any regular function $f$ on $X$ determines a function $f^{\mathrm{an}}$ on $X^{\mathrm{an}}$ valued in $\mathbb C$ (in the case where $k=\mathbb R$, we identify $(\mathbb A^1_k)^{\mathrm{an}}$ with the upper half-plane in $\mathbb C$). By Proposition \ref{Pro:fonctorialiteberko}, the map $f^{\mathrm{an}}$ is a continuous complex function on $X^{\mathrm{an}}$.
\end{rema}

{
\begin{prop}\label{Pro: Berkovich space of reduced scheme}
Let $X$ be a $k$-scheme, $X_{\mathrm{red}}$ be the reduced scheme associated with $X$, and $i:X_{\mathrm{red}}\rightarrow X$ be the canonical morphism. Then the associated continuous map of Berkovich spaces $i^{\mathrm{an}}:X_{\mathrm{red}}^{\mathrm{an}}\rightarrow X^{\mathrm{an}}$ is a homeomorphism. 
\end{prop}
\begin{proof}
Note that the restrictions of the functors $F_{X}$ and $F_{X_{\mathrm{red}}}$ {to} $\mathbf{E}_k$ are the same. Therefore $i^{\mathrm{an}}$ is a bijection of sets. Moreover, it is an homeomorphism if we equip $X_{\mathrm{red}}^{\mathrm{an}}$ and $X^{\mathrm{an}}$ with the Zariski topologies. Let $U$ be a Zariski open subset of $X$. By definition $\mathcal O_{X_{\mathrm{red}}}(U)$ is the reduced ring associated with $\mathcal O_X(U)$. For any nilpotent element $s$ in $\mathcal O_X(U)$ one has $s(x)=0$ for any $x\in X$. As a consequence, if $f$ is a regular function of $X$ on $U$ and if $\overline f$ is its canonical image in $\mathcal O_{X_{\mathrm{red}}}(U)$, then one has $|\overline f|=|f|$ on $U^{\mathrm{an}}$. Thus the Berkovich topologies on $X^{\mathrm{an}}$ and $X_{\mathrm{red}}^{\mathrm{an}}$ are the same.
\end{proof}
}

\begin{rema}\label{remark:reduction:point}
We assume that the absolute value $|\ndot|$ is non-Archimedean and non-trivial. Let $\mathfrak o_k$ be the valuation ring of $(k,|\ndot|)$. Let $\mathscr A$ be a finitely generated $\mathfrak o_k$-algebra, which contains $\mathfrak o_k$ as a subring. Let $A=\mathscr A\otimes_{\mathfrak o_k}k$, which identifies with the localisation of $\mathscr A$ with respect to $\mathfrak o_k\setminus\{0\}$. Note that the Berkovich space $(\Spec A)^{\mathrm{an}}$ identifies with the set of all multiplicative seminorms on $A$ extending the absolute value $|\ndot|$ on $k$. If $x$ is a point of $(\Spec A)^{\mathrm{an}}$, we denote by $\widehat{\kappa}(x)$ the completion of the residue field $\kappa(x)$ with respect to the absolute value $|\ndot|_x$, and we let $p_x$ be the $k$-morphism from $\Spec\widehat{\kappa}(x)$ to $\Spec A$ corresponding to the point $j(x)\in\Spec A$ (see Definition \ref{Def:specification} for the specification map $j$). Then $p_x$ extends to an $\mathfrak o_k$-morphism $\mathscr P_x$ from $\Spec\mathfrak o_x$ to $\Spec\mathscr A$ if and only if $|a|_x\leqslant 1$ for any $a\in\mathscr A$. In this case the image of the maximal ideal $\mathfrak m_x$ of $\mathfrak o_x$ by the morphism $\mathscr P_x:\Spec\mathfrak o_x\rightarrow \Spec\mathscr A$ identifies with the prime ideal \[(\mathscr A,|\ndot|_x)_{<1}:=\{a\in\mathscr A\,:\,|a|_x<1\}\text{ of $\mathscr A$},\] which lies in the fibre $(\Spec\mathscr A)_\circ$ of $\Spec\mathscr A$ over the maximal ideal of $\mathfrak o_k$. We denote by $(\Spec A)^{\mathrm{an}}_{\mathscr A}$ the subset of $(\Spec A)^{\mathrm{an}}$ of points $x$ such that \[\sup_{a\in\mathscr A}|a|_x\leqslant 1\] and by $r_{\mathscr A}:(\Spec A)^{\mathrm{an}}_{\mathscr A}\rightarrow(\Spec\mathscr A)_\circ$ the map sending $x\in(\Spec A)^{\mathrm{an}}_{\mathscr A}$ to $(\mathscr A,|\ndot|_x)_{<1}$, called the \emph{reduction map}\index{reduction map}. Note that the reduction map is always surjective (cf. \cite[Proposition~2.4.4]{Berkovich90} or \cite[4.13 and Proposition~4.14]{Gubler13}). 
\end{rema}

\begin{prop}\label{Pro:integralclosure} We assume that the absolute value $|\ndot|$ is non-trivial and non-Archimedean and let $\mathfrak o_k$ be the valuation ring of $(k,|\ndot|)$.
Let $\mathscr A$ be a finitely generated $\mathfrak o_k$-algebra and $A$ be the localisation of $\mathscr A$ with respect to $\mathfrak o_k\setminus\{0\}$. Then the integral closure of $\mathscr A$ in $A$ identifies with \[\bigcap_{x\in(\Spec A)_{\mathscr A}^{\mathrm{an}}}(A,|\ndot|_x)_{\leqslant 1},\] where \[(A,|\ndot|_x)_{\leqslant 1}=\{a\in A\,:\,|a|_x\leqslant 1\}.\]
In particular, if $(k,|\ndot|)$ is discrete, $\mathscr A$ is flat over $\mathfrak o_k$ and $\mathscr A/\varpi \mathscr A$ is reduced,
then \[\mathscr A = \bigcap_{x\in(\Spec A)_{\mathscr A}^{\mathrm{an}}}(A,|\ndot|_x)_{\leqslant 1},\] where $\varpi$ is a uniformizing parameter of
$(k,|\ndot|)$. 
\end{prop}
\begin{proof}
Let $\mathscr B$ be the integral closure of $\mathscr A$ in $A$. We first show that  
$\mathscr B$ is contained in  $(A, |\ndot|_x)_{\leqslant 1}$
for any
$x \in (\Spec A)^{\mathrm{an}}_{\mathscr A}$.
If $a \in \mathscr B$, then there are $a_1, \ldots, a_n \in \mathscr A$ such that
$a^n + a_1 a^{n-1} + \cdots + a_n = 0$. 
Therefore
\[
|a|_x^n  = |a^n|_x = | a_1 a^{n-1} + \cdots + a_n |_x \leqslant
\max_{i\in\{1, \ldots, n\}}  |a_i|_x \cdot|a|_x^{n-i}  \leqslant \max_{i\in\{1, \ldots, n\}}  |a|_x^{n-i} ,
\]
which implies that $|a|_x \leqslant 1$.

Let $a \in A$ such that $a$ is not integral over $\mathscr A$. Since $A$ is a $k$-algebra of finite type, it is a Noetherian ring which is non-zero (since $a\in A$). In particular, it admits only finitely many minimal prime ideals $S^{-1}\mathfrak p_1,\ldots,S^{-1}\mathfrak p_n$, where $\mathfrak p_1,\ldots,\mathfrak p_n$ are prime ideals of $\mathscr A$ which do not intersect $S:=\mathfrak o_k\setminus\{0\}$. We show that there exists $j\in\{1,\ldots,n\}$ such that the canonical image of $a$ in $A/S^{-1}\mathfrak p_j$ is not integral over $\mathscr A/\mathfrak p_j$. Assume that, for any $i\in\{1,\ldots,n\}$, $f_i$ is a monic polynomial in $(\mathscr A/\mathfrak p_i)[T]$ such that $f_i(\lambda_i)=0$, where $\lambda_i$ is the class of $a$ in $A/S^{-1}\mathfrak p_i$. Let $F_i$ be a monic polynomial in $\mathscr A[T]$ whose reduction modulo $\mathfrak p_i[T]$ {coincides} with $f_i$. One has 
$F_i({a})\in S^{-1}\mathfrak p_i$ for any $i\in\{1,\ldots,n\}$. Let $F$ be the product of the polynomials $F_1,\ldots,F_n$. Then $F({a})$ belongs to the intersection $\bigcap_{i=1}^nS^{-1}\mathfrak p_i$, hence is nilpotent, which implies that $a$ is integral over $\mathscr A$. To show that there exists $x\in(\Spec A)^{\mathrm{an}}_{\mathscr A}$ such that $|a|_x>1$ we may replace $\mathscr A$ (resp. $A$) by $\mathscr A/\mathfrak p_j$ (resp. $A/S^{-1}\mathfrak p_j$) and hence assume that $\mathscr A$ is an integral domain without loss of generality. 

Let $b = a^{-1}$ in the fraction field of $A$.
We assert that 
\[
b\mathscr A[b] \cap \mathfrak o_k \not= \{ 0 \}\quad\text{and}\quad
1 \not\in b\mathscr A[b].
\]
We set $a = a'/s$ for some $a' \in \mathscr A$ and 
$s \in {S}$.
Then $s = b a' \in b\mathscr A[b] \cap \mathfrak o_k$, so that $b\mathscr A[b] \cap \mathfrak o_k \not= \{ 0 \}$.
Next we assume that $1 \in b \mathscr A[b]$. Then 
there exist $m\in\mathbb N_{\geqslant 1}$ and $(a_1',\ldots,a_m')\in\mathscr A^m$ such that
\[
1 = a'_1 b + a'_2 b^2 + \cdots + a'_{m} b^{m},
\]
or equivalently
$a^{m} = a'_1 a^{m-1} + \cdots + a'_{m}$, which is a contradiction.

Let $\mathfrak p$ be a maximal ideal of $\mathscr A[b]$ such that $b\mathscr A[b] \subseteq \mathfrak p$.
As $\mathfrak p \cap \mathfrak o_k \not= \{ 0 \}$ and
$\mathfrak p \cap \mathfrak o_k \subseteq \mathfrak m_k$ (where $\mathfrak m_k$ is the maximal ideal of $\mathfrak o_k$), we have
$\mathfrak p \cap \mathfrak o_k = \mathfrak m_k$ (since the Krull dimension of $\mathfrak o_k$ is  
$1$\footnote{It suffices to see $\mathfrak m_k \subseteq \mathfrak p$
for a non-zero prime ideal $\mathfrak p$ of $\mathfrak o_k$. Fix $e \in \mathfrak p \setminus \{ 0 \}$. If $x \in \mathfrak m_k$, then
$x^n e^{-1} \in \mathfrak o_k$ for some positive integer $n$ because $| x | < 1$,
so that $x^n \in \mathfrak o_k e \subseteq \mathfrak p$, and
hence $x \in \mathfrak p$.}), and hence $\mathfrak p$ lies in the  fibre $(\Spec \mathscr A[b])_{\circ}$ of $\Spec\mathscr A[b]$ over $\mathfrak m_k$.
Note that $\mathscr A[b]$ is finitely generated over $\mathfrak o_k$ and $\mathscr A[b] \otimes_{\mathfrak o_k} k = A[b]$. 
Thus, since the reduction map 
\[
r_{\mathscr A[b]} : (\Spec A[b])^{\mathrm{an}}_{\mathscr A[b]} \longrightarrow (\Spec \mathscr A[b])_{\circ}
\]
is surjective,
there is $x \in  (\Spec A[b])^{\mathrm{an}}_{\mathscr A[b]}$ such that
$r_{\mathscr A[b]}(x) = \mathfrak p$. Clearly $x \in (\Spec A)^{\mathrm{an}}_{\mathscr A}$.
As $b \in \mathfrak p$, we have $|b|_x < 1$, so that $|a|_x > 1$ because $ab = 1$.
Therefore,
\[
a \not\in  \bigcap_{ x \in (\Spec A)^{\mathrm{an}}_{\mathscr A}} (A, |\ndot|_x)_{\leqslant 1},
\]
as required.

Finally we consider the last assertion. We assume that there is 
\[
a \in \bigcap_{x\in(\Spec A)_{\mathscr A}^{\mathrm{an}}}(A,|\ndot|_x)_{\leqslant 1} \ \setminus\ \mathscr A.
\]
By the previous result, there are $a_1, \ldots, a_n \in \mathscr A$ such that
\[
a^n + a_1 a^{n-1} + \cdots + a_{n-1} a + a_n = 0.
\]
One can choose a positive integer $e$ such that
$\varpi^e a \in \mathscr A$ and $\varpi^{e-1} a \not\in \mathscr A$. As
\[
(\varpi^e a)^n + \varpi^e a_1 (\varpi^e a)^{n-1} + \cdots + \varpi^{e(n-1)}a_{n-1} (\varpi^e a) + \varpi^{en}a_n = 0,
\]
$(\varpi^e a)^n = 0$ in $\mathscr A/\varpi \mathscr A$, so that $\varpi^e a = 0$ in $\mathscr A/\varpi \mathscr A$ because
$\mathscr A/\varpi \mathscr A$ is reduced. Therefore there is $a' \in \mathscr A$ such that $\varpi^e a = \varpi a'$, and hence
$\varpi^{e-1} a = a' \in \mathscr A$ because $\mathscr A$ is flat over $\mathfrak o_k$. This is a contradiction.

\end{proof}

\subsection{Metric on a vector bundle} Let $X$ be a scheme over $\Spec k$. We denote by $\mathcal F_{X^{\mathrm{an}}}$ the sheaf of all real-valued  functions on the Berkovich {topological} space $X^{\mathrm{an}}$. Let  $\mathcal C^0_{X^{\mathrm{an}}}$ be the subsheaf of $\mathcal F_{X^{\mathrm{an}}}$ of continuous functions.

\begin{defi}
Let $E$ be  a locally free $\mathcal O_X$-module of finite rank. We call \emph{metric}\index{metric} on $E$ any family $\varphi=\{|\ndot|_{\varphi}(x)\}_{x\in X^{\mathrm{an}}}$, where each $|\ndot|_{\varphi}(x)$
 is a norm on $E(x):=E\otimes\widehat{\kappa}(x)$, $\widehat{\kappa}(x)$ being the completion of $\kappa(x)$ with respect to the absolute value $|\ndot|_x$.
We use the symbol $|\ndot|_{\varphi}$ instead of the usual double bar symbol 
in order to distinguish a local norm from a global seminorm (cf. Definition~\ref{def:sup:norm}).

Note that the family $\varphi$ actually defines a morphism of sheaves (of sets) from $E$ to $j_*(\mathcal F_{X^{\mathrm{an}}})$, which sends each section $s$ of $E$ over a Zariski open subset $U$ of $X$ to the function $|s|_{\varphi}:U^{\mathrm{an}}\rightarrow\mathbb R_{\geqslant 0}$ sending $x\in U^{\mathrm{an}}$ to \[|s|_\varphi(x):=|s(x)|_{\varphi}(x),\] where $s(x)$ denotes the reduction of $s$ in $E(x)$. If this morphism of sheaves takes values in $j_*(\mathcal C^0_{X^{\mathrm{an}}})$ (namely,  for any section $s$ of $E$ on a Zariski open subset of $X$, the function $|s|_\varphi$ is continuous with respect to the Berkovich topology), we say that the metric $\varphi$ is \emph{continuous}\index{metric!continuous ---}.
\end{defi}

\begin{rema}
Let $E$ be a locally free $\mathcal O_X$-module of finite rank,
equipped with a continuous metric $\varphi$. Let $F$ be a locally free sub-$\mathcal O_X$-module of $E$. For any $x\in X^{\mathrm{an}}$, the restriction of the norm $|\ndot|_\varphi(x)$ {to} $F(x)$ defines a norm on $F(x)$. These norms actually define a continuous metric on $F$. However, we don't know if, for any quotient vector bundle of $E$, the quotient norms of $|\ndot|_\varphi(x)$ ($x\in X^{\mathrm{an}}$) define a continuous metric on the quotient bundle.
\end{rema}

The following lemma is used in the proof of Proposition~\ref{Pro:pullback}.

\begin{lemm}\label{Lem:encadrementcontinue}
Let $M$ be a topological space and $f$ be a non-negative function on $M$. Suppose that, for any $\alpha\in\intervalle]01[$, there exists a continuous function $f_\alpha$ on $M$ such that $\alpha f_\alpha\leqslant f\leqslant f_\alpha$. Then the function $f$ is continuous.
\end{lemm}
\begin{proof}
Let $x_0$ be a point of $M$. From the inequalities $\alpha f_\alpha\leqslant f\leqslant f_\alpha$, we deduce
\[\liminf_{x\rightarrow x_0}\alpha f_\alpha(x)\leqslant\liminf_{x\rightarrow x_0} f(x)\leqslant\limsup_{x\rightarrow x_0}f(x)\leqslant\limsup_{x\rightarrow x_0}f_\alpha(x).\]
Since the function $f_\alpha$ is continuous, we obtain
\[\alpha f_\alpha(x_0)\leqslant\liminf_{x\rightarrow x_0} f(x)\leqslant\limsup_{x\rightarrow x_0}f(x)\leqslant f_\alpha(x_0).\]
Moreover, one has $\alpha f_\alpha(x_0)\leqslant f(x_0)\leqslant f_\alpha(x_0)$. Hence
\[\liminf_{x\rightarrow x_0}f(x)\leqslant\limsup_{x\rightarrow x_0}f(x)\leqslant\alpha^{-1}f(x_0)\leqslant\alpha^{-1}f_\alpha(x_0)\leqslant\alpha^{-2}\liminf_{x\rightarrow x_0}f(x).\]
Since $\alpha\in\intervalle]01[$ is arbitrary and $\liminf_{x\rightarrow x_0}f(x)$ is finite, we obtain
\[\liminf_{x\rightarrow x_0}f(x)=\limsup_{x\rightarrow x_0}f(x)=f(x_0).\]
The proposition is thus proved.
\end{proof}

\begin{defi}
Let $\pi:X\rightarrow\Spec(k)$ be a $k$-scheme and $\overline{V} = (V,\|\ndot\|)$ be a finite-dimensional normed vector space over $k$. For any $x\in X^{\mathrm{an}}$, let $|\ndot|_{\overline{V},\varepsilon}(x)$ be the norm on $V\otimes_k\widehat{\kappa}(x)$ induced by $\|\ndot\|$ by $\varepsilon$-extension of scalars, and $|\ndot|_{\overline V,\pi}$ be the norm on $V\otimes_k\widehat{\kappa}(x)$ induced by $\|\ndot\|$ by $\pi$-extension of scalars (see \S\ref{Subsec:extensionofscalars}). If $|\ndot|$ is Archimedean and if the norm $\|\ndot\|$ is induced by an inner product, for any $x\in X^{\mathrm{an}}$, we denote by $|\ndot|_{\overline{V},\mathrm{HS}}(x)$ the norm on  $V\otimes_k\widehat{\kappa}(x)$ induced by $\|\ndot\|$ by orthogonal extension of scalars.
For simplicity, the norms $|\ndot|_{\overline{V},\varepsilon}(x)$, $|\ndot|_{\overline V,\pi}(x)$ and $|\ndot|_{\overline V,\mathrm{HS}}(x)$ are often denoted by $|\ndot|_{\varepsilon}(x)$, $|\ndot|_{\pi}(x)$ and $|\ndot|_{\mathrm{HS}}(x)$, respectively.
\end{defi}

\begin{prop}\label{Pro:pullback}
The norms $|\ndot|_{\overline{V},\varepsilon}(x)$, $x\in X^{\mathrm{an}}$ define a continuous metric on the locally free $\mathcal O_X$-module  
$\pi^*(V)$.
\end{prop}
\begin{proof} Let $U$ be a Zariski open subset of $X$ and $s$ be a section of $\pi^*(V)$ over $U$. It {suffices} to prove that the function $(x\in U^{\mathrm{an}})\mapsto |s|_{\overline{V},\varepsilon}(x)$ is continuous on $U^{\mathrm{an}}$.

We first treat the non-Archimedean case. By Proposition \ref{Pro:existenceoforthogonal}, for any $\alpha\in\intervalle]01[$, there exists an $\alpha$-orthogonal basis $\{e_i\}_{i=1}^n$ of $V$. By Proposition \ref{Pro:alpha-orthgonalextension}, for any $x\in X^{\mathrm{an}}$, $\{e_i\}_{i=1}^n$ is also an $\alpha$-orthogonal basis of $(V\otimes_k\widehat{\kappa}(x),|\ndot|_{\overline{V},\varepsilon}(x))$. We can write $s$ into the form $s=f_1e_1+\cdots+f_ne_n$, where $f_1,\ldots,f_n$ are regular functions on $U$. Since $\{e_i\}_{i=1}^n$ is an $\alpha$-orthogonal basis and the norm $|\ndot|_{\overline{V},\varepsilon}(x)$ is ultrametric for any $x\in U^{\mathrm{an}}$, one has
\[\forall\,x\in U^{\mathrm{an}},\quad\alpha\max_{i\in\{1,\ldots,n\}}|f_i|_x\cdot|e_i|_{\overline{V},\varepsilon}(x)\leqslant|s|_{\overline{V},\varepsilon}(x)\leqslant\max_{i\in\{1,\ldots,n\}}|f_i|_x\cdot|e_i|_{\overline{V},\varepsilon}(x).\]
By Proposition \ref{Pro:extensiondecorps}, one has $|e_i|_{\overline{V},\varepsilon}(x)=\|e_i\|_{**}$ for any $x\in X^{\mathrm{an}}$. Hence
\[\forall\,x\in U^{\mathrm{an}},\quad\alpha\max_{i\in\{1,\ldots,n\}}|f_i|_x\cdot\|e_i\|_{**}\leqslant|s|_{\overline{V},\varepsilon}(x)\leqslant\max_{i\in\{1,\ldots,n\}}|f_i|_x\cdot\|e_i\|_{**}.\] 
Note that the function $(x\in U^{\mathrm{an}})\mapsto |f_i|_x$ is continuous for any $i$. Hence the function
\[(x\in U^{\mathrm{an}})\longmapsto\max_{i\in\{1,\ldots,n\}}|f_i|_x\cdot\|e_i\|_{**}\]
is also continuous. Since $\alpha\in \intervalle]01[$ is arbitrary, by Lemma \ref{Lem:encadrementcontinue} we obtain that the function $(x\in U^{\mathrm{an}})\mapsto |s|_{\overline{V},\varepsilon}(x)$ is continuous.

We now consider the Archimedean case. Let $\{e_i\}_{i=1}^n$ be a basis of $V$. We write the section $s$ in the form $f_1e_1+\cdots+f_ne_n$, where $f_1,\ldots,f_n$ are regular functions on $U$. Note that $f_1^{\mathrm{an}},\ldots,f_n^{\mathrm{an}}$ are continuous complex functions on $U^{\mathrm{an}}$. Since the norm $\|\ndot\|_{\mathbb C,\varepsilon}$ is a continuous function on $V_{\mathbb C}$, we obtain that the map (see Remark \ref{Rem: Berkovich function})
\[(x\in X^{\mathrm{an}})\longmapsto |s|_{\overline{V},\varepsilon}(x)=\|f_1^{\mathrm{an}}(x)e_1+\cdots+f_n^{\mathrm{an}}(x)e_n\|_{\mathbb C,\varepsilon}\]
is a continuous function on $U^{\mathrm{an}}$. The proposition is thus proved.
\end{proof}

\begin{prop}
We assume that the absolute value $|\ndot|$ is Archimedean. Let $\overline V=(V,\|\ndot\|)$ be a finite-dimensional normed vector space over $k$.
\begin{enumerate}
\renewcommand{\labelenumi}{{\rm(\arabic{enumi})}}
\item The norms $|\ndot|_{\overline{V},\pi}(x)$, $x\in X^{\mathrm{an}}$ define a continuous metric on the locally free $\mathcal O_X$-module $\pi^*(V)$.
\item If the norm $\|\ndot\|$ is induced by an inner product, then the norms $|\ndot|_{\overline V,\mathrm{HS}}(x)$, $x\in X^{\mathrm{an}}$ define a continuous metric on the locally free $\mathcal O_X$-module $\pi^*(V)$. 
\end{enumerate}
\end{prop}
\begin{proof}
The proof is quite similar to the second part of the proof of Proposition \ref{Pro:pullback}, where we use the continuity of the norms $\|\ndot\|_{\mathbb C,\pi}$ and $\|\ndot\|_{\mathbb C,\mathrm{HS}}$ (in the case where the norm $\|\ndot\|$ is induced by an inner product) on the topological space $V_{\mathbb C}$.
\end{proof}

\begin{prop}\label{Pro:dualnormondualspace}
We assume that the field $k$ is Archimedean. Let $\pi:X\rightarrow\Spec(k)$ be a $k$-scheme and $\overline{V} = (V,\|\ndot\|)$ be a finite-dimensional normed vector space over $\mathbb R$. For any $x\in X^{\mathrm{an}}$, let $|\ndot|_{\overline{V}}(x)$ be the norm on $V\otimes_k\widehat{\kappa}(x)$ induced by $\|\ndot\|$ by $\natural$-extension of scalars, where $\natural=\varepsilon$, $\pi$ or $\mathrm{HS}$ (in the case where $\|\ndot\|$ is induced by an inner product) and let $|\ndot|_{\overline{V}}(x)_*$ be the dual norm on $V^\vee\otimes_k\widehat{\kappa}(x)$ of $|\ndot|_{\overline{V}}(x)$. Then the norms $|\ndot|_{\overline{V}}(x)_*$, $x\in X^{\mathrm{an}}$ define a continuous metric on the locally free $\mathcal O_X$-module $\pi^*(V^\vee)$.
\end{prop}
\begin{proof}
Let $\{\alpha_i\}_{i=1}^n$ be a basis of $V^\vee$. Locally on a Zariski open subset $U$ of $X$, any element $s\in H^0(U,\pi^*(V^\vee))$ can be written in the form $s=f_1\alpha_1+\cdots+f_n\alpha_n$, where $f_1,\ldots,f_n$ are regular functions on $U$. Let $\|\ndot\|_{\mathbb C,\natural, *}$ be the dual norm of $\|\ndot\|_{\mathbb C,\natural}$ (the norm on $V_{\mathbb C}$ induced by $\|\ndot\|$ by $\natural$-extension of scalars). We claim that \begin{equation}\label{Equ:norm1des}|s|_{\overline{V}}(x)_*=\|f_1^{\mathrm{an}}(x)\alpha_1+\cdots+f_n^{\mathrm{an}}(x)\alpha_n\|_{\mathbb C,\natural, *}.\end{equation} The equality follows from the definition of $|\ndot|_{\overline{V}}(x)_*$ when $\kappa(x)=\mathbb C$. In the case where $\kappa(x)=\mathbb R$, the norm $|\ndot|_{\overline{V}}(x)_*$ is the dual norm of $\|\ndot\|$. Hence it coincides with the restriction of $\|\ndot\|_{\mathbb C,\natural, *}$ {to} $V^\vee$ (see Proposition~\ref{Pro:comparisonofdualnormes:scalar:extension} \ref{Pro:comparisonofdualnormes:scalar:extension:epsilon}, \ref{Pro:comparisonofdualnormes:scalar:extension:epsilon:pi} and Remark \ref{Rem:extensiondoubledual} \ref{Item:extensiondoubledual 2} for the cases $\natural=\varepsilon$, $\pi$ and $\mathrm{HS}$, respectively). Thus the equality \eqref{Equ:norm1des} also holds in this case. Since the norm $\|\ndot\|_{\mathbb C,\natural, *}$ is a continuous function on $V^\vee_{\mathbb C}$, we obtain that the function $(x \in U^{\mathrm{an}}) \mapsto |s|_{\overline{V}}(x)_*$ is continuous. 
\end{proof}

\begin{defi}\label{def:sup:norm}

If the $k$-scheme $X$ is proper, then the associated Berkovich space $X^{\mathrm{an}}$ is compact (see \cite{Berkovich90} Proposition 3.4.8). In particular, if $E$ is a locally free $\mathcal O_X$-module 
equipped with a continuous metric $\varphi$, for any global section $s\in H^0(X,E)$, the number
\[\|s\|_{\varphi}:=\sup_{x\in X^{\mathrm{an}}}|s|_\varphi(x)\]
is finite. Thus we obtain a map $\|\ndot\|_{\varphi}: H^0(X,E)\rightarrow\mathbb R_+$, which is actually a seminorm 
on the $k$-vector space $ H^0(X,E)$.

Let $X_{\mathrm{red}}$ be the reduced scheme associated with $X$ and $E_{\mathrm{red}} := E \otimes_{\mathcal O_X} \mathcal O_{X_{\mathrm{red}}}$. Note that the natural morphism $ X^{\mathrm{an}}_{\mathrm{red}}\to X^{\mathrm{an}} $ is
a homeomorphism (see Proposition \ref{Pro: Berkovich space of reduced scheme}), so that to give a continuous metric of $E$ on $X^{\mathrm{an}}$ is equivalent to give
a continuous metric of $E_{\mathrm{red}}$ on $X^{\mathrm{an}}_{\mathrm{red}}$.
The corresponding metric of $E_{\mathrm{red}}$ is denoted by $\varphi_{\mathrm{red}}$. 
Moreover, if $X$ is proper and we denote the natural homomorphism $H^0(X, E) \to H^0(X_{\mathrm{red}}, E_{\mathrm{red}})$ by $\gamma_{E}$,
then it is easy to see that $\| s \|_{\varphi} = \| \gamma_E (s) \|_{\varphi_{\mathrm{red}}}$ for any $s\in H^0(X,E)$.
By (1) in the following proposition, the null space of $\|\ndot\|_{\varphi}$ coincides with
the kernel of $\gamma_E$, which is denoted by $\mathcal N(X, E)$. The induced norm on $H^0(X, E)/\mathcal N(X, E)$
is denoted by $\|\ndot\|^{\sim}_{\varphi}$.
\end{defi}

\begin{prop}\phantomsection\label{prop:def:sup:norm}
\begin{enumerate}[label=\rm(\arabic*)]
\item\label{Item: reduce case is a norm}
If $X$ is reduced, then $\|\ndot\|_{\varphi}$ is actually a norm.

\item\label{Item: quotient by null space is surjective}
For any $x \in X^{\mathrm{an}}$, the image of $\mathcal N(X, E) \otimes_{k} \widehat{\kappa}(x)$ by the natural homomorphism $H^0(X, E) \otimes_{k} \widehat{\kappa}(x) \to E \otimes_{\mathcal O_X} \widehat{\kappa}(x)$ is zero, so that
one has the induced homomorphism $(H^0(X, E)/\mathcal N(X, E)) \otimes_{k} \widehat{\kappa}(x) \to E \otimes_{\mathcal O_X} \widehat{\kappa}(x)$.
Moreover, if $E$ is generated by global sections, then
$(H^0(X, E)/\mathcal N(X, E)) \otimes_{k} \widehat{\kappa}(x) \to E \otimes_{\mathcal O_X} \widehat{\kappa}(x)$
is surjective for all $x \in X^{\mathrm{an}}$.

\item\label{Item: null space nilpotent}
If $E$ is invertible and $s \in \mathcal N(X, E)$, then there is a positive integer $n_0$ such that
$s^{\otimes n} = 0$ for all $n \geqslant n_0$.
\end{enumerate}
\end{prop}

\begin{proof}
\ref{Item: reduce case is a norm} It is sufficient to see that if $\|s\|_{\varphi}=0$, then $s = 0$.
Let $\eta_1, \ldots, \eta_r$ be the generic points of the irreducible components of $X$.
Let $\tilde{\eta}_i$ be a point of $X^{\mathrm{an}}$ such that $j(\tilde{\eta}_i) = \eta_i$.
By our assumption, $|s|_{\varphi}(\tilde{\eta}_i) = 0$, so that $s(\eta_i) = 0$ for all $i$.
Therefore, one has the assertion because $X$ is reduced.

\medskip\ref{Item: quotient by null space is surjective} We denote the natural homeomorphism $X^{\mathrm{an}} \to X^{\mathrm{an}}_{\mathrm{red}}$ by $p$.
Then we have the following commutative diagram:
\[
\xymatrix{
H^0(X, E) \otimes_k \widehat{\kappa}(x) \ar[r]^(0.45){\sim} \ar[d] & 
H^0(X, E) \otimes_k \widehat{\kappa}(p(x)) \ar[r]^(0.44){\gamma_E \otimes \operatorname{id}} \ar[d] & 
H^0(X_{\mathrm{red}}, E_{\mathrm{red}}) \otimes_k \widehat{\kappa}(p(x)) \ar[dl] \\
E \otimes_{\mathcal O_X} \widehat{\kappa}(x)  \ar[r]^(0.42){\sim} & 
E_{\mathrm{red}} \otimes_{\mathcal O_{X_{\mathrm{red}}}} \widehat{\kappa}(p(x)) & \\
}
\]
because $E_{\mathrm{red}} =
E \otimes_{\mathcal O_X} \mathcal O_{X_{\mathrm{red}}}$.
Therefore one has (2).

\medskip\ref{Item: null space nilpotent} Let $X = \bigcup_{i=1}^N \Spec(A_i)$ be an affine open covering of $X$ such that, for each $i$,
there is a local basis $\omega_i$ of $E$ over $\Spec(A_i)$. For any $i\in\{1,\ldots,N\}$, let $a_i\in A_i$ such that $s=a_i\omega_i$ on $\Spec A_i$. As $\rest{s}{X_{\mathrm{red}}} = 0$, $a_i$ is a nilpotent element of $A_i$, so that
we can find a positive integer $n_0$ such that $s^{\otimes n_0} = 0$. Therefore $s^{\otimes n} = 0$
for $n \geqslant n_0$ because $s^{\otimes n} = s^{\otimes n_0} s^{\otimes n - n_0}$.
\end{proof}

\subsection{Base change}\label{SubSec:basechange} Let  $k'/k$ be a field extension equipped with an absolute value $|\ndot|'$ which extends $|\ndot|$ on $k$.  We assume that $k'$ is complete with respect to this absolute value.

Let $X$ be a scheme over $\Spec k$, $X'$ be the fibre product $X\times_{\Spec k}\Spec k'$ and $p:X'\rightarrow X$ be the morphism of projection. If $(K,|\ndot|_K)$ is a valued extension of $k'$ and $f:\Spec K\rightarrow X'$ is a $k'$-point of $X'$ valued in $K$, then the composition morphism $\pi\circ f$ is a $k$-point of $X$ valued in $K$. This construction is functorial and thus determines by passing to colimit a surjective map between Berkovich spaces from ${X'}^{\mathrm{an}}$ to $X^{\mathrm{an}}$ which we denote by $p^{\natural}$. We emphasise that ${X'}^{\mathrm{an}}$ is constructed from the projective $k'$-scheme $X'$. Thus $p^{\natural}$ differs from the map between Berkovich spaces associated with $p$ (considered as a $k$-morphism of schemes) as in Definition \ref{Def:mapassociatetoamor}.

\begin{prop}\label{Pro:continuouspnatural}
The map $p^{\natural}:{X'}^{\mathrm{an}}\rightarrow X^{\mathrm{an}}$ defined above is continuous with respect to the Berkovich topology.
\end{prop}
\begin{proof}
Let $U$ be a Zariski open subset of $X$ and $g$ be a regular function on $U$. We denote by $g'$ the pull-back of $g$ by $p$, which is a regular function on $p^{-1}(U)$. For any $y\in{X'}^{\mathrm{an}}$ one has $|g'|(y)=|g|(p^{\natural}(y))$. Hence $|g|\circ p^{\natural}$ is a continuous function on $p^{-1}(U)^{\mathrm{an}}$. Therefore the map $p^{\natural}$ is continuous.
\end{proof}

\begin{defi}\label{def:metric:base:change}
Let $E$ be a locally free $\mathcal O_X$-module of finite rank, equipped with a metric $\varphi$ and let $E_{k'}$ be the pull-back of $E$ by the projection morphism $p:X'\rightarrow X$. If $y$ is a point of ${X'}^{\mathrm{an}}$ and if $x=p^{\natural}(y)$, then the norm $|\ndot|_{\varphi}(x)$ on $E(x)=E\otimes\widehat{\kappa}(x)$ induces by $\varepsilon$-extension (resp. $\pi$-extension)  of scalars a norm on $E_{k'}(y)\cong E(x)\otimes_{\widehat{\kappa}(x)}\widehat{\kappa}(y)$, denoted by $|\ndot|_{\varphi_{k',\varepsilon}}(y)$ (resp. $|\ndot|_{\varphi_{k',\pi}}(y)$). These norms define a metric on $E_{k'}$, denoted by $\varphi_{k',\varepsilon}$ (resp. $\varphi_{k',\pi}$), called the metric \emph{induced by $\varphi$ by $\varepsilon$-extension\index{metric!induced by $\varepsilon$-extension of scalars} (resp. $\pi$-extension)} of scalars\index{metric!induced by $\varepsilon$-extension of scalars}.

Assume that the norm $|\ndot|_{\varphi}(x)$ on $E(x)$ is induced by an inner product. For any point $y\in {X'}^{\mathrm{an}}$ such that $x=p^\natural(y)$, we let $|\ndot|_{\varphi_{k',\mathrm{HS}}}(y)$ be the norm on $E_{k'}(y)$ induced by $|\ndot|_{\varphi}(x)$ by orthogonal extension of scalars. These norms define a metric on $E_{k'}$, denoted by $\varphi_{k',\mathrm{HS}}$ and called the metric \emph{induced by $\varphi$ by orthogonal extension of scalars.}\index{metric!induced by orthogonal extension of scalars}

In the case where $E$ is an invertible $\mathcal O_X$-module, the three metrics $\varphi_{k',\varepsilon}$, $\varphi_{k',\pi}$ and $\varphi_{k,\mathrm{HS}}$ are the same (see Remark \ref{Rem:extensiondoubledual} \ref{item: dimension 1 ext}), and are just denoted by $\varphi_{k'}$.
\end{defi}

\begin{prop}\label{prop:cont:inv:pullback:base:change}
Let $L$ be an invertible $\mathcal O_X$-module equipped with a continuous metric $\varphi$. 
Then one has the following:
\begin{enumerate}[label=\rm(\arabic*)]
\item\label{Item: varphi k' continuous}
The metric $\varphi_{k'}$ is continuous.

\item\label{Item: computation norm s varphi}
For all $s \in H^0(X, L)$, $\| s \|_{\varphi} = \| p^*(s) \|_{\varphi_{k'}}$ .
\end{enumerate}
\end{prop}

\begin{proof}
Let $U$ be a Zariski open subset of $X$ and $s \in H^0(U, L)$. Then one has
\begin{equation}\label{eqn:prop:cont:inv:pullback:base:change:01}
|p^*(s)|_{\varphi_{k'}}=|s|_{\varphi}\circ p^{\natural}|_{p^{-1}(U)^{\mathrm{an}}},
\end{equation}
so that the assertion \ref{Item: computation norm s varphi} follows.
If we assume that $L$ is trivialised by $s$ over $U$, then
$p^*(s)$ is a section in $H^0(p^{-1}(U),L_{k'})$ which trivialises $L_{k'}$ on $p^{-1}(U)$.
Moreover, by the above \eqref{eqn:prop:cont:inv:pullback:base:change:01} together with Proposition \ref{Pro:continuouspnatural},  we obtain that $|p^*(s)|_{\varphi_{k'}}$ is a continuous function on $p^{-1}(U)^{\mathrm{an}}$.
Therefore one has the assertion \ref{Item: varphi k' continuous}.
\end{proof}

\section{Metrics on  invertible sheaves}\label{Subsec:metricsoninvertiblesheaves}

Let $X$ be a scheme over $\Spec k$. In this section, we discuss constructions and properties of metrics on invertible $\mathcal O_X$-modules.

\subsection{Dual metric and tensor product metric}
Let $L$ be an invertible $\mathcal O_X$-module
and $\varphi$ be a metric on $L$. Note that for any $x\in X^{\mathrm{an}}$, the norm $|\ndot|_{\varphi}(x)$ is determined by its value on any non-zero element of $L\otimes\widehat{\kappa}(x)$. In particular, to verify that the metric $\varphi$ is continuous, it suffices to prove that there exists a covering $\{U_i\}_{i\in I}$ of $X$ by affine open subsets and for each $i\in I$ there exists a section $s_i \in H^0(U_i,L)$ which trivialises the invertible sheaf $L$ on $U_i$ such that the function $|s_i|_\varphi$ is continuous on the topological space $U_i^{\mathrm{an}}$.

\begin{defi}\label{Def:dualmetric}Let $L$ be an invertible $\mathcal O_X$-module. If $\varphi$ is a metric on $L$, then the dual $\mathcal O_X$-module $L^\vee$ is naturally equipped with a metric $\varphi^\vee$ such that, for sections $\alpha$ and $s$ of $L^\vee$ and $L$ over a Zariski open subset $U$ of $X$ respectively, one has 
\[\forall\,x\in U^{\mathrm{an}},
\quad|\alpha(s)|(x)=|s|_{\varphi}(x) \cdot|\alpha|_{\varphi^{\vee}}(x).\]
We call $\varphi^\vee$ the \emph{dual metric}\index{dual metric}\index{metric!dual ---} of $\varphi$ and we also use the expression $-\varphi$ to denote the metric $\varphi^\vee$. 
\end{defi}

\begin{prop}
\label{Pro:dualcontinuousmetric}Let $X$ be a $k$-scheme, $L$ be an invertible $\mathcal O_X$-module and $\varphi$ be a metric on $L$.
If $\varphi$ is a continuous metric, then $\varphi^\vee$ is also continuous.
\end{prop}
\begin{proof}
Let $U$ be a Zariski open subset of $X$ on which the invertible sheaves $L$ and $L^\vee$ are trivialised by sections $s\in\Gamma(U,L)$ and $\alpha\in\Gamma(U,L^\vee)$ respectively. Then $\alpha(s)$ is a regular function, and
\[|\alpha|_{\varphi^\vee}=\frac{|\alpha(s)|}{|s|_\varphi}\]
on $U^{\mathrm{an}}$. Since the functions $|\alpha(s)|$ and $|s|_\varphi$ are all continuous, also is $|\alpha|_{\varphi^\vee}$. Since $U$ is arbitrary, we obtain that $\varphi^\vee$ is a continuous metric.
\end{proof}

\begin{defi}\label{Def:rootmetric}Let $L$ be an invertible $\mathcal O_X$-module and $n$ be a positive integer.
Suppose given a metric $\varphi$ on $L^{\otimes n}$. Then the maps \[(s\in H^0(U,L))\longmapsto |s^n|_{\varphi}^{1/n},\]
with $U$ running over the set of all Zariski open subsets of $X$, define a metric on $L$, denoted by $\frac 1n\varphi$. If the metric $\varphi$ is continuous, then also is $\frac 1n\varphi$.
\end{defi}

\begin{defi}
\label{Def:sumoftwometrics}
Suppose given two invertible $\mathcal O_X$-modules $L_1$ and $L_2$, equipped with  metrics $\varphi_1$ and $\varphi_2$ respectively. We denote by $\varphi_1+\varphi_2$ the  metric on $L_1\otimes L_2$ such that, for any Zariski open subset $U$ of $X$ and all sections $s_1\in H^0(U,L_1)$, $s_2\in H^0(U,L_2)$, one has
\[\forall\,x\in U^{\mathrm{an}},
\quad |s_1\cdot s_2|_{\varphi_1+\varphi_2}(x)=|s_1|_{\varphi_1}(x)\cdot|s_2|_{\varphi_2}(x).\] The metric $\varphi_1+\varphi_2$ is called \emph{tensor product}\index{tensor product metric}\index{metric!tensor product ---} of $\varphi_1$ and $\varphi_2$. Note that, if the metrics $\varphi_1$ and $\varphi_2$ are continuous, then also is $\varphi_1+\varphi_2$.
We also use the expression $\varphi_1-\varphi_2$ to denote the  metric $\varphi_1+\varphi_2^\vee$ on $L_1\otimes L_2^\vee$. If $L$ is an invertible $\mathcal O_X$-module equipped with a  metric $\varphi$, for any integer $n\in\mathbb N_{\geqslant 1}$, we use the expression $n\varphi$ to denote the metric $\varphi+\cdots+\varphi$ ($n$ copies) on $L^{\otimes n}$.
\end{defi}

\begin{prop}\label{Pro:sousmultiplicativesec}
Let $X$ be a scheme over $\Spec k$, $L_1$ and $L_2$ be invertible $\mathcal O_X$-modules,  
and $\varphi_1$ and $\varphi_2$ be continuous metrics on $L_1$ and $L_2$, respectively.
Then the canonical $k$-linear homomorphism $H^0(X,L_1)\otimes_kH^0(X,L_2)\rightarrow H^0(X,L_1\otimes L_2)$, sending $s_1\otimes s_2\in H^0(X,L_1)\otimes_kH^0(X,L_2)$ to $s_1\cdot s_2$, has operator norm $\leqslant 1$, where we consider the $\pi$-tensor product of $\|\ndot\|_{\varphi_1}$ and $\|\ndot\|_{\varphi_2}$ on the tensor product space, and the norm $\|\ndot\|_{\varphi_1+\varphi_2}$ on $H^0(X,L_1\otimes L_2)$. In particular, if 
$s_1$ and $s_2$ are elements in $H^0(X,L_1)$ and $H^0(X,L_2)$, respectively,
then the following inequality holds
\begin{equation}\label{Equ:sousmultiplicativesec}
\|s_1\cdot s_2\|_{\varphi_1 + \varphi_2}\leqslant\|s_1\|_{\varphi_1}\cdot\|s_2\|_{\varphi_2}.
\end{equation} 
\end{prop}

\begin{proof} Let $\eta$ be an element of $H^0(X,L_1)\otimes_kH^0(X,L_2)$, which is written as
\[\eta=\sum_{i=1}^Ns_1^{(i)}\otimes s_2^{(i)},\] where $s_1^{(1)},\ldots,s_1^{(N)}$ are elements in $H^0(X,L_1)$, $s_2^{(1)},\ldots,s_2^{(N)}$ are elements in $H^0(X,L_2)$. Let $s$ be the element
\[\sum_{i=1}^{N}s_1^{(i)}\cdot s_2^{(i)}\]
in $H^0(X,L_1\otimes L_2)$, which is the image of $\eta$ by the canonical homomorphism \[H^0(X,L_1)\otimes_kH^0(X,L_2)\longrightarrow H^0(X,L_1\otimes L_2).\]
For any $x\in X^{\mathrm{an}}$ one has
\[\begin{split}&\quad\;|s|_{\varphi_1+\varphi_2}(x)=\Big|\sum_{i=1}^Ns_1^{(i)}\cdot s_2^{(i)}\Big|_{\varphi_1+\varphi_2}(x)\leqslant\sum_{i=1}^N|s_1^{(i)}\cdot s_2^{(i)}|_{\varphi_1+\varphi_2}(x)\\&=\sum_{i=1}^{N}|s_1^{(i)}|_{\varphi_1}(x)\cdot|s_2^{(i)}|_{\varphi_2}(x)\leqslant\sum_{i=1}^{N}\|s_1^{(i)}\|_{\varphi_1}\cdot\|s_2^{(i)}\|_{\varphi_2}.
\end{split}\]
Since $x\in X^{\mathrm{an}}$ is arbitrary, we obtain
\[\|s\|_{\varphi_1+\varphi_2}
\leqslant\sum_{i=1}^{N}\|s_1^{(i)}\|_{\varphi_1}\cdot\|s_2^{(i)}\|_{\varphi_2}.\]
Therefore $\|s\|_{\varphi_1+\varphi_2}\leqslant \|\eta\|_\pi$, where $\|\ndot\|_\pi$ denotes the $\pi$-tensor product of $\|\ndot\|_{\varphi_1}$ and $\|\ndot\|_{\varphi_2}$. The first assertion is thus proved.

If 
$s_1$ and $s_2$ are elements in $H^0(X,L_1)$ and $H^0(X,L_2)$ respectively, then one has\[\|s_1 \cdot s_2 \|_{\varphi_1 + \varphi_2}=\sup_{x\in X^{\mathrm{an}}} |s_1 \cdot s_2|_{\varphi_1 + \varphi_2}(x)=\sup_{x\in X^{\mathrm{an}}}|s_1|_{\varphi_1}(x)  \cdot |s_2|_{\varphi_2}(x)
\leqslant  \|s_1\|_{\varphi_1} \cdot \|s_2 \|_{\varphi_2},\]
as required.
\end{proof}

\begin{rema}
Assume that the absolute value $|\ndot|$ is non-Archimedean. The statement of Proposition \ref{Pro:sousmultiplicativesec} remains true if we consider the $\varepsilon$-tensor product $\norm{\ndot}_{\varepsilon}$ of $\norm{\ndot}_{\varphi_1}$ and $\norm{\ndot}_{\varphi_2}$ on the tensor product space  $H^0(X,L_1)\otimes_kH^0(X,L_2)$. In fact, by Proposition \ref{Pro:alphatenso}, if $\{e_i\}_{i=1}^n$ and $\{f_j\}_{j=1}^m$ are $\alpha$-orthogonal basis  of $(H^0(X,L_1),\norm{\ndot}_{\varphi_1})$ and $(H^0(X,L_2),\norm{\ndot}_{\varphi_2})$ respectively, then \[\{e_i\otimes f_j\}_{(i,j)\in\{1,\ldots,n\}\times\{1,\ldots,m\}}\]
is an $\alpha^2$-orthogonal basis with respect to $\norm{\ndot}_{\varepsilon}$, where $\alpha\in\mathopen{]}0,1\mathclose{[}$. For any
\[\eta=\sum_{i=1}^n\sum_{j=1}^ma_{ij}e_i\otimes e_j\in H^0(X,L_1)\otimes_kH^0(X,L_2),\]
one has
\[\norm{\eta}_{\varepsilon}\geqslant \alpha^2\max_{(i,j)\in\{1,\ldots,n\}\times\{1,\ldots,m\}}|a_{ij}|\cdot\norm{e_i}\cdot\norm{f_j},\]
which is bounded from below by $\alpha^2$ times the norm of the canonical image of $\eta$ in $H^0(X,L_1\otimes L_2)$ since the norm $\norm{\ndot}_{\varphi_1+\varphi_2}$ is ultrametric.
\end{rema}

\subsection{Distance between metrics}\label{subsec:distance:metric}
Let $\varphi$ be a metric on $\mathcal O_X$. Then $-\ln|\boldsymbol{1}|_{\varphi}$ is a  function on $X^{\mathrm{an}}$, where $\boldsymbol{1}$ denotes the section of unity of $\mathcal O_X$. If $\varphi$ is a continuous metric, then $-\ln|\boldsymbol{1}|_{\varphi}$ is a continuous function. Conversely, any real-valued function $g$ on $X^{\mathrm{an}}$ determines a metric on $\mathcal O_X$ such that the norm at $x\in X^{\mathrm{an}}$ of the section of unity of $\mathcal O_X$ is $\mathrm{e}^{-g(x)}$. The metric is continuous if and only if the function $g$ is continuous. Therefore the set of all metrics on $\mathcal O_X$ is canonically in bijection with the set of all real-valued function on $X^{\mathrm{an}}$. This correspondance also maps bijectively the set of all continuous metrics on $\mathcal O_X$ to the set $C^0(X^{\mathrm{an}})$ of all continuous real-valued functions on $X^{\mathrm{an}}$.

\begin{defi}\label{Def:localdistancemetric}Let $L$ be an invertible $\mathcal O_X$-module. If $\varphi$ and $\varphi'$ are two  metrics on $L$, then $\varphi'-\varphi$ is a  metric on $L\otimes L^\vee\cong\mathcal O_X$, hence corresponds to a real valued function on $X^{\mathrm{an}}
$. By abuse of notation, we use the expression $\varphi'-\varphi$ to denote this function. We say that the metric $\varphi'$ is \emph{larger}\index{metric!larger than another metric} than $\varphi$ if $\varphi'-\varphi$ is a non-negative function and we use the expressions $\varphi'\geqslant\varphi$ or $\varphi\leqslant\varphi'$ to denote the relation ``\emph{$\varphi'$ is larger than $\varphi$}''. If $\varphi$ and $\varphi'$ are  metrics on $L$, we denote by $d(\varphi,\varphi')$ the element
\[\sup_{x\in X^{\mathrm{an}}}|\varphi'-\varphi|(x)\in\mathbb R_{\geqslant 0}\cup\{+\infty\},\]
called the \emph{distance}\index{distance between metrics}\index{metric!distance between ---}\label{Page:distance} between $\varphi$ and $\varphi'$. Note that one has
\begin{equation}
d(\varphi,\varphi')=\sup_{x\in X^{\mathrm{an}}}\big|\ln|\ndot|_{\varphi}(x)-\ln|\ndot|_{\varphi'}(x)\big|.
\end{equation}

\end{defi}

\begin{prop}\label{Prop:Equ:comparaisondistancesup}
If the $k$-scheme $X$ is proper (so that the sup seminorms are defined), then
\begin{equation}\label{Equ:comparaisondistancesup}
d(\|\ndot\|^{\sim}_{\varphi},\|\ndot\|^{\sim}_{\varphi'})\leqslant d(\varphi,\varphi')
\end{equation}
(see \S\ref{Subsec:Distance} for the notion of distance between two norms and \S\ref{Sec: Norms} for the notion of norm associated with a seminorm).
\end{prop}

\begin{proof}
Fix $s \in H^0(X, L) \setminus \mathcal N(X, L)$.
For $\epsilon > 0$, one can choose $x \in X^{\mathrm{an}}$ such that
$\mathrm{e}^{-\epsilon} \| s \|_{\varphi_1} \leqslant |s|_{\varphi_1}(x)$. Then
\[
\ln \|s \|_{\varphi_1} - \ln \|s \|_{\varphi_2} \leqslant
\ln |s|_{\varphi_1}(x) - \ln |s|_{\varphi_2}(x) + \epsilon \leqslant d(\varphi_1, \varphi_2) + \epsilon,
\]
so that one has $\ln \|s \|_{\varphi_1} - \ln \|s \|_{\varphi_2} \leqslant d(\varphi_1, \varphi_2)$
because $\epsilon$ is an arbitrary positive number. In the same way, $\ln \|s \|_{\varphi_2} - \ln \|s \|_{\varphi_1} \leqslant d(\varphi_1, \varphi_2)$, and hence one obtains
\[
\big| \ln \|s \|_{\varphi_1} - \ln \|s \|_{\varphi_2} \big| \leqslant d(\varphi_1, \varphi_2),
\]
which implies the assertion of the proposition.
\end{proof}

For any integer $n\in\mathbb Z$ one has 
$n\varphi'-n\varphi=n(\varphi'-\varphi)$ and hence
\begin{equation}\label{Equ:distancen}d(n\varphi',n\varphi)=|n|d(\varphi',\varphi).\end{equation}
The distance function verifies the triangle inequality: if $\varphi_1,\varphi_2$ and $\varphi_3$ are three continuous metrics on $L$, then one has 
\begin{equation}\label{Equ:triangleinequality}d(\varphi_1,\varphi_3)\leqslant d(\varphi_1,\varphi_2)+d(\varphi_2,\varphi_3)\end{equation}
because \[\big|\ln |\ndot|_{\varphi_1}(x) - \ln |\ndot|_{\varphi_3}(x)\big| \leqslant \big|\ln |\ndot|_{\varphi_1}(x) -\ln |\ndot|_{\varphi_2}(x)\big| + \big|\ln |\ndot|_{\varphi_2}(x) - \ln |\ndot|_{\varphi_3}(x)\big|\]
for any $x\in X^{\mathrm{an}}$.

\begin{defi}\label{Def:pull-back}
Let $Y$ and $X$ be two schemes over $\Spec k$, and $f:Y\rightarrow X$ be a $k$-morphism. Suppose given an invertible $\mathcal O_X$-module $L$, equipped with a metric $\varphi$. Then the metric $\varphi$ induces by pull-back a metric $f^*(\varphi)$ on $Y$ such that, for any $y\in Y^{\mathrm{an}}$, the norm $|\ndot|_{f^*(\varphi)}(y)$ is induced by $|\ndot|_{\varphi}(f(y))$ by extension of scalars. The metric $f^*(\varphi)$ is called the \emph{pull-back}\index{metric!pull-back}\index{pull-back of a metric} of $\varphi$ by $f$. For any section $s$ of $L$ on a Zariski open subset $U$ of $X$, one has \begin{equation}\label{Equ:pullback}|f^*(s)|_{f^*(\varphi)}=|s|_\varphi\circ f^{\mathrm{an}}|_{f^{-1}(U)^{\mathrm{an}}}. \end{equation}
In particular, if the metric $\varphi$ is continuous, then also is $f^*(\varphi)$.
\end{defi}

\begin{prop}
Let $Y$ and $X$ be two schemes over $\Spec k$, $f:Y\rightarrow X$ be a $k$-morphism, $L$ be an invertible $\mathcal O_X$-module, and $\varphi$ and $\varphi'$ be two metrics on $L$. Then one has
\[d(f^*(\varphi),f^*(\varphi'))\leqslant d(\varphi,\varphi').\]
Moreover, the equality holds if $f : Y \to X$ is surjective.
\end{prop} 
\begin{proof}
By \eqref{Equ:pullback}, one has $f^*(\varphi)-f^*(\varphi')=(\varphi-\varphi')\circ f^{\mathrm{an}}$.
Hence
\[d(f^*(\varphi),f^*(\varphi'))=\sup_{y\in Y^{\mathrm{an}}}|f^*(\varphi)-f^*(\varphi')|(y)\leqslant\sup_{y\in Y^{\mathrm{an}}}|\varphi-\varphi'|(y)=d(\varphi,\varphi').\]
If $f : Y \to X$ is surjective, then $f^{\mathrm{an}} : Y^{\mathrm{an}} \to X^{\mathrm{an}}$ is also
surjective, so that the last assertion follows.
\end{proof}

\subsection{Fubini-Study metric} Let $V$ be a finite-dimensional vector space over $k$. We denote by $\pi:\mathbb P(V)\rightarrow\Spec(k)$ the projective space of $V$. Note that the functor $F_{\mathbb P(V)}$ from the category $\mathbf{A}_k$ of $k$-algebras to the category of sets corresponding to $\mathbb P(V)$ (see \S\ref{Sec: Metrised vector bundles}) sends any $k$-algebra $A$ to the set of all projective quotient $A$-modules of $V\otimes_kA$ which are of rank $1$. By gluing morphisms of schemes, we obtain that, for any $k$-scheme $f:X\rightarrow\Spec(k)$, the set of all $k$-morphisms from $X$ to $\mathbb P(V)$ is in functorial bijection with the set of all invertible quotient $\mathcal O_X$-module of $f^*(V)$. In the case where $X$ is the projective space $\mathbb P(V)$, the invertible quotient $\mathcal O_X$-module of $\pi^*(V)$ corresponding to the identity map $\mathbb P(V)\rightarrow\mathbb P(V)$ is called the \emph{universal invertible sheaf}\index{universal invertible sheaf}, denoted by $\mathcal O_V(1)$. It verifies the following universal property: for any $k$-scheme $f:X\rightarrow\Spec k$, a $k$-morphism $g:X\rightarrow\mathbb P(V)$ corresponds to the invertible quotient
\[g^*(p):g^*(\pi^*(V))\cong f^*(V)\longrightarrow g^*(\mathcal O_V(1)),\]
where $p:\pi^*(V)\rightarrow\mathcal O_V(1)$ is the quotient homomorphism defining the universal invertible sheaf.

Let $\overline{V} = (V, \|\ndot\|)$ be a normed vector space of finite rank over $k$.
For any point $x$ in the Berkovich space $\mathbb P(V)^{\mathrm{an}}$, if the absolute value $|\ndot|$ is non-Archimedean, we denote by  $|\ndot|_{\overline{V}}(x)$ the norm on $V\otimes_k\widehat{\kappa}(x)$ induced by $\|\ndot\|$ by $\varepsilon$-extension of scalars; if the absolute value $|\ndot|$ is Archimedean, we denote by  $|\ndot|_{\overline{V}}(x)$ the norm on $V\otimes_k\widehat{\kappa}(x)$ induced by $\|\ndot\|$ by $\pi$-extension of scalars. We emphasise that, in the case where $\widehat{\kappa}(x)=k$ (namely $x$ corresponds to a rational point of $\mathbb P(V)$), the vector space $V\otimes_k\widehat{\kappa}(x)$ is canonically isomorphic to $V$ and the norm $|\ndot|_{\overline V}(x)$ identifies with the double dual norm of $\|\ndot\|$. We denote by $|\ndot|_{\overline{V},\mathrm{FS}}(x)$ the quotient norm on $\mathcal O_V(1)(x)=\mathcal O_V(1)\otimes_{\mathcal O_{\mathbb P(V)}}\widehat{\kappa}(x)$ induced by the norm $|\ndot|_{\overline{V}}(x)$ on $V\otimes_k\widehat{\kappa}(x)$, called the \emph{Fubini-Study norm}\index{Fubini-Study norm} on $\mathcal O_V(1)(x)$ induced by $\|\ndot\|$.
For simplicity, the norm $|\ndot|_{\overline{V},\mathrm{FS}}(x)$ is often denoted by
$|\ndot|_{\mathrm{FS}}(x)$.

\begin{rema}
It is a natural question to determine if the Fubini-Study metric can be defined in a uniform way (for non-Archimedean and Archimedean cases). Let $\overline V=(V,\|\ndot\|)$ be a finite-dimensional vector space over $k$. For any point $x\in X^{\mathrm{an}}$, we let $|\ndot|_{\overline V,\varepsilon}(x)$ and $|\ndot|_{\overline V,\pi}(x)$ be the norms on $V\otimes_k\widehat{\kappa}(x)$ induced by $\|\ndot\|$ by $\varepsilon$-extension and $\pi$-extension of scalars respectively. If the absolute value $|\ndot|$ is non-Archimedean, then both norms $|\ndot|_{\overline V,\varepsilon}(x)$ and $|\ndot|_{\overline V,\pi}(x)$ induce the same quotient norm on $\mathcal O_V(1)(x)$. In fact, by Proposition~\ref{Pro:comparisonofdualnormes:scalar:extension} \ref{Pro:comparisonofdualnormes:scalar:extension:epsilon}, \ref{Pro:comparisonofdualnormes:scalar:extension:epsilon:pi}, the dual norms of both $|\ndot|_{\overline V,\varepsilon}$ and $|\ndot|_{\overline V,\pi}$ identify with $\|\ndot\|_{*,\widehat{\kappa}(x),\varepsilon}$, and hence induce the same restricted norm on $\mathcal O_V(1)(x)^\vee$. By Proposition \ref{Pro:dualquotient},  the dual norms of the quotient norms on $\mathcal O_V(1)(x)$ of $|\ndot|_{\overline V,\varepsilon}$ and $|\ndot|_{\overline V,\pi}$ are the same. Since $\mathcal O_V(1)(x)$ is a vector space of rank $1$ on $\widehat{\kappa}(x)$, we obtain that these quotient norms are the same. In other words, in both the Archimedean and non-Archimedean cases, we may use the $\pi$-extension of scalars to define the Fubini-Study metric. However, for the reason of applications in the study of adelic vector bundles, it is more convenient to consider the $\varepsilon$-extension of scalars for the non-Archimedean case. We emphasis however that, in the Archimedean case, if we apply  the $\varepsilon$-extension of scalars instead of the $\pi$-extension of scalars, in general we obtain a different metric from the Fubini-Study metric.
\end{rema}

\begin{prop}\label{prop:FS:metric:continuous}
Let $\overline{V} = (V,\|\ndot\|)$ be a finite-dimensional normed vector space over $k$. Then the norms $|\ndot|_{\overline{V},\mathrm{FS}}(x)$, $x\in \mathbb P(V)^{\mathrm{an}}$ described above define a continuous metric on the universal invertible sheaf $\mathcal O_V(1)$.
\end{prop} 
\begin{proof} By Proposition \ref{Pro:doubedualandquotient} (see also Remark \ref{Rem:extensiondoubledual}), for any $x\in \mathbb P(V)^{\mathrm{an}}$, the norms $\|\ndot\|$ and $\|\ndot\|_{**}$ induce the same Fubini-Study norm on $\mathcal O_V(1)(x)$. Hence we may assume without loss of generality that the norm $\|\ndot\|$ is ultrametric when $(k,|\ndot|)$ is non-Archimedean.

For any $x\in \mathbb P(V)^{\mathrm{an}}$, let $|\ndot|_{\overline{V}}(x)$ be the norm on $V\otimes_k\widehat{\kappa}(x)$ induced by $\|\ndot\|$ by $\pi$-extension of scalars, and let $|\ndot|_{\overline{V}}(x)_*$ be the dual norm of $|\ndot|_{\overline{V}}(x)$. The norms $|\ndot|_{\overline{V}}(x)_*$ define a metric $\varphi$ on $\pi^*(V^\vee)$. By Proposition~\ref{Pro:comparisonofdualnormes:scalar:extension} \ref{Pro:comparisonofdualnormes:scalar:extension:epsilon}, \ref{Pro:comparisonofdualnormes:scalar:extension:epsilon:pi}, the norm $|\ndot|_{\overline{V}}(x)_*$ coincides with the norm induced by $\|\ndot\|_*$ by $\varepsilon$-extension of scalars. Therefore, by Proposition \ref{Pro:pullback}, we obtain that the metric $\varphi$ is continuous. 

The dual norm of the Fubini-Study norm $|\ndot|_{\overline{V},\mathrm{FS}}(x)$ then coincides with the restriction of $|\ndot|_{\overline{V}}(x)_*$ {to} $\mathcal O_V(1)^\vee\otimes\widehat{\kappa}(x)$ by Proposition~\ref{Pro:dualquotient}. Hence these dual norms (for $x\in \mathbb P(V)^{\mathrm{an}}$) form a continuous metric on $\mathcal O_V(1)^\vee$. Therefore the Fubini-Study norms $|\ndot|_{\overline{V},\mathrm{FS}}(x)$, $x\in \mathbb P(V)^{\mathrm{an}}$ define a continuous metric on $\mathcal O_V(1)$ (see Proposition \ref{Pro:dualcontinuousmetric}).
\end{proof}

\begin{defi}\label{Def:Fubini-Studymetric}
Let $\overline{V} = (V,\|\ndot\|)$ be a finite-dimensional normed vector space over $k$. The continuous metric on $\mathcal O_V(1)$ formed by the Fubini-Study norms $|\ndot|_{\overline{V},\mathrm{FS}}(x)$ with $x\in \mathbb P(V)^{\mathrm{an}}$ is called the \emph{Fubini-Study metric}\index{metric!Fubini-Study metric}\index{Fubini-Study metric} on $\mathcal O_V(1)$ associated with the norm $\|\ndot\|$ on $V$.
\end{defi}

{
\begin{rema}\label{Rem: sup norm of Fubini-Study}
Let $\overline V=(V,\|\ndot\|)$ be a normed vector space over $k$ and $s$ be an element of $V$. For any $x\in\mathbb P(V)^{\mathrm{an}}$ such that $s(x)\neq 0$, by definition one has 
\[|s|_{\overline V,\mathrm{FS}}(x)=\inf_{
\begin{subarray}{c}t\in V\otimes_k\widehat{\kappa}(x)\\
t(x)=s(x)
\end{subarray}}\|t\|_{\widehat{\kappa}(x),\natural}\]
with $\natural=\varepsilon$ if $|\ndot|$ is non-Archimedean, and $\natural=\pi$ if $|\ndot|$ is Archimedean.
In particular, one has (see Proposition \ref{Pro:extensiondecorps}) 
\begin{equation}\label{Equ: norm s(x) bound}|s|_{\overline V,\mathrm{FS}}(x)\leqslant \|s\|_{\widehat{\kappa}(x),\natural}=\|s\|_{**}.\end{equation}
Moreover, any rational point $y\in\mathbb P(V)(k)$ corresponds to a non-zero element $\beta_y:V\rightarrow k$ in the dual vector space $V^\vee$. The dual norm of $\beta_y$ identifies with the inverse of the quotient norm of $1\in k$. Therefore one has
\[|s|_{\overline V,\mathrm{FS}}(y)=\frac{|\beta_y(s)|}{\|\beta_y\|_*}\]
and hence
\[\sup_{y\in\mathbb P(V)(k)}|s(y)|_{\overline V,\mathrm{FS}}(y)=\sup_{y\in\mathbb P(V)(k)}\frac{|\beta_y(s)|}{\|\beta_y\|_*}=\|s\|_{**}.\]
Combing with \eqref{Equ: norm s(x) bound}, we obtain $\|\ndot\|_{\overline V,\mathrm{FS}}=\|\ndot\|_{**}$.

Let $(k',|\ndot|)$ be a valued extension of  $(k,|\ndot|)$. Note that the fibre product $\mathbb P(V)\times_{\Spec k}\Spec k'$ identifies with the projective space of $V':=V\otimes_kk'$. The Fubini-Study metric on $\mathcal O_V(1)$ induces by base change a continuous metric on $\mathcal O_{V'}(1)$ which we denote by $\{|\ndot|_{\overline V,\mathrm{FS},k'}(x)\}_{x\in\mathbb P(V')^{\mathrm{an}}}$. By Corollary \ref{Cor:extsucc} and Proposition \ref{Pro:quotientnormextensionscalar}, we obtain that this metric coincides with the Fubini-Study metric associated with the norm $\|\ndot\|_{k',\natural}$ on $V'$, where $\natural=\varepsilon$ if $|\ndot|$ is non-Archimedean and $\natural=\pi$ if $|\ndot|$ is Archimedean. In particular, one has
\begin{equation}\|\ndot\|_{\overline V,\mathrm{FS},k'}=\|\ndot\|_{\overline{V}',\mathrm{FS}}.\end{equation}  
\end{rema}}

\begin{defi}\label{Def: Quotient metric}Let $f:X\rightarrow\Spec k$ be a  $k$-scheme and $L$ be an invertible $\mathcal O_X$-module. Suppose given a finite-dimensional vector space $V$ over $k$ and a surjective {$\mathcal O_X$-homomorphism} $\beta:f^*(V)\rightarrow L$. Then the homomorphism $\beta$ corresponds to a $k$-morphism of schemes $g:X\rightarrow\mathbb P(V)$ such that $g^*(\mathcal O_V(1))$ is canonically isomorphic to $L$. If $V$ is equipped with a norm $\|\ndot\|$, then the Fubini-Study metric on $\mathcal O_V(1)$ induces by pull-back a continuous metric on $L$, called the \emph{quotient metric}\index{quotient metric}\index{metric!quotient ---} induced by the normed vector space $(V,\|\ndot\|)$ and the surjective homomorphism $\beta$.
\end{defi}  

\begin{defi}\label{Def: rightarrow ell}
Let $\ell$ be a section of $L$ over a Zariski open set $U$, which trivialises the invertible sheaf $L$.
The section $\ell$ yields the isomorphism $\iota : \mathcal O_U \to L|_U$ given by
$a \mapsto a\ell$. We define $\vec{\ell} : f^*(V)|_U \to \mathcal O_U$
by $\iota^{-1} \circ \beta_U$, that is, the following diagram is commutative:
\[
\xymatrix{ f^{*}(V)|_U \ar[r]^{\vec{\ell}} \ar[rd]_{\beta_U} & \mathcal O_U \ar[d]^{\iota} \\
& L|_U}
\]
\end{defi}
If $\{e_i\}_{i=1}^r$ is a basis of $V$, $\{e_i^{\vee}\}_{i=1}^r$ is the dual basis
of $\{e_i\}_{i=1}^r$ and $\beta_U(e_i) = a_i \ell$ for $i\in\{1, \ldots, n\}$, 
then $\vec{\ell}$ is given by
\[
\vec{\ell} = a_{1} e_1^{\vee} + \cdots + a_{r} e_r^{\vee}.
\]
For each $x \in U^{\mathrm{an}}$, the evaluation of $\vec{\ell}$ at $x$ is denoted by $\vec{\ell}_x$,
that is, \[\vec{\ell}_x = a_{1}(x) e_1^{\vee} + \cdots + a_{r}(x) e_r^{\vee}\in \mathrm{Hom}_{\widehat{\kappa}(x)}(V_{\widehat{\kappa}(x)},\widehat{\kappa}(x)).\]

\begin{prop}\label{prop:FS:metric:evaluation}
Let $f:X\rightarrow\Spec k$ be a k-scheme and $L$ be an invertible $\mathcal O_X$-module. Suppose given a finite-dimensional normed vector space $(V,\|\ndot\|)$ over $k$ and a surjective $\mathcal O_X$-homomorphism $\beta:f^*(V)\rightarrow L$. Let $\varphi$ be the quotient metric induced by $(V,\|\ndot\|)$ and $\beta$. For any section $\ell$ of $L$ on a Zariski open subset $U$ of $X$ which trivialises $L$ on $U$, one has
\begin{equation}\label{eqn:prop:FS:metric:evaluation:01}\forall\,x\in U^{\mathrm{an}},\quad 
|\ell|_{\varphi}(x)=\| \vec{\ell}_x \|_{\widehat{\kappa}(x),\natural,*}^{-1} = \| \vec{\ell}_x \|_{*,\widehat{\kappa}(x),\varepsilon}^{-1},
\end{equation}
where $\natural=\varepsilon$ if $|\ndot|$ is non-Archimedean and $\natural=\pi$ if $|\ndot|$ is Archimedean.
\end{prop}

\begin{proof}
For each $x \in U^{\mathrm{an}}$, one has the following commutative diagram:
\[
\xymatrix{ V_{\widehat{\kappa}(x)} \ar[r]^{\vec{\ell}_x} \ar[rd]_{\beta_x} & \widehat{\kappa}(x) \ar[d]^{\iota_x} \\
& L \otimes_{\mathcal O_U} \widehat{\kappa}(x)}
\]
By (1) in Lemma~\ref{lem:quotient:dim:1:comm}, the operator norm of $ \beta_x$ is $1$.
Moreover, the operator norm of $ \vec{\ell}_x $ is $\| \vec{\ell}_x \|_{\widehat{\kappa}(x),\natural,*}$ and that of 
$ \iota_x $ is $|\ell|_{\varphi}(x)$. As the operator norm of $\beta_x$ is the product of the operator norms of $\vec{\ell}_x$ and $\iota_x$
by (2) in Proposition~\ref{lem:quotient:dim:1:comm}, we obtain $|\ell|_\varphi(x)=\|\vec{\ell}_x\|_{\widehat{\kappa}(x),\natural,*}$.
The equality $\| \vec{\ell}_x \|_{\widehat{\kappa}(x),\natural,*}^{-1}=\|\vec{\ell}_x\|_{*,\widehat{\kappa}(x),\varepsilon}$  follows from Proposition~\ref{Pro:comparisonofdualnormes:scalar:extension} \ref{Pro:comparisonofdualnormes:scalar:extension:epsilon}, \ref{Pro:comparisonofdualnormes:scalar:extension:epsilon:pi} for the non-Archimedean and Archimedean cases, respectively.
\end{proof}

\begin{rema}\label{Rem:existenceofcontmetr}
Let $X$ be a quasi-projective $k$-scheme. The above construction shows that any ample invertible $\mathcal O_X$-module admits a continuous metric. By Proposition \ref{Pro:dualcontinuousmetric}, we deduce that, more generally, any invertible $\mathcal O_X$-module admits a continuous metric.
\end{rema}

{
\begin{rema}\label{Rem: powers of quotient}
Let $f:X\rightarrow\Spec k$ be a $k$-scheme, $L$ be an invertible $\mathcal O_X$-module, $V$ be a finite-dimensional vector space over $k$, and $\beta:f^*(V)\rightarrow L$ be a surjective $\mathcal O_X$-homomorphism. Let $n$ be an integer, $n\geqslant 1$. Then $\beta$ induces a surjective $\mathcal O_X$-homomorphism $\beta^{\otimes n}:f^*(V^{\otimes n})\rightarrow L^{\otimes n}$. Let $\|\ndot\|$ be a norm on $V$ and $\varphi$ be the quotient metric induced by $(V,\|\ndot\|)$ and $\beta$. We claim that $n\varphi$ is the quotient metric induced by $(V^{\otimes n},\|\ndot\|_{\natural})$ and $\beta^{\otimes n}$, where $\|\ndot\|_{\natural}$ denotes the $\natural$-tensor power of the norm $\|\ndot\|$, $\natural=\varepsilon$ if $|\ndot|$ is non-Archimedean and $\natural=\pi$ if $|\ndot|$ is Archimedean. In fact, if $x$ is an point in $X^{\mathrm{an}}$ and $\ell$ is a non-zero element of $L\otimes\widehat{\kappa}(x)$, then $\ell^{\otimes n}$ is a non-zero element of $L^{\otimes n}\otimes\widehat{\kappa}(x)$. By Propositions \ref{Pro:epstensoextcorps} and \ref{Pro: produit pi tensoriel extension}, the norm $\|\ndot\|_{\natural,\widehat{\kappa}(x),\natural}$ on $V^{\otimes n}\otimes\widehat{\kappa}(x)$ coincides with the $\natural$-tensor power of $\|\ndot\|_{\widehat{\kappa}(x),\natural}$.
Consider the dual homomorphism 
\[\beta_x^{\vee\otimes n}:(L^{\vee}\otimes\widehat{\kappa}(x))^{\otimes n}\longrightarrow (V^{\vee}\otimes\widehat{\kappa}(x))^{\otimes n}.\]
By Proposition \ref{Pro:dualitypiepsilon} and Corollary \ref{Cor: dual tensor product}, the dual norm $\|\ndot\|_{\natural,\widehat{\kappa}(x),\natural,*}$ coincides with the $\varepsilon$-tensor power of $\|\ndot\|_{\widehat{\kappa}(x),\natural,*}=\|\ndot\|_{*,\widehat{\kappa}(x),\varepsilon}$ (see Proposition~\ref{Pro:comparisonofdualnormes:scalar:extension} \ref{Pro:comparisonofdualnormes:scalar:extension:epsilon},\ref{Pro:comparisonofdualnormes:scalar:extension:epsilon:pi}). In both cases one has
\[\|\beta_x^{\vee\otimes n}(\ell^{\vee\otimes n})\|_{\natural,\widehat{\kappa}(x),\natural,*}=\|\beta_x^{\vee}(\ell^\vee)\|_{\widehat{\kappa}(x),\natural,*}^n=|\ell^{\otimes n}|_{n\varphi}(x)^{-1},\]
where the first equality comes from Remark \ref{Rem:produittenrk1}, and the second comes from Proposition \ref{Pro:dualquotient}.
\end{rema}}

\begin{prop}\label{Pro:distrancequot}
Let $f:X\rightarrow\Spec k$ be a scheme over $\Spec k$ and $L$ be an invertible $\mathcal O_X$-module. Let $V$ be a finite-dimensional vector space and $\beta:f^*(V)\rightarrow L$ be a surjective homomorphism. If $\|\ndot\|$ and $\|\ndot\|'$ are two norms on $V$ and if $\varphi$ and $\varphi'$ are quotient metrics on $L$ induced by $\overline{V} = (V,\|\ndot\|)$ and $\overline{V}{}' = (V,\|\ndot\|')$ (and the surjective homomorphism $\beta$) respectively, then one has $d(\varphi,\varphi')\leqslant d(\|\ndot\|,\|\ndot\|')$.
\end{prop}
\begin{proof}
Let $x$ be a point of $X^{\mathrm{an}}$, $|\ndot|_{\overline{V}}(x)$ and $|\ndot|_{\overline{V}'}(x)$ be the norms on $V\otimes\widehat{\kappa}(x)$ induced by $\|\ndot\|$ and $\|\ndot\|'$ by extension of scalars. Proposition \ref{Pro:distanceextension} leads to \[d(|\ndot|_{\overline{V}}(x),|\ndot|_{\overline{V}'}(x))\leqslant d(\|\ndot\|,\|\ndot\|').\] Since $|\ndot|_{\varphi}(x)$ and $|\ndot|_{\varphi'}(x)$ are respectively the quotient norms of $|\ndot|_{\overline{V}}(x)$ and $|\ndot|_{\overline{V}'}(x)$, by Proposition \ref{Pro:distancequotientandsub} one has $d(|\ndot|_{\varphi}(x),|\ndot|_{\varphi'}(x))\leqslant d(|\ndot|_{\overline{V}}(x),|\ndot|_{\overline{V}'}(x))$. Therefore
\[d(\varphi,\varphi')=\sup_{x\in X^{\mathrm{an}}}d(|\ndot|_{\varphi}(x),|\ndot|_{\varphi'}(x)) \leqslant \sup_{x\in X^{\mathrm{an}}}d(|\ndot|_{\overline{V}}(x),|\ndot|_{\overline{V}'}(x))\leqslant 
d(\|\ndot\|,\|\ndot\|'),\]
as required.
\end{proof}

\begin{defi}\label{Def:sequenceofmetrics}
Let $\pi:X\rightarrow\Spec(k)$ be a projective $k$-scheme, $L$ be an invertible $\mathcal O_X$-module, which is generated by global sections, and $\varphi$ be a continuous metric on $L$.
By Proposition~\ref{prop:def:sup:norm},
for each $x \in X^{\mathrm{an}}$, the homomorphism 
$(H^0(X, L)/\mathcal N(X, L)) \otimes \widehat{\kappa}(x) \to L \otimes \widehat{\kappa}(x)$
induced by $H^0(X, L)  \otimes \widehat{\kappa}(x) \to L \otimes \widehat{\kappa}(x)$
is surjective, so that one has a quotient norm $|\ndot|'(x)$ on $L \otimes \widehat{\kappa}(x)$ induced by
$\|\ndot\|^{\sim}_{\varphi}$.
The family $\{ |\ndot|'(x) \}_{x \in X^{\mathrm{an}}}$ of metrics is denoted by $\varphi_{\mathrm{FS}}$, that is,
$|\ndot|_{\varphi_{\mathrm{FS}}}(x) := |\ndot|'(x)$ for $x \in X^{\mathrm{an}}$, called the \emph{Fubini-Study metric associated with $\varphi$}\index{Fubini-Study metric}\index{metric!Fubini-Study ---}.
Let $X_{\mathrm{red}}$ be the reduced scheme associated with $X$ and $L_{\mathrm{red}} := L \otimes_{\mathcal O_X} \mathcal O_{X_{\mathrm{red}}}$. Let $V$ be the image of $H^0(X, L) \to H^0(X_{\mathrm{red}}, L_{\mathrm{red}})$.
Then $\|\ndot\|^{\sim}_{\varphi}$ is a norm of $V$ and $\varphi_{\mathrm{FS}}$ is the quotient metric
induced by the surjection $V \otimes \mathcal O_{X_{\mathrm{red}}} \to L_{\mathrm{red}}$ and
$\|\ndot\|^{\sim}_{\varphi}$. In particular, $\varphi_{\mathrm{FS}}$ is continuous.
For an integer $n \geqslant 1$, we set $\varphi_n = (n\varphi)_{\mathrm{FS}}$.

By Propositions \ref{Pro:quotientderang1norm} and \ref{Pro: quotient metric pi ext2}, for any point $x\in X^{\mathrm{an}}$ and any non-zero element $\ell\in L\otimes_{\mathcal O_X}\widehat{\kappa}(x)$, one has 
\begin{equation}\label{Equ: local norm by the global one}|\ell|_{\varphi_{\mathrm{FS}}}(x) =\inf_{\begin{subarray}{c}
s\in H^0(X,L),\,\lambda\in\widehat{\kappa}(x)^{\times}\\
s(x)=\lambda\ell
\end{subarray}}|\lambda|_x^{-1}\cdot\|s\|_{\varphi}.\end{equation}
This equality is fundamental in the study of quotient metrics.
\end{defi}

\begin{prop}\label{Pro:metricinduitparnormsup}
Let $\pi:X\rightarrow\Spec k$ be a projective $k$-scheme and $L$ be an invertible $\mathcal O_X$-module generated by global sections, equipped with a continuous metric $\varphi$. Then the following assertions hold.
\begin{enumerate}[label=\rm(\arabic*)]
\item\label{item: varphi n bounded from below} For any integer $n\in\mathbb N_{\geqslant 1}$, one has $\varphi_n\geqslant n\varphi$, where $\varphi_n$ denotes the Fubini-Study metric associated with $n\varphi$.
\item\label{item: same sup norm} The sup seminorm $\|\ndot\|_{\varphi_n}$ on $V_n:=H^0(X,L^{\otimes n})$ induced by $\varphi_n$ coincides with $\|\ndot\|_{n\varphi}$.
\item\label{item: subadditivity} Let $M$ be another invertible $\mathcal O_X$-module generated by global sections, equipped with a continuous metric $\psi$.
Then $(\varphi + \psi)_{\mathrm{FS}} \leqslant \varphi_{\mathrm{FS}} + \psi_{\mathrm{FS}}$. In particular, for any pair $(m,n)$ of positive integers one has $\varphi_{n+m}\leqslant\varphi_n+\varphi_m$.
\item\label{Item: distance varphi n}  For any integer $n\geqslant 1$, one has $d(\varphi_n,n\varphi)\leqslant nd(\varphi_1,\varphi)$. In particular, if $\varphi_1=\varphi$ then $\varphi_n=n\varphi$ for any $n\in\mathbb N_{\geqslant 1}$.

\item\label{item: change of metrics} Let $\varphi'$ be another continuous metric on $L$. Then one has $d(\varphi_n,\varphi_n')\leqslant nd(\varphi,\varphi')$ for any $n\in\mathbb N_{\geqslant 1}$.
\end{enumerate}
\end{prop}
\begin{proof} For any $n\in\mathbb N_{\geqslant 1}$, we denote by $V_n$ the vector space $H^0(X,L^{\otimes n})$ over $k$.

\ref{item: varphi n bounded from below} Let $x$ be a point of $X^{\mathrm{an}}$ and $\ell$ be an element of $L^{\otimes n}\otimes\widehat{\kappa}(x)$. 
Note that, for any $s \in V_n$, one has $|s|_{n\varphi}(x) \leqslant \| s \|_{n\varphi}$, so that,
by \eqref{Equ: local norm by the global one}, one obtains
\[|\ell|_{\varphi_n}(x)=\inf_{\begin{subarray}{c}s\in V_n,\,\lambda\in\widehat{\kappa}(x)^{\times},\\ 
s(x)=\lambda\ell\end{subarray}}|\lambda|_x^{-1}\cdot\|s\|_{n\varphi} \geqslant
\inf_{\begin{subarray}{c}s\in V_n,\,\lambda\in\widehat{\kappa}(x)^{\times},\\ 
s(x)=\lambda\ell\end{subarray}}|\lambda|_x^{-1}\cdot|s|_{n\varphi}(x) = |\ell|_{n\varphi}(x),
\] 
as desired.

\medskip\ref{item: same sup norm} By \ref{item: varphi n bounded from below}, one has $\|\ndot\|_{\varphi_n}\geqslant\|\ndot\|_{n\varphi}$. In the following, we prove the converse inequality. If $s$ is a global section of $L^{\otimes n}$, for any $x\in X^{\mathrm{an}}$, by \eqref{Equ: local norm by the global one}, one has 
\[|s|_{\varphi_n}(x)=\inf_{\begin{subarray}{c}t\in V_n,\,\lambda\in\widehat{\kappa}(x)^\times\\
t(x)=\lambda s(x)\end{subarray}}|\lambda|_x^{-1}\cdot\|t\|_{n\varphi}\leqslant\|s\|_{n\varphi}.\]
Hence
$\|s\|_{\varphi_n}=\sup_{x\in X^{\mathrm{an}}}|s|_{\varphi_n}(x)\leqslant\|s\|_{n\varphi}$.

\medskip\ref{item: subadditivity} Let $x$ be a point of $X^{\mathrm{an}}$, $\ell$ and $\ell'$ be elements of $L\otimes\widehat{\kappa}(x)$ and $M\otimes\widehat{\kappa}(x)$ respectively. By \eqref{Equ: local norm by the global one} together with
\eqref{Equ:sousmultiplicativesec} one has
{\allowdisplaybreaks
\begin{align*}
|\ell \cdot\ell'|_{(\varphi + \psi)_{\mathrm{FS}}}(x) & =\inf_{\begin{subarray}{c}s\in H^0(X,L\otimes M),\,\lambda\in\widehat{\kappa}(x)^{\times}\\ s(x)=\lambda\ell\cdot\ell'\end{subarray}}|\lambda|_x^{-1}\cdot\|s\|_{\varphi+\psi}\\
&\leqslant\inf_{\begin{subarray}{c}(t,t')\in H^0(X,L)\times H^0(X,M)\\
(\mu,\eta)\in(\widehat{\kappa}(x)^{\times})^2\\
t(x)=\mu\ell,\,t'(x)=\eta\ell'
\end{subarray}}|\mu\eta|_x^{-1}\cdot\|t \cdot t'\|_{\varphi+\psi}\\
&\leqslant\inf_{\begin{subarray}{c}(t,t')\in H^0(X,L)\times H^0(X,M)\\
(\mu,\eta)\in(\widehat{\kappa}(x)^{\times})^2\\
t(x)=\mu\ell,\,t'(x)=\eta\ell'
\end{subarray}}\left(|\mu|_x^{-1}\cdot\|t\|_{\varphi}\right)\left(|\eta|_x^{-1}\cdot\|t'\|_{\psi}\right)\\
& = |\ell|_{\varphi_{\mathrm{FS}}}(x)\cdot|\ell'|_{\psi_{\mathrm{FS}}}(x),
\end{align*}}
For the last assertion, note that 
\[
\varphi_{n+m} = (n\varphi + m\varphi)_{\mathrm{FS}} \leqslant (n\varphi)_{\mathrm{FS}} + (m\varphi)_{\mathrm{FS}} = \varphi_n + \varphi_m.
\]

\medskip\ref{Item: distance varphi n} By \ref{item: subadditivity}, one has $\varphi_n\leqslant n\varphi_1$. Moreover, by \ref{item: varphi n bounded from below}, one has $\varphi_n\geqslant n\varphi$ and $\varphi_1\geqslant \varphi$. Hence
$0\leqslant\varphi_n-n\varphi\leqslant n\varphi_1-n\varphi=n(\varphi_1-\varphi)$, 
which implies
\[d(\varphi_n,n\varphi)=\sup_{x\in X^{\mathrm{an}}}(\varphi_n-n\varphi)(x)\leqslant n\sup_{x\in X^{\mathrm{an}}}(\varphi_1-\varphi)(x)=nd(\varphi_1,\varphi).\]

\medskip\ref{item: change of metrics} 
By Proposition \ref{Pro:distrancequot} together with \eqref{Equ:comparaisondistancesup}, one has 
\[d(\varphi_n,\varphi_n')\leqslant d(\|\ndot\|^{\sim}_{n\varphi},\|\ndot\|^{\sim}_{n\varphi'})
\leqslant d(n\varphi, n\varphi') = nd(\varphi, \varphi'),\]
where the equality comes from \eqref{Equ:distancen}.
\end{proof}

\begin{prop}\label{Pro:positivityofquotientmetric}
Let $\pi:X\rightarrow\Spec k$ be a projective $k$-scheme and $L$ be an invertible $\mathcal O_X$-module. Suppose given a normed vector space $(V,\|\ndot\|)$ and a surjective homomorphisms $\beta:\pi^*(V)\rightarrow L$. Let $\varphi$ be the quotient metric on $L$ induced by $(V,\|\ndot\|)$ and $\beta$. Then,  one has the following:
\begin{enumerate}[label=\rm(\arabic*)]
\item\label{Item: norm of f v bounded by norm of v} Let $f:V\rightarrow H^0(X,L)$ be the adjoint homomorphism of $\beta:\pi^*(V)\rightarrow L$.
Then, $\|f(v)\|_{\varphi}\leqslant\|v\|$ for any $v \in V$.
\item\label{Item: for any n varphi n = n varphi} For any integer $n\geqslant 1$, $\varphi_n=n\varphi$.
\end{enumerate}
\end{prop}
\begin{proof}

\ref{Item: norm of f v bounded by norm of v} By Propositions \ref{Pro:quotientderang1norm} and \ref{Pro: quotient metric pi ext2},  
for $x \in X^{\mathrm{an}}$,
\[|f(v)|_{\varphi}(x)=\inf_{\begin{subarray}{c}t\in V,\,\lambda\in\widehat{\kappa}(x)^{\times}\\ f(t)(x)=\lambda f(v)(x)
\end{subarray}}|\lambda|_x^{-1}\cdot\|t\|\leqslant\|v\|,\]
so that one has \ref{Item: norm of f v bounded by norm of v}.

\medskip
\ref{Item: for any n varphi n = n varphi}
By Proposition \ref{Pro:metricinduitparnormsup} \ref{Item: distance varphi n}, it suffices to verify that $\varphi_1=\varphi$. Note that $\varphi_1\geqslant\varphi$ (by Proposition \ref{Pro:metricinduitparnormsup} \ref{item: varphi n bounded from below}). In the following, we prove the converse inequality. 
Let $x$ be a point of $X^{\mathrm{an}}$ and $\ell$ be an element of $L\otimes\widehat{\kappa}(x)$. 
By \ref{Item: norm of f v bounded by norm of v} and \eqref{Equ: local norm by the global one} together with
Propositions \ref{Pro:quotientderang1norm} and \ref{Pro: quotient metric pi ext2}, 
one has
\begin{align*}|\ell|_{\varphi_1}(x) & =\inf_{\begin{subarray}{c}
s\in H^0(X,L)\\
\lambda\in\widehat{\kappa}(x)^{\times}\\
s(x) =\lambda\ell
\end{subarray}}|\lambda|_x^{-1}\cdot\|s\|_{\varphi}\leqslant\inf_{\begin{subarray}{c}
s'\in V,\,\lambda\in\widehat{\kappa}(x)^{\times}\\
f(s')(x)=\lambda\ell
\end{subarray}}|\lambda|_x^{-1}\cdot\|f(s')\|_{\varphi} \\
& \leqslant \inf_{\begin{subarray}{c}
s'\in V,\,\lambda\in\widehat{\kappa}(x)^{\times}\\
f(s')(x)=\lambda\ell\end{subarray}}|\lambda|_x^{-1}\cdot\|s'\|=|\ell|_{\varphi}(x). \end{align*}
Therefore one has $\varphi_1=\varphi$.
\end{proof}

\section{Semi-positive metrics}

Let $(k,|\ndot|)$ be a complete valued field and $\pi:X\rightarrow\Spec k$  be a projective $k$-scheme. In this section, we discuss positivity conditions of continuous metrics on invertible $\mathcal O_X$-modules.

\subsection{Definition and basic properties} Let $L$ be an invertible $\mathcal O_X$-module equipped with a continuous metric $\varphi$. We assume that $L$ is generated by global sections. We have constructed in Definition \ref{Def:sequenceofmetrics} a sequence of quotient metrics $\{\varphi_n\}_{n\in\mathbb N_{\geqslant 1}}$. By Proposition \ref{Pro:metricinduitparnormsup} (1) and (4), we obtain that  $\{d(\varphi_n,n\varphi)\}_{n\in\mathbb N_{\geqslant 1}}$ is a sub-additive non-negative sequence and the normalised sequence $\{d(\varphi_n,n\varphi)/n\}_{n\in\mathbb N_{\geqslant 1}}$ is bounded from above. Hence the sequence $\{d(\varphi_n,n\varphi)/n\}_{n\in\mathbb N_{\geqslant 1}}$ converges in $\mathbb R_+$. We denote by $\defp(\varphi)$ the limit
\begin{equation}
\defp(\varphi):=\lim_{n\rightarrow+\infty}\frac 1nd(\varphi_n,n\varphi),\end{equation}
called the \emph{default of positivity}\index{default of positivity} of the metric $\varphi$.
By definition, for any integer $m\in\mathbb N_{\geqslant 1}$, one has \begin{equation}\label{Equ:dphomogene}
\defp(m\varphi)=m\defp(\varphi).\end{equation}
We say that the metric $\varphi$ is \emph{semipositive}\index{semipositive}\index{metric!semipositive ---} if one has $\defp(\varphi)=0$.
Clearly, if the metric $\varphi$ is semipositive, then for any integer $m\geqslant 1$, the metric $m\varphi$ on $L^{\otimes m}$ is also semipositive. Conversely, if there is an integer $m\geqslant 1$ such that $m\varphi$ is semipositive, then the metric $\varphi$ is also semipositive.

More generally, we assume that $L$ is semiample (namely a positive tensor power of $L$ is generated by global sections).
Let $n$ be a positive integer such that $L^{\otimes n}$ is generated by global sections.
The quantity $\defp(n\varphi)/n$ does not depend on the choice of $n$ by \eqref{Equ:dphomogene},
so that we define $\defp(\varphi)$ to be $\defp(n\varphi)/n$.
It is easy to see that \eqref{Equ:dphomogene} still holds under the assumption that $L$ is semiample.
We say that $\varphi$ is \emph{semipositive}\index{semipositive}\index{metric!semipositive ---} if $\defp(\varphi) = 0$.

\begin{rema}
Let $(V,\|\ndot\|)$ be a normed vector space of finite rank over $k$ and $\beta:\pi^*(V)\rightarrow L$ be
a surjective homomorphism. Let $\varphi$ be the quotient metric on $L$ induced by $(V,\|\ndot\|)$ and $\beta$. Then,
by Proposition~\ref{Pro:positivityofquotientmetric}, $\varphi$ is semipositive.
\end{rema}

\begin{prop}\label{Pro:puissancepositve}
Let $L$ be a semiample invertible $\mathcal O_X$-module, equipped with a continuous metric $\varphi$. If $\varphi$ is semipositive, then $n\varphi$ is semipositive for any $n\in\mathbb N_{\geqslant 1}$. Conversely, if there exists an integer $n\in\mathbb N_{\geqslant 1}$ such that $n\varphi$ is semipositive, then the metric $\varphi$ is also semipositive.
\end{prop}

The following proposition shows that semipositive metrics form a closed subset in the topological space of continuous metrics.

\begin{prop}\label{Pro:limitsemipostive}
Let $L$ be a semiample invertible $\mathcal O_X$-module, equipped with a continuous metric $\varphi$. Suppose that there is a sequence of semipositive metrics $\{\varphi^{(m)}\}_{n\in\mathbb N}$ on $L$ such that \[\lim_{m\rightarrow+\infty}d(\varphi^{(m)},\varphi)=0.\] Then the metric $\varphi$ is also semipositive.
\end{prop}
\begin{proof}
For any integer $p\geqslant 1$, one has \[d(p\varphi^{(m)},p\varphi)=pd(\varphi^{(m)},\varphi).\] Therefore, by replacing $L$ by a certain tensor power $L^{\otimes p}$ and $\varphi$ by $p\varphi$, we may assume without loss of generality that $L$ is generated by global sections. Thus the metrics $\varphi^{(m)}_n$ and $\varphi_n$ are well defined for any $m\in\mathbb N$ and any $n\in\mathbb N_{\geqslant 1}$. Moreover, by Proposition \ref{Pro:metricinduitparnormsup} (5) we obtain that \[d(\varphi^{(m)}_n,\varphi_n)\leqslant nd(\varphi^{(m)},\varphi).\] Note that for $m\in\mathbb N$ and $n\in\mathbb N_{\geqslant 1}$ one has
\[\begin{split}d(\varphi_n,n\varphi)&\leqslant d(\varphi_n,\varphi_n^{(m)})+d(\varphi^{(m)}_n,n\varphi^{(m)})+d(n\varphi^{(m)},n\varphi)\\
&\leqslant
d(\varphi^{(m)}_n,n\varphi^{(m)})+2d(\varphi^{(m)},\varphi).
\end{split}\] 
By taking the limit when $n$ tends to the infinity, we obtain
\[\defp(\varphi)\leqslant 2d(\varphi^{(m)},\varphi)+\defp(\varphi^{(m)})=2d(\varphi^{(m)},\varphi),\]
where the equality comes from the hypothesis that the metrics $\varphi^{(m)}$ are semipositive. By taking the limit when $m$ tends to the infinity, we obtain the semipositivity of the metric $\varphi$.
\end{proof}

\begin{rema}
Proposition \ref{Pro:positivityofquotientmetric} shows that quotient metrics on an invertible $\mathcal O_X$-module are semipositive. Let $L$ be a semiample invertible $\mathcal O_X$-module and $\varphi$ be a continuous metric on $L$. For any $n\in\mathbb N_{\geqslant 1}$ such that $L^{\otimes n}$ is generated by global sections, let $\varphi^{(n)}$ be a continuous metric on $L$ such that $n\varphi^{(n)}$ is a quotient metric. If
$\lim_{n\rightarrow+\infty}d(\varphi^{(n)},\varphi)=0$, then 
the metric $\varphi$ is semipositive.
This is a consequence of Propositions \ref{Pro:puissancepositve} and \ref{Pro:limitsemipostive}. 
\end{rema}

\begin{prop}\label{Pro:subadditivedelta}
Let $X$ be a projective $k$-scheme, $L$ and $L'$ be semiample invertible $\mathcal O_X$-modules, equipped with continuous metrics $\varphi$ and $\varphi'$, respectively. One has $\defp(\varphi+\varphi')\leqslant\defp(\varphi)+\defp(\varphi')$. In particular, if both metrics $\varphi$ and $\varphi'$ are semipositive, then the  metric $\varphi+\varphi'$ on the tensor product $L\otimes L'$ is also semipositive.
\end{prop}
\begin{proof}
By \eqref{Equ:dphomogene}, we may assume that $L$ and $L'$ are generated by global sections.
For any integer $n\geqslant 1$, one has a natural $k$-linear homomorphism \[H^0(X,L^{\otimes n})\otimes H^0(X,L'{}^{\otimes n})\longrightarrow H^0(X,(L\otimes L')^{\otimes n})\]
given by the tensor product.
Moreover, by Proposition~\ref{Pro:sousmultiplicativesec},
for $s\in H^0(X,L^{\otimes n})$ and $s'\in H^0(X,L'{}^{\otimes n})$ one has
\[\|ss'\|_{n(\varphi+\varphi')}\leqslant\|s\|_{n\varphi}\cdot\|s'\|_{n\varphi'}.\]
By 
Proposition \ref{Pro:metricinduitparnormsup} (3), we obtain $(\varphi+\varphi')_n\leqslant\varphi_n+\varphi'_n$ and hence 
\[d((\varphi+\varphi')_n,n(\varphi+\varphi'))\leqslant d(\varphi_n,n\varphi)+d(\varphi_n',n\varphi').\]
Dividing the two sides of the inequality by $n$, by passing to limit when $n$ tends to the infinity we obtain $\defp(\varphi+\varphi')\leqslant\defp(\varphi)+\defp(\varphi')$. 
\end{proof}

\begin{prop}\label{prop:equivalent:semiample:metrized}
Let $L$ be a semiample invertible $\mathcal O_X$-module, equipped with a continuous metric $\varphi$. 
Then the following are equivalent:
\begin{enumerate}[label=\rm(\arabic*)]
\item
The metric $\varphi$ is semipositive.

\item
For any $\epsilon > 0$, there is a positive integer $n$ such that, for all $x \in X^{\mathrm{an}}$,
we can find $s \in H^0(X, L^{\otimes n})_{\widehat{\kappa}(x)} \setminus \{ 0 \}$ with
$\| s \|_{n\varphi,\widehat{\kappa}(x)} \leqslant \mathrm{e}^{n \epsilon} |s|_{n\varphi}(x)$.
\end{enumerate}
\end{prop}

\begin{proof}
(1) $\Longrightarrow$ (2):
By our assumption, there is a positive integer $n$ such that
\[
|\ndot|_{n\varphi}(x) \leqslant |\ndot|_{\varphi_n}(x) \leqslant \mathrm{e}^{n\epsilon/2} |\ndot|_{n\varphi}(x)
\]
for all $x \in X^{\mathrm{an}}$. Moreover, there is an $s \in H^0(X, L^{\otimes n})_{\widehat{\kappa}(x)} \setminus \{ 0 \}$
such that $\|s \|_{n\varphi, \widehat{\kappa}(x)} \leqslant \mathrm{e}^{n\epsilon/2} |s|_{\varphi_n}(x)$. Therefore,
\[
\|s \|_{n\varphi, \widehat{\kappa}(x)} \leqslant \mathrm{e}^{n\epsilon/2} |s|_{\varphi_n}(x) \leqslant \mathrm{e}^{n\epsilon}|s|_{n\varphi}(x).
\]

\medskip
(2) $\Longrightarrow$ (1): For a positive integer $m$, there is a positive integer $a_m$ such that,
for any $x \in X^{\mathrm{an}}$, we can find $s \in H^0(X, L^{\otimes a_m})_{\widehat{\kappa}(x)} \setminus \{ 0 \}$ with
$\| s \|_{a_m \varphi, \widehat{\kappa}(x)} \leqslant \mathrm{e}^{a_m/m} |s|_{a_m \varphi}(x)$. Note that
\[
|s|_{a_m \varphi}(x) \leqslant |s|_{\varphi_{a_m}}(x) \leqslant  \mathrm{e}^{a_m/m} |s|_{a_m \varphi}(x),
\]
which implies that
\[
0 \leqslant \frac{1}{a_m} \left( \ln |\ndot|_{\varphi_{a_m}}(x) - \ln |\ndot|_{a_m \varphi}(x) \right) \leqslant \frac{1}{m}
\]
for all $x \in X^{\mathrm{an}}$, so that $\defp(\varphi)=0$.
\end{proof}

\begin{theo}\label{thm:semipositive:positive:current}
Let $X$ be an irreducible and reduced projective scheme over $\Spec\mathbb C$, $L$ be a semiample invertible $\mathcal O_X$-module and
$\varphi$ be a continuous metric of $L$. Then the following are equivalent:
\begin{enumerate}[label=\rm(\arabic*)]
\item\label{Item: first chern current positive}
The first Chern current $c_1(L, \varphi)$ is positive.

\item\label{Item: uniform extension on points}
For any positive number $\epsilon > 0$, there is a positive integer $n$ such that,
for all $x \in X$, we can find $s \in H^0(X, L^{\otimes n}) \setminus \{ 0 \}$ with
$\| s\|_{n \varphi} \leqslant \mathrm{e}^{n\epsilon} |s|_{n\varphi}(x)$.

\item\label{Item: metric semipositive}
The metric $\varphi$ is semipositive.
\end{enumerate}
\end{theo}

\begin{proof}
The proof of
``\ref{Item: first chern current positive} $\Longrightarrow$ \ref{Item: uniform extension on points}'' is very technical. For the proof, we refer to the papers
\cite{Zhang95} and \cite[Theorem~0.2]{Moriwaki2015}.
``\ref{Item: uniform extension on points} $\Longrightarrow$ \ref{Item: metric semipositive}'' is nothing more than  Proposition~\ref{prop:equivalent:semiample:metrized}.
Here let us consider the following claim:

\begin{enonce}{Claim}\label{claim:thm:semipositive:positive:current:01}
Let $M$ be an invertible $\mathcal O_X$-module, $\overline{V} = (V,\|\ndot\|)$ be a finite-dimensional normed vector
space over $\mathbb C$ and $V \otimes_{\mathbb C} \mathcal O_X \to M$ be a surjective homomorphism.
We assume that there is a basis $\{e_i\}_{i=1}^r$ of $V$ such that
\[
\forall\, (a_1, \ldots, a_r) \in \mathbb C^r, \quad\| a_1 e_1 + \cdots + a_r e_r \| = \max \{ |a_1|, \ldots, |a_r| \}.
\] 
Let $\psi$ be the quotient metric of $M$ induced by $\overline{V}$ and $V \otimes_{\mathbb C} \mathcal O_X \to M$.
Then the first Chern current $c_1(M, \psi)$ is semipositive.
\end{enonce}

\begin{proof}
Let $\{e_i^{\vee}\}_{i=1}^r$ be the dual basis of $V$. Then it is easy to see that
the dual norm $\|\ndot\|_*$ of $\|\ndot\|$ is given by
\[
\forall a_1, \ldots, a_r \in \mathbb C, \quad\| a_1 e^{\vee}_1 + \cdots + a_r e^{\vee}_r \|_* = |a_1| + \cdots + |a_r|.
\]
For $v \in V$, the induced global section of $M$ over $X$ is denoted by $\tilde{v}$.
Let $s$ be a local basis of $M$ over a Zariski open set $U$.
We set $\tilde{e}_i = a_i s$ for some holomorphic function $a_i$ on $U$. Then, by Proposition~\ref{prop:FS:metric:evaluation}, the function
\[
x\longmapsto - \ln |s|_{\psi}(x) = \ln ( |a_1|(x) + \cdots + |a_r|(x) )
\] 
is plurisubharmonic on $U^{\mathrm{an}}$ because $\ln ( |a_1|(\ndot) + \cdots + |a_r|(\ndot) )$ is plurisubharmonic on $U^{\mathrm{an}}$. 
\end{proof}

Let us see that \ref{Item: metric semipositive} $\Longrightarrow$ \ref{Item: first chern current positive}. Clearly we may assume that $L$ is generated by global sections.
For each $n \geqslant 1$, let $r_n := \dim_{\mathbb C} H^0(X, L^{\otimes n})$ and
$\{e_{n,i}\}_{i=1}^{r_n}$ be an orthonormal basis of $H^0(X, L^{\otimes n})$ with respect to $\|\ndot\|_{n \varphi}$.
If we set 
\[
\| a_1 e_{n,1} + \cdots + a_{r_n} e_{n, r}\|'_n := \max\{ |a_1|, \ldots, |a_{r_n}| \}
\]
for $a_1, \ldots, a_{r_n} \in \mathbb C$, then $\|\ndot\|'_n \leqslant \|\ndot\|_{n\varphi} \leqslant r_n \|\ndot\|'_n$. Let $\psi_n$ be the quotient metric of $L^{\otimes n}$ by $\|\ndot\|'_n$.
Then $d(\varphi_n, \psi_n) \leqslant \ln (r_n)$ because $\psi_n \leqslant \varphi_n \leqslant r_n \psi_n$.
Therefore, as
\[
d({\textstyle\frac 1n} \psi_n, \varphi) \leqslant d(\textstyle\frac 1n 
\psi_n, {\textstyle\frac 1n} \varphi_n) + d({\textstyle\frac 1n}\varphi_n, \varphi)
\leqslant \frac{1}{n}\ln (r_n) + d({\textstyle\frac 1n}\varphi_n, \varphi),
\]
one has $\lim_{n\to\infty} d(\frac 1n \psi_n, \varphi) = 0$ by our assumption.
This means that, for a local basis $s$ of $L$ over an open set $U$,
the sequence $\{ -\frac 1n \ln |s^{\otimes n}|_{\psi_n} \}_{n=1}^{\infty}$ converges to 
$-\ln |s|_{\varphi}$ uniformly on any compact set in $U$. 
As $-\frac 1n \ln |s^{\otimes n}|_{\psi_n}$ is plurisubharmonic by the above claim,
$-\ln |s|_{\varphi}$ is also plurisubharmonic, as required.
\end{proof}

\begin{coro}
Let $T$ be a reduced complex analytic space and $\|\ndot\|$ be a norm of $\mathbb C^n$.
If $f_1, \ldots, f_n$ are holomorphic functions on $T$, then
$\log \|( f_1, \ldots, f_n) \|$ is plurisubharmonic on $T$.
\end{coro}

\begin{proof}
First of all, recall the following fact (cf. \cite[Corollary~2.9.5]{Klimek1991}):
\begin{quote}
\it\quad If $u$ is a plurisubharmonic function on $\mathbb C^n$ and
$f_1, \ldots, f_n$ are holomorphic functions on $T$, then
$u(f_1, \ldots, f_n)$ is plurisubharmonic on $T$.
\end{quote}
Thus it is sufficient to see that $f(z_1, \ldots, z_n) := \log \| (z_1, \ldots, z_n) \|$ is plurisubharmonic on $\mathbb C^n \setminus \{ (0,\ldots, 0) \}$. 
Let $\mathbb C^n \otimes_{\mathbb C} \mathcal O_{\mathbb P_{\mathbb C}^{n-1}} \to \mathcal O_{\mathbb P_{\mathbb C}^{n-1}}(1)$ be the surjective homomorphism given by
$e_i \mapsto X_i$, where $\{e_i\}_{i=1}^n$ is the standard basis of $\mathbb C^n$ and
$(X_1 : \cdots : X_n)$ is a homogeneous coordinate of $\mathbb P^{n-1}_{\mathbb C}$.
Let us consider the dual norm $\|\ndot\|_*$ of $\|\ndot\|$ on $\mathbb C^n$, that is,
we identify the dual space $(\mathbb C^n)^{\vee}$ with $\mathbb C^n$ in the natural way.
Let $\varphi$ be the quotient metric of $\mathcal O_{\mathbb P_{\mathbb C}^{n-1}}(1)$ induced
by $\|\ndot\|_*$ and $\mathbb C^n \otimes_{\mathbb C} \mathcal O_{\mathbb P_{\mathbb C}^{n-1}} \to \mathcal O_{\mathbb P_{\mathbb C}^{n-1}}(1)$.
Note that $X_i$ gives a local basis of $\mathcal O_{\mathbb P_{\mathbb C}^{n-1}}(1)$ over $\{ X_i \not= 0 \}$ and
$\vec{X}_i = \sum_{j=1}^n (X_j/X_i) e_j$ (see Definition \ref{Def: rightarrow ell}), so that
by Proposition~\ref{prop:FS:metric:evaluation} together with the fact $\|\ndot\|_{**} = \|\ndot\|$, one has
$-\log |X_i|_{\varphi} = \log \| \vec{X}_i \|$
on $\{ X_i \not= 0 \}$. Therefore, by Theorem~\ref{thm:semipositive:positive:current} together with the previous fact, the function
\[(z_1,\ldots,z_n)\longmapsto
\log \|\textstyle (\frac{z_1}{z_i}, \ldots, \frac{z_{i-1}}{z_i}, 1, \frac{z_{i+1}}{z_i}, \ldots, \frac{z_n}{z_i})\|
\]
is plurisubharmonic on $\mathbb C^n \setminus \{ z_i = 0 \}$.
Note that 
\[
f(z_1, \ldots, z_n) = \log |z_i| + \log \| \textstyle (\frac{z_1}{z_i}, \ldots, \frac{z_{i-1}}{z_i}, 1, \frac{z_{i+1}}{z_i}, \ldots, \frac{z_n}{z_i})\|,
\]
so that $f$ is plurisubharmonic on $\mathbb C^n \setminus \{ z_i = 0 \}$ for all $i$, and hence
$f$ is plurisubharmonic on $\mathbb C^n \setminus \{ (0,\ldots, 0) \}$.
\end{proof}

\subsection{Model metrics}
In this subsection, we assume that the absolute value $|\ndot|$ is non-Archimedean and \emph{non-trivial}. We denote by $\mathfrak o_k$ the valuation ring of $(k,|\ndot|)$. Let $X\rightarrow\Spec(k)$ be a projective $k$-scheme and $L$ be an invertible $\mathcal O_X$-module. By \emph{model}\index{model} of $(X,L)$, we refer to a projective and flat $\mathfrak o_k$-scheme $\mathscr X$ equipped with an invertible $\mathcal O_{\mathscr X}$-module $\mathscr L$ such that the generic fibre of $\mathscr X$ coincides with $X$ and that the restriction of $\mathscr L$ {to} $X$ coincides with $L$. As in Definition \ref{Def:specification}, we denote by $j:X^{\mathrm{an}}\rightarrow X$ the specification map. 

Let $x$ be a point in $X^{\mathrm{an}}$ and let $p_x:\Spec\widehat{\kappa}(x)\rightarrow X$ be the $k$-morphism of schemes defined by $x$, where $\widehat{\kappa}(x)$ is the completion of the residue field $\kappa(x)$ of $j(x)$ with respect to the absolute value $|\ndot|_x$. The composition of $p_x$ with the inclusion morphism $X\rightarrow\mathscr X$ then defines a $\mathfrak o_k$-morphism from $\Spec\widehat{\kappa}(x)$ to $\mathscr X$. By definition $L\otimes\widehat{\kappa}(x)$ is the pull-back sheaf $p_x^{*}(L)$. By the valuative criterion of properness (see \cite{EGAII} Chapter II, Theorem 7.3.8), there exists a unique $\mathfrak o_k$-morphism $\mathscr P_x$ from $\Spec(\mathfrak o_x)$ to $\mathscr X$
which {identifies} with $p_x$ on the generic fibre, where $\mathfrak o_x$ is the valuation ring of $\widehat{\kappa}(x)$. 
The image of the maximal ideal of $\mathfrak o_x$ by $\mathscr P_x$, denoted by $r_{\mathscr X}(x)$, is
called the \emph{reduction point}\index{reduction point} of $x$. Note that $r_{\mathscr X}(x)$ {belongs to} 
the special fibre of $\mathscr X \to \Spec(\mathfrak o_k)$. Furthermore 
$\mathscr P_x^*(\mathscr L)$ is a lattice in $L\otimes\widehat{\kappa}(x)$ (see \S\ref{Subsec:latticesandnorms}). We denote by 
$|\ndot|_{\mathscr L}(x)$ the norm on $L\otimes\widehat{\kappa}(x)$ defined by this lattice, namely
\[\forall\,\ell\in L\otimes\widehat{\kappa}(x),\quad |\ell|_{\mathscr L}(x):=\inf\{|a|_x
\,:\,a\in\widehat{\kappa}(x)^{\times},\,a^{-1}\ell\in\mathscr P_x^*(\mathscr L)\}.\]
The family of norms $\{|\ndot|_{\mathscr L}(x)\}_{x\in X^{\mathrm{an}}}$ forms a metric on $L$ which we denote by $\varphi_{\mathscr L}$, called the metric \emph{induced} by the model $(\mathscr X,\mathscr L)$\index{metric!induced by a model}.

\begin{rema}
Let $(\mathscr X, \mathscr L)$ be a model of $(X, L)$.
Then $H^0(\mathscr X, \mathscr L)$ is a lattice of $H^0(X, L)$.
Indeed, since $X$ is a proper $k$-scheme, $H^0(\mathscr X,\mathscr L)$ is an $\mathfrak o_k$-module of finite type such that  (see \cite{NouveauEGA1}, Chapter I, Proposition 9.3.2) \[H^0(X,L)=H^0(\mathscr X,\mathscr L)\otimes_{\mathfrak o_k}k.\] Moreover, since $\mathscr X$ is a flat $\mathfrak o_k$-scheme, the $\mathfrak o_k$-module $H^0(\mathscr X,\mathscr L)$ is torsion free. Therefore the canonical map $H^0(\mathscr X,\mathscr L)\rightarrow H^0(X,L)$ is injective and hence it is a lattice in $H^0(X,L)$ (see Definition \ref{Def:lattice}). 
\end{rema}

\begin{rema}
Let $\mathscr X$ and $\mathscr Y$ be projective $\mathfrak o_k$-schemes and $f:\mathscr Y\rightarrow\mathscr X$ an $\mathfrak o_k$-morphisme. Let $X$ and $Y$ be the generic fibres of $\mathscr X$ and $\mathscr Y$ respectively, and $f_k:X\rightarrow Y$ be the morphism induced by $f$. Let $\mathscr L$ be an invertible sheaf on $\mathscr X$ and $L$ be the restriction of $\mathscr L$ {to} the generic fibre $X$. Then the model $(\mathscr X,\mathscr L)$ induces a metric $|\ndot|_{\mathscr L}$ on the invertible sheaf $L$. The couple $(\mathscr Y,f^*(\mathscr L))$ forms a model of $(Y,f_k^*(L))$. Note that the model metric $|\ndot|_{f^*(\mathscr L)}$ on $f_k^*(L)$ coincides with the pull-back of the metric $|\ndot|_{\mathscr L}$ by $f_k$ (see Definition \ref{Def:pull-back}).
\end{rema}

\begin{prop}\label{Pro:modelFubiniStudy}
Let $L$ be an invertible $\mathcal O_{X}$-module which is generated by global sections.
Let $\varphi_{\mathscr L}$ be the metric induced by a model $(\mathscr X, \mathscr L)$ of $(X, L)$. 
Let $\varphi$ be a continuous metric of $L$ and \[\mathscr H = \{ s \in H^0(\mathscr X, \mathscr L) \,:\, \| s\|_{\varphi} \leqslant 1 \}.\]
Then one has the following:
\begin{enumerate}[label=\rm(\arabic*)]
\item\label{Item: if generated by model sections the varphi smaller than varphi l} 
If $\mathscr H \otimes_{\mathfrak o_k} \mathcal O_{\mathscr X} \to \mathscr L$ is surjective, then
$\varphi \leqslant \varphi_{\mathscr L}$.

\item\label{Item: quotient by global sections varphi larger}
If $\varphi$ is the quotient metric on $L$ induced by $\|\ndot\|_{H^0(\mathscr X, \mathscr L)}$ (see Definition \ref{Def:inducedbylattice} for the norm induced by a lattice),then $\varphi \geqslant \varphi_{\mathscr L}$.

\item\label{Item: Quotient induced by the lattice metric equal}
If $\varphi$ is the quotient metric on $L$ induced by $\|\ndot\|_{H^0(\mathscr X, \mathscr L)}$,  
and the natural homomorphism $H^0(\mathscr X, \mathscr L) \otimes_{\mathfrak o_k} \mathcal O_{\mathscr X} \to \mathscr L$ is surjective, then $\varphi = \varphi_{\mathscr L}$.
\end{enumerate}
\end{prop}

\begin{proof}
For $x \in X^{\mathrm{an}}$, let $p_x:\Spec\widehat{\kappa}(x)\rightarrow X$ be the $k$-morphism of schemes defined by $x$, and $\mathscr P_x:\Spec\mathfrak o_x\rightarrow\mathscr X$ be the $\mathfrak o_k$-morphism extending $p_x$.
Moreover, let $\pi_x:H^0(X,L)\otimes_k\widehat{\kappa}(x)\rightarrow p_x^*(L)$ be the natural homomorphism.

\medskip
\ref{Item: if generated by model sections the varphi smaller than varphi l} By our assumption, $\mathscr H \otimes_{\mathfrak o_k} \mathfrak o_x \to \mathscr P_x^*(\mathscr L)$
is surjective, so that, 
if $\ell$ lies in $\mathscr P_x^*(\mathscr L)$, then there exist $s_1,\ldots,s_n$ in $\mathscr H$ and $a_1,\ldots,a_n$ in $\mathfrak o_x$ such that $\pi_x(a_1s_1+\cdots+a_ns_n)=\ell$. For any $i\in\{1,\ldots,n\}$, let $\ell_i=\pi_x(s_i)$. Then one has $\ell=a_1\ell_1+\cdots+a_n\ell_n$. As $s_i \in \mathscr H$,
one has $|\ell_i|_{\varphi}(x)\leqslant 1$ for any $i$, which leads to $|\ell|_{\varphi}(x)\leqslant 1$. By Proposition \ref{Pro:comparisonbyunitball}, one obtains $|\ndot|_{\varphi}(x)\leqslant|\ndot|_{\varphi_{\mathscr L}}(x)$. The assertion \ref{Item: if generated by model sections the varphi smaller than varphi l} is thus proved.   

\medskip
\ref{Item: quotient by global sections varphi larger} Note that $p_x^*(L)$ is a quotient vector space of rank $1$ of $H^0(X, L) \otimes_k\widehat{\kappa}(x)$ and 
$|\ndot|_{\varphi}(x)$ is the quotient norm on $p_x^*(L)$ induced by $\|\ndot\|_{H^0(\mathscr X, \mathscr L),\widehat{\kappa}(x)}$.
By Proposition \ref{Pro:quotientderang1norm}, for $\ell\in p_x^*(L)\setminus\{0\}$ one has
\begin{equation}\label{Equ:normofellx}|\ell|_{\varphi}(x)=\inf_{\begin{subarray}{c}s\in H^0(X,L),\,\lambda\in\widehat{\kappa}(x)^{\times}\\
\pi_x(s)=\lambda \ell\end{subarray}}|\lambda|^{-1}\|s\|_{H^0(\mathscr X, \mathscr L)}.\end{equation}
Let $s\in H^0(X,L)$ and $\lambda\in\widehat{\kappa}(x)^{\times}$ such that $\pi_x(s)=\lambda\ell$. By definition one has 
\[\|s\|_{H^0(\mathscr X, \mathscr L)}=\inf\{|a|\,:\,a\in k^{\times},\;a^{-1}s\in H^0(\mathscr X, \mathscr L)\}.\]
If $a$ is an element in $k^{\times}$ such that $a^{-1}s\in H^0(\mathscr X, \mathscr L)$, then $a^{-1}\lambda\ell\in\mathscr P_x^*(\mathscr L)$ because $\mathscr P_x^*(\mathscr L)$ contains the image of $H^0(\mathscr X, \mathscr L)\otimes_{\mathfrak o_k}\mathfrak o_x$ in $p_x^*(L)$ by $\pi_x$. Hence \[|a^{-1}\lambda\ell|_{\varphi_{\mathscr L}}(x)=|a|^{-1}|\lambda|\cdot|\ell|_{\varphi_{\mathscr L}}(x)\leqslant 1,\]
which implies that $|\ell|_{\varphi_{\mathscr L}}(x)\leqslant |\lambda|^{-1}|a|$. Since $a$ is arbitrary with $a^{-1} s \in H^0(\mathscr X, \mathscr L)$, one obtains $|\ell|_{\varphi_{\mathscr L}}(x)\leqslant|\lambda|^{-1}\|s\|_{H^0(\mathscr X, \mathscr L)}$, which leads to $|\ell|_{\varphi_{\mathscr L}}(x)\leqslant|\ell|_{\varphi}(x)$.

\medskip
\ref{Item: Quotient induced by the lattice metric equal} By \ref{Item: quotient by global sections varphi larger}, it is sufficient to see $\varphi \leqslant \varphi_{\mathscr L}$.
Note that for $s \in H^0(\mathscr X, \mathscr L)$,
one has $\| s \|_{\varphi} \leqslant 1$, so that $\mathscr H = H^0(\mathscr X, \mathscr L)$. Thus, by (1), one obtains $\varphi \leqslant \varphi_{\mathscr L}$.
\end{proof}

\begin{coro}\label{Coro:modelFubiniStudy}
Let $E$ be a {finite-dimensional} vector space over $k$, $\mathcal E$ be a lattice in $E$ and $\|\ndot\|_{\mathcal E}$ be the norm on $E$ induced by the lattice $\mathcal E$ (see Definition \ref{Def:inducedbylattice}). Then the Fubini-Study  metric (see Definition \ref{Def:Fubini-Studymetric}) on the invertible $\mathcal O_{\mathbb P(E)}$-module $\mathcal O_{E}(1)$ induced by $\|\ndot\|_{\mathcal E}$ coincides with the metric induced by the model $(\mathbb P(\mathcal E),\mathcal O_{\mathcal E}(1))$ of $(\mathbb P(E),\mathcal O_E(1))$.
\end{coro}

\begin{rema}\phantomsection
\label{Rem:restrictioofmodel}
\begin{enumerate}[label=\rm(\arabic*)]
\item\label{Item: Prufer domain} Recall that a valuation ring is a Pr\"ufer domain, which is a generalisation of Dedekind domain (non-necessarily Noetherian). In particular, an $\mathfrak o_k$-module is flat if and only if it is torsion free (see \cite{Bourbaki64} Chapter VII, \S2, Exercices 12 and 14). 

\item \label{Item: Zariski closure}
Let $\mathscr A$ be a flat $\mathfrak o_k$-algebra, and $A=S^{-1}\mathscr A$, where $S=\mathfrak o_k\setminus\{0\}$. Note that the canonical homomorphism $\mathscr A\rightarrow A$ is injective. If $B$ is a quotient ring of $A$ with respect to the ideal $I$ and if $\mathscr I=I\cap\mathscr A$, then the quotient ring $\mathscr B=\mathscr A/\mathscr I$ is a torsion-free $\mathfrak o_k$-module, and hence is flat. 
The scheme $\Spec(\mathscr B)$ is called the \emph{Zariski closure of $\Spec(B)$ in $\Spec(\mathscr A)$}\index{Zariski closure}.

\item \label{Item:model Zariski closure}
Let $X$ be a projective scheme over $\Spec k$, equipped with an invertible $\mathcal O_X$-module $L$. Let $(\mathscr X,\mathscr L)$ be a model of $(X,L)$. Then the global section space $H^0(\mathscr X,\mathscr L)$ is a torsion-free $\mathfrak o_k$-module. Moreover, if $Y$ is a closed subscheme of $X$ and if $\mathscr Y$ is the Zariski closure of $Y$ in $\mathscr Y$, then $(\mathscr Y,\mathscr L|_{\mathscr Y})$ is a model of $(Y,L|_Y)$.  Moreover, the metric on $L|_Y$ induced by $\mathscr L|_{\mathscr Y}$ coincides with the restriction of the metric $\varphi_{\mathscr L}$ {to} $L|_Y$. 

\item\label{Item: model and quotient}
Let $\pi:X\rightarrow\Spec k$ be a projective scheme over $\Spec k$, $E$ be a finite-dimensional vector space over $k$ and $\beta:\pi^*(E)\rightarrow L$ be a surjective homomorphism, which defines a $k$-morphism $f:X\rightarrow\mathbb P(E)$. Suppose that $f$ is a closed immersion. Let $\mathcal E$ be a lattice in $E$, $\mathscr X$ be the Zariski closure of $X$ in $\mathbb P(\mathcal E)$ and $\mathscr L$ be the restriction of $\mathcal O_{\mathcal E}(1)$ {to} $\mathscr X$. Then Proposition \ref{Pro:modelFubiniStudy} shows that the quotient metric on $L$ induced by the norm $\|\ndot\|_{\mathcal E}$ coincides with the metric induced by the model $(\mathscr X,\mathscr L)$.
\end{enumerate}
\end{rema}

\begin{prop}\label{Pro:tensormodel}
Let $L$ and $M$ are two invertible $\mathcal O_X$-modules. Suppose that $\mathscr X$ is a proper and flat $\mathfrak o_k$-scheme such that $\mathscr X_k=X$. If $\mathscr L$ and $\mathscr M$ are invertible $\mathcal O_{\mathscr X}$-modules such that $(\mathscr X,\mathscr L)$ and $(\mathscr X,\mathscr M)$ are models of $(X,L)$ and $(X,M)$ respectively, then one has
$\varphi_{\mathscr L\otimes\mathscr M}=\varphi_{\mathscr L}+\varphi_{\mathscr M}$. 
\end{prop}
\begin{proof}
Let $x$ be a point of $X^{\mathrm{an}}$, $s$ and $t$ are element in $\mathscr L_{r_{\mathscr X}(x)}$ and $\mathscr M_{r_{\mathscr X}(x)}$ which trivialise the invertible sheaves $\mathscr L$ and $\mathscr M$ around $r_{\mathscr X}(x)$ respectively, and $\ell=\mathscr P_x^*(s)$, $m=\mathscr P_x^*(t)$, where $\mathscr P_x:\Spec(\mathfrak o_x)\rightarrow\mathscr X$ is the unique $\mathfrak o_k$-morphism extending the $k$-morphism $\Spec\widehat{\kappa}(x)\rightarrow X$ corresponding to the point $x$. Since $\ell$ and $m$ are generators of the free $\mathfrak o_{x}$-modules (of rank $1$) $\mathscr P_x^*(\mathscr L)$ and $\mathscr P_x^*(\mathscr M)$, one has $|\ell|_{\mathscr L}(x)=|m|_{\mathscr M}(x)=1$. Moreover, $\mathscr P_x^*(s\otimes t)$ is a generator of the free $\mathfrak o_x$-module $\mathscr P_x^*(\mathscr L\otimes\mathscr M)$, we obtain that $|\ell\otimes m|_{\mathscr L\otimes\mathscr M}(x)=1$. Hence one has $\varphi_{\mathscr L\otimes\mathscr M}=\varphi_{\mathscr L}+\varphi_{\mathscr M}$.
\end{proof}

\begin{prop}\label{Pro:latticeby model}
Let $X$ be a proper $k$-scheme and $L$ be an invertible $\mathcal O_X$-module. Let $(\mathscr X,\mathscr L)$ be a model of $(X,L)$. 
\begin{enumerate}[label=\rm(\arabic*)]
\item\label{Item: model metric continous} The metric $\varphi_{\mathscr L}$ on $L$ is continuous.
\item\label{Item: global section in the unit ball} The $\mathfrak o_k$-module $H^0(\mathscr X,\mathscr L)$ is a lattice in $H^0(X,L)$, which is contained in the unit ball of $H^0(X,L)$ with respect to the seminorm $\|\ndot\|_{\varphi_{\mathscr L}}$. 
\item\label{Item: global section is unit ball} If the valuation $|\ndot|$ is discrete, $X$ is reduced and the central fibre of $\mathscr X \to \Spec(\mathfrak o_k)$ is reduced, then
\[\{ s \in H^0(X, L) \,:\, \| s \|_{\varphi_{\mathscr L}} \leqslant 1 \} = H^0(\mathscr X, \mathscr L)\] and
$ \| s \|_{\varphi_{\mathscr L}} = \|\ndot\|_{H^0(\mathscr X, \mathscr L)}$. 
\end{enumerate}
\end{prop}
\begin{proof}
\ref{Item: model metric continous}
We choose an ample invertible $\mathcal O_{\mathscr X}$-module $\mathscr A$ such that $\mathscr L \otimes \mathscr A$ and $\mathscr A$ are generated by global sections.
Let $\varphi_1$ and $\varphi_2$ be the quotient metrics of $L \otimes A$ and $A$ induced by $\|\ndot\|_{H^0(\mathscr X, \mathscr L \otimes \mathscr A)}$ and $\|\ndot\|_{H^0(\mathscr X, \mathscr A)}$, respectively, where $A = \rest{\mathscr A}{X}$. Then $\varphi_1$ and $\varphi_2$ are continuous by Prposition~\ref{prop:FS:metric:continuous} and Definition~\ref{Def: Quotient metric}. Moreover, by Proposition~\ref{Pro:modelFubiniStudy}, $\varphi_{\mathscr L \otimes \mathscr A} = \varphi_1$ and $\varphi_{\mathscr A} = \varphi_2$. Therefore one has the assertion because $\varphi_{\mathscr L} = \varphi_{\mathscr L \otimes \mathscr A} -  \varphi_{\mathscr A} = \varphi_1 - \varphi_2$ by Proposition~\ref{Pro:tensormodel}.
\if01
Let $\mathscr U$ be a Zariski open subset of $\mathscr X$ on which the sheaf $\mathscr L$ is trivial and $U$ be the intersection of $\mathscr U$ with the generic fibre $X$. Let $s$ be a section of $\mathscr L$ on $\mathscr U$ which trivialises $\mathscr L$. By abuse of notation, we still denote by $s$ its restriction {to} $U$. Then by definition one has $|s|_{\mathscr L}(x)=1$ for any $x\in U^{\mathrm{an}}$. Therefore the metric $\varphi_{\mathscr L}$ is continuous. 
\fi

\medskip\ref{Item: global section in the unit ball} 
Let $s$ be a section in $ H^0(\mathscr X,\mathscr L)$, viewed as an element in $H^0(X,L)$, by definition one has
\[\forall\,x\in X^{\mathrm{an}},\quad\|s\|_{\mathscr L}(x)\leqslant 1.\]
Hence $H^0(\mathscr X,\mathscr L)$ is contained in the closed unit ball of $(H^0(X,L),\|\ndot\|_{\varphi_{\mathscr L}})$.

\medskip\ref{Item: global section is unit ball} Let $\mathscr E$ be the unit ball of $H^0(X, L)$ with respect to $\|\ndot\|_{\varphi_{\mathscr L}}$.
First let us see that 
\begin{equation}\label{eqn:Pro:latticeby model:00}
H^0(\mathscr X, \mathscr L) = \mathscr E.
\end{equation}
By \ref{Item: global section in the unit ball}, $H^0(\mathscr X, \mathscr L) \subseteq \mathscr E$.
Let $\mathscr X = \bigcup_{i=1}^N \Spec(\mathscr A_i)$ be an affine open covering of $\mathscr X$ such that $\mathscr A_i$ is of finite type over $\mathfrak o_{k}$ and
$\mathscr L$ has a local basis $\ell_i$ over $\Spec(\mathscr A_i)$. For $s \in \mathscr E$,
we set $s = f_i \ell_i$ and $f_i \in A_i := S^{-1} \mathscr A_i$, where $S = \mathfrak o_{k} \setminus \{ 0 \}$. Then, for $x\in(\Spec A_i)_{\mathscr A_i}^{\mathrm{an}}$ (cf. Remark~\ref{remark:reduction:point}),
$|s|_{\varphi_{\mathscr L}}(x) = |f_i|_x \leqslant 1$. Therefore, as the central fibre
$\mathscr X_{\circ}$ of $\mathscr X \to \Spec(\mathfrak o_k)$ is reduced,
by the last assertion of Proposition~\ref{Pro:integralclosure}, one has $f_i \in \mathscr A_i$, and hence
$s \in H^0(\mathscr X, \mathscr L )$.

Next we need to see that
\begin{equation}\label{eqn:Pro:latticeby model:01}
\quad
\|\ndot\|_{\varphi_{\mathscr L}} = \|\ndot\|_{\mathscr E}.
\end{equation}
Let $\varpi$ be a uniformising parameter of $\mathfrak o_{k}$, $\mathscr X_{\circ}$ be the fibre of $\mathscr X$ over the maximal ideal of $\mathfrak o_k$
and $\mathscr L_{\circ}$ be the restriction of $\mathscr L$ to 
$\mathscr X_{\circ}$, that is, $\mathscr L_{\circ} = \mathscr L/\varpi \mathscr L$.
{The} short exact sequence
\[
\xymatrix{\relax 0 \ar[r]& \mathscr L \ar[r]^-{\varpi\cdot}& \mathscr L \ar[r]&
\mathscr L_{\circ} \ar[r]& 0}
\]
gives rise to an exact sequence:
\[
\xymatrix{\relax 0 \ar[r]& H^0(\mathscr X, \mathscr L)\ar[r]^- {\varpi\cdot} &
 H^0(\mathscr X, \mathscr L) \ar[r]& H^0(\mathscr X_{\circ}, \mathscr L_{\circ}),}
\]
that is, the natural homomorphism 
\begin{equation}\label{eqn:Pro:latticeby model:02}
H^0(\mathscr X, \mathscr L)/\varpi H^0(\mathscr X, \mathscr L) \longrightarrow H^0(\mathscr X_{\circ}, \mathscr L_{\circ})
\end{equation}
is injective. Moreover,
by Proposition~\ref{Pro:normetreausauxdisc}, one has
\begin{equation}\label{eqn:Pro:latticeby model:03}
\|\ndot\|_{\varphi_{\mathscr L}} \leqslant \|\ndot\|_{\mathscr E}.
\end{equation}

Here we claim that if $\| s\|_{\mathscr E} = 1$ for $s \in H^0(X, L)$, then
$\|s \|_{\varphi_{\mathscr L}} = 1$. Obviously $\|s \|_{\varphi_{\mathscr L}} \leqslant 1$ by \eqref{eqn:Pro:latticeby model:03}.
As $H^0(\mathscr X, \mathscr L ) = \mathscr E$ by \eqref{eqn:Pro:latticeby model:00},
one has $s \in H^0(\mathscr X, \mathscr L )$ and $s$ is not zero in
$H^0(\mathscr X, \mathscr L)/\varpi H^0(\mathscr X, \mathscr L)$,
so that by the injectivity of \eqref{eqn:Pro:latticeby model:02}, $s$ is not zero in $H^0(\mathscr X_{\circ}, \mathscr L_{\circ})$.
Let $\xi$ be a closed point $\mathscr X_{\circ}$ with $s(\xi) \not= 0$.
Let $\ell_{\xi}$ be a local basis of $\mathscr L$ around $\xi$. Then $s = f \ell_{\xi}$ and $f \in \mathcal O_{\mathscr X,\xi}^{\times}$. On the other hand,
since the reduction map $r : X^{\mathrm{an}} \to \mathscr X_{\circ}$ is surjective, one can find $x \in X^{\mathrm{an}}$
with $r(x) = \xi$. Then $| s |_{\varphi_{\mathscr L}}(x) = |f|_x = 1$, so that $\|s \|_{\varphi_{\mathscr L}} = 1$, as desired.

In general, for $s \in H^0(X, L) \setminus \{ 0 \}$, there is an integer $e$ such that $\| \varpi^e s \|_{\mathscr E} = 1$, so that
$\| \varpi^e s \|_{\varphi_{\mathscr L}} = 1$, and hence
$\| s \|_{\mathscr E} = \|  s \|_{\varphi_{\mathscr L}} = |\varpi|^{-e}$.
\end{proof}

\begin{prop}
\label{prop:vanishing:mu:nef:big}
Let $X$ be a projective $k$-scheme and $L$ be an ample invertible $\mathcal O_X$-module. If $(\mathscr X,\mathscr L)$ is a model of $(X,L)$ such that $L$ is ample and $\mathscr L$ is nef,
then the metric $\varphi_{\mathscr L}$ is semipositive.
\end{prop}
\begin{proof} Let $\pi:\mathscr X\rightarrow\Spec\mathfrak o_k$ be the structural morphism.
First we assume that 
$\mathscr L$ is ample. We choose a positive integer $n$ such that
$\mathscr L^{\otimes n}$ is very ample.
Then we have a closed embedding $\iota : \mathscr X \to \mathbb P(\mathcal E_n)$ with  $\mathcal E_n:=H^0(\mathscr X,\mathscr L^{\otimes n})$, which is induced by the canonical (surjective) homomorphism $\pi^*(\pi_*(\mathscr L^{\otimes n}))=\pi^*(\mathcal E_n)\rightarrow\mathscr L^{\otimes n}$. Note that one has
$\mathscr L^{\otimes n} = \iota^{*}(\mathcal O_{\mathcal E_n}(1))$. Moreover, $\mathcal E_n$  is a lattice in $E_n:=H^0(X,L^{\otimes n})$. By Proposition \ref{Pro:modelFubiniStudy}, the metric $\varphi_{\mathscr L^{\otimes n}}=n\varphi_{\mathscr L}$ (see Proposition \ref{Pro:tensormodel} for this equality) coincides with the quotient metric on $ L^{\otimes n}$ induced by the norm $\|\ndot\|_{\mathcal E_n}$, hence is semipositive (see Proposition \ref{Pro:positivityofquotientmetric}). By Proposition \ref{Pro:puissancepositve}, we obtain that the metric $\varphi_{\mathscr L}$ is also semipositive.

In the following, we treat the general nef case. Let $\mathscr M$ be an ample invertible sheaf on $\mathscr X$
and $M$ be the restriction of $\mathscr M$ {to} $X$. Since $\mathscr L$ is nef, we obtain that $\mathscr L^{\otimes n}\otimes\mathscr M$ is ample for any $n\in\mathbb N_{n\geqslant 1}$. Moreover, since $L$ is ample, for a sufficiently positive integer $n_0\geqslant 1$, the invertible $\mathcal O_X$-module $L^{\otimes n_0}\otimes M^\vee$ is generated by global sections. Thus for any integer $n>n_0$ one has (by Proposition \ref{Pro:subadditivedelta})
\[n\defp(\varphi_{\mathscr L})=\defp(\varphi_{\mathscr L^{\otimes n}})\leqslant\defp(\varphi_{\mathscr L^{\otimes n_0}\otimes\mathscr M^{\vee}})+\defp(\varphi_{\mathscr L^{\otimes(n-n_0)}\otimes\mathscr M})=\defp(\varphi_{\mathscr L^{\otimes n_0}\otimes\mathscr M^{\vee}}),\]
where the first equality comes from \eqref{Equ:dphomogene}, and the second equality comes from the semi-positivity of the metric 
$\varphi_{\mathscr L^{\otimes(n-n_0)}\otimes\mathscr M}$. Since $n\geqslant n_0$ is arbitrary, we obtain that $\defp(\varphi_{\mathscr L})=0$, namely $\varphi_{\mathscr L}$ is a semipositive metric.
\end{proof}

\begin{prop}\label{Cor:convergencemodel}
Let $\pi:X\rightarrow\Spec k$ be a projective scheme over $\Spec k$ and $L$ be an ample invertible $\mathcal O_X$-module, equipped with a continuous metric $\varphi$ which is semipositive. Then, for sufficiently positive integer $n$, there exists a sequence {$\{(\mathscr X_n,\mathscr L_n)\}_{n\in\mathbb N,\,n\geqslant 1}$} of models of $(X,L^{\otimes n})$ such that each $\mathscr L_n$ is ample, and that 
\begin{equation}\label{Equ:convergencemodel}\lim_{n\rightarrow+\infty}\frac 1n d(\varphi_{\mathscr L_n},n\varphi)=0.\end{equation}
\end{prop}
\begin{proof} Let $\lambda\in \intervalle]01[$ be a number such that \[\lambda<\sup\{|a|\,:\,a\in k^{\times},\,|a|<1\}.\] For any $n\in\mathbb N_{\geqslant 1}$, let $V_n=H^0(X,L^{\otimes n})$.
Since $L$ is an ample, for sufficiently positive integer $n$, the canonical homomorphism $\pi^*(V_n)\rightarrow L^{\otimes n}$ is surjective, and the corresponding $k$-morphism $X\rightarrow\mathbb P(V_n)$ is a closed embedding. By Proposition \ref{Pro:approximation}, there exists a lattice of finite type $\mathcal V_n$ of $V_n$ such that \[d(\|\ndot\|_{\mathcal V_n},\|\ndot\|_{n\varphi})\leqslant \ln(\lambda^{-1}).\] Let $\mathscr X_n$ be the Zariski closure of $X$ in $\mathbb P(\mathcal V_n)$ and $\mathscr L_n$ be the restriction of $\mathcal O_{\mathbb P(\mathcal V_n)}(1)$ {to} $\mathscr X$, then $(\mathscr X,\mathscr L_n)$ is a model of $(X,L)$ with $\mathscr L_n$ being  ample. Moreover, the metric on $L^{\otimes n}$ induced by $\mathscr L_n$ coincides with the quotient metric on $L^{\otimes n}$ induced by $(V_n,\|\ndot\|_{\mathcal V_n})$ and the canonical quotient homomorphism $\pi^*(V_n)\rightarrow L^{\otimes n}$. Therefore by Propositon \ref{Pro:distrancequot} one has \[d(\varphi_{\mathscr L_n},n\varphi)\leqslant d(\|\ndot\|_{\mathcal V_n},\|\ndot\|_{n\varphi})\leqslant\ln(\lambda^{-1})\]
as required.
\end{proof}

{\subsection{Purity}

Let $X\rightarrow\Spec k$ be a projective $k$-scheme, $L$ be an invertible $\mathcal O_X$-module and $\varphi$ be a continuous metric on $X$. If the norm $\norm{\ndot}_{\varphi}$ on $H^0(X,L)$ is pure, we say that the metric $\varphi$ is \emph{pure}\index{pure}\index{metric!pure ---}. If $n\varphi$ is pure for all $n\in\mathbb N_{\geqslant 1}$, we say that $\varphi$ is \emph{stably pure}\index{stably pure}\index{metric!stably pure}. Note that, if the absolute value $|\ndot|$ is not discrete, then any continuous metric on $L$ is stably pure
(cf. Proposition~\ref{Pro:notdiscret}).

\begin{prop}
We assume that the absolute value $|\ndot|$ is discrete. Let $X$ be a projective $k$-scheme and $L$ be an invertible $\mathcal O_X$-module. If $(\mathscr X,\mathscr L)$ is a model of $(X,L)$ such that the central fibre of $\mathscr X\rightarrow\Spec\mathfrak o_k$ is reduced, then the metric $\varphi_{\mathscr L}$ is stably pure.  
\end{prop}
\begin{proof}
Note that for any $n\in\mathbb N_{\geqslant 1}$ one has $n\varphi_{\mathscr L}=\varphi_{\mathscr L^{\otimes n}}$ and $(\mathscr X,\mathscr L^{\otimes n})$ is a model of $(X,L^{\otimes n})$. Therefore it suffices to show that $\varphi_{\mathscr L}$ is a pure metric. By Proposition \ref{Pro:latticeby model}, the norm $\norm{\ndot}_{\varphi_{\mathscr L}}$ is induced by the lattice $H^0(\mathscr X,\mathscr L)$, hence it is pure. 
\end{proof}}

\subsection{Extension property}

In this subsection, we introduce the extension property of an ample invertible module with a semipositive continuous metric, that is,
an extension of a section with a control on the norm.

Throughout this subsection,
let $\pi:X\rightarrow\Spec(k)$ be a projective $k$-scheme and $L$ be an invertible $\mathcal O_X$-module, 
equipped with a continuous metric $\varphi$. Let us begin with the following lemma.

\begin{lemm}\label{lemm:properties:a:varphi:n:l}
Let $Y$ be a closed subscheme of $X$.  For $n \in \mathbb N$, $n \geqslant 1$,
let $\gamma_n : H^0(X, L^{\otimes n}) \to H^0(Y, \rest{L^{\otimes n}}{Y})$ be the restriction map.
For any  
element $\ell$ of $H^0(Y,L|_Y) \setminus \mathcal N(Y, L|_Y)$, we define $a_{\varphi, n}(\ell) \in [0, \infty]$ to be
\[
a_{\varphi, n}(\ell) := \begin{cases}
\infty & \text{if $\gamma_n^{-1}(\{\ell^{\otimes n}\}) = \varnothing$},\\
{\displaystyle \inf_{\begin{subarray}{c}
s\in H^0(X,L^{\otimes n})\\
s|_Y=\ell^{\otimes n}
\end{subarray}}\bigg(\ln\|s\|_{n\varphi}-\ln\|\ell\|^n_{\varphi|_Y}\bigg)} & \text{otherwise},
\end{cases}
\]
where $\varphi|_Y$ denotes the \emph{restriction}\index{restriction of a metric}\index{metric!restriction} of $\varphi$ {to} $L|_Y$, defined as the pull-back of $\varphi$ by the inclusion morphisme $Y\rightarrow X$ (see Definition \ref{Def:pull-back}). 
Then we have the following:
\begin{enumerate}[label=\rm(\arabic*)]
\item\label{Item: sequence an is subadditive}
The sequence $\{ a_{\varphi, n}(\ell) \}_{n \in \mathbb N}$ is subadditive, that is,
$a_{\varphi, n+n'}(\ell) \leqslant a_{\varphi, n}(\ell) + a_{\varphi, n'}(\ell)$ for $n, n' \in \mathbb N$.

\item\label{Item: continuity of an}
Let $\varphi'$ be another continuous metric of $L$. If $\gamma_n^{-1}(\{\ell^{\otimes n}\}) \neq \varnothing$, then
\[
| a_{\varphi, n}(\ell) - a_{\varphi', n}(\ell)| \leqslant 2n\, d(\varphi, \varphi').
\]
\end{enumerate}
\end{lemm}

\begin{proof}
\ref{Item: sequence an is subadditive} Clearly we may assume that $\gamma_n^{-1}(\{\ell^{\otimes n}\}) \not= \varnothing$ and $\gamma_{n'}^{-1}(\{\ell^{\otimes n'}\}) \not= \varnothing$.
Then $\gamma_{n+n'}^{-1}(\{\ell^{\otimes n+n'}\}) \not= \varnothing$, so that
{\allowdisplaybreaks
\begin{align*}
a_{\varphi, n+n'}(\ell) & = \inf_{\begin{subarray}{c}
s''\in H^0(L^{\otimes n+n'})\\
s''|_Y=\ell^{\otimes n + n'}
\end{subarray}}\bigg(\ln\|s''\|_{(n+n')\varphi}-\ln\|\ell\|^{n+n'}_{\varphi|_Y}\bigg) \\
& \leqslant \inf_{\begin{subarray}{c}
(s,s') \in H^0(L^{\otimes n}) \times H^0(L^{\otimes n'})\\
s|_Y=\ell^{\otimes n}, s'|_Y=\ell^{\otimes n'}
\end{subarray}}\bigg(\ln\|s \otimes s'\|_{(n+n')\varphi}-\ln\|\ell\|^{n+n'}_{\varphi|_Y}\bigg) \\
& \leqslant \inf_{\begin{subarray}{c}
(s,s') \in H^0(L^{\otimes n}) \times H^0(L^{\otimes n'})\\
s|_Y=\ell^{\otimes n}, s'|_Y=\ell^{\otimes n'}
\end{subarray}}\bigg(\ln\|s\|_{n\varphi} + \ln\|s'\|_{n'\varphi} - \ln\|\ell\|^{n+n'}_{\varphi|_Y}\bigg) \\
& = a_{\varphi, n}(\ell) + a_{\varphi, n'}(\ell),
\end{align*}}
as required.

\medskip
\ref{Item: continuity of an} 
Clearly we may assume that $a_{\varphi, n}(\ell) \geqslant a_{\varphi', n}(\ell)$.
For any $\epsilon > 0$, choose $s \in H^0(X, L^{\otimes n})$ such that $\rest{s}{Y} = \ell^{\otimes n}$ and
\[
\ln\|s\|_{n\varphi'}-\ln\|\ell\|^n_{\varphi'|_Y} \leqslant a_{\varphi', n}(\ell) + \epsilon.
\]
Then, by using \eqref{Equ:comparaisondistancesup} and \eqref{Equ:distancen},
{\allowdisplaybreaks
\begin{align*}
a_{\varphi, n}(\ell) - a_{\varphi', n}(\ell) & \leqslant (\ln\|s\|_{n\varphi}-\ln\|\ell\|^n_{\varphi|_Y})  - (\ln\|s\|_{n\varphi'}-\ln\|\ell\|^n_{\varphi'|_Y}) + \epsilon \\
& \leqslant \big| \ln\|s\|_{n\varphi} - \ln\|s\|_{n\varphi'} \big| + \big| \ln\|\ell\|^n_{\varphi|_Y} - \ln\|\ell\|^n_{\varphi'|_Y}\big| + \epsilon \\
& \leqslant d(\|\ndot\|_{n\varphi}, \|\ndot\|_{n\varphi'}) + n d(\|\ndot\|_{\varphi|_Y}, \|\ndot\|_{\varphi'|_Y}) + \epsilon \\
& \leqslant d(n\varphi, n\varphi') + n d(\varphi|_Y, \varphi'|_Y) + \epsilon \leqslant 2n\, d(\varphi, \varphi') + \epsilon,
\end{align*}}
so that the assertion follows because $\epsilon$ is an arbitrary positive number.
\end{proof}

\begin{defi} Let $Y$ be a closed subscheme of $X$. 
For $\ell \in H^0(Y, \rest{L}{Y})$,
we say that $\ell$ has the \emph{extension property}\index{extension property} for the metric $\varphi$ if,
for any $\epsilon>0$, there exists $n_0\in\mathbb N$, $n_0\geqslant 1$, such that for any integer $n$, $n\geqslant n_0$, there exists a section $s\in H^0(X,L^{\otimes n})$ satisfying 
\begin{equation}\label{eqn:extension:prop}
s|_Y=\ell^{\otimes n}\quad\text{and}\quad
\|s\|_{n\varphi}\leqslant\mathrm{e}^{\epsilon n}\|\ell\|_{\varphi|_Y}^{n}.
\end{equation}
If $\ell \in \mathcal N(Y, \rest{L}{Y})$, then $\ell$ has the extension property for the metric $\varphi$.
Indeed, by Proposition~\ref{prop:def:sup:norm} \ref{item: subadditivity},
there is a positive integer $n_0$ such that $\ell^{\otimes n} = 0$ for all integer $n \geqslant n_0$, so that if we choose
$s = 0 \in H^0(X, L^{\otimes n})$, then the above properties \eqref{eqn:extension:prop} hold.
In this sense, in order to check the extension property, we may assume that $\ell \not\in \mathcal N(Y, \rest{L}{Y})$.

For any non-zero element $\ell$ of $H^0(Y,L|_Y) {\setminus \mathcal N(Y, \rest{L}{Y})}$, we let
\begin{equation}\label{obstruction index}
\lambda_{\varphi}(\ell)=\limsup_{n\rightarrow+\infty} \frac{a_{\varphi, n}(\ell)}{n} \in[0,+\infty].
\end{equation} 
We call $\lambda_{\varphi}(\ell)$ the \emph{extension obstruction index}\index{extension obstruction index} of $\ell$.
\end{defi}

\begin{defi}
We assume that $H^0(X, L^{\otimes n}) \to H^0(Y, \rest{L}{Y}^{\otimes n})$
is surjective for all $n \geqslant 1$. Let $\|\ndot\|_{n\varphi,\mathrm{quot}}$
be the quotient seminorm of $H^0(Y, \rest{L}{Y}^{\otimes n})$ induced by
$\|\ndot\|_{n\varphi}$ and the surjective homomorphism
$H^0(X, L^{\otimes n}) \to H^0(Y, \rest{L}{Y}^{\otimes n})$.
For $\ell \in H^0(Y, \rest{L}{Y}^{\otimes n})$, 
we define $\|\ell\|_{\varphi,\mathrm{quot}}^{(n)}$ to be
\[
\|\ell\|_{\varphi,\mathrm{quot}}^{(n)} := \big(\|\ell^{\otimes n}\|_{n\varphi,\mathrm{quot}}\big)^{1/n}.
\]
It is easy to see that 
\begin{equation}
\| \ell\|_{\varphi,\mathrm{quot}}^{(\infty)} := \lim\limits_{n\to\infty} \| \ell\|_{\varphi,\mathrm{quot}}^{(n)} = \inf\limits_{n > 0} \| \ell\|_{\varphi,\mathrm{quot}}^{(n)} \in \mathbb R_{\geqslant 0}
\end{equation}
and
\begin{equation}
\| \ell \|_{\varphi|_Y} \leqslant \| \ell\|_{\varphi,\mathrm{quot}}^{(\infty)}
\end{equation}
because
\[
\|\ell_n \otimes \ell_{n'} \|_{(n+n')\varphi,\mathrm{quot}} \leqslant
\|\ell_n \|_{n\varphi,\mathrm{quot}} \|\ell_{n'} \|_{n'\varphi,\mathrm{quot}}\quad\text{and}\quad
\| \ell_n \|_{n \varphi|_Y} \leqslant \|\ell_n \|_{n\varphi,\mathrm{quot}}
\]
for all $\ell_n \in H^0(Y, \rest{L}{Y}^{\otimes n})$ and $\ell_{n'} \in H^0(Y, \rest{L}{Y}^{\otimes n'})$.
\if01
Therefore one can see that
\begin{enumerate}[label=(\alph*)]
\item $\lim\limits_{n\to\infty} \| \ell\|_{\varphi,\mathrm{quot}}^{(n)}$ exists and
$\lim\limits_{n\to\infty} \| \ell\|_{\varphi,\mathrm{quot}}^{(n)} = \inf\limits_{n > 0} \| \ell\|_{\varphi,\mathrm{quot}}^{(n)}$. We denote $\lim\limits_{n\to\infty} \| \ell\|_{\varphi,\mathrm{quot}}^{(n)}$ by $\| \ell\|_{\varphi,\mathrm{quot}}^{(\infty)}$.

\item $\| \ell \|_{\varphi|_Y} \leqslant \| \ell\|_{\varphi,\mathrm{quot}}^{(\infty)}$.
\end{enumerate}
\fi
\end{defi}

\begin{prop}\label{prop:extension:prop:equiv:cond}
We assume that $\ell \not\in \mathcal N(Y, \rest{L}{Y})$ and
there exists a positive integer $n_1$ such that, for all $n \geqslant n_1$,
$\ell^{\otimes n}$ lies in the image of the restriction map $H^0(X,L^{\otimes n})\rightarrow H^0(Y,\rest{L}{Y}^{\otimes n})$.
Then one has the following:
\begin{enumerate}[label=\rm(\arabic*)]
\item\label{Item: computation of lambda} ${\displaystyle \lambda_{\varphi}(\ell) = \lim_{n\rightarrow+\infty} \frac{a_{\varphi, n}(\ell)}{n} = \inf_{n \geqslant 1} \frac{a_{\varphi, n}(\ell)}{n}}$.

\item The following are equivalent:
\begin{enumerate}[label=(\arabic*)]
\renewcommand{\labelenumii}{(\arabic{enumi}.\alph{enumii})}
\item\label{Item: elle has ext property 1} $\ell$ has the extension property for $\varphi$.

\item\label{Item: elle has ext property 2} For any $\epsilon > 0$, there are a positive integer $n$ and a section $s \in H^0(X, L^{\otimes n})$ such that
$\rest{s}{Y} = \ell^{\otimes n}$ and $\|s\|_{n\varphi}\leqslant\mathrm{e}^{\epsilon n}\|\ell\|_{\varphi|_Y}^{n}$.

\item\label{Item: elle has ext property 3} $\lambda_{\varphi}(\ell) = 0$.
\end{enumerate}
\item\label{Item: equivalent assume surjectivite} We assume that $H^0(X, L^{\otimes n}) \to H^0(Y, \rest{L}{Y}^{\otimes n})$
is surjective for all $n \geqslant 1$. Then, 
the above equivalent properties are also equivalent to 
$\| \ell \|_{\varphi|_Y} = \| \ell\|_{\varphi,\mathrm{quot}}^{(\infty)}$. 
\end{enumerate}
\end{prop}

\begin{proof}
\ref{Item: computation of lambda} is a consequence of Fekete's lemma because
the sequence $\{a_{\varphi,n}\}_{n\in\mathbb N}$ is subadditive by Lemma~\ref{lemm:properties:a:varphi:n:l}. 

\medskip``\ref{Item: elle has ext property 1} $\Longrightarrow$ \ref{Item: elle has ext property 2}'' is obvious.

\medskip``\ref{Item: elle has ext property 2} $\Longrightarrow$ \ref{Item: elle has ext property 3}'':
For any $\epsilon > 0$, there is a positive integer $n$
such that $a_{\varphi, n}(\ell) \leqslant n \epsilon$, so that, by \ref{Item: computation of lambda},
\[
0 \leqslant \lambda_{\varphi}(\ell) = \inf_{n \geqslant 1} \frac{a_{\varphi, n}(\ell)}{n} \leqslant \epsilon,
\]
and hence one has \ref{Item: elle has ext property 3}. 

\medskip``\ref{Item: elle has ext property 3} $\Longrightarrow$ \ref{Item: elle has ext property 1}'':
Since ${\displaystyle \lambda_{\varphi}(\ell) = \lim_{n\rightarrow+\infty} \frac{a_{\varphi, n}(\ell)}{n}}$ by \ref{Item: computation of lambda}, 
we can see \ref{Item: elle has ext property 1}.

\medskip 
\ref{Item: equivalent assume surjectivite} Note that $a_{\varphi, n}(\ell)/n = \ln \|\ell\|_{\varphi,\mathrm{qout}}^{(n)} - \ln \|\ell\|_{\varphi|_Y}$, so that $\lambda_{\varphi}(\ell) = \ln \|\ell\|_{\varphi,\mathrm{qout}}^{(\infty)} - \ln \|\ell\|_{\varphi|_Y}$. Thus the assertion follows.
\end{proof}

\begin{rema}\phantomsection\label{Rem:extension}
\begin{enumerate}[label=\rm(\arabic*)] 
\item\label{Item: continuity of a n}
Let $\varphi'$ be another metric on $L$. By Lemma~\ref{lemm:properties:a:varphi:n:l}  
one has
\[|a_{\varphi,n}(\ell)-a_{\varphi',n}(\ell)|\leqslant 2n\,d(\varphi,\varphi'),\]
provided that $a_{\varphi,n}(\ell)$ or $a_{\varphi',n}(\ell)$ is finite. We deduce from this inequality that, $\lambda_{\varphi}(\ell)$ is finite if and only if $\lambda_{\varphi'}(\ell)$ is finite. Moreover, when these numbers are finite, one has
\begin{equation}\label{distancelambda}|\lambda_{\varphi}(\ell)-\lambda_{\varphi'}(\ell)|\leqslant 2 d(\varphi,\varphi').\end{equation}

\item\label{Item: dialation lambda}
We assume that $\lambda_{\varphi}(\ell) < \infty$. Then one has
\begin{equation}\label{lambda:linearity}
\lambda_{n\varphi}(\ell^{\otimes n}) = n \lambda_{\varphi}(\ell)
\end{equation}
for all $n > 0$. Indeed,
\[
\lambda_{n\varphi}(\ell^{\otimes n}) = \lim_{m\to\infty} \frac{a_{n\varphi, m}(\ell^{\otimes n})}{m}
= \lim_{m\to\infty} \frac{a_{\varphi, nm}(\ell)}{m} = n \lim_{m\to\infty} \frac{a_{\varphi, nm}(\ell)}{nm} = n \lambda_{\varphi}(\ell).
\]

\item\label{Item: comparison lambda restriction} Let $X'$ be a closed subscheme of $X$ such that $Y \subseteq X'$.
We assume that there is a positive integer $n_0$ such that, for all $n \geqslant n_0$,
$H^0(X, L^{\otimes n}) \to H^0(X', \rest{L^{\otimes n}}{X'})$ is surjective.
Then 
\begin{equation}\label{lambda:comp:subscheme}
\lambda_{\rest{\varphi}{X'}}(\ell) \leqslant \lambda_{\varphi}(\ell).
\end{equation}
Indeed, as $\| s|_{X'} \|_{\rest{\varphi}{X'}} \leqslant \| s \|_{\varphi}$
for all $s \in H^0(X, L^{\otimes n})$ and $H^0(X, L^{\otimes n}) \to H^0(X', \rest{L^{\otimes n}}{X'})$ is surjective for all $n \geqslant n_0$,
one has $a_{\rest{\varphi}{X'}, n}(\ell) \leqslant a_{\varphi, n}(\ell)$ for all $n \geqslant n_0$, so that
the assertion follows.
\end{enumerate}
\end{rema}

\subsubsection{A generalisation of a result in \cite{Zhang95} and \cite{Moriwaki2014}}

Let $X$ be a $d$-dimensional integral smooth scheme over $\mathbb C$. Let $Y$ be a closed and reduced subscheme of $X$ defined by an ideal sheaf $I$ on $X$, that is,
$I = \sqrt{I}$ and $Y = \Spec (\mathcal O_X/I)$.
Let $\mu_I : X_I \to X$ be the blowing-up along $I$, that is,
$X_I = \operatorname{Proj}\left( \bigoplus_{m=0}^{\infty} I^m \right)$.
Let $\tilde{\mu} : \widetilde{X}_I \to X_I$ be the normalisation of $X_I$.
Furthermore, let $\mu' : X' \to \widetilde{X}_I$ be a desingularisation of $\widetilde{X}_I$ such that
$\mu'$ yields an isomorphism 
\[
X' \setminus {\mu'}^{-1}(\operatorname{Sing}(\widetilde{X}_I))
\overset{\sim}{\longrightarrow} \widetilde{X}_I \setminus \operatorname{Sing}(\widetilde{X}_I).
\]
We denote the compositions 
\[
X' \overset{\mu'}{\longrightarrow} \widetilde{X}_I \overset{\tilde{\mu}}{\longrightarrow} X_I \overset{\mu_I}{\longrightarrow} X.
\]
by $\mu$, that is, $\mu := \mu_I \circ \tilde{\mu} \circ \mu'$.
Note that $X' \setminus \mu^{-1}(Y) \overset{\sim}{\longrightarrow} X \setminus Y$ via $\mu$.

\begin{lemm}\label{lemm:mult:blowing:normalization:desingularization}
There are positive integers $m_0$ and $c$ such that 
$\mu_*(I^m \mathcal O_{X'}) \subseteq I^{m -c}$ for all $m \geqslant m_0 + c$.
\end{lemm}

\begin{proof}
Let us consider the following claim:

\begin{enonce}{Claim}
\begin{enumerate}[label=\rm(\alph*)]
\item\label{Item: blow up lower star} 
$\mu'_*(I^m \mathcal O_{X'}) = I^m \mathcal O_{\widetilde{X}_I}$ for all integer $m \geqslant 0$.

\item\label{Item: exists m grater than c such that}
There is a positive integer $c$ such that
$\tilde{\mu}_{*}(I^m \mathcal O_{\widetilde{X}_I}) \cap \mathcal O_{X_I} \subseteq I^{m-c}\mathcal O_{X_I}$ for all integer $m \geqslant c$.

\item\label{Item: exists m0 such that push forward}
There is a positive integer $m_0$ such that
$\mu_{I, *}(I^m \mathcal O_{X_I}) = I^m$ for all integer $m \geqslant m_0$.
\end{enumerate}
\end{enonce}

\begin{proof}
\ref{Item: blow up lower star} Note that $I^m \mathcal O_{\widetilde{X}_I}$ is invertible and $I^m \mathcal O_{X'} = {\mu'}^*(I^m \mathcal O_{\widetilde{X}_I})$.
Moreover as $\widetilde{X}_I$ is normal, $\mu'_*(\mathcal O_{X'}) = \mathcal O_{\widetilde{X}_I}$, so that
the assertion follows from the projection formula.

\medskip\ref{Item: exists m grater than c such that} We choose an affine open covering $X_I = \bigcup_{i=1}^N \Spec(A_i)$.
Let $\widetilde{A}_i$ be the normalisation of $A_i$. Then
$\widetilde{X}_I = \bigcup_{i=1}^N \Spec(\widetilde{A}_i)$ is an affine open covering.
Note that $\widetilde{A}_i$ is a finitely generated $A_i$-module, so that, by Artin-Lees lemma (cf. \cite[Corollary~10.10]{Atiyah_Macdonald1964}),
there is a positive constant $c_i$ such that $I^m \widetilde{A}_i \cap A_i = I^{m-c_i} (I^c \widetilde{A}_i \cap A_i)$
for all $m \geqslant c_i$, which implies $I^m \widetilde{A}_i \cap A_i \subseteq  I^{m-c_i} A_i$. Therefore, if we set
$c = \max \{ c_1, \ldots, c_N \}$, then one has the assertion.

\medskip\ref{Item: exists m0 such that push forward} This is essentially proved in \cite[Chapter~II, Theorem~5.19]{Hart77}. Indeed, at the final line in the proof of the above reference, 
it says that ``$S'_n = S_n$ for all sufficiently large $n$'', which is nothing more than the assertion of \ref{Item: exists m0 such that push forward}
because $\mathcal O_{X_I}(m) = I^m \mathcal O_{X_I}$.
\end{proof}

Let us go back to the proof of the lemma.
This is a local question, so that we may assume that $X = \operatorname{Spec}(A)$ for some
finitely generated regular $\mathbb C$-algebra $A$.
Note that $\mu_*(I^m \mathcal O_{X'}) \subseteq \mathcal O_X$. Therefore, it is sufficient to see that,
if $f \in I^m \mathcal O_{X'}$ for $f \in A$, then $f \in I^{m-c}$. First of all, by \ref{Item: blow up lower star},
$f \in I^m \mathcal O_{\widetilde{X}_I}$, so that $f \in \tilde{\mu}_{*}(I^m \mathcal O_{\widetilde{X}_I}) \cap \mathcal O_{X_I}$, and
hence, by \ref{Item: exists m grater than c such that}, $f \in I^{m-c} \mathcal O_{X_I}$. Note that $m -c \geqslant m_0$.
Therefore, one has $f \in I^{m-c}$, as required.
\end{proof}

We assume that $X$ is projective. 
Let $L$ be an ample invertible $\mathcal O_X$-module and $\varphi$ be a $C^{\infty}$-metric of $L$ such that $c_1(L, \varphi)$ is positive.
Let $U$ be an open set (in the sense of the classical topology) of $X$
such that $Y \subseteq U$. 
The proof of \cite[Theorem~7.6]{Moriwaki2014} works well even if we change the exponent $d$ of $\rho$ by
a positive number $\delta$ except (3) in Claim~2, which should be
\[
\text{``If $\delta \geqslant d$, then $\rho^{-\delta}$ is not integrable on any neighborhood of $Y$''.}
\]
At page 231, line 6 from the bottom, one constructs a $C^{\infty}$-section $l'$ of $L^{\otimes n}$ over $X$, which is holomorphic on $U'$ and 
satisfies the integrability condition
\begin{equation}\label{eqn:lem:delta:mult:Y:00}
\int_X |l'|^2 \rho^{-\delta} \Phi < \infty.
\end{equation}

Let $E_1, \ldots, E_r$ be irreducible components of $\mu^{-1}(Y)$.
We set 
\[
I \mathcal O_{X'} = -(a_1 E_1 + \cdots + a_r E_r)\quad\text{and}\quad
K_{X'} = \mu^*(K_X) + b_1 E_1 + \cdots + b_r E_r.
\]
Note that $a_i > 0$ and $b_i > 0$ for all $i$. 

\begin{lemm}\label{lemm:mult:along:Ei:lower:bound}
If $e_i$ is the multiplicity of $l'$ along $E_i$, then
$e_i > a_i \delta - b_i - 1$ for $i=1, \ldots, r$.
\end{lemm}

\begin{proof}
Let $\eta$ be a closed point of $E_i \setminus \operatorname{Sing}(E_1 + \cdots + E_r)$ and $\xi = \mu(\eta)$.
Let $y_1, \ldots, y_d$ be a local coordinate of $X'$ on an open neighborhood $W$ of $\eta$ such that $E_i$ is defined by $y_1 = 0$.
Let $x_1, \ldots, x_d$ be a local coordinate of $X'$ on an open neighborhood $V$ of $\xi$.
In the following, if it is necessary, we will shrink $V$ and $W$ freely.
First of all, we may assume that
\begin{equation}\label{eqn:lem:delta:mult:Y:01}
\mu(W) \subseteq V.
\end{equation}
Moreover, we can find a positive constant $C$ such that 
\begin{equation}\label{eqn:lem:delta:mult:Y:02}
\Phi \geqslant C (\sqrt{-1})^d (dx_1 \wedge d\bar{x}_1) \wedge \cdots \wedge (dx_d \wedge d\bar{x}_d)
\end{equation}
on $V$. Let $\omega$ be a local basis of $L$ on $V$.
Note that $\rho$ can be writen by
\[
\rho = (|f_1|^2 + \cdots + |f_N|^2)|\omega^m|
\]
on $V$, where $f_1, \ldots, f_N$ are generators of $I$ on $V$, so that
there is a positive constant $C'$ such that 
\begin{equation}\label{eqn:lem:delta:mult:Y:03}
\mu^*(\rho) \leqslant C' |y_1|^{2a_i}
\end{equation}
on $W$. Further one has
\begin{multline*}
\mu^*((\sqrt{-1})^d (dz_1 \wedge d\bar{z}_1) \wedge \cdots \wedge (dz_d \wedge d\bar{z}_d)) \\
= |y_1|^{2b_i}|u|^2(\sqrt{-1})^d (dy_1 \wedge d\bar{y}_1) \wedge \cdots \wedge (dy_d \wedge d\bar{y}_d)
\end{multline*}
on $W$, where $u$ is a nowhere vanishing holomorphic function $W$, so that
\begin{multline}\label{eqn:lem:delta:mult:Y:04}
\mu^*((\sqrt{-1})^d (dz_1 \wedge d\bar{z}_1) \wedge \cdots \wedge (dz_d \wedge d\bar{z}_d)) \\
\geqslant C'' |y_1|^{2b_i}(\sqrt{-1})^d (dy_1 \wedge d\bar{y}_1) \wedge \cdots \wedge (dy_d \wedge d\bar{y}_d)
\end{multline}
holds on $W$ for some positive constant $C''$.
If we set $l' = f' \omega^n$, then $f' = y_1^{e_i} g$ on $W$ such that $g$ is not identically zero on $\rest{E_i}{W}$, so that
one can find $(0, \alpha_2, \ldots, \alpha_n) \in \rest{E_i}{W}$ and a positive number $r$ such that
$g \not= 0$ on 
\[
W_r = \{ ({y_1, \ldots, y_r}) \in W \,:\, 
\text{$|{y}_j - \alpha_j| \leqslant r$ for all $j=1, \ldots, d$} \},
\]
where $\alpha_1 = 0$.
Therefore one can find a positive constant $C'''$ such that 
\begin{equation}\label{eqn:lem:delta:mult:Y:05}
\mu^*(|l'|^2) \geqslant C''' |y_1|^{2e_i}
\end{equation}
on $W_r$.
Thus, if we set $y_j - \alpha_j = r_j \exp(\sqrt{-1}\theta_j)$ for $j=1, \ldots, d$, then,
by \eqref{eqn:lem:delta:mult:Y:00}, \eqref{eqn:lem:delta:mult:Y:01}, \eqref{eqn:lem:delta:mult:Y:02}, \eqref{eqn:lem:delta:mult:Y:03}, \eqref{eqn:lem:delta:mult:Y:04} and \eqref{eqn:lem:delta:mult:Y:05},
{\allowdisplaybreaks
\begin{align*}
\infty & > \int_X |l'|^2 \rho^{-\delta} \Phi \geqslant  \int_{V} |l'|^2 \rho^{-\delta} \Phi \\
& \geqslant
\int_{V} C |l'|^2 \rho^{-\delta} (\sqrt{-1})^d (dx_1 \wedge d\bar{x}_1) \wedge \cdots \wedge (dx_d \wedge d\bar{x}_d) \\
& \geqslant C \int_{W_r} \mu^*(|l'|^2 \rho^{-\delta} (\sqrt{-1})^d (dx_1 \wedge d\bar{x}_1) \wedge \cdots \wedge (dx_d \wedge d\bar{x}_d)) \\
& \geqslant C {C'}^{-\delta} C'' C'''\int_{W_r} |y_1|^{2 e_i + 2b_i -2a_i \delta}(\sqrt{-1})^d (dy_1 \wedge d\bar{y}_1) \wedge \cdots \wedge (dy_d \wedge d\bar{y}_d) \\
& \geqslant C {C'}^{-\delta} C'' C'''\int_{[0,r]^d \times [0, 2\pi]^d} r_1^{2 e_i + 2b_i -2a_i \delta + 1} r_2 \cdots r_d dr_1 \cdots dr_d d\theta_1 \cdots d \theta_d \\
& \geqslant C {C'}^{-\delta} C'' C'''2 \pi^d r^{2(d-1)} \int_{[0,r]} r_1^{2 e_i + 2b_i -2a_i \delta + 1} dr_1,
\end{align*}}
so that $2 e_i + 2b_i -2a_i \delta + 1 > -1$, as required.
\end{proof}

By virtue of Lemma~\ref{lemm:mult:along:Ei:lower:bound} together with Lemma~\ref{lemm:mult:blowing:normalization:desingularization}, one has the following generalisation of
\cite[Lemma~2.6]{Zhang95} and \cite[Theorem~7.6]{Moriwaki2014}.

\begin{theo}\label{thm:extension:norm:control:non:reduced:Y}
Let $Y'$ be a closed subscheme of $X$ defined by an ideal sheaf $J$ on $X$, that is, $Y' = \Spec (\mathcal O_X/J)$. Let $U$ be a Zariski open subset of $X$ containing $Y'$.
Then there are a positive integer $n_0$ and a positive constant $C$ such that,
for any integer $n \geqslant n_0$ and $l_U \in H^0(U, L^{\otimes n})$ with
$\|l_U\|_{n \varphi_U} < \infty$,
there is $l \in H^0(X, L^{\otimes n})$ such that $\rest{l}{Y'} = \rest{l_U}{Y'}$ and
$\| l \|_{n\varphi} \leqslant C n^{2d} \|l_U\|_{n \varphi_U}$, where
$\|l_U\|_{n \varphi_U} :=\sup \{ |l_U|_{n \varphi}(x) \,:\, x \in U^{\mathrm{an}}\}$.
\end{theo}

\begin{proof}
Let $I := \sqrt{J}$ and $Y := \Spec(\mathcal O_X/I)$. 
One can find $a \in \mathbb Z_{\geqslant 1}$ such that $I^{a} \subseteq J$.
We fix a positive number $\delta$ with 
\[
\delta \geqslant \max_{i = 1, \ldots, r} \left\{ \frac{b_i+1}{a_i} \right\} + c + m_0 + a,
\]
where $c$ and $m_0$ are the positive integers in Lemma~\ref{lemm:mult:blowing:normalization:desingularization}. 
The proof of \cite[Theorem~7.6]{Moriwaki2014} is carried out by using the exponent $\delta$ 
instead of $d$.
The point is to show that $\rest{l}{Y'} = \rest{l_U}{Y'}$.
By Lemma~\ref{lemm:mult:along:Ei:lower:bound},
\[
e_i > a_i \delta - b_i - 1 \geqslant a_i(c + m_1 + a),
\]
so that, there is a Zariski closed set $Z$ of $\mu^{-1}(Y)$ such that
$\dim Z \leqslant d-2$ and
\[
\mu^*(l') \in \rest{I^{c + m_1 + a} \mu^*(L^{\otimes n})}{\mu^{-1}(U') \setminus Z}.
\]
As $I^{c + m_1 + a} \mu^*(L^{\otimes n})$ is invertible,
one can see that $\mu^*(l') \in \rest{I^{c + m_1 + a} \mu^*(L^{\otimes n})}{\mu^{-1}(U')}$, and hence,
by Lemma~\ref{lemm:mult:blowing:normalization:desingularization},
$l' \in I^{m_1 + a} L^{\otimes n} \subseteq I^{a} L^{\otimes n}$.
Therefore the class of $l'$ in  $L^{\otimes n}/JL^{\otimes n}$ is zero over $Y'$,
and hence $\rest{l}{Y'} = \rest{l_U}{Y'}$.
The remaining estimates are same as the proof of \cite[Theorem~7.6]{Moriwaki2014}.
\end{proof}

\begin{coro}\label{cor:extension:norm:control:non:reduced:Y}
Let $X$, $L$ and $Y'$ be the same as in Theorem~\ref{thm:extension:norm:control:non:reduced:Y}.
Let $\varphi$ a continuous metric of $L$ such that the first Chern current 
$c_1(L, \varphi)$ is positive. Then, for $\ell \in H^0(Y', \rest{L}{Y'})$,
$\ell$ has the extension property for $\varphi$.
\end{coro}

\begin{proof}
Clearly we may assume that $\ell \not\in \mathcal N(Y', \rest{L}{Y'})$.
As $L$ is ample, there is a $C^{\infty}$-metric $\psi$ on $L$ such that
$c_1(L, \psi)$ is a positive form.

\begin{enonce}{Claim}
If the corollary holds for any $C^{\infty}$-metric of $L$ with the semipositive Chern form,
then the corollary holds in general.
\end{enonce}

\begin{proof}
Let $\phi$ be a continuous function such that $\psi - \varphi = \phi$.
It is well known that there is a sequence $\{ \phi_{n} \}_{n=1}^{\infty}$
of $C^{\infty}$-functions on $X^{\mathrm{an}}$ such that
$\varphi_n := \psi - \phi_n$ is a $C^{\infty}$-metric of $L$ with the semipositive Chern form and
$\{ \phi_n \}_{n=1}^{\infty}$ converges uniformly to $\phi$ (for example, see \cite[Theorem~1]{BlockiKolodziej2007} or \cite[Lemma~4.2]{Moriwaki2012}).
Thus $\lim_{n\to\infty} d(\varphi, \varphi_n) = 0$.
By our assumption, $\lambda_{\varphi_n}(\ell) = 0$, so that
$\lambda_{\varphi}(\ell) \leqslant 2 d(\varphi, \varphi_n)$, and hence
the assertion follows.
\end{proof}

We fix a positive number $\epsilon$.
By the above claim, we may assume that $\varphi$ is $C^{\infty}$, so that
if we set $f = \psi - \varphi$, then $f$ is a $C^{\infty}$ function.
Note that for $\lambda \in \intervalle]01[$, $\varphi + \lambda f$ gives rise to a positive Chern form
because $\varphi + \lambda f = (1-\lambda)\varphi + \lambda \psi$.
We choose $\lambda_0 \in \intervalle]01[$ such that 
\[
\lambda_0 \sup \{ |f(x)| \,:\, x \in X^{\mathrm{an}} \} \leqslant \epsilon.
\]
We set $\varphi' = \varphi + \lambda_0 f$. Then
\begin{equation}\label{eqn:cor:extension:norm:control:non:reduced:Y:01}
\mathrm{e}^{-\epsilon}|\ndot|_{\varphi}(x) \leqslant |\ndot|_{\varphi'}(x) \leqslant \mathrm{e}^{\epsilon} |\ndot|_{\varphi}(x)
\end{equation}
for all $x \in X^{\mathrm{an}}$.
We choose a positive integer $a$ such that 
\[
H^0(X, L^{\otimes a}) \to H^0(Y', \rest{L^{\otimes a}}{Y'})
\]
is surjective, so that one can find $t$ such that $t \in H^0(X, L^{\otimes a})$ and
$\rest{t}{Y'} = \ell^{\otimes a}$.
We also choose an open set $U$ of $X$ such that $Y' \subseteq U$ and
\begin{equation}\label{eqn:cor:extension:norm:control:non:reduced:Y:02}
\| t \|_{a\varphi'_U} \leqslant \mathrm{e}^{a\epsilon} \| l^{\otimes a} \|_{\rest{a\varphi'}{Y'}}.
\end{equation}
By the above theorem, there are a positive integer $n_1$ and a positive constant $C$ such that,
for any $n \geqslant n_1$, one can find $s \in H^0(X, L^{\otimes a n})$ such that
$\rest{s}{Y'} = \ell^{\otimes a n}$ and 
\begin{equation}\label{eqn:cor:extension:norm:control:non:reduced:Y:03}
\| s \|_{na \varphi'} \leqslant C n^{2d} \| t^{\otimes n}\|_{na\varphi'_U}.
\end{equation}
Let $n_2$ be a positive integer such that $n_2 \geqslant n_1$ and
\begin{equation}\label{eqn:cor:extension:norm:control:non:reduced:Y:04}
C n^{2d} \leqslant \mathrm{e}^{\epsilon an}
\end{equation}
for $n \geqslant n_2$. Therefore, using \eqref{eqn:cor:extension:norm:control:non:reduced:Y:01},
\eqref{eqn:cor:extension:norm:control:non:reduced:Y:02},
\eqref{eqn:cor:extension:norm:control:non:reduced:Y:03} and
\eqref{eqn:cor:extension:norm:control:non:reduced:Y:04},
one has
\begin{align*}
\| s \|_{na\varphi} & \leqslant \mathrm{e}^{na\epsilon} \| s \|_{na\varphi'}  \leqslant 
\mathrm{e}^{na\epsilon} \left( C n^{2d} \| t\|_{a\varphi'_U}^n \right)  \\
& \leqslant
\mathrm{e}^{2na\epsilon} \| t\|_{a\varphi'_U}^n
\leqslant \mathrm{e}^{3na\epsilon} \| l \|_{\rest{\varphi'}{Y'}}^{an} \leqslant \mathrm{e}^{4na\epsilon} \| l \|_{\rest{\varphi}{Y'}}^{an},
\end{align*}
which means that $\lambda_{a\varphi}(l^{\otimes a}) \leqslant 4a \epsilon$, so that
$\lambda_{\varphi}(l) \leqslant 4 \epsilon$. Therefore one has $\lambda_{\varphi}(l) = 0$ because
$\epsilon$ is an arbitrary positive number.
\end{proof}

\subsubsection{Extension property over an Archimedean field}
\label{Extension property over an Archimedean field}

We assume that $k$ is either $\mathbb R$ or $\mathbb C$ and
the absolute value of $k$ is the standard absolute value.

\begin{theo}\label{thm:extension:property:over:Archimedean}
Let $X$ be a projective scheme over $k$, $L$ be an ample invertible $\mathcal O_X$-module and
$\varphi$ be a semipositive continuous metric metric of $L$.
For any closed subscheme $Y$ of $X$ and any 
$\ell \in H^0(Y, \rest{L}{Y})$,
$\ell$ has the extension property for $\varphi$.
\end{theo}

\begin{proof}
Clearly we may assume that $\ell \not\in \mathcal N(Y, \rest{L}{Y})$.
Let us see the following claim:

\begin{enonce}{Claim}\label{claim:thm:extension:property:over:Archimedean:01}
\begin{enumerate}[label=\rm(\arabic*)]
\item\label{Item: FS has extension property}
We assume that $k = \mathbb C$, $X = \mathbb P^n$, $L = \mathcal O(1)$ 
and $\varphi$ is the Fubini-Study metric arising from
a norm $\|\ndot\|$ on $H^0(\mathbb P^n, \mathcal O(1))$.
Then the assertion of the theorem holds.

\item\label{Item: FS extension R}
We assume that $k = \mathbb R$, $X = \mathbb P^n$, $L = \mathcal O(1)$ 
and $\varphi$ is the Fubini-Study metric arising from
a norm $\|\ndot\|$ on $H^0(\mathbb P^n, \mathcal O(1))$.
Then the assertion of the theorem holds.

\end{enumerate}
\end{enonce}

\begin{proof}
\ref{Item: FS has extension property} By Theorem~\ref{thm:semipositive:positive:current}, the first Chern current $c_1(L, \varphi)$ is positive, so that (1) is a consequence of Corollary~\ref{cor:extension:norm:control:non:reduced:Y}.

\medskip
\ref{Item: FS extension R} We consider $X_{\mathbb C}$, $L_{\mathbb C}$ and $\varphi_{\mathbb C}$.
Then $\varphi_{\mathbb C}$ is the Fubini-Study metric induced by the norm $\|\ndot\|_{\mathbb C}$
on $H^0(X, L) \otimes_{\mathbb R} \mathbb C$ by Proposition~\ref{Pro:quotientnormextensionscalar}. Thus,
by using the case (1),
for any $\epsilon > 0$, there is a positive integer $n_0$ such that, for any $n \geqslant n_0$,
we can find $s \in H^0(X_{\mathbb C}, L_{\mathbb C}^{\otimes n})$ with
$\rest{s}{Y_{\mathbb C}} = \ell^{\otimes n}$ and $\| s\|_{n\varphi_{\mathbb C}} \leqslant \mathrm{e}^{n\epsilon} \| \ell \|_{\rest{\varphi_{\mathbb C}}{Y}}^n$. 
First of all, note that $\| \ell \|_{\rest{\varphi_{\mathbb C}}{Y}} = \| \ell \|_{\rest{\varphi}{Y}}$. If we set $s = \sigma + i \tau$ ($\sigma, \tau \in H^0(X, L^{\otimes n})$),
then $(\rest{\sigma}{Y}) + i (\rest{\tau}{Y}) = \ell^{\otimes n}$, so that $\rest{\tau}{Y} = 0$, and hence
$\rest{\sigma}{Y} = \ell^{\otimes n}$. Moreover, for any $x\in X^{\mathrm{an}}_{\mathbb C}$ one has $\|\sigma+i\tau\|_{n\varphi_{\mathbb C}}(x)=\|\sigma-i\tau\|_{n\varphi_{\mathbb C}}(\overline x)$, so that $\|\sigma+i\tau\|_{n\varphi_{\mathbb C}}=\|\sigma-i\tau\|_{n\varphi_{\mathbb C}}$. We then deduce that \[2\|\sigma\|_{n\varphi}=2\|\sigma_{n\varphi_{\mathbb C}}\|\leqslant\|\sigma+i\tau\|_{n\varphi_{\mathbb C}}+\|\sigma-i\tau\|_{n\varphi_{\mathbb C}}=2\|s\|_{n\varphi_{\mathbb C}}.\]
Thus one has the assertion in this case.
\end{proof}

We choose $n_1$ such that, for all $n \geqslant n_1$,
$L^{\otimes n}$ is very ample.
Then $\varphi_n$ is the restriction of the Fubini-Study metric $\varphi_{\|\ndot\|_{n\varphi}}$
of $\mathcal O(1)$ {to} $\mathbb P(H^0(X, L^{\otimes n}))$ induced by the norm $\|\ndot\|_{n\varphi}$.
Thus, by the above claim together with \eqref{lambda:comp:subscheme}, 
\[
0 \leqslant \lambda_{\varphi_n}(\ell^{\otimes n}) \leqslant \lambda_{\varphi_{\|\ndot\|_{n\varphi}}}(\ell^{\otimes n}) = 0,
\]
and hence
$\lambda_{n\varphi}(\ell^{\otimes n}) \leqslant 2d(n \varphi, \varphi_n)$ by \eqref{distancelambda}.
Since $\lambda_{n\varphi}(\ell^{\otimes n}) = n \lambda_{\varphi}(\ell)$ by \eqref{lambda:linearity}, one has
\[
0 \leqslant \lambda_{\varphi}(\ell) \leqslant 2d(\varphi, \textstyle{\frac 1n}\varphi_n).
\]
Therefore, the assertion follows.
\end{proof}

\subsubsection{Extension property over a non-Archimedean field}
\label{Extension property over a non-Archimedean field}
In this subsection, we fix a field $k$ equipped with a \emph{non-Archimedean} absolute value $|\ndot|$, under which the field $k$ is complete. 

\begin{prop}\label{Pro:relaemvement} 
We assume that $|\ndot|$ is non-trivial.
Let $X$ be a projective $k$-scheme, $L$ be an invertible $\mathcal O_X$-module and $(\mathscr X,\mathscr L)$ be a model of $(X,L)$. Let $s$ be a global section of $L$ such that $\|s\|_{\varphi_{\mathscr L}}\leqslant 1$. Then there exists an element $a\in\mathfrak o_k\setminus\{0\}$ such that $as^n$ belongs to $H^0(\mathscr X,\mathscr L^{\otimes n})$ for all integers $n\geqslant 1$. 
\end{prop}
\begin{proof}
Let $(\mathscr U_i)_{i=1}^N$ be a covering of $\mathscr X$ by affine open subsets, such that $\mathscr U_i = \Spec(\mathscr A_i)$ and the invertible sheaf $\mathscr L$ trivialises on each $\mathscr U_i$, that is,
$\rest{\mathscr L}{\mathscr U_i} = \mathscr A_i s_i$ for some $s_i \in \mathscr L(\mathscr U_i)$. 
Then the restriction of $s$ {to} $U_i:=\mathscr U_i\cap X$ can be written in the form $\lambda_is_i$, where $\lambda_i\in A_i = S^{-1}\mathscr A_i$, 
$S=\mathfrak o_k\setminus\{0\}$. 
Note that for $x\in (U_i)^{\mathrm{an}}_{\mathscr A_i}$, the reduction point of $x$ is in $\mathscr U_i$
(cf. Remark~\ref{remark:reduction:point}), so that
since $\|s\|_{\varphi_{\mathscr L}}\leqslant 1$, we obtain that  
\[
|\lambda_i|_x = |\lambda_i|_x\cdot|s_i|_{\varphi_{\mathscr L}}(x) = |s|_{\varphi_{\mathscr L}}(x) \leqslant \|s\|_{\varphi_{\mathscr L}}\leqslant 1
\]
for any $x\in (U_i)^{\mathrm{an}}_{\mathscr A_i}$.  
By Proposition \ref{Pro:integralclosure}, $\lambda_i$ is integral over the ring 
$\mathscr A_i$,  
namely $\mathscr A_i[\lambda_i]$ is an $\mathscr A_i$-module of finite type. In particular, there exists an integer $d_i\geqslant 1$ such that, for any integer $n\geqslant 1$,
\[\lambda_i^n\in\mathscr A_i +\mathscr A_i \lambda_i+\cdots+\mathscr A_i \lambda_i^{d_i}.\] 
Moreover,
there exists $a_i\in \mathfrak o_k\setminus\{0\}$ such that \[\{a_i\lambda_i,\ldots,a_i\lambda_i^{d_i}\}\subset\mathscr A_i.\] 
We then obtain that $a_i\lambda_i^n\in\mathscr A_i$ for any integer $n\geqslant 1$. Finally, let $a=\prod_{i=1}^Na_i\in\mathfrak o_k\setminus\{0\}$. For any integer $n\geqslant 1$ and any $i\in\{1,\ldots,N\}$, one has \[(as^n)|_{U_i}=(a\lambda_i^n)s_i^n\in H^0(\mathscr U_i,\mathscr L^{\otimes n}).\]
Since $\mathscr X$ is flat over $\Spec(\mathfrak o_k)$, these sections glue together to be a global section of $\mathscr L^{\otimes n}$. The proposition is thus proved. 
\end{proof}

\begin{prop}\phantomsection\label{Pro:liftingformodel}
We assume that $|\ndot|$ is non-trivial.
Let $X$ be a projective $k$-scheme, $L$ be an ample invertible $\mathcal O_X$-module and $Y$ be a closed subscheme of $X$. Let $u\geqslant 1$ be an integer and $(\mathscr X,\mathscr L)$ be a model of $(X,L^{\otimes u})$ such that $\mathscr L$ is ample. Assume that $\varphi=\frac 1u\varphi_{\mathscr L}$. Let $\varphi_Y$ be the restriction of the metric $\varphi$ {to} $L|_Y$. For any positive number $\epsilon$ and any $\ell\in H^0(Y,L|_Y)$, 
there exists an integer $n\geqslant 1$ and a section $s\in H^0(X,L^{\otimes n})$ such that  $s|_Y=\ell^n$ and 
\[\|s\|_{n\varphi}\leqslant\mathrm{e}^{n\epsilon}\|\ell\|_{\varphi_Y}^n.\]
In other words, one has $\lambda_{\varphi}(\ell) =0$ if $\ell \not\in  N(Y, \rest{L}{Y})$. 
\end{prop}
\begin{proof}
 We choose a positive integer $m$ such that 
\[
\mathrm{e}^{-m\epsilon/2} < \sup \{ |a| \,:\, a \in k^{\times},\ |a| < 1 \}.
\]
By Proposition \ref{Pro:dilatation}, on $H^0(Y, \rest{L}{Y}^{\otimes m})/\mathcal N(Y, \rest{L}{Y}^{\otimes m})$,
there is $\alpha \in k^{\times}$ such that 
\begin{equation}\label{Equ:encadreentdialat}\mathrm{e}^{-m\epsilon/2}\leqslant\|\alpha\ell^m\|_{m\varphi_Y}\leqslant 1.\end{equation}
Let $\mathscr Y$ be the Zariski closure of $Y$ in $\mathscr X$. By Proposition \ref{Pro:relaemvement}, there exists an element $\beta\in\mathfrak o_k\setminus\{0\}$ such that \[\beta(\alpha\ell^m)^{pu}\in H^0(\mathscr Y,\mathscr L^{\otimes mp}|_{\mathscr Y})\]
for any integer $p\geqslant 1$. Moreover, since the invertible sheaf $\mathscr L$ is ample, for sufficiently positive integer $p$,  the restriction map \[H^0(\mathscr X,\mathscr L^{\otimes mp})\longrightarrow H^0(\mathscr Y,\mathscr L^{\otimes mp}|_{\mathscr Y})\] is surjective. Hence we can choose $p\in\mathbb N_{\geqslant 1}$ such that 
\begin{equation}\label{eqn:Pro:liftingformodel:01}
|\beta|^{-1}\leqslant\mathrm{e}^{mpu\epsilon/2}
\end{equation}
and that there exists $t\in H^0(\mathscr X,\mathscr L^{\otimes mp})$ verifying $t|_{\mathscr Y}=\beta(\alpha\ell^m)^{pu}$. We then take $n=mpu$ and $s=\beta^{-1}\alpha^{-pu}t\in H^0(X,L^{\otimes n})$. One has $s|_{\mathscr Y}=\ell^{n}$ and
\[\begin{split}
\|s\|_{n\varphi} & =|\beta|^{-1}\cdot|\alpha|^{-pu}\cdot\|t\|_{n\varphi} \leqslant|\beta|^{-1}\cdot|\alpha|^{-pu}  \\
& \leqslant \mathrm{e}^{mpu\epsilon/2} \cdot \left( \mathrm{e}^{m\epsilon/2} \|\ell^m\|_{m\varphi_Y}\right)^{pu}   = \mathrm{e}^{n\varepsilon} \cdot \|\ell^m\|_{m\varphi_Y}^{pu}  \leqslant\mathrm{e}^{n\epsilon}\|\ell\|_{\varphi_Y}^n 
\end{split},\]
where the first inequality comes from Proposition \ref{Pro:latticeby model}(2), the second one from \eqref{eqn:Pro:liftingformodel:01} and \eqref{Equ:encadreentdialat}, and the last one from \eqref{Equ:sousmultiplicativesec}. 
The first part of the proposition is thus proved,
so that the last assertion follows from Proposition~\ref{prop:extension:prop:equiv:cond}.
\end{proof}

\begin{theo}\label{Thm:extensionpropertynontrivial}
Let $X$ be a projective $k$-scheme and $L$ be an ample invertible $\mathcal O_X$-module, equipped with a 
semipositive continuous metric $\varphi$. Let $Y$ be a closed subscheme of $X$ and $\ell$ be 
an element in $H^0(Y,L|_Y)$. Then $\ell$ has the extension property for $\varphi$.
\end{theo}
\begin{proof}
Clearly we may assume that $\ell \not\in \mathcal N(Y, \rest{L}{Y})$.  
First we assume that the absolute value $|\ndot|$ on $k$ is non-trivial. 
By Proposition \ref{Cor:convergencemodel}, for sufficiently positive integer $n$, there exists a sequence $(\mathscr X_n,\mathscr L_n)$  
of models of $(X,L^{\otimes n})$ such that each $\mathscr L_n$ is ample, and that 
\begin{equation}\label{Equ:convergencemodel1}\lim_{n\rightarrow+\infty}\frac 1n d(\varphi_{\mathscr L_n},n\varphi)=0.\end{equation}
For any $n\in\mathbb N_{\geqslant 1}$, let $\varphi^{(n)}=\frac1n\varphi_{\mathscr L_n}$.
By Proposition \ref{Pro:liftingformodel} (see also Remark \ref{Rem:extension}), one has $\lambda_{\varphi^{(n)}}(\ell)=0$. Therefore, by the relations \eqref{distancelambda} and \eqref{Equ:convergencemodel1}, we obtain that $\lambda_{\varphi}(\ell)=0$.

\bigskip
In the following, we treat the trivial valuation case. The main idea is to introduce the field of formal Laurent series over $k$ in order to reduce the problem to the non-trivial valuation case. We assume that the absolute value $|\ndot|$ is trivial. We denote by $k'$ the field $k(\!(T)\!)$ of formal Laurent series
over $k$, namely $k'$ is the fraction field of the ring $k\lbr T\rbr$ of formal series over $k$. Note that $k\lbr T\rbr$ is a discrete valuation ring.

\begin{enonce}{Claim}\label{Claim:Thm:extensionpropertynontrivial:01}
The field extension $k \subseteq k'$ is separable.
\end{enonce}

\begin{proof}
We may assume that the characteristic $p$ of $k$ is  positive.
First let us see the following claim:

\begin{enonce}{SubClaim}
Let $E$ be a finite extension of $k$ and
$\{\omega_i\}_{i=1}^e$ be a basis of $E$ over $k$. 
Then we have the following:
\begin{enumerate}[label=\rm(\roman*)]
\item\label{Item: if combination is zero in ET then gi are same}
Let $(g_1, \ldots, g_e) \in (k')^e$. If $\omega_1 g_1 + \cdots + \omega_e g_e = 0$ in $E(\!(T)\!)$, then $g_1 = \cdots = g_e = 0$.

\item\label{Item: if combination is zero in ET then ci are same}
Let $\{f_i\}_{i=1}^s$ be a family of elements in $k'$ which is linearly
independent over $k$ and $(c'_1, \ldots, c'_s)$ be an element of $E^s$. If
$c'_1 f_1 + \cdots + c'_s f_s = 0$ in $E(\!(T)\!)$, then $c'_1 = \cdots = c'_s = 0$.
\end{enumerate}
\end{enonce}

\begin{proof}
\ref{Item: if combination is zero in ET then gi are same} is trivial if we consider the coefficients of $g_1, \ldots, g_e$.

\medskip\ref{Item: if combination is zero in ET then ci are same}
We set $c'_i = \sum_{j=1}^e  c_{ij} \omega_j$ for some $c_{ij} \in k$. Then
\[
\sum_{i=1}^s c'_i f_i = \sum_{j=1}^e \left(\sum_{i=1}^s  c_{ij} f_i \right) \omega_j = 0,
\]
so that, by (i), $\sum_{i=1}^s  c_{ij} f_i = 0$ for all $j$.
Therefore $c_{ij} = 0$ for all $i, j$, as desired.
\end{proof}

By \cite[Th\'{e}or\`{e}me~2 in Chapter V, \S25, n${}^{\circ}$4]{Bourbaki_A4-7},
it is sufficient to see that
if $f_1, \ldots, f_s \in k(\!(T)\!)$ are linearly independent over $k$,
then $f_1^p, \ldots, f_s^p$ are  linearly independent over $k$.
We assume that $c_1 f_1^p + \cdots + c_s f_s^p = 0$ for some $c_1, \ldots, c_s \in k$.
Let $E$ be a finite extension field of $k$
such that we can find $c'_{i} \in E$ with $c_{i} = (c'_{i})^p$ for all $i$.
Then $\sum_{i=1}^s c'_{i}f_i = 0$ because
\[
0 = \sum_{i=1}^s c_{i}f_i^p  = \left( \sum_{i=1}^s c'_{i}f_i \right)^p.
\]
Thus, by (ii), one has $c'_{i} = 0$, as requested.
\end{proof}

Let us consider a subset $\Sigma$ of $\mathbb R$ given by
\[
\Sigma = \bigcup_{n=0}^{\infty} \ \bigcup_{(v, v') \in (H^0(X, L^{\otimes n}) \setminus \mathcal N(X, L^{\otimes n}))^2
} \mathbb Q(\ln \| v \|_{n\varphi} - \ln \| v' \|_{n\varphi}).
\]
Since 
\[
\{ \|v\|_{n\varphi} \,:\, v \in H^0(X, L^{\otimes n}) \} 
\]
is a finite set by Corollary~\ref{Cor:finitevlaue},
one has $\#(\Sigma) \leqslant \aleph_0$, so that one can choose $\alpha \in \mathbb R_{>0} \setminus \Sigma$.
We denote by $v_T(\ndot)$ the corresponding valuation on $k'$, and by $|\ndot|'$ the absolute value on $k'$ defined as 
\[\forall\,a\in k',\;|a|'=\mathrm{e}^{-\alpha v_T(a)}.\] 
Note that this absolute value extends the trivial absolute value on $k$. We denote by $X_{k'}$ and
$Y_{k'}$ the fibre products $X\times_{\Spec k}\Spec k'$ and $Y\times_{\Spec k}\Spec k'$, respectively,
and by $p:X_{k'}\rightarrow X$ and $p_Y : Y_{k'}\rightarrow Y$ the morphism of projections.

As explained in \S\ref{SubSec:basechange}, the morphism $p$ corresponds to a map $p^{\natural}:X_{k'}^{\mathrm{an}}\rightarrow X^{\mathrm{an}}$. This map is actually surjective. 
In fact, if $K$ is a field extension of $k$, equipped with an absolute value extending the trivial absolute value on $k$, then we can equip the field $K(T)$ of rational functions of one variable with the absolute value such that   
\[\forall\,F=a_0+a_1T+\cdots+a_nT^n\in K[T],\quad |F|=\max_{i\in\{0,\ldots,n\}}|a_i|\cdot \mathrm{e}^{-\alpha i}.\]
This absolue value extends the restriction of $|\ndot|'$ {to} $k(T)$. Hence the completion $\widehat{K(T)}$ of $K(T)$ with this absoute value is a valued extension of $k'$. If $f:\Spec K\rightarrow X$ is a $k$-morphism defining a point $x$ in $X^{\mathrm{an}}$, then it gives rise to a $k'$-morphism from $\Spec\widehat{K(T)}$ to $X_{k'}$, which defines a point $y$ in ${X'}^{\mathrm{an}}$ such that $p^{\natural}(y)=x$.

The surjectivity of $p^{\natural}$ implies that the restriction of the seminorm $\|\ndot\|_{n\varphi_{k'}}$ {to} $H^0(X,L)$ coincides with $\|\ndot\|_{n\varphi}$.
In fact, if $s$ is a section in $H^0(X,L^{\otimes n})$, then one has $\|p^*(s)\|_{n\varphi_{k'}}=\|s\|_{n\varphi}\circ p^{\natural}$. Therefore
\[\|p^*(s)\|_{n\varphi_{k'}}=\sup_{y\in {X'}^{\mathrm{an}}}\|p^*(s)\|_{n\varphi_{k'}}(y)=\sup_{y\in {X'}^{\mathrm{an}}}\|s\|_{n\varphi}(p^{\natural}(y))=\|s\|_{n\varphi},\]
where the last equality comes from the surjectivity of the map $p^{\natural}$. For $(v, v') \in (H^0(X,L^{\otimes n}) \setminus \mathcal N(X,L^{\otimes n}))^2$, if $\|v\|_{n\varphi}/\| v' \|_{n\varphi} \in |{k'}^{\times}|$,
then \[ \ln \|v\|_{n\varphi} - \ln \| v' \|_{n\varphi} = -\alpha\, v_T(a(T)) \] for some $a(T) \in {k'}^{\times}$. As
\[\alpha\not\in\bigcup_{(v,v')\in (H^0(X,L^{\otimes n}) \setminus \mathcal N(X,L^{\otimes n}))^2}\mathbb Q(\ln\|v\|_{n\varphi}-\ln\|v'\|_{n\varphi}),\]
we obtain $v_T(a(T)) = 0$, so that $\|v\|_{n\varphi} = \| v' \|_{n\varphi}$. By Proposition~\ref{prop:extension:trivial:Laurent}, the seminorm $\|\ndot\|_{n\varphi_{k'}}$ identifies with $\norm{\ndot}_{n\varphi,k',\varepsilon}$, the $\varepsilon$-extension of scalars of $\norm{\ndot}_{n\varphi}$.

Let $X_{\mathrm{red}}$ and $Y_{\mathrm{red}}$ be the reduced schemes associated with $X$ and $Y$, respectively. By Claim~\ref{Claim:Thm:extensionpropertynontrivial:01},
$X_{\mathrm{red}, k'} := X_{\mathrm{red}} \times_{\Spec k} \Spec k'$ and $Y_{\mathrm{red}, k'} := Y_{\mathrm{red}} \times_{\Spec k} \Spec k'$ are reduced (see \cite[Proposition IV.(4.6.1)]{EGAIV_2}), so that
\[
\mathcal N\big(X_{k'}, L_{k'}^{\otimes n}\big) = \mathcal N(X, L^{\otimes n}) \otimes_k k'
\quad\text{and}\quad
\mathcal N\big(Y_{k'}, \rest{L_{k'}}{Y_{k'}}^{\otimes n}\big) = \mathcal N(Y, \rest{L}{Y}) \otimes_k k',
\]
where $L_{k'} = L \otimes_k k'$.

By \eqref{distancelambda}, without 
loss of generality, we may assume that $L$ is very ample and that $\varphi$ is the quotient metric on $L$ induced by a ultrametric norm $\|\ndot\|$ on $V=H^0(X,L)/\mathcal N(X, L)$ and the natural surjection
$\beta : V \otimes \mathcal O_{X_{\mathrm{red}}} \to \rest{L}{X_{\mathrm{red}}}$. 
We may also assume that the restriction map $H^0(X,L^{\otimes n})\rightarrow H^0(Y,L|_Y^{\otimes n})$ is surjective for all $n \geqslant 1$.

For $n \geqslant 1$, let 
\[
V_n := H^0(X, L^{\otimes n})/\mathcal N(X, L^{\otimes n})\quad\text{and}\quad
V_{Y, n} = H^0(Y, \rest{L}{Y}^{\otimes n})/\mathcal N(Y, \rest{L}{Y}^{\otimes n}).
\]
Note that $V_1 = V$, and $V_n$ and $V_{Y, n}$ are isomorphic to the images
of 
\[
H^0(X, L^{\otimes n}) \to H^0(X_{\mathrm{red}}, \rest{L}{X_{\mathrm{red}}}^{\otimes n})
\quad\text{and}\quad
H^0(Y, \rest{L}{Y}^{\otimes n}) \to H^0(Y_{\mathrm{red}}, \rest{L}{Y_{\mathrm{red}}}^{\otimes n}),
\]
respectively. 
Let $V_{k'}$ be the vector space $V\otimes_kk'$ and let $\|\ndot\|_{k',\varepsilon}$ be the norm on $V_{k'}$ induced by $V$ by $\varepsilon$-extension of scalars. 
Then  
the surjective homomorphism $\beta:V \otimes \mathcal O_{X_{\mathrm{red}}} \rightarrow \rest{L}{X_{\mathrm{red}}}$ 
induces a surjective homomorphism 
$\beta_{k'}:V_{k'} \otimes \mathcal O_{X_{\mathrm{red}, k'}} \rightarrow \rest{L_{k'}}{X_{\mathrm{red}, k'}}
$. By Proposition \ref{Pro:quotientnormextensionscalar}, the metric $\varphi_{k'}$ of $L_{k'}$ obtained by $\varphi$ by extension of scalars coincides with the quotient metric on $L_{k'}$ induced by $(V_{k'},\|\ndot\|_{k',\varepsilon})$ and $\beta_{k'}$. Therefore, by Theorem \ref{Thm:extensionpropertynontrivial}, 
for any $\epsilon > 0$,
there exist an integer $n\geqslant 1$ 
and a section $s'\in H^0(X_{k'},L_{k'}^{\otimes n})$ such that $s'|_{Y_{k'}}=p_Y^*(\ell)^n$ and
\[\|s'\|_{n\varphi_{k'}}\leqslant \mathrm{e}^{n\epsilon}\|p_Y^*({\ell})\|_{\varphi_{Y,k'}}^n\] where  
$\varphi_{Y,k'}$ is the metric on $(L|_Y) \otimes_k k'\cong L_{k'}|_{Y_{k'}}$ induced by $\varphi_Y$ by extension of scalars, which equals the restriction of $\varphi_{k'}$ {to} $Y_{k'}$. 

Let $\{e_1,\ldots,e_{\alpha_n},f_1,\ldots,f_{\beta_n}\}$ be an orthogonal basis of $H^0(X,L^{\otimes n})$ with respect to $\norm{\ndot}_{n\varphi}$ such that $\{f_1,\ldots,f_{\beta_n}\}$ form a basis of the kernel of the restriction map $H^0(X,L^{\otimes n})\rightarrow H^0(Y,L|_Y^{\otimes n})$ (see Proposition \ref{Pro:existenceepsorth} for the existence of such an orthogonal basis).
We set \[s' = \sum_{i=1}^{\alpha_n} a_i(T) e_i+\sum_{j=1}^{\beta_n}b_j(T)f_j,\]
where 
$(a_1(T), \ldots, a_{\alpha_n}(T),b_1(T),\ldots,b_{\beta_n}(T)) \in (k')^{\alpha_n+\beta_n}$.
Note that \[\rest{s'}{Y_{k'}} = a_1(T) \rest{e_1}{Y} + \cdots + a_{\alpha_n}(T) \rest{e_{\alpha_n}}{Y}\] and
$\{ \rest{e_1}{Y}, \ldots, \rest{e_{\alpha_n}}{Y} \}$ forms a basis of $H^0(Y, \rest{L}{Y}^{\otimes n})$.
Since the restriction of $s'$ {to} $Y_{k'}$ can be written as the pull-back of an element in $H^0(Y,L|_Y^{\otimes n})$, 
one can see that $a_1(T), \ldots, a_{\alpha_n}(T) \in k$, so that $a_1(T), \ldots, a_{\alpha_n}(T)$ are denoted by $a_1, \ldots, a_{\alpha_n}$. 
Therefore if we set $s=a_1e_1+\cdots+a_{\alpha_n} e_{\alpha_n}$,
then $s\in H^0(X,L^{\otimes n})$, $s|_Y=\ell^{\otimes n}$, and
\begin{align*}
\|s\|_{n\varphi} & \leqslant \max_{i\in\{1,\ldots, \alpha_n\}}\{ |a_i|\cdot\|e_i\|_{n\varphi} \}\\
&\leqslant\max\Big\{\max_{i\in\{1,\ldots,\alpha_n\}}|a_i|\cdot\|e_i\|_{n\varphi},\max_{j\in\{1,\ldots,\beta_n\}}|b_{j}(T)|'\cdot\|f_{j}\|_{n\varphi}\Big\}\\
& = \|s'\|_{n\varphi_{k'}}\leqslant \mathrm{e}^{\epsilon n}\|p_Y^*(\ell)\|_{\varphi_{Y,k'}}^n = \mathrm{e}^{\epsilon n}\|\ell\|_{\varphi_Y}^n,
\end{align*}
where the first equality comes from the fact that $\{e_1,\ldots,e_{\alpha_n},f_1,\ldots,f_{\beta_n}\}$ forms an orthogonal basis of $H^0(X_{k'},L_{k'}^{\otimes n})$ with respect to $\norm{\ndot}_{n\varphi_{k'}}$ (see Proposition \ref{Pro:alpha-orthgonalextension}). The theorem is thus proved.
\end{proof}

\begin{rema}\label{remark:reduced:cohomology:group}
Let $X$ be a scheme of finite type over $k$ and $\mathcal F$ be a coherent $\mathcal O_X$-module.
Let $X_{\mathrm{red}}$ be the reduced structure of $X$. The \emph{reduced $i$-th cohomology group of $\mathcal F$}\index{reduced cohomology group},
denoted by $H^i_{\mathrm{red}}(X, \mathcal F)$, is
defined to be the image of the homomorphism \[
H^i(X, \mathcal F) \to H^i\big(X_{\mathrm{red}}, \rest{F}{X_{\mathrm{red}}}\big).
\]
Using the reduced cohomology group, one has the following variant of Theorem~\ref{Thm:extensionpropertynontrivial}:

\begin{quote}
\hskip\parindent\emph{We assume that $X$ is projective.
Let $L$ be an ample invertible $\mathcal O_X$-module, equipped with
a semipositive continuous metric $\varphi$. Let $Y$ be a closed subscheme of $X$ and
$\ell$ be an element of $H^0_{\mathrm{red}}(Y, \rest{L}{Y})$. Then, for any $\epsilon > 0$,
there exist a positive integer $n$ and $s \in H^0_{\mathrm{red}}(X, L^{\otimes n})$ such that
$\rest{s}{Y} = l^{\otimes n}$ and $\|s\|'_{n\varphi} \leqslant \mathrm{e}^{n\epsilon} (\|l\|'_{\varphi_Y})^n$,
where $\rest{s}{Y}$ is the image of $s$ by the homomorphism $H^0_{\mathrm{red}}(X, L^{\otimes n}) \to
H^0_{\mathrm{red}}(Y, \rest{L}{Y}^{\otimes n})$.}
\end{quote}
Since $H^0_{\mathrm{red}}(X, L^{\otimes n})$ and $H^0_{\mathrm{red}}(Y, \rest{L}{Y}^{\otimes n})$ are
subgroups of $H^0(X_{\mathrm{red}}, \rest{L}{X_{\mathrm{red}}}^{\otimes n})$ and $H^0(Y_{\mathrm{red}}, \rest{L}{Y_{\mathrm{red}}}^{\otimes n})$, the proof of the above result can be done in the similar way as Theorem~\ref{Thm:extensionpropertynontrivial}.
\end{rema}

\if00
\section{Cartier divisors}
In this section, we recall Cartier divisors and linear systems.
Further we introduce $\mathbb Q$-Cartier and $\mathbb R$-Cartier divisors on
an integral scheme and study their basic properties.

\subsection{Reminder on Cartier divisors} 
Let $X$ be a scheme. For any open subset $U$ of $X$, we denote by $S_X(U)$ the set of all elements $a\in\mathcal O_X(U)$ such that the homomorphisme of $\mathcal O_U$-modules $\mathcal O_U\rightarrow\mathcal O_U$ defined as the homothety by $a$ is injective. This is a multiplicative subset of $\mathcal O_X(U)$. Moreover, $S_X$ is a subsheaf of sets of $\mathcal O_X$. We denote by  $\mathscr M_X$ the sheaf of rings associated with the presheaf
\[U\longmapsto \mathcal O_X(U)[S_X(U)^{-1}],\]
called the sheaf of \emph{meromorphic functions}\footnote{The definition in \cite[IV.20]{EGAIV_4} is not adequate, see \cite{Kleiman79} for details.} on $X$. The canonical homomorphisms $\mathcal O_X(U)\rightarrow\mathcal O_X(U)[S_X(U)^{-1}]$ of rings induce a homomorphism of sheaves of rings $\mathcal O_X\rightarrow\mathscr M_X$.  The sections of $\mathscr M_X$ on an open subset $U$ of $X$ are  
called \emph{meromorphic fonctions}\index{meromorphic functions} on $U$.

\begin{rema}
Let $U$ be an open subset of $X$. Any element $a$ in $S_X(U)$ is a regular element (namely the homothety $\mathcal O_X(U)\rightarrow\mathcal O_X(U)$ defined by $a$ is injective). of $\mathcal O_X(U)$. It is not true in general that $S_X(U)$ contains all regular elements of $\mathcal O_X(U)$. However, if $U$ is an affine open subset of $X$, $S_X(U)$ identifies with the set of all regular elements in  $\mathcal O_X(U)$. In fact, an element $a\in\mathcal O_X(U)$ is in $S_X(U)$ if and only if its image in the local ring $\mathcal O_{X,x}$ is a regular element for any $x\in U$. The announced property thus results from the faithful flatness of
\[\mathcal O_X(U)\longrightarrow\bigoplus_{x\in U}\mathcal O_{X,x},\]
provided that $U$ is an affine open subset. In particular, if $X$ is an integral scheme, then $\mathscr M_X$ is the constant sheaf associated to the local ring of the generic point of $X$ (which is a field). We use the notation $R(X)$ to denote the field of all meromorphic functions on $X$. 
\end{rema}

Denote by $\mathcal O_X^{\times}$ the sheaf of abelian groupes described as follows. For any open subset $U$ of $X$, $\mathcal O_X^{\times}(U)$ is the set of elements  $a\in\mathcal O_X(U)$ such that the homothety $\mathcal O_U\rightarrow\mathcal O_U$ defined by $a$ is an isomorphism. Similarly, denote by $\mathscr M_X^{\times}$ the sheaf of abelian groupes on $X$ whose section space over any open subset $U\subseteq X$ is the set of all meromorphic functions $\varphi$ on $U$ such that the homothety $\mathscr M_U\rightarrow\mathscr M_U$ defined by $\varphi$ is an isomorphism. This is a subsheaf of multiplicative monoids of $\mathscr M_X$.

\begin{defi}\label{Def:Cartierdivisor}
We call \emph{Cartier divisor}\index{Cartier divisor} on $X$ any global section of the quotient sheaf $\mathscr M_X^\times/\mathcal O_X^{\times}$, that is,
a data of a Zariski open covering $X = \bigcup_{\alpha} U_{\alpha}$ of $X$ and a section $s_{\alpha} \in \mathscr M_X^{\times}(U_{\alpha})$ over $U_{\alpha}$ for each $\alpha$, which is called a \emph{local equation}\index{local equation} over $U_{\alpha}$, 
such that $s_{\alpha}s^{-1}_{\beta} \in \mathcal O_X^{\times}(U_{\alpha} \cap U_{\beta})$ for all $\alpha$ and $\beta$. The group of all Cartier divisors on $X$ is denoted by $\Div(X)$, where the group law is written additively. We say that a Cartier divisor $D$ is \emph{effective}\index{effective divisor}\index{Cartier divisor!effective ---} if it is a section of $(\mathscr M_X^\times\cap\mathcal O_X)/\mathcal O_X^{\times}$. We use the expression $D\geqslant 0$ to denote the condition ``$D$ is effective''. 
Moreover, for $D_1, D_2 \in \Div(X)$, an expression $D_1 \geqslant D_2$ is defined by $D_1 - D_2 \geqslant 0$.
\end{defi}

The Cartier divisors are closely related to invertible sheaves.
Let $D$ be a Cartier divisor on $X$. Denote by $\mathcal O_X(D)$ the sub-$\mathcal O_X$-module of $\mathscr M_X$ generated by $-D$. Namely, if $X = \bigcup_{\alpha} U_{\alpha}$ is an open covering of $X$ and $s_{\alpha}$ is a local equation of $D$ over $U_{\alpha}$, then $\rest{\mathcal O_X(D)}{U_{\alpha}} =
\mathcal O_{U_{\alpha}} s_{\alpha}^{-1}$. Note that $\mathcal O_X(D)$ is an invertible $\mathcal O_X$-module since it is locally generated by a regular element. We say that the Cartier divisor $D$ is \emph{ample}\index{ample Cartier divisor}\index{Cartier divisor!ample ---} (resp. \emph{very ample}\index{very ample Cartier divisor}\index{Cartier divisor!very ample ---}) if the invertible sheaf $\mathcal O_X(D)$ is ample (resp. very ample).
 By definition, one has $\mathcal O_X(-D)\cong\mathcal O_X(D)^\vee$. Moreover, if $D_1$ and $D_2$ are two Cartier divisors, then $\mathcal O_X(D_1+D_2)\cong\mathcal O_X(D_1)\otimes\mathcal O_X(D_2)$. Thus the map sending a divisor $D$ to the isomorphism class of $\mathcal O_X(D)$ defines a homomorphism from $\Div(X)$ to the Picard group $\mathrm{Pic}(X)$ (namely the group of isomorphism classes of invertible $\mathcal O_X$-modules). This homomorphism is surjective notably when $X$ is a reduced scheme with locally finite irreducible components, or a quasi-projective scheme over a {Noetherian} ring (cf. \cite[IV.21.3.4-5]{EGAIV_4}). We recall a simple proof of this result for the particular case where $X$ is an integral scheme.

\begin{prop}\label{Pro:surjdivi}
Let $X$ be an integral scheme. Then the homomorphism $\Div(X)\rightarrow\Pic(X)$ constructed above is surjective.
\end{prop}
\begin{proof}
Let $\eta$ be the generic point of $X$ and $R(X)$ be the field of all meromorphic functions on $X$. Let $L$ be an invertible sheaf and $s$ be a non-zero element in $L_\eta$. Then the maps $ H^0(U,L)\rightarrow R(X)$, $t\mapsto t_\eta/s$ (where $U$ denotes an open subset of $X$) define a $\mathcal O_X$-linear homomorphisme from $L$ to $\mathscr M_X$. The images of local trivialisations of $L$ by this homomorphism {define}  a global section of $\mathscr M_X^{\times}/\mathcal O_X^{\times}$, whose opposite $D$ is a Cartier divisor such that $\mathcal O_X(D)\cong L$. 
\end{proof}

\begin{rema}
Let $X$ be an integral scheme and $L$ be an invertible $\mathcal O_X$-module. Let $\eta$ be the generic point of $X$. We call \emph{rational section}\index{rational section} of $L$ any element in $L_\eta$. Note that for any non-empty open subset $U$ of $X$, the restriction map $ H^0(U,L)\rightarrow L_\eta$ is injective. By abuse of language, we also call a section of $L$ on a non-empty open subset of $X$ a rational section of $L$. The proof of the above proposition shows that any non-zero rational section $s$ of $L$ defines a Cartier divisor of $X$, which we denote by $\mathrm{div}(s)$. One can verify that, if $L$ and $L'$ are two invertible $\mathcal O_X$-modules and if $s$ and $s'$ are respectively non-zero rational sections of $L$ and $L'$, then one has $\mathrm{div}(ss')=\mathrm{div}(s)+\mathrm{div}(s')$.
\end{rema}

The exact sequence of abelian sheaves 
\[\xymatrix{1\ar[r]&\mathcal O_X^{\times}\ar[r]&\mathscr M_X^{\times}\ar[r]&\mathscr M_X^\times/\mathcal O_X^{\times}\ar[r]&0}\]
induces a long exact sequence of cohomology groups 
\begin{equation}\label{Equ:grouppicard}\xymatrix{\relax 1\ar[r]& H^0(X,\mathcal O_X^{\times})\ar[r]&H^0(X,\mathscr M_X^{\times})\ar[r]^-{\mathrm{div}}&\Div(X)\ar[r]^-\theta& H^1(X,\mathcal O_X^{\times})}.\end{equation}
Note that the cohomology group $H^1(X,\mathcal O_X^{\times})$ identifies with the Picard group $\Pic(X)$ of $X$ (cf. \cite[0.5.6.3]{NouveauEGA1}), and $\theta$ is just the group homomorphism sending any Cartier divisor $D$ to the isomorphism class of the invertible $\mathcal O_X$-module $\mathcal O_X(D)$. The image of the group homomorphism $\mathrm{div}(\ndot)$ is denoted by $\mathrm{PDiv}(X)$. The divisors in $\mathrm{PDiv}(X)$ are called \emph{principal divisors}\index{principal divisor}\index{Cartier divisor!principal ---}. The quotient group $\Div(X)/\mathrm{PDiv}(X)$ is called the \emph{divisor class group} of $X$, denoted by $\mathrm{Cl}(X)$. The exactness of the sequence \eqref{Equ:grouppicard} shows that the homomorphism from $\mathrm{Cl}(X)$ to $\Pic(X)$, sending the equivalent class of a {Cartier} divisor $D$ to the isomorphism class of the invertible sheaf $\mathcal O_X(D)$ is injective. It is an isomorphism once $X$ is a reduced scheme with locally finite irreducible components, or a quasi-projective scheme over a Noetherian ring. We write this result as a corollary of Proposition \ref{Pro:surjdivi} in the particular case where $X$ is an integral scheme.

\begin{coro}
Let $X$ be an integral scheme. The homomorphism  from $\mathrm{Cl}(X)$ to $\Pic(X)$ sending the equivalence class of a divisor class $D$ to the isomorphism class of $\mathcal O_X(D)$ is an isomorphism. 
\end{coro}

If two Cartier divisors $D$ and $D'$ of $X$ differ by a principal divisor, namely lie in the same class in $\mathrm{Cl}(X)$, we say that they are \emph{linearly equivalent}\index{linearly equivalent}\index{Cartier divisor!linearly equivalent}.

\begin{prop}\label{prop:Effective:Weil:Effective}
We assume that $X$ is locally Noetherian and normal. Let $D$ be a Cartier divisor on $X$ and $D = \sum_{\Gamma \in X^{(1)}} a_{\Gamma} \Gamma$
be the expansion as a Weil divisor, where $X^{(1)}$ is the set of all codimension one points of $X$. Then
$D \geqslant 0$ if and only if $a_{\Gamma} \geqslant 0$ for all $\Gamma \in X^{(1)}$.
\end{prop}

\begin{proof}
It is suuffuent to show that if $a_{\Gamma} \geqslant 0$ for all $\Gamma \in X^{(1)}$, then $D \geqslant 0$.
Let $f_x$ be a local equation of $D$ at $x \in X$. By our assumption, $f_x \in \mathcal O_{X,\Gamma}$ for all $\Gamma \in X^{(1)}$ and $x \in \overline{\{\Gamma\}}$, so that, by virtue of \cite[THEOREM~38]{Matsumura1979}, 
\[
f_x \in \bigcap\limits_{x \in \overline{\{\Gamma\}},\, \Gamma \in X^{(1)}} \mathcal O_{X,\Gamma} = \mathcal O_{X,x},
\]
and hence the assertion follows.
\end{proof}

\subsection{Linear system of a divisor}

In this subsection, we fix an integral scheme $X$, and denote by $R(X)$ the field of all rational functions on $X$. 
\begin{defi}\label{linear:system:Cartier:div}
Let $D$ be a Cartier divisor of $X$. We define 
\[H^0(D):=\{f\in R(X)^{\times}\,:\,\mathrm{div}(f)+D\geqslant 0\}\cup\{0\},\]
called the \emph{complete linear system}\index{complete linear system} of the divisor $D$. It forms
a subgroup of $R(X)$ with respect to the additive composition law and is invariant by the multiplication by a scalar in $K$. Hence it is a $K$-vector subspace of $R(X)$.
\end{defi}
We obtain from the definition that, if $D$ and $D'$ are two Cartier divisors which are linearly equivalent, and $g$ is a non-zero rational function such that $D'=D+\mathrm{div}(g)$. Then the map $f\mapsto fg$ defines a bijection from $H^0(D)$ to $H^0(D')$.

Let $D$ be a Cartier divisor of $X$. Being a sub-$\mathcal O_X$-module of $\mathscr M_X$, the invertible sheaf $\mathcal O_X(D)$ shares the same generic fibre with $\mathscr M_X$, which is also canonically isomorphic to the field $R(X)$. Therefore the unit element in $R(X)$ defines a rational section of $D$ which we denote by $s_D$. One can verify that $\mathrm{div}(s_D)=D$ and $s_{D+D'}=s_Ds_{D'}$ for any couple $(D,D')$ of Cartier divisors of $X$.

\begin{lemm}\label{Lem:effectivdivisor}
Let $D$ be a Cartier divisor of $X$ and $s_D$ be the meromorphic section of $\mathcal O_X(D)$ constructed above. Then $D$ is an effective divisor if and only if $s_D$ extends to a global section of $\mathcal O_X(D)$.
\end{lemm}
\begin{proof}
Assume that $D$ is an effective divisor. Then the invertible $\mathcal O_X$-module $\mathcal O_X(-D)$ is actually an invertible ideal sheaf of $\mathcal O_X$ since it is generated by $D$. Let $s:\mathcal O_X\rightarrow\mathcal O_X(D)$ be the homomorphism of $\mathcal O_X$-modules which is dual to the inclusion map $\mathcal O_X(-D)\rightarrow\mathcal O_X$. It defines a global section of $\mathcal O_X(D)$ whose value at the generic point coincides with $s_D$.

Conversely, if $L$ is an invertible sheaf on $X$ and if $s$ is a non-zero global section of $L$, then $\mathrm{div}(s)$ is an effective Cartier divisor. In particular, if $s_D$ extends to a global section of $\mathcal O_X(D)$, then $D=\mathrm{div}(s_D)$ is an effective divisor.  
\end{proof}

\begin{prop}\label{Pro:H0Dcriterion}
Let $D$ be a Cartier divisor of $X$. A rational function $f$ lies in $H^0(D)$ if and only if $fs_D$ extends to a global section of $\mathcal O_X(D)$. 
\end{prop}
\begin{proof}
By definition, for any non-zero meromorphic function $f\in K$, the relation $fs_D=s_{\mathrm{div}(f)+D}$ holds.
The Lemma \ref{Lem:effectivdivisor} shows that $fs_D$ extends to a global section of $\mathcal O_X\otimes\mathcal O_X(D)\cong\mathcal O_X(D)$ if and only if $\mathrm{div}(f)+D$ is an effective divisor. The proposition is thus proved.
\end{proof}

\begin{rema}
\begin{enumerate}[label=\rm(\arabic*)]
\item\label{Item: cartier divisor and invertible sheaf} Let $L$ be an invertible $\mathcal O_X$-module.  The Proposition \ref{Pro:H0Dcriterion} shows that, if $s$ is a non-zero meromorphic section of $L$, then the relation $t\mapsto t/s$ defines an isomorphism between the groups $H^0(X,L)$ and $H^0(\mathrm{div}(s))$. In particular, if $D$ is a Cartier divisor, then $H^0(D)$ is canonically isomorphic to $H^0(X,\mathcal O_X(D))$.
\item\label{Item: global section of Cartier divisor} Assume that the scheme $X$ is defined over a ground field $k$, then the field of rational functions $R(X)$ is an extension of $k$. Moreover, for any Cartier divisor $D$ of $X$, $H^0(D)$ is a $k$-vector subspace of $R(X)$.  
\end{enumerate}
\end{rema}

\subsection{$\mathbb Q$-Cartier and $\mathbb R$-Cartier divisors}

{As in the previous subsection, $X$ denotes an integral scheme.}
Let $\mathbb K$ be either $\mathbb Z$, $\mathbb Q$ or $\mathbb R$.
An element of $\Div_{\mathbb K}(X) := \Div(X) \otimes_{\mathbb Z} \mathbb K$ is called a
\emph{$\mathbb K$-Cartier divisor}\index{Cartier divisor!K-Cartier divisor@$K$-Cartier divisor} on $X$. Note that a $\mathbb Z$-Cartier divisor is 
a usual Cartier divisor.
A $\mathbb K$-Cartier divisor can be regarded as
an element of 
\[
H^0(X, (\mathscr M_X^{\times}/\mathcal O_X^{\times}) \otimes_{\mathbb Z} \mathbb K) =
H^0(X, (\mathscr M_X^{\times} \otimes_{\mathbb Z} \mathbb K)/(\mathcal O_X^{\times} \otimes_{\mathbb Z} \mathbb K) ),
\]
so that, for any point $x \in X$, there are an open neighborhood $U$ of $x$ and
$f \in \rest{(\mathscr M_X^{\times} \otimes_{\mathbb Z} \mathbb K)}{U}$ such that
$D$ is defined by $f$ over $U$. Note that if $f' \in \rest{(\mathscr M_X^{\times} \otimes_{\mathbb Z} \mathbb K)}{U}$ also defines $D$ over $U$, then $f/f' \in \rest{(\mathcal O_X^{\times} \otimes_{\mathbb Z} \mathbb K)}{U}$. The element $f$ is called a \emph{local equation}\index{local equation} of $D$. Moreover, the morphism of groups $K(X)^{\times}\rightarrow\mathrm{Div}(X)$ induces by extension of scalars a $\mathbb K$-linear map $K(X)^{\times}\otimes_{\mathbb Z}\mathbb K\rightarrow\mathrm{Div}_{\mathbb K}(X)$ which we denote by $\mathrm{div}_{\mathbb K}(\ndot)$.

Let $D$ be a $\mathbb K$-Cartier divisor on $X$.
Let $f_x$ be a local equation of $D$ around $x$.
Note that the condition $f_x \in \mathcal O_X^{\times} \otimes_{\mathbb Z} \mathbb K$
does not depend on the choice of the local equation of $D$ around $x$, so that
we define $\operatorname{Supp}_{\mathbb K}(D)$ to be
\[
\operatorname{Supp}_{\mathbb K}(D) = \{ x \in X \,:\, f_x \not\in \mathcal O_X^{\times} \otimes_{\mathbb Z} \mathbb K \}.
\]

\begin{prop}
\begin{enumerate}[label=\rm(\arabic*)]
\item
$\operatorname{Supp}_{\mathbb K}(D)$ is a closed subset of $X$.

\item
If $D$ is a Cartier divisor, then $\operatorname{Supp}_{\mathbb Q}(D) =
\bigcap_{n=1}^{\infty} \operatorname{Supp}_{\mathbb Z}(nD)$.
In particular, $\operatorname{Supp}_{\mathbb Q}(D) \subseteq
\operatorname{Supp}_{\mathbb Z}(D)$. Moreover, if $X$ is normal, then
$\operatorname{Supp}_{\mathbb Q}(D) =
\operatorname{Supp}_{\mathbb Z}(D)$.

\item
If $D$ is a $\mathbb Q$-Cartier divisor, then
$\operatorname{Supp}_{\mathbb Q}(D) =  \operatorname{Supp}_{\mathbb R}(D)$.
\end{enumerate}
\end{prop}

\begin{proof}
The proof can be found in \cite[Section~1.2]{Moriwaki16}.
\end{proof}

\begin{defi}
Let $D$ be a $\mathbb K$-Cartier divisor on $X$.
We say that $D$ is $\mathbb K$-effective, denoted by $D \geqslant_{\mathbb K} 0$, if, for every $x \in X$,
a local equation of $D$ can be expressed by $f_1^{a_1} \cdots f_r^{a_r}$, where
$f_1, \ldots, f_r \in \mathcal O_{X, x} \setminus \{ 0 \}$ and
$a_1, \ldots, a_r \in \mathbb R_{>0}$.
Similarly as Definition~\ref{linear:system:Cartier:div}, we define $H^0_{\mathbb K}(D)$  to be
\[
H^0_{\mathbb K}(X, D) := \{ \varphi \in R(X)^{\times}\,:\,\mathrm{div}(\varphi)+D\geqslant_{\mathbb K} 0\} \cup \{ 0 \}.
\]
Note that in the case where $\mathbb K = \mathbb Z$,
$H^0_{\mathbb Z}(X, D)$ coincides with $H^0(D)$ in Definition~\ref{linear:system:Cartier:div}.
\end{defi}

\begin{prop}\label{prop:basic:prop:K:Cartier:div}
Let $D$ be a $\mathbb K$-Cartier divisor on $X$. Then we have the following:
\begin{enumerate}[label=\rm(\arabic*)]
\item\label{Item: effective Q Cartier}
We assume that $\mathbb K = \mathbb Q$.
Then $D \geqslant_{\mathbb Q} 0$ if and only if $D \geqslant_{\mathbb R} 0$.

\item\label{Item: bijective of H0Q and H0R}
We assume that $\mathbb K = \mathbb Q$.
Then the natural map $H^0_{\mathbb Q}(X, D) \to H^0_{\mathbb R}(X, D)$ is bijective.

\item\label{Item: effective Cartier Normal}
We assume that $\mathbb K = \mathbb Z$ and $X$ is locally Noetherian and normal.
Then $D \geqslant_{\mathbb Z} 0$ if and only if $D \geqslant_{\mathbb Q} 0$.

\item\label{Item: bijective of H0Z and H0Q Normal}
We assume that $\mathbb K = \mathbb Z$ and $X$ is locally Noetherian and normal.
Then the natural map $H^0_{\mathbb Z}(X, D) \to H^0_{\mathbb Q}(X, D)$ is bijective.

\item\label{Item: product of global sections} If $a \in H^0(X, \mathcal O_X)$ and $\varphi \in H^0_{\mathbb K}(X, D)$,
then $a \varphi \in H^0_{\mathbb K}(X, D)$.

\item\label{Item: normal case global section submodule}
If $X$ is locally Noetherian and normal, then $H^0_{\mathbb K}(X, D)$ forms a $H^0(X, \mathcal O_X)$-submodule of $R(X)$.
\end{enumerate}
\end{prop}

\begin{proof}
\ref{Item: effective Q Cartier} Obviously $D \geqslant_{\mathbb Q} 0$ implies $D \geqslant_{\mathbb R} 0$.
Conversely we assume that $D \geqslant_{\mathbb R} 0$.
Let $f_x$ be a local equation of $D$ at $x \in X$.
Then there are $f_1, \ldots, f_r \in \mathcal O_{X, x} \setminus \{ 0 \}$ and $a_1, \ldots, a_r \in \mathbb R_{>0}$ such that
$f_x = f_1^{a_1} \cdots f_{r}^{a_r}$. Note that $f_x \in R(X)^{\times} \otimes_{\mathbb Z} \mathbb Q$, so that,
by Lemma~\ref{lemm:approximation:Real:by:Rational} as below, there are $a'_1, \ldots, a'_r \in \mathbb Q_{>0}$ such that
$f_x = f_1^{a'_1} \cdots f_{r}^{a'_r}$.
Therefore, $D$ is $\mathbb Q$-effective.

\medskip
\ref{Item: bijective of H0Q and H0R} Let $\varphi \in H^0_{\mathbb R}(X, D) \setminus \{ 0 \}$. 
Then $D + \mathrm{div}(\varphi) \geqslant_{\mathbb R} 0$.
Note that $D + \mathrm{div}(\varphi)$ is a $\mathbb Q$-Cartier divisor, so that, by (1),
$D + \mathrm{div}(\varphi) \geqslant_{\mathbb Q} 0$, which means $\varphi \in H^0_{\mathbb Q}(X, D) \setminus \{ 0 \}$.

\medskip
\ref{Item: effective Cartier Normal}
We assume that $D \geqslant_{\mathbb Q} 0$. Let $D = \sum_{\Gamma} a_{\Gamma} \Gamma$ be the expansion as a Weil divisor.
Then $a_{\Gamma} \geqslant 0$ for all $\Gamma$. Thus $D \geqslant 0$ by Proposition~\ref{prop:Effective:Weil:Effective}.

\medskip
\ref{Item: bijective of H0Z and H0Q Normal} is a consequence of \ref{Item: effective Cartier Normal}.
 
\ref{Item: product of global sections} is obvious.

\medskip
\ref{Item: normal case global section submodule} By \ref{Item: product of global sections}, it is sufficient to show that if
$\varphi, \psi \in H^0_{\mathbb K}(X, D)$, then $\varphi + \psi \in H^0_{\mathbb K}(X, D)$.
If we set $D = \sum_{\Gamma} \alpha_{\Gamma} \Gamma$ ($\alpha_{\Gamma} \in \mathbb K$) and
$\mathrm{div}(\varphi) = \sum_{\Gamma} \ord_{\Gamma} (\varphi) \Gamma$ as a Weil divisor
for $\varphi \in R(X)^{\times}$, then 
\[
\mathrm{div}(\varphi)+D\geqslant_{\mathbb K} 0
\quad\Longleftrightarrow\quad \forall\ \Gamma, \ 
\ord_{\Gamma}(\varphi) + \alpha_{\Gamma} \geqslant 0
\]
by \cite[Lemma~1.2.4]{Moriwaki16} together with \ref{Item: effective Q Cartier}.
Moreover, for $\varphi, \psi \in R(X)$, 
\[
\ord_{\Gamma}(\varphi + \psi) \geqslant \min \{ \ord_{\Gamma}(\varphi), \ord_{\Gamma}(\psi) \}.
\]
Therefore \ref{Item: normal case global section submodule} follows.
\end{proof}

\begin{lemm}
\label{lemm:approximation:Real:by:Rational}
Let $V$ be a vector space over $\mathbb Q$. Then we have the following:
\begin{enumerate}[label=\rm(\arabic*)]
\item\label{Item: intersection WR and V}
$W_{\mathbb R} \cap V = W$ for any vector subspace $W$ of $V$.

\item\label{Item: approximation of R-linear combination}
Let $x, x_1, \ldots, x_r \in V$ such that
$x = a_1 x_1 + \cdots + a_r x_r$
for some $a_1, \ldots, a_r \in \mathbb R$. Then, for any $\epsilon > 0$,
there are $a'_1, \ldots, a'_r \in \mathbb Q$ such that
$x = a'_1 x_1 + \cdots + a'_r x_r$ and $| a'_i - a_i | \leqslant \epsilon$ for all $i$.
\end{enumerate}
\end{lemm}

\begin{proof}
\ref{Item: intersection WR and V}
is obvious because $V/W \to (V/W)_{\mathbb R}$ is injective and $(V/W)_{\mathbb R} = V_{\mathbb R}/W_{\mathbb R}$.

\medskip
\ref{Item: approximation of R-linear combination}
We consider the homomorphism $\psi:\mathbb Q^r\to V$ sending $(t_1,\ldots,t_r)\in\mathbb Q^r$ to $t_1x_1+\cdots+t_rx_r$. Denote by $W$ the image of $\psi$. By (1), the point $x$ belongs to $W$. We pick an element $b$ in $\psi^{-1}(\{x\})$. Let $\psi_{\mathbb R}:\mathbb R^r \to V_{\mathbb R}$ be the scalar extension of $\psi$, that is,
$\psi_{\mathbb R}(\alpha_1, \ldots, \alpha_r) = \alpha_1 x_1 + \cdots + \alpha_r x_r$ for any $(\alpha_1,\ldots,\alpha_r)\in\mathbb R^r$, whose image is $W_{\mathbb R}$. As $\Ker(\psi_{\mathbb R}) = \Ker(\psi)_{\mathbb R}$, $\Ker(\psi)$ is dense in $\Ker(\psi_{\mathbb R})$. Therefore, $\psi^{-1}(\{x\})=b+\Ker(\psi)$ is dense in $\psi_{\mathbb R}^{-1}(\{x\})=b+\Ker(\psi_{\mathbb R})$, which implies the assertion of \ref{Item: approximation of R-linear combination}.
\end{proof}

\begin{exem}
\label{example:R:Cartier:vs:Cartier} The study of effective $\mathbb Q$-Cartier or $\mathbb R$-Cartier divisors on non-normal schemes is more subtle than that in the normal case. This phenomenon can be shown by the following examples, which have been discussed in \cite{Chen_Moriwaki2017b}.
Let $X := \operatorname{Proj}(k[T_0, T_1, T_2]/(T_0T_2^2 - T_1^3)$ over a field $k$, 
$U_i := \{ T_i \not= 0 \}\cap X$ 
($i=0,1,2$) and $x := T_1/T_0, y := T_2/T_0$ on $U_0$.
Then $U_0 = X \setminus \{ (0:0:1) \}$ and $U_2 = X \setminus \{ (1:0:0) \}$, so that $X = U_0 \cup U_2$.
Note that $y/x$ is not regular at $(1:0:0)$ and
$y/x \in \mathcal O_{X, \zeta}^{\times}$ for all $\zeta \in U_0 \cap U_2$.
Let $D$ be a Cartier divisor on $X$ given by
\[
D = \begin{cases}
(y/x) & \text{on $U_0$}, \\
(1) & \text{on $U_2$}.
\end{cases}
\]

\begin{enumerate}[label=\rm(\arabic*)]
\item As $y/x$ is not regular at $(1:0:0)$, $D$ is not effective as a Cartier divisor.
On the other hand, since
\[
2D = \begin{cases}
(x) & \text{on $U_0$}, \\
(1) & \text{on $U_2$},
\end{cases}
\]
$D$ is effective as a $\mathbb Q$-Cartier divisor.
As a consequence, $1 \not\in H^0(X, D)$ and $1 \in H^0_{\mathbb Q}(X, D)$, that is,
$H^0(X, D) \to H^0_{\mathbb Q}(X, D)$ is not surjective.

\item
We assume that $\mathrm{char}(k) = 0$.
We set $\varphi := x/y$. 
As $\varphi = T_1/T_2$ is regular on $U_2$, $\varphi \in H^0(X, D)$.
Here let us see $1 + \varphi \not\in H^0_{\mathbb Q}(X, D)$.
We assume the contrary, that is, $1 + \varphi \in H^0_{\mathbb Q}(X, D)$. Then
\[
(1 + \varphi)(y/x) = 1 + y/x
\]
is $\mathbb Q$-effective on $U_0$, so that there is a positive integer $N$ such that $(1 + y/x)^N$ is regular on $U_0$.
Here we claim that $(y/x)^i$ is regular over $U_0$ for an integer $i \geqslant 2$.
Indeed, we set $i = 2j + \epsilon$, where $j \geqslant 1$ and $\epsilon \in \{ 0, 1 \}$.
Then as
\[
(y/x)^i = (y/x)^{2j + \epsilon} = (y^2)^jy^{\epsilon}x^{-2j - \epsilon} = x^{j-\epsilon} y^{\epsilon},
\]
the assertion follows. 
Note that
\[
y/x = (1/N) \left( (1 + y/x)^N - 1 - \sum_{i=2}^N \binom{N}{i} (y/x)^i \right),
\]
so that $y/x$ is regular on $U_0$. This is a contradiction because $y/x$ is not regular on $U_0$.

\item
Next we assume that $\mathrm{char}(k) = 2$. We set $U'_0 := U_0 \setminus \{ (1 : 1 : 1) \}$.
Note that $X = U'_0 \cup U_2$ and
$1 + y/x \in \mathcal O^{\times}_{X, \zeta}$ for all $\zeta \in U'_0 \cap U_2$, so that
we set
\[
D' := \begin{cases}
(1 + y/x) & \text{on $U'_0$}, \\
(1) & \text{on $U_2$}.
\end{cases}
\]
Since $y/x$ is not regular at $(1:0:0)$, we have $D' \not= 0$.
Moreover, as $(1 + y/x)^2 = 1 + x$, we have
\[
2D' = \begin{cases}
(1+x) & \text{on $U'_0$}, \\
(1) & \text{on $U_2$},
\end{cases}
\]
and hence $2D' = 0$ because $1 + x \in \mathcal O_{X, \zeta}^{\times}$ for all $\zeta \in U'_0$.
Therefore, the natural homomorphism $\Div(X) \to \Div_{\mathbb K }(X)$ is not injective.
Furthermore $\operatorname{Supp}_{\mathbb K }(D') = \varnothing$, but $\operatorname{Supp}_{\mathbb Z}(D') = \{ (1:0:0) \}$.
\end{enumerate}
\end{exem}

\begin{prop}\label{prop:normal:R:effective:sum}
We assume that $X$ is locally Noetherian and normal. Let $D$ be an $\mathbb R$-effective $\mathbb R$-Cartier divisor on $X$.
{Then} there are effective Cartier divisors $D_1, \ldots, D_n$ and {positive real numbers} $a_1, \ldots, a_n$
such that $D = a_1 D_1 + \cdots + a_n D_n$.
\end{prop}

\begin{proof}
If $D = 0$, then the assertion is obvious, so that we may assume that $D \not= 0$.
We choose prime divisors $\Gamma_1, \ldots, \Gamma_n$ and $a_1, \ldots, a_n \in \mathbb R_{>0}$ such that
$D = a_1 \Gamma_1 + \cdots + a_n \Gamma_n$ as a Weil divisor. We set
\[
\begin{cases}
V = \{ E = c_1  \Gamma_1 + \cdots + c_n \Gamma_n \,:\, \text{$(c_1, \ldots, c_n) \in \mathbb Q^n$ and $E$ is a $\mathbb Q$-Cartier divisor} \},\\
V_{\mathbb R} := V \otimes_{\mathbb Q} \mathbb R,\quad P = V_{\mathbb R} \cap (\mathbb R_{>0}  \Gamma_1 + \cdots + \mathbb R_{>0} \Gamma_n).
\end{cases}
\]
Then $P$ is an open cone in $V_{\mathbb R}$ and $D \in P$.  Thus the assertion follows.
\end{proof}

We assume that $X$ is projective over a field $k$.
An $\mathbb K$-Cartier divisor $D$ on $X$ is said to be \emph{ample}\index{ample K-Cartier divisor@ample $K$-Cartier divisor} if there are
ample Cartier divisors $D_1, \ldots, D_n$ and $(a_1, \ldots, a_n) \in \mathbb K_{>0}^n$ such that
$D = a_1 D_1 + \cdots + a_n D_n$. 

\begin{prop}\label{prop:ample:R:Cartier:small:perturbation}
Let $A$ be an ample $\mathbb R$-Cartier divisor on $X$ and $D_1, \ldots, D_m$ be Cartier divisors on $X$.
Then there is a positive number $\delta$ such that
$A + \sum_{j=1}^m \delta_j D_j$ is ample for all $\delta_1, \ldots, \delta_m \in \mathbb R$ with
$|\delta_1| + \cdots + |\delta_m| < \delta$. In particular, the ampleness of $\mathbb R$-Cartier divisors is
an open condition.
\end{prop}

\begin{proof}
We choose ample Cartier {divisors} $A_1, \ldots, A_n$ and $(a_1, \ldots, a_n) \in \mathbb R_{>0}^n$ such that
$A = a_1 A_1 + \cdots + a_n A_n$. Let $l$ be a positive rational number such that $lA_1 \pm D_j$ is ample
for all $j=1, \ldots, m$. Note that
\begin{multline*}
A + \sum_{j=1}^m \delta_j D_j = \sum_{j=1}^m |\delta_j|\left( l A_1 + \mathrm{sign}(\delta_j) D_j\right) + 
\left( a_1  - l(|\delta_1| + \cdots + |\delta_m|) \right)A_1 + \sum_{i=2}^m a_i A_i,
\end{multline*}
where
\[
\mathrm{sign}(a) = \begin{cases} 1 & \text{if $a \geqslant 0$}, \\
-1 & \text{if $a < 0$}.
\end{cases}
\]
Therefore, if we choose $\delta = a_1/l$, then $A + \sum_{j=1}^m \delta_j D_j$ is ample.
\end{proof}

\begin{prop}\label{prop:approximation:R:sec:by:Q:sec}
We assume that $X$ is locally Noetherian and normal. Let $D$ be an $\mathbb R$-effective $\mathbb R$-Cartier divisor on $X$.
Let $s_1, \ldots, s_n \in \mathrm{Rat}(X)^{\times} \otimes_{\mathbb Z} \mathbb Q$ and $(a_1, \ldots, a_n) \in \mathbb R^n$ such that 
$a_1, \ldots, a_n$ are linearly independent over $\mathbb Q$ and 
$D + (s_1^{a_1} \cdots s_n^{a_n})$ is $\mathbb R$-effective. Then, 
%
%
for any $\epsilon > 0$, there is a positive number $\delta$ such that
if 
$|a'_1 - a_1| + \cdots + |a'_n - a_n| \leqslant \delta$,
then $(1+\epsilon)D + (s_1^{a'_1} \cdots s_n^{a'_n})$ is $\mathbb R$-effective.

\end{prop}

\begin{proof}
We set $\phi = s_1^{a_1} \cdots s_n^{a_n}$.
Let us see that $\mathrm{Supp}((s_i)) \subseteq \mathrm{Supp}((\phi))$ for all $i$. Otherwise there is a prime divisor $\Gamma$ such that
$\ord_{\Gamma}(s_i) \not= 0$ and $\ord_{\Gamma}(\phi) = 0$, so that
$\sum_{j=1}^n a_j \ord_{\Gamma}(s_j) = 0$, which contradicts to the linear independency of $a_1, \ldots, a_n$ over $\mathbb Q$.
If $\mathrm{Supp}((\phi)) = \emptyset$, then $\mathrm{Supp}((s_i)) = \emptyset$ for all $i$, and hence the assertion is obvious, so that
we may assume that $\mathrm{Supp}((\phi)) \not= \emptyset$.
Let $\Gamma_1, \ldots, \Gamma_m$ be distinct prime divisors such that $\mathrm{Supp}((\phi)) = \Gamma_1 \cup \cdots \cup \Gamma_m$.
Then we can set $(s_i) = \sum_{l=1}^m h_{il} \Gamma_l$ for some $h_{il} \in \mathbb Q$. If we set $\gamma_l = \sum_{i=1}^{n} a_i h_{il}$, then
$((\phi)) = \sum_{l=1}^m \gamma_l \Gamma_l$. As $\mathrm{Supp}((\phi)) = \Gamma_1 \cup \cdots \cup \Gamma_m$, one has
$\gamma_l \not= 0$ for all $l$. We set
\[
L_+ = \{ l \in \{ 1, \ldots, m \} \mid \gamma_l > 0 \}\quad\text{and}\quad L_- = \{ l \in \{ 1, \ldots, m \} \mid \gamma_l < 0 \}.
\]
As $\ord_{\Gamma_l}(D) + \gamma_l \geqslant 0$, one has $\ord_{\Gamma_l}(D) > 0$ for all $l \in L_-$.
We set $C = \max_{i,l} \{ |h_{il}| \}$ and choose $\delta > 0$ such that
\[
C \delta < \min \{ |\gamma_1|, \ldots, |\gamma_m| \}\quad\text{and}\quad C \delta \leqslant \epsilon \min_{l \in L_-} \{ \ord_{\Gamma_l}(D) \}.
\]
Let $a'_1, \ldots, a'_n \in \mathbb R$ such that $|a'_1 - a_1| + \cdots + |a'_n - a_n | \leqslant \delta$. If we set $\gamma'_l = \sum_{i=1}^{n} a'_i h_{il}$,
then $|\gamma'_l - \gamma_l| \leqslant C\delta$, so that $\{ l \mid \gamma'_l > 0 \} = L_+$ and $\{ l \mid \gamma'_l < 0 \} = L_-$.
Further, for $l \in L_-$, we have
\[
(1+\epsilon) \ord_{\Gamma_l}(D) + \gamma'_l = (\ord_{\Gamma_l}(D) + \gamma_l) + (\epsilon \ord_{\Gamma_l}(D) + (\gamma'_l - \gamma_l)) \geqslant \epsilon \ord_{\Gamma_l}(D) - C\delta
\geqslant 0.
\]
Therefore $(1+\epsilon)D + (s_1^{a'_1} \cdots s_n^{a'_n})$ is $\mathbb R$-effective because 
\[
\ord_{\Gamma_l}((1+\epsilon)D + (s_1^{a'_1} \cdots s_n^{a'_n})) = (1+\epsilon) \ord_{\Gamma_l}(D)
+ \gamma'_l
\]
for $l=1, \ldots, m$.
\end{proof}

\section{Green functions}
Let $k$ be a field equipped with an absolute value $|\ndot|$, which is complete. If $|\ndot|$ is Archimedean, we assume that it is the usual absolute value on $\mathbb R$ or $\mathbb C$. Let $X$ be an integral projective scheme over $\Spec k$.

\subsection{Green functions of Cartier divisors}\label{Subsec: Green function}
Let $X^{\mathrm{an}}$ be the Berkovich {topological} space associated with $X$. 
We denote by $C^0_{\mathrm{gen}}(X^{\mathrm{an}})$ the set of all continuous functions on a non-empty Zariski open subset of $X^{\mathrm{an}}$, modulo the following equivalence relation
\[\text{$f\sim g$}
\quad\overset{\mathrm{def}}{\Longleftrightarrow}\quad
\text{$f$ and $g$ coincide on a non-empty Zariski open subset}.\]
Note that the addition and the multiplication of functions induce a structure of $\mathbb R$-algebra on the set $C^0_{\mathrm{gen}}(X^{\mathrm{an}})$. Moreover, for any non-empty Zariski open subset $U$ of $X$, we have a natural $\mathbb R$-algebra homomorphism from $C^0(U^{\mathrm{an}})$ to $C^0_{\mathrm{gen}}(X^{\mathrm{an}})$. Since $U^{\mathrm{an}}$ is dense in $X^{\mathrm{an}}$ (see \cite{Berkovich90} Corollary 3.4.5), we obtain that this homomorphism is injective. Moreover, the $\mathbb R$-algebra $C^0_{\mathrm{gen}}(X^{\mathrm{an}})$\index{C0gen@$C^0_{\mathrm{gen}}$} is actually the colimit of the system $C^0(U^{\mathrm{an}})$ in the category of $\mathbb R$-algebras, where $U$ runs over the set of all non-empty Zariski open subsets of $X$. We say that an element of $C^0_{\mathrm{gen}}(X^{\mathrm{an}})$ \emph{extends to a continuous function on $U^{\mathrm{an}}$}\index{extends to a continuous function} if it belongs to the image of the canonical homomorphism $C^0(U^{\mathrm{an}})\rightarrow C^0_{\mathrm{gen}}(X^{\mathrm{an}})$.

\begin{rema}
Let $f$ be an element of $C^0_{\mathrm{gen}}(X^{\mathrm{an}})$. If $U$ is a non-empty Zariski open subset of $X$ such that $f$ extends to a continuous function on $U^{\mathrm{an}}$, then, by the injectivity of the canonical homomorphism $C^0(U^{\mathrm{an}})\rightarrow C^0_{\mathrm{gen}}(X^{\mathrm{an}})$ there exists a unique continuous function on $U^{\mathrm{an}}$ whose canonical image in $C^0_{\mathrm{gen}}(X^{\mathrm{an}})$ is $f$. Therefore, by gluing of continuous functions we obtain the existence of a largest Zariski open subset $U_f$ {of $X$} such that $f$ extends to a continuous function on $U_f^{\mathrm{an}}$. The set $U^{\mathrm{an}}$ is called the \emph{domain of definition}\index{domain of definition} of the element $f$. By abuse of notation, we still use the expression $f$ to denote the continuous function on $U_f^{\mathrm{an}}$ corresponding to the element $f\in C^0_{\mathrm{gen}}(X^{\mathrm{an}})$.
\end{rema}

If $f$ is a non-zero rational function on $X$, then it is an invertible regular function on a non-empty Zariski open subset $U$ of $X$. Therefore $\ln|f|$ is a continuous function on $U^{\mathrm{an}}$, which determines an element of $C^0_{\mathrm{gen}}(X^{\mathrm{an}})$. Note that this element does not depend on the choice of the non-empty Zariski open subset $U$. We still denote by $\ln|f|$ this element by abuse of notation.

\begin{defi}
Let $D$ be a Cartier divisor on $X$. We call \emph{Green function}\index{Green function} of $D$ any element $g$ of $C^0_{\mathrm{gen}}(X^{\mathrm{an}})$ such that, 
for any local equation $f$ of $D$ over a non-empty Zariski open subset $U$,
the element $g+\ln|f|$ of $C^0_{\mathrm{gen}}(X^{\mathrm{an}})$ extends to a continuous function on $U^{\mathrm{an}}$.
\end{defi}

\begin{exem}
Let $f$ be a non-zero rational function on $X$. Then $\mathrm{div}(f)$ is a Cartier divisor. By definition, $-\ln|f|$ is a Green function of $\mathrm{div}(f)$.
More generally, let $L$ be an invertible $\mathcal O_X$-module, equipped with a continuous metric $\varphi$. Let $s$ a non-zero rational section of $L$. Then the function $-\ln|s|_\varphi$, which is well defined outside of the zero {points and poles} of the section $s$ and is continuous, determines an element of $C^0_{\mathrm{gen}}(X^{\mathrm{an}})$. It is actually a Green function of the divisor $\mathrm{div}(s)$. In particular, we deduce from Remark \ref{Rem:existenceofcontmetr} that, for any Cartier divisor $D$ on $X$, there exists a Green function of $D$.
\end{exem}

\begin{rema}
One can also construct a metrized invertible sheaf from a Cartier divisor equipped with a Green function. Let $D$ be a Cartier divisor on $X$ and $g$ be a Green function of $D$. If $f$ is a rational function of $X$ which defines the divisor $D$ on a non-empty Zariski open subset, then the element $f^{-1}s_D$ is a section of the invertible sheaf $\mathcal O_X(D)$ which trivialises the latter on $U$. Note that the element $-(g+\ln|f|)$ of $C^0_{\mathrm{gen}}(X^{\mathrm{an}})$ extends to a continuous function on $U^{\mathrm{an}}$. We denote by $|f^{-1}s_D|_g$ the exponential of this function, which defines a continuous metric on the restriction of $L$ {to} $U$. By gluing we obtain a continuous metric on $L$ which we denote by $\varphi_g$. By definition one has $g=-\ln|s_D|_g$ in $C^0_{\mathrm{gen}}(X^{\mathrm{an}})$.
\end{rema}

\begin{prop}\phantomsection\label{prop:linear:Green:function}
\begin{enumerate}[label=\rm(\arabic*)]
\item\label{Item: Green function of trivial Cartier divisor} An element $g$ in $C^0_{\mathrm{gen}}(X^{\mathrm{an}})$ is a Green function of the trivial Cartier divisor if and only if it extends to a continuous function on $X^{\mathrm{an}}$.

\item
Let $D$ and $D'$ be Cartier divisors on $X$ and $g, g'$ be Green functions of $D$ and $D'$,
respectively. Then, for $(a, a') \in \mathbb Z^2$, $a g + a'g'$ is a Green function of $aD + a'D'$.
\end{enumerate}
\end{prop}

\begin{proof}
(1) is obvious.

\medskip
(2) Let $f$ and $f'$ be local equations of $D$ and $D'$, respectively.
Then $f^a{f'}^{a'}$ is a local equation of $aD + a'D'$.
As $g + \ln |f|$ and $g' + \ln |f'|$ {extend to continuous functions} locally,
\[
a(g + \ln |f|) + a'(g' + \ln |f'|) = a g + a' g' + \ln \left|f^a\cdot {f'}^{a'}\right|
\]
is locally continuous, as required.
\end{proof}

We denote by $\widehat{\mathrm{Div}}(X)$ the set of all pairs of the form $(D,g)$, where $D$ is a Cartier divisor on $X$ and $g$ is a Green function of $D$. The above proposition shows that $\widehat{\mathrm{Div}}(X)$ forms a commutative group with the composition law \[((D_1,g_1),(D_2,g_2))\longmapsto (D_1+D_2,g_1+g_2).\] One has a natural homomorphism of groups $\widehat{\mathrm{Div}}(X)\rightarrow\mathrm{Div}(X)$ sending $(D,g)$ to $D$. The kernel of this homomorphism is $C^0(X^{\mathrm{an}})$.

\subsection{Green functions for $\mathbb Q$-Cartier and $\mathbb R$-Cartier divisors}
Let $\mathbb K$ be either $\mathbb Q$ or $\mathbb R$.
Let $f$ be an element of $R(X)^{\times} \otimes_{\mathbb Z} \mathbb K$, that is,
\[f = f_1^{a_1} \cdots f_r^{a_r},\quad (f_1, \ldots, f_r) \in (R(X)^{\times})^r\text{ and
} (a_1, \ldots, a_r) \in \mathbb K^r.\] Then one can consider an element of
$C^0_{\mathrm{gen}}(X^{\mathrm{an}})$ given by
$a_1 \ln |f_1| + \cdots + a_r \ln | f_r|$, which dose not depend on the choice of
the expression $f = f_1^{a_1} \cdots f_r^{a_r}$. Indeed, let $f = g_1^{b_1} \cdots g_{l}^{b_{l}}$
be another expression of $f$.
Let us choose an affine open set $U = \Spec(A)$ such that $f_1, \ldots, f_r, g_1, \ldots, g_l$ {belong to  $A^{\times}$.}
For $x \in U^{\mathrm{an}}$, as the seminorm $|\ndot|_x$ is multiplicative, 
$|\ndot|_x$ naturally extends to a map $|\ndot|_x : A^{\times} \otimes_{\mathbb Z} \mathbb K \to
\mathbb R$, so that $|f_1|_x^{a_1} \cdots |f_r|_x^{a_r} = |g_1|_x^{b_1} \cdots |g_{l}|_x^{b_{l}}$.
Therefore,
\[
a_1 \ln |f_1| + \cdots + a_r \ln | f_r| = b_1 \ln |g_1| + \cdots + b_l \ln |g_l|
\]
on $U^{\mathrm{an}}$, which shows the assertion. We denote the above function by $\ln |f|$.

\begin{defi}
Let $D$ be {a} $\mathbb K$-Cartier divisor on $X$.
We say an element $g\in C^0_{\mathrm{gen}}(X^{\mathrm{an}})$ is a \emph{$D$-Green function}\index{Green function} or \emph{Green function of $D$} if,
for any point $x \in X^{\mathrm{an}}$ and {any} local equation $f$ of $D$ on a Zariski neighbourhood of $j(x)$, $g + \log |f|$ extends {to}
a continuous function around $x$.
We denote by $\widehat{\mathrm{Div}}_{\mathbb K}(X)$ the set of all pairs of the form $(D,g)$, where $D$ is a $\mathbb K$-Cartier divisor on $X$ and $g$ is a Green function of $D$. Note that $\widehat{\mathrm{Div}}_{\mathbb K}(X)$ is actually a vector space over $\mathbb K$, which is the quotient of $\widehat{\mathrm{Div}}(X)\otimes_{\mathbb Z}\mathbb K$  by the vector subspace generated by elements of the form $\lambda(D,g)-(\lambda D,\lambda g)$, where $(D,g)\in\widehat{\mathrm{Div}}(X)$ and $\lambda\in\mathbb K$.
\end{defi}

\begin{prop}
\begin{enumerate}[label=\rm(\arabic*)]
\item\label{Item: Green function of the trivial divisor}
Let $g$ be a Green function of the trivial $\mathbb K$-Cartier divisor. Then $g$ extends to a continuous function on $X^{\mathrm{an}}$.

\item\label{Item: additivity of Green functions}
Let $D$ and $D'$ be $\mathbb K$-Cartier divisors on $X$ and $g, g'$ be Green functions of $D$ and $D'$,
respectively. Then, for $a, a' \in \mathbb K$, $a g + a'g'$ is a Green function of $aD + a'D'$.
\end{enumerate}
\end{prop}

\begin{proof}
It can be proven in the same way as Proposition~\ref{prop:linear:Green:function}.
\end{proof}

\begin{prop}\label{Pro:e-gextension}
Let $\mathbb K$ be either $\mathbb Z$ or $\mathbb Q$ or $\mathbb R$.
Let $D$ be an effective $\mathbb K$-Cartier divisor on $X$ and $g$ be a Green function of $D$. Then the element $\mathrm{e}^{-g}$ of $C^0_{\mathrm{gen}}(X)$ extends to a non-negative continuous function on $X^{\mathrm{an}}$.  
\end{prop}
\begin{proof}
Locally on a Zariski open subset $U = \Spec(A)$ of $X$, the divisor $D$ is defined by 
$f_1^{a_1} \cdots f_r^{a_r}$ (where $f_1, \ldots, f_r$ are elements of $A \setminus \{ 0 \}$ and $a_1, \ldots, a_r$ are elements of $\mathbb K_{>0}$) and  
the element $g+\ln|f|$ of $C^0_{\mathrm{gen}}(X^{\mathrm{an}})$ extends to a continuous function on $U^{\mathrm{an}}$. Hence $\mathrm{e}^{-g}=|f|\cdot \mathrm{e}^{-(g+\ln|f|)}$ extends to a continuous function on $X^{\mathrm{an}}$, which is non-negative. 
\end{proof}

\begin{defi}\label{Def: H0 metric}
Let $\mathbb K$ be either $\mathbb Z$ or $\mathbb Q$ or $\mathbb R$.
\begin{enumerate}[label=\rm(\arabic*)]
\item
For $f \in H_{\mathbb K}^0(D)$, $|f|\exp(-g)$ extends to a continuous function. Indeed,
as $D + (f)$ is effective and $g - \ln |f|$ is a Green function of $D + (f)$, by the above proposition,
$|f|\exp(-g) = \exp(-(g - \ln |f|))$ is a continuous function.
We denote the function $|f|\exp(-g)$ by $|f|_g$. Moreover, $\sup \{ |f|_g(x) \,:\, x \in X^{\mathrm{an}} \}$
is denoted by $\| f \|_{g}$.

\item
Let $D$ be an effective $\mathbb K$-Cartier divisor on $X$ and $g$ be a Green function of $X$. By abuse of notation, we use the expression $g$ to denote the map $-\ln(\mathrm{e}^{-g}):X^{\mathrm{an}}\rightarrow \mathbb R\cup\{+\infty\}$, where $\mathrm{e}^{-g}$ is the non-negative continuous function on $X^{\mathrm{an}}$ described in Proposition \ref{Pro:e-gextension}. We say that an element $(D,g)$ of $\widehat{\mathrm{Div}}(X)$ or $\widehat{\mathrm{Div}}_{\mathbb R}(X)$ is \emph{effective}\index{effective} if $D$ is effective and 
the map $g$ takes non-negative values.

\item Let $\overline{D} = (D, g)$ be an element of $\widehat{\mathrm{Div}}_{\mathbb K}(X)$.
We define $\widehat{H}^0_{\mathbb K}(\overline{D})$ to be 
\[
\widehat{H}^0_{\mathbb K}(\overline{D}) := \{ f \in R(X)^{\times} \,:\,
\text{$\overline{D}+\widehat{(f)}$ is effective}\}\cup\{0\}.
\]
Note that $\widehat{H}^0_{\mathbb K}(\overline{D}) = \{ f \in H^0_{\mathbb K}(D) \,:\, \|f \|_g \leqslant 1 \}$.
\end{enumerate}
\end{defi}

\begin{rema}\label{Rem: linear bijection isometry}
Let $(D,g)$ be an element of $\widehat{\mathrm{Div}}_{\mathbb K}(X)$ and $s$ be an element of $R(X)^{\times}$. Let $(D',g')=(D+\mathrm{div}(s),g-\ln|s|)$. Then the map $H^0_{\mathbb K}(D')\rightarrow H^0_{\mathbb K}(D)$ sending $f\in H^0_{\mathbb K}(D')$ to $fs$ is a bijection. Moreover, {for any $f\in H^0_{\mathbb K}(D')$ one has $\norm{f}_{g'}=\norm{fs}_g$}.
\end{rema}

%% file: ch3_2019_03_23.tex

\chapter{Adelic curves}

\IfChapVersion
\ChapVersion{Version of Chapter 3 : \\ \StrSubstitute{\DateChapThree}{_}{\_}}
\fi

The theory of ad\`eles in the study of global fields was firstly introduced by Chevalley \cite[Chapitre III]{Chevalley51} in the function field setting and by Weil \cite{Weil51} in the number field setting. This theory allows to consider all places of a global field in a unified way. It also leads to a uniform approach in the geometry of numbers in global fields, either via the adelic version of Minkowski's theorems and Siegel's lemma developed by McFeat \cite{McFeat71}, Bombieri-Vaaler \cite{Bombieri-Vaaler83}, Thunder \cite{Thunder96}, Roy-Thunder \cite{Roy_Thunder96}, or via the study of adelic vector bundles developed by Gaudron \cite{Gaudron08}, generalising the slope theory introduced by {Stuhler \cite{Stuhler76}, Grayson \cite{Grayson84} and} Bost \cite{BostBour96,Bost2001}. The adelic point of view is also closely related to the Arakelov geometry approach to the height theory in arithmetic geometry. Recall that the Arakelov height theory has been developed by Arakelov \cite{Arakelov74,Arakelov75}, Szpiro \cite{Szprio85}, Faltings \cite{Faltings91}, Bost-Gillet-Soul\'e \cite{BGS94},  (compare to the approach of Philippon \cite{Philippon86}, see also \cite{Soule91} for the comparison of these approaches). We refer the readers to \cite{Zhang95a} for an application of the Arakelov height theory in the adelic setting to the Bogomolov problem. The Arakelov height theory has been generalized by Moriwaki \cite{Moriwaki00} to setting of finitely generated field over a number field (see also \cite{Moriwaki04} for a panoramic view of this theory).

The purpose of this chapter is to develop a formal setting of adelic curves, which permits to include the above examples of global fields and finitely generated extensions of global fields, as well as less standard examples such as the trivial absolute value, polarised projective varieties and arithmetic varieties, and the combination of different adelic structures. More concretely, we consider a field equipped with a family of absolute values on the field, indexed by a measure space, which verifies a ``\emph{product formula}''\index{product formula@product formula} (see Section \ref{SubSec:adelicarithcurve} below). This construction is similar to that of $M$-field introduced by Gubler \cite{Gubler97} (see \cite{BPS14} for the height theory of toric varieties in this setting, and the work of Ben Yaakov and Hrushovski in the model theory framework). Moreover, Gaudron and R\'emond \cite{Gaudron_Remond14} have studied Siegel's lemma for fields of algebraic numbers with a similar point of view. However, our main concern is to establish a general setting with which we can develop not only the height theory but also the geometry of adelic vector bundles and birational Arakelov theory. Therefore our choice is  different from the previous works. In particular, we require that the absolute values are well defined for all places (same as the setting of \emph{globally valued field}\index{globally valued field@globally valued field} of Ben Yaakov and Hrushovski, compare to \cite[\S2]{Gubler97}) and we pay a particular attention to the algebraic coverings of adelic curves and the measurability properties (see Sections \ref{Sec:separableext}-\ref{Sec:algebracexte}). We prove that, for any adelic curve $S$ with underlying field $K$ and any algebraic extension $L$ of $K$, there exists a natural structure of adelic curve on $L$ whose measure space is fibred over that of $S$ with a disintegration kernel (compare to \cite{Gaudron_Remond14}). Curiously, even in the simplest case of finite separable extensions, this result is far from simple (see for example Theorem \ref{Thm:mesurabilite}). The main subtleties appear in the proof of the measurability of the fibre integral, which is neither classic in the theory of disintegration of measures nor in the extension of absolute values in algebraic number theory. We combine the monotone class theorem (in a functional form) in measure theory with divers technics in algebra and number theory such as symmetric polynomials and Vandermonde matrix to resolve this problem.

The chapter is organized as follows. In Section \ref{SubSec:adelicarithcurve} we give the definition of adelic curves and discuss several basic measurability properties concerning Archimedean absolute values. Various examples of adelic curves are presented in Section 1.2. In the subsequent two sections, we discuss algebraic extensions of adelic curves. The finite separable extension case is treated in Section 1.3, where we establish the measurability of fibre integrals (Theorem \ref{Thm:mesurabilite}) and the construction of the extended adelic curve (Theorem \ref{Thm:constructionofextension}). Section \ref{Sec:algebracexte} is devoted to the generalisation of these results to arbitrary algebraic extensions case, where the compatibility of the construction in the situation of successive extensions (Theorem \ref{Thm:successivealgextension}) is proved. These results will serve as the fundament for the geometry of adelic vector bundles and birational Arakelov geometry over adelic curves developed in further chapters.

\section{Definition of Adelic curves}\label{SubSec:adelicarithcurve}
Let $K$ be a commutative field 
and $M_K$ be the set of all absolute values on $K$.
We call \emph{adelic structure on $K$}\index{adelic structure on K@adelic structure on $K$} a measure space $(\Omega,\mathcal A,\nu)$ equipped with a map $\phi:\omega\mapsto |\ndot|_\omega$ from $\Omega$ to $M_K$
satisfying the following properties:
\begin{enumerate}[label=\rm(\roman*)]
\item
$\mathcal A$ is a $\sigma$-algebra on $\Omega$ and $\nu$ is a measure on $(\Omega,\mathcal A)$;
\item
for any $a\in K^{\times}:=K\setminus\{0\}$, the function $\omega\mapsto\ln|a|_\omega$ is $\mathcal A$-measurable, integrable with respect to the measure $\nu$.
\end{enumerate}
The data $(K,(\Omega,\mathcal A,\nu),\phi)$ is called an \emph{adelic curve}\index{adelic curve@adelic curve}.
Moreover, the space $\Omega$ and the map $\phi : \Omega \to M_K$ are called
a \emph{parameter space of $M_K$}\index{parameter space of MK@parameter space of $M_K$} and a \emph{parameter map}\index{parameter map@parameter map}, respectively.
We do not require neither the injectivity nor the surjectivity of $\phi$.
Further, if {the equality}
\begin{equation}\label{Equ:productformula}\int_{\Omega}\ln|a|_\omega\,\nu(\mathrm{d}\omega)=0\end{equation}
holds for each $a \in K^{\times}$, then
the adelic curve $(K,(\Omega,\mathcal A,\nu),\phi)$ is said to be \emph{proper}\index{proper@proper}.
The equation \eqref{Equ:productformula} is called
a \emph{product formula}\index{product formula@product formula}.

The set of all $\omega\in\Omega$ such that $|\ndot|_\omega$ is Archimedean (resp. non-Archimedean) is written as $\Omega_\infty$ (resp. $\Omega_{\mathrm{fin}}$).
For any element $\omega\in\Omega$, let $K_\omega$ be the completion of $K$ with respect to the absolute value $|\ndot|_\omega$. Note that $|\ndot|_\omega$ extends by continuity to an absolute value on $K_\omega$ which we still denote by $|\ndot|_\omega$.

\begin{prop}\label{Pro:mesurabilite}
Let $(K,(\Omega,\mathcal A,\nu),\phi)$ be an adelic curve. The set $\Omega_\infty$ of all $\omega\in\Omega$ such that the absolute value $|\ndot|_\omega$ is Archimedean belongs to $\mathcal A$. 
\end{prop}
\begin{proof}
The result of the proposition is trivial if the characteristic of $K$ is positive since in this case the set $\Omega_\infty$ is empty. In the following, we assume that the characteristic of $K$ is zero. Let $f$ be the function on $\Omega$ defined as $f(\omega)=\ln|2|_\omega$. Then $\Omega_\infty=\{\omega\in\Omega\,:\,f(\omega)>0\}$. Hence $\Omega_\infty$ is a measurable set.
\end{proof}

In the case where $|\ndot|_\omega$ is Archimedean, the field $K_\omega$ is equal to $\mathbb R$ or $\mathbb C$. However, the absolute value $|\ndot|_\omega$ does not necessarily identify with the usual absolute value on $K_{\omega}$. By Ostrowski's theorem (see \cite{Neukirch} Chapter II, Theorem 4.2),  there exists a number $\kappa(\omega)\in \intervalle]01]$ such that $|\ndot|_\omega$ equals $|\ndot|^{\kappa(\omega)}$ on $\mathbb Q$, where $|\ndot|$ denotes the usual absolute value on $\mathbb C$. Therefore one has $|\ndot|_\omega=|\ndot|^{\kappa(\omega)}$ on $K_\omega=\mathbb R$ or $\mathbb C$. 

\begin{prop}\label{Pro:mesurabilitekappa} If we extend the domain of definition of the function $\kappa$ to $\Omega$ by taking the value $0$ on $\Omega\setminus\Omega_\infty$,
then the function $\kappa$ is $\mathcal A$-measurable and integrable with respect to $\nu$. In particular, if the function $\kappa$ is bounded from below on $\Omega_\infty$ by a positive number, then one has $\nu(\Omega_{\infty})<+\infty$.
\end{prop}
\begin{proof}
The result of the proposition is trivial if $\Omega_\infty$ is empty. In the following, we assume that $\Omega_\infty$ is non-empty. In this case the field $K$ is of characteristic zero. One has
\[\forall\,\omega\in\Omega_\infty,\quad\ln|2|_\omega=\kappa(\omega)\ln(2),\]
so that
\[
\forall\,\omega\in\Omega,\quad\kappa(\omega) = \frac{\max \{ 0, \ln |2|_{\omega} \}}{\ln(2)}.
\]
Therefore the $\mathcal A$-measurability and $\nu$-integrability of the function $\omega\mapsto\ln|2|_\omega$ imply the results of the proposition.
\end{proof}

\begin{prop}
Let $S=(K,(\Omega,\mathcal A,\nu),\phi)$ be an adelic curve. We assume that the field $K$ is countable. Let $\Omega_0$ be the set of points $\omega\in\Omega$ such that the absolute value $|\ndot|_\omega$ on $K$ is trivial. Then $\Omega_0$ belongs to $\mathcal A$.
\end{prop}
\begin{proof}
By definition,
\[\Omega_0=\bigcap_{a\in K^{\times}}\big\{\omega\in\Omega\,:\,|a|_\omega=1\big\}.\]
Since the function $(\omega\in\Omega)\mapsto |a|_\omega$ is $\mathcal A$-measurable, the set $\{\omega\in\Omega\,:\,|a|_\omega=1\}$ belongs to $\mathcal A$. Since $K^{\times}$ is a countable set, we obtain that $\Omega_0$ also belongs to $\mathcal A$.
\end{proof}

\section{Examples}\label{Sec: Examples of adelic curves} We introduce several fundamental examples of proper adelic curves. Some of them are very classic objects in algebraic geometry or in arithmetic geometry. 

\subsection{\bf Function fields}\label{Subsec: function field}
Let $k$ be a field, $C$ be a regular projective curve over $\Spec k$ and $K$ be the field of all rational functions on $C$. We denote by $\Omega$ the set of all closed points of the curve $C$, equipped with the discrete $\sigma$-algebra $\mathcal A$ (namely $\mathcal A$ is the $\sigma$-algebra of all subsets of $\Omega$). For any closed point $x$ of $C$, the local ring $\mathcal O_{C,x}$ is a discrete valuation ring whose fraction field is $K$. We let $\mathrm{ord}_x(\cdot):K\rightarrow\mathbb Z\cup\{+\infty\}$ be the discrete valuation on $K$ of valuation ring $\mathcal O_{C,x}$ and $|\ndot|_x$ be the absolute value on $K$ defined as \[\forall\,a\in K^\times,\quad |a|_x=\mathrm{e}^{-\mathrm{ord}_x(a)}.\]
Let $n_x:=[k(x):k]$ be the degree of the residue field of $x$. Thus we obtain a map $\phi:\Omega\rightarrow M_K$ sending $x\in\Omega$ to $|\ndot|_x$. We equip the measurable space $(\Omega,\mathcal A)$ with the measure $\nu$ such that $\nu(\{x\})=n_x$. The relation
\[\forall\,a\in K^\times,\quad\sum_{x\in \Omega}n_x\,\mathrm{ord}_x(a)=0\]
shows that the equality 
\[\int_{\Omega}\ln|a|_x\,\nu(\mathrm{d}x)=0\]
holds for any $a\in K^{\times}$. Therefore $(K,(\Omega,\mathcal A,\nu),\phi)$ is a proper adelic curve. 

\subsection{\bf Number fields}\label{Subsec:Numberfields}Let $K$ be a number field. Denote by $\Omega$ the set of all places of $K$, equipped with the discrete $\sigma$-algebra. For any $\omega\in \Omega$, let $|\ndot|_\omega$ be the absolute value on $K$ in the equivalence class $\omega$, which extends either the usual Archimedean absolute value on $\mathbb Q$ or one of the $p$-adic absolute values (such that the absolute value of $p$ is $1/p$). Thus we obtain a map $\phi:\Omega\rightarrow M_K$ sending $\omega\in\Omega$ to $|\ndot|_\omega$. For each 
$\omega\in {\Omega}$, let $n_\omega$ be the local degree $[K_\omega:\mathbb Q_\omega]$, where $K_\omega$ and $\mathbb Q_\omega$ denote respectively the completion of $K$ and $\mathbb Q$ with respect to the absolute value $|\ndot|_\omega$. Let $\nu$ be the measure on $(\Omega,\mathcal A)$ such that $\nu(\{\omega\})=n_\omega$ for any $\omega\in\Omega$. Note that the usual product formula (cf. \cite{Neukirch} Chapter III, Proposition 1.3) asserts that
\[\forall\,a\in K^{\times},\quad\prod_{\omega\in\Omega}|a|_\omega^{[K_\omega:\mathbb Q_\omega]}=1,\]
which can also be written in the form
\[\forall\,a\in K^{\times},\quad\int_{\Omega}\ln|a|_\omega\,\nu(\mathrm{d}\omega)=0.\]
Therefore $(K,(\Omega,\mathcal A,\nu),\phi)$ is a proper adelic curve. 

\subsection{\bf Copies of the trivial absolute value}\label{Subsec: copies of trivial absoulte value}Let $K$ be any field and $(\Omega,\mathcal A,\nu)$ be an arbitrary measure space. For each $\omega\in \Omega$, let $|\ndot|_\omega$ be the trivial absolute value on $K$, namely one has $|a|_\omega=1$ for any $a\in K^{\times}$. We denote by $\phi:\Omega\rightarrow M_K$ the map sending all elements of $\Omega$ to the trivial absolute value on $K$. Then the equality
\[\forall\,a\in K^{\times},\quad\int_{\Omega}\ln|a|_\omega\,\nu(\mathrm{d}\omega)=0\] 
is trivially satisfied. Therefore the data $(K,(\Omega,\mathcal A,\nu),\phi)$ form a proper adelic curve.

\subsection{Polarised varieties}\label{Subsec:polarizedvar} Let $k$ be a field and $X$ be {an integral and} normal projective scheme of dimension $d\geqslant 1$ over $\Spec k$. Let $K=k(X)$ be the field of rational functions on $X$ and $\Omega=X^{(1)}$ be the set of all prime divisors in $X$, equipped with the discrete $\sigma$-algebra $\mathcal A$. We also fix a  family $\{D_i\}_{i=1}^{d-1}$ of ample divisors on $X$. Let $\nu$ be the measure on $(\Omega,\mathcal A)$ such that
\[\forall\,Y\in\Omega=X^{(1)},\quad\nu(\{Y\})=\deg(D_1\cdots D_{d-1}\cap[Y]).\]
Thus we obtain a measure space $(\Omega,\mathcal A,\nu)$.

For each $Y\in\Omega$, let $\mathcal O_{X,Y}$ be the local ring of $X$ on the generic point of $Y$. It is a discrete valuation ring since it is a Noetherian normal domain of Krull dimension $1$. Moreover, its fraction field is $K$. We denote by $\mathrm{ord}_Y(\ndot)$ the corresponding valuation on $K$ and by $|\ndot|_Y$ the absolute value on $K$ with $|\ndot|_Y:=\mathrm{e}^{-\mathrm{ord}_Y(\ndot)}$. Thus we obtain a map $\phi$ from $\Omega$ to the set of all absolute values on $K$, sending $Y\in\Omega$ to $|\ndot|_Y$. 

For any rational function $f\in K^{\times}$, let $(f)$ be the principal divisor associated with $f$, which is
\[(f):=\sum_{Y\in \Omega}\mathrm{ord}_Y(f)\cdot Y.\]
Therefore, the relation $\deg(D_1\cdots D_{d-1}\cdot(f))=0$ can be written as
\[\int_{\Omega}\ln|f|_Y\,\nu(\mathrm{d}Y)=0.\]
Hence $(K,(\Omega,\mathcal A,\nu),\phi)$ is a proper adelic curve.

\subsection{\bf Function field over $\mathbb Q$}\label{subsec:fun:field:Q} Let $K=\mathbb Q(T)$ be the field of rational functions of one variable 
{$T$} and with coefficients in $\mathbb Q$. We consider $K$ as the field of all rational functions on $\mathbb P^1_{\mathbb Q}$. Any closed point $x\in\mathbb P^1_{\mathbb Q}$ defines a discrete valuation on $K$ which we denote by $\mathrm{ord}_x(\ndot)$. Let $\infty$ be the rational point of $\mathbb P^1_{\mathbb Q}$ such that \[\mathrm{ord}_\infty(f/g)=\deg(g)-\deg(f)\]
for polynomials $f$ and $g$ in $\mathbb Q[T]$ such that $g\neq 0$. Then the open subscheme $\mathbb P^1_{\mathbb Q}\setminus\{\infty\}$ is isomorphic to $\mathbb A^1_{\mathbb Q}$. Therefore any closed point $x$ of $\mathbb P^1_{\mathbb Q}$ different from $\infty$ corresponds to an irreducible polynomial $F_x$ in $\mathbb Q[T]$ (up to dilation by a scalar in $\mathbb Q^{\times}$). By convention we assume that $F_x\in\mathbb Z[T]$ and that the coefficients of $F_x$ are coprime. Let $H(x)$ be the Mahler measure of the polynomial $F_x$, defined as 
\[H(x):=\exp\bigg(\int_0^1\ln|F_x(\mathrm{e}^{2\pi it})|\,\mathrm{d}t\bigg).\]
Note that if the polynomial $F_x$ is written in the form 
\[F_x(T)=a_dT^d+\cdots+a_1T+a_0=a_d(T-\alpha_1)\cdots(T-\alpha_d),\]
then one has {(by Jensen's formula, see \cite{MR1554908})}
\[H(x)=|a_d|\prod_{j=1}^d\max\{1,|\alpha_j|\}\geqslant 1.\]
Let $|\ndot|_x$ be the absolute value on $\mathbb Q(T)$ such that 
\[\forall\,\varphi\in\mathbb Q(T),\quad|\varphi|_x=H(x)^{-\ord_x(\varphi)}.\]

For any prime number $p$, let $|\ndot|_p$ be the natural extension to $\mathbb Q(T)$ of the $p$-adic absolute value on $\mathbb Q$ constructed as follows. For any  
\[f=a_dT^d+\cdots+a_1T+a_0\in\mathbb Q[T]\]
let \[|f|_p:=\max_{j\in\{0,\ldots,d\}}|a_j|_p.\]
Note that one has $|fg|_p=|f|_p\cdot|g|_p$ for $f$ and $g$ in $\mathbb Q[T]$ (see \cite{Bombieri_Gubler} Lemma 1.6.3 for example) and thus the function $|\ndot|_p$ on $\mathbb Q[T]$ extends in a unique way to a multiplicative function on $\mathbb Q(T)$. Moreover, the extended function satisfies the strong triangle inequality and therefore defines a non-Archimedean absolute value on $\mathbb Q(T)$.

Denote by $[0,1]_*$ the set of $t\in[0,1]$ such that $\mathrm{e}^{2\pi it}$ is transcendental. For any $t\in [0,1]_*$, let $|\ndot|_t$ be the absolute value on $\mathbb Q(T)$ such that 
\[\forall\,\varphi\in\mathbb Q(T),\quad|\varphi|_t:=
|\varphi({\mathrm{e}^{2\pi it}})|,\]
where $|\ndot|$ denotes the usual absolute value {of $\mathbb C$}. The absolute value $|\ndot|_t$ is Archimedean.

Denote by $\Omega$ the disjoint union $\Omega_{h}\amalg\mathcal P\amalg [0,1]_*$, where $\Omega_{h}$ is the set of all closed points of {$\mathbb P_{\mathbb Q}^1\setminus\{\infty\}$}, and $\mathcal P$ denotes the set of all prime numbers. Let $\phi:\Omega\rightarrow M_K$ be the map sending $\omega\in\Omega$ to $|\ndot|_\omega$. We equip $\Omega_{h}$ and $\mathcal P$ with the discrete $\sigma$-algebras, and $[0,1]_*$ with the restriction of the Borel $\sigma$-algebra on $[0,1]$. Let $\mathcal A$ be the $\sigma$-algebra on $\Omega$ generated by the above $\sigma$-algebras on $\Omega_{h}$, $\mathcal P$ and $[0,1]_*$ respectively. Let $\nu$ be the measure on $\Omega$ such that $\nu(\{x\})=1$ for $x\in{\Omega}_{h}$, that ${\nu}(\{p\})=1$ for any prime number $p$ and that the restriction of $\nu$ on $[0,1]_*$ coincides with the Lebesgue measure. Then for any $f\in K[T]\setminus\{0\}$ one has 
\[\int_{\Omega}\ln|f|_\omega\,\nu(\mathrm{d}\omega)=\sum_{x\in\Omega_{h}}\ln|f|_{x}+\sum_{p\in\mathcal P}\ln|f|_p+\int_{[0,1]_*}\ln|f(\mathrm{e}^{2\pi it})|\,\mathrm{d}t.\]
Since $[0,1]\setminus [0,1]_*$ is negligible with respect to the Lebesgue measure, we obtain that 
\[\int_{[0,1]_*}\ln|f(\mathrm{e}^{2\pi it})|\,\mathrm{d}t=\int_{0}^1\ln|f(\mathrm{e}^{2\pi it})|\,\mathrm{d}t\]
is equal to the logarithm of the Mahler measure of the polynomial $f$. In particular, if we write the polynomial $f$ in the form 
\[f=aF_{x_1}^{r_1}\cdots F_{x_n}^{r_n},\]
where $x_1,\ldots,x_n$ are distinct closed points of $\mathbb P^1_{\mathbb Q}\setminus\{\infty\}$, and $a\in\mathbb Q^{\times}$. Then one has
\[\int_{[0,1]_*}\ln|f(\mathrm{e}^{2\pi it})|\,\mathrm{d}t=\ln|a|+\sum_{j=1}^nr_j\ln H(x_j).\]
Therefore one has  
\[\int_{\Omega}\ln|f|_\omega\,\nu(\mathrm{d}\omega)=\sum_{j=1}^n(-r_j)\ln H(x_j)+\sum_{p\in\mathcal P}\ln|a|_p+\int_{[0,1]_*}\ln|f(\mathrm{e}^{2\pi it})|\,\mathrm{d}t=0.\]
Hence $(K,(\Omega,\mathcal A,\nu),\phi)$ is a proper adelic curve.

This example of proper adelic curve is much less classic and may looks artificial. However, it is actually very natural from the Arakelov geometry point of view. In fact, one can consider $\mathbb Q(T)$ as the field of the rational functions on the arithmetic variety 
{$\mathbb P^1_{\mathbb Z} := \operatorname{Proj}(\mathbb Z[X, Y])$ with $T = X/Y$.} Then the relation
\[\forall\,\varphi\in K^{\times},\quad\int_{\Omega}\ln|\varphi|_{\omega}\,\nu(\mathrm{d}\omega)=0\] 
can be interpreted as 
\[\widehat{\deg}(\widehat{c}_1({\mathcal O_{\mathbb P_{\mathbb Z}^1}(1), \|\ndot\|})\cdot\widehat{(\varphi)})=0,\]
where {$\|\ndot\|$ is the continuous Hermitian metric of $\mathcal O_{\mathbb P_{\mathbb Z}^1}(1)$ given by
\[\| a X + b Y \|(\xi_1 : \xi_2) = \frac{|a \xi_1 + b \xi_2|}{\max \{ |\xi_1|, |\xi_2| \}}\] and}
$\widehat{(\varphi)}$ is the arithmetic divisor associated to the rational function $\varphi$. The integrals of $\ln|\varphi|_\omega$ on $\Omega_{h}$, $\mathcal P$ and $I_*$ correspond to the horizontal, vertical and Archimedean contributions in the arithmetic intersection product, respectively. 

\if01
{\color{blue}\subsection{Function field over a finitely generated subfield of $\mathbb R$}
Let $K_0$ be a subfield of $\mathbb R$ which is finitely generated over $\mathbb Q$. Consider the field $K=K_0(T)$ of rational functions with coefficients in $K_0$. Denote by $[0,1]_{K_0}$ the subset of $[0,1]$ of points $t$ such that $\mathrm{e}^{2\pi it}$ is transcendental over $K_0$. Note that $[0,1]\setminus[0,1]_{K_0}$ is a countable set since $K_0$ is finitely generated over $\mathbb Q$. Let $D^\circ=\{z\in\mathbb C\,:\,|z|<1\}$ be the open unit disc of $\mathbb C$ and $\Omega=D^\circ\amalg [0,1]_{K_0}$, equipped with the $\sigma$-algebra $\mathcal A$ generated by the discrete $\sigma$-algebra on $D^\circ$ and the Borel $\sigma$-algebra on $[0,1]_{K_0}$. Denote by $\nu$ the measure on $(\Omega,\mathcal A)$ such that $\nu(\{x\})=1$ for any $x\in D^\circ$ and that the restriction of $\nu$ on $[0,1]_{K_0}$ identifies with the restriction of the Lebesgue measure. We construction a map $\phi:\Omega\rightarrow M_K$, $(\omega\in\Omega)\mapsto |\ndot|_\omega$ as follows. For any $t\in[0,1]_{K_0}$, let $|f|_t:=|f(\mathrm{e}^{2\pi it})|$ for any $f\in K$, where $|\ndot|$ denotes the usual absolute value on $\mathbb C$. For any $x\in D^\circ$ and any $f\in K$, we consider $f$ as a meromorphic function on $\mathbb C$ and let $|f|_x$ be $\mathrm{e}^{-\mathrm{ord}_x(f)}$, where $\mathrm{ord}_x(f)$ denotes the multiplicity of $x$ as point of zero or pole of $f$ (if $x$ is neither a point of zero nor a pole of $f$, then $\mathrm{ord}_x(f)=0$). Thus $(K,(\Omega,\mathcal A,\nu),\phi)$ forms an adelic curve. This adelic curve is proper since the argument principal shows that, for any $f\in K^{\times}$ one has
\[\int_0^1\ln|f(\mathrm{e}^{2\pi it})|\,\mathrm{d}t=\] 
}  
\fi

\subsection{Polarised arithmetic variety}\label{subsec:polarized;arith:var}
The previous example treated in Subsection~\ref{subsec:fun:field:Q} can be considered as a very particular case of adelic structures arising from polarised arithmetic varieties.
Let $K$ be a finitely generated field over $\mathbb Q$ {and $d$ be its transcendental degree over $\mathbb Q$}. 
Let $k$ be the set of all algebraic elements of $K$ over $\mathbb Q$.
Note that $k$ is a finite extension over $\mathbb Q$.
A normal model of $K$ over $\mathbb Q$ means
{an integral and normal projective scheme} $X$ over $\mathbb Q$ such that
the rational function field of $X$ is $K$.

For simplicity, the set of all $\mathbb C$-valued points of $\Spec(K)$ is denoted by $K(\mathbb C)$,
that is, $K(\mathbb C)$ is the set of all embeddings of $K$ into $\mathbb C$. 
Let $X$ be a \emph{normal {projective} model} of $K$ over $\mathbb Q$, {namely $X$ is an integral normal projective $\mathbb Q$-scheme, whose field of rational functions identifies with $K$}. 
Let $\Spec(K) \to X$ be the canonical morphism.
Considering the composition 
\[
\Spec(\mathbb C) \longrightarrow \Spec(K) \longrightarrow X,
\]
we may treat $K(\mathbb C)$ as a subset of $X(\mathbb C)$.
Note that
\[
K(\mathbb C) = X(\mathbb C) \setminus \bigcup_{Y \subsetneq X} Y(\mathbb C),
\]
where $Y$ runs over all prime divisors on $X$.
Indeed, ``$\subseteq$'' is obvious.
Conversely, let 
$x \in X(\mathbb C) \setminus \bigcup_{Y \subsetneq X} Y(\mathbb C)$. 
Then, for any $f \in K^{\times}$, $f$ has no zero and pole at $x$ as a rational function on 
$X(\mathbb C)$, so that we have a homomorphism $K \to \mathbb C$ given by $f \mapsto f(x)$,
as required.
Note that {the restriction to $K(\mathbb C)$ of the Zariski topology on $X(\mathbb C)$} does not depend on the choice of $X$.
In fact, for any non-empty Zariski open set $U$ of $X$, $K(\mathbb C)$ is a subset of $U(\mathbb C)$, so that
if $X'$ is another normal model of $K$ over $\mathbb Q$
and $U$ is a common open set of $X$ and $X'$, 
then $K(\mathbb C)$ is a subset of $U(\mathbb C)$.
For $x \in K(\mathbb C)$, we set $|\ndot|_x :=|\sigma_x(\ndot)|$, where $\sigma_x$ is the
corresponding embedding $K \hookrightarrow \mathbb C$.

Let $O_k$ be the ring of integers in $k$.
For {any} maximal ideal $\mathfrak p$ of $O_k$, let $v_{\mathfrak p}$ be the absolute value of $k$ given by
\[
v_{\mathfrak p}(\ndot) = \#(O_k/\mathfrak p)^{-\mathrm{ord}_{\mathfrak p}(\ndot)}.
\]
Let $k_{\mathfrak p}$ be the completion of $k$ with respect to $v_{\mathfrak p}$.
By abuse of notation, the natural extension of $v_{\mathfrak p}$ to $k_{\mathfrak p}$ is
also denoted by $v_{\mathfrak p}$. Let $X$ be a normal {projective} model of $K$, 
$X_{\mathfrak p} := X \times_{\Spec(k)} \Spec(k_{\mathfrak p})$
and $h : X_{\mathfrak p} \to X$ the natural projection.
Let $X_{\mathfrak p}^{\mathrm{an}}$ be the analytification of $X_{\mathfrak p}$ in the sense of
Berkovich \cite{Berkovich90}.
For $x \in X_{\mathfrak p}^{\mathrm{an}}$, the associated scheme point of $X_{\mathfrak p}$ is denoted by $p_x$.
We say that $x$ is a {\em generic point} of $X_{\mathfrak p}^{\mathrm{an}}$ if
$h(p_x)$ is the generic point of $X$.
The set of all generic points of $X_{\mathfrak p}^{\mathrm{an}}$ is denoted by 
$K_{\mathfrak p}^{\mathrm{an}}$.
If $U$ is a non-empty Zariski open set of $X$, then
$K_{\mathfrak p}^{\mathrm{an}} \subseteq
U_{\mathfrak p}^{\mathrm{an}}$, so that $K_{\mathfrak p}^{\mathrm{an}}$ and
the {Berkovich} topology of $K_{\mathfrak p}^{\mathrm{an}}$ 
do not depend on the choice of the model $X$.
Moreover, as before, we can see that
\[
K_{\mathfrak p}^{\mathrm{an}} = X_{\mathfrak p}^{\mathrm{an}} \setminus
\bigcup_{Y \subsetneq X} Y^{\mathrm{an}}_{\mathfrak p},
\]
where $Y$ runs over all prime divisors on $X$.
For $x \in K_{\mathfrak p}^{\mathrm{an}}$, the corresponding seminorm $|\ndot|_x$ induces
an absolute value of $K$
because $K$ is contained in the residue field of $X_{\mathfrak p}$ at $p_x$. 
By abuse of notation, it is also denoted by $|\ndot|_x$.

The Zariski-Riemann space 
$\mathrm{ZR}(K/k)$ of $K$ over $k$ is
defined by the set of all discrete valuation rings $\mathcal O$ such that
$k \subseteq \mathcal O \subseteq K$ and the fraction field of $\mathcal O$ is $K$.
For $\mathcal O \in \mathrm{ZR}(K/k)$, the associated valuation of $K$ is denoted by 
$\operatorname{ord}_{\mathcal O}$.
We set $|\ndot|_{\mathcal O} := \exp(-\operatorname{ord}_{\mathcal O}(\ndot))$.

Let $\Omega^{\mathrm{fin}}_{\mathrm{geom}} := \mathrm{ZR}(K/k)$,
$\Omega^{\mathrm{fin}}_{\mathfrak p} := K_{\mathfrak p}^{\mathrm{an}}$,
$\Omega^{\infty} := K(\mathbb C)$ and
\[
\Omega:= \Omega^{\mathrm{fin}}_{\mathrm{geom}} \amalg  \coprod_{\mathfrak p \in \mathrm{Max}(O_k)}
\Omega^{\mathrm{fin}}_{\mathfrak p} \amalg \Omega^{\infty},
\]
where $\mathrm{Max}(O_k)$ is the set of all maximal ideals of $O_k$.
{Let $\phi:\Omega \to M_K$ be the map $\omega \mapsto |\ndot|_{\omega}$.}
Here we consider the $\sigma$-algebra $\mathcal A$ on $\Omega$ generated by
the discrete $\sigma$-algebra on $\Omega^{\mathrm{fin}}$,
the Borel $\sigma$-algebra on $\Omega^{\mathrm{fin}}_{\mathfrak p}$ 
with respect to the topology of $K_{\mathfrak p}^{\mathrm{an}}$ for each 
$\mathfrak p \in \mathrm{Max}(O_k)$, and
the Borel $\sigma$-algebra on $\Omega^{\infty}$ with respect to the topology of $K(\mathbb C)$.
In order to introduce a measure on $(\Omega, \mathcal A)$,
let us fix a normal model $X$ of $K$ and
nef adelic arithmetic $\mathbb R$-Cartier divisors
\[
\overline{D}_1 = (D_1, g_1), \ldots, \overline{D}_d = (D_d, g_d)
\]
of $C^0$-type on $X$ (for details of adelic arithmetic $\mathbb R$-Cartier divisors, see \cite{Moriwaki16}).
The collection $(X; \overline{D}_1, \ldots, \overline{D}_d)$ is called
a {\em polarisation of $K$}.
Let $X^{(1)}$ be the set of all prime divisors on $X$.
The Radon measure on $X^{\mathrm{an}}_{\mathfrak p}$ given by
\[\varphi \longmapsto \widehat{\deg}_{\mathfrak p}((D_1, g_{1,\mathfrak p}) \cdots (D_d, g_{d,\mathfrak p}); \varphi)\]
is denoted by $\mu_{(D_1, g_{1,\mathfrak p}), \ldots, (D_d, g_{d,\mathfrak p})}$.
A measure $\nu$ on $\Omega$ is defined as follows:
$\nu$ on $\Omega^{\mathrm{fin}}_{\mathrm{geom}}$ is a discrete measure given by 
\[
\nu(\{ \mathcal O \}) =
\begin{cases}
\widehat{\deg}\left(\overline{D}_1 \cdots \overline{D}_d \cdot (\Gamma, 0)\right) & \text{if $\mathcal O = \mathcal O_{X, \Gamma}$
for some $\Gamma \in X^{(1)}$}, \\
0 & \text{otherwise},
\end{cases}
\]
$\nu$ on $\Omega^{\mathrm{fin}}_{\mathfrak p}$ is the restriction of 
$2\mu_{(D_1, g_{1,\mathfrak p}), \ldots, (D_d, g_{d,\mathfrak p})}$ 
to $K^{\mathrm{an}}_{\mathfrak p}$,
and
$\nu$ on $\Omega^{\infty}$ is given by $2 c_1(D_1, g_{1,\infty}) \wedge \ldots \wedge
c_1(D_d, g_{d,\infty})$.
Then $(K, (\Omega, \mathcal A, \nu))$ yields a proper adelic structure of $K$. Indeed,
for each $f \in K^{\times}$, the product formula can be checked as follows:
\begin{multline*}
\int_{\Omega} \ln |f|_{\omega} \,\nu(d\omega) = \sum_{\Gamma \in X^{(1)}} 
- \operatorname{ord}_{\Gamma}(f) \widehat{\deg}\left(\overline{D}_1 \cdots \overline{D}_d \cdot (\Gamma,0) \right) \\
\qquad\qquad\qquad+ \sum_{\mathfrak p \in \mathrm{Max}(O_k)}
\int_{K_{\mathfrak p}^{\mathrm{an}}} \ln |f|^2 d \mu_{(D_1, g_{1,\mathfrak p}), \ldots, (D_d, g_{d,\mathfrak p})} \\
+ \int_{K(\mathbb C)} \ln |f|^2c_1(D_1, g_{1,\infty}) \wedge \ldots \wedge
c_1(D_d, g_{d,\infty}).
\end{multline*}
For a proper subvariety of $Y$ of $X$, 
$Y_{\mathfrak p}^{\mathrm{an}}$ and $Y(\mathbb C)$ are null sets with respect to the measures
$\mu_{(D_1, g_{1,\mathfrak p}), \ldots, (D_d, g_{d,\mathfrak p})}$ and
$c_1(D_1, g_{1,\infty}) \wedge \ldots \wedge c_1(D_d, g_{d,\infty})$, respectively.
In addition, we have only countably many prime divisors on $X$. Therefore, 
the above equation implies 
\begin{multline*}
\int_{\Omega} \ln |f|_{\omega} \,\nu(d\omega) = 
- \widehat{\deg}\left(\overline{D}_1 \cdots \overline{D}_d \cdot ((f),0) \right) \\
\qquad\qquad\qquad + \sum_{\mathfrak p \in \mathrm{Max}(O_k)}
\int_{X_{\mathfrak p}^{\mathrm{an}}} \ln |f|^2 d \mu_{(D_1, g_{1,\mathfrak p}), \ldots, (D_d, g_{d,\mathfrak p})} \\
+ \int_{X(\mathbb C)} \ln |f|^2c_1(D_1, g_{1,\infty}) \wedge \ldots \wedge
c_1(D_d, g_{d,\infty}).
\end{multline*}
On the other hand,
\begin{align*}
0 & = \widehat{\deg}\left(\overline{D}_1 \cdots \overline{D}_d \cdot \widehat{(f)} \right) \\
& =
\widehat{\deg}\left(\overline{D}_1 \cdots \overline{D}_d \cdot \left((f), \sum_{\mathfrak p \in \mathrm{Max}(O_k)} -
\ln |f|^2[\mathfrak p] - \ln |f|^2[\infty] \right)\right)\\
& = \widehat{\deg}\left(\overline{D}_1 \cdots \overline{D}_d \cdot ((f),0) \right) \\
& \qquad - \sum_{\mathfrak p \in \mathrm{Max}(O_k)}
\int_{X_{\mathfrak p}^{\mathrm{an}}} \ln |f|^2 d\mu_{(D_1, g_{1,\mathfrak p}), \ldots, (D_d, g_{d,\mathfrak p})} \\
& \qquad\qquad -\int_{X(\mathbb C)} \ln |f|^2c_1(D_1, g_{1,\infty}) \wedge \ldots \wedge
c_1(D_d, g_{d,\infty}),
\end{align*}
as desired.

This proper adelic structure is denoted by $S(X; \overline{D}_1, \ldots, \overline{D}_d)$.

\subsection{Amalgamation of adelic structures}
Let $K$ be a field, \[\big((\Omega,\mathcal A,\nu),\phi:\Omega\rightarrow M_K\big)\;\text{ and }\;\big((\Omega',\mathcal A',\nu'),\phi':\Omega'\rightarrow M_K\big)\] be two adelic structures on $K$. Then the disjoint union of measure spaces \[(\Omega,\mathcal A,\nu)\amalg (\Omega',\mathcal A',\nu')\] together with the map  $\Phi:\Omega\amalg\Omega'\rightarrow M_K$ extending both $\phi$ and $\phi'$ form also an adelic structure on $K$. 
If $S$ and $S'$ denote the adelic curves $(K,(\Omega,\mathcal A,\nu),\phi)$ and $(K,(\Omega',\mathcal A',\nu'),\phi')$ respectively, then  we use the expression $S\amalg S'$ to denote the adelic curve \[\big(K,(\Omega,\mathcal A,\nu)\amalg(\Omega',\mathcal A',\nu'),\Phi\big),\] called the \emph{amalgamation}\index{amalgamation of adelic curves@amalgamation of adelic curves} of the adelic curves $S$ and $S'$. Similarly, one can define the amalgamation for any finite family of adelic arithmetic structures on the field $K$.
Note that if $S$ and $S'$ are proper, then {$S \amalg S'$} is also proper. In fact, for any $a\in K^\times$ one has
\[\int_{\Omega}\ln|a|_\omega\,\nu(\mathrm{d}\omega)+\int_{\Omega'}\ln|a|_\omega\,\nu'(\mathrm{d}\omega)=0.\]

\subsection{Restriction of adelic structure to a subfield}\label{Subsec: adelic curve restriction}
Let $S=(K,(\Omega,\mathcal A,\nu),\phi)$ be an adelic curve and let $K_0$ be a subfield of $K$. Let $\phi_0:\Omega\rightarrow M_{K_0}$ be the map sending $\omega\in\Omega$ to the restriction of $|\ndot|_\omega$ to $K_0$. Then $\phi_0$ defines an adelic structure on $K_0$, called the \emph{restriction}\index{restriction of the adelic structure@restriction of the adelic structure} of the adelic structure of $S$ {to} $K_0$. If $S$ is proper, then its restriction to $K_0$ is also proper.

\subsection{Restriction of adelic structure to a measurable subset}

Let $S=(K,(\Omega,\mathcal A,\nu),\phi)$ be an adelic curve and $\Omega_0$ be an element of $\mathcal A$. Let $\mathcal A_0$ be the restriction of the $\sigma$-algebra $\mathcal A$ {to} $\Omega_0$ and $\nu_0$ be the restriction of the measure {to} $\Omega_0$. Then $(K,(\Omega_0,\mathcal A_0,\nu_0),\phi|_{\Omega_0})$ is an adelic curve, called the \emph{restriction of $S$ {to} $\Omega_0$}\index{restriction of S to Omega0@restriction of $S$ to $\Omega_0$}. Note that this adelic curve is not necessarily proper, even if the adelic curve $S$ is proper.

\section{Finite separable extensions}\label{Sec:separableext}Let $S=(K,(\Omega_K,\mathcal A_K,\nu_K),\phi_K)$ be an adelic curve.  Let $K'/K$ be a finite and separable  extension. For each $\omega\in\Omega_K$, let $M_{K',\omega}$ be the set of all absolute values on $K'$ which extend the absolute value $|\ndot|_\omega$ on $K$. Let $\Omega_{K'}$ be the disjoint union
\[\coprod_{\omega\in\Omega_K}M_{K',\omega}.\]
One has a natural projection $\pi_{K'/K}:\Omega_{K'}\rightarrow\Omega_{K}$ which sends the elements of $M_{K',\omega}$ to $\omega$. Let $\phi_{K'}:\Omega_{K'}\rightarrow M_{K'}$ be the map induced by the inclusion maps $M_{K',\omega}\rightarrow M_{K'}$. If $x$ is an element {of} $\Omega_{K'}$, we also use the expression $|\ndot|_x$ to denote the corresponding absolute value. Note that the following diagram is commutative
\[\xymatrix{\relax \Omega_{K'}\ar[r]^-{\pi_{K'/K}}\ar[d]_-{\phi_{K'}}&\Omega_K\ar[d]^-{\phi_K}\\
M_{K'}\ar[r]_-{\varpi_{K'/K}}&M_K}\] 
and identifies $\Omega_{K'}$ with the fibre product of $M_{K'}$ and $\Omega_K$ over $M_K$ in the category of sets, where $\varpi_{K'/K}$ sends any absolute value on $K'$ to its restriction {to} $K$.
We equip the set $\Omega_{K'}$ with the $\sigma$-algebra $\mathcal A_{K'}$ generated by $\pi_{K'/K}$ and all real-valued functions of the form $(x\in \Omega_{K'})\mapsto |a|_x$, where $a$ runs over $K'$. Namely it is the smallest $\sigma$-algebra on $\Omega_{K'}$ which makes these maps measurable\footnote{A map $f : X' \to X$ of measurable spaces $(X', \mathcal A')$ and $(X, \mathcal A)$ is said to be measurable if $f^{-1}(B) \in \mathcal A'$ for all $B \in \mathcal A$.}, where we consider the $\sigma$-algebra $\mathcal A_K$ on $\Omega_K$ and the Borel $\sigma$-algebra on $\mathbb R$.

We aim to construct a measure $\nu_{K'}$ on the measurable space $(\Omega_{K'},\mathcal A_{K'})$ such that the direct image of $\nu_{K'}$ by $\pi_{K'/K}$ coincides with $\nu_K$. Note that on each fibre $M_{K',\omega}$ of $\pi_{K'/K}$ there is a natural probability measure $\mathbb P_{K',\omega}$ such that
\begin{equation}\label{Equ:Pomega}\forall\,x\in M_{K',\omega},\quad\mathbb P_{K',\omega}(\{x\})=\frac{[K_x':K_\omega]}{[K':K]}.\end{equation}
We refer to \cite{Neukirch} Chapter II, Corollary 8.4 for a proof of the equality
\begin{equation}\label{Equ:average}\sum_{x\in M_{K',\omega}}\frac{[K_x':K_w]}{[K':K]}=1.\end{equation}
Intuitively the family of probability measures $\{\mathbb P_{K',\omega}\}_{\omega\in\Omega_K}$ should form the disintegration of the measure $\nu_{K'}$ with respect to $\nu_K$. However, as we will show below, the construction of the measure $\nu_{K'}$ relies actually on a subtil application of the monotone class theorem and the properties of extensions of absolute values.

\subsection{ Integration along fibres}
If $f$ is a function on $\Omega_{K'}$ valued in $\mathbb R$, we define $I_{K'/K}(f)$ to be the function on $\Omega_K$ which sends $\omega\in\Omega_K$ to
\[\sum_{x\in M_{K',\omega}}\frac{[K'_x:K_\omega]}{[K':K]}f(x).\]
This is an $\mathbb R$-linear operator from the vector space of all real-valued functions on $\Omega_{K'}$ to that of all real-valued functions on $\Omega_K$. The equality \eqref{Equ:average} shows that $I_{K'/K}$ sends the constant function $1$ on $\Omega_{K'}$ to that on $\Omega_K$. The following properties of the linear operator $I_{K'/K}$ are straightforward.

\begin{prop}\label{Pro:coefficent}
Let $f$ be a real-valued function on $\Omega_{K'}$ and $\varphi$ be a real-valued function on $\Omega_K$. Let $\widetilde{\varphi}=\varphi\circ\pi_{K'/K}$. Then 
\[I_{K'/K}(\widetilde{\varphi}f)=\varphi I_{K'/K}(f).\]
\end{prop}
\begin{proof}
By definition, for any $\omega\in\Omega_K$ one has
\[\big({I_{K'/K}}(\widetilde{\varphi}f)\big)(\omega)=\sum_{x\in M_{K',\omega}}\frac{[K_x':K_\omega]}{[K':K]}\widetilde{\varphi}(x)f(x)=\sum_{x\in M_{K',\omega}}\frac{[K_x':K_\omega]}{[K':K]}\varphi(\omega)f(x).\]
\end{proof}

\begin{prop}\label{Pro:integrationsucc}
Let $K''/K'/K$ be successive finite separable extensions of fields. Let $f$ be a real-valued function on $\Omega_{K'}$ and $\widetilde f=f\circ\pi_{K''/K'}$, where $\pi_{K''/K'}:\Omega_{K''}\rightarrow\Omega_{K'}$ sends any absolute value in $\Omega_{K'',\omega}$ to its restriction {to} $K'$, viewed as an element in $\Omega_{K',\omega}$ ($\omega\in\Omega_K$). Then one has
\begin{equation}\label{Equ:intsucc}I_{K''/K}(\widetilde f)=I_{K'/K}(f).\end{equation}
\end{prop}
\begin{proof}
For any $\omega\in\Omega_K$, one has
\[\begin{split}&\quad\;\big(I_{K''/K}(\widetilde f)\big)(\omega)=\sum_{y\in M_{K'',\omega}}\frac{[K''_y:K_\omega]}{[K'':K]}\widetilde f(y)\\&=\sum_{x\in M_{K',\omega}}\sum_{\begin{subarray}{c}y\in M_{K'',\omega}\\y|_{K'}=x\end{subarray}}\frac{[K_y'':K'_x]}{[K'':K']}\cdot\frac{[K'_x:K_\omega]}{[K':K]}f(x),\end{split}\]
which can also be written as
\[\sum_{x\in M_{K',\omega}}\frac{[K'_x:K_\omega]}{[K':K]}f(x)\sum_{\begin{subarray}{c}y\in M_{K'',\omega}\\
y|_{K'}=x
\end{subarray}}\frac{[K_y'':K_x']}{[K'':K']}.\]
Therefore the desired equality follows from the relation
\[\sum_{\begin{subarray}{c}y\in M_{K'',\omega}\\
y|_{K'}=x
\end{subarray}}\frac{[K_y'':K_x']}{[K'':K']}=1.\]
\end{proof}

\begin{coro}\label{Cor:projection of adelic curves}
Let $\varphi$ be an $\mathcal A_K$-measurable function on $\Omega_K$. One has
\[I_{K'/K}(\varphi\circ\pi_{K'/K})=\varphi.\]
\end{coro}
\begin{proof}
It suffices to apply the previous proposition to the successive extensions $K'/K/K$ and then use the fact that $I_{K/K}$ is the identity map to obtain the result.
\end{proof}

\subsection{Measurability of fibre integrals}
The following theorem shows that the operator $I_{K'/K}$ sends an $\mathcal A_{K'}$-measurable function to an $\mathcal A_K$-measurable function. This result is fundamental in the construction of a suitable measure on the measurable space $(\Omega_{K'},\mathcal A_{K'})$.
\begin{theo}\label{Thm:mesurabilite}
For any real-valued $\mathcal A_{K'}$-measurable function $f$, the function $I_{K'/K}(f)$ is $\mathcal A_K$-measurable.
\end{theo}
\begin{proof}
{\it Step 1:} We first prove that, if $a$ is a primitive element of the finite separable extension $K'/K$ (namely $K'=K(a)$) and if $f_a$ is the function on $\Omega_{K'}$ sending $x\in\Omega_{K'}$ to $|a|_x$, then the function $I_{K'/K}(f_a)$ is $\mathcal A_K$-measurable. 

Let $K^{\mathrm{ac}}$ be an algebraic closure of $K$ containing $K'$. For each $\omega\in\Omega_K$, we extend the absolute value $|\ndot|_\omega$ to $K^{\mathrm{ac}}$ via an embedding of $K^{\mathrm{ac}}$ into an algebraic closure $K_\omega^{\mathrm{ac}}$ of $K_\omega$. We still denote by $|\ndot|_\omega$ the extended absolute value on $K^{\mathrm{ac}}$ by abuse of notation.

\begin{lemm}\label{Lem:mesurabilites1}
Let {$d\in\mathbb N_{\geqslant 1}$ and} $\{\alpha_1,\ldots,\alpha_d\}$ be a finite family of {distinct} elements in $K^{\mathrm{ac}}$. For any $\omega\in\Omega_K$, one has
\begin{equation}\label{Equ:product}\max_{j\in\{1,\ldots,d\}}|\alpha_j|_\omega=\limsup_{N\rightarrow+\infty}\bigg|\sum_{i=1}^d\alpha_i^N\bigg|_\omega^{\frac{1}{N}}.\end{equation}
{Moreover, for any $c\in K^{\mathrm{ac}}$, the function 
\[(\omega\in\Omega_K)\longmapsto\max_{\tau\in\operatorname{Aut}_K(K^{\mathrm{ac}})}|\tau(c)|_\omega\]
is $\mathcal A_K$-measurable.}
\end{lemm}
\begin{proof} 
First of all, by the triangle inequality one has
\[\bigg|\sum_{i=1}^d\alpha_i^N\bigg|_\omega^{\frac{1}{N}}\leqslant d^{1/N}{\max\{|\alpha_1|_\omega,\ldots,|\alpha_d|_\omega\}}.\]
Therefore
\[{\max\{|\alpha_1|_\omega,\ldots,|\alpha_d|_\omega\}}\geqslant\limsup_{N\rightarrow+\infty}\bigg|\sum_{i=1}^d\alpha_i^N\bigg|_\omega^{\frac{1}{N}}.\]

Without loss of generality, we can assume that 
\[|\alpha_1|_\omega=\ldots=|\alpha_\ell|_\omega>|\alpha_{\ell+1}|_\omega\geqslant\ldots\geqslant|\alpha_d|_\omega,\]
where $\ell\in\{1,\ldots,d\}$.
For $i\in\{1,\ldots,\ell\}$, let $\beta_i=\alpha_i/\alpha_1$. One has $\beta_1=1$ and 
\[|\beta_1|_\omega=\ldots=|\beta_\ell|_\omega=1.\]
For any integer $N\geqslant 1$, one has
\[\begin{pmatrix}
\beta_1^N+\cdots+\beta_{\ell}^N\\
\vdots\\
\beta_1^{N+\ell-1}+\cdots+\beta_\ell^{N+\ell-1}
\end{pmatrix}=\begin{pmatrix}
1&\ldots&1\\
\beta_1^{1}&\ldots&\beta_\ell^{1}\\
\vdots&\ddots&\vdots\\
\beta_1^{\ell-1}&\ldots&\beta_\ell^{\ell-1}
\end{pmatrix}\begin{pmatrix}\beta_1^N\\
\vdots\\
\beta_\ell^N
\end{pmatrix}\]
Let $\widehat{K^{\mathrm{ac}}}$ be the completion of $K^{\mathrm{ac}}$ with respect to $|\ndot|_{\omega}$. We equip the vector space {$(\widehat{K^{\mathrm{ac}}})^{\ell}$} with the  following norm
\[\forall\,(z_1,\ldots,z_\ell)\in {(\widehat{K^{\mathrm{ac}}}})^\ell,\quad \|(z_1,\ldots,z_\ell)\|:=\max_{i\in\{1,\ldots,\ell\}}|z_i|_\omega. \]
Then the vector $(\beta_1^N,\ldots,\beta_\ell^N)$ has norm $1$ with respect to $\|\ndot\|$. Moreover, the Vandermonde matrix above is invertible since $\beta_1,\ldots,\beta_\ell$ are distinct. Therefore the norm of the vector
\[\bigg(\sum_{i=1}^\ell\beta_i^N,\ldots,\sum_{i=1}^\ell\beta_{i}^{N+\ell-1}\bigg)\]
is bounded from below by a positive constant which does not depend on $N$, which shows that the sequence 
\[\bigg|\sum_{i=1}^d\Big(\frac{\alpha_i}{\alpha_1}\Big)^N\bigg|_{\omega},\quad N\in\mathbb N,\,N\geqslant 1\]
does not converge to zero when $N\rightarrow+\infty$. This implies that
\[\limsup_{N\rightarrow+\infty}\bigg|\sum_{i=1}^d\alpha_i^N\bigg|_\omega^{\frac{1}{N}}\geqslant|\alpha_1|_\omega={\max\{|\alpha_1|_\omega,\ldots,|\alpha_d|_\omega\}}.\]
The equality \eqref{Equ:product} is thus proved.   

{We now proceed with the proof of the second statement.} Let \[T^d-\lambda_1T^{d-1}+\cdots+(-1)^d\lambda_d\in K[T]\] 
be the minimal polynomial of {$c$ over $K$}, and $\alpha_1,\ldots,\alpha_d$ be its roots in $K^{\mathrm{ac}}$. Since the extension $K'/K$ is separable, these roots are distinct. By definition, for $k\in\{1,\ldots,d\}$ one has
\[\lambda_k=\sum_{\begin{subarray}{c}(i_1,\ldots,i_k)\in\{1,\ldots,d\}^k\\ i_1<\ldots<i_k\end{subarray}}\alpha_{i_1}\cdots\alpha_{i_k}.\]
By the fundamental theorem on symmetric polynomials (see for example \cite[\S10-11]{Edwards84}), if $F$ is a polynomial in $K[X_1,\ldots,X_d]$ which is invariant by the action of the symmetric group $\mathfrak S_d$ by permuting the variables, then there exists a polynomial $G\in K[T_1,\ldots,T_d]$ such that
\[F(\alpha_1,\ldots,\alpha_d)=G(\lambda_1,\ldots,\lambda_d).\]
In particular, one has $F(\alpha_1,\ldots,\alpha_d)\in K$ and hence the function
\[(\omega\in\Omega_K)\longmapsto |F(\alpha_1,\ldots,\alpha_d)|_\omega\]
is $\mathcal A_K$-measurable.  For any $N\in\mathbb N$, $N\geqslant 1$, the sum $\sum_{i=1}^d\alpha_i^N$ can be written as a symmetric polynomial evaluated at $(\alpha_1,\ldots,\alpha_d)$, thus the function 
\[(\omega\in\Omega_K)\longmapsto \bigg|\sum_{i=1}^d\alpha_i^N\bigg|_\omega\]
is $\mathcal A_K$-measurable. Combining this observation with the equality \eqref{Equ:product}, we obtain that the function 
\[{(\omega\in\Omega_K)\longmapsto
\max_{\tau\in\operatorname{Aut}_K(K^{\mathrm{ac}})}}|\tau(c)|_\omega=\max_{i\in\{1,\ldots,d\}}|\alpha_i|_\omega\]
is $\mathcal A_K$-measurable.
\end{proof}

{We now continue with the proof of the statement that the function $I_{K'/K}$ is $\mathcal A_K$-measurable. Let $\{\gamma_1,\ldots,\gamma_n\}$ be the orbit of $a$ under the action of $\operatorname{Aut}_K(K^{\mathrm{ac}})$. For any $\omega\in\Omega_K$, let $(s_1(\omega),\ldots,s_n(\omega))$ be the array $(|\gamma_1|_\omega,\ldots,|\gamma_n|_\omega)$ sorted in the decreasing order. Let $k$ be an arbitrary element of $\{1,\ldots,n\}$. For any $\omega\in\Omega_K$, one has
\[s_1(\omega)\cdots s_k(\omega)=\max_{\begin{subarray}{c}(i_1,\ldots,i_k)\in\{1,\ldots,n\}^k\\
i_1<\ldots<i_k\end{subarray}}\max_{\tau\in\operatorname{Aut}_K(K^{\mathrm{ac}})}|\tau(\gamma_{i_1}\cdots\gamma_{i_k})|_\omega.\]
By Lemma \ref{Lem:mesurabilites1}, we obtain that the function $s_1\cdots s_k$ is $\mathcal A_K$-measurable. Therefore all the functions $s_1,\ldots,s_n$ on $\Omega_K$ are $\mathcal A_K$-measurable. } In particular, if $f_a:\Omega_{K'}\rightarrow\mathbb R$ is the function sending $x\in\Omega_{K'}$ to $|a|_x$, where $a$ is the primitive element of the finite separable extension $K'/K$ fixed in the beginning of the step, then for any $\omega\in\Omega_K$ one has (we refer the readers to \cite[page 163]{Neukirch} for the second equality)
\[{\big(I_{K'/K}(f_a)\big)(\omega)=\sum_{x\in M_{K',\omega}}\frac{[K'_x:K_\omega]}{[K':K]}f_a(x)=\frac{1}{n}
\sum_{i=1}^n|\gamma_i|_\omega=\frac{1}{n}\sum_{i=1}^ns_i(\omega),}\]
which implies that $I_{K'/K}(f_a)$ is $\mathcal A_K$-measurable.

{\it Step 2:} We then prove that, for any element $b\in K'$, the function $I_{K'/K}(f_b)$ on $\Omega_K$ is $\mathcal A_K$-measurable, where $f_b$ denotes the function on $\Omega_{K'}$ sending $x\in\Omega_{K'}$ to $|b|_x$.

We consider the sub-extension $K(b)/K$ of $K'/K$. It is a finite and separable extension of $K$ and $b$ is a primitive element. Let $g$ be the function on $\Omega_{K(b)}$ sending $y\in\Omega_{K(b)}$ to $|b|_y$. One has
\[f_b=g\circ\pi_{K'/K(b)},\]
where the map $\pi_{K'/K(b)}:\Omega_{K'}\rightarrow\Omega_{K(b)}$ is defined as in Proposition \ref{Pro:integrationsucc}. By \eqref{Equ:intsucc}, one obtains
\[I_{K'/K}(f_b)=I_{K(b)/K}(g).\]
By the result obtained in Step 1, the function $I_{K(b)/K}(g)$ on $\Omega_K$ is $\mathcal A_K$-measurable. This proves the measurability of the function $I_{K'/K}(f_b)$.

{\it Step 3:} We are now able to apply the monotone class theorem to prove the announced measurability property. 

Let $\mathcal H$ be the set of all non-negative and bounded functions $f$ on $\Omega_{K'}$ such that the function $I_{K'/K}(f)$ on $\Omega_K$ is $\mathcal A_K$-measurable. Note that the constant function $1$ on $\Omega_{K'}$ belongs to $\mathcal H$ since $I_{K'/K}(1)$ coincides with the constant function $1$ on $\Omega_K$. If $f$ and $g$ are two functions in $\mathcal H$ such that $f\geqslant g$, then $f-g\in\mathcal H$ since $I_{K'/K}(f-g)=I_{K'/K}(f)-I_{K'/K}(g)$ is $\mathcal A_K$-measurable. Moreover, Proposition \ref{Pro:coefficent} shows that, if $f$ and $g$ are two functions in $\mathcal H$, and $\varphi$ and $\psi$ are two non-negative and bounded $\mathcal A_K$-measurable functions on $\Omega_K$, $\widetilde\varphi=\varphi\circ\pi_{K'/K}$ and $\widetilde\psi=\psi\circ\pi_{K'/K}$, then the function $\widetilde\varphi f+\widetilde\psi g$ belongs to $\mathcal H$ since \[I_{K'/K}(\widetilde\varphi f+\widetilde\psi g)=\varphi I_{K'/K}(f)+\psi I_{K'/K}(g)\]
is $\mathcal A_K$-measurable. Finally, the operator $I_{K'/K}$ preserves pointwise limit. Therefore, if $\{f_n\}_{n\in\mathbb N}$ is a uniformly bounded sequence of functions in $\mathcal H$ which converges pointwisely to a function $f$, then one has $f\in\mathcal H$. These properties show that $\mathcal H$ is a $\lambda$-family (see Definition \ref{Def:lambda-family}) on $\Omega_{K'}$.

Let $\mathcal C$ be the set of all non-negative and bounded functions on $\Omega_{K'}$ which can be written in the form $f_b\widetilde{\varphi}$, where $b$ is an element of $K'$, $\varphi$ is a non-negative and bounded $\mathcal A_K$-measurable function on $\Omega_K$, and $\widetilde{\varphi}=\varphi\circ\pi_{K'/K}$. Note that  if $b_1$ and $b_2$ are two elements of $K'$ then one has
\[f_{b_1b_2}=f_{b_1}f_{b_2}.\]
Therefore the family $\mathcal C$ is stable by multiplication. By the result obtained in Step 2 and Proposition \ref{Pro:coefficent}, we obtain that $\mathcal C\subseteq\mathcal H$. The monotone class theorem (Theorem \ref{Thm:monotoneclass}) then implies that the family $\mathcal H$ contains all non-negative and bounded $\sigma(\mathcal C)$-measurable functions. Finally, any non-negative $\sigma(\mathcal C)$-measurable function on $\Omega_{K'}$ can be written as the limit of an increasing sequence of non-negative and \emph{bounded} $\sigma(\mathcal C)$-measurable functions, and any real-valued $\sigma(\mathcal C)$-measurable function is the difference {of} two non-negative $\sigma(\mathcal C)$-measurable functions. Therefore, for any real-valued $\sigma(\mathcal C)$-measurable function $f$, the function $I_{K'/K}(f)$ is $\mathcal A_K$-measurable.

{\it Step 4:} It remains to prove that the $\sigma$-algebras $\sigma(\mathcal C)$ and $\mathcal A_{K'}$ are the same. Clearly one has $\sigma(\mathcal C)\subseteq\mathcal A_{K'}$ since any function in $\mathcal C$ is $\mathcal A_{K'}$-measurable. To prove the equality it suffices to show that any function of the form $f_b$ with $b\in K'$ is $\sigma(\mathcal C)$-measurable. Let \[T^m+\mu_1T^{m-1}+\cdots+\mu_m\in K[T]\] be the minimal polynomial of $b$ over $K$. Then for any $\omega\in\Omega_K$ and any $x\in M_{K',\omega}$, one has
\begin{equation}\label{Equ:troncation}|b|_x\leqslant m\cdot\max\{1,|\mu_1|_{\omega},\ldots,|\mu_m|_{\omega}\}\end{equation}
since otherwise one should have
\[1>\frac{|\mu_1|_{\omega}}{|b|_x}+\cdots+\frac{|\mu_m|_\omega}{|b|_x}\geqslant\sum_{i=1}^m\frac{|\mu_i|_\omega}{|b|_x^i},\]
which contradicts the equality
\[b^m=-\mu_1b^{m-1}-\cdots-\mu_m.\]
For any $N\in\mathbb N$, let $A_N$ be the set
\[\{\omega\in\Omega_K\,:\,\max\{|\mu_1|_\omega,\ldots,|\mu_m|_\omega\}\leqslant N\}\in\mathcal A_K.\]
The relation \eqref{Equ:troncation} shows that the function $f_b\cdot(\indic_{A_N}\circ\pi_{K'/K})$ is non-negative and bounded. Hence it belongs to $\mathcal C$. Finally, since 
\[f_b=\lim_{N\rightarrow+\infty}f_b\cdot(\indic_{A_N}\circ\pi_{K'/K}),\]
we obtain that the function $f_b$ is $\sigma(\mathcal C)$-integrable. The theorem is thus proved.
\end{proof}

\subsection{Construction of the measure}
In this subsection, we describe the construction of a suitable measure on the measurable space $(\Omega_{K'},\mathcal A_{K'})$ to form an adelic structure on $K'$ and prove some compatibility results.

\begin{defi}
We denote by $\nu_{K'}:\mathcal A_{K'}\rightarrow \mathbb R_+\cup\{+\infty\}$ the map defined as follows:
\[\forall\,A\in\mathcal A_{K'},\quad \nu_{K'}(A):=\int_{\Omega_K}I_{K'/K}(\indic_A)\,\mathrm{d}\nu_K.\]
By Theorem \ref{Thm:mesurabilite}, this map is well defined.
\end{defi}

\begin{theo}\phantomsection\label{Thm:constructionofextension}\begin{enumerate}[label=\rm(\arabic*)]\item\label{Item: nu K' is a measure}
The map $\nu_{K'}$ is a measure on the measurable space $(\Omega_{K'},\mathcal A_{K'})$ such that, for any non-negative $\mathcal A_{K'}$-measurable function $f$ on $\Omega_{K'}$ one has 
\begin{equation}\label{Equ:convergence}\int_{\Omega_{K'}}f\,\mathrm{d}\nu_{K'}=\int_{\Omega_K}I_{K'/K}(f)\,\mathrm{d}\nu_K.\end{equation}
\item\label{Item: criterion of measurability} A real-valued $\mathcal A_{K'}$-measurable function $f$ on $\Omega_{K'}$ is integrable with respect to $\nu_{K'}$ if and only if $I_{K'/K}(|f|)$ is integrable with respect to $\nu_{K}$. Moreover, the equality \eqref{Equ:convergence} also holds for all real-valued $\mathcal A_K'$-measurable functions on $\Omega_{K'}$ which are integrable with respect to $\nu_{K'}$.
\item\label{Item: direct image of the measure} The direct image of the measure $\nu_{K'}$ by the measurable map $\pi_{K'/K}$ coincides with $\nu_K$, namely for any real-valued $\mathcal A_K$-measurable function $\varphi$ on $\Omega_K$ which is non-negative (resp. integrable with respect to $\nu_K$), the function $\varphi\circ\pi_{K'/K}$ is non-negative (resp. integrable with respect to $\nu_{K'}$), and one has
\begin{equation}\label{Equ:directimage}\int_{\Omega_{K'}}(\varphi\circ\pi_{K'/K})\,\mathrm{d}\nu_{K'}=\int_{\Omega_K}\varphi\,\mathrm{d}\nu_K.\end{equation} 
\item\label{Item: Adelic curve base change} $S' = (K',(\Omega_{K'},\mathcal A_{K'},\nu_{K'}),\phi_{K'})$ is an adelic curve. 
\item\label{Item: productformula:norm}
For $b \in K' \setminus \{ 0 \}$, one has
\begin{equation}\label{Equ:productformula:norm}
[K': K] \int_{\Omega_{K'}} \ln |b|_x\, \mathrm{d} \nu_{K'} = \int_{\Omega_{K}} \ln | N_{K'/K}(b) |_{\omega}\, \mathrm{d} \nu_K,
\end{equation}
where $N_{K'/K}(b)$ is the norm of $b$ with respect to the extension $K'/K$.
In particular, if $S$ is proper, then $S'$ is also proper.
\end{enumerate}
\end{theo}
\begin{proof}
\ref{Item: nu K' is a measure} The operator $I_{K'/K}$ preserve pointwise limits. Therefore, if $\{A_n\}_{n\in\mathbb N}$ is a countable family of disjoint sets in $\mathcal A_{K'}$ and if $A=\bigcup_{n\in\mathbb N}A_n$, one has
\[\nu_{K'}(A)=\int_{\Omega_K}I_{K'/K}(\indic_A)\,\mathrm{d}\nu_K=\int_{\Omega_K}\sum_{n\in\mathbb N}I_{K'/K}(\indic_{A_n})\,\mathrm{d}\nu_K=\sum_{n\in\mathbb N}\nu_{K'}(A_n),\]
where the last equality comes from the monotone convergence theorem.

The set of all non-negative and bounded $\mathcal A_{K'}$-measurable functions $f$ which verify the equality \eqref{Equ:convergence} forms a $\lambda$-family. Moreover, this $\lambda$-family contains the set of all functions of the form $\indic_{A}$ ($A\in\mathcal A_{K'}$), which is stable by multiplication. By Theorem \ref{Thm:monotoneclass}, we obtain that the equality \eqref{Equ:convergence} actually holds for all non-negative and bounded $\mathcal A_{K'}$-measurable functions, and hence holds for general non-negative $\mathcal A_{K'}$-measurable functions by the monotone convergence theorem again. 

\ref{Item: criterion of measurability} The equality \eqref{Equ:convergence} clearly implies that a real-valued $\mathcal A_{K'}$-measurable function $f$ on $\Omega_{K'}$ is integrable with respect to $\nu_{K'}$ if and only if $I_{K'/K}(|f|)$ is integrable with respect to $\nu_{K}$. Moreover, if $f$ is a real-valued $\mathcal A_{K'}$-measurable function on $\Omega_{K'}$ which is integrable with respect to $\nu_{K'}$, then the equality \eqref{Equ:convergence} applied to $\max(f,0)$ and $-\min(f,0)$ shows that 
\[\int_{\Omega_{K'}}\max(f,0)\,\mathrm{d}\nu_{K'}=\int_{\Omega_K}I_{K'/K}(\max(f,0))\,\mathrm{d}\nu_K\]
and
\[\int_{\Omega_{K'}}(-\min(f,0))\,\mathrm{d}\nu_{K'}=\int_{\Omega_K}I_{K'/K}(-\min(f,0))\,\mathrm{d}\nu_K.\]
Since these numbers are finite, the difference of the above two equalities leads to
\[\int_{\Omega_{K'}}f\,\mathrm{d}\nu_{K'}=\int_{\Omega_K}I_{K'/K}(f)\,\mathrm{d}\nu_K.\]

\ref{Item: direct image of the measure} By the two assertions  proved above, one has
\[\int_{\Omega_{K'}}(\varphi\circ\pi_{K'/K})\,\mathrm{d}\nu_{K'}=\int_{\Omega_K}I_{K'/K}(\varphi\circ\pi_{K'/K})\,\mathrm{d}\nu_K.\] 
By Corollary \ref{Cor:projection}, one has
\[I_{K'/K}(\varphi\circ\pi_{K'/K})=\varphi.\]
Thus we obtain \eqref{Equ:directimage}.

\ref{Item: Adelic curve base change} 
Let $b$ be an element in $K'\setminus\{0\}$ and $f_b$ be the function on $\Omega_{K'}$ sending $x\in\Omega_{K'}$ to $|b|_x$. Let $\lambda=N_{K'/K}(b)$ be the norm of $b$ with respect to the extension $K'/K$. For any $\omega\in\Omega_K$ one has (see \cite{Neukirch} Chapter II, Corollary 8.4 and page 161)
\[\prod_{x\in M_{K',\omega}}|b|_x^{[K'_x:K_\omega]}=|\lambda|_\omega,\]
which implies that 
\begin{equation}\label{Equ:fibreintf}I_{K/K'}(\ln f_b)=\frac{1}{[K':K]}\ln f_\lambda,\end{equation}
where 
$f_\lambda$ is the function  on $\Omega_K$ which sends $\omega\in\Omega_K$ to $|\lambda|_\omega$.

Let $\Omega_{K,\infty}$ be the set of all $\omega\in\Omega_{K}$ such that $|\ndot|_\omega$ is an Archimedean absolute value. By Proposition \ref{Pro:mesurabilite}, this is an element of $\mathcal A_K$. Similarly, let $\Omega_{K',\infty}$ be the set of all $x\in \Omega_{K'}$ such that the absolute value $|\ndot|_x$ is Archimedean. One has
$\Omega_{K',\infty}=\pi_{K'/K}^{-1}(\Omega_{K,\infty})$.
We will prove the integrability of $\ln f_b$ on $\Omega_{K'}\setminus\Omega_{K',\infty}$ and on $\Omega_{K',\infty}$, respectively. 
For this purpose we use a refinement of the method in the Step 4 of the proof of Theorem \ref{Thm:mesurabilite}. 

Let \[T^m+\mu_1T^{m-1}+\cdots+\mu_m\in K[T]\] be the minimal polynomial of $b$ over $K$. Then for any $\omega\in\Omega_{K}\setminus\Omega_{K,\infty}$ and any $x\in M_{K',\omega}$, one has
\begin{equation*}|b|_x\leqslant \max\{1,|\mu_1|_{\omega},\ldots,|\mu_m|_{\omega}\}.\end{equation*}
Otherwise one should have
\[1>\max_{i\in\{1,\ldots,m\}}\frac{|\mu_i|_\omega}{|b|_x}\geqslant\max_{i\in\{1,\ldots,m\}}\frac{|\mu_i|_\omega}{|b|_x^i}.\]
However, the equality
\[b^m=-\mu_1b^{m-1}-\cdots-\mu_m.\] 
implies that
\[|b|_x^m\leqslant\max_{i\in\{1,\ldots,m\}}|\mu_i|_\omega\cdot|b|_x^{m-i},\]
which leads to a contradiction.  Therefore, if we denote by $g$ the function
\begin{equation}\label{Equ:foncg}\omega\longmapsto\max\{0,\ln|\mu_1|_\omega,\ldots\ln|\mu_m|_\omega\}\end{equation} on $\Omega_K$, then $\ln f_b$ is bounded from above by $g\circ\pi_{K'/K}$ on $\Omega_K\setminus\Omega_{K,\infty}$. Moreover, by the definition of adelic curves, the functions $\omega\mapsto\ln|\mu_i|_\omega$ is integrable for any $i\in\{1,\ldots,m\}$, hence also is the function $g$. 
The function \[(g\circ\pi_{K'/K}-\ln f_b)\indic_{\Omega_{K'}\setminus\Omega_{K',\infty}}\] is  non-negative, and 
\[\begin{split}&\quad\;I_{K'/K}((g\circ\pi_{K'/K}-\ln f_b)\indic_{\Omega_{K'}\setminus\Omega_{K',\infty}})\\&=\big(I_{K'/K}(g\circ\pi_{K'/K})-I_{K'/K}(\ln f_b)\big)\indic_{\Omega_K\setminus\Omega_{K,\infty}}
\\&=\Big(g-\frac{1}{[K':K]}\ln f_\lambda\Big)\indic_{\Omega_K\setminus\Omega_{K,\infty}}\end{split}\] 
is an integrable function with respect to $\nu_K$, where the first equality comes from Proposition \ref{Pro:coefficent} and the fact that $I_{K'/K}$ is a linear operator, and the second equality comes from Corollary \ref{Cor:projection} and \eqref{Equ:fibreintf}. By the second assertion of the theorem, the function $(g\circ\pi_{K'/K}-\ln f_b)\indic_{\Omega_{K'}\setminus\Omega_{K',\infty}}$ is integrable, and hence also is $(\ln f_b)\indic_{\Omega_{K'}\setminus\Omega_{K',\infty}}$.

We now consider the Archimedean case. We assume that $\Omega_\infty$ is non-empty. Then the characteristic of the field $K$ is zero. In particular, it contains $\mathbb Q$ as its prime field.  Moreover, for any $x\in M_{K',\omega}$, one has
\begin{equation}|b|_x\leqslant m^{\kappa(\omega)}\cdot \max\{1,|\mu_1|_{\omega},\ldots,|\mu_m|_{\omega}\},\end{equation}
where $\kappa(\omega)$ is the exponent of $|\ndot|_\omega$ as a power of the usual absolute value on $K_\omega=\mathbb R$ or $\mathbb C$.
Otherwise one should have
\[1>\frac{|\mu_1|_{\omega}^{1/\kappa(\omega)}}{|b|_x^{1/\kappa(\omega)}}+\cdots+\frac{|\mu_m|_\omega^{1/\kappa(\omega)}}{|b|_x^{1/\kappa(\omega)}}\geqslant\sum_{i=1}^m\frac{|\mu_i|_\omega^{1/\kappa(\omega)}}{|b|_x^{i/\kappa(\omega)}}.\]
However, the equality
\[b^m=-\mu_1b^{m-1}-\cdots-\mu_m.\] 
implies that
\[|b|_x^{m/\kappa(\omega)}\leqslant |\mu_1|^{1/\kappa(\omega)}_{\omega}|b|_x^{(m-1)/\kappa(\omega)}+\cdots+|\mu_m|_{\omega}^{1/\kappa(\omega)}\]
since $|\ndot|_\omega^{1/\kappa(\omega)}$ and $|\ndot|_{x}^{1/\kappa(\omega)}$ are absolute values on $K'$ and $K$ respectively (which extend the usual absolute value on $\mathbb Q$). Therefore, the function $\ln f_b$ is bounded from above by $(\ln(m)\kappa+g)\circ\pi_{K'/K}$ on $\Omega_{K,\infty}$, where $g$ is the function defined in \eqref{Equ:foncg}, and we have extended the function $\kappa$ on $\Omega_{K}$ by taking the value $0$ on $\Omega_K\setminus\Omega_{K,\infty}$. Since the function $\ln(m)\kappa+g$ is integrable, by the same argument as in the non-Archimedean case, we obtain the integrability of the function $\ln f_b$ on $\Omega_{K',\infty}$.

\ref{Item: productformula:norm}
follows from  \eqref{Equ:convergence} and \eqref{Equ:fibreintf}.

\end{proof}

\section{General algebraic extensions}\label{Sec:algebracexte}

Let $S=(K,(\Omega_K,\mathcal A_K,\nu_K),\phi_K)$ be an adelic curve. In this section, we consider the construction of adelic curves from $S$ whose underlying field are general algebraic extensions of $K$. 

\subsection{Finite extension}
Let $K''$ be a finite extension of $K$ and $K'$ be the separable closure of $K$ in the field $K''$. By the result of the previous section, one can construct an adelic structure on the field $K'$ which we denote by $((\Omega_{K'},\mathcal A_{K'},\nu_{K'}),\phi_{K'})$.

Note that $K''$ is a purely inseparable extension of $K'$ (see \cite{Bourbaki_A4-7} Chapter V, \S7, no.7, Proposition 13.a). If $q$ is the degree of the extension $K''/K'$, then for any $\alpha\in K''$ one has $\alpha^q\in K'$  (see \cite{Bourbaki_A4-7} Chapter V, \S5, no.1, Proposition 1). In particular, any absolute value $|\ndot|$ on $K'$ extends in a unique way to $K''$ and one has
\begin{equation}\label{Equ: purely inseparable ext}\forall\,\alpha\in K'',\quad |\alpha|=|\alpha^q|^{1/q},\end{equation}
where $|\alpha|$ denotes the extended absolute value on $K''$ evaluated on $\alpha$, and $|\alpha^q|$ denotes the initial absolute value on $K'$ evaluated on $\alpha^q$. In other words, the sets $M_{K'}$ and $M_{K''}$ are in canonical bijection. This observation permits to construct, for any $\alpha\in K''\setminus\{0\}$ the function 
\[\Omega_{K'}\rightarrow\mathbb R,\quad (x\in \Omega_{K'})\mapsto\ln|\alpha|_x\]
This function is clearly $\mathcal A_{K'}$-measurable since one has
\[\forall\,x\in\Omega_{K'},\quad \ln|\alpha|_x=\frac{1}{q}\ln|\alpha^q|_x.\]
Moreover, it is also integrable with respect to $\nu_{K'}$ and one has
\begin{equation}\label{Eqn:product:formula:pure:inseparable}
\int_{\Omega_{K'}}\ln|\alpha|_x\,\nu_{K'}(\mathrm{d}x)=\frac{1}{q}\int_{\Omega_{K'}}\ln|\alpha^q|_x\,\nu_{K'}(dx). 
\end{equation}
This fact shows that 
\[(K'',(\Omega_{K'},\mathcal A_{K'},\nu_{K'}),\phi_{K'})\]
is actually an adelic curve, where we identify $M_{K'}$ with $M_{K''}$. Note that the relation \eqref{Equ: purely inseparable ext} shows that $\mathcal A_{K'}$ is also the smallest $\sigma$-algebra on $\Omega_{K'}$ which makes the canonical projection map $\Omega_{K'}\rightarrow\Omega_K$ and the functions $(x\in\Omega_{K'})\mapsto |\alpha|_x$ measurable, where $\alpha\in K''$.

\begin{defi}\label{Equ:finiteextension}
We denote by $S\otimes_KK''$ the adelic curve \[(K'',(\Omega_{K'},\mathcal A_{K'},\nu_{K'}),\phi_{K'})\]
constructed as above, called the \emph{finite extension}\index{finite extension@finite extension} of $S$ induced by the extension of fields $K''/K$. We also use the expression $\pi_{K''/K}$ to denote the projection map $\pi_{K'/K}:\Omega_{K'}\rightarrow\Omega_K$ described in the previous section. Similarly, we also use the expression $I_{K''/K}$ to denote the operator $I_{K'/K}$. Note that $\Omega_K'$  identifies also with the fibre product of $\Omega_K$ and $M_{K''}$ over $M_K$ in the category of sets since $M_{K'}$ and $M_{K''}$ are the same. Similarly, $\phi_{K'}$ identifies with the  projection map from $\Omega_{K'}=\Omega_K\times_{M_K}M_{K''}$ to $M_{K''}=M_{K'}$.
Note that if $S$ is proper, then $S \otimes_K K''$ is also proper (cf.
\eqref{Eqn:product:formula:pure:inseparable} and Theorem~\ref{Thm:constructionofextension}, \ref{Item: productformula:norm}).
\end{defi}

In the following, we will prove that the above construction of adelic curves is compatible with successive finite extensions. The lemma below is important for the proof.

\begin{lemm}\label{Lem:completion-separable}
Let $L/K$ be a finite extension of fields, $|\ndot|_v$ be an absolute value on $K$, and $|\ndot|_w$ be an absolute value on $L$ extending $|\ndot|_v$. Let $K^{\mathrm{sc}}$ be the separable closure of $K$ in $L$. Then the completion $K^{\mathrm{sc}}_{w}$ of $K^{\mathrm{sc}}$ with respect to the absolute value $|\ndot|_w$ identifies with the separable closure of $K_v$ in $L_w$.
\end{lemm}
\begin{proof}
The case where $|\ndot|_v$ is Archimedean is trivial since the characteristic of the field $K$ is then zero and hence $K^{\mathrm{sc}}=L$. In the following, we assume that $|\ndot|_v$ is non-Archimedean. We first prove that the extension $K_w^{\mathrm{sc}}/K_v$ is separable. Let $\alpha\in K^{\mathrm{sc}}$ be a primitive element (see \cite{Bourbaki_A4-7} Chapter V, \S7, no.4 Theorem 1 for its existence) of the separable extension $K^{\mathrm{sc}}/K$  and let $F$ be its minimal polynomial. Assume that $F$ is decomposed in $K^{\mathrm{sc}}_v[T]$ into the product of distinct irreducible polynomials as $F=F_1\cdots F_m$. Since $F$ is a separable polynomial, the same are the polynomials $F_1,\ldots,F_m$. For any extension $|\ndot|_w$ of the absolute value $|\ndot|_v$ to $K^{\mathrm{sc}}$, there exists an index $i\in\{1,\ldots,m\}$ such that $K_w^{\mathrm{sc}}\cong K_{v}[T]/(F_i)$ (see \cite{Neukirch} Chapter II, Propositions 8.2 and 8.3). Therefore the extension $K_w^{\mathrm{sc}}/K_v$ is separable.

In the following, we prove that the extension $L_w/K_w^{\mathrm{sc}}$ is purely inseparable. Note that $L_w=K_v(L)=K_w^{\mathrm{sc}}(L)$. Since the extension $L/K^{\mathrm{sc}}$ is purely inseparable, we obtain that the extension $L_w/K_w^{\mathrm{sc}}$ is also purely inseparable since it is generated by purely inseparable elements (see \cite{Bourbaki_A4-7} Chapter V, \S7, no.2, the corollary of Proposition 2). By \cite{Bourbaki_A4-7} Chapter V, \S7, no.7, Proposition 13.c, we obtain that $K_w^{\mathrm{sc}}$ is the separable closure of $K_v$ in $L_w$ since $K_w^{\mathrm{sc}}/K_v$ is separable and $L_w/K_w^{\mathrm{sc}}$ is purely inseparable. 
\end{proof}

\begin{rema}\label{Rem:remark 1.4.3}
Let $K''/K$ be a finite extension of fields and denote by $(K'',(\Omega'',\mathcal A'',\nu''),\phi'')$ the adelic curve $S\otimes_KK''$. The above lemma allows to write the operator $I_{K''/K}$ in the following form: for any real-valued $\mathcal A''$-measurable function $f$ on $\Omega''$
\[\forall\,\omega\in\Omega_K,\quad(I_{K''/K}(f))(\omega)=\sum_{x\in M_{K'',\omega}}\frac{[K_x'':K_\omega]_s}{[K'':K]_s}f(x),\]
where for any finite extension $L_2/L_1$ of fields, the expression $[L_2:L_1]_s$ denotes the separable degree of the extension $L_2/L_1$.
\end{rema}

\begin{prop}\label{Pro:successiveextension}
Let $K_2/K_1/K$ be successive finite extensions of fields. If $((\Omega_K,\mathcal A_K,\nu_K),\phi_K)$ is an adelic structure on $K$ and $S$ is the corresponding adelic curve, then one has \[(S\otimes_KK_1)\otimes_{K_1}K_2=S\otimes_KK_2.\]
Moreover, one has
\begin{equation}\label{Equ:compatibility}\pi_{K_1/K}\circ\pi_{K_2/K_1}=\pi_{K_2/K},\quad I_{K_1/K}\circ I_{K_2/K_1}=I_{K_2/K},\end{equation}
where we have used the conventions of notation described in Definition \ref{Equ:finiteextension}.
\end{prop}
\begin{proof}
We denote by $(K_1,(\Omega_1,\mathcal A_1,\nu_1),\phi_1)$ and $(K_2,(\Omega_2,\mathcal A_2,\nu_2),\phi_2)$ the adelic curves $S\otimes_KK_1$ and $S\otimes_KK_2$ respectively. First of all, set-theoretically one has
\[\Omega_2=\coprod_{\omega\in\Omega_K}M_{K_2,\omega}=\coprod_{\omega\in\Omega_K}\coprod_{x\in M_{K_1,\omega}}M_{K_2,x}=\coprod_{x\in\Omega_1}M_{K_2,x},\] and hence  
\[\pi_{K_2/K}=\pi_{K_1/K}\circ\pi_{K_2/K_1}\]
Moreover, if $f$ is a real-valued function on $\Omega_2$, by Lemma \ref{Lem:completion-separable} and Remark \ref{Rem:remark 1.4.3}, for any $\omega\in\Omega_K$ one has
\[\begin{split}&\quad\;(I_{K_2/K}(f))(\omega)=\sum_{y\in M_{K_2,\omega}}\frac{[K_{2,y}:K_\omega]_s}{[K_2:K]_s}f(y)\\&=\sum_{x\in M_{K_1,\omega}}\frac{[K_{1,x}:K_\omega]_s}{[K_1:K]_s}\sum_{y\in M_{K_2,x}}\frac{[K_{2,y}:K_{1,x}]_s}{[K_2:K_1]_s}f(y)=I_{K_1/K}(I_{K_2/K_1}(f))(\omega),
\end{split}\]
where in the second equality we have used the multiplicativity of the separable degree (see \cite{Bourbaki_A4-7} Chapter V, \S6, no.5). 

We then show that the map $\pi_{K_2/K_1}$ is $\mathcal A_2$-measurable. Since the $\sigma$-algebra $\mathcal A_1$ is generated by $\pi_{K_1/K}$ and functions of the form $f_a:(x\in\Omega_1)\mapsto |a|_x$ with $a\in K_1$, it suffices to prove that the maps $\pi_{K_1/K}\circ\pi_{K_2/K_1}$ and $f_a\circ\pi_{K_2/K_1}$ are $\mathcal A_2$-measurable. We have shown that $\pi_{K_1/K}\circ\pi_{K_2/K_1}=\pi_{K_2/K}$, which is clearly $\mathcal A_2$-measurable by the definition of the adelic curve $S\otimes_KK_2$.  Moreover, if $a$ is an element in $K_1$, then the map $f_a\circ\pi_{K_2/K_1}$ sends $y\in\Omega_2$ to $|a|_y$. Hence it is also $\mathcal A_2$-measurable. In particular, the $\sigma$-algebra $\mathcal A_2$ contains the  $\sigma$-algebra $\mathcal A_2'$ in the adelic structure of $(S\otimes_KK_1)\otimes_{K_1}K_2$, namely the smallest $\sigma$-algebra which makes the map $\pi_{K_2/K_1}$ and all functions of the form $(y\in\Omega_K)\mapsto |\alpha|_y$ measurable, where $\alpha$ runs over $K_2$.

From the equality $\pi_{K_1/K}\circ\pi_{K_2/K_1}=\pi_{K_2/K}$ we also obtain that $\pi_{K_2/K}^{-1}(\mathcal A_K)$ is contained in $\pi_{K_2/K_1}^{-1}(\mathcal A_1)$, and thus is contained in $\mathcal A_2'$ since $\pi_{K_2/K_1}$ is $\mathcal A_2'$-measurable. Since $\mathcal A_2$ is the smallest $\sigma$-algebra on $\Omega_2$ which makes the map $\pi_{K_2/K}$ and all functions of the form $(y\in\Omega_K)\mapsto |\alpha|_y$ measurable, we obtain that $\mathcal A_2\subseteq\mathcal A_2'$. Combining with the result obtained above, we obtain that the $\sigma$-algebras $\mathcal A_2$ and $\mathcal A_2'$ coincide.

Finally, the relation $I_{K_1/K}\circ I_{K_2/K_1}=I_{K_2/K}$ shows that the measure $\nu_2$ coincides with the measure in the adelic structure of the adelic curve $(S\otimes_KK_1)\otimes_{K_1}K_2$. The proposition is thus proved.
\end{proof}

\subsection{General algebraic extensions} We now consider an algebraic extension $L$ of $K$ which is not necessarily finite. Let $\mathscr E_{L/K}$ be the set of all finite extensions of $K$ which are contained in $L$. This set is ordered by the relation of inclusion. Moreover, it is also filtered in the sense that, if $K_1$ and $K_2$ are two finite extensions of $K$ which are contained in $L$, then there exists a finite extension $K_3\in\mathscr E_{L/K}$ such that $K_3\supseteq K_1\cup K_2$.

By the result obtained in the previous subsection, for each element $K''$ in $\mathscr E_{L/K}$, we can equipped $K''$ with a natural adelic structure induced from the adelic structure of $S$, as described in Definition \ref{Equ:finiteextension}. We denote by $((\Omega_{K''},\mathcal A_{K''},\nu_{K''}),\phi_{K''})$ this adelic structure. Moreover, Proposition \ref{Pro:successiveextension} shows that, for successive finite extensions $K_2/K_1/K$ of the field $K$ which are contained in $L$, there exist a natural projection \[\pi_{K_2/K_1}:(\Omega_{K_2},\mathcal A_{K_2})\longrightarrow(\Omega_{K_1},\mathcal A_{K_1})\]
together with a disintegration operator $I_{K_2/K_1}$ from the vector space of all real-valued $\mathcal A_{K_2}$-measurable functions on $\Omega_{K_2}$ to that of all real-valued $\mathcal A_{K_1}$-measurable functions on $\Omega_{K_1}$, which sends $\nu_{K_2}$-integrable functions to $\nu_{K_1}$-integrable functions. These data actually define a functor from a filtered ordered set to the category of measure spaces. Intuitively one can define an adelic structure on $L$ whose measure space part is the projective limit of this functor. However, the projective limite in the category of measure spaces does not exist in general  (the product of infinitely many measures need not make sense). Therefore more careful treatment is needed for our setting of projective system of finite extensions of adelic curves. {Our strategy is to construct  the fibres as projective limits of probability spaces, which always exist.}

For any $\omega\in\Omega_K$, let $M_{L,\omega}$ be the set of all absolute values on $L$ which extend $|\ndot|_\omega$. Let $\Omega_L$ be the disjoint union of all $M_{L,\omega}$ with $\omega$ runs through $\Omega_K$. In other words, $\Omega_L$ is the fibre product of $\Omega_K$ and $M_L$ over $M_K$. The inclusion maps $M_{L,\omega}\rightarrow M_{L}$ define a map from $\Omega_L$ to $M_L$ which we denote by $\phi_L$. Moreover, for any extension $K''\in\mathscr E_{L/K}$ and any $\omega\in\Omega_K$, one has a natural map from $M_{L,\omega}$ to $M_{K'',\omega}$ defined by restriction of absolute values. These maps induce a map from $\Omega_L$ to $\Omega_{K''}$ which we denote by $\pi_{L/K''}$. If $x$ is an element of $\Omega_{K''}$, we denote by $M_{L,x}$ the set of absolute values on $L$ which extends $|\ndot|_x$. It identifies with the inverse image of $\{x\}$ by $\pi_{L/K''}$. If $K_1\subseteq K_2$ are extensions in $\mathscr E_{L/K}$, then one has
\begin{equation}\label{Equ:compatibilityproj}\pi_{L/K_1}=\pi_{K_2/K_1}\circ\pi_{L/K_2}.\end{equation}

\begin{prop}
For any $K''\in\mathscr E_{L/K}$, the map $\pi_{L/K''}:\Omega_L\rightarrow\Omega_{K''}$ is surjective. 
\end{prop}
\begin{proof}
Let $x$ be an element in $\Omega_{K''}$. The absolute value $|\ndot|_\omega$ extends in a unique way to the algebraic closure $(K_{x}'')^{\mathrm{ac}}$ of $K_{x}''$ (see \cite{Neukirch} Chapter II, Theorem 4.8). Therefore, if we choose an embedding of $L$ into $(K_{x}'')^{\mathrm{ac}}$, then we obtain an absolute value on $L$ which extends $|\ndot|_x$. 
\end{proof}

Similarly to \cite[Lemma 2.1]{Gaudron_Remond14}, the set $\Omega_L$ described above gives an explicit construction of the projective limit of the projective system $\{\Omega_{K''}\}_{K''\in\mathscr E_{L/K}}$ in the category of sets, where $\{\pi_{L/K''}\}_{K''\in\mathscr E_{L/K}}$ are universal maps. In fact, any absolute value on $L$ is uniquely determined by its restrictions on the subfields in $\mathscr E_{L/K}$. We equip $\Omega_L$ with the $\sigma$-algebra $\mathcal A_L$ generated by the maps $\pi_{L/K''}$ (namely the smallest $\sigma$-algebra which makes all maps $\pi_{L/K''}$ measurable) where $K''$ runs over $\mathscr E_{L/K}$. Thus $(\Omega_L,\mathcal A_L)$ identifies with the projective limit of the projective system $\big\{(\Omega_{K''},\mathcal A_{K''})\big\}_{K''\in\mathscr E_{L/K}}$ in the category of mesurable spaces.

Let $\omega$ be an element in $\Omega_K$. We equip $M_{L,\omega}$ with the smallest $\sigma$-algebra $\mathcal A_{L,\omega}$ such that the restriction of $\pi_{L/K''}$ on $M_{L,\omega}$ is measurable for any $K''\in\mathscr E_{L/K}$, where we consider the discrete $\sigma$-algebra on $M_{K'',\omega}=\pi_{K''/K}^{-1}(\{\omega\})$. Let $V_{L,\omega}$ be the set of all real-valued functions on $\Omega_L$ which can be written in the form $f\circ(\pi_{L/K''}|_{M_{L,\omega}})$, where $K''$ is an element of $\mathscr E_{L/K}$, and $f$ is a function on $M_{K'',\omega}$. Let \[I_{L/K,\omega}:V_{L,\omega}\longrightarrow\mathbb R\]
be the map which sends any function of the form $f\circ(\pi_{L/K''}|_{M_{L,\omega}})$ to the integral
\begin{equation}\label{Equ:fibreint}\int_{M_{K'',\omega}}f\,\mathrm{d}\mathbb P_{K',\omega},\end{equation}
where $K'$ is the separable closure of $K$ in $K''$, and $\mathbb P_{K',\omega}$ is the probability measure on $M_{K',\omega}$ defined in \eqref{Equ:Pomega}. Similarly to \eqref{Equ:compatibility}, the fibre integral is compatible with successive finite extensions of the field $K$ and the map $I_{L/K,\omega}$ is well defined since the value of the integral \eqref{Equ:fibreint} does not depend on the choice of the field $K''$ upon which we write the function in $V_{L,\omega}$ as the composition of a function on $M_{K'',\omega}$ with $\pi_{L/K''}|_{M_{L,\omega}}$.

\begin{prop}
The set $V_{L,\omega}$ forms an algebra over $\mathbb R$ with respect to the composition laws of addition and multiplication of functions, and the map $I_{L/K,\omega}:V_{L,\omega}\rightarrow\mathbb R$ is an $\mathbb R$-linear operator. Moreover, it induces a probability measure on the measurable space $(M_{L,\omega},\mathcal A_{L,\omega})$.
\end{prop}
\begin{proof}
The first assertion comes from the fact that the set $\mathscr E_{L/K}$ is filtered, which implies that any finite collection of functions in $V_{L,\omega}$ descend on the same space $M_{K'',\omega}$, where $K''\in\mathscr E_{L/K}$. In particular, the family $\mathcal D$ of subsets $A\subseteq M_{L,\omega}$ such that $\indic_A\in V_{L,\omega}$ is an algebra (of sets), which generates $\mathcal A_{L,\omega}$ as a $\sigma$-algebra. Moreover, the map $\mathbb P_{L,\omega}:\mathcal D\rightarrow\mathbb R_+$ which sends $A\in\mathcal D$ to $I_{L/K,\omega}(\indic_A)$ is an additive functional. Clearly it sends $M_{L,\omega}$ to $1$.

The $\sigma$-algebra $\mathcal A_{L,\omega}$ is actually the Borel algebra of the projective limit topology on $M_{L,\omega}$ (namely the most coarse topology on $M_{L,\omega}$ which makes all maps $\pi_{L/K''}$ continuous, where $K''\in\mathscr E_{L/K}$). This topology also identifies with the induced topology on $M_{L,\omega}$ viewed as a subset of $\prod_{F\in\mathscr E_{L/K}}M_{F,\omega}$ (equipped with the product topology), where on each set $M_{F,\omega}$ we consider the discrete topology. Note that $M_{L,\omega}$ is actually a closed subset of this product space since it is the intersection of closed subsets of the form
\[W_{K''}:=\bigg\{(x_{F})_{F\in\mathscr E_{L/K}}\in\prod_{F\in\mathscr E_{L/K}}M_{F,\omega}\,:\,\pi_{K''/F}(x_{K''})=x_F\text{ for }F\subseteq K''\bigg\}.\]
Therefore, by Tychonoff's theorem, we obtain that $M_{L,\omega}$ is actually a compact topological space. Moreover, any set in $\mathcal D$ is open and closed since it is the inverse image of a discrete set by a continuous map. Therefore, the sets in $\mathcal D$ are open and compact. As a consequence, if $\{A_n\}_{n\in\mathbb N}$ is a sequence of disjoint sets in $\mathcal D$ whose union also lies in $\mathcal D$, then for sufficiently large $n$ one has $A_n=\varnothing$. Hence the function $\mathbb P_{L,\omega}:\mathcal D\rightarrow\mathbb R_+$ is actually $\sigma$-additive. By Carath\'eodory's extension theorem, the function $\mathbb P_{L,\omega}$ extends to a Borel probability measure on $(M_{L,\omega},\mathcal A_{L,\omega})$ such that 
\[I_{L/K,\omega}(f)=\int_{M_{L,\omega}}f\,\mathrm{d}\mathbb P_{L,\omega}.\]
The proposition is thus proved.
\end{proof}

\begin{rema}\label{Rem:operateurILK}
Let $V_L$ be the vector space of all real-valued functions $f$ on $\Omega_L$ which can be written as $g\circ \pi_{L/K''}$, where $K''/K$ is a finite extension which is contained in $L$. Then the above construction leads to a linear operator $I_{L/K}$ from $V_L$ to the vector space of all real-valued functions on $\Omega$, sending $f\in V_L$ to the function \[(\omega\in\Omega)\longmapsto I_{L/K,\omega}(f|_{M_{L,\omega}}).\]
Clearly, if $g$ is real-valued function on $\Omega$, then $I_{L/K}(g\circ\pi_{L/K})=g$.
\end{rema}

The above proposition allows to define the fibre integration for  non-negative $\mathcal A_L$-measurable functions on $\Omega_{L}$.

\begin{prop}\label{Pro:fibreint}
Let $f$ be a non-negative $\mathcal A_L$-measurable function on $\Omega_{L}$. For any $\omega\in\Omega_K$, the restriction of $f$ {to} $M_{L,\omega}$ is $\mathcal A_{L,\omega}$-measurable. Moreover,  the map $I_{L/K}(f)$ from $\Omega_{K}$ to $[0,+\infty]$ which sends $\omega\in\Omega_K$ to $\int_{M_{L,\omega}}f(x)\,\mathbb P_{L,\omega}(\mathrm{d}x)$ is $\mathcal A_K$-measurable.
\end{prop}
\begin{proof}
Let $\mathcal H$ be the set of all bounded non-negative functions $f$ on $\Omega_K$ such that $f|_{M_{L,\omega}}$ is $\mathcal A_{L,\omega}$-measurable for any $\omega\in\Omega_K$ and that the map
\[(\omega\in\Omega_K)\longmapsto\int_{M_{L,\omega}}f(x)\,\mathbb P_{L,\omega}(\mathrm{d}x)\]
is $\mathcal A_K$-measurable. Then $\mathcal H$ is a $\lambda$-family of non-negative functions on $\Omega_{L}$ (see Definition \ref{Def:lambda-family}). Moreover, the set $\mathcal H$ contains the subset $\mathcal C$ of all bounded non-negative functions of the form $g\circ\pi_{L/K''}$, where $K''$ is an element in $\mathscr E_{L/K}$ and $g$ is an $\mathcal A_{K''}$-measurable function on $\Omega_{K''}$. In fact, for any $\omega\in\Omega_K$, one has
\[\int_{M_{L,\omega}}g(\pi_{L/K''}(x))\,\mathbb P_{L,\omega}(\mathrm{d}x)=I_{K''/K}(g)(\omega).\]
Since $\mathcal C$ is stable under multiplication, by the monotone class theorem \ref{Thm:monotoneclass} we obtain that the family $\mathcal H$ actually contains all non-negative, bounded and $\sigma(\mathcal C)$-measurable functions. By definition, $\mathcal A_L$ is the $\sigma$-algebra generated by the maps $\pi_{L/K''}$ with $K''\in\mathscr E_{L/K}$. Therefore one has $\mathcal A_L=\sigma(\mathcal C)$. Thus we obtain the result of the proposition for bounded non-negative $\mathcal A_L$-measurable functions. For general non-negative $\mathcal A_L$-measurable function $f$, we can apply the assertion of the proposition to the functions $\{\min(f,n)\}_{n\in\mathbb N}$ which form an increasing sequence converging to $f$. Passing to limit when $n$ goes to the infinity, we obtain the result for $f$.
\end{proof}

The above proposition allows to construct a measure $\nu_L$ on the measurable space $(\Omega_L,\mathcal A_L)$ such that, for any subset $A$ of $\mathcal A_L$, one has
\[\nu_L(A)=\int_{\Omega_K}\Big(\int_{M_{L,\omega}}\indic_A(x)\,\mathbb P_{L,\omega}(\mathrm{d}x)\Big)\nu_K(\mathrm{d}\omega).\] 
For any non-negative $\mathcal A_L$-measurable function $f$ on $\Omega_L$, one has
\begin{equation}\label{Equ:fubini}\int_{\Omega_L}f(x)\,\nu_L(\mathrm{d}x)=\int_{\Omega_K}\Big(\int_{M_{L,\omega}}f(x)\,\mathbb P_{L,\omega}(\mathrm{d}x)\Big)\nu_K(\mathrm{d}\omega).\end{equation}
We denote by $I_{L/K}(f)$ the map from $\Omega_K$ to $[0,+\infty]$ which sends $\omega\in\Omega_K$ to $\int_{M_{L,\omega}}f(x)\,\mathbb P_{L,\omega}(\mathrm{d}x)$. More generally, for any $\mathcal A_L$-measurable function $f$ such that $I_{L/K}(|f|)$ is a real-valued function, we define $I_{L/K}(f)$ as the real-valued function
\[I_{L/K}(\max(f,0))-I_{L/K}(-\min(f,0)).\]
Note that, if $f$ is of the form $g\circ\pi_{L/K''}$ where $K''\in\mathscr E_{L/K}$ and $g$ is an $\mathcal A_{K''}$-measurable function, then $I_{L/K}(f)$ is always well defined, and one has
\begin{equation}\label{Equ:ILoverK}
I_{L/K}(g\circ \pi_{L/K''})=I_{K''/K}(g).
\end{equation}
With this notation, the equality \eqref{Equ:fubini} can also be written as
\begin{equation}\label{Equ:fubinibis}\int_{\Omega_L}f\,\mathrm{d}\nu_L=\int_{\Omega_K}I_{L/K}(f)\,\mathrm{d}\nu_K.\end{equation}
Thus we obtain the following result.

\begin{prop}\label{Pro:fibreintegralint}
An $\mathcal A_L$-measurable function $f$ is $\nu_L$-integrable if and only if $I_{L/K}(|f|)$ is $\nu_K$-integrable. Moreover, $I_{L/K}$ defines a continuous linear operator from $\mathscr L^1(\Omega_L,\mathcal A_L,\nu_L)$ to $\mathscr L^1(\Omega_K,\mathcal A_K,\nu_K)$, and  the equality \eqref{Equ:fubinibis} also holds for $\nu_L$-integrable functions. 
\end{prop}

Note that the relations \eqref{Equ:ILoverK} and \eqref{Equ:fubinibis} also imply that, if $g$ is a function in $\mathscr L^1(\Omega_K,\mathcal A_K,\nu_K)$, then $g\circ\pi_{L/K}$ belongs to $\mathscr L^1(\Omega_L,\mathcal A_L,\nu_L)$, and one has
\begin{equation}
\label{Equ:integrationofcomposition}
\int_{\Omega_L}(g\circ\pi_{L/K})\,\mathrm{d}\nu_L=\int_{\Omega_K}g\,\mathrm{d}\nu_K.
\end{equation}

The following proposition shows that $(L,(\Omega_L,\mathcal A_L,\nu_L),\phi_L)$ forms an adelic curve.

\begin{prop}\label{Pro:measurablitliAL}
For any non-zero element $a\in L$, the function \[(z\in\Omega_L)\longmapsto \ln |a|_z\]
is $\mathcal A_L$-measurable. Moreover, if $S = (K, (\Omega_K, \mathcal A_K,\nu_K), \phi_K)$ is proper, then $(L, (\Omega_L, \mathcal A_L,\nu_L), \phi_L)$
is also proper.
\end{prop}
\begin{proof}
Denote by $g$ the function on $\Omega_L$ such that $g(z)=\ln|a|_z$. We choose a finite extension $K''\in\mathscr E_{L/K}$ which contains $a$. Let $f:\Omega_{K''}\rightarrow\mathbb R$ be the function which sends $x\in\Omega_{K''}$ to $\ln|a|_x$. Then $f$ is an $\mathcal A_{K''}$-measurable function on $\Omega_{K''}$. Since the function $g$ identifies with the composition $f\circ\pi_{L/K''}$, we obtain that $g$ is $\mathcal A_L$-measurable. 
\if01
Moreover, for any $\omega\in\Omega_K$, one has
\[I_{L/K}(g)(\omega)=I_{L/K,\omega}(g)=I_{K''/K}(f).\]
Therefore we obtain 
\[\int_{\Omega_L}g\,\mathrm{d}\nu_L=\int_{\Omega_K}I_{K''/K}(f)\,\mathrm{d}\nu_K=0,\]
where the second equality comes from \eqref{Equ:convergence} and the product formula for the adelic curve $(K,(\Omega_K,\mathcal A_K,\nu_K),\phi_K)$. The proposition is thus proved.
\fi

We assume that $S$ is proper.
For any $\omega\in\Omega_K$, one has
\[I_{L/K}(g)(\omega)=I_{L/K,\omega}(g)=I_{K''/K}(f)(\omega).\]
Therefore we obtain 
\[\int_{\Omega_L}g\,\mathrm{d}\nu_L=\int_{\Omega_K}I_{K''/K}(f)\,\mathrm{d}\nu_K=0,\]
where the second equality comes from Theorem~\ref{Thm:constructionofextension},
\ref{Item: nu K' is a measure} and \ref{Item: productformula:norm}.
\end{proof}

\begin{defi}
Let $S=(K,(\Omega_K,\mathcal A_K,\nu_K),\phi_K)$ be an adelic curve and $L/K$ be an algebraic extension. The adelic curve $(L,(\Omega_L,\mathcal A_L,\nu_L),\phi_L)$ is called an \emph{algebraic extension}\index{algebraic extension@algebraic extension} of $S$, denoted by $S\otimes_KL$.
\end{defi}

The following result, which is similar to Proposition \ref{Pro:successiveextension}, shows the compatibility property of the algebraic extensions of adelic curves.

\begin{theo}\label{Thm:successivealgextension}
Let $S=(K,(\Omega_K,\mathcal A_K,\nu_K),\phi_K)$ be an adelic curve and $L_2/L_1/K$ be successive algebraic extensions of fields. Then one has \begin{equation}\label{Equ:extensiongeneral}(S\otimes_KL_1)\otimes_{L_1}L_2=S\otimes_KL_2.\end{equation}
Moreover, the following relations hold
\begin{equation}\label{Equ:compatibilityL}\pi_{L_1/K}\circ\pi_{L_2/L_1}=\pi_{L_2/K},\quad I_{L_1/K}\circ I_{L_2/L_1}=I_{L_2/K}.\end{equation}
\end{theo}
\begin{proof}
Let $(L_1,(\Omega_1,\mathcal A_1,\nu_1),\phi_1)$ and $(L_2,(\Omega_2,\mathcal A_2,\nu_2),\phi_2)$ be the adelic curves $S\otimes_KL_1$ and $S\otimes_KL_2$ respectively. First of all, set-theoretically one has
\[\Omega_2=\coprod_{\omega\in\Omega_K}M_{L_2,\omega}=\coprod_{\omega\in\Omega_K}\coprod_{x\in M_{L_1,\omega}}M_{L_2,x}=\coprod_{x\in\Omega_1}M_{L_2,x},\] and hence  
\[\pi_{L_2/K}=\pi_{L_1/K}\circ\pi_{L_2/L_1}.\]
Moreover, for any extension $K_1\in\mathscr E_{L_1/K}$ the map $\pi_{L_2/K_1}$ is $\mathcal A_2$-measurable. Therefore $\pi_{L_2/L_1}:\Omega_2\rightarrow\Omega_1$ is an $\mathcal A_2$-measurable map since the $\sigma$-algebra $\mathcal A_1$ is generated by the maps $\pi_{L_1/K_1}$ with $K_1\in\mathscr E_{L_1/K}$.

We now proceed with the proof of the equalities \eqref{Equ:extensiongeneral} and \eqref{Equ:compatibilityL}  with the supplementary assumption that the extension $L_2/L_1$ is \emph{finite}. We first show that the $\sigma$-algebra $\mathcal A_2$ coincides with that in the adelic structure of $(S\otimes_KL_1)\otimes L_2$, namely the smallest $\sigma$-algebra $\mathcal A_2'$ on $\Omega_2$ such that $\pi_{L_2/L_1}$ and all functions of the form $(y\in\Omega_2)\mapsto |a|_y$ are $\mathcal A_2'$-measurable, where $a\in L_2$. We have already shown that the map $\pi_{L_2/L_1}$ is measurable. Hence by proposition \ref{Pro:measurablitliAL}, we obtain that $\mathcal A_2'\subseteq\mathcal A_2$. Conversely, for any extension $K_1\in\mathscr E_{L_1/K}$ one has 
\[\pi_{L_2/K_1}=\pi_{L_1/K_1}\circ\pi_{L_2/L_1},\]
and hence
\[\pi_{L_2/K_1}^{-1}(\mathcal A_{K_1})=\pi_{L_2/L_1}^{-1}(\pi_{L_1/K_1}^{-1}(\mathcal A_{K_1}))\subseteq \pi_{L_2/L_1}^{-1}(\mathcal A_1)\subseteq\mathcal A_2'. \]
If $K_2$ is an extension in $\mathscr E_{L_2/K}$, then $K_1:=K_2\cap L_1\in\mathscr E_{L_1/K}$. Moreover, the $\sigma$-algebra $\mathcal A_{K_2}$ is generated by $\pi_{K_2/K_1}$ and the functions of the form $x\mapsto|a|_x$ on $\Omega_{K_2}$, where $a\in K_2$. Note that $\pi_{L_2/K_1}=\pi_{L_1/K_1}\circ\pi_{L_2/L_1}$ is $\mathcal A_2'$-measurable, and for any $a\in K_2$, the composition of the function $x\mapsto|a|_x$ on $\Omega_{K_2}$ with $\pi_{L_2/K_2}$, which identifies with the function $y\mapsto |a|_y$ on $\Omega_{L_2}$, is also $\mathcal A_2'$-measurable, we obtain that the map $\pi_{L_2/K_2}$ is actually $\mathcal A_2'$-measurable. Since $\mathcal A_2$ is the smallest $\sigma$-algebra which makes all $\pi_{L_2/K_2}$ measurable, where $K_2\in\mathscr E_{L_2/K}$, we obtain $\mathcal A_2\subseteq\mathcal A_2'$. Therefore one has $\mathcal A_2=\mathcal A_2'$. 

It remains to establish the relation $I_{L_1/K}(I_{L_2/L_1}(f))=I_{L_2/K}(f)$ for any non-negative $\mathcal A_{2}$-measurable function on $\Omega_2$. By induction it suffices to treat the case where $[L_2:L_1]$ is a prime number. Moreover, similarly to the proof of Proposition \ref{Pro:fibreint}, by using the monotone class theorem, we only need to verify the equality $I_{L_1/K}(I_{L_2/L_1}(f))=I_{L_2/K}(f)$ for functions $f$ of the form $g\circ\pi_{L_2/K_2}$, where $K_2\in\mathscr E_{L_2/K}$ and $g$ is a non-negative $\mathcal A_{K_2}$-measurable function on $\Omega_{K_2}$. If $K_2$ belongs to $\mathscr E_{L_1/K}$, one has $I_{L_2/L_1}(f)=g\circ\pi_{L_1/K_2}$, and therefore
\[I_{L_1/K}(I_{L_2/L_1}(f))=I_{K_2/K}(g)=I_{L_2/K}(f),\]
where the second equality comes from \eqref{Equ:ILoverK}.
Otherwise one has $[K_2:K_1]=[L_2:L_1]$ since $[K_2:K_1]$ divides $[L_2:L_1]$ which is a prime number, where $K_1$ denotes the intersection of $K_2$ with $L_1$. Moreover, there exists an element $a\in K_2$ such that $K_2=K_1(a)$ and $L_2=L_1(a)$. If $a$ is totally inseparable over $K_1$, then it is also totally inseparable over $L_1$. In this case $I_{L_2/L_1}(f)=f=g\circ\pi_{L_1/K_1}$ and therefore the equality $I_{L_1/K}(I_{L_2/L_1}(f))=I_{L_2/K}(f)$ also holds in this case.

In the following, we assume that the element $a$ is separable over $K_1$. Let \[P(T)=T^p+b_1T^{p-1}+\cdots+b_p\in K_1[T]\]
be the minimal polynomial of $a$ over $K_1$. Since $p$ is a prime number, and $a\not\in L_1$, we obtain that it is also the minimal polynomial of $a$ over $L_1$. In particular, the element $a$ is also separable over $L_1$. Let $y$ be an element in $\Omega_1$ and $x=\pi_{L_1/K_1}(y)$. Assume that
\[P=P_1\cdots P_r\]
is the splitting of the polynomial $P$ in the ring $K_{1,x}[T]$ into the product of distinct irreducible polynomials. Then the polynomials $P_1,\ldots,P_r$ correspond to points $x_1,\ldots,x_r$ which form the set $\pi_{K_2/K_1}^{-1}(\{x\})$. Moreover, one has (see \eqref{Equ:Pomega} for the definition of $\mathbb P_{K_2,x}$)
\[\mathbb P_{K_2,x}(\{x_i\})=\frac{\deg(P_i)}{p}.\]
Assume that each $P_i$ splits in $L_{1,y}$ into the product of distinct irreducible polynomials as 
\[P_i=Q_{i,1}\cdots Q_{i,n_i}.\]
Then each factor $Q_{i,j}$ corresponds to a point $y_{i,j}$ in $\pi_{L_2/L_1}^{-1}(\{y\})$ and one has
\[\mathbb P_{L_2,y}(x_{i,j})=\frac{\deg(Q_{i,j})}{p}.\]
Therefore, if a non-negative function $f$ on $\Omega_2$ is of the form $g\circ\pi_{L_2/K_2}$, where $g$ is a $\mathcal A_{K_2}$-measurable function on $\Omega_{K_2}$, then one has
\[\big(I_{L_2/L_1}(f)\big)(y)=\sum_{i=1}^r\sum_{j=1}^{n_i}\frac{\deg(Q_{i,j})}{p}g(x)=\frac{\deg(P_i)}{p}g(x)=\big(I_{K_2/K_1}(g)\big)(x),\]
which shows that 
\[I_{L_2/L_1}(f)=I_{K_2/K_1}(g)\circ\pi_{L_1/K_1}.\]
Therefore one has
\[\begin{split}&\quad\;I_{L_1/K}(I_{L_2/L_1}(f))=I_{L_1/K}(I_{K_2/K_1}(g)\circ\pi_{L_1/K_1})\\
&=I_{K_1/K}(I_{K_2/K_1}(g))=I_{K_2/K}(g)=I_{L_2/K}(f),
\end{split}\]
where the second and the last equalities come from \eqref{Equ:ILoverK}. Thus we have established the second equality in \eqref{Equ:compatibilityL}, which implies that the measure in the adelic structure of $(S\otimes_KL_1)\otimes_{L_1}L_2$ coincides with the measure $\nu_2$ in the adelic structure of $S\otimes_KL_2$. The theorem is then established in the particular case where $[L_2:L_1]$ is finite.

In the following, we will prove the general case of the theorem. Note that the previously proved case actually implies that, for any finite extension $L''$ of $L_1$, the map $\pi_{L_2/L''}$ is $\mathcal A_2$-measurable since the $\sigma$-algebra $\mathcal A_{L''}$ in the adelic structure of $(S\otimes_KL_1)\otimes L''$ coincides with that in the adelic structure of $S\otimes_KL''$, which is generated by the maps $\pi_{L''/K''}$ with $K''\in\mathscr E_{L''/K}$. In particular, if we denote by $\mathcal A_2'$ the $\sigma$-algebra in the adelic structure of $(S\otimes_KL_1)\otimes_{L_1}L_2$, then one has $\mathcal A_2'\subseteq\mathcal A_2$. Conversely, for any $K_2\in\mathscr E_{L_2/K}$, one has $\pi_{L_2/K_2}=\pi_{L''/K_2}\circ\pi_{L_2/L''}$, where $L''=L_1K_2$ is an element in $\mathscr E_{L_2/L_1}$. Hence $\pi_{L_2/K_2}$ is $\mathcal A_2'$-measurable. Since $K_2\in\mathscr E_{L_2/L}$ is arbitrary, we obtain that $\mathcal A_2\subseteq\mathcal A_2'$ and hence $\mathcal A_2=\mathcal A_2'$.

Again it remains to establish the equality $I_{L_1/K}(I_{L_2/L_1}(f))=I_{L_2/K}(f)$ for any non-negative $\mathcal A_{2}$-measurable function $f$ on $\Omega_2$, which can be written in the form $g\circ\pi_{L_2/K_2}$, where $K_2\in\mathscr E_{L_2/K}$. Let $L''=K_2L_1$. One has $L''\in\mathscr E_{L_2/L_1}$ and $g\circ\pi_{L_2/K_2}=g\circ\pi_{L''/K_2}\circ\pi_{L_2/L''}$. Therefore \eqref{Equ:ILoverK} implies
\[I_{L_2/K}(g\circ\pi_{L_2/K_2})=I_{L''/K}(g\circ\pi_{L''/K_2})=I_{K_2/K}(g).\]
Moreover, also by \eqref{Equ:ILoverK} one obtains
\[\begin{split}&\quad\;I_{L_1/K}(I_{L_2/L_1}(g\circ\pi_{L_2/K_2}))=I_{L_1/K}(I_{L_2/L_1}(g\circ\pi_{L''/K_2}\circ\pi_{L_2/L''})))\\
&=I_{L_1/K}(I_{L''/L_1}(g\circ\pi_{L''/K_2}))=I_{L''/K}(g\circ\pi_{L''/K_2})=I_{K_2/K}(g),\end{split}\]
where the third equality comes from the proved case of finite extensions. Thus we establish the relation $I_{L_1/K}\circ I_{L_2/L_1}=I_{L_2/K}$ which implies that the measure $\nu_2$ identifies with that in the adelic structure of $(S\otimes_KL_1)\otimes_{L_1}L_2$. The theorem is thus proved.
\end{proof}

\begin{prop}\label{Pro:sigma algebra on extension}
Let $S=(K,(\Omega_K,\mathcal A_K,\nu_K),\phi_K)$ be an adelic curve and $L/K$ be an algebraic extension. Let $(L,(\Omega_L,\mathcal A_L,\nu_L),\phi_L)$ be the adelic curve $S\otimes_KL$. Then $\mathcal A_L$ is the smallest $\sigma$-algebra making the canonical projection map $\pi_{L/K}:\Omega_L\rightarrow\Omega_K$ and the functions $(x\in\Omega_L)\mapsto |a|_x$ measurable for all $a \in L$.
\end{prop}
\begin{proof}
By definition the projection map $\pi_{L/K}$ is $\mathcal A_L$-measurable. Moreover, by Proposition \ref{Pro:measurablitliAL}, for any $a\in L$, the function $(x\in\Omega_L)\mapsto |a|_x$ is $\mathcal A_L$-measurable. 

Suppose that $F$ is a  map from a measurable space $(E,\mathcal E)$ to $\Omega_L$ such that the composed map $\pi_{L/K}\circ F$ and the functions $(y\in E)\mapsto |a|_{F(y)}$ are measurable, where $a\in L$. We will show that $F$ is measurable if we consider the $\sigma$-algebra $\mathcal A_L$ on $\Omega_L$. This implies that $\mathcal A_L$ is contained in the smallest $\sigma$-algebra making the canonical projection $\Omega_L\rightarrow\Omega_K$ and the functions $(x\in\Omega_L)\mapsto |a|_x$ measurable, where $a\in L$.

Recall that $\mathcal A_L$ is the smallest $\sigma$-algebra making the projection maps $\pi_{L/K''}:\Omega_L\rightarrow\Omega_{K''}$ measurable, where $K''/K$ runs over the set of finite extensions contained in $L$. To show the measurability of $F$ it suffices to verify the measurability of $\pi_{L/K''}\circ F$ for any finite extension $K''/K$ contained in $L$. Moreover, since $\mathcal A_{K''}$ is the smallest $\sigma$-algebra making the projection map $\pi_{K''/K}:\Omega_{K''}\rightarrow\Omega_K$ and the functions $(x\in\Omega_{K''})\mapsto|a|_x$ measurable, where $a\in K''$, we are reduced to verify the measurability of $\pi_{K''/K}\circ\pi_{L/K''}\circ F=\pi_{L/K}\circ F$ and 
\begin{equation}\label{Equ: function composed y}(y\in E)\longmapsto |a|_{\pi_{L/K''}(F(y))},\quad\text{where $a\in K''$.}\end{equation}
By the assumption on $F$, the map $\pi_{L/K}\circ F$ is measurable. Moreover, since $a\in K''$, one has
\[|a|_{\pi_{L/K''}(F(y))}=|a|_{F(y)}.\]
Hence the function in  \eqref{Equ: function composed y} is also measurable. The proposition is thus proved.
\end{proof}

\section{Height function and Northcott property}
Let $S = (K, (\Omega, \mathcal{A}, \nu), \phi)$ be a \emph{proper} 
adelic curve and
$K^{\mathrm{ac}}$ be an algebraic closure of $K$.
Let $S \otimes_K K^{\mathrm{ac}} = (K^{\mathrm{ac}}, (\Omega_{K^{\mathrm{ac}}}, \mathcal{A}_{K^{\mathrm{ac}}}, \nu_{K^{\mathrm{ac}}}), \phi_{K^{\mathrm{ac}}})$ be the algebraic extension of $S$ by $K^{\mathrm{ac}}$.

\begin{defi}\label{def:height:adelic:curve}
For $(a_0, a_1, \ldots, a_n) \in (K^{\mathrm{ac}})^{n+1} \setminus \{ (0,\ldots, 0) \}$, we define $h_S(a_0, \ldots, a_n)$ to be
\[
h_S(a_0, \ldots, a_n) := \int_{\Omega_{K^{\mathrm{ac}}}} \ln \left( \max \{ |a_0|_{\chi}, \ldots, |a_n|_{\chi} \} \right) \nu_{K^{\mathrm{ac}}}(d\chi).
\]
By the product formula, $h_S(\lambda a_0, \ldots, \lambda a_n) = h_S(a_0, \ldots, a_n)$ 
for all $\lambda \in K^{\mathrm{ac}} \setminus \{ 0 \}$, so that
there is a map $\mathbb P^n(K^{\mathrm{ac}}) \to \mathbb R$ such that
the following diagram is commutative:
\[
\xymatrix{ (K^{\mathrm{ac}})^{n+1} \setminus \{ (0,\ldots, 0) \} \ar[d] \ar[r]^-{h_S} & \mathbb R \\
 \mathbb P^n(K^{\mathrm{ac}}) \ar[ur] }
\]
By abuse of notation, the map $\mathbb P^n(K^{\mathrm{ac}}) \to \mathbb R$ is also denoted by $h_{S}$. For $x \in \mathbb P^n(K^{\mathrm{ac}})$, the value $h_S(x)$ is called the \emph{height}\index{height@height} of $x$
with respect to the adelic curve $S$.
\end{defi}

\begin{defi}\label{def:northcott:property}
We say that $S$ has the Northcott property if the set
$\left\{ a \in K \,:\, h_S(1 : a) \leqslant C \right\}$
is finite for any $C \geqslant 0$.
In the cases of Example~\ref{Subsec:Numberfields}, Example~\ref{subsec:fun:field:Q} and Example~\ref{subsec:polarized;arith:var},
the Northcott property holds (for details, see \cite{Moriwaki00} and \cite{Moriwaki04}).
\end{defi}

The purpose of this section is to prove the following theorem:

\begin{theo}[Northcott's theorem]\label{thm:Northcott:thm}
If $S$ has the Northcott property, then
the set $\{ x \in \mathbb P^n(K^{\mathrm{ac}}) \,:\, h_S(x) \leqslant C,\ [K(x) : K] \leqslant \delta \}$ is finite
for any $C$ and $\delta$.
\end{theo}

Before starting the proof of Theorem~\ref{thm:Northcott:thm},
we need to prepare two lemmas.

\begin{lemm}\label{lem:invariant:auto:integral}
Let $K'$ be a finite normal extension of $K$.
Then
$h_S(1 : \sigma(\alpha)) = h_S(1 : \alpha)$
for all $\alpha \in K'$ and $\sigma \in \operatorname{Aut}_K(K')$. 
\end{lemm}

\begin{proof}
Let $K''$ be the separable closure of $K$ in $K'$ and $q = [K' : K'']$.
Then $\rest{\sigma}{K''} \in \operatorname{Aut}_K(K'')$ and $a^q \in K''$.
If the assertion holds for the extension $K''/K$, then
\[
q h_S(1 : \sigma(\alpha)) = h_S(1 : \sigma(\alpha^q)) = h_S(1 : \alpha^q) = q h_S(1 : \alpha),
\]
so that we may assume that the extension $K'/K$ is separable.

For $\chi \in \pi^{-1}_{K'/K}(\{ \omega \})$ and $\tau \in \mathrm{Gal}(K'/K)$, 
a map $(\beta \in K') \mapsto |\tau(\beta)|_{\chi}$ gives rise to
an element of $\pi^{-1}_{K'/K}(\{ \omega \})$, which is denoted by ${\chi}^{\tau}$. In this way,
one has an action $\mathrm{Gal}(K'/K) \times \pi^{-1}_{K'/K}(\{ \omega \}) \to \pi^{-1}_{K'/K}(\{ \omega \})$
given by $(\tau, \chi) \mapsto {\chi}^{\tau}$. Note that the action is transitive (cf. \cite[Chapter~II, Proposition~9.1]{Neukirch}) 
and $\mathrm{Gal}(K'_{\chi}/K_{\omega}) = \mathrm{Stab}_{\mathrm{Gal}(K'/K)}(\chi)$
(cf. \cite[Chapter~II, Proposition~9.6]{Neukirch}). In particular,
$[K'_{\chi}: K_{\omega}] = [K'_{\chi'}: K_{\omega}]$ for all $\chi, \chi' \in \pi^{-1}_{K'/K}(\{ \omega \})$.
Therefore,
{\allowdisplaybreaks
\begin{align*}
I_{K'/K}(\ln(\max \{1, |\sigma(\alpha)| \}))(\omega) & = \sum_{\chi \in \pi^{-1}_{K'/K}(\{ \omega \})}
\frac{[K'_{\chi}: K_{\omega}]}{[K' : K]} \ln(\max \{1,  |\sigma(\alpha)|_{\chi}\}) \\
& = \sum_{\chi \in \pi^{-1}_{K'/K}(\{ \omega \})}
\frac{[K'_{{\chi}^{\sigma}}: K_{\omega}]}{[K' : K]} \ln(\max \{1,  |\alpha|_{{\chi}^{\sigma}} \}) \\
& = I_{K'/K}(\ln(\max \{1, |\alpha| \}))(\omega),
\end{align*}}
and hence the assertion follows.
\end{proof}

For a polynomial $F = a_n X^n + \cdots + a_1 X + a_0 \in K^{\mathrm{ac}}[X] \setminus \{ 0 \}$, we define $h_S(F)$ to be
$h_S(F) := h_S(a_n : \cdots : a_1 : a_0)$.
If we set $\| F \|_{\chi} := \max \{ |a_n|_{\chi}, \ldots, |a_0|_{\chi} \}$ for $\chi \in \Omega_{K^{\mathrm{ac}}}$ as in Subsection~\ref{Subsec:Norm:polynomial}, then
\[
h_S(F) = \int_{\Omega_{K^{\mathrm{ac}}}} \ln (\|F\|_{\chi})\, \nu_{K^{\mathrm{ac}}}(d\chi).
\]

\begin{lemm}\label{lem:height:product:polynomial}
For $F, G \in K^{\mathrm{ac}}[X] \setminus \{ 0 \}$, one has
\[
h_S(FG) \leqslant h_S(F) + h_S(G) + \ln \min \{\deg(F)+1,\deg(G)+1 \} \int_{\Omega_{\infty}} \nu(d\omega).
\]
\end{lemm}

\begin{proof}
By Proposition~\ref{prop:norm:polynomial},
{\allowdisplaybreaks
\begin{align*}
h_S(FG) & = \int_{\Omega_{K^{\mathrm{ac}}}} \ln (\|FG\|_{\chi}) \nu_{K^{\mathrm{ac}}}(d\chi) \\
& \leqslant \int_{\Omega_{K^{\mathrm{ac}}}} \left( \ln (\|F\|_{\chi}) +\ln(\|G\|_{\chi}) \right) \nu_{K^{\mathrm{ac}}}(d\chi) \\
& \hskip5em +
\int_{\Omega_{\infty}} \ln (\min \{\deg(F) + 1, \deg(G) + 1 \})\nu(d\omega),
\end{align*}}
so that the assertion follows.
\end{proof}

\begin{proof}[Proof of Theorem~\ref{thm:Northcott:thm}]
Clearly we may assume that $C \geqslant 0$ and $\delta \geqslant 1$.
Let us begin with the following special case:

\begin{enonce}{Claim}
The set
$\left\{ \alpha \in K^{\mathrm{ac}} \,:\,
h_S(1 : \alpha) \leqslant C,\
[K(\alpha) : K] \leqslant \delta \right\}$
is finite.
\end{enonce}

\begin{proof}
Let $F$ be the minimal monic polynomial of $\alpha$ over $K$. We set
$F = X^n + a_{n-1} X^{n-1} + \cdots + a_1 X + a_0 = (X -\alpha_1) \cdots (X - \alpha_n)$
and $K' = K(\alpha_1, \ldots, \alpha_n)$, where $\alpha_1 = \alpha$.
Then, by Lemma~\ref{lem:invariant:auto:integral} and Lemma~\ref{lem:height:product:polynomial},
\begin{align*}
h_S(F) & \leqslant \sum_{i=1}^n h_{S}(X-\alpha_i) + (n-1)\ln(2) \int_{\Omega_{\infty}} \nu(d\omega) \\
& = \sum_{i=1}^n h_S(1 : \alpha_i) + (n-1)\ln(2) \int_{\Omega_{\infty}} \nu(d\omega) \\
& = n h_S(1 : \alpha)  + (n-1)\ln(2) \int_{\Omega_{\infty}} \nu(d\omega) 
\leqslant \delta C + (\delta -1) \ln(2) \int_{\Omega_{\infty}} \nu(d\omega).
\end{align*}
Note that 
$h_S(1 : a_i) \leqslant h_S(F)$ and $a_i \in K$
for all $i=0, \ldots, n-1$,
so that one can see that there are finitely many possibilities of $F$ because $S$ has the Northcott property.
Therefore the assertion of the claim follows.
\end{proof}

Let us go back to the proof of Theorem~\ref{thm:Northcott:thm}.
For $i=0,\ldots, n$, let
\[
\Upsilon_i := \{ x = (x_0 : \cdots : x_n) \in \mathbb P^n(K^{ac}) \,:\, h_S(x) \leqslant C,\  [K(x):K] \leqslant \delta,\ x_i \not= 0 \}.
\]
It is sufficient to show that $\#(\Upsilon_i) < \infty$ for all $i$. 
Without of loss of generality, we may assume that $i=0$.
Then
\[
\Upsilon_0 = \{ (a_1, \ldots, a_n) \in (K^{ac})^n \,:\, h_S(1, a_1, \ldots, a_n) \leqslant C,\ 
[K(a_1, \ldots, a_n) : K] \leqslant \delta  \}.
\]
Note that $[K(a_i) : K] \leqslant [K(a_1, \ldots, a_n) : K]$ and
$h_S(1 : a_i) \leqslant h_S(1, a_1, \ldots, a_n)$
for all $i = 1, \ldots, n$. Thus the assertion is a consequence of the above special case.
\end{proof}

\begin{coro}
We assume that $S$ has the Northcott property. 
Let $K'$ be a finite extension of $K$. Then $S \otimes_K {K'}$ has also the Northcott property.
\end{coro}

\begin{rema}
Theorem~\ref{thm:Northcott:thm} can be generalized to the case of an adelic vector bundle. For details,
see Proposition~\ref{Pro:Northcott for arithmetic varieties}.
\end{rema}

\section{Measurability of  automorphism actions}

Let $S=(K,(\Omega,\mathcal A,\nu),\phi)$ be an adelic curve and $L/K$ be an algebraic  extension. We denote by $\mathrm{Aut}_K(L)$ the group of field automorphisms of $L$ which are $K$-linear. The group $\mathrm{Aut}_K(L)$ acts on $M_{L}$ as follows: for any $\tau\in \mathrm{Aut}_K(L)$ and any $x\in M_{L}$, one has 
\[\forall\, a\in L,\quad |a|_{\tau(x)}=|\tau(a)|_{x}.\] 
Moreover, by definition the restrictions of the absolute values $|\ndot|_x$ and $|\ndot|_{\tau(x)}$ on $K$ are the same. Therefore we obtain an action of the $K$-linear automorphism group $\mathrm{Aut}_K(L)$ on the set $\Omega_{L}=\Omega\times_{M_K}M_{L}$ (where we consider trivial actions of $\mathrm{Aut}_K(L)$ on $\Omega$ and on $M_K$).

\begin{prop}\label{Pro: measurability of Galois action}
Let $S=(K,(\Omega_K,\mathcal A_K,\nu_K),\phi_K)$ be an adelic curve and $L/K$ be an algebraic extension. For any $\tau\in\mathrm{Aut}_K(L)$, the action of $\tau$ on $\Omega_{L}$ is measurable, where on $\Omega_L$ we consider the $\sigma$-algebra $\mathcal A_L$ in the adelic structure of $S\otimes_K{L}$. 
\end{prop}
\begin{proof}
By Proposition \ref{Pro:sigma algebra on extension} the $\sigma$-algebra $\mathcal A_L$ is the smallest $\sigma$-algebra which makes measurable the canonical projection map $\pi_{L/K}:\Omega_{L}\rightarrow\Omega_K$ and the functions $(x\in\Omega_{L})\mapsto|a|_x$, where $a\in L$. Let $\tau\in\mathrm{Aut}_K(L)$. To show the measurability of the action of $\tau$ on $\Omega_{L}$, it suffices to verify the measurability of the map $\pi_{L/K}\circ\tau$ and the functions $(x\in\Omega_{L})\mapsto|a|_{\tau(x)}$ with $a\in L$. Note that by definition $\pi_{L/K}\circ\tau=\pi_{L/K}$ and $|\alpha|_{\tau(x)}=|\tau(\alpha)|_x$. The proposition is thus proved.
\end{proof}

\begin{prop}\label{Pro: max of measurable function}
Let $S=(K,(\Omega_K,\mathcal A_K,\nu_K),\phi_K)$ be an adelic curve and $L/K$ be a finite extension. Let $F:\Omega_L\rightarrow\mathbb R$ be an $\mathcal A_L$-measurable function. For any $\omega\in\Omega$, let \[f(\omega)=\max_{x\in\pi_{L/K}^{-1}(\{\omega\})}F(x).\]
Then the function $f:\Omega_K\rightarrow\mathbb R$ is $\mathcal A_K$-measurable.
\end{prop}
\begin{proof} We first assume that the extension $L/K$ is normal. 
By \cite{Bourbaki64}, Chapitre VI, \S 8, $\text{n}^{\circ} 6$, Proposition 7, for any $\omega\in\Omega_K$, the action of the $K$-linear automorphism group $\mathrm{Aut}_K(L)$ on 
$M_{L,\omega}$ is transitive. As a consequence, if we denote by $\widetilde F$ the function \[\displaystyle\max_{\tau\in\mathrm{Aut}_K(L)}F\circ\tau,\] then for each $\omega\in\Omega_K$, the restriction of $\widetilde F$ on $\pi_{L/K}^{-1}(\{\omega\})$ is constante, the value of which is equal to $f(\omega)$. By Proposition \ref{Pro: measurability of Galois action}, for any $\tau\in\mathrm{Aut}_K(L)$, the action of $\tau$ on $\Omega_L$ is measurable and hence the function $F\circ\tau$ is $\mathcal A_L$-measurable. Since $\mathrm{Aut}_K(L)$ is a finite set, we deduce that the function $\widetilde F$ is also $\mathcal A_L$-measurable. By Proposition \ref{Pro:fibreint}, the function $f=I_{L/K}(\widetilde F)$ is $\mathcal A_K$-measurable. 

In the general case, we pick a finite normal extension $L_1/K$ which contains $L$. By applying the proved result to the function $F\circ\pi_{L_1/L}$, we still obtain the measurability of the function $f$. The proposition is thus proved.  
\end{proof}

\section{Morphisms of adelic curves}\label{Sec: Morphism of adelic curve}

In this section, we consider morphism of adelic curves. 
\begin{defi}
Let $S=(K,(\Omega,\mathcal A,\nu),\phi)$ and $S'=(K',(\Omega',\mathcal A',\nu'),\phi')$ be two adelic curves, we call \emph{morphism from $S'$ to $S$}\index{morphism from S' to S@morphism from $S'$ to $S$} any triplet $\alpha=(\alpha^{\#},\alpha_{\#},I_\alpha)$, where 
\begin{enumerate}[label=\rm(\alph*)]
\item\label{Item: alpha upper sharp homomorphism} $\alpha^\#:K\rightarrow K'$ is a field homomorphism,
\item\label{Item: alpha lower sharp} $\alpha_\#:(\Omega',\mathcal A')\rightarrow(\Omega,\mathcal A)$ is a measurable map such that the following diagram is commutative
\[\xymatrix{\Omega'\ar[d]_-{\phi'}\ar[r]^-{\alpha_{\#}}&\Omega\ar[d]^-{\phi}\\
M_{K'}\ar[r]_-{-\circ \alpha^{\#}}&M_K
}\]
and that the direct image of $\nu'$ by $\alpha_{\#}$ coincides with $\nu$, namely, for any function  $f\in\mathscr L^1(\Omega,\mathcal A,\nu)$, one has
\[\int_{\Omega'}f\circ \alpha_{\#}\,\mathrm{d}\nu'=\int_\Omega f\,\mathrm{d}\nu.\]
\item\label{Item: morphism of adelic curve integral} $I_\alpha: L^1(\Omega',\mathcal A',\nu')\rightarrow L^1(\Omega,\mathcal A,\nu)$ is a disintegration kernel of $\alpha_\#$, namely $I_\alpha$ is a linear map such that, for any element $g\in L^1(\Omega',\mathcal A',\nu')$, one has
\[\int_{\Omega}I_\alpha(g)\,\mathrm{d}\nu=\int_{\Omega'}g\,\mathrm{d}\nu',\]
and for any function $f\in\mathscr L^1(\Omega,\mathcal A,\nu)$ one has $I_{\alpha}$ sends the equivalence class of $f\circ\alpha_\#$ to that of $f$.
\end{enumerate}
Naturally, if $S$, $S'$ and $S''$ are adelic curves and if $\alpha=(\alpha^\#,\alpha_{\#},I_\alpha):S'\rightarrow S$ and $\beta=(\beta^\#,\beta_\#,I_\beta):S''\rightarrow S$ are morphisms of adelic curves, then $\alpha\circ\beta:=(\beta^\#\circ\alpha^\#,\alpha_\#\circ\beta_\#,I_\alpha\circ I_\beta)$ forms a morphism of adelic curves from $S''$ to $S$. Thus the adelic curves and their morphisms form a category.
\end{defi}  

\begin{exem}
Let $S=(K,(\Omega,\mathcal A,\nu),\phi)$ be an adelic curve. 
\begin{enumerate}[label=\rm(\arabic*)]
\item If $S\otimes_KK'=(K',(\Omega_{K'},\mathcal A_{K'},\nu_{K'}),\phi_{K'})$ is an algebraic extension of $S$, then the triplet $(K\hookrightarrow K',\pi_{K'/K},\overline I_{K'/K})$ is a morphism of adelic curves from $S\otimes_KK'$ to $S$, where $\overline{I}_{K'/K}:L^1(\Omega_{K'},\mathcal A_{K'},\nu_{K'})\rightarrow L^1(\Omega,\mathcal A,\nu)$ is the linear map induced by $I_{K'/K}:\mathscr L^1(\Omega_{K'},\mathcal A_{K'},\nu_{K'})\rightarrow\mathscr L^1(\Omega,\mathcal A,\nu)$.
\item Assume that $K_0$ is a subfield of $K$. Let $S_0$ be the field $K_0$ equipped with the restriction to $K_0$ of the adelic structure of $S$ (see Subsection \ref{Subsec: adelic curve restriction}). Then the triplet $(K_0\hookrightarrow K,\Id_\Omega,\Id_{ L^1(\Omega,\mathcal A,\nu)})$ forms a morphism of adelic curves from $S$ to $S_0$.
\item Let $K=\mathbb Q(T)$ be the field of rational functions of one variable $T$ with coefficients in $\mathbb Q$. Let $S=(K,(\Omega,\mathcal A,\nu),\phi)$ constructed in Subsection \ref{subsec:fun:field:Q}. Recall that $(\Omega,\mathcal A,\nu)$ is written as a disjoint union $\Omega_h\coprod\mathcal P\coprod [0,1]_*$, where $\Omega_h$ is the set of closed points of $\mathbb P^1_{\mathbb Q}$, $\mathcal P$ is the set of prime numbers, and $[0,1]_*$ is the subset of $[0,1]$ of $t$ such that $\mathrm{e}^{2\pi it}$ is transcendental. Let $S_{\mathbb Q}=(\mathbb Q,(\Omega_{\mathbb Q},\mathcal A_{\mathbb Q},\nu_{\mathbb Q}),\phi_{\mathbb Q})$ be the adelic curve defined in Subsection \ref{Subsec:Numberfields} and $\widetilde{S}_{\mathbb Q}=(\mathbb Q,(\widetilde{\Omega}_{\mathbb Q},\widetilde{\mathcal A}_{\mathbb Q},\widetilde{\nu}_{\mathbb Q}),\widetilde{\phi}_{\mathbb Q})$ be the adelic curve consisting of the filed $\mathbb Q$ equipped with the amalgamation of the adelic structure of $S_{\mathbb Q}$ and a family of copies of the trivial absolute value on $\mathbb Q$ indexed by $\Omega_h$. We can also write $\widetilde{\Omega}$ as the disjoint union of three subsets $\Omega_h\coprod\mathcal P\coprod\{\infty\}$, where $0$ denotes the trivial absolute value on $\mathbb Q$ and $\infty$ denotes the infinite place of $\mathbb Q$. Let $\alpha^{\#}:\mathbb Q\rightarrow\mathbb Q(T)$ be the inclusion map. Let $\alpha_{\#}:\Omega\rightarrow\widetilde{\Omega}_{\mathbb Q}$ be the map which sends any element of $\Omega_h\coprod\mathcal P$ to itself and send any element of $[0,1]_*$ identically to $\infty$. Finally, let $I_\alpha: L^1(\Omega,\mathcal A,\nu)\rightarrow  L^1(\widetilde{\Omega}_{\mathbb Q},\widetilde{\mathcal A}_{\mathbb Q},\widetilde{\nu}_{\mathbb Q})$ be the linear map sending the equivalence class of any function $f\in\mathscr L^1(\Omega,\mathcal A,\nu)$ to that of the function $I_\alpha(f)$ sending $\omega\in\Omega_h\coprod\mathcal P$ to $f(\omega)$ and $\infty$ to $\int_{[0,1]_*}f(t)\,\mathrm{d}t$. Then the triplet $(\alpha^\#,\alpha_\#,I_\alpha)$ forms a morphism of adelic curves from $S$ to $\widetilde{S}_{\mathbb Q}$.
\end{enumerate}
\end{exem}

%% file: ch4_2019_03_23.tex

\chapter{Vector bundles on adelic curves: global theory}\label{Chp:globaltheory}

\IfChapVersion
\ChapVersion{Version of Chapter 4 : \\ \StrSubstitute{\DateChapFour}{_}{\_}}
\fi

The purpose of this chapter is to study the geometry of adelic curves, notably the divisors and vector bundles.

\section{Norm families}
Let $S=(K,(\Omega,\mathcal A,\nu),\phi)$ be an adelic curve (see \S\ref{SubSec:adelicarithcurve}). Recall that for any $\omega\in\Omega$, we denote by $|\ndot|_\omega$ the absolute value of $K$ indexed by $\omega$. Note that in the case where $|\ndot|_\omega$ is Archimedean, there exists a constant $\kappa(\omega)$, $0<\kappa(\omega)\leqslant 1$, such that $|\ndot|_\omega=|\ndot|^{\kappa(\omega)}$, where $|\ndot|$ denotes the usual absolute value on $\mathbb R$ or $\mathbb C$. For simplicity, we assume that $\kappa(\omega)=1$ for any $\omega\in\Omega_\infty$, namely $|\ndot|_\omega$ identifies with the usual absolute value on $\mathbb R$ or $\mathbb C$. Note that this assumption is harmless for the generality of the theory since in general case we can replace the absolute values $\{|\ndot|_\omega\}_{\omega\in\Omega_\infty}$ by the usual ones and consider the measure $\mathrm{d}\widetilde{\nu}=(\indic_{\Omega\setminus\Omega_\infty}+\kappa\indic_{\Omega_\infty})\,\mathrm{d}\nu$
instead. 

\subsection{Definition and algebraic constructions}\label{Subsec:Norm families}

Let $E$ be a vector space of finite rank over $K$. We denote by $\mathcal N_E$ the set of norm families\index{norm family} $\{\|\ndot\|_\omega\}_{\omega\in \Omega}$, where each $\|\ndot\|_\omega$ is a norm on $E_{K_\omega}:=E\otimes_KK_\omega$. We say that a norm family $\xi=\{\|\ndot\|_\omega\}_{\omega\in \Omega}$ is \emph{ultrametric on $\Omega\setminus\Omega_\infty$}\index{ultrametric}\index{norm family!ultrametric ---} if the norm $\|\ndot\|_\omega$ is ultrametric for any $\omega\in\Omega\setminus\Omega_\infty$.
We say that a norm family $\{\|\ndot\|_\omega\}_{\omega\in \Omega}$ in $\mathcal N_E$ is \emph{Hermitian}\index{Hermitian}\index{norm family!Hermitian ---} if the following conditions are satisfied:
\begin{enumerate}[label=\rm(\alph*)]
\item for any $\omega\in\Omega\setminus\Omega_\infty$, the norm $\|\ndot\|_\omega$ is ultrametric (namely the norm family is ultrametric on $\Omega\setminus\Omega_\infty$);
\item for any $\omega\in\Omega_\infty$, the norm $\|\ndot\|_{\omega}$ is induced by an inner product (see \S\ref{Subsec:innerproduct}), namely there exists an inner product $\langle\;,\;\rangle_{\omega}$ on $E_{K_\omega}$ such that 
$\|\ell\|_\omega=\langle\ell,\ell\rangle_\omega^{1/2}$
for any $\ell\in E_{K_\omega}$.
\end{enumerate} We denote by $\mathcal H_E$ the subset of $\mathcal N_E$ consisting of all Hermitian norm families. In the following, we describe some algebraic constructions of norm families.

{
\subsubsection{Multiplication by a numerical function}\label{Subsubsec:Multiplication by a numerical function} Let $E$ be a finite-dimensional vector space over $K$, $\xi=\{\norm{\ndot}_\omega\}_{\omega\in\Omega}\in\mathcal N_E$ and $f:\Omega\rightarrow\intervalle{]}{0}{+\infty}{[}$ be a positive function on $\Omega$. We denote by $f\xi$ the norm family $\{f(\omega)\norm{\ndot}_{\omega}\}_{\omega\in\Omega}$ in $\mathcal N_E$.
}

\subsubsection{Restrict and quotient norm families}\label{page:restricquot}

Let $E$ be a finite-dimensional vector space over $K$ and $\xi=\{\|\ndot\|_{E,\omega}\}_{\omega\in\Omega}$ be a norm family in $\mathcal N_E$. Let $F$ be a vector subspace of $E$. For any $\omega\in\Omega$, let $\|\ndot\|_{F,\omega}$ be the restriction of the norm $\|\ndot\|_{E,\omega}$ {to} $F_{K_\omega}$ (see Definition \ref{Def:restriction}). Then $\{\|\ndot\|_{F,\omega}\}_{\omega\in\Omega}$ forms a norm family in $\mathcal N_F$, called the \emph{restriction}\index{restriction}\index{norm family!restriction} of $\xi$ to $F$. Similarly, if $G$ is a quotient vector space of $E$, then each norm $\|\ndot\|_{E,\omega}$ induces by quotient a norm $\|\ndot\|_{G,\omega}$ on $G_{K_\omega}$ (see \S\ref{Subsec:Quotientnorm}). Thus we obtain a norm family $\{\|\ndot\|_{G,\omega}\}_{\omega\in\Omega}$ in $\mathcal N_{G}$, called the \emph{quotient}\index{quotient}\index{norm family!quotient ---} of $\xi$ on $G$. Note that, if the norm family $\xi$ is ultrametric on $\Omega\setminus\Omega_\infty$ (resp. Hermitian), then all its restrictions and quotients are also ultrametric on $\Omega\setminus\Omega_\infty$ (resp. Hermitian).

\subsubsection{Direct sums} Let $\mathscr S$ be the set of all convex and continuous functions $f:[0,1]\rightarrow [0,1]$ such that $\max\{t,1-t\}\leqslant f(t)$ for any $t\in[0,1]$. If $E$ and $F$ are finite-dimensional vector spaces over $K$ and if $\xi_E=\{\|\ndot\|_{E,\omega}\}_{\omega\in\Omega}$ and $\xi_F=\{\|\ndot\|_{F,\omega}\}_{\omega\in\Omega}$ are respectively norm families in $\mathcal N_E$ and $\mathcal N_F$, for any family $\psi=\{\psi_\omega\}_{\omega\in\Omega}$ of elements in $\mathscr S$ we define a norm family $\xi_{E}\oplus_{\psi}\xi_F=\{\|\ndot\|_\omega\}_{\omega\in\Omega}$ in $\mathcal N_{E\oplus F}$ (see Subsection \ref{Subsec:directsums}) such that, for any $(x,y)\in E_{K_\omega}\oplus F_{K_\omega}$,
\[\|(x,y)\|_\omega:=
(\|x\|_{E,\omega}+\|y\|_{F,\omega})\psi_\omega\Big(\frac{\|x\|_{E,\omega}}{\|x\|_{E,\omega}+\|y\|_{F,\omega}}\Big).\]
We call $\xi_{E}\oplus_{\psi}\xi_F$ the $\psi$-\emph{direct sum}\index{direct sum}\index{norm family!direct sum} of $\xi_E$ and $\xi_F$. 
If both norm families $\xi_E$ and $\xi_F$ are Hermitian, and if
\[\psi_\omega(a,b)=\begin{cases}
\max\{a,b\},&\omega\in\Omega\setminus\Omega_\infty,\\
(a^2+b^2)^{1/2},&\omega\in\Omega_\infty,
\end{cases}\]
then the direct sum $\xi_E\oplus_\psi\xi_F$ belongs to $\mathcal H_{E\oplus F}$. We call it the \emph{orthogonal direct sum}\index{orthogonal direct sum}\index{norm family!orthogonal direct sum} of $\xi_E$ and $\xi_F$. 

\subsubsection{Dual norm family} Let $E$ be a vector space of finite rank over $K$ and $\xi=\{\|\ndot\|_{\omega}\}_{\omega\in\Omega}$ be a norm family in $\mathcal N_E$. The dual norms (see \S\ref{Subsec:dualnorm}) $\{\|\ndot\|_{\omega,*}\}_{\omega\in\Omega}$ form a norm family in $\mathcal N_{E^\vee}$, called the \emph{dual}\index{dual}\index{normal family!dual ---} of $\xi$, denoted by $\xi^\vee$. Note that $\xi^\vee$ is always ultrametric on $\Omega\setminus\Omega_\infty$, and it is Hermitian if $\xi$ is. Moreover, if for any $\omega\in\Omega\setminus\Omega_\infty$ the norm $\|\ndot\|_\omega$ is ultrametric, then one has $(\xi^{\vee})^\vee=\xi$, where we identify $E$ with its double dual space $E^{\vee\vee}$ (see Proposition \ref{Pro:doubledualarch} and Corollary \ref{Cor:doubledual}). In particular, if $E$ is a $K$-vector space of dimension $1$, then one has $(\xi^\vee)^\vee=\xi$.\label{doubledual}

\subsubsection{Tensor products} 
Let $\{E_i\}_{i=1}^n$ be a family of finite-dimensional vector spaces over $K$. For any $i\in\{1,\ldots,n\}$, let $\xi_i=\{\|\ndot\|_{i,\omega}\}_{\omega\in\Omega}$ be an element in $\mathcal N_{E_i}$. We denote by $\xi_1\otimes_\pi\cdots\otimes_\pi\xi_n$ the norm family $\{\|\ndot\|_{\omega,\pi}\}_{\omega\in\Omega}$ in $\mathcal N_{E_1\otimes\cdots\otimes E_n}$, where $\|\ndot\|_{\omega,\pi}$ is the $\pi$-tensor product of the norms $\|\ndot\|_{i,\omega}$, $i\in\{1,\ldots,n\}$ (see \S\ref{Subsec:tensorproduct}). The norm family $\xi_1\otimes_\pi\cdots\otimes_\pi\xi_n$ is called the $\pi$-\emph{tensor product}\index{pi tensor product@$\pi$-tensor product}\index{norm family!pi tensor product@$\pi$-tensor product} of $\xi_1,\ldots,\xi_n$. Similarly, we denote by $\xi_1\otimes_\varepsilon\cdots\otimes_{\varepsilon}\xi_n$ the norm family in $\mathcal N_{E_1\otimes\cdots\otimes E_n}$ consisting of $\varepsilon$-tensor products (see \S\ref{Subsec:tensorproduct}) of the norms $\|\ndot\|_{1,\omega},\ldots,\|\ndot\|_{n,\omega}$. We call $\xi_1\otimes_\varepsilon\cdots\otimes_{\varepsilon}\xi_n$ the $\varepsilon$-\emph{tensor product}\index{epsilon tensor product@$\varepsilon$-tensor product}\index{norm family!epsilon tensor product@$\varepsilon$-tensor product} of $\xi_1,\ldots,\xi_n$. We also introduce the following mixed version of $\varepsilon$-tensor product and $\pi$-tensor product. We denote by $\xi_1\otimes_{\varepsilon,\pi}\cdots\otimes_{\varepsilon,\pi}\xi_n$ consisting of the norms $\|\ndot\|_{\omega,\varepsilon}$ if $\omega\in\Omega\setminus\Omega_\infty$ and $\|\ndot\|_{\omega,\pi}$ if $\omega\in\Omega_\infty$. This norm family is called the $\varepsilon,\pi$-tensor product of $\xi_1,\ldots,\xi_n$. Note that the $\varepsilon$-tensor product and the $\varepsilon,\pi$-tensor product are both ultrametric on $\Omega\setminus\Omega_\infty$.

If all norm families $\xi_1,\ldots,\xi_n$ are Hermitian, we denote by $\xi_1\otimes\cdots\otimes\xi_n$ the norm family $\{\|\ndot\|_{\omega}\}_{\omega\in\Omega}$ in $\mathcal H_{E_1\otimes\cdots\otimes E_n}$, where for each $\omega\in\Omega\setminus\Omega_\infty$, the norm $\|\ndot\|_{\omega}$ is the $\varepsilon$-tensor product of $\{\|\ndot\|_{i,\omega}\}_{i\in\{1,\ldots,n\}}$, and for each $\omega\in\Omega_\infty$, the norm $\|\ndot\|_{\omega}$ is the orthogonal tensor product (see \S\ref{Subsec:Hilbert-Schmidt}) of $\{\|\ndot\|_{i,\omega}\}_{i\in\{1,\ldots,n\}}$. We call $\xi_1\otimes\cdots\otimes\xi_n$ the \emph{orthogonal tensor product}\index{orthogonal tensor product}\index{norm family!orthogonal tensor product} of $\xi_1,\ldots,\xi_n$.\label{Page:orthogonaltp}

\subsubsection{Exterior powers}
Let $E$ be a finite-dimensional vector space over $K$ and $\xi=\{\|\ndot\|_{E,\omega}\}_{\omega\in\Omega}$ be a norm family in $\mathcal N_E$. Let $i$ be a non-negative integer. We equip $E^{\otimes i}$ with the  $\varepsilon,\pi$-tensor power of the norm family $\xi$, which induces by  quotient a norm family on the exterior power $\Lambda^i(E)$ which we denote by $\Lambda^i\xi$.

\subsubsection{Determinant}
Let $E$ be a finite-dimensional vector space over $K$ and $\xi=\{\|\ndot\|_{E,\omega}\}_{\omega\in\Omega}$ be a norm family in $\mathcal N_E$. Each norm $\|\ndot\|_{E,\omega}$ induces a determinant norm $\|\ndot\|_{\det(E),\omega}$ on $\det(E)\otimes_KK_\omega\cong\det(E_{K_{\omega}})$ (see \S\ref{SubSec:determinant norm}). The norm family $\{\|\ndot\|_{\det(E),\omega}\}_{\omega\in\Omega}$ is called the \emph{determinant}\index{determinant}\index{norm family!determinant} of $\xi$, denoted by $\det(\xi)$. By Proposition \ref{Pro:doubledualdet} we obtain that $\det(\xi)$ coincides with $\Lambda^r\xi$, where $r$ is the rank of $E$ over $K$.

\label{Page:extensionofscalars}
\subsubsection{Extension of scalars} Let $E$ be a vector space of finite rank over $K$ and $\xi=\{\|\ndot\|_{\omega}\}_{\omega\in\Omega}$ be a norm family in $\mathcal N_E$. Let $L/K$ be an algebraic extension of the field $K$ and $((\Omega_L,\mathcal A_L,\nu_L),\phi_L)$ be the adelic structure of the adelic curve $S\otimes_KL$. We construct a norm family $\xi_L=\{\|\ndot\|_x\}_{x\in\Omega_L}\in\mathcal N_{E\otimes_KL}$ as follows: for any $x\in\Omega_L$ whose canonical image in $\Omega$ is $\omega$, if $|\ndot|_\omega$ is non-Archimedean, $\|\ndot\|_{x}$ is the norm $\|\ndot\|_{\omega,L_x,\varepsilon}$ on $(E_{K_\omega})\otimes_{K_\omega}L_x$ induced by $\|\ndot\|_\omega$ by $\varepsilon$-extension of scalars; otherwise $\|\ndot\|_{x}$ is the norm $\|\ndot\|_{\omega,L_x,\pi}$ on $(E_{K_\omega})\otimes_{K_\omega}L_x$ induced by $\|\ndot\|_\omega$ by $\pi$-extension of scalars (see \S\ref{Subsec:extensionofscalars}). By Proposition~\ref{Pro:comparisonofdualnormes:scalar:extension} \ref{Pro:comparisonofdualnormes:scalar:extension:epsilon}, 
if the rank of $E$ over $K$ is $1$, then the norm family $(\xi^\vee)_L$ identifies with the dual norm family of $\xi_L$. Moreover, by Corollary \ref{Cor:extsucc}, if $L_2/L_1/K$ are successive algebraic extensions, then one has $(\xi_{L_1})_{L_2}=\xi_{L_2}$, where we identify $E\otimes_KL_2$ with $(E\otimes_KL_1)\otimes_{L_1}L_2$.

Let $\xi=\{\|\ndot\|_{\omega}\}_{\omega\in\Omega}\in\mathcal H_E$ be a Hermitian norm family, where for $\omega\in\Omega$, the norm $\|\ndot\|_\omega$ is induced by an inner product $\langle\kern.15em,\kern.02em\rangle_\omega$. We denote by $\xi^H_{L}=\{\|\ndot\|_x\}_{x\in\Omega_L}\in\mathcal H_{E\otimes_KL}$ the Hermitian norm family defined as follows. For any  $x\in\Omega_{L}\setminus\Omega_{L,\infty}$ over $\omega\in\Omega\setminus\Omega_\infty$, one has $\|\ndot\|_x=\|\ndot\|_{\omega,L_{x},\varepsilon}$; for any $x\in\Omega_{L,\infty}$ over $\omega\in\Omega_\infty$, $\|\ndot\|_x$ is the norm $\|\ndot\|_{\omega,L_x,\mathrm{HS}}$ induced by the inner product $\langle\kern.15em,\kern.02em\rangle_{\omega,L_x}$ on $E_{L_x}$ which extends $\langle\kern.15em,\kern.02em\rangle_\omega$ on $E_{K_\omega}$ (namely $\|\ndot\|_x$ is the orthogonal tensor product norm of $\|\ndot\|_\omega$ and $|\ndot|_x$ if we identify $E_{L_x}$ with $E_{K_\omega}\otimes_{K_\omega}L_x$, see Remark \ref{Rem:extensiondoubledual}). One has $(\xi^\vee)^H_L=(\xi^H_L)^\vee$ (see Proposition~\ref{Pro:comparisonofdualnormes:scalar:extension} \ref{Pro:comparisonofdualnormes:scalar:extension:epsilon}
for the ultrametric part and Remark \ref{Rem:extensiondoubledual} for the inner product part).

{
\subsubsection{Comparison of norm families}
Let $E$ be a finite-dimensional vector space over $K$, and $\xi=\{\norm{\ndot}_\omega\}_{\omega\in\Omega}$ and $\xi'=\{\norm{\ndot}_\omega'\}_{\omega\in\Omega}$ be two elements of $\mathcal N_E$. We say that $\xi$ is \emph{smaller than}
\index{norm family!smaller than another norm family} $\xi'$ if for any $\omega\in\Omega$ one has $\norm{\ndot}_\omega\leqslant\norm{\ndot}_\omega'$. The condition ``$\xi$ is smaller than $\xi$'' is denoted by $\xi\leqslant \xi'$ or $\xi'\geqslant \xi$.
}

\subsubsection{Local distance} Let $E$ be a finite-dimensional vector space over $K$, and $\xi=\{\|\ndot\|_\omega\}_{\omega\in\Omega}$ and $\xi'=\{\|\ndot\|_\omega'\}_{\omega\in\Omega}$ be two norm families in $\mathcal N_E$. For any $\omega\in\Omega$, let
\[d_\omega(\xi,\xi'):=\sup_{s\in E_{K_\omega}\setminus \{ 0 \}}
\Big|\ln\|s\|_\omega-\ln\|s\|_{\omega}'\Big|.\]
We call $d_\omega(\xi,\xi')$ the \emph{local distance}\index{local distance}\index{norm family!local distance}\label{Page:localdist} on $\omega$ of the norm families $\xi$ and $\xi'$. By Proposition \ref{Pro:distanceofoperatornorms}, one has
\begin{equation}\label{Equ:distdualnorm}d_\omega(\xi^\vee,(\xi')^{\vee})\leqslant d_\omega(\xi,\xi'),\end{equation}
and the equality holds if $\omega\in\Omega_\infty$ or if $\|\ndot\|_\omega$ and $\|\ndot\|_{\omega}'$ are both ultrametric.

\subsection{Dominated norm families}

Let $E$ be a finite-dimensional vector space over $K$ and $\xi=\{\|\ndot\|_\omega\}_{\omega\in\Omega}$ be a norm family in $\mathcal N_E$.

\begin{defi}
We say that the norm family $\xi$ is \emph{upper dominated}\index{upper dominated}\index{norm family!upper dominated} if, 
for any non-zero element $s\in E$, there exists a $\nu$-integrable function $A(\ndot)$ on $\Omega$ such that $\ln\|s\|_\omega\leqslant A(\omega)$ $\nu$-almost everywhere. Note that the upper dominancy is equivalent to 
\[
\forall s \in E \setminus \{ 0 \},\quad
\upint_{\Omega} \ln\|s\|_{\omega}\, \nu(d\omega) <+\infty 
\]
with the notation of Definition \ref{Def:upper and lower integral}. Similarly, we say that the norm family $\xi$ is \emph{lower dominated}\index{lower dominated}\index{norm family!lower dominated} if, 
for any non-zero element $s\in E$, there exists a $\nu$-integrable function $B(\ndot)$ on $\Omega$ such that $B(\omega) \leqslant \ln\|s\|_\omega$ $\nu$-almost everywhere. Note that the lower dominancy is equivalent to 
\[
\forall s \in E \setminus \{ 0 \},\quad
\lowint_{\Omega} \ln\|s\|_{\omega}\, \nu(d\omega) > -\infty. 
\]
\end{defi}

\begin{defi}\label{Def:dominatednormfamily}
We say that $\xi$ is \emph{dominated}\index{dominated}\index{norm family!dominated} if $\xi$ and $\xi^{\vee}$ are both
upper dominated. Note that the upper dominancy of $\xi$ and $\xi^{\vee}$ implies the lower dominancy of $\xi^{\vee}$ and $\xi$, respectively, because (see Proposition \ref{Pro:upperlownega})
$\ln \| \alpha \|_{\omega, *} + \ln \|s \|_{\omega} \geqslant 0$
for all $s \in E$ and $\alpha \in E^{\vee}$ with $\alpha(s) = 1$, so that
if $\xi$ is dominated, then $\xi$ and $\xi^{\vee}$ are upper and lower dominated.\end{defi}

\begin{rema}\label{Rem: dominancy of double dual}
If $\xi$ is a dominated norm family, then also is $\xi^\vee$. In fact, for any $\omega\in\Omega$ one has $\|\ndot\|_{\omega,**}\leqslant\|\ndot\|_{\omega}$ (see \eqref{Equ:doubledual}). Therefore the upper dominancy of $\xi$ implies that of $\xi^{\vee\vee}$. The converse is true when $\|\ndot\|_\omega$ is ultrametric for $\omega\in\Omega\setminus\Omega_\infty$ since in this case one has $\|\ndot\|_{\omega,**}=\|\ndot\|_\omega$ for any $\omega\in\Omega$ (see Proposition \ref{Pro:doubledualarch} and Corollary \ref{Cor:doubledual}).
\end{rema}

\begin{rema}
It is \emph{not} true that if $\xi$ is upper and lower dominated then it is dominated. Consider the following example. Let $K$ be an infinite field. We equip $K$ with the discrete $\sigma$-algebra $\mathcal A$ and let $\nu$ be the atomic measure on $K$ such that $\nu(\{a\})=1$ for any $a\in K$. For any $a\in K$, let $|\ndot|_a$ be the trivial absolute value on $K$. Then $S=(K,(K,\mathcal A,\nu),\{|\ndot|_a\}_{a\in K})$ forms an adelic curve. Consider now the vector space $E=K^2$ over $K$. For any $a\in K$ let $\|\ndot\|_a$ be the norm on $K^2$ such that
\[\|(x,y)\|_a=\begin{cases}
0,&\text{if }x=y=0,\\
1/2,&\text{if }y=ax,\,x\neq 0,\\
1,&\text{else}.
\end{cases}\] 
Then for any vector $s\in K^2$, $s\neq 0$, one has $\|s\|_a=1$ for all except at most one  $a\in K$. Therefore the function $(a\in K)\mapsto\ln\|s\|_a$ on $K$ is integrable, and in particular dominated. Therefore the norm family $\xi=\{\|\ndot\|_a\}_{a\in K}$ is upper dominated and lower dominated. However, for any $a\in K$, the dual norm $\|\ndot\|_{a,*}$ on $K^2$ (we identify $K^2$ with the dual vector space of itself in the canonical way) satisfies
\[\|(x,y)\|_{a,*}=\begin{cases}
0,&\text{if }x=y=0,\\
1,&\text{if }x=-ay,\,y\neq 0\\
2,&\text{else}.
\end{cases}\]
Therefore, for any non-zero element $s\in K^2$, one has $\ln\|s\|_{a,*}=\ln(2)$ for all except at most one element $a\in K$. The dual norm family $\xi^\vee$ is thus not upper dominated.
\end{rema}

\begin{exem}\label{Exa:adelicbasis} A fundamental example of dominated norm family is that arising from a basis.
Let $E$ be a vector space of finite rank $r$ over $K$, and $\boldsymbol{e}=\{e_1,\ldots,e_r\}$ be a basis of $E$ over $K$. For any algebraic extension $L/K$ and any $x\in\Omega_{L}$,  let $\|\ndot\|_{\boldsymbol{e},x}$ be the norm on $E\otimes_K L_x$ such that, for any $(\lambda_1,\ldots,\lambda_r)
\in L_x^r$,
\[\|\lambda_1e_1+\cdots+\lambda_re_r\|_{\boldsymbol{e},x}:=\begin{cases}\max_{i\in\{1,\ldots,r\}}|\lambda_i|_x,&\text{if $x\in\Omega_{L}\setminus\Omega_{L,\infty}$}\\[1ex]
|\lambda_1|_x+\cdots+|\lambda_r|_x,&\text{if $x\in\Omega_{L,\infty}$,}
\end{cases}\]
where $\Omega_{L,\infty}$ denotes the set of  all $x\in\Omega_L$ such that the absolute value $|\ndot|_x$ is Archimedean. Let $\xi_{\boldsymbol{e}}$ be the norm family $\{\|\ndot\|_{\boldsymbol{e},\omega}\}_{\omega\in\Omega}$. Note that one has $\xi_{\boldsymbol{e},L}=\{\|\ndot\|_{\boldsymbol{e},x}\}_{x\in\Omega_L}$ for any algebraic extension $L/K$. Moreover, for any  non-zero vector $s=a_1e_1+\cdots+a_re_r\in E\otimes_KL$, with $(a_1,\ldots,a_r)\in L^r$, one has
\[\forall\,x\in\Omega_L,\quad\ln\|s\|_{\boldsymbol{e},x} \leqslant \max_{\begin{subarray}{c}i\in\{1,\ldots,r\}\\
a_i\neq 0\end{subarray}}\ln|a_i|_{x}+\ln(r)\indic_{\Omega_{L,\infty}}(x).\]
Since the functions $x\mapsto\ln|a|_x$ are $\nu_L$-integrable for all $a\in L\setminus\{0\}$ and since $\nu_L(\Omega_{L,\infty})<+\infty$ (see Proposition \ref{Pro:mesurabilitekappa}), we obtain that the function $(x\in\Omega_L)\mapsto\ln\|s\|_{\boldsymbol{e},x}$ is  $\nu_L$-integrable. If we denote by $\{e_1^\vee,\ldots,e_r^\vee\}$ the dual basis of $\boldsymbol{e}$, then for any $\alpha=a_1e_1^\vee+\cdots+a_re_r^\vee\in E^\vee\otimes_KL$ with $(a_1,\ldots,a_r)\in L^r$ and any $x\in\Omega_L$ one has 
\[
\|\alpha\|_{\boldsymbol{e},x,*}= 
\max\{|a_1|_x,\ldots,|a_r|_x\}.
\]
 Therefore the function $(x\in\Omega_L)\mapsto\|\alpha\|_{\boldsymbol{e},x,*}$ is $\mathcal A$-measurable. If $\alpha\neq 0$, then the function $(x\in\Omega_L)\mapsto\ln\|\alpha\|_{\mathbf{e},x,*}$ is $\nu_L$-integrable.
Hence the norm family $\xi_{\boldsymbol{e},L}$ is dominated. Note that (see page \pageref{Page:localdist} for the definition of the local distance function)
\begin{equation}\label{Equ:distnormbasis}\forall\,x\in\Omega_L,\quad d_x((\xi_{\boldsymbol{e},L})^\vee,\xi_{\boldsymbol{e}^\vee,L})\leqslant\ln(r)\indic_{\Omega_{L,\infty}}(x),
\end{equation}
where $\boldsymbol{e}^\vee$ denotes the dual basis of $\boldsymbol{e}$.
\end{exem}

\begin{prop}\label{Pro:adelicvectorbundledist}
Let $E$ be a finite-dimensional vector space over $K$, $\xi_1$ and $\xi_2$ be norm families in $\mathcal N_E$. We assume that $\xi_1$ is dominated. If the local distance function $(\omega\in\Omega)\mapsto d_\omega(\xi_{1},\xi_2)$ is $\nu$-dominated (see Definition \ref{Def:dominancefunc}), then the norm family $\xi_2$ is dominated. In particular, if there exists a basis $\boldsymbol{e}$ of $E$ over $K$ such that the function $(\omega\in\Omega)\mapsto d_{\omega}(\xi_{\boldsymbol{e}},\xi_2)$ is $\nu$-dominated, then the norm family $\xi_2$ is dominated.
\end{prop}
\begin{proof} 
Assume that $\xi_i$ is of the form $\{\|\ndot\|_{i,\omega}\}_{\omega\in\Omega}$, $i\in\{1,2\}$.
Let $s$ be a non-zero element in $E$. For any $\omega\in\Omega$, one has,
\[\ln\|s\|_{2,\omega}-\ln\|s\|_{1,\omega}\leqslant d_\omega(\xi_1,\xi_2)\quad\text{$\nu$-almost everywhere}.\]
Moreover, since the norm family $\xi_1$ is dominated, one has $\upint_{\Omega} \ln\|s\|_{1,\omega}\, \nu(d\omega)<+\infty$;
since the local distance function $d(\xi_1,\xi_2)$ is dominated, one has 
$\upint_{\Omega} d_{\omega}(\xi_1,\xi_2)\,\nu(d\omega) <+\infty$. 
Therefore by Proposition \ref{Pro:subsuperadditive} one has
\[\upint_{\Omega} \ln\|s\|_{2,\omega}\, \nu(d\omega) \leqslant \upint_{\Omega} \ln\|s\|_{1,\omega} \,\nu(d\omega) +\upint_{\Omega} d_{\omega} (\xi_1,\xi_2)\, \nu(d\omega) <+\infty.\]

By \eqref{Equ:distdualnorm}, one has $d_\omega(\xi_1^\vee,\xi_2^\vee)\leqslant d_\omega(\xi_1,\xi_2)$ for any $\omega$. Hence the same argument as above applied to the dual norm families shows that 
\[\forall\,\alpha\in E^\vee\setminus\{0\}, \quad \upint_{\Omega} \ln\|\alpha\|_{2,\omega,*}\,\nu(d\omega) <+\infty.\]
Therefore, the norm family $\xi_2$ is dominated. To establish the last assertion, it suffices to apply the obtained result to the case where $\xi_1=\xi_{\boldsymbol{e}}$ (see Example \ref{Exa:adelicbasis} for the fact that the norm family $\xi_{\boldsymbol{e}}$ is dominated).
\end{proof}

\begin{prop}\label{Pro:comparaisonavecbasenorm}
Let $E$ be a vector space of finite rank over $K$ and $\xi=\{\|\ndot\|_{\omega}\}_{\omega\in\Omega}$ be an element of $\mathcal N_E$ which is dominated.
Then for any basis $\boldsymbol{e}=\{e_1,\ldots,e_r\}$ of $E$, there exists a $\nu$-integrable function $A_{\boldsymbol{e}}$ on $\Omega$ such that, for any algebraic extension $L/K$ and any $x\in\Omega_L$ one has (note that $\xi_K=\xi^{\vee\vee}$ in the case where $L=K$)
\begin{equation}\label{Equ:encadrement}
d_x(\xi_L,\xi_{\boldsymbol{e},L})\leqslant A_{\boldsymbol{e}}(\pi_{L/K}(x)).\end{equation}
In particular, the local distance function $(x\in\Omega_L)\mapsto d_x(\xi_L,\xi_{\boldsymbol{e},L})$ is $\nu_L$-dominated.
\end{prop}
\begin{proof} Let $x$ be an element in $\Omega_L$.
Assume that $s$ is a non-zero vector of $E_{L_x}$ which is written  as $s=a_1e_1+\cdots+a_re_r$ with $(a_1,\ldots,a_r)\in L_x^r$. Then one has
\[
\|s\|_x \leqslant \begin{cases}
\displaystyle\max_{i\in\{1,\ldots,r\}}(|a_i|_x\cdot\|e_i\|_x), & \text{if $x \in \Omega_L\setminus\Omega_{L,\infty}$},\vspace{5pt}\\[1ex]
\displaystyle \sum_{i=1}^r|a_i|_x\cdot\|e_i\|_x, & \text{if $x \in
\Omega_{L,\infty}$}.
\end{cases}
\]
Thus $\ln\|s\|_x\leqslant\ln\|s\|_{\boldsymbol{e},x}+\max_{i\in\{1,\ldots,r\}}\ln\norm{e_i}_x$.
Moreover, one can also interpret
\[\sup_{0\neq s\in E_{L_x}}\frac{\|s\|_{\boldsymbol{e},x}}{\|s\|_x}\]
as the operator norm of the $L_x$-linear map 
\[\mathrm{Id}_{E_{L_x}}:(E_{L_x},\|\ndot\|_x)\longrightarrow (E_{L_x},\|\ndot\|_{\boldsymbol{e},x}),\]
which is equal to the operator norm of the dual $L_x$-linear {map}
\[\mathrm{Id}_{E_{L_x}^\vee}:(E_{L_x}^\vee,\|\ndot\|_{\boldsymbol{e},x,*})\longrightarrow(E_{L_x}^\vee,\|\ndot\|_{x,*})\]
since the norm $\|\ndot\|_{\boldsymbol{e},x}$ is reflexive (see Proposition \ref{Pro:normdualopt}).
Let $\{e_i^\vee\}_{i=1}^r$ be the dual basis of $\boldsymbol{e}$. For any $\alpha=b_1e_1^\vee+\cdots+b_re_r^\vee\in E_{L_x}^\vee$, one has
\[
\|\alpha\|_{x,*} \leqslant
\begin{cases}
\displaystyle\max_{i\in\{1,\ldots,r\}}(|b_i|_x\cdot\|e_i^\vee\|_{x,*}), & \text{if $x \in \Omega_L\setminus\Omega_{L,\infty}$},\vspace{5pt}\\[1ex]
{\displaystyle\sum_{i=1}^r|b_i|_x\cdot\|e_i^\vee\|_{x,*}},
& \text{if $x \in \Omega_{L,\infty}$}.
\end{cases}
\]
Thus we obtain 
\[\ln\|\alpha\|_{x,*}-\ln\|\alpha\|_{\boldsymbol{e},x,*}\leqslant\max_{i\in\{1,\ldots,r\}}\ln\|e_i^\vee\|_{x,*}+\ln(r)\indic_{\Omega_{L,\infty}}(x).\]
Therefore, for any $s\in E_L$, one has
\begin{equation}\label{Equ:majnorme2}-\max_{i\in\{1,\ldots,r\}}\ln\norm{e_i}_x\leqslant \ln\|s\|_{\boldsymbol{e},x}-\ln\|s\|_x\leqslant\max_{i\in\{1,\ldots,r\}}\ln\|e_i^\vee\|_{x,*}+\ln(r)\indic_{\Omega_{L,\infty}}(x).\end{equation}

Note that, if $\omega=\pi_{L/K}(x)$, then one has (see Proposition \ref{Pro:extensiondecorps})
\[\forall\,i\in\{1,\ldots,r\},\quad\|e_i\|_x=\|e_i\|_{\omega,**}.\]
Moreover, if $\omega=\pi_{L/K}(x)$ belongs to $\Omega\setminus\Omega_\infty$, then
\[\|e_i^\vee\|_{x,*}=\|e_i^\vee\|_{\omega,L_x,\varepsilon,*}=\|e_i^\vee\|_{\omega,*,L_x,\varepsilon}=\|e_i^\vee\|_{\omega,*},\]
where the second equality comes from Proposition~\ref{Pro:comparisonofdualnormes:scalar:extension} \ref{Pro:comparisonofdualnormes:scalar:extension:epsilon} 
and the last one comes from Proposition \ref{Pro:extensiondecorps}. If $\omega=\pi_{L/K}(x)\in\Omega_\infty$, then
\[\|e_i^\vee\|_{x,*}=\|e_i^\vee\|_{\omega,L_x,\pi,*}=\|e_i^\vee\|_{\omega,*,L_x,\varepsilon}=\|e_i^\vee\|_{\omega,*},\] 
where the second equality comes from 
Proposition~\ref{Pro:comparisonofdualnormes:scalar:extension} \ref{Pro:comparisonofdualnormes:scalar:extension:epsilon:pi}, 
and the last one comes from Proposition \ref{Pro:extensiondecorps}.
Since the norm family $\xi$ is dominated, there exists a $\nu$-integrable function $A$ on $\Omega$ such that (see Remark \ref{Rem: dominancy of double dual})
\[\max_{i\in \{1,\ldots,r\}}\max\{\ln\|e_i\|_{\omega,**},\ln\|e_i^\vee\|_{\omega,*}\}\leqslant A(\omega)\quad\text{$\nu$-almost everywhere.}\]
Therefore, by  \eqref{Equ:majnorme2}, we obtain
\[d_x(\xi_L,\xi_{\boldsymbol{e},L})\leqslant A_{\boldsymbol{e}}(\pi_{L/K}(x))\quad\text{$\nu$-{almost everywhere}},\]
with 
\[\forall\,\omega\in\Omega,\quad A_{\boldsymbol{e}}(\omega):=A(\omega)+\ln(r)\indic_{\Omega_\infty}(\omega).\]
Note that the function $A_{\boldsymbol{e}}$ is $\nu$-integrable on $(\Omega,\mathcal A)$. The proposition is thus proved. 
\end{proof}

\begin{coro}
\label{Cor: distance dominant}
Let $E$ be a vector space of finite rank over $K$, $\xi_1$ and $\xi_2$ be norm families in $\mathcal N_E$ which are dominated and ultrametric on $\Omega\setminus\Omega_\infty$. Then the local distance function $(\omega\in\Omega)\mapsto d_\omega(\xi_1,\xi_2)$ is $\nu$-dominated.
\end{coro}
\begin{proof}
Let $\boldsymbol{e}$ be a basis of $E$ over $K$. By Proposition \ref{Pro:comparaisonavecbasenorm}, the local distance functions $(\omega\in\Omega)\mapsto d_\omega(\xi_1,\xi_{\boldsymbol{e}})$ and $(\omega\in\Omega)\mapsto d_\omega(\xi_2,\xi_{\boldsymbol{e}})$ are both $\nu$-dominated. Since for any $\omega\in\Omega$ one has
\[d_\omega(\xi_1,\xi_2)\leqslant d_\omega(\xi_1,\xi_{\boldsymbol{e}})+d_\omega(\xi_2,\xi_{\boldsymbol{e}}),\] by Propositions \ref{Pro: preserve order int sup and inf} and \ref{Pro:subsuperadditive}
the function $(\omega\in\Omega)\mapsto d_\omega(\xi_1,\xi_2)$ is $\nu$-dominated.
\end{proof}

{
\begin{rema}
The assertion of Corollary \ref{Cor: distance dominant} does not necessarily hold without the condition that the norm families are ultrametric on $\Omega\setminus\Omega_\infty$. Consider the following counter-example. Let $K$ be an infinite field, $\mathcal A$ be the discrete $\sigma$-algebra on $K$ and $\nu$ be the atomic measure on $K$ such that $\nu(\{a\})=1$ for any $a\in K$. For any $a\in K$, let $|\ndot|_a$ be the trivial absolute value on $K$. We consider the adelic curve $S=(K,(K,\mathcal A,\nu),\{|\ndot|_a\}_{a\in K})$. For any $a\in K$, let $\norm{\ndot}_a$ be the norm on $K^2$ such that 
\[\norm{(x,y)}_a=\begin{cases}
0,&\text{if $x=y=0$},\\
2,&\text{if $y=ax$, $x\neq 0$},\\
1,&\text{else}.
\end{cases}\]
Note that for any $s\in K^2\setminus\{(0,0)\}$, one has $\norm{s}_a=1$ for all except at most one $a\in K$. Therefore the norm family $\xi=\{\norm{\ndot}_a\}_{a\in K}$ is upper dominated. If we identify $K^2$ with the dual vector space of itself in the canonical way, then for any $a\in K$ one has 
\[\forall\,(x,y)\in K^2,\quad\norm{(x,y)}_{a,*}=\begin{cases}
0,&\text{if $x=y=0$},\\
1,&\text{else}.
\end{cases}\]
Hence the dual norm family $\xi^\vee$ is also upper dominated. Now let $\boldsymbol{e}=\{(1,0),(0,1)\}$ be the canonical basis of $K^2$. For any $a\in K$ one has
\[\forall\,(x,y)\in K^2,\quad\norm{(x,y)}_{\boldsymbol{e},a}=\begin{cases}
0,&\text{if $x=y=0$},\\
1,&\text{else}.
\end{cases}\]
Therefore one has $d_a(\xi,\xi_{\boldsymbol{e}})=\ln(2)$ for any $a\in K$. Since $K$ is an infinite set, the local distance function $(a\in K)\mapsto d_a(\xi,\xi_{\boldsymbol{e}})$ is clearly not upper dominated.
\end{rema}}

\begin{coro}\label{Cor:dominatedanddist}
Let $E$ be a finite-dimensional vector space over $K$ and $\xi=\{\|\ndot\|_\omega\}_{\omega\in\Omega}$ be a norm family in $\mathcal N_E$. The following assertions are equivalent:
\begin{enumerate}[label=\rm(\arabic*)]
\item\label{Item: xi is dominated Cor:dominatedanddist} the norm family $\xi$ is dominated and the local distance function $(\omega\in\Omega)\mapsto d_\omega(\xi,\xi^{\vee\vee})$ is $\nu$-dominated;
\item\label{Item: local distance to any basis is dominated} for any basis $\boldsymbol{e}$ of $E$, the local distance function $(\omega\in
\Omega)\mapsto d_\omega(\xi,\xi_{\boldsymbol{e}})$ is $\nu$-dominated;
\item\label{Item: local distance to some basis is dominated} there exists a basis $\boldsymbol{e}$ of $E$ such that the local distance function $(\omega\in\Omega)\mapsto d_\omega(\xi,\xi_{\boldsymbol{e}})$ is $\nu$-dominated.
\end{enumerate}
\end{coro}
\begin{proof}
``\ref{Item: xi is dominated Cor:dominatedanddist}$\Longrightarrow$\ref{Item: local distance to any basis is dominated}'': Note that the norm family $\xi^{\vee\vee}$ is ultrametric on $\Omega\setminus\Omega_\infty$. Moreover it is dominated since $\xi$ is dominated (see Remark \ref{Rem: dominancy of double dual}). By Proposition \ref{Pro:comparaisonavecbasenorm}, we obtain that, for any basis $\boldsymbol{e}$ of $E$, the local distance function $(\omega\in\Omega)\mapsto d_\omega(\xi^{\vee\vee},\xi_{\boldsymbol{e}})$ is $\nu$-dominated. By the assumption that the function $(\omega\in\Omega)\mapsto d_\omega(\xi,\xi^{\vee\vee})$ is $\nu$-dominated, we deduce that the function $(\omega\in\Omega)\mapsto d_\omega(\xi,\xi_{\boldsymbol{e}})$ is also $\nu$-dominated. 

``\ref{Item: local distance to any basis is dominated}$\Longrightarrow$\ref{Item: local distance to some basis is dominated}'' is trivial.

``\ref{Item: local distance to some basis is dominated}$\Longrightarrow$\ref{Item: xi is dominated Cor:dominatedanddist}'': By  Proposition \ref{Pro:adelicvectorbundledist}, the norm family $\xi$ is $\nu$-dominated. By Proposition \ref{Pro:comparaisonavecbasenorm}, the function $(\omega\in\Omega)\mapsto d_\omega(\xi^{\vee\vee},\xi_{\boldsymbol{e}})$ is $\nu$-dominated. Since $d_\omega(\xi,\xi^{\vee\vee})\leqslant d_\omega(\xi,\xi_{\boldsymbol{e}})+d_\omega(\xi^{\vee\vee},\xi_{\boldsymbol{e}})$, we obtain that the function $(\omega\in\Omega)\mapsto d_\omega(\xi,\xi^{\vee\vee}) $ is $\nu$-dominated. 
\end{proof}

\begin{defi}\label{Def: strongly dominated}
Let $E$ be a finite-dimensional vector space over $K$ and $\xi$ {be} a norm family on $E$. We say that $\xi$ is \emph{strongly dominated}\index{strongly dominated}\index{norm family!strongly dominated} if it is dominated and if the function $(\omega\in\Omega)\mapsto d_\omega(\xi,\xi^{\vee\vee})$ is $\nu$-dominated (or equivalently, $(E,\xi)$ satisfies any of the assertions in Corollary \ref{Cor:dominatedanddist}). 
\end{defi}

\begin{rema}\label{Rem: strongly dominated} Let $E$ be a finite-dimensional vector space over $K$ and $\xi$ {be} a norm family on $E$.
If the norm family $\xi$ is ultrametric on $\Omega\setminus\Omega_\infty$, then the function $(\omega\in\Omega)\mapsto d_\omega(\xi,\xi^{\vee\vee})$ is identically zero (and hence is $\nu$-dominated). Therefore, in this case $\xi$  is dominated if and only if it is strongly dominated. In particular, in the case where $E$ is of dimension $1$ over $K$, the norm family $\xi$ is dominated if and only if it is strongly dominated. Moreover, if $\xi$ is a dominated norm family on a finite-dimensional vector space, then the dual norm family $\xi$ is strongly dominated since it is dominated (see Remark \ref{Rem: dominancy of double dual}) and ultrametric on $\Omega\setminus\Omega_\infty$.
\end{rema}

\begin{coro}\label{Cor:dominanceext}
Let $E$ be a vector space of finite rank over $K$ and $\xi=\{\|\ndot\|_{\omega}\}_{\omega\in\Omega}$ be a norm family in $\mathcal N_E$. If $\xi$ is dominated, then for any algebraic extension $L/K$ the norm family $\xi_L$ is strongly dominated. Conversely, if there exists an algebraic extension $L/K$ such that the norm family $\xi_L$ is dominated and if the function $(\omega\in\Omega)\mapsto d_\omega(\xi,\xi^{\vee\vee})$ is $\nu$-dominated, then the norm family $\xi$ is also dominated. 
\end{coro}
\begin{proof} Assume that the norm family $\xi$ is dominated. 
By Proposition \ref{Pro:comparaisonavecbasenorm}, for any basis $\boldsymbol{e}$ of $E$, the local distance function $(x\in\Omega_L)\mapsto d_x(\xi_{L},\xi_{\boldsymbol{e},L})$ is $\nu_L$-dominated. Therefore, by Proposition \ref{Pro:adelicvectorbundledist} we obtain that the norm family $\xi_{L}$ is strongly dominated.

Conversely, we assume that $L/K$ is an algebraic extension and $\xi_L$ is dominated. Let $\boldsymbol{e}$ be a basis of $E$ over $K$. Since $\xi_{L}$ is dominated and since the norms in the family $\xi_L$ corresponding to non-Archimedean absolute values are ultrametric, the function $f:(x\in\Omega_L)\mapsto d_x(\xi_{L},\xi_{\boldsymbol{e},L})$ is $\nu_L$-dominated. Moreover, by Proposition \ref{Pro:distanceextension}, one has $f=g\circ \pi_{L/K}$, where $g$ sends $\omega\in\Omega$ to $d_\omega(\xi^{\vee\vee},\xi_{\boldsymbol{e}})$ (one has $\xi_{\boldsymbol{e}}=\xi_{\boldsymbol{e}}^{\vee\vee}$ since it is ultrametric on $\Omega\setminus\Omega_\infty$). Since the function $f$ is $\nu_L$-dominated, there exists an $\nu_L$-integrable function $A$ on $\Omega_L$ such that $|f|\leqslant A$ almost everywhere. By Proposition \ref{Pro:fibreintegralint}, the function $I_{L/K}(A)$ is $\nu$-integrable. Moreover, one has $I_{L/K}(|f|)=|g|$ since $|f|=|g|\circ\pi_{L/K}$ (see Remark \ref{Rem:operateurILK}). Therefore, the function $g$ is $\nu$-dominated by the $\nu$-integrable function $I_{L/K}(A)$, which implies that the norm family $\xi^{\vee\vee}$ is dominated. Finally, by Proposition \ref{Pro:adelicvectorbundledist} and the assumption that the function $(\omega\in\Omega)\mapsto d_{\omega}(\xi,\xi^{\vee\vee})$ is $\nu$-dominated, we obtain that the norm family $\xi$ is dominated.
\end{proof}

\begin{rema}\label{Rem:existencehermitindomnormf} Let $E$ be a finite-dimensional vector space over $K$.
Corollary \ref{Cor:dominatedanddist} implies that there exist Hermitian norm families  on $E$ which are dominated. In fact, let $\boldsymbol{e}=\{e_i\}_{i=1}^r$ be a basis of $E$ over $K$. Consider the following norm family $\xi=\{\|\ndot\|_\omega\}_{\omega\in\Omega}$ with
\[\forall\,(a_1,\ldots,a_r)\in K_\omega^r,\quad\|a_1e_1+\cdots+a_re_r\|_\omega=\begin{cases}
\max_{i\in\{1,\ldots,r\}}|a_i|_\omega,&\text{if $\omega\in\Omega\setminus\Omega_\infty$,}\\
\big(\sum_{i=1}^r|a_i|_\omega^2\big)^{1/2},&\text{if $\omega\in\Omega_\infty$.}
\end{cases}\] It is a Hermitian norm family on $E$.
Note that one has $d_\omega(\xi,\xi_{\boldsymbol{e}})\leqslant\frac 12\ln(r)\indic_{\Omega_\infty}(\omega)$. Therefore $(\omega\in\Omega)\mapsto d_\omega(\xi,\xi_{\boldsymbol{e}})$ is a $\nu$-dominated function on $\Omega$. By Corollary \ref{Cor:dominatedanddist}, we obtain that $\xi$ is a dominated norm family.

Let $\xi$ be a Hermitian norm family on $E$. If $\xi$ is dominated, then for any algebraic extension $L/K$, the norm family $\xi^H_L$ is dominated. In fact, by Corollary \ref{Cor:dominatedanddist}, the norm family $\xi_L$ is dominated. By Proposition \ref{Pro:extensiondecorps} \ref{item: comparison}, the local distance function $(x\in \Omega_{L})\mapsto d_x(\xi_L,\xi_L^H)$ is bounded from above by $\frac 12\ln(2)\indic_{\Omega_{L,\infty}}$. By Proposition \ref{Pro:adelicvectorbundledist}, we obtain that the norm family $\xi_L^H$ is dominated. 
\end{rema}

\begin{prop}\label{Pro:normfamilyfromabasis}
Let $E$ be a finite-dimensional vector space over $K$ and $L/K$ be an algebraic extension of fields. For any $\omega\in\Omega$, we fix an extension $|\ndot|_{L,\omega}$ on $L$ of the absolute value $|\ndot|_\omega$ and denote by $L_\omega$ the completion of $L$ with respect to the extended absolute value. Let $\boldsymbol{e}=\{e_i\}_{i=1}^r$ be a basis of $E\otimes_KL$. For any $\omega\in\Omega$, let $\|\ndot\|_{\omega}'$ be the norm on $E\otimes_KL_\omega$ defined as
\[\forall\,(\lambda_1,\ldots,\lambda_r)\in L_\omega^r,\quad
\|\lambda_1e_1+\cdots+\lambda_re_r\|_{\omega}'=\max_{i\in\{1,\ldots,r\}}|\lambda_i|_{L,\omega}\] 
and let $\|\ndot\|_\omega$ be the restriction of $\|\ndot\|_\omega'$ {to} $E\otimes_KK_\omega$. Then the norm family $\xi=\{\|\ndot\|_\omega\}_{\omega\in\Omega}$ in $\mathcal N_E$ is strongly dominated.
\end{prop}
\begin{proof}
We first prove that, for any element $b\in L$, the function \[(\omega\in\Omega)\longmapsto\ln|b|_{L,\omega}\]
is bounded from above by a $\nu$-integrable function. Let 
\[F(X)=X^n+a_{n-1}X^{n-1}+\cdots+a_0\in K[X]\]
be the minimal polynomial
of $b$. By the same argument as in the proof of Theorem \ref{Thm:constructionofextension} \ref{Item: Adelic curve base change}, we obtain that 
\[\ln|b|_\omega\leqslant \indic_{\Omega_\infty}(\omega)\ln(n)+\max\{0,\ln|a_0|_\omega,\ldots,\ln|a_{n-1}|_{\omega}\}.\]
By Proposition \ref{Pro:mesurabilitekappa}, the function $\indic_{\Omega_{\infty}}$ is $\nu$-integrable. Moreover, by the definition of adelic curves, for any $i\in\{0,\ldots,n-1\}$ such that $a_i\neq 0$, the function $(\omega\in\Omega)\mapsto\ln|a_i|_\omega$ is also $\nu$-integrable, we thus obtain the assertion. 

Let $\boldsymbol{f}=\{f_i\}_{i=1}^r$ be a basis of $E$ over $K$ and $A=(a_{ij})_{(i,j)\in\{1,\ldots,r\}^2}\in M_{r\times r}(L)$ be the transition matrix between $\boldsymbol{e}$ and $\boldsymbol{f}$, namely
\[\forall\,i\in\{1,\ldots,r\},\quad f_i=\sum_{j=1}^r a_{ij}e_j.\]
Let $(b_{ij})_{(i,j)\in\{1,\ldots,r\}^2}\in M_{r\times r}(L)$ be the inverse matrix of $A$. Then one has
\[\forall\,i\in\{1,\ldots,r\},\quad e_i=\sum_{j=1}^r b_{ij}f_j.\]
By the above assertion, there exists a $\nu$-integrable function $g$ on $\Omega$ such that
\[\forall\,\omega\in\Omega,\; \max_{(i,j)\in\{1,\ldots,r\}^2}\max\{\ln|a_{ij}|_\omega,\ln|b_{ij}|_\omega\}\leqslant g(\omega).\]

We will prove that the local distance function $d(\xi,\xi_{\boldsymbol{f}})$ is $\nu$-dominated. Let $\omega\in\Omega$ and $x=\lambda_1f_1+\cdots+\lambda_rf_r$ be an element of $E\otimes_{K}K_\omega$. One has
\[x=\sum_{i=1}^r\lambda_i\sum_{j=1}^r a_{ij}e_j=\sum_{j=1}^r\bigg(\sum_{i=1}^ra_{ij}\lambda_i\bigg)e_j.\]
Therefore 
\[\begin{split}\ln\|x\|_\omega&=\max_{j\in\{1,\ldots,r\}}\ln\bigg|\sum_{i=1}^ra_{ij}\lambda_i\bigg|_{L,\omega}\leqslant \max_{i\in\{1,\ldots,r\}}\ln|\lambda_i|_\omega+g(\omega)+\ln(r)\indic_{\Omega_\infty}{(\omega)}\\
&\leqslant\ln\|x\|_{\boldsymbol{f},\omega}+g(\omega)+\ln(r)\indic_{\Omega_\infty}{(\omega)}.
\end{split}\]
Similarly, if we write $x$ as $x=\mu_1e_1+\cdots+\mu_re_r$, with $(\mu_1,\ldots,\mu_r)\in L_\omega^r$, one has
\[x=\sum_{i=1}^r\mu_i\sum_{j=1}^rb_{ij} f_j=\sum_{j=1}^r\bigg(\sum_{i=1}^rb_{ij}\mu_i\bigg)f_j.\]
Namely $\lambda_j=\sum_{i=1}^rb_{ij}\mu_i$ for any $j\in\{1,\ldots,r\}$. If $\omega\in\Omega\setminus\Omega_\infty$ then
\[\begin{split}\ln\|x\|_{\boldsymbol{f},\omega}&=\ln\Big(\max_{j\in\{1,\ldots,r\}}|\lambda_j|_\omega\Big)=\ln\bigg(\max_{j\in\{1,\ldots,r\}}\bigg|\sum_{i=1}^rb_{ij}\mu_i\bigg|_{\omega}\bigg)\\
&\leqslant \ln\Big(\max_{i\in\{1,\ldots,r\}}|\mu_i|_{L,\omega}\Big)+g(\omega)=\ln\|x\|_{\omega}+g(\omega).
\end{split}\]
If $\omega\in\Omega_\infty$, then
\[\begin{split}\ln\|x\|_{\boldsymbol{f},\omega}&=\ln\Big(\sum_{j=1}^r|\lambda_j|_\omega\Big)=\ln\bigg(\sum_{j=1}^r\bigg|\sum_{i=1}^rb_{ij}\mu_i\bigg|_{\omega}\bigg)\\
&\leqslant \ln\Big(\max_{i\in\{1,\ldots,r\}}|\mu_i|_{L,\omega}\Big)+g(\omega)+\ln(r^2)=\ln\|x\|_{\omega}+g(\omega)+\ln(r^2).
\end{split}\]
Therefore, one has
\[\forall\,\omega\in\Omega,\quad d_\omega(\xi,\xi_{\boldsymbol{f}})\leqslant g(\omega)+2\ln(r)\indic_{\Omega_\infty}{(\omega)},\]
which implies that the local distance function $d(\xi,\xi_{\boldsymbol{f}})$ is $\nu$-dominated. By Corollary \ref{Cor:dominatedanddist}, we obtain that the norm family $\xi$ is dominated.
\end{proof}

The following proposition is a criterion of the dominance property in the case where the vector space is of dimension $1$.
\begin{prop}\label{Pro:line bundle}
Let $E$ be a vector space of rank $1$ over $K$ and $\xi=\{\|\ndot\|_\omega\}_{\omega\in\Omega}$ be a norm family on $E$. Then the following conditions are equivalent:
\begin{enumerate}[label=\rm(\arabic*)]
\item\label{Item: xi is dominated Pro:line bundle} the norm family $\xi$ is dominated;
\item\label{Item: for any s , function is dominated} for any non-zero element $s\in E$, the function $(\omega\in\Omega)\mapsto\ln\|s\|_\omega$ is $\nu$-dominated;
\item\label{Item: exists s function is dominated}  there exists a non-zero element $s\in E$ such that the function $(\omega\in\Omega)\mapsto\ln\|s\|_{\omega}$ is $\nu$-dominated.
\end{enumerate} 
\end{prop}
\begin{proof}
``\ref{Item: xi is dominated Pro:line bundle}$\Longrightarrow \ref{Item: for any s , function is dominated}\Longrightarrow\ref{Item: exists s function is dominated}$'' are trivial. In the following, we prove ``\ref{Item: exists s function is dominated}$\Longrightarrow$\ref{Item: xi is dominated Pro:line bundle}''. If
$s'$ is a non-zero element in $E$, then we can write it in the form $s'=as$, where $a$ is a non-zero element in $K$. Then one has
\[\forall\,\omega\in\Omega,\quad \ln\|s'\|_\omega=\ln|a|_\omega+\ln\|s\|_{\omega}=\ln|a|_\omega+\ln\|s\|_{\omega}.\]  
Since the function $(\omega\in\Omega)\mapsto \ln\|s\|_{\omega}$ is $\nu$-dominated and the function $(\omega\in\Omega)\mapsto \ln|a|_\omega$ is $\nu$-integrable, we obtain that the function $(\omega\in\Omega)\mapsto \ln\|s'\|_{\omega}$ is $\nu$-dominated. Moreover, if we denote by $s^\vee$ the dual element of $s$ in $E^\vee$, then one has
\begin{equation}\label{Equ:normdudual}\ln\|s^\vee\|_{\omega,*}=-\ln\|s\|_\omega\end{equation}
for any $\omega\in\Omega$. By the same argument as above we obtain that, for any non-zero element $\alpha\in E^\vee$, the function $(\omega\in\Omega)\mapsto\ln\|\alpha\|_{\omega,*}$ is $\nu$-dominated. Therefore the norm family $\xi$ is dominated.
\end{proof}

Let $K'$ be a finite extension field of $K$ and let 
\[
S_{K'} = S \otimes_K K'=(K',(\Omega_{K'},\mathcal A_{K'},\nu_{K'}),\phi_{K'})
\]
be the algebraic extension of $S$ by $K'$.
Let $E$ be a finite-dimensional vector space over $K$ and $E_{K'} =  E \otimes_K K'$.
Note that, for any $\omega\in\Omega$ and any $\omega'\in\Omega_{K'}$ such that $\pi_{K'/K}(\omega')=\omega$, the vector space $E\otimes_KK_\omega$ can be naturally considered as a $K_\omega$-vector subspace of $E_{K'}\otimes_{K'}K_{\omega'}'$.

\begin{prop}\label{prop:dominated:extensin:field}
Let $\xi = \{ \|\ndot\|_{\omega} \}_{\omega \in \Omega}$ and
$\xi' = \{ \|\ndot\|'_{\omega'} \}_{\omega' \in \Omega_{K'}}$ be norm familes of $E$ and $E_{K'}$, respectively, such that
\begin{equation}\label{Equ: restriction of norm prime is the initial norm}
\forall\, \omega \in \Omega,\ \forall\, \omega' \in \pi^{-1}_{K'/K}(\{\omega\}),\ 
\forall s \in E \otimes_K K_{\omega},\quad \| s \|_{\omega} = \| s \|'_{\omega'}.
\end{equation}
If $\xi'$ is dominated (resp. strongly dominated), then $\xi$ is also dominated (resp. strongly dominated).
\end{prop}
\begin{proof}
Assume that $\xi$ is dominated. Let $s$ be a non-zero element in $E$. By the assumption \eqref{Equ: restriction of norm prime is the initial norm} 
one has 
\[\upint_{\Omega}\ln\norm{s}_\omega\,\nu(\mathrm{d}\omega)=\upint_{\Omega_{K'}}\ln\norm{s}'_{\omega'}\,\nu_{K'}(\mathrm{d}\omega')<+\infty.\]
Hence the norm family $\xi$ is upper dominated. Let $\alpha$ be a non-zero element in $E^\vee$. For any $\omega\in\Omega$ and any $\omega'\in\pi_{K'/K}^{-1}(\{\omega\})$, one has 
\[\norm{\alpha}_{\omega,*}=\sup_{s\in (E\otimes_K{K_\omega})\setminus\{0\}}\frac{|\alpha(s)|_\omega}{\norm{s}_\omega}=\sup_{s\in (E\otimes_K{K_\omega})\setminus\{0\}}\frac{|\alpha(s)|_\omega}{\norm{s}'_{\omega'}}\leqslant\norm{\alpha}'_{\omega',*}.\]
Since $(\xi')^\vee$ is upper dominated, we deduce that $\xi^\vee$ is also upper dominated.

Assume that $\xi'$ is strongly dominated. Let $\pmb e = \{e_i\}_{i=1}^n$ be a basis of $E$. Then, for $\omega \in \Omega$ and $\omega' \in \pi_{K'/K}^{-1}(\omega)$,
\begin{align*}
d_{\omega'}(\xi', \xi_{\pmb e, K'}) & =
\sup_{ s' \in (E_{K'} \otimes_{K'} K'_{\omega'}) \setminus \{ 0 \}} \Big| \ln \|s'\|'_{\omega'} - \ln \| s'\|_{\pmb e, K', \omega'} \Big| \\
& \geqslant 
\sup_{ s \in (E \otimes_K K_{\omega}) \setminus \{ 0 \}} \Big| \ln \|s\|'_{\omega'} - \ln \| s\|_{\pmb e, K',\omega'} \Big| \\
& = \sup_{ s \in (E \otimes _K K_{\omega}) \setminus \{ 0 \}} \Big| \ln \|s\|_{\omega} - \ln \| s\|_{\pmb e,\omega} \Big| = d_{\omega}(\xi, \xi_{\pmb e}).
\end{align*}
By our assumption together with Corollary~\ref{Cor:dominatedanddist},
the function $(\omega' \in \Omega_{K'}) \mapsto d_{\omega'}(\xi', \xi_{\pmb e, K'})$ is $\nu_{K'}$-dominated, that is,
there is an integrable function $A'$ on $\Omega_{K'}$ such that
$d_{\omega'}(\xi', \xi_{\pmb e, K'}) \leqslant A'(\omega')$ for all $\omega' \in \Omega_{K'}$, so that the above estimate implies
that $d_{\omega}(\xi, \xi_{\pmb e}) \leqslant I_{K'/K}(A')(\omega)$ for all $\omega \in \Omega$.
By Proposition~\ref{Pro:fibreintegralint}, one has
\[
\int_{\Omega_{K'}} A'(\omega')\, \nu_{K'}(\mathrm{d} \omega') = \int_{\Omega} I_{K'/K}(A')(\omega)\, \nu(\mathrm{d}\omega)
\]
and hence $\xi$ is strongly dominated by Corollary~\ref{Cor:dominatedanddist} again.
\end{proof}

\begin{coro}\label{coro:dominated:pullback:dominated:orginal}
Let $f:X\rightarrow\Spec K$ be a geometrically reduced  
projective $K$-scheme and $L$ be an invertible $\mathcal O_X$-module.
Let $X_{K'} := X \times_{\Spec K} \Spec K'$ and $L_{K'} := L \otimes_K K'$.
For each $\omega \in \Omega$ and $\omega' \in \Omega_{K'}$, 
$X_{\omega}$, $L_{\omega}$, $X_{K', \omega'}$ and $L_{K', \omega'}$ are defined by
\[
\begin{cases}
X_{\omega} := X \times_{\Spec(K)} \Spec(K_{\omega}), & L_{\omega} := L \otimes_K K_{\omega},\\
X_{K', \omega'} := X_{K'} \times_{\Spec(K')} \Spec(K'_{\omega'}), & L_{K', \omega'} := L_{K'} \otimes_{K'} K'_{\omega'}.
\end{cases}
\]
Moreover, for each $\omega \in \Omega$ and $\omega' \in \pi_{K'/K}^{-1}(\omega)$, 
let $\varphi_{\omega}$ be a metric of $L_{\omega}$ on $X_{\omega}$, and
$\varphi_{K', \omega'}$ be the metric of $L_{K'}$ obtained by $\varphi_{\omega}$
by the extension of scalars (cf. Definition~\ref{def:metric:base:change}).
Let $\|\ndot\|_{\varphi_\omega}$ and $\|\ndot\|_{\varphi_{K', \omega'}}$
be the sup norms on $H^0(X_{\omega},L_{\omega})$ and $H^0(X_{K', \omega'},L_{K', \omega'})$
obtained the metrics $\varphi_\omega$ and $\varphi_{K', \omega'}$, respectively. 
If $\xi_{K'} = \big\{ \|\ndot\|_{\varphi_{K', \omega'}} \big\}_{\omega' \in \Omega_{K'}}$ on $H^0(X_{K'}, L_{K'})$
is dominated, then $\xi=\{\|\ndot\|_{\varphi_\omega}\}_{\omega\in\Omega}$ on $H^0(X, L)$ is also dominated.
\end{coro}

\begin{proof}
For $\omega \in \Omega$, $\omega' \in \pi_{K'/K}^{-1}(\omega)$ and $s \in H(X_{\omega}, L_{\omega})$,
one has $\| s \|_{\varphi_{K', \omega'}} = \|s \|_{\varphi_\omega}$ (see Proposition \ref{prop:cont:inv:pullback:base:change}),
so that
the assertion follows from Proposition~\ref{prop:dominated:extensin:field}.
\end{proof}

The following proposition shows that the dominance property is actually preserved by most of the algebraic constructions on norm families.

\begin{prop}\phantomsection\label{Pro:dominancealgebraic}
\begin{enumerate}[label=\rm(\arabic*)]
\item\label{Item:subdom} Let $E$ be a finite-dimensional vector space over $K$ and $\xi$ be a dominated (resp. strongly dominated) norm family on $E$. The the restriction of $\xi$ {to} any vector subspace of $E$ is a dominated {(resp. strongly dominated)} norm family.
\item\label{Item:quotdom}Let $E$ be a finite-dimensional vector space over $K$ and $\xi$ be a dominated {(resp. strongly dominated)} norm family on $E$. Then the quotient norm family of $\xi$ on any quotient vector space of $E$ is a dominated (resp. strongly dominated) norm family.
\item\label{Item:dualdom}Let $E$ be a finite-dimensional vector space over $K$ and $\xi$ be an element in $\mathcal N_E$. If $\xi$ is dominated, then the norm family $\xi^\vee$ is strongly dominated.
\item\label{Item:sommedom}Let $E$ and $F$ be finite-dimensional vector spaces over $K$, and $\xi_E$ and $\xi_F$ be elements in $\mathcal N_E$ and $\mathcal N_F$, respectively. Let $\psi:\Omega\rightarrow\mathscr S$ be a map
such that $\psi=\psi_0$ outside of a measurable subset $\Omega'$ of $\Omega$ with $\nu(\Omega')<+\infty$, where $\psi_0$ denotes the function in $\mathscr S$ sending $t\in[0,1]$ to $\max\{t,1-t\}$. If both norm families $\xi_E$ and $\xi_F$ are dominated (resp. strongly dominated), then  the $\psi$-direct sum $\xi_E\oplus_{\psi}\xi_F$ is also dominated (resp. strongly dominated).
\item\label{Item:tensordom}Let $E$ and $F$ be finite-dimensional vector spaces over $K$, and $\xi_E$ and $\xi_F$ be elements in $\mathcal N_E$ and $\mathcal N_F$, respectively. Assume that both norm families $\xi_E$ and $\xi_F$ are dominated. Then the $\varepsilon$-tensor product $\xi_E\otimes_{\varepsilon}\xi_F$ and the  $\varepsilon,\pi$-tensor product $\xi_E\otimes_{\varepsilon,\pi}\xi_F$ are strongly dominated. If in addition both norm families $\xi_E$ and $\xi_F$ are Hermitian, then the orthogonal tensor product $\xi_E\otimes\xi_F$ is strongly dominated.
\item\label{Item:detdom}Let $E$ be a finite-dimensional vector space over $K$ and $\xi$ be an element in $\mathcal N_E$. Assume that $\xi$ is dominated. Then, for any $i\in\mathbb N$, the exterior power norm family $\Lambda^i\xi$ is strongly dominated. In particular, the determinant norm family $\det(\xi)$ is strongly dominated.
\end{enumerate}
\end{prop}
\begin{proof}
{\ref{Item:subdom} and \ref{Item:quotdom} in the dominated case: We first show the following claim: if $\xi=\{\norm{\ndot}_\omega\}_{\omega\in\Omega}$ is an upper dominated norm family, then all its restrictions and quotients are also upper dominated. Let $F$ be a vector space of $E$ and $\xi_F=\{\norm{\ndot}_\omega\}_{\omega\in\Omega}$ be the restriction of $\xi$ {to} $F$. For any $s\in F\setminus\{0\}$ and any $\omega\in\Omega$ one has $\norm{s}_{F,\omega}=\norm{s}_{\omega}$. Since the norm family $\xi$ is upper dominated, the function $(\omega\in\Omega)\mapsto \ln\norm{s}_{F,\omega}$ is upper dominated. Let $G$ be a quotient vector space of $E$ and $\xi_G=\{\norm{\ndot}_{G,\omega}\}_{\omega\in\Omega}$ be the quotient of $\xi$ on $G$. For any $t\in G\setminus\{0\}$ and any $s\in E$ which represents the class $t$ in $G$, one has $\norm{t}_{G,\omega}\leqslant\norm{s}_{\omega}$ for any $\omega\in\Omega$. Therefore the function $(\omega\in\Omega)\mapsto\ln\norm{t}_{G,\omega}$ is upper dominated.

Let $\xi=\{\norm{\ndot}_{\omega}\}_{\omega\in\Omega}$ be a norm family on $E$ such that the dual norm family $\xi^\vee$ is upper dominated. Let $G$ be a quotient vector space of $E$. We identify $G^\vee$ with a vector subspace of $E^\vee$. By Proposition \ref{Pro:dualquotient}, if $\xi_G$ denotes the quotient norm family of $\xi$ on $G$, then $\xi_G^\vee$ identifies with the restriction of $\xi^\vee$ {to} $G^\vee$. By the claim proved above, we obtain that $\xi_G^\vee$ is upper dominated. Similarly, if $F$ is a vector subspace of $E$ and if $\xi_F=\{\norm{\ndot}_{F,\omega}\}_{\omega\in\Omega}$ is the restriction of $\xi$ {to} $F$, then, for any $\omega\in\Omega$, $\norm{\ndot}_{F,\omega,*}$ is bounded from above by the quotient of the norm $\norm{\ndot}_{\omega,*}$ on $F^\vee$ (viewed as a quotient vector space of $E^\vee$). Therefore, by the claim proved above, we obtain that the norm family $\xi_F^\vee$ is upper dominated.

\ref{Item:subdom} in the strongly dominated case: Let $F$ be a vector subspace of $E$.
Let $\boldsymbol{f}$ be a basis of $F$. We complete it into a basis $\boldsymbol{e}$ of $E$. For any $\omega\in\Omega$ one has $d_\omega(\xi_F,\xi_{\boldsymbol{f}})\leqslant d_\omega(\xi,\xi_{\boldsymbol{e}})$. Since the function $(\omega\in\Omega)\mapsto d_\omega(\xi,\xi_{\boldsymbol{e}})$ is $\nu$-dominated, also is the function $(\omega\in\Omega)\mapsto d_\omega(\xi_F,\xi_{\boldsymbol{f}})$. By Corollary \ref{Cor:dominatedanddist}, we obtain that the norm family $\xi_F$ is {strongly} dominated.

\ref{Item:quotdom} in the strongly dominated case: Let $\boldsymbol{g}=\{g_i\}_{i=1}^m$ 
be a basis of $G$. For any $i\in\{1,\ldots,m\}$, we choose a vector $e_i$ in $E$ such that the canonical image of $e_i$ in $G$ is $g_i$. We complete the family $\{e_i\}_{i=1}^m$ into a basis $\boldsymbol{e}$ of $E$. Then for any $\omega\in\Omega$ one has $d_\omega(\xi_G,\xi_{\boldsymbol{g}})\leqslant d_\omega(\xi,\xi_{\boldsymbol{e}})$. Therefore, the function $(\omega\in\Omega)\mapsto d_\omega(\xi_G,\xi_{\boldsymbol{g}})$ is $\nu$-dominated, which implies that $\xi_G$ is {strongly} dominated.
}

\ref{Item:dualdom} has already been shown in Remark \ref{Rem: dominancy of double dual}, see also Remark \ref{Rem: strongly dominated} for the strong dominancy.

\ref{Item:sommedom} {in the dominated case: We first show the following claim: if both norm families $\xi_E$ and $\xi_F$ are upper dominated, then also is the direct sum $\xi_E\oplus_\psi\xi_F$. In fact, if $(s,t)$ is an element in $E\oplus F$, for $\omega\in\Omega\setminus\Omega'$ one has
\[\|(s,t)\|_\omega=\max\{\|s\|_{E,\omega},\|t\|_{F,\omega}\},\]
and for $\omega\in\Omega'$, one has
\[\|(s,t)\|_\omega\leqslant \|s\|_{E,\omega}+\|t\|_{F,\omega}\leqslant 2\max\{\|s\|_{E,\omega},\|t\|_{F,\omega}\},\]
where $\|\ndot\|_\omega$ denotes the norm indexed by $\omega$ in $\xi_E\oplus_\psi\xi_F$. Therefore the function $(\omega\in\Omega)\mapsto\norm{(s,t)}_{\omega}$ is upper dominated.

Let $\psi'$ be the map from $\Omega$ to $\mathscr S$ sending any $\omega\in\Omega\setminus\Omega_\infty$ to $\psi_0$ and any $\omega\in\Omega_\infty$ to $\psi(\omega)_*$ (see Definition \ref{Def:dualdirecsum}). By Proposition \ref{Pro:dualdirectsum}
we obtain that $(\xi_E\oplus_\psi\xi_F)^\vee$ identifies with $\xi_E^\vee\oplus_{\psi'}\xi_F^\vee$. By the claim proved above, if $\xi_E^\vee$ and $\xi_F^\vee$ are upper dominated, then also is $(\xi_E\oplus_{\psi}\xi_F)^\vee$.

\ref{Item:sommedom} in the strongly dominated case:} Let $\boldsymbol{e}'$ and $\boldsymbol{e}''$ be bases of $E$ and $F$ respectively, and $\boldsymbol{e}$ be the disjoint union of $\boldsymbol{e}'$ and $\boldsymbol{e}''$, viewed as a basis of $E\oplus F$. Since $\xi_E$ and $\xi_F$ are both dominated, by Corollary \ref{Cor:dominatedanddist} there exist $\nu$-integrable functions $A'$ and $A''$ such that 
\[d_\omega(\xi_E,\xi_{\boldsymbol{e}'})\leqslant A'(\omega),\quad d_\omega(\xi_F,\xi_{\boldsymbol{e}''})\leqslant A''(\omega)\quad\text{$\nu$-almost everywhere}.
\]
Moreover, if $(s,t)$ is an element in $E\oplus F$, for $\omega\in\Omega\setminus\Omega'$ one has
\[\|(s,t)\|_\omega=\max\{\|s\|_{E,\omega},\|t\|_{F,\omega}\},\]
and for $\omega\in\Omega'$, one has
\[\max\{\|s\|_{E,\omega},\|t\|_{F,\omega}\}\leqslant\|(s,t)\|_\omega\leqslant \|s\|_{E,\omega}+\|t\|_{F,\omega}\leqslant 2\max\{\|s\|_{E,\omega},\|t\|_{F,\omega}\},\]
where $\|\ndot\|_\omega$ denotes the norm indexed by $\omega$ in $\xi_E\oplus_\psi\xi_F$. Therefore
\[d_\omega(\xi_E\oplus_\psi\xi_F,\xi_{\boldsymbol{e}})\leqslant \max\{A'(\omega),A''(\omega)\}+\ln(2)\indic_{\Omega'}(\omega)\quad\text{$\nu$-almost everywhere}.\]
Note that the function $(\omega\in\Omega)\mapsto \max\{A'(\omega),A''(\omega)\}+\ln(2)\indic_{\Omega'}(\omega)$ is $\nu$-integrable. Hence the norm family $\xi_E\oplus_\psi\xi_F$ is strongly dominated (by Corollary \ref{Cor:dominatedanddist}).

\ref{Item:tensordom} By \ref{Item:dualdom}, the norm families $\xi_E^{\vee\vee}$ and $\xi_F^{\vee\vee}$ are both dominated. Therefore, without loss of generality, we may assume that $\|\ndot\|_{E,\omega}$ and $\|\ndot\|_{F,\omega}$ are ultrametric norms for $\omega\in\Omega\setminus\Omega_\infty$ (see Definition \ref{Def:tensorproducts}
, see also Proposition \ref{Pro:doubedualandquotient}). Let $\boldsymbol{e}=\{e_i\}_{i=1}^n$ and $\boldsymbol{f}=\{f_j\}_{j=1}^m$ be bases of $E$ and $F$ over $K$, and let $\boldsymbol{e}\otimes\boldsymbol{f}=\{e_i\otimes f_j\}_{(i,j)\in\{1,\ldots,n\}\times\{1,\ldots,m\}}$. Note that $\boldsymbol{e}\otimes\boldsymbol{f}$ is a basis of $E\otimes F$. Moreover, for any $\omega\in\Omega$, the norm $\|\ndot\|_{\boldsymbol{e}\otimes\boldsymbol{f},\omega}$ identifies with the $\varepsilon$-tensor product of the norms $\|\ndot\|_{E,\omega}$ and $\|\ndot\|_{F,\omega}$. Since the norm families $\xi_E$ and $\xi_F$ are dominated, there exist $\nu$-integrable functions $A_E$ and $A_F$ on $\Omega$ such that\begin{equation}\label{Equ:majorationdedistance}d_\omega(\xi_E,\xi_{\boldsymbol{e}})=\sup_{0\neq s\in E}\Big|\ln\|s\|_{E,\omega}-\ln\|s\|_{\boldsymbol{e},\omega}\Big|\leqslant A_E(\omega)\quad\text{$\nu$-almost everywhere},\end{equation}
and
\[d_\omega(\xi_F,\xi_{\boldsymbol{f}})=\sup_{0\neq t\in F}\Big|\ln\|t\|_{F,\omega}-\ln\|t\|_{\boldsymbol{f},\omega}\Big|\leqslant A_F(\omega)\quad\text{$\nu$-almost everywhere}.\]
By \eqref{Equ:distdualnorm}, we obtain
\[d_\omega(\xi_E^\vee,\xi_{\boldsymbol{e}}^\vee)=\sup_{0\neq \alpha\in E^\vee}\Big|\ln\|\alpha\|_{E,\omega,*}-\ln\|\alpha\|_{\boldsymbol{e},\omega,*}\Big|\leqslant A_E(\omega)\quad\text{$\nu$-almost everywhere},\]
which implies (see \eqref{Equ:distnormbasis})
\[d_\omega(\xi_E^\vee,\xi_{\boldsymbol{e}^\vee})\leqslant A_E(\omega)+\ln(n)\indic_{\Omega_\infty}(\omega)\quad\text{$\nu$-almost everywhere}.\]
Therefore, for any $\omega\in\Omega$ and any non-zero tensor $\varphi\in\Hom_K(E^\vee,F)\cong E\otimes_KF$, one has (see Remark \ref{Rem:operateureps})
\[\begin{split}\Big|\ln\|\varphi\|_{\varepsilon,\omega}-\ln\|\varphi\|_{\boldsymbol{e}\otimes\boldsymbol{f},\omega}\Big|&=\bigg|\sup_{0\neq\alpha\in E^\vee}\ln\frac{\|\varphi(\alpha)\|_{F,\omega}}{\|\alpha\|_{E,\omega,*}}-\sup_{0\neq\alpha\in E^\vee}\ln\frac{\|\varphi(\alpha)\|_{\boldsymbol{f},\omega}}{\|\alpha\|_{\boldsymbol{e}^\vee,\omega}}\bigg|\\
&\leqslant A_E(\omega)+A_F(\omega)+\ln(n)\indic_{\Omega_\infty}(\omega)\end{split}\]
$\nu$-almost everywhere,
where $\|\ndot\|_{\varepsilon,\omega}$ denotes the $\varepsilon$-tensor product of $\|\ndot\|_{E,\omega}$ and $\|\ndot\|_{F,\omega}$. Therefore the $\varepsilon$-tensor product norm family $\xi_E\otimes_{\varepsilon}\xi_F$ is dominated (and hence is strongly dominated since it is ultrametric on $\Omega\setminus\Omega_\infty$, see Remark \ref{Rem: strongly dominated}). By using the fact that $\xi_E\otimes_{\varepsilon,\pi}\xi_F=(\xi_E^\vee\otimes_{\varepsilon}\xi_F^\vee)^\vee$ (see Corollary \ref{Cor: dual tensor product} and Proposition \ref{Pro:dualitypiepsilon} for the non-Archimedean and the Archimedean cases respectively), we deduce the dominance property of $\xi_E\otimes_{\varepsilon,\pi}\xi_F$ from the above result and the assertion \ref{Item:dualdom} of the Proposition.  Finally, by Propositions \ref{Pro:comparisonofHSandepsnorms} and \ref{Pro:adelicvectorbundledist}, we deduce that the orthogonal tensor product norm family $\xi_E\otimes\xi_F$ is also strongly dominated, provided that $\xi_E$ and $\xi_F$ are both Hermitian.

\ref{Item:detdom} is a direct consequence of \ref{Item:quotdom} and \ref{Item:tensordom} since $\Lambda^i\xi$ is a quotient norm family of the $i$-th $\varepsilon,\pi$-tensor power of $\xi$.
\end{proof}

\begin{rema}\label{Rem:metricDE}Let $E$ be a finite-dimensional vector space over $K$. We denote by $\mathcal D_E$ the subset of $\mathcal N_E$ of all strongly dominated norm families $\xi=\{\|\ndot\|_\omega\}_{\omega\in\Omega}$ on $E$. 
By Corollary \ref{Cor: distance dominant}, we obtain that, for any pair $(\xi,\xi')$ of norm families in $\mathcal D_E$, the local distance function $(\omega\in\Omega)\mapsto d_\omega(\xi,\xi')$ is $\nu$-dominated. This observation allows us to construct a function $\mathrm{dist}(\ndot,\ndot)$ on $\mathcal D_E\times\mathcal D_E$, defined as (see \S\ref{Sec:ulint} for the definition of the upper integral  
$\upint_{\Omega} h(\omega) \,\nu(d\omega)$) 
\[\mathrm{dist}(\xi,\xi'):=\upint_{\Omega} d_{\omega} (\xi,\xi') \,\nu(\mathrm{d}\omega).\]Clearly this function is symmetric with respect to its two variables, and verifies the triangle inequality, where the latter {assertion} follows from the triangle inequality of the local distance function and Proposition \ref{Pro:subsuperadditive}. Therefore, $\mathrm{dist}(\ndot,\ndot)$ is actually a pseudometric on $\mathcal D_E$. Moreover, for any pair $(\xi,\xi')$ of elements of $\mathcal D_E$, $\mathrm{dist}(\xi,\xi')=0$ if and only if $\xi_\omega=\xi'_\omega$ $\nu$-almost everywhere (see Proposition \ref{Pro:seminormed}). Therefore, the pseudometric $\mathrm{dist}(\ndot,\ndot)$ induces a metric on the quotient space of $\mathcal D_E$ modulo the equivalence relation
\[\xi\sim \xi'\Longleftrightarrow \xi_\omega=\xi_\omega'\quad\nu\text{-almost everywhere}.\]
This quotient metric space is actually complete. In fact, assume that $\{\xi_n\}_{n\in\mathbb N}$ is a Cauchy sequence in $\mathcal D_E$. Then we can pick a subsequence $\{\xi_{n_k}\}_{k\in\mathbb N}$ such that
\[\forall\,k\in\mathbb N,\quad \upint_{\Omega} \Big(\indic_{\{d(\xi_{n_k},\xi_{n_{k+1}})\geqslant 2^{-k}\}}\Big)\, \nu(d\omega) \leqslant 2^{-k}.\]
The set of $\omega\in\Omega$ such that $\{\xi_{n_k,\omega}\}_{k\in\mathbb N}$ is not a Cauchy sequence (with respect the the metric defined in \S\ref{Subsec:Distance}) is a $\nu$-negligible set. Let $\xi$ be a norm family such that $\{\xi_{n_k,\omega}\}_{k\in\mathbb N}$ converges to $\xi_\omega$ $\nu$-almost everywhere (see Remark \ref{Rem:localcompleteness} for the local completeness). Then, by the same argument as in the proof of Proposition \ref{Pro:completeness}, we obtain that $\mathrm{dist}(\xi_n,\xi)$ converges to $0$ when $n$ goes to the infinity. 
\end{rema}

\subsection{Measurability of norm families} Let $E$ be a vector space of finite rank over $K$ and $\xi=\{\|\ndot\|_\omega\}_{\omega\in\Omega}$ be a norm family in $\mathcal N_E$. We say that the norm family $\xi$ is $\mathcal A$-\emph{measurable}\index{measurable norm family}\index{norm family!measurable} (or simply \emph{measurable} when there is no ambiguity on the $\sigma$-algebra $\mathcal A$) if for any $s\in E$ the function $(\omega\in\Omega)\mapsto \|s\|_\omega$ is $\mathcal A$-measurable.  By definition, if the norm family $\xi$ is $\mathcal A$-measurable on $\Omega$, then also is its restriction {to} a vector subspace of $E$. The following proposition shows that measurable direct sums preserve the measurability of norm families.

\begin{prop}\label{Pro:mesurabilitedesomme}
Let $E$ and $F$ be finite-dimensional vector spaces over $K$ and $\xi_E=\{\|\ndot\|_{E,\omega}\}_{\omega\in\Omega}$, $\xi_F=\{\|\ndot\|_{F,\omega}\}_{\omega\in\Omega}$ be respectively norm families in $\mathcal N_E$ and $\mathcal N_F$, which are both $\mathcal A$-measurable. For any map $\psi:\Omega\rightarrow\mathscr S$ which is $\mathcal A$-measurable, where we consider the Borel $\sigma$-algebra on $\mathscr S$ induced by the topology of uniform convergence, the direct sum $\xi_E\oplus_\psi\xi_F=\{\|\ndot\|_{\psi,\omega}\}_{\omega\in\Omega}$ is $\mathcal A$-measurable.
\end{prop}
\begin{proof}
Consider the map 
$g:\mathscr S\times \intervalle[0{+\infty}[^2\rightarrow\mathbb R$
\[g(\eta,a,b)\longrightarrow\begin{cases}
0,&a+b=0,\\
(a+b)\eta(a/(a+b)),&a+b\neq 0.
\end{cases}\]
We claim that the map $g$ is continuous. Let $\{(\eta_n,a_n,b_n)\}_{n\in\mathbb N}$ be a sequence in $\mathscr S\times\intervalle[0{+\infty}[^2$ which converges to $(\eta,a,b)\in\mathscr S\times\intervalle[0{+\infty}{[}^2$. If $a+b\neq 0$, then $a_n/(a_n+b_n)$ converges to $a/(a+b)$, and therefore
\[\begin{split}&\quad\;\big|\eta_n(a_n/(a_n+b_n))-\eta(a/(a+b))\big|\\&\leqslant\big|\eta_n(a_n/(a_n+b_n))-\eta(a_n/(a_n+b_n))\big|+\big|\eta(a_n/(a_n+b_n))-\eta(a/(a+b))\big|\\
&\leqslant\|\eta_n-\eta\|_{\sup}+\big|\eta(a_n/(a_n+b_n))-\eta(a/(a+b))\big|
\end{split}\]  
converges to $0$ when $n$ tends to the infinity. We then deduce that \[\lim_{n\rightarrow+\infty}g(\eta_n,a_n,b_n)=g(\eta,a,b).\]
If $a+b=0$, then 
\[\lim_{n\rightarrow+\infty}g(\eta_n,a_n,b_n)=0=g(\eta,a,b)\]
since the sequence of functions $\{\eta_n\}_{n\in\mathbb N}$ is uniformly bounded and $a_n+b_n$ converges to $0$ when $n$ tends to the infinity.

Note that $\mathscr S$ is a closed subset of $C^0([0,1])$, the space of all continuous real functions on $[0,1]$. Since $C^0([0,1])$ admits a countable topological basis (see \cite{Bourbaki06} Chapter X, \S3.3, Theorem 1), also is $\mathscr S$. Being a metric space, the topological space $\mathscr S$ is thus separable. Therefore, the Borel $\sigma$-algebra of the product topological space $\mathscr S\times\intervalle[0{+\infty}{[}^2$ coincides with the product $\sigma$-algebra of the Borel $\sigma$-algebras of $\mathscr S$ and $\intervalle[0{+\infty}[^2$  (see \cite{Kallenberg} Lemma 1.2). In particular, the function $F$ is measurable with respect to the product $\sigma$-algebra. If $(s,t)$ is an element in $E\oplus F$, then one has
\[\|(s,t)\|_{\psi,\omega}=g(\psi(\omega),\|s\|_{E,\omega},\|t\|_{F,\omega}),\]
which is an $\mathcal A$-measurable function since the maps $\psi$, $\omega\mapsto\|s\|_{E,\omega}$ and $\omega\mapsto\|t\|_{F,\omega}$ are all $\mathcal A$-measurable.
\end{proof}

\begin{prop}\phantomsection\label{Pro:mesurability}
\begin{enumerate}[label=\rm(\arabic*)]
\item\label{Item:criterion of measurability} Let $E$ be a vector space of dimension $1$ over $K$ and $\xi$ be a norm family in $\mathcal N_E$. Then $\xi$ is $\mathcal A$-measurable if and only if there exists an element $s\in E\setminus\{0\}$ such that the function $(\omega\in\Omega)\mapsto \|s\|_{\omega}$ is $\mathcal A$-measurable.
\item\label{Item:dual measurability} Let $E$ be a vector space of dimension $1$ over $K$ and $\xi$ be a norm family in $\mathcal N_E$ which is $\mathcal A$-measurable. Then the dual norm family $\xi^\vee$ is  also $\mathcal A$-measurable.
\item\label{Item:tensor measurability} Let $E_1$ and $E_2$ be vector spaces of dimension $1$ over $K$, and $\xi_1$ and $\xi_2$ be norm families in $\mathcal N_{E_1}$ and $\mathcal N_{E_2}$ respectively. We assume that {both} norm families $\xi_1$ and $\xi_2$ are $\mathcal A$-measurable. Then the tensor product $\xi_1\otimes\xi_2$ (which is also equal to $\xi_1\otimes_{\varepsilon}\xi_2$ and $\xi_1\otimes_{\pi}\xi_2$) is also $\mathcal A$-measurable.
\item\label{Item:power measurability} Let $E$ be a vector space of dimension $1$ over $K$ and $\xi$ be a norm family in $\mathcal N_E$. Assume that there exists an integer $n\geqslant 1$ such that $\xi^{\otimes n}$ is $\mathcal A$-measurable, then the norm family $\xi$ is also measurable.
\end{enumerate}
\end{prop}
\begin{proof}
\ref{Item:criterion of measurability} The necessity follows from the definition. For the sufficiency, we assume that there exists $s\in E\setminus\{0\}$ such that the function $(\omega\in\Omega)\mapsto\|s\|_\omega$ is $\mathcal A$-measurable. If $s'$ is a general element in $E$, there exists $a\in K$ such that $s'=as$. Note that for any $\omega\in\Omega$ one has
\[\|s'\|_\omega=|a|_\omega\cdot\|s\|_\omega.\]
Since the function $(\omega\in\Omega)\mapsto |a|_\omega$ is $\mathcal A$-measurable, we obtain that the function $(\omega\in\Omega)\rightarrow\|s'\|_{\omega}$ is $\mathcal A$-measurable.

\ref{Item:dual measurability} Let $s$ be a non-zero element in $E$ and $\alpha$ be the element in $E^\vee$ such that $\alpha(s)=1$. For any $\omega\in\Omega$ one has $\|\alpha\|_\omega=\|s\|_\omega^{-1}$. Since the function $(\omega\in\Omega)\mapsto \|s\|_{\omega}$ is $\mathcal A$-measurable, we obtain that the function $(\omega\in\Omega)\mapsto\|\alpha\|_{\omega}$ is also $\mathcal A$-measurable. Therefore, by \ref{Item:criterion of measurability} we obtain that the norm family $\xi^\vee$ is $\mathcal A$-measurable on $\Omega$.

\ref{Item:tensor measurability} Let $s_1$ and $s_2$ be non-zero elements of $E_1$ and $E_2$ respectively. Then $s_1\otimes s_2$ is a non-zero element of $E_1\otimes E_2$. Moreover, for any $\omega\in\Omega$ one has $\|s_1\otimes s_2\|_{\omega}=\|s_1\|_{\omega}\cdot\|s_2\|_{\omega}$. Since the functions $(\omega\in\Omega)\mapsto \|s_1\|_\omega$ and $(\omega\in\Omega)\mapsto \|s_2\|_\omega$ are $\mathcal A$-measurable, we obtain that the function $(\omega\in\Omega)\mapsto \|s_1\otimes s_2\|_\omega$ is $\mathcal A$-measurable. By \ref{Item:criterion of measurability}, the norm family $\xi_1\otimes\xi_2$ is $\mathcal A$-measurable.

\ref{Item:power measurability} Let $s$ be a non-zero element of $E$. For any $\omega\in\Omega$, one has $\|s^{\otimes n}\|_\omega=\|s\|_\omega^{n}$. Hence $\|s\|_\omega=\|s^{\otimes n}\|_\omega^{1/n}$. Since $\xi^{\otimes n}$ is $\mathcal A$-measurable, the function $(\omega\in\Omega)\mapsto \|s^{\otimes n}\|_\omega$ is $\mathcal A$-measurable. As a consequence, the function $(\omega\in\Omega)\mapsto \|s\|_\omega$ is also $\mathcal A$-measurable. By \ref{Item:criterion of measurability}, we obtain that the norm family $\xi$ is $\mathcal A$-measurable.
\end{proof}

\begin{rema}
It is not clear that other algebraic constructions of norm families preserve the $\mathcal A$-measurability. We consider the following counter-example. Let $K=\mathbb R$ and $(\Omega,\mathcal A,\nu)$ be the set $\mathbb R$ equipped with the Borel $\sigma$-algebra and the Lebesgue measure. Let $\phi:\Omega\rightarrow M_{\mathbb R}$ be the constant map which sends any point of $\Omega$ to the trivial absolute value. Then $(K,(\Omega,\mathcal A,\nu),\phi)$ is an adelic curve. Let $f:\mathbb R\rightarrow \intervalle]01]$ be a map which is not Borel measurable. Let $E$ be a vector space of dimension $2$ over $\mathbb R$ and $\{e_1,e_2\}$ be a basis of $E$. For each $t\in\Omega=\mathbb R$, let $\|\ndot\|_t$ be the norm on $E$ such that $\|\lambda(e_1+te_2)\|_t=f(t)$ for $\lambda\in\mathbb R\setminus\{0\}$, and $\|s\|_t=1$ if $s\in E\setminus\mathbb R(e_1+te_2)$. Then $
\xi=\{\|\ndot\|_t\}_{t\in\mathbb R}$ is an element in $\mathcal N_E$. Note that, for any $s\in E$ the function $t\mapsto\|s\|_t$ is Borel measurable on $\mathbb R$ since it is constant except at most one point of $\mathbb R$. However, if we denote by $G$ the quotient space $E/\mathbb Re_2$ and by $\xi_G=\{\|\ndot\|_{G,t}\}_{t\in\mathbb R}$ the quotient norm family of $\xi$. The one has
\[\forall\,t\in\mathbb R,\;\big\|[e_1]\big\|_{G,t}=f(t).\]
Therefore the quotient norm family $\xi_G$ is not $\mathcal A$-measurable on $\Omega$.
\end{rema}

The following results show that, at least in the particular case where $K$ is a countable set, the algebraic constructions of norm families defined in the previous subsection preserve the $\mathcal A$-measurability of norm families.

\begin{prop}\label{Pro:mesurabilityofquotient}
We assume that the field $K$ admits a countable subfield $K_0$ which is dense in $K_\omega$ for any $\omega\in\Omega$.
\begin{enumerate}[label=\rm(\arabic*)]
\item\label{Item: measruble norm family 1} Let $E$ be a vector space of finite rank over $K$ and $\xi=\{\|\ndot\|_{\omega}\}_{\omega\in\Omega}$ be a norm family in $\mathcal N_E$ which is $\mathcal A$-measurable. Then
\begin{enumerate}[label=\rm(\arabic*)]\renewcommand{\labelenumii}{\rm(\arabic{enumi}.\alph{enumii})}
\item\label{Item: measurability of quotient norm family} for any quotient space $G$ of $E$, the quotient norm family $\xi_G=\{\|\ndot\|_{G,\omega}\}_{\omega\in\Omega}$ on $G$ of $\xi$ is $\mathcal A$-measurable;
\item\label{Item: dual norm family is measurable} the dual norm family $\xi^\vee=\{\|\ndot\|_{\omega,*}\}_{\omega\in\Omega}$ is $\mathcal A$-measurable;
\item\label{Item: extension of scalars measurable} for any algebraic extension $L/K$, the norm family $\xi_L=\{\|\ndot\|_x\}_{x\in\Omega_L}$ on $E_L:=E\otimes_KL$ is $\mathcal A_L$-measurable. If in addition $\xi$ is Hermitian, then the norm family $\xi_L^H$ is $\mathcal A_L$-mesurable.
\end{enumerate}
\item\label{Item: mearsurable norm family tensor} Let $E$ and $F$ be finite-dimensional vector spaces over $K$, $\xi_E=\{\|\ndot\|_{E,\omega}\}_{\omega\in\Omega}$ and $\xi_F=\{\|\ndot\|_{F,\omega}\}_{\omega\in\Omega}$ be respectively norm families in $\mathcal N_E$ and $\mathcal N_F$. We assume that $\xi_E$ and $\xi_F$ are both $\mathcal A$-measurable. Then
\begin{enumerate}[label=\rm{(\arabic*)}]\renewcommand{\labelenumii}{\rm(\arabic{enumi}.\alph{enumii})}
\item\label{Item: measurability of tensor product} the $\pi$-tensor product $\xi_E\otimes_\pi\xi_F$, the $\varepsilon$-tensor product $\xi_E\otimes_\varepsilon\xi_F$ and the $\varepsilon,\pi$-tensor product $\xi_E\otimes_{\varepsilon,\pi}\xi_F$ are all $\mathcal A$-measurable;
\item\label{Item: measurability of Hermitian tensor product} if in addition the norm families $\xi_E$ and $\xi_F$ are Hermitian, the orthogonal tensor product $\xi_E\otimes\xi_F=\{\|\ndot\|_{E\otimes F,\omega}\}_{\omega\in\Omega}$ is $\mathcal A$-measurable.
\end{enumerate}
\item Let $E$ be a finite-dimensional vector space over $K$ and $\xi_E=\{\|\ndot\|_{E,\omega}\}_{\omega\in\Omega}$ be a norm family in $\mathcal N_E$ which is $\mathcal A$-measurable. Then the exterior norm family $\Lambda^i\xi$ is $\mathcal A$-measurable for any $i\in\mathbb N$. In particular, the determinant norm family $\det(\xi)$ is $\mathcal A$-measurable.
\item Let $E$ be a finite-dimensional vector space over $K$, and $\xi=\{\|\ndot\|_\omega\}_{\omega\in\Omega}$ and $\xi'=\{\|\ndot\|'_{\omega}\}_{\omega\in\Omega}$ be two $\mathcal A$-measurable norm families in $\mathcal N_E$. Then the local distance function $\omega\mapsto d_\omega(\xi,\xi')$ is $\mathcal A$-measurable.
\end{enumerate}
\end{prop}
\begin{proof}\ref{Item: measurability of quotient norm family} Let $p:E\rightarrow G$ be the projection map and $F$ be its kernel. Let $F_0$ be a finite-dimensional $K_0$-vector subspace of $F$ which generates $F$ as a vector space over $K$. Note that for any $\omega\in\Omega$ the set $F_0$ is dense in $F_{K_\omega}$. For any $\ell\in G$ and any $\omega\in\Omega$, one has
\[\|\ell\|_{G,\omega}=\inf_{s\in E,\,p(s)=\ell}\|s\|_{\omega}=\inf_{t\in F_0}\|s_0+t\|_{\omega},\]
where $s_0$ is an element in $E$ such that $p(s_0)=\ell$.
As the norm family $\xi_E$ is $\mathcal A$-measurable on $\Omega$, the function $(\omega\in\Omega)\mapsto\|s\|_{\omega}$ is $\mathcal A$-measurable for any $s\in E$. Hence the function $\omega\mapsto\|\ell\|_{G,\omega}$ is also $\mathcal A$-measurable since it is the infimum of a countable family of $\mathcal A$-measurable functions.

\ref{Item: dual norm family is measurable} Let $E_0$ be a finite-dimensional $K_0$-vector subspace of $E$ which generates $E$ as a vector space over $K$. For any $\alpha\in E^\vee$ and any $\omega\in\Omega$, one has
\[\|\alpha\|_{\omega,*}=\sup_{s\in E\setminus\{0\}}\frac{|\alpha(s)|_\omega}{\|s\|_{\omega}}=\sup_{s\in E_0\setminus\{0\}}\frac{|\alpha(s)|_\omega}{\|s\|_{\omega}}\]
since $E_0\setminus\{0\}$ is dense in $E_{K_\omega}\setminus\{0\}$. As the norm family $\xi_E$ is $\mathcal A$-measurable on $\Omega$, the function $(\omega\in\Omega)\mapsto\|s\|_{\omega}$ is $\mathcal A$-measurable for any $s\in E$. Moreover, $\alpha(s)$ belongs to $K$, and thus the function $\omega\mapsto|\alpha(s)|_\omega$ is $\mathcal A$-measurable. Hence the function $(\omega\in\Omega)\mapsto\|\alpha\|_{\omega,*}$ is $\mathcal A$-measurable since it is the supremum of a countable family of $\mathcal A$-measurable functions.

\ref{Item: extension of scalars measurable} Let $H_0$ be a finite-dimensional $K_0$-vector subspace of $E^\vee$ which generates $E^\vee$ as a vector space over $K$. Then $H_0\setminus\{0\}$ is dense in $E^\vee_{K_\omega}\setminus\{0\}$ for any $\omega\in\Omega$.  Let $s$ be an element in $E_L$. For any  $x\in\Omega_L$, let
\[\|s\|_{\omega,L_x,\varepsilon}=\sup_{\varphi\in E^\vee_{K_\omega}\setminus\{0\}}\frac{|\varphi(s)|_x}{\|\varphi\|_{\omega,*}}=\sup_{\varphi\in H_0\setminus\{0\}}\frac{|\varphi(s)|_x}{\|\varphi\|_{\omega,*}},\]
where $\omega=\pi_{L/K}(x)$, and $\|\ndot\|_{\omega,*}$ denotes the dual norm of $\|\ndot\|_\omega$. We have seen in \ref{Item: dual norm family is measurable} that the dual norm family $\xi^\vee=\{\|\ndot\|_{\omega,*}\}_{\omega\in\Omega}$ is $\mathcal A$-mesurable on $\Omega$. Therefore the function \[(x\in \Omega_L)\longmapsto\|\varphi\|_{\pi_{L/K}(x),*},\] which is the composition of the $\mathcal A$-measurable function $\omega\mapsto\|\varphi\|_{\omega,*}$ with $\pi_{L/K}$, is $\mathcal A_L$-measurable. Moreover, the function $x\mapsto|\varphi(s)|_x$ on $\Omega_L$ is  $\mathcal A_L$-mesurable. Therefore, the function $x\mapsto \|s\|_{\omega,L_x,\varepsilon}$ on $\Omega_L$, which is the supremum of a countable family of measurable functions, is also measurable. Therefore the norm family $\{\|\ndot\|_x'\}_{x\in\Omega_L}$ is measurable. This result applied to $\xi^\vee$ shows that the norm family $\{\|\ndot\|_{\omega,*,L_x,\varepsilon}\}_{x\in\Omega_L}$ is measurable. By 
Proposition~\ref{Pro:comparisonofdualnormes:scalar:extension} \ref{Pro:comparisonofdualnormes:scalar:extension:epsilon}, \ref{Pro:comparisonofdualnormes:scalar:extension:epsilon:pi}, 
the norm family $\xi_L$ identifies with the dual of $\{\|\ndot\|_{\omega,*,L_x,\varepsilon}\}_{x\in\Omega_L}$, and hence is measurable.

Assume that the norm family $\xi$ is Hermitian. Let $s$ be an element of $E_L$, which is written as $s_1\otimes\lambda_1+\cdots+s_n\otimes\lambda_n$, where $(s_1,\ldots,s_n)\in E^n$ and $(\lambda_1,\ldots,\lambda_n)\in L^n$. For any $x\in\Omega_L$ one has 
\[\|s\|_x^2=\sum_{i=1}^n\sum_{j=1}^n\langle s_i,s_j\rangle_{\pi_{L/K}(x)}\langle\lambda_i,\lambda_j\rangle_x.\]
Note that the function \[(\omega\in\Omega_\infty)\mapsto\langle s_i,s_j\rangle_\omega=\frac 12\big(\|s_i+s_j\|_\omega^2-\|s_i\|_\omega^2-\|s_j\|_\omega^2\big)\]
is $\mathcal A|_{\Omega_\infty}$-measurable (hence its composition with $\pi_{L/K}|_{\Omega_{L,\infty}}$ is $\mathcal A_L|_{\Omega_{L,\infty}}$-measurable) and the function
\[(x\in\Omega_{L,\infty})\longmapsto\langle \lambda_i,\lambda_j\rangle_x=\frac 12\big(|\lambda_i+\lambda_j|_x^2-|\lambda_i|^2_x-|\lambda_j|_x^2  \big)\]
is $\mathcal A_L|_{\Omega_{L,\infty}}$-measurable. Therefore the function $(x\in\Omega_L)\mapsto\|s\|_x^2$ is $\mathcal A_L$-measurable on $\Omega_{L,\infty}$. Moreover, by the measurability of $\xi_L$ proved above, this function is also $\mathcal A_L$-measurable on $\Omega_L\setminus\Omega_{L,\infty}$. Hence it is $\mathcal A_L$-measurable.

\ref{Item: measurability of tensor product} 
Let $s$ be an element in $E\otimes F$, $\varphi$ be the $K$-linear map from $E^\vee$ to $F$ which corresponds to $s$, and $r$ be the rank of $\varphi$. Let $\{\varphi_i\}_{i=1}^{n}$ be a basis of $E^\vee$ such that $\varphi_{r+1},\ldots,\varphi_n$ belong to the kernel of $f$ and let $\{e_i\}_{i=1}^n$ be the dual basis of $\{\varphi_i\}_{i=1}^n$. For $i\in\{1,\ldots,r\}$, let $f_i$ be the image of $\alpha_i$ by $\varphi$. We complet the family $\{f_i\}_{i=1}^r$ to a basis $\{f_j\}_{j=1}^m$ of $F$. One has
\[s=e_1\otimes f_1+\cdots+e_r\otimes f_r.\]
Let $E_0$ and $F_0$ be $K_0$-vector subspaces of $E$ and $F$ generated by $\{e_i\}_{i=1}^n$ and $\{f_j\}_{j=1}^m$ respectively.

By definition, for $\omega\in\Omega$ one has 
\[\|s\|_{\pi,\omega}=\inf\left\{\sum_{i=1}^N\|x_i\|_{E,\omega}\cdot\|y_i\|_{F,\omega}\,:\,{\footnotesize\begin{array}{l}
\text{$s=x_1\otimes y_1+\cdots+x_N\otimes y_N$ for some $N\in\mathbb N$,}\\ (x_1,\ldots,x_N)\in E_{K_\omega}^N\text{ and }(y_1,\ldots,y_N)\in F_{K_\omega}^N
\end{array}}\right\}.\]
We claim that $\|s\|_{\pi,\omega}$ is eqal to 
\begin{equation}\label{Equ:infotimespi}\|s\|_{\pi,\omega}':=\inf\left\{\sum_{i=1}^N\|x_i\|_{E,\omega}\cdot\|y_i\|_{F,\omega}\,:\,{\footnotesize \begin{array}{l}
\text{$s=x_1\otimes y_1+\cdots+x_N\otimes y_N$ for some $N\in\mathbb N$}\\
(x_1,\ldots,x_N)\in E_0^N,\text{ and }(y_1,\ldots,y_N)\in F_0^N
\end{array}}\right\}.\end{equation}
Clearly $\|s\|_{\pi,\omega}$ is bounded from above by $\|s\|_{\pi,\omega}'$.
We will show that $\|s\|_{\pi,\omega}\geqslant\|s\|_{\pi,\omega}'$. By Proposition \ref{Pro:topologicalnormedspace}, there exists $\alpha\in \intervalle]01]$ such that the bases  $\{e_i\}_{i=1}^n$ and $\{f_j\}_{j=1}^m$ of $E_{K_\omega}$ and $F_{K_\omega}$ are both $\alpha$-orthogonal (see Definition \ref{Def:alphaorthogonal}).  Assume that $s$ is written in the form 
\[s=x_1\otimes y_1+\cdots+x_N\otimes y_N,\]
where $(x_1,\ldots,x_N)\in E_{K_\omega}^N$ and $(y_1,\ldots,y_N)\in F_{K_\omega}^N$. For any $\epsilon>0$, there exists $(x_1',\ldots,x_N')\in E_0^N$ and $(y_1',\ldots,y_N')\in F_0^N$ such that 
\begin{equation}\label{Equ:densityapprox}\sup_{\ell\in\{1,\ldots,N\}}\max\{\|x_\ell-x_\ell'\|_{E,\omega},\|y_\ell-y_\ell'\|_{F,\omega}\}\leqslant\epsilon.\end{equation} 
We write $x_\ell-x_{\ell}'$ into the form
\[x_{\ell}-x_{\ell}'=\sum_{i=1}^n a_{\ell,i}e_i,\quad\text{where } (a_{\ell,1},\ldots,a_{\ell,n})\in K_\omega^n.\]
Since the basis $\{e_i\}_{i=1}^n$ is $\alpha$-orthogonal, one has
\begin{equation}\label{Equ:measurablitypi1}\epsilon\geqslant\|x_{\ell}-x_{\ell}'\|_{E,\omega}\geqslant\alpha\max_{i\in\{1,\ldots,n\}}|a_{\ell,i}|_{\omega}\cdot\\|e_i\|_{E,\omega}.\end{equation}
Similarly, if we write $y_{\ell}-y_{\ell}'$ as
\[\sum_{j=1}^mb_{\ell,j}f_j, \quad (a_{\ell,1},\ldots,a_{\ell,m})\in K_\omega^m,\]
one has
\begin{equation}\label{Equ:measurablitypi2}\epsilon\geqslant\|y_{\ell}-y_{\ell}'\|_{F,\omega}\geqslant\alpha\max_{j\in\{1,\ldots,m\}}|b_{\ell,j}|_\omega\cdot\|f_j\|_{F,\omega}.\end{equation}
Let 
\[M=\sup_{\ell\in\{1,\ldots,N\}}\max\{\|x_\ell'\|_{E,\omega},\|y_\ell\|_{F,\omega}\}.\]
If we write $x_\ell'$ and $y_{\ell}$ into linear combination of $\{e_i\}_{i=1}^n$ and $\{f_j\}_{j=1}^m$ respectively:
\[x_{\ell}'=\sum_{i=1}^nc_{\ell,i}e_i,\quad y_{\ell}=\sum_{j=1}^md_{\ell,j}f_j,\]
one has
\begin{equation}\label{Equ:measurablitypi3}M\geqslant\alpha\max_{i\in\{1,\ldots,n\}}|c_{\ell,i}|_\omega\cdot\|e_i\|_{E,\omega}\text{ and } M\geqslant\alpha\max_{j\in\{1,\ldots,m\}}|d_{\ell,j}|_\omega\cdot\|f_j\|_{F,\omega}.\end{equation}
Note that
\[\begin{split}&\quad\;s-x_1'\otimes y_1'-\cdots-x_N'\otimes y_N'=\sum_{\ell=1}^N(x_\ell\otimes y_\ell-x_\ell'\otimes y_\ell')\\&=\sum_{\ell=1}^N((x_\ell-x_{\ell}')\otimes y_\ell+x_{\ell}'\otimes(y_\ell-y_{\ell}'))\\
&=\sum_{i=1}^n\sum_{j=1}^m\Big(\sum_{\ell=1}^Na_{\ell,i}d_{\ell,j}+c_{\ell,i}b_{\ell,j}\Big)e_i\otimes f_j.
\end{split}\]
For $(i,j)\in \{1,\ldots,n\}\times\{1,\ldots,m\}$, let
\[A_{i,j}=\sum_{\ell=1}^Na_{\ell,i}d_{\ell,j}+c_{\ell,i}b_{\ell,j}.\]
Since $\{e_i\otimes f_j\}_{(i,j)\in\{1,\ldots,n\}\times\{1,\ldots,m\}}$ is a basis of $E_0\otimes_{K_0}F_0$ over $K_0$ and a basis of $E_{K_\omega}\otimes_{K_\omega}F_{K_\omega}$ over $K_\omega$, one obtains that $A_{i,j}\in K_0$ for any $(i,j)\in \{1,\ldots,n\}\times\{1,\ldots,m\}$. Therefore one has
\[\begin{split}\|s\|_{\pi,\omega}'&\leqslant\sum_{\ell=1}^N\|x_\ell'\|_{E,\omega}\cdot\|y_{\ell}'\|_{F,\omega}+\sum_{i=1}^n\sum_{j=1}^m|A_{i,j}|_\omega\cdot\|e_i\|_{E,\omega}\cdot\|f_i\|_{F,\omega}\\
&\leqslant\sum_{\ell=1}^N(\|x_\ell\|_{E,\omega}+\epsilon)(\|y_{\ell}\|_{F,\omega}+\epsilon)+\sum_{i=1}^n\sum_{j=1}^m|A_{i,j}|_\omega\cdot\|e_i\|_{E,\omega}\cdot\|f_i\|_{F,\omega}\\
&\leqslant \sum_{\ell=1}^N(\|x_\ell\|_{E,\omega}+\epsilon)(\|y_{\ell}\|_{F,\omega}+\epsilon)+\alpha^{-2}\epsilon mnMN,\end{split}\]
where the first inequality comes from the definition \eqref{Equ:infotimespi} of $\|\ndot\|_{\pi,\omega}'$, the second inequality results from \eqref{Equ:densityapprox}, and the third inequality comes from \eqref{Equ:measurablitypi1}, \eqref{Equ:measurablitypi2} and \eqref{Equ:measurablitypi3}. Since $\epsilon$ is arbitrary, we obtain 
$\|s\|_{\pi,\omega}'\leqslant\sum_{\ell=1}^N\|x_\ell\|\cdot\|y_{\ell}\|$,
which leads to $\|s\|_{\pi,\omega}'\leqslant\|s\|_{\pi,\omega}$ since the writing $s=x_1\otimes y_1+\cdots+x_N\otimes y_N$ is arbitrary. 

As the set 
\[\bigcup_{N\in\mathbb N}\{(x_1,\ldots,x_N,y_1,\ldots,y_N)\in E_0^N\times F_0^N\,:\,s=x_1\otimes y_1+\cdots+x_N\otimes y_N\}\]
is countable, we obtain that the function $(\omega\in\Omega)\mapsto\|s\|_{\pi,\omega}=\|s\|_{\pi',\omega}$ is $\mathcal A$-measurable. Therefore the norm family $\xi_E\otimes_\pi\xi_F$ is $\mathcal A$-measurable. 

By Proposition \ref{Pro:dualitypiepsilon} one has
$\xi_E\otimes_\varepsilon\xi_F=(\xi_E^\vee\otimes_{\pi}\xi_F^\vee)^\vee$.
Hence by the \ref{Item: dual norm family is measurable} of the proposition established above, we obtain that the norm family $\xi_E\otimes_\epsilon\xi_F$ is also $\mathcal A$-measurable. Finally, by Corollary \ref{Cor: dual tensor product} and Proposition \ref{Pro:dualitypiepsilon} one has $\xi_E\otimes_{\varepsilon,\pi}\xi_F=(\xi_E^\vee\otimes_{\varepsilon}\xi_F^\vee)^\vee$. Therefore the norm family $\xi_E\otimes_{\varepsilon,\pi}\xi_F$ is also $\mathcal A$-measurable.

\ref{Item: measurability of Hermitian tensor product} We proceed with the measurability of the orthogonal tensor product norm family in assuming that both norm families $\xi_E$ and $\xi_F$ are Hermitian. In the first step, we treat a particular case where $E=F^\vee$ and $\xi_E=\xi_F^\vee$. In this case the tensor product space $E\otimes_KF$ is isomorphic to the space $\mathrm{End}_K(F)$ of $K$-linear endomorphisms of $F$, and the $\varepsilon$-tensor product norm family $\xi_F^\vee\otimes_\varepsilon\xi_F$ identifies with the family of operator norms. Let $\{x_i\}_{i=1}^n$ be a basis of $F$ over $K$ and $F_0$ be the $K_0$-vector subspace of $F$ generated by $\{x_i\}_{i=1}^n$. By using the basis $\{x_i\}_{i=1}^n$ one can identify $\mathrm{End}_{K_0}(F_0)$ with $M_{n\times n}(K_0)$, the space of all matrices of size $n\times n$ with coefficients in $K_0$. Similarly, for any $\omega\in\Omega$, one can identify $\mathrm{End}_{K_\omega}(F_{K_\omega})$ with the space $M_{n\times n}(K_\omega)$ of all matrices of size $n\times n$ with coefficients in $K_\omega$. In particular, $\mathrm{End}_{K_0}(F_0)$ is dense in $\mathrm{End}_{K_\omega}(F_{K_\omega})$.

For any $f\in\mathrm{End}_K(F)$ and any $\omega\in\Omega_\infty$, let $\|f\|_{\mathrm{HS},\omega}$ be the Hilbert-Schmidt norm of $f$. By Proposition \ref{Pro:Hilbert-Schmidt}, one has
\[\|f\|_{\mathrm{HS},\omega}=\bigg(\sum_{i=1}^r\inf_{\begin{subarray}{c}g\in\mathrm{End}_{K_\omega}(F_{K_\omega})\\
\rang(g)\leqslant i-1
\end{subarray}}\|f-g\|_\omega^{2}\bigg)^{1/2},\]
where $\|\ndot\|_\omega$ denotes the operator norm on $\mathrm{End}_{K_\omega}(F_{K_\omega})$.
Note that the set 
\[\{g\in\mathrm{End}_{K_0}(F_0)\,:\,\rang(g)\leqslant i-1\}\] is dense in \[\{g\in\mathrm{End}_{K_\omega}(F_{K_\omega})\,:\,\rang(g)\leqslant i-1\}.\] Hence
\[\|f\|_{\mathrm{HS},\omega}=\bigg(\sum_{i=1}^r\inf_{\begin{subarray}{c}g\in\mathrm{End}_{K_0}(F_0)\\
\rang(g)\leqslant i-1
\end{subarray}}\|f-g\|_\omega^{2}\bigg)^{1/2}.\]
By the result of \ref{Item: measurability of tensor product} on the  measurability of $\varepsilon$-tensor product, the function \[(\omega\in\Omega_\infty)\mapsto \|f-g\|_{\omega}\] is $\mathcal A|_{\Omega_\infty}$-measurable. Hence we deduce the $\mathcal A|_{\Omega_\infty}$-measurability of the function $(\omega\in\Omega_\infty)\mapsto \|f\|_{\mathrm{HS},\omega}$. Finally, since the norm $\|\ndot\|_{E\otimes F,\omega}$ identifies with the $\varepsilon$-tensor product $\|\ndot\|_{\varepsilon,\omega}$ for $\omega\in\Omega\setminus\Omega_\infty$, by the result of \ref{Item: measurability of tensor product}, the function $(\omega\in\Omega\setminus\Omega_\infty)\mapsto\|f\|_{E\otimes F,\omega}$ is  $\mathcal A|_{\Omega\setminus\Omega_\infty}$-measurable. Thus we obtain the $\mathcal A$-measurability of the norm family $\xi_E\otimes\xi_F$.

We now consider the general case. Let $T$ be a tensor vector in $E\otimes_KF$, viewed as a linear map from $E^\vee$ to $F$. Let $G$ be the direct sum $E^\vee\oplus F$ and $\xi_G$ be the orthogonal direct sum of $\xi_E^\vee$ and $\xi_F$. By Proposition \ref{Pro:mesurabilitedesomme} and the result obtained in (1.b), the norm family $\xi_G$ is $\mathcal A$-measurable. Moreover, the linear map $T:E^\vee\rightarrow F$ induces a $K$-linear endomorphism $f=\begin{pmatrix}0&0\\T&0
\end{pmatrix}$ of $E^\vee\oplus F$. For any $\omega\in\Omega_\infty$, the Hilbert-Schmidt norm of $T$ with respect to $\|\ndot\|_{E,\omega}$ and $\|\ndot\|_{F,\omega}$ identifies with the Hilbert-Schmidt norm of $f$ with respect to the orthogonal direct sum norm $\|\ndot\|_{E^\vee\oplus F,\omega}$. By the particular case proved above, we obtain the measurability of the function $(\omega\in\Omega_\infty)\mapsto\|T\|_{E\otimes F,\omega}$. Combined with the measurability of $\xi_E\otimes_{\varepsilon}\xi_F$ proved in \ref{Item: measurability of tensor product}, which implies the measurability of the function $(\omega\in\Omega\setminus\Omega_\infty)\mapsto \|T\|_{E\otimes F,\omega}$, we obtain that  the function $(\omega\in\Omega)\mapsto \|T\|_{E\otimes F,\omega}$ is $\mathcal A$-measurable. The assertion is thus proved.

(3) We equip $E^{\otimes i}$ with the $i$-th $\varepsilon,\pi$-tensor power of $\xi$. By \ref{Item: measurability of tensor product}, this norm family is $\mathcal A$-measurable. The exterior power norm family is its quotient. By \ref{Item: measurability of quotient norm family} we obtain that the norm family $\Lambda^i\xi$ is $\mathcal A$-measurable.

(4) Let $E_0$ be a finite-dimensional $K_0$-vector subspace of $E$, which generates $E$ as a $K$-vector space. For any $\omega\in\Omega$, $E_0$ is dense in $E_{K_\omega}$. Therefore one has
\[d_\omega(\xi,\xi')=\sup_{0\neq s\in E_0}\Big|\ln\|s\|_\omega-\ln\|s\|_{\omega}'\Big|.\]
Since $E_0$ is a countable set and since the functions $(\omega\in\Omega)\mapsto\|s\|_\omega$ and $(\omega\in\Omega)\mapsto\|s\|_{\omega}'$ are both measurable, we obtain that the function $(\omega\in\Omega)\mapsto d_\omega(\xi,\xi')$ is also measurable.
\end{proof}

\begin{rema}\label{remark:number:field:measurable}
We assume that $K$ is a number field. Let $S=(K,(\Omega,\mathcal A,\nu),\phi)$ be the standard adelic curve 
as in Subsection~\ref{Subsec:Numberfields}.
Let $E$ be a finite-dimensional vector space over $K$ and $\xi = \{ \|\ndot\|_{\omega} \}_{\omega \in \Omega}$
be a norm family of $E$. Since $\mathcal A$ is the disctere $\sigma$-algebra,
every function on $\Omega$ is measurable, so that $(E,\xi)$ is measurable.
\end{rema}

\begin{theo}\label{Thm: Hermitian approximation via measurable selection}We assume that $K$ admits a countable subfield $\widetilde K$ which is dense in  $K_\omega$ for any $\omega\in\Omega_\infty$. Let $E$ be a finite-dimensional vector space over $K$ and $\xi=\{\norm{\ndot}_\omega\}_{\omega\in\Omega}$ be a measurable norm family on $E$ which is ultrametric on $\Omega\setminus\Omega_\infty$. For any $\epsilon>0$, there exists a \emph{measurable} Hermitian norm family $\xi^H=\{\norm{\ndot}_\omega^H\}_{\omega\in\Omega}$ on $E$ such that $\norm{\ndot}_\omega^H=\norm{\ndot}_\omega$ for any $\omega\in\Omega\setminus\Omega_\infty$ and that
\[\norm{\ndot}_{\omega}\leqslant\norm{\ndot}_\omega^H\leqslant (r+\epsilon)^{1/2}\norm{\ndot}_\omega\]
for any $\omega\in\Omega_\infty$, where $r$ denotes the rank of $E$ over $K$.
\end{theo}
\begin{proof}
The assertion is trivial when $\Omega_\infty$ is empty. In what follows, we assume that $\Omega_\infty$ is non-empty. In this case the field $K$ is of characteristic $0$.
We divide the proof of the theorem into two steps.

\medskip
{\bf Step 1:} In this step, we  show that
there is a family $\{ \varphi_{\omega} \}_{\omega \in \Omega_{\infty}}$ of embeddings $\varphi_{\omega} : K \to
\mathbb C$ such that $|\ndot|_{\omega} = |\varphi_{\omega}(\ndot)|$ for all $\omega \in \Omega_{\infty}$
and, for any $a \in \widetilde K$, the function
$(\omega \in \Omega_{\infty}) \mapsto (\varphi_{\omega}(a) \in \mathbb C)$ is $\mathcal A|_{\Omega_\infty}$-measurable, where $|\ndot|$ denotes the usual absolute value on $\mathbb C$.

For each $\omega\in\Omega_\infty$, we fix an embedding $\widetilde\varphi_\omega:K\rightarrow\mathbb C$ such that $|\ndot|_\omega=|\widetilde\varphi_\omega(\ndot)|$. 
We denote by $f_\omega(a)$
and $g_\omega(a)$ the real part and the imaginary part of $\widetilde\varphi_\omega(a)$, respectively, 
that is, $\widetilde\varphi_\omega(a) = f_\omega(a) + \sqrt{-1} g_\omega(a)$ for $a \in K$. 
We claim that, for any $a\in K$, the function $(\omega\in\Omega_\infty)\mapsto f_\omega(a)$ is $\mathcal A|_{\Omega_\infty}$-measurable. In fact, we can write $f_\omega(a)$ as $|a+\textstyle{\frac 12}|^2_\omega-|a|^2_\omega-\frac{1}{4}$. Therefore the claim follows from the definition of adelic curve. Moreover, for any $a\in K$, the function $(\omega\in\Omega_\infty)\mapsto g_\omega(a)^2$ is measurable since we can write $g_\omega(a)^2$ as $|a|_\omega^2-f_\omega(a)^2$. In particular, the function $(\omega \in \Omega_{\infty}) \mapsto |g_{\omega}(a)|$ is measurable.
As a consequence, for any couple of elements $(a,b)$ in $K$, the function $(\omega\in\Omega_\infty)\mapsto g_\omega(a)g_\omega(b)$ is measurable because 
\[g_\omega(a)g_\omega(b)=\textstyle{\frac 12}(g_\omega(a+b)^2-g_\omega(a)^2-g_\omega(b)^2).\]

\begin{enonce}{Claim}\label{Claim: measurable criterion}
Let $a$ be an element of $K$. Assume that $s:\Omega_\infty\rightarrow\{1,-1\}$ is a map such that the function $(\omega\in\Omega_\infty)\mapsto s(\omega)g_\omega(a)$ is $\mathcal A|_{\Omega_\infty}$-measurable.  Then for any function $\eta:\Omega_\infty\rightarrow\mathbb R$ such that the function $(\omega\in\Omega_\infty)\mapsto g_\omega(a)\eta(\omega)$ is $\mathcal A|_{\Omega_\infty}$-measurable, the function $(\omega\in\Omega_\infty)\mapsto s(\omega)\eta(\omega)\indic_{g_\omega(a)\neq 0}$ is $\mathcal A|_{\Omega_\infty}$-measurable. In particular, for any $b\in K$, the function $(\omega\in\Omega_\infty)\mapsto s(\omega)g_\omega(b)\indic_{g_\omega(a)\neq 0}$ is $\mathcal A|_{\Omega_\infty}$-measurable.
\end{enonce}
\begin{proof}
Let $h_{a,s}:\Omega_\infty\rightarrow\mathbb R$ be the function defined by
\[h_{a,s}(\omega):=\begin{cases}
s(\omega)g_\omega(a)^{-1}=(s(\omega)g_\omega(a))^{-1},&\text{if $g_\omega(a)\neq 0$},\\
0,&\text{if $g_{\omega}(a)=0$}.
\end{cases}\] 
As the function $(\omega\in\Omega_\infty)\mapsto s(\omega)g_\omega(a)$ is $\mathcal A|_{\Omega_\infty}$-measurable, also is $h_{a,s}$. Since the function $(\omega\in\Omega_\infty)\mapsto g_{\omega}(a)\eta(\omega)$ is measurable, we deduce that the function 
\[(\omega\in\Omega_\infty)\longmapsto h_{a,s}(\omega)g_\omega(a)\eta(\omega)=s(\omega)\eta(\omega)\indic_{g_\omega(a)\neq 0}\]
is $\mathcal A|_{\Omega_\infty}$-measurable, which proves the first assertion. The second assertion follows from the first one and the fact that the function $(\omega\in\Omega_\infty)\mapsto g_\omega(a)g_\omega(b)$ is $\mathcal A|_{\Omega_\infty}$-measurable. 
\end{proof}

In the following, we show that there exists a map $s:\Omega_\infty\rightarrow\{1,-1\}$ such that, for any $a\in\widetilde K$, the function $(\omega\in\Omega_\infty)\mapsto s(\omega)g_\omega(a)$ is $\mathcal A|_{\Omega_\infty}$-measurable. Since $\widetilde K$ is countable, we can write its elements in a sequence $\{a_n\}_{n\in\mathbb N}$. We will construct by induction a decreasing sequence of functions $s_n:\Omega_\infty\rightarrow\{1,-1\}$ ($n\in\mathbb N$) which satisfy the following conditions:
\begin{enumerate}[label=\rm(\arabic*)]
\item for any $n\in\mathbb N$ and any $i\in\{0,\ldots,n\}$, the function $(\omega\in\Omega_\infty)\mapsto s_n(\omega)g_\omega(a_i)$ is  $\mathcal A|_{\Omega_\infty}$-measurable,
\item for any $n\in\mathbb N$, 
one has \[\{\omega\in\Omega_\infty\,:\,s_n(\omega)=1\}\supseteq \{\omega\in\Omega_\infty\,:\,g_\omega(a_0)=\cdots=g_\omega(a_n)=0\}.\]
\end{enumerate}
In the case where $n=0$, we just choose 
\[s_0(\omega)=\begin{cases}
1&\text{if $g_\omega(a_0)\geqslant 0$,}\\
-1&\text{if $g_\omega(a_0)<0$}.
\end{cases}\]
Then one has $s_0(\omega) g_\omega(a_0)=|g_\omega(a_0)|$. Therefore the function \[(\omega\in\Omega_\infty)\mapsto s_0(\omega)g_\omega(a_i)\] is $\mathcal A|_{\Omega_\infty}$-measurable. Moreover, by definition one has 
\[\{\omega\in\Omega_\infty\,:\,s_0(\omega)=1\}\supseteq \{\omega\in\Omega_\infty\,:\,g_\omega(a_0)=0\}.\] Assume that the functions $s_0,\ldots,s_n$ have been constructed, which satisfy the conditions above. We choose
\[s_{n+1}(\omega)=\begin{cases}
-s_n(\omega),&\text{if $g_\omega(a_0)=\cdots= g_{\omega}(a_n)=0$ and $g_{\omega}(a_{n+1})< 0$},\\
s_n(\omega),&\text{otherwise}.
\end{cases}\]
Clearly, if $g_\omega(a_0)=\cdots=g_\omega(a_n)=g_\omega(a_{n+1})=0$, then $s_{n+1}(\omega)=s_n(\omega)=1$,
so that the above condition (2) for $s_{n+1}$ is satisfied. Moreover, if $g_\omega(a_0)=\cdots=g_{\omega}(a_n)=0$ and 
$g_\omega(a_{n+1})<0$, then $s_{n}(\omega)=1$ and $s_{n+1}(\omega)=-s_{n}(\omega)=-1$. Hence we always have $s_{n+1}(\omega)\leqslant s_n(\omega)$ for any $\omega\in\Omega_\infty$. For any $i\in\{0,\ldots,n\}$ and any $\omega\in\Omega_\infty$ one has 
\[s_{n+1}(\omega)g_\omega(a_i)=s_n(\omega)g_\omega(a_i)\]
since $s_{n+1}(\omega)=s_n(\omega)$ once $g_\omega(a_i)\neq 0$. Hence by the induction hypothesis the function $(\omega\in\Omega_\infty)\mapsto s_{n+1}(\omega)g_\omega(a_i)$ is $\mathcal A|_{\Omega_\infty}$-measurable. 

In what follows, we show that the function $(\omega\in\Omega_\infty)\mapsto s_{n+1}(\omega)g_\omega(a_{n+1})$ is also $\mathcal A|_{\Omega_\infty}$-measurable. For any $i\in\{0,\ldots,n\}$, the set \[\{\omega\in\Omega_\infty\,:\,g_\omega(a_i)=0\}=\{\omega\in\Omega_\infty\,:\,g_\omega(a_i)^2=0\}\] belongs to $\mathcal A$ since the function $(\omega\in\Omega_\infty)\mapsto g_\omega(a_i)^2$ is $\mathcal A|_{\Omega_\infty}$-measurable. 
By Claim \ref{Claim: measurable criterion},  the function \begin{equation}\label{Equ: sn+1 gan+1}(\omega\in\Omega_\infty)\longmapsto s_{n+1}(\omega)g_{\omega}(a_{n+1})\indic_{g_\omega(a_0)\neq 0}\end{equation} is $\mathcal A|_{\Omega_\infty}$-measurable. Similarly, for any $i\in\{1,\ldots,n\}$ we deduce from the $\mathcal A|_{\Omega_\infty}$-measurability of the function 
\[(\omega\in\Omega)\longmapsto g_\omega(a_{i})g_\omega(a_{n+1})\indic_{g_\omega(a_0)=0}\cdots\indic_{g_\omega(a_{i-1})=0}\]
that the function 
\begin{equation}\label{Equ: sn+1 gan+1i}(\omega\in\Omega)\longmapsto s_{n+1}(\omega)g_\omega(a_{n+1})\indic_{g_\omega(a_0)=0}\cdots\indic_{g_\omega(a_{i-1})=0}\indic_{g_\omega(a_i)\neq 0} \end{equation}
is $\mathcal A|_{\Omega_\infty}$-measurable. Moreover, the function 
\begin{equation}\label{Equ: sn+1gan+1indic}(\omega\in\Omega_\infty)\longmapsto s_{n+1}(\omega)g_\omega(a_{n+1})\indic_{g_\omega(a_0)=\cdots=g_\omega(a_n)=0}\end{equation}
is $\mathcal A|_{\Omega_\infty}$-measurable since $s_{n+1}(\omega)g_{\omega}(a_{n+1})=|g_{\omega}(a_{n+1})|$ once the condition $g_\omega(a_0)=\cdots=g_\omega(a_n)=0$ holds and since the set 
\[\{\omega\in\Omega_\infty\,:\,g_\omega(a_0)=\cdots=g_\omega(a_n)=0\}=\{\omega\in\Omega_\infty\,:\,g_\omega(a_0)^2=\cdots=g_\omega(a_n)^2=0\}\]
belongs to $\mathcal A$. We then obtain the $\mathcal A|_{\Omega_\infty}$-measurability of the function $(\omega\in\Omega_\infty)\mapsto s_{n+1}(\omega)g_\omega(a_{n+1})$ since we can write the function as the sum of \eqref{Equ: sn+1 gan+1}, \eqref{Equ: sn+1gan+1indic} and \eqref{Equ: sn+1 gan+1i} for $i\in\{1,\ldots,n\}$.  
Let $s$ be the limit of the decreasing sequence of functions $\{s_n\}_{n\in\mathbb N}$. For any $n\in\mathbb N$, the function $(\omega\in\Omega_\infty)\mapsto s_m(\omega)g_\omega(a_n)$ is $\mathcal A|_{\Omega_\infty}$-measurable for any integer $m\geqslant n$. By passing to limit when $m\rightarrow+\infty$, we obtain that the function $(\omega\in\Omega_\infty)\mapsto s(\omega)g_\omega(a_n)$ is $\mathcal A|_{\Omega_\infty}$-measurable.

Here we define $\varphi_{\omega} : K \to \mathbb C$ to be
\[
\varphi_{\omega} := \begin{cases}
\widetilde{\varphi}_{\omega} & \text{if $s(\omega) = 1$,} \\
\text{the complex conjugation of $\widetilde{\varphi}_{\omega}$} &  \text{if $s(\omega) = -1$.}
\end{cases}
\]
By the measurability result proved above, we obtain that, for any $a\in \widetilde K$, the function $(\omega\in\Omega_\infty)\mapsto\varphi_\omega(a)$ is an $\mathcal A|_{\Omega_\infty}$-measurable.

\medskip
{\bf Step 2: }
Let  $\Omega_{\infty,c}$ be the set of $\omega\in\Omega_\infty$ such that $K_\omega=\mathbb C$. Then one has
\[\Omega_{\infty}\setminus\Omega_{\infty,c}=\bigcap_{a\in\widetilde K}\{\omega\in\Omega_\infty\,:\,\varphi_\omega(a)\in\mathbb R\}.\]
Therefore, the sets $\Omega_{\infty,c}$ and $\Omega_{\infty}\setminus\Omega_{\infty,c}$ belong to $\mathcal A$ (by Proposition \ref{Pro:mesurabilite}, the set $\Omega_\infty$ belongs to $\mathcal A$).

Let $\{e_i\}_{i=1}^r$ be a basis of $E$ over $K$. We consider $\mathbb C^{r\times r}$ as the space of complex matrices of size $r\times r$ and equip it with the Euclidean topology and Borel $\sigma$-algebra. 

{
Let $r'=r+\varepsilon$. For any $\omega\in\Omega_{\infty,c}$, let $G(\omega)$ be the set of all matrices $A\in\mathbb C^{r\times r}$ such that, for any $(\lambda_1,\ldots,\lambda_r)\in\widetilde K^r$, if we note $(x_1,\ldots,x_r)=(\varphi_\omega(\lambda_1),\ldots,\varphi_\omega(\lambda_r))$, then one has
\[\norm{x_1e_1+\cdots+x_re_r}_{\omega}^2\leqslant (\overline x_1,\ldots,\overline x_r)A^*A\begin{pmatrix}
x_1\\
\vdots\\
x_r
\end{pmatrix}\leqslant r'\norm{x_1e_1+\cdots+x_re_r}_\omega^2,\]
where $A^* = {}^t\overline{A}$.
For any $\omega\in\Omega_\infty\setminus\Omega_{\infty,c}$, let $G(\omega)$ be the set of all matrices $A\in \mathbb C^{r\times r}$ such that, for any $(\lambda_1,\ldots,\lambda_r)\in \widetilde K^r$ and any $(\lambda_1',\ldots,\lambda_r')\in\widetilde K^r$, if we note \[(z_1,\ldots,z_r):=(\varphi_\omega(\lambda_1)+\sqrt{-1}\varphi_\omega(\lambda_1'),\ldots,\varphi_\omega(\lambda_r)+\sqrt{-1}\varphi_\omega(\lambda_r')),\]
then
\[ \norm{z_1e_1+\cdots+z_re_r}^2_{\omega,\mathbb C}\leqslant(\overline z_1,\ldots,\overline z_r)A^*A\begin{pmatrix}z_1\\
\vdots\\
z_r
\end{pmatrix}\leqslant r'\norm{z_1e_1+\cdots+z_re_r}^2_{\omega,\mathbb C},\]
where  the norm $\norm{\ndot}_{\omega,\mathbb C}$ is defined as follows
\[\forall\,(x,y)\in E_{K_\omega}^2,\quad\norm{x+\sqrt{-1}y}_{\omega,\mathbb C}:=(\norm{x}_\omega^2+\norm{y}_\omega^2)^{-1/2}.\]
For any $A\in\mathbb C^{r\times r}$, let $E(A)$ be the subset of $\mathbb C^r$ defined as 
\[\{x\in\mathbb C^r\,:\,{}^t \overline x A^*Ax\leqslant 1 \}.\]
For any $\omega\in\Omega_{\infty}$, let
\[B_\omega:=\{(x_1,\ldots,x_r)\in\mathbb C^r\,:\,\norm{x_1e_1+\cdots+x_re_r}_{\omega,\mathbb C}\leqslant 1\}.\]
Then a matrix $A$ belongs to $G(\omega)$ if and only if the following conditions hold
\[(r')^{-1/2}B_\omega\subseteq E(A)\subseteq B_\omega.\]
In fact, this relation is equivalent to, for any $(z_1,\ldots,z_r)\in\mathbb C^r$,
\begin{equation}\label{Equ: containing unit ball}\norm{z_1e_1+\cdots+z_re_r}^2_{\omega,\mathbb C}\leqslant (\overline z_1,\ldots,\overline z_r)A^*A\begin{pmatrix}z_1\\
\vdots\\
z_r
\end{pmatrix}\leqslant r'\norm{z_1e_1+\cdots+z_re_r}^2_{\omega,\mathbb C}.\end{equation}
Hence it implies that $A\in G(\omega)$. The converse implication is also true since $\varphi_\omega(\widetilde K)+\sqrt{-1}\varphi_\omega(\widetilde K)$ is dense in $\mathbb C$ if $K_\omega=\mathbb R$, and $\varphi_\omega(\widetilde K)$ is dense in $\mathbb C$ if $K_\omega=\mathbb C$.}

Let $\delta>0$ such that $(1+\delta)\sqrt{r}<\sqrt{r'}$.
By Theorem \ref{Thm:John} (see also Remark \ref{rema:Thm:John}), there exists a Hermitian matrix $A_0$ such that 
\[(1+\delta)^{-1} r^{-1/2}B_\omega\subseteq E(A_0)\subseteq (1+\delta)^{-1}B_\omega.\]
This shows that $A_0$ belongs to the interior of $G(\omega)$ and hence the interior $G(\omega)^\circ$ of $G(\omega)$ is not empty. We denote by $F(\omega)$ the closure of $G(\omega)^\circ$.  
If $U$ is an open subset of $\mathbb C^{r\times r}$, one has
\[\begin{split}\{\omega\in\Omega_\infty\,:\,U&\cap F(\omega)\neq \varnothing\}=\{\omega\in\Omega_\infty\,:\,U\cap G(\omega)^\circ\neq \varnothing\}\\&=\bigcup_{A\in\overline{\mathbb Q}^{r\times r}\cap U}\{\omega\in\Omega_\infty\,:\,A\in G(\omega)^\circ\}.\end{split}\]
Note that for any matrix $A\in\mathbb C^{r\times r}$
\[\{\omega\in\Omega_{\infty,c}\,:\,A\in G(\omega)^\circ\}=\bigcup_{\mu\in\mathbb Q\cap (0,1)}\bigcap_{\lambda=(\lambda_1,\ldots,\lambda_r)\in\widetilde K^r}\Omega_{\infty,c}(\mu,\lambda),\]
where $\Omega_{\infty,c}(\mu,\lambda_1,\ldots,\lambda_r)$ is the set of $\omega\in\Omega_{\infty,c}$ such that
\[\begin{split} (1+\mu)\norm{\lambda_1e_1+\cdots+\lambda_re_r}^2_\omega\leqslant(\overline{\varphi_\omega(\lambda_1)},&\ldots,\overline{\varphi_\omega(\lambda_r)})A^*A\begin{pmatrix}\varphi_\omega(\lambda_1)\\
\vdots\\
\varphi_\omega(\lambda_r)
\end{pmatrix}\\
&\leqslant (1-\mu)r'\norm{\lambda_1e_1+\cdots+\lambda_re_r}^2_\omega.\end{split}\]
Note that $\Omega_\infty(\mu,\lambda_1,\ldots,\lambda_r)$ belongs to $\mathcal A$ since the functions \[(\omega\in\Omega_\infty)\longmapsto \norm{\lambda_1e_1+\cdots+\lambda_re_r}_\omega^2\] and \[(\omega\in\Omega_\omega)\longmapsto \varphi_\omega(\lambda_i) \quad\text{($i\in\{1,\ldots,r\}$)}\] are $\mathcal A|_{\Omega_\infty}$-measurable. {We then deduce that $\{\omega\in\Omega_{\infty,c}\,:\,A\in G(\omega)^\circ\}$ belongs to $\mathcal A$. Similarly, \[\{\omega\in\Omega_\infty\setminus\Omega_{\infty,c}\,:\,A\in G(\omega)^\circ\}=\bigcup_{\mu\in\mathbb Q\cap (0,1)}\bigcap_{\begin{subarray}{c}\lambda=(\lambda_1,\ldots,\lambda_r)\in\widetilde K^r\\
\lambda'=(\lambda_1',\ldots,\lambda_r')\in\widetilde K^r
\end{subarray}}(\Omega_\infty\setminus\Omega_{\infty,c})(\mu,\lambda,\lambda'),\]
where $(\Omega_\infty\setminus\Omega_{\infty,c})(\mu,\lambda,\lambda')$ is the set of $\omega\in\Omega_{\infty,c}$ such that
\[\begin{split} &(1+\mu)\Big(\norm{\lambda_1e_1+\cdots+\lambda_re_r}^2_\omega+\norm{\lambda_1e_1+\cdots+\lambda_re_r}^2_\omega\Big)\\
&\qquad\leqslant(\overline z_1(\omega),\ldots,\overline z_r(\omega))A^*A\begin{pmatrix}z_1(\omega)\\
\vdots\\
z_n(\omega)
\end{pmatrix}\\
&\qquad\qquad\leqslant (1-\mu)r'\Big(\norm{\lambda_1e_1+\cdots+\lambda_re_r}^2_\omega+\norm{\lambda_1'e_1+\cdots+\lambda_r'e_r}^2_\omega\Big).\end{split}\]
Hence $\{\omega\in\Omega_\infty\setminus\Omega_{\infty,c}\,:\,A\in G(\omega)^\circ\}$ belongs to $\mathcal A$.
}

Gathering the results we obtained, one can conclude that
\[\{ \omega \in \Omega_{\infty} \,:\, U \cap F({\omega}) \not= \varnothing \}\] belongs to $\mathcal A$, so that by 
the measurable selection theorem of Kuratowski and Ryll-Nardzweski (see \ref{Thm: measurable selection}), we obtain that there exists an $\mathcal A|_{\Omega_\infty}$-measurable map $\alpha:\Omega_\infty\rightarrow M$ such that $\alpha(\omega)$ belongs to $F(\omega)$ for any $(\omega\in\Omega)$. Finally, for any $\omega\in\Omega_\infty$ and any $(\lambda_1,\ldots,\lambda_r)\in K_\omega^r$ and  we let (where we extend $\varphi_\omega$ by continuity to $K_\omega\rightarrow\mathbb C$)
\[\norm{\lambda_1e_1+\cdots+\lambda_re_r}_\omega^H:=
\left[(\overline{\varphi_\omega(\lambda_1)},\ldots,{\overline{\varphi_\omega(\lambda_r)}})\alpha(\omega)^*\alpha(\omega)\begin{pmatrix}\varphi_\omega(\lambda_1)\\
\vdots\\
\varphi_\omega(\lambda_r)
\end{pmatrix}\right]^{1/2}.\]
Then $\norm{\ndot}_\omega^H$ is a Hermitian norm which satisfies
\[\norm{\ndot}_\omega\leqslant\norm{\ndot}_\omega^H\leqslant (r+\varepsilon)^{1/2}\norm{\ndot}_\omega.\]
For $\omega\in\Omega\setminus\Omega_\infty$, let $\norm{\ndot}_{\omega}^H:=\norm{\ndot}_{\omega}$. Then by the measurability of the map $\alpha(\ndot)$ we obtain that the norm family $\xi^H:=\{\norm{\ndot}_\omega^H\}_{\omega\in\Omega}$ is measurable. The theorem is thus proved. 
\end{proof}

\subsection{Adelic vector bundles}

In this section, we {introduce} the notion of adelic vector bundles on an adelic curve $S=(K,(\Omega,\mathcal A,\nu),\phi)$. An adelic vector bundle is a finite-dimensional vector space $E$ over $K$ equipped with a family of norms indexed by $\Omega$, which satisfies some measurability and dominance conditions so that the height of non-zero vectors is well defined (see Definition \ref{Def:adelicvecbun}). In the classic setting of global fields, the notion of adelic vector bundles was defined differently in the literature (see for example \cite{Gaudron08}): one requires that almost all norms come from a common integral model of $E$. However, in our setting it is not relevant to consider integral models. The readers will discover the link between our definition and the classic one via the dominance property described in Proposition \ref{Pro:comparaisonavecbasenorm}.

\begin{defi}
\label{Def:adelicvecbun} 
Let $E$ be a finite-dimensional vector space over $K$, and $\xi$ be a norm family in $\mathcal N_E$.
If both norm families $\xi$ and $\xi^\vee$ are $\mathcal A$-measurable on $\Omega$ and if $\xi$ is dominated (resp. strongly dominated), we say that the couple $(E,\xi)$ is \emph{an adelic vector bundle}\index{adelic vector bundle} (resp. a \emph{strongly adelic vector bundle})\index{adelic vector bundle!strongly ---}\index{strongly adelic vector bundle} on $S$.

If  the norm family $\xi$ is Hermitian (in this case $(E,\xi)$ is necessarily a strongly adelic vector bundle), we say that $(E,\xi)$ is a \emph{Hermitian adelic vector bundle}\index{Hermitian adelic vector bundle}\index{adelic vector bundle!Hermitian ---} on $S$. If the rank of $E$ is $1$ (in this case $\xi$ is necessarily Hermitian), we say that $(E,\xi)$ is an \emph{adelic line bundle}\index{adelic line bundle} on $S$. 
\end{defi}

\begin{prop}\label{Pro:criterion of adelic line bundle}
Let $E$ be a vector space of rank $1$ over $K$ and $\xi$ be a norm family in $\mathcal N_E$. If $\xi$ is $\mathcal A$-measurable and dominated, then $(E,\xi)$ is an adelic line bundle on $S$. 
\end{prop}
\begin{proof}
Since $E$ is of rank $1$ over $K$, any dominated norm family is strongly dominated (see Remark \ref{Rem: strongly dominated}). Moreover, by Proposition \ref {Pro:mesurability}, if $\xi$ is $\mathcal A$-measurable, then also is $\xi^\vee$.
\end{proof}

\begin{prop}\label{Pro:dual adelic line bundle}
Let $(E,\xi)$ be an adelic line bundle on $S$. Then $(E^\vee,\xi^\vee)$ is  an adelic line bundle on $S$.
\end{prop}
\begin{proof}
By definition, the norm family $\xi^\vee$ is $\mathcal A$-measurable on $\Omega$.  Moreover, by Proposition \ref{Pro:dominancealgebraic} \ref{Item:dualdom}, the norm family $\xi^\vee$ is dominated. By Proposition \ref{Pro:criterion of adelic line bundle}, we obtain that $(E^\vee,\xi^\vee)$ is an adelic line bundle on $S$.
\end{proof}

The following proposition is fundamental in the height theory of rational points in a projective space over an adelic curve.
\begin{prop}\label{Pro:quotientandsubofrank1}
Let $(E,\xi)$ be an adelic vector bundle on $S$. 
\begin{enumerate}[label=\rm(\arabic*)]
\item\label{Item: subspace of rank one adelic line bundle} Any vector subspace of rank $1$ of $E$ equipped with the restriction of the norm family $\xi$ forms an adelic line bundle on $S$.
\item\label{Item: quotient space of rank one adelic line bundle} Any quotient vector space of rank $1$ of $E$ equipped with the quotient norm family of $\xi$ forms an adelic line bundle on $S$.
\end{enumerate}
\end{prop}
\begin{proof}
\ref{Item: subspace of rank one adelic line bundle} Let $F$ be a vector subspace of rank $1$ of $E$ and $\xi_F$ be the restriction of $\xi$ {to} $F$. Clearly $\xi_F$ is $\mathcal A$-measurable. Moreover, by Proposition \ref{Pro:dominancealgebraic} \ref{Item:subdom}, the norm family $\xi_F$ is dominated. By Proposition \ref{Pro:criterion of adelic line bundle}, $(F,\xi_F)$ is an adelic line bundle on $S$.

\ref{Item: quotient space of rank one adelic line bundle} Let $G$ be a quotient space of rank $1$ of $E$ and $\xi_G$ be the quotient of the norm family $\xi$ on $G$. Then $G^\vee$ identifies with a vector subspace of rank $1$ of $E^\vee$ and $\xi_G^\vee$ identifies with the restriction of $\xi^\vee$ {to} $G^\vee$ (see Proposition \ref{Pro:dualquotient}). Therefore $(G^\vee,\xi_G^\vee)$ is an adelic line bundle on $S$. Finally, by Proposition \ref{Pro:dual adelic line bundle} and the fact that $\xi_G=\xi_G^{\vee\vee}$ (where we identify $G$ with $G^{\vee\vee}$), we obtain that  $(G,\xi_G)$ is an adelic line bundle on $S$.
\end{proof}

\begin{prop}\phantomsection\label{Pro:critereadelicvb}\begin{enumerate}[label=\rm(\arabic*)]
\item\label{Item: dual adelic vector bundle 01} Let $(E,\xi)$ be an adelic vector bundle (resp. a strongly adelic vector bundle) on $S$, $F$ be a vector subspace of $E$ and $\xi_F$ be the restriction of $\xi$ {to} $F$. If the norm family  $\xi_F^\vee$ is $\mathcal A$-measurable, then $(F,\xi_F)$ is an adelic vector bundle (resp. a strongly adelic vector bundle) on $S$.
\item\label{Item: quotient adelic vector bundle} Let $(E,\xi)$ be an adelic vector bundle (resp. a strongly adelic vector bundle) on $S$, $G$ be a quotient vector space of $E$ and $\xi_G$ be the quotient norm family of $\xi$. If the norm {family} $\xi_G$ is $\mathcal A$-measurable, then $(G,\xi_G)$ is an adelic vector bundle (resp. a strongly adelic vector bundle) on $S$.
\item\label{Item: dual adelic vector bundle 02} Let $(E,\xi)$ be an adelic vector bundle on $S$. Assume that the norm family $\xi^{\vee\vee}$ is measurable. Then $(E^\vee,\xi^{\vee})$ is a strongly adelic vector bundle on $S$.
\item\label{Item: tensor adelic vector bundle} Let $(E,\xi_E)$ and $(F,\xi_F)$ be adelic vector bundles on $S$. If the norm families $\xi_E\otimes_\varepsilon \xi_F$ and $(\xi_E\otimes_\varepsilon \xi_F)^\vee$ are $\mathcal A$-measurable, then $(E\otimes F,\xi_E\otimes_{\varepsilon}\xi_F)$ is a strongly adelic vector bundle on $S$. Similarly,  $(E\otimes F,\xi_E\otimes_{\varepsilon,\pi}\xi_F)$ is a strongly adelic vector bundle on $S$ provided that the both norm families $\xi_E\otimes_{\varepsilon,\pi}\xi_F$ and $(\xi_E\otimes_{\varepsilon,\pi}\xi_F)^\vee$ are measurable. If in addition $\xi_E$ and $\xi_F$ are both Hermitian, and if both norm families $\xi_E\otimes\xi_F$ and $(\xi_E\otimes\xi_F)^\vee$ are $\mathcal A$-measurable, then the orthogonal tensor product $(E\otimes F,\xi_E\otimes\xi_F)$ is a Hermitian adelic vector bundle on $S$.
\item\label{Item:detadelic} Let $(E,\xi)$ be an adelic vector bundle on $S$. If $\det(\xi)$ is $\mathcal A$-measurable then $(\det(E),\det(\xi))$ is an adelic line bundle on $S$.
\end{enumerate}
\end{prop}
\begin{proof}
These assertions are direct consequences of Proposition \ref{Pro:dominancealgebraic}. We just mention below some particular points. For the assertion \ref{Item: dual adelic vector bundle 01}, since $\xi$ is $\mathcal A$-measurable, by definition $\xi_F$ is also measurable. For the assertion \ref{Item: quotient adelic vector bundle}, since $\xi^\vee$ is $\mathcal A$-measurable, and $\xi_G^\vee$ identifies with the restriction of $\xi^\vee$ {to} $G^\vee$, it is also $\mathcal A$-measurable. For the last assertion, since $\det(E)$ is of rank $1$ on $K$, the $\mathcal A$-measurability of $\det(\xi)$ implies that of its dual.
\end{proof}

\begin{coro}
Let $(E_1,\xi_1)$ and $(E_2,\xi_2)$ be adelic line bundles on $S$. Then the tensor product $(E_1\otimes E_2,\xi_1\otimes\xi_2)$ is also an adelic line bundle on $S$.
\end{coro}
\begin{proof}
By Proposition \ref{Pro:mesurability} \ref{Item:tensor measurability}, the tensor product norm family $\xi_1\otimes\xi_2$ is $\mathcal A$-measurable. By Proposition \ref{Pro:critereadelicvb} \ref{Item: tensor adelic vector bundle}, we obtain that $(E_1\otimes E_2,\xi_1\otimes\xi_2)$ is an adelic line bundle on $S$.
\end{proof}

\begin{rema}\label{Rem: measuarbility under supplementay condition}
By using the measurability results obtained in the previous subsection (notably Proposition \ref{Pro:mesurabilityofquotient}), we obtain that the assertions of Proposition \ref{Pro:critereadelicvb} remain true without measurability assumptions, if the $\sigma$-algebra $\mathcal A$ is discrete, or if the field $K$ admits a countable subfield which is dense in all completions $K_\omega$, $\omega\in\Omega$.
\end{rema}

In the case of direct sums, the measurability result in Proposition \ref{Pro:mesurabilitedesomme} leads to the following criterion (without any condition on $K$).

\begin{prop}
Let $(E,\xi_E)$ and $(F,\xi_F)$ be adelic vector bundles (resp. a strongly adelic vector bundle) on $S$, and $\psi:(\omega\in\Omega)\mapsto \psi_\omega\in\mathscr S$ be an $\mathcal A$-measurable map. We assume that there exists a measurable subset $\Omega'$ of $\Omega$ such that $\nu(\Omega')<+\infty$ and that $\psi_\omega=\psi_0$ on $\Omega\setminus\Omega'$, where $\psi_0$ denotes the function in $\mathscr S$ sending $t\in[0,1]$ to $\max\{t,1-t\}$. Then $(E\oplus F,\xi_E\oplus_{\psi}\xi_F)$ is an adelic vector bundle (resp. a strongly adelic vector bundle) on $S$.
\end{prop}
\begin{proof}
Since $(E,\xi_E)$ and $(F,\xi_F)$ are adelic vector bundles (resp. strongly adelic vector bundles) on $S$, the norm families $\xi_E$, $\xi_F$, $\xi_E^\vee$ and $\xi_F^\vee$ are all $\mathcal A$-measurable, and the norm families $\xi_E$ and $\xi_F$ are dominated (resp. strongly dominated). 

By Proposition \ref{Pro:mesurabilitedesomme}, the $\psi$-direct sum $\xi_E\oplus_\psi\xi_F$ is also $\mathcal A$-measurable. Let $\psi'=\{\psi'_\omega\}_{\omega\in\Omega}$ be the family in $\mathscr S$ such that $\psi_\omega=\psi_0$ on $\Omega\setminus\Omega_\infty$ and $\psi_\omega'=\psi_{\omega,*}$ (see Definition \ref{Def:dualdirecsum}) on $\Omega_\infty$, then one has
\[(\xi_E\oplus_\psi\xi_F)^\vee=\xi_E^\vee\oplus_{\psi'}\xi_F^\vee.\]
Note that the map from $\mathscr S$ to itself sending $\varphi\in\mathscr S$ to $\varphi_*$ is continuous. This is a consequence of \eqref{Equ:fonctioncorrespondsant} and Proposition \ref{Pro:distanceofoperatornorms}. Therefore, the map $\psi'$ is also $\mathcal A$-measurable. Still by Proposition \ref{Pro:mesurabilitedesomme}, we obtain that the norm family $(\xi_E\oplus\xi_F)^\vee$ is $\mathcal A$-measurable.

By Proposition \ref{Pro:dominancealgebraic} \ref{Item:sommedom}, the norm family $\xi_E\oplus_\psi\xi_F$ is dominated (resp. strongly dominated). Therefore $(E\oplus F,\xi_E\oplus_{\psi}\xi_F)$ is an adelic vector bundle (resp. strongly adelic vector bundle) on $S$.
\end{proof}

\subsection{Examples} In this subsection, we present several examples of adelic vector bundles, which include most classical constructions.

\subsubsection*{Torsion free coherent sheaves}
Let $k$ be a field and $X$ be a normal projective scheme of dimension $d\geqslant 1$ over $\Spec k$, equipped with a family $\{D_i\}_{i=1}^{d-1}$ of ample divisors on $X$. Let $K=k(X)$ be the field of rational functions on $X$ and $\Omega=X^{(1)}$, equipped with the discrete $\sigma$-algebra. We have seen in \S\ref{Subsec:polarizedvar} that $S=(K,(\Omega,\mathcal A,\nu),\phi)$ is an adelic curve, where the measure $\nu$ is defined as
\[\forall\,Y\in\Omega=X^{(1)},\quad\nu(\{Y\})=\deg(D_1\cdots D_{d-1}\cap[Y])\]
and the map $\phi:\Omega\rightarrow M_K$ sends $Y\in\Omega$ to $|\ndot|_Y=\mathrm{e}^{-\ord_Y(\ndot)}$.

Let $\mathcal E$ be a torsion-free (namely the canonical homomorphism $\mathcal E\rightarrow\mathcal E^{\vee\vee}$ is injective) coherent sheaf on $X$ and $E:=\mathcal E\otimes_{\mathcal O_X}K$. The latter is a finite-dimensional vector space over $K$. Moreover, for any $Y\in\Omega$, $\mathcal E\otimes_{\mathcal O_X}\mathcal O_{X,Y}$ is a torsion-free module of finite type over the discrete valuation ring $\mathcal O_{X,Y}$ (the local ring of $X$ at the generic point of $Y$), hence is a free $\mathcal O_{X,Y}$-module of finite rank. We define a norm $\|\ndot\|_Y$ on $E\otimes_K{K}_Y$ as follows
\[\forall\,s\in E\otimes_KK_Y,\quad\|s\|_Y:=\inf\{|a|_Y\,:\,a\in K_Y^{\times},\;a^{-1}s\in\mathcal E\otimes_{\mathcal O_X}\widehat{\mathcal O}_{X,Y}\},\]
where $\widehat{\mathcal O}_{X,Y}$ is the completion of $\mathcal O_{X,Y}$, which identifies with the valuation ring of $K_Y$. This norm is clearly ultrametric. Thus we obtain a Hermitian norm family  in $\mathcal N_E$, which we denote by $\xi_{\mathcal E}$. Note that the dual norm family $\xi_{\mathcal E}^\vee$ identifies with $\xi_{\mathcal E^\vee}$, where $\mathcal E^{\vee}$ denotes the dual $\mathcal O_X$-module of $\mathcal E$.

Since torsion-free coherent sheaves are locally free on codimension $1$, we obtain that, for any basis $\boldsymbol{e}$ of $E$, the norms $\|\ndot\|_Y$ and $\|\ndot\|_{\boldsymbol{e},Y}$ are identical for all but a finite number of $Y\in\Omega$. Therefore, the couple $(E,\xi_{\mathcal E})$ is a strongly adelic vector bundle on $S$.  

\subsubsection*{Hermitian vector bundles on an arithmetic curve}
Let $K$ be a number field and $\mathcal O_K$ be the ring of algebraic integers in $K$. Recall that a Hermitian vector bundle on $\Spec\mathcal O_K$ is by definition a couple $(\mathcal E,\{\|\ndot\|_\sigma\}_{\sigma:K\rightarrow\mathbb C})$, where $\mathcal E$ is a projective $\mathcal O_K$-module of finite rank, and for any embedding $\sigma:K\rightarrow\mathbb C$, $\|\ndot\|_\sigma$ is a Hermitian norm on $\mathcal E\otimes_{\mathcal O_K,\sigma}\mathbb C$. We also require that the norms $\{\|\ndot\|_\sigma\}_{\sigma:K\rightarrow\mathbb C}$ are invariant under the complex conjugation, namely for $s_1,\ldots,s_n$ in $\mathcal E$, $\lambda_1,\ldots,\lambda_n$ in $\mathbb C$, and $\sigma:K\rightarrow\mathbb C$, one has
\[\|\lambda_1\otimes s_1+\cdots+\lambda_n\otimes s_n\|_{\sigma}=\|\overline{\lambda}_1\otimes s_1+\cdots+\overline{\lambda}_n\otimes s_n\|_{\overline{\sigma}}.\] We let $E:=\mathcal E\otimes_{\mathcal O_K}K$.

Let $S=(K,(\Omega,\mathcal A,\nu),\phi)$ be the adelic curve associated with the number field $K$, as described in \S\ref{Subsec:Numberfields}. Recall that $\Omega$ is the set of all places of $K$, $\mathcal A$ is the discrete $\sigma$-algebra on $\Omega$ and $\nu(\{\omega\})=[K_\omega:\mathbb Q_\omega]$. 

Recall that any finite place of $K$ is determined by a maximal ideal $\mathfrak p$ of $\mathcal O_K$. Let $\widehat{\mathcal O}_{K,\mathfrak p}$ be the completion of the local ring $\mathcal O_{K,\mathfrak p}$, which is also the valuation ring of $K_{\mathfrak p}$. We construct a norm $\|\ndot\|_{\mathfrak p}$ as follows
\[\forall\,s\in E\otimes_KK_{\mathfrak p},\quad\|s\|_{\mathfrak p}:=\inf\{|a|_{\mathfrak p}\,:\,a\in K_{\mathfrak p}^{\times},\;a^{-1}s\in\mathcal E\otimes_{\mathcal O_K}\widehat{\mathcal O}_{K,\mathfrak p}\}.\]
Let $v$ be an Archimedean place of $K$. Then $v$ corresponds to an embedding $\sigma$ of $K$ into $\mathbb C$, we let $\|\ndot\|_v$ be the restriction of $\|\ndot\|_{\sigma}$ {to} $E\otimes_KK_v$. Note that the condition that the norms $\{\|\ndot\|_\sigma\}_{\sigma:K\rightarrow\mathbb C}$ are invariant under the complex conjugation ensures that the norm $\|\ndot\|_v$ does not depend on the choice of the embedding $\sigma:K\rightarrow\mathbb C$ corresponding to $v$. Thus we obtain a norm family $\xi=\{\|\ndot\|_v\}_{v\in\Omega}$ in $\mathcal N_E$. Since $\mathcal E$ is a locally free sheaf, we obtain that, for any basis $\boldsymbol{e}$ of $E$ over $K$, one has $\|\ndot\|_v=\|\ndot\|_{\boldsymbol{e},v}$ for all but a finite number of $v$.

\subsubsection*{Ultrametrically normed vector space over a trivially valued field}
Let $K$ be an arbitrary field and $\Omega$ be the one point set $\{\omega\}$. Let $|\ndot|_{\omega}$ be the trivial absolute value on $K$. We then obtain an adelic curve $S$ by taking the discrete $\sigma$-algebra $\mathcal A$ on $\Omega$ and the measure $\nu$ on $(\Omega,\mathcal A)$ such that $\nu(\{\omega\})=1$. Then any ultrametrically normed finite-dimensional vector space over $K$ is a strongly adelic vector bundle on $S$.

\begin{rema}
Let $K$ be a number field and $E$ be a finite-dimensional vector space over $K$. In \cite[\S3]{Gaudron08}, a structure of adelic vector bundle on $E$ has been defined as a norm family $\xi\in\mathcal N_E$ such that, for all but finitely many $\omega\in\Omega$ (where $\Omega$ denotes the set of all places of $K$), the norm $\xi_\omega$ is induced by a projective $\mathcal O_K$-module of finite type $\mathcal E$. Clearly such a structure of adelic vector bundle is a dominated norm family. We denote by $\mathcal D_E^\circ$ the subset of $\mathcal D_E$ consisting of all structures of adelic vector bundles in the sense of \cite{Gaudron08}. We claim that $\mathcal D_E^\circ$ is dense in $\mathcal D_E$ (with respect to the metric $\mathrm{dist}(\ndot,\ndot)$ defined in Remark \ref{Rem:metricDE}). In other words, given a dominated norm family $\xi$ on $E$, there exists a sequence $\{\xi_n\}_{n\in\mathbb N}$ in $\mathcal D_E^\circ$ which converges to $\xi$. In fact, we can choose an arbitrary element $\xi_0$ in $\mathcal D_E^\circ$. The main point is that, if we modify finitely many norms in the family $\xi_0$, we still obtain a norm family in $\mathcal D_E^\circ$. Since the local distance function $d(\xi,\xi_0)$ is $\nu$-dominated, we can construct a sequence $\{\Omega_n\}_{n\geqslant 1}$ of subsets of $\Omega$, such that $\Omega\setminus\Omega_n$ is a finite set and that 
\begin{equation}\label{Equ:convergestoxi}\lim_{n\rightarrow+\infty}\upint_{\Omega}\indic_{\Omega_n}(\omega)d_\omega(\xi,\xi_0)\,\nu(\mathrm{d}\omega)=0.\end{equation}
We then let $\xi_n$ be the norm family such that
\[\xi_{n,\omega}=\begin{cases}
\xi_{0,\omega},&{\omega\in\Omega_n}\\
\xi_\omega,&\omega\in\Omega\setminus\Omega_n.
\end{cases}\]
Then the sequence $\{\xi_n\}_{n\in\mathbb N}$ is contained in $\mathcal D_E^\circ$ and converges to $\xi$ (see \eqref{Equ:convergestoxi}). Combined with the completeness of the space $\mathcal D_E$ explained in Remark \ref{Rem:metricDE}, we obtain that $\mathcal D_E$ is actually the completion of the metric space $\mathcal D_E^\circ$. 
\end{rema}

\section{Adelic divisors} 
Let $S=(K,(\Omega,\mathcal A,\nu),\phi)$ be an adelic curve. 
We call \emph{adelic divisor}\index{adelic divisor} on $S$ any element $\zeta$ in the vector space $L^1(\Omega,\mathcal A,\nu)$
(see Section~\ref{Subsec: L1 space}). For the reason of customs of arithmetic geometry, we use the notation $\widehat{\mathrm{Div}}_{\mathbb R}(S)$ to denote the vector space $L^1(\Omega,\mathcal A,\nu)$.

If $\zeta$ is an adelic divisor on $S$, we define its \emph{Arakelov degree}\index{Arakelov degree} as
\begin{equation}\label{Equ:degre}{\deg}(\zeta):=\int_{\Omega}\zeta(\omega)\,\nu(d\omega)\in\mathbb R.\end{equation} The function $\deg$ is a continuous linear form on $\widehat{\mathrm{Div}}_{\mathbb R}(S)$. If $a$ is an element of $K^{\times}$, we denote by $\widehat{(a)}$ the adelic divisor represented by the function which sends $\omega\in\Omega$ to $-\ln|a|_\omega$, called the \emph{adelic divisor associated with $a$}.
The map \[\widehat{(\ndot)}:K^{\times}\longrightarrow\widehat{\mathrm{Div}}_{\mathbb R}(S)\] is additive and hence extends to an $\mathbb R$-linear homomorphism from $K^{\times}\otimes_{\mathbb Z}\mathbb R$ to $\widehat{\mathrm{Div}}_{\mathbb R}(S)$, which we denote by $\widehat{(\ndot)}_{\mathbb R}$. The closure of the image of this map is denoted by $\widehat{\mathrm{PDiv}}_{\mathbb R}(S)$ and the elements of this vector space are called \emph{principal adelic divisors}\index{principal adelic divisor}\index{adelic divisor!principal ---}. We denote by $\widehat{\mathrm{Cl}}_{\mathbb R}(S)$ the quotient space $\widehat{\mathrm{Div}}_{\mathbb R}(S)/\widehat{\mathrm{PDiv}}_{\mathbb R}(S)$. Note that it forms actually a Banach space with respect to the quotient norm.
Two adelic divisors lying in the same equivalent class in $\widehat{\mathrm{Cl}}_{\mathbb R}(S)$ are said to be \emph{$\mathbb R$-linearly equivalent}\index{R-linearly equivalent@$\mathbb R$-linearly equivalent}. 

We say that an adelic divisor $\zeta$ on $S$ is \emph{effective}\index{effective}\index{adelic divisor!effective ---} if $\zeta$ is 
$\nu$-almost everywhere non-negative. Denote by $\widehat{\mathrm{Div}}_{\mathbb R}(S)_+$ the cone of all effective adelic divisors on $S$. Clearly, if $\zeta$ is effective, then ${\deg}(\zeta)\geqslant 0$.

Let $S'=(K',(\Omega',\mathcal A',\nu'),\phi')$ be another adelic curve and $\alpha=(\alpha^\#,\alpha_\#,I_\alpha):S'\rightarrow S$ be a morphism of adelic curves (see Section \ref{Sec: Morphism of adelic curve}). If $\zeta$ is an adelic divisor on $S$, which is represented by an element $f\in\mathscr L^1(\Omega,\mathcal A,\nu)$, we denote by $\alpha^*(\zeta)$ the adelic divisor on $S'$ represented by the function $f\circ\alpha_\#$ (the equivalence class of $f\circ\alpha_\#$ does not depend on the choice of the representative $f$ since $\nu$ identifies with the direct image of $\nu'$ by $\alpha_{\#}$). If $\zeta'$ is an adelic divisor on $S'$, we denote by $\alpha_*(\zeta')$ the adelic divisor $I_\alpha(\zeta')$ on $S$. Since $I_\alpha$ is a disintegration kernel of $\nu'$ over $\nu$, one has $\alpha_*(\alpha^*(\zeta))=\zeta$ for any adelic divisor $\zeta$ on $S$.

From now on we assume that $S$ is proper. Then one has
${\deg}(\zeta)=0$ if $\zeta$ is a principal adelic divisor. This is a direct consequence of the product formula and the fact that ${\deg}(\ndot)$ is a continuous linear operator. Therefore the $\mathbb R$-linear map ${\deg}(\ndot)$ induces by passing to quotient a continuous $\mathbb R$-linear map from $\widehat{\mathrm{Cl}}_{\mathbb R}(S)$ to $\mathbb R$ which sends any class $[\zeta]$ to ${\deg}(\zeta)$. We still denote this linear map by ${\deg}(\ndot)$ by abuse of notation.

\section{Arakelov degree and slopes}
The purpose of this section is to generalise the theory of Arakelov degree and slopes to the setting of adelic vector bundles over adelic curves. Throughout the section, let $S=(K,(\Omega,\mathcal A,\nu),\phi)$ be a proper adelic curve. For all subsections except the first one, we assume in addition that, either the $\sigma$-algebra $\mathcal A$ is discrete, or there exists a countable subfield $K_0$ of $K$ which is dense in all $K_\omega$, $\omega\in\Omega$.

\subsection{Arakelov degree of adelic line bundles}
In this subsection, we fix an adelic vector bundle $\overline E=(E,\xi)$ on $S$.

\begin{defi}
If $s$ is a non-zero vector in $E$, by Proposition \ref{Pro:quotientandsubofrank1}, the vector space $Ks$ equipped with the induced norm family forms an adelic line bundle on $S$. In particular, the function $\ln\|s\|:(\omega\in\Omega)\mapsto\ln\|s\|_\omega$ is $\nu$-dominated and $\mathcal A$-measurable, and hence is $\nu$-integrable. 
The adelic divisor given by $(\omega\in\Omega)\mapsto -\ln\|s\|_\omega$ is denoted by $\widehat{\mathrm{div}}_{\xi}(s)$, which is called
the \emph{adelic divisor}\index{adelic divisor} of $s$ with respect to $\xi$.
We define the \emph{Arakelov degree}\index{Arakelov degree} of $s$ with respect to $\xi$ as the Arakelov degree of $\widehat{\mathrm{div}}_{\xi}(s)$, that is,
\[\widehat{\deg}_{\xi}(s):= \deg\left(\widehat{\mathrm{div}}_{\xi}(s)\right) = -\int_{\Omega}\ln\|s\|_\omega\,\nu(\mathrm{d}\omega). \]
Moreover, by the product formula \eqref{Equ:productformula} we obtain that, for any $a\in K^{\times}$, one has
\begin{equation}\label{Equ:upperandlowerdegree}\widehat{\deg}_\xi(as)=\widehat{\deg}_\xi(s).\end{equation}
\end{defi}

\begin{rema}
We assume $E = K^{n}$ and $\xi = \{ \|\ndot\|_{\omega} \}_{\omega \in \Omega}$ is given by
\[
\|(a_1, \ldots, a_n)\|_{\omega} = \max \{ |a_1|_{\omega}, \ldots, |a_n|_{\omega} \}
\]
for each $\omega \in \Omega$.
Then $h_S(s) = -\widehat{\deg}_{\xi}(s)$ for all $s \in E \setminus \{ 0 \}$ (for the definition of $h_S$, see
Defintion~\ref{def:height:adelic:curve}).
\end{rema}

\begin{defi}
Let $(E,\xi)$ be an adelic line bundle on $S$. We call \emph{Arakelov degree}\index{Arakelov degree} of $(E,\xi)$ the number $\widehat{\deg}_\xi(s)$, where $s$ is a non-zero element of $E$. Note that the relation \eqref{Equ:upperandlowerdegree} shows that the definition does not depend on the choice of the non-zero element $s$ in $E$. We denote the Arakelov degrees of $(E,\xi)$ by $\widehat{\deg}(E,\xi)$.
\end{defi}

\begin{prop}\label{Pro:degreedualline}
Let $(E,\xi)$ be an adelic line bundle on $S$. Then $(E^\vee,\xi^\vee)$ is also an adelic line bundle on $S$. Moreover, one has
\begin{equation}\label{Equ:dualdegree}
\widehat{\deg}(E^\vee,\xi^\vee)=-\widehat{\deg}(E,\xi).\end{equation}
\end{prop}
\begin{proof}By Proposition \ref{Pro:dual adelic line bundle}, the couple $(E^\vee,\xi^\vee)$ is also an adelic line bundle, so that the Arakelov degree $\widehat{\deg}(E^\vee,\xi^\vee)$ is well defined.
If $\alpha$ is a non-zero element of $E^\vee$ and $s$ is a non-zero element of $E$ then one has
\[\forall\,\omega\in\Omega,\quad |\alpha(s)|_\omega=\|\alpha\|_{\omega,*}\cdot\|s\|_{\omega}.\]
By the product formula
\[\int_{\Omega}\ln|\alpha(s)|_\omega\,\nu(\mathrm{d}\omega)=0,\]
we obtain the equality \eqref{Equ:dualdegree}.
\end{proof}

\begin{prop}\label{Pro:degreduproduit}
Let $(E_1,\xi_1)$ and $(E_2,\xi_2)$ be adelic line bundles on $S$. Let $E=E_1\otimes_KE_2$ and $\xi=\xi_1\otimes_{\varepsilon}\xi_2$ (which is also equal to $\xi_1\otimes\xi_2$ and $\xi_1\otimes_{\pi}\xi_2$). Then one has
\begin{equation}\label{Equ: degree tensor product}
\widehat{\deg}(E_1\otimes E_2,\xi_1\otimes \xi_2)=\widehat{\deg}(E_1,\xi_1)+\widehat{\deg}(E_2,\xi_2).
\end{equation}
\end{prop}
\begin{proof}
Let $s_1$ and $s_2$ be non-zero elements of $E_1$ and $E_2$, respectively. For any $\omega\in\Omega$, one has
\begin{equation}\label{Equ:additivite}\ln\|s_1\otimes s_2\|_\omega=\ln\|s_1\|_\omega+\ln\|s_2\|_\omega.\end{equation}
By taking the integral with respect to $\nu$, we obtain the equality \eqref{Equ: degree tensor product}.
\end{proof}

\subsection{Arakelov degree of adelic vector bundles}

From now on and until the end of the section, we assume that, either the $\sigma$-algebra $\mathcal A$ is discrete, or there exists a countable subfield $K_0$ of $K$ which is dense in each $K_\omega$, where $\omega\in\Omega$.

\begin{defi}
Let $(E,\xi)$ be an adelic vector bundle on $S$. By Proposition \ref{Pro:critereadelicvb}, $(\det(E),\det(\xi))$ is an adelic line bundle on $S$. We define the \emph{Arakelov degree}\index{Arakelov degree}\index{adelic vector bundle!Arakelov degree} of $(E,\xi)$ as \[\widehat{\deg}(E,\xi):=\widehat{\deg}(\det(E),\det(\xi)).\]  Note that, the Arakelov degree of the zero adelic vector bundle is $0$. By Proposition \ref{Pro:doubledualdet}, one has $\det(\xi)=\det(\xi^{\vee\vee})$. Therefore
\begin{equation}\label{Equ:doublddual}
\widehat{\deg}(E,\xi)=\widehat{\deg}(E,\xi^{\vee\vee}).\end{equation}
\end{defi}

\begin{prop}
Let $(E,\xi)$ be a \emph{Hermitian} adelic vector bundle on $S$. One has
\begin{equation}\label{Equ:degreeduale} \widehat{\deg}(E,\xi)=-\widehat{\deg}(E^\vee,\xi^\vee).
\end{equation}
\end{prop}
\begin{proof}
The determinant of $E^\vee$ is canonically isomorphic to the dual vector space of $\det(E)$, and the norm family $\det(\xi^\vee)$ identifies with $\det(\xi)^\vee$ under this isomorphism (see Proposition \ref{Pro:delta=1}), provided that $\xi$ is Hermitian. Therefore, by Proposition \ref{Pro:degreedualline} we obtain the equalities.
\end{proof}

\begin{defi}
Let $(E,\xi)$ be an adelic vector bundle on $S$. We denote by $\delta(\xi)$ the function on $\Omega$ sending $\omega\in\Omega$ to $\delta_\omega(\xi):=\delta(E_{K_\omega},\|\ndot\|_{\omega})$ (see \S\ref{Subsec:dualdet}). Recall that the function $\delta(\xi)$ is identically $1$ on $\Omega\setminus\Omega_\infty$ (see Proposition \ref{Pro:delta=1}), and takes value in $[1,r^{r/2}]$ on $\Omega_\infty$ (see Proposition \ref{Pro:minorationdelta} and the inequalities \eqref{Equ:majorationdelta} and \eqref{Equ:majorationdelta2}), where $r$ is the rank of $E$ over $K$. In particular, the function $\ln\delta(\xi)$ is $\nu$-dominated since it is bounded and vanishes outside a set of finite measure.

Similarly, we denote by $\Delta(\xi)$ the function on $\Omega$ sending $\omega\in\Omega$ to $\Delta_\omega(\xi):=\Delta(E_{K_\omega},\norm{\ndot}_\omega)$ (see Definition \ref{Def:Delta:V:norm}). This function is bounded from below by the constant function $1$. Moreover, it identifies with the constant function $1$ if $\xi$ is Hermitian.
\end{defi}

\begin{prop}\label{Pro:degdual}
Let $(E,\xi=\{\norm{\ndot}_\omega\}_{\omega\in\Omega})$ be an adelic vector bundle on $S$. Then the function $(\omega\in\Omega)\mapsto\delta_\omega(\xi)$ is $\mathcal A$-measurable and its logarithm is integrable with respect to $\nu$. It is a constant function of value $1$ when $\xi$ is Hermitian. Moreover, the following relations hold
\begin{multline}\label{Equ:encadrementdelta}\qquad
0\leqslant\widehat{\deg}(E,\xi)+\widehat{\deg}(E^\vee,\xi^\vee)\\
\hskip-2em =\int_\Omega \ln(\delta_\omega(\xi))\,\nu(d\omega)\leqslant \frac{1}{2}\rang(E)\ln(\rang(E))\nu(\Omega_\infty).
\end{multline}
In particular, one has $\widehat{\deg}(E,\xi)+\widehat{\deg}(E^\vee,\xi^\vee)=0$ if for any $\omega\in\Omega_\infty$ the norm $\norm{\ndot}_\omega$ is induced by an inner product.
\end{prop}
\begin{proof} By Proposition \ref{Pro:critereadelicvb}, we obtain that both couples $(\det(E),\det(\xi))$ and $(\det(E^\vee),\det(\xi^\vee))$ are adelic line bundles on $S$.
Let $\eta$ be a non-zero element in $\det(E)$ and $\eta^\vee$ be its dual element in $\det(\xi)$. By \eqref{Equ:definitiondedelta}, one has
\[(-\ln\|\eta\|_{\omega,\det})+(-\ln\|\eta^\vee\|_{\omega,*,\det})=\ln\delta_\omega(\xi).\] Therefore the function $(\omega\in\Omega)\mapsto\delta_\omega(\xi)$ is $\mathcal A$-measurable. Moreover, by Proposition \ref{Pro:delta=1} we obtain that $\delta_\omega(\xi)=1$ if $\omega\in\Omega\setminus\Omega_\infty$ or if $\omega\in\Omega_\infty$ and the norm $\|\ndot\|_\omega$ is induced by an inner product. Therefore \eqref{Equ:encadrementdelta} follows from the inequalities \[0\leqslant\ln\delta(\xi)\leqslant \frac{1}{2}\rang(E)\ln(\rang(E))\indic_{\Omega_\infty},\]
which also implies the $\nu$-integrability of $\ln\delta(\xi)$.
\end{proof}

{
\begin{prop}
\label{Pro: integrability of Delta}Let $(E,\xi)$ be a strongly adelic vector bundle on $S$. The function $\ln\Delta(\xi)$ is $\nu$-dominated. 
\end{prop}
\begin{proof}
Let $\boldsymbol{e}=\{e_i\}_{i=1}^r$ be a basis of $E$ over $K$. By Corollary \ref{Cor:dominatedanddist}, the local distance function $(\omega\in\Omega)\mapsto d_\omega(\xi,\xi_{\boldsymbol{e}})$ is $\nu$-dominated. We write the norm families $\xi$ and $\xi_{\boldsymbol{e}}$ in the form of $\xi=\{\norm{\ndot}_\omega\}_{\omega\in\Omega}$ and $\xi_{\boldsymbol{e}}=\{\norm{\ndot}_{\boldsymbol{e},\omega}\}_{\omega\in\Omega}$, respectively. Let $\omega\in\Omega$. If $\norm{\ndot}_{h,\omega}$ is a norm on $E_{K_\omega}$ bounded from below by $\norm{\ndot}_{\boldsymbol{e},\omega}$, which is either ultrametric or induced by an inner product, then the norm $\mathrm{e}^{d_\omega(\xi,\xi_{\boldsymbol{e}})}\norm{\ndot}_{h,\omega}$ is bounded from below by $\norm{\ndot}_{\omega}$. This norm is also ultrametric or induced by an inner product. Therefore we obtain that \[\ln\Delta_\omega(\xi)\leqslant \ln\Delta_\omega(\xi_{\boldsymbol{e}})+d_\omega(\xi,\xi_{\boldsymbol{e}})\leqslant \frac 12 r\ln(r)\indic_{\Omega_\infty}(\omega)+d_\omega(\xi,\xi_{\boldsymbol{e}}),\]
where the second inequality comes from \eqref{Equ:majorationdelta}.
Since $\nu(\Omega_\infty)<+\infty$ (see Proposition \ref{Pro:mesurabilitekappa}), we obtain that the function $\ln\Delta(\xi)$ is $\nu$-dominated.  
\end{proof}

}

\begin{defi}
Let $\overline E=(E,\xi)$ be an adelic vector bundle on $S$. We denote by $\delta(\overline E)$ the integral $\int_\Omega\ln(\delta_\omega(\xi))\,\nu(\mathrm{d}\omega)$, which is also equal to $\widehat{\deg}(E,\xi)+\widehat{\deg}(E^\vee,\xi^\vee)$ (see Proposition \ref{Pro:degdual}). We denote by $\Delta(\overline E)$ the lower integral $\lowint_{\Omega}\ln(\Delta_\omega(\xi))\,\nu(\mathrm{d}\omega)$, which takes value in $[0,+\infty]$. It is finite once the function $(\omega\in\Omega)\rightarrow d_\omega(\xi,\xi^{\vee\vee})$ is $\nu$-dominated (see Proposition \ref{Pro: integrability of Delta}), namely $\overline E$ is a strongly adelic vector bundle.
\end{defi}

Let $E$, $F$ and $G$ be  vector spaces of finite rank over $K$, and $\xi_E$, $\xi_F$ and $\xi_G$ be norm families in $\mathcal N_E$, $\mathcal N_F$ and $\mathcal N_G$ respectively. We say that a diagram
\[\xymatrix{\relax 0\ar[r]&(F,\xi_F)\ar[r]^-f&(E,\xi_E)\ar[r]^-g& (G,\xi_G)\ar[r]&0}\]
is an \emph{exact sequence}\index{exact sequence} if 
\[\xymatrix{\relax 0\ar[r]&F\ar[r]^-f& E\ar[r]^-g& G\ar[r]&0}\] 
is an exact sequence of vector spaces over $K$
and if the norm family $\xi_F$ is the restriction of $\xi_E$ {to} $F$, and the norm family $\xi_G$ is the quotient of the norm family of $\xi_E$ on $G$ (see \S\ref{Subsec:Norm families}, page \pageref{page:restricquot}). Note that if $\xi_E$ is Hermitian, then both norm families $\xi_F$ and $\xi_G$ are Hermitian.

\begin{prop}
\label{Pro:suiteexactedeg}
Let $(E,\xi)$ be an adelic vector bundle over $S$ and
\[0=E_0\subseteq E_1\subseteq\ldots\subseteq E_n=E\]
be a flag of vector subspaces of $E$. For any $i\in\{1,\ldots,n\}$, let $\xi_i$ be the restriction of $\xi$ {to} $E_i$ and let $\eta_i$ be the quotient norm family of $\xi_i$ on $E_i/E_{i-1}$. Then one has
\begin{equation}\label{Equ:deglowerbounde}\sum_{i=1}^n\widehat{\deg}(E_i/E_{i-1},\eta_i)\leqslant\widehat{\deg}(E,\xi)\end{equation}
and
\begin{equation}
\label{Equ:degupperbound Delta}
\widehat{\deg}(E,\xi)-\Delta(E,\xi)\leqslant\sum_{i=1}^n\Big(\widehat{\deg}(E_i/E_{i-1},\eta_i)-\Delta(E_i/E_{i-1},\eta_i)\Big).
\end{equation}
If in addition $\xi$ is ultrametric on $\Omega\setminus\Omega_\infty$, then
\begin{equation}
\label{Equ:degupperbounde}
\widehat{\deg}(E,\xi)-\delta(E,\xi)\leqslant\sum_{i=1}^n\Big(\widehat{\deg}(E_i/E_{i-1},\eta_i)-\delta(E_i/E_{i-1},\eta_i)\Big).\end{equation}
In particular, if $\xi$ is Hermitian, then one has
\begin{equation}\label{Equ:additivity degree}\widehat{\deg}(E,\xi)=\sum_{i=1}^n\widehat{\deg}(E_i/E_{i-1},\eta_i).\end{equation}
\end{prop}
\begin{proof}
For any $i\in\{1,\ldots,n\}$, we have an exact sequence
\[\xymatrix{0\ar[r]&(E_{i-1},\xi_{i-1})\ar[r]&(E_i,\xi_i)\ar[r]&(E_i/E_{i-1},\eta_i)\ar[r]&0}.\]
In particular, one has a canonical isomorphism (see \cite{Bourbaki70} Chapter III, \S 7, no.7, Proposition 10) 
\[\det(E_i)\cong \det(E_{i-1})\otimes\det(E_i/E_{i-1}).\] 
For any $i\in\{1,\ldots,n\}$, we pick a non-zero element $\alpha_i\in\det(E_i/E_{i-1})$ and let $\beta_i=\alpha_1\otimes\cdots\otimes\alpha_i$, viewed as an element in $\det(E_i)$. By convention, let $\beta_0$ be the element $1\in\det(E_0)\cong K$. {By 
Corollary \ref{Cor:exactsequencenorm} and Proposition \ref{Pro:minormationdenorme}, one has
\[\ln\|\alpha_i\|_\omega+\ln\|\beta_{i-1}\|_\omega+\ln\frac{\Delta_\omega(\xi_{i-1})\Delta_\omega(\eta_i)}{\Delta_\omega(\xi_i)}\leqslant\ln\|\beta_i\|_{\omega}\leqslant\ln\|\alpha_i\|_\omega+\ln\|\beta_{i-1}\|_{\omega}.\]
If $\xi$ is ultrametric on $\Omega\setminus\Omega_\infty$, then by Proposition \ref{Pro:exactesequenceanddelta}, one has
\[\ln\|\alpha_i\|_\omega+\ln\|\beta_{i-1}\|_\omega+\ln\frac{\delta_\omega(\xi_{i-1})\delta_\omega(\eta_i)}{\delta_\omega(\xi_i)}\leqslant\ln\|\beta_i\|_{\omega},\]
Taking the sum with respect to $i$, we obtain 
\begin{equation}\label{Equ: encadrement de ln betan}\sum_{i=1}^n-\ln\|\alpha_i\|_\omega\leqslant -\ln\|\beta_n\|_\omega\leqslant \Big(\sum_{i=1}^n-\ln\|\alpha_i\|_\omega\Big)+\ln\Delta_\omega(\xi)- \Big(\sum_{i=1}^n\ln\Delta_\omega(\eta_i)\Big)\end{equation}
and, in the case where $\xi$ is ultrametric on $\Omega\setminus\Omega_\infty$,  \[-\ln\|\beta_n\|_\omega\leqslant \Big(\sum_{i=1}^n-\ln\|\alpha_i\|_\omega\Big)+\ln\delta_\omega(\xi)- \Big(\sum_{i=1}^n\ln\delta_\omega(\eta_i)\Big).\]
}
By taking the  integrals with respect to $\nu$, we obtain the inequalities \eqref{Equ:deglowerbounde}  and \eqref{Equ:degupperbounde}. Moreover, \eqref{Equ: encadrement de ln betan} leads to 
\[-\ln\norm{\beta_n}_\omega+\sum_{i=1}^n\ln\Delta_\omega(\eta_i)\leqslant\Big(\sum_{i=1}^n-\ln\|\alpha_i\|_\omega\Big)+\ln\Delta_\omega(\xi).\]
Taking the lower integrals, by Proposition \ref{Pro:translationbyintegrable} and the inequality \eqref{Equ: superadditive} we obtain \eqref{Equ:degupperbound Delta}.

If $\xi$ is a Hermitian norm family, then each $\eta_i$ is a Hermitian norm family, $i\in\{1,\ldots,n\}$. By Proposition \ref{Pro:degdual}, all functions $\ln\delta(\xi)$ and $\ln\delta(\eta_i)$ vanish. Therefore the equality \eqref{Equ:additivity degree} holds.
\end{proof}

\begin{prop}\label{Pro: Arakelov degree preserved by extension of scalars}
Let $(E,\xi)$ be an adelic vector bundle on $S$. If $L/K$ is an algebraic extension of fields, then one has $\det(\xi_L)=\det(\xi)_L$. In particular, 
$\widehat{\deg}(E,\xi)=\widehat{\deg}(E_L,\xi_L)$.
\end{prop}
\begin{proof}
The relation $\det(\xi_L)=(\det\xi)_L$ comes from 
\ref{Pro:extensionofdet:epsilon} (for the non-Archimedean case) and \ref{Pro:extensionofdet:pi} (for the Archimedean case)
in Proposition~\ref{Pro:extensionofdet}.

Let $\alpha$ be a non-zero element of $\det(E)$. For any $x\in\Omega_L$ and $\omega=\pi_{L/K}(x)$,  one has
$\ln\|\alpha\|_x=\ln\|\alpha\|_\omega$.
By taking the integral with respect to $\nu_L$, we obtain $\widehat{\deg}(E,\xi)=\widehat{\deg}(E_L,\xi_L)$.
\end{proof}

\begin{defi}\label{Def:heightofmap}Let $(E,\xi_E)$ and $(F,\xi_F)$ be adelic vector bundles on $S$. Let $f:E\rightarrow F$ be a $K$-linear map. If $f$ is non-zero, we define
\emph{local height function}\index{local height function} of $f$ as the real-valued function on $\Omega$ sending $\omega\in\Omega$ to $\ln\|f_{K_\omega}\|$, where $f_{K_\omega}$ is the $K_\omega$-linear map $E_{K_\omega}\rightarrow F_{K_\omega}$ induced by $f$, and $\|f_{K_\omega}\|$ is its operator norm.
\end{defi}

{
\begin{prop}\label{Pro: integrability of local height}
Let $(E,\xi_E)$ and $(F,\xi_F)$ be adelic vector bundles, and $f:E\rightarrow F$ be a $K$-linear map. The local height function of $f$ is $\mathcal A$-measurable. If $(E,\xi_E)$ and $(F,\xi_F)$ are strongly adelic vector bundles, then the local height function of $f$ is $\nu$-dominated.
\end{prop}
\begin{proof}
We first prove the measurability of the local height function. The statement is trivial when the $\sigma$-algebra $\mathcal A$ is discrete. In the following, we prove the measurability under the hypothesis that there exists a countable subfield $K_0$ of $K$ which is dense in each $K_\omega$. In this case there exists a countable sub-$K_0$-module $E_0$ of $E$ which generates $E$ as a vector space over $K$. For any $\omega\in\Omega$, viewed as a subset of $E_{K_\omega}$, $E_0$ is dense. Therefore one has
\[\norm{f_{K_\omega}}=\inf_{x\in E_0\setminus\{0\}}\frac{\norm{f(x)}_{F,\omega}}{\norm{x}_{E,\omega}.}\]
Since $(E,\xi_E)$ and $(F,\xi_F)$ are adelic vector bundles, for any $x\in E_0$, the functions $(\omega\in\Omega)\mapsto\norm{x}_{E,\omega}$ and $(\omega\in\Omega)\mapsto \norm{f(x)}_{F,\omega}$ are $\mathcal A$-measurable. Therefore the function $(\omega\in\Omega)\mapsto \ln\norm{f_{K_\omega}}$ is $\mathcal A$-measurable.

We now proceed with the proof of the dominancy of the function. We consider $f$ as an element of $E^\vee\otimes F$ and equip this vector space with the norm family $\xi_E^\vee\otimes_\varepsilon \xi_F$ denoted by $\{\norm{\ndot}_{\omega,\varepsilon}\}_{\omega\in\Omega}$. By Remark \ref{Rem:operateureps}, the norm of \[f_{K_\omega}:(E_{K_\omega},\norm{\ndot}_{E,\omega,**})\rightarrow (F_{K_\omega},\norm{\ndot}_{F,\omega,**})\] identifies with $\norm{f}_{\omega,\varepsilon}$. By Proposition \ref{Pro:dominancealgebraic} \ref{Item:dualdom} and \ref{Item:tensordom}, the norm family $\xi_E^\vee\otimes_\varepsilon \xi_F$ is dominated. Hence the function $(\omega\in\Omega)\mapsto \ln\norm{f}_{\omega,\varepsilon}$ is $\nu$-dominated.  Note that one has
\[\big|\ln\norm{f_{K_\omega}}-\ln\norm{f}_{\omega,\varepsilon}\big|\leqslant d_\omega(\xi_E,\xi_E^{\vee\vee}) +d_\omega(\xi_F,\xi_F^{\vee\vee}),\]
where $\norm{f_{K_\omega}}$ denotes the operator norm of $f_{K_\omega}:(E_{K_\omega},\norm{\ndot}_{E,\omega})\rightarrow (F_{K_\omega},\norm{\ndot}_{F,\omega})$. As the local distance functions $(\omega\in\Omega)\mapsto d_\omega(\xi_E,\xi_E^{\vee\vee})$ and $(\omega\in\Omega)\mapsto d_\omega(\xi_F,\xi_F^{\vee\vee})$ are $\nu$-dominated, we obtain that the function $(\omega\in\Omega)\mapsto\ln\norm{f_{K_\omega}}$ is $\nu$-dominated. The proposition is thus proved. 
\end{proof}

\begin{defi}
Let $(E,\xi_E)$ and $(F,\xi_F)$ be adelic vector bundles, and $f:E\rightarrow F$ be a $K$-linear map. We define the \emph{height}\index{height} $h(f)$ of $f$ as the lower integral
\[\lowint_{\Omega}\ln\|f_{K_\omega}\|\,\nu(\mathrm{d}\omega).\] 
By Remark \ref{Rem:operateureps}, in the case where $\xi_E$ and $\xi_F$ are ultrametric on $\Omega\setminus\Omega_\infty$, one has
\[{h}(f)=-\widehat{\deg}_{\xi_E^\vee\otimes_\varepsilon\xi_F}(f).\]
\end{defi}
}

\begin{prop}\label{Pro:slopeinequality1}
Let $(E_1,\xi_1)$ and $(E_2,\xi_2)$ be adelic vector bundles on $S$ and $f:E_1\rightarrow E_2$ be a $K$-linear isomorphism. One has
\begin{equation}\label{Equ:relationdeg}
\widehat{\deg}(E_1,\xi_1)-\widehat{\deg}(E_2,\xi_2)={h}(\det(f)).
\end{equation}
In particular,
\begin{equation}\label{Equ: inequality of degrees}\widehat{\deg}(E_1,\xi_1)\leqslant\widehat{\deg}(E_2,\xi_2)+rh(f).\end{equation}
\end{prop} 
\begin{proof}
By definition one has
\begin{equation}\begin{split}{h}(\det(f))&=-\widehat{\deg}\big(\det(E_1)^\vee\otimes\det(E_2),\det(\xi_1)^\vee\otimes\det(\xi_2)\big)\\
&=-\widehat{\deg}(E_1,\xi_1)+\widehat{\deg}(E_2,\xi_2),
\end{split}\end{equation}
where the second equality comes from Propositions \ref{Pro:degreedualline} and \ref{Pro:degreduproduit}.
Finally the inequality \eqref{Equ: inequality of degrees} is a consequence of \eqref{Equ:relationdeg} and Proposition~\ref{Pro:Hadamard}.
\end{proof}

\subsection{Arakelov degree of tensor adelic vector bundles}

Let $\overline E=(E,\xi_E)$ and $\overline F=(F,\xi_F)$ be adelic vector bundles on $S$. We denote by $\overline E\otimes_{\varepsilon,\pi}\overline F$ the couple $(E\otimes_KF,\xi_E\otimes_{\varepsilon,\pi}\xi_F)$, called the $\varepsilon,\pi$-tensor product of $\overline E$ and $\overline F$. By Proposition \ref{Pro:critereadelicvb} (see also Remark \ref{Rem: measuarbility under supplementay condition}), $\overline E\otimes_{\varepsilon,\pi}\overline F$ is an adelic vector bundle on $S$. If both $\overline E$ and $\overline F$ are Hermitian  adelic vector bundles, we denote by $\overline E\otimes\overline F$ the couple $(E\otimes_KF,\xi_E\otimes\xi_F)$, called the \emph{orthogonal tensor product}\index{orthogonal tensor product} of $\overline E$ and $\overline F$.  By Proposition \ref{Pro:critereadelicvb}, $\overline E\otimes\overline F$ is a Hermitian adelic vector bundle on $S$.

\begin{prop}\label{Pro: degre of tensor}
Let $\overline E=(E,\xi_E)$ and $\overline F=(F,\xi_F)$ be adelic vector bundles on $S$. One has
\begin{equation}\label{Equ: degree of tensors 1}\widehat{\deg}(\overline E\otimes_{\varepsilon,\pi}\overline F)=\rang(F)\,\widehat{\deg}(\overline E)+\rang(E)\,\widehat{\deg}(\overline F).\end{equation}
If $\overline E$ and $\overline F$ are Hermitian  adelic vector bundles, then one has
\begin{equation}\label{Equ: degree of tensors 2}\widehat{\deg}(\overline E\otimes \overline F)=\rang(F)\,\widehat{\deg}(\overline E)+\rang(E)\,\widehat{\deg}(\overline F).\end{equation}
\end{prop}
\begin{proof}Let $n$ and $m$ be the ranks of $E$ and $F$ over $K$ respectively. 
By Propositions \ref{Pro: tensor product and deteminant pi} and \ref{Pro: tensor product and deteminant epsilon}, under the canonical isomorphism \[\det(E)^{\otimes m}\otimes\det(F)^{\otimes n}\cong \det(E\otimes_KF),\] the norm family $\det(\xi_E)^{\otimes m}\otimes\det(\xi_F)^{\otimes n}$ identifies with $\det(\xi_E\otimes_{\varepsilon,\pi}\xi_F)$. Therefore the equality \eqref{Equ: degree of tensors 1} results from Proposition \ref{Pro:degreduproduit}. 

The equality \eqref{Equ: degree of tensors 2} can be proved in a similar way by using Propositions \ref{Pro: tensor product and deteminant HS} and \ref{Pro: tensor product and deteminant epsilon}.
\end{proof}

\subsection{Positive degree} Let $(E,\xi)$  be an adelic vector bundle on $S$. We define the \emph{positive degree}\index{positive degree}\index{adelic vector bundle!positive degree} of $(E,\xi)$ as
\[\widehat{\deg}_+(E,\xi):=\sup_{F\subseteq E}\widehat{\deg}(F,\xi_F),\]
where $F$ runs over the set of all vector subspaces of $E$, and $\xi_F$ denotes the restriction of $\xi$ {to} $F$. Clearly one has $\widehat{\deg}(E,\xi)\geqslant 0$.

\begin{prop}\label{prop:positive:degree:basic:properties}
Let $(E,\xi_E)$ be an adelic vector bundle on $S$, $F$ be a vector subspace of $E$ and $G$ be the quotient space of $E$ by $F$. Let $\xi_F$ be the restriction of $\xi_E$ {to} $F$ and $\xi_G$ be the quotient of $\xi_E$ on $G$. Then one has
\begin{equation}\label{Equ:degree+Delta}
\widehat{\deg}_+(F,\xi_F)\leqslant\widehat{\deg}_+(E,\xi_E)\leqslant\widehat{\deg}_+(F,\xi_F)+\widehat{\deg}_+(G,\xi_G)+\Delta(E,\xi_E).
\end{equation}
If in addition $(E,\xi_E)$ is ultrametric on $\Omega\setminus\Omega_\infty$, then 
\begin{equation}\label{Equ:degree+}
\widehat{\deg}_+(F,\xi_F)\leqslant\widehat{\deg}_+(E,\xi_E)\leqslant\widehat{\deg}_+(F,\xi_F)+\widehat{\deg}_+(G,\xi_G)+\delta(E,\xi_E).
\end{equation}
\end{prop}
\begin{proof}
The first inequality of \eqref{Equ:degree+Delta} follows directly from the definition of positive degree. In the following, we prove the second inequality of \eqref{Equ:degree+Delta}.

Let $E_1$ be a vector subspace of $E$, $F_1=F\cap E_1$ and $G_1$ be the canonical image of $E_1$ in $G$. Then we obtain the following short exact sequence of adelic vector bundles:
\[\xymatrix{0\ar[r]&(F_1,\xi_{F_1})\ar[r]&(E_1,\xi_{E_1})\ar[r]&(G_1,\xi_{G_1})\ar[r]&0},\]
where $\xi_{E_1}$ is the restriction of the norm family $\xi_E$, $\xi_{F_1}$ is the restriction of $\xi_{E_1}$ {to} $F_1$ and $\xi_{G_1}$ is the quotient norm family of $\xi_{E_1}$ on $G_1$. Note that the norm family $\xi_{F_1}$ coincides with the restricted norm family of $\xi_F$ induced by the inclusion map $F_1\rightarrow F$. Moreover, if we denote by $\xi_{G_1}'$ the restricted norm family induced by the inclusion map $G_1\rightarrow G$, then the identity map $(G_1,\xi_{G_1})\rightarrow (G_1,\xi_{G_1}')$ has norm $\leqslant 1$ (see Proposition \ref{prop:quotient:norm:linear:map} \ref{SubItem: diagram of quotient seminorm2}) on any $\omega\in\Omega$. In particular, by Proposition \ref{Pro:slopeinequality1} one has
\[\widehat{\deg}(G_1,\xi_{G_1})\leqslant\widehat{\deg}(G_1,\xi_{G_1}')\leqslant\widehat{\deg}_+(G,\xi_G).\]
 Therefore, by Proposition \ref{Pro:suiteexactedeg}, one has
\[\begin{split}\widehat{\deg}(E_1,\xi_{E_1})&\leqslant\widehat{\deg}(F_1,\xi_{F_1})+\widehat{\deg}( G_1,\xi_{G_1})+\Delta(E_1,\xi_{E_1})\\
&\leqslant\widehat{\deg}_+(F,\xi_F)+\widehat{\deg}_+(G,\xi_G)+\Delta(E_1,\xi_{E_1})\\
&\leqslant\widehat{\deg}_+(F,\xi_F)+\widehat{\deg}_+(G,\xi_G)+\Delta(E,\xi_{E}),\end{split}\]
where the last inequality comes from Corollary \ref{Cor: comparaison de Delta}. Similarly, by Propositions \ref{Pro:suiteexactedeg}  and  \ref{Pro:exactesequenceanddelta}, in the case where $\xi_E$ is ultrametric on $\Omega\setminus\Omega_\infty$ we have
\[\widehat{\deg}(E_1,\xi_{E_1})\leqslant\widehat{\deg}_+(F,\xi_F)+\widehat{\deg}_+(G,\xi_G)+\delta(E,\xi_E).\]
Since $E_1\subseteq E$ is arbitrary, we obtain \eqref{Equ:degree+}.
\end{proof}

\begin{prop}\label{prop:deg:plus:comp}
Let $(E, \xi = \{ \|\ndot\|_{\omega} \}_{\omega \in \Omega})$ be an adelic vector bundle on $S$.  
Then we have the following:
\begin{enumerate}[label=\rm(\arabic*)]
\item\label{Item: comparaison degre plus sous fibre}
Let $(F, \eta)$ be an adelic vector bundle on $S$ such that $F$ is a vector subspace of $E$ and $\eta \geqslant \xi_F$ on $F$.
Then $\widehat{\deg}_+(F, \eta) \leqslant \widehat{\deg}_+(F, \xi_F) \leqslant \widehat{\deg}_+(E, \xi)$.

\item\label{Item: degre plus tordu}
Let $\varphi$ be an integrable function on $\Omega$.
Then 
\[
\begin{cases}
{\displaystyle \widehat{\deg}_+(E, \exp(-\varphi) \xi) \leqslant \widehat{\deg}_+(E, \xi) + \dim_K( E) \left| \int_{S} \varphi(\omega) \,\nu(\mathrm{d}\omega) \right|},\\[1.5ex]
{\displaystyle \widehat{\deg}_+(E, \exp(\varphi) \xi) \geqslant \widehat{\deg}_+(E, \xi) - \dim_K (E) \left| \int_{S} \varphi(\omega) \,\nu(\mathrm{d}\omega) \right|.}\\
\end{cases}
\]
Moreover, if 
${\displaystyle \int_S \varphi\, \nu(\mathrm{d}\omega) \geqslant 0}$,
then
\[
\widehat{\deg}_+(E, \xi) \leqslant \widehat{\deg}_+(E, \exp(-\varphi) \xi)\quad\text{and}\quad
\widehat{\deg}_+(E, \exp(\varphi) \xi) \leqslant \widehat{\deg}_+(E, \xi).
\]
\end{enumerate}
\end{prop}

\begin{proof}
(1) The inequality $\widehat{\deg}_+(F, \xi_F) \leqslant \widehat{\deg}_+(E, \xi)$ has been proved in Proposition~\ref{prop:positive:degree:basic:properties}.
For $\epsilon > 0$, there is a vector subspace $W$ of $F$ such that
\[\widehat{\deg}_+(F, \eta) - \epsilon \leqslant \widehat{\deg}(W, \eta_W),\] so that
\[
\widehat{\deg}_+(F, \eta) - \varepsilon \leqslant \widehat{\deg}(W, \eta_W) \leqslant \widehat{\deg} (W, \xi_W) \leqslant  \widehat{\deg}_+(F, \xi_F),
\]
as required.

(2) 
Let $F$ be a vector subspace of $E$ over $K$. Then
\begin{equation}\label{eqn:prop:deg:plus:comp:01}
\widehat{\deg}(F, \exp(-\varphi)\xi_F) = \widehat{\deg}(F, \xi_{F}) + \dim_K (F) \int_{\Omega} \varphi(\omega)\,\nu(\mathrm{d}\omega),
\end{equation}
so that if
${\displaystyle \int_S \varphi \nu(d\omega) \geqslant 0}$,
then $\widehat{\deg}(F, \exp(-\varphi)\xi_F) \geqslant \widehat{\deg}(F, \xi_{F})$,
which {leads} to the third inequality. Moreover, by \eqref{eqn:prop:deg:plus:comp:01},
\[
\widehat{\deg}(F, \exp(-\varphi)\xi_F) \leqslant \widehat{\deg}_+(E, \exp(-\varphi) \xi) + \dim_K (E) \left| \int_{\Omega} \varphi(\omega)\, \nu(\mathrm{d}\omega)\right|,
\]
and hence the first inequality follows.

If we set $\xi' = \exp(\varphi)\xi$, then the first and third inequalities imply
the second and fourth inequalities, respectively.
\end{proof}

\subsection{Riemann-Roch formula}
Here we consider a Riemann-Roch formula of an adelic vector bundle on $S$.

\begin{prop}
\label{prop:Riemann:Roch:adelic:curve}
Let $\overline{V} = (V, \xi)$ be an adelic vector bundle on $S$. 
Then one has
\begin{equation}\label{Equ: Riemann-Roch}
0 \leqslant \widehat{\deg}(\overline{V}) - \left(
\widehat{\deg}_+(\overline{V}{}^{\vee\vee}) - \widehat{\deg}_+(\overline{V}{}^{\vee}) \right)  \leqslant \delta(\overline V).
\end{equation}
If $\xi$ is ultrametric on $\Omega\setminus\Omega_\infty$, then one has
\begin{equation}\label{Equ: Riemann-Roch bis}
0 \leqslant \widehat{\deg}(\overline{V}) - \left(
\widehat{\deg}_+(\overline{V}) - \widehat{\deg}_+(\overline{V}{}^{\vee}) \right)  \leqslant \delta(\overline V).
\end{equation}
\end{prop}

\begin{proof}
Let $\epsilon > 0$. We choose a vector subspace $W$ of $V$ such that
$\widehat{\deg}(W,\xi_W) \geqslant \widehat{\deg}_+(\overline{V}) - \epsilon$, where $\xi_W$ is the restriction of $\xi$ {to} $W$.
Let $\xi^{\vee}$ be the dual norm {family} of $\xi$ on $V^{\vee}$, $\xi_{V/W}$ be the quotient norm {family} of $\xi$ on $V/W$,
and $\xi_{V/W}^{\vee}$ be the dual norm family of $\xi_{V/W}$ on $(V/W)^{\vee}$. If we consider $(V/W)^{\vee}$ as a vector subspace of $V^{\vee}$, then
$\xi_{V/W}^{\vee}$ identifies with the restriction of $\xi^{\vee}$ {to} $(V/W)^{\vee}$ by Proposition~\ref{Pro:dualquotient}, so that
\[
\widehat{\deg}((V/W)^{\vee},\xi_{V/W}^{\vee}) \leqslant \widehat{\deg}_+(V^{\vee},\xi^{\vee}).
\]
On the other hand, one has
\[
\widehat{\deg}((V/W)^{\vee},\xi_{V/W}^{\vee}) + \widehat{\deg}(V/W, \xi_{V/W}) \geqslant 0
\]
by Proposition~\ref{Pro:degdual} and
\[
\widehat{\deg}(W,\xi_W) + \widehat{\deg}(V/W,\xi_{V/W}) \leqslant \widehat{\deg}(V,\xi)
\]
by Proposition~\ref{Pro:suiteexactedeg}. Therefore,
\begin{align*}
\widehat{\deg}_+(V^{\vee},\xi^{\vee}) & \geqslant \widehat{\deg}((V/W)^{\vee}, \xi_{V/W}^{\vee}) \geqslant -\operatorname{\widehat{\deg}}(V/W,\xi_{V/W}) \\
& \geqslant \widehat{\deg}(W,\xi_W) - \widehat{\deg}(V,\xi)
\geqslant \widehat{\deg}_+(V,\xi) - \widehat{\deg}(V,\xi) - \epsilon 
\end{align*}
and hence
\begin{equation}
\label{eqn:prop:Riemann:Roch:adelic:curve:01}
\widehat{\deg}_+(\overline{V}) - \widehat{\deg}_+(\overline{V}{}^{\vee}) \leqslant \widehat{\deg}(\overline{V}).
\end{equation}
Replacing $\overline{V}$ by $\overline{V}{}^{\vee}$ in \eqref{eqn:prop:Riemann:Roch:adelic:curve:01}, we obtain
\begin{equation}\label{eqn:prop:Riemann:Roch:adelic:curve:02}
\widehat{\deg}_+(\overline{V}{}^{\vee}) - \widehat{\deg}_+(\overline{V}{}^{\vee\vee}) \leqslant \widehat{\deg}(\overline{V}{}^{\vee}),
\end{equation}
which, by Proposition \ref{Pro:degdual}, implies the second inequality of \eqref{Equ: Riemann-Roch}.
Replacing $\overline{V}$ by $\overline{V}{}^{\vee\vee}$ in \eqref{eqn:prop:Riemann:Roch:adelic:curve:01}, by the fact that $\norm{\ndot}_{\omega,**,*}=\norm{\ndot}_{\omega,*}$ for any $\omega\in\Omega$ (see Proposition \ref{Pro:doubedualandquotient} \ref{Item: seminorm and double dual induce the same dual norm}) and the equality \eqref{Equ:doublddual} we obtain
\begin{equation}\label{eqn:prop:Riemann:Roch:adelic:curve:03}
\widehat{\deg}_+(\overline{V}{}^{\vee\vee}) - \widehat{\deg}_+(\overline{V}{}^{\vee}) \leqslant \widehat{\deg}(\overline{V}),
\end{equation}
which leads to the first inequality of \eqref{Equ: Riemann-Roch}.

In the case where $\xi$ is ultrametric on $\Omega\setminus\Omega_\infty$, one has $\|\ndot\|_{\omega} = \|\ndot\|_{\omega, **}$ for all $\omega \in \Omega$
by Proposition~\ref{Pro:doubledualarch} (Archimedean case) and Corollary~\ref{Cor:doubledual}. Thus \eqref{Equ: Riemann-Roch bis} follows from \eqref{Equ: Riemann-Roch}.
\end{proof}

\subsection{Comparison of $\widehat{\deg}_+$ and $\widehat{h}^0$ in the classic setting}
In this subsection, we compare $\widehat{\deg}_+$ with the invariant $\widehat{h}^0$ in the classic settings of vector bundles on a regular projective curve and  Hermitian vector bundles on an arithmetic curve. 

\subsubsection{Function field case}
Let $k$ be a field, $C$ be a regular projective curve over $\Spec k$ and $K=k(C)$ be the field of rational functions on $C$. Let $\Omega$ be the set of all closed points of the curve $C$, equipped with the discrete $\sigma$-algebra $\mathcal A$ and the measure $\nu$ such that $\nu(\{x\})=[k(x):k]$ for any $x\in\Omega$. Let $\phi:\Omega\rightarrow M_K$ be the map sending $x$ to $|\ndot|_x=\mathrm{e}^{-\mathrm{ord}_x(\ndot)}$. We have seen in \S\ref{Subsec: function field} that $S=(K,(\Omega,\mathcal A,\nu),\phi)$ is an adelic curve.

Let $\mathcal E$ be a locally free $\mathcal O_C$-module of finite type and $E=\mathcal E_K:=\mathcal E\otimes_{\mathcal O_C}K$ be its generic fibre. For any $x\in\Omega$, let $\|\ndot\|_x$ be the norm on $E\otimes_KK_\omega$ defined as
\[\forall\,s\in E\otimes_KK_\omega,\quad\|s\|_x=\inf\{|a|_x\,:\,a\in K_\omega^{\times},\;a^{-1}s\in\widehat{\mathcal O}_{C,x},\}\]
where $\widehat{\mathcal O}_{C,x}$ is the completion of $\mathcal O_{C,\omega}$, which identifies with the valuation ring of $K_\omega$. Then $\xi_{\mathcal E}=\{\|\ndot\|_x\}_{x\in \Omega}$ forms a Hermitian norm family on $E$ and $(E,\xi_{\mathcal E})$ is an adelic vector bundle on $\Omega$. Note that the Arakelov degree of $(E,\xi_{\mathcal E})$ identifies with the degree of the locally free $\mathcal O_C$-module $\mathcal E$, namely
\[\widehat{\deg}(E,\xi_{\mathcal E})=\deg(c_1(\mathcal E)\cap[C]).\]
The Harder-Narasimhan flag of $(E,\xi_{\mathcal E})$ is also related to the classic construction of Harder-Narasimhan filtration of $\mathcal E$. In fact, there exists a unique flag of locally free $\mathcal O_C$-modules
\[0=\mathcal E_0\subsetneq\mathcal E_1\subsetneq \ldots\subsetneq\mathcal E_n=\mathcal E\]
such that each sub-quotient $\mathcal E_i/\mathcal E_{i-1}$ is a locally free $\mathcal O_C$-module which is semistable and that
\[\mu(\mathcal E_1/\mathcal E_0)>\mu(\mathcal E_2/\mathcal E_1)>\ldots>\mu(\mathcal E_n/\mathcal E_{n-1}).\] 
Then the Harder-Narasimhan flag of the Hermitian adelic vector bundle $(E,\xi_{\mathcal E})$ is given by
\[0=\mathcal E_{0,K}\subsetneq\mathcal E_{1,K}\subsetneq\ldots\subsetneq\mathcal E_{n,K}=E.\]

The notion of positive degree for locally free $\mathcal O_C$-modules of finite rank has been proposed in \cite{Chen15} and compared with the rank (over $k$) of the vector space of global sections, by using the Riemann-Roch formula on curves. For any locally free $\mathcal O_C$-module of finite rank $\mathcal E$, we denote by $h^0(\mathcal E)$ the rank of $H^0(C,\mathcal E)$ over $k$.

\begin{theo}
Let $g(C)$ be the genus of $C$ relatively to the field $k$ (namely $g(C)=h^0(\omega_{C/k})$, $\omega_{C/k}$ being the relative dualising sheaf of $C$ over $\Spec k$). For any locally free $\mathcal O_C$-module of finite rank $\mathcal E$, one has
\begin{equation}
|h^0(\mathcal E)-\widehat{\deg}_+(E,\xi_{\mathcal E})|\leqslant \rang_K(E)\max(g(C)-1,1).
\end{equation}
\end{theo}
We refer the readers to \cite[\S2]{Chen15} for a proof.

\subsubsection{Number field case}

Let $K$ be a number field and $\Omega$ be the set of all places of $K$, equipped with the discrete $\sigma$-algebra $\mathcal A$. For each $\omega\in\Omega$, we denote by $|\ndot|_\omega$ the absolute value on $K$ extending either the usual absolute value on $\mathbb Q$ or one of the $p$-adic absolute values (with $|p|_\omega=p^{-1}$ in the latter case). We let $K_\omega$ (resp. $\mathbb Q_\omega$) be the completion of $K$ (resp. $\mathbb Q$) with respect to the absolute value $|\ndot|_\omega$. Let $\nu$ be the measure on the measurable space $(\Omega,\mathcal A)$ such that $\nu(\{\omega\})=[K_\omega:\mathbb Q_\omega]$. Then $S=(K,(\Omega,\mathcal A,\nu),\phi:\omega\mapsto |\ndot|_{\omega})$ forms an adelic curve.

Let $\mathcal O_K$ be the ring of algebraic integers in $K$. Recall that a Hermitian vector bundle on $\Spec\mathcal O_K$ is by definition the data $\overline{\mathcal E}=(\mathcal E,\{\|\ndot\|_\sigma\}_{\sigma\in\Omega_\infty})$ of a projective $\mathcal O_K$-module of finite type $\mathcal E$ together with a family of norms, where $\|\ndot\|_\sigma$ is a norm on the vector space $\mathcal E\otimes_{\mathcal O_K}K_\omega$ which is induced by an inner product. Similarly to the function field case, the $\mathcal O_K$-module structure of $\mathcal E$ induces, for each non-Archimedean place $\mathfrak p\in\Omega\setminus\Omega_\infty$, a ultrametric norm $\|\ndot\|_{\mathfrak p}$ on $\mathcal E\otimes_{\mathcal O_K}K_{\mathfrak p}$ as follows~:
\[\forall\,s\in \mathcal E\otimes_{\mathcal O_K}K_{\mathfrak p},\quad \|s\|_{\mathfrak p}=\inf\{|a|_{\mathfrak p}\,:\,a\in K_{\mathfrak p}^{\times},\;a^{-1}s\in\mathcal O_{\mathfrak p}\},\]
where $\mathcal O_{\mathfrak p}$ is the valuation ring of $K_{\mathfrak p}$. Let $E$ be $\mathcal E\otimes_{\mathcal O_K}K$ and let $\xi_{\overline{\mathcal E}}$ be the norm family $\{\|\ndot\|_\omega\}_{\omega\in\Omega}$. Then the couple $(E,\xi_{\overline{\mathcal E}})$ form an adelic vector bundle on $S$, which is said to be \emph{induced by $\overline{\mathcal E}$}. 

Recall that the space  $\widehat{H}^0(\overline{\mathcal E})$ of ``global sections'' of $\overline{\mathcal E}$ is defined as
\[\widehat{H}^0(\overline{\mathcal E}):=\{s\in\mathcal E\,:\,\sup\nolimits_{\sigma\in\Omega_\infty}\|s\|_\sigma\leqslant 1\}=\{s\in E\,:\,\sup\nolimits_{\omega\in\Omega}\|s\|_\omega\leqslant 1\}.\]
This is a finite set. However it does not possess a natural vector space structure over a base field. We define (compare to the case of function field of a regular projective curve over a finite field) $\widehat{h}^0(\overline{\mathcal E})$ to be $\ln\mathrm{card}(\widehat{H}^0(\overline{\mathcal E}))$. The invariants $\widehat{h}^0(\overline{\mathcal E})$ and $\widehat{\deg}_+(E,\xi_{\overline{\mathcal E}})$ have been compared in \cite[\S6]{Chen15} (see also \cite{Chen08}). We denote by 
\begin{enumerate}[label=--]
\item $B_n$ the unit ball in $\mathbb R^n$, where $n\in\mathbb N$,
\item $\mathrm{vol}(B_n)$ the Lebesgue measure of $B_n$, which is equal to $\pi^{n/2}/\Gamma(n/2+1)$,
\item $r_1(K)$ the number of real places of $K$,
\item $r_2(K)$ the number of complex places of $K$.
\end{enumerate}

\begin{theo}
Let $\overline{\mathcal E}$ be a Hermitian vector bundle on $\Spec\mathcal O_K$, and $(E,\xi_{\overline{\mathcal C}})$ be the adelic vector bundle on $S$ induced by $\overline{\mathcal E}$. then
\[|\widehat{h}^0(\overline{\mathcal E})-\widehat{\deg}_+(E,\xi_{\overline{\mathcal E}})|\leqslant \rang_K(E)\ln|\mathfrak D_K|+C(K,\rang_K(E)),\]
where $\mathfrak D_K$ is the discriminant of $K$ over $\mathbb Q$, and for any  
$n\in\mathbb N$,
\[\begin{split}C(K,n):=n[K:\mathbb Q]&\ln(3)+n(r_1(K)+r_2(K))\ln(2)+\frac n2\ln|\mathfrak D_K|-r_1\ln(\mathrm{vol}(B_n)n!)\\
&-r_2\ln(V(B_{2n})(2n)!)+\ln(([K:\mathbb Q]n)!)\end{split}\]
\end{theo}

\subsection{Slopes and slope inequalities}
Let $(E,\xi)$ be an adelic vector bundle on $S$ such that $E\neq\{0\}$. We define the \emph{slope}\index{slope} {of} $(E,\xi)$ as
\[\widehat{\mu}(E,\xi):=\frac{\widehat{\deg}(E,\xi)}{\rang(E)}.\]
We define the \emph{maximal slope} of $(E,\xi)$ as 
\[\widehat{\mu}_{\max}(E,\xi):=\sup_{0\neq F\subseteq E}\widehat{\mu}(F,\xi_F),\]
where $F$ runs over the set of all non-zero vector subspaces of $E$ and $\xi_F$ denotes the restriction of $\xi$ to $F$.   
Similarly, we define the \emph{minimal slope}\index{minimal slope}\index{slope!minimal ---}
of $(E,\xi)$ as 
\begin{equation}\label{Equ: definition of mu min}\widehat{\mu}_{\min}(E,\xi)=\inf_{ E \twoheadrightarrow G \not=\{0\}}\widehat{\mu}_{\max}(G,\xi_G),\end{equation}
where $G$ runs over the set of all non-zero quotient spaces of $E$, and $\xi_G$ denotes the quotient norm family of $\xi$. By definition one has
$\widehat{\mu}_{\min}(\overline E)\leqslant\widehat{\mu}_{\max}(\overline E)$ and $\widehat{\mu}(\overline E)\leqslant\widehat{\mu}_{\max}(\overline E)$ (note that here the vector space $E$ has been assumed to be non-zero).
If $\overline E=\overline{\boldsymbol{0}}$ is the zero adelic vector bundle, we define by convention
\[\widehat{\mu}_{\max}(\overline {\boldsymbol{0}}):=-\infty,\quad\widehat{\mu}(\overline{\boldsymbol{0}}):=0,\quad\widehat{\mu}_{\min}(\overline {\boldsymbol{0}}):=+\infty.\]

\begin{prop}\label{Pro: mu min plus mu max dual}
Let $\overline E=(E,\xi_E)$ be a non-zero  
adelic vector bundle on $S$. One has
\begin{equation}\label{Equ: sum of mu min and mu max dual}\widehat{\mu}_{\min}(\overline E)+\widehat{\mu}_{\max}(\overline E{}^\vee)\geqslant 0,\end{equation}
provided that $\widehat{\mu}_{\max}(\overline E{}^\vee)<+\infty$ (we will show in Proposition \ref{Pro: finiteness of slopes} that this condition is always satisfied, {and, as a consequence of the current proposition, that one has $\widehat{\mu}_{\min}(\overline E)>-\infty$}).
\end{prop}
\begin{proof}
Let $G$ be a non-zero quotient vector space of $E$ and $\xi_G$ be the quotient norm family of $\xi_E$. Note that $G^\vee$ identifies with a vector subspace of $E^\vee$ and by Proposition \ref{Pro:dualquotient}, the dual norm family $\xi_G^\vee$ identifies with the restriction of $\xi_E^\vee$ {to} $G^\vee$. By Proposition \ref{Pro:degdual}, one has
\[0\leqslant\widehat{\mu}(G,\xi_G)+\widehat{\mu}(G^\vee,\xi_G^\vee)\leqslant\widehat{\mu}_{\max}(G,\xi_G)+\widehat{\mu}_{\max}(E^\vee,\xi_E^\vee).\]
Since $G$ is arbitrary, we obtain the inequality \eqref{Equ: sum of mu min and mu max dual}.  
\end{proof}

Classically in the setting of vector bundles over a regular projective curve or in that of Hermitian vector bundle over an arithmetic curve, the minimal slope is rather defined as the minimal value of slopes of quotient bundles. A direct analogue would replace $\widehat{\mu}_{\max}$ by $\widehat{\mu}$ in \eqref{Equ: definition of mu min} for the definition of the minimal slope. However, it can be shown that the two definitions are actually equivalent.

\begin{prop}\label{Pro: mu min equivalent to classiical one}
Let $\overline E$ be a non-zero adelic vector bundle on $S$. One has
\[\widehat{\mu}_{\min}(\overline E)=\inf_{ E \twoheadrightarrow G \not=\{0\}}\widehat{\mu}(\overline G),\]
where $G$ runs over the set of non-zero quotient vector space of $E$, and in $\overline G$ we consider the quotient norm family of that in $\overline E$.
\end{prop}
\begin{proof}
For any non-zero quotient vector space $G$ of $E$, one has $\widehat{\mu}(\overline G)\leqslant\widehat{\mu}_{\max}(\overline G)$. Therefore the case where $\widehat{\mu}_{\min}(\overline E)=-\infty$ is trivial. In the following, we assume that $\widehat{\mu}_{\min}(\overline E)>-\infty$.

Let $\epsilon>0$. Among the non-zero quotient vector spaces of $E$ of maximal slope bounded from above by $\widehat{\mu}_{\min}(\overline E)+\epsilon$, we choose a $G$ having the least rank. We claim that the maximal slope of $\overline G$ does not exceed $\widehat{\mu}(\overline G)+\epsilon$ and we will prove this assertion by contradiction. Assume the contrary, namely $\widehat{\mu}_{\max}(\overline G)>\widehat{\mu}(\overline G)+\epsilon$. Let $G'$ be a non-zero vector subspace of $G$ such that 
$\widehat{\mu}(\overline{G'})> \widehat{\mu}(\overline G)$.
We suppose in addition that $G'$ is of maximal rank among the non-zero vector subspaces verifying this condition.
Note that by definition $G'\subsetneq G$.
 Let $H=G/G'$. Since $G'$ is not equal to $G$, $H$ is a non-zero quotient vector space of $E$. Let $H''$ be a non-zero vector subspace of $H$ and $G''$ be the inverse image of $H''$ by the quotient map $G\rightarrow H=G/G'$. We have a short exact sequence
\[\xymatrix{0\ar[r]& G'\ar[r]&G''\ar[r]&H''\ar[r]&0},\]
which leads to (by Proposition \ref{Pro:suiteexactedeg})
\begin{equation}\label{Equ: deg G' plus deg}\widehat{\deg}(\overline{G'})+\widehat{\deg}(\overline {H''})\leqslant\widehat{\deg}(\overline{G''}).\end{equation}
Since $H''$ is non-zero, the rank of $G''$ is greater than that of $G'$, and hence \[\widehat{\mu}(\overline{G''})\leqslant\widehat{\mu}(\overline G)<\widehat{\mu}(\overline{G'})\]
by the maximality assumption of $\rang_K(G')$. Therefore \eqref{Equ: deg G' plus deg} leads to
$\widehat{\mu}(\overline{H''})\leqslant\widehat{\mu}(\overline{G'})$. Since $H''$ is arbitrary, we obtain
\[\widehat{\mu}_{\max}(\overline H)\leqslant\widehat{\mu}(\overline{G'})\leqslant\widehat{\mu}_{\max}(\overline G)\leqslant\widehat{\mu}_{\min}(\overline E)+\epsilon.\]
This contradicts the minimality assumption of $\rang_K(G)$, which proves the claim. Since $\epsilon$ is arbitrary, the proposition is proved.
\end{proof}

\begin{coro}\label{Cor: majoration of mu min plus mu max}
Let $\overline E=(E,\xi_E)$ be a non-zero adelic vector bundle on $S$. One has
\begin{equation}\label{Equ: majoration of mu min plus mu max}
\widehat{\mu}_{\min}(\overline E{}^\vee)+\widehat{\mu}_{\max}(\overline E)\leqslant\frac 12\ln(\rang_K(E))\,\nu(\Omega_\infty).
\end{equation}
Moreover, one has $\widehat{\mu}_{\min}(\overline E{}^\vee)+\widehat{\mu}_{\max}(\overline E)=0$ if $\overline E$ is Hermitian.
\end{coro}
\begin{proof}
Let $F$ be a non-zero vector subspace of $E$ and $\xi_F=\{\norm{\ndot}_{F,\omega}\}_{\omega\in\Omega}$ be the restriction of $\xi_E=\{\norm{\ndot}_{E,\omega}\}_{\omega\in\Omega}$ {to} $F$. For any $\omega\in\Omega$, $\norm{\ndot}_{F,\omega,*}$ is bounded from above by the quotient norm of $\norm{\ndot}_{E,\omega,*}$ by the canonical surjective map $E^\vee_{K_\omega}\rightarrow F^\vee_{K_\omega}$. 
 Hence by Proposition \ref{Pro: mu min equivalent to classiical one}, one has $\widehat{\mu}(\overline{F}{}^\vee)\geqslant\widehat{\mu}_{\min}(\overline E{}^\vee)$. Therefore,
\[\begin{split}\widehat{\mu}(\overline F)+\widehat{\mu}_{\min}(\overline E{}^\vee)&\leqslant\widehat{\mu}(\overline F)+\widehat{\mu}(\overline{F}{}^\vee)=\frac{1}{\rang_K(F)}\int_{\Omega}\ln(\delta_\omega(\xi_F))\,\nu(d\omega)\\
&\leqslant\frac 12\ln(\rang_K(F))\nu(\Omega_\infty),\end{split}\]
where the equality follows from Proposition \ref{Pro:degdual} and the last inequality comes from Remark \ref{rema:Thm:John}. Since $F$ is arbitrary, we obtain \eqref{Equ: majoration of mu min plus mu max}. 

If $\overline E$ is Hermitian, then for any non-zero vector subspace $F$ of $E$ one has 
\[\widehat{\mu}(\overline F)+\widehat{\mu}_{\min}(\overline E{}^\vee)\leqslant\widehat{\mu}(\overline F)+\widehat{\mu}(\overline{F}{}^\vee)=0,\] 
which leads to $\widehat{\mu}_{\min}(\overline E{}^\vee)+\widehat{\mu}_{\max}(\overline E)\leqslant 0$. As we have seen that $\widehat{\mu}_{\min}(\overline E{}^\vee)+\widehat{\mu}_{\max}(\overline E)\geqslant 0$ in Proposition \ref{Pro: mu min plus mu max dual} {(note that $\overline{E}{}^{\vee\vee}=\overline E$ when $\overline E$ is Hermitian)}, the equality $\widehat{\mu}_{\min}(\overline E{}^\vee)+\widehat{\mu}_{\max}(\overline E)=0$ holds.
\end{proof}

\subsection{Finiteness of slopes}

Let $\overline E=(E,\xi)$ be an adelic vector bundle on $S$ such that $\xi$ is Hermitian. We assume that the vector space $E$ is non-zero and we denote by $\Theta(E)$ the set of all $K$-vector subspaces of $E$. For any $F\in\Theta(E)$, the vector subspace $F$ equipped with the restricted norm family forms a Hermitian adelic vector bundle on $S$ (see Proposition \ref{Pro:critereadelicvb}). We denote by $\overline F$ this Hermitian adelic vector bundle. Note that the rank and the Arakelov degree defines two functions on $\Theta(E)$, which satisfy the following relations: for any pair $(E_1,E_2)$ of elements in $\Theta(E)$, one has

\begin{gather}
\rang_K(E_1\cap E_2)+\rang_K(E_1+E_2)=\rang_K(E_1)+\rang_K(E_2),\\\label{Equ:convexitededeg}
\widehat{\deg}(\overline{E_1\cap E_2})+\widehat{\deg}(\overline{E_1+E_2})\geqslant \widehat{\deg}(\overline E_1)+\widehat{\deg}(\overline E_2),
\end{gather}
where the inequality \eqref{Equ:convexitededeg} comes from Corollary \ref{Cor:convexitenorme}.

\begin{prop}\label{Pro: existence of maximal slope}
Let $E$ be a non-zero vector space of finite rank over $K$ and $\Theta(E)$ be the set of all vector subspaces of $E$. Assume given two functions $r:\Theta(E)\rightarrow\mathbb R_+$ and $d:\Theta(E)\rightarrow\mathbb R$ which verify the following conditions:
\begin{enumerate}[label=\rm(\arabic*)]
\item the function $r(\ndot)$ takes value $0$ on the zero vector subspace of $E$ and takes positive values on non-zero vector subspaces of $E$;
\item for any couple $(E_1,E_2)$ of elements in $\Theta(E)$ one has \[r(E_1\cap E_2)+r(E_1+E_2)=r(E_1)+r(E_2)\]
and 
\[d(E_1\cap E_2)+d(E_1+E_2)\geqslant d(E_1)+d(E_2);\]
\item {$d(\{0\})\leqslant 0$}.
\end{enumerate}
Then the function $\mu=d/r$ attains its maximal value $\mu_{\max}$ on the set $\Theta^*(E)$ of all non-zero vector subspaces of $E$. Moreover, there exists a non-zero vector subspace $E_{\mathrm{des}}$ of $E$ such that $\mu(E_{\mathrm{des}})=\mu_{\max}$ and which contains all non-zero vector subspaces $F$ of $E$ such that $\mu(F)=\mu_{\max}$.
\end{prop}
\begin{proof}
The first relation in the condition (2) implies that, if $L_1,\ldots,L_n$ are vector subspaces of rank $1$ of $E$, which are linearly independent, then one has
\[r(L_1+\cdots+L_n)=r(L_1)+\cdots+r(L_n).\]
In particular, if $L$ and $L'$ are different vector subspaces of rank $1$ in $E$ then one has $r(L)=r(L')$. In fact, let $s$ and $s'$ be non-zero vectors of $L$ and $L'$ respectively, and let $L''$ be the vector subspace of $E$ generated by $s+s'$ (which is of rank $1$). Then one has \[r(L)+r(L'')=r(L+L')=r(L')+r(L'').\]  Therefore the function $r(\ndot)$ is proportional to the rank function. Without loss of generality, we may assume that the function $r(\ndot)$ identifies with the rank function.

We prove the proposition by induction on the rank of the vector space $E$. The case where the $r(E)=1$ is trivial. In the following, we assume that $r(E)\geqslant 2$ and that the proposition has been proved for vector spaces of rank $<r(E)$. If for any non-zero vector subspace $F$ of $E$ one has $\mu(F)\leqslant\mu(E)$, then there is nothing to prove since $\mu(E)=\mu_{\max}$ and $E=E_{\mathrm{des}}$. Otherwise there exists a non-zero vector subspace $E'$ of $E$ such that $\mu(E')>\mu(E)$. Moreover, we can choose $E'$ such that $r(E')$ is maximal (among the non-zero vector subspaces of $E$ having this property). Clearly one has $r(E')<r(E)$. Hence by the induction hypothesis the {restriction of the} function $\mu(\ndot)$ {to $\Theta^*(E')$} attains its maximum, and among the non-zero vector subspaces of $E'$ on which the restriction of the function $\mu(\ndot)$ {to} $\Theta(E')$ attains the maximal value there exists a greatest one $E'_{\mathrm{des}}$ with respect to the relation of inclusion. Let $E_{\mathrm{des}}:=E'_{\mathrm{des}}$ be this vector space. We claim that $E_{\mathrm{des}}$ verifies the properties announced in the proposition.

Let $F$ be a non-zero vector subspace of $E$. If $F\subseteq E'$, then clearly one has $\mu(F)\leqslant\mu(E_{\mathrm{des}})$. Otherwise the rank of $F\cap E'$ is smaller than $r(F)$ and the rank of $F+E'$ is greater than $r(E')$. Moreover, since $F\cap E'\subseteq E'$, one has {(here we use the condition that $d(\{0\})\leqslant 0$ to treat the case where $F\cap E'=\{0\}$)} \[{d(F\cap E')\leqslant\mu(E_{\mathrm{des}})r(F\cap E');}\] since $F+E'\supsetneq E'$, one has $\mu(F+E')\leqslant\mu(E)<\mu(E')$. Therefore
\[\begin{split}d(F\cap E')+d(F+E')&=\mu(F\cap E')r(F\cap E')+\mu(F+E')r(F+E')\\
&<\mu(E_{\mathrm{des}})r(F\cap E')+\mu(E')r(F+E').
\end{split}\] 
Combining this relation with the inequality in the condition (2) of the proposition, we obtain 
\[\mu(E_{\mathrm{des}})r(F\cap E')+\mu(E')r(F+E')> \mu(E')r(E')+\mu(F)r(F).\]
By the equality in the condition (2), we deduce
\[\mu(F)r(F)<\mu(E_{\mathrm{des}})r(F\cap E')+\mu(E')(r(F)-r(F\cap E'))\leqslant\mu(E_{\mathrm{des}})r(F).\]
Therefore, the function $\mu(\ndot)$ attains its maximal value $\mu_{\max}$ at $E_{\mathrm{des}}$. Moreover, if $F$ is a non-zero vector subspace of $E$ such that $\mu(F)=\mu(E_{\mathrm{des}})$, then one should have $F\subseteq E'$, and hence $F\subseteq E_{\mathrm{des}}$ by the induction hypothesis. The proposition is thus proved.
\end{proof}

\begin{defi}
Let $\overline E$ be a non-zero Hermitian adelic vector bundle on $S$. We can apply the above proposition to the functions of rank and of Arakelov degree to obtain the existence of a (unique) non-zero vector subspace $E_{\mathrm{des}}$ of $E$ such that 
\[\widehat{\mu}(\overline E_{\mathrm{des}})=\widehat{\mu}_{\max}(\overline E)=\sup_{0\neq F\in\Theta(E)}\widehat{\mu}(\overline F)\]
and containing all non-zero vector subspaces of $E$ on which the function $\widehat{\mu}$ attains the maximal slope of $\overline E$. The vector subspace $E_{\mathrm{des}}$ is called the \emph{destabilising vector subspace}\index{destabilising vector subspace} of the Hermitian adelic vector bundle $\overline E$.
If $\overline E_{\mathrm{des}}=\overline E$, we say that the Hermitian adelic vector bundle $\overline E$ is \emph{semistable}\index{semistable}\index{adelic vector bundle!semistable ---}. In particular, for any non-zero Hermitian adelic vector bundle $\overline E$ on $S$, the Hermitian adelic vector bundle $\overline E_{\mathrm{des}}$ is always semistable.
\end{defi}

{
\begin{prop}\label{Pro: finiteness of slopes}
Let $(E,\xi)$ be a non-zero adelic vector bundle on $S$. Then one has $\widehat{\mu}_{\max}(E,\xi)<+\infty$ and $\widehat{\mu}_{\min}(E,\xi)>-\infty$.
\end{prop}
\begin{proof}
{Let $r$ be the rank of $E$ over $K$.} We first prove the inequality   $\widehat{\mu}_{\max}(E,\xi)<+\infty$ in the particular case where $\xi$ is ultrametric on $\Omega\setminus\Omega_\infty$. By Theorem \ref{Thm: Hermitian approximation via measurable selection}, there exists a measurable Hermitian norm family $\xi^H$ on $E$ such that 
\[\forall\,\omega\in\Omega,\quad d_\omega(\xi,\xi^H)\leqslant \frac 12\ln(r+1)\indic_{\Omega_\infty}(\omega).\]
Therefore, for any non-zero vector subspace $F$ of $E$ one has
\[\big|\widehat{\mu}(F,\xi_F)-\widehat{\mu}(F,\xi^H_F)\big|\leqslant\frac 12\ln(r+1)\nu(\Omega_\infty). \]
Moreover, by Proposition \ref{Pro: existence of maximal slope}, the maximal slope $\widehat{\mu}_{\max}(E,\xi^H)$ is finite. Therefore one has $\widehat{\mu}_{\max}(E,\xi)<+\infty$.

We now proceed with the proof of the relation $\widehat{\mu}_{\max}(E,\xi)<+\infty$ in the general case. We write $\xi$ in the form $\{\norm{\ndot}_\omega\}_{\omega\in\Omega}$. Note that for any $\omega\in\Omega$ one has $\norm{\ndot}_{\omega,**}\leqslant\norm{\ndot}_{\omega}$. Therefore $\widehat{\mu}_{\max}(E,\xi^{\vee\vee})\geqslant\widehat{\mu}_{\max}(E,\xi)$. Note that the norm family $\xi^{\vee\vee}$ is ultrametric on $\Omega\setminus\Omega_\infty$, and $(E,\xi^{\vee\vee})$ is an adelic vector bundle (see Proposition \ref{Pro:critereadelicvb} \ref{Item: dual adelic vector bundle 01}). 
By the particular case proved above, one has $\widehat{\mu}_{\max}(E,\xi^{\vee\vee})<+\infty$. Thus we obtain $\widehat{\mu}_{\max}(E,\xi)<+\infty$.

Applying the above proved result to $(E^\vee,\xi^\vee)$ we obtain $\widehat{\mu}_{\max}(E^\vee,\xi^\vee)<+\infty$. Therefore, by Proposition \ref{Pro: mu min plus mu max dual} we obtain $\widehat{\mu}_{\min}(E,\xi)>-\infty$.
 
\end{proof}}

\subsection{Some slope estimates}

The following proposition is a natural generalisation of the slope inequalities to the setting of adelic curves. We refer the readers to \cite[\S4.1]{Bost2001} for this theory in the classic setting of Hermitian vector bundles over an algebraic integer ring.

\begin{prop}\label{Pro:inegalitdepente}
Let $(E_1,\xi_1)$ and $(E_2,\xi_2)$ be 
adelic vector bundles on $S$, and $f:E_1\rightarrow E_2$ be a $K$-linear map.
\begin{enumerate}[label=\rm(\arabic*)]
\item\label{Item: slope inequality mu max} If $f$ is injective, then one has $\widehat{\mu}_{\max}( E_1,\xi_1)\leqslant\widehat{\mu}_{\max}(E_2,\xi_2)+h(f)$.
\item\label{Item: slope inequality mu min} If $f$ is surjective, then one has $\widehat{\mu}_{\min}( E_1,\xi_1)\leqslant\widehat{\mu}_{\min}( E_2,\xi_2)+h(f)$.
\item\label{Item: slope inequality mu mix} If $f$ is non-zero, then one has $\widehat{\mu}_{\min}( E_1,\xi_1)\leqslant\widehat{\mu}_{\max}(E_2,\xi_2)+h(f)$.
\end{enumerate}
\end{prop}
\begin{proof}
\ref{Item: slope inequality mu max} The assertion is trivial if $f$ is the zero map since in this case $E_1=\{0\}$ and $\widehat{\mu}_{\max}(E_1,\xi_1)=-\infty$ by convention. In the following, we assume that the linear map $f$ is non-zero. Let $F_1$ be a non-zero vector subspace of $E_1$ and $F_2$ be its image in $E_2$. Let $g:F_1\rightarrow F_2$ be the restriction of $f$ {to} $F_1$. It is an isomorphism of vector spaces. Moreover, if we equip $F_1$ and $F_2$ with induced norm families, by Proposition \ref{Pro:slopeinequality1} one has
\[\widehat{\mu}( F_1,\xi_{F_1})\leqslant\widehat{\mu}(F_2,\xi_{F_2})+h(g)\leqslant\widehat{\mu}_{\max}(E_2,\xi_2)+h(f),\]
where $\xi_{F_1}$ and $\xi_{F_2}$ are restrictions of $\xi_1$ and $\xi_2$ {to} $F_1$ and $F_2$, respectively.
Since $F_1$ is arbitrary, we obtain the inequality $\widehat{\mu}_{\max}( E_1,\xi_1)\leqslant\widehat{\mu}_{\max}(E_2,\xi_2)+h(f)$.

\ref{Item: slope inequality mu min} The assertion is trivial if $f$ is the zero map since in this case $E_2=\{0\}$ and $\widehat{\mu}_{\min}(E_2,\xi_2)=+\infty$ by convention. In the following, we assume that the linear map $f$ is non-zero. Let $G_2$ be a non-zero quotient vector space of $E_2$ and $\widetilde f$ be the composition of $f$ with the quotient map $E_2\rightarrow G_2$. Let $F_1$ be the kernel of $\widetilde f$, $G_1$ be the quotient space $E_1/F_1$, and $g:G_1\rightarrow G_2$ be the $K$-linear map induced by $\widetilde{f}$. It is a $K$-linear isomorphism. By \ref{Item: slope inequality mu max}, one has
\[\widehat{\mu}_{\min}(E_1,\xi_1)\leqslant\widehat{\mu}_{\max}( G_1,\xi_{G_1})\leqslant\widehat{\mu}_{\max}({G}_2,\xi_{G_2})+h(g)\leqslant\widehat{\mu}_{\max}( G_2,\xi_{G_2})+h(f),\]
where $\xi_{G_1}$ and $\xi_{G_2}$ are the quotient norm family of $\xi_1$ and $\xi_2$, respectively.
Since $G_2$ is arbitrary, one obtains $\widehat{\mu}_{\min}( E_1,\xi_1)\leqslant\widehat{\mu}_{\min}(E_2,\xi_2)+h(f)$.

\ref{Item: slope inequality mu mix} Let $G$ be the image of $E_1$ by $f$, which is non-zero since $f$ is non-zero. We equip $G$ with the restriction $\xi_G$ of $\xi_2$ {to} $G$. As $G$ is non-zero, one has $\widehat{\mu}_{\min}(G,\xi_G)\leqslant\widehat{\mu}_{\max}(G,\xi_G)\leqslant\widehat{\mu}_{\max}( E_2,\xi_2)$. By (2), one has $\widehat{\mu}_{\min}(E_1,\xi_1)\leqslant\widehat{\mu}_{\min}(G,\xi_G)+h(f)$. Hence $\widehat{\mu}_{\min}(E_1,\xi_1)\leqslant\widehat{\mu}_{\max}(E_2,\xi_2)+h(f)$.
\end{proof}

\begin{prop}\label{Pro: exact sequence mu min}
Let $(E',\xi')$ and $(E,\xi)$ be adelic vector bundles on $S$, and $f:E'\rightarrow E$ be an injective $K$-linear map. Let $E''$ be the quotient vector space $E/f(E')$ and $\xi''$ be the quotient norm family of $\xi$ on $E''$. Then the following inequality holds
\begin{equation}\label{Equ: estimation de mu min pour suite exacte}\widehat{\mu}_{\min}(E,\xi)\geqslant \min(\widehat{\mu}_{\min}(E',\xi')-h(f),\widehat{\mu}_{\min}(E'',\xi'')).\end{equation}
If in addition $\widehat{\mu}_{\min}(E',\xi')-h(f)\geqslant\widehat{\mu}_{\min}(E'',\xi'')$, then the equality $\widehat{\mu}_{\min}(E,\xi)=\widehat{\mu}_{\min}(E'',\xi'')$ holds.
\end{prop}
\begin{proof}
The inequality \eqref{Equ: estimation de mu min pour suite exacte} is trivial if $E=\{0\}$ since in this case one has $\widehat{\mu}_{\min}(E,\xi)=+\infty$ by convention. Moreover, one has $E''=\{0\}$ since $E''$ is a quotient vector space of $E$. Therefore the equality $\widehat{\mu}_{\min}(E,\xi)=\widehat{\mu}_{\min}(E'',\xi'')$ holds.

In the following, we assume that $E\neq\{0\}$. Let $Q$ be a quotient vector space of $E$ and $\xi_Q$ be the quotient norm family of $\xi$. Let $\pi:E\rightarrow Q$ be the quotient map. If the composed map $\pi f$ is non-zero, by Proposition \ref{Pro:inegalitdepente} \ref{Item: slope inequality mu mix}, one has 
\[\widehat{\mu}_{\min}(E',\xi')\leqslant\widehat{\mu}_{\max}(Q,\xi_Q)+h(\pi f)\leqslant\widehat{\mu}_{\max}(Q,\xi_Q)+h(f).\]
Otherwise the quotient map $\pi:E\rightarrow Q$ factorises through $E''$ and by Proposition \ref{Pro:inegalitdepente} \ref{Item: slope inequality mu min} one has
\[\widehat{\mu}_{\min}(E'',\xi'')\leqslant\widehat{\mu}_{\min}(Q,\xi_Q)\leqslant\widehat{\mu}_{\max}(Q,\xi_Q).\]
Since $Q$ is arbitrary, the inequality \eqref{Equ: estimation de mu min pour suite exacte} is true.

If $\widehat{\mu}_{\min}(E',\xi')-h(f)\geqslant\widehat{\mu}_{\min}(E'',\xi'')$, then \eqref{Equ: estimation de mu min pour suite exacte} implies $\widehat{\mu}_{\min}(E,\xi)\geqslant\widehat{\mu}_{\min}(E'',\xi'')$. Moreover, by  Proposition \ref{Pro:inegalitdepente} \ref{Item: slope inequality mu min} one has $\widehat{\mu}_{\min}(E,\xi)\leqslant\widehat{\mu}_{\min}(E'',\xi'')$. Hence the equality $\widehat{\mu}_{\min}(E,\xi)=\widehat{\mu}_{\min}(E'',\xi'')$ holds.
\end{proof}

\begin{prop}\label{Pro:successive embeddings}
Let $\{(E_i,\xi_{E_i})\}_{i=1}^n$ be a family of adelic vector bundles, where $n\in\mathbb N$, $n\geqslant 2$. Assume that 
\begin{equation}\label{Equ: sequence of alphas}\xymatrix{E_0:=\{0\}\ar[r]^-{\alpha_1}&E_1\ar[r]^-{\alpha_2}&E_2\ar[r]^-{\alpha_3}&E_3\ar[r]&\cdots\ar[r]^-{\alpha_{n-1}}&E_{n-1}\ar[r]^-{\alpha_n}&E_n}\end{equation}
is a sequence of injective $K$-linear maps. For any $i\in\{1,\ldots,n\}$, let $\beta_i=\alpha_n \scirc \cdots \scirc \alpha_{i+1}$, where by convention $\beta_n:=\mathrm{Id}_{E_n}$. For any $i\in\{1,\ldots,n\}$, let $Q_i$ be the quotient space $E_i/\alpha_i(E_{i-1})$ and $\xi_{Q_i}$ be the quotient norm family of $\xi_{E_i}$. Then one has
\[\widehat{\mu}_{\min}(E_n,\xi_{E_n})\geqslant\min_{i\in\{1,\ldots,n\}} \big(\widehat{\mu}_{\min}(Q_i,\xi_{Q_i})-h(\beta_i)\big)\]
\end{prop}
\begin{proof}
The case where $n=2$ was proved in Proposition \ref{Pro: exact sequence mu min}. In the following, we assume that $n\geqslant 3$ and that the proposition has been proved for the case of $n-1$ adelic vector bundles. For any $i\in\{2,\ldots,n\}$, let $E_i'$ be the cokernel of the composed linear map
\[\alpha_i\scirc \cdots \scirc\alpha_2:E_1\longrightarrow E_i\]
and let $\xi_{E_i'}$ be the quotient norm family of $\xi_{E_i}$ on $E_i'$. Let $E_1'=\{0\}$. Then the sequence \eqref{Equ: sequence of alphas} induces a sequence of $K$-linear maps 
\begin{equation}\label{Equ: sequence of alpha prime}\xymatrix{E_1':=\{0\}\ar[r]^-{\alpha_2'}&E_2'\ar[r]^-{\alpha_3'}&E_3'\ar[r]&\cdots\ar[r]^-{\alpha_{n-1}'}&E_{n-1}'\ar[r]^-{\alpha_n'}&E_n'}.\end{equation}
For any $i\in\{2,\ldots,n\}$, let $\beta_i'=\alpha_n'\scirc \cdots \scirc\alpha_{i+1}'$, where by convention {$\beta_n'=\mathrm{Id}_{E_n'}$}. For any $i\in\{2,\ldots,n\}$, let $Q_i'$ be the quotient space $E_i'/\alpha_i'(E_{i-1}')$ and $\xi_{Q_i'}$ be the quotient norm family of $\xi_{E_i'}$. Note that $Q_i'$ is canonically isomorphic to $Q_i$, and under the canonical isomorphism $Q_i\cong Q_i'$, the norm family $\xi_{Q_i}$ identifies with $\xi_{Q_i'}$ (see Proposition~\ref{prop:quotient:norm:linear:map}). Therefore one has $\widehat{\mu}_{\min}(Q_i,\xi_{Q_i})=\widehat{\mu}_{\min}(Q_i',\xi_{Q_i'})$ for any $i\in\{2,\ldots,n\}$. Applying the induction hypothesis to \eqref{Equ: sequence of alpha prime} we obtain
\[\widehat{\mu}_{\min}(E_n',\xi_{E_n'})\geqslant\min_{i\in\{2,\ldots,n\}}(\widehat{\mu}_{\min}(Q_i,\xi_{Q_i})-h(\beta_i'))\geqslant \min_{i\in\{2,\ldots,n\}}(\widehat{\mu}_{\min}(Q_i,\xi_{Q_i})-h(\beta_i)),\]
where the second inequality comes from Proposition~\ref{prop:quotient:norm:linear:map}. Finally, by Proposition \ref{Pro: exact sequence mu min} one has 
\[\widehat{\mu}_{\min}(E_n,\xi_{E_n})\geqslant \min\big\{\widehat{\mu}_{\min}(E_n',\xi_{E_n'}),\widehat{\mu}_{\min}(E_1,\xi_{E_1})-h(\beta_1)\big\}.\]
The proposition is thus proved.
\end{proof}

\begin{prop}\label{Pro: slope of tensor}
Let $\overline E$ and $\overline F$ be adelic vector bundles on $S$. One has
\[\widehat{\mu}(\overline E\otimes_{\varepsilon,\pi}\overline F)=\widehat{\mu}(\overline E)+\widehat{\mu}(\overline F).\]
If $\overline E$ and $\overline F$ are both Hermitian, then
\[\widehat{\mu}(\overline E\otimes\overline F)=\widehat{\mu}(\overline E)+\widehat{\mu}(\overline F).\]
\end{prop}
\begin{proof}
These equalities are direct consequences of Proposition \ref{Pro: degre of tensor}. 
\end{proof}

\begin{prop}
Let $\overline E=(E,\xi_E)$ and $\overline F=(F,\xi_F)$ be adelic vector bundles on $S$. One has
\begin{equation}\label{Equ: mu max tensor lower bound}\widehat{\mu}_{\max}(\overline E\otimes_{\varepsilon,\pi}\overline F)\geqslant\widehat{\mu}_{\max}(\overline E)+\widehat{\mu}_{\max}(\overline F).\end{equation}
If $\overline E$ and $\overline F$ are  Hermitian adelic vector bundles, then 
\begin{equation}\label{Equ: mu max tensor lower bound2}\widehat{\mu}_{\max}(\overline E\otimes\overline F)\geqslant\widehat{\mu}_{\max}(\overline E)+\widehat{\mu}_{\max}(\overline F).\end{equation}
\end{prop}
\begin{proof} 
Let $E_1$ and $F_1$ be vector subspaces of $E$ and $F$ respectively. Let $\xi_{E_1}$ and $\xi_{F_1}$ be the restrictions of $\xi_E$ and $\xi_F$ {to} $E_1$ and $F_1$ respectively. By Proposition \ref{Pro: slope of tensor}, one has
\[\widehat{\mu}(\overline E_1\otimes_{\varepsilon,\pi}\overline F_1)=\widehat{\mu}(\overline E_1)+\widehat{\mu}(\overline F_1).\]
If $\xi_E$ and $\xi_F$ are both Hermitian, then 
\[\widehat{\mu}(\overline E_1\otimes\overline F_1)=\widehat{\mu}(\overline E_1)+\widehat{\mu}(\overline F_1).\] By Proposition \ref{Pro: restriction and tensors}, if we denote by $\xi$ the restriction of $\xi_E\otimes_{\varepsilon,\pi}\xi_F$ {to} $E_1\otimes F_1$, then the identity map from $\overline E_1\otimes_{\varepsilon,\pi}\overline F_1$ to $(E_1\otimes_kF_1,\xi)$   has  height $\leqslant 0$ and therefore \[\widehat{\mu}(\overline E_1)+\widehat{\mu}(\overline F_1)=\widehat{\mu}(\overline E_1\otimes_{\varepsilon,\pi}\overline F_1)\leqslant\widehat{\mu}(E_1\otimes_kF_1,\xi)\leqslant\widehat{\mu}_{\max}(\overline E\otimes_{\varepsilon,\pi}\overline F).\]
Similarly, if both norm families $\xi_E$ and $\xi_F$ are Hermitian, then by Proposition \ref{Pro: pi tensor sub Archimedean} the restriction of $\xi_E\otimes\xi_F$ {to} $E_1\otimes_KF_1$ identifies with $\xi_{E_1}\otimes\xi_{F_1}$. Hence
\[\widehat{\mu}(\overline E_1)+\widehat{\mu}(\overline F_1)=\widehat{\mu}(\overline E_1\otimes\overline F_1)\leqslant\widehat{\mu}_{\max}(\overline E\otimes\overline F).\]
Since $E_1$ and $F_1$ are arbitrary, we obtain the inequalities \eqref{Equ: mu max tensor lower bound} and \eqref{Equ: mu max tensor lower bound2}.
\end{proof}

\begin{lemm}\label{lem:mu:max:deg:plus}
Let $(E, \xi)$ be an adelic vector bundle over $S$. Then we have the following:
\begin{enumerate}[label=\rm(\arabic*)]
\item Let $\psi$ be an integrable function on $\Omega$. Then
\[
\widehat{\mu}_{\max}(E, \mathrm{e}^{\psi}\xi) = \widehat{\mu}_{\max}(E, \xi) - \int_{\Omega} \psi\, \nu(\mathrm{d}\omega).
\]

\item If $\widehat{\mu}_{\max}(E, \xi) \leqslant 0$, then $\widehat{\deg}_+(E, \xi) = 0$.
\end{enumerate}
\end{lemm}

\begin{proof}
(1) Let $F$ be a non-zero vector subspace of $E$. Then, as
\[
\widehat{\mu}(F, \mathrm{e}^{\psi}\xi_{F}) = \widehat{\mu}(F, \xi_{F}) - \int_{\Omega} \psi \,\nu(\mathrm{d}\omega),
\] 
we obtain (1).

\medskip
(2) Let $F$ be a non-zero vector subspace of $E$. By our assumption, $\widehat{\mu}(F, \xi_F) \leqslant 0$, that is,
$\widehat{\deg}(F, \xi_F) \leqslant 0$, so that the assertion follows.
\end{proof}

\subsection{Harder-Narasimhan filtration: Hermitian case} 

It had been discovered by Stuhler \cite{Stuhler76} (generalised by Grayson \cite{Grayson84}) that the Euclidean lattices and vector bundles on projective algebraic curves share some common constructions and properties such as slopes and Harder-Narasimhan filtration etc. Later Bost has developed the slope theory of Hermitian vector bundles over spectra of algebraic integer rings, see \cite[Appendice]{BostBour96} (see also \cite{Soule97} and \cite[\S4.1]{Bost2001} for more details, and \cite{Gaudron08,Bost_Chen,Gaudron_Remond14} for {further} generalisations).

The Hermitian adelic vector bundles on $S$ form a category in which a theory of Hader-Narasimhan filtration can be developed in a functorial way. We refer the readers to \cite{Andre09,Chen10a} for more details. In this subsection, we adopt a more direct approach as in \cite{Bost_Chen}.

Let $\overline E$ be a non-zero  Hermitian adelic vector bundle on $S$. We can construct in a recursive way a flag
\[0=E_0\subsetneq E_1\subsetneq \ldots\subsetneq E_n=E\]
of vector subspaces of $E$ such that $\overline{E_i/E_{i-1}}=\overline{(E/E_{i-1})}_{\mathrm{des}}$, called the \emph{Harder-Narasimhan flag}\index{Harder-Narasimhan!flag} of $\overline E$, where $E/E_{i-1}$ is equipped with the quotient norm family, and $E_{i}/E_{i-1}$ is equipped with the subquotient norm family (namely the restriction of the norm family of $\overline{E/E_{i-1}}$ {to} $E_i/E_{i-1}$).

\begin{prop}
Let $\overline E$ be a non-zero  Hermitian adelic vector bundle on $S$ and
\[0=E_0\subsetneq E_1\subsetneq \ldots\subsetneq E_n=E\]
be the Harder-Narasimhan flag of $\overline E$. Then each subquotient $\overline{E_i/E_{i-1}}$ is a semistable Hermitian adelic vector bundle. Moreover, if we let $\mu_i=\widehat{\mu}(\overline{E_i/E_{i-1}})$ for $i\in\{1,\ldots,n\}$, then one has
$\mu_1>\ldots>\mu_n$. 
\end{prop}
\begin{proof}
We reason by induction on the length $n$ of the Harder-Narasimhan flag. When $n=1$, the assertion is trivial. In the following, we suppose that $n\geqslant 2$. By definition
\[0=E_1/E_1\subsetneq E_2/E_1\subsetneq\ldots\subsetneq E_n/E_1\] is the Harder-Narasimhan flag of $\overline{E/E_1}$. Therefore the induction hypothesis leads to 
$\mu_2>\ldots>\mu_n$. It remains to establish $\mu_1>\mu_2$. Since $E_1$ is the destabilising vector subspace of $E$ and $E_2$ contains strictly $E_1$, one has \begin{equation}\label{Equ:mu1>mue2}\mu_1=\widehat{\mu}(\overline E_1)>\widehat{\mu}(\overline E_2).\end{equation}
Moreover, 
\[\xymatrix{0\ar[r]&\overline E_1\ar[r]&\overline E_2\ar[r]&\overline{E_2/E_1}\ar[r]&0}\]
forms an exact sequence of adelic vector bundles on $S$. Therefore one has
\[\widehat{\deg}(\overline E_2)=\widehat{\deg}(\overline E_1)+\widehat{\deg}(\overline{E_2/E_1})=\mu_1\rang_K(E_1)+\mu_2\rang_K(E_2/E_1).\]
By \eqref{Equ:mu1>mue2} we obtain
\[\mu_1\rang_K(E_1)+\mu_2\rang_K(E_2/E_1)<\mu_1\rang(E_2)\]
and hence $\mu_1>\mu_2$. The proposition is thus proved.
\end{proof}

The Harder-Narasimhan flag and the slopes of the successive subquotients in the previous proposition permit to construct an $\mathbb R$-filtration $\mathcal F_{\mathrm{hn}}$ on the vector space $E$, called the \emph{Harder-Narasimhan $\mathbb R$-filtration}\index{Harder-Narasimhan!R-filtration@$\mathbb R$-filtration} as follows:
\begin{equation}\label{Equ:HNfiltration}\forall\,t\in\mathbb R,\quad\mathcal F^{t}_{\mathrm{hn}}(\overline E):=E_i\;\text{ if $\mu_{i+1}<t\leqslant\mu_i$},\end{equation}
where by convention $\mu_0=+\infty$ and $\mu_{n+1}=-\infty$. If $\overline E$ is the zero Hermitian adelic vector bundle, by convention its Harder-Narasimhan $\mathbb R$-filtration is defined as the only $\mathbb R$-filtration of the zero vector space: for any $t\in\mathbb R$ one has $\mathcal F^t(\overline E)=\{0\}$.  Note that the Harder-Narasimhan $\mathbb R$-filtration is locally constant on the left, namely $\mathcal F^{t-\varepsilon}_{\mathrm{hn}}(\overline E)=\mathcal F^{t}_{\mathrm{hn}}(\overline E)$ if $\varepsilon>0$ is sufficiently small. Moreover, each subquotient
\[\mathrm{Sq}_{\mathrm{hn}}^t(\overline E):=\mathcal F_{\mathrm{hn}}^t(\overline E)/\mathcal F_{\mathrm{hn}}^{t+}(\overline E)\]
with $\mathcal F_{\mathrm{hn}}^{t+}(\overline E):=\bigcup_{\varepsilon>0}\mathcal F_{\mathrm{hn}}^{t+\varepsilon}(\overline E)$, viewed as a Hermitian adelic vector bundle in considering the \emph{subqutient norm family}\index{subquotient}\index{norm family!subquotient}, namely the quotient of the restricted norm family on $\mathcal F_{\mathrm{hn}}^t(\overline E)$, is either zero or a semistable Hermitian adelic vector bundle of slope $t$.
The following proposition shows that the Harder-Narasimhan $\mathbb R$-filtration is actually characterized by these properties.

\begin{prop}\label{Pro:characterisationdehnfiltration}
Let $\overline E$ be a non-zero Hermitian adelic vector bundle on $S$ and $\mathcal F$ be a decreasing $\mathbb R$-filtration of $E$ which is separated, exhaustive\footnote{Let $E$ be a vector space over $K$ and $(\mathcal F^t(E))_{t\in\mathbb R}$ be a decreasing $\mathbb R$-filtration of $E$. We say that the filtration $\mathcal F$ is \emph{separated}\index{separated} if $\mathcal F^t(E)=\{0\}$ for sufficiently positive $t$. We say that the filtration $\mathcal F$ is \emph{exhaustive}\index{exhaustive} if $\mathcal F^t(E)=E$ for sufficiently negative $t$. } and locally constant on the left. Assume that each subquotient $\mathrm{Sq}(\overline E):=\mathcal F^t(E)/\mathcal F^{t+}(E)$ equipped with the subquotient norm family, is either zero or a semistable Hermitian adelic vector bundle of slope $t$. Then the $\mathbb R$-filtration $\mathcal F$ coincides with the Harder-Narasimhan $\mathbb R$-filtration of $\overline E$.
\end{prop}
\begin{proof}
We will prove an alternative statement as follows. Let 
\begin{equation}\label{Equ:flagF}0=F_0\subsetneq F_1\subsetneq\ldots\subsetneq F_m\subsetneq E\end{equation}
be a flag of vector subspaces of $E$. We will prove that, if each subquotient $\overline{F_i/F_{i-1}}$ ($i\in\{1,\ldots,m\}$) is a semistable Hermitian adelic vector bundle and if the relations
\[\widehat{\mu}(\overline{F_1/F_0})>\ldots>\widehat{\mu}(\overline{F_{m}/F_{m-1}})\]
hold, then \eqref{Equ:flagF} identifies with the Harder-Narasimhan flag of $\overline E$. This alternative statement is actually equivalent to the form announced in the proposition. In fact, the data of an $\mathbb R$-filtration of $E$ is equivalent to that of a flag (of vector subspaces figuring in the $\mathbb R$-filtration) and a decreasing sequenc of real numbers indicating the indices where the $\mathbb R$-filtration has jumps see \S\ref{Subsection: Trivial valuation}, notably Remark \ref{Rem: R-filtration as flag plus slopes}. 

We will prove the statement by induction on the rank of $E$. The case where $\rang_K(E)=1$ is trivial. In the following, we assume that $\rang_K(E)\geqslant2$ and that the alternative assertion has been proved for any Hermitian adelic vector bundle of rank $<\rang_K(E)$. Let
\[0=E_0\subsetneq E_1\subsetneq\ldots\subsetneq E_n=E\]
be the Harder-Narasimhan flag of $\overline E$. We begin by showing that $F_1=E_1$. Since $E_1$ is the destabilizing vector subspace of $\overline E$, one has $\widehat{\mu}(\overline F_1)\leqslant\widehat{\mu}(\overline E_1)$. Moreover one has
\[0=F_0\cap E_1\subseteq F_1\cap E_1\subseteq F_2\cap E_1\subseteq\ldots F_m\cap E_1=E_1.\]
Note that each subquotient $(F_{i}\cap E_1)/(F_{i-1}\cap E_1)$ identifies with a vector subspace of $F_i/F_{i-1}$. Since $\overline{F_i/F_{i-1}}$ is semistable, one has
\[\begin{split}\widehat{\deg}\big(\overline{(F_{i}\cap E_1)/(F_{i-1}\cap E_1)}\big)&\leqslant\widehat{\mu}(\overline{F_i/F_{i-1}})\rang_K\big({(F_{i}\cap E_1)/(F_{i-1}\cap E_1)}\big)\\
&\leqslant\widehat{\mu}(\overline F_1)\rang_K\big({(F_{i}\cap E_1)/(F_{i-1}\cap E_1)}\big),
\end{split}\]
where the second inequality is strict if $i>1$ and if ${(F_{i}\cap E_1)/(F_{i-1}\cap E_1)}$ is non-zero. Therefore we obtain
\begin{equation}\label{Equ:sommededegre}\begin{split}\widehat{\deg}(\overline E_1)&=\sum_{i=1}^m\widehat{\deg}\big(\overline{(F_{i}\cap E_1)/(F_{i-1}\cap E_1)}\big)\leqslant\widehat{\mu}(\overline F_1)\rang(E_1).
\end{split}\end{equation}
Combining with the inequality {$\widehat{\mu}(\overline F_1)\leqslant\widehat{\mu}(\overline E_1)=\widehat{\mu}_{\max}(\overline E)$}, we deduce that the inequality \eqref{Equ:sommededegre} is actually an equality, which also implies that ${(F_{i}\cap E_1)/(F_{i-1}\cap E_1)}=\{0\}$ once $i>1$. Therefore one has $F_1=E_1$, which also leads to the alternative assertion in the particular case where $\overline E$ is semistable.

In the case where $\overline E$ is not semistable, namely $n\geqslant 2$, note that 
\begin{equation}\label{Equ:hdfiltrationofeE}0=E_1/E_1\subsetneq E_2/E_1\subsetneq\ldots\subsetneq E_{n}/E_{1}=E/E_1\end{equation} is the Harder-Narasimhan flag of $\overline{E/E_1}=\overline{E/F_1}$. By the induction hypothesis applied to $\overline{E/F_1}$, we obtain that the flag
\[0=F_1/F_1\subsetneq F_2/F_1\subsetneq\ldots\subsetneq F_m/F_1=E/F_1\]
coincides with \eqref{Equ:hdfiltrationofeE}. The proposition is thus proved. 
\end{proof}

\begin{defi}
Let $\overline E$ be a non-zero Hermitian adelic vector bundle on $S$, and 
\[0=E_0\subsetneq E_1\subsetneq \ldots\subsetneq E_n=E\]
be its Harder-Narasimhan flag.  For any $i\in\{1,\ldots,\rang_K(E)\}$, there exists a unique $j\in\{1,\ldots,n\}$ such that $\rang_K(E_{j-1})<i\leqslant \rang_K(E_j)$. We let $\widehat{\mu}_i(\overline E)$ be 
the slope $\widehat{\mu}(\overline{E_j/E_{j-1}})$, called the \emph{$i$-th slope}\index{ith slope@$i$-th slope} of $\overline E$. Clearly one has
\[\widehat{\mu}_1(\overline E)\geqslant\ldots\geqslant\widehat{\mu}_r(\overline E),\]
where $r$ is the rank of $E$ over $K$. Moreover, by definition $\widehat{\mu}_1(\overline E)$ coincides with the maximal slope of $\overline E$.
\end{defi}

{
\begin{rema}
As in the classic case of vector bundles on projective curves or Hermitian vector bundles over algebraic integer rings, one can naturally construction Harder-Narasimhan polygones associated with Hermitian adelic vector bundles on  adelic curves.  Let $\overline E$ be a non-zero Hermitian adelic vector bundle on the adelic curve $S$. We consider the convex hull $C_{\overline E}$ in $\mathbb R^2$ of the points 
$(\rang_K(F),\widehat{\deg}(\overline F))$, where $F$ runs over the set of all vector subspaces of $E$. The upper boundary of this convex set identifies with the graph of a concave function $P_{\overline E}$ on $[0,\rang_K(E)]$ which is affine on each interval $[i-1,i]$ with $i\in\{1,\ldots,\rang_K(E)\}$. This function is called \emph{Harder-Narasimhan polygon}\index{Harder-Narasimhan!polygon} of $\overline E$.  If
\[0=E_0\subsetneq E_1\subsetneq \ldots\subsetneq E_n=E\] is the Harder-Narasimhan flag of $\overline E$, then the abscissae on which the Harder-Narasimhan polygon $P_{\overline E}$ changes slopes are exactly $\rang_K(E_i)$ for $i\in\{0,\ldots,n\}$. Moreover, the value of $P_{\overline E}$ on $\rang_K(E_i)$ is $\widehat{\deg}(\overline E_i)$.  
\end{rema}}

\begin{prop}\label{Pro:mumindudual}
Let $\overline E$ be a non-zero Hermitian adelic vector bundle on $S$ and $r$ be the rank of $E$ over $K$. One has
\begin{equation}\label{Equ:mumindudual}\widehat{\mu}_{r}(\overline E)=-\widehat{\mu}_{1}(\overline E{}^\vee).\end{equation}
In particular, $\widehat{\mu}_{r}(\overline E)$ is equal to $\widehat{\mu}_{\min}(\overline E)$.
\end{prop}
\begin{proof}
Let 
\[0=E_0\subsetneq E_1\subsetneq \ldots\subsetneq E_n=E\]
be the Harder-Narasimhan flag of $\overline E$. Note that
\begin{equation}\label{Equ:hddeEvee}0=(E/E_n)^\vee\subsetneq (E/E_{n-1})^{\vee}\subsetneq\ldots\subsetneq (E/E_0)^{\vee}=E^\vee\end{equation}
is a flag of vector subspaces of $E^\vee$, and for $i\in\{1,\ldots,n\}$ one has
\[(E/E_{i-1})^\vee/(E/E_i)^\vee\cong (E_i/E_{i-1}).\]
By Proposition \ref{Pro:characterisationdehnfiltration} (notably the alternative form stated in the proof), we obtain that \eqref{Equ:hddeEvee} is actually the Harder-Narasimhan flag of $\overline E{}^\vee$. Therefore
\[\widehat{\mu}_{1}(\overline E{}^\vee)=\widehat{\mu}((\overline{E/E_{n-1}})^\vee)=-\widehat{\mu}(\overline{E/E_{n-1}})=-\widehat{\mu}_{r}(\overline E),\]
where the second equality comes from Proposition \ref{Pro:degdual}.

Note that $\overline{E/E_{n-1}}$ is a non-zero quotient Hermitian adelic bundle of $\overline E$ which is semistable (so that $\widehat{\mu}(\overline{E/E_{n-1}})=\widehat{\mu}_{\max}(\overline{E/E_{n-1}})$). Therefore one has $\widehat{\mu}_{\min}(\overline E)\leqslant\widehat{\mu}_r(\overline E)$.  
Conversely, if $F$ is a vector subspace of $E$ such that $F\subsetneq E$ and $G$ is the quotient space $E/F$. Then $G^\vee$ identifies with a non-zero vector subspace of $E^\vee$. Moreover, by Proposition \ref{Pro:dualquotient} the dual of the quotient norm family of $\overline G$ identifies with the restriction of the dual norm family in {the adelic vector bundle structure of} $\overline E{}^\vee$. Hence one has
\[\widehat{\mu}(\overline G{}^\vee)\leqslant\widehat{\mu}_{1}(\overline E{}^\vee)=-\widehat{\mu}_{r}(\overline E).\]
Still by Proposition  \ref{Pro:degdual}, one obtains
\[\widehat{\mu}_{\max}(\overline G)\geqslant \widehat{\mu}(\overline G)\geqslant\widehat{\mu}_{r}(\overline E).\]
The equality $\widehat{\mu}_{\min}(\overline E)=\widehat{\mu}_r(\overline E)$ is thus proved. 
\end{proof}

The following proposition, which results from the slope inequalities, provides the functoriality of Harder-Narasimhan $\mathbb R$-filtration (see \cite{Chen10a} for the meaning of the functoriality of Harder-Narasimhan $\mathbb R$-filtration).

\begin{prop}\label{Pro:fonctorialite}
Let $\overline E$ and $\overline F$ be two Hermitian adelic vector bundles on $S$, and $f:E\rightarrow F$ be a non-zero $K$-linear map. For any $t\in\mathbb R$ one has 
\[f(\mathcal F^t_{\mathrm{hn}}(\overline E))\subseteq\mathcal F_{\mathrm{hn}}^{t-h(f)}(\overline F).\]
\end{prop}
\begin{proof}
We will actually show by contradiction that the composition of maps
\[\xymatrix{\relax\mathcal F_{\mathrm{hn}}^t(\overline E)\ar[r]^-f& F\ar[r]&F/\mathcal F_{\mathrm{hn}}^{t-h(f)}(\overline F)}\]
is zero. If this map is not zero, then by Proposition \ref{Pro:inegalitdepente} \ref{Item: slope inequality mu mix} one obtains
\[\widehat{\mu}_{\min}(\mathcal F_{\mathrm{hn}}^t(\overline E))\leqslant\widehat{\mu}_{\max}(\overline{F}/\mathcal F_{\mathrm{hn}}^{t-h(f)}(\overline F))+h(f).\]
By \eqref{Equ:HNfiltration} we obtain
\[t\leqslant\widehat{\mu}_{\min}(\mathcal F_{\mathrm{hn}}^t(\overline E))\leqslant\widehat{\mu}_{\max}(\overline{F}/\mathcal F_{\mathrm{hn}}^{t-h(f)}(\overline F))+h(f)<t-h(f)+h(f)=t,\]
which leads to a contradiction.
\end{proof}

\begin{coro}\label{Cor:coonstructionofFHN}
Let $\overline E$ be a non-zero Hermitian adelic vector bundle on $S$. One has
\[\mathcal F_{\mathrm{hn}}^t(\overline E)=\sum_{\begin{subarray}{c}0\neq F\in\Theta(E)\\
\widehat{\mu}_{\min}(\overline F)\geqslant t\end{subarray}}F,\]
where $F$ runs over the set $\Theta(E)$ of all non-zero vector subspaces of $E$ with minimal slope $\geqslant t$. In other words, $\mathcal F^t_{\mathrm{hn}}(\overline E)$ is the largest vector subspace of $E$ whose minimal slope is bounded from below by $t$. 
\end{coro}
\begin{proof}
By the definition of the Harder-Narasimhan $\mathbb R$-filtration (see \eqref{Equ:HNfiltration}), for any $t\in\mathbb R$ one has
$\widehat{\mu}_{\min}(\mathcal F_{\mathrm{hn}}^t(\overline E))\geqslant t$. Moreover, if $F$ is a non-zero vector subspace of $E$, then one has $\mathcal F^t_{\mathrm{hn}}(\overline F)=F$ provided that $t\leqslant\widehat{\mu}_{\min}(\overline F)$. Therefore the proposition \ref{Pro:fonctorialite} applied to the inclusion map $F\rightarrow E$ leads to $F\subseteq\mathcal F_{\mathrm{hn}}^t(E)$.
\end{proof}

\begin{prop}
Let $\overline E$ be a non-zero Hermitian adelic vector bundle on $S$ and $r$ be the rank of $E$ over $K$. The following equalities hold:
\begin{gather}\label{Equ: degree sum}
\widehat{\deg}(\overline E)=\sum_{i=1}^r\widehat{\mu}_i(\overline E)=-\int_{\mathbb R}t\,\mathrm{d}\rang(\mathcal F_{\mathrm{hn}}^t(\overline E)),\\
\label{Equ: degree sum+}\widehat{\deg}_+(\overline E)=\sum_{i=1}^r\max\{\widehat{\mu}_i(\overline E),0\}=\int_0^{+\infty}\rang(\mathcal F^t_{\mathrm{hn}}(\overline E))\,\mathrm{d}t.
\end{gather}
\end{prop}
\begin{proof}
By definition the sum of Dirac measures 
\[\sum_{i=1}^r\delta_{\widehat{\mu}_i(\overline E)}\]
identifies with the derivative $-\mathrm{d}\rang(\mathcal F^t_{\mathrm{hn}}(\overline E))$ in the sense of distribution. Therefore, the second equality in \eqref{Equ: degree sum} is true, and the second equality in \eqref{Equ: degree sum+} follows from the relation
\[\sum_{i=1}^r\max\{\widehat{\mu}_i(\overline E),0\}=-\int_0^{+\infty}t\,\mathrm{d}\rang(\mathcal F_{\mathrm{hn}}^t(\overline E))\] 
and the formula of integration by part.

Let \[0=E_0\subsetneq E_1\subsetneq \ldots\subsetneq E_n=E\]
be the Harder-Narasimhan flag of $\overline E$. By Proposition \ref{Pro:suiteexactedeg} one has
\[\widehat{\deg}(\overline E)=\sum_{j=1}^n\widehat{\deg}(\overline{E_j/E_{j-1}})=\sum_{j=1}^n\widehat{\mu}(\overline{E_j/E_{j-1}})\rang_K(E_j/E_{j-1})=\sum_{i=1}^r\widehat{\mu}_i(\overline E),\]
which proves \eqref{Equ: degree sum}. 

Let $\ell$ be the largest element in $\{1,\ldots,n\}$ such that $\widehat{\mu}(\overline E_{\ell}/\overline E_{\ell-1})\geqslant 0$. If $\widehat{\mu}(\overline {E_j/E_{j-1}})<0$ for any $j\in\{1,\ldots,n\}$, by convention we let $\ell=0$. Then by \eqref{Equ: degree sum} one has
\[\widehat{\deg}(\overline E_{\ell})=\sum_{i=1}^r\max\{\widehat{\mu}_{i}(\overline E),0\}.\]
Hence we obtain
\[\sum_{i=1}^r\max\{\widehat{\mu}_{i}(\overline E),0\}\leqslant \widehat{\deg}_+(\overline E).\]
Conversely, if $F$ is a non-zero vector subspace of $E$ and $m$ is its rank over $K$, by Proposition \ref{Pro:fonctorialite} one has $\widehat{\mu}_i(\overline F)\leqslant\widehat{\mu}_i(\overline E)$ for any $i\in\{1,\ldots,m\}$. Therefore by \eqref{Equ: degree sum} one obtains
\[\widehat{\deg}(\overline F)=\sum_{i=1}^m\widehat{\mu}_i(\overline F)\leqslant\sum_{i=1}^m\widehat{\mu}_i(\overline E)\leqslant\sum_{i=1}^r\max\{\widehat{\mu}_i(\overline E),0\}.\]
\end{proof}

\subsection{Harder-Narasimhan filtration: general case}
Inspired by Corollary \ref{Cor:coonstructionofFHN}, we extend the definition of Harder-Narasimhan $\mathbb R$-filtration to the setting of general adelic vector bundles.

\begin{defi}
Let $\overline E$ be a non-zero adelic vector bundle on $S$. For any $t\in\mathbb R$, let
\begin{equation}\label{Equ: Harder-Narasimhan filtration in general}\mathcal F_{\mathrm{hn}}^t(\overline E):=\bigcap_{\epsilon>0}\sum_{\begin{subarray}{c}{\{0\}\neq F \in\Theta( E)}\\
\widehat{\mu}_{\min}(\overline F)\geqslant t-\epsilon\end{subarray}}F,\end{equation}
{where $\Theta(E)$ denotes the set of vector subspaces of $E$.}
 By the finiteness of maximal and minimal slopes proved in Proposition \ref{Pro: finiteness of slopes}, we obtain that $\mathcal F^{t}_{\mathrm{hn}}(\overline E)=E$ when $t$ is sufficiently negative, and $\mathcal F^{t}_{\mathrm{hn}}(\overline E)=\{0\}$ is sufficiently positive. By convention we let $\mathcal F_{\mathrm{hn}}^{+\infty}(\overline E)=\{0\}$ and $\mathcal F_{\mathrm{hn}}^{-\infty}=E$.
\end{defi}

\begin{prop}\label{Pro: general hn filtration}
Let $\overline E$ be a non-zero adelic vector bundle on $S$. For any $t\in\mathbb R$, the vector space $\mathcal F_{\mathrm{hn}}^t(\overline E)$ equipped with the induced norm family has a minimal slope $\geqslant t$. In particular, one has
\[\forall\,t\in\mathbb R,\quad\mathcal F_{\mathrm{hn}}^t(\overline E)=\sum_{\begin{subarray}{c}{\{0\}\neq F \in\Theta( E)}\\
\widehat{\mu}_{\min}(\overline F)\geqslant t
\end{subarray}}F\]and
\[\widehat{\mu}_{\min}(\overline E)=\max\{t\in\mathbb R\,:\,\mathcal F_{\mathrm{hn}}^t(\overline E)=E\}.\]
\end{prop}
\begin{proof} {Let $t\in\mathbb R$. For sufficiently small $\varepsilon>0$, one has
\[F_{\mathrm{hn}}^t(\overline E)=\sum_{\begin{subarray}{c}{\{0\}\neq F \in\Theta( E)}\\
\widehat{\mu}_{\min}(\overline F)\geqslant t-\epsilon\end{subarray}}F\]}
Let $M$ be a non-zero quotient vector space of $\mathcal F^t_{\mathrm{hn}}(\overline E)$. By definition, for any $\epsilon>0$ there exists a vector subspace $F_\epsilon$ of $E$ such that $\widehat{\mu}_{\min}(\overline{ F_\epsilon})\geqslant t-\epsilon$ and that the composed map $F_\epsilon\rightarrow\mathcal F_{\mathrm{hn}}^t(\overline E)\rightarrow M$ is non-zero. By the slope inequality (see Proposition \ref{Pro:inegalitdepente}) we have $t-\epsilon\leqslant\widehat{\mu}_{\min}(\overline{F_\epsilon})\leqslant\widehat{\mu}_{\max}(\overline M)$, which leads to $\widehat{\mu}_{\max}(\overline M)\geqslant t$ since $\epsilon>0$ is arbitrary. As $M$ is arbitrary, we obtain the first statement.

By the first statement of the proposition, for any $t\in\mathbb R$ such that $\mathcal F_{\mathrm{hn}}^t(\overline E)=E$, one has $\widehat{\mu}_{\min}(\overline E)\geqslant t$. Conversely, by definition, if $t$ is a real number such that $\widehat{\mu}_{\min}(\overline E)\geqslant t$, then $E\subseteq\mathcal F_{\mathrm{hn}}^t(\overline E)$ and hence $E=\mathcal F_{\mathrm{hn}}^t(\overline E)$. Therefore, the equality $\widehat{\mu}_{\min}(\overline E)=\max\{t\in\mathbb R\,:\,\mathcal F_{\mathrm{hn}}^t(\overline E)=E\}$ holds.
\end{proof}

\begin{defi}
Let $\overline E$ be a non-zero adelic vector bundle on $S$. For any $i\in\{1,\ldots,\rang_K(E)\}$, we let
\[\widehat{\mu}_i(\overline E):=\sup\{t\in\mathbb R\,:\,\rang_K(\mathcal F_{\mathrm{hn}}^t(\overline E))\geqslant i\}.\]
The number $\widehat{\mu}_i(\overline E)$ is called the \emph{$i$-th slope}\index{ith slope@$i$-th slope} of $\overline E$. Proposition \ref{Pro: general hn filtration} shows that the last slope of $\overline E$ identifies with the minimal slope $\widehat{\mu}_{\min}(\overline E)$ of $\overline E$. 
\end{defi}

\begin{rema}\label{Rem: comparison between mu max and mu 1}
Let $\overline E$ be a non-zero adelic vector bundle. In general the first slope $\widehat{\mu}_1(\overline E)$ does not coincide with $\widehat{\mu}_{\max}(\overline E)$ and we only have an inequality $\widehat{\mu}_1(\overline E)\leqslant\widehat{\mu}_{\max}(\overline E)$. Moreover, if the norm family of $\overline E$ is ultrametric on $\Omega\setminus\Omega_\infty$, then one has $ \widehat{\mu}_{\max}(\overline E)\leqslant\widehat{\mu}_1(\overline E)+\frac 12\ln(\rang_K(E))\nu(\Omega_\infty)$.
This follows from \eqref{Equ: degree bounded from above by sum of slopes delta} and \eqref{Equ:encadrementdelta}.
\end{rema}

With the extended definition, the statement of Proposition \ref{Pro:fonctorialite} still holds for general adelic vector bundles.

\begin{prop}
\label{Pro:fonctorialite general}
Let $\overline E$ and $\overline F$ be adelic vector bundles on $S$, and $f:E\rightarrow F$ be a non-zero $K$-linear map. For any $t\in\mathbb R$ one has 
\[f(\mathcal F^t_{\mathrm{hn}}(\overline E))\subseteq\mathcal F_{\mathrm{hn}}^{t-h(f)}(\overline F).\]
\end{prop}
\begin{proof}
Let $M$ be a non-zero vector subspace of $E$ such that $\widehat{\mu}_{\min}(\overline M)\geqslant t$. By Proposition \ref{Pro:inegalitdepente} \ref{Item: slope inequality mu min}, one has 
\[\widehat{\mu}_{\min}(\overline M)\leqslant\widehat{\mu}_{\min}(\overline{f(M)})+h(f|_M)\leqslant\widehat{\mu}_{\min}(\overline{f(M)})+h(f).\]
Therefore $f(M)\subseteq\mathcal F_{\mathrm{hn}}^{t-h(f)}(F)$. 
\end{proof}

\begin{prop}\label{Pro: degree bounded by integrals}
Let $\overline E$ be a non-zero adelic vector bundle on $S$. Let $r$ be the rank of $E$ over $K$. Then the following inequalities hold:
\begin{gather}\label{Equ: degree bounded from below by sum of slopes}
\widehat{\deg}(\overline E)\geqslant\sum_{i=1}^r\widehat{\mu}_i(\overline E)=-\int_{\mathbb R}t\,\mathrm{d}\,\rang(\mathcal F_{\mathrm{hn}}^t(\overline E)),\\ \label{Equ: positive degree bounded by sum of positive slopes}
\widehat{\deg}_+(\overline E)\geqslant\sum_{i=1}^r\max\{\widehat{\mu}_i(\overline E),0\}=\int_0^{+\infty}\rang(\mathcal F_{\mathrm{hn}}^t(\overline E))\,\mathrm{d}t.
\end{gather}
\end{prop}
\begin{proof}
For $i\in\{1,\ldots,r\}$, let $E_i$ be $\mathcal F_{\mathrm{hn}}^{\widehat{\mu}_i(\overline E)}(\overline E)$. Let $E_0=\{0\}$. Then for each $i\in\{1,\ldots,r\}$, such that $E_i\supsetneq E_{i-1}$ one has
\[\widehat{\mu}(\overline{E_i/E_{i-1}})\geqslant\widehat{\mu}_{\min}(\overline E_i)=\widehat{\mu}_i(\overline E),\]
where the last equality comes from Proposition \ref{Pro: general hn filtration} and the fact that the restriction of the $\mathbb R$-filtration $\mathcal F_{\mathrm{hn}}$ {to} $E_i$ coincides with the Harder-Narasimhan $\mathbb R$-filtration of $\overline E_i$. Therefore, by Proposition \ref{Pro:suiteexactedeg} one has
\[\widehat{\deg}(\overline E)\geqslant\sum_{\begin{subarray}{c}
i\in\{1,\ldots,r\}\\
E_i\supsetneq E_{i-1}
\end{subarray}}\widehat{\deg}(\overline{E_i/E_{i-1}})\geqslant \sum_{\begin{subarray}{c}
i\in\{1,\ldots,r\}\\
E_i\supsetneq E_{i-1}
\end{subarray}}\rang(E_i/E_{i-1})\widehat{\mu}_i(\overline E)=\sum_{i=1}^r\widehat{\mu}_i(\overline E),\]
which proves \eqref{Equ: degree bounded from below by sum of slopes}. Finally, if we let $j$ be the largest index in $\{1,\ldots,r\}$ such that $\widehat{\mu}_j(\overline E)\geqslant 0$. Then by what we have proved 
\[\widehat{\deg}(\overline E_j)\geqslant\sum_{i=1}^j\widehat{\mu}_i(\overline E)=\sum_{i=1}^r\max(\widehat{\mu}_i(\overline E),0).\]
Therefore, the inequality \eqref{Equ: positive degree bounded by sum of positive slopes} holds.
\end{proof}

\begin{prop}\label{Pro: degree bounded above by integrals}
Let $\overline E=(E,\xi)$ be a non-zero adelic vector bundle on $S$. Let $r$ be the rank of $E$ over $K$. Then the following inequalities hold:
\begin{equation}\label{Equ: degree bounded from above by sum of slopes}
\widehat{\deg}(\overline E)\leqslant\sum_{i=1}^r\widehat{\mu}_i(\overline E)+\Delta(\overline E).\end{equation}
{If in addition $\xi$ is ultrametric on $\Omega\setminus\Omega_\infty$, then one has
\begin{equation}\label{Equ: degree bounded from above by sum of slopes delta}
\widehat{\deg}(\overline E)\leqslant\sum_{i=1}^r\widehat{\mu}_i(\overline E)+\delta(\overline E).\end{equation}}
\end{prop}
\begin{proof}
We reason by induction on the rank of $E$ over $K$. The case where $\rang_K(E)=1$ is trivial since in this case $\overline E$ is Hermitian. In the following, we assume that $\rang_K(E)>1$ and that the proposition has been proved for adelic vector bundles of rank $<\rang_K(E)$. 

The Harder-Narasimhan $\mathbb R$-filtration corresponds to an increasing flag
\[0=E_1\subsetneq E_2\subsetneq\ldots\subsetneq E_n=E\]
and a decreasing sequence of numbers $\mu_1>\ldots>\mu_n$ corresponding to the points of jump of the $\mathbb R$-filtration. By Proposition \ref{Pro: general hn filtration}, the minimal slope of $\overline E$ is equal to $\mu_n$.

Let $\epsilon$ be a positive number such that $\epsilon<\mu_{n-1}-\mu_{n}$ and $E'$ be a vector subspace of $E$ such that $E'\subsetneq E$ and $\widehat{\mu}_{\max}(\overline{E/E'})\leqslant\widehat{\mu}_{\min}(\overline E)+\epsilon=\mu_n+\epsilon$. By Proposition \ref{Pro: mu min equivalent to classiical one}, one has
\[\widehat{\mu}(\overline{E/E'})\geqslant\widehat{\mu}_{\min}(\overline E)=\mu_n.\]
Therefore, one has
\[\mu_n\leqslant\widehat{\mu}(\overline{E/E'})\leqslant\widehat{\mu}_{\max}(\overline{E/E'})\leqslant\mu_n+\epsilon.\]
Moreover, by Proposition \ref{Pro: general hn filtration}, one has
\[\widehat{\mu}_{\min}(\overline E_{n-1})\geqslant\mu_{n-1}>\mu_n+\epsilon.\]
By Proposition \ref{Pro:inegalitdepente}, we obtain that the composed map $E_{n-1}\rightarrow E\rightarrow E/E'$ is zero, or equivalently, $E_{n-1}$ is contained in $E'$. Note that for any vector subspace $F$ of $E'$ such that $F\supsetneq E_{n-1}$ one has $\widehat{\mu}_{\min}(\overline F)\leqslant\mu_n$, otherwise the Harder-Narasimhan $\mathbb R$-filtration of $\overline E$ could not correspond to the flag $0=E_1\subsetneq E_2\subsetneq\ldots\subsetneq E_n=E$ and the decreasing sequence $\mu_1>\ldots>\mu_n$. Therefore, the Harder-Narasimhan $\mathbb R$-filtration of $\overline E{}'$ corresponds to a flag of the form ($\ell\in\mathbb N$)
\[0=E_1\subsetneq \ldots\subsetneq E_{n-1}\subsetneq E_n'\subsetneq \ldots\subsetneq E_{n-1+\ell}'=E'\]
together with a decreasing sequence $\mu_1>\ldots>\mu_{n-1}>\mu_n'>\ldots>\mu_{n-1+\ell}'$, where $\mu_n'\leqslant\mu_n$ whenever $\ell\geqslant 1$. We apply the induction hypothesis to $\overline E{}'$ and obtain
\[\widehat{\deg}(\overline E{}')\leqslant\sum_{i=1}^{n-1}\mu_i\rang_K(E_i/E_{i-1})+\sum_{i=n}^{n-1+\ell}\mu_i'\rang_K(E_{i}'/E_{i-1'})+\Delta(\overline E{}'),\]
with the convention $E_{n-1}=E_{n-1}'$. By the condition that $\mu_n'\leqslant\mu_n$ whenever $\ell\geqslant 1$ we obtain
\[\widehat{\deg}(\overline E{}')\leqslant\sum_{i=1}^{n-1}\mu_i\rang_K(E_i/E_{i-1})+\mu_n\rang_K(E'/E_{n-1})+\Delta(\overline E{}').\]
Finally, by Proposition \ref{Pro:suiteexactedeg} (notably the inequality \eqref{Equ:degupperbound Delta}) one obtains
\[\begin{split}\widehat{\deg}(\overline E)&\leqslant\widehat{\deg}(\overline E{}')+\widehat{\deg}(\overline{E/E'})+\Delta(\overline E)-\Delta(\overline E{}')-\Delta(\overline{E/E'})\\
&\leqslant\sum_{i=1}^{n-1}\mu_i\rang_K(E_i/E_{i-1})+\mu_n\rang_K(E'/E_{n-1})+(\mu_n+\epsilon)\rang_K(E/E')+\Delta(\overline E)\\
&=\sum_{j=1}^r\widehat{\mu}_j(\overline E)+\epsilon\rang(E/E')+\Delta(\overline E)\leqslant\sum_{j=1}^r\widehat{\mu}_j(\overline E)+\epsilon\rang_K(E)+\Delta(\overline E).
\end{split}\]
Since $\epsilon$ is arbitrary, we obtain the inequality \eqref{Equ: degree bounded from above by sum of slopes}. In the case where $\xi$ is ultrametric on $\Omega\setminus\Omega_\infty$ as above allows to deduce \eqref{Equ: degree bounded from above by sum of slopes delta} from \eqref{Equ:degupperbounde}. The proposition is thus proved.
\end{proof}

\begin{coro}
Let $\overline E=(E,\xi)$ be a non-zero adelic vector bundle on $S$ and $r$ be the rank of $E$ over $K$. One has
\begin{equation}\label{Equ: bound of deg plus E adelic vector bundle}
\widehat{\deg}_+(\overline E)\leqslant \sum_{i=1}^{r}\max\{\widehat{\mu}_i(\overline E),0\}+\Delta(\overline E).
\end{equation}
If in addition $\xi$ is ultrametric on $\Omega\setminus\Omega_\infty$, then one has
\begin{equation}\label{Equ: bound of deg plus E adelic vector bundle bis}
\widehat{\deg}_+(\overline E)\leqslant \sum_{i=1}^{r}\max\{\widehat{\mu}_i(\overline E),0\}+\delta(\overline E).
\end{equation}
\end{coro}
\begin{proof}
Let $F$ be a non-zero vector subspace of $E$ and $m$ be the rank of $F$ over $K$. By \eqref{Equ: degree bounded from above by sum of slopes} one has
\[\widehat{\deg}(\overline F)\leqslant\sum_{j=1}^m\widehat{\mu}_j(\overline F)+\Delta(\overline F)\leqslant\sum_{j=1}^m\max\{\widehat{\mu}_j(\overline F),0\}+\Delta(\overline F).\]
Note that 
\[\sum_{j=1}^m\max\{\widehat{\mu}_j(\overline F),0\}=-\int_{\mathbb R}\max\{t,0\}\,\mathrm{d}(\rang_K(\mathcal F_{\mathrm{hn}}^t(\overline F)))=\int_0^{+\infty}\rang_K(\mathcal F_{\mathrm{hn}}^t(\overline F))\,\mathrm{d}t.\]
Moreover, by Proposition \ref{Pro:fonctorialite general}, for any $t\in\mathbb R$, one has
\[\rang_K(\mathcal F_{\mathrm{hn}}^t(\overline F))\leqslant\rang_K(\mathcal F_{\mathrm{hn}}^t(\overline E)).\]
Therefore,
\[\widehat{\deg}(\overline F)\leqslant\int_0^{+\infty}\rang_K(\mathcal F_{\mathrm{hn}}^t(\overline E))\,\mathrm{d}t+\Delta(\overline F)=\sum_{i=1}^{r}\max\{\widehat{\mu}_i(\overline E),0\}+\Delta(\overline F).\]
Note that $\Delta(\overline F)\leqslant\Delta(\overline E)$ (see Corollary \ref{Cor: comparaison de Delta}). By taking the supremum with respect to $F$, we obtain the inequality \eqref{Equ: bound of deg plus E adelic vector bundle}.

The proof of the inequality \eqref{Equ: bound of deg plus E adelic vector bundle bis} is quite similar, where we combine the above argument with the inequality \eqref{Equ: degree bounded from above by sum of slopes delta}.
\end{proof}

\begin{rema}\label{remark:deg:plus:mu:max}
Let $(E, \xi)$ be an adelic vector bundle on $S$. Then one has the following inequality:
if $\widehat{\deg}_+(E,\xi) > 0$, then
\begin{equation}\label{eqn:remark:deg:plus:mu:max:01}
\widehat{\deg}_+(E,\xi) \leqslant \mathrm{rk}_K(E)\,\widehat{\mu}_{\max}(E, \xi).
\end{equation}
As a consequence, we {obtain}
\begin{equation}\label{eqn:remark:deg:plus:mu:max:02}
\widehat{\deg}_+(E,\xi) \leqslant \mathrm{rk}_K(E)\max\{\widehat{\mu}_{\max}(E, \xi), 0 \}
\end{equation}
in general. The inequality \eqref{eqn:remark:deg:plus:mu:max:02} is weaker than \eqref{Equ: bound of deg plus E adelic vector bundle}
and \eqref{Equ: bound of deg plus E adelic vector bundle bis}, 
but it holds without {an error term}.
Moreover, the inequality \eqref{eqn:remark:deg:plus:mu:max:01} can be proved as follow:
for any $\epsilon \in \intervalle{]}{0}{\ \widehat{\deg}_+(E,\xi)}{[}$,
one can find a non-zero vector subspace $F$ of $E$ such that
$0 \leqslant \widehat{\deg}_+(E,\xi) - \epsilon \leqslant \widehat{\deg}(F, \xi_F)$,
so that
\[
0 < \frac{\widehat{\deg}_+(E, \xi) - \epsilon}{\mathrm{rk}_K(E)} \leqslant
\frac{\widehat{\deg}(F, \xi_F)}{\mathrm{rk}_K(E)} \leqslant \widehat{\mu}(F, \xi_F) \leqslant \widehat{\mu}_{\max}(E, \xi),
\]
which implies \eqref{eqn:remark:deg:plus:mu:max:01}.
\end{rema}

\begin{defi}\label{Def: tilde deg}
Let $\overline E$ be an adelic vector bundle on $S$ and $r$ be the rank of $E$ over $K$. We denote by $\widetilde{\deg}(\overline E)$ the sum $\widehat{\mu}_1(\overline E)+\cdots+\widehat{\mu}_r(\overline E)$. If $\overline E$ is the zero adelic vector bundle on  $S$, then by convention $\widetilde{\deg}(\overline E)$ is defined to be $0$. If $\overline E$ is non-zero, we define $\widetilde{\mu}(\overline E)$ to be the quotient $\widetilde{\deg}(\overline E)/\rang_K(E)$.
\end{defi}

{
\begin{prop}
Let $\overline E$ and $\overline F$ be non-zero adelic vector bundles on $S$ and $f:E\rightarrow F$ be a $K$-linear map. 
\begin{enumerate}[label=\rm(\arabic*)]
\item\label{Item: degree tilde E bounded by degree tilde F plus height f} Suppose that $f$ is a bijection. Then one has
\begin{equation}\label{Equ: comparison of sum of successive slopes}\widetilde{\deg}(\overline E)\leqslant\widetilde{\deg}(\overline F)+\rang_K(F)\cdot h(f).\end{equation}
\item\label{Item: slope inequality for non Hermitian case} Suppose that $f$ is injective. Then $\widehat{\mu}_1(\overline E)\leqslant\widehat{\mu}_1(\overline F)+h(f)$.
\end{enumerate}
\end{prop}
\begin{proof}
\ref{Item: degree tilde E bounded by degree tilde F plus height f} By Proposition \ref{Pro:fonctorialite general}, for any $t\in\mathbb R$ one has
\[f(\mathcal F_{\mathrm{hn}}^t(\overline E))\subseteq \mathcal F_{\mathrm{hn}}^{t-h(f)}(\overline F).\]
Therefore the inequality \eqref{Equ: comparison of sum of successive slopes} follows from Proposition \ref{Pro: comparaison des minima}.

\ref{Item: slope inequality for non Hermitian case} Let $\lambda=\widehat{\mu}_1(\overline E)$. Then $\mathcal F_{\mathrm{hn}}^\lambda(\overline E)\neq \{0\}$. Since $f$ is injective, by Proposition \ref{Pro:fonctorialite general}, this implies that $\mathcal F_{\mathrm{hn}}^{\lambda-h(f)}(\overline E)\neq\{0\}$ and hence $\lambda-h(f)\leqslant\widehat{\mu}_1(\overline F)$.
\end{proof}
}

\begin{prop}\label{Pro: tilde degree under subquotient}
Let $\overline E=(E,\xi)$ be an adelic vector bundle on $S$ and 
\[0=E_0\subseteq E_1\subseteq\ldots\subseteq E_n\]
be a flag of vector subspaces of $E$. One has
\begin{equation}\label{Equ: bound of deg E minus Delta E 01}\widehat{\deg}(\overline E)-\Delta(\overline E)\leqslant\sum_{i=1}^n\widetilde{\deg}(\overline {E_i/E_{i-1}})\leqslant\widehat{\deg}(\overline E)\end{equation}
If in addition $\xi$ is ultrametric on $\Omega\setminus\Omega_\infty$, one has 
\begin{equation}\label{Equ: bound of deg E minus Delta E 02}\widehat{\deg}(\overline E)-\delta(\overline E)\leqslant\sum_{i=1}^n\widetilde{\deg}(\overline {E_i/E_{i-1}})\leqslant\widehat{\deg}(\overline E)\end{equation}
\end{prop}
\begin{proof}
By Propositions \ref{Pro: degree bounded by integrals} and \ref{Pro: degree bounded above by integrals} (notably the inequality \eqref{Equ: degree bounded from above by sum of slopes}), for any $i\in\{1,\ldots,n\}$, one has
\[\widehat{\deg}(\overline{E_i/E_{i-1}})-\Delta(\overline{E_i/E_{i-1}})\leqslant\widetilde{\deg}(\overline{E_i/E_{i-1}})\leqslant\widehat{\deg}(\overline{E_i/E_{i-1}}).\]
Taking the sum with respect to $i$, by Proposition \ref{Pro:suiteexactedeg}, we obtain 
\[\widehat{\deg}(\overline E)-\Delta(\overline E)\leqslant\sum_{i=1}^n\widetilde{\deg}(\overline{E_i/E_{i-1}})\leqslant\widehat{\deg}(\overline E).\] 
In the case where $\xi$ is ultrametric on $\Omega\setminus\Omega_\infty$, the above argument combined with \eqref{Equ: degree bounded from above by sum of slopes delta} leads to \eqref{Equ: bound of deg E minus Delta E 02}.
\end{proof}

\begin{defi}
Let $\overline E$ be a non-zero adelic vector bundle on $S$. We say that $\overline E$ is \emph{semistable}\index{semistable}\index{adelic vector bundle!semistable} if its Harder-Narasimhan $\mathbb R$-filtration only have one jump point, namely one has $\widehat{\mu}_1(\overline E)=\cdots=\widehat{\mu}_r(\overline E)$ with $r=\rang_K(E)$. By definition, {the following conditions are equivalent:
\begin{enumerate}[label=\rm(\arabic*)]
\item $\overline E$ is semistable;
\item  for any non-zero vector subspace $F$ of $E$, one has $\widehat{\mu}_{\min}(\overline F)\leqslant\widehat{\mu}_{\min}(\overline E)$;
\item $\widetilde{\mu}(\overline E)=\widehat{\mu}_{\min}(\overline E)$.
\end{enumerate}}
\end{defi}

\begin{theo}\label{Thm: HN subquotient semistable}
Let $\overline E$ be a non-zero adelic vector bundle on $S$. We assume that the Harder-Narasimhan $\mathbb R$-filtration corresponds to the flag 
\begin{equation}\label{Equ: Harder-Narasimhan filtration general case}0=E_0\subsetneq E_1\subsetneq \ldots\subsetneq E_n=E\end{equation}
and the decreasing sequence $\mu_1>\ldots>\mu_n$ of real numbers. Then each subquotient $\overline{E_i/E_{i-1}}$ is semistable and $\widetilde{\mu}(\overline{E_i/E_{i-1}})=\mu_i$, $i\in\{1,\ldots,n\}$. Moreover, \eqref{Equ: Harder-Narasimhan filtration general case} is the only flag of vector subspaces of $E$ such that each subquotient $\overline{E_i/E_{i-1}}$ is semistable and 
\[\widetilde{\mu}(\overline{E_1/E_0})>\ldots>\widetilde{\mu}(\overline{E_n/E_{n-1}}).\]
\end{theo}
\begin{proof}We begin with showing that each subquotient $\overline{E_i/E_{i-1}}$ is semistable and that $\widehat{\mu}(\overline{E_i/E_{i-1}})=\mu_i$.
The case where $i=1$ results from the definition of Harder-Narasimhan $\mathbb R$-filtration. In what follows, we suppose that $i\geqslant 2$. 

By Proposition \ref{Pro: general hn filtration}, for any $j\in\{1,\ldots,n\}$ one has $\widehat{\mu}_{\min}(\overline E_j)\geqslant \mu_j$. Moreover, by definition of the Harder-Narasimhan $\mathbb R$-filtration, one has $\widehat{\mu}_{\min}(\overline E_j)\leqslant\mu_j$. Hence we obtain the equality  $\widehat{\mu}_{\min}(\overline E_j)=\mu_j$.

We claim that any vector subspace $G'$ of $E_i/E_{i-1}$ has a minimal slope $\leqslant\mu_i$. Let $\pi:E_i\rightarrow E_{i}/E_{i-1}$ be the canonical quotient map and $E_i'$ be the preimage of $G'$ by the quotient map $\pi$. Since $E_i'$ contains strictly $E_{i-1}$, one has $\widehat{\mu}_{\min}(\overline {E_i'})\leqslant\mu_i$. For any $\epsilon>0$ there exists a quotient vector space $H'$ of $E_i'$ such that {$\widehat{\mu}_{\max}(\overline{H}{}')\leqslant\mu_i+\epsilon$}. If $\epsilon<\mu_{i-1}-\mu_i$, then $\widehat{\mu}_{\min}(\overline{E}_{i-1})=\mu_{i-1}>\mu_i+\epsilon$. By Proposition \ref{Pro:inegalitdepente} {\ref{Item: slope inequality mu mix}}, we obtain that the composed map $E_{i-1}\rightarrow E_i'\rightarrow H'$ is zero, or equivalently, $H'$ is actually a quotient vector space of $E_i'/E_{i-1}=G'$. Hence we obtain $\widehat{\mu}_{\min}(\overline{G}{}')\leqslant\mu_i$. Therefore one has $\widehat{\mu}_1(\overline{E_i/E_{i-1}})\leqslant\mu_i\leqslant\mu_{\min}(\overline{E_i/E_{i-1}})$, which implies that $\overline {E_i/E_{i-1}}$ is semistable and $\widetilde{\mu}(
\overline{E_i/E_{i-1}})=\mu_i$.

We now proceed with the proof of the uniqueness by induction on the rank of $E$ over $K$. The case where $\rang_K(E)=1$ is trivial. In the following, we assume that the assertion has been proved for non-zero adelic vector bundles of rank $<\rang_K(E)$. We still denote by 
\[0=E_0\subsetneq E_1\subsetneq\ldots\subsetneq E_n=E\]
the Harder-Narasimhan flag of $\overline E$.  
Let 
\[0=F_0\subsetneq F_1\subsetneq\ldots\subsetneq F_m=E\]
be a flag of vector subspaces of $E$ such that each subquotient $\overline{F_j/F_{j-1}}$ is semistable and that
\[\widetilde{\mu}(\overline{F_1/F_0})>\ldots>\widetilde{\mu}(\overline{F_m/F_{m-1}}).\]
Since the subquotients $\overline{F_j/F_{j-1}}$ are semistable, we can rewrite these inequalities as
\begin{equation}\label{Equ: decreasing minimal slopes}\widehat{\mu}_{\min}(\overline{F_1/F_0})>\ldots>\widehat{\mu}_{\min}(\overline{F_m/F_{m-1}}).\end{equation}
We claim that $E_1$ is actually contained in $F_1$. Assume that $i$ is the smallest index in $\{1
,\ldots,m\}$ such that $E_1\subset F_i$. We identifie $E_1/(E_1\cap F_{i-1})$ with a vector subspace of $F_i/F_{i-1}$. Since $\overline{F_i/F_{i-1}}$ is semistable, one has
\[\widehat{\mu}_{\min}(\overline{E_1/(E_1\cap F_{i-1})})\leqslant\widehat{\mu}_{\min}(\overline{F_i/F_{i-1}}).\]
If $i>1$, then by \eqref{Equ: decreasing minimal slopes} one has $\widehat{\mu}_{\min}(\overline{F_i/F_{i-1}})<\widehat{\mu}_{\min}(\overline{F_1})\leqslant\widehat{\mu}_{\min}(\overline{E_1})$,
which leads to a contradiction since $\widehat{\mu}_{\min}(\overline{E_1/(E_1\cap F_{i-1})})\geqslant\widehat{\mu}_{\min}(\overline E_1)$. Therefore one has $E_1\subseteq F_1$. If the inclusion is strict, then by the definition of Harder-Narasimhan filtration  one has $\widehat{\mu}_{\min}(\overline{F_1})<\widehat{\mu}_{\min}(\overline{E_1})$. This contradicts the semi-stability of $\overline{F_1}$. Therefore we have $E_1=F_1$. Moreover, for any vector subspace $M$ of $E$ which contains strictly $E_1$, one has $\widehat{\mu}_{\min}(\overline M)<\widehat{\mu}_{\min}(\overline E_1)$. Hence, by Proposition \ref{Pro: exact sequence mu min} one has $\widehat{\mu}_{\min}(\overline M)=\widehat{\mu}_{\min}(\overline{M/E_1})$. Therefore, if $E/E_1$ is non-zero, then
\[0=E_1/E_1\subsetneq\ldots\subsetneq E_n/E_1=E/E_1\]
is the Harder-Narasimhan flag of $E/E_1$. By the induction hypothesis one has $n=m$ and $E_i=F_i$ for any $i\in\{2,\ldots,n\}$. The uniqueness is thus proved.
\end{proof}

\begin{prop}
Let $\overline E$ be a non-zero adelic vector bundle on $S$. The following assertions are equivalent:
\begin{enumerate}[label=\rm(\arabic*)]
\item $\overline E$ is semistable,
\item for any non-zero vector subspace $F$ of $E$, one has $\widetilde{\mu}(\overline F)\leqslant\widetilde{\mu}(\overline E)$
\item for any non-zero quotient vector space $G$ of $E$, one has $\widetilde{\mu}(\overline G)\geqslant\widetilde{\mu}(\overline E)$.
\end{enumerate}
\end{prop}
\begin{proof} Let $r$ be the rank of $E$ over $K$.
Assume that $\overline E$ is semistable. Then one has $\widehat{\mu}_1(\overline E)=\cdots=\widehat{\mu}_r(\overline E)=\widetilde{\mu}(\overline E)$. If $F$ is a non-zero vector subspace of $E$, then by Proposition \ref{Pro:fonctorialite general} we obtain that, for any $t\in\mathbb R$, one has $\mathcal F_{\mathrm{hn}}^t(\overline F)\subseteq\mathcal F_{\mathrm{hn}}^t(\overline E)$. Therefore, for any $i\in\{1,\ldots,\rang_K(F)\}$ one has $\widehat{\mu}_i(\overline F)\leqslant\widetilde{\mu}(\overline E)$, which implies $\widetilde{\mu}(\overline F)\leqslant\widetilde{\mu}(\overline E)$. Similarly, if $G$ is a non-zero quotient vector space of $E$ and $\pi:E\rightarrow G$ is the quotient map, then, by Proposition \ref{Pro:fonctorialite general}, one has $\pi(\mathcal F^{t}_{\mathrm{hn}}(\overline E))\subseteq\mathcal F^t_{\mathrm{hn}}(\overline G)$ for any $t\in\mathbb R$. Therefor for any $i\in\{1,\ldots,\rang_K(G)\}$ one has $\widehat{\mu}_i(\overline G)\geqslant\widetilde{\mu}(\overline E)$, which implies that $\widetilde{\mu}(\overline G)\geqslant\widetilde{\mu}(\overline E)$. Hence we have proved the implications  (1)$\Rightarrow$(2) and (1)$\Rightarrow$(3). 

We will prove the converse implications by contraposition. Suppose that $\overline E$ is not semistable and its Harder-Narasimhan $\mathbb R$-filtration corresponds to the flag
\[0=E_0\subsetneq E_1\subsetneq\ldots\subsetneq E_n=E,\]
and the successive jump points $\mu_1<\ldots<\mu_n$,
where $n\in\mathbb N$, $n\geqslant 1$. Then one has
\[\widetilde{\mu}(\overline E)=\frac{1}{\rang_K(E)}\sum_{i=1}^n\mu_i\rang_K(E_i/E_{i-1}).\]
By Theorem \ref{Thm: HN subquotient semistable} we obtain that $\widetilde{\mu}(\overline E_1)=\mu_1>\widetilde{\mu}(\overline E)$ and $\widetilde{\mu}(\overline{E_n/E_{n-1}})=\mu_n<\widetilde{\mu}(\overline E)$. The proposition is thus proved.
\end{proof}

\begin{rema}\label{Rem: adelic vector bundle on trivial valuation}
Consider the particular case where the adelic curve consists of exactly one copy of the trivial absolute value on $K$ (of measure $1$ with respect to $\nu$). In this case an adelic vector bundle on $S$ is just a finite-dimensional vector space $E$ over $K$ equipped with a norm $\norm{\ndot}$ (which is not necessary ultrametric), where we consider the trivial absolute value on $K$. We have shown in \S\ref{Subsection: Trivial valuation} that ultrametric norms on a finite-dimensional vector space over $K$ correspond bijectively to $\mathbb R$-filtrations on the same vector space. In particular, if $(E,\norm{\ndot})$ is a Hermitian adelic vector bundle on $S$, then the $\mathbb R$-filtration on $E$ corresponding to $\norm{\ndot}$ identifies with the Harder-Narasimhan $\mathbb R$-filtration of $(E,\norm{\ndot})$.
\end{rema}

{
\begin{prop}\label{Pro: condition of semistablity}
We equip $K$ with the trivial absolute value. Let $(E,\norm{\ndot})$ be a finite-dimensional normed vector space over $K$, which is also considered as an adelic vector bundle as in Remark \ref{Rem: adelic vector bundle on trivial valuation}. 
The adelic vector bundle $(E,\norm{\ndot})$ is semistable if and only if the double dual norm $\norm{\ndot}_{**}$ is constant on $E\setminus\{0\}$. {Moreover, in this case one has \[-\ln\norm{x}_{**}=\widehat{\mu}(E,\norm{\ndot})=\widehat{\mu}_{\min}(E,\norm{\ndot})\] for any $x\in E\setminus\{0\}$.}
\end{prop}
\begin{proof} 
First we assume that $(E,\norm{\ndot})$ is semistable.
 Let $\{e_i\}_{i=1}^r$ be an $\alpha$-orthogonal basis of $(E,\norm{\ndot})$, where $\alpha\in\intervalle{]}{0}{1}{[}$. Without loss of generality, we assume that $\norm{e_1}\leqslant\ldots\leqslant\norm{e_r}$. Moreover, by Proposition  \ref{Pro:orthogonalesthadamard} one has
\[\widehat{\deg}(E,\norm{\ndot})\leqslant -r\ln(\alpha)-\sum_{i=1}^r\ln\norm{e_i}.\]
In particular, if 
${\norm{e_r}}/{\norm{e_1}}>\alpha^{-r}$, that is, \[\displaystyle -\frac{1}{r} \ln \| e_1 \| > -\frac{1}{r} \ln \| e_r \| - \ln (\alpha),\] then
\begin{align*}
-\ln\norm{e_1} &= -\frac{r-1}{r}\ln\norm{e_1} - \frac{1}{r}\ln\norm{e_1} > -\frac{1}{r}\sum_{i=1}^{r-1} \ln\norm{e_i} -\frac{1}{r} \ln \| e_r \| - \ln (\alpha) \\
& \geqslant\widehat{\mu}(E,\norm{\ndot})\geqslant\widehat{\mu}_{\min}(E,\norm{\ndot}),
\end{align*}
which shows that $(E,\norm{\ndot})$ is not semistable, so that
${\norm{e_r}}/{\norm{e_1}}\leqslant \alpha^{-r}$. 
This observation
shows that, for any $\alpha$-orthogonal basis $\{e_i\}_{i=1}^r$ of $E$, one has
\begin{equation}\label{Equ: difference of ln norm ei eij}\max_{(i,j)\in\{1,\ldots,r\}^2}\Big|\ln\norm{e_i}-\ln\norm{e_j}\Big|\leqslant -r\ln(\alpha).\end{equation}
Note that $\{e_i\}_{i=1}^r$ is also an $\alpha$-orthogonal basis of $(E,\norm{\ndot}_{**})$ (see Proposition \ref{Pro:alphaorthogonale}). Moreover, we deduce from \eqref{Equ: difference of ln norm ei eij} and \eqref{Equ:majoration par doubledual} that 
\begin{equation}\label{Equ: estimate of distance of ei and ej}
\max_{(i,j)\in\{1,\ldots,r\}^2}\Big|\ln\norm{e_i}_{**}-\ln\norm{e_j}_{**}\Big|\leqslant -(r+1)\ln(\alpha).\end{equation}
Note that one has $\widehat{\deg}(E,\norm{\ndot})=\widehat{\deg}(E,\norm{\ndot}_{**})$ (see Proposition \ref{Pro:doubledualdet}). Moreover, by Propositions \ref{Pro:hadamard} and \ref{Pro:orthogonalesthadamard} one has
\[-\sum_{i=1}^r\ln\norm{e_i}_{**}\leqslant\widehat{\deg}(E,\norm{\ndot})\leqslant -r\ln(\alpha)-\sum_{i=1}^r\ln\norm{e_i}_{**}.\]
Combining this estimate with  \eqref{Equ: estimate of distance of ei and ej} we obtain
\[\max_{i\in\{1,\ldots,r\}}\Big|\ln\norm{e_i}_{**}-\widehat{\mu}(E,\norm{\ndot})\Big|\leqslant-(r+2)\ln(\alpha).\]
In particular, for any $(\lambda_1,\ldots,\lambda_r)\in K^r\setminus\{(0,\ldots,0)\}$, one has
\[\Big|\ln\norm{\lambda_1e_1+\cdots+\lambda_re_r}_{**}-\widehat{\mu}(E,\norm{\ndot})\Big|\leqslant -(r+3)\ln(\alpha)\] Since $(E,\norm{\ndot})$ admits an $\alpha$-orthogonal basis for any $\alpha\in\intervalle{]}{0}{1}{[}$ (see Corollary \ref{Cor: existence of alpha orthogonal}), we obtain that the restriction of {$\ln\norm{\ndot}_{**}$} {to} $E\setminus\{0\}$ is constant (which is equal to {$-\operatorname{\widehat{\mu}}(E,\norm{\ndot})$}).
 
Assume now that the double dual norm $\norm{\ndot}_{**}$ is constant on $E\setminus\{0\}$. Since $\norm{\ndot}_{**}$ and $\norm{\ndot}$ induce the same dual norm on $E^\vee$, we obtain that the restriction {of $\ln\norm{\ndot}_{*}$ to} $E^\vee\setminus\{0\}$ is constant and takes {$-\operatorname{\widehat{\mu}}(E^\vee,\norm{\ndot}_*)$} as its value. Note that one has {$-\operatorname{\widehat{\mu}}(E^\vee,\norm{\ndot}_*)=\widehat{\mu}(E,\norm{\ndot})$} by Proposition \ref{Pro:degdual}. We will show that $(E,\norm{\ndot})$ is semistable. First we show that $\widehat{\mu}_{\min}(E,\norm{\ndot})=\widehat{\mu}(E,\norm{\ndot})$. Let $G$ be a non-zero quotient vector space of $E$ and $\norm{\ndot}_G$ be the quotient norm of $\norm{\ndot}$ on $G$. By Proposition \ref{Pro:dualquotient}, $\norm{\ndot}_{G,*}$ coincides with the restriction of $\norm{\ndot}_*$ {to} $G^\vee$. Since the {function $\ln\norm{\ndot}_*$} takes constant value {$\widehat{\mu}(E,\norm{\ndot})$} on $E^\vee\setminus\{0\}$ we obtain that \[\widehat{\mu}(G,\norm{\ndot}_G)=-\operatorname{\widehat{\mu}}(G^\vee,\norm{\ndot}_{G,*})=\widehat{\mu}(E,\norm{\ndot}).\] Therefore $\widehat{\mu}_{\min}(E,\norm{\ndot})=\widehat{\mu}(E,\norm{\ndot})$. Now for any non-zero vector subspace $F$ of $E$ one has \[\widehat{\mu}_{\min}(F,\norm{\ndot}_F)\leqslant\widehat{\mu}(F,\norm{\ndot}_F)\leqslant\widehat{\mu}(E,\norm{\ndot}_E),\] where $\norm{\ndot}_F$ denotes the restriction of $\norm{\ndot}$ {to} $F$. In fact, {$\ln\norm{\ndot}_F$} is bounded from below by the restriction of {$\ln\norm{\ndot}_{**}$} {to} $F$, which is constant on $F\setminus\{0\}$ of value {$-\operatorname{\widehat{\mu}}(E,\norm{\ndot})$}. Therefore $(E,\norm{\ndot})$ is semistable.
\end{proof}

{

\begin{rema}\label{Rem: cretiereno on semistabilit}
We keep the notation of the previous proposition. Note that the norms $\norm{\ndot}$ and $\norm{\ndot}_{**}$ induce the same dual norm on $E^\vee$ (see Proposition \ref{Pro:doubedualandquotient} \ref{Item: seminorm and double dual induce the same dual norm}), so that we obtain that the adelic vector bundle $(E,\norm{\ndot})$ is semistable if and only if the restriction of the function $\norm{\ndot}_*$ on $E^\vee\setminus\{0\}$ is constant. Moreover, in this case one has (see Proposition \ref{Pro:delta=1})
\[\forall\,\varphi\in E^\vee\setminus\{0\},\quad -\ln\norm{\varphi}_*=-\operatorname{\widehat{\deg}}(E,\norm{\ndot}).\]
\end{rema}

\begin{rema}
In the case where $\norm{\ndot}$ is ultrmetric, the normed vector space $(E,\norm{\ndot})$ corresponds to a sequence 
\[0=E_0\subsetneq E_1\subsetneq\ldots\subsetneq E_n=E\]
of vector subspaces of $E$ and a decreasing sequence $\mu_1>\ldots>\mu_n$ of real numbers (see Remark \ref{Rem: R-filtration as flag plus slopes}). Note that, for any $i\in\{1,\ldots, n\}$ the restriction of the subquotient norm $\norm{\ndot}_{E_i/E_{i-1}}$ to $(E_i/E_{i-1})\setminus\{0\}$ is constant and takes $\mathrm{e}^{-\mu_i}$ as its value. Therefore Proposition \ref{Pro: condition of semistablity} implies that $(E_i/E_{i-1},\norm{\ndot}_{E_i/E_{i-1}})$ is semistable and admits $\mu_i$ as its minimal slope. Therefore Theorem \ref{Thm: HN subquotient semistable} shows that 
\[0=E_0\subsetneq E_1\subsetneq\ldots\subsetneq E_n=E\]
is the Harder-Narasimhan flag of the adelic vector bundle $(E,\norm{\ndot})$.
\end{rema}

\begin{prop}
We equip $K$ with the trivial absolute value. Consider a finite-dimensional non-zero normed vector space $(E,\norm{\ndot})$ over $K$. The Harder-Narasimhan flags of $(E,\norm{\ndot})$ and $(E,\norm{\ndot}_{**})$ are the same. Moreover, for any $i\in\{1,\ldots,\rang_K(E)\}$ one has $\widehat{\mu}_i(E,\norm{\ndot})=\widehat{\mu}_i(E,\norm{\ndot}_{**})$.
\end{prop}
\begin{proof}
Let $n$ be the rank of $E$ over $K$. We reason by induction on $n$. First of all, if $(E,\norm{\ndot})$ is semistable, then by Proposition \ref{Pro: condition of semistablity} (see also its proof), the function $-\ln\norm{\ndot}_{**}$ is constant on $E\setminus\{0\}$ and takes $\widehat{\mu}(E,\norm{\ndot})=\widehat{\mu}_{\min}(\norm{\ndot})$ as its value. Therefore the assertion of the proposition holds in this case, and in particular the assertion is true when $n=1$. In the following we suppose that $(E,\norm{\ndot})$ is not semistable (hence $n\geqslant 2$) and that the proposition has been proved for normed vector spaces of dimension $\leqslant n-1$. 

For any $i\in\{1,\ldots,n\}$, let $\mu_i=\widehat{\mu}_i(E,\norm{\ndot}_{**})$ and $E_i$ be the ball of radius $\mathrm{e}^{-\mu_i}$ in $(E,\norm{\ndot}_{**})$ centered at the origin.  Let
\[\beta=\min\{\mu_i-\mu_{i-1}\,:\,i\in\{2,\ldots,n\},\;\mu_i>\mu_{i-1}\}\]
and $\alpha$ be an element of $\mathopen{]}0,1\mathclose{[}$ such that $\alpha>\mathrm{e}^{-\beta/n}$. Let $\{e_i\}_{i=1}^n$ be an $\alpha$-orthogonal bases of $(E,\norm{\ndot})$. By Proposition \ref{Pro:alphaorthogonale}, it is also an $\alpha$-orthgonal basis of $(E,\norm{\ndot}_{**})$. By Proposition \ref{Pro: alpha orth is orth}, $\{e_i\}_{i=1}^n$ is an orthogonal basis of $(E,\norm{\ndot}_{**})$. Without loss of generality, we may assume that \[\{e_i\}_{i=1}^n\cap E_1=\{e_1,\ldots,e_m\},\]
where $m$ is the rank of $E_1$ over $K$ (see Proposition \ref{Pro: alpha orth is orth} \ref{Item: criterion of orthgonal basis}). Let $\{e_i^\vee\}_{i=1}^n$ be the dual basis of $\{e_i\}_{i=1}^n$ and $\norm{\ndot}_1$ be the restriction of the norm $\norm{\ndot}$ to $E_1$. For any $i\in\{1,\ldots,m\}$, let $\varphi_i$ be the restriction of $e_i^\vee$ on $E_1$. Then $\{\varphi_i\}_{i=1}^m$ forms a basis of $E_1^\vee$, which is the dual basis of $\{e_i\}_{i=1}^m$. By lemma \ref{Lem:normofdualbasis}, $\{\varphi_i\}_{i=1}^m$ is an $\alpha$-orthogonal basis of $E_1^\vee$, and one has
\[\forall\,i\in\{1,\ldots,m\},\quad \norm{e_i}^{-1}\leqslant\norm{\varphi_i}_{1,*}\leqslant\alpha^{-1}\norm{e_i}^{-1}.\] 
By Proposition \ref{Pro:alphaorthogonale}, one has
\[\forall\,i\in\{1,\ldots,m\},\quad \alpha\norm{e_i}\leqslant \norm{e_i}_{**}=\mathrm{e}^{-\mu_1}\leqslant\norm{e_i}.\]
Therefore one obtains
\[\forall\,i\in\{1,\ldots,m\},\quad \alpha\mathrm{e}^{\mu_1}\leqslant\norm{\varphi_i}_{1,*}\leqslant\alpha^{-1}\mathrm{e}^{\mu_1}.\]
As a consequence, for a general non-zero element $\varphi$ of $E_1^\vee$, which is written in the form $\lambda_1\varphi_1+\cdots+\lambda_m\varphi_m$, one has
\[\norm{\varphi}_{1,*}\leqslant\max_{\begin{subarray}{c}
i\in\{1,\ldots,m\}\\
\lambda_i\neq 0
\end{subarray}}\norm{\varphi_i}_{1,*}\leqslant\alpha^{-1}\mathrm{e}^{\mu_1}\]
and
\[\norm{\varphi}_{1,*}\geqslant \alpha\max_{\begin{subarray}{c}
i\in\{1,\ldots,m\}\\
\lambda_i\neq 0
\end{subarray}}\norm{\varphi_i}_{1,*}\geqslant\alpha^2\mathrm{e}^{\mu_1}.\]
Since $(E,\norm{\ndot})$ admits an $\alpha$-orthogonal basis for any $\alpha\in\mathopen{]}0,1\mathclose{[}$, we obtain that the restriction of $\norm{\ndot}_{1,*}$ on $E_1\setminus\{0\}$ is constant and takes $\mathrm{e}^{\mu_1}$ as its value. Therefore, Proposition \ref{Pro: condition of semistablity} (see also Remark \ref{Rem: cretiereno on semistabilit}), we obtain that $(E_1,\norm{\ndot})$ is semistable and admits $\mu_1$ as its minimal slope.

By Proposition \ref{Pro:dualquotient}, one has (see Definition \ref{Def:restriction} and Subsection \ref{Subsec:Quotientnorm} for notation)
\[\norm{\ndot}_{*,(E/E_1)^\vee\hookrightarrow E^\vee}=\norm{\ndot}_{E\twoheadrightarrow E/E_1,*}.\]
Moreover, since $\norm{\ndot}_{*}$ is ultrametric, by Proposition \ref{Pro:quotientdualnonarch} one has
\[\norm{\ndot}_{**,E\twoheadrightarrow E/E_1}=\norm{\ndot}_{*,(E/E_1)^\vee\hookrightarrow E^\vee,*}=\norm{\ndot}_{E\twoheadrightarrow E/E_1,**}.\]
Applying the induction hypothesis to $(E/E_1,\norm{\ndot}_{E\twoheadrightarrow E/E_1})$ we obtain that the Harder-Narasimhan flags and the successive slopes of $(E/E_1,\norm{\ndot}_{E\twoheadrightarrow E/E_1})$  and $(E/E_1,\norm{\ndot}_{**,E\twoheadrightarrow E/E_1})$ are the same. Therefore, by Theorem \ref{Thm: HN subquotient semistable} we obtain that the Harder-Narasimhan flag and the successive slopes of $(E,\norm{\ndot})$ coincides with those of $(E,\norm{\ndot}_{**})$. The proposition is thus proved.
\end{proof}
}

\begin{rema}
In the framework of linear code, Randriambololona \cite{Randriam17} has proposed a Harder-Narasimhan theory based on semimodular degree functions on the modular lattice of vector subspaces. Note that our approach, which relies on the Arakelov degree function of quotient vector spaces, has a very different nature from the classic method (due to the fact that the equality \eqref{Equ:additivity degree} and the inequality \eqref{Equ:convexitededeg} fail in general for non-Hermitian adelic vector bundles). It is an intriguing question to compare the Harder-Narasimhan filtrations constructed in our setting and in \cite{Randriam17}.
\end{rema}
}

{
\subsection{Absolute positive degree and absolute maximal slope}

We have seen in Proposition \ref{Pro: Arakelov degree preserved by extension of scalars} that the Arakelov degree is preserved by extension of scalars. In this subsection, we discuss the behaviour of the maximal slope and the positive degree under extension of scalars to the algebraic closure of $K$. We denote by $K^{\mathrm{ac}}$ the algebraic closure of the field $K$.

\begin{defi}
Let $(E,\xi)$ be an adelic vector bundle on $S$. We denote by $\sposdeg(E,\xi)$ the positive degree of $(E_{K^{\mathrm{ac}}},\xi_{K^{\mathrm{ac}}})$, called the \emph{absolute positive degree}\index{absolute positive degree} of $(E,\xi)$. If $E$ is non-zero, we denote by $\smaxslop(E,\xi)$ the maximal slope of $(E_{K^{\mathrm{ac}}},\xi_{K^{\mathrm{ac}}})$, called the \emph{absolute maximal slope}\index{absolute maximal slope} of $(E,\xi)$. 
\end{defi}

\begin{prop}\label{Pro: comparaison of degree plus and mu max}
Let $(E,\xi)$ be an adelic vector bundle on $S$.
One has \[\widehat{\deg}_+(E,\xi)\leqslant\sposdeg(E,\xi)\text{ and }\widehat{\mu}_{\max}(E,\xi)\leqslant\smaxslop(E,\xi).\] Moreover, for any algebraic extension $L$ of $K$, one has \[\sposdeg(E,\xi)=\sposdeg(E_L,\xi_L)\text{
and }\smaxslop(E,\xi)=\smaxslop(E_L,\xi_L).\]
\end{prop}
\begin{proof}
Let $F$ be a vector subspace of $E$ and $\xi_F$ be the restriction of $\xi$ {to} $F$. Let $\xi_{F_{K^{\mathrm{ac}}}}$ be the restriction of $\xi_{K^{\mathrm{ac}}}$ {to} $F_{K^{\mathrm{ac}}}$.
By Proposition \ref{Pro: restriction of epsion tensor} \ref{Item: restriction epsilon tensor}, \ref{Item: pi restriction extension}, the identity map $(F_{K^{\mathrm{ac}}},\xi_{F,{K^{\mathrm{ac}}}})\rightarrow (F_{K^{\mathrm{ac}}},\xi_{F_{K^{\mathrm{ac}}}})$ has norm $\leqslant 1$ on any $\omega\in\Omega$. By Proposition \ref{Pro:slopeinequality1}, one has
\[\widehat{\deg}(F,\xi_F)=\widehat{\deg}(F_{K^{\mathrm{ac}}},\xi_{F,{K^{\mathrm{ac}}}})\leqslant\widehat{\deg}(F_{K^{\mathrm{ac}}},\xi_{F_{K^{\mathrm{ac}}}})\leqslant \widehat{\deg}_+(E_{K^{\mathrm{ac}}},\xi_{K^{\mathrm{ac}}})=\sposdeg(E,\xi),\]
where the first equality comes from Proposition \ref{Pro: Arakelov degree preserved by extension of scalars}. Similarly, if $F$ is non-zero, one has
\[\widehat{\mu}(F,\xi_F)\leqslant\widehat{\mu}(F_{K^{\mathrm{ac}}},\xi_{F,K^{\mathrm{ac}}})\leqslant\widehat{\mu}(F_{K^{\mathrm{ac}}},\xi_{F_{K^{\mathrm{ac}}}})\leqslant\widehat{\mu}_{\max}(F_{K^{\mathrm{ac}}},\xi_{F_{K^{\mathrm{ac}}}})=\smaxslop(F,\xi_F).\] 
Since $F$ is arbitrary, we obtain \[\widehat{\deg}_+(E,\xi)\leqslant\sposdeg(E,\xi)\quad\text{ and }\quad\widehat{\mu}_{\max}(E,\xi)\leqslant\smaxslop(E,\xi).\]

By Corollary \ref{Cor:extsucc}, if $L$ is an algebraic extension of $K$, then one has $(\xi_L)_{K^{\mathrm{ac}}}=\xi_{K^{\mathrm{ac}}}$. Therefore $\sposdeg(E,\xi)=\sposdeg(E_L,\xi_L)$, and $\smaxslop(E,\xi)=\smaxslop(E_L,\xi_L)$.
\end{proof}

\begin{prop}\label{Pro: stable by extension of scalars mu max}
Assume that the field $K$ is perfect. Let $(E,\xi)$ be a \emph{Hermitian} adelic vector bundle on $S$. Then one has 
\[\widehat{\deg}_+(E,\xi)=\sposdeg(E,\xi)\quad\text{ and }\quad\widehat{\mu}_{\max}(E,\xi)=\smaxslop(E,\xi).\]
\end{prop} 
\begin{proof}
Without loss of generality, we may assume that the vector space $E$ is non-zero. Let
\[\{0\}= \widetilde E_0\subsetneq \widetilde E_1\subsetneq\ldots\subsetneq \widetilde E_n=E_{K^{\mathrm{a}}}\]
be the Harder-Narasimhan flag of $(E_{K^{\mathrm{ac}}},\xi_{K^{\mathrm{ac}}})$. By the uniqueness of Harder-Narasimhan filtration (see Proposition \ref{Pro:characterisationdehnfiltration}), for any $K$-automorphism $\tau$ of $K^{\mathrm{ac}}$ and any $i\in\{1,\ldots,n\}$, the vector space $\widetilde E_i$ is stable by $\tau$. Since the filed $K$ is perfect, by Galois descent (see \cite{Bourbaki_A4-7}, Chapter V, \S10, no.4, Corollary of Proposition 6), there exists a flag
\begin{equation}\label{Equ: sequence e0 etc}\{0\}=E_0\subsetneq E_1\subsetneq\ldots\subsetneq E_n=E\end{equation}
such that $\widetilde E_i=E_{i,K^{\mathrm{a}}}$ for any $i\in\{1,\ldots,n\}$. Moreover, by Propositions \ref{Pro: restriction of epsion tensor} \ref{Item: eps restriction extension}, \ref{Item: pi restriction extension} (here we use the hypothesis that $\xi$ is Hermitian), if we denote by $\xi_i$ the restriction of $\xi$ {to} $E_i$, the $\xi_{i,K^{\mathrm{ac}}}$ coincides with the restriction $\widetilde{\xi}_i$ of $\xi_{K^{\mathrm{ac}}}$ {to} $\widetilde E_i$. Therefore by Proposition \ref{Pro: Arakelov degree preserved by extension of scalars} one has $\widehat{\deg}(E_i,\xi_i)=\widehat{\deg}(\widetilde E_i,\widetilde{\xi}_i)$. We then deduce that the slopes of $E_i/E_{i-1}$ and $\widetilde E_i/\widetilde E_{i-1}$ (equipped with subquotient norm families) are the same. Hence by Proposition \ref{Pro:characterisationdehnfiltration} we obtain that \eqref{Equ: sequence e0 etc} is the Harder-Narasimhan flag of $(E,\xi)$. Therefore,
\[\widehat{\mu}_{\max}(E,\xi)=\widehat{\mu}(E_1,\xi_1)=\widehat{\mu}(\widetilde E_1,\widetilde \xi_1)=\smaxslop(E,\xi)\]
and
\[\widehat{\deg}_+(E,\xi)=\max_{i\in\{0,\ldots,n\}}\widehat{\deg}(E_i,\xi_i)=\max_{i\in\{0,\ldots,n\}}\widehat{\deg}(\widetilde E_i,\widetilde \xi_i)=\sposdeg(E,\xi).\]
\end{proof}

\subsection{Successive minima}

The successive minima are classic invariants of Hermitian vector bundles on an arithmetic curve. In this subsection, we extend their construction (more precisely, the construction of successive minima of Roy-Thunder \cite{Roy_Thunder96}) to the setting of adelic vector bundles on an adelic curve.

\begin{defi}
Let $(E,\xi)$ be an adelic vector bundle on $S$ and $r$ be the rank of $E$ over $K$.
For any $i\in\{1,\ldots,r\}$, let 
\[\nu_i(E,\xi):=\sup\{t\in\mathbb R\,:\,\rang_K(\mathrm{Vect}_K(\{s\in E_K\,:\,\widehat{\deg}_\xi(s)\geqslant t\}))\geqslant i\},\]
called the \emph{$i^{\text{th}}$ (logarithmic) minimum}\index{ith logarithmic minimum@$i$-th logarithmic minimum} of $(E,\xi)$.
In other words, $\nu_i(E,\xi)$ is the supremum of the set of real numbers $t$ such that there exist at list $i$ linearly independent vectors of Arakelov degree $\geqslant t$. Clearly one has
\[\nu_1(E,\xi)\geqslant\ldots\geqslant\nu_r(E,\xi).\]
The first minimum $\nu_1(E,\xi)$ is also denoted by $\nu_{\max}(E,\xi)$, and the last minimum $\nu_r(E,\xi)$ is also denoted by $\nu_{\min}(E,\xi)$.  For any $t\in\mathbb R$, let
\[{\mathcal F_{\mathrm{m}}^{t}(E,\xi):=\bigcap_{\varepsilon>0}\mathrm{Vect}_K(\{s\in E\setminus\{0\}\,:\,\widehat{\deg}_\xi(s)\geqslant t-\varepsilon\}).}\] By definition \[\nu_i(E,\xi)=\sup\{t\in\mathbb R\,:\,\rang_K(\mathcal F_{\mathrm{m}}^t(\overline E))\geqslant i\}.\]
If $E$ is the zero vector space, then by convention we define \[\nu_{\max}(E,\xi):=-\infty\quad\text{ and }\quad\nu_{\min}( E,\xi):=+\infty.\]

We also define the absolute version of the successive minima as follows. For any $i\in\{1,\ldots,r\}$, let
$\nu_i^{\mathrm{a}}(E,\xi):=\nu_i(E_{K^{\mathrm{a}}},\xi_{K^{\mathrm{a}}})$, where $K^{\mathrm{a}}$ denotes the algebraic closure of $(E,\xi)$. Similarly, we let \[\nu_{\max}^{\mathrm{a}}(E,\xi):=\nu_{\max}(E_{K^{\mathrm{a}}},\xi_{K^{\mathrm{a}}})\quad\text{ and }\quad\nu_{\min}^{\mathrm{a}}(E,\xi):=\nu_{\min}(E_{K^{\mathrm{a}}},\xi_{K^{\mathrm{a}}}).\]
\end{defi}

\begin{prop}
Let $\overline E=(E,\xi)$ be a non-zero adelic vector bundle on $S$. For any $t\in\mathbb R$ one has
\begin{equation}\label{Equ: filtration by minima}{\mathcal F_{\mathrm{m}}^t(\overline E)=\bigcap_{\varepsilon>0}\sum_{\begin{subarray}{c}
0\neq F\subseteq E\\
\nu_{\min}(\overline F)\geqslant t-\varepsilon
\end{subarray}}F,}\end{equation}
where $F$ runs over the set of all non-zero vector subspaces of $E$, and in the structure of adelic vector bundle of $\overline F$ we consider the restricted norm family.
\end{prop}
\begin{proof}
Let $\varepsilon>0$ and  $F$ be a non-zero vector subspace of $E$ such that $\nu_{\min}(\overline F)>t-\varepsilon$. There exists a basis $\{s_i\}_{i=1}^n$ of $F$ over $K$ such that 
\[\min_{i\in\{1,\ldots,n\}}\widehat{\deg}_\xi(s_i)\geqslant t-2\varepsilon.\]
Therefore one has
\[\sum_{\begin{subarray}{c}0\neq F\subseteq E\\
\widehat{\mu}_{\min}(F)\geqslant t-\varepsilon\end{subarray}}F\subseteq \operatorname{Vect}_K(\{s\in E\setminus\{0\}\,:\,\widehat{\deg}_{\xi}(s)\geqslant t-2\varepsilon\}).\]

Conversely, for any $t\in\mathbb R$ such that $\mathcal F_{\mathrm{m}}^t(\overline E)\neq\{0\}$ {and any $\varepsilon>0$}, there exist  elements  $u_1,\ldots,u_r$ in $E$ which generates $\mathcal F_{\mathrm{m}}^t(\overline E)$ as vector space over $K$ and such that 
\[\forall\,j\in\{1,\ldots,r\},\quad \widehat{\deg}_{\xi}(u_j)\geqslant {t-\varepsilon}.\]
Therefore $\nu_{\min}(\mathcal F_{\mathrm{m}}^t(\overline E))\geqslant t$.
\end{proof}

\begin{prop}\label{Pro: comparison of nu under extension of fields}
Let $(E,\xi)$ be an adelic vector bundle on $S$ and $r$ be the rank of $E$ over $K$. For any $i\in\{1,\ldots,r\}$, one has $\nu_i(E,\xi)\leqslant\nu_i^{\mathrm{a}}(E,\xi)$. Moreover, for any algebraic extension $L$ of $K$ and any $i\in\{1,\ldots,r\}$ one has $\nu_{i}^{\mathrm{a}}(E,\xi)=\nu_i^{\mathrm{a}}(E_L,\xi_L)$.
\end{prop}
\begin{proof}
Let $\{s_j\}_{j=1}^i$ be a linearly independent family in $E$. Then it is also a linearly independent family in $E_{K^{\mathrm{a}}}$. Moreover, by the same argument as in the proof of Proposition \ref{Pro: comparaison of degree plus and mu max}, for any $j\in\{1,\ldots,i\}$ one has $\widehat{\deg}_{\xi}(s_j)\leqslant\widehat{\deg}_{\xi_{K^{\mathrm{a}}}}(s_j)$. Therefore $\nu_i(E,\xi)\leqslant\nu_i(E_{K^{\mathrm{a}}},\xi_{K^{\mathrm{a}}})=\nu_i^{\mathrm{a}}(E,\xi)$. The equality $\nu_i^{\mathrm{a}}(E_L,\xi_L)=\nu_i^{\mathrm{a}}(E,\xi)$ comes from the relation $(\xi_L)_{K^{\mathrm{a}}}=\xi_{K^{\mathrm{a}}}$, which is a consequence of Corollary \ref{Cor:extsucc}.
\end{proof}

The following proposition is straightforward from the definition of the (absolute) fist minimum and the (absolute) maximal slope.

\begin{prop}\label{Pro: first minimum bounded from above by maximal slope}
If $(E,\xi)$ is an adelic vector bundle on $S$, then one has 
\begin{equation}\label{Equ: udeg bounded by mu max}
\nu_1(E,\xi)\leqslant\widehat{\mu}_{\max}(E,\xi)\quad\text{ and }\quad \nu_1^{\mathrm{a}}(E,\xi)\leqslant\smaxslop(E,\xi).
\end{equation}
\end{prop}

\subsection{Minkowski property}

\begin{defi}\label{Def: Minkowski property}
Let $S=(K,(\Omega,\mathcal A,\nu),\phi)$ be an adelic curve. Let $C$ be a non-negative real number. We say that the adelic curve $S$ satisfies the \emph{Minkowski property of level $\geqslant C$}\index{Minkowski property} if, for any adelic vector bundle $(E,\xi)$ on $S$ such that $\xi$ is ultrametric on $\Omega\setminus\Omega_\infty$, one has
\[\nu_1(E,\xi)\geqslant\widehat{\mu}_{\max}(E,\xi)-C\ln(\rang_K(E)).\]
We say that $S$ satisfies the \emph{absolute Minkowski property of level $\geqslant C$}\index{absolute Minkowski property} if, for any adelic vector bundle $(E,\xi)$ on $S$, one has
\[\nu_1^{\mathrm{a}}(E,\xi)\geqslant\smaxslop(E,\xi)- C\ln(\rang_K(E)).\] 
\end{defi} 

\begin{rema}
Let $\overline V$ be an Euclidean lattice. The first theorem of Minkowski can be stated as (see \cite[\S3.2]{Bost_Kunnemann} for more details)
\[\nu_1(\overline V)\geqslant\widehat{\mu}_{\max}(\overline V)-\frac 12\ln(\rang(V)).\]
Hence the Minkowski property should be considered as an analogue in the general setting of adelic curve of the statement of the first theorem of Minkowski. For general number fields, it has been shown in \cite[\S5]{MR3604914} that, for any adelic vector bundle $\overline E$ of rank $n$ over  a number field $K$, one has
\[\nu_1(\overline E)\geqslant\widehat{\mu}_{\max}(\overline E)-\frac 12\ln(n)-\frac 12\ln|\mathfrak D_{K/\mathbb Q}| \]
and
\[\nu_1^{\mathrm{a}}(\overline E)\geqslant\widehat{\mu}_{\max}(\overline E)-\frac 12\sum_{\ell=2}^n\frac{1}{\ell}=\widehat{\mu}_{\max}^{\mathrm{a}}(\overline E)-\frac 12\sum_{\ell=2}^n\frac{1}{\ell}\]
where $\mathfrak D_{K/\mathbb Q}$ is the discriminant of $K$ over $\mathbb Q$. Therefore the adelic curve corresponding to a number field satisfies the Minkowski property of level $\geqslant \frac 12+\frac 12\ln|\mathfrak D_{K/\mathbb Q}|$ and the absolute Minkowski property of level $\geqslant \frac 12$.

In the function field case, given a regular projective curve (over a base field), by Riemann-Roch formula there exists a constant $A>0$ which only depends on the  curve, such that 
\[\nu_1(E)\geqslant\mu_{\max}(E)-A\]
for any vector bundle $E$ on the curve (see \cite[Remark 8.3]{Chen15}). Therefore the Minkowski property of level $\geqslant A/\ln(2)$ is satisfied in this case. Moreover, if the base field is of characteristic zero, then it has been shown in \cite{Bost_Chen} that the absolute Minkowski property of level $\geqslant 0$ is satisfied.

The Minkowski property may fail for general adelic curves. Consider the adelic curve $S=(\mathbb Q,(\mathbb Q,\mathcal A,\nu),\phi)$ consisting of the field of rational numbers, the measure space of $\mathbb Q$ equipped with the discrete $\sigma$-algebra and the atomic measure such that $\nu(\{\omega\})=1$ for any $\omega\in\mathbb Q$, together with the map $\phi$ sending any $\omega\in\mathbb Q$ to the trivial absolute value on $\mathbb Q$. We write the rational numbers into a sequence $\{q_n\}_{n\in\mathbb N}$. For any $n\in\mathbb N_{\geqslant 2}$, consider the following adelic vector bundle $\overline E_{n}$ on $S$. Let $E_{n}=K^2$. For $m\in\mathbb N$ such that $m< n$ let $\norm{\ndot}_{q_m}$ be the norm on $E_{n}$ defined as
\[\norm{(x,y)}_{q_m}=\begin{cases}
\mathrm{e}^{-1},&\text{if there exists $a\in K^{\times}$ such that $(x,y)=a(1,q_m)$},\\
0,&\text{if $(x,y)=(0,0)$,}\\
1,&\text{else}.
\end{cases}\]
For $m\in\mathbb N$ such that $m\geqslant n$, let $\norm{\ndot}_{q_m}$ be the norm on $K^2$ such that $\norm{(x,y)}=1$ for any $(x,y)\in K^2\setminus\{(0,0)\}$.
Then by definition one has $\widehat{\deg}(\overline E_n)=n$ and hence $\widehat{\mu}(\overline E_n)=n/2$. Moreover, for any vector subspace $F$ of $E$, either there exists $m\in\{0,\ldots,n-1\}$ such that $F=K(1,q_m)$ and thus $\widehat{\mu}(\overline F)=1$, or one has $\widehat{\mu}(\overline F)=0$. Since $n\geqslant 2$, we obtain that the adelic vector bundle $\overline E_n$ is semistable and of slope $n/2$. Moreover the first minimum of $\overline E_n$ is $1$. Therefore it is not possible to find a constant $C$ only depending on $S$ such that $\widehat{\mu}_{\max}(\overline E_n)$ is bounded from above by $\nu_1(\overline E_n)+C\ln(2)$.
\end{rema}

In the literature, the (absolute) Minkowski property is closely related to the semistability of tensor vector bundles and the estimation of the maximal slope of them. We refer the readers to \cite{Andre11,Bost_Chen,Gaudron_Remond13} for more detailed discussions. In the following, we prove several slope estimates in assuming the Minkowski property.

\begin{prop}
Let $\overline E=(E,\xi_E)$ and $\overline F=(F,\xi_F)$ be adelic vector bundles on $S$. We assume that $\xi_E$ and $\xi_F$ are ultrametric on $\Omega\setminus\Omega_\infty$. One has
\begin{equation}\label{Equ: nu 1 tensor bounded}\nu_1(\overline E\otimes_{\varepsilon}\overline F)\leqslant\widehat{\mu}_{\max}(\overline E)+\widehat{\mu}_{\max}(\overline F).\end{equation}
\end{prop}
\begin{proof}
Let $f$ be a non-zero element of $E\otimes_KF$, viewed as a $K$-linear map from $E^\vee$ to $F$. By Proposition \ref{Pro:inegalitdepente} \ref{Item: slope inequality mu mix}, one has
\[\widehat{\mu}_{\min}(\overline E{}^\vee)\leqslant\widehat{\mu}_{\max}(\overline F)+h(f)=\widehat{\mu}_{\max}(\overline F)-\widehat{\deg}_{\xi_E\otimes_{\varepsilon}\xi_F}(f).\]
By Proposition \ref{Pro: mu min plus mu max dual}, we obtain
\[0\leqslant\widehat{\mu}_{\max}(\overline E)+\widehat{\mu}_{\max}(\overline F)-\widehat{\deg}_{\xi_E\otimes_{\varepsilon}\xi_F}(f).\]
Since $f$ is arbitrary, we obtain the inequality \ref{Equ: nu 1 tensor bounded}.
\end{proof}

\begin{coro}Let $C$ be a non-negative real number.
We assume that the adelic curve $S$ satisfies the Minkowski property of level $\geqslant C$. {Let $\overline E=(E,\xi_E)$ and $\overline F=(F,\xi_F)$ be adelic vector bundles on $S$.
\begin{enumerate}[label=\rm(\arabic*)]
\item Assume that $\xi_E$ and $\xi_F$ are ultrametric on $\Omega\setminus\Omega_\infty$. Then
\begin{equation}
\label{Equ: upper bound mu max tensor}\widehat{\mu}_{\max}(\overline E\otimes_{\varepsilon}\overline F)\leqslant \widehat{\mu}_{\max}(\overline E)+\widehat{\mu}_{\max}(\overline F)+C\ln(\rang_K(E)\cdot\rang_K(F)),
\end{equation}
\item One has
\begin{equation}
\label{Equ: lower bound mu min tensor}
\widehat{\mu}_{\min}(\overline E\otimes_{\varepsilon,\pi}\overline F)\geqslant\widehat{\mu}_{\min}(\overline E)+\widehat{\mu}_{\min}(\overline F)-(C+{\textstyle\frac 12}\nu(\Omega_\infty))\ln(\rang_K(E)\cdot\rang_K(F)).
\end{equation}
\end{enumerate}}
\end{coro}
\begin{proof}
By the assumption of Minkowski property, we have
\[\nu_1(\overline E\otimes_{\varepsilon}\overline F)\geqslant\widehat{\mu}_{\max}(\overline E\otimes_{\varepsilon}\overline F)-C\ln(\rang_K(E)\cdot\rang_K(F)).\]
Hence \eqref{Equ: upper bound mu max tensor} follows from \eqref{Equ: nu 1 tensor bounded}.

If we apply the inequality \eqref{Equ: upper bound mu max tensor} to $\overline E{}^\vee$ and $\overline F{}^\vee$ (note that $\xi_E^\vee$ and $\xi_F^\vee$ are always ultrametric on $\Omega\setminus\Omega_\infty$), we obtain
\[\widehat{\mu}_{\max}(\overline E{}^\vee\otimes_{\varepsilon}\overline F{}^\vee)\leqslant\widehat{\mu}_{\max}(\overline E{}^\vee)+\widehat{\mu}_{\max}(\overline F{}^\vee)+C\ln(\rang_K(E)\cdot\rang_K(F)).\]
By Proposition \ref{Pro: mu min plus mu max dual} we deduce that
\[\widehat{\mu}_{\min}(\overline E\otimes_{\varepsilon,\pi}\overline F)\geqslant-\widehat{\mu}_{\max}(\overline E{}^\vee)-\widehat{\mu}_{\max}(\overline F{}^\vee)-C\ln(\rang_K(E)\cdot\rang_K(F)).\]
Finally, by Corollary \ref{Cor: majoration of mu min plus mu max} we obtain \eqref{Equ: lower bound mu min tensor}.
\end{proof}

\begin{prop}\label{Pro: comparaison mu min et nu min}
Let $\overline E$ be a non-zero adelic vector bundle on $S$. One has
\[\nu_{\min}(\overline E)\leqslant\widehat{\mu}_{\min}(\overline E).\]
\end{prop}
\begin{proof}
Let $r$ be the rank of $E$ over $K$. Assume that $\overline E$ is of the form $\overline E=(E,\xi)$, with $\xi=\{\|\ndot\|_{\omega}\}_{\omega\in\Omega}$. Let $t$ be a real number and $\{s_i\}_{i=1}^r$ be a basis of $E$ over $K$ such that $\widehat{\deg}_\xi(s_i)\geqslant t$ for any $i\in\{1,\ldots,r\}$. 

Let $G$ be a quotient vector space of $E$ and $\xi_G=\{\|\ndot\|_{G,\omega}\}_{\omega\in\Omega}$ be the quotient norm family of $\xi$ on $G$. For any $i\in\{1,\ldots,r\}$, let $\alpha_i$ be the canonical image of $s_i$ in $G$. Without loss of generality, we may assume that $\{\alpha_1,\ldots,\alpha_n\}$ form a basis of $G$ over $K$. For any $\omega\in\Omega$, one has
\[\|\alpha_1\wedge\cdots\wedge \alpha_n\|_{G,\omega,\det}\leqslant\prod_{i=1}^n\|\alpha_i\|_{G,\omega}\leqslant\prod_{i=1}^n\|s_i\|_\omega,\]
where the first inequality comes from Proposition \ref{Pro:Hadamard's inequality}. Therefore one has
\[\widehat{\deg}(\overline G)\geqslant\sum_{i=1}^n\widehat{\deg}_{\xi}(s_i)\geqslant nt,\]
which implies $\widehat{\mu}(\overline G)\geqslant t$. Therefore we obtain $\nu_{\min}(\overline E)\leqslant\widehat{\mu}_{\min}(\overline E)$.
\end{proof}

\begin{coro}\label{Cor: comparison between mu i and nu i}
Let $\overline E=(E,\xi)$ be a non-zero adelic vector bundle on $S$. For any $i\in\{1,\ldots,\rang_K(E)\}$ one has $\nu_{i}(\overline E)\leqslant\widehat{\mu}_i(\overline E)$.
\end{coro}
\begin{proof}By the relations \eqref{Equ: filtration by minima} and \eqref{Equ: Harder-Narasimhan filtration in general},  Proposition \ref{Pro: comparaison mu min et nu min} leads to $\mathcal F_{\mathrm{m}}^t(\overline E)\subseteq \mathcal F^t_{\mathrm{hn}}(\overline E)$ for any $t\in\mathbb R$. Therefore, by Proposition \ref{Pro: comparaison des minima} we obtain that $\nu_i(\overline E)\leqslant\widehat{\mu}_i(\overline E)$ for any $i\in\{1,\ldots,\rang_K(E)\}$.
\end{proof}

\begin{defi}
Let $S=(K,(\Omega,\mathcal A,\nu),\phi)$ be an adelic curve. Let $C$ be a non-negative real number. We say that $S$ satisfies the \emph{strong Minkowski property of level $\geqslant C$}\index{strong Minkowski property} if for any adelic vector bundle $(E,\xi)$ on $S$ such that $\xi$ is ultrametric on $\Omega\setminus\Omega_\infty$ one has 
\begin{equation}\label{Equ: majoration de mu min par ny min}\nu_{\min}(E,\xi)\geqslant\widehat{\mu}_{\min}(E,\xi)-C\ln(\rang_K(E)).\end{equation}\end{defi}

\begin{prop}\label{Pro: strong Minkowski}
Assume that the adelic curve $S$ satisfies the strong Minkowski property of  level $\geqslant C$. For any non-zero adelic vector bundle $\overline E=(E,\xi)$ on $S$ one has
\begin{equation}\label{Equ: bound of successive slope by successive minima}\widehat{\mu}_i(\overline E)\leqslant \nu_i(\overline E)+C\ln(\rang_K(E)).\end{equation}
\end{prop}
\begin{proof}
Since the adelic curve $S$ satisfies the strong Minkowski property of  level $\geqslant C$, by the relation \eqref{Equ: majoration de mu min par ny min} we obtain that
\[\forall\,t\in\mathbb R,\quad\mathcal F_{\mathrm{hn}}^t(\overline E)\subseteq\mathcal F_{\mathrm{m}}^{t-C\ln(\rang_K(E))}(\overline E).\]
Therefore, by Proposition \ref{Pro: comparaison des minima} we obtain the inequality \eqref{Equ: bound of successive slope by successive minima}.
\end{proof}

\begin{rema}\label{Rem: number field is strong Minkowski}
Proposition \ref{Pro: strong Minkowski} shows that, if the adelic curve $S$ satisfies the strong Minkowski property of level $\geqslant C$, then it also satisfies Minkowski property of level $\geqslant C$. Moreover, the transference theorem of Gaudron \cite[Theorem 36]{Gaudron18} shows that, for any \emph{Hermitian} adelic vector bundle $\overline E$ of rank $n$ over a number field $K$, one has 
\[\nu_{\min}(\overline E)-\widehat{\mu}_{\min}(\overline E)=\nu_{\min}(\overline E)+\widehat{\mu}_{\max}(\overline E{}^\vee)\geqslant \nu_{\min}(\overline E)+\widehat{\nu}_{\max}(\overline E{}^\vee)\geqslant\ln(n)+\ln|\mathfrak D_{K/\mathbb Q}|,\]
where $\mathfrak D_{K/\mathbb Q}$ is the discriminant of $K$ over $\mathbb Q$. We then deduce that, if $\overline E$ is a general adelic vector bundle of rank $n$ over $K$, which is not necessarily Hermitian, one has (by Theorem \ref{Thm: Hermitian approximation via measurable selection})
\[\nu_{\min}(\overline E)\geqslant\widehat{\mu}_{\min}(\overline E)-\Big(1+\frac 12[K:\mathbb Q]\Big)\ln(n)+\ln|\mathfrak D_{K/\mathbb Q}|.\]
Therefore the adelic curve corresponding to a number field $K$ satisfies the strong Minkowski property of level $1+\frac{1}{2}[K:\mathbb Q]\ln|\mathfrak D_{K/\mathbb Q}|$.
\end{rema}

\section{Adelic vector bundles over number fields}
Throughout this section, we fix a number field $K$ and the standard adelic curve $S =(K,(\Omega,\mathcal A,\nu),\phi)$ of $K$ as in Subsection~\ref{Subsec:Numberfields}.
Note that $S$ is proper.
Denote by $\Omega_{\mathrm{fin}}$ the set $\Omega\setminus\Omega_\infty$ of finite places of $K$, and by $\mathfrak o_K$ the ring of algebraic integers  in $K$. 
Note that the absolute value $|\ndot|_{\omega}$ at $\omega$ is given by
\[
\forall\, x \in K_{\omega},\quad
|x|_{\omega} = \begin{cases}
\text{the standard absolute value of $x$ in either $\mathbb R$ or $\mathbb C$} & \text{if $\omega \in \Omega_{\infty}$}, \\[1ex]
{\displaystyle \exp\left( \frac{-\log p_{\omega}\ord_{\omega}(x)}{\ord_{\omega}(p_{\omega})} \right)} & \text{if $\omega \in \Omega_{\mathrm{fin}}$},
\end{cases}
\]
where 
$p_{\omega}$ is the characteristic of the residue
field of the valuation ring of
$K_{\omega}$. 
Further, for $\omega \in \Omega_{\mathrm{fin}}$, let $\mathfrak o_{K,\omega}$ be the localisation of $\mathfrak o_K$ at $\omega$
and $\mathfrak o_{\omega}$ be the valuation ring of 
the completion $K_{\omega}$ of $K$ with respect to $\omega$, that is,
\[
\mathfrak o_{K,\omega} = \{ a \in K \mid |a|_{\omega} \leqslant 1 \}\quad\text{and}\quad
\mathfrak o_{\omega} = \{ a \in K_{\omega} \mid |a|_{\omega} \leqslant 1 \}.
\]
Moreover,
$\nu(\{\omega\}) = [K_{\omega}: \mathbb Q_{\omega}]$ for $\omega \in \Omega$ and
$\sum_{\omega \in \Omega_{\infty}} \nu (\{ \omega \}) = [K : \mathbb Q]$.

\medskip
Let $E$ be a finite-dimensional vector space over $K$ and
$\xi = \{ \|\ndot\|_{\omega} \}_{\omega \in \Omega}$ be a norm family of $E$ over $S$.
In this section, we always assume that $\|\ndot\|_{\omega}$ is ultrametric for every $\omega \in \Omega_{\mathrm{fin}}$. 
For $\omega \in \Omega_{\mathrm{fin}}$, we set
\[
E_{\omega} := E \otimes_K K_{\omega}\quad\text{and}\quad
\mathscr E_{\omega} := \{ x \in E_{\omega} \mid \| x \|_{\omega} \leqslant 1 \}.
\]
By Proposition~\ref{prop:ultrametric:ball:lattice} and \ref{Pro:normetreausauxdisc},
$\mathscr E_{\omega}$ is a free $\mathfrak o_{\omega}$-module of rank $\dim_K E$ and $\mathscr E_{\omega} \otimes_{\mathfrak o_{\omega}} K_{\omega} = E_{\omega}$.

\begin{rema}\label{remark:local:E:complete:E}
As in the next subsection, let $(E,\xi)^{\omega}_{\leqslant 1} := \{ x \in E \mid \| x \|_{\omega} \leqslant 1 \}$. Then one can see the following:
\begin{enumerate}[label=\rm(\arabic*)]
\item
$(E,\xi)^{\omega}_{\leqslant 1}$ is a free $\mathfrak o_{K,\omega}$-module.

\item
$(E,\xi)^{\omega}_{\leqslant 1} \otimes_{\mathfrak o_{K,\omega}} \mathfrak o_{\omega} =
\mathscr E_{\omega}$.
\end{enumerate}

(1) Fix a basis of $(x_i)_{i=1}^r$ of $E$.
We consider a norm $\|\ndot\|'_{\omega}$ on $E_{\omega}$ given by
\[
\forall\, \lambda_1, \ldots, \lambda_r \in K_{\omega},\quad
\| \lambda_1 x_1 + \cdots + \lambda_r x_r \|'_{\omega} = \max \{ |\lambda_1|_{\omega}, \ldots, |\lambda_r|_{\omega} \}.
\]
By Proposition~\ref{Pro:topologicalnormedspace}, there is a positive integer $n$ such that $|\varpi_{\omega}|^n \|\ndot\|'_{\omega}
\leqslant \|\ndot\|_{\omega}$, where $\varpi_{\omega}$ is a uniformizing parameter of
$\mathfrak o_{K, \omega}$. Therefore,
\[
(E,\xi)^{\omega}_{\leqslant 1} \subseteq 
\mathfrak o_{K,\omega} e_1 \varpi_{\omega}^{-n} + \cdots + \mathfrak o_{K,\omega} e_r \varpi_{\omega}^{-n},
\]
as required.

(2)
Obviously $(E, \xi)^{\omega}_{\leqslant 1} \otimes_{\mathfrak o_{K, \omega}} \mathfrak o_{\omega}\subseteq \mathscr E_{\omega}$.
Let $(e_i)_{i=1}^r$ be
a free basis of $(E, \xi)^{\omega}_{\leqslant 1}$ over $\mathfrak o_{K,\omega}$.
For $x \in \mathscr E_{\omega}$, 
we choose $a_1, \ldots, a_r \in K_{\omega}$ such that
$x = a_1 e_1 + \cdots + a_r e_r$.
One can find $a'_1, \ldots, a'_r \in K$ such that
\[
| a_i - a'_i |_{\omega} \leqslant \frac{1}{2} |a_i|_{\omega}\ (\forall\, i)\quad\text{and}\quad
\| x - (a'_1 e_1 + \cdots + a'_r e_r) \|_{\omega} \leqslant \frac{1}{2} \| x \|_{\omega}.
\]
If we set $x' = a'_1 e_1 + \cdots + a'_r e_r$, then
$\|x \|_{\omega} = \|x'\|_{\omega}$ and $|a_i|_{\omega} = |a'_i|_{\omega}$ for all $i$.
In particular $x' \in (E, \xi)^{\omega}_{\leqslant 1}$, so that $a'_i \in \mathfrak o_{K,\omega}$, and hence $|a_i|_{\omega} = |a'_i|_{\omega}\leqslant 1$.
Therefore, $x \in (E, \xi)^{\omega}_{\leqslant 1} \otimes_{\mathfrak o_{K, \omega}} \mathfrak o_{\omega}$.
\end{rema}

\subsection{Coherency for a norm family}
Let $E$ be a finite-dimensional vector space over $K$.
Let $\xi = \{ \|\ndot\|_{\omega} \}_{\omega \in \Omega}$ be a norm family of $E$ over $S$. 
We define $(E, \xi)^{\mathrm{fin}}_{\leqslant 1}$ and
$(E,\xi)^{\omega}_{\leqslant 1}$ ($\omega \in \Omega_{\mathrm{fin}}$) to be
\[
\begin{cases}
(E, \xi)^{\mathrm{fin}}_{\leqslant 1} := \{ x \in E \mid \| \text{$x \|_{\omega} \leqslant 1$ for all $\omega \in \Omega_{\mathrm{fin}}$} \}, \\
(E,\xi)^{\omega}_{\leqslant 1} := \{ x \in E \mid \| x \|_{\omega} \leqslant 1 \}.
\end{cases}
\]
Note that $(E, \xi)^{\mathrm{fin}}_{\leqslant 1}$ and $(E,\xi)^{\omega}_{\leqslant 1}$ are an $\mathfrak o_K$-module and
an $\mathfrak o_{K,\omega}$-module, respectively.
Furthermore, by Remark~\ref{remark:local:E:complete:E}, $(E,\xi)^{\omega}_{\leqslant 1}$ is
a free $\mathfrak o_{K,\omega}$-module and $(E,\xi)^{\omega}_{\leqslant 1} \otimes_{\mathfrak o_{K,\omega}} \mathfrak o_{\omega} = \mathscr E_{\omega}$.
Let us begin with the following proposition.

\begin{prop}\label{prop:finiteness:cond:classical:setting}
The following are equivalent:
\begin{enumerate}[label=\rm(\arabic*)]
\item
For any $v \in E$, $\| v \|_{\omega} \leqslant 1$ except finitely many $\omega \in \Omega_{\mathrm{fin}}$.

\item
$(E,\xi)^{\mathrm{fin}}_{\leqslant 1} \otimes_{\mathfrak o_K} \mathfrak o_{K,\omega} = (E,\xi)^{\omega}_{\leqslant 1}$
for all $\omega \in \Omega_{\mathrm{fin}}$.

\item
$(E, \xi)^{\mathrm{fin}}_{\leqslant 1} \otimes_{\mathfrak o_K} \mathfrak o_{K,\omega} = (E,\xi)^{\omega}_{\leqslant 1}$
for some $\omega \in \Omega_{\mathrm{fin}}$.

\item
$(E, \xi)^{\mathrm{fin}}_{\leqslant 1} \otimes_{\mathfrak o_K} K = E$.
\end{enumerate}
Moreover, under the above equivalent conditions,
$(E, \xi)^{\mathrm{fin}}_{\leqslant 1} \otimes_{\mathbb Z} \mathbb Q = E$.
\end{prop}

\begin{proof}
First of all, let us see the following claim:

\begin{enonce}{Claim}\label{claim:prop:finiteness:cond:classical:setting:01}
\begin{enumerate}[label=\rm(\alph*)]
\item
Let $S$ be a finite subset of $\Omega_{\mathrm{fin}}$. Then there is $f \in \mathfrak o_K \setminus \{ 0 \}$ such that
\[
\ord_{\omega}(f) \begin{cases}
> 0 & \text{if $\omega \in S$}, \\
= 0 & \text{if $\omega \not\in S$}.
\end{cases}
\]

\item
$(E, \xi)^{\omega}_{\leqslant 1} \otimes_{\mathfrak o_{K,\omega}} K = E$ for all $\omega \in \Omega_{\mathrm{fin}}$.
\end{enumerate}
\end{enonce}

\begin{proof}
(a) Let us consider the ideal given by $\mathfrak I = \prod_{\mathfrak p \in S} \mathfrak p$.
As the class group of $K$ is finite,
there are a positive integer $a$ and $f \in \mathfrak o_K$
such that $f \mathfrak o_K = \mathfrak I^{a}$, as required.

\medskip
(b) Obviously $(E, \xi)^{\omega}_{\leqslant 1} \otimes_{\mathfrak o_{K,\omega}} K \subseteq E$. For $v \in E$, there is
$a \in \mathfrak o_{K,\omega}\setminus \{0\}$ such that $av \in (E, \xi)^{\omega}_{\leqslant 1}$, which shows the converse
inclusion.
\end{proof}

(1) $\Longrightarrow$ (2): 
Clearly $(E, \xi)^{\mathrm{fin}}_{\leqslant 1} \otimes_{\mathfrak o_K} \mathfrak o_{K,\omega} \subseteq (E, \xi)^{\omega}_{\leqslant 1}$.
Conversely, for $v \in (E, \xi)^{\omega}_{\leqslant 1}$, as $S = \{ \omega' \in \Omega_{\mathrm{fin}} \mid \| v \|_{\omega'} > 1 \}$ is finite,
there is $f \in \mathfrak o_K \setminus \{ 0 \}$ such that $| f |_{\omega'} < 1$ for $\omega' \in S$ and
$| f |_{\omega'} = 1$ for $\omega' \in \Omega_{\mathrm{fin}} \setminus S$ by the above claim (a). 
Thus, there is a positive integer $n$ such that
$f^{n} v \in (E, \xi)^{\mathrm{fin}}_{\leqslant 1}$. Note that $f \in \mathfrak o_{K,\omega}^{\times}$. Thus the converse inclusion holds.

\medskip
``(2) $\Longrightarrow$ (3)'' is obvious and ``(3) $\Longrightarrow$ (4)'' follows from (b) in the claim.
Let us see that ``(4) $\Longrightarrow$ (1)''.
For $v \in E$, there is $a \in \mathfrak o_K \setminus \{ 0 \}$ such that $a v \in (E, \xi)^{\mathrm{fin}}_{\leqslant 1}$, that is,
$|a|_{\omega} \| v \|_{\omega} \leq 1$ for all $\omega \in \Omega_{\mathrm{fin}}$. Note that $|a|_{\omega} = 1$ except finitely many $\omega$,
so that one has (1).

\medskip
Note that $(E, \xi)^{\mathrm{fin}}_{\leqslant 1} \otimes_{\mathfrak o_K} K$ and $(E,\xi)^{\mathrm{fin}}_{\leqslant 1} \otimes_{\mathbb Z} \mathbb Q$
are the localizations of $(E, \|\ndot\|)^{\mathrm{fin}}_{\leqslant 1}$ with respect to $\mathfrak o_K \setminus \{0\}$ and $\mathbb Z \setminus \{ 0\}$, respectively.
Therefore, for the last assertion, it is sufficient to show that, for $\alpha \in \mathfrak o_K \setminus \{0\}$, there is
$\alpha' \in \mathfrak o_K \setminus \{0\}$ such that $\alpha \alpha' \in \mathbb Z \setminus \{ 0 \}$.
Indeed, one can find $a_1, \ldots, a_n \in \mathbb Z$ such that $\alpha^n + a_1 \alpha^{n-1} + \cdots + a_{n-1} \alpha  + a_n =0$.
We may assume that $a_n \not= 0$. Thus
$\alpha (\alpha^{n-1} + a_1 \alpha^{n-2} + \cdots + a_{n-1}) = -a_n \in \mathbb Z \setminus \{0\}$.
\end{proof}

\begin{defi}\label{def:finiteness:cond:classical:setting}
We say that $(E, \xi)$ is \emph{coherent}\index{coherent} if the equivalent conditions of Proposition~\ref{prop:finiteness:cond:classical:setting}
are satisfied.
\end{defi}

\begin{prop}\label{prop:strongly:cohenent}
If there are a non-empty open set $U$ of $\Spec(\mathfrak o_K)$ and
a locally free $\mathfrak o_U$-module $\mathscr E$ such that
$\mathscr E \otimes_{\mathfrak o_U} K = E$ and
$\|\ndot\|_{\omega} = \|\ndot\|_{\mathscr E \otimes_{\mathfrak o_U} \mathfrak o_{\omega}}$ for all $\omega \in U \cap \Omega_{\mathrm{fin}}$,
then $(E, \xi)$ is coherent and dominated,
where $\mathfrak o_U$ is the ring of regular functions on the open set $U$.
\end{prop}

\begin{proof}
For $s \in E \setminus \{ 0 \}$, we can find a non-empty open set $U' \subseteq U$
such that $s \in \mathscr E \otimes_{\mathfrak o_U} \mathfrak o_{\omega}$ and $\mathscr E \otimes_{\mathfrak o_U} \mathfrak o_{\omega}/\mathfrak o_{\omega} s$ is torsion free
for all $\omega \in U' \cap \Omega_{\mathrm{fin}}$, so that
\[
\| s \|_{\omega} = \| s\|_{\mathscr E \otimes_{\mathfrak o_U} \mathfrak o_{\omega}} = 1.
\]
In particular, $(E,\xi)$ is upper-dominated and coherent.
Let $\mathscr E^{\vee}$ be the dual of $\mathscr E$ over $U$.
Note that $\mathscr E^{\vee} \otimes_{\mathfrak o_U} \mathfrak o_{\omega} = (\mathscr E \otimes_{\mathfrak o_U} \mathfrak o_{\omega})^{\vee}$, so that
by Propsotion~\ref{coro:dual:lattice:norm},
$\|\ndot\|_{\omega,*} = \|\ndot\|_{\mathscr E^{\vee} \otimes_{\mathfrak o_U} \mathfrak o_{\omega}}$ for all $\omega \in U \cap \Omega_{\mathrm{fin}}$.
Therefore, in the same way as above, one can see that
$(E^{\vee},\xi^{\vee})$ is upper-dominated, and hence $(E, \xi)$ is dominated.
\end{proof}

\subsection{Finite generation of a dominated vector bundle over $S$}
Let $E$ be a finite-dimensional vector space over $K$ and $\xi = \{ \|\ndot\|_{\omega} \}_{\omega \in \Omega}$ be a norm family of $E$ over $S$.
The purpose of this subsection is to prove the following theorem.

\begin{prop}\label{prop:finite:generation:dominated}
If $\xi$ is dominated, then
\[
(E, \xi)^{\mathrm{fin}}_{\leqslant 1} := \{ x \in E \,:\, \text{$\| x \|_{\omega} \leqslant 1$ for all $\omega \in \Omega_{\mathrm{fin}}$} \}
\]
is a finitely generated $\mathfrak o_{K}$-module.
\end{prop}

\begin{proof}
First we assume that $\dim_K E = 1$. Fix $x \in E \setminus \{ 0 \}$.  
For each $\omega \in \Omega_{\mathrm{fin}}$, 
let $a_{\omega}$ be the smallest integer $a$ with $a \geqslant - \ln \|x\|_{\omega}/\ln |\varpi_{\omega}|_{\omega}$,
where $\varpi_{\omega}$ is a local parameter of $\mathfrak o_{K,\omega}$.
 
As $\xi$ is lower dominated, there is an integrable function $A(\omega)$ on $\Omega$
such that 
\[
\forall\, \omega \in \Omega,\quad
- \ln \| x \|_{\omega} \leqslant A(\omega).
\]

Here we assume that there are infinitely many $\omega \in \Omega_{\mathrm{fin}}$ with
$a_{\omega} \leqslant -1$. 
As 
$a_{\omega} \leqslant -1$ implies
$-\ln |\varpi_{\omega}|_{\omega} \leqslant - \ln \|x\|_{\omega}$, one has
\begin{align*}
A(\omega) \nu(\{ \omega \}) & \geqslant -\ln \| x \|_{\omega}\nu(\{ \omega \}) \\
& \geqslant -\ln |\varpi_{\omega}|_{\omega} \nu(\{ \omega \})
= \ln \#(\mathfrak o_K/\mathfrak p_{\omega}) \geqslant \ln p_{\omega},
\end{align*}
where $\mathfrak p_\omega$ is the maximal ideal of $\mathfrak o_K$ and $p_\omega$ is the characteristic of the residue field $\mathfrak o_K/\mathfrak p_\omega$, which gives a contradiction to the integrability of the function $A(\ndot)$. Therefore,
$a_{\omega} \geqslant 0$ except finitely many $\omega \in \Omega_{\mathrm{fin}}$.

Note that
\begin{align*}
\| a x \|_{\omega} = |a|_{\omega} \| x \|_{\omega} \leqslant 1 & \Longleftrightarrow\ |\varpi_{\omega}|_{\omega}^{\ord_{\omega}(a) + \ln \|x\|_{\omega}/\ln |\varpi_{\omega}|_{\omega}} \leqslant 1 \\
& \Longleftrightarrow\ \ord_{\omega}(a) \geqslant -\ln \| x\|_{\omega}/\ln |\varpi_{\omega}|_{\omega} \\
& \Longleftrightarrow\ \ord_{\omega}(a) \geqslant a_{\omega}.
\end{align*}
Therefore 
\[
\{ a \in K \,:\, \text{$\| a x \|_{\omega} \leqslant 1$ for all $\omega \in \Omega_{\mathrm{fin}}$} \} 
= \{ a \in K \,:\, \text{$\ord_{\omega}(a) - a_{\omega} \geqslant 0$ for all $\omega \in \Omega_{\mathrm{fin}}$} \}
\]
is finitely generated over $\mathfrak o_K$ by Lemma~\ref{lem:finite:gen:number:field}. Thus one has the assertion in the case where $\dim_K E = 1$.

\medskip
In general, we prove the theorem by induction on $\dim_K E$.
By the previous observation, we may assume $\dim_K E \geqslant 2$.
Fix $x \in E \setminus \{ 0 \}$. We set $E' = Kx$ and $E'' = E/E'$.
Let $\xi'$ be the norm family on $E'$ given by the restriction of $\xi$, and
$\xi''$ be the norm family on $E''$ given by the quotient of $\xi$. 
Then $\xi'$ and $\xi''$ are dominated by Proposition~\ref{Pro:dominancealgebraic}, so that,
by the hypothesis of induction, $(E', \xi')^{\mathrm{fin}}_{\leqslant 1}$ and
$(E'', \xi'')^{\mathrm{fin}}_{\leqslant 1}$ are finitely generated over $\mathfrak o_K$.
Note that $\beta((E, \xi)^{\mathrm{fin}}_{\leqslant 1}) \subseteq (E'', \xi'')^{\mathrm{fin}}_{\leqslant 1}$,
where $\beta$ is the canonical homomorphism $E \to E''$.
In particular, 
$\beta((E, \xi)^{\mathrm{fin}}_{\leqslant 1})$ is finitely generated over $\mathfrak o_K$ because $\mathfrak o_K$ is Noetherian.
Therefore, one has the exact sequence 
\[
0 \longrightarrow (E', \xi')^{\mathrm{fin}}_{\leqslant 1} \longrightarrow (E, \xi)^{\mathrm{fin}}_{\leqslant 1} \longrightarrow \beta((E, \xi)^{\mathrm{fin}}_{\leqslant 1}) \longrightarrow 0,
\] 
and hence the assertion follows.
\end{proof}

\begin{lemm}\label{lem:finite:gen:number:field}
Let $\{ b_{\omega} \}_{\omega \in \Omega_{\mathrm{fin}}}$ be a family of integers indexed by $\Omega_{\mathrm{fin}}$.
Then
\[
\mathfrak o_K(\{ b_{\omega} \}_{\omega \in \Omega_{\mathrm{fin}}}) := \{ a \in K \,:\, \text{$\ord_{\omega}(a) + b_{\omega} \geqslant 0$ for all $\omega \in \Omega_{\mathrm{fin}}$} \}
\]
is finitely generated over $\mathfrak o_K$ if and only if either $b_{\omega} \leqslant 0$ except finitely many $\omega$, or
$b_{\omega} < 0$ for infinitely many $\omega$.
\end{lemm}

\begin{proof}
We set $S = \{ \omega \in \Omega_{\mathrm{fin}} \,:\, b_{\omega} \geqslant 1 \}$ and 
$T = \{ \omega \in \Omega_{\mathrm{fin}} \,:\, b_{\omega} \leqslant -1 \}$

First we assume that $b_{\omega} \leqslant 0$ except finitely many $\omega$, that is, $\#(S) < \infty$. Then
one can choose $f \in \mathfrak o_K \setminus \{ 0 \}$ such that
$\ord_{\omega}(f) \geqslant b_{\omega}$ for all $\omega \in \Omega_{\mathrm{fin}}$.
Note that $\{ a \in K \,:\, \text{$\ord_{\omega}(a) \geqslant 0$ for all $\omega \in \Omega_{\mathrm{fin}}$} \} = \mathfrak o_K$, so that
$\mathfrak o_K(\{ b_{\omega} \}_{\omega \in \Omega_{\mathrm{fin}}}) f \subseteq \mathfrak o_K$.
Thus $\mathfrak o_K(\{ b_{\omega} \}_{\omega \in \Omega_{\mathrm{fin}}}) f$ is finitely generated over $\mathfrak o_K$ because
$\mathfrak o_K$ is Noetherian. Therefore $\mathfrak o_K(\{ b_{\omega} \}_{\omega \in \Omega_{\mathrm{fin}}})$ 
is also finitely generated over $\mathfrak o_K$.

Next we assume that $b_{\omega} < 0$ for infinitely many $\omega$, that is, $\#(T) = \infty$.
In this case, $\mathfrak o_K(\{ b_{\omega} \}_{\omega \in \Omega_{\mathrm{fin}}}) = \{ 0 \}$.
Indeed, if $x \in \mathfrak o_K(\{ b_{\omega} \}_{\omega \in \Omega_{\mathrm{fin}}}) \setminus \{ 0 \}$,
then $\ord_{\omega}(x) \geq 1$ for all $\omega \in T$, which is a contradiction.

Finally we assume that $\#(S) = \infty$ and $\#(T) < \infty$. In this case, we need to show that
$\mathfrak o_K(\{ b_{\omega} \}_{\omega \in \Omega_{\mathrm{fin}}})$ is not finitely generated over $\mathfrak o_K$.
We set $S = \{ \omega_1, \omega_2, \ldots, \omega_n, \ldots \}$.
For each positive integer $N$, let us consider a family $\{ b_{N, \omega} \}_{\omega \in \Omega_{\mathrm{fin}}}$ of integers given by
\[
b_{N, \omega} = \begin{cases}
0 & \text{if $\omega \in \{\omega_n \,:\, n \geqslant N+1 \}$},\\
b_{\omega} & \text{otherwise}.
\end{cases}
\]
Then one has a strictly increasing sequence of finitely generated $\mathfrak o_K$-modules:
\[
\mathfrak o_K(\{ b_{1, \omega} \}_{\omega \in \Omega_{\mathrm{fin}}}) \subsetneq 
\mathfrak o_K(\{ b_{2, \omega} \}_{\omega \in \Omega_{\mathrm{fin}}}) \subsetneq \cdots \subsetneq
\mathfrak o_K(\{ b_{N, \omega} \}_{\omega \in \Omega_{\mathrm{fin}}}) \subsetneq \cdots
\]
such that $\bigcup_{N=1}^{\infty} \mathfrak o_K(\{ b_{N, \omega} \}_{\omega \in \Omega_{\mathrm{fin}}}) = \mathfrak o_K(\{ b_{\omega} \}_{\omega \in \Omega_{\mathrm{fin}}})$.
Therefore $\mathfrak o_{K}(\{ b_{\omega} \}_{\omega \in \Omega_{\mathrm{fin}}})$ is not finitely generated over $\mathfrak o_K$.
\end{proof}

\begin{enonce}[remark]{Example}
Let $\{ b_{\omega} \}_{\omega \in \Omega_{\mathrm{fin}}}$ be a family of integers indexed by $\Omega_{\mathrm{fin}}$.
To each $\omega \in \Omega_{\mathrm{fin}}$, we assign a norm $\|\ndot\|_{\omega}$ of $K_{\omega}$ given by
\[
\| x \|_{\omega} = \exp\left(\frac{-b_{\omega} \log p_{\omega}}{\ord_{\omega}(p_{\omega})}\right) |x|_{\omega}
\]
for $x \in K_{\omega}$.
Moreover, for ${\omega} \in \Omega_{\infty}$ let $\|\ndot\|_{\omega}$ be the standard absolute value of either $\mathbb R$ or $\mathbb C$. 
Then $\xi := \{ \|\ndot\|_{\omega} \}_{\omega \in \Omega}$ yields a norm family on $K$.
Note that, for $\omega \in \Omega_{\mathrm{fin}}$, $\| x \|_{\omega} \leqslant 1$ if and only if
$\ord_{\omega}(x) + b_{\omega} \geqslant 0$ for $x \in K$, that is,
\[
(K, \xi)^{\mathrm{fin}}_{\leqslant 1} = \mathfrak o_{K}(\{ b_ {\omega} \}_{\omega \in \Omega_{\mathrm{fin}}}).
\]
For example, if we set $b_{\omega} = 1$ for all $\omega \in \Omega_{\mathrm{fin}}$,
then $(K, \xi)^{\mathrm{fin}}_{\leqslant 1}$ is not finitely generated over $\mathfrak o_K$ by Lemma~\ref{lem:finite:gen:number:field}.
\end{enonce}

\subsection{Invariants $\bblambda$ and $\bbsigma$}\label{subsection:inv:lambda:classical:settings}
Let $\left(E, \xi = \{ \|\ndot\|_{\omega} \}_{\omega \in \Omega}\right)$ be an 
adelic vector bundle on $S$. 
Let
$\mathscr E := (E, \xi)^{\mathrm{fin}}_{\leqslant 1} = \{ x \in E \,:\, \text{$\| x \|_{\omega} \leqslant 1$ for
all $\omega \in \Omega_{\mathrm{fin}}$} \}$.
If $\xi$ is coherent and dominated, then,
by Proposition~\ref{prop:finiteness:cond:classical:setting}, \ref{prop:finite:generation:dominated} and Remark~\ref{remark:local:E:complete:E},
$\mathscr E$ is finitely generated $\mathfrak o_K$-module,
$\mathscr E \otimes_{\mathfrak o_K} K = E$ and
$\mathscr E \otimes_{\mathfrak o_K} \mathfrak o_{\omega} 
= \mathscr E_{\omega}$
for all $\omega \in \Omega_{\mathrm{fin}}$.

We define $\|\ndot\|_{\infty}$ to be
\[
\forall\, x \in E, \quad
\| x \|_{\infty} := \max_{\omega \in \Omega_{\infty}} \{ \| \iota_{\omega}(x) \|_{\sigma} \},
\]
where $\iota_{\omega}$ is the canonical homomorphism
$E \to E_{\omega}$. 
Under the assumption that $\xi$ is coherent and dominated, the invariant $\bblambda(E,\xi)$ is defined to be
\[
\bblambda(E, \xi) := 
\begin{cases}
\infty & \text{if $E = \{ 0 \}$}, \\[1ex]
\sup \left\{ \lambda \in \mathbb R \, :\,  
\begin{array}{l} \text{There is a basis $e_1, \ldots, e_r$ of $E$} \\
\text{over $K$ such that $e_1, \ldots, e_r \in \mathscr E$} \\
\text{and $\max \{ \| e_1 \|_{\infty}, \ldots,  \| e_r \|_{\infty} \} \leqslant e^{-\lambda}$} 
\end{array} \right\} & \text{otherwise}.
\end{cases}
\]

\medskip
By Proposition~\ref{Pro:normetreausauxdisc},
\begin{equation}\label{subsection:inv:lambda:classical:settings:eqn:02}
0 \leqslant \sup_{x \in E_{\omega} \setminus \{ 0 \}} \ln \left(\frac{\| x \|_{\mathscr E_{\omega}}}{\| x \|_{\omega}}\right) \leqslant -\log |\varpi_{\omega}|_{\omega}
\end{equation}
for any $\omega \in \Omega_{\mathrm{fin}}$, where $\|\ndot\|_{\mathscr E_{\omega}}$ is the norm arising from $\mathscr E_{\omega}$ 
(cf. Subsection~\ref{Subsec:latticesandnorms}). The \emph{impurity}\index{impurity} $\bbsigma(E,\xi)$ of $(E,\xi)$ is defined to be
\[
\bbsigma(E,\xi) := 
\sum_{\omega \in \Omega_{\mathrm{fin}}} \sup_{x \in E_{\omega} \setminus \{ 0 \}} \ln \left(\frac{\| x \|_{\mathscr E_{\omega}}}{\| x \|_{\omega}}\right) \nu(\{\omega \}) \in [0, \infty].
\]
Note that $\bbsigma(E,\xi) = 0$
if and only if $\|\ndot\|_{\omega} = \|\ndot\|_{\mathscr E_{\omega}}$ for all $\omega \in \Omega_{\mathrm{fin}}$.
Moreover, if $\xi$ is coherent and  dominated, then, by Proposition~\ref{prop:strongly:cohenent},
$\xi' = \{ \|\ndot\|_{\mathscr E_{\omega}} \}_{\omega \in \Omega_{\mathrm{fin}}} \cup \{ \|\ndot\|_{\omega} \}_{\omega \in \Omega_{\infty}}$ is also coherent and dominated, so that $\bbsigma(E, \xi) < \infty$ by Corollary~\ref{Cor:dominatedanddist}.

\begin{prop}
We assume that $\xi$ is coherent. Then the following are equivalent:
\begin{enumerate}[label=\rm(\arabic*)]
\item
$\xi$ is dominated.

\item
$\mathscr E$ is finitely generated over $\mathfrak o_K$ and $\bbsigma(E, \xi) < \infty$.
\end{enumerate}
\end{prop}

\begin{proof}
It is sufficient to see that (2) $\Longrightarrow$ (1).
If we set 
\[
\xi' = \{ \|\ndot\|_{\mathscr E_{\omega}} \}_{\omega \in \Omega_{\mathrm{fin}}} \cup \{ \|\ndot\|_{\omega} \}_{\omega \in \Omega_{\infty}},
\] 
then $\xi'$ is dominated by Proposition~\ref{prop:strongly:cohenent} together with Proposition~\ref{prop:finiteness:cond:classical:setting} and Remark~\ref{remark:local:E:complete:E}.
Therefore the assertion follows from the assumption $\bbsigma(E, \xi) < \infty$.
\end{proof}

\begin{prop}\label{prop:comp:lambda:nu:number:field}
We assume that $\xi$ is coherent and dominated.
There is a constants $c_K$
depending only on $K$ 
such that
\[
[K :\mathbb Q]\bblambda(E,\xi) \leqslant \nu_{\min}(E,\xi) \leqslant [K :\mathbb Q]\bblambda(E,\xi) + 
\bbsigma(E,\xi) + c_K.
\]
\end{prop}

\begin{proof}
First we consider the inequality $[K :\mathbb Q]\bblambda(E,\xi) \leqslant \nu_{\min}(E,\xi)$.
We set $\lambda = \bblambda(E,\xi)$. For $\epsilon > 0$, 
there is a basis $\{ e_1, \ldots, e_r \}$ of $E$ over $K$
such that $e_i \in \mathscr E$ and
$\| e_i \|_{\infty} \leqslant e^{-\lambda + \epsilon}$ for all $i$.
On the other hand,
\begin{align*}
\widehat{\deg}_{\xi}(e_i) & = \sum_{\omega \in \Omega} - \ln \| e_i \|_{\omega}\nu(\{\omega\}) \geqslant
\sum_{\omega \in \Omega_{\infty}} - \ln \| e_i \|_{\omega} \nu(\{\omega\}) \geqslant 
\sum_{\omega \in \Omega_{\infty}} (\lambda - \epsilon)\nu(\{\omega\}) \\
& = [K : \mathbb Q] (\lambda - \epsilon),
\end{align*}
so that the assertion follows.

\medskip
(2) Next let us see the second inequality
\[
\nu_{\min}(E,\xi) \leqslant [K :\mathbb Q]\bblambda(E,\xi) + 
\bbsigma(E,\xi) + c_K.
\]
For $\epsilon > 0$, 
there is a basis $\{ e'_1, \ldots, e'_{r} \}$ of $E$ over $K$
such that $\widehat{\deg}(e'_i) \geq \nu_{\min}(E,\xi) - \epsilon$ for all $i$.
We set $E_i = K e'_i$ and $\mathscr E_i = \mathscr E \cap E_i$.
By Lemma~\ref{lemm:finite:invertible:OK:module} below,
there is an $e''_i \in \mathscr E_i$ such that $\#(\mathscr E_i / \mathfrak o_K e''_i ) \leqslant C'_K$,
where $C'_K$ is a constant depending only on the number field $K$.
Therefore $\nu_{\min}(E,\xi)-\epsilon$ is bounded from above by
\begin{align*}
 & \quad\; \widehat{\deg}_{\xi}(e'_i) =  \widehat{\deg}_{\xi}(e''_i) 
= \sum_{\omega \in \Omega_{\mathrm{fin}}} - \ln \| e''_i \|_{\omega} \nu(\{\omega\})
+ \sum_{\omega \in \Omega_{\infty}} - \ln \| e''_i \|_{\omega} \nu(\{\omega\}) \\
& = \sum_{\omega \in \Omega_{\mathrm{fin}}} - \ln \| e''_i \|_{\mathscr E_\omega}\nu(\{\omega\}) +
\sum_{\omega \in \Omega_{\mathrm{fin}}} \ln \left(\frac{\| e''_i \|_{\mathscr E_{\omega}}}{\| e''_i \|_{\omega}}\right) \nu(\{\omega\}) \\
& \hskip21em + \sum_{\omega \in \Omega_{\infty}} - \ln \| e''_i \|_{\omega} \nu(\{\omega\})\\
& \leqslant \ln \#(\mathscr E_i /\mathfrak o_K e''_i) +
\sum_{\omega \in \Omega_{\mathrm{fin}}} \ln (\|\mathrm{Id}_{E_{\omega}}\|^{\mathrm{op}}_{\omega}) \nu(\{\omega\})
+ \sum_{\omega \in \Omega_{\infty}} - \ln \| e''_i \|_{\omega}\nu(\{\omega\}) \\
& \leqslant \ln C'_K + \bbsigma(E,\xi) + \sum_{\omega \in \Omega_{\infty}} -\ln \| e''_i \|_{\omega} \nu(\{\omega\}).
\end{align*}
If we set 
\[
A = \frac{1}{[K : \mathbb Q]} \sum_{\omega \in \Omega_{\infty}} \ln \| e''_i \|_{\omega}\nu(\{\omega\}),
\]
then $\sum_{\omega \in \Omega_{\infty}} (\ln \| e''_i \|_{\omega} - A) \nu(\{\omega\})= 0$.
Let $\{ u_1, \ldots, u_s \}$ be a free basis of $\mathfrak o_K^{\times}$ modulo the torsion subgroup.
Then, by Dirichlet's unit theorem, there are $a'_{i1}, \ldots, a'_{is} \in \mathbb R$ such that
\[
\ln \| e''_i \|_{\omega} - A  = \sum_{j=1}^s a'_{ij} \ln |u_j|_{\omega}
\]
for all $\omega \in \Omega_{\infty}$. 
Let $a_{ij}$ be the round-up of $a'_{ij}$.
Then
\[
\sum_{j=1}^s (a'_{ij} - a_{ij})  \ln |u_j|_{\omega} \leqslant \sum_{j=1}^s |a'_{ij} - a_{ij}|\cdot \big|\ln |u_j|_{\omega}\big|
\leqslant \sum_{j=1}^s \big|\ln |u_j|_{\omega}\big| \leqslant C''_K,
\]
where $C''_K =  \sum_{\omega \in \Omega_{\infty}} \sum_{j=1}^s |\ln |u_j|_{\omega}|$.
Therefore,
\begin{align*}
-A & = \sum_{j=1}^s a'_{ij} \ln |u_j|_{\omega} - \ln \| e''_i \|_{\omega} \leqslant C''_K  + \sum_{j=1}^s a_{ij} \ln |u_j|_{\omega} - \ln \| e''_i \|_{\omega} \\
& = C''_K - \ln \| v_i e''_i\|_{\omega},
\end{align*}
where $v_i = \prod_{j=1}^s u_j^{-a_{ij}}$, and hence, if we set $e_i = v_i e''_i$, then $e_i \in \mathscr E$ and
\[
\nu_{\min}(E,\xi) - \epsilon \leqslant \ln C'_K +\bbsigma(E,\xi) + [K : \mathbb Q]C''_K - [K : \mathbb Q] \ln \| e_i\|_{\omega},
\]
that is, there is a constant $c_K$ depending only on $K$ such that 
\[
\nu_{\min}(E,\xi) - \epsilon \leqslant c_K +\bbsigma(E,\xi) - [K : \mathbb Q] \ln \| e_i\|_{\omega}
\]
for all $i$ and $\omega \in \Omega_{\infty}$. We choose $i$ and $\omega$ such that
$\max \{ \|e_1\|_{\infty}, \ldots, \|e_r\|_{\infty} \} = \| e_i \|_{\omega}$.
Then, as $e^{-\bblambda(E,\xi)} \leqslant \| e_i \|_{\omega}$, that is,
$-\ln \| e_i \|_{\omega} \leqslant \bblambda(E, \xi)$, 
\[
\nu_{\min}(E,\xi) - \epsilon \leqslant c_K +\bbsigma(E,\xi) + [K : \mathbb Q] \bblambda(E, \xi),
\]
and hence the assertion follows.
\end{proof}

\begin{lemm}
\label{lemm:finite:invertible:OK:module}
There is a constant $e_K$ depending only on $K$ such that, for any invertible $\mathfrak o_K$-module $\mathscr L$,
we can find $l \in \mathscr L \setminus \{ 0 \}$ such that $\#(\mathscr L/\mathfrak o_K l) \leqslant e_K$.
\end{lemm}

\begin{proof}
Since the class group is finite, there are finitely many 
invertible $\mathfrak o_K$-modules $\mathscr L_1, \ldots, \mathscr L_h$ such that,
for any invertible $\mathfrak o_K$-module $\mathscr L$, there is $\mathscr L_i$ such that $\mathscr L_i \simeq \mathscr L$.
For each $i = 1, \ldots, h$, fix $l_i \in \mathscr L_i \setminus \{ 0 \}$.
Let $\mathscr L$ be an invertible $\mathfrak o_K$-module. Then there are $\mathscr L_i$ and an isomorphism
$\varphi : \mathscr L_i \to \mathscr L$. If we set $l = \varphi(l_i)$, then
$\mathscr L_i/\mathfrak o_K l_i \simeq \mathscr L /\mathfrak o_K l$, as required.
\end{proof}

%% file: ch5_2019_03_23.tex

\chapter{Slopes of tensor product}
\label{Chap: Slsemistableopes of tensor product}

\IfChapVersion
\ChapVersion{Version of Chapter 5 : \\ \StrSubstitute{\DateChapFive}{_}{\_}}
\fi
The purpose of this chapter is to study the minimal slope of the tensor product of a finite of adelic vector bundles on an adelic curve. More precisely, give a family $\overline E_1,\ldots,\overline E_d$ of adelic vector bundles over an adelic curve $S$, we give a lower bound of $\widehat{\mu}_{\min}(\overline E_1\otimes_{\varepsilon,\pi}\cdots\otimes_{\varepsilon,\pi}\overline{E}_d)$ in terms of the sum of the minimal slopes of $\overline E_i$ minus a term which is the product of three half of the measure of the infinite places and the sum of $\ln(\rang(E_i))$, see Corollary \ref{Cor: tensorial minimal slope property} for details. This result, whose form is similar to the main results of \cite{Gaudron_Remond13,Bost_Chen,Chen_pm}, does not rely on the comparison of successive minima and the height proved in \cite{Zhang95}, which des not hold for general adelic curves. Our method inspires the work of Totaro \cite{Totaro96} on $p$-adic Hodge theory and relies on the geometry invariant theory on Grassmannian. The chapter is organised as follows. In the first section, we regroup several fundamental property of $\mathbb R$-filtrations. We then recall in the second section some basic notions and results of the geometric invariant theory, in particular the Hilbert-Mumford criterion of the semistability. In the third section we give an estimate for the slope of a quotient adelic vector bundle of the tensor product adelic vector bundle, under the assumption that the underlying quotient space, viewed as a rational point of the Grassmannian (with the Pl\"ucker coordinates), is semistable in the sense of geometric invariant theory. In the fifth section, we prove a non-stability criterion which generalises \cite[Proposition 1]{Totaro96}. Finally, we prove in the sixth section the lower bound of the minimal slope of the tensor product adelic vector bundle in the general case.

\section{Reminder on $\mathbb R$-filtrations}
Let $K$ be a field. We equip $K$ with the trivial absolute value $|\ndot|$ such that $|a|=1$ for any $a\in K\setminus\{0\}$. Note that $K$ equipped with the trivial absolute value forms an adelic curve whose underlying measure space is a one point set equipped with the counting measure (which is a probability measure), see \S\ref{Subsec: copies of trivial absoulte value}. Moreover, any finite-dimensional normed vector space over $(K,|\ndot|)$ can be considered as an adelic vector bundle on $S$. In fact, if $V$ is a finite-dimensional vector space over $K$, any norm on $V$ can be considered as a norm family indexed by the one point set. This norm family is clearly measurable. It is also dominated since all norms on $V$ are equivalent (see Corollaries \ref{Coro:equivalenceofnrom} and \ref{Cor:dominatedanddist}). 

Let $V$ be a finite-dimensional vector space over $K$. Recall that the set of ultrametric norms on $V$ are canonically in bijection with the set of $\mathbb R$-filtrations on $V$ (see Remark~\ref{Rem: R-filtration as flag plus slopes}).
If $\|\ndot\|$ is an ultrametric norm on $V$, then the balls centered at  the origin are vector subspaces of $V$, and $\{(V,\|\ndot\|)_{\leqslant\mathrm{e}^{-t}}\}_{t\in\mathbb R}$ is the corresponding $\mathbb R$-filtration. Conversely, given an $\mathbb R$-filtration $\mathcal F$ on $V$, we define a function $\lambda_{\mathcal F}:V\rightarrow\mathbb R\cup\{+\infty\}$ as follows
\[\forall\,x\in V,\quad\lambda_{\mathcal F}(x):=\sup\{t\in\mathbb R\,:\,x\in\mathcal F^t(V)\}.\] 
Then the ultrametric norm $\|\ndot\|_{\mathcal F}$ corresponding to the $\mathbb R$-filtration $\mathcal F$ is given by
\[\forall\,x\in V,\quad \|x\|_{\mathcal F}=\mathrm{e}^{-\lambda_{\mathcal F}(x)}.\]

{
\begin{defi}
Let $V$ be a finite-dimensional vector space over $K$ and $\mathcal F$ be an $\mathbb R$-filtration on $V$. For any $t\in\mathbb R$, we denote by $\operatorname{sq}^t_{\mathcal F}(V)$ the quotient vector space
\[\mathcal F^t(V)\Big/\bigcup_{\varepsilon>0}\mathcal F^{t+\varepsilon}(V).\]
Clearly, if $\mathcal F$ corresponds to the flag 
\[0=V_0\subsetneq V_1\subsetneq\ldots\subsetneq V_n=V\]
of vector subspaces of $V$ together with the sequence
\[\mu_1>\ldots>\mu_n\]
in $\mathbb R$, then 
\[\forall\,i\in\{1,\ldots,n\},\quad\operatorname{sq}_{\mathcal F}^{\mu_i}(V)=V_i/V_{i-1},\]
and $\operatorname{sq}_{\mathcal F}^t(V)=\{0\}$ if $t\not\in\{\mu_1,\ldots,\mu_n\}$. 
\end{defi}
}

\begin{prop}\label{Pro: property of normed vector space over trivial valued field}
Let $(V,\|\ndot\|)$ be a finite-dimensional ultrametrically normed vector space over $K$. The following assertions hold. 
\begin{enumerate}[label=\rm(\arabic*)]
\item\label{Item: admit an orthogonal basis} The normed vector space $(V,\|\ndot\|)$ admits an orthogonal basis. 
\item\label{Item: compute of Arakelov degree} If $\boldsymbol{e}=\{e_i\}_{i=1}^r$ is an orthogonal basis of $(V,\|\ndot\|)$, then the Arakelov degree of $(V,\|\ndot\|)$ is equal to $\lambda_{\mathcal F}(e_1)+\cdots+\lambda_{\mathcal F}(e_r)$.
\item\label{Item: criterion of orthgonal} A basis $\boldsymbol{e}=\{e_i\}_{i=1}^r$ of $V$ is orthogonal if and only if it is compatible with the $\mathbb R$-filtration $\mathcal F$, namely $\#(\mathcal F^t(V)\cap\boldsymbol{e})=\rang(\mathcal F^t(V))$ for any $t\in\mathbb R$.
\item\label{Item: semistability condition} Assume that the vector space $V$ is non-zero. The adelic vector bundle $(V,\|\ndot\|)$ is semistable if and only if the function $\|\ndot\|$ is constant on $V\setminus\{0\}$.
\item \label{Item: Harder-Narasimhan filtration of V}The Harder-Narasimhan $\mathbb R$-filtration of $(V,\|\ndot\|)$ identifies with $\mathcal F$.
\item\label{Item: orthogonal basis and slopes} Let $\boldsymbol{e}=\{e_i\}_{i=1}^r$ be an orthogonal basis of $(V,\|\ndot\|)$. Then the sequence of successive slopes of $(V,\|\ndot\|)$ identifies with the sorted sequence of $\{\lambda_{\mathcal F}(e_i)\}_{i=1}^r$.
\item\label{Item: HN filtration subquotient in trivial valuation case} Let 
\[0=V_0\subsetneq V_1\subsetneq \ldots\subsetneq V_r=V\]
be a \emph{complete} flag of vector subspaces of $V$. For any $i\in\{1,\ldots,r\}$, let $\|\ndot\|_i$ be the subquotient norm of $\|\ndot\|$ on the vector space $V_i/V_{i-1}$. Then the sequence of successive slopes of $(V,\|\ndot\|)$ identifies with the sorted sequence of \[\left\{\hdeg(V_i/V_{i-1},\|\ndot\|_i)\right\}_{i=1}^r.\]
\end{enumerate}
\end{prop}
\begin{proof}
\ref{Item: admit an orthogonal basis} Note that the valued field $(K,|\ndot|)$ is locally compact. By Proposition \ref{Pro:existenceepsorth}, there exists an orthogonal basis of $(V,\|\ndot\|)$.

\ref{Item: compute of Arakelov degree} Let $\boldsymbol{e}=\{e_i\}_{i=1}^r$ be an orthogonal basis of $(V,\|\ndot\|)$. By Proposition \ref{Pro:orthogonalesthadamard}, it is an Hadamard basis, namely 
\[\|e_1\wedge\cdots\wedge e_r\|=\prod_{i=1}^r\|e_i\|.\] 
Therefore one has
\[\hdeg(V,\|\ndot\|)=-\ln\|e_1\wedge\cdots\wedge e_r\|=-\sum_{i=1}^r\ln\|e_i\|=\sum_{i=1}^r\lambda_{\mathcal F}(e_i)\]

\ref{Item: criterion of orthgonal}  Assume that the $\mathbb R$-filtration $\mathcal F$ corresponds to the flag 
\[0=V_0\subsetneq V_1\subsetneq \ldots\subsetneq V_n=V\]
together with the sequence \[\mu_1>\ldots>\mu_n\]
(cf. Remark~\ref{Rem: R-filtration as flag plus slopes}).
Let $\boldsymbol{e}=\{e_i\}_{i=1}^r$ be a basis of $V$. Then $\boldsymbol{e}$ is compatible with the $\mathbb R$-filtration $\mathcal F$ if and only if $\#(\boldsymbol{e}\cap V_j)=\rang(V_j)$ for any $j\in\{1,\ldots,n\}$. By Proposition \ref{Pro: alpha orth is orth} \ref{Item: criterion of orthgonal basis}, this condition is also equivalent to the orthogonality of the basis $\boldsymbol{e}$.

\ref{Item: semistability condition} follows directly from Proposition \ref{Pro: condition of semistablity} since $\norm{\ndot}=\norm{\ndot}_{**}$ (see Corollary \ref{Cor:doubledual}).
 
\ref{Item: Harder-Narasimhan filtration of V} The $\mathbb R$-filtration corresponds to an increasing flag
\[0=V_0\subsetneq V_1\subsetneq\ldots\subsetneq V_n=V\]
of vector subspaces of $V$, together with a decreasing sequence of real numbers
\[\mu_1>\ldots>\mu_n.\]
Note that for any $i\in\{1,\ldots,n\}$ and any $x\in V_{i}\setminus V_{i-1}$ one has $\lambda_{\mathcal F}(x)=\mu_i$. In particular, the subquotient norm $\|\ndot\|_i$ on $V_i/V_{i-1}$ induced by $\|\ndot\|$ takes constant value $\mathrm{e}^{-\mu_i}$ on $(V_i/V_{i-1})\setminus\{0\}$. Therefore, by \ref{Item: semistability condition} the adelic vector bundle $(V_i/V_{i-1},\|\ndot\|_i)$ is semistable of slope $\mu_i$. By Proposition \ref{Pro:characterisationdehnfiltration}, we obtain that $\mathcal F$ is the Harder-Narasimhan $\mathbb R$-filtration of $(V,\|\ndot\|)$.

\ref{Item: orthogonal basis and slopes} Assume that the $\mathbb R$-filtration $\mathcal F$ corresponds to the flag
\[0=V_0\subsetneq V_1\subsetneq\ldots\subsetneq V_n=V\]
and the sequence
\[\mu_1>\ldots>\mu_n.\]
By definition, the value $\mu_i$ appears exactly $\rang(V_i/V_{i-1})$ times in the successive slopes of $(V,\|\ndot\|)$. Moreover, a basis $\boldsymbol{e}$ is orthogonal if and only if it is compatible with the flag \[0=V_0\subsetneq V_1\subsetneq\ldots\subsetneq V_n=V,\]
or equivalently, for any $i\in\{1,\ldots,n\}$, the set $\boldsymbol{e}\cap (V_i\setminus V_{i-1})$ contains exactly $\rang(V_i/V_{i-1})$ elements. Since the function $\lambda_{\mathcal F}(\ndot)$ takes the constant value $\mu_i$ on $V_i\setminus V_{i-1}$, we obtain the assertion.

\ref{Item: HN filtration subquotient in trivial valuation case} By Proposition \ref{Pro:existenceepsorth}
, there exists an orthogonal basis $\boldsymbol{e}=\{e_i\}_{i=1}^r$ which is compatible with the flag
\[0=V_0\subsetneq V_1\subsetneq\ldots\subsetneq V_{r}=V.\] Without loss of generality, we may assume that $e_i\in V_i\setminus V_{i-1}$ for any $i\in\{1,\ldots,r\}$ 
Since the basis $\boldsymbol{e}=\{e_i\}_{i=1}^r$ is orthogonal, the image of $e_i$ in $V_i/V_{i-1}$ has norm $\|e_i\|$. In fact, any element $x$ in $e_i+V_{i-1}$ can be written in the form 
\[e_i+\sum_{j=1}^{i-1}a_je_j\]
and hence $\|x\|_i\geqslant \|e_i\|$. Therefore one has \[\hdeg(V_i/V_{i-1},\|\ndot\|_i)=-\ln\|e_i\|=\lambda_{\mathcal F}(e_i).\]  
\end{proof}

\begin{prop}\label{Pro: simultaneous orthogonal trivial valuation case}
Let $V$ be a finite-dimensional vector space over $K$ and $\|\ndot\|$ and $\|\ndot\|'$ be two ultrametric norms on $V$. Then there exists a basis $\boldsymbol{e}$ of $V$ which is orthogonal with respect to $\|\ndot\|$ and $\|\ndot\|'$ simultaneously.
\end{prop}
\begin{proof}
Let $\mathcal F$ be the $\mathbb R$-filtration on $V$ associated with the norm $\|\ndot\|$, which corresponds to a flag
\[0=V_0\subsetneq V_1\subsetneq \ldots\subsetneq V_n=V\]
together with a sequence $\mu_1>\ldots>\mu_n$.
By  Proposition \ref{Pro:existenceepsorth}, there exists an orthogonal basis $\boldsymbol{e}$ of $(V,\|\ndot\|')$ which is compatible with the the flag
\[0=V_0\subsetneq V_1\subsetneq \ldots\subsetneq V_n=V.\]
By Proposition \ref{Pro: property of normed vector space over trivial valued field} \ref{Item: criterion of orthgonal}, we obtain that $\boldsymbol{e}$ is also orthogonal with respect to $\|\ndot\|$.
\end{proof}

\begin{coro}
Let $V$ be a finite-dimensional vector space over $K$, and
\[0=V_0\subsetneq V_1\subsetneq\ldots\subsetneq V_n = V\quad\text{and}\quad 0=W_0\subsetneq W_1\subsetneq\ldots\subsetneq W_m=V\]
be two flags of vector subspaces of $V$. There exists a basis $\boldsymbol{e}$ of $V$ which is compatible with the two flags simultaneously.
\end{coro}
\begin{proof}
By choosing two decreasing sequences of real numbers $\mu_1>\ldots>\mu_n$ and $\lambda_1>\ldots>\lambda_m$ we obtain two $\mathbb R$-filtrations on $V$, which correspond to two ultrametric norms on $V$. By Proposition \ref{Pro: simultaneous orthogonal trivial valuation case}, there exists a basis of $V$ which is orthogonal with respect to the two norms simultaneously. By Proposition \ref{Pro: property of normed vector space over trivial valued field} \ref{Item: criterion of orthgonal}, this basis is compatible with the two flags simultaneously.
\end{proof}

{
\begin{defi}\label{Def: tensor product of filtrations}
Let $d\in\mathbb N_{\geqslant 2}$ and $(E_j)_{j=1}^d$ be a family of finite-dimensional vector spaces over $K$. For any $j\in\{1,\ldots,d\}$, let $\mathcal F_j$ be an $\mathbb R$-filtrations on $E_j$, which corresponds to an ultrametric norm $\norm{\ndot}_j$ on $E_j$. Let $\norm{\ndot}_{\varepsilon}$ be the $\varepsilon$-tensor product of the norms $\|\ndot\|_1,\ldots,\norm{\ndot}_d$. The $\mathbb R$-filtration on the tensor product space \[E_1\otimes_K\cdots\otimes_K E_d\] corresponding to $\norm{\ndot}_{\varepsilon}$ is called the \emph{tensor product}\index{tensor product@tensor product $\mathbb R$-filtration} of the $\mathbb R$-filtrations $\mathcal F_1,\ldots,\mathcal F_d$, which is denoted by $\mathcal F_1\otimes\cdots\otimes\mathcal F_d$.
\end{defi}

\begin{rema}\label{Rem: tensor product filtrations}
We keep the notation of Definition \ref{Def: tensor product of filtrations}. For any $j\in\{1,\ldots,d\}$, let $\boldsymbol{e}^{(j)}=\{e_i^{(j)}\}_{i=1}^{n_j}$ be an orthogonal base of $(E_j,\|\ndot\|_j)$. 
By Proposition \ref{Pro:alphatenso} together with Remark \ref{Rem:produittenrk1}, one has the following:
\begin{enumerate}[label=\rm(\roman*)]
\item $\big\{e_{i_1}^{(1)}\otimes\cdots\otimes e_{i_d}^{(d)}\big\}_{(i_1,\ldots,i_d)\in\prod_{j=1}^d\{1,\ldots,n_j\}}$ forms an
orthogonal basis of the vector space $E_1\otimes_K\cdots\otimes_KE_d$ with respect to $\norm{\ndot}_{\varepsilon}$.
\item $\big\|e_{i_1}^{(1)}\otimes\cdots\otimes e_{i_d}^{(d)}\big\|_{\varepsilon}=\prod_{j=1}^d\big\|e_{i_j}^{(j)}\big\|_j$ for any $(i_1,\ldots,i_d)\in\prod_{j=1}^d\{1,\ldots,n_j\}$.
\end{enumerate}
Therefore, if we denote by $\mathcal F$ the tensor product $\mathbb R$-filtration $\mathcal F_1\otimes\cdots\otimes\mathcal F_d$, then the vector space $\mathcal F^t(E_1\otimes_K\cdots\otimes_KE_d)$ is generated by the vectors $e_{i_1}^{(1)}\otimes\cdots\otimes e_{i_d}^{(d)}$ such that 
\[\lambda_{\mathcal F_1}(e_{i_1}^{(1)})+\cdots+\lambda_{\mathcal F_d}(e_{i_d}^{(d)})\geqslant t.\]
Therefore, one has
\[\begin{split}\mathcal F^t(E_1\otimes_K\cdots\otimes_K E_d)&=\sum_{t_1+\cdots+t_d\geqslant t}\mathcal F_1^{t_1}(E_1)\otimes_K\cdots\otimes_K\mathcal F_d^{t_d}(E_d)\\
&=\sum_{t_1+\cdots+t_d= t}\mathcal F_1^{t_1}(E_1)\otimes_K\cdots\otimes_K\mathcal F_d^{t_d}(E_d).
\end{split}\]
Furthermore, by (i), (ii) and Proposition~\ref{Pro: property of normed vector space over trivial valued field} \ref{Item: orthogonal basis and slopes},
if $(E_j, \|\ndot\|_j)$ are all semistable, where $j\in\{1,\ldots,d\}$, then
$(E_1\otimes_K\cdots\otimes_K E_d,\|\ndot\|_{\varepsilon})$ is also semistable.
\end{rema}}

\section{Reminder on geometric invariant theory}

Let $K$ be a field. By \emph{group scheme over $\Spec K$}\index{group scheme over Spec K@group scheme over $\Spec K$} or by \emph{$K$-group scheme}\index{K-group scheme@$K$-group scheme}, we mean a $K$-scheme $\pi:G\rightarrow\Spec K$ equipped with a $K$-morphism $m_G:G\times_K G\rightarrow G$ (called the \emph{group scheme structure of $G$}\index{group scheme structure of G@group scheme structure of $G$}) such that, for any $K$-scheme $f:S\rightarrow\Spec K$, the set $G(S)$ of $K$-morphisms from $S$ to $G$ equipped with the composition law $m_G(S):G(S)\times G(S)\rightarrow G(S)$ forms a group. Note that the maps of inverse 
$\iota_G(S):G(S)\rightarrow G(S)$ and the maps of unity \[e_G(S):\Spec K(S)=\{S\stackrel{f}{\rightarrow}\Spec K\}\longrightarrow G(S)\] actually define $K$-morphisms $\iota_G:G\rightarrow G$ and $e_G:\Spec K\rightarrow G$, which make the following diagrams commutative:
\[\xymatrix{G\times_KG\times_KG\ar[rr]^-{m_G\times \Id_G}\ar[d]_-{\Id_G\times m_G}&&G\times_KG\ar[d]^-{m_G}\\
G\times_KG\ar[rr]_-{m_G}&&G}\]

\[\xymatrix{G\ar[rrd]_-{\Id_G}\ar[rr]^-{(e_G\pi,\Id_G)}&&G\times_KG\ar[d]^-{m_G}&&G\ar[lld]^-{\Id_G}\ar[ll]_-{(\Id_G,e_G\pi)}\\
&&G}\]

\[\xymatrix{G\ar[r]^-{(\Id_G,\iota_G)}\ar[d]_-{(\iota_G,\Id_G)}\ar[rd]|-{\parbox[c][1em][c]{20 pt}{\hspace{2pt}\small$e_G\pi$}}&G\times_KG\ar[d]^-{m_G}\\G\times_KG\ar[r]_-{m_G}&G}\]

Let $G$ and $H$ be group schemes over $\Spec K$. We call \emph{morphism of $K$-group schemes}\index{morphism of K-group schemes@morphism of $K$-group schemes} from $G$ to $H$ any $K$-morphism $f:G\rightarrow H$ such that, for any $K$-scheme $S$, le morphism $f$ induces a morphism of groups $f(S):G(S)\rightarrow H(S)$.

\begin{exem}
Let $V$ be a finite-dimensional vector space over $K$. We denote by $\mathbb{GL}(V)$ the open subscheme of the affine $K$-scheme $\mathbb A(\mathrm{End}(V)^\vee)$
\if01
\lquery{90pt}{{It seems that a dual is needed here. Could you please help me to check this point?} {To define the affine space $\mathbb A(V)$, one has two ways: either $\mathbb A(V) := \mathrm{Spec}(S(V))$ or $\mathbb A(V) := \mathrm{Spec}(S(V^{\vee}))$, where $S(V) = \bigoplus_{m=0}^{\infty} \mathrm{Sym}^m(V)$. If we think $-$ of $\mathrm{Spec}(-)$ is the space of functions on $V$, then it is better to choose $S(V^{\vee})$. In our case, it does not matter, so that
the choice of $S(V)$ seems to be good.}} \fi
defined by the non-vanishing of the determinant. For any $K$-scheme $\pi:S\rightarrow\Spec K$, one has
\[\mathbb{GL}(V)(S)=\mathrm{Aut}_{\mathcal O_S}(\pi^*(V)).\]
The set $\mathbb{GL}(V)(S)$ is canonically equipped with a structure of group. The group structures of $\mathbb{GL}(V)(S)$ where $S$ runs over the set of $K$-schemes define a $K$-morphism $\mathbb{GL}(V)\times_K\mathbb{GL}(V)\rightarrow \mathbb{GL}(V)$, which makes $\mathbb{GL}(V)$ a group scheme over $K$. The group scheme $\mathbb{GL}(V)$ is called the \emph{general linear group scheme associated with $V$}\index{general linear group scheme associated with V@general linear group scheme associated with $V$}.
\end{exem}

\begin{defi}
Let $G$ be a group scheme over $\Spec K$ and $X$ be a scheme over $\Spec K$. As \emph{action}\index{action@action} of $G$ on $X$, we refer to a $K$-morphism $f:G\times_KX\rightarrow S$ such that, for any $K$-scheme $S$, the map
\[f(S):G(S)\times X(S)\longrightarrow X(S)\]
defines an action of the group $G(S)$ on $X(S)$. 
\end{defi}

\begin{exem}\label{Exa: action of GL on proj}
Let $V$ be a finite-dimensional vector space over $K$ and $X$ be the projective space $\mathbb P(V)$. Recall that, for any $K$-scheme $p:S\rightarrow\Spec K$, $\mathbb P(V)(S)$ identifies with the set of all invertible quotient modules of $p^*(V)$. Note that the automorphism of $p^*(V)$ acts naturally on the set $\mathbb P(V)(S)$ of invertible quotient modules of $p^*(V)$. Hence we obtain an action of the general linear group scheme $\mathbb{GL}(V)$ on the projective space $\mathbb P(V)$.  

More generally, let $G$ be a group scheme over $\Spec K$. By (finite-dimensional) \emph{linear representation}\index{linear representation@linear representation} of $G$ we refer to a morphism of $K$-group schemes from $G$ to certain $\mathbb{GL}(V)$, where $V$ is a finite-dimensional vector space over $K$. Note that such a morphism induces an action of $G$ on the projective space $\mathbb P(V)$. This action is said to be \emph{linear}\index{linear action@linear action}.
\end{exem}

Let $G$ be a group scheme over $\Spec K$ which acts on a $K$-scheme $X$. We denote by $f:G\times_K X\rightarrow X$ the action of $G$ on $X$ and by $\pr_2:G\times_KX\rightarrow X$ the projection to the second coordinate. Let $L$ be an invertible $\mathcal O_X$-module. We call \emph{$G$-linear structure on $L$}\index{G-linear structure on L@$G$-linear structure on $L$} any isomorphism $\eta$ of $\mathcal O_{G\times_KX}$-modules from $f^*(L)$ to $\pr_2^*(L)$ such that the following diagram commutes
\begin{equation}\label{Equ: compatibility condition}\begin{gathered}\xymatrix{(\pr_2\circ (\Id_G\times f))^*(L)\ar@{=}[rr]&&(f\circ\pr_{23})^*(L)\ar[d]^-{\pr_{23}^*(\eta)}\\
(f\circ(\Id_G\times f))^*(L)\ar[u]^-{(\Id_G\times f)^*(\eta)}\ar@{=}[d]&&(\pr_2\circ\pr_{23})^*(L)\ar@{=}[d]\\
(f\circ(m_G\times\Id_X))^*(L)\ar[rr]_-{(m_G\times\Id_X)^*(\eta)}&&(\pr_2\circ(m\times \Id_X))^*(L)}
\end{gathered}\end{equation}
where $\pr_{23}:G\times_KG\times_K X\rightarrow G\times_KX$ is the projection to the second and the third coordinates, and $m_G:G\times_{K}G
\rightarrow G$ is the group scheme structure on $G$. The couple $(L,\eta)$ is called a \emph{$G$-linearised invertible $\mathcal O_X$-module}\index{G-linearised invertible OX-module@$G$-linearised invertible $\mathcal O_X$-module}.

Note that, if $\eta:f^*(L)\rightarrow\pr_2^*(L)$ is a $G$-linear structure on $L$, then $\eta^\vee:f^*(L^\vee)\rightarrow\pr_2^*(L^\vee)$ is a $G$-linear structure on $L^\vee$. Moreover, if $(L_1,\eta_1)$ and $(L_2,\eta_2)$ are $G$-linearised invertible $\mathcal O_X$-modules, then \[\eta_1\otimes\eta_2:f^*(L_1)\otimes f^*(L_2)\cong f^*(L_1\otimes L_2)\longrightarrow \pr_2^*(L_1\otimes L_2)\cong \pr_2^*(L_1)\otimes\pr_2^*(L_2)\]
is a $G$-linear structure on $L_1\otimes L_2$.

\begin{exem}
Let $G$ a group scheme over $K$ and $V$ be a finite-dimensional vector space over $K$. A linear action of $G$ on $\mathbb P(V)$ defines canonically a $G$-linear structure on the universal invertible sheaf $\mathcal O_V(1)$. Let $f:G\times_K\mathbb P(V)\rightarrow \mathbb P(V)$ be a linear action of $G$ on $\mathbb P(V)$. Let $\pi:\mathbb P(V)\rightarrow\Spec K$ be the structural morphism and  $\beta:\pi^*(V)\rightarrow\mathcal O_V(1)$ be the tautological invertible quotient sheaf of $\pi^*(V)$. Note that the morphism
\[\xymatrix{G\times_K\mathbb P(V)\ar[rr]^-{(\pr_1,f)}&&G\times_K\mathbb P(V)}\]  
is an isomorphism of $K$-schemes, the inverse of which is given by the composed $K$-morphism
\[\xymatrix{G\times_K\mathbb P(V)\ar[rr]^-{(\pr_1,\iota_G\pr_1,\pr_2)}&&G\times_KG\times_K\mathbb P(V)\ar[rr]^-{(\pr_1,f\pr_{23})}&&G\times_K\mathbb P(V)}.\]
Moreover, one has $\pr_2\circ(\pr_1,f)=f$. Therefore $((\pr_1,f)^*)^{-1}$ defines an isomorphism from $f^*(\mathcal O_V(1))$ to $\pr_2^*(\mathcal O_V(1))$ as invertible quotient modules of $f^*(\pi^*(V))\cong\pr_2^*(\pi^*(V))$. The fact that the action of $G$ on $\mathbb P(V)$ is linear shows that this isomorphism actually defines a $G$-linear structure on $\mathcal O_V(1)$. 
\end{exem}

\begin{defi}
We denote by $\mathbb{G}_{\mathrm{m},K}=\Spec K[T,T^{-1}]$ the multiplicative group scheme over $\Spec K$ (recall that for any $K$-scheme $S$ one has $\mathbb G_{\mathrm{m},K}(S)=\mathcal O_S(S)^{\times}$).  If $G$ is a group scheme over $\Spec K$, by \emph{one-parameter subgroup}\index{one-parameter subgroup@one-parameter subgroup} of $G$ any morphism of $K$-group schemes from $\mathbb{G}_{\mathrm{m},K}$ to $G$.
\end{defi}

Let $G$ be a group scheme over $K$, which acts on a $K$-scheme $X$. Denote by $f:G\times_KX\rightarrow X$ the action. If 
$\varphi:\mathbb{G}_{\mathrm{m},K}\rightarrow G$
is a one-parameter subgroup of $G$, then $f$ and $\varphi$ induce an action of $\mathbb{G}_{\mathrm{m},K}$ on $X$, denoted by $f_{\varphi}$. Note that  $f_{\varphi}$ is the composed morphism
\[\xymatrix{\mathbb{G}_{\mathrm{m},K}\times_K X\ar[rr]^-{\varphi\times\Id_X}&&G\times_KX\ar[r]^-f&X}.\]

Let $g:\mathbb{G}_{\mathrm{m},K}\times_KX\rightarrow X$ be an action of the multiplicative group $\mathbb{G}_{\mathrm{m},K}$ on a \emph{proper} $K$-scheme. Suppose that $x:\Spec K\rightarrow X$ is a rational point of $X$. The \emph{orbit}\index{orbit@orbit} of $x$ by the action of $\mathbb{G}_{\mathrm{m},K}$ is by definition the following composed morphism $\mathrm{orb}_x$
\[\xymatrix{\mathbb{G}_{\mathrm{m},K}\cong\mathbb{G}_{\mathrm{m},K}\times_K\Spec K\ar[rr]^-{\Id_{\mathbb{G}_{\mathrm{m},K}}\times x}&&\mathbb{G}_{\mathrm{m},K}\times_KX\ar[r]^-g&X.}\]
Since $X$ is proper over $\Spec K$, by the valuative criterion of properness, the morphism $\mathrm{orb}_x:\mathbb{G}_{\mathrm{m},K}\rightarrow X$ extends in a unique way to a $K$-morphism $\widetilde{\mathrm{orb}}_x:\mathbb A^1_K=\Spec K[T]\rightarrow X$. The image by $\widetilde{\mathrm{orb}}_x$ of the rational point of $\mathbb A_K^1$ corresponding to the prime ideal $(T)$ is denoted by $\widetilde x_g$. Note that $\widetilde x_g$ is a rational point of $X$ which is invariant by the action of $\mathbb{G}_{\mathrm{m},K}$.

Assume that $L$ is a $\mathbb{G}_{\mathrm{m},K}$-linearised invertible $\mathcal O_X$-module. Since $\widetilde x_g$ is a fixed rational point of the action $g$, the $\mathbb{G}_{\mathrm{m},K}$-linear structure corresponds to an action of $\mathbb{G}_{\mathrm{m},K}$ on $\widetilde{x}_g^*(L)$, which is induced by an endomorphism of the $K$-group scheme $\mathbb{G}_{\mathrm{m},K}$. Note that any endomorphism of the $K$-group scheme $\mathbb{G}_{\mathrm{m},K}$ is of the form $t\mapsto t^n$, where the exponent $n$ is an integer. We denote by $\mu(x,L)$ the \emph{opposite}\index{opposite@opposite} of the exponent corresponding to the action of $\mathbb{G}_{\mathrm{m},K}$ on $\widetilde{x}_g^*(L)$. Note that our choice of the constant $\mu(x,L)$ conforms with that of the book \cite{Mumford94}.

More generally, if $G$ is a $K$-group scheme, $f:G\times_KX\rightarrow X$ is an action of $G$ on a proper $K$-scheme $X$ and if $\varphi:\mathbb{G}_{\mathrm{m},K}\rightarrow G$ is a one-parameter subgroup of $G$, for any $x\in X(K)$ we denote by $\mu(x,\varphi, L)$ the exponent corresponding to the action of $\mathbb{G}_{\mathrm{m},K}$ on $\widetilde x_{f_\lambda}^*(L)$.

\begin{exem}\label{Exa: coefficient mu of one-parameter group}
Consider the one-parameter subgroups of the general linear group. Let $E$ be a finite-dimensional vector space over $K$ and $\lambda:\mathbb G_{m,K}\rightarrow\mathbb{GL}(E)$ be a one-parameter subgroup. By \cite[II.\S2, $\text{n}^{\circ}$2, 2.5]{Demazure_Gabriel70}, we can decompose the vector space $E$ as a direct sum of $K$-vector subspaces  $E_1,\ldots,E_n$ which are invariant by the action of $\mathbb G_{m,K}$, and integers $(a_1,\ldots,a_n)$ such that the action of $\mathbb G_{m,K}$ on $E_i$ is given by $t\mapsto t^{a_i}\mathrm{Id}_{E_i}$. Therefore the one-parameter subgroup $\varphi$ determines an $\mathbb R$-filtration $\mathcal F_\varphi$ on $E$ such that 
\[\mathcal F_\varphi^t(E)=\bigoplus_{\begin{subarray}{c}i\in\{1,\ldots,n\}\\
-a_i\geqslant t
\end{subarray}}E_i.\]
We now consider the canonical action of $\mathbb{GL}(E)$ on the projective space $\mathbb P(E)$ (see Example \ref{Exa: action of GL on proj}). Let $x$ be a rational point of $\mathbb P(E)$ and $\pi_x:E\rightarrow K$ be the one dimensional quotient vector space of $E$ corresponding to $x$. Then the morphism $\mathrm{orb}_x:\mathbb G_{m,K}=\Spec K[T,T^{-1}]\rightarrow\mathbb P(E)$ is represented by the surjective $K[T,T^{-1}]$-linear map
$p_x:E\otimes_KK[T,T^{-1}]\longrightarrow K[T,T^{-1}]$ sending $v_i\otimes 1$ to $\pi_x(v_i)T^{a_i}$ for any $v_i\in E_i$. 
The extended morphism $\widetilde{\mathrm{ord}}_x:\mathbb A_K^1\rightarrow\mathbb P(E)$ corresponds to the surjective $K[T]$-linear map
\[\widetilde p_x:E\otimes_KK[T]\longrightarrow K[T]\cdot T^{-\mu(x,\mathcal O_E(1))}\]
given by the restriction of $p_x$ on $E\otimes_KK[T]$, where $\mathcal O_E(1)$ denotes the universal invertible sheaf. In particular, one has
\[\mu(x,\varphi,\mathcal O_E(1))=-\min\{a_i\,:\,i\in\{1,\ldots,n\},\;\pi_x(E_i)\neq \{0\}\}.\] 
Therefore, we can interpret the constant $\mu(x,\varphi,\mathcal O_E(1))$ via the $\mathbb R$-filtration $\mathcal F_\varphi$. In fact, the $\mathbb R$-filtration $\mathcal F_\varphi$ induces by the surjective map $\pi_x:E\rightarrow K$ an $\mathbb R$-filtration on $K$ (viewed as a one-dimensional vector space over $K$), which corresponds to the quotient norm of $\norm{\ndot}_{\mathcal F_\varphi}$. Then the number $\mu(x,\varphi,\mathcal O_E(1))$ is equal to the jump point of this quotient $\mathbb R$-filtration.
\end{exem}

The following theorem of Hilbert-Mumford relates the positivity of $\mu(x,\lambda,L)$ to the non-vanishing of a global section invariant by the  action of the $K$-group scheme.

\begin{theo}
\label{Thm:Hilbert-Mumford}
We assume that the field $K$ is perfect.
Let $G$ be a reductive $K$-group scheme acting on a projective $K$-scheme $X$, $L$ be a $G$-linearised ample invertible $\mathcal O_X$-module. For any rational point $x\in X(K)$, the following two conditions are equivalent:
\begin{enumerate}[label=\rm(\arabic*)]
\item for any one-parameter subgroup $\lambda:\mathbb{G}_{\mathrm{m},K}\rightarrow G$ of $G$, $\mu(x,\lambda,L)\geqslant 0$;
\item\label{Item: HM criterion} there exists an integer $n\in\mathbb N_{\geqslant 1}$ and a section $s\in H^0(X,L^{\otimes n})$ which is invariant under the action of $G(K)$ such that $x$ lies outside of the vanishing locus of $s$.
\end{enumerate}
\end{theo}

{We just explain why the condition \ref{Item: HM criterion} implies the positivity of $\mu(x,\lambda,L)$ for any one-parameter group. Let \[\widetilde{\mathrm{orb}}_x:\mathbb A^1_K=\Spec K[T]\longrightarrow X\] be the extended orbit of the rational point $x$ by the action of $\mathbb{G}_{\mathrm{m},K}$ via $\lambda$. Then the pull-back of $L$ by $\widetilde{\mathrm{orb}}_x$ corresponds to a free $K[T]$-module of rank $1$, which is equipped with a linear action of $\mathbb{G}_{\mathrm{m},K}$. This action corresponds to an invertible element of the tensorial algebra 
\[K[t,t^{-1}]\otimes_KK[T]\cong K[t,t^{-1},T],\]
where $t$ and $T$ are variables. Moreover, the compatibility condition \eqref{Equ: compatibility condition} shows that $\eta(t,T)$ satisfies the following relation
\[\eta(t,T)\eta(u,T)=\eta(tu,T) \text{ in $K[t,t^{-1},u,u^{-1},T]$},\]
where $t$, $u$, and $T$ are variables. Therefore $\eta(t,T)$ is of the form $t^{-\mu(x,\lambda,L)}$.
We fix a section $m$ of $\widetilde{\mathrm{orb}}_x^{\raisebox{-10pt}{\scriptsize{*}}}(L)$ which trivialises this invertible sheaf.
Note that the pull-back of the section $s$ is an element of this free $K[T]$-module which is invariant by the action of $\mathbb{G}_{\mathrm{m},K}(K)=K^{\times}$. We write $s$ in the form $P(T)m$, where $P\in k[T]$. Then the action of an element $a\in K^{\times}$ on $s$ gives the section $P(aT)a^{-\mu(x,\lambda,L)}m$. Hence $P$ is a homogeneous polynomial and $\mu(x,\lambda,L)$ is equal to the degree of $P$, which is a non-negative integer.} We refer the readers to \cite[\S2.1]{Mumford94} for a proof of the above theorem. See also \cite{Kempf78} and \cite{Rousseau78}.

\begin{defi}
Under the assumption and with the notation of Theorem \ref{Thm:Hilbert-Mumford}, if $x\in X(K)$ satisfies the equivalent conditions of the theorem, we say that the point $x$ is \emph{semistable}\index{semistable@semistable} with respect to the $G$-linearised invertible $\mathcal O_X$-module $L$.
\end{defi}

{
\begin{rema}\label{Rem: action of gm on produt}
Let $d\in\mathbb N_{\geqslant 2}$ and $\{E_j\}_{j=1}^d$ be a family of finite-dimensional non-zero vector spaces over $K$. Any one-parameter subgroup \[\lambda:\mathbb{G}_{\mathrm{m},K}\longrightarrow \mathbb{GL}(E_1)\times_K\cdots\times_K\mathbb{GL}(E_d)\] can be written in the form $(\lambda_1,\cdots,\lambda_d)$, where $\lambda_j:\mathbb{G}_{\mathrm{m},K}\rightarrow\mathbb{GL}(E_j)$ is a one-parameter subgroup of $\mathbb{GL}(E_j)$, $j\in\{1,\ldots,d\}$. We can then decompose the vector space $E_j$ into the direct sum of eigenspaces of the action $\lambda_j$ as follows:
\[E_j=E_{j,1}\oplus\cdots\oplus E_{j,n_j},\]
where each $E_{j,i}$ is stable by the action of $\lambda_j$, and on $E_{j,i}$ the action of $\mathbb{G}_{\mathrm{m},K}$ is given by $t\mapsto t^{a_{j,i}}\mathrm{Id}_{E_{j,i}}$, $i\in\{1,\ldots,n_j\}$. Note that the one-parameter subgroup $\lambda$ induces an action of $\mathbb{G}_{\mathrm{m},K}$ on the tensor product space $E_1\otimes_K\cdots\otimes_K E_d$ via the canonical morphisme of $K$-group schemes \[\mathbb{GL}(E_1)\times_K\cdots\times_K\mathbb{GL}(E_d)\longrightarrow\mathbb{GL}(E_1\otimes_K\cdots\otimes_KE_d).\] For any $(i_1,\ldots,i_d)\in\prod_{j=1}^d\{1,\ldots,n_j\}$, the vector subspace $E_{1,i_1}\otimes_K\cdots\otimes_KE_{d,i_d}$ of $E_1\otimes_K\cdots\otimes_KE_d$ is invariant by the action of $\mathbb{G}_{\mathrm{m},K}$, and on $E_{1,i_1}\otimes_K\cdots\otimes_KE_{d,i_d}$ the action of $\mathbb{G}_{\mathrm{m},K}$ is given by \[t\longmapsto t^{a_{1,i_1}+\cdots+a_{d,i_d}}\mathrm{Id}_{E_{1,i_1}\otimes_K\cdots\otimes_K E_{d,i_d}}.\] We construct an $\mathbb R$-filtration $\mathcal F_\lambda$ on $E_1\otimes_K\cdots\otimes_K E_d$ as follows
\begin{equation}\label{Equ:filtration sur le produit}
\mathcal F^t_\lambda(E_1\otimes_K\cdots\otimes_K E_d) :=\sum_{\begin{subarray}{c}
(i_1,\ldots,i_d)\in\prod_{j=1}^d\{1,\ldots,n_j\}\\
-a_{1,i_1}-\cdots-a_{d,i_d}\geqslant t
\end{subarray}}E_{1,i_1}\otimes_K\cdots\otimes_K E_{d,i_d}.\end{equation} 
Moreover, if we denote by $\mathcal F_{\lambda_j}$ the $\mathbb R$-filtrations on $E_j$ defined by 
\[\mathcal F^{t}_{\lambda_j}(E_j)=\sum_{\begin{subarray}{c}
i\in\{1,\ldots,n_j\}\\
-a_i\geqslant t
\end{subarray}}E_{j,i},\]
then the $\mathbb R$-filtration $\mathcal F_\lambda$ defined in \eqref{Equ:filtration sur le produit} identifies with the tensor product of $\mathcal F_{\lambda_1},\ldots,\mathcal F_{\lambda_d}$ (see Definition \ref{Def: tensor product of filtrations}, see also Remark \ref{Rem: tensor product filtrations}). Conversely, for any $\mathbb R$-filtration $\mathcal F_j$ with integral jump points on the vector spaces $E_j$, then there exists a one-parameter subgroup $\lambda_j:\mathbb{G}_{m,K}\rightarrow\mathbb{GL}(E_j)$ such that $\mathcal F_{\lambda_j}=\mathcal F_{j}$. This comes from Proposition \ref{Pro:existenceepsorth} which allows us to construct the actions of $\mathbb G_{m,K}$ on $E_{j}$ diagonally with respect to an orthogonal basis.  

More generally, for any integer $r\geqslant 1$, any one-parameter subgroup \[\lambda=(\lambda_1,\ldots,\lambda_d):\mathbb{G}_{\mathrm{m,K}}\rightarrow\mathbb{GL}(E_1)\times_K\cdots\times_K\mathbb{GL}(E_d)\] induces an action of $\mathbb{G}_{\mathrm{m},K}$ on the $K$-vector space \[(E_1\otimes_K\cdots\otimes_KE_d)^{\otimes r}.\] Again the $\mathbb R$-filtration on $(E_1\otimes_K\cdots\otimes_KE_d)^{\otimes r}$ corresponding to the eigenspace decomposition of the action of $\mathbb{G}_{\mathrm{m},K}$ identifies with \[(\mathcal F_{\lambda_1}\otimes\cdots\otimes\mathcal F_{\lambda_d})^{\otimes r}.\]

For any $j\in\{1,\ldots,d\}$, let $a_j$ be the rank of $E_j$. Consider a non-zero quotient vector space $V$ of $E_1\otimes_K\cdots\otimes_K E_d$. Let $r$ be the rank of $V$ over $K$. The canonical surjective map $E_1\otimes_K\cdots\otimes_K E_d\rightarrow V$ determines a rational point $x$ of \[\pi:P=\mathbb P((E_1\otimes_K\cdots\otimes E_d)^{\otimes r})\longrightarrow\Spec K,\] which corresponds the composed map 
\[(E_1\otimes_K\cdots\otimes_K E_d)^{\otimes r}\longrightarrow V^{\otimes r}\longrightarrow\det(V).\]
We consider the semistability of the point $x$ with respect to the $\mathbb{GL}(E_1)\times_K\cdots\times_K\mathbb{GL}(E_d)$-invertible sheaf \[L:=\mathcal O_P(a_1\cdots a_d)\otimes\pi^*(\det(E_1^{\vee})^{\otimes rb_1}\otimes\cdots\otimes\det(E_d^{\vee})^{\otimes rb_d}),\]
where for any $j\in\{1,\ldots,d\}$,
\[b_j:=\frac{a_1\cdots a_d}{a_j}.\] Let 
\[\lambda=(\lambda_1,\ldots,\lambda_d):\mathbb G_{m,K}\longrightarrow\mathbb{GL}(E_1)\times_K\cdots\times_K\mathbb{GL}(E_d)\] be a one-parameter subgroup, which determine, for each $j\in\{1,\ldots,d\}$, an $\mathbb R$-filtration $\mathcal F_{\lambda_j}$ on $E_j$. We let $\norm{\ndot}_{j}$ be the ultrametric norm on $E_{j}$ corresponding to $\mathcal F_{\lambda_j}$, where we consider the trivial absolute value on $K$. We equip $E_1\otimes_K\cdots\otimes_K E_d$ with the $\varepsilon$-tensor product of the norms $\norm{\ndot}_j$ and equip $V$ with the quotient norm. By Example \ref{Exa: coefficient mu of one-parameter group} and Proposition \ref{Pro:doubledualdet}, we obtain that
\[\mu(x,\lambda,\mathcal O_P(1))=\hdeg(\overline V).\]
Moreover, by definition, for any $j\in\{1,\ldots,d\}$ one has
\[\mu(x,\lambda_j,\pi^*(\det(E_j^{\vee})))=-\hdeg(E_j,\norm{\ndot}_j)\]
Therefore we obtain
\begin{equation}\label{Equ: computation of mu x varphi L}\mu(x,\lambda,L)=a_1\cdots a_d\hdeg(\overline V)-r\sum_{j=1}^{d}b_j\hdeg(\overline E_j).\end{equation}
Therefore we deduce from the Hilbert-Mumford criterion the following result.
\end{rema}

\begin{coro}\label{Cor: Hilbert-Mumford} We equip $K$ with the trivial absolute value. 
Let $\{E_j\}_{j=1}^d$ be a finite family of finite-dimensional non-zero vector spaces over $K$, and $V$ be a non-zero quotient vector space of $E_1\otimes_K\cdots\otimes_KE_d$. Let $r$ be the rank of $V$, and, for any $j\in\{1,\ldots,d\}$, let $a_j$ the rank of $E_j$. Let \[P=\mathbb P((E_1\otimes_K\cdots\otimes_K E_d)^{\otimes r}),\] 
$\pi:P\rightarrow\Spec K$ be the canonical morphism and \[L=\mathcal O_P(a_1\cdots a_d)\otimes\pi^*(\det(E_1^{\vee})^{\otimes rb_1}\otimes\cdots\otimes\det(E_d^{\vee})^{\otimes rb_d}),\]
where 
\[\forall\,j\in\{1,\ldots,d\},\quad b_j=\frac{a_1\cdots a_d}{a_j}.\]
Then the composed surjective map \begin{equation}\label{Equ: det comme un point pluker}(E_1\otimes_K\cdots\otimes_KE_d)^{\otimes r}\longrightarrow V^{\otimes r}\longrightarrow\det(V),\end{equation} viewed as a rational point $x$ of $P$, is semistable with respect to the $\mathbb{GL}(E_1)\times_K\cdots\times_K\mathbb{GL}(E_d)$-linearised invertible sheaf $L$ if and only if, for all ultrametric norms $\norm{\ndot}_j$ on $E_j$, $j\in\{1,\ldots,d\}$, if we equip $V$ with the quotient norm of the $\varepsilon$-tensor product of $\norm{\ndot}_1,\ldots,\norm{\ndot}_d$, then one has
\begin{equation}\label{Equ: minoration of mu V}\widehat{\mu}(\overline V)\geqslant\sum_{j=1}^d\widehat{\mu}(\overline E_j).\end{equation}
\end{coro}
\begin{proof} Assume that the inequality \eqref{Equ: minoration of mu V} holds for any choice of norms $\norm{\ndot}_j$. 
By \eqref{Equ: computation of mu x varphi L}, for any one-parameter subgroup $\lambda:\mathbb{G}_{m,K}\rightarrow\mathbb{GL}(E_1)\times_K\cdots\times_K\mathbb{GL}(E_d)$, one has
\[\mu(x,\lambda,L)=a_1\cdots a_dr\Big(\widehat{\mu}(\overline V)-\sum_{j=1}^d\widehat{\mu}(\overline E_j)\Big)\geqslant 0.\]
Hence the rational point $x$ of $P$ defined by \eqref{Equ: det comme un point pluker} is semistable with respect to $L$. 

Conversely, by Remark \ref{Rem: action of gm on produt} the semi-stability of the rational points $x$ implies that the inequality \eqref{Equ: minoration of mu V} holds for any choice of ultrametric norms $\norm{\ndot}_j$ such that $\ln\norm{E_j\setminus\{0\}}_j\subseteq\mathbb Z$. As a consequence the inequality \eqref{Equ: minoration of mu V} holds for any choice of ultrametric norms $\norm{\ndot}_j$ such that $\ln\norm{E_j\setminus\{0\}}_j\subseteq\mathbb Q$. In fact, in this case there exists $n\in\mathbb N_{>0}$ such that the (finite) set $\ln\norm{E_j\setminus\{0\}}_j$ is contained in $\frac 1n\mathbb Z$ for any $j\in\{1,\ldots,d\}$. Note that the $n^{\text{th}}$ power of the function $\norm{\ndot}_j$ forms a norm on $E_j$. If we denote by $\norm{\ndot}_{V}$ the quotient norm of the $\varepsilon$-tensor product of $\norm{\ndot}_1,\ldots,\norm{\ndot}_d$, then the quotient norm of the $\varepsilon$-tensor product of $\norm{\ndot}_1^{n},\ldots,\norm{\ndot}_d^{n}$ is $\norm{\ndot}_V^n$. Note that \[\forall\,j\in\{1,\ldots,d\},\quad \ln\norm{E_j\setminus\{0\}}_j^{n}=n\ln\norm{E_j\setminus\{0\}}_j\subseteq\mathbb Z\] and hence
\[n\widehat{\mu}(V,\norm{\ndot}_V)=\widehat{\mu}(V,\norm{\ndot}_V^n)\geqslant\sum_{j=1}^d\widehat{\mu}(E_j,\norm{\ndot}_j^{n})=n\sum_{j=1}^d\widehat{\mu}(E_j,\norm{\ndot}_j).\]
Finally the general case follows from a limite procedure by using Proposition \ref{Pro:slopeinequality1}. 
\end{proof}}

\section{Estimate for the minimal slope under semi-stability assumption}

In this section, we fix an adelic curve $S=(K,(\Omega,\mathcal A,\nu),\phi)$ and assume in addition that the field $K$ is perfect. {We fix an integer $d\geqslant 2$ and we let $\{\overline E_j=(E_j,\xi_j)\}_{j=1}^d$ be a family of non-zero adelic vector bundles on $S$. Let $V$ be a quotient vector space of $E_1\otimes_K\cdots\otimes_KE_d$ and $r$ be the rank of $V$ over $K$. For any $j\in\{1,\ldots,d\}$, let $a_j$ be the rank of $E_j$. We equip $V$ with the quotient norm family of $\xi_1\otimes_{\varepsilon,\pi}\cdots\otimes_{\varepsilon,\pi}\xi_d$. Note that the quotient map $E_1\otimes_K\cdots\otimes_K E_d\rightarrow V$ induces a surjective map \[\Lambda^r(E_1\otimes_K\cdots\otimes_K E_d)\longrightarrow\Lambda^rV=\det(V).\] Consider the composed map
\[(E_1\otimes_K\cdots\otimes_K E_d)^{\otimes r}\longrightarrow \Lambda^r(E_1\otimes_K\cdots\otimes_K E_d)\longrightarrow \det(V),\]
which permits to identify $\det(V)$ as a rational point of $P=\mathbb P((E_1\otimes_K\cdots\otimes_K E_d)^{\otimes r})$. Denote by $\pi:P\rightarrow\Spec K$ the structural morphism and by $L$ the invertible sheaf \begin{equation}\label{Equ: invertible sheaf L}\mathcal O_P(a_1\cdots a_d)\otimes\pi^*(\det(E_1^{\vee})^{\otimes rb_1}\otimes\cdots\otimes\det(E_d^{\vee})^{r\otimes b_d}),\end{equation}
where 
\[\forall\,j\in\{1,\ldots,d\},\quad b_j=\frac{a_1\cdots a_d}{a_j}.\]
We equip $L$ with its natural $\mathbb{GL}(E_1)\times_K\cdots\times_K\mathbb{GL}(E_d)$-linear structure. Note that $L$ and $\mathcal O_P(a_1\cdots a_d)$ are isomorphic as invertible $\mathcal O_P$-modules, however the natural $\mathbb{GL}(E_1)\times_K\cdots\times_K\mathbb{GL}(E_d)$-linear structure on these two invertible sheaves are different.

Our purpose is to estimate $\widehat{\mu}(\overline V)$ under the additional assumption that, as a rational point of $P=\mathbb P((E_1\otimes_K\cdots\otimes_K E_d)^{\otimes r})$, the determinant line $\det(V)$ is semistable with respect to the $\mathbb{GL}(E_1)\times_K\cdots\times_K\mathbb{GL}(E_d)$-linearised  invertible sheaf $L$. 

\begin{prop}\label{Pro: estimate in semistable case} We equip $V$ with the quotient norm family $\xi_V$ of $\xi_1\otimes_{\varepsilon,\pi}\cdots\otimes_{\varepsilon,\pi}\xi_d$.
Assume that, as a rational point of the $K$-scheme $P$, $\det(V)$ is semistable with respect to the $\mathbb{GL}(E_1)\times_K\cdots\times_K\mathbb{GL}(E_d)$-linearised invertible sheaf 
$L$ defined in \eqref{Equ: invertible sheaf L}.
Then the following inequality holds:
\begin{equation}
\label{Equ: estimate in the semistable case}
\widehat{\mu}(V,\xi_V)\geqslant \sum_{j=1}^d\Big(\widehat{\mu}(E_j,\xi_j)-\nu(\Omega_\infty)\ln(\rang(E_j))\Big).
\end{equation}
\end{prop}
\begin{proof} For any integer $m\in\mathbb N_{\geqslant 1}$, let $\mathfrak S_m$ be the symmetric group of $\{1,\ldots,m\}$. Let $A=a_1\cdots a_d$.
By the first principal theorem of the classic invariant theory (see \cite[Chapter III]{Weyl} and \cite[Appendix 1]{Atiyah_Bott_Patodi}, see also \cite[Theorem 3.3]{Chen_pm}), there exist an integer $n\geqslant 1$ and an element $(\sigma_1,\ldots,\sigma_d)\in \mathfrak S_{nrA}^d$ such that the composed map
\begin{equation}\label{Equ: evaluation diagram}
\begin{gathered}
\xymatrix{\relax 
\det(V)^{\vee\otimes nA}\ar[d]\\
E_1^{\vee\otimes nrA}\otimes\cdots\otimes E_d^{\vee\otimes nrA}\ar[d]_-{\sigma_1\otimes\cdots\otimes\sigma_d}\\ E_1^{\vee\otimes nrA}\otimes \cdots\otimes E_d^{\vee\otimes nrA}\ar[d]_-{\varpi}\\
\det(E_1^{\vee})^{\otimes nrb_1}\otimes\cdots\otimes\det(E_d^{\vee})^{\otimes nrb_d}}
\end{gathered}
\end{equation}
is non-zero. Since $\xi_V$ is the quotient norm family of $\xi_1\otimes_{\varepsilon,\pi}\cdots\otimes_{\varepsilon,\pi}\xi_d$, the determinant norm family $\det(\xi_V)$ is the quotient norm family of the $\varepsilon,\pi$-tensor product norm family $\xi_1^{\otimes_{\varepsilon,\pi}r}\otimes_{\varepsilon,\pi}\cdots\otimes_{\varepsilon,\pi}\xi_d^{\otimes_{\varepsilon,\pi}r}$ by the following composed map (this is a consequence of Propositions \ref{prop:quotient:norm:linear:map} \ref{Item: successive quotient seminorm}, \ref{Pro:quotientavecpitensor} and \ref{Pro: tensor product and deteminant epsilon})
\[E_1^{\otimes r}\otimes_K \cdots\otimes_K E_d^{\otimes r}\longrightarrow \Lambda^r(E_1\otimes_K\cdots\otimes_K E_d)\longrightarrow\Lambda^r(V)=\det(V).\] 
By passing to the dual vector space, we obtain that the dual of the determinant norm family $\det(\xi_V)^\vee$ identifies with the restrict norm family of $\xi_1^{\vee\otimes_{\varepsilon}r}\otimes_{\varepsilon}\cdots\otimes_{\varepsilon}\xi_d^{\vee\otimes_\varepsilon r}$. This is a consequence of  Proposition \ref{Pro:dualitypiepsilon}, Corollary \ref{Cor: dual tensor product} and Proposition \ref{Pro:dualquotient}.

By Proposition \ref{Pro: varepsilon determinant}, the height of the $K$-linear map $\varpi$ in \eqref{Equ: evaluation diagram} is bounded from above by \[\sum_{j=1}^d nrb_j\ln(a_j!),\] where we consider the norm family $\xi_1^{\vee\otimes_{\varepsilon}nrA}\otimes_{\varepsilon}\cdots\otimes_{\varepsilon}\xi_d^{\vee\otimes_{\varepsilon}nrA}$ on $E_1^{\vee\otimes{nrA}}\otimes_K\cdots\otimes_K E_d^{\vee\otimes{nrA}}$, and the norm family $\det(\xi_1)^{\vee\otimes nrb_1}\otimes\cdots\otimes\det(\xi_d)^{\vee\otimes nrb_d}$ on $\det(E_1^\vee)^{\otimes nrb_1}\otimes_K\cdots\otimes_K\det(E_d^\vee)^{\otimes nrb_d}$.

Therefore by the slope inequality we obtain
\[-nA\hdeg(V,\xi_V)\leqslant -\sum_{j=1}^d \Big(nrb_j\hdeg(E_j,\xi_j)-\nu(\Omega_\infty)nrb_j\ln(a_j!)\Big),\]
which leads to
\[\begin{split}\widehat{\mu}(V,\xi_V)&\geqslant\sum_{j=1}^d
\Big(\widehat{\mu}(E_j,\xi_j)-\frac{1}{a_j}\ln(a_j!)\nu(\Omega_\infty)\Big)\\
&\geqslant\sum_{j=1}^d
\Big(\widehat{\mu}(E_j,\xi_j)-\ln(a_j)\nu(\Omega_\infty)\Big).
\end{split}\]

\end{proof}

\begin{rema}
Assume that the norm families $\xi_1,\ldots,\xi_d$ are Hermitian. If we equip $V$ with the quotient norm family $\xi_V^{H}$ of the orthogonal tensor product $\xi_1\otimes\cdots\otimes\xi_d$, then a similar argument as above leads to the following inequality (where we use Proposition \ref{Pro:determinant HS} to compute the height of $\varpi$)
\begin{equation}
\label{Equ: estimate in the semistable Hermitian case}
\widehat{\mu}(V,\xi_V^H)\geqslant \sum_{j=1}^d\Big(\widehat{\mu}(E_j,\xi_j)-\frac 12\nu(\Omega_\infty)\ln(\rang(E_j))\Big).
\end{equation}
\end{rema}
}

\section{An interpretation of the geometric semistability} Let $E$ be a finite-dimensional non-zero vector space over $K$ and $r$ be the rank of $E$ over $K$. 
We denote by $\mathbf{Fil}(E)$ the set of $\mathbb R$-filtrations $\mathcal E$ on $E$. Let $\mathbf{Fil}_0(E)$ be the subset of $\mathbf{Fil}(E)$ of $\mathbb R$-filtrations $\mathcal E$ such that $\hdeg(E,\|\ndot\|_{\mathcal E})=0$, where $\|\ndot\|_{\mathcal E}$ is 
the norm on $E$ associated with the $\mathbb R$-filtration $\mathcal E$
(here we consider the trivial absolute value on $K$), in other words,
\[\forall\,x\in E,\quad \|x\|_{\mathcal E}=\exp(-\sup\{t\in\mathbb R\,:\,x\in\mathcal E^t(E)\}).\]

Let $\mathcal E_1$ and $\mathcal E_2$ be two elements in $\mathbf{Fil}(E)$. 
By Proposition \ref{Pro: simultaneous orthogonal trivial valuation case}, there exists a basis $\boldsymbol{e}=\{e_i\}_{i=1}^r$ of $E$ which is orthogonal with respect to the norms $\|\ndot\|_{\mathcal E_1}$ and $\|\ndot\|_{\mathcal E_2}$ simultaneously.
We denote by $\langle\mathcal E_1,\mathcal E_2\rangle$ the number
\begin{equation}\label{Equ: produit scalaire}\frac{1}{r}\sum_{i=1}^r(-\ln\|e_i\|_{\mathcal E_1})(-\ln\|e_i\|_{\mathcal E_2}).\end{equation}
As shown by the following proposition, this number actually does not depend on the choice of the basis $\boldsymbol{e}$.

\begin{prop}Let $E$ be a finite-dimensional non-zero vector space over $K$, and $\mathcal E_1$ and $\mathcal E_2$ be $\mathbb R$-filtrations on $E$. If $\boldsymbol{e}=\{e_i\}_{i=1}^r$ is a basis of $E$ which is compatible with the $\mathbb R$-filtrations $\mathcal E_1$ and $\mathcal E_2$ simultaneously, then the following equality holds
\begin{equation}\label{Equ: coupling of filtrations}\sum_{i=1}^r\lambda_{\mathcal E_1}(e_i)\lambda_{\mathcal E_2}(e_i)=\sum_{t\in\mathbb R}t\,\hdeg(\mathrm{sq}^t_{\mathcal E_1}(E),\|\ndot\|_{\mathcal E_2,\mathrm{sq}^t_{\mathcal E_1}(E)}),\end{equation}
where $\|\ndot\|_{\mathcal E_2,\mathrm{sq}^t_{\mathcal E_1}(E)}$ denotes the subquotient norm induced by $\|\ndot\|_{\mathcal E_2}$ on the vector space $\mathrm{sq}^t_{\mathcal E_1}(E)$.
\end{prop}
{
\begin{proof}
We assume that the $\mathbb R$-filtration $\mathcal E_1$ corresponds to the flag 
\[0=V_0\subsetneq V_1\subsetneq \ldots\subsetneq V_n=V\] together with the sequence
\[\mu_1>\ldots>\mu_n.\]
Then the right hand side of the formula can be written as
\[\sum_{j=1}^n \mu_j\,\hdeg(V_j/V_{j-1},\|\ndot\|_{\mathcal E_2,j}),\]
where $\|\ndot\|_{\mathcal E_2,j}$ is the subquotient norm on $V_j/V_{j-1}$ induced by $\|\ndot\|_{\mathcal E_2}$. By Proposition \ref{Pro: property of normed vector space over trivial valued field} \ref{Item: criterion of orthgonal} the basis $\boldsymbol{e}$ is compatible with respect to the flag \[0=V_0\subsetneq V_1\subsetneq \ldots\subsetneq V_n=V.\]
By Proposition \ref{Pro:orthogonal basis in sub and quotient spaces}, the canonical image of $\boldsymbol{e}\cap (V_j\setminus V_{j-1})$ in $V_j/V_{j-1}$ forms an orthogonal basis of $(V_j/V_{j-1},\|\ndot\|_{\mathcal E_2,j})$. Moreover, for any $x\in \boldsymbol{e}\cap (V_j\setminus V_{j-1})$ one has
\[\|x\|_{\mathcal E_2}=\|[x]\|_{\mathcal E_2,j}\]
and $\|x\|_{\mathcal E_1}=\mathrm{e}^{-\mu_j}$
Therefore,
\[\begin{split}&\quad\;\sum_{i=1}^r\lambda_{\mathcal E_1}(e_i)\lambda_{\mathcal E_2}(e_i)=\sum_{j=1}^n\sum_{x\in\boldsymbol{e}\cap(V_j\setminus V_{j-1})}\lambda_{\mathcal E_1}(x)\lambda_{\mathcal E_2}(x)\\
&=
\sum_{j=1}^n\mu_j\sum_{x\in\boldsymbol{e}\cap(V_j\setminus V_{j-1})}(-\ln\|[x]\|_{\mathcal E_2,j})=\sum_{j=1}^n\mu_j\hdeg(V_j/V_{j-1},\|\ndot\|_{\mathcal E_2,j}),
\end{split}\]
where the last equality comes from Proposition \ref{Pro: property of normed vector space over trivial valued field} (\ref{Item: compute of Arakelov degree}). The equality \eqref{Equ: coupling of filtrations} is thus proved.
\end{proof}
}

We say that an $\mathbb R$-filtration $\mathcal E\in\mathbf{Fil}(E)$ is \emph{trivial}\index{trivial@trivial} if the function $Z_{\mathcal E}$ is constantly zero, or equivalently, $\langle\mathcal E,\mathcal E\rangle=0$. 

{
\begin{lemm}\label{Lem: minimal value of a concave function}
Let $V$ be a finite-dimensional non-zero vector space over $\mathbb R$, equipped with an inner product $\emptyinnprod$. Let $\{\ell_i\}_{i=1}^n$ be a finite family of linear forms on $V$, where $n\in\mathbb N$, $n\geqslant 1$. Let $\theta:V\setminus\{0\}\rightarrow\mathbb R$ be the function defined by
\[\theta(x)=\max_{i\in\{1,\ldots,n\}}\frac{\ell_i(x)}{\norm{x}},\]
where $\norm{\ndot}$ is the norm induced by the inner product $\emptyinnprod$. 
Then the function $\theta$ attains its minimal value on $V\setminus\{0\}$. Moreover, if $c$ is the minimal value of $\theta$ and if $x_0$ is a point of $V\setminus\{0\}$ minimising the function $\theta$, then for any $x\in\mathbb R^n$ one has
\begin{equation}\label{Equ: minoration of theta}\theta(x)\geqslant c\frac{\langle x,x_0\rangle}{\norm{x}\cdot\norm{x_0}}.\end{equation}
\end{lemm}
\begin{proof} 
Note that the function $\theta$ is invariant by positive dilatations, namely for any $x\in V\setminus\{0\}$ and any $\lambda>0$ one has $\theta(\lambda x)=\theta(x)$. Moreover, the function $\theta$ is clearly continuous. Hence it attains its minimal value, which is equal to $\min_{x\in V,\,\norm{x}=1}\theta(x)$.

To show the inequality \ref{Equ: minoration of theta}, we may assume without loss of generality that $\norm{x}=\norm{x_0}=1$. Note that for any $t\in[0,1]$ one has
\[\begin{split}c\norm{tx+(1-t)x_0}&\leqslant \norm{tx+(1-t)x_0}\theta(tx+(1-t)x_0)\\
&=\max_{i\in\{1,\ldots,n\}}\ell_i(tx+(1-t)x_0)\leqslant t\theta(x)+(1-t)c. \end{split}\]
Note that when $t=0$ one has
\[c\norm{tx+(1-t)x_0}=c=t\theta(x)+(1-t)c.\]
Therefore the right derivative at $t=0$ of the convex function \[t\in[0,1]\longmapsto c\norm{tx+(1-t)x_0}\] is bounded from above by $\theta(x)-c$, which leads to 
\[c\frac{t\norm{x}^2-(1-t)\norm{x}^2+(1-2t)\langle x,x_0\rangle}{\norm{tx+(1-t)x_0}}\Big|_{t=0}=c(\langle x,x_0\rangle -1)\leqslant\theta(x)-c,\]
namely $\theta(x)\geqslant c\langle x,x_0\rangle$.
\end{proof}
}

{
\begin{theo}\label{Thm: criterion of semistability}
Let $d\in\mathbb N_{\geqslant 2}$,  $\{E_j\}_{j=1}^d$ be a family of finite-dimensional non-zero vector spaces over $K$, and $V$ be a quotient vector space of $E_1\otimes_K\cdots\otimes_KE_d$. Let $r$ be the rank of $V$ over $K$. For any $j\in\{1,\ldots,d\}$, let $a_j$ be the rank of $E_j$ over $K$. Then the following conditions are equivalent.
\begin{enumerate}[label=\rm(\arabic*)]
\item\label{Item: point not semistable} The rational point $x$ of \[\pi:P=\mathbb P((E_1\otimes_K\cdots\otimes_KE_d)^{\otimes r})\longrightarrow\Spec K\] corresponding to $\det(V)$ is \emph{not} semistable with respect to the $\mathbb{GL}(E_1)\times\cdots\times\mathbb{GL}(E_d)$-linearised invertible sheaf \[L:=\mathcal O_P(a_1\cdots a_d)\otimes\pi^*(\det(E_1^\vee)^{\otimes rb_1}\otimes\cdots\det(E_d^\vee)^{\otimes rb_d}),\] where $b_j=a_1\cdots a_d/a_j$ for any $j\in\{1,\ldots,d\}$.
\item\label{Item: function theta negative} Let $\mathbf{S}$ be the subset of $\mathbf{Fil}_0(E_1)\times\cdots\times\mathbf{Fil}_0(E_d)$ consisting of the filtrations $(\mathcal F_1,\ldots,\mathcal F_d)$ which are not simultaneously trivial. For each \[{\mathcal F}=(\mathcal F_1,\ldots,\mathcal F_d)\in\mathbf{Fil}(E_1)\times\cdots\times\mathbf{Fil}(E_d),\] let $\|\ndot\|_{{\mathcal F},V}$ be the quotient norm on $V$ of the $\varepsilon$-tensor product of $\norm{\ndot}_{\mathcal F_1},\ldots,\norm{\ndot}_{\mathcal F_d}$. Then the function $\Theta:\mathbf{S}\rightarrow\mathbb R$ defined as
\[\forall\,{\mathcal F}=(\mathcal F_1,\ldots,\mathcal F_d)\in\mathbf{S},\quad \Theta({\mathcal F})=\frac{\widehat{\mu}(V,\|\ndot\|_{{\mathcal F},V})}{(\langle\mathcal F_1,\mathcal F_1\rangle+\cdots+\langle\mathcal F_d,\mathcal F_d\rangle )^{1/2}}\] 
attains its minimal value, which is negative.
\end{enumerate}
Moreover, if the above conditions are satisfied and if $\mathcal E=(\mathcal E_1,\ldots,\mathcal E_d)$ is a minimal point of the function $\Theta$, then for any ${\mathcal F}=(\mathcal F_1,\ldots,\mathcal F_d)\in\mathbf{Fil}(E_1)\times\cdots\times\mathbf{Fil}(E_d)$ one has
\begin{equation}\label{Equ: minoration de mu par minimal value}\widehat{\mu}(V,\|\ndot\|_{{\mathcal F},V})\geqslant \sum_{j=1}^d\widehat{\mu}(E_j,\|\ndot\|_{\mathcal F_j})+c\frac{\langle\mathcal E_1,\mathcal F_1\rangle+\cdots+\langle\mathcal E_d,\mathcal F_d\rangle}{(\langle\mathcal E_1,\mathcal E_1\rangle+\cdots+\langle\mathcal E_d,\mathcal E_d\rangle)^{1/2}},\end{equation}
where $c$ is the minimal value of $\Theta$.
\end{theo}
\begin{proof}
Assume that the condition \ref{Item: function theta negative} holds, then there exists $\mathcal F=(\mathcal F_1,\ldots,\mathcal F_d)\in\mathbf{S}$ such that 
\[\widehat{\mu}(V,\norm{\ndot}_{\mathcal F,V})<0=\widehat{\mu}(E_1,\norm{\ndot}_{\mathcal F_1})+\cdots+\widehat{\mu}(E_d,\norm{\ndot}_{\mathcal F_d}).\]
By Corollary \ref{Cor: Hilbert-Mumford}, the point $x$ is not semistable with respect to $L$. Conversely, if the condition \ref{Item: point not semistable} holds, then there exist \[\mathcal F'=(\mathcal F_1',\ldots,\mathcal F_d')\in\mathbf{Fil}(E_1)\times\cdots\times\mathbf{Fil}(E_d)\]  
such that 
\[\widehat{\mu}(V,\norm{\ndot}_{\mathcal F',V})<\sum_{j=1}^d\widehat{\mu}(E,\norm{\ndot}_{\mathcal F_j'}).\]
For $j\in\{1,\ldots,d\}$, let $\mathcal F_j$ be $\mathbb R$-filtrations on $E_j$ such that 
\[\norm{\ndot}_{\mathcal F_j}=\exp(\widehat{\mu}(E,\norm{\ndot}_{\mathcal F_j'}))\norm{\ndot}_{\mathcal F_j'}.\] Then one has
$\mathcal F_j\in\mathbf{Fil}_0(E_j)$ for any $j\in\{1,\ldots,d\}$.   Moreover, if we denote by $\mathcal F$ the vector $(\mathcal F_1,\ldots,\mathcal F_d)$, then one has 
\[\widehat{\mu}(V,\norm{\ndot}_{\mathcal F,V})=\widehat{\mu}(V,\norm{\ndot}_{\mathcal F',V})-\sum_{j=1}^d\widehat{\mu}(E,\norm{\ndot}_{\mathcal F_j'})<0.\]
In particular, the $\mathbb R$-filtrations $\mathcal F_j$ are not simultaneously trivial since otherwise we should have $\widehat{\mu}(V,\norm{\ndot}_{\mathcal F,V})=0$. Therefore one has $\mathcal F\in\mathbf{S}$, which implies that the function $\Theta$ takes at least a negative value.

In the following, we show that the function $\Theta$ attains its minimal value. For any $j\in\{1,\ldots,d\}$, let $\mathscr B_j$ be the set of bases of $E_j$. For $n\in\mathbb N$, let $\Delta_n$ be the vector subspace of $\mathbb R^n$ of vectors $(z_1,\ldots,z_n)$ such that $z_1+\cdots+z_n=0$.  For any 
\[\boldsymbol{e}=(\boldsymbol{e}^{(1)},\ldots,\boldsymbol{e}^{(d)})\in\mathscr B_1\times\cdots\times\mathscr B_d,\] let \[\Psi_{\boldsymbol{e}}:\Delta_{a_1}\times\cdots\times\Delta_{a_d}\rightarrow\mathbf{Fil}_0(E_1)\times\cdots\times\mathbf{Fil}_0(E_d)\] be the map sending $(\boldsymbol{y}^{(1)},\ldots,\boldsymbol{y}^{(d)})$ to the vector of $\mathbb R$-filtrations $(\mathcal F_1,\ldots,\mathcal F_d)$ such that, for any $j\in\{1,\ldots,d\}$, $\boldsymbol{e}^{(j)}=\{e^{(j)}_i\}_{i=1}^{a_j}$ forms an orthogonal basis of $\norm{\ndot}_{\mathcal F_j}$ with 
\[\boldsymbol{y}^{(j)}=(\lambda_{\mathcal F_j}(e_1^{(j)}),\ldots,\lambda_{\mathcal F_j}(e_{a_j}^{(j)})).\] 
For $j\in\{1,\ldots,d\}$, let $\boldsymbol{y}^{(j)}=(y^{(j)}_1,\ldots,y^{(j)}_{a_j})$ be an element of $\Delta_{a_j}$.
If \[\mathcal F=(\mathcal F_1,\ldots,\mathcal F_d)=\Psi_{\boldsymbol{e}}(\boldsymbol{y}^{(1)},\ldots,\boldsymbol{y}^{(d)}),\] then one has  
$\widehat{\mu}(V,\norm{\ndot}_{\mathcal F,V})$ is equal to the maximal value of 
\[\sum_{j=1}^dy_{i_{1}^{(j)}}^{(j)}+\cdots+y_{i^{(j)}_{r}}^{(j)}\]
for those $(i_1^{(j)},\ldots,i_r^{(j)})\in\{1,\ldots,a_j\}^r$ such that the image of \[(e_{i_1^{(1)}}^{(1)}\otimes\cdots\otimes e_{i_r^{(1)}}^{(1)})\otimes\cdots\otimes(e_{i_1^{(d)}}^{(d)}\otimes\cdots\otimes e_{i_r^{(d)}}^{(d)}) \] by the canonical composed surjective map
\[E_1^{\otimes r}\otimes_K\cdots\otimes_K E_d^{\otimes r}\cong (E_1\otimes_K\cdots\otimes_KE_d)^{\otimes r}\longrightarrow V^{\otimes r}\longrightarrow\det(V)\]
is non-zero. Moreover, for any $j\in\{1,\ldots,d\}$ one has
\[\langle\mathcal F_j,\mathcal F_j\rangle=\sum_{i=1}^{a_j}(y^{(j)}_i)^2.\]  
Therefore, the composition of $\Theta$ with the restriction of $\Psi_{\boldsymbol{e}}$ on \[(\Delta_{a_1}\times\ldots\times\Delta_{a_d})\setminus\{(0,\ldots,0)\}\] defines a continuous function on  $(\Delta_{a_1}\times\ldots\times\Delta_{a_d})\setminus\{(0,\ldots,0)\}$ which is invariant by dilatation by elements in $\mathbb R_{>0}$. It hence attains its minimal value. Moreover, although $\mathscr B_1\times\cdots\times\mathscr B_d$ may contain infinitely many elements, from the expression of the value $\widehat{\mu}(V,\norm{\ndot}_{\mathcal E\otimes_\varepsilon\mathcal F,V})$ as above we obtain that there are only finitely many (at most $2^{(a_1\cdots a_d)^r}$) possibility for the composed function \[\Theta\circ(\Psi_{\boldsymbol{e}}|_{\Delta_{a_1}\times\cdots\times\Delta_{a_d}\setminus\{(0,\ldots,0)\}}).\] Therefore, the function $\Theta$ attains its minimal value, which is negative since $\Theta$ takes at least one negative value.

In the following, we prove the inequality \eqref{Equ: minoration de mu par minimal value}. Let $\mathcal E=(\mathcal E_1,\ldots,\mathcal E_d)$ be an element of $\mathbf{Fil}_0(E_1)\times\cdots\times\mathbf{Fil}_0(E_d)$ which minimise the function $\Theta$. Let $\mathcal F=(\mathcal F_1,\ldots,\mathcal F_d)$ be an element of $\mathbf{Fil}(E_1)\times\cdots\times\mathbf{Fil}(E_d)$. Note that, if  $\mathcal F_j'$ is the $\mathbb R$-filtration of $E_j$ such that
\[\norm{\ndot}_{\mathcal F_j'}=\exp(\hdeg(E_j,\norm{\ndot}_{\mathcal F_j}))\norm{\ndot}_{\mathcal F_j},\]
then one has $\mathcal F_j'\in\mathbf{Fil}_0(E_j)$ and $\langle\mathcal E_j,\mathcal F_j'\rangle=\langle\mathcal E_j,\mathcal F_j\rangle$. Moreover, if we denote by $\mathcal F'$ the vector $(\mathcal F_1',\ldots,\mathcal F_d')$, then
\[\widehat{\mu}(V,\norm{\ndot}_{\mathcal F',V})=\widehat{\mu}(V,\norm{\ndot}_{\mathcal F,V})-\sum_{j=1}^d\hdeg(E_j,\mathcal F_j).\]
Therefore, to show the inequality \eqref{Equ: minoration de mu par minimal value}, it suffices to treat the case where $\mathcal F\in\mathbf{Fil}_0(E_1)\times\cdots\times\mathbf{Fil}_0(E_d)$.

By Proposition \ref{Pro: simultaneous orthogonal trivial valuation case}, for any $j\in\{1,\ldots,d\}$, there exists a basis $\boldsymbol{e}^{(j)}$ of $E$ which is orthogonal with respect to the norms $\|\ndot\|_{\mathcal E_j}$ and $\|\ndot\|_{\mathcal F_j}$ simultaneously.  Therefore the inequality \eqref{Equ: minoration de mu par minimal value} follows from Lemma \ref{Lem: minimal value of a concave function}.
\end{proof}
}

\section{Lifting and refinement of filtrations}

Let $V$ be a finite-dimensional vector space over $K$ and 
\begin{equation}\label{Equ: flag of sub spaces of V}0=V_0\subsetneq V_1\subsetneq\ldots\subsetneq V_n=V\end{equation}
be a flag of vector subspaces of $V$. Suppose given, for any $i\in\{1,\ldots,n\}$, an $\mathbb R$-filtration $\mathcal F_i$ of the sub-quotient vector space $V_i/V_{i-1}$. We will construct an $\mathbb R$-filtration on $V$ from the data of $\{\mathcal F_i\}_{i=1}^n$ as follows. For any $i\in\{1,\ldots,n\}$, let $\widetilde{\boldsymbol{e}}^{(i)}$ be a basis of $V_i/V_{i-1}$ which is orthogonal with respect to the norm $\|\ndot\|_{\mathcal F_i}$, where we consider the trivial valuation on $K$. The basis $\widetilde{\boldsymbol{e}}^{(i)}$ gives rise to a linearly independent family $\boldsymbol{e}^{(i)}$ in $V_i$. Let $\boldsymbol{e}=\bigcup_{i=1}^n\boldsymbol{e}^{(i)}$ be the (disjoint) union of $\boldsymbol{e}^{(i)}$, $i\in\{1,\ldots,n\}$. Note that $\boldsymbol{e}$ forms actually a basis of $V$ over $K$. We define an ultrametric norm $\|\ndot\|$ on $V$ such that $\boldsymbol{e}$ is an orthogonal basis under this norm and that, for any $i\in\{1,\ldots, n\}$ and any $x\in\boldsymbol{e}^{(i)}$, the norm of $x$ is $\|\widetilde x\|_{\mathcal F_i}$, where $\widetilde x$ denotes the class of $x$ in $V_i/V_{i-1}$.

\begin{rema}\begin{enumerate}[label=\rm(\arabic*)]
\item For any $i\in\{1,\ldots,n\}$, the subquotient norm on $V_i/V_{i-1}$ induced by $\|\ndot\|$ coincides with $\|\ndot\|_{\mathcal F_i}$. In particular, one has
\[\hdeg(V,\|\ndot\|)=\sum_{i=1}^n\hdeg(V_i/V_{i-1},\|\ndot\|_{\mathcal F_i}).\]
In particular, the Arakelov degree of $(V,\|\ndot\|)$ doesn't depend on the choice of $\boldsymbol{e}$.
\item Let $\mathcal F$ be the $\mathbb R$-filtration corresponding to the ultrametric norm $\norm{\ndot}$. Assume that \eqref{Equ: flag of sub spaces of V} is the flag of vector subspaces of $V$ defined by an $\mathbb R$-filtration $\mathcal G$ on $V$ and that $\mu_1>\ldots>\mu_n$ are jump points of the $\mathbb R$-filtration $\mathcal G$. Then we can compute $\langle\mathcal F,\mathcal G\rangle$ as follows:
\[\langle\mathcal F,\mathcal G\rangle=\frac{1}{\rang(V)}\sum_{i=1}^n\mu_i\hdeg(V_i/V_{i-1},\|\ndot\|_{\mathcal F_i}).\]
\end{enumerate}
\end{rema}

\begin{defi}
The $\mathbb R$-filtration on $V$ corresponding to the norm $\|\ndot\|$ constructed above is called a \emph{lifting}\index{lifting@lifting} of the family $\{\mathcal F_i\}_{i=1}^n$ (relatively to the basis $\boldsymbol{e}$). We emphasis that the lifting depends on the choice of the basis $\boldsymbol{e}$.
\end{defi}

\begin{defi}
Let $V$ be a finite-dimensional vector space over $k$ and $\mathcal F$ be an $\mathbb R$-filtration on $V$. We call \emph{refinement}\index{refinement@refinement} of $\mathcal F$ any flag
\[0=V_0\subsetneq V_1\subsetneq\ldots\subsetneq V_n=V\]
of vector subspaces of $V$ together with a non-increasing sequence 
\[t_1\geqslant\ldots\geqslant t_n\]
such that, for any $i\in\{1,\ldots,n\}$ and any $x\in V_i\setminus V_{i-1}$, one has $\|x\|_{\mathcal F}=\mathrm{e}^{-t_i}$.
\end{defi}

\begin{rema}\label{Rem: construction of refinement}
Let $V$ be a finite-dimensional vector space over $K$ and $\mathcal F$ be an $\mathbb R$-filtration on $V$. Recall that the $\mathbb R$-filtration $\mathcal F$ corresponds to a flag
\[0=W_0\subsetneq W_1\subsetneq \ldots\subsetneq W_m=V\]
together with a decreasing sequence 
\[\lambda_1>\ldots>\lambda_m.\]
To choose a refinement of $\mathcal F$ is equivalent to specify, for any $j\in\{1,\ldots,m\}$, a flag
\[0=V_j^{(0)}/W_{j-1}\subsetneq V_j^{(1)}/W_{j-1}\subsetneq\ldots\subsetneq V_j^{(n_j)}/W_{j-1}=W_j/W_{j-1}\] 
of $W_j/W_{j-1}$. The corresponding refinement is given by the flag
\[0=V_0\subsetneq V_1^{(1)}\subsetneq\ldots\subsetneq V_1^{(n_1)}\subsetneq\ldots\subsetneq V_m^{(1)}\subsetneq\ldots\subsetneq V_m^{(n_m)}=V\]
and the non-increasing sequence
\[\underbrace{\lambda_1=\cdots=\lambda_1}_{n_1\text{ copies}}>\underbrace{\lambda_2=\cdots=\lambda_2}_{n_2\text{ copies}}>\ldots>\underbrace{\lambda_m=\cdots=\lambda_m}_{n_m\text{ copies}}.\]
\end{rema}

{
\begin{prop}\label{Pro: existence of refinement of tensor filtration}
Let $d\in\mathbb N_{\geqslant 2}$, $\{E_j\}_{j=1}^d$ be a family of finite-dimensional vector spaces over $K$, and $(\mathcal F_1,\ldots,\mathcal F_d)\in\mathbf{Fil}(E_1)\times\cdots\mathbf{Fil}(E_d)$. Let $G=E_1\otimes_K\cdots\otimes_K E_d$ be the tensor product space and let $\mathcal G$ be the tensor product $\mathbb R$-filtration of $\mathcal F_1,\ldots,\mathcal F_d$ (namely the $\mathbb R$-filtration on $G$ corresponding to the $\varepsilon$-tensor product of $\|\ndot\|_{\mathcal F_1},\ldots,\|\ndot\|_{\mathcal F_d}$). Then there exists a refinement 
\[0=G_0\subsetneq G_1\subsetneq\ldots\subsetneq G_n=G,\quad t_1\geqslant\ldots\geqslant t_n\]
of the $\mathbb R$-filtration $\mathcal G$, such that, for any $i\in\{1,\ldots,n\}$, the subquotient $G_i/G_{i-1}$ is canonically isomorphic to a tensor product of subquotients of the form \[\mathrm{sq}_{\mathcal F_1}^{\lambda_{i,1}}(E_1)\otimes_K\cdots\otimes_K \mathrm{sq}_{\mathcal F_d}^{\lambda_{i,d}}(E_d)\] with $\lambda_{i,1}+\cdots+\lambda_{i,d}=t_i$.
\end{prop}
\begin{proof}
Since $\mathcal G$ is the tensor product $\mathbb R$-filtration of $\mathcal F_1,\ldots,\mathcal F_d$, one has
\[\mathcal G^t(G)=\sum_{\mu_1+\cdots+\mu_d\geqslant t}\mathcal F_1^{\mu_1}(E_1)\otimes_K\cdots\otimes_K\mathcal F_d^{\mu_d}(E_d).\]
Therefore,
\[\mathrm{sq}_{\mathcal G}^t(G)=\bigoplus_{\mu_1+\cdots+\mu_d=t}\mathrm{sq}_{\mathcal F_1}^{\mu_1}(E_1)\otimes_K\cdots\otimes_K\mathrm{sq}^{\mu_d}_{\mathcal F_d}(E_d).\]
For any $t\in\mathbb R$ there exists clearly a flag of $\mathrm{sq}_{\mathcal G}^t(G)$ whose successive subquotient is of the form $\mathrm{sq}_{\mathcal F_1}^{\mu_1}(E_1)\otimes_K\cdots\otimes_K\mathrm{sq}_{\mathcal F_d}^{\mu_d}(E_d)$ with $\mu_1+\cdots+\mu_d=t$. Hence we can construct a refinement of the $\mathbb R$-filtration $\mathcal G$ by using the construction in Remark \ref{Rem: construction of refinement}. \end{proof}

\begin{rema}\label{Rem:comparison subquotients}
We keep the notation of Proposition \ref{Pro: existence of refinement of tensor filtration} and suppose that each $E_j$ is equipped with a norm family $\xi_j$ such that $(E_j,\xi_j)$ forms an adelic vector bundle, and we equip $G$ with the $\varepsilon,\pi$-tensor product norm family $\xi_G=\xi_1\otimes_{\varepsilon,\pi}\cdots\otimes_{\varepsilon,\pi}\xi_d$. For any $t\in\mathbb R$ and any $j\in\{1,\ldots,d\}$, let $\xi_j^t$ be the induced norm families of $\xi_j$ on $\mathcal F_j^t(E_j)$ and $\xi_{j,\mathrm{sq}}^t$ be quotient norm family of $\xi_j^t$ on $\mathrm{sq}_{\mathcal F_j}^t(E_j)$. By Proposition \ref{Pro:quotientavecpitensor} and \ref{Pro:quotientr1eps}, for any  $(\mu_1,\ldots,\mu_d)\in\mathbb R^d$, the quotient norm family of $\xi_1^{\mu_1}\otimes_{\varepsilon,\pi}\cdots\otimes_{\varepsilon,\pi}\xi_d^{\mu_d}$ on $\mathrm{sq}_{\mathcal F_1}^{\mu_1}(E_1)\otimes_K\cdots\otimes_K\mathrm{sq}_{\mathcal F_d}^{\mu_d}(E_d)$ identifies with $\xi_{1,\mathrm{sq}}^{\mu_1}\otimes_{\varepsilon,\pi}\cdots\otimes_{\varepsilon,\pi}\xi_{d,\mathrm{sq}}^{\mu_d}$.

We consider a refinement
\[0=G_0\subsetneq G_1\subsetneq\ldots\subsetneq G_n=G,\quad t_1\geqslant\ldots\geqslant t_n\]
of the $\mathbb R$-filtration $\mathcal G$ such that each subquotient $G_i/G_{i-1}$ is canonically isomorphic to a tensor product of the form $\mathrm{sq}_{\mathcal F_1}^{\lambda_{i,1}}(E_1)\otimes_K\cdots\otimes_K\mathrm{sq}_{\mathcal F_d}^{\lambda_{i,d}}(E_i)$ with $\lambda_{i,1}+\cdots+\lambda_{i,d}=t_i$. Note that the canonicity of the isomorphism means that the vector space $G_i$ contains $\mathcal F^{\lambda_{i,1}}(E_1)\otimes_K\cdots\otimes_K\mathcal F^{\lambda_{i,d}}(E_d)$ and the composition
\[\mathcal F^{\lambda_{i,1}}(E_1)\otimes_K\cdots\otimes_K\mathcal F^{\lambda_{i,d}}(E_d)\longrightarrow G_i\longrightarrow G_{i}/G_{i-1}\]
of the inclusion map $\mathcal F^{\lambda_{i,1}}(E_1)\otimes_K\cdots\otimes_K\mathcal F^{\lambda_{i,d}}(E_d)\rightarrow G_i$ with the quotient map $G_i\rightarrow G_i/G_{i-1}$ induces an isomorphism  
\[\varphi_i:\mathrm{sq}_{\mathcal F_1}^{\lambda_{i,1}}(E_1)\otimes_K\cdots\otimes_K\mathrm{sq}_{\mathcal F_d}^{\lambda_{i,d}}(E_d)\longrightarrow G_i/G_{i-1}.\] We are interested in the comparison between $\xi_{1,\mathrm{sq}}^{\lambda_{i,1}}\otimes_{\varepsilon,\pi}\cdots\otimes_{\varepsilon,\pi}\xi_{d,\mathrm{sq}}^{\lambda_{i,d}}$ and the subquotient norm family of $\xi_G$ on $G_i/G_{i-1}$. By Propositions \ref{Pro: restriction and tensors}, 
the restriction of $\xi_G$ on $\mathcal F_1^{\lambda_{i,1}}(E_1)\otimes_K\cdots\otimes_K\mathcal F_d^{\lambda_{i,d}}(E_d)$ is bounded from above by $\xi_1^{\lambda_{i,1}}\otimes_{\varepsilon,\pi}\cdots\otimes_{\varepsilon,\pi}\xi_{d}^{\lambda_{i,d}}$. Therefore, for any $\omega\in\Omega$, the isomorphism $\varphi_{i,\omega}$ has an operator norm $\leqslant 1$.

Assume that the norm families $\xi_1,\ldots,\xi_d$ are Hermitian. Let $\widetilde{\xi}_G$ be the orthogonal tensor product of $\xi_1,\ldots,\xi_d$. If we equip $\mathrm{sq}_{\mathcal F_1}^{\lambda_{i,1}}(E_1)\otimes_K\cdots\otimes_K\mathrm{sq}_{\mathcal F_d}^{\lambda_{i,d}}(E_d)$  with the orthogonal product norm family $\xi_{1,\mathrm{sq}}^{\lambda_{i,1}}\otimes\cdots\otimes\xi_{d,\mathrm{sq}}^{\lambda_{i,d}}$ and $G_i/G_{i-1}$ with the subquotient norm family of $\widetilde{\xi}_G$ on $G_i/G_{i-1}$, then, for any $\omega\in\Omega$, the operator norm of $\varphi_{i,\omega}$ is bounded from above by $1$. This follows from the fact that the restriction of $\xi_G$ on $\mathcal F_1^{\lambda_{i,1}}(E_1)\otimes_K\cdots\otimes_K\mathcal F_d^{\lambda_{i,d}}(E_d)$ identifies with $\xi_1^{\lambda_{i,1}}\otimes\cdots\otimes\xi_{d}^{\lambda_{i,d}}$ (see Proposition \ref{Pro: pi tensor sub Archimedean}).
\end{rema}
}

\section{Estimation in general case}

In this section, we establish the following result.
{
\begin{theo}\label{Thm: estimation of hn filtrations for tensor product}
Let $d\in\mathbb N_{\geqslant 2}$, $\{(E_j,\xi_{j})\}_{j=1}^d$  be  a family of non-zero \emph{Hermitian} adelic vector bundles on $S$, and $V$ be a non-zero quotient vector space of $E_1\otimes_K\cdots\otimes_K E_d$. For any $j\in\{1,\ldots,d\}$ Let $\mathcal H_j$ be the Harder-Narasimhan $\mathbb R$-filtrations of $(E_j,\xi_j)$, and $\|\ndot\|_V$ be the quotient norm of the $\varepsilon$-tensor product of $\norm{\ndot}_{\mathcal H_1},\ldots,\norm{\ndot}_{\mathcal H_d}$, $\xi_V$ be the quotient norm family of the $\varepsilon,\pi$-tensor product $\xi_1\otimes_{\varepsilon,\pi}\cdots\otimes_{\varepsilon,\pi}\xi_d$, and $\widetilde{\xi}_V$ be the quotient norm family of the orthogonal tensor product $\xi_1\otimes\cdots\otimes\xi_d$. Then the following inequality holds
\begin{gather}\label{Equ: minoration de mu hat}
\widehat{\mu}(V,\xi_V)\geqslant\widehat{\mu}(V,\|\ndot\|_V)-\nu(\Omega_\infty)\sum_{j=1}^d\ln(\rang(E_j)),\\
\widehat{\mu}(V,\widetilde{\xi}_V)\geqslant\widehat{\mu}(V,\|\ndot\|_V)-\frac 12\nu(\Omega_\infty)\sum_{j=1}^d\ln(\rang(E_j)).
\end{gather}
In particular, if all adelic vector bundles $(E_j,\xi_j)$ and $\overline F$ are semistable, then one has
\begin{gather}\label{Equ: minoration mu min tensor }
\widehat{\mu}(V,\widetilde{\xi}_V)\geqslant\sum_{j=1}^d\Big(\widehat{\mu}(E_j,\xi_j)-\nu(\Omega_\infty)\ln(\rang(E_j))\Big),\\
\label{Equ: minoration mu min tensor 2}\widehat{\mu}(V,\xi_V)\geqslant\sum_{j=1}^d\Big(\widehat{\mu}(E_j,\xi_j)-\frac 12\nu(\Omega_\infty)\ln(\rang(E_j))\Big).
\end{gather}
\end{theo}
\begin{proof} For any $j\in\{1,\ldots,d\}$, let $a_j$ be the rank of $E_j$ over $K$.
We reason by induction on $A=a_1+\cdots+a_d$. The theorem is clearly true when $A=d$ (namely $\rang(E_j)=1$ for any $j$). In the following, we assume that the theorem has been proved for any family of adelic vector bundles whose ranks have the sum $<A$.

{\it Step 1:} In this step, we assume that the adelic vector bundles $(E_j,\xi_j)$ are not simultaneously semistable, or equivalently, at least one of the $\mathbb R$-filtrations $\mathcal H_j$ has more than one jump point. Let $G$ be the tensor product space $E_1\otimes_K\cdots\otimes_KE_d$ and $\mathcal G\in\mathbf{Fil}(G)$ be the tensor product of the $\mathbb R$-filtrations $\mathcal H_1,\ldots,\mathcal H_d$, which corresponds to the $\varepsilon$-tensor product of the norms $\norm{\ndot}_{\mathcal F_1},\ldots\norm{\ndot}_{\mathcal F_d}$. We choose a refinement
\[0=G_0\subsetneq G_1\subsetneq\ldots\subsetneq G_n=G,\quad t_1\geqslant\ldots\geqslant t_n\]
of the $\mathbb R$-filtration $\mathcal G$ such that, for any $i\in\{1,\ldots,n\}$, the subquotient $G_i/G_{i-1}$ is canonically isomorphic to a tensor product of the form \[\mathrm{sq}^{\lambda_{i,1}}_{\mathcal H_1}(E_1)\otimes_K\otimes\cdots\otimes_K\mathrm{sq}^{\lambda_{i,d}}_{\mathcal H_d}(E_d)\] with $\lambda_{i,1}+\cdots+\lambda_{i,d}=t_i$ (see Proposition \ref{Pro: existence of refinement of tensor filtration}). The assumption that at least one of the $\mathbb R$-filtrations $\mathcal F_{j}$ has more than one jump point implies that \[\rang(\mathrm{sq}^{\lambda_{i,1}}_{\mathcal H_1}(E_1))+\cdots+\rang(\mathrm{sq}^{\lambda_{i,d}}_{\mathcal H_d}(E_d))<a_1+\cdots+a_d.\]

For any $i\in\{1,\ldots,n\}$ and any $j\in\{1,\ldots,d\}$, denoted by $\xi_{j,\mathrm{sq}}^{\lambda_{i,j}}$ the subquotient norm families of $\xi_j$ on $\mathrm{sq}^{\lambda_{i,j}}_{\mathcal H_j}(E_j)$. For any $i\in\{1,\ldots,n\}$, let $\xi_{G_i/G_{i-1}}$ be the subquotient norm family of $\xi_G$ on $G_i/G_{i-1}$,  $\xi_{V_i/V_{i-1}}'$ be the quotient norm family of $\xi_{G_i/G_{i-1}}$ on $V_i/V_{i-1}$, $\xi_{V_i/V_{i-1}}$ be the subquotient norm family of $\xi_V$ on $V_i/V_{i-1}$, and ${\xi}_{V_i/V_{i-1}}''$ be the quotient norm family of $\xi_{1,\mathrm{sq}}^{\lambda_{i,1}}\otimes_{\varepsilon,\pi}\cdots\otimes_{\varepsilon,\pi}\xi_{d,\mathrm{sq}}^{\lambda_{i,d}}$
on $V_i/V_{i-1}$, where we identify $G_i/G_{i-1}$ with \[\mathrm{sq}^{\lambda_{i,1}}_{\mathcal H_1}(E_1)\otimes_K\cdots\otimes_K\mathrm{sq}^{\lambda_{i,d}}_{\mathcal H_d}(E_d).\] By Proposition \ref{prop:quotient:norm:linear:map} \ref{Item: diagram of quotient seminorm}, one has 
\begin{equation}\label{Equ: minormation degree Vi}\hdeg(V_i/V_{i-1},\xi_{V_i/V_{i-1}})\geqslant \hdeg(V_i/V_{i-1},\xi_{V_i/V_{i-1}}').\end{equation}
By Remark \ref{Rem:comparison subquotients}, one has
\begin{equation}\label{Equ: minormation degree Vi2}\hdeg(V_i/V_{i-1},{\xi}_{V_i/V_{i-1}}')\geqslant\hdeg(V_i/V_{i-1},{\xi}_{V_i/V_{i-1}}'').\end{equation}
Moreover, for any $(i,j)\in\{1,\ldots,n\}\times\{1,\ldots,d\}$, the Hermitian adelic vector bundle $(\mathrm{sq}^{\lambda_{i,j}}_{\mathcal H_j}(E_j),\xi_{j,\mathrm{sq}}^{\lambda_{i,j}})$ is semistable of slope $\lambda_{i,j}$. Therefore, by the induction hypothesis one has
\[\begin{split}&\quad\;\hdeg(V_i/V_{i-1},\xi_{V_i/V_{i-1}})\geqslant\hdeg(V_i/V_{i-1},\xi_{V_i/V_{i-1}}'')\\
&\geqslant\rang(V_i/V_{i-1})
\sum_{j=1}^d\big(\lambda_{i,j}-\nu(\Omega_\infty)\ln(a_j)\big).
\end{split}\]
Taking the sum with respect to $i\in\{1,\ldots,n\}$, we obtain 
\[\begin{split}\hdeg(V,\xi_V)&=\sum_{i=1}^n\hdeg(V_i/V_{i-1},\xi_{V_i/V_{i-1}})\\
&\geqslant\sum_{i=1}^n\rang(V_i/V_{i-1})\sum_{j=1}^d\big(\lambda_{i,j}-\nu(\Omega_\infty)\ln(a_j)\big)\\
&=\sum_{i=1}^n\rang(V_i/V_{i-1})t_i-\rang(V)\nu(\Omega_\infty)\sum_{j=1}^d\ln(a_j)\\
&=\hdeg(V,\|\ndot\|_V)-\rang(V)\nu(\Omega_\infty)\sum_{j=1}^d\ln(a_j),
\end{split}\]
which leads to
\[\widehat{\mu}(V,\xi_V)\geqslant\widehat{\mu}(V,\|\ndot\|_V)-\nu(\Omega_\infty)\sum_{j=1}^d\ln(a_j).\]

Similarly, if we denote by $\widetilde{\xi}_{V_i/V_{i-1}}$ the subquotient norm family of $\widetilde{\xi}_V$, then the induction hypothesis gives 
\[\hdeg(V_i/V_{i-1},\widetilde{\xi}_{V_i/V_{i-1}})\geqslant\rang(V_i/V_{i-1})\sum_{j=1}^d\Big(\lambda_{i,j}-\frac 12\nu(\Omega_\infty)\ln(a_j)\Big),\]
which leads to 
\[\widehat{\mu}(V,\xi_V)\geqslant\widehat{\mu}(V,\|\ndot\|_V)-\frac 12\nu(\Omega_\infty)\sum_{j=1}^d\ln(a_j).\]

{\it Step 2:} In this step, we assume that all adelic vector bundles $(E_j,\xi_j)$ are semistable. Note that the Harder-Narasimhan $\mathbb R$-filtrations of $(E_j,\xi_j)$ has then only one jump point. Therefore, it suffices to prove \eqref{Equ: minoration mu min tensor } and \eqref{Equ: minoration mu min tensor 2}. Note that the case where $\det(V)$ is semistable as a rational point of \[P:=\mathbb P((E_1\otimes_K\cdots\otimes_KE_d)^{\otimes\rang(V)})\] has been proved in Proposition \ref{Pro: estimate in semistable case}. In the following, we assume that $\det(V)$ is not semistable as a rational point of $P$.

For each \[\mathcal F=(\mathcal F_1,\ldots,\mathcal F_d)\in\mathbf{Fil}(E_1)\times\cdots\times\mathbf{Fil}(E_d),\] let $\|\ndot\|_{\mathcal F,V}$ be the quotient norm on $V$ of the $\varepsilon$-tensor product of $\|\ndot\|_{\mathcal F_1},\ldots,\norm{\ndot}_{\mathcal F_d}$. Let $\mathbf{S}$ be the subset of $\mathbf{Fil}_0(E_1)\times\cdots\times\mathbf{Fil}_0(E_d)$ consisting of vectors $(\mathcal F_1,\ldots,\mathcal F_d)$ such that the filtrations $\mathcal F_1,\ldots,\mathcal F_d$ are not simultaneously trivial.  Then, by Theorem \ref{Thm: criterion of semistability}, the function $\Theta:\mathbf{S}\rightarrow\mathbb R$
\[\forall\,\mathcal F=(\mathcal F_1,\ldots,\mathcal F_d)\in\mathbf{S},\quad \Theta(\mathcal F):=\frac{\widehat{\mu}(V,\|\ndot\|_{\mathcal F,V})}{(\langle\mathcal F_1,\mathcal F_1\rangle+\cdots+\langle\mathcal F_d,\mathcal F_d\rangle )^{1/2}}\] 
attains its minimal value $c$, which is negative. In the following, we denote by $\mathcal E=(\mathcal E_1,\ldots,\mathcal E_d)$ a minimal point of the function $\Theta$. Then, for any $\mathcal F=(\mathcal F_1,\ldots,\mathcal F_d)\in\mathbf{Fil}(E_1)\times\cdots\times\mathbf{Fil}(E_d)$, one has
\[\widehat{\mu}(V,\|\ndot\|_{\mathcal F,V})\geqslant \sum_{j=1}^d\widehat{\mu}(E_j,\|\ndot\|_{\mathcal F_j})+c\frac{\langle\mathcal E_1,\mathcal F_1\rangle+\cdots+\langle\mathcal E_d,\mathcal F_d\rangle}{(\langle\mathcal E_1,\mathcal E_1\rangle+\cdots+\langle\mathcal E_d,\mathcal E_d\rangle)^{1/2}}.\]

In the following, for each $j\in\{1,\ldots,d\}$, we denote by $\mathcal F_j$ the $\mathbb R$-filtration on $E$ which induces on each subquotient $\mathrm{sq}^t_{\mathcal E_j}(E_j)$ the Harder-Narasimhan filtration of this vector space equipped with the subquotient norm family $\xi_{j,\mathrm{sq},\mathcal E_j}^t$ of $\xi_j$. By \eqref{Equ: degree sum}, for any $t\in\mathbb R$, if we denote by $\norm{\ndot}_{\mathcal F_j,\mathrm{sq},t}$ the subquotient norm of $\norm{\ndot}_{\mathcal F_j}$ on $\mathrm{sq}_{\mathcal E_j}^t(E_j)$, then one has
\[\hdeg(\mathrm{sq}^t_{\mathcal E_j}(E_j),\|\ndot\|_{\mathcal F_j,\mathrm{sq},t})=\hdeg(\mathrm{sq}^t_{\mathcal E_j}(E_j),\xi_{j,\mathrm{sq},\mathcal E_j}^t).\]
Taking the sum with respect to $t$, by Proposition \ref{Pro:suiteexactedeg} and the assumption that $\xi_j$ is Hermitian we obtain that
\[\hdeg(E_j,\|\ndot\|_{\mathcal E_j})=\hdeg(E_j,\xi_j).\]
Moreover, by \eqref{Equ: degree sum} one has
\[\langle\mathcal E_j,\mathcal F_j\rangle=\frac{1}{a_j}\sum_{t\in\mathbb R}t\sum_{i=1}^{r_j(t)}\widehat{\mu}_i(\mathrm{sq}^t_{\mathcal E_j}(E_j),\xi_{j,\mathrm{sq},\mathcal E_j}^t)=\frac{1}{a_j}\sum_{t\in\mathbb R}t\,\hdeg(\mathrm{sq}^t_{\mathcal E_j}(E_j),\xi_{j,\mathrm{sq},\mathcal E_j}^t),\]
where $r_j(t)=\rang_K(\mathrm{sq}_{\mathcal E_j}^t(E_j))$, and the second equality comes from \eqref{Equ: degree sum}. For any $j\in\{1,\ldots,d\}$ and any $u\in\mathbb R$, let 
\[\Psi_j(u)=\sum_{t< u}\hdeg(\mathrm{sq}^t_{\mathcal E_j}(E_j),\xi_{j,\mathrm{sq},\mathcal E_j}^t)=\hdeg(\overline{E/\mathcal E_j^u(E_j)}),\]
where we consider the quotient norm family on $E/\mathcal E_j^u(E_j)$.
Since $(E_j,\xi_j)$ is semistable, one has
\[\widehat{\mu}(\overline{E/\mathcal E_j^u(E_j)})\geqslant \widehat{\mu}_{\min}(E_j,\xi_j)=\widehat{\mu}(E_j,\xi_j)\]
and hence 
\[\Psi_j(u)\geqslant\widehat{\mu}(E_j,\xi_j)\rang(E/\mathcal E_j^u(E_j)).  \]
By Abel's summation formula we obtain
\[\langle\mathcal E_j,\mathcal F_j\rangle=\frac{1}{a_j}\int_{\mathbb R}t\,\mathrm{d}\Psi_j(t)=M_j\frac{\Psi_j(M_j)}{a_j}-\frac{1}{a_j}\int_{-\infty}^{M_j}\Psi_j(t)\,\mathrm{d}t,\]
where $M_j$ is a sufficiently positive number such that $\mathcal E_j^{M_j}(E_j)=\{0\}$. Therefore one has
\[\begin{split}\langle\mathcal E_j,\mathcal F_j\rangle&\leqslant M_j\frac{\hdeg(E_j,\xi_j)}{a_j}-\frac{1}{a_j}\widehat{\mu}(E_j,\xi_j)\int_{-\infty}^{M_j}\rang(E_j/\mathcal E_j^t(E_j))\,\mathrm{d}t\\
&=\frac{\widehat{\mu}(E_j,\xi_j)}{a_j}\int_{-\infty}^{M_j} t\,\mathrm{d}\rang(E_j/\mathcal E_j^t(E_j))=\frac{\widehat{\mu}(E_j,\xi_j)}{a_j}\hdeg(E_j,\|\ndot\|_{\mathcal E_j})=0.\end{split}\]
Therefore we obtain 
\[\widehat{\mu}(V,\|\ndot\|_{\mathcal F,V})\geqslant\sum_{j=1}^d\widehat{\mu}(E_j,\xi_j).\]

It remains to compare $\widehat{\mu}(V,\|\ndot\|_{\mathcal F,V})$ with the slopes of $(V,\xi_V)$ and $(V,\widetilde{\xi}_V)$. We choose a refinement
\[0=G_0\subsetneq G_1\subsetneq\ldots\subsetneq G_n=E_1\otimes_K\cdots\otimes_KE_d,\quad t_1\geqslant\ldots\geqslant t_n\]
of the $\mathbb R$-filtration $\mathcal E_1\otimes\cdots\otimes\mathcal E_d$ such that, for any $i\in\{1,\ldots,n\}$, the subquotient $G_i/G_{i-1}$ is canonically isomorphic to a tensor product of the form \[\mathrm{sq}_{\mathcal E_1}^{\lambda_{i,1}}(E_1)\otimes_K\cdots\otimes_K\mathrm{sq}_{\mathcal E_d}^{\lambda_{i,d}}(E_d)\] with $\lambda_{i,1}+\cdots+\lambda_{i,d}=t_i$. For any $i\in\{1,\ldots,n\}$, let $\|\ndot\|_{\mathcal F,G_i/G_{i-1}}$ be the subquotient norm on $G_i/G_{i-1}$ of the $\varepsilon$-tensor product of $\norm{\ndot}_{\mathcal F_1},\ldots,\norm{\ndot}_{\mathcal F_d}$. By the construction of $\mathcal F_1,\ldots,\mathcal F_d$, the subquotient norm  $\|\ndot\|_{\mathcal F,G_i/G_{i-1}}$ on $G_i/G_{i-1}$ corresponds to the tensor product of the Harder-Narasimhan $\mathbb R$-filtrations of $(\mathrm{sq}_{\mathcal E_1}^{\lambda_{i,1}}(E_1),\xi_{1,\mathrm{sq},\mathcal E_1}^t),\ldots,(\mathrm{sq}_{\mathcal E_d}^{\lambda_{i,d}}(E_d),\xi_{d,\mathrm{sq},\mathcal E_d}^t)$. By the induction hypothesis, and the same argument showing \eqref{Equ: minormation degree Vi} and \eqref{Equ: minormation degree Vi2}, we obtain
\begin{equation*}\hdeg({V_i/V_{i-1}},\xi_{V_i/V_{i-1}})\geqslant\hdeg(V_i/V_{i-1},\|\ndot\|_{\mathcal F,V_i/V_{i-1}})-\rang(V_i/V_{i-1})\nu(\Omega_\infty)\sum_{j=1}^d\ln(a_j),\end{equation*}
where $\xi_{V_i/V_{i-1}}$ is the subquotient norm family of $\xi_V$ on $V_i/V_{i-1}$, $\|\ndot\|_{\mathcal F,V_i/V_{i-1}}$ is the quotient norm of $\|\ndot\|_{\mathcal F,G_i/G_{i-1}}$, which identifies with the subquotient norm of $\|\ndot\|_{F,V}$ since the flag $0=G_0\subsetneq G_1\subsetneq\ldots\subsetneq G_n$ is compatible with the $\mathbb R$-filtration $\mathcal F_1\otimes\cdots\otimes\mathcal F_d$. Taking the sum of the above formula with respect to $i\in\{1,\ldots,n\}$, we obtain
\[\hdeg(V,\xi_V)\geqslant\hdeg(V,\|\ndot\|_{\mathcal F,V})-\rang(V)\nu(\Omega_\infty)\sum_{j=1}^d\ln(a_j),\]
which leads to 
\[\widehat{\mu}(V,\xi_V)\geqslant\widehat{\mu}(V,\|\ndot\|_{\mathcal F,V})-\nu(\Omega_\infty)\sum_{j=1}^d\ln(a_j).\]
Similarly, one has 
\[\widehat{\mu}(V,\widetilde{\xi}_V)\geqslant\widehat{\mu}(V,\norm{\ndot}_{\mathcal F,V})-\frac 12\nu(\Omega_\infty)\sum_{j=1}^d\ln(a_j).\]
The theorem is thus proved.
\end{proof}

\begin{coro}\label{Cor: tensorial minimal slope property}
Let $d\in\mathbb N_{\geqslant 2}$, $\{(E,\xi_j)\}_{j=1}^d$ be a family of  adelic vector bundles on $S$ and $V$ be a non-zero quotient vector space of $E_1\otimes_K\cdots\otimes_KE_d$. Let $\xi_V$ be the quotient norm families of $\xi_1\otimes_{\varepsilon,\pi}\cdots\otimes_{\varepsilon,\pi}\xi_d$.
Then one has
\begin{equation}\label{Equ: minoration de mu min dans le cas non hermitien}\widehat{\mu}_{\min}(V,\xi_V)\geqslant\sum_{j=1}^d\Big(\widehat{\mu}_{\min}(E_j,\xi_j)-\frac 32\nu(\Omega_\infty)\ln(\rang(E_j))\Big).\end{equation}
If all norm families $\xi_1,\ldots,\xi_d$ are Hermitian, then one has
\begin{gather}\label{Equ: minoration de mu min dans le cas hermitien}
\widehat{\mu}_{\min}(V,\xi_V)\geqslant\sum_{j=1}^d\Big(\widehat{\mu}_{\min}(E_j,\xi_j)-\nu(\Omega_\infty)\ln(\rang(E_j))\Big),\\\label{Equ: minoration de mu min dans le cas hermitien orth}
\widehat{\mu}_{\min}(V,\widetilde\xi_V)\geqslant\sum_{j=1}^d\Big(\widehat{\mu}_{\min}(E_j,\xi_j)-\frac 12\nu(\Omega_\infty)\ln(\rang(E_j))\Big),
\end{gather}
where $\widetilde{\xi}_V$ is the quotient norm family of the orthogonal tensor product $\xi_1\otimes\cdots\otimes\xi_d$.
\end{coro}
\begin{proof}
We begin with the proof of the Hermitian case. To establish \eqref{Equ: minoration de mu min dans le cas hermitien} it suffices to prove weaker inequalities
\begin{gather}\label{Equ: minoration de mu W par la somme de mu min}\widehat{\mu}(V,\xi_V)\geqslant\sum_{j=1}^d\Big(\widehat{\mu}_{\min}(E_j,\xi_j)-\nu(\Omega_\infty)\ln(\rang(E_j))\Big),\\
\label{Equ: minoration de mu W par la somme de mu min orth}
\widehat{\mu}(V,\widetilde\xi_V)\geqslant\sum_{j=1}^d\Big(\widehat{\mu}_{\min}(E_j,\xi_j)-\frac 12\nu(\Omega_\infty)\ln(\rang(E_j))\Big)
\end{gather}
for all non-zero quotient vector space $V$ of $E_1\otimes_K\cdots\otimes_KE_d$. For any $j\in\{1,\ldots,d\}$, let $\mathcal H_j$ be the Harder-Narasimhan $\mathbb R$-filtration of $(E_j,\xi_j)$. Let $\|\ndot\|_V$ be the quotient norm of the $\varepsilon$-tensor product of $\norm{\ndot}_{\mathcal H_1},\ldots,\norm{\ndot}_{\mathcal H_d}$. By Theorem \ref{Thm: estimation of hn filtrations for tensor product}, one has
\begin{gather}
\label{Equ: minoration de mu W par la somme de mu min2}\widehat{\mu}(V,\xi_V)\geqslant\widehat{\mu}(V,\|\ndot\|_V)-\nu(\Omega_\infty)\sum_{j=1}^d\ln(\rang(E_j)),\\
\label{Equ: minoration de mu W par la somme de mu min2 orth}
\widehat{\mu}(V,\widetilde\xi_V)\geqslant\widehat{\mu}(V,\|\ndot\|_V)-\frac 12\nu(\Omega_\infty)\sum_{j=1}^d\ln(\rang(E_j)).\end{gather}
Moreover, since $\|\ndot\|_V$ is the quotient norm of the $\varepsilon$-tensor product of $\|\ndot\|_{\mathcal H_1},\ldots,\norm{\ndot}_{\mathcal H_d}$, one has (see Remark \ref{Rem: tensor product filtrations}) \[\widehat{\mu}(V,\|\ndot\|_V)\geqslant\sum_{j=1}^d\widehat{\mu}_{\min}(E_j,\xi_j).\] Therefore \eqref{Equ: minoration de mu W par la somme de mu min} follows from \eqref{Equ: minoration de mu W par la somme de mu min2} and \eqref{Equ: minoration de mu W par la somme de mu min orth} follows from \eqref{Equ: minoration de mu W par la somme de mu min2 orth}.

In the following, we proceed with the proof of \eqref{Equ: minoration de mu min dans le cas non hermitien} in the general (non-necessarily Hermitian) case. Note that one has 
\[\xi_1^{\vee\vee}\otimes_{\varepsilon,\pi}\cdots\otimes_{\varepsilon,\pi}\xi_d^{\vee\vee}=\xi_1\otimes_{\varepsilon,\pi}\cdots\otimes_{\varepsilon,\pi}\xi_d.\]
Moreover, since $\xi_j^{\vee\vee}\leqslant\xi_j$, one has
\[\widehat{\mu}_{\min}(E_j,\xi_j^{\vee\vee})\geqslant\widehat{\mu}_{\min}(E_j,\xi_j)\]
for any $j\in\{1,\ldots,d\}$.
Therefore, by replacing $\xi_j$ by $\xi_j^{\vee\vee}$ we may suppose without loss of generality that all $\xi_j$ are non-Archimedean on $\Omega\setminus\Omega_\infty$.

Assume that $\xi_j$ is of the form $\{\norm{\ndot}_{j,\omega}\}_{j=1}^d$. By Theorem \ref{Thm: Hermitian approximation via measurable selection}, for any $\epsilon>0$ and any $j\in\{1,\ldots,d\}$ there exist measurable Hermitian norm families $\xi_j^H=\{\norm{\ndot}_{j,\omega}^H\}_{\omega\in\Omega}$ of $E_j$ such that $\norm{\ndot}_{j,\omega}^H=\norm{\ndot}_{j,\omega}$ for any $\omega\in\Omega\setminus\Omega_\infty$and 
\[ \norm{\ndot}_{j,\omega}\leqslant\norm{\ndot}_{j,\omega}^H\leqslant (\rang(E_j)+\epsilon)^{1/2}\norm{\ndot}_{j,\omega}\]
for any $\omega\in\Omega_\infty$. 
By the slope inequality (see Proposition \ref{Pro:inegalitdepente}) one has 
\begin{gather}\label{Equ: minoration Lowner}
\begin{gathered}\widehat{\mu}_{\min}(E_j,\xi_j^H)\geqslant\widehat{\mu}(E_j,\xi_j)-\frac 12\nu(\Omega_\infty)\ln(\rang(E_j)+\epsilon).
\end{gathered}\end{gather}
Moreover, if we denote by $\xi_V'$ the quotient norm family of $\xi_1^H\otimes_{\varepsilon,\pi}\cdots\otimes_{\varepsilon,\pi}\xi_d^H$ on $V$, one has
\begin{equation}\label{Equ: minoration Lowner quotient}\widehat{\mu}_{\min}(V,\xi_V')\leqslant\widehat{\mu}_{\min}(V,\xi_V)\end{equation}
by the slope inequality.
Applying the Hermitian case of the corollary to $(E_j,\xi_{j}^H)$ ($j\in\{1,\ldots,d\}$) and $(V,\xi_V')$, we obtain
\[\widehat{\mu}_{\min}(V,\xi_V')\geqslant\sum_{j=1}^d\Big(\widehat{\mu}_{\min}(E_j,\xi_j^H)-\nu(\Omega_\infty)\ln(\rang(E_j))\Big).\]
Combining this inequality with \eqref{Equ: minoration Lowner} and \eqref{Equ: minoration Lowner quotient}, by passing to limite when $\epsilon$ tend to $0+$ we obtain \eqref{Equ: minoration de mu min dans le cas non hermitien}. The corollary is thus proved.
\end{proof}

\begin{rema}
In the case where the adelic curve $S$ comes from an arithmetic curve. The inequality \eqref{Equ: minoration de mu min dans le cas hermitien orth} recovers essentially the second inequality of \cite[Corollary 5.4]{Gaudron_Remond13}, which strengthen \cite[Theorem 1]{Chen_pm}. From the methodological point of view, the arguments in this chapter rely on the geometric invariant theory without using the theorem of successive minima of Zhang, which is a key argument in \cite{Gaudron_Remond13,Bost_Chen} (see \cite[Theorem 5.2]{Zhang95}, see also \cite[\S3]{Gaudron_Remond13}).
\end{rema}
}

%% file: ch6_2019_03_23.tex

\chapter{Adelic line bundles on arithmetic varieties}

\IfChapVersion
\ChapVersion{Version of Chapter 6 : \\ \StrSubstitute{\DateChapSix}{_}{\_}}
\fi

In this chapter, we fix a proper adelic curve $S=(K,(\Omega,\mathcal A,\nu),\phi)$. We assume that, either the $\sigma$-algebra $\mathcal A$ is discrete, or there exists a countable subfield $K_0$ of $K$ which is dense in the completion $K_{\omega}$ of $K$ with respect to any $\omega \in \Omega$.

\section{Metrised line bundles on an arithmetic variety}
Let $X$ be a projective scheme over $\Spec K$ and $L$ be an invertible $\mathcal O_X$-module. For any $\omega\in\Omega$, we let $X_\omega$ be the fibre product $X\times_{\Spec K}\Spec K_\omega$ (recall that $K_\omega$ is the completion of $K$ with respect to $|\ndot|_\omega$) and $L_\omega$ be the pull-back of $L$ by the canonical projection morphism $X_\omega\rightarrow X$. By \emph{metric family}\index{metric family} on $L$, we refer to a family of continuous metrics $\varphi=\{\varphi_{\omega}\}_{\omega\in\Omega}$, where $\varphi_{\omega}$ is a continuous metric on $L_\omega$. If $\varphi=\{\varphi_\omega\}_{\omega\in\Omega}$ and $\varphi'=\{\varphi_\omega'\}_{\omega\in\Omega}$ are two metric families on $L$, the \emph{local distance}\index{local distance}\index{metric family!local distance} of $\varphi$ and $\varphi'$ at $\omega\in\Omega$ is defined as (see Definition \ref{Def:localdistancemetric})
\[d_\omega(\varphi,\varphi'):=d(\varphi_\omega,\varphi_\omega').\]
The global distance between $\varphi$ and $\varphi'$ is defined as the upper integral
\begin{equation}\label{eqn:def:global:distance}\operatorname{dist}(\varphi,\varphi'):=\upint_{\Omega}d_\omega(\varphi,\varphi')\,\nu(\mathrm{d}\omega).\end{equation}

If $\varphi=\{\varphi_\omega\}_{\omega\in\Omega}$ is a metric family on $L$, then the dual metrics $\{-\varphi_{\omega}\}_{\omega\in\Omega}$ form a metric family on $L^\vee$, denoted by $-\varphi$. If $L$ and $L'$ are invertible $\mathcal O_X$-modules, and $\varphi=\{\varphi_\omega\}_{\omega\in\Omega}$ and $\varphi'=\{\varphi'_\omega\}_{\omega\in\Omega}$ are metric families on $L$ and $L'$ respectively, then $\{\varphi_\omega+\varphi'_\omega\}_{\omega\in\Omega}$ is a metric family on $L\otimes L'$, denoted by $\varphi+\varphi'$. The metric family $\varphi+(-\varphi')$ on $L\otimes L'{}^\vee$ is also denoted by $\varphi-\varphi'$. Similarly,  for any integer $n\geqslant 0$, $\{n\varphi_\omega\}_{\omega\in\Omega}$ is a metric family on $L^{\otimes n}$, denoted by $n\varphi$.

{
\begin{defi}\label{Def: pull-back metric family}
Let $Y$ and $X$ be projective schemes over $\Spec K$ and $f:Y\rightarrow X$ be a projective $K$-morphism. Let $L$ be an invertible $\mathcal O_X$-module equipped with a metric family $\varphi=\{\varphi_\omega\}_{\omega\in\Omega}$. We denote by $f^*(\varphi)$ the metric family $\{f_\omega^*(\varphi_\omega)\}_{\omega\in\Omega}$ on $f^*(L)$, where for any $\omega\in\Omega$, $f_\omega:Y_\omega\rightarrow X_\omega$ is the $K_\omega$-morphisme induced by $f$, and $f_\omega^*(\varphi_\omega)$ is defined in Definition \ref{Def:pull-back}. The norm family $f^*(\varphi)$ is called the \emph{pull-back}\index{pull-back}\index{metric family!pull-back} of $\varphi$ by $f$. In the particular case where $Y$ is a closed subscheme of $X$ and $f$ is the canonical immersion, the norm family $f^*(\varphi)$ is also denoted by $\varphi|_Y$ and called the \emph{restriction}\index{restriction}\index{metric family!restriction} of $\varphi$ on $Y$.
\end{defi}

The following properties are straightforward from the definition.

\begin{prop}\label{Pro: pull back by point}
Let $X$ be a projective scheme over $\Spec K$ and $f:Y\rightarrow X$ be a projective morphism of $K$-schemes. 
\begin{enumerate}[label=\rm(\arabic*)]
\item If $L_1$ and $L_2$ are two invertible $\mathcal O_X$-modules, and $\varphi_1$ and $\varphi_2$ are metric families on $L_1$ and $L_2$ respectively, then one has $f^*(\varphi_1+\varphi_2)=f^*(\varphi_1)+ f^*(\varphi_2)$ on $f^*(L_1\otimes L_2)\cong f^*(L_1) \otimes f^*(L_2)$.
\item For any invertible $\mathcal O_X$-module $L$ and any metric family $\varphi$ on $L$, one has $f^*(-\varphi)=-f^*(\varphi)$ on $f^*(L^{\vee})\cong f^*(L)^\vee$.
\end{enumerate}
\end{prop} 

\begin{rema}\label{Rem: metric family on one point}
Let us consider the particular case where $X$ is the spectrum of a finite extension $K'$ of the field $K$. For any $\omega\in\Omega$, the Berkovich space of $X_{\omega}$ 
identifies with the discrete set of absolute values on $K'$ extending $|\ndot|_\omega$ on $K$. Moreover, any invertible $\mathcal O_X$-module $L$ could be considered as a vector space of rank $1$ over $K'$, and any metric family on $L$ is just a norm family with respect to the adelic curve $S\otimes_KK'$ (cf. Definition \ref{Equ:finiteextension}) if we consider $L$ as a vector space over $K'$. 
\end{rema}
}

\subsection{Quotient metric families}

Let $E$ be a finite-dimensional vector space over $K$ and $\xi=\{\|\ndot\|_{\omega}\}_{\omega\in\Omega}$ be a norm family on $E$. Let $f:X\rightarrow\Spec K$ be a projective $K$-scheme and $L$ be an invertible $\mathcal O_X$-module. Suppose given a surjective homomorphism $\beta:f^*(E)\rightarrow L$. For any $\omega\in\Omega$, the morphism $f:X\rightarrow\Spec K$ induces by base change a morphism $f_\omega$ from $X_\omega:=X\times_{\Spec K}\Spec K_\omega$ to $\Spec K_\omega$. We denote by $L_\omega$ the pull-back of $L$ on $X_\omega$. The homomorphism $\beta$ induces a surjective homomorphism $\beta_\omega:f_\omega^*(E)\rightarrow L_\omega$. Therefore, the norm $\|\ndot\|_\omega$ induces a quotient metric $\varphi_\omega$ on $L_\omega$ (see Definition \ref{Def: Quotient metric}). The family $\varphi=\{\varphi_\omega\}_{\omega\in\Omega}$ is called the \emph{quotient metric family induced by $(E,\xi)$ and $\beta$}\index{quotient metric family}\index{metric family!quotient ---}.

Let $\pmb{e} = \{e_i\}_{i=1}^r$ be a basis of $E$.
For each $\omega \in \Omega$, 
let $\|\ndot\|_{\pmb{e},\omega}$ be the norm on $E_{\omega} := E \otimes_K K_{\omega}$ given by
\[
\forall\, a_1, \ldots, a_r \in K_{\omega},\quad
\| a_1 e_1 + \cdots + a_r e_r \|_{\pmb{e},\omega} :=
\begin{cases}
\max \{ |a_1|_{\omega}, \ldots, |a_r|_{\omega} \} & \text{if $\omega \in \Omega\setminus\Omega_{\infty}$},\\
|a_1|_{\omega} + \cdots + |a_r|_{\omega} & \text{if $\omega \in \Omega_{\infty}$},
\end{cases}
\]
and let $\varphi_{\pmb{e},\omega}$ be
the metric of $L_{\omega}$ induced by $\|\ndot\|_{\pmb{e},\omega}$ and
the surjective homomorphism $E_{\omega} \otimes_{K_{\omega}} \mathcal O_{X_{\omega}}
\to L_{\omega}$.
Let $\xi_{\pmb{e}} := \{ \|\ndot\|_{\pmb{e},\omega} \}_{\omega \in \Omega}$
and let $\varphi_{\pmb{e}} := \{ \varphi_{\pmb{e},\omega} \}_{\omega \in \Omega}$.
The metric family $\varphi_{\pmb{e}}$ is called the \emph{quotient metric family of $L$ induced by 
$\beta$ and $\pmb{e}$}.

\begin{rema}\label{Rem: double dual induce the same metric family}
We keep the above notation. Given a fixed surjective homomorphism $\beta:f^*(E)\rightarrow L$,
by Proposition \ref{Pro:doubedualandquotient} (see also Remark \ref{Rem:extensiondoubledual}), the norm family $\xi$ and the double dual norm family $\xi^{\vee\vee}$ induce the same quotient metric family on $L$. 
\end{rema}

\begin{prop}\label{Pro:distance of quotient metric families}
Let $(E,\xi)$ and $(E',\xi')$ be finite-dimensional vector spaces equipped with \emph{dominated}\index{dominated}\index{metric family!dominated} norm families. Let $f:X\rightarrow \Spec K$ be a projective scheme over $\Spec K$, $\beta:f^*(E)\rightarrow L$ and $\beta':{f'}^*(E')\rightarrow L$ be two surjective homomorphisms inducing closed immersions $i:X\rightarrow\mathbb P(E)$ and $i':X\rightarrow\mathbb P(E')$, and $\varphi$ and $\varphi'$ be quotient metric families induced by $(E,\xi)$ and $\beta$, and by $(E',\xi')$ and $\beta'$, respectively. Then the local distance function $(\omega\in\Omega)\mapsto d_\omega(\varphi,\varphi')$ is $\nu$-dominated.
\end{prop}
\begin{proof}
We begin with the particular case where $E=E'$ and $\beta=\beta'$. By Proposition \ref{Pro:distrancequot}, for any $\omega\in\Omega$ one has
\[d_\omega(\varphi,\varphi')\leqslant d_\omega(\xi^{\vee\vee},\xi'{}^{\vee\vee}).\]
Note that the norm families $\xi^{\vee\vee}$ and $\xi'{}^{\vee\vee}$ are strongly dominated (see Remark \ref{Rem: strongly dominated})
By Corollary \ref{Cor:dominatedanddist} and the triangle inequality of the local distance function, we obtain that the function $(\omega\in\Omega)\mapsto d_\omega(\xi^{\vee\vee},\xi'{}^{\vee\vee})$ is $\nu$-dominated and then deduce that the function $(\omega\in\Omega)\mapsto d_\omega(\varphi,\varphi')$ is $\nu$-dominated.

In the general case, since $i$ and $i'$ are closed immersions, 
by Serre's vanishing therorem (cf. \cite[Theorem~5.2, Chapter~III]{Hart77}), there exists an integer $n\geqslant 1$ such that
\[
\Gamma(\mathbb P(E),\mathcal O_{\mathbb P(E)}(n))\longrightarrow\Gamma(X,L^{\otimes n})
\quad\text{and}\quad
\Gamma(\mathbb P(E'),\mathcal O_{\mathbb P(E')}(n))\longrightarrow\Gamma(X,L^{\otimes n})
\]
are surjective, so that the natural homomorphisms
\[E^{\otimes n}\longrightarrow\mathrm{Sym}^n(E)=\Gamma(\mathbb P(E),\mathcal O_{\mathbb P(E)}(n))\longrightarrow\Gamma(X,L^{\otimes n})\]
and 
\[{E'}{}^{\otimes n}\longrightarrow\mathrm{Sym}^n(E')=\Gamma(\mathbb P(E'),\mathcal O_{\mathbb P(E')}(n))\longrightarrow\Gamma(X,L^{\otimes n})\]
are both surjective. Therefore both surjective homomorphisms $\beta$ and $\beta'$ factorise through $f^*\Gamma(X,L^{\otimes n})$. Moreover, by Remark \ref{Rem: powers of quotient}, if we equip $E^{\otimes n}$ and ${E'}{}^{\otimes n}$ with the $\varepsilon,\pi$-tensor power norm families  (see \S\ref{Subsec:Norm families}) of $\xi$ and $\xi'$ respectively, then the corresponding quotient metric families are $n\varphi$ and $n\varphi'$ respectively. Note that the $\varepsilon,\pi$-tensor powers of $\xi$ and $\xi'$ are dominated (see Proposition \ref{Pro:dominancealgebraic}
 \ref{Item:tensordom}). 
Therefore, by the special case proved above, we obtain that the function 
\[(\omega\in\Omega)\longmapsto d_\omega(n\varphi,n\varphi')=nd_\omega(\varphi,\varphi')\]
is $\nu$-dominated. The proposition is thus proved.
\end{proof}

\subsection{Dominated metric families}

Throughout this subsection, let $f:X\rightarrow \Spec K$ be a projective $K$-scheme.

\begin{defi}\label{Def: dominated metric family:very ample}
Let $L$ be a very ample invertible $\mathcal O_X$-module. We say that a metric family $\varphi$ on $L$ is \emph{dominated}\index{dominated}\index{metric family!dominated ---} if there exist a finite-dimensional vector space $E$ over $K$, a dominated norm family $\xi$ on $E$, and a surjective homomorphism $\beta:f^*(E)\rightarrow L$ inducing a closed immersion $X\rightarrow\mathbb P(E)$, such that the quotient metric family $\varphi'$ induced by $(E,\xi)$ and $\beta$ satisfies the following condition: \begin{quote}\emph{the local distance function $(\omega\in\Omega)\mapsto d_\omega(\varphi,\varphi')$ is $\nu$-dominated.}\end{quote}
\end{defi}

\begin{rema}
With the above definition, Proposition \ref{Pro:distance of quotient metric families} implies the 
following assertions. Let $E$ be a finite-dimensional vector space over $K$ equipped with a dominated norm family $\xi$. 
Let $L$ be an invertible $\mathcal O_X$-module and $\beta:f^*(E)\rightarrow L$ be a surjective homomorphism inducing a closed immersion $X\rightarrow\mathbb P(E)$. Then the quotient metric family induced by $(E,\xi)$ and $\beta$ is dominated. Moreover, if $\varphi_1$ and $\varphi_2$ are two metric families on $L$ which are dominated, then the local distance function $(\omega\in\Omega)\mapsto d_\omega(\varphi_1,\varphi_2)$ is $\nu$-dominated.
\end{rema}

\begin{prop}\label{Pro: dominancy by tensor product}
Let $L_1$ and $L_2$ be very ample invertible $\mathcal O_X$-modules. Assume that $\varphi_1$ and $\varphi_2$ are dominated metric families on $L_1$ and $L_2$ respectively. Then $\varphi_1+\varphi_2$ is a dominated metric family on $L_1\otimes L_2$.
\end{prop}
\begin{proof}
Since the metric families $\varphi_1$ and $\varphi_2$ are dominated, there exist finite-dimensional vector spaces $E_1$ and $E_2$ over $K$, dominated norm families $\xi_1$ and $\xi_2$ on $E_1$ and $E_2$ respectively, and surjective homomorphisms $\beta_1:f^*(E_1)\rightarrow L_1$ and $\beta_2:f^*(E_2)\rightarrow L_2$ inducing closed immersions $X\rightarrow\mathbb P(E_1)$ and $X\rightarrow\mathbb P(E_2)$ respectively, such that if we denote by $\widetilde\varphi_1$ and $\widetilde\varphi_2$ the quotient metric families induced by $(E_1,\xi_1)$ and $\beta_1$ and by $(E_2,\xi_2)$ and $\beta_2$ respectively, then the local distance functions
\[(\omega\in\Omega)\longmapsto d_\omega(\varphi_1,\widetilde\varphi_1)\text{\quad and\quad}(\omega\in\Omega)\longmapsto d_\omega(\varphi_2,\widetilde\varphi_2)\]are $\nu$-dominated. Consider now the composition morphism
\[\iota:\xymatrix{\relax X\ar[r]^-{(\iota_1,\iota_2)}&\mathbb P(E_1)\times_K\mathbb P(E_2)\ar[r]^-{\varsigma}&\mathbb P(E_1\otimes_KE_2)},\]
where $\iota_1$ and $\iota_2$ are closed immersions corresponding to $\beta_1$ and $\beta_2$, and $\varsigma$ is the Segre embedding. Note that $\iota$ is the closed immersion corresponding to the surjective homomorphism \[\beta_1\otimes\beta_2:f^*(E_1\otimes_KE_2)\cong f^*(E_1)\otimes_{\mathcal O_X}f^*(E_2)\longrightarrow L_1\otimes_{\mathcal O_X}L_2. \]
Moreover, if we equip $E_1\otimes_KE_2$ with the $\varepsilon,\pi$-tensor product norm family of $\xi_1$ and $\xi_2$, then the quotient metric family on $L_1\otimes L_2$ induced by $(E_1\otimes_KE_2,\xi_1\otimes_{\varepsilon,\pi}\xi_2)$ and $\beta_1\otimes\beta_2$ identifies with $\widetilde{\varphi}_1+\widetilde{\varphi}_2$. This is a consequence of Proposition \ref{Pro:quotientr1eps} (for the non-Archimedean case) and  Proposition \ref{Pro:quotientavecpitensor} (for the Archimedean case). Since
\[\forall\,\omega\in\Omega,\quad d_\omega(\varphi_1+\varphi_2,\widetilde\varphi_1+\widetilde\varphi_2)\leqslant d_\omega(\varphi_1,\widetilde{\varphi}_1)+d_\omega(\varphi_2,\widetilde\varphi_2),\]
we obtain that the function $(\omega\in\Omega)\mapsto d_\omega(\varphi_1+\varphi_2,\widetilde{\varphi}_1+\widetilde{\varphi}_2)$ is $\nu$-dominated. Therefore the metric family $\varphi_1+\varphi_2$ is dominated.
\end{proof}

\begin{defi}\label{Def: dominated metric family}
Let $L$ be an invertible $\mathcal O_X$-module and $\varphi$ be a metric family on $L$. We say that $\varphi$ is \emph{dominated}\index{dominated}\index{metric family!dominated} if there exist two very ample invertible $\mathcal O_X$-modules $L_1$ and $L_2$ together with dominated metric families $\varphi_1$ and $\varphi_2$ on $L_1$ and $L_2$ respectively, such that 
$L= L_2\otimes L_1^\vee$ and $\varphi=\varphi_2-\varphi_1$. 
\end{defi}

\begin{rema}
In the case where the invertible $\mathcal O_X$-module $L$ is very ample, the condition of dominancy in Definition \ref{Def: dominated metric family} is actually equivalent to that in Definition \ref{Def: dominated metric family:very ample}. In order to explain this fact (in avoiding confusions), in this remark we temporary say that a metric family $\varphi$ on a very ample invertible $\mathcal O_X$-module $L$ is \emph{strictly dominated}\index{strictly dominated}\index{metric family!strictly dominated} if it satisfies the condition in Definition \ref{Def: dominated metric family:very ample}. Clearly, if $\varphi$ is strictly dominated, then it is dominated (namely satisfies the condition in Definition \ref{Def: dominated metric family}) since we can write $L$ as $L^{\otimes 2}\otimes L^\vee$ and $\varphi$ as $2\varphi-\varphi$. Conversely, if $\varphi$ is dominated, then there exist very ample invertible $\mathcal O_X$-modules $L_1$ and $L_2$ such that $L\cong L_2\otimes L_1^\vee$ and strictly dominated metric families $\varphi_1$ and $\varphi_2$ on $L_1$ and $L_2$ such that $\varphi=\varphi_2-\varphi_1$. We pick an arbitrary strictly dominated metric family $\varphi'$ on $L$. By Proposition \ref{Pro: dominancy by tensor product}, we obtain that $\varphi'+\varphi_1$ is a strictly dominated metric family on $L_2$. Hence the local distance function \[(\omega\in\Omega)\longmapsto d_\omega(\varphi_2,\varphi'+\varphi_1)=d_\omega(\varphi+\varphi_1,\varphi'+\varphi_1)=d_\omega(\varphi,\varphi')\] is $\nu$-dominated. Therefore the metric family $\varphi$ is strictly dominated. 
\end{rema}

{
\begin{prop}
\label{Pro: Domination of quotient metric family}
Let $E$ be a finite-dimensional vector space over $K$ equipped with a norm family $\xi$, $\beta:f^*(E)\rightarrow L$ be a surjective homomorphism (we \emph{do not} assume that $\beta$ induces a closed immersion), and $\varphi$ be the quotient metric family induced by $(E,\xi)$ and $\beta$. Suppose that $\xi$ is a dominated norm family. Then $\varphi$ is a dominated metric family.
\end{prop}
\begin{proof}
Since $X$ is a projective $K$-scheme, there exists a very ample invertible $\mathcal O_X$-module $L'$. Let $E'$ be a finite dimensional vector space over $K$ and $\beta':f^*(E')\rightarrow L'$ be a surjective homomorphism, which induces a closed embedding of $X$ in $\mathbb P(E')$, which we denote by $\lambda'$. Let $\lambda:X\rightarrow\mathbb P(E)$ be the $K$-morphism induced by $\beta$. Then the tensor product homomorphism \[\beta\otimes\beta':f^*(E)\otimes_{\mathcal O_X} f^*(E')\cong f^*(E\otimes_KE')\longrightarrow L\otimes_{\mathcal O_X}L'\] corresponds to the composed $K$-morphisme
\[\xymatrix{\relax X\ar[r]^-{(\lambda,\lambda')}&\mathbb P(E)\times_K\mathbb P(E')\ar[r]^-{\varsigma}&\mathbb P(E\otimes_KE')},\]
where $\varsigma$ is the Segre embedding. Since $X$ is separated over $\Spec K$ and $\lambda'$ is a closed immersion, the morphism $(\lambda,\lambda')$ is a closed immersion. Therefore, the morphism from $X$ to $\mathbb P(E\otimes_K E')$ induced by $\beta\otimes\beta'$ is a closed embedding. 

Let $\xi'$ be a dominated norm family on $E'$ and $\varphi'$ be the metric family on $L'$ induced by $(E',\xi')$ and $\beta'$. By definition the metric family $\varphi'$ is dominated (see Proposition \ref{Pro:dominancealgebraic} \ref{Item:tensordom}). Moreover, the metric family $\varphi+\varphi'$ is induced by $(E\otimes E',\xi\otimes_{\varepsilon,\pi}\xi')$ and $\beta\otimes\beta'$ (see the proof of Proposition \ref{Pro: dominancy by tensor product}). As $\xi$ and $\xi'$ are dominated, we obtain that $\xi\otimes_{\varepsilon,\pi}\xi'$ is dominated. Therefore the metric families $\varphi+\varphi'$ is dominated. Hence the metric family $\varphi$ is also dominated.
\end{proof}
}

\begin{prop}\label{Pro: dominancy preserved by operators}
Let $L$ and $L'$ be invertible $\mathcal O_X$-modules, and $\varphi$ and $\varphi'$ be
metric families on $L$ and $L'$, respectively.
\begin{enumerate}[label=\rm(\arabic*)]
\item\label{Item: dual dom}  
If $\varphi$ is dominated, then the dual metric family $-\varphi$ on $L^\vee$ is dominated.
\item\label{Item: tensor dom} 
If $\varphi$ and $\varphi'$ are dominated, then
the tensor product metric family $\varphi+\varphi'$ on $L\otimes L'$ is dominated.
\item\label{Item: local distance of dominated metric families} 
If $L = L'$ and $\varphi$ and $\varphi'$ are dominated, then
the local distance function $(\omega\in\Omega)\mapsto d_\omega(\varphi,\varphi')$ is $\nu$-dominated.

\item\label{Item: dom plus dist dom implies dom} If $L = L'$, $\varphi'$ is dominated and
the local distance function $(\omega\in\Omega)\mapsto d_\omega(\varphi,\varphi')$ is $\nu$-dominated,
then $\varphi$ is dominated.

\item\label{Item: dilatation dominated} If $r \varphi$ is dominated for some non-zero integer $r$, then
$\varphi$ is dominated.

\item\label{Item: restriction dominated} Let $g:Y\rightarrow X$ be a projective morphism of $K$-schemes. If $\varphi$ is dominated, then
$g^*(\varphi)$ is also dominated. 

\end{enumerate}
\end{prop}
\begin{proof}
\ref{Item: dual dom} Let $L_1$ and $L_2$ be very ample invertible $\mathcal O_X$-modules and $\varphi_1$ and $\varphi_2$ be dominated metric families on $L_1$ and $L_2$ respectively, such that $L\cong L_2\otimes L_1^\vee$ and that $\varphi=\varphi_2-\varphi_1$. Then one has $L^\vee\cong L_1\otimes L_2^\vee$ and $-\varphi=\varphi_1-\varphi_2$. Hence the metric family $-\varphi$ is dominated.

\ref{Item: tensor dom} Let $L_1$, $L_2$, $L_1'$ and $L_2'$ be very ample invertible $\mathcal O_X$-modules, $\varphi_1$, $\varphi_2$, $\varphi_1'$ and $\varphi_2'$ be dominated metric families on $L_1$, $L_2$, $L_1'$ and $L_2'$ respectively, such that $L\cong L_2\otimes L_1^\vee$, $L'\cong L_2'\otimes L_1'{}^\vee$, $\varphi=\varphi_2-\varphi_1$ and $\varphi'=\varphi_2'-\varphi_1'$. Note that $L\otimes L'\cong (L_2\otimes L_2')\otimes (L_1\otimes L_1')^\vee$, and $\varphi+\varphi'=(\varphi_2+\varphi_2')-(\varphi_1+\varphi_1')$. By Proposition \ref{Pro: dominancy by tensor product}, the metric families $\varphi_2+\varphi_2'$ and $\varphi_1+\varphi_1'$ are dominated. Hence $\varphi+\varphi'$ is dominated.

\ref{Item: local distance of dominated metric families} Let $L_1$ be a very ample invertible $\mathcal O_X$-module and $\varphi_1$ be dominated metric family on $L_1$. Let $\varphi_2=\varphi+\varphi_1$ and $\varphi_2'=\varphi'+\varphi_1$. By \ref{Item: tensor dom}, the metric families $\varphi_2$ and $\varphi_2'$ are dominated. Since the invertible $\mathcal O_X$-module $L_2$ is very ample, by Proposition \ref{Pro:distance of quotient metric families} we obtain that the local distance function $(\omega\in\Omega)\mapsto d_\omega(\varphi_2,\varphi_2')=d_\omega(\varphi,\varphi')$ is $\nu$-dominated.

\ref{Item: dom plus dist dom implies dom} First we assume that $L$ is very ample.
As $\varphi'$ is dominated, there exist a finite-dimensional vector space $E$ over $K$, a dominated norm family $\xi$ on $E$, and a surjective homomorphism $\beta:f^*(E)\rightarrow L$ inducing a closed immersion $X\rightarrow\mathbb P(E)$ such that, if $\psi$ is the quotient metric family induced by $(E,\xi)$ and $\beta$, then the local distance function $(\omega\in\Omega)\mapsto d_\omega(\varphi',\psi)$ is $\nu$-dominated.
Note that \[
\forall\, \omega \in \Omega,\quad
d_{\omega}(\varphi, \psi) \leqslant d_{\omega}(\varphi, \varphi') + d_{\omega}(\varphi', \psi).
\]
Thus $(\omega\in\Omega)\mapsto d_\omega(\varphi,\psi)$ is $\nu$-dominated, as required.

In general, there are very ample invertible $\mathcal O_{X}$-modules $L_1$ and $L_2$,
and dominated metric families $\varphi'_1$ and $\varphi'_2$ on $L_1$ and $L_2$, respectively,
such that $L = L_1 \otimes L_2^{\vee}$ and $\varphi' = \varphi'_1 - \varphi'_2$.
We set $\varphi_1 = \varphi + \varphi'_2$ and $\varphi_2 = \varphi'_2$. Then $\varphi = \varphi_1 - \varphi_2$, and $\varphi_1$ and $\varphi_2$ are metric families of $L_1$ and $L_2$, respectively.
As
\[
d_{\omega}(\varphi_1, \varphi'_1) = d_{\omega}(\varphi_1 - \varphi_2, \varphi'_1-\varphi'_2)
= d_{\omega}(\varphi, \varphi'),
\]
$(\omega\in\Omega)\mapsto d_\omega(\varphi_1,\varphi'_1)$ is $\nu$-dominated, so that
$\varphi_1$ is dominated by the previous observation.
Therefore, $\varphi$ is also dominated.

\ref{Item: dilatation dominated} As $-r\varphi$ is dominated by \ref{Item: dual dom}, we may assume that $r > 0$.
We choose a dominated metric family $\psi$ on $L$. By \ref{Item: tensor dom}, $r \psi$ is dominated, so that, by \ref{Item: local distance of dominated metric families},
\[
(\omega\in\Omega)\mapsto d_\omega(r\varphi,r\psi) = r d_\omega(\varphi,\psi)
\]
is $\nu$-dominated. Therefore $(\omega\in\Omega)\mapsto d_\omega(\varphi,\psi)$
is also $\nu$-dominated, and hence the assertion follows from \ref{Item: dom plus dist dom implies dom}.

\ref{Item: restriction dominated} First we assume that $L$ is very ample. Then
there exist a finite-dimensional vector space $E$ over $K$, a dominated norm family $\xi$ on $E$, a surjective homomorphism $\beta:f^*(E)\rightarrow L$ inducing a closed immersion $X\rightarrow\mathbb P(E)$ such that, if $\psi$ denotes the quotient metric family induced by $(E,\xi)$ and $\beta$, then the local distance function $(\omega\in\Omega)\mapsto d_\omega(\varphi,\psi)$ is $\nu$-dominated. {Note that $\beta : E \otimes_K \mathcal O_{X} \to L$ yields
the surjective homomorphism $g^*(\beta) : E \otimes_K \mathcal O_{Y} \to g^*(L)$.
Moreover, $g^*(\psi)$ coincides with quotient metric family induced by $(E,\xi)$ and $g^*(\beta)$. By Proposition \ref{Pro: Domination of quotient metric family}, $g^*(\psi)$ is a dominated metric family. Moreover, for any $\omega\in\Omega$ one has
\[
d_\omega(g^*(\varphi),g^*(\psi))\leqslant d_\omega(\psi, \varphi).
\]
By \ref{Item: dom plus dist dom implies dom} we obtain that $g^*(\varphi)$ is a dominated metric family.

In general, there are very ample invertible $\mathcal O_{X}$-modules $L_1$ and $L_2$
such that $L = L_1 \otimes L_2^{\vee}$.
Let $\varphi_1$ be dominated metric family on $L_1$. If we set $\varphi_2 = \varphi_1 - \varphi$,
then $\varphi_2$ is dominated by \ref{Item: dual dom} and \ref{Item: tensor dom}. By the previous case,
$g^*(\varphi_1)$ and $g^*(\varphi_2)$ are dominated, so that,
by \ref{Item: dual dom} and \ref{Item: tensor dom} again, $g^*(\varphi) = g^*(\varphi_1) - g^*(\varphi_2)$ is
also dominated.}
\end{proof}

\begin{theo}\label{Thm: Fubini-Study dominated}
Let $f:X\rightarrow\Spec K$ be a geometrically reduced projective $K$-scheme and $L$ be an invertible $\mathcal O_X$-module, equipped with a dominated metric family $\varphi=\{\varphi_\omega\}_{\omega\in\Omega}$. For any $\omega\in\Omega$, let $\|\ndot\|_{\varphi_\omega}$ be the sup norm on $H^0(X,L)\otimes_KK_\omega$ corresponding to the metric $\varphi_\omega$. Then the norm family $\xi=\{\|\ndot\|_{\varphi_\omega}\}_{\omega\in\Omega}$ on $H^0(X,L)$ is strongly dominated.
\end{theo}

\begin{proof}
Let us begin with the following claim:

\begin{enonce}{Claim}\label{Claim:Thm: Fubini-Study dominated:00}
If the assertion of the theorem holds under the assumption that $X$ is geometrically integral,
then it holds in general.
\end{enonce}

\begin{proof}
One can find a finite extension $K'$ of $K$ and the irreducible decomposition
$X_1 \cup \cdots \cup X_n$ of $X_{K'}$ such that $X_1, \ldots, X_n$ are geometrically integral.
We use the same notation as in Corollary~\ref{coro:dominated:pullback:dominated:orginal},
which says that it is sufficient to see that $\xi_{K'}$ is dominated.


Let $\psi = \{ \psi_{\omega} \}_{\omega \in \Omega}$ be another metric family of $L$.
Then $d_{\omega'}(\varphi_{K'}, \psi_{K'}) = d_{\omega}(\varphi, \psi)$ for all $\omega \in \Omega$ and
$\omega' \in \Omega_{K'}$ with $\pi_{K'/K}(\omega') = \omega$.
Therefore, one can see that $\varphi_{K'} = \{\varphi_{K',\omega'}\}_{\omega'\in\Omega_{K'}}$
is dominated.

Let $\varphi_{K',i}$ be the restriction of $\varphi_{K'}$ to $X_i$.
Then, by Proposition~\ref{Pro: dominancy preserved by operators}, \ref{Item: restriction dominated}, $\varphi_{K',i}$ is also dominated for all $i$.
On the other hand, one has
the natural injective homomorphism
\[
H^0(X_{K',\omega'}, L_{K', \omega'}) \to H^0(X_{1,\omega'}, L_{K', \omega'}) \oplus \cdots \oplus 
H^0(X_{n,\omega'}, L_{K',\omega'}).
\]
Here we give a norm $\|\ndot\|_{\omega'}$ on $H^0(X_{1,\omega'}, L_{K',\omega'}) \oplus \cdots \oplus 
H^0(X_{n,\omega'}, L_{K',\omega'})$ given by 
\[
\|(x_1, \ldots, x_n)\|_{\omega'} = \max_{i\in \{ 1, \ldots, n\}} \big\{ \| x_i \|_{\varphi_{K',i,\omega'}} \big\}.
\]
By our assumption, $\xi_{K',i}$ is dominated for all $i$, so that $\{ \|\ndot\|_{\omega'} \}_{\omega\in \Omega_{K'}}$
is also dominated by Proposition~\ref{Pro:dominancealgebraic}, \ref{Item:sommedom}.
Note that $\|\ndot\|_{\varphi_{K',\omega'}}$ is the restriction norm of $\|\ndot\|_{\omega'}$, and hence
$\xi_{K'} = \{ \|\ndot\|_{\varphi_{K',\omega'}}\}_{\omega' \in \Omega_{K'}}$ is
dominated by Proposition~\ref{Pro:dominancealgebraic}, \ref{Item:subdom}, as desired.
\end{proof}

From now on, we assume that $X$ is geometrically integral.

\begin{enonce}{Claim}\label{Claim:Thm: Fubini-Study dominated:01}
If $L$ is very ample, the assertion of the theorem holds.
\end{enonce}

\begin{proof}

Let $E=H^0(X,L)$, $r = \dim_K(E)$ and $\beta:f^*(E)\rightarrow L$ be the canonical surjective homomorphism which induces a closed immersion of $X$ in $\mathbb P(E)$. 
Note that any non-zero element of $E$ can not identically vanish on $X$.
Hence there exist a finite extension $K'$ of $K$ together with closed points
$P_1,\ldots,P_r$ of $X$ such that
the residue filed $\kappa(P_i)$ at $P_i$ is contained in $K'$, and that we have a strictly decreasing sequence of $K'$-vector spaces
\[
E_0 \supsetneq E_1 \supsetneq E_2 \supsetneq \cdots \supsetneq E_{r-1} \supsetneq E_r = \{ 0 \},
\]
where $E_0 = E \otimes_K K'$ and \[E_i = \{ s \in E\otimes_K K' \,:\, s(P_1) = \cdots = s(P_i) = 0 \}\]
for $i\in\{1, \ldots, r\}$. 
In order to prove Claim~\ref{Claim:Thm: Fubini-Study dominated:01},
by virtue of Corollary~\ref{coro:dominated:pullback:dominated:orginal}, 
we may assume that $K' = K$.

Let $\omega_1, \ldots, \omega_r$ be local bases of $L$ around $P_1, \ldots, P_r$, respectively.
For each $i = 1, \ldots, r$, we define $\theta_{i} \in E^{\vee}$ to be
\[
\forall\, s \in E,\quad
\theta_i(s) = f_s(P_i)\qquad(\text{$s = f_s \omega_i$ around $P_i$}).
\]
Note that $\theta_1, \ldots, \theta_r$ are linearly independent over $K$.
Let $\boldsymbol{e}=\{e_i\}_{i=1}^r$ be the dual basis of $\{\theta_i\}_{i=1}^r$ over $K$, that is,
$(e_1, \ldots, e_r) \in E^r$ and $\theta_i(e_j) = \delta_{ij}$.
Here we define a norm $\|\ndot\|'_{\omega}$ as follows:
for any element $s\in E_{\omega}$ written as $s=\lambda_1e_1+\cdots+\lambda_re_r$ with $(\lambda_1,\ldots,\lambda_r)\in K_\omega^r$, 
\[\|s\|'_\omega=\max_{i\in\{1,\ldots,r\}}|\lambda_i|_\omega.\]
Let $\xi'$ be the norm family of $E$ given by $\{ \|\ndot\|'_{\omega} \}_{\omega \in \Omega}$. 
Since the volume of $\Omega_{\infty}$ is finite, 
$\xi'$ is dominated by Corollary~\ref{Cor:dominatedanddist}.

Let $\varphi'$ be the metric family of $L$ induced by $(E,\xi')$ and $\beta$.
Note that $\varphi'$ is dominated.
Let us see that $\|\ndot\|'_\omega=\|\ndot\|_{\varphi'_\omega}$ for any $\omega\in\Omega$.
First of all, by Proposition~\ref{Pro:positivityofquotientmetric},
$\|\ndot\|_{\varphi'_\omega} \leqslant \|\ndot\|_\omega'$.
For any $(i, j)\in\{1,\ldots,r\}^2$,  
one has 
\[
|e_j|_{\varphi'_\omega}(P_{i})=\begin{cases} 1, & \text{if $i=j$}, \\
0, & \text{if $i\not=j$}.
\end{cases}
\]
Therefore, for any $(\lambda_1,\ldots,\lambda_r)\in K_\omega^r$,
\[
\|\lambda_1e_1+\cdots+\lambda_re_r\|_{\varphi'_\omega} \geqslant\max_{i\in\{1,\ldots,r\}}|\lambda_1e_1+\cdots+\lambda_re_r|_{\varphi'_\omega}(P_{i})
=\max_{i\in\{1,\ldots,r\}}|\lambda_i|_\omega,
\]
as required.

By the inequality \eqref{Equ:comparaisondistancesup} in Subsection~\ref{subsec:distance:metric}, one has
\[
d_{\omega}(\|\ndot\|_{\varphi_{\omega}}, \|\ndot\|'_{\omega})
= d_{\omega}(\|\ndot\|_{\varphi_{\omega}}, \|\ndot\|_{\varphi'_{\omega}})
\leqslant d_{\omega}(\varphi_{\omega}, \varphi'_{\omega}).
\]
Therefore, $\omega \mapsto d_{\omega}(\|\ndot\|_{\varphi_{\omega}}, \|\ndot\|'_{\omega})$
is $\nu$-dominaited, and hence $\xi$ is dominaited by Proposition \ref{Pro:adelicvectorbundledist}.
\end{proof}

Finally let us consider the following claim:

\begin{enonce}{Claim}\label{Claim:Thm: Fubini-Study dominated:02}
For any non-zero element $s\in H^0(X,L)$, 
the function 
$(\omega\in\Omega)\longmapsto \ln\|s\|_{\varphi_\omega}$
is $\nu$-dominated. 
\end{enonce}

\begin{proof}
Fix a non-zero element $s\in H^0(X,L)$. 

\medskip
Let us construct a dominated metric family $\varphi'$ of $L$ such that the function 
$(\omega\in\Omega)\longmapsto \ln\|s\|_{\varphi'_\omega}$ is $\nu$-dominated.
Let $L_1$ be a very ample invertible $\mathcal O_X$-module such that $L_2:=L\otimes L_1$ is also very ample. The multiplication by the non-zero section $s$ defines an injective $K$-linear map from $H^0(X,L_1)$ to $H^0(X,L_2)$. We choose a dominated norm family $\xi_2'=\{\|\ndot\|_{2,\omega}'\}_{\omega,\in\Omega}$ on $H^0(X,L_2)$ and let $\xi_1'$ be the restriction of $\xi_2'$ on $H^0(X,L_1)$ via this injective map. By Proposition \ref{Pro:dominancealgebraic} \ref{Item:subdom}, the norm family $\xi_1'$ is also dominated. Let $\varphi_1'$ and $\varphi_2'$ be the quotient metric family induced by $\xi_1'$ and $\xi_2'$ respectively, where we consider the canonical surjective homomorphisms $f^*(H^0(X,L_1))\rightarrow L_1$ and $f^*(H^0(X,L_2))\rightarrow L_2$.  
We set $\varphi'=\varphi_2'-\varphi_1'$. By Propositions \ref{Pro:quotientderang1norm} and \ref{Pro: quotient metric pi ext2}, for all $\omega\in\Omega$, $x\in X_\omega^{\mathrm{an}}$ such that $s(x)\neq 0$ and $\ell\in L_{1,\omega}\otimes\widehat{\kappa}(x)\setminus\{0\}$, one has
\[\begin{split}|\ell|_{\varphi_{1,\omega}'}(x)&=\inf_{\begin{subarray}{c}u\in H^0(X,L_1),\,\lambda\in\widehat{\kappa}(x)^{\times}\\
u(x)=\lambda\ell
\end{subarray}}|\lambda|_x^{-1}\cdot\|su\|_{2,\omega}'\\
&\geqslant\inf_{\begin{subarray}{c}
v\in H^0(X,L_2),\,\lambda\in\widehat{\kappa}(x)^{\times}\\
v(x)=\lambda s(x)\ell
\end{subarray}}|\lambda|_x^{-1}\cdot\|v\|_{2,\omega}'\\&=|s(x)\ell|_{\varphi_{2,\omega}'}(x)=|s(x)|_{\varphi'_{\omega}}(x)\cdot|\ell|_{\varphi_{1,\omega}'}(x),
\end{split}\]
which leads to the inequality $\|s\|_{\varphi'_\omega,\sup}\leqslant 1$. Moreover, by Proposition \ref{Pro:sousmultiplicativesec}, for any non-zero section $u\in H^0(X,L_1)$, one has
\[\|su\|_{\varphi_{2,\omega}'}\leqslant\|s\|_{\varphi'_{\omega}}\cdot\|u\|_{\varphi_{1,\omega}'}.
\]
Therefore, the function $(\omega\in\Omega)\mapsto \ln\|s\|_{\varphi'_\omega}$ is non-positive and bounded from below by a $\nu$-dominated function. 

\medskip
By the inequality \eqref{Equ:comparaisondistancesup} in Subsection~\ref{subsec:distance:metric}, one has
\[
d_{\omega}(\|\ndot\|_{\varphi_{\omega}}, \|\ndot\|_{\varphi'_{\omega}})
\leqslant d_{\omega}(\varphi_{\omega}, \varphi'_{\omega}).
\]
On the other hand, by Proposition~\ref{Pro: dominancy preserved by operators} \ref{Item: local distance of dominated metric families}, the function $(\omega\in\Omega)\longmapsto d_{\omega}(\varphi_{\omega}, \varphi'_{\omega})$
is $\nu$-dominated, so that the assertion follows.
\end{proof}

We now proceed with the proof of the theorem. Let $L_1$ be a very ample invertible $\mathcal O_X$-module such that $L_2:=L\otimes L_1$ is also very ample. We fix a non-zero global section $t$ of $L_1$, which defines an injective $K$-linear map from $H^0(X,L)$ to $H^0(X,L_2)$. We choose a dominated metric family $\varphi_1=\{\varphi_{1,\omega}\}_{\omega\in\Omega}$ on $L_1$ such that $\|t\|_{\varphi_{1,\omega}}\leqslant 1$ for any $\omega\in\Omega$, which is possible if we take a strongly
dominated norm family $\xi_1=\{\|\ndot\|_{1,\omega}\}_{\omega\in\Omega}$ on $H^0(X,L_1)$ such that $\|t\|_{1,\omega}\leqslant 1$ for any $\omega\in\Omega$, and choose $\varphi_1$ as the quotient metric family induced by $(H^0(X,L_1),\xi_1)$ and the canonical surjective homomorphism $f^*(H^0(X,L_1))\rightarrow L_1$. Let $\varphi_2=\{\varphi_{2,\omega}\}_{\omega\in\Omega}$ be the metric family $\varphi+\varphi_1$ on $L_2$. By Proposition \ref{Pro: dominancy preserved by operators} \ref{Item: tensor dom}, the metric family $\varphi_2$ is dominated. Let $\xi_2$ be the norm family $\{\|\ndot\|_{\varphi_{2,\omega}}\}$ on $H^0(X,L_2)$. By 
Claim~\ref{Claim:Thm: Fubini-Study dominated:01}, the norm family $\xi_2$ is strongly dominated. 

Let $\{s_1,\ldots,s_m\}$ be a basis of $H^0(X,L)$. For any $i\in\{1,\ldots,m\}$, let $t_i=ts_i\in H^0(X,L_2)$. We choose sections $t_{m+1},\ldots,t_n$ in $H^0(X,L_2)$ such that $\{t_1,\ldots,t_n\} $ forms a basis of $H^0(X,L_2)$. Let ${\xi}_2^\circ=\{\|\ndot\|_{2,\omega}^\circ\}_{\omega\in\Omega}$ be the norm family on $H^0(X,L_2)$ such that, for any $\omega\in\Omega$ and any $(\lambda_1,\ldots,\lambda_n)\in K_\omega^n$,
\[\|\lambda_1t_1+\cdots+\lambda_nt_n\|_{2,\omega}^\circ=\begin{cases}\max_{i\in\{1,\ldots,n\}}|\lambda_i|_\omega,&\text{if $\omega\in\Omega\setminus\Omega_\infty$,}\\
|\lambda_1|_\omega+\cdots+|\lambda_n|_\omega,&\text{if $\omega\in\Omega_\infty$}.
\end{cases}\]

For any $\omega\in\Omega$ and any $(\lambda_1,\ldots,\lambda_\omega)\in K_\omega^m$, one has
\begin{equation*}\|\lambda_1s_1+\cdots+\lambda_ms_m\|_{\varphi_{\omega}}\geqslant\|\lambda_1t_1+\cdots+\lambda_mt_m\|_{\varphi_{2,\omega}}
\end{equation*}
since $\|t\|_{\varphi_{1,\omega},\sup}\leqslant 1$. As the norm family $\xi_2$ is strongly dominated, the local distance function $(\omega\in\Omega)\mapsto d_\omega(\xi_2,\xi_2^\circ)$ is $\nu$-dominated (see Corollary \ref{Cor:dominatedanddist}).  In particular, there exists a $\nu$-dominated function $A$ on $\Omega$ such that, for any $\omega\in\Omega$ and any $(\lambda_1,\ldots,\lambda_m)\in K_\omega^m$, one has
\begin{equation}\label{Equ: minoration norm lambda1s1 etc}\ln\|\lambda_1s_1+\cdots+\lambda_ms_m\|_{\varphi_{\omega}}\geqslant\ln\|\lambda_1t_1+\cdots+\lambda_mt_m\|_{2,\omega}^\circ-A(\omega).
\end{equation}
Moreover,
\begin{equation*}\|\lambda_1s_1+\cdots+\lambda_ms_m\|_{\varphi_\omega}\leqslant\|\lambda_1t_1+\cdots+\lambda_mt_m\|_{2,\omega}^\circ\cdot\max_{i\in\{1,\ldots,m\}}\|s_i\|_{\varphi_\omega}\end{equation*}
By Claim~\ref{Claim:Thm: Fubini-Study dominated:02}, for any $i\in\{1,\ldots,m\}$, the function $(\omega\in\Omega)\mapsto\ln\|s_i\|_{\varphi_\omega}$ is $\nu$-dominated. Therefore, there exists a $\nu$-dominated function $B$ on $\Omega$ such that, for any $\omega\in\Omega$ and any $(\lambda_1,\ldots,\lambda_m)\in K_\omega^m$, one has
\begin{equation}
\label{Equ: majoration norm lambda1s1 etc}
\ln\|\lambda_1s_1+\cdots+\lambda_ms_m\|_{\varphi_{\omega}}\leqslant\ln\|\lambda_1t_1+\cdots+\lambda_mt_m\|_{2,\omega}^\circ+B(\omega).
\end{equation}
The inequalities \eqref{Equ: minoration norm lambda1s1 etc} and \eqref{Equ: majoration norm lambda1s1 etc} imply that the local distance function $(\omega\in\Omega)\mapsto d_\omega(\xi,\xi^\circ)$ is $\nu$-dominated, where $\xi^\circ$ is the restricted norm family of $\xi_2^\circ$ on $H^0(X,L)$. The strong dominancy of $\xi$ then follows from Corollary \ref{Cor:dominatedanddist}. The theorem is thus proved.
\end{proof}

{
\begin{prop}\label{Pro: dominated on each point}
Let $K'$ be a finite extension of $K$ and $X=\Spec K'$. Let $L$ be an invertible $\mathcal O_X$-module. Then a metric family $\varphi$ on $L$ is dominated if and only if the corresponding norm family $\xi_L$ on $L$ relatively to the adelic curve $S\otimes_KK'$ (cf. Remark \ref{Rem: metric family on one point}) is dominated.
\end{prop}
\begin{proof}
First we assume that there exists a finite-dimensional vector space $E$ over $K$, a surjective $K'$-linear map $\beta:E\otimes_KK'\rightarrow L$, and a dominated norm family $\xi_E$ on $E$ such that the local distance function $(\omega\in\Omega)\mapsto d_\omega(\varphi,\varphi')$ is $\nu$-dominated, where $\varphi'$ denotes the metric family on $L$ induced by $(E,\xi_E)$ and $\beta$. Denote by  $\xi_L'$ the norm families on $L$ relatively to the adelic curve $S\otimes_KK'$, which correspond to the metric families  $\varphi'$. By definition $\xi_L'$ identifies with the quotient norm family of $\xi_{K'}$ induced by $\beta$. Since $\xi$ is dominated, also is $\xi_{K'}$ (cf. Corollary \ref{Cor:dominanceext}). Hence the norm family $\xi_L'$ is also dominated (cf. Proposition \ref{Pro:dominancealgebraic} \ref{Item:quotdom}). By Proposition \ref{Pro:adelicvectorbundledist}, we obtain that the norm family $\xi_L$ is dominated.

Conversely, we assume that the norm family $\xi_L$ is dominated. Let $e$ be a generator of $L$ as vector space over $K'$ and let $\xi_{e}=\{\norm{\ndot}_{e,\omega}\}_{\omega\in\Omega}$ be the norm family on $Ke$ such that $\norm{\lambda e}_{e,\omega}=|\lambda|_\omega$ for any $\omega\in\Omega$. Note that the norm family $\xi_{e,K'}$ on $(Ke)\otimes_KK'\cong L$ is dominated. By Corollary \ref{Cor:dominatedanddist}, the local distance function $(x\in\Omega_{K'})\mapsto d_x(\xi_{e,K'},\xi_L)$ is $\nu_{K'}$-dominated. Let $\varphi'$ be the metric family induced by $(Ke,\xi_{e})$ and the canonical isomorphism $(Ke)\otimes_KK'\cong L$. For any $x\in \Omega_{K'}$ one has $d_x(\xi_{e,K'},\xi_L)=d_x(\varphi',\varphi)$. By Proposition \ref{Pro: dominancy preserved by operators} \ref{Item: dom plus dist dom implies dom}, the metric family $\varphi$ is dominated.
\end{proof}
}

\subsection{Universally dense point families}

In this subsection, we consider universally dense point families (cf. Lemma~\ref{Lem: density of closed pts}) and their consequences.

\begin{lemm}\label{Lem: density of closed pts}
Let $K$ be a field, $X$ be a scheme locally of finite type over $\Spec K$ and $K'/K$ be a field extension. Let $X_{K'}$ be the fibre product $X\times_{\Spec K}\Spec K'$ and $\pi:X_{K'}\rightarrow X$ be the morphism of projection. For any closed point $P$ of $X$, the set $\pi^{-1}(\{P\})$ is finite and consists of closed points of $X_{K'}$. Moreover, if $F$ is a set consisting of  
closed points of $X$, which is Zariski dense in $X$, then the subset $\pi^{-1}(F)$ of $X_{K'}$ is Zariski dense.
\end{lemm}
\begin{proof}
Let $P'$ be a point of $X_{K'}$ such that $P=\pi(P')$ is a closed point. Since $X$ is locally of finite type over $\Spec K$, the residue field of $P$ is a finite extension of $K$ (this is a consequence of Zariski's lemma, see \cite{Zariski47}). As the residue field of $P'$ is a quotient ring of $\kappa(P)\otimes_KK'$, we obtain that it is a finite extension of $K'$. Moreover, since $\kappa(P)\otimes_KK'$ is an Artinian $K'$-algebra, the set $\pi^{-1}(\{P\})$ is finite.

In the following, we fix an algebraic closure $K'{}^{\mathrm{ac}}$ of the field $K'$ and we denote by $K^{\mathrm{ac}}$ the algebraic closure of $K$ in $K'{}^{\mathrm{ac}}$. For any closed point $P$ of $X$, we choose an arbitrary embedding of $\kappa(P)$ in $K^{\mathrm{ac}}$ so that we can consider the residue field $\kappa(P)$ as a subfield of $K^{\mathrm{ac}}$. Similarly, for any $P'\in\pi^{-1}(\{P\})$, we choose an embedding of the residue field $\kappa(P')$ in $K'{}^{\mathrm{ac}}$ which extends the embedding $\kappa(P)\rightarrow K^{\mathrm{ac}}$.

To prove the lemma it suffices to verify that, for any affine open subset $U$ of $X$, $\pi^{-1}(U\cap F)$ is Zariski dense in $U_{K'}$. Therefore we may assume without loss of generality that $X$ is an affine scheme of finite type over $K$.  We let $A$ be the coordinate ring of $X$. Thus the coordinate ring of $X_{K'} $ is $A\otimes_KK'$. Let $f$ be an element of $A\otimes_KK'$ and $\widetilde f$ be the canonical image of $f$ in $A\otimes_KK'{}^{\mathrm{ac}}$. We write $\widetilde f$ as a linear combination
\[\widetilde f=a_1g_1+\cdots+a_ng_n,\]
where $g_1,\ldots,g_n$ are elements in $A\otimes_KK^{\mathrm{ac}}$, and $a_1,\ldots,a_n$ are elements of $K'{}^{\mathrm{ac}}$ which are linearly independent over $K^{\mathrm{ac}}$. Let $P$ be a closed point of $X$ and $\mathfrak m_P$ be the maximal ideal of $A$ corresponding to $P$. Assume that for any $P'\in\pi^{-1}(\{P\})$ one has $f(P')=0$, where $f(P')$ denotes the image of $f$ by the projection map $A\otimes_KK'\rightarrow(A\otimes_KK')/\mathfrak m_{P'}$, with $\mathfrak m_{P'}$ being the (maximal) ideal of $A\otimes_KK'$ corresponding to $P'$. Then the canonical image of $f$ in $(A/\mathfrak m_P)\otimes_KK'$ is nilpotent, which implies that the canonical image of  $\widetilde f$ in $(A/\mathfrak m_P)\otimes_KK'{}^{\mathrm{ac}}$ is nilpotent. In particular, the canonical image of $\widetilde f$ by the composed map
\[A\otimes_KK'{}^{\mathrm{ac}}\longrightarrow (A/\mathfrak m_P)\otimes_KK'{}^{\mathrm{ac}}=\kappa(P)\otimes_KK'{}^{\mathrm{ac}}\longrightarrow K'{}^{\mathrm{ac}}\]
is zero, where the last map in the above diagram is given by $\lambda\otimes \mu\mapsto \lambda\mu$ for any $\lambda\in K(P)\subseteq K^{\mathrm{ac}}$ and $\mu\in K'{}^{\mathrm{ac}}$. In other words, one has
\[a_1g_1(P)+\cdots+a_ng_n(P)=0,\]
where for each $i\in\{1,\ldots,n\}$, $g_n(P)$ denotes the image of $g_n$ by the composed map
\[A\otimes_KK^{\mathrm{ac}}\longrightarrow (A/\mathfrak m_P)\otimes_KK^{\mathrm{ac}}=\kappa(P)\otimes_KK^{\mathrm{ac}}\longrightarrow K^{\mathrm{ac}}.\]
Since $a_1,\ldots,a_n$ are linearly independent over $K^{\mathrm{a}}$, we obtain that $g_1(P)=\cdots=g_n(P)=0$. Since this holds for any $P\in F$ and since $F$ is Zariski dense in $X$, we obtain that $g_1,\ldots,g_n$ are nilpotent elements in $A\otimes_KK^{\mathrm{ac}}$. Therefore $\widetilde f$ is a nilpotent element of $A\otimes_KK'{}^{\mathrm{ac}}$. Since the extension $K'{}^{\mathrm{ac}}/K^{\mathrm{ac}}$ equips $K'{}^{\mathrm{ac}}$ with a structure of $K^{\mathrm{ac}}$-algebra which is faithfully flat, the canonical map $A\otimes_KK^{\mathrm{ac}}\rightarrow A\otimes_KK'{}^{\mathrm{ac}}$ is injective (see \cite{Demazure_Gabriel70}, Chapitre I, \S1, $\text{n}^\circ$2, Lemme 2.7).  Therefore, $f$ is a nilpotent element in $A\otimes_KK'$. This shows that $\pi^{-1}(F)$ is Zariski dense in $X_{K'}$. 
\end{proof}

\begin{prop}\label{Pro: mesurability out of omega0}
Let $S=(K,(\Omega,\mathcal A,\nu),\phi)$ be an adelic curve. Let $X$ be a projective $K$-scheme, $L$ be an invertible $\mathcal O_X$-module, equipped with a metric family $\varphi$. We assume that, for any closed point $P$ in $X$, the norm family $P^*(\varphi)$ (cf. Remark \ref{Rem: metric family on one point}) on $P^*(L)$ is measurable. Then, for any $s\in H^0(X,L)$, the function \[(\omega\in\Omega\setminus\Omega_0)\longmapsto \|s\|_{\varphi_\omega}\]
is measurable, where $\Omega_0$ denotes the set of $\omega\in\Omega$ such that $|\ndot|_\omega$ is the trivial absolute value, and we consider the restriction of the $\sigma$-algebra $\mathcal A$ on $\Omega\setminus\Omega_0$.
\end{prop}
\begin{proof}
As $X$ is projective over $K$, considering the coefficients of defining homogeneous polynomials of $X$, we can find a subfield $K_0$ of $K$ which is finitely generated over the prime field of $K$ (and hence $K_0$ is countable) and a projective scheme $X_0$ over $\Spec K_0$ such that $X\cong X_0\times_{\Spec K_0}\Spec K$. Let $\mathscr P$ be the set of closed points $P$ in $X$ whose canonical image in $X_0$ is a closed point. By Lemma \ref{Lem: density of closed pts}, $\mathscr P$ is a Zariski dense and countable subset of $X$.

By the assumption of the proposition, for any closed point $P$ of $X$, the function
\[(x\in\Omega_{\kappa(P)})\longmapsto |s|_{\varphi_{\omega}}(P_x), \qquad(\omega = \pi_{\kappa(P)/K}(x))\]
is $\mathcal A_{\kappa(P)}$-measurable. By Proposition \ref{Pro: max of measurable function}, we obtain that the function
\[(\omega\in\Omega)\longmapsto \max_{x\in\pi_{\kappa(P)/K}^{-1}(\{\omega\})}|s|_{\varphi_\omega}(P_x)\]
is $\mathcal A$-measurable.  Therefore, the function
\begin{equation}\label{Equ: max along closed pts}(\omega\in\Omega)\longmapsto\sup_{P\in\mathscr P}\max_{
x\in
\pi_{\kappa(P)/K}^{-1}(\{\omega\})
}|s|_{\varphi_\omega}(x)\end{equation}
is $\mathcal A$-measurable since $\mathscr P$ is countable.

To obtain the conclusion of the proposition, it remains to show that the function coincides with $\omega\mapsto \|s\|_{\varphi_\omega,\sup}$ on $\Omega\setminus\Omega_0$. For this purpose it suffices to verify that, for any $\omega\in\Omega\setminus\Omega_0$, the set 
\[F_\omega=\left\{P_x\,:\,\text{$P\in\mathscr P$ and $x\in\pi_{\kappa(P)/K}^{-1}(\{\omega\})$} \right\}\]
is dense in $X_\omega^{\mathrm{an}}$, where $X_\omega:=X\times_{\Spec K}\Spec K_\omega$. Let $j_\omega: X_\omega^{\mathrm{an}}\rightarrow X_\omega$ be the specification map. By Lemma \ref{Lem: density of closed pts}, $j_\omega(F_\omega)$ is Zariski dense in $X_\omega$ and hence $F_\omega$ is dense in $X_\omega^{\mathrm{an}}$ with respect to the Berkovich topology (see \cite[Corollary 3.4.5]{Berkovich90}). The proposition is thus proved.
\end{proof}

We assume that $K$ is equipped with the trivial absolute value $|\ndot|_0$.
Let $F$ be a finitely generated field over $K$ such that
the transcendence degree of $F$ over $K$ is one.
Let $C_F$ be a regular projective curve over $K$ such that
the function field of $C_F$ is $F$, that is, $C_F$ is the unique regular model of $F$ over $K$.
It is well-known that,
for any absolute value $|\ndot|$ of $F$ over $K$ (i.e. the restriction of $|\ndot|$ to $K$ is trivial), 
there are a closed point $\xi$ of $C_F$ and $q \in \mathbb R_{\geqslant 0}$ such that
$|\varphi| = \exp(-q \ord_{\xi}(\varphi))$ for all $\varphi \in F^{\times}$ (see \cite[\S I.6]{Hart77} and \cite[Proposition II.(3.3)]{Neukirch}).
Note that in the case where $q = 0$, the absolute value is trivial. 
We say that $q$ is the \emph{exponent}\index{exponent} of $|\ndot|$.
The absolute value given by $\exp(-q \ord_{\xi}(\ndot))$ is denoted by $|\ndot|_{(\xi, q)}$.
Let $X$ be a projective scheme over $\Spec K$. The Berkovich space associated with $X$ is denote by $X^{\mathrm{an}}$ (see Definition \ref{Def:mapassociatetoamor}). Let $j : X^{\mathrm{an}} \to X$ be the specification map.
Let us consider the following subsets $X_0^{\mathrm{an}}$, $X_{1,\mathbb Q}^{\mathrm{an}}$ and $X_{\leqslant 1,\mathbb Q}^{\mathrm{an}}$ in $X^{\mathrm{an}}$:
\[
\begin{cases}
X_0^{\mathrm{an}} := \{ x \in X^{\mathrm{an}} \,:\, \text{$j(x)$ is closed} \},\\[1ex]
X_{1,\mathbb Q}^{\mathrm{an}} := 
\left\{ x \in X^{\mathrm{an}} \left| \begin{array}{l} \text{the Zariski closure of $\{ j(x) \}$ has dimension one} \\
\text{and the exponent of the corresponding absolute value} \\
\text{is rational} \end{array} \kern-.5em\right.\right\},\\[1ex]
X_{\leqslant 1,\mathbb Q}^{\mathrm{an}} := X_0^{\mathrm{an}} \cup X_{1,\mathbb Q}^{\mathrm{an}}.
\end{cases}
\]

\begin{lemm}\label{Lem: density in Berkovich space}
$X_{\leqslant 1,\mathbb Q}^{\mathrm{an}}$ is dense in $X^{\mathrm{an}}$ with respect to the Berkovich topology.
\end{lemm}

\begin{proof}
To prove the lemma, 
we need to show that, for any regular function $f$ over an affine open subset $U = \Spec A$ of $X$ and for any point $x\in U^{\mathrm{an}}$, the value $|f|(x)$ belongs to the  closure $T$ of the set $\{|f|(z) \,:\, z\in X_{\leqslant 1,\mathbb Q}^{\mathrm{an}} \cap U^{\mathrm{an}}\}$ in $\mathbb R$. 
First let us see the following claim:

\begin{enonce}{Claim}
\begin{enumerate}
\renewcommand{\labelenumi}{(\alph{enumi})}
\item If $f$ is not a nilpotent element, then $1 \in T$.

\item If $f$ has a zero point in $U$, then $0 \in T$.
\end{enumerate}
\end{enonce}

\begin{proof}
(a) As $f$ is not a nilpotent element, there is a closed point $z$ of $U$ such that $f(z) \not= 0$, so that
$1 = |f|(z) \in T$.

(b) In this case one can find a closed point $z'$ with $f(z') = 0$.
Therefore, $0 = |f|(z') \in T$.
\end{proof}

Let us go back to the proof of the lemma.
Let $X' =\Spec A'$ be the Zariki closure of $\{ j(x) \}$ in $U$, where
$A'$ is the quotient domain of $A$ by the prime ideal corresponding to $j(x)$. 
Let $|\ndot|_x$ be the absolute value on the field of fractions of $A'$ corresponding to $x$. 
If $\dim X' = 0$, then
$j(x)$ is closed, so that the assertion is obvious. Moreover, if $|f'|_x$ is either $0$ or $1$, then the assertion is also obvious by the above claim. In particular, if $f' = \rest{f}{X'}$ is algebraic over $K$,
then $|f|(x) = |f'|_x$ is either $0$ or $1$, and hence the assertion is true. Therefore we may assume that $\dim X' \geqslant 1$, $f'$ is transcendental over $K$ and $|f'|_x \in \mathbb R_{\geqslant 0} \setminus \{ 0, 1 \}$.

Consider the ring $A' \otimes_{K[f']} K(f')$, where $K(f')$ is the fraction field of $K[f']$.
This is a localisation of the ring $A'$ 
with respect to the multiplicatively closed subset $K[f'] \setminus \{ 0 \}$. 
We pick a closed point $\zeta'\in \Spec(A' \otimes_{K[f']} K(f'))$ and let $\zeta$ be the canonical image of $\zeta'$ in $U$. Then the point $\zeta\in U$ has dimension $1$ and the canonical image $f''$ of $f'$ in the residue field $\kappa(\zeta)$ is transcendental over $K$
because $f'$ is an element of the constant field $K(f')$ of the variety $\Spec(A' \otimes_{K[f']} K(f'))$. 
In particular, the natural homomorphism $K[f'] \to K[f'']$ is an isomorphism, which yields an isomorphism
$K(f') \overset{\sim}{\longrightarrow} K(f'')$.
Let $|\ndot|'_x$ be the restriction of $|\ndot|_x$ to $K(f')$, and $|\ndot|''_x$ be the absolute value of $K(f'')$ such that
the above isomorphism gives rise to an isometry 
\[
(K(f'),|\ndot|_x') \overset{\sim}{\longrightarrow} (K(f''), |\ndot|''_x).
\]
Then $|f''|''_x = |f'|'_x = |f|(x)$.
Let $|\ndot|_{\zeta}$ be an extension of $|\ndot|_x''$ to the residue field $\kappa(\zeta)$ and
$C_{\zeta}$ be a regular and projective model of $\kappa(\zeta)$ over $K$.
Then there are a closed point $\xi$ of $C_{\zeta}$ and $q \in \mathbb R_{>0}$ such that
$|\ndot|_{\zeta} = |\ndot|_{(\xi,q)}$.
Thus the assertion follows if we consider a sequence $\{q_n\}_{n=1}^{\infty}$ of rational numbers 
such that $\lim_{n\to\infty} q_n = q$.
\end{proof}

\begin{rema}
In the case where the absolute value of $K$ is non-trivial, 
$X_{0}^{\mathrm{an}}$ is dense in $X^{\mathrm{an}}$ (cf. Lemma \ref{Lem: density of closed pts} and \cite[Corollary 3.4.5]{Berkovich90}).
However we need one more layer $X_{1,\mathbb Q}^{\mathrm{an}}$ if the absolute value of $K$ is trivial.
Moreover, if the dimension of every irreducible component of $X$ is greater than or equal to one, then
$X_{1,\mathbb Q}^{\mathrm{an}}$ is dense in $X^{\mathrm{an}}$ with respect to the Berkovich topology.
Indeed, it is sufficient to show that
$X_{0}^{\mathrm{an}}$ is contained in the closure of $X_{1,\mathbb Q}^{\mathrm{an}}$. 
Let $x \in X_{0}^{\mathrm{an}}$ and choose a subvariety $C'$ of $X$ such that $\dim C' = 1$ and $j(x) \in C'$.
Let $\mu : C \to C'$ be the normalisation of $C'$ and $\xi \in C$ with $\mu(\xi) = j(x)$.
Note that $\lim_{\substack{q \to \infty \\ q \in \mathbb Q_{>0}}} |\ndot|_{(\xi, q)}$ gives rise to the trivial valuation of the residue field $\kappa(\xi)$,
which means that $x$ belongs to the closure of $X_{1,\mathbb Q}^{\mathrm{an}}$. 
\end{rema}

\begin{rema}\label{Rem: countability}
If $K$ is countable, then $X_{\leqslant 1,\mathbb Q}^{\mathrm{an}}$ is also countable. In fact, the set of all closed points of a projective scheme over $K$ is countable. Therefore, $X_0^{\mathrm{an}}$ is countable. Moreover, if we fix an increasing sequence \[K_1\subseteq K_2\subseteq\ldots\subseteq K_n\subseteq K_{n+1}\subseteq\ldots\] of finite extensions of the field $K(T)$ of rational functions such that $\bigcup_{n\in\mathbb N,\,n\geqslant 1}K_n$ is the algebraic closure of $K(T)$, then any point $z\in X_1^{\mathrm{an}}$ is represented by a point of $X$ valued in certain $K_n$ equipped with an absolute value over $K$ of rational exponent. Suppose that $K_n$ identifies with the rational function field of the projective curve $C_n$ over $K$, there are only countably many such absolute values since $K$ is assumed to be countable. Hence $X_{1,\mathbb Q}^{\mathrm{an}}$ is also countable.
\end{rema}

\begin{rema}
We assume that $K$ is uncountable and the absolute value of $K$ is trivial.
In this case, we can not expect a dense countable subset of $X^{\mathrm{an}}$.
Indeed, let $S$ be a countable subset of $\mathbb P_K^{1,\mathrm{an}}$.
Let $r : \mathbb P_K^{1,\mathrm{an}} \to \mathbb P_K^{1}$ be the reduction map.
As $r(S)$ is countable and $K$ is uncountable, there is a closed point $\xi \in \mathbb P_K^{1}$ such that
$\xi \not\in r(S)$. We set $I := \{ \exp(-q \ord_{\xi}(\ndot)) \,:\, q \in \intervalle{]}0{\infty}{[} \}$.
Then $I$ is an open set of $\mathbb P_K^{1,\mathrm{an}}$ and $r(I) = \{ \xi \}$, so that
$I \cap S = \emptyset$, which shows that $S$ is not dense in $\mathbb P_K^{1,\mathrm{an}}$.
\end{rema}

Let $S=(K,(\Omega,\mathcal A,\nu),\phi)$ be an adelic curve. Denote by $\Omega_0$ the set of $\omega\in\Omega$ such that the absolute value $|\ndot|_\omega$ on $K$ is trivial. Let $\mathcal A_0$ be the restriction of the $\sigma$-algebra on $\Omega_0$. Let $X$ be a projective scheme over $\Spec K$. We equip $K$ with the trivial absolute value and denote by $X^{\mathrm{an}}$ be the Berkovich space associated with $X$ (see Definition \ref{Def:mapassociatetoamor}). Suppose given an invertible $\mathcal O_X$-module $L$ equipped with a metric family $\varphi=\{\varphi_\omega\}_{\omega\in\Omega}$.  For any point $x\in X^{\mathrm{an}}$, the metric $\varphi$ induces, for any $\omega\in\Omega$, a norm $|\ndot|_{\varphi_\omega}(x)$ on the one-dimensional vector space $L\otimes_{\mathcal O_X}\widehat{\kappa}(x)$.

\begin{prop}\label{Pro: measurability on omega0}
Let $S=(K,(\Omega,\mathcal A,\nu),\phi)$ be an adelic curve. We assume that the field $K$ is countable. Let $X$ be a projective scheme over $\Spec K$, $L$ be an invertible $\mathcal O_X$-module and $\varphi$ be 
metric family on $L$. 
Suppose that, for any $x\in X_{\leqslant 1,\mathbb Q}^{\mathrm{an}}$ and any non-zero element $\ell$ in $L\otimes_{\mathcal O_X}\widehat{\kappa}(x)$, the function $(\omega\in\Omega_0)\mapsto |\ell|_{\varphi_\omega}(x)$ is $\mathcal A_0$-measurable.  Then for any $s\in H^0(X,L)$, the function
\[(\omega\in\Omega_0)\longmapsto \|s\|_{\varphi_\omega}\]
is measurable on $(\Omega_0,\mathcal A_0)$.
\end{prop}
\begin{proof}
By Lemma \ref{Lem: density in Berkovich space}, we can write $\|s\|_{\varphi_\omega}$ as
\[\|s\|_{\varphi_\omega}=\sup_{z\in X_{\leqslant 1,\mathbb Q}^{\mathrm{an}}}|s|_{\varphi_\omega}(z).\]
By the assumption of the proposition, each function $(\omega\in\Omega_0)\mapsto|s|_{\varphi_\omega}(z)$ is $\mathcal A_0$-measurable. Since $X_{\leqslant 1,\mathbb Q}^{\mathrm{an}}$ is a countable set (see Remark \ref{Rem: countability}), we deduce that the function $(\omega\in\Omega_0)\mapsto\|s\|_{\varphi_\omega}$ is also $\mathcal A_0$-measurable. The proposition is thus proved.
\end{proof}

\subsection{Measurable metric families}

\begin{defi}\label{Def: measurability of metric family}
Let $S=(K,(\Omega,\mathcal A,\nu),\phi)$ be an adelic curve, $X$ be a projective scheme over $\Spec K$ and $L$ be an invertible $\mathcal O_X$-module. We say that a metric family $\varphi=\{\varphi_\omega\}_{\omega\in\Omega}$ on $L$ is \emph{measurable}\index{measurable}\index{metric family!measurable} if the following conditions are satisfied:
\begin{enumerate}[label=\rm(\alph*)]
\item for any closed point $P$ in $X$, the norm family $P^*(\varphi)$ on $P^*(L)$ is measurable,
\item for any point $x\in X^{\mathrm{an}}_{\leqslant 1,\mathbb Q}$ (where we consider the trivial absolute value on $K$ in the construction of the Berkovich space $X^{\mathrm{an}}$) and any element $\ell$ in $L\otimes_{\mathcal O_X}\widehat{\kappa}(x)$, the function $(\omega\in\Omega_0)\mapsto|\ell|_{\varphi_\omega}(x)$ is $\mathcal A_0$-measurable, where $\Omega_0$ the set of $\omega\in\Omega$ such that the absolute value $|\ndot|_\omega$ on $K$ is trivial, and $\mathcal A_0$ is the restriction of the $\sigma$-algebra on $\Omega_0$.
\end{enumerate}
\end{defi}

\begin{prop}\label{Pro: measurability preserved by operators}
Let $S=(K,(\Omega,\mathcal A,\nu),\phi)$ be an adelic curve and $X$ be a projective scheme over $\Spec K$. 
\begin{enumerate}
\renewcommand{\labelenumi}{{\rm(\arabic{enumi})}}
\item\label{Item: measurability dual} Let $L$ be an invertible $\mathcal O_X$-module equipped with a measurable metric family $\varphi$, then the dual metric family $-\varphi$ on $L^\vee$ is measurable.
\item\label{Item: measurability tensor} Let $L_1$ and $L_2$ be two invertible $\mathcal O_X$-modules. If $\varphi_1$ and $\varphi_2$ are measurable metric families on $L_1$ and $L_2$ respectively, then the metric family $\varphi_1+\varphi_2$ on $L_1\otimes L_2$ is measurable.
\item \label{Item: measurability power} Let $L$ be an invertible $\mathcal O_X$-module equipped with a metric family $\varphi$. Suppose that there exists an integer $n\geqslant 1$ such that $n\varphi$ is measurable, then the metric family $\varphi$ is also measurable.

\item Let $L$ be an invertible $\mathcal O_X$-module equipped with a measurable metric family $\varphi$, and $f:Y\rightarrow X$ be a projective morphism of $K$-schemes. Then $f^*({\varphi})$ is measurable.
\end{enumerate}
\end{prop}
\begin{proof}
(\ref{Item: measurability dual}) Let $P$ be a closed point of $X$. One has $P^*(-\varphi)=P^*(\varphi)^\vee$. Since $\varphi$ is measurable, the norm family $P^*(\varphi)$ is measurable. By Proposition \ref{Pro:mesurability} \ref{Item:dual measurability}, we obtain that $P^*(-\varphi)$ is also measurable.

Let $x$ be a point of $X_{\leqslant 1,\mathbb Q}^{\mathrm{an}}$, $\ell$ be a non-zero element of $L\otimes_{\mathcal O_X}\widehat{\kappa}(x)$ and $\ell^\vee$ be the dual element of $\ell$ in $L^\vee\otimes_{\mathcal O_X}\widehat{\kappa}(x)$. Then for any $\omega\in\Omega_0$ one has 
\[|\ell|_{\varphi_\omega}(x)\cdot|\ell^\vee|_{-\varphi_\omega}(x)=1.\]
Since the function $(\omega\in\Omega_0)\mapsto |\ell|_{\varphi_\omega}(x)$ is $\mathcal A_0$-measurable, also is the function $(\omega\in\Omega_0)\mapsto |\ell^\vee|_{-\varphi_\omega}(x)$. Therefore the metric family $-\varphi$ is measurable.

(\ref{Item: measurability tensor}) Let $P$ be a closed point of $X$. One has $P^*(\varphi_1+\varphi_2)=P^*(\varphi_1)\otimes P^*(\varphi_2)$. Since $\varphi_1$ and $\varphi_2$ are both measurable, the norm families $P^*(\varphi_1)$ and $P^*(\varphi_2)$ are measurable. By Proposition \ref{Pro:mesurability} \ref{Item:tensor measurability}, the tensor norm family $P^*(\varphi_1)\otimes P^*(\varphi_2)$ is also measurable.

Let $x$ be a point of $X_{\leqslant 1,\mathbb Q}^{\mathrm{an}}$, $\ell_1$ and $\ell_2$ be non-zero elements of $L_1\otimes_{\mathcal O_X}\widehat{\kappa}(x)$ and $L_2\otimes_{\mathcal O_X}\widehat{\kappa}(x)$, respectively. For any $\omega\in\Omega_0$ one has
\[|\ell_1|_{\varphi_{1,\omega}}(x)\cdot|\ell_2|_{\varphi_{2,\omega}}(x)=|\ell_1\otimes\ell_2|_{(\varphi_1+\varphi_2)_\omega}(x).\]
Since the metric families $\varphi_1$ and $\varphi_2$ are measurable, the functions $(\omega\in\Omega_0)\mapsto |\ell_1|_{\varphi_{1,\omega}}$ and $(\omega\in\Omega_0)\mapsto |\ell_2|_{\varphi_{2,\omega}}$ are both $\mathcal A_0$-measurable. Hence the function $(\omega\in\Omega_0)\mapsto |\ell_1\otimes\ell_2|_{(\varphi_1+\varphi_2)_\omega}(x)$ is $\mathcal A_0$-measurable.

(\ref{Item: measurability power}) Let $P$ be a closed point of $X$. One has $P^*(n\varphi)=P^*(\varphi)^{\otimes n}$. Since $n\varphi$ is measurable, the norm family $P^*(n\varphi)$ is measurable. By Proposition \ref{Pro:mesurability} \ref{Item:power measurability}, we obtain that the norm family $P^*(\varphi)$ is also measurable.

Let $x$ be a point of $X_{\leqslant 1,\mathbb Q}^{\mathrm{an}}$ and $\ell$ be a non-zero element of $L\otimes_{\mathcal O_X}\widehat{\kappa}(x)$. For any $\omega\in\Omega_0$ one has
$|\ell^{\otimes n}|_{n\varphi_\omega}(x)=|\ell|_{\varphi_\omega}(x)^n$ and hence $|\ell|_{\varphi_\omega}(x)=|\ell|_{n\varphi_\omega}(x)^{1/n}$.
Since the metric family $n\varphi$ is measurable, we obtain that the function $(\omega\in\Omega_0)\mapsto |\ell^{\otimes n}|_{n\varphi_\omega}(x)$ is $\mathcal A_0$-measurable. Hence the function $(\omega\in\Omega_0)\mapsto |\ell|_{\varphi_\omega}(x)$ is also $\mathcal A_0$-measurable. Therefore the metric family $\varphi$ is measurable.

(4) This is obvious by the definition of the measurability of $\varphi$.
\end{proof}

The following proposition shows that the measurability of metric family is a property stable by pointwise limit.

\begin{prop}
Let $S=(K,(\Omega,\mathcal A,\nu),\phi)$ be an adelic curve, $X$ be a projective scheme over $\Spec K$ and $L$ be an invertible $\mathcal O_X$-module. Let $\varphi$ and $\{\varphi_n\}_{n\in\mathbb N}$ be metric families on $L$ such that, for any $\omega\in\Omega$ and any $x\in X_\omega^{\mathrm{an}}$ (where $X_\omega^{\mathrm{an}}$ is the Berkovich space associated with $X_\omega:=X\times_{\Spec K}\Spec K_\omega$) one has 
\[\lim_{n\rightarrow+\infty}d(|\ndot|_{\varphi_{n,\omega}}(x),|\ndot|_{\varphi_\omega}(x))=0.\]
Assume that the metric families $\varphi_n$, $n\in\mathbb N$, are measurable. Then the limite metric family $\varphi$ is also measurable.
\end{prop}
\begin{proof}
Let $P$ be a closed point of $X$ and $s$ be a non-zero element in $P^*(L)$. By the assumption of the proposition, for any $\omega'\in\Omega_{K(P)}$ over $\omega\in\Omega$, the norm $\|\ndot\|_{\omega'}$ indexed by $\omega'$ in the norm family $P^*(\varphi)$ is given by $|\ndot|_{\varphi_\omega}(P_{\omega'})$, where $P_{\omega'}$ is the algebraic point of $X_\omega$ determined by $P$  and $|\ndot|_{\omega'}$. As well, for any $n\in\mathbb N$, the norm $\|\ndot\|_{n,\omega'}$ indexed by $\omega'$  in $P^*(\varphi_n)$ is given by $|\ndot|_{\varphi_{n,\omega}}(P_{\omega'})$. Therefore, by the limit  assumption of the proposition, the sequence of functions
\[(\omega'\in\Omega_{K(P)})\longmapsto \|s\|_{n,\omega'}\]
converges pointwisely to $(\omega'\in\Omega_{K(P)})\mapsto\|s\|_{\omega'}$, which implies that the latter function is $\mathcal A_{K(P)}$-measurable by the measurability assumption of the proposition. 

Similarly, for any point $x\in X^{\mathrm{an}}_{\leqslant 1,\mathbb Q}$ (where we consider the trivial absolute value on $K$ in the construction of $X^{\mathrm{an}}$), and any element $\ell\in L\otimes_{\mathcal O_X}\widehat{\kappa}(x)$, the sequence of functions
\[(\omega\in\Omega_0)\longmapsto |\ell|_{\varphi_{n,\omega}}(x),\quad n\in\mathbb N\]
converges pointwisely to $(\omega\in\Omega_0)\mapsto|\ell|_{\varphi_\omega}(x)$. Since each function in the sequence is $\mathcal A_0$-measurable, also is the limit function. The proposition is thus proved.
\end{proof}

\begin{prop}
Let $S=(K,(\Omega,\mathcal A,\nu),\phi)$ be an adelic curve, $f:X\rightarrow\Spec K$ be a projective $K$-scheme and $L$ be an invertible $\mathcal O_X$-module. We assume that $K$ admits a countable subfield which is dense in every $K_\omega$, $\omega\in\Omega$. Let $E$ be a finite-dimensional vector space over $K$, equipped with a measurable norm family $\xi=\{\|\ndot\|_\omega\}_{\omega\in\Omega}$. Suppose give a surjective homomorphism $\beta:f^*(E)\rightarrow L$ and let $\varphi$ be the quotient metric family on $L$ induced by $(E,\xi)$ and $\beta$. Then the metric family $\varphi$ is measurable.  
\end{prop}
\begin{proof} 
By Proposition \ref{Item: dual norm family is measurable}, the double dual norm family $\xi^{\vee\vee}$ is measurable. By
Remark \ref{Rem: double dual induce the same metric family}, we can replace $\xi$ by $\xi^{\vee\vee}$ without changing the corresponding quotient metric family. Hence we may assume without loss of generality that the norm $\|\ndot\|_\omega$ is ultrametric for any $\omega\in\Omega\setminus\Omega_\infty$.

Let $P$ be a closed point of $X$. Then the norm family $P^*(\varphi)$ identifies with the quotient norm family of $\xi\otimes K(P)$ induced by the surjective homomorphism \[P^*(\beta):P^*(f^*(E))\cong E\otimes_KK(P)\longrightarrow P^*(L).\]
By Proposition \ref{Pro:mesurabilityofquotient} \ref{Item: measurability of tensor product} and \ref{Item: measurability of quotient norm family}, the norm family $P^*(\varphi)$ is measurable.

Assume that $\Omega_0$ is not empty. In this case the field $K$ itself is countable. Let $x$ be a point of $X_{\leqslant 1,\mathbb Q}^{\mathrm{an}}$ and $\ell$ be a non-zero element of $L\otimes_{\mathcal O_X}\widehat{\kappa}(x)$. By Proposition \ref{Pro:quotientderang1norm}, for any $\omega\in\Omega_0$ one has
\[|\ell|_{\varphi_\omega}(x)=\inf_{\begin{subarray}{c}s\in E,\,\lambda\in\widehat{\kappa}(x)^{\times}\\
s(x)=\lambda\ell\end{subarray}}|\lambda|_x^{-1}\cdot\|s\|_{\omega},\]
where $s(x)$ denotes the image of $s$ by the quotient map \[\beta_x:E\otimes_K\widehat{\kappa}(x)\rightarrow L\otimes_{\mathcal O_X}\widehat{\kappa}(x).\]
Since the norm family $\xi$ is measurable, the function $(\omega\in\Omega_0)\mapsto\|s\|_\omega$ is $\mathcal A_0$-measurable. Moreover, since $K$ is a countable field, the vector space $E$ is a countable set. Hence we obtain that the function $(\omega\in\Omega_0)\mapsto |\ell|_{\varphi_\omega}(x)$ is $\mathcal A_0$-measurable. Therefore, the metric family $\varphi$ is measurable.
\end{proof}

\begin{defi}
Let $S=(K,(\Omega,\mathcal A,\nu),\phi)$ be an adelic curve, $\pi:X\rightarrow\Spec K$ be a projective $K$-scheme, $L$ be an invertible $\mathcal O_X$-module and $\varphi$ be a metric family on $L$. We denote by $\pi_*(\varphi)$ the norm family $\{\|\ndot\|_{\varphi_\omega}\}_{\omega\in\Omega}$ on $\pi_*(L)$.
\end{defi}

The propositions \ref{Pro: mesurability out of omega0} and \ref{Pro: measurability on omega0} can be resumed as follows.

\begin{theo}\label{Thm: measurability of linear series}
Let $S=(K,(\Omega,\mathcal A,\nu),\phi)$ be a adelic curve, $X$ be a projective $K$-scheme and $L$ be an invertible $\mathcal O_X$-module. We assume that, either $\Omega_0$ is empty, or the field $K$ is countable. For any measurable metric family $\varphi$ on $L$, the norm family $\pi_*(\varphi)$ is measurable.
\end{theo}

\section{Adelic line bundle and Adelic divisors}
In this section, we fix an adelic curve $S=(K,(\Omega,\mathcal A,\nu),\phi)$ be an adelic curve. 

\begin{defi}
Let $X$ be a projective scheme over $\Spec K$. We call \emph{adelic line bundle}\index{adelic line bundle} on $X$ any invertible $\mathcal O_X$-module $L$ equipped with a metric family $\varphi$ which is dominated an measurable.

\subsection{Height function}\label{Subsec: height}
Let $(L,\varphi)$ be an adelic line bundle on $X$. For any closed point $P$ of $X$, the norm family $P^*(\varphi)$ on $P^*(L)$ is measurable and dominated (see Propositions \ref{Pro: dominancy preserved by operators} \ref{Item: restriction dominated} and \ref{Pro: dominated on each point} for the dominance of $P^*(\varphi)$). Therefore $(P^*(L),P^*(\varphi))$ is an adelic line bundle on $S$ (cf. Proposition \ref{Pro:criterion of adelic line bundle}). We denote by $h_{(L,\varphi)}(P)$ the Arakelov degree of this adelic line bundle, called the \emph{height}\index{height} of $P$ with respect to the adelic line bundle $(L,\varphi)$. By abuse of notation, we also denote by $h_{(L,\varphi)}(\ndot)$ the function on the set $X(K^{\mathrm{ac}})$ of algebraic points of $X$ sending any $K$-morphism $\Spec K^{\mathrm{ac}}\rightarrow X$ to the height of the image of the $K$-morphism.
\end{defi}

{

\begin{prop}\label{Pro: additivity of height function}
Let $X$ be a projective $K$-scheme, $(L_1,\varphi_1)$ and $(L_2,\varphi_2)$ be adelic line bundles on $X$.
\begin{enumerate}[label=\rm(\arabic*)]
\item\label{Item: additivity of height function} For any closed point $P$ of $X$, one has
\[h_{(L_1\otimes L_2,\varphi_1+\varphi_2)}(P)=h_{(L_1,\varphi_1)}(P)+h_{(L_2,\varphi_2)}(P).\]
\item\label{Item: distance height function} Assume that $L_1$ and $L_2$ are the same invertible $\mathcal O_X$-module $L$. Then, for any closed point $P$ of $X$, one has
\begin{equation}\label{Equ: difference height}|h_{(L,\varphi_1)}(P)-h_{(L,\varphi_2)}(P)|\leqslant \operatorname{dist}(\varphi_1,\varphi_2).\end{equation}
\end{enumerate}
\end{prop}
\begin{proof}
\ref{Item: additivity of height function} This follows directly from Proposition \ref{Pro: pull back by point} and \ref{Pro:degreduproduit}.

\ref{Item: distance height function} Let $P$ be a closed point of $X$. By definition, for any $\omega\in\Omega$ and any $x\in M_{K(P),\omega}$ one has $d_x(P^*(\varphi_1),P^*(\varphi_2))\leqslant d_\omega(\varphi_1,\varphi_2)$. Therefore, by taking the integral with respect to $x$, we obtain the inequality \eqref{Equ: difference height}.
\end{proof}

The following proposition shows that, if the adelic curve $S$ has the Northcott property, then the height function associated with an adelic line bundle with ample underlying invertible sheaf has a finiteness property of Northcott type.
\begin{prop}\label{Pro:Northcott for arithmetic varieties}
Assume that the adelic curve $S$ has the Northcott property (cf. Definition~\ref{def:northcott:property}). Let $X$ be a projective $K$-scheme and $(L,\varphi)$ be an adelic line bundle on $X$ such that $L$ is ample. For all positive real numbers $\delta$ and $C$, the set
\begin{equation}\label{Equ: Northcott property}\{P\in X(K^{\mathrm{ac}})\,:\,h_{(L,\varphi)}(P)\leqslant C,\;[K(P):K]\leqslant \delta\}\end{equation}
is finite.
\end{prop}
\begin{proof}
By Proposition \ref{Pro: additivity of height function}, for any integer $n\geqslant 1$, one has $h_{(nL,n\varphi)}=nh_{(L,\varphi)}$ as function on $X(K^{\mathrm{ac}})$. Therefore, without loss of generality we may assume that $L$ is very ample. Moreover, by Proposition \ref{Pro: additivity of height function} \ref{Item: distance height function}, to prove the proposition it suffice to check the finiteness of \eqref{Equ: Northcott property} for a particular choice of metric family $\varphi$. Thus we may assume without loss of generality that there exist a finite dimensional vector space $E$ over $K$ and a surjective homomorphism $\beta:E\otimes_K\mathcal O_K\rightarrow L$ inducing a closed immersion $X\rightarrow\mathbb P(E)$, together with a basis $\boldsymbol{e}=\{e_i\}_{i=0}^r$ of $E$ over $K$, such that $\varphi$ identifies with the quotient metric family induced by $(E,\xi_{\boldsymbol{e}})$ and $\beta$ (see Example \ref{Exa:adelicbasis} for the definition of $\boldsymbol{e}$). Let $P$ be a closed point of $X$. Then $(P^*(L)\otimes_{K(P)}K^{\mathrm{ac}},P^*(\varphi)_{K^{\mathrm{ac}}})$ is a quotient adelic line bundle of $(E_{K^{\mathrm{ac}}},\xi_{\boldsymbol{e},K^{\mathrm{ac}}})$ and hence $(L_{K^{\mathrm{ac}}}^{\vee},P^*(\varphi)_{K^{\mathrm{ac}}}^\vee)$ identifies with an adelic line subbundle of $(E_{K^{\mathrm{ac}}}^\vee,\xi_{\boldsymbol{e},K^{\mathrm{ac}}}^\vee)$.  Suppose that, as a vector subspace of rank $1$ of $E_{K^{\mathrm{ac}}}^\vee$, $L_{K^{\mathrm{ac}}}^\vee$ is generated by the vector $a_0e_0^\vee+\cdots+a_re_r^\vee$, where \[[a_0:\cdots:a_r]\in\mathbb P^r(K^{\mathrm{ac}}),\] then one has
\[\begin{split}&\quad\;\widehat{\deg}(P^*(L),P^*(\varphi))=\widehat{\deg}(L_{K^{\mathrm{ac}}},P^*(\varphi)_{K^{\mathrm{ac}}})=-\widehat{\deg}(L_{K^{\mathrm{ac}}}^\vee,P^*(\varphi)_{K^{\mathrm{ac}}}^\vee)\\
&=\int_{\Omega_{K^{\mathrm{ac}}}} \ln \left( \max \{ |a_0|_{\chi}, \ldots, |a_r|_{\chi} \} \right) \nu_{K^{\mathrm{ac}}}(d\chi).\end{split}\]
Therefore the finiteness of \eqref{Equ: Northcott property} follows from Theorem \ref{thm:Northcott:thm}.
\end{proof}}

\subsection{Essential minimum}
In this subsection, we fix an \emph{integral} projective scheme  $X$ over $\Spec K$. For any adelic line bundle $(L,\varphi)$ on $X$, we define the \emph{essential minimum} of $(L,\varphi)$ as
\[\widehat{\mu}_{\mathrm{ess}}(L,\varphi):=\sup_{Z\subsetneq X}\inf_{P\in(X\setminus Z)(K^{\mathrm{ac}})}h_{(L,\varphi)}(P),\]
where $Z$ runs over the set of all strict Zariski closed subsets of $X$, and $P$ runs over the set of closed points of the open subscheme $X\setminus Z$ of $X$. By Proposition \ref{Pro: additivity of height function} \ref{Item: additivity of height function}, for any integer $n\geqslant 1$, one has $\widehat{\mu}_{\mathrm{ess}}(L^{\otimes n},n\varphi)=n\operatorname{\widehat{\mu}_{\mathrm{ess}}}(L,\varphi)$.

\begin{prop}\label{Pro: superadditivity of ess min}
\label{Pro: super-additivity of ess min}Let $(L_1,\varphi_1)$ and $(L_2,\varphi_2)$ be adelic line bundles on $X$. Then
\[\widehat{\mu}_{\mathrm{ess}}(L_1\otimes L_2,\varphi_1+\varphi_2)\geqslant\widehat{\mu}_{\mathrm{ess}}(L_1,\varphi_1)+\widehat{\mu}_{\mathrm{ess}}(L_2,\varphi_2).\]
\end{prop}
\begin{proof}
Let $Z_1$ and $Z_2$ be strict Zariski closed subsets of $X$. Then $Z=Z_1\cup Z_2$ is also a strict Zariski closed subset of $X$. If $P$ is an element of $(X\setminus Z)(K^{\mathrm{ac}})$, one has 
\[\begin{split}&\quad\;h_{(L_1\otimes L_2,\varphi_1+\varphi_2)}(P)=h_{(L_1,\varphi_1)}(P)+h_{(L_2,\varphi_2)}(P)\\&\geqslant\inf_{Q_1\in(X\setminus Z_1)(K^{\mathrm{ac}})}h_{(L_1,\varphi_1)}(Q_1)+\inf_{Q_2\in(X\setminus Z_2)(K^{\mathrm{ac}})}h_{(L_2,\varphi_2)}(Q_2).
\end{split}\]
Taking the infimum with respect to $P\in (X\setminus Z)(K^{\mathrm{ac}})$ and then the supremum with respect to $(Z_1,Z_2)$, we obtain that 
\[\widehat{\mu}_{\mathrm{ess}}(L_1\otimes L_2,\varphi_1+\varphi_2)\geqslant \widehat{\mu}_{\mathrm{ess}}(L_1,\varphi_1)+\widehat{\mu}_{\mathrm{ess}}(L_2,\varphi_2).\] 
\end{proof}

\begin{prop}\label{Pro: properties of essential minimum}
Let $(L,\varphi)$ be an adelic line bundle on $X$. 
\begin{enumerate}[label=\rm(\arabic*)]
\item\label{Item: essential minimum as inf} The essential minimum of $(L,\varphi)$ identifies with the infimum of the set of real numbers $C$ such that $\{P\in X(K^{\mathrm{ac}})\,:\,h_{(L,\varphi)}(P)\leqslant C\}$ is Zariski dense in $X$.

\item\label{Item: essential minimum birational} If $X'$ is an integral projective $K$-scheme and $f:X'\rightarrow X$ is a birational projective morphism, then one has $\widehat{\mu}_{\mathrm{ess}}(L,\varphi)=\widehat{\mu}_{\mathrm{ess}}(f^*(L),f^*(\varphi))$.  
\end{enumerate}
\end{prop}
\begin{proof}
\ref{Item: essential minimum as inf} Let $C$ be a real number such that the set \[M_C:=\{P\in X(K^{\mathrm{ac}})\,:\,h_{(L,\varphi)}(P)\leqslant C\}\] is Zariski dense. Then, for any Zariski closed subset $Z$ of $X$ such that $Z\subsetneq X$, the set $M_C$ is not contained in $Z$, namely $(X\setminus Z)(K^{\mathrm{ac}})$ contains at least one element of $M_C$. Therefore one has
\[\inf_{P\in(X\setminus Z)(K^{\mathrm{ac}})}h_{(L,\varphi)}(P)\leqslant C.\]
We then obtain that $\widehat{\mu}_{\mathrm{ess}}(L,\varphi)$ is bounded from above by the infimum of the set of real numbers $C$ such that $M_C$ is Zariski dense in $X$.  

Conversely, if $C$ is a real number such that $M_C$ is not Zariski dense in $X$, then one has
\[\widehat{\mu}_{\mathrm{ess}}(L,\varphi)\geqslant\inf_{P\in(X\setminus \overline M_C^{\mathrm{Zar}})(K^{\mathrm{ac}})}h_{(L,\varphi)}(P)\geqslant C.\]
Since $C$ is arbitrary, we obtain that $\widehat{\mu}_{\mathrm{ess}}(L,\varphi)$ is bounded from below by the infimum of the set of real numbers $C$ such that $M_C$ is Zariski dense in $X$.

\ref{Item: essential minimum birational}  Denote by $(L',\varphi')$ the adelic line bundle $(f^*(L),f^*(\varphi))$, and by  $Y$ the exceptional locus of $f$. Note that the restriction of $f$ to $Z'\setminus f^{-1}(Y)$ is an isomorphism between $X'\setminus f^{-1}(Y)$ and $X\setminus Y$.  
Let $Z$ be a Zariski closed subset of $X$ such that $Z\subsetneq X$. Let $Z'=f^{-1}(Z)$. It is a Zariski closed subset of $X'$ such that $Z'\subsetneq X'$. Therefore, 
\[\widehat{\mu}_{\mathrm{ess}}(L',\varphi')\geqslant
\inf_{P\in(X'\setminus(Z'\cup f^{-1}(Y)))(K^{\mathrm{ac}})}h_{(L',\varphi')}(P)\geqslant\inf_{Q\in (X\setminus Y)}h_{(L,\varphi)}(Q).\]
Since $Y$ is arbitrary, we obtain that $\widehat{\mu}_{\mathrm{ess}}(L',\varphi')\geqslant\widehat{\mu}_{\mathrm{ess}}(L,\varphi)$.

Let $C$ be a real number such that the set $M_C$ of points $Q\in X(K^{\mathrm{ac}})$ with $h_{(L,\varphi)}(Q)\leqslant C$ is Zariski dense. Then the set $M_C\cap(X\setminus Y)(K^{\mathrm{ac}})$ is also Zariski dense in $X$. This implies that $f^{-1}(M_C\cap(X\setminus Y)(K^{\mathrm{ac}}))$ is Zariski dense in $X'$. Note that for any $P\in f^{-1}(M_C\cap(X\setminus Y)(K^{\mathrm{ac}}))$ one has $h_{(L',\varphi')}(P)=h_{(L,\varphi)}(f(P))$. Therefore the set of $P\in X'(K^{\mathrm{ac}})$ with $h_{(L',\varphi')}(P)\leqslant C$ is Zariski dense, which implies that $\widehat{\mu}_{\mathrm{ess}}(L',\varphi')\leqslant C$. Since $C>\widehat{\mu}_{\mathrm{ess}}(L,\varphi)$ is arbitrary, we obtain that $\widehat{\mu}_{\mathrm{ess}}(L',\varphi')\geqslant\widehat{\mu}_{\mathrm{ess}}(L,\varphi)$.   
\end{proof}

\begin{prop}
\label{Pro: lower bound of mu ess} We assume that, either the $\sigma$-algebra $\mathcal A$ is discrete, or the field $K$ admits a countable subfield which is dense in each $K_\omega$, $\omega\in\Omega$. Let $f:X\rightarrow\Spec K$ be an integral projective scheme and $(L,\varphi)$ be an adelic line bundle on $X$. If $s$ is a non-zero global section of $L$, then 
\[\widehat{\mu}_{\mathrm{ess}}(L,\varphi)\geqslant\widehat{\deg}_{f_*(\varphi)}(s).\]
In particular, $\widehat{\mu}_{\mathrm{ess}}(L,\varphi)>-\infty$ once $L$ admits a non-zero global section.
\end{prop}
\begin{proof}
For any closed point $P$ outside of the zero locus of $s$, one has 
\[\begin{split}h_{(L,\varphi)}(P)&=-\int_{\chi\in\Omega_{K^{\mathrm{ac}}}}\ln|s|_{\varphi_{\pi_{K^{\mathrm{ac}}/K}}(\chi)}(\sigma_\chi(P))\,\nu_{K^{\mathrm{ac}}}(\mathrm{d}\chi)\\&\geqslant-\int_{\chi\in\Omega_{K^{\mathrm{ac}}}}\ln\norm{s}_{\varphi_{\pi_{K^{\mathrm{ac}}/K}}(\chi)}\,\nu_{K^{\mathrm{ac}}}(\mathrm{d}\chi)\\&=-\int_{\omega\in\Omega}\ln\norm{s}_{\varphi_\omega}\,\nu(\mathrm{d}\omega)=\widehat{\deg}_{f_*(\varphi)}(s).
\end{split}\]
This leads to the inequality $\widehat{\mu}_{\mathrm{ess}}(L,\varphi)\geqslant\widehat{\deg}_{f_*(\varphi)}(s)$. 
\end{proof}

\begin{prop}\label{Pro: upper bound of mu ess}
Let $(L,\varphi)$ be an adelic line bundle on $X$. One has $\widehat{\mu}_{\mathrm{ess}}(L,\varphi)<+\infty$.
\end{prop}
\begin{proof}
If $\varphi$ and $\varphi'$ are metric families on $L$ such that $(L,\varphi)$ and $(L,\varphi')$ are adelic line bundles on $X$, then for any $P\in X(K^{\mathrm{ac}})$ one has
\[|h_{(L,\varphi)}(P)-h_{(L,\varphi')}(P)|\leqslant \operatorname{dist}(\varphi,\varphi')\]
(see \eqref{eqn:def:global:distance} for the definition of $\operatorname{dist}(\varphi,\varphi')$).
Therefore, to show the proposition, it suffice to prove the assertion for a particular choice of the metric family $\varphi$. This observation allows us to change the metric family whenever necessary in the proof.

Let $M$ be a very ample invertible $\mathcal O_X$-module such that $M\otimes L$ is also very ample. Let $\varphi_M$ be a metric family on $M$ such that $(M,\varphi_M)$ forms an adelic line bundle on $X$. By Proposition \ref{Pro: lower bound of mu ess}, one has $\widehat{\mu}_{\mathrm{ess}}(M,\varphi_M)>-\infty$. Moreover, by  Proposition \ref{Pro: superadditivity of ess min} one has $\widehat{\mu}_{\mathrm{ess}}(L\otimes M,\varphi+\varphi_M)\geqslant\widehat{\mu}_{\mathrm{ess}}(L,\varphi)+\widehat{\mu}_{\mathrm{ess}}(M,\varphi_M)$.
Therefore, by replacing $L$ by $L\otimes M$ we may assume without loss of generality that $L$ is a very ample invertible $\mathcal O_X$-module.

By Noetherian normalisation we obtain that there existe a positive integer $n$, an integral projective $K$-scheme $X'$, a birational projective $K$-morphism $f:X'\rightarrow X$, together with a generically finite projective $K$-morphism $g:X'\rightarrow\mathbb P^r_K$ (where $r$ is the Krull dimension of $X$) such that $g^*(\mathcal O(1))\cong f^*(L^{\otimes n})$, where $\mathcal O(1)$ denotes the universal invertible sheaf on $\mathbb P_K^r$. We can for example construct first a rational morphism from $X$ to $\mathbb P_K^r$ corresponding to an injective finite homogeneous homomorphism from the polynomial algebra to the Cox ring of some power of $L$. This step is guaranteed by the fact that the Cox ring $\bigoplus_{m\in\mathbb N}H^0(X,L^{\otimes m})$ is finitely generated, by using Noether normalisation, see \cite[\S13.1]{Eise}. 
Then we can take $X'$ as the blowing-up of $X$ along the locus where the rational morphism is not defined. By Proposition \ref{Pro: properties of essential minimum} \ref{Item: essential minimum birational}, one has $n\widehat{\mu}_{\mathrm{ess}}(L,\varphi)=n\widehat{\mu}_{\mathrm{ess}}(f^*(L),f^*(\varphi))=\widehat{\mu}_{\mathrm{ess}}(f^*(L^{\otimes n}),nf^*(\varphi))$. Therefore we can reduce the problem to the case where there exists a generically finite projective $K$-morphism $g:X\rightarrow\mathbb P_K^r$ such that $L\cong g^*(\mathcal O(1))$. 

We identify $\mathbb P_K^r$ with $\mathbb P(K^{r+1})$ and equip $K^{r+1}$ with the norm family $\xi$ associated with the canonical basis (see Example \ref{Exa:adelicbasis}). Let $\varphi_0$ be the quotient metric family on $\mathcal O(1)$ induced by $(K^{r+1},\xi)$ and the canonical surjective homomorphism $K^{r+1}\otimes_K\mathcal O_{\mathbb P_K^r}\rightarrow\mathcal O(1)$. As explained above, we may assume without loss of generality that $\varphi=g^*(\varphi_0)$. In particular, for any closed point $P$ of $X$, one has $h_{(L,\varphi)}(P)=h_{(\mathcal O(1),\varphi_0)}(g(P))$. Moreover, similarly as in the proof of Proposition \ref{Pro:Northcott for arithmetic varieties}, for any element $[a_0:\ldots:a_r]\in \mathbb P_K^r(K^{\mathrm{ac}})$, one has
\[h_{(\mathcal O(1),\varphi_0)}([a_0:\ldots:a_r])=\int_{\Omega_{K^{\mathrm{ac}}}} \ln \left( \max \{ |a_0|_{\chi}, \ldots, |a_r|_{\chi} \} \right) \nu_{K^{\mathrm{ac}}}(d\chi).\]
In particular, if $a_0,\ldots,a_r$ are all roots of the unity, then one has $h_{(\mathcal O(1),\varphi_0)}([a_0:\ldots:a_r])=0$. This implies that the set of closed points in $X$ having non-positive height (with respect to $(L,\varphi)$) is Zariski dense. Therefore $\widehat{\mu}_{\mathrm{ess}}(L,\varphi)\leqslant 0$. The proposition is thus proved.
\end{proof}

\subsection{Adelic divisors}
In this subsection we fix a \emph{geometrically integral} projective scheme over $\Spec K$. If $D$ is a Cartier divisor on $X$, for any $\omega\in\Omega$, $D$ induces by base change a Cartier divisor on $X_{\omega}$, which we denote by $D_\omega$. 

Let $D$ be a Cartier divisor on $X$. We call \emph{Green function family}\index{Green function family} of $D$ any family $g=\{g_\omega\}_{\omega\in\Omega}$ parametrised by $\Omega$ such that each $g_\omega$ is a Green function of $D_\omega$ (cf. Subsection \ref{Subsec: Green function}). Note that each Green function $g_\omega$ determines a continuous metric on the invertible sheaf $\mathcal O_{X_{\omega}}(D_\omega)\cong\mathcal O_X(D)\otimes_{\mathcal O_X}\mathcal O_{X_{\omega}}$, which we denote by $\varphi_{g_\omega}$. Thus the collection $\{\varphi_{g_\omega}\}_{\omega\in\Omega}$ forms a metric family on $\mathcal O_X(D)$ which we denote by $\varphi_g$ and called the metric family associated with $g$. We say that the Green function family $g$ is \emph{dominated}\index{dominated}\index{Green function family!dominated} (resp. \emph{measurable}) if the associated metric family $\varphi_g$ is dominated (resp. measurable). In the case where $g$ is dominated and measurable, we say that the couple $(D,g)$ is an \emph{adelic Cartier divisor}\index{adelic Cartier divisor} on $X$. Note that this condition is equivalent to the assertion that $(\mathcal O_X(D),\varphi_g)$ is an adelic line bundle on $X$. In this case we denote by $h_{(D,g)}$ the height function $h_{(\mathcal O_X(g),\varphi_g)}$ on $X(K^{\mathrm{ac}})$.

Let $D$ and $D'$ be Cartier divisors on $X$, $g=\{g_\omega\}_{\omega\in\Omega}$ and $g'=\{g_\omega'\}_{\omega\in\Omega}$ be Green function families of $D$ and $D'$, respectively. We denote by $g+g'$ the Green function family $\{g_\omega+g_\omega'\}_{\omega\in\Omega}$ of $D+D'$. Moreover, we denote by $-g$ the Green function family $\{-g_\omega\}_{\omega\in\Omega}$ of $-D$. Note that, if $(D,g)$ and $(D',g')$ are adelic Cartier divisors then $(D+D',g+g')$ and $(-D,-g)$ are also adelic Cartier divisors. This follows from Propositions \ref{Pro: dominancy preserved by operators} and \ref{Pro: measurability preserved by operators}. Therefore, the set of adelic Cartier divisors forms an abelian group, which we denote by $\widehat{\mathrm{Div}}(X)$.   

\begin{rema}
In the case where the Cartier  divisor $D$ is trivial, a Green function family on $D$ can be considered as a family $\{g_\omega\}_{\omega\in\Omega}$ of continuous real-valued functions, where $g_\omega$ is a continuous function on $X_{\omega}^{\mathrm{an}}$. It is dominated if and only if the function $(\omega\in\Omega)\mapsto\sup_{x\in X_{\omega}^{\mathrm{an}}}|g_\omega|(x)$ is $\nu$-dominated. It is measurable if the following two conditions are satisfied (cf. Definition \ref{Def: measurability of metric family}):
\begin{enumerate}[label=\rm(\alph*)]
\item for any closed point $P$ of $X$, the function $(\omega\in\Omega)\mapsto g_\omega(P)$ is $\mathcal A$-measurable,
\item for any point $x\in X_{\leqslant 1,\mathbb Q}^{\mathrm{an}}$ (where we consider the trivial absolute value on $K$ in the construction of $X^{\mathrm{an}}$), the function $(\omega\in\Omega_0)\mapsto g_\omega(x)$ is $\mathcal A_0$-measurable, where $\Omega_0$ is the set of $\omega\in\Omega$ such that $|\ndot|_\omega$ is trivial. 
\end{enumerate}
The set of all dominated and measurable Green function families on the trivial Cartier divisor forms actually a vector space over $\mathbb R$, which we denote by $\widehat{C}^0(X)$.
\end{rema}

\begin{defi}
Let $\mathbb K$ be either $\mathbb Q$ or $\mathbb R$. We denote by $\widehat{\mathrm{Div}}_{\mathrm{\mathbb K}}(X)$ the $\mathbb K$-vector space $\widehat{\mathrm{Div}}(X)\otimes_{\mathbb Z}\mathbb K$ modulo the vector subspace generated by elements of the form
\[(0,g_1) \otimes \lambda_1 +\cdots+(0,g_n) \otimes \lambda_n -(0,\lambda_1g_1+\cdots+\lambda_ng_n),\]
where $\{g_i\}_{i=1}^n$ is a finite family of elements in $\widehat{C}^0(X)$, and $(\lambda_1,\ldots,\lambda_n)\in\mathbb K^n$. 
In the other words, $\widehat{\mathrm{Div}}_{\mathrm{\mathbb K}}(X)$ consists of
pairs (see \S\ref{Subsec: Green function} for the notation of $C^0_{\mathrm{gen}}(X_\omega^{\mathrm{an}})$) 
\[
(D, \{ g_{\omega} \}_{\omega \in \Omega}) \in \Div_{\mathbb K}(X) \times \prod_{\omega \in \Omega} C^0_{\mathrm{gen}}(X_\omega^{\mathrm{an}})
\]
such that $D = a_1 D_1 + \cdots + a_n D_n$ and $g_{\omega} = a_1 g_{1, \omega} + \cdots + a_n g_{n,\omega}$
for some $(D_1, g_1), \ldots, (D_n, g_n) \in \widehat{\mathrm{Div}}(X)$ and $a_1, \ldots, a_n \in \mathbb K$.
For $\lambda_1, \lambda_2 \in \mathbb K$ and $(D_1, g_1), (D_2, g_2) \in \widehat{\mathrm{Div}}_{\mathrm{\mathbb K}}(X)$,
$\lambda_1 (D_1, g_1) + \lambda_2(D_2, g_2)$ is defined by $(\lambda_1 D_1 + \lambda_2 D_2, \lambda_1 g_1 + \lambda_2 g_2)$.
Note that $\lambda_1 (D_1, g_1) + \lambda_2(D_2, g_2) \in \widehat{\mathrm{Div}}_{\mathbb K}(X)$.
In this sense, $\widehat{\mathrm{Div}}_{\mathbb K}(X)$ forms a vector space over $\mathbb K$.

The elements in $\widehat{\mathrm{Div}}_{\mathbb K}(X)$ are called adelic $\mathbb K$-Cartier divisors on $X$. For any element $\overline D$ written in the form $\lambda_1\overline D_1+\cdots+\lambda_n\overline D_n$ with $(\overline D_1,\ldots,\overline D_n)\in\widehat{\mathrm{Div}}(X)$ and $(\lambda_1,\ldots,\lambda_n)\in\mathbb K^n$, we define a function $h_{\overline D}:X(K^{\mathrm{ac}})\rightarrow\mathbb R$ such that for any $P\in X(K^{\mathrm{ac}})$, 
\[h_{\overline D}(P):=\sum_{i=1}^n\lambda_ih_{\overline D_i}(P).\]
Note that the Proposition \ref{Pro: additivity of height function} \ref{Item: additivity of height function} shows that this map is actually well defined. 
\end{defi}

\begin{rema}\label{Rem: decomposition of adelic divisor}
Let $\overline D$ be an element of $\widehat{\mathrm{Div}}_{\mathbb K}(X)$, which is written in the form $\lambda_1(D_1,g_1)+\cdots+\lambda_n(D_n,g_n)$, where $(\lambda_1,\ldots,\lambda_n)\in\mathbb K^n$, and for any $i\in\{1,\ldots,n\}$, $(D_i,g_i)$ is an element of $\widehat{\mathrm{Div}}(X)$. Then, for any $\omega\in\Omega$, the element $\lambda_1D_{1,\omega}+\cdots+\lambda_nD_{n,\omega}$ of $\mathrm{Div}_{\mathbb K}(X)$ is equal to $D_\omega$, where $D=\lambda_1D_1+\cdots+\lambda_nD_n\in\mathrm{Div}_{\mathbb K}(X)$. Moreover, assume that $g_i$ is written in the form $\{g_{i,\omega}\}_{\omega\in\Omega}$, where $g_{i,\omega}$ is a Green function of $D_{i,\omega}$. Then for any $\omega\in\Omega$, the element $\lambda_1g_{1,\omega}+\cdots+\lambda_ng_{n,\omega}$ is a Green function of the $\mathbb K$-Cartier divisor $D_\omega$, which does not depend on the choice of the decomposition $\overline D=\lambda_1(D_1,g_1)+\cdots+\lambda_n(D_n,g_n)$. Thus we can write $\overline D$ in the form $(D,g)$, where $D$ is a $\mathbb K$-Cartier divisor of $X$ and $g$ is a family of Green functions of the form $\{g_\omega\}_{\omega\in\Omega}$, with $g_\omega$ being a Green function of $D_\omega$. Note that the measurability of the Green function families $g_1,\ldots,g_n$ implies the following statements:
\begin{enumerate}[label=\rm(\alph*)]
\item for any closed point $P$ of $X$ outside of the support of $D$, the function $(\omega\in\Omega)\mapsto g_\omega(P)$ is well defined and is $\mathcal A$-measurable,
\item for any point $x\in X_{\leqslant 1,\mathbb Q}^{\mathrm{an}}$ outside of the analytification of the support of $D$, the function $(\omega\in\Omega_0)\mapsto g_\omega(x)$ is well defined and is $\mathcal A_0$-measurable.
\end{enumerate} 
Moreover, if $D$ belongs to $\mathrm{Div}(X)$, then $g$ is a dominated Green function family of $D$. This statement results directly from the following proposition.
\end{rema}

\begin{exem}
Let $s$ be a non-zero rational function on $X$. For any $\omega\in\Omega$, we consider $s$ as a non-zero rational function on $X_{\omega}$. Note that $-\ln|s|_\omega$ is  a Green function of the principal Cartier divisor $\mathrm{div}(s)$. Note that the Green function family $\{-\ln|s|_\omega\}_{\omega\in\Omega}$ is measurable and dominated since the corresponding metric family on $\mathcal O_X(\mathrm{div}(f))\cong\mathcal O_X$ is trivial. Thus \[(s\in K(X)^{\times})\longmapsto \widehat{\mathrm{div}}(s):=(\mathrm{div}(s),\{-\ln|s|_\omega\}_{\omega\in\Omega})\]
defines a morphism of groups from $(K(X)^{\times},\times)$ to $\widehat{\mathrm{Div}}(X)$. The adelic Cartier divisors belonging to the image of this morphism are called \emph{principal adelic Cartier divisors}\index{principal adelic Cartier divisor}\index{adelic Cartier divisor!principal ---}. Moreover, for $\mathbb K\in\{\mathbb Q,\mathbb R\}$ this morphism induces a $\mathbb K$-linear map $\widehat{\mathrm{div}}_{\mathbb K}:K(X)^{\times}\otimes_{\mathbb Z}\mathbb K\rightarrow\widehat{\mathrm{Div}}_{\mathbb K}(X)$ sending to $s_1^{\lambda_1}\cdots s_n^{\lambda_n}$ to $\lambda_1\widehat{\mathrm{div}}(s_1)+\cdots+\lambda_n\widehat{\mathrm{div}}(s_n)$. The adelic $\mathbb K$-Cartier divisors belonging to the image of this $\mathbb K$-linear map are said to be \emph{principal}. 
\end{exem}

Let $(D,g)$ be an adelic $\mathbb K$-Cartier divisor on $S$.
For $\phi \in H^0_{\mathbb K}(X, D)$, $|\phi|_{\omega} \exp(-g_{\omega})$ extends to a continuous function on $X_{\omega}^{\mathrm{an}}$ by
Proposition~\ref{Pro:e-gextension}. We denote by 
\[\|\phi\|_{g_{\omega}}:=\sup_{x \in X_{\omega}^{\mathrm{an}}} \left\{ \left(|\phi|_{\omega} \exp(-g_{\omega})\right)(x) \right\}.\]

\begin{prop}\label{prop:adelic:R:Cartier:div:integrable} We assume that, either the $\sigma$-algebra $\mathcal A$ is discrete,
or the field $K$ admits a countable subfield which is dense in every $K_\omega$, $\omega\in\Omega$. Let $(D,g)$ be an adelic $\mathbb K$-Cartier divisor on $S$ and $\phi\in H^0_{\mathbb K}(X,D)$. The function on $\Omega$ given by
\[
(\omega \in \Omega) \mapsto \ln \|\phi\|_{g_{\omega}} = \sup_{x \in X_{\omega}^{\mathrm{an}}} \left\{ (- g_{\omega} + \log |\phi|_{\omega})(x) \right\}
\]
is $\nu$-integrable.
\end{prop}

\begin{proof}
Note that $D' := D + (\phi) \geqslant_{\mathbb K} 0$,
$g'_{\omega} := g_{\omega} - \log|\phi|_{\omega}$ is a Green function of $D'_{\omega}$
and 
$|\phi|_{g_{\omega}} = |1|_{g'_{\omega}}$ on $X_{\omega}^{\mathrm{an}}$, so that we may assume that
$D $ is $\mathbb K$-effective and $\phi = 1$. 

Let $X'$ be the normalisation of $X$. Since $X$ and $X'$ have the same function field,
$X'$ is also geometrically integral over $K$. Moreover, let $D'$ (resp. $g'_{\omega}$) be the pull-back of $D$ by $X' \to X$ 
(resp. $X'_{\omega} \to X_{\omega}$). Then $g' = \{ g'_{\omega} \}_{\omega \in \Omega}$ is a family of Green functions of $D'$ over $S$.
Note that $\| 1 \|_{g_{\omega}} = \| 1 \|_{g'_{\omega}}$, so that we may further assume that $X$ is normal.

First we consider the case $\mathbb K = \mathbb Q$.
Then there is a positive integer $N$ such that $ND$ is a Cartier divisor.
Then $\phi^N \in H^0(X, ND)$ and $\omega \mapsto \ln \|\phi^N \|_{Ng_{\omega}}$ is integrable on $\Omega$
by Theorem~\ref{Thm: Fubini-Study dominated} and Theorem~\ref{Thm: measurability of linear series}.
Note that $\ln \|\phi^N \|_{Ng_{\omega}} = N \ln \|\phi \|_{g_{\omega}}$, so that
$\omega \mapsto \ln \|\phi \|_{g_{\omega}}$ is also integrable on $\Omega$.

Next we consider the case $\mathbb K = \mathbb R$.
By Proposition~\ref{prop:normal:R:effective:sum}, 
there are effective Cartier divisors $D_1, \ldots, D_r$ and $a_1, \ldots, a_r \in \mathbb R_{\geqslant 0}$
such that $D = a_1 D_1 + \cdots + a_r D_r$. We choose a family of Green functions $g_i = \{ g_{i,\omega} \}_{\omega \in \Omega}$
of $D_i$ over $S$ such that $(D_i, g_i)$ is an adelic Cartier divisor over $S$ for each $i$ and
\[
(D, g) = (a_1D_1 + \cdots + a_rD_r, a_1g_1 + \cdots + a_r g_r).
\]
If we set $\psi_{i}(\omega) = \ln \| 1 \|_{g_{i,\omega}}$ and
$g'_{i,\omega} := g_{i,\omega} + \psi_i(\omega)$
for $i=1, \ldots, r$,
then $\psi_{i}$ is integrable on $\Omega$ and 
\[
\| 1 \|_{g'_{i,\omega}} = \| 1 \|_{g_{i,\omega}}\exp(-\psi_i(\omega)) = 1,
\] 
so that $g'_{i,\omega} \geqslant 0$ for all $i$ and $\omega$.
Note that if we set $g' = a_1 g'_1 + \cdots +  a_n g'_n$, then
\[
\ln \| 1 \|_{g',\omega} = \ln \| 1 \|_{g,\omega} - (a_1 \psi_{1}(\omega) +
\cdots + a_n \psi_{n}(\omega)).
\]
Therefore, we may assume that $g_{i,\omega} \geqslant 0$ for all $i$ and $\omega$.

For each $i$, we choose a sequence $\{ a_{i,n} \}_{n=1}^{\infty}$ of non-negative rational numbers such that
\[
0 \leqslant a_{i,n} - a_i \leqslant \frac{a_i}{n}\quad\text{and}\quad a_{i, n+1} \leqslant a_{i, n}
\]
for all $n$. We set 
\[
(D_n, h_n) := (a_{1,n}D_1 + \cdots + a_{r,n}D_r, a_{1,n} g_1 + \cdots + a_{r,n} g_r).
\]
Then $D_n$ is effective and
\[
-h_{n, \omega} \leqslant -g \leqslant \frac{n}{n+1}(- h_{n, \omega}) \leqslant 0\quad\text{and}\quad
-h_{n,\omega} \leqslant -h_{n+1,\omega}
\]
for all $n$ and $\omega$. If we set \[A(\omega) = \sup_{x \in X_{\omega}^{\mathrm{an}}} \{ -g_{\omega}(x) \}\,\text{ and }\,
A_n(\omega) = \sup_{x \in X_{\omega}^{\mathrm{an}}} \{ -h_{n,\omega}(x) \}, \]
then
\[
A_n(\omega) \leqslant A(\omega) \leqslant \frac{n}{n+1} A_n(\omega) \leqslant 0 \quad\text{and}\quad
A_{n}(\omega) \leqslant A_{n+1}(\omega)
\]
for all $n$ and $\omega$. Thus $\lim_{n\to\infty} A_n(\omega) = A(\omega)$ and $A_n(\omega) \leqslant A(\omega) \leqslant 0$.
Note that $\omega \mapsto A_n(\omega)$ is integrable for all $n$.
Therefore, by monotone convergence theorem,
$A(\omega)$ is integrable.
\end{proof}

\begin{coro}\label{coro:adelic:R:Cartier:div:integrable}We keep the hypothesis of Proposition \ref{prop:adelic:R:Cartier:div:integrable}.
Let $(D, g)$ be an adelic $\mathbb K$-Cartier divisor on $X$. Let $\phi \in K(X)^{\times} \otimes_{\mathbb Z} \mathbb K$ such that
$D + (\phi) \geqslant_{\mathbb K} 0$. Then
the function \[(\omega\in\Omega)\longmapsto \ln \| \phi \|_{g_{\omega}} = \sup_{x \in X_{\omega}^{\mathrm{an}}} \{ (-g_{\omega} + \log |\phi|_{\omega})(x) \}\] is $\nu$-integrable.
\end{coro}

\begin{proof}
If we set $D' = (D) + (\phi)$ and $g'_{\omega} = g_{\omega} - \ln |\phi|_{\omega}$, then $(D', g' = \{ g'_{\omega} \}_{\omega \in \Omega})$ is
an adelic $\mathbb K$-Cartier divisor on $X$. Thus the assertion follows from Proposition~\ref{prop:adelic:R:Cartier:div:integrable}.
\end{proof}

\begin{coro}\label{Coro: Green family on the trivial divisor} We keep the hypothesis of Proposition \ref{prop:adelic:R:Cartier:div:integrable}.
Let $(\mathbf{0},g)$ be an adelic $\mathbb K$-Cartier divisor on $X$ whose underlying $\mathbb K$-Cartier divisor is trivial. Assume that $g$ is written in the form $\{g_\omega\}_{\omega\in\Omega}$, where $g_\omega$ is considered as a continuous function on $X_{\omega}^{\mathrm{an}}$. Then the function \[(\omega\in\Omega)\longmapsto\sup_{x\in X_{\omega}^{\mathrm{an}}}|g_\omega(x)|\] is $\nu$-integrable.
\end{coro}

For any $\overline D\in\widehat{\mathrm{Div}}_{\mathbb K}(X)$, we define the \emph{essential minimum}\index{essential minimum} of $\overline D$ as
\[\widehat{\mu}_{\mathrm{ess}}(\overline D):=\sup_{Z\subsetneq X}\inf_{P\in(X\setminus Z)(K^{\mathrm{ac}})}h_{\overline D}(P),\]
where $Z$ runs over the set of all strict Zariski closed subsets of $X$, and $P$ runs over the set of closed points of the open subscheme $X\setminus Z$ of $X$. It turns out that the analogue of Proposition \ref{Pro: superadditivity of ess min} and  Proposition \ref{Pro: properties of essential minimum} \ref{Item: essential minimum as inf} holds for adelic $\mathbb K$-Cartier divisors (with essentially the same proof). We resume these statements as follows.

\begin{prop}\label{Pro: super additivity mu ess}
Let $\overline D$ be an adelic $\mathbb K$-Cartier divisor on $X$. Then $\widehat{\mu}(\overline D)$ identifies with the infimum of the set of real numbers $C$ such that $\{P\in X(K^{\mathrm{ac}})\,:\,h_{\overline D}(P)\leqslant C\}$ is Zariski dense in $X$. Moreover, if $\overline D_1$ and $\overline D_2$ are adelic $\mathbb K$-Cartier divisors on $X$, then $\widehat{\mu}_{\mathrm{ess}}(\overline D_1+\overline D_2)\geqslant\widehat{\mu}_{\mathrm{ess}}(\overline D_1)+\widehat{\mu}_{\mathrm{ess}}(\overline D_2)$.  
\end{prop} 

Similarly as in the case of  adelic line bundles, the essential minimum of adelic $\mathbb K$-Cartier divisors never takes the value of $+\infty$.
\begin{prop}\label{Pro: upper bound ess min 2}
Let $\overline D$ be an adelic $\mathbb K$-Cartier divisor on $X$. One has $\widehat{\mu}_{\mathrm{ess}}(\overline D)<+\infty$.
\end{prop}
\begin{proof}
Assume that $\overline D$ is written in the form $\overline D=\lambda_1\overline D_1+\cdots+\lambda_n\overline D_n$, where $D_1,\ldots,D_n$ are very ample Cartier divisors on $X$ and $(\lambda_1,\ldots,\lambda_n)\in\mathbb K^n$. By Proposition \ref{Pro: lower bound of mu ess}, for any $i\in\{1,\ldots,n\}$ one has $\widehat{\mu}_{\mathrm{ess}}(\overline D_i)>-\infty$. We choose $(\lambda_1',\ldots,\lambda_n')\in (\mathbb K\cap\mathbb R_{>0})^n$ such that $\lambda_i+\mu_i\in\mathbb Z$ for any $i\in\{1,\ldots,n\}$. Let 
\[\overline E:=\sum_{i=1}^n(\lambda_i+\lambda_i')\overline D_i.\] By Proposition \ref{Pro: upper bound of mu ess}, one has
$\widehat{\mu}_{\mathrm{ess}}(\overline E)<+\infty$. Moreover, by Proposition \ref{Pro: super additivity mu ess}, one has
\[\widehat{\mu}_{\mathrm{ess}}(\overline E)\geqslant\widehat{\mu}_{\mathrm{ess}}(\overline D)+\sum_{i=1}^n\lambda_i'\,\widehat{\mu}_{\mathrm{ess}}(\overline D_i).\]
Since $\widehat{\mu}_{\mathrm{ess}}(\overline D_i)>-\infty$ and $\lambda_i'>0$ for any $i\in\{1,\ldots,n\}$, we deduce that $\widehat{\mu}_{\mathrm{ess}}(\overline D)<+\infty$.
\end{proof}

\begin{defi}\label{def:norm:family:adelic:R:div}
We assume that $X$ is normal. 
Let $\overline D=(D,g)$ be an adelic $\mathbb K$-Cartier divisor, where the Green function  family $g$ is written in the form $\{g_\omega\}_{\omega\in\Omega}$. 
For $\omega \in \Omega$, let $X'_{\omega}$ be the normalization of $X_{\omega}$ and
$D'_{\omega}$ (resp. $g'_{\omega}$) be the pull-back of $D$ by $X'_{\omega} \to X_{\omega}$ (resp. the pull-back of $g_{\omega}$ by ${X'_{\omega}}^{\mathrm{an}} \to X_{\omega}^{\mathrm{an}}$). By using the natural injective homomorphism $H^0_{\mathbb K}(X, D) \otimes_{K} K_{\omega} \to H^0_{\mathbb K}(X'_{\omega}, D'_{\omega})$ and
$g'_{\omega}$, one has a norm $\|\ndot\|_{g_{\omega}}$ on $H^0_{\mathbb K}(X, D) \otimes_{K} K_{\omega}$ (cf. Definition \ref{Def: H0 metric}).
The norm family $\{\norm{\ndot}_{g_\omega}\}_{\omega\in\Omega}$ is denoted by $\xi_g$.
\end{defi}

\begin{theo} We assume that, either the $\sigma$-algebra $\mathcal A$ is discrete, 
or the field $K$ admits a countable subfield which is dense in every $K_\omega$, $\omega\in\Omega$.
Suppose that $X$ is normal. Then
the couple $(H^0_{\mathbb K}(D),\xi_g)$ is a strictly adelic vector bundle on $S$.
\end{theo}
\begin{proof}
The measurability of $\xi_g$ is a consequence of Proposition~\ref{prop:adelic:R:Cartier:div:integrable}.
Let us consider the dominancy of $\xi_g$.
By using \cite[Lemma~5.2.3]{Moriwaki2012},
$\overline D$ is written in the form 
\[
\lambda_1(D_1,g_1)+\cdots+\lambda_n(D_n,g_n),
\]
where $(D_i,g_i)$'s are elements of $\widehat{\mathrm{Div}}(X)$ such that $D_1,\ldots,D_n$ are effective, and $(\lambda_1,\ldots,\lambda_n)\in\mathbb K^n$. 
Let $(\lambda'_1,\ldots,\lambda'_n)$ be an element of $(\mathbb K\cap\mathbb R_{>0})^n$ such that  $\lambda_i+\lambda_i'\in\mathbb Z_{>0}$ for any $i\in\{1,\ldots,n\}$. 
Let 
\[
(D',g'):=(\lambda_1+\lambda_1')(D_1,g_1)+\cdots+(\lambda_n+\lambda_n')(D_n,g_n),
\]
which is viewed as an adelic Cartier divisor on $X$. Since $D_i$ is effective, we obtain that $1$ belongs to $H^0(D_i)$. Moreover, by Proposition~\ref{prop:adelic:R:Cartier:div:integrable}, the function $(\omega\in\Omega)\mapsto \ln\norm{1}_{g_{i,\omega}}$ is $\nu$-integrable.

Let $\boldsymbol{e}=\{e_i\}_{i=1}^m$ be a basis of $H^0_{\mathbb K}(D)$. We complete it into a basis $\boldsymbol{e}'=\{e_i\}_{i=1}^r$ of $H^0_{\mathbb K}(D')$. By Theorem \ref{Thm: Fubini-Study dominated}, the norm family $\xi_{g'}:=\{\norm{\ndot}_{g_\omega'}\}_{\omega\in\Omega}$ is strongly dominated,
so that, by Corollary \ref{Cor:dominatedanddist}, the local distance function $(\omega\in\Omega)\mapsto d_\omega(\xi_{g'},\xi_{\boldsymbol{e}'})$ is $\nu$-dominated. 
Further, by Proposition~\ref{prop:adelic:R:Cartier:div:integrable}, the function $(\omega\in\Omega)\mapsto \ln\norm{e_i}_{g_{\omega}}$ is $\nu$-integrable for each $i$.

Let $\omega\in\Omega$ and $(a_1,\ldots,a_m)\in K_\omega^m$, one has 
\begin{multline*}
\ln\norm{a_1e_1+\cdots+a_me_m}_{g_\omega} \leqslant \max_{i\in\{1,\ldots,m\}}\{ \ln |a_i| + \ln \norm{e_i}_{g_\omega}\}+\indic_{\Omega_\infty}(\omega)\ln(m) \\
\leqslant \ln \| a_1 e_1 + \cdots + a_m e_m \|_{\xi_{\boldsymbol{e}}} + \max_{i\in\{1,\ldots,m\}}\{ \ln \norm{e_i}_{g_\omega}\} +\indic_{\Omega_\infty}(\omega)\ln(m).
\end{multline*}
Moreover,
\begin{multline*}
\ln\norm{a_1e_1+\cdots+a_me_m}_{g_\omega}\geqslant\ln\norm{a_1e_1+\cdots+a_me_m}_{g_\omega'}-\sum_{i=1}^n\lambda_i'\ln\norm{1}_{g_{i,\omega}}\\
\geqslant\ln\norm{a_1e_1+\cdots+a_me_m}_{\xi_{\boldsymbol{e}}}-d_\omega(\xi_{g'},\xi_{\boldsymbol{e}'})-\sum_{i=1}^n\lambda_i'\ln\norm{1}_{g_{i,\omega}},
\end{multline*}
and hence one obtains
\begin{multline*}
d_{\omega}(\xi_g, \xi_{\boldsymbol{e}}) \\
\leqslant \max \left\{ \max_{i\in\{1,\ldots,m\}}\left\{ \big|\ln \norm{e_i}_{g_\omega}\big|\right\} +\indic_{\Omega_\infty}(\omega)\ln(m), d_\omega(\xi_{g'},\ \xi_{\boldsymbol{e}'})+\sum_{i=1}^n\lambda_i' \big|\ln\norm{1}_{g_{i,\omega}}\big|\right\}
\end{multline*}
Therefore the local distance function $(\omega\in\Omega)\mapsto d_\omega(\xi_g,\xi_{\boldsymbol{e}})$ is $\nu$-dominated, which implies that the norm family $\xi_g$ is strongly dominated (cf. Corollary \ref{Cor:dominatedanddist}).
\end{proof}

\section{Okounkov bodies and concave transform}

\subsection{Reminder on some facts about convex sets}
In this subsection, we recall some basic facts about convex sets in finite-dimensional vector spaces, which will be used in the subsequening subsections.

\begin{prop}\label{Pro:interiorofconve}
Let $V$ be a finite-dimensional vector space over $\mathbb R$. Suppose that $C_1$ and $C_2$ are two convex subsets of $V$ which have the same closure in $V$, then the interiors $C_1^\circ$ and $C_2^\circ$ are also the same. 
\end{prop}
\begin{proof}
It suffices to prove that, if $C$ is a convex subset of $V$, then the interior of the closure $\overline C$ coincides with the interior $C^\circ$ of $C$. Let $x$ be an element of $V$. If $x$ does not lie in  $C^\circ$, then (by Hahn-Banach theorem, see \cite[Theorem 3.4]{Rudin73}) there exists an affine function $q$ on $V$ such that the restriction of $q$ on $C^\circ$ is non-negative but $q(x)\leqslant 0$. Since the set $\{y\in V\,:\,q(y)\geqslant 0\}$ is closed, it contains $\overline C$. Moreover, the interior of this set is $\{y\in V\,:\,q(y)>0\}$, which contains $(\overline C)^\circ$. Hence $x$ cannot lie in $(\overline C)^\circ$. Therefore one has $C^\circ\supseteq(\overline{C})^\circ$, which actually implies the equality of these two sets.  
\end{proof} 

\begin{prop}
Let $V$ be a finite-dimensional vector space over $\mathbb R$ and $(C_i)_{i\in I}$ be a family of convex subsets of $W$. Suppose that the family $(C_i)_{i\in I}$ is filtered, namely, for any couple $(i_1,i_2)$ of indices in $I$, there exists $j\in I$ such that $C_{i_1}\cup C_{i_2}\subseteq C_j$. Let $C$ be the union of $C_i$, $i\in I$. Then the interior of $C$ identifies with the union of $C_i^\circ$, $i\in I$. 
\end{prop}
\begin{proof}
Since the family $(C_i)_{i\in I}$ is filtered, for any couple of points $(x,y)$ in $C$, there exists an index $i\in I$ such that $\{x,y\}\subseteq C_i$. Therefore $C$ is a convex subset of $V$. As a consequence, for any point $x$ of the interior $C^\circ$, there exists points $x_1,\ldots,x_n$ in $C$ such that the point $x$ is contained in the interior of the convex hull of $x_1,\ldots,x_n$. Still by the assumption that the family $(C_i)_{i\in I}$ is filtered, there exists $j\in I$ such that $\{x_1,\ldots,x_n\}\subseteq C_j$. Hence one has $x\in C_j^\circ$. 

\end{proof}

\subsection{Graded semigroups}

Let $V$ be a finite-dimensional vector space over $\mathbb R$. By \emph{graded semigroup}\index{graded semigroup} in $V$ we refer to a non-empty subset $\Gamma$ of $\mathbb N_{\geqslant 1}\times V$ which is stable by addition. If $\Gamma$ is a graded semigroup in $V$, for any $n\in\mathbb N_{\geqslant 1}$ we denote by $\Gamma_n$ the projection of $\Gamma\cap(\{n\}\times V)$ in $V$. Let $\mathbb N(\Gamma)$ be the set of all $n\in\mathbb N_{\geqslant 1}$ such that $\Gamma_n$ is non-empty. This is a non-empty sub-semigroup of $\mathbb N_{\geqslant 1}$. We denote by $\mathbb Z(\Gamma)$ the subgroup of $\mathbb Z$ generated by $\mathbb N(\Gamma)$.

\begin{prop}
Let $\Gamma$ be a graded semigroup in $V$. Then there exist at most finitely many positive elements of 
$\mathbb Z(\Gamma) \setminus \mathbb N(\Gamma)$.
\end{prop}
\begin{proof}
The group $\mathbb Z(\Gamma)$ is  non-zero since $\Gamma$ is not empty. Hence there exists a positive integer $m$ such that $\mathbb Z(\Gamma)=m\mathbb Z$. Assume that $m$ is written in the form \[m=a_1n_1+\cdots+a_{\ell}n_{\ell},\] where $n_1,\ldots,n_{\ell}$ are elements in $\mathbb N(\Gamma)$ and $a_1,\ldots,a_{\ell}$ are integers. Since $\mathbb N(\Gamma)\subseteq \mathbb Z(\Gamma)$, there exists a positive integer $N$ such that $n_1+\cdots+n_\ell=mN$. Let \[b=N\cdot\max_{i\in\{1,\ldots,\ell\}}|a_i|.\]
We claim that $mn\in\mathbb N(\Gamma)$ for any $n\geqslant Nb$. In fact, we can write such $n$ in the form $n=cN+r$ where $c\in\mathbb N_{\geqslant b}$ and $r\in\{0,\ldots,N-1\}$. Thus
\[mn=cmN+mr=c(n_1+\cdots+n_\ell)+r(a_1n_1+\cdots+a_{\ell}n_{\ell})=(c+ra_1)n_1+\cdots+(c+ra_{\ell})n_{\ell}.\]
Since $c\geqslant b$ and $r<N$, we obtain that $c+ra_i\geqslant 0$ for any $i\in\{1,\ldots,\ell\}$. Hence $mn\in\mathbb N(\Gamma)$.
\end{proof}

\begin{defi}
Let $\Gamma$ be a graded semigroup in $V$. We denote by $\Delta(\Gamma)$ the closure of the set
\[\bigcup_{n\in\mathbb N,\,n\geqslant 1}\{n^{-1}\alpha\,:\,\alpha\in\Gamma_n\}\subset V.\]
\end{defi}

\begin{prop}\label{Pro:Okounkov body of graded semigroup}
Let $\Gamma$ be a graded semigroup in $V$. The set $\Delta(\Gamma)$ is a closed convex subset of $V$.
\end{prop}
\begin{proof}
It suffices to prove the convexity of the set $\Delta(\Gamma)$. Observe that, if $n$ and $m$ are two positive integers, $\alpha$ and $\beta$ are elements of $\Gamma_n$ and $\Gamma_m$, respectively. We show that, for any $\epsilon\in [0,1]\cap\mathbb Q$, one has $\epsilon n^{-1}\alpha+(1-\epsilon)m^{-1}\beta\in\Delta(\Gamma)$. Let $\epsilon=p/q$ be a rational number in $[0,1]$, where $q\in\mathbb N_{\geqslant 1}$. One has
\[\epsilon n^{-1}\alpha+(1-\epsilon)m^{-1}\beta=\frac{p}{qn}\alpha+\frac{q-p}{qm}\beta=(qmn)^{-1}(pm\alpha+(q-p)n\beta).\]
Since $\alpha\in\Gamma_n$ and $\beta\in\Gamma_m$, one has $pm\alpha+(q-p)n\beta\in\Gamma_{qmn}$. Therefore \[\epsilon n^{-1}\alpha+(1-\epsilon)m^{-1}\beta\in\Delta(\Gamma).\]

Let $H$ be the set \[\bigcup_{n\in\mathbb N,\,n\geqslant 1}\{n^{-1}\alpha\,:\,\alpha\in\Gamma_n\}.\]
Let $x$ and $y$ be two points in $\Delta(\Gamma)$, and $\epsilon\in [0,1]$. By definition, there exists two sequences $\{x_n\}_{n\in\mathbb N}$ and $\{y_n\}_{n\in\mathbb N}$ in $H$ such that 
\[\lim_{n\rightarrow+\infty}x_n=x,
\quad\lim_{n\rightarrow+\infty}y_n=y.\]
Let $\{\epsilon_n\}_{n\in\mathbb N}$ be a sequence in $[0,1]\cap\mathbb Q$ which converges to $\epsilon$. By what we have shown above, for any $n\in\mathbb N$ one has $\epsilon_nx_n+(1-\epsilon_n)y_n\in H$. Moreover, one has
\[\lim_{n\rightarrow+\infty}\epsilon_nx_n+(1-\epsilon_n)y_n=\epsilon x+(1-\epsilon)y.\]
Therefore $\epsilon x+(1-\epsilon)y\in\Delta(\Gamma)$.
\end{proof}

Let $\Gamma$ be a graded semigroup in $V$. We denote by $\Gamma_{\mathbb R}$ the $\mathbb R$-vector subspace of $\mathbb R\times V$ generated by $\Gamma$. 
For $n \in \mathbb Z$, let $A(\Gamma)_n$ be the projection of $\Gamma_{\mathbb R}\cap(\{n\}\times V)$ to $V$.
Especially, $A(\Gamma)_1$ is denoted by $A(\Gamma)$. Note that $A(\Gamma)_0$ is a vector subspace of $V$, which is a translation of the affine subspace $A(\Gamma)$.
Since $A(\Gamma)_n$ 
is the image of an affine subspace of $\mathbb R\times V$ by a linear map, it is an affine subspace in $V$. Note that any element in $A(\Gamma)=A(\Gamma)_1$ can be written in the form
\[\lambda_1\gamma_1+\cdots+\lambda_{\ell}\gamma_{\ell},\]
where for $i\in\{1,\ldots,\ell\}$, $\gamma_i\in\Gamma_{n_i}$, $n_i\in\mathbb N$, $n_i\geqslant 1$, and $(\lambda_1,\ldots,\lambda_\ell)$ is an element in $\mathbb R^{\ell}$ such that $\lambda_1n_1+\cdots+\lambda_\ell n_\ell=1$.
We denote by $\Gamma_{\mathbb Z}$ the subgroup of $\mathbb R\times V$ generated by $\Gamma$. For any $n\in\mathbb Z$, $n\geqslant 1$, let $\Gamma_{\mathbb Z,n}$ be the image of $\Gamma_{\mathbb Z}\cap (\{n\}\times V)$ in $V$ by the canonical projection. Note that $\Gamma_{\mathbb Z,n}$ is non-empty if and only if $n\in\mathbb Z(\Gamma)$.

\begin{prop}\label{Pro:proprietegradedsemigropu}
Let $\Gamma$ be a graded semigroup in $V$. We assume that $\Gamma_{\mathbb Z}$ is a discrete subset of $\mathbb R\times V$.
\begin{enumerate}[label=\rm(\arabic*)]
\item\label{Item: Gamma Z 0 a lattice} The set $\Gamma_{\mathbb Z,0}$ is a lattice in $A(\Gamma)_0$. 
\item\label{Item: Gamma n to Gamma n' bijection} For any $n,n'\in\mathbb Z(\Gamma)$ and any $\gamma_0\in\Gamma_{\mathbb Z,n}$, the map from $\Gamma_{\mathbb Z,n'}$ to $\Gamma_{\mathbb Z,n+n'}$, sending $\gamma\in\Gamma_{\mathbb Z,n'}$ to $\gamma+\gamma_0$, is a bijection. 
\item\label{Item:intersection K gamma n} For any convex and compact subset $K$ of $A(\Gamma)$ which is contained in the relative interior of $\Delta(\Gamma)$, one has
\begin{equation}\label{Equ:intersectionK}K\cap \{n^{-1}\gamma\,:\,\gamma\in\Gamma_n\}=K\cap\{n^{-1}\gamma\,:\,\gamma\in\Gamma_{\mathbb Z,n}\}\end{equation}
for sufficiently positive $n\in\mathbb N(\Gamma)$.
\end{enumerate}
\end{prop}
\begin{proof}
\ref{Item: Gamma Z 0 a lattice} Let $n$ be an element in $\mathbb Z(\Gamma)$ and $\gamma_0\in\Gamma_{\mathbb Z,n}$. By definition, an element $x\in V$ lies in $A(\Gamma)_0$ if and only if $x+n^{-1}\gamma_0\in A(\Gamma)$. In other words, $A(\Gamma)_0$ is precisely the vector subspace of $V$ of all vectors $\gamma$ which can be written in the form
\begin{equation}\label{Equ:gammaaslineacomb}\gamma=\lambda_1\gamma_1+\cdots+\lambda_\ell\gamma_\ell,\end{equation}
where for any $i\in\{1,\ldots,\ell\}$, $\gamma_i\in\Gamma_{n_i}$ with $n_i\in\mathbb N(\Gamma)$, and $(\lambda_1,\ldots,\lambda_\ell)$ is an element in $\mathbb R^\ell$ such that $\lambda_1n_1+\cdots+\lambda_\ell n_\ell=0$. Note that the set $\Gamma_{\mathbb Z,0}$ is characterized by the same condition, except that $(\lambda_1,\ldots,\lambda_\ell)$ is required to be in $\mathbb Z^\ell$ 
Therefore $\Gamma_{\mathbb Z,0}$ is a subset (and hence a subgroup) of $A(\Gamma)_0$. Moreover, we can also rewrite \eqref{Equ:gammaaslineacomb} as
\[\gamma=\frac{\lambda_1}{n}(n\gamma_1-n_1\gamma_0)+\cdots+\frac{\lambda_\ell}{n}(n\gamma_\ell-n_\ell\gamma_0).\] 
Since $n\gamma_i-n_i\gamma_0$ belongs to $\Gamma_{\mathbb Z,0}$ for $i\in\{1,\ldots,\ell\}$, we obtain that $A(\Gamma)_0$ is generated by $\Gamma_{\mathbb Z,0}$ as a vector space over $\mathbb R$. Moreover, since $\Gamma_{\mathbb Z}$ is a discrete subspace of $\mathbb R\times V$, the set $\Gamma_{\mathbb Z,0}\subseteq V$ is also discrete. Hence it forms a lattice in $A(\Gamma)_0$.

\ref{Item: Gamma n to Gamma n' bijection} This comes from the definition of $\Gamma_{\mathbb Z}$. In particular, the inverse map is given by $(\gamma'\in\Gamma_{\mathbb Z,n+n'})\mapsto\gamma'-\gamma_0$.

\ref{Item:intersection K gamma n} Let $\Theta$ be the family of all sub-semigroups of $\Gamma$ which are finitely generated. The family of convex sets $\{\Delta(\Gamma')\}_{\Gamma'\in\Theta}$ is filtered. Let $C$ be the union of all $\Delta(\Gamma')$, $\Gamma'\in\Theta$. By definition, the closure of $C$ coincides with $\Delta(\Gamma)$. Therefore (by Proposition \ref{Pro:interiorofconve}), the  interior of $\Delta(\Gamma)$ relatively to $A(\Gamma)$  identifies with that of $C$, which is equal to $\bigcup_{\Gamma'\in\Theta}\Delta(\Gamma')^\circ$, where $\Delta(\Gamma')^\circ$ denotes the relative interior of $\Delta(\Gamma')$ in $A(\Gamma)$. Since $K$ is a compact subset of $\Delta(\Gamma)^\circ$ and since the family $(\Delta(\Gamma')^\circ)_{\Gamma'\in\Theta}$ is filtered, there exists $\Gamma'\in\Theta$ such that $K\subseteq\Delta(\Gamma')^\circ$. Moreover, since $\Gamma_{\mathbb Z}$ is a discrete subgroup of $\mathbb R\times V$, it is actually finitely generated. Hence by possibly enlarging $\Gamma'$ we may assume that $\Gamma'_{\mathbb Z}=\Gamma_{\mathbb Z}$. Therefore, without loss of generality, we may assume that the semigroup $\Gamma$ is finitely generated.

Let $\{x_i\}_{i=1}^\ell$ be a system of generators of $\Gamma$, where $x_i=(n_i,\gamma_i)$. Then $\Delta(\Gamma)$ is just the convex hull of $n_i^{-1}\gamma_i$ ($i\in\{1,\ldots,\ell\}$).
The set 
\[F=\{\lambda_1x_1+\cdots+\lambda_\ell x_\ell\,|\,(\lambda_1,\ldots,\lambda_\ell)\in[0,1]^\ell\}\]
is a compact subset of $\mathbb R\times V$. Therefore the intersection of $F$ with $\Gamma_{\mathbb Z}$ is finite since $\Gamma_{\mathbb Z}$ is supposed to be discrete. In particular, there exists $x_0=(n_0,\gamma_0)\in\Gamma$ such that $x_0+y\in\Gamma$ for any $y\in F\cap\Gamma_{\mathbb Z}$. Let $n$ be an element of $\mathbb Z(\Gamma)$, $n\geqslant 1$, and let $\gamma\in\Gamma_{\mathbb Z,n}$. If $n^{-1}\gamma$ belongs to $\Delta(\Gamma)$, then there exists $(a_1,\ldots,a_\ell)\in\mathbb R_+^\ell$ such that $a_1n_1+\cdots+a_\ell n_\ell=n$ and that $\gamma=a_1\gamma_1+\cdots+a_\ell\gamma_\ell$. Let $b_i=\lfloor a_i\rfloor$ and $\lambda_i=a_i-b_i$
 for any $i\in\{1,\ldots,\ell\}$. We write $x=(n,\gamma)$ in the form $x=x'+y$ with
\[x'=b_1x_1+\cdots+b_\ell x_\ell\in\Gamma,\quad y=\lambda_1x_1+\cdots+\lambda_\ell x_\ell\in F.\]
Since $x\in\Gamma_{\mathbb Z}$, also is $y$. Hence $y\in F\cap\Gamma_{\mathbb Z}$. Thus $x+x_0=x'+(y+x_0)\in\Gamma$. In particular, one has \[\gamma+\gamma_0\in\Gamma_{n+n_0}.\] Now we introduce a norm $\|\ndot\|$ on $V$. Since $K$ is a compact subset of the relative interior of $\Delta(\Gamma)$, there exists $\epsilon>0$ such that, for any $u\in K$, the ball \[B(u,\epsilon)=\{u'\in W\,:\,\|u-u'\|\leqslant \epsilon\}\]
is contained in $\Delta(\Gamma)$. Moreover, the set $K$ is bounded. Therefore, for sufficiently positive integer $n\in\mathbb N(\Gamma)$, if $\beta$ is an element in $\Gamma_{\mathbb Z,n}\cap nK$, then one has \[(n-n_0)^{-1}(\beta-\gamma_0)\in\Delta(\Gamma),\] which implies that $\beta\in\Gamma_n$ by the above argument. The equality \eqref{Equ:intersectionK} is thus proved.
\end{proof}

\begin{defi}\label{Def: normalisation of Lebesgue measure}
Let $\Gamma$ be a graded  semigroup in $V$ such that $\Gamma_{\mathbb Z}$ is discrete. Let $A(\Gamma)_0$ be the vector subspace of $V$ which is the translation of the affine subspace $A(\Gamma)$. We equip $A(\Gamma)_0$ with the normalised Lebesgue measure such that the mass of a fundamental domain of the lattice $\Gamma_{\mathbb Z,0}$ in $A(\Gamma)_0$ is $1$. This measure induces by translation a Borel measure on $A(\Gamma)$. We denote by $\eta_{\Gamma}$ the restriction of this Borel measure to 
the closed convex set $\Delta(\Gamma)$, that is, for 
any function $f\in C_c(A(\Gamma))$ (namely $f$ is continuous on $A(\Gamma)$ and of compact support), one has

\[\int_{A(\Gamma)}f(x)\,\eta_{\Gamma}(\mathrm{d}x)=\int_{\Delta(\Gamma)}f(\gamma)\,\mathrm{d}\gamma,\]
where $\mathrm{d}\gamma$ denotes the normalised Lebesgue measure.
\end{defi}

The following theorem is the key point of the Newton-Okounkov body approach to the study of graded linear series \cite{Okounkov96,Kaveh_Khovanskii,Lazarsfeld_Mustata08}. Here we adopte the form presented in the Bourbaki seminar lecture of Boucksom \cite{Boucksom14}.

\begin{theo}\label{Thm:measureOkounkov}
Let $\Gamma$ be a graded semigroup in $V$ such that $\Gamma_{\mathbb Z}$ is discrete. For any integer $n\in\mathbb N(\Gamma)$, we denote by $\eta_{\Gamma,n}$ the Radon measure on $A(\Gamma)$ such that, for any function $f\in C_c(A(\Gamma))$ one has
\[\int_{A(\Gamma)}f(x)\,\eta_{\Gamma,n}(\mathrm{d}x)=\frac{1}{n^{\kappa}}\sum_{\gamma\in\Gamma_n}f(n^{-1}\gamma),\]
where $\kappa$ is the dimension of the affine space $A(\Gamma)$. Then the sequence of measures $\{\eta_{\Gamma,n}\}_{n\in\mathbb N(\Gamma)}$ 
converges vaguely (see \S\ref{Sec: vague convergence}) to the Radon measure $\eta_{\Gamma}$.
\end{theo}
\begin{proof}
Recall that the vague convergence in the statement of the theorem signifies that the sequence $\{\eta_{\Gamma,n}\}_{n\in\mathbb N(\Gamma)}$, viewed as a sequence of positive linear functionals on $C_c(A(\Gamma))$, converges pointwisely to $\eta_{\Gamma}$. In other words, for any continuous function $f$ on $A(\Gamma)$ of compact support, one has
\begin{equation}\label{Equ:vaguelimit}\lim_{n\in\mathbb N(\Gamma),\,n\rightarrow+\infty}\frac{1}{n^{\kappa}}\sum_{\gamma\in\Gamma_n}f(n^{-1}\gamma)=\int_{\Delta(\Gamma)}f(\gamma)\,\mathrm{d}\gamma.\end{equation} 
Note that the direct image preserves the vague convergence. Therefore, it suffices to prove that, for any non-negative continuous function $f$ on $\Delta(\Gamma)$ which is of compact support, the equality \eqref{Equ:vaguelimit} holds.

For any $n\in\mathbb N(\Gamma)$ one has
\[\frac{1}{n^\kappa}\sum_{\gamma\in\Gamma_n}f(n^{-1}\gamma)\leqslant\frac{1}{n^\kappa}\sum_{\gamma\in\Gamma_{\mathbb Z,n}\cap n\Delta(\Gamma)}f(n^{-1}\gamma).\]
Note that the right had side of the inequality is the $n^{\mathrm{th}}$ Riemann sum of the function $f$ on the convex set $\Delta(\Gamma)$. Therefore one has
\[\lim_{n\in\mathbb N(\Gamma),\,n\rightarrow+\infty}\frac{1}{n^\kappa}\sum_{\gamma\in\Gamma_{\mathbb Z,n}\cap n\Delta(\Gamma)}f(n^{-1}\gamma)=\int_{\Delta(\Gamma)}f(\gamma)\,\mathrm{d}\gamma,\] 
which implies
\[\limsup_{n\in\mathbb N(\Gamma),\,n\rightarrow+\infty}\frac{1}{n^{\kappa}}\sum_{\gamma\in\Gamma_n}f(n^{-1}\gamma)\leqslant\int_{\Delta(\Gamma)}f(\gamma)\,\mathrm{d}\gamma.\]
Moreover, if $g$ is a continuous function on $\Delta(\Gamma)$ whose support is contained in $\Delta(\Gamma)^\circ$ (the relative interior of $\Delta(\Gamma)$ in $A(\Gamma)$) and which is bounded from above by $f$, by Proposition \ref{Pro:proprietegradedsemigropu} \ref{Item:intersection K gamma n}, for sufficiently positive $n$ one has
\[\sum_{\gamma\in\Gamma_{\mathbb Z,n}\cap n\Delta(\Gamma)}g(n^{-1}\gamma)=\sum_{\gamma\in\Gamma_n}g(n^{-1}\gamma).\]
Hence one has
\[\liminf_{n\in\mathbb N(\Gamma),\,n\rightarrow+\infty}\frac{1}{n^{\kappa}}\sum_{\gamma\in\Gamma_n}f(n^{-1}\gamma)\geqslant\lim_{n\in\mathbb N(\Gamma),\,n\rightarrow+\infty}\frac{1}{n^{\kappa}}\sum_{\gamma\in\Gamma_n}g(n^{-1}\gamma)=\int_{\Delta(\Gamma)^\circ}g(\gamma)\,\mathrm{d}\gamma.\]
Since the restriction of the function $f$ on $\Delta(\Gamma)^\circ$ can be written as the limite of an increasing sequence of continous functions with support contained in $\Delta(\Gamma)^\circ$, by the monotone convergence theorem of Levi, one has
\[\liminf_{n\in\mathbb N(\Gamma),\,n\rightarrow+\infty}\frac{1}{n^{\kappa}}\sum_{\gamma\in\Gamma_n}f(n^{-1}\gamma)\geqslant\int_{\Delta(\Gamma)^\circ}f(\gamma)\,\mathrm{d}\gamma.\] 
Finally, since the border of $\Delta(\Gamma)$ has Lebesgue measure $0$, we obtain the desired result.
\end{proof}

\begin{defi}
Let $\Gamma$ be a graded semigroup in $V$. The dimension of the affine space $A(\Gamma)$ is called the \emph{Kodaira dimension}\index{Kodaira dimension} of $\Gamma$.
\end{defi}

\begin{coro}\label{Cor:convergencemesure}
We keep the notation and the hypotheses of Theorem \ref{Thm:measureOkounkov}. For any convex subset $C$ of $\Delta(\Gamma)$  one has
\begin{equation}\label{Equ:masslimite}\lim_{n\in\mathbb N(\Gamma),\,n\rightarrow+\infty}\frac{\#(\Gamma_n\cap nC)}{n^\kappa}=\eta_{\Gamma}(C),\end{equation}
where $\kappa$ is the Kodaira dimension of $\Gamma$.
\end{coro}
\begin{proof}
Let $C^\circ$ be the relative interior of $C$ in $A(\Gamma)$. If $C^\circ$ is empty, then one has $\eta_{\Gamma}(C)=0$. Moreover, for $n\in\mathbb N_{\geqslant 1}$, one has $\#(\Gamma_{n,\mathbb Z}\cap nC)=o(n^{\kappa})$ since $\Gamma_{n,\mathbb Z}$ is a translation of lattice (see Proposition \ref{Pro:proprietegradedsemigropu}). Therefore, one has
\[\lim_{n\in\mathbb N(\Gamma),\,n\rightarrow+\infty}\frac{\#(\Gamma_n\cap nC)}{n^{\kappa}}=0.\]

In the following, we assume that $C^\circ$ is not empty. Let $K$ be a compact convex subset of $C^{\circ}$. We can find a function $f\in C_c(A(\Gamma))$ with $0\leqslant f\leqslant \indic_C$,  $f|_K\equiv 1$. Then one has
\[\forall\,n\in\mathbb N(\Gamma),\;\frac{\#(\Gamma_n\cap nC)}{n^\kappa}\geqslant\int_{A(\Gamma)} f(x)\,\eta_{\Gamma,n}(\mathrm{d}x),\]
which leads to (by Theorem \ref{Thm:measureOkounkov})
\[\liminf_{n\in\mathbb N(\Gamma),\,n\rightarrow+\infty}\frac{\#(\Gamma_n\cap nC)}{n^\kappa}\geqslant\int_{A(\Gamma)} f(x)\,\eta_{\Gamma}(\mathrm{d}x)\geqslant\eta_{\Gamma}(K).\]
Since $K$ is arbitrary, we obtain
\[\liminf_{n\in\mathbb N(\Gamma),\,n\rightarrow+\infty}\frac{\#(\Gamma_n\cap nC)}{n^{\kappa}}\geqslant\eta_{\Gamma}(C^\circ)=\eta_{\Gamma}(C).\]
In particular, if $C$ is not bounded, then 
\[\lim_{n\in\mathbb N(\Gamma),\,n\rightarrow+\infty}\frac{\#(\Gamma_n\cap nC)}{n^{\kappa}}=\eta_{\Gamma}(C)=+\infty.\]

In the following, we assume in addition that the convex set $C$ is bounded. Denote by $\overline{C}$ the closure of the convex set $C$. It is a conex and compact subset of $A(\Gamma)$. Let $K$ be a compact subset of $A(\Gamma)$ such that the relative interior of $K$ contains $\overline C$.
For any non-negative function $g\in C_c(A(\Gamma))$ with support contained $K$ and such that $0\leqslant g\leqslant 1$, $g|_C\equiv 1$, one has
\[\forall\,n\in\mathbb N(C),\;\frac{\#(\Gamma_n\cap nC)}{n^\kappa}\leqslant\int_{A(\Gamma)} g(x)\,\eta_{\Gamma,n}(\mathrm{d}x).\]
By Theorem \ref{Thm:measureOkounkov}, we obtain
\[\limsup_{n\in\mathbb N(\Gamma),\,n\rightarrow+\infty}\frac{\#(\Gamma_n\cap nC)}{n^\kappa}\leqslant\int_{A(\Gamma)}g(x)\,\eta_{\Gamma}(\mathrm{d}x)\leqslant\eta_{\Gamma}(\mathrm{d}x)\leqslant\eta_{\Gamma}(K).\] 
Since $K$ is arbitrary, we obtain
\[\limsup_{n\rightarrow+\infty}\frac{\#(\Gamma_n\cap nC)}{n^\kappa}\leqslant\int_{A(\Gamma)}g(x)\,\eta_{\Gamma}(\mathrm{d}x)\leqslant\eta_{\Gamma}(C).\]
\end{proof}

\subsection{Concave transform} Let $V$ be a finite-dimensional vector space over $\mathbb R$ and $\Gamma$ be a graded semigroup in $V$ such that $\Gamma_{\mathbb Z}$ is discrete. We suppose given a map $\delta:\mathbb N_{\geqslant 1}\rightarrow\mathbb R$ such that $\delta(n)/n$ tends to $0$ when $n\rightarrow+\infty$.

\begin{defi}
Let $g:\Gamma\rightarrow\mathbb R$ be a function. We say that the function $g$ is 
\emph{strongly $\delta$-superadditive}\index{strongly superadditive} if for any $\ell\in\mathbb N_{\geqslant 2}$ and for all elements $(n_1,\gamma_1),\ldots,(n_\ell,\gamma_\ell)$ in $\Gamma$, one has
\begin{equation}\label{Equ:almostsuperadditive}g(n_1+\cdots+n_\ell,\gamma_1+\cdots+\gamma_\ell)\geqslant\sum_{i=1}^\ell(g(n_i,\gamma_i)-\delta(n_i)).\end{equation}
\end{defi}
The purpose of this subsection is to prove the following result.

\begin{theo}\label{Thm:concavetransform}
Let $\Gamma$ be a graded semigroup in $V$. We assume that $\Gamma_{\mathbb Z}$ is discrete and that $\Delta(\Gamma)$ is compact. Suppose given a function $g$ on $\Gamma$ which is strongly $\delta$-superadditive for certain function $\delta:\mathbb N_{\geqslant 1}\rightarrow\mathbb R$ such that \[\lim_{n\rightarrow+\infty}\frac{\delta(n)}{n}=0.\] For any $n\in\mathbb N(\Gamma)$, let $\nu_n$ be the Borel probability measure on $\mathbb R$ 
given by \[\forall\, f \in C_c(\mathbb R),\quad \int_{\mathbb R} f(t)\,\nu_n(\mathrm{d}t)=\frac{1}{\#\Gamma_n}\sum_{\gamma\in\Gamma_n}f\left({\textstyle \frac 1n}g(n,\gamma)\right).\]
The the sequence of measures $\{\nu_n\}_{n\in\mathbb N(\Gamma)}$ converges vaguely to a Borel measure $\nu_{\Gamma}$ on $\mathbb R$. Moreover, $\nu_{\Gamma}$ is either the zero measure or a probability measure, and in the latter case the sequence $\{\nu_n\}_{n\in\mathbb N(\Gamma)}$ actually converges weakly to $\nu_{\Gamma}$ (see Theorem \ref{Thm: criterion of weak convergence}) and there exists a concave function $G_{\Gamma}:\Delta(\Gamma)^\circ\rightarrow\mathbb R$ such that $\nu_{\Gamma}$ identifies with the direct image of \[\frac{1}{\eta_{\Gamma}(\Delta(\Gamma))}\eta_{\Gamma}\] by the map $G_\Gamma$.
\end{theo}
\begin{proof} 
We introduce an auxiliary function $\widetilde{g}$ on $\Gamma$ taking values in $\mathbb R\cup\{+\infty\}$ as follows:
\begin{equation}\label{Equ:gtilde}\forall\,u\in\Gamma,\quad\widetilde{g}(u)=\limsup_{n\rightarrow+\infty}\frac{g(nu)}{n}.\end{equation}
Note that the sequence defining $\widetilde{g}(u)$ is bounded from below and hence the sup limit does not take the value $-\infty$. 
The proof of the theorem is decomposed into the following steps.

{\it Step 1: The sup limit in the formula \eqref{Equ:gtilde} is actually a limit.} This follows from the following 
generalisation of Fekete's lemma (the case where $\delta(n) = 0$ for all $n$):
let $\{a_n\}_{n\geqslant 1}$ be a sequence in $\mathbb R$ such that, for any $\ell\in\mathbb N_{\geqslant 2}$ and for all $n_1,\ldots,n_\ell$ in $\mathbb N_{\geqslant 1}$ one has
\[a_{n_1+\cdots+n_\ell}\geqslant\sum_{i=1}^\ell(a_{n_i}-\delta(n_i)),\]
then the sequence $\{a_n/n\}_{n\geqslant 1}$ converges in $\mathbb R\cup\{+\infty\}$. In fact, if $p$ is an integer, $p\geqslant 1$ and if $m\in\mathbb N$, $r\in\{1,\ldots,p\}$ one has
\[a_{mp+r}\geqslant ma_{p}+a_r-m\delta(p)-\delta(r),\]
and hence
\[\frac{a_{mp+r}}{mp+r}\geqslant\frac{m}{mp+r}a_p+\frac{a_r}{mp+r}-\frac{m\delta(p)+\delta(r)}{mp+r}.\]
Therefore
\[\liminf_{n\rightarrow+\infty}\frac{a_n}{n}\geqslant\frac{a_p}{p}-\frac{\delta(p)}{p}.\]
In particular, $\liminf_{n\rightarrow+\infty}a_n/n\geqslant a_1-\delta(1)>-\infty$. Moreover, this inequality also implies that
\[\liminf_{n\rightarrow+\infty}\frac{a_n}{n}\geqslant\limsup_{p\rightarrow+\infty}\Big(\frac{a_p}{p}-\frac{\delta(p)}{p}\Big)=\limsup_{p\rightarrow+\infty}\frac{a_p}{p},\]
which leads to the convergence of the sequence $\{a_n/n\}_{n\geqslant 1}$.

{\it Step 2: Some properties of the function $\widetilde{g}$.} Let $u_1=(n_1,\gamma_1)$ and $u_2=(n_2,\gamma_2)$ be two elements in $\Gamma$. For any $n\in\mathbb N_{\geqslant 1}$ one has
\[g(n(u_1+u_2))\geqslant g(nu_1)+g(nu_2)-\delta(nn_1)-\delta(nn_2)\]
and hence
\[\frac{g(n(u_1+u_2))}{n}\geqslant\frac{g(nu_1)}{n}+\frac{g(nu_2)}{n}-\frac{\delta(nn_1)+\delta(nn_2)}{n}.\]
By taking the limit when $n\rightarrow+\infty$, we obtain $\widetilde{g}(u_1+u_2)\geqslant\widetilde{g}(u_1)+\widetilde{g}(u_2)$. In other words, the function $\widetilde{g}$ is superadditive.

Let $(n,\gamma)$ be an element of $\Gamma$. Note that for any $N\in\mathbb N_{\geqslant 1}$ one has
\[\frac{{g}(Nn,N\gamma)}{N}\geqslant g(n,\gamma)-\delta(n).\]
By taking the limit when $N\rightarrow+\infty$, we obtain
\begin{equation}\label{Equ:minorationgtilde}\widetilde{g}(n,\gamma)\geqslant g(n,\gamma)-\delta(n).\end{equation}

{\it Step 3: Construction of the function $G_{\Gamma}$.} For any $t\in\mathbb R$, let $\Gamma^t$ be the set of all $(n,\gamma)\in\Gamma$ such that $\widetilde g(n,\gamma)\geqslant nt$. It is actually a sub-semigroup of $\Gamma$ since $\widetilde{g}$ is super-additive. Note that $\{\Gamma^t\}_{t\in\mathbb R}$ is a decreasing family of sub-semigroups of $\Gamma$ and hence $\{\Delta(\Gamma^t)\}_{t\in\mathbb R}$ is a decreasing family of closed convex subsets of $\Delta(\Gamma)$. We define the function $G_{\Gamma}:\Delta(\Gamma)\rightarrow\mathbb R\cup\{+\infty\}$ as follows:
\[\forall\,x\in\Delta(\Gamma),\quad G_{\Gamma}(x)=\sup\{t\in\mathbb R\,:\,x\in\Delta(\Gamma^t)\}.\]
By definition, if $t$ is a real number, then $G_\Gamma(x)\geqslant t$ if and only if $x\in\bigcap_{s<t}\Delta(\Gamma^s)$.
We claim that the function $G_{\Gamma}$ is concave. In fact, since the function $\widetilde{g}$ is super-additive, we obtain that, if $s$ and $t$ are two real numbers and if $\epsilon\in [0,1]\cap\mathbb Q$, for $u\in\Gamma^s$ and $v\in\Gamma^t$ one has
\[N(\epsilon u+(1-\epsilon)v)\in\Gamma^{\epsilon s+(1-\epsilon)t}, \]
where $N$ is an element in $\mathbb N_{\geqslant 1}$ such that $N\epsilon\in\mathbb N$. Therefore one has \[\epsilon\Delta(\Gamma^s)+(1-\epsilon)\Delta(\Gamma^t)\subseteq\Delta(\Gamma^{\epsilon s+(1- \epsilon)t}).\]
In general, if we choose a sequence $\{ \epsilon_n \}$ of rational numbers such that $\lim_{n\to\infty} \epsilon_n = \epsilon$ and $\epsilon_n s + (1 -\epsilon_n) t \geqslant \epsilon s + (1 -\epsilon) t$ for all $n$,
then
\[
\epsilon_n \Delta(\Gamma^s)+(1-\epsilon_n )\Delta(\Gamma^t)\subseteq\Delta(\Gamma^{\epsilon_n s+(1-\epsilon_n) t})
\subseteq\Delta(\Gamma^{\epsilon s+(1-\epsilon) t}),
\]
and hence $\epsilon \Delta(\Gamma^s)+(1-\epsilon )\Delta(\Gamma^t)\subseteq\Delta(\Gamma^{\epsilon s+(1-\epsilon) t})$.
Combining with the definition of the function $G_{\Gamma}$, we obtain the concavity of $G_{\Gamma}$. In particular, the restriction of the function $G_{\Gamma}$ on $\Delta(\Gamma)^\circ$ is either finite or identically $+\infty$, and it is a continuous function on $\Delta(\Gamma)^\circ$ when it is finite.

{\it Step 4. Abundance of $\Gamma^t_{\mathbb Z}$.} Let $t$ be an element of $\mathbb R$ such that $t<\sup_{x\in\Delta(\Gamma)}G_{\Gamma}(x)$. We will prove that $\Gamma_{\mathbb Z}^t=\Gamma_{\mathbb Z}$ (and hence $A(\Gamma^t)=A(\Gamma)$). 
Note that $\Gamma_{\mathbb Z}$ is finitely generated because $\Gamma_{\mathbb Z}$ is discrete.
Let $u_i=(n_i,\gamma_i)$, $i\in\{1,\ldots,\ell\}$ be a family of elements in $\Gamma$ which forms a system of generators in $\Gamma_{\mathbb Z}$. Since $t<\sup_{x\in\Delta(\Gamma)}G_{\Gamma}(x)$, there exists $\epsilon>0$ such that $\Gamma^{t+\epsilon}$ is not empty. Let $u_0=(n_0,\gamma_0)$ be an element in $\Gamma^{t+\epsilon}$. By definition, one has $\widetilde g(u_0)\geqslant n_0(t +\epsilon)$. Therefore, for sufficiently positive integer $p$, one has
\[\forall\,i\in\{1,\ldots,\ell\},\quad\widetilde{g}(pu_0+u_i)\geqslant p\widetilde{g}(u_0)+\widetilde{g}(u_i)\geqslant (pn_0+n_i)t,\]
namely $pu_0+u_i\in\Gamma^t$ for any $i\in\{1,\ldots,\ell\}$, which leads to $\Gamma^t_{\mathbb Z}=\Gamma_{\mathbb Z}$.

{\it Step 5: Lower bound of the function $g$.} We fix a (closed) fundamental domain $F$ of the lattice $\Gamma_{\mathbb Z,0}$ (see Proposition \ref{Pro:proprietegradedsemigropu}
\ref{Item: Gamma Z 0 a lattice}). For $n\in\mathbb N(\Gamma)$, we call an $n$-\emph{cell}\index{n-cell@$n$-cell} in $A(\Gamma)_n$ any closed convex subset of $A(\Gamma)_n$ of the form $\gamma_0+F$, where $\gamma_0$ is an element in $\Gamma_{\mathbb Z,n}$. 
We say that a compact subset $K$ of $A(\Gamma)$ is $n$-\emph{tileable}\index{tileable} if it can be written as a union of $n$-cells in $A(\Gamma)_n$. Note that, if $K$ is $n$-tileable, then, for any integer $p\geqslant 1$, the set $pK$ is $pn$-tileable since $pF$ can be written as the union of $p^{\kappa}$ $0$-celles. 

Let $t$ be a real number such that $t<\sup_{x\in\Delta(\Gamma)}G_{\Gamma}(x)$, and $\epsilon$ be a positive real number. Let $m\geqslant 1$ be the generator of the group $\mathbb Z(\Gamma)$. 
Suppose given a compact subset $K$ of $\Delta(\Gamma^t)^\circ$. We assume that there exists an integer $n\in\mathbb N(\Gamma)$ such that $nK$ is $n$-tileable. 

By Proposition \ref{Pro:proprietegradedsemigropu} \ref{Item:intersection K gamma n}, there exists an integer $n_0\in\mathbb N_{\geqslant 1}$ which verifies the following conditions (in the condition \ref{Item: chap6 condition 2}  we also use the result of Step 4 to identify $\Gamma_{mn_0,\mathbb Z}$ with $\Gamma_{mn_0,\mathbb Z}^t$): 
\begin{enumerate}[label=\rm(\arabic*)]
\item\label{Item: chap6 condition 1} $mn_0K$ is $mn_0$-tileable;
\item\label{Item: chap6 condition 2} $mnK\cap\Gamma_{mn,\mathbb Z}\subseteq\Gamma_{mn}^t$ for any $n\in\mathbb N$, $n\geqslant n_0$;
\item\label{Item: chap6 condition 3} for any integer $q\geqslant mn_0$, $\delta(q)/q<\epsilon/3$.
\end{enumerate}
For simplifying the notation, in the following we denote by $\Theta$ the set $mn_0K\cap\Gamma_{mn_0,\mathbb Z}$. Note that the condition \ref{Item: chap6 condition 2} implies that $\widetilde{g}(mn_0,\gamma)\geqslant t$ for any $\gamma\in\Theta$. Therefore by the definition of the function $\widetilde{g}$ and the finiteness of the set $
\Theta$, we obtain that there exists an integer $N_0$ divisible by $n_0$ such that 
\begin{equation}\label{Equ:estimationggamma}\frac{1}{mN_0}g(mN_0,(N_0/n_0)\gamma)\geqslant t-\frac{\epsilon}{3}\end{equation}
for any $\gamma\in\Theta$. 

Let $N$ be an integer, $N\geqslant n_0$. Let $\alpha$ be an element in $mNK\cap\Gamma_{mN}$ and $x=(n_0/N)\alpha$. Since $mn_0K$ is $mn_0$-tileable, there exists an $mn_0$-cell $C$ such that 
$x$
belongs to $C$. We write $C$ as $\gamma_0+F$ with $\gamma_0\in\Gamma_{mn_0}$. Let $\{e_1,\ldots,e_{\kappa}\}$ be the basis of $\Gamma_{\mathbb Z}$ defining the fundamental domain $F$. Then the point $x$ can be written in a unique way as
\[x=\gamma_0+\sum_{i=1}^\kappa\lambda_i e_i,\]
where \[\forall\,i\in\{1,\ldots,\kappa\},\quad\lambda_i\in [0,1].\] Moreover, since $N(x-\gamma_0)=n_0\alpha-N\gamma_0\in \Gamma_{0,\mathbb Z}$, we obtain that $N\lambda_i\in\mathbb Z$ for any $i\in\{1,\ldots,\kappa\}$. Without loss of generality, we may assume that $\lambda_1\geqslant\ldots\geqslant\lambda_{\kappa}$. Then we can rewrite $x$ as
\[x=\sum_{i=0}^\kappa(\lambda_{i}-\lambda_{i+1})\gamma_{i},\]
where by convention $\lambda_0=1$, $\lambda_{\kappa+1}=0$, and for $i\in\{1,\ldots,\kappa\}$, $\gamma_i=\gamma_0+e_1+\cdots+e_i$. Note that $\gamma_0,\ldots\gamma_{\kappa}$ are vertices of the $mn_0$-cell $C$, hence belong to $\Theta$. For any $i\in\{0,\ldots,\kappa\}$, let $b_i$ be the integral part of \[\frac{N}{n_0}(\lambda_i-\lambda_{i+1}).\]
One has
\[N-(\kappa+1)n_0+1\leqslant n_0\sum_{i=0}^{\kappa}b_i\leqslant N.\]
Therefore, we can write $\alpha$ as
\[\alpha=\sum_{i=0}^\kappa b_i\gamma_i+\beta', \]
where $\beta\in\Gamma_{mr',\mathbb Z}\cap mr'K$, with
\[r'=N-n_0\sum_{i=0}^\kappa b_i\in\{0,\ldots,(\kappa+1)n_0-1\}.\]
Note that we have assumed that $N\geqslant n_0$. Therefore, if $r'\leqslant n_0-1$, then there exists at least an indice $b_i$ which is $>0$. In this case, we replace $\beta'$ by $\beta'+\gamma_i$ and $b_i$ by $b_i-1$. Thus we obtain the existence of a decomposition of $\alpha$ into the form
\[\alpha=\sum_{i=0}^\kappa a_i\gamma_i+\beta\]
with $a_i\in\mathbb N$ for $i\in\{0,\ldots,\kappa\}$, and $\beta\in\Gamma_{mr,\mathbb Z}\cap mrK$ with
\[r\in\{n_0,\ldots,(\kappa+1)n_0-1\}.\] The advantage of the new decomposition is that $\beta$ actually belongs to $\Gamma_{mr}$ (see the condition \ref{Item: chap6 condition 2} above).
Finally, we write each $a_i$ in the form $a_i=p_iN_0/n_0+r_i$ with $p_i\in\mathbb N$ and $r_i\in\{0,\ldots,N_0/n_0-1\}$. Then we can decompose $\alpha$ as
\[\alpha=\sum_{i=0}^\kappa p_i(N_0/n_0)\gamma_i+\omega,\]
where 
\[\omega=\beta+\sum_{i=0}^\kappa r_i\gamma_i.\]
The element $\omega$ belongs to certain $\Gamma_{ms}\cap msK$ with
\[s\in\{n_0,\ldots,(\kappa+1)N_0-1\}.\]
Hence by \eqref{Equ:almostsuperadditive} and \eqref{Equ:estimationggamma} one obtains
\begin{equation}\label{Equ: gmNalpha over mN}\begin{split}&\quad\;\frac{g(mN,\alpha)}{mN}\\&\geqslant\frac{1}{mN}\bigg(\sum_{i=0}^\kappa p_ig(mN_0,(N_0/n_0)\gamma_i)+g(ms,\omega)-\delta(mN_0)\sum_{i=0}^\kappa p_i-\delta(ms)\bigg)\\
&\geqslant\frac{N_0P}{N}(t-\epsilon/3)+\frac{g(ms,\omega)}{mN}-\frac{P}{mN}\delta(mN_0)-\frac{\delta(ms)}{N},
\end{split}\end{equation}
where 
\[P=p_0+\cdots+p_\kappa=\frac{N-s}{N_0}.\]
Therefore we obtain
\[\liminf_{N\rightarrow+\infty}\inf_{\alpha\in mNK\cap\Gamma_{mN}}\frac{g(mN,\alpha)}{mN}\geqslant t-\frac{2\epsilon}{3},\]
where we have used the condition (3) above to obtain
\[\frac{P\delta(mN_0)}{mN}=\frac{N-s}{N}\cdot\frac{\delta(mN_0)}{mN_0}\leqslant\frac{N-s}{N}\cdot\frac{\epsilon}{3}.\]
Therefore, there exists an integer $N'$ depending on $t$, $\epsilon$ and $K$ such that $g(mN,\alpha)\geqslant mN(t-\epsilon)$ for any $N\geqslant N'$ and any $\alpha\in\Gamma_{mN}\cap mNK$.

{\it Step 6: Convergence of measures.} We now proceed with the proof of the convergence of the measures. We first consider the case where $G_{\Gamma}$ is identically $+\infty$ on the interior of $\Delta(\Gamma)$. Let $f$ be a non-negative continuous function with compact support on $\mathbb R$ and $t_0\in\mathbb R$ be a real number which is larger than the supremum of the support of the function $f$.  Let $K$ be a compact subset of $\Delta(\Gamma)^\circ$. By the results in Step 5, we obtain that, there exists $n_0\in\mathbb N$ such that, for any $n\in\mathbb N(\Gamma)$, $n\geqslant n_0$ and any $\alpha\in\Gamma_n\cap nK$, one has $g(N,\alpha)\geqslant nt_0$. Hence 
\[\int_{\mathbb R}f(t)\,\nu_n(\mathrm{d}t)\leqslant\Big(\frac{\#(\Gamma_n\setminus nK)}{\#\Gamma_n}\Big)M=\Big(1-\frac{\#(\Gamma_n\cap nK)}{\#\Gamma_n}\Big)M,\]
where $M=\sup_{t\in\mathbb R}f(t)$.
By Corollary \ref{Cor:convergencemesure}, one has
\[\lim_{n\in\mathbb N(\Gamma),\,n\rightarrow+\infty}\frac{\#(\Gamma_n\cap nK)}{\#\Gamma_n}=\frac{\eta_{\Gamma}(K)}{\eta_{\Gamma}(\Delta(\Gamma))}.\]
Since $K$ is arbitrary, we obtain 
\[\lim_{n\in\mathbb N(\Gamma),\,n\rightarrow+\infty}\int_{\mathbb R}f(t)\,\nu_n(\mathrm{d}t)=0.\]

In the following, we assume that the function $G_{\Gamma}$ is finite. In this case, the direct image $\nu_{\Gamma}$ of $\eta_{\Gamma}(\Delta(\Gamma))^{-1}\eta_{\Gamma}$ by $G_\Gamma$ is a Borel provability measure on $\mathbb R$. We denote by $F$ its probability distribution function, namely
\[\forall\,t\in\mathbb R,\;F(t)=\nu_{\Gamma}(\intervalle{]}{-\infty}{t}{
]})=1-\frac{\eta_{\Gamma}(\Delta(\Gamma^t))}{\eta_{\Gamma}(\Delta(\Gamma))}.\]
By Corollary \ref{Cor:convergencemesure}, one has
\begin{equation}\label{Equ:F(t)}F(t)=1-\lim_{n\in\mathbb N(\Gamma),\,n\rightarrow+\infty}\frac{\#(\Gamma_n\cap n\Delta(\Gamma^t))}{\#\Gamma_n}.\end{equation}
The function $F$ is continuous on $\mathbb R$, except possibly at the point $\sup_{x\in\Delta(\Gamma)}G_{\Gamma}(x)$ (the discontinuity of the function $F$ happens precisely when the function $G_\Gamma$ is constant on $\Delta(\Gamma)^\circ$). For any $n\in\mathbb N(\Gamma)$, let $F_n$ be the probability distribution function of $\nu_n$.

If $(n,\gamma)$ is an element of $\Gamma$, then one has \[G(n^{-1}\gamma)\geqslant\frac{\widetilde{g}(n,\gamma)}{n}\geqslant \frac{g(n,\gamma)}{n}-\frac{\delta(n)}{n},\]
where the second inequality comes from \eqref{Equ:minorationgtilde}. Therefore we obtain
\[\forall\,t\in\mathbb R,\;\{(n,\gamma)\in\Gamma\,:\,G(n^{-1}\gamma)> t-\delta(n)/n\}\supseteq\{(n,\gamma)\in\Gamma\,:\,g(n,\gamma)/n> t\},\]
which implies (by \eqref{Equ:F(t)})
\begin{equation}\label{Equ:convergenceFn(t)1}\forall\,\epsilon>0,\quad 1-\liminf_{n\in\mathbb N(\Gamma),\,n\rightarrow+\infty}F_n(t)\leqslant 1-F(t-\epsilon).\end{equation}

Conversely, for any $t\in\mathbb R$, any $\epsilon>0$ and any compact subset $K$ of $\Delta(\Gamma^{t+\epsilon})^\circ$ (the relative interior of $\Delta(\Gamma^{t+\epsilon})$ with respect to $A(\Gamma)$). By the result obtained in Step 5, we obtain that, there exists $N_0\in\mathbb N$ such that, for any $n\in\mathbb N(\Gamma)$, $n\geqslant N_0$, one has 
\[\forall\,\gamma\in\Gamma_n\cap nK,\quad g(n,\gamma)\geqslant nt.\]
Therefore, we obtain
\begin{equation}\label{Equ:convergenceFn(t)2}\forall\,\epsilon>0,\quad 1-\limsup_{n\in\mathbb N(\Gamma),\,n\rightarrow+\infty}F_n(t)\geqslant 1-F(t+\epsilon).\end{equation}  
The estimates \eqref{Equ:convergenceFn(t)1} and \eqref{Equ:convergenceFn(t)2} leads to the convergence of $\{F_n(t)\}_{n\in\mathbb N}$ to $F(t)$ if $t\in\mathbb R$ is a point of continuity of the function $F$, which implies the weak convergence of the sequence $\{\nu_{n}\}_{n\in\mathbb N(
\Gamma)}$ to $\nu_\Gamma$ (see \cite[\S I.4]{Petrov75} for more details about weak convergence of Borel probability measures on $\mathbb R$).
\end{proof}

\begin{defi}
Let $\Gamma$ be a graded semigroup in $V$ such that $\Gamma_{\mathbb Z}$ is discrete, and $g:\Gamma\rightarrow\mathbb R$ and $\delta:\mathbb N_{\geqslant 1}\rightarrow\mathbb R$ be functions. We say that the function $g$ is $\delta$-\emph{superadditive}\index{superadditive} if for all elements $(n_1,\gamma_1)$ and $(n_2,\gamma_2)$ in $\Gamma$, one has
\begin{equation}\label{Equ:almostsuperadditive2}g(n_1+n_2,\gamma_1+\gamma_2)\geqslant g(n_1,\gamma_1)+g(n_2,\gamma_2)-\delta(n_1)-\delta(n_2).\end{equation}
\end{defi}

\begin{lemm}
Let $\delta:\mathbb N_{\geqslant 1}\rightarrow\mathbb R_{\geqslant 0}$ be an increasing function such that
\[\sum_{a\in\mathbb N}\frac{\delta(2^a)}{2^a}<+\infty.\]
Then one has
\begin{equation}\label{Equ: limit of delta n over n}
\lim_{n\rightarrow+\infty}\frac{\delta(n)}{n}=0
\end{equation}
and
\begin{equation}\label{Equ: sum lim 2alpha}
\lim_{a\rightarrow+\infty}\frac{1}{2^a}\sum_{i=0}^{a} \delta(2^i)=0.
\end{equation}
\end{lemm}
\begin{proof}
For $n\in\mathbb N_{\geqslant 1}$, let $a(n)=\lfloor\log_2n\rfloor$. One has $2^{a(n)}\leqslant n< 2^{a(n)+1}$. Hence
\[\frac{\delta(n)}{n}\leqslant \frac{\delta(2^{a(n)+1})}{2^{a(n)}}.\]
By the hypothesis of the lemma, one has
\[\lim_{n\rightarrow+\infty}\frac{\delta(2^{a(n)+1})}{2^{a(n)+1}}=0,\]
which implies \eqref{Equ: limit of delta n over n}.

For any $a\in\mathbb N$, let 
\[S_a:=\sum_{i\in\mathbb N,\,{i\geqslant a}}\frac{\delta(2^i)}{2^i}.\]
By Abel's summation formula, one has
\[\sum_{i=0}^a\delta(2^i)=\sum_{i=0}^a(S_i-S_{i+1})2^i=S_0-S_{a+1}2^a+\sum_{i=1}^a S_i2^{i-1}.\]
As the sequence $\{S_a\}_{a\in\mathbb N}$ converges to $0$, one has
\[\lim_{a\rightarrow+\infty}\frac 1{2^a}\sum_{i=1}^a S_i2^{i-1}=0,\]
which implies the relation \eqref{Equ: sum lim 2alpha}.
\end{proof}

\begin{prop}\label{Pro: weak superadditive conv}
Let $\delta:\mathbb N_{\geqslant 1}\rightarrow\mathbb R_{\geqslant 0}$ be an increasing function such that 
\begin{equation}\label{Equ: convergence of sum delta}\sum_{a\in\mathbb N_{\geqslant1}}\frac{\delta(2^a)}{2^a}<+\infty.\end{equation}
Let $\{b_n\}_{n\in\mathbb N}$ be a sequence of real numbers. We assume that there exists an integer $n_0>0$ such that, for any couple $(n,m)$ of  integers which are $\geqslant n_0$, one has 
\begin{equation}\label{Equ: delta sup additivite faible}b_{n+m}\geqslant b_n+b_m-\delta(n)-\delta(m).\end{equation}
Then the sequence $\{b_n/n\}_{n\in\mathbb N_{\geqslant 1}}$ converges in $\mathbb R\cup\{+\infty\}$.
\end{prop}
\begin{proof}
We first treat the case where $n_0=1$.  For any $n\in\mathbb N_{\geqslant 1}$ one has
\[b_{2n}\geqslant 2b_{n}-2\delta(n),\]
and hence by induction we obtain that 
\begin{equation}\label{Equ: minoration b 2 aplpha n}b_{2^a n}\geqslant 2^a \bigg(b_n-\sum_{i=0}^{a-1} \frac{\delta(2^{i}n)}{2^i}\bigg).\end{equation}
In particular, one has
\[\frac{b_{2^a}}{2^a}\geqslant b_1-\sum_{i=0}^{a-1}\frac{\delta(2^i)}{2^i},\]
which implies that 
\[\limsup_{n\rightarrow+\infty}\frac{b_n}{n}>-\infty.\]

For any $a\in\mathbb N$, let
\[S_a=\sum_{i\in\mathbb N,\, i\geqslant a}\frac{\delta(2^i)}{2^i}.\]
By the hypothesis \eqref{Equ: convergence of sum delta}, we have
\begin{equation}\label{Equ: limite Salpha 0}\lim_{a\rightarrow+\infty}S_a=0.\end{equation}
Let $n\in\mathbb N_{\geqslant 1}$ and let $a(n)$ be the unique natural number such that $2^{a(n)}\leqslant n<2^{a(n)+1}$, namely $a(n)=\lfloor \log_2n\rfloor$. 
Let $p$ be an element in $\mathbb N_{\geqslant 1}$, which is written in $2$-adic basis as
\[p=\sum_{i=0}^\kappa\epsilon_i2^i\]
with $\epsilon_i\in\{0,1\}$ for $i\in\{0,\ldots,\kappa\}$ and $\epsilon_\kappa=1$. For any $r\in\{0,\ldots,n-1\}$, by \eqref{Equ: delta sup additivite faible} one has
\begin{equation}\label{Equ:minoartion bnp+r}b_{np+r}\geqslant b_{np}+b_r-\delta(np)-\delta(r)\geqslant b_{np}+b_r-2\delta(np).\end{equation}
Moreover, by induction on $\kappa$ one has
\begin{equation*}
\begin{split}
b_{np}&\geqslant\sum_{i=0}^{\kappa}\epsilon_ib_{2^in}-2\sum_{i=1}^\kappa\epsilon_i\delta(2^in)\\&\geqslant\sum_{i=0}^\kappa\epsilon_i2^ib_{n}-\sum_{i=1}^\kappa\epsilon_i\bigg(2^i\sum_{j=0}^{i-1}
\frac{\delta(2^jn)}{2^j}+2\delta(2^in)\bigg)\\
&\geqslant pb_n-2p\sum_{j=0}^{\kappa}\frac{\delta(2^jn)}{2^j}.
\end{split}
\end{equation*}
Since $n\geqslant 2^{a(n)}$ we deduce that
\begin{equation}\label{Equ: minormation de bnp}
b_{np}\geqslant pb_n-p2^{\beta+1}S_{a(n)}\geqslant pb_n-pnS_{a(n)}
\end{equation}
Combining \eqref{Equ:minoartion bnp+r} and \eqref{Equ: minormation de bnp}, we obtain
\[\frac{b_{np+r}}{np+r}\geqslant\frac{pb_n+b_r}{np+r}-\frac{np}{np+r}S_{a(n)}-2\frac{\delta(np)}{np+r}.\]
Taking the infimum limit when $p\rightarrow+\infty$, by \eqref{Equ: limit of delta n over n} we obtain 
\[\liminf_{m\rightarrow+\infty} \frac{b_m}{m}\geqslant \frac{b_n}{n}-S_{a(n)},\]
which implies, by \eqref{Equ: limite Salpha 0}, that
\[\liminf_{m\rightarrow+\infty} \frac{b_m}{m}\geqslant\limsup_{n\rightarrow+\infty}\frac{b_n}{n}.\]
Therefore the sequence $\{b_n/n\}_{n\in\mathbb N_{\geqslant 1}}$ converges in $\mathbb R\cup\{+\infty\}$.

For the general case, we apply the obtained result on the sequence $\{b_{n_0k}\}_{k\in\mathbb N_{\geqslant 1}}$ and obtain the convergence of the sequence $\{b_{n_0k}/k\}_{k\in\mathbb N_{\geqslant 1}}$. Moreover, if $\ell$ is an element in $\{n_0,\ldots,2n_0-1\}$, then for any $k\in\mathbb N_{\geqslant 1}$ one has
\[b_{n_0(k+2)}-b_{2n_0-\ell}+\delta(n_0k+\ell)+\delta(2n_0-\ell)\geqslant b_{n_0k+\ell}\geqslant
b_{n_0k}+b_{\ell}-\delta(n_0k)-\delta(\ell).\]
Dividing this formula by $n_0k+\ell$ and taking the limit when $k\rightarrow+\infty$, we obtain
\[\lim_{k\rightarrow+\infty}\frac{b_{n_0k+\ell}}{n_0k+\ell}=\lim_{k\rightarrow+\infty}\frac{b_{n_0k}}{n_0k}.\]
Since $\ell$ is arbitrary, we obtain the statement announced in the proposition.
\end{proof}

\begin{theo}\label{Thm:concavetransform 2}
Let $\Gamma$ be a graded semigroup in $V$. We assume that $\Gamma_{\mathbb Z}$ is discrete and that $\Delta(\Gamma)$ is compact. Suppose given a function $g$ on $\Gamma$ which is $\delta$-superadditive for certain increasing function $\delta:\mathbb N_{\geqslant 1}\rightarrow\mathbb R$ such that 
\begin{equation}\label{Equ: somme delta 2alpha sur 2alpha fini}\sum_{a\in\mathbb N}\frac{\delta(2^a)}{2^{a}}<+\infty.\end{equation} For any $n\in\mathbb N(\Gamma)$, let $\nu_n$ be the Borel probability measure on $\mathbb R$ such that 
\[\int_{\mathbb R} f(t)\,\nu_n(\mathrm{d}t)=\frac{1}{\#\Gamma_n}\sum_{\gamma\in\Gamma_n}f(n^{-1}g(n,\gamma)).\]
The the sequence of measures $\{\nu_n\}_{n\in\mathbb N(\Gamma)}$ converges vaguely to a Borel measure $\nu_{\Gamma}$ on $\mathbb R$. Moreover, $\nu_{\Gamma}$ is either the zero measure or a probability measure, and in the latter case the sequence $\{\nu_n\}_{n\in\mathbb N(\Gamma)}$ actually converges weakly to $\nu_{\Gamma}$ and there exists a concave function $G:\Delta(\Gamma)^\circ\rightarrow\mathbb R$, called \emph{concave transform}\index{concave transform} of $g$, such that $\nu_{\Gamma}$ identifies with the direct image of \[\frac{1}{\eta_{\Gamma}(\Delta(\Gamma))}\eta_{\Gamma}\] by the map $G_\Gamma$.
\end{theo}
\begin{proof}
The proof is very similar to that of Theorem \ref{Thm:concavetransform}. We will sketch it in emphasising the difference. Let $u=(\ell,\gamma)$ be an element in $\Gamma$, where $\ell\geqslant 1$. Since the function $g$ is $\delta$-superadditive, for any pair $(n,m)\in\mathbb N_{\geqslant 1}$ one has
\[g((n+m)u)\geqslant g(nu)+g(mu)-\delta(n\ell)-\delta(m\ell).\]  Moreover, if we let $b$ be an integer such that $\ell\leqslant 2^b$, then, by the increasing property of the function $\delta$, one has
\[\sum_{a\in\mathbb N}\frac{\delta(2^a\ell)}{2^a}\leqslant 2^b\sum_{a\in\mathbb N}\frac{\delta(2^{a+b})}{2^{a+b}}<+\infty.\] 
By Proposition \ref{Pro: weak superadditive conv}, for any $u\in\Gamma$, the sequence $\{g(nu)/n\}_{n\in\mathbb N_{\geqslant 1}}$ converges in $\mathbb R\cup\{+\infty\}$. We denote by $\widetilde g(u)$ the limite of the sequence. Moreover, the convergence of the series $\sum_{a\in\mathbb N}{\delta(2^a)}/{2^{a}}$ implies that
\[\lim_{a\rightarrow+\infty}\frac{\delta(2^a)}{2^a}=0.\]
Still by the hypothesis that the function $\delta(\ndot)$ is increasing, we deduce that 
\[\lim_{n\rightarrow+\infty}\frac{\delta(n)}{n}=0.\]
Therefore, by the same argument as in the Step 2 of the proof of Theorem \ref{Thm:concavetransform}, we obtain that the function $\widetilde g$ is superadditive, namely, for any pair $u_1,u_2$ of elements in $\Gamma$, one has
\[\widetilde g(u_1+u_2)\geqslant\widetilde g(u_1)+\widetilde g(u_2).\]
Moreover, for any $(n,\gamma)\in \Gamma$ and any $a\in\mathbb N_{\geqslant 1}$ one has
\begin{equation}\label{Equ: minormation of g 2alpha n}g(2^a n,2^a\gamma)\geqslant 2^a g(n,\gamma)-\sum_{i=0}^{a-1}2^{a-i}\delta(2^in).\end{equation}
Let $b(n)=\lfloor \log_2n\rfloor+1$. One has $2^{b(n)-1}\leqslant n< 2^{b(n)}$. Let 
\[R(n)=2^{b(n)}\sum_{i=b(n)}^{+\infty}\frac{\delta(2^i)}{2^i}.\]
Note that one has
\[\lim_{n\rightarrow+\infty}\frac{R(n)}{n}=0\]
by the hypothesis \eqref{Equ: somme delta 2alpha sur 2alpha fini}.  
By the increasing property of the function $\delta$ one has
\[\sum_{i=0}^{a-1}2^{a-i}\delta(2^in)\leqslant 2^{a+b(n)}\sum_{i=0}^{a-1}\frac{\delta(2^{i+b(n)})}{2^{i+b(n)}}\leqslant 2^a R(n).\]
Therefore the inequality \eqref{Equ: minormation of g 2alpha n} leads to 
\[\widetilde g(n,\gamma)\geqslant g(n,\gamma)-R(n).\]
We then proceed as in the Steps 3-6 of the proof of Theorem \ref{Thm:concavetransform}, except that in the counterpart of the minoration \eqref{Equ: gmNalpha over mN} we need more elaborated estimate as in \eqref{Equ: minormation de bnp}. 
\end{proof}

{
\begin{rema}\label{remark:sup:G:Boucksom_Chen}
We keep the notations in the proof of Theorems \ref{Thm:concavetransform} and \ref{Thm:concavetransform 2}. 
By virtue of \cite[Lemma 1.6]{Boucksom_Chen} (see also \cite{Chen_JTNB}), we obtain that, for any real number $t$ such that 
\[t<\lim_{n\in\mathbb N(\Gamma),\,n\rightarrow+\infty}\max_{\gamma\in\Gamma_n}\frac{\widetilde g(n,\gamma)}{n}=\lim_{n\in\mathbb N(\Gamma),\,n\rightarrow+\infty}\frac{g(n,\gamma)}{n},\]
the set $\{x\in\Delta(\Gamma)^\circ\,:\,G(x)\geqslant t\}$ has a positive measure with respect to $\eta_\Gamma$ (and hence is not empty). In particular, we obtain 
\begin{equation}\label{Equ: sup of concave transform}\sup_{x\in\Delta(\Gamma)^\circ}G(x)=\lim_{n\in\mathbb N(\Gamma),\,n\rightarrow+\infty}\frac{g(n,\gamma)}{n}\end{equation}
\end{rema}
}

\subsection{Applications to the study of graded algebras}\label{Subsec: applications to the study of graded algebras}

Let $d\geqslant 1$ be an integer. We call \emph{monomial order}\index{monomial order} on $\mathbb Z^d$ any total order $\leqslant$ on $\mathbb Z^d$ such that $0\leqslant\alpha$ for any $\alpha\in\mathbb N^d$
and that $\alpha\leqslant\alpha'$ implies $\alpha+\beta\leqslant\alpha'+\beta$ for all $\alpha$, $\alpha'$ and $\beta$ in $\mathbb Z^d$. For example, the lexicographic order on $\mathbb Z^d$ is a monomial order.

Given a monomial order $\leqslant$ on $
\mathbb Z^d$, we construction a $\mathbb Z^d$-valuation \[v:k\lbr T_1,\ldots,T_d\rbr\rightarrow\mathbb Z^d\cup\{\infty\}\] as follows. For any $\alpha=(a_1,\ldots,a_d)\in\mathbb N^d$ we denote by $T^\alpha$ the monomial $T_1^{a_1}\cdots T_d^{a_d}$. For any formal series $F$ written as
\[F(T_1,\ldots,T_d)=\lambda_\alpha T^\alpha+\sum_{\alpha<\beta}\lambda_\beta T^\beta,\quad \lambda_\alpha\neq 0, \]
we let $v_{\geqslant }(F):=\alpha$. If $F=0$ is the zero formal series, let $v(0)=\infty$. It is easy to check that the map $v$ satisfies the following axioms of valuation
: for any $(F,G)\in k\lbr T_1,\ldots,T_d\rbr^2$, one has $v(FG)=v(F)+v(G)$ and $v(F+G)\geqslant\min(v(F),v(G))$, and the equality $v(F+G)=\min(v(F),v(G))$ holds when $v(F)\neq v(G)$. In particular, if we denote by $R$ the fraction field of $k\lbr T_1,\ldots,T_d\rbr$, then the map $v:k\lbr T_1,\ldots,T_d\rbr\rightarrow\mathbb Z^d\cup\{\infty\}$ extends to a map $v:R\rightarrow \mathbb Z^d\cup\{\infty\}$ such that, for any $(F,G)\in k\lbr T_1,\ldots,T_d\rbr$ with $G\neq 0$, one has $v(F/G)=v(F)-v(G)$. The valuation map $v:R\rightarrow\mathbb Z^d\cup\{\infty\}$ allows to define a $\mathbb Z^d$-filtration $\mathcal G$ of $R$ as follows
\[\forall\,\alpha\in\mathbb Z^d,\quad \mathcal G_{\geqslant\alpha}(R):=\{f\in R\,:\,v(f)\geqslant\alpha\}.\]
Note that for $(f,g)\in R^2$ one has $v(fg)=v(f)+v(g)$ and $v(f+g)\geqslant\min(v(f),v(g))$. Therefore, for $(\alpha,\beta)\in\mathbb Z^d\times\mathbb Z^d$ one has
\begin{equation}\label{Equ: multiplicativity of g}\mathcal G_{\geqslant\alpha}(R)\cdot\mathcal G_{\geqslant\beta}(R)\subseteq\mathcal G_{\geqslant\alpha+\beta}(R).\end{equation}
For any $\alpha\in\mathbb Z^d$, we let
\[\mathcal G_{>\alpha}(R):=\{f\in R\,:\,v(f)>\alpha\}\text{ and }\mathrm{gr}_\alpha(R):=\mathcal G_{\geqslant\alpha}(R)/\mathcal G_{>\alpha}(R).\]
The relation \eqref{Equ: multiplicativity of g} shows that the $k$-algebra structure on $R$ induces by passing to graduation a $k$-algebra structure on 
\[\mathrm{gr}(R):=\bigoplus_{\alpha\in\mathbb Z}\mathrm{gr}_\alpha(R)\]
so that $\mathrm{gr}(R)$ is isomorphic to the group algebra $k[\mathbb Z^d]$.

Let $V_\sbullet=\bigoplus_{n\in\mathbb N}V_n$ be a graded sub-$k$-algebra of the polynomial ring $R[T]$ (viewed as a graded $k$-algebra with the grading by the degree on $T$). The filtration $\mathcal G$ on $R$ induces an $\mathbb Z^d$-filtration on each homogeneous component $V_n$. The direct sum of  subquotients of $V_n$ form an $\mathbb N\times\mathbb Z^d$-graded sub-$k$-algebra \[\mathrm{gr}(V_\sbullet)=\bigoplus_{(n,\alpha)\in\mathbb N\times\mathbb Z^d}\mathrm{gr}_{(n,\alpha)}(V_\sbullet)\] of $\mathrm{gr}(R)[T]\cong k[\mathbb N\times\mathbb Z^d]$. In particular, $\mathrm{gr}(V_\sbullet)$ is an integral ring, and each homogeneous component $\mathrm{gr}_{(n,\alpha)}(V_\sbullet)$ is either zero or $k$-vector space of dimension $1$. In particular, the set 
\[\Gamma(V_\sbullet):=\{(n,\alpha)\in\mathbb N_{\geqslant 1}\times\mathbb Z^d\,:\,\mathrm{gr}_{(n,\alpha)}(V_\sbullet)\neq\{0\}\}\]
is a sub-semigroup of $\mathbb N\times\mathbb Z^d$, called the \emph{Newton-Okounkov semigroup}\index{Newton-Okounkov semigroup} of $V$. The algebra $\mathrm{gr}(V_\sbullet)$ is canonically isomorphic to the semigroup $k$-algebra associated with $\Gamma(V_\sbullet)$.   Denote by $\Delta(V_\sbullet)$ the closure of the subset
\[\{n^{-1}\alpha\,:\,(n,\alpha)\in\Gamma(V_\sbullet)\}\]
of $\mathbb R^d$, called the \emph{Newton-Okounkov body}\index{Newton-Okounkov body} of $V_\sbullet$. Let $A(V_\sbullet)$ be the affine subspace of $\mathbb R^d$ the canonical projection of $\Gamma(V_\sbullet)\cap(\{1\}\times\mathbb R^d)$ in $\mathbb R^d$. By Proposition \ref{Pro:Okounkov body of graded semigroup}, $\Delta(V_\sbullet)$ is a closed convex subset of $A(V_\sbullet)$. Moreover, the relative interior of $\Delta(V_\sbullet)$ in $A(V_\sbullet)$ is not empty. The dimension of the affine space $A(V_\sbullet)$ is called the \emph{Kodaira dimension}\index{Kodaira dimension} of the graded linear series $V_\sbullet$.

\begin{prop}\label{Pro:converges to volume of Okounkov body}
Let $V_\sbullet=\bigoplus_{n\in\mathbb N}V_n$ be a graded sub-$k$-algebra of the polynomial ring $R[T]$. One has 
\[\lim_{n\in\mathbb N(V_\sbullet),\,n\rightarrow+\infty}\frac{\rang_k(V_n)}{n^\kappa}=\vol(\Delta(V_\sbullet)),\]
where $\mathbb N(V_\sbullet)$ is the set of $n\in\mathbb N$ such that $V_n\neq\{0\}$, $\kappa$ is the Kodaira dimension of $V_\sbullet$, and $\vol(\ndot)$ is the Lebesgue measure which is normalised with respect to the semi-group $\Gamma(V_\sbullet)$ as in Definition \ref{Def: normalisation of Lebesgue measure}.
\end{prop}
\begin{proof}
It is a direct consequence of Corollary \ref{Cor:convergencemesure}. 
\end{proof}

\begin{defi} Let $V_\sbullet$ be a graded sub-$k$-algebra of $R[T]$ such that $V_n$ is of finite rank for any $n\in\mathbb N$. 
\begin{enumerate}[label=\rm(\alph*)]
\item We say that $V_\sbullet$ is \emph{of subfinite type}\index{subfinite type} if it is contained in a graded sub-$k$-algebra of $R[T]$ which is of finite type (over $k$).
\item We call \emph{$\mathbb R$-filtration on $V_\sbullet$} any collection $\mathcal F_\sbullet=\{\mathcal F_n\}_{n\in\mathbb N}$, where $\mathcal F_n$ is an $\mathbb R$-filtration on $V_n$. \item Let $\delta:\mathbb N_{\geqslant 1}\rightarrow\mathbb R_{\geqslant 0}$ be a function. We say that an $\mathbb R$-filtration $\mathcal F_\sbullet$ on $V_\sbullet$ is \emph{strongly $\delta$-superadditive}\index{strongly superadditive} if for any $\ell\in\mathbb N_{\geqslant 1}$ and all $(n_1,\ldots,n_\ell)\in\mathbb N_{\geqslant 1}^\ell$ and $(t_1,\ldots,t_\ell)\in\mathbb R^\ell$, one has
\[\mathcal F_{n_1}^{t_1}(V_{n_1})\cdots\mathcal F_{n_\ell}^{t_\ell}(V_{n_\ell})\subseteq\mathcal F_{n_1+\cdots+n_\ell}^{t_1+\cdots+t_\ell-\delta(n_1)-\cdots-\delta(n_\ell)}V_{n_1+\cdots+n_\ell}.\]
We say that the $\mathbb R$-filtration $\mathcal F_\sbullet$ is \emph{$\delta$-superadditive}\index{superadditive} if the above relation  holds in the particular case where $\ell=2$, namely, for any $(n_1,n_2)\in\mathbb N_{\geqslant 1}^2$ and any $(t_1,t_2)\in\mathbb R^2$  
\[ \mathcal F_{n_1}^{t_1}(V_{n_1})\mathcal F_{n_2}^{t_2}(V_{n_2})\subseteq\mathcal F_{n_1+n_2}^{t_1+t_2-\delta(n_1)-\delta(n_2)}V_{n_1+n_2}.\]
\end{enumerate}
\end{defi}

In the following theorem, we fix a graded sub-$k$-algebra  $V_\sbullet$ of subfinite type of $R[T]$, which is equipped with an $\mathbb R$-filtration $\mathcal F_\sbullet$. We suppose in addition that $\mathbb N(V_\sbullet):=\{n\in\mathbb N\,:\,n\geqslant1,\;V_n\neq\{0\}\}$ is not empty. For each $n\in\mathbb N(V_\sbullet)$, let $\mathbb P_n$ be the Borel probability measure on $\mathbb R$ such that, for any positive Borel function $f$ on $\mathbb R$, one has 
\[\int_{\mathbb R}f(t)\,\mathbb P_n(\mathrm{d}t)=\frac{1}{\rang(V_n)}\sum_{i=1}^{\rang(V_n)}f\left({\textstyle \frac 1n}\widehat{\mu}_i(V_n,\|\ndot\|_{\mathcal F_n})\right),\]
where $\|\ndot\|_{\mathcal F_n}$ is the norm on $V_n$ associated with the $\mathbb R$-filtration $\mathcal F_n$.

\begin{theo}\label{Thm: convergence de mesure pour le cas filtraetions}
Let $\delta:\mathbb N_{\geqslant 1}\rightarrow\mathbb R_{\geqslant 0}$ be an increasing function. We suppose that, either  $\mathcal F_\sbullet$ is strongly $\delta$-superadditive and 
$\lim_{n\rightarrow+\infty}{\delta(n)}/{n}=0$,
or $\mathcal F_\sbullet$ is $\delta$-superadditive and 
$\sum_{a\in\mathbb N}{\delta(2^a)}/{2^a}<+\infty$.   Then the sequence of measures $\{\mathbb P_n\}_{n\in\mathbb N(V_\sbullet)}$ converges vaguely to a limite Borel measure $\mathbb P_{\mathcal F_\sbullet}$ on $\mathbb R$, which is the direct image of the uniform probability measure on $\Delta(V_\sbullet)^\circ$ by a concave function $G_{\mathcal F_\sbullet}:\Delta(V_\sbullet)^\circ\rightarrow\mathbb R\cup\{+\infty\}$, called the \emph{concave transform}\index{concave transform} of $\mathcal F_\sbullet$. Moreover, $\mathbb P_{\mathcal F_\sbullet}$ is either the zero measure or a probability measure, and, in the case where it is a probability measure, $\{\mathbb P_n\}_{n\in\mathbb N(V_\sbullet)}$ also converges weakly to $\mathbb P_{\mathcal F_\sbullet}$.
\end{theo}
\begin{proof}
Let $\Gamma(V_\sbullet)$ be the Newton-Oknounkov semigroup of $V_\sbullet$. Since $V_\sbullet$ is contained in an $\mathbb N$-graded sub-algebra of finite type of $R[T]$, the group $\Gamma(V_\sbullet)_{\mathbb Z}$ is discrete and the Newton-Okounkov body $\Delta(V_\sbullet)$ is compact.  For any $\gamma=(n,\alpha)\in\Gamma(V_\sbullet)$, let $\|\ndot\|_{u}$ be the subquotient norm on $\mathrm{gr}_\gamma(V_\sbullet)$ induced by $\|\ndot\|_{\mathcal F_n}$ and let $g_{\mathcal F_\sbullet}(\gamma)$ be the Arakelov degree of $(\mathrm{gr}_{\gamma}(V_\sbullet),\|\ndot\|_{\gamma})$. Since the $\mathbb R$-filtration $\mathcal F_\sbullet$ is strongly $\delta$-superadditive (resp. $\delta$-superadditive), the function $g_{\mathcal F_\sbullet}$ on $\Gamma(V_\sbullet)$ is strongly $\delta$-superadditive (resp. $\delta$-superadditive). Moreover, by Proposition \ref{Pro: property of normed vector space over trivial valued field} \ref{Item: HN filtration subquotient in trivial valuation case}, the sequence of successive slopes of $(V_n,\|\ndot\|_{\mathcal F_n})$ identifies with the sorted sequence of $\{g_{\mathcal F_\sbullet}(n,\alpha)\}_{\alpha\in\Gamma(V_\sbullet\}_n}$. Therefore the Borel probability measure $\mathbb P_n$ verifies
\[\int_{\mathbb R}f(t)\,\mathbb P_n(\mathrm{d}t)=\frac{1}{\rang(V_n)}\sum_{\alpha\in\Gamma(V_\sbullet)_n}f\left({\textstyle\frac 1n}g_{\mathcal F_\sbullet}(n,\alpha)\right).\] 
Therefore the assertion follows from Theorem \ref{Thm:concavetransform} (resp. Theorem \ref{Thm:concavetransform 2}).
\end{proof}

{
\begin{rema}
\label{Rem: dilation} We keep the notation and the hypothesis of Theorem \ref{Thm: convergence de mesure pour le cas filtraetions}.
\begin{enumerate}[label=\rm{(\arabic*)}]
\item By \eqref{Equ: sup of concave transform} we obtain that 
\begin{equation}\label{Equ: maxmal value of concave transfor is asymptotic first slope}
\sup_{x\in\Delta(V_\sbullet)^\circ}G_{\mathcal F_\sbullet}(x)=\lim_{n\in\mathbb N(V_\sbullet),\,n\rightarrow+\infty}\frac1n\widehat{\mu}_1(V_n,\norm{\ndot}_{\mathcal F_n}).
\end{equation}
\item Let $m$ be an integer, we denote by $V_\sbullet^{(m)}$ the graded sub-$k$-algebra of $R[T]$ such that $V_{n}^{(m)}=V_{nm}$ for any $n\in\mathbb N$. Then one has
\[\Gamma(V_\sbullet^{(m)})=\{(n,\alpha)\in\mathbb N_{\geqslant 1}\times\mathbb Z^d\,:\,\mathrm{gr}_{(nm,\alpha)}(V_\sbullet)\neq\{0\}\}.\]
Therefore one has $\Delta(V_\sbullet^{(m)})=m\Delta(V_\sbullet)$. Denote by $\mathcal F_\sbullet^{(m)}$ the family of filtrations $\{\mathcal F_{mn}\}_{n\in\mathbb N}$ on $V_\sbullet^{(m)}$, then one has \[\forall\,x\in\Delta(V_\sbullet)^\circ,\quad G_{\mathcal F_\sbullet^{(m)}}(mx)=mG_{\mathcal F_\sbullet}(x).\]
In particular, $\mathbb P_{\mathcal F_\sbullet^{(m)}}$ identifies with the direct image of $\mathbb P_{\mathcal F_\sbullet}$ by the dilatation map $(t\in\mathbb R)\mapsto mt$.
\end{enumerate}
\end{rema}

\begin{rema}\label{Rem: super additivity of filtered graded linear series}
Let $U_\sbullet$, $V_\sbullet$ and $W_\sbullet$ be graded sub-$k$-algebras of subfinite type of $R[T]$. Suppose that, for any $n\in\mathbb N$, one has \[U_n+V_n:=\{x+y\,:\,x\in U_n,\; y\in V_n\}\subseteq W_n.\] Then, for all $n\in\mathbb N$ and $(\alpha,\beta)\in\mathbb N^d$ such that $(n,\alpha)\in\Gamma(U_\sbullet)$ and $(n,\beta)\in\Gamma(V_\sbullet)$, one has $(n,\alpha+\beta)\in\Gamma(W_\sbullet)$. Therefore, $\Delta(U_\sbullet)+\Delta(V_\sbullet)\subseteq\Delta(W_\sbullet)$.

Assume that the graded sub-$k$-algebras $U_\sbullet$, $V_\sbullet$ and $W_\sbullet$ are equipped with $\mathbb R$-filtrations $\mathcal F^U_{\sbullet}$, $\mathcal F^V_{\sbullet}$ and $\mathcal F^W_{\sbullet}$ respectively. Let $\delta:\mathbb N_{\geqslant 1}\rightarrow\mathbb R_{\geqslant 0}$. We suppose that, either  $\mathcal F^U_\sbullet$, $\mathcal F^V_\sbullet$ and $\mathcal F^W_\sbullet$ are $\delta$-superadditive and 
$\lim_{n\rightarrow+\infty}{\delta(n)}/{n}=0$,
or $\mathcal F^U_\sbullet$, $\mathcal F^V_\sbullet$ and $\mathcal F^W_\sbullet$ are weakly $\delta$-superadditive and 
$\sum_{a\in\mathbb N}{\delta(2^a)}/{2^a}<+\infty$. Let $\epsilon:\mathbb N_{\geqslant 1}\rightarrow\mathbb R_{\geqslant 0}$ be a map such that $\lim_{n\rightarrow+\infty}\epsilon(n)/n=0$. Suppose that, for any $n\in\mathbb N_{\geqslant 1}$ and any $(t_1,t_2)\in\mathbb R^2$, one has
\[\mathcal F_n^{U,t_1}(U_n)\cdot\mathcal F_n^{V,t_2}(V_n)\subseteq\mathcal F_n^{W,t_1+t_2-\epsilon(n)}(W_n).\]
Then, for all $(n,m)\in\mathbb N_{\geqslant 1}^2$ and $(\alpha,\beta)\in\mathbb N^d$ such that $(n,\alpha)\in\Gamma(U_\sbullet)$ and $(n,\beta)\in\Gamma(V_\sbullet)$, one has 
\[g_{\mathcal F^U_\sbullet}(mn,m\alpha)+g_{\mathcal F^V_\sbullet}(mn,m\beta)\leqslant g_{\mathcal F^W_\sbullet}(mn,m\alpha+m\beta)+\epsilon(mn).\]
Dividing the two sides of the inequality by $m$, by passing to limit when $m\rightarrow+\infty$, we obtain that 
\[\widetilde g_{\mathcal F^U_\sbullet}(n,\alpha)+\widetilde g_{\mathcal F^V_\sbullet}(n,\beta)\leqslant \widetilde g_{\mathcal F^W_\sbullet}(n,\alpha+\beta).\]
Therefore, for $(x,y)\in\Delta(U_\sbullet)\times\Delta(V_\sbullet)$, one has
\[G_{\mathcal F_\sbullet^U}(x)+G_{\mathcal F_\sbullet^V}(y)\leqslant G_{\mathcal F_\sbullet^W}(x+y).\]
\end{rema}
}

\subsection{Applications to the study of the volume function}

\begin{defi}
Let $C_0$ be a non-negative real number. We say that the adelic curve $S$ \emph{satisfies the tensorial minimal slope property of level $\geqslant C_0$}\index{tensorial minimal slope property of level C0@tensorial minimal slope property of level $\geqslant C_0$} if, for any couple $(\overline E,\overline F)$ of adelic vector bundles on $S$, the following inequality holds
\begin{equation}\label{Equ: minimal slope property}\widehat{\mu}_{\min}(\overline E\otimes_{\varepsilon,\pi}\overline F)\geqslant\widehat{\mu}_{\min}(\overline E)+\widehat{\mu}_{\min}(\overline F)-C_0\ln(\rang_K(E)\cdot\rang_K(F)).\end{equation}
Recall that we have proved in Corollary \ref{Cor: tensorial minimal slope property} that, if the field $K$ is perfect, then the adelic curve $S$ satisfies the tensorial minimal slope property of level $\geqslant \frac 32\nu(\Omega_\infty)$.
\end{defi}

We let $R=\mathrm{Frac}(K\lbr T_1,\ldots,T_d\rbr)$ be the fraction field of the $K$-algebra of formal series of $d$ variables $T_1,\ldots,T_d$ and we equip $\mathbb Z^d$ with a monomial order $\leqslant$ and $R$ with the corresponding $\mathbb Z^d$-filtration as explained in Subsection \ref{Subsec: applications to the study of graded algebras}.

\begin{defi}\label{Def: graded algebra of adelic vector bundles}
We call \emph{graded $K$-algebra of adelic vector bundles with respect to $R$}\index{graded algebra of adelic vector bundles} any family $\overline E_\sbullet=\{(E_n,\xi_n)\}_{n\in\mathbb N}$ of adelic vector bundles on $S$ such that the following conditions are satisfied:
\begin{enumerate}[label=\rm(\alph*)]\item $E_\sbullet=\bigoplus_{n\in\mathbb N} E_nT^n$ forms a graded sub-$K$-algebra of subfinite type of the polynomial ring $R[T]$;
\item for any $n\in\mathbb N$, the norm family $\xi_n$ is ultrametric on $\Omega\setminus\Omega_\infty$;
\item\label{Item: submultiplicative} assume that $\xi_n$ is of the form $\{\norm{\ndot}_{n,\omega}\}_{\omega\in\Omega}$, then, for all $\omega\in\Omega$, $(n_1,n_2)\in\mathbb N_{\geqslant 1}^2$, and $(s_1,s_2)\in E_{n_1,K_\omega}\times E_{n_2,K_\omega}$, one has
$\norm{s_1\cdot s_2}_{n_1+n_2,\omega}\leqslant\norm{s_1}_{n_1,\omega}\cdot\norm{s_2}_{n_2,\omega}$.
\end{enumerate}
\end{defi}

\begin{prop}\label{Pro: convergence theorem}
Let $\overline E_\sbullet=\{(E_n,\xi_n)\}_{n\in\mathbb N}$ be a graded $K$-algebra of adelic vector bundles with respect to $R$.  For any $n\in\mathbb N$, we equip $E_n$ with the Harder-Narasimhan $\mathbb R$-filtration $\mathcal F_n$. Then the collection $\mathcal F_\sbullet=\{\mathcal F_n\}_{n\in\mathbb N}$ forms an $\mathbb R$-filtration on $E_\sbullet$ which is $\delta$-superadditive, where $\delta$ denotes the function $\mathbb N_{\geqslant 1}\rightarrow\mathbb R_{\geqslant 0}$ sending $n\in\mathbb N_{\geqslant 1}$ to $C\ln(\rang_K(E_n))$.
\end{prop}
\begin{proof}
Let $n_1$ and $n_2$ be elements of $\mathbb N_{\geqslant 1}$. By the condition \ref{Item: submultiplicative} in Definition \ref{Def: graded algebra of adelic vector bundles}, for any $\omega\in\Omega$ and $s_1^{(1)}\otimes s_2^{(1)}+\cdots+s_1^{(N)}\otimes s_2^{(N)}\in E_{n_1,K_\omega}\otimes_{K_\omega} E_{n_2,K_\omega}$, one has
\[\norm{s_1^{(1)}s_2^{(1)}+\cdots+s_1^{(N)}s_2^{(N)}}_{n_1+n_2,\omega}\leqslant\begin{cases}
\max_{i\in\{1,\ldots,N\}}\norm{s_1^{(i)}}_{n_1,\omega}\cdot\norm{s_2^{(i)}}_{n_2,\omega},&\omega\in\Omega\setminus\Omega_\infty,\\
\sum_{i=1}^N\norm{s_1^{(i)}}_{n_1,\omega}\cdot\norm{s_2^{(i)}}_{n_2,\omega},&\omega\in\Omega_\infty.
\end{cases}\] 
Therefore, the canonical $K_\omega$-linear map $E_{n_1,K_\omega}\otimes_{K_\omega}E_{n_2,K_\omega}\rightarrow E_{n_1+n_2,K_\omega}$ is of operator norm $\leqslant 1$.
Let $F_1$ and $F_2$ be non-zero vector subspace of $E_{n_1}$ and $E_{n_2}$, respectively, and let $G$ be the image of $F_1\otimes_K F_2$ by the canonical $K$-linear map $E_{n_1}\otimes_K E_{n_2}\rightarrow E_{n_1+n_2}$. By Proposition \ref{Pro:inegalitdepente} \ref{Item: slope inequality mu min}, one has 
\[\begin{split}&\quad\;\widehat{\mu}_{\min}(\overline G)\geqslant\widehat{\mu}_{\min}(\overline F_1\otimes_{\varepsilon,\pi}\overline F_2)\\&\geqslant\widehat{\mu}_{\min}(\overline F_1)+\widehat{\mu}_{\min}(F_2)-C(\ln(\rang_K(F_1)))-C(\ln(\rang_K(F_2)))\\&\geqslant\widehat{\mu}_{\min}(\overline F_1)+\widehat{\mu}_{\min}(F_2)-C(\ln(\rang_K(E_{n_1})))-C(\ln(\rang_K(E_{n_2}))),
\end{split} \]
where the second inequality comes from \eqref{Equ: tensorial minimal slope superadditivity}.
By Proposition \ref{Pro: general hn filtration}, we obtain that, for any $(t_1,t_2)\in\mathbb R^2$, one has
\[\mathcal F_{n_1}^{t_1}(E_{n_1})\cdot\mathcal F_{n_2}^{t_2}(E_{n_2})\subseteq\mathcal F_{n_1+n_2}^{t_1+t_2-\delta(n_1)-\delta(n_2)}(E_{n_1+n_2}).\] 
The proposition is thus proved.
\end{proof}

\begin{coro}\label{Cor: convergence theorem for graded adelic vector bundles}
Let $\overline E_\sbullet=\{(E_n,\xi_n)\}_{n\in\mathbb N}$ be a graded $K$-algebra of adelic vector bundles with respect to $R$. We assume that $\mathbb N(E_\sbullet)$ does not reduce to $\{0\}$ and we denote by $q\in\mathbb N$ a generator of the group $\mathbb Z(E_\sbullet)$. Suppose in addition that  
\begin{equation}\label{Equ: somme de log rang est fini}\sum_{
a\in\mathbb N,\,
2^aq\in\mathbb N(E_\sbullet)
}\frac{\ln(\rang_K(E_{2^aq}))}{2^a}<+\infty.\end{equation}
For each $n\in\mathbb N(V_\sbullet)$, let $\mathbb P_n$ be the Borel probability measure on $\mathbb R$ such that, for any positive Borel function $f$ on $\mathbb R$, one has
\[\int_{\mathbb R}f(t)\,\mathbb P_n(\mathrm{d}t)=\frac{1}{\rang(E_n)}\sum_{\rang(E_n)}f(\textstyle{\frac 1n}\widehat{\mu}_i(\overline E_n)).\]
Then the sequence of measures $\{\mathbb P_n\}_{n\in\mathbb N(E_\sbullet)}$ converges vaguely to a limite Borel measure $\mathbb P_{\overline E_\sbullet}$, which is the direct image of the uniform distribution on $\Delta(E_\sbullet)$ by a concave function $G_{\overline E_\sbullet}:\Delta(\overline E_\sbullet)\rightarrow\mathbb R\cup\{+
\infty\}$. Moreover, the limite measure is either zero or a Borel probability measure, and in the latter case the sequence $\{\mathbb P_n\}_{n\in\mathbb N(E_\sbullet)}$ also converges weakly to $\mathbb P_{\overline{E}_\sbullet}$.
\end{coro}
\begin{proof}
This is a direct consequence of Proposition \ref{Pro: convergence theorem} and Theorem \ref{Thm: convergence de mesure pour le cas filtraetions}.
\end{proof}

\begin{rema} We keep the notation and the conditions of Corollary \ref{Cor: convergence theorem for graded adelic vector bundles}. We suppose that  $(\frac1n\widehat{\mu}_{1}(\overline E_n))_{n\in\mathbb N(V_\sbullet),\,n\geqslant 1}$ is bounded from above. Then the limit measure $\mathbb P_{\overline{E}_\sbullet}$ is a probability measure. The weak convergence of $\{\mathbb P_n\}_{n\in\mathbb N(E_\sbullet)}$ to $\mathbb P_{\overline{E}_\sbullet}$ implies that 
\begin{equation}\label{Equ: convergence de volume arithm}\begin{split}\int_{\mathbb R}\max\{t,0\}\,\mathbb P_{\overline{E}_\sbullet}(\mathrm{d}t)&=\lim_{n\in\mathbb N(E_\sbullet),\,n\rightarrow+\infty}
\frac{1}{n\rang(E_n)}\sum_{i=1}^{\rang(E_n)}\max\{\widehat{\mu_i}(\overline E_n),0\}\\
&=\lim_{n\in\mathbb N(E_\sbullet),\,n\rightarrow+\infty}\frac{\widehat{\deg}_+(\overline E_n)}{n\rang(E_n)}.
\end{split}\end{equation}
If in addition the sequence $(\frac1n\widehat{\mu}_{\min}(\overline E_n))_{n\in\mathbb N(V_\sbullet),\,n\geqslant 1}$ is bounded from below, then one has
\begin{equation}\label{Equ: convergence de mu normalise}\begin{split}\int_{\mathbb R}t\,\mathbb P_{\overline{E}_\sbullet}(\mathrm{d}t)
&=\lim_{n\in\mathbb N(E_\sbullet),\,n\rightarrow+\infty}
\frac{1}{n\rang(E_n)}\sum_{i=1}^{\rang(E_n)}\widehat{\mu_i}(\overline E_n)
\\
&=\lim_{n\in\mathbb N(E_\sbullet),\,n\rightarrow+\infty}\frac{\widehat{\mu}(\overline E_n)}{n}.
\end{split}\end{equation}
\end{rema}

\begin{prop}\label{Pro: super additivity of G}
Let 
\[
\begin{cases}
\overline U_\sbullet=\{(U_n,\{\norm{\ndot}_{U_n,\omega}\}_{\omega\in\Omega})\}_{n\in\mathbb N},\\
\overline V_\sbullet=\{(V_n,\{\norm{\ndot}_{V_n,\omega}\}_{\omega\in\Omega})\}_{n\in\mathbb N},\\
\overline W_\sbullet=\{(W_n,\{\norm{\ndot}_{W_n,\omega}\}_{\omega\in\Omega})\}_{n\in\mathbb N}
\end{cases}
\]
be graded $K$-algebras of adelic vector bundles with respect to $R$. We assume that 
\begin{equation}
\sum_{a\in\mathbb N,\;2^aq\in\mathbb N(W_\sbullet)}\frac{\ln(\rang(W_{2^aq}))}{2^a}<+\infty,\end{equation}
where $q\in\mathbb N$ is a generator of the group $\mathbb Z(W_\sbullet)$. Suppose that, for any $n\in\mathbb N$ one has $U_n\cdot V_n\subseteq W_n$, and 
\begin{equation}\label{Equ: submultiplicative section norm}\forall\,\omega\in\Omega,\;\forall\,(s,s')\in U_{n,K_\omega}\times V_{n,K_\omega},\quad \norm{ss'}_{W_n,\omega}\leqslant\norm{s}_{U_n,\omega}\cdot\norm{s'}_{V_n,\omega}.\end{equation}
Then, for any $(x,y)\in\Delta(U_\sbullet)\times\Delta(V_\sbullet)$, one has
\begin{equation}\label{Equ: super additivity of G}G_{\overline W_\sbullet}(x+y)\geqslant G_{\overline U_\sbullet}(x)+G_{\overline V_\sbullet}(y).\end{equation}
\end{prop}
\begin{proof} Denote by $\delta:\mathbb N_{\geqslant 1}\rightarrow\mathbb R_{\geqslant 0}$ the function sending $n\in\mathbb N_{\geqslant 1}$ to
\[C\ln(\rang_K(U_n))+C\ln(\rang_K(V_n))\]
Let $n\in\mathbb N$, $n\geqslant 1$. Suppose that $E$ is a non-zero vector subspace of $U_n$ and $F$ is a non-zero vector subspace of $V_n$. Since the adelic curve $S$ satisfies the tensorial minimal slope superadditivity of level $\geqslant C$, one has 
\[\widehat{\mu}_{\min}(\overline E\otimes_{\varepsilon,\pi}\overline F)\geqslant\widehat{\mu}_{\min}(E)+\widehat{\mu}_{\min}(F)-\delta(n).\]
Moreover, by \eqref{Equ: submultiplicative section norm} the canonical $K$-linear map $E\otimes F\rightarrow W_n$ has height $\leqslant 0$ if we consider the adelic vector bundles $\overline E\otimes_{\varepsilon,\pi}\overline F$ and $(W_n,\{\norm{\ndot}_{W_n,\omega}\}_{\omega\in\Omega}$. Therefore, if we denote by $\mathcal F_n^U$, $\mathcal F_n^V$ and $\mathcal F_n^W$ the Harder-Narasimhan $\mathbb R$-filtrations of $\overline U_n$, $\overline V_n$ and $\overline W_n$ respectively, then, for any $(t,t')\in\mathbb R^2$,
\[\mathcal F^{U,t}_n(U_n)\cdot\mathcal F^{V,t'}_n(V_n)\subseteq\mathcal F_n^{W,t+t'-\delta(n)}.\]
By Remark \ref{Rem: super additivity of filtered graded linear series}, we obtain the inequality \eqref{Equ: super additivity of G}.
\end{proof}

\begin{rema}
Let $V_\sbullet$ be a graded $k$-algebra. We say that $V_\sbullet$ is \emph{of subfinite type}\index{subfinite type} if it is contained in a graded $k$-algebra of finite type. It is not true that any integral graded $k$-algebra of subfinite type can be identifies as a graded sub-$k$-algebra of the ring of polynomials (of one variable) with coefficients in the fraction field of the formal series ring (with finitely many variables) over $k$ since the latter condition implies that $V_\sbullet$ admits a valuation of one-dimensional leaves and in particular $V_\sbullet$ is geometrically integral. We refer to \cite[Remark 5.3]{Chen_Ikoma} for more details. Moreover, the combination of the methods in \cite[\S4]{Chen_JTNB} and \cite[\S5]{Chen_Ikoma} 
allows to obtain a generalisation of Corollary \ref{Cor: convergence theorem for graded adelic vector bundles} and Proposition \ref{Pro: super additivity of G} to the case of graded algebras of adelic vector bundles whose underlying graded $k$-algebras are integral domain and of subfinite type over $k$. Note that the $\mathbb R$-filtration by slopes of a graded algebra of adelic vector bundles is not necessarily superadditive and we need an argument similar to the Step 2 in the proof of Theorem \ref{Thm:concavetransform} in order to replace the $\mathbb R$-filtration by slopes by a superadditive $\mathbb R$-filtration while keeping the asymptotic behaviour of the distribution of average on the jump points of the $\mathbb R$-filtrations. The approach can serve to remove the hypothesis that the scheme $X$ admits a regular rational points in Theorem \ref{Thm: limit theorem for adelic line bundles} and Theorem \ref{Thm: Brunn-Minkowski for adelic line bundles}.
\end{rema}

\begin{defi}\label{Def: strong tensorial minimal slope property}
Let $C_1$ be a non-negative real number. We say that the adelic curve $S$ satisfies \emph{the strong tensorial minimal slope property of level $\geqslant C_1$} if, for any integer $n\in\mathbb N_{\geqslant 2}$ and any family $\{\overline E_i\}_{i=1}^n$  of $d$ non-zero adelic vector bundles on $S$, the following inequality holds
\[\widehat{\mu}(\overline E_{1}\otimes_{\varepsilon,\pi}\cdots\otimes_{\varepsilon,\pi}\overline E_n)\geqslant\sum_{i=1}^n\big(\widehat{\mu}_{\min}(\overline E_i)-C_1\ln(\rang_K(E_i))\big).\]
Note Corollary \ref{Cor: tensorial minimal slope property} shows that the adelic curve $S$ satisfies the strong tensorial minimal slope property of level $\geqslant\frac32\nu(\Omega_\infty)$, provided that the field $K$ is perfect.
\end{defi}

\begin{rema}
Let $C$ be a non-negative real number. We say that the adelic curve $S$ satisfies the strong minimal slope property of level $\geqslant C$. Let $\overline E_\sbullet=\{(E_n,\xi_n)\}_{n\in\mathbb N}$ be a graded $K$-algebra of adelic vector bundles with respect to $R$.  For any $n\in\mathbb N$, we equip $E_n$ with the Harder-Narasimhan $\mathbb R$-filtration $\mathcal F_n$. Then the collection $\mathcal F_\sbullet=\{\mathcal F_n\}_{n\in\mathbb N}$ forms an $\mathbb R$-filtration on $E_\sbullet$ which is strongly $\delta$-superadditive, where $\delta$ denotes the function $\mathbb N_{\geqslant 1}\rightarrow\mathbb R_{\geqslant 0}$ sending $n\in\mathbb N_{\geqslant 1}$ to $C\ln(\rang_K(E_n))$. Therefore, by Theorem \ref{Thm: convergence de mesure pour le cas filtraetions}, if the condition 
\[\lim_{n\rightarrow+\infty}\frac{\ln(\rang(E_n))}{n}=0\]
is satisfied, then the sequence of Borel probability measures $\{\mathbb P_n\}_{n\in\mathbb N(E_\sbullet)}$, defined by
\[\int_{\mathbb R}f(t)\,\mathbb P_n(\mathrm{d}t)=\frac{1}{\rang(E_n)}\sum_{i=1}^{\rang(E_n)}\delta_{\frac 1n\widehat{\mu}_{i}(\overline E_n)},\]
converges vaguely to a Borel measure $\mathbb P_{\overline{E}_\sbullet}$, which is either the zero measure or a Borel probability measure. In the latter case, the sequence $\{\mathbb P_n\}_{n\in\mathbb N(E_\sbullet)}$ converges weakly to $\mathbb P_{\overline{E}_{\sbullet}}$. Similarly, the assertion of Proposition \ref{Pro: super additivity of G}
 holds under the condition
\[\lim_{n\rightarrow+\infty}\frac{\ln(\rang(W_n))}{n}=0.\]
\end{rema}

\begin{rema}
In the number field setting, Yuan \cite{Yuan09,MR3369367}
has proposed another method to associate to each adelic line bundle a convex body which computes the arithmetic volume of the adelic line bundle. His method relies on multiplicity estimates of arithmetic global sections with respect to a flag of subvarieties of the fibre of the arithmetic variety over a finite place, which is similar to \cite{Lazarsfeld_Mustata08}. Note that in the general setting of adelic curves the set of ``arithmetic global sections'' is not necessarily finite and the classic formula relating the  volume function and the asymptotic behaviour of ``arithmetic global sections'' does not hold in general.
\end{rema}

\section{Asymptotic invariants of graded linear series}

In this section, we fix an integral projective $K$-scheme $X$ and denote by $\pi:X\rightarrow\Spec K$ the structural morphism. Let $d$ be the Krull dimension of the scheme $X$.

\subsection{Asymptotic maximal slope}
Let $(L,\varphi)$ be an adelic line bundle on an integral projective $K$-scheme $\pi:X\rightarrow\Spec K$. For any $n\in\mathbb N_{\geqslant 1}$, the metric family $n\varphi$ on $L^{\otimes n}$ induces a norm family $\{\|\ndot\|_{n\varphi_\omega}\}_{\omega\in\Omega}$ on the linear series $\pi_*(L^{\otimes n})=H^0(X,L^{\otimes n})$ which we denote by $\pi_*(n\varphi)$. By Theorem \ref{Thm: Fubini-Study dominated} and \ref{Thm: measurability of linear series}, the pair $(\pi_*(L^{\otimes n}),\pi_*(n\varphi))$ forms an adelic vector bundle on $S$, which we denote by $\pi_*(L^{\otimes n},n\varphi)$. Note that $\pi_*(n\varphi)$ is ultrametric on $\Omega\setminus\Omega_\infty$. We can then compute divers arithmetic invariants of these adelic vector bundles. The asymptotic behaviour of these arithmetic invariants describes the positivity of the adelic line bundle $(L,\varphi)$.

Let $(L,\varphi)$ be an adelic line bundle on $X$. We define
\[\nu_1^{\mathrm{asy}}(L,\varphi):=\limsup_{n\rightarrow+\infty}\frac{\nu_1(\pi_*(L^{\otimes n},n\varphi))}{n},\]
called the \emph{asymptotic first minimum}\index{asymptotic first minimum} of $(L,\varphi)$. Similarly, we define
\begin{equation}\label{Equ: asymptotic maximal slope}\widehat{\mu}_{\max}^{\mathrm{asy}}(L,\varphi):=\limsup_{n\rightarrow+\infty}\frac{\widehat{\mu}_{\max}(\pi_*(L^{\otimes n},n\varphi))}{n},\end{equation}
called the \emph{asymptotic maximal slope}\index{asymptotic maximal slope} of $(L,\varphi)$. Note that all adelic vector bundles $\pi_*(L^{\otimes n},n\varphi)$ are ultrametric on $\Omega\setminus\Omega_\infty$. Therefore, by Remark \ref{Rem: comparison between mu max and mu 1} and the fact that 
\[\lim_{n\rightarrow+\infty}\frac{\ln(\rang_K(H^0(X,L^{\otimes n})))}{n}=0,\]
we obtain that
\begin{equation}\widehat{\mu}_{\max}^{\mathrm{asy}}(L,\varphi)=\limsup_{n\rightarrow+\infty}\frac{\widehat{\mu}_1(\pi_*(L^{\otimes n},n\varphi))}{n}.\end{equation}

Let $\mathbb K$ be either $\mathbb Z$ or $\mathbb Q$ or $\mathbb R$. From now on we assume that $\mathbb K = \mathbb Z$ or $X$ is normal.
Let $(D,g)$ be an adelic $\mathbb K$-Cartier divisor on $X$. Note that if $X$ is not normal and $\mathbb K$ is either $\mathbb Q$ or $\mathbb R$,
then $H^0_{\mathbb K}(X, D)$ is not necessarily a vector space over $K$ (cf. Example~\ref{example:R:Cartier:vs:Cartier}). Similarly we can define
$\nu_1^{\mathrm{asy}}(D, g)$ and $\widehat{\mu}_{\max}^{\mathrm{asy}}(D, g)$ as follows:
\[
\begin{cases}
{\displaystyle \nu_1^{\mathrm{asy}}(D, g):=\limsup_{n\rightarrow+\infty}\frac{\nu_1(H^0_{\mathbb K}(X, nD), \xi_{ng})}{n}}, \\
{\displaystyle \widehat{\mu}_{\max}^{\mathrm{asy}}(D, g) := \limsup_{n\rightarrow+\infty}\frac{\widehat{\mu}_{\max}(H^0_{\mathbb K}(X, nD), \xi_{ng})}{n}.}
\end{cases}
\]
In the case where $\mathbb K = \mathbb Z$,
\[
\nu_1^{\mathrm{asy}}(D, g) = \nu_1^{\mathrm{asy}}(\mathcal O_X(D),\varphi_g)\quad\text{and}\quad
\widehat{\mu}_{\max}^{\mathrm{asy}}(D, g) = \widehat{\mu}_{\max}^{\mathrm{asy}}(\mathcal O_X(D), \varphi_g),
\]
where $\varphi_g$ is the metric family of $L$ defined by $g$.

\begin{prop}\label{prop:Minkowski:nu:1:mu:max}
Let $(D, g)$ be an adelic $\mathbb K$-Cartier divisor on $S$. Then one has $\nu_1^{\mathrm{asy}}(D, g)\leqslant\widehat{\mu}_{\max}^{\mathrm{asy}}(D, g)$. The equality holds when $S$ satisfies the Minkowski property of certain level (see Definition \ref{Def: Minkowski property}). 
\end{prop}
\begin{proof}
By \eqref{Equ: udeg bounded by mu max}, for any $n\in\mathbb N_{n\geqslant 1}$ one has
\[\nu_1(H^0_{\mathbb K}(X, nD), \xi_{ng})\leqslant\widehat{\mu}_{\max}(H^0_{\mathbb K}(X, nD), \xi_{ng}).\]
Therefore $\nu_1^{\mathrm{asy}}(D,g)\leqslant\widehat{\mu}_{\max}^{\mathrm{asy}}(D,g)$. 

If the adelic curve $S$ satisfies the Minkowski property of level $\geqslant C$, where $C\geqslant 0$, then  
\[\nu_1(H^0_{\mathbb K}(X, nD), \xi_{ng})\geqslant\widehat{\mu}_{\max}(H^0_{\mathbb K}(X, nD), \xi_{ng})-C\ln(\rang_K(H^0_{\mathbb K}(X, nD))).\]
Since $X$ is a projective scheme, one has 
\[\rang_K(H^0_{\mathbb K}(X, nD))=O(n^{\dim(X)}).\]
Hence
\[\lim_{n\rightarrow+\infty}\frac{\ln(\rang_K(H^0_{\mathbb K}(X, nD)))}{n}=0.\]
Therefore $\nu_1^{\mathrm{asy}}(D, g)\geqslant\widehat{\mu}_{\max}^{\mathrm{asy}}(D, g)$.
\end{proof}

Let $(D, g)$ be an adelic $\mathbb R$-Cartier divisor on $X$.
Let $\mathbb K$ be either $\mathbb Q$ or $\mathbb R$.
We set 
\[
\begin{cases}
\Gamma^{\times}_{\mathbb K}(D) = \{ \phi \in K(X)^{\times} \otimes_{\mathbb Z} \mathbb K \,:\, D + (\phi) \geqslant_{\mathbb K} 0 \},\\[1ex]
\nu^{\mathrm{asy}}_{1,\mathbb K}(D,g) = \begin{cases}
\sup \left\{ \widehat{\deg}_{\xi_g}(s) \,:\, s \in \Gamma^{\times}_{\mathbb K}(D) \right\} & \text{if $\Gamma^{\times}_{\mathbb K}(D) \not= \emptyset$},\\
-\infty & \text{if $\Gamma^{\times}_{\mathbb K}(D) = \emptyset$},
\end{cases}
\end{cases}
\]
where (cf. Corollary~\ref{coro:adelic:R:Cartier:div:integrable})
\[ \hdeg_{\xi_g}(s) = -\int_{\Omega}  \ln \| s \|_{g_\omega}\, \nu(d\omega).\] 
Note that  
$\nu^{\mathrm{asy}}_{1}(D,g) = \nu^{\mathrm{asy}}_{1,\mathbb Q}(D,g) \leqslant \nu^{\mathrm{asy}}_{1,\mathbb R}(D,g)$.

\begin{prop}\label{prop:nu1:concave}
We assume that $X$ is normal. Let $(D, g)$ and $(D', g')$ be adelic $\mathbb R$-Cartier divisors on $X$.
Then one has the following:
\begin{enumerate}[label=\rm(\arabic*)]
\item
$\nu^{\mathrm{asy}}_{1, \mathbb K}((D, g) + (D', g')) \geqslant \nu^{\mathrm{asy}}_{1,\mathbb K}(D,g) + \nu^{\mathrm{asy}}_{1,\mathbb K}(D',g')$.

\item
$\nu^{\mathrm{asy}}_{1,\mathbb K}(a(D,g)) = a \nu^{\mathrm{asy}}_{1,\mathbb K}(D,g)$ for all $a \in \mathbb K_{\geqslant 0}$.  
%
\end{enumerate}
\end{prop}

\begin{proof}
(1) Clearly we may assume that $\Gamma^{\times}_{\mathbb K}(D) \not= \emptyset$ and $\Gamma^{\times}_{\mathbb K}(D') \not= \emptyset$.
If $\phi \in \Gamma^{\times}_{\mathbb K}(D)$ and $\phi' \in \Gamma^{\times}_{\mathbb K}(D')$,
then $\| \phi \phi'\|_{g+g'} \leqslant \| \phi \|_{g} \| \phi' \|_{g'}$, so that
\[
-\log \| \phi \|_{g} -\log \| \phi' \|_{g'} \leqslant -\log \| \phi\phi' \|_{g+g'} \leqslant \nu_{1,\mathbb K}^{\mathrm{asy}}(D+D', g+g').
\]
Therefore one has (1).

(2) Clearly we may assume that $a > 0$ and $\Gamma^{\times}_{\mathbb K}(D) \not= \emptyset$.
Then one has a bijective correspondence $(s \in \Gamma^{\times}_{\mathbb K}(D)) \mapsto (s^a \in \Gamma^{\times}_{\mathbb K}(aD))$.
Moreover, $\widehat{\deg}_{\xi_{ag}}(s^a) = a \widehat{\deg}_{\xi_g}(s)$ for $s \in \Gamma^{\times}_{\mathbb K}(D)$. Thus the assertion follows.
\end{proof}

\begin{theo}\label{thm:nu1:Q:R:comp}
We assume that $X$ is normal. Let $(D, g)$ be an adelic $\mathbb R$-Cartier divisor on $X$. If $\Gamma^{\times}_{\mathbb Q}(D) \not= \emptyset$,
then $\nu^{\mathrm{asy}}_{1,\mathbb Q}(D,g) = \nu^{\mathrm{asy}}_{1,\mathbb R}(D,g)$. In particular, if $D$ is big, then
$\nu^{\mathrm{asy}}_{1}(a(D,g)) = a \nu^{\mathrm{asy}}_{1}(D,g)$ for all $a \in \mathbb R_{\geqslant 0}$. 
\end{theo}

\begin{proof}
By our assumption, we can find $\psi \in \Gamma^{\times}_{\mathbb Q}(D)$. Then the map
\[\alpha_{\psi} : \Gamma^{\times}_{\mathbb K}(D) \to \Gamma^{\times}_{\mathbb K}(D + (\psi))\] given by $\phi \mapsto \phi\psi^{-1}$
is bijective and, for $\phi \in \Gamma^{\times}_{\mathbb K}(D)$,
$\| \phi \|_g = \| \alpha_{\psi}(\phi) \|_{g - \log |\psi|}$, so that
\[
\nu_{1,\mathbb K}^{\mathrm{asy}}(D, g) = \nu_{1,\mathbb K}^{\mathrm{asy}}(D + (\psi), g - \log|\psi|).
\]
Therefore we may assume that $D$ is effective.
Moreover, for an integrable function $\varphi$ on $\Omega$, 
\[
\nu_{1,\mathbb K}^{\mathrm{asy}}(D, g+\varphi) = \nu_{1,\mathbb K}^{\mathrm{asy}}(D, g) + \int_{\Omega} \varphi \nu(d\omega),\\
\]
so that we may further assume that \[{\displaystyle \int_{\Omega} - \log \| 1 \|_{g_{\omega}} \nu(d\omega) \geqslant 0}.\]

For $\phi \in \Gamma^{\times}_{\mathbb R}(D)$, we choose $s_1, \ldots, s_r \in K(X)^{\times} \otimes_{\mathbb Z} \mathbb Q$
and $a_1, \ldots, a_r \in \mathbb R$ such that $\phi = s_1^{a_1} \cdots s_r^{a_r}$ and
$a_1, \ldots, a_r$ are linearly independent over $\mathbb Q$. 
We set $\| x \|_0 = |x_1| + \cdots + |x_{\textcolor{mgreen}{r}}|$ and $s^x = s_1^{x_1} \cdots s_r^{\textcolor{mgreen}{x_r}}$ for $x = (x_1, \ldots, x_r) \in \mathbb R^r$, so that
if we denote $(a_1, \ldots, a_r)$ by $\alpha$, then $\phi = s^{\alpha}$. By Proposition~\ref{prop:approximation:R:sec:by:Q:sec},
for any a positive rational number $\varepsilon$,
there is a positive number $\delta$ such that if $\| \alpha' - \alpha \|_0 \leqslant \delta$ for $\alpha' \in \mathbb R^r$,
then $(1 + \varepsilon)D + (s^{\alpha'})$ is effective.
We choose a basis $\{ \omega_1, \ldots, \omega_r \}$ of $\mathbb Q^r$ such that
$\| \omega_j - a \|_0 \leqslant \delta$ for all $j \in \{1, \ldots, r \}$, so that $(1+\varepsilon)D + (s^{\omega_j}) \geqslant 0$ for all $j$.
Here we set $\alpha = \lambda_1 \omega_1 + \cdots + \lambda_r \omega_r$.
Further, if we define a norm $\|\ndot\|_{\omega}$ by $\| x_1 \omega_1 + \cdots + x_r \omega_r\|_{\omega} = |x_1| + \cdots + |x_r|$ for $x_1, \ldots, x_r \in \mathbb R$,
then there is a positive constant $C$ such that $C\|\ndot\|_0 \leqslant \|\ndot\|_{\omega}$.
Note that for any $t > 0$, there is $\alpha' \in \mathbb Q^r$ such that if we set $\alpha' = \lambda'_1 \omega_1 + \cdots + \lambda'_r \omega_r$, then
$\lambda'_j \geqslant \lambda_j$ ($\forall\, j$) and $\| \alpha' - \alpha \|_{\omega} \leqslant t$. Indeed, for each $j$, one can find
$\lambda'_j \in \mathbb Q$ such that $0 \leqslant \lambda'_j - \lambda_j \leqslant t/r$, and hence $\| \alpha' - \alpha \|_{\omega} \leqslant t$.
Therefore, we can also choose a sequence $\{ \alpha_n \}_{n=1}^{\infty}$ of $\mathbb Q^r$ 
with the following properties:
\begin{enumerate}[label=\rm(\roman*)]
\item $\alpha_n \in \alpha + \mathbb R_{\geqslant 0} \omega_1 + \cdots + \mathbb R_{\geqslant 0} \omega_r$ for all $n  \geqslant 1$.
\item $\| \alpha_n - \alpha \|_{\omega} \leqslant \min \{ \varepsilon/((1+\varepsilon)n), C\delta \}$ for all $n  \geqslant 1$.
\end{enumerate}
Since $\| \alpha_n - \alpha \|_{0} \leqslant \delta$ by (ii),  $(1 + \varepsilon) D + (s^{\alpha_n}) \geqslant 0$ for all $n \geqslant 1$. Moreover,
if we set $\alpha_n - \alpha = \sum_{j=1}^r \lambda_j^{(n)} \omega_j$, then 
$\lambda_j^{(n)} \geqslant 0\ (\forall\, j)$ and
$\sum_{j=1}^r \lambda_j^{(n)} \leqslant \varepsilon/((1+\varepsilon)n)$.
Therefore, if we denote $\varepsilon - (1+\varepsilon) \sum_{j=1}^r \lambda_j^{(n)}$
by $\kappa_n$, then $\kappa_n \geqslant 0$ and, for each $\omega \in \Omega$,
\[
|s^{\alpha_n}|_{(1+\varepsilon)g_{\omega}} = | s^{\alpha}|_{g_{\omega}} |s^{\omega_1\lambda_1^{(n)}} \cdots s^{\omega_r\lambda_r^{(n)}} |_{\varepsilon g_{\omega}}
= |\phi|_{g_{\omega}} | s^{\omega_1} |_{(1+\varepsilon) g_{\omega}}^{\lambda_1^{(n)}} \cdots | s^{\omega_r} |_{(1+\varepsilon) g_{\omega}}^{\lambda_r^{(n)}}
| 1 |_{\kappa_n g_{\omega}},
\]
so that $\| s^{\alpha_n}\|_{(1+\varepsilon)g_{\omega}} \leqslant \| \phi\|_{g_{\omega}} \| s^{\omega_1} \|_{(1+\varepsilon)g_\omega}^{\lambda_1^{(n)}} \cdots \| s_{\omega_r} \|_{(1+\varepsilon)g_{\omega}}^{\lambda_r^{(n)}}\| 1 \|_{g_{\omega}}^{\kappa_n}$.
Therefore, since
\[
\int_{\Omega}- \log \| s^{\alpha_n} \|_{(1+\varepsilon)g_{\omega}}\nu(d\omega) \leqslant \nu_{1,\mathbb Q}^{\mathrm{asy}}((1+\varepsilon)(D,g))
\]
and ${\displaystyle \kappa_n \int_{\Omega} -\log \| 1 \|_{g_{\omega}} \nu(d\omega) \geqslant 0}$,
one has
\[
 \int_{\Omega}  - \log \| \phi \|_{g_{\omega}} \nu(d\omega)+ \sum_{j=1}^r\lambda_j^{(n)} \int_{\Omega} - \log \| s^{\omega_j} \|_{(1+\varepsilon) g_\omega} \nu(d\omega)\leqslant \nu_{1,\mathbb Q}^{\mathrm{asy}}((1+\varepsilon)(D, g)),
\]
so that taking $n\to\infty$, we obtain
\[
\int_{\Omega} - \log \| \phi \|_{g_\omega}  \nu(d\omega)\leqslant \nu_{1,\mathbb Q}^{\mathrm{asy}}((1+\varepsilon)(D,g)),
\]
and hence, as 
$\varepsilon$ is a rational number, by Proposition~\ref{prop:nu1:concave}, (2), one can see
\[
\nu_{1,\mathbb R}^{\mathrm{asy}}(D, g) \leqslant \nu_{1,\mathbb Q}^{\mathrm{asy}}((1+\varepsilon)(D, g)) = (1+\varepsilon) \nu_{1,\mathbb Q}^{\mathrm{asy}}(D, g),
\]
which implies $\nu_{1,\mathbb R}^{\mathrm{asy}}(D, g) \leqslant \nu_{1,\mathbb Q}^{\mathrm{asy}}(D, g)$, as required.
\end{proof}

\begin{prop}\label{Pro: comparison ess minimum and asy max slope}
Let $(D, g)$ be an adelic $\mathbb K$-Cartier divisor on $X$.
We assume that either $\mathbb K = \mathbb Z$ or $X$ is normal.
Then one has
\[
\widehat{\mu}_{\mathrm{ess}}(D,g)\geqslant\widehat{\mu}_{\max}(H^0_{\mathbb K}(X, D), \xi_g).
\]
In particular,
\begin{equation}\label{Equ: bound of mu max asy by mu ess}
\widehat{\mu}_{\mathrm{ess}}(D,g)\geqslant\widehat{\mu}_{\max}^{\mathrm{asy}}(D,g).
\end{equation}
\end{prop}

\begin{proof}
The second inequality is a consequence of the first inequality because $\widehat{\mu}_{\mathrm{ess}}(nD,ng) = n\widehat{\mu}_{\mathrm{ess}}(D,g)$.

Let $U$ be a non-empty Zariski open set of $X$ given by
\[
\{ x \in X \,:\, \text{$X \to \Spec K$ is smooth at $x$ and $x \not\in \mathrm{Supp}_{\mathbb K}(D)$} \}.
\]
Note that $\exp(-g_{\omega})$ is a positive continuous function on $U_{\omega}^{\mathrm{an}}$ for each $\omega \in \Omega$ and that,
for $\phi \in H^0_{\mathbb K}(X, D)$ and $x \in U$, one has $\phi \in \mathcal O_{X,x}$.
Let $t$ be a real number such that $\widehat{\mu}_{\mathrm{ess}}(D,g) < t$. 
Then there is an infinite subset $\Lambda$ of $U(K^{\mathrm{ac}})$ such that
$\Lambda$ is Zariski dense in $U$ and $h_{(D,g)}(P) \leqslant t$ for all $P \in \Lambda$.

Let $F$ be a non-zero vector subspace of $H^0_{\mathbb K}(X,D)$. Then
there exist $P_1,\ldots,P_{\dim F} \in \Lambda$ such that the evaluation map
\[f:F\otimes_{K}K^{\mathrm{ac}}\longrightarrow\bigoplus_{i=1}^{\dim F} \kappa(P_i) \]
is a bijection. 
For $\chi \in \Omega_{K^{\mathrm{an}}}$, let $P_{i, \chi}$ be the unique extension of $P_i \in X_{K^{\mathrm{an}}}$ to $(X_{K^{\mathrm{an}}})_{\chi}^{\mathrm{an}}$.
Let $\|\ndot\|_{P_{i, \chi}}$ be a norm of $\kappa(P_i)_{\chi}$ given by $\| 1 \|_{P_{i,\chi}} = \exp(-g_{\pi(\chi)}(\mu(P_{i, \chi})))$,
where $\pi$ is the canonical map $\Omega_{K^{\mathrm{an}}} \to \Omega$ and
$\mu : (X_{K^{\mathrm{an}}})_{\chi}^{\mathrm{an}} \to X_{\pi(\chi)}^{\mathrm{an}}$ is also the canonical morphism as analytic spaces.
We set $\xi_i := \{ \|\ndot\|_{P_{i,\chi}} \}_{\chi \in \Omega_{K^{\mathrm{an}}}}$.
We equip $\bigoplus_{i=1}^{\dim F}\kappa(P_i)$ with the $\psi_0$-direct sum $\xi = \{ \|\ndot\|_{\chi} \}_{\chi \in \Omega_{K^{\mathrm{an}}}}$ of $\xi_1, \ldots, \xi_{\dim F}$, where $\psi_0$ denotes the function from $[0,1]$ to $[0,1]$ sending $x\in[0,1]$ to $\max(x,1-x)$ (see Subsections \ref{Subsec:directsums} and \ref{Subsec:Norm families}). Note that, if we denote by $\{e_i\}_{i=1}^{\dim F}$ a basis of $\bigoplus_{i=1}^{\dim F}\kappa(P_i)$ such that $e_i\in \kappa(P_i)$, then this basis is orthogonal with respect to $\|\ndot\|_\chi$ for any $\chi\in\Omega_{K^{\mathrm{an}}}$. By Proposition \ref{Pro:orthogonalesthadamard}, this basis is also an Hadamard basis with respect to $\|\ndot\|_\chi$ for any $\chi\in\Omega_{K^{\mathrm{an}}}$. In particular, one has
\[\widehat{\deg}\Big(\bigoplus_{i=1}^{\dim F}(\kappa(P_i) ,\xi_i)\Big)=\sum_{i=1}^{\dim F} h_{(D, g)}(P_i)\leqslant (\dim F) t.\]
Moreover, for any $\chi\in \Omega_{K^{\mathrm{ac}}}$ the operator norm of $f_{\chi}$ is $\leqslant 1$. Therefore, by Proposition \ref{Pro:slopeinequality1}, one has 
\[
\widehat{\mu}(\overline F) = \widehat{\mu}(\overline F \otimes K^{\mathrm{an}}) \leqslant \frac{1}{\dim F} \widehat{\deg}\Big(\bigoplus_{i=1}^{\dim F}(\kappa(P_i) ,\xi_i)\Big) \leqslant t.
\]
Since $F$ is arbitrary, we obtain $\widehat{\mu}_{\max}(H^0_{\mathbb K}(X, D), \xi_g)\leqslant t$. 
Therefore \eqref{Equ: bound of mu max asy by mu ess} follows because 
$t$ is an arbitrary real number with $t>\widehat{\mu}_{\mathrm{ess}}(L,\varphi)$.
\end{proof}

\subsection{Volume function}
We assume that there exists $C\geqslant 0$ such that the adelic curve $S$ verifies the tensorial minimal slope superadditivity.
\begin{defi}
Let $(L,\varphi)$ be an adelic line bundle on $X$. We define the \emph{arithmetic volume}\index{arithmetic volume} of $(L,\varphi)$ as
\[\widehat{\vol}(L,\varphi):=\limsup_{n\rightarrow+\infty}\frac{\widehat{\deg}_+(\pi_*(L^{\otimes n},n\varphi))}{n^{d+1}/(d+1)!}.\]
We say that $(L,\varphi)$ is \emph{big}\index{big} if $\widehat{\vol}(L,\varphi)>0$.
\end{defi}

{
Assume that the $K$-scheme $X$ admit a regular rational point of $X$. Then the local ring $\mathcal O_{X,P}$ is a regular local ring. By Cohen's structure theorem of complete regular local rings \cite[Proposition 10.16]{Eise}, the formal completion of $\mathcal O_{X,P}$ is isomorphic to the algebra of formal series $K\lbr T_1,\ldots,T_d\rbr$, where $d$ is the Krull dimension of $X$. If $L$ is an invertible $\mathcal O_X$-module, by choosing a local generator of the $\mathcal O_{X,P}$-module $L_P$, we can identify the graded linear series $\bigoplus_{n\in\mathbb N}H^0(X,L^{\otimes n})$ as a graded sub-$K$-algebra (of subfinite type) of $K\lbr T_1,\ldots,T_d\rbr[T]$. We denote by $\Delta(L)$ the Newton-Okounknov body of this graded algebra (see \S\ref{Subsec: applications to the study of graded algebras} for the construction of $\Delta(L)$). For any $n\in\mathbb N$, let $r_n:=\rang_K(H^0(X,L^{\otimes n}))>0$. By Proposition \ref{Pro:converges to volume of Okounkov body} one has
\[\int_{\Delta(L)}1\,\mathrm{d}x=\lim_{r_n>0,\,n\rightarrow+\infty}\frac{r_n}{n^{\kappa}},\]
where $\kappa$ is the Kodaira-Iitaka dimension of the graded linear series $\bigoplus_{n\in\mathbb N}H^0(X,L^{\otimes n})$ (which is also called the \emph{Kodaira-Iitaka dimension}\index{Kodaira-Iitaka dimension} of $L$). In particular, if $L$ is a big line bundle, namely 
\[\mathrm{vol}(L):=\limsup_{n\rightarrow+\infty}\frac{r_n}{n^d/d!}>0,\]
or equivalently, $\kappa=d$, one has
\[\mathrm{vol}(L)=d!\int_{\Delta(L)}1\,\mathrm{d}x.\]

If $(L,\varphi)$ is an adelic line bundle of $X$, then the family $\{(\pi_*(L^{\otimes n},n\varphi))\}_{n\in\mathbb N}$ forms a graded $K$-algebra of adelic vector bundles with respect to $\mathrm{Frac}(K\lbr T_1,\ldots,T_d\rbr)$ (see Definition \ref{Def: graded algebra of adelic vector bundles}). For any $n\in\mathbb N_{\geqslant 1}$ such that $r_n:=\rang_K(H^0(X,L^{\otimes n}))>0$, let $\mathbb P_{(L,\varphi),n}$ be the Borel probability measure on $\mathbb R$ such that, for any positive Borel function $f$ on $\mathbb R$, one has
\[\int_{\mathbb R}f(t)\,\mathbb P_{(L,\varphi),n}(\mathrm{d}t)=\frac{1}{r_n}\sum_{i=1}^{r_n}f(\textstyle{\frac 1n}\widehat{\mu}_i(\pi_*(L^{\otimes n},n\varphi))).\] }

\begin{theo}\label{Thm: limit theorem for adelic line bundles}
Assume that the scheme $X$ admits a regular rational point.  Let $(L,\varphi)$ be an adelic line bundle on $X$. For any $n\in\mathbb N$, let $r_n=\rang_K(H^0(X,L^{\otimes n}))$. Assume that there exists $n\in\mathbb N_{\geqslant 1}$ such that $r_n>0$. 
Then the sequence of measures $\{\mathbb P_{(L,\varphi),n}\}_{n\in\mathbb N,\,r_n>0}$ converges weakly to a Borel probability measure $\mathbb P_{(L,\varphi)}$, which is the direct image of a concave real-valued function $G_{(L,\varphi)}$ on $\Delta(L)^\circ$.  In particular, if $(L,\varphi)$ is big, then the invertible $\mathcal O_X$-module is big. Moreover, in the case where $L$ is big, the sequence 
\begin{equation}\label{Equ: sequence defining volume hat}\frac{\widehat{\deg}_+(\pi_*(L^{\otimes n},n\varphi))}{n^{d+1}/(d+1)!},\quad n\in\mathbb N,\;r_n>0\end{equation}
converges to $\widehat{\mathrm{vol}}(L,\varphi)$, which is also equal to
\begin{equation}\label{Equ: formula for arithmetic volume}(d+1)\mathrm{vol}(L)\int_{\intervalle{[}{0}{+\infty}{[}}t\,\mathbb P_{(L,\varphi)}(\mathrm{d}t)=(d+1)\int_{\Delta(L)^\circ}\max(G_{(L,\varphi)}(x),0)\,\mathrm{d}x.\end{equation}
\end{theo}
{
\begin{proof}
We deduce from Corollary \ref{Cor: convergence theorem for graded adelic vector bundles} that the sequence $\{\mathbb P_{(L,\varphi),n}\}_{n\in\mathbb N,\,r_n>0}$ converges vaguely to a Borel measure $\mathbb P_{(L,\varphi)}$ on $\mathbb R$, which is the direct image of the uniform probability measure on $\Delta(L)^\circ$ by a concave function $G_{(L,\varphi)}:\Delta(L)^\circ\rightarrow\mathbb R\cup\{+\infty\}$. Moreover, by Proposition \ref{Pro: comparison ess minimum and asy max slope} and \ref{Pro: upper bound of mu ess}, we obtain that the supports of the Borel probability measures $\mathbb P_{(L,\varphi),n}$ are uniformly bounded from above. The function $G_{(L,\varphi)}$ is then bounded from above and hence limit measure $\mathbb P_{(L,\varphi)}$ is a Borel probability measure and the sequence $\{\mathbb P_n\}_{n\in\mathbb N,\,r_n>0}$ converges weakly to $\mathbb P_{(L,\varphi)}$. In particular, the sequence $\{\frac{1}{nr_n}\widehat{\deg}_+(L^{\otimes n},n\varphi)\}_{n\in\mathbb N,\,r_n>0}$ converges to \[\int_{\intervalle{[}{0}{+\infty}{[}}t\,\mathbb P_{(L,\varphi)}(\mathrm{d}t)=\frac{1}{\mathrm{vol}(\Delta(L))}\int_{\Delta(L)^\circ}\max(G_{(L,\varphi)}(x),0)\,\mathrm{d}x\] since
\[\int_{\intervalle{[}{0}{+\infty}{[}}t\,\mathbb P_n(\mathrm{d}t)=\int_{\intervalle{[}{0}{\widehat{\mu}_{\max}^{\mathrm{asy}}(L,\varphi)}{[}}t\,\mathbb P_n(\mathrm{d}t)=\frac{1}{nr_n}\sum_{i=1}^{r_n}\max\{\widehat{\mu}_i(\pi_*(L^{\otimes n},n\varphi)),0\},\]
and by \eqref{Equ: bound of deg plus E adelic vector bundle bis} and \eqref{Equ: positive degree bounded by sum of positive slopes},
\[\bigg|\widehat{\deg}_+(\pi_*(L^{\otimes n},n\varphi))-\sum_{i=1}^{r_n}\max\{\widehat{\mu}_i(\pi_*(L^{\otimes n},n\varphi)),0\}\bigg|\leqslant \frac 12\ln(r_n)\nu(\Omega_\infty).\]
In particular, in the case where $\widehat{\vol}(L,\varphi)>0$, one has 
\[\limsup_{n\rightarrow+\infty}\frac{r_n}{n^{d}}>0,\]
namely the invertible $\mathcal O_X$-module $L$ is big. Moreover, in the case where $L$ is big, one has
\[\lim_{n\rightarrow+\infty}\frac{r_n}{n^d/d!}=\vol(L)>0.\]
Therefore the sequence \eqref{Equ: sequence defining volume hat} converges (to $\widehat{\vol}(L,\varphi)$ by definition), which is equal to \eqref{Equ: formula for arithmetic volume}.
\end{proof}

}

\begin{theo}\label{Thm: Brunn-Minkowski for adelic line bundles}
Assume that the scheme $X$ admits a regular rational point. Let $(L_1,\varphi_1)$ and $(L_2,\varphi_2)$ be big adelic line bundles on $X$. Then the following inequality of Brunn-Minkowski type holds
\begin{equation}\label{Equ: Brunn Minkowski L varphi}
\widehat{\mathrm{vol}}(L_1\otimes L_2,\varphi_1+\varphi_2)^{1/(d+1)}\geqslant\widehat{\mathrm{vol}}(L_1,\varphi_1)^{1/(d+1)}+\widehat{\mathrm{vol}}(L_2,\varphi_2)^{1/(d+1)}.
\end{equation}
\end{theo}
\begin{proof}
For any adelic line bundle $(L,\varphi)$ on $X$ such that $L$ is big, we denote by $\widehat{\Delta}(L,\varphi)$ the closure of the convex set $\{(x,t)\in\Delta(L)^\circ\times\mathbb R\,:\,0\leqslant t\leqslant G_{(L,\varphi)}(x)\}$. Then Theorem \ref{Thm: limit theorem for adelic line bundles} implies that
\[\widehat{\mathrm{vol}}(L,\varphi)=(d+1)\int_{\widehat{\Delta}(L,\varphi)} 1\,\mathrm{d}(x,t)\] 
By Proposition \ref{Pro: super additivity of G}, one has 
\[\widehat{\Delta}(L_1,\varphi_1)+\widehat{\Delta}(L_2,\varphi_2)\leqslant\widehat{\Delta}(L_1\otimes L_2,\varphi_1+\varphi_2).\]
Therefore the relation \eqref{Equ: Brunn Minkowski L varphi} follows from the classic Brunn-Minkowski inequality.
\end{proof}

\subsection{Volume of adelic $\mathbb R$-Cartier divisors}
We assume that $X$ is normal and
geometrically integral and admits a regular rational point $P$. We identify the formal completion of $\mathcal O_{X,P}$ with $K\lbr T_1,\ldots,T_d\rbr$, which allow us to embed the rational function field $K(X)$ into the fraction field $R=\mathrm{Frac}(K\lbr T_1,\ldots,T_d\rbr)$. We also suppose that  there exists $C\geqslant 0$ such that the adelic curve $S$ verifies the tensorial minimal slope superadditivity of level $\geqslant C$. In the following, the symbol $\mathbb K$ denotes $\mathbb Z$, $\mathbb Q$ or $\mathbb R$.

Let $D$ be a $\mathbb K$-Cartier divisor. We identify $\bigoplus_{n\in\mathbb N}H^0(nD)$ with a graded sub-$K$-algebra of subfinite type of $R[T]$. We denote by $\Delta(D)$ the Newton-Okounkov body of this graded algebra (see \S\ref{Subsec: applications to the study of graded algebras} for its construction).

\begin{defi}
Let $(D,g)$ be an adelic $\mathbb K$-Cartier divisor. 
We define the \emph{arithmetic volume}\index{arithmetic volume} of $(D,g)$ as
\[\widehat{\vol}(D,g):=\limsup_{n\rightarrow+\infty}\frac{\widehat{\deg}_+(H^0_{\mathbb K}(nD),\xi_{ng})}{n^{d+1}/(d+1)!}\]
(for the definition of the norm family $\xi_{ng}$, see Definition~\ref{def:norm:family:adelic:R:div}).
We say that $(D,g)$ is \emph{big}\index{big} if $\widehat{\mathrm{vol}}(D,g)>0$. Note that for any $s\in K(X)^{\times}$ one has  (see Remark \ref{Rem: linear bijection isometry})\begin{equation}\label{Equ: invariance of vol by principal divisor}\widehat{\vol}((D,g)+\widehat{\mathrm{div}}(s))=\widehat{\vol}(D,g).\end{equation}
Moreover, $(D,g)$ is said to be \emph{arithmetically $\mathbb K$-effective}\index{arithmetically effective}, which is denoted by $(D, g) \geqslant_{\mathbb K} (0,0)$, if $D$ is $\mathbb K$-effective and $g_{\omega} \geqslant 0$ for all
$\omega \in \Omega$. For adelic $\mathbb K$-Cartier divisors $(D_1, g_1)$ and $(D_2, g_2)$ on $X$,
\[
(D_1, g_1) \geqslant_{\mathbb K} (D_2, g_2) \quad\overset{\mathrm{def}}{\Longleftrightarrow}\quad (D_1, g_1) - (D_2, g_1) \geqslant_{\mathbb K} (0,0).
\]
Note that if $(D_1, g_1) \geqslant_{\mathbb K} (D_2, g_2)$, then $\widehat{\deg}_+(H^0_{\mathbb K}(D_1), \xi_{g_1}) \geqslant
\widehat{\deg}_+(H^0_{\mathbb K}(D_2), \xi_{g_2})$. In particular, $\widehat{\mathrm{vol}}(D_1, g_1) \geqslant \widehat{\mathrm{vol}}(D_2, g_2)$.
By using Prposition~\ref{prop:adelic:R:Cartier:div:integrable}, if $D \geqslant_{\mathbb K} 0$, then there is a family of $D$-Green functions $g$ over $S$ such that $(D, g) \geqslant_{\mathbb K} (0,0)$.
\end{defi}

Let $(D,g)$ be an adelic $\mathbb K$-Cartier divisor. The family $\{(H^0_{\mathbb K}(nD),\xi_{ng})\}_{n\in\mathbb N}$ forms a graded $K$-algebra of adelic vector bundles with respect to $R=\mathrm{Frac}(K\lbr T_1,\ldots,T_d\rbr)$. For any $n\in\mathbb N_{\geqslant 1}$ such that $r_n:=\rang_K(H^0(nD))>0$, we let $\mathbb P_{(D,g),n}$ be the Borel probability measure on $\mathbb R$ such that, for any positive Borel function on $\mathbb R$, one has
\[\int_{\mathbb R}f(t)\,\mathbb P_{(D,g),n}(\mathrm{d}t)=\frac{1}{r_n}\sum_{i=1}^{r_n}f(\textstyle{\frac 1n\widehat{\mu}_i(H^0_{\mathbb K}(nD),\xi_{ng})}).\]

\begin{theo}\label{Thm: volume of arithmetic divisor}
Let $(D,g)$ be an adelic $\mathbb K$-Cartier divisor. For any $n\in\mathbb N$, let $r_n=\rang_K(H^0_{\mathbb K}(nD))$. Assume that there exists $n\in\mathbb N_{n\geqslant 1}$ such that $r_n>0$. Then the sequence of measures $\{\mathbb P_{(D,g),n}\}_{n\in\mathbb N,r_n>0}$ converges weakly to a Borel probability measure $\mathbb P_{(D,g)}$, which is the direct image of a concave real-valued function $G_{(D,g)}$ on $\Delta(D)^\circ$. In particular, if $(D,g)$ is big, then the $\mathbb K$-Cartier divisor $D$ is big. Moreover, in the case where $D$ is big, the sequence 
\[\frac{\widehat{\deg}_+(H^0_{\mathbb K}(nD),\xi_{ng})}{n^{d+1}/(d+1)!},\quad n\in\mathbb N,\;r_n>0\]
converges to $\widehat{\mathrm{vol}}(D,g)$, which is also equal to
\begin{equation}\label{Equ:volume arithmetique}
(d+1)\mathrm{vol}(L)\int_{[0,+\infty[}t\,\mathbb P_{(D,g)}(\mathrm{d}t)=(d+1)\int_{\Delta(D)^\circ}\max(G_{(D,g)}(x),0)\,\mathrm{d}x.
\end{equation}
\end{theo}
\begin{proof}
We omit the proof since it is quite similar to that of Theorem \ref{Thm: limit theorem for adelic line bundles}. 
\end{proof}

\begin{coro}\label{coro:Hilbert:Samuel:with:real:case}
Let $(D, g)$ and $(A, h)$ be adelic $\mathbb R$-Cartier divisors on $X$. We assume that $D$ is big.
Then
\[
\lim_{t\to\infty} \frac{\widehat{\deg}_+(H^0_{\mathbb R}(X, tD + A), \xi_{tg + h})}{t^{d+1}/(d+1)!} = \widehat{\mathrm{vol}}(D, g),
\]
where $t$ is a positive real number.
\end{coro}

\begin{proof}
Let us begin with the following claim:

\begin{enonce}{Claim}\label{claim:coro:Hilbert:Samuel:with:real:case:01}
$\widehat{\mathrm{vol}}(aD, ag) = a^{d+1}\widehat{\mathrm{vol}}(D, g)$ for a positive integer $a$.
\end{enonce}

\begin{proof}
By Theorem~\ref{Thm: volume of arithmetic divisor}, 
\begin{align*}
\widehat{\mathrm{vol}}(aD, ag) & = \lim_{n\to\infty} \frac{\widehat{\deg}_+(H^0_{\mathbb K}(naD),\xi_{nag})}{n^{d+1}/(d+1)!} \\
& = a^{d+1}\lim_{n\to\infty} \frac{\widehat{\deg}_+(H^0_{\mathbb K}(naD),\xi_{nag})}{(na)^{d+1}/(d+1)!} = a^{d+1} \widehat{\mathrm{vol}}(D, g),
\end{align*}
as required.
\end{proof}

\begin{enonce}{Claim}\label{claim:coro:Hilbert:Samuel:with:real:case:02}
If $D$ is $\mathbb R$-effective, then the assertion of the corollary holds.
\end{enonce}

\begin{proof}
Choose positive integers $n_0$ and $n_1$ such that $H^0_{\mathbb R}(X, n_0D + A) \not= \{ 0 \}$ and
$H^0_{\mathbb R}(X, n_1D - A) \not= \{ 0 \}$, so that one can take $\phi \in H^0_{\mathbb R}(X, n_0D + A) \setminus \{ 0 \}$ and
$\psi \in H^0_{\mathbb R}(X, n_1D - A) \setminus \{ 0 \}$. 
Let us consider the following injective homomorphisms
\[
\alpha_t : H^0_{\mathbb R}((\lfloor t\rfloor-n_0)D) \to H^0_{\mathbb R}(tD + A)
\quad\text{and}\quad
\beta_t : H^0_{\mathbb R}(tD + A) \to H^0_{\mathbb R}((\lceil t \rceil + n_1)D)
\]
given by $f \mapsto f\phi$ and $f \mapsto f\psi$, respectively. Note that
\[
\Vert \alpha_t(f) \Vert_{tg_{\omega} + h_{\omega}} \leqslant \| f \|_{(\lfloor t\rfloor-n_0)g_{\omega}} \|\phi\|_{(t-\lfloor t\rfloor + n_0)g_{\omega} + h_{\omega}} \leqslant
\| f \|_{(\lfloor t\rfloor-n_0)g_{\omega}} \|\phi\|_{n_0g_{\omega} + h_{\omega}} \| 1 \|_{g_{\omega}}^{t - \lfloor t\rfloor},
\]
so that,  
by Proposition~\ref{prop:deg:plus:comp}, \ref{Item: comparaison degre plus sous fibre} and \ref{Item: degre plus tordu},
\begin{multline}\label{eqn:coro:Hilbert:Samuel:with:real:case:01}
\widehat{\deg}_+(H^0_{\mathbb R}((\lfloor t\rfloor-n_0)D), \xi_{([t] - n_0)g}) \leqslant \widehat{\deg}_+(H^0_{\mathbb R}(tD + A), \xi_{tg+h}) \\
+ (\dim_K H^0_{\mathbb R}((\lfloor t\rfloor-n_0)D)) \int_{\Omega} \left( \big|\ln \|\phi\|_{n_0g_{\omega} + h_{\omega}} \big| +  \big | \ln \| 1 \|_{g_{\omega}} \big| \right) \nu(d\omega).
\end{multline}
In the same way, one has
\begin{multline}\label{eqn:coro:Hilbert:Samuel:with:real:case:02}
\widehat{\deg}_+(H^0_{\mathbb R}(tD + A), \xi_{tg+h}) \leqslant \widehat{\deg}_+(H^0_{\mathbb R}((\lceil t \rceil+n_1)D), \xi_{(\lceil t \rceil+n_1)g}) \\
+ (\dim_K H^0_{\mathbb R}(tD + A)) \int_{\Omega} \left( \big|\ln \|\psi\|_{n_1g_{\omega} - h_{\omega}} \big| +  \big | \ln \| 1 \|_{g_{\omega}} \big| \right) \nu(d\omega).
\end{multline}
Note that 
\[
\lim_{t \to \infty}\frac{\widehat{\deg}_+(H^0_{\mathbb R}((\lfloor t\rfloor -n_0)D), \xi_{(\lfloor t\rfloor - n_0)g})}{t^{d+1}/(d+1)!} = \widehat{\mathrm{vol}}(D, g),\ 
\lim_{t \to \infty}\frac{\dim_K H^0_{\mathbb R}((\lfloor t\rfloor -n_0)D)}{t^{d+1}/(d+1)!} = 0,
\]
so that, by \eqref{eqn:coro:Hilbert:Samuel:with:real:case:01}, one has
\[
\widehat{\mathrm{vol}}(D,g)  \leqslant \liminf_{n\to\infty} \frac{\widehat{\deg}_+(H^0_{\mathbb R} (tD + A), \xi_{tg+h})}{t^{d}/(d+1)!}.
\]
Similarly, by using \eqref{eqn:coro:Hilbert:Samuel:with:real:case:02},
\[
  \limsup_{n\to\infty} \frac{\widehat{\deg}_+(H^0_{\mathbb R} (tD + A), \xi_{tg+h})}{t^{d}/(d+1)!} \leqslant \widehat{\mathrm{vol}}(D,g).
\]
Thus the assertion of the claim follows.
\end{proof}

\begin{enonce}{Claim}\label{claim:coro:Hilbert:Samuel:with:real:case:03}
If there is $\phi \in K(X)^{\times}$ such that $D' := D + (\phi)$ is $\mathbb R$-effective,
then the assertion of the corollary holds.
\end{enonce}

\begin{proof}
We set $(D', g') = (D, g) + \widehat{(\phi)}$.
We choose an arithmetically $\mathbb R$-effective adelic Cartier divisor $(B, k)$ on $X$ such that $(B, k) \pm \widehat{(\phi)}$ are arithmetically $\mathbb R$-effective.
Then, as $(B, k) \pm (t- \lfloor t\rfloor)\widehat{(\phi)}$ are arithmetically $\mathbb R$-effective, one has
\[
t(D', g') + (A, h) - (B, k) - \lfloor t\rfloor \widehat{(\phi)} \leqslant t(D, g) + (A, h) \leqslant t(D', g') + (A, h) + (B, k) + \lfloor t\rfloor \widehat{(\phi)}.
\]
Thus
\begin{multline*}
\widehat{\deg}_+(H^0_{\mathbb R}(tD'+A - B), \xi_{tg'+h-k}) \leqslant
\widehat{\deg}_+(H^0_{\mathbb R}(tD+A), \xi_{tg+h}) \\
\leqslant
\widehat{\deg}_+(H^0_{\mathbb R}(tD'+A + B), \xi_{tg'+h+k}),
\end{multline*}
so that, by using Claim~\ref{claim:coro:Hilbert:Samuel:with:real:case:02},
\begin{multline*}
\widehat{\mathrm{vol}}(D', g') = \lim_{t\to\infty} \frac{\widehat{\deg}_+(H^0_{\mathbb R}(tD'+A - B), \xi_{tg'+h-k})}{t^{d+1}/(d+1)!} \\
\leqslant
\liminf_{t\to\infty}\frac{\widehat{\deg}_+(H^0_{\mathbb R}(tD+A), \xi_{tg+h})}{t^{d+1}/(d+1)!} \leqslant
\limsup_{t\to\infty}\frac{\widehat{\deg}_+(H^0_{\mathbb R}(tD+A), \xi_{tg+h})}{t^{d+1}/(d+1)!} \\
\leqslant
\lim_{t\to\infty} \frac{\widehat{\deg}_+(H^0_{\mathbb R}(tD'+A + B), \xi_{tg'+h+k})}{t^{d+1}/(d+1)!} = \widehat{\mathrm{vol}}(D', g').
\end{multline*}
Therefore one has the claim because $\widehat{\mathrm{vol}}(D, g) = \widehat{\mathrm{vol}}(D', g')$
\end{proof}

\medskip
In general, there are a positive integer $a$ and $f \in K(X)^{\times}$ such that $a D + (f)$ is $\mathbb R$-effective, so that, by using 
Claim~\ref{claim:coro:Hilbert:Samuel:with:real:case:01} and Claim~\ref{claim:coro:Hilbert:Samuel:with:real:case:03},
one has
\begin{align*}
\widehat{\mathrm{vol}}(D, g) & = \frac{1}{a^{d+1}} \widehat{\mathrm{vol}}(aD, ag) = \frac{1}{a^{d+1}}
\lim_{t\to\infty}\frac{\widehat{\deg}_+(H^0_{\mathbb R}(taD+A), \xi_{tag+h})}{t^{d+1}/(d+1)!} \\
& = \lim_{t\to\infty}\frac{\widehat{\deg}_+(H^0_{\mathbb R}(taD+A), \xi_{tag+h})}{(ta)^{d+1}/(d+1)!} 
= \lim_{t\to\infty}\frac{\widehat{\deg}_+(H^0_{\mathbb R}(tD+A), \xi_{tg+h})}{t^{d+1}/(d+1)!},
\end{align*}
as required.
\end{proof}

\begin{coro}\label{coro:homogeneous:volume:adelic:div}
Let $(D, g)$ be an adelic $\mathbb R$-Cartier divisor on $X$. Then, for any $a \in \mathbb R_{\geqslant 0}$,
$\widehat{\mathrm{vol}}(aD, ag) = a^{d+1}\widehat{\mathrm{vol}}(D, g)$.
\end{coro}

\begin{proof}
Clearly we may assume that $a > 0$. If $D$ is not big, then $aD$ is also not big, so that
$\widehat{\mathrm{vol}}(D, g) = 0$ and $\widehat{\mathrm{vol}}(aD, ag) = 0$.
Thus the assertion follows in this case.
If $D$ is big, then by Corollary~\ref{coro:Hilbert:Samuel:with:real:case},
\begin{align*}
\widehat{\mathrm{vol}}(aD, ag) & = \lim_{t\to\infty} \frac{\widehat{\deg}_+(H^0_{\mathbb K}(taD),\xi_{tag})}{t^{d+1}/(d+1)!} \\
& = a^{d+1}\lim_{t\to\infty} \frac{\widehat{\deg}_+(H^0_{\mathbb K}(taD),\xi_{tag})}{(ta)^{d+1}/(d+1)!} = a^{d+1} \widehat{\mathrm{vol}}(D, g),
\end{align*}
as required. 
\end{proof}

\begin{rema} Let $(D_1,g_1)$ and $(D_2,g_2)$ be adelic $\mathbb K$-Cartier divisors. 
Proposition \ref{Pro: super additivity of G} shows that, if $(x,y)\in\Delta(D_1)^\circ\times\Delta(D_2)^\circ$, one has 
\[G_{(D_1+D_2,g_1+g_2)}(x+y)\geqslant G_{(D_1,g_1)}(x)+G_{(D_2,g_2)}(y).\]
\end{rema}

Similar to Theorem \ref{Thm: Brunn-Minkowski for adelic line bundles}, an analogue of Brunn-Minkowski inequality holds for adelic $\mathbb K$-Cartier divisors.

\begin{theo}\label{Thm: Brunn-Minkowski for adelic divisor}
Let $(D_1,g_1)$ and $(D_2,g_2)$ be big adelic $\mathbb K$-Cartier divisors on $X$. 
Then the following inequality holds
\begin{equation}\label{Equ: Brunn Minkowski divisors}
\widehat{\mathrm{vol}}(D_1+D_2,g_1+g_2)^{1/(d+1)}\geqslant\widehat{\mathrm{vol}}(D_1,g_1)^{1/(d+1)}+\widehat{\mathrm{vol}}(D_2,g_2)^{1/(d+1)}.
\end{equation} 
\end{theo}
\begin{proof}
The proof of \eqref{Equ: Brunn Minkowski divisors} is similar to that of \eqref{Equ: Brunn Minkowski L varphi}, which relies on the inequality
\[\forall\,(x,y)\in\Delta(D_1)\times\Delta(D_2),\quad G_{(D_1+D_2,g_1+g_2)}(x+y)\geqslant G_{(D_1,g_1)}(x)+G_{(D_2,g_2)}(y).\]
\end{proof}

Let us consider a criterion for the bigness of adelic $\mathbb K$-Cartier divisors.

\begin{lemm}\label{Lem: maximal value of G}
Let $(D,g)$ be an adelic $\mathbb K$-Cartier divisor on $X$ such that $D$ is big. Then $ \sup_{x\in\Delta(D)^\circ}G_{(D,g)}(x)$ is equal to $\widehat{\mu}_{\max}^{\mathrm{asy}}(D,g)$.
\end{lemm}
\begin{proof}
By \eqref{Equ: maxmal value of concave transfor is asymptotic first slope}, we obtain that the maximal value of $G_{(D,g)}$ is equal to $\lim_{n\rightarrow+\infty}\frac 1n\widehat{\mu}_1(H^0_{\mathbb K}(nD),\xi_{ng})$. Note that all norm families $\xi_{ng}$ are ultrametric on $\Omega\setminus\Omega_\infty$. 

By Remark \ref{Rem: comparison between mu max and mu 1} and the relation \[\lim_{n\rightarrow+\infty}\frac 1n\ln\rang_K(H^0_{\mathbb K}(nD))=0,\] we obtain the equality $\sup\limits_{x \in \Delta(D)^\circ} G_{(D,g)}(x) = \widehat{\mu}_{\mathrm{\max}}^{\mathrm{asy}}(D, g)$.
\end{proof}

\begin{prop}\label{prop:big:small:sec}
Let $(D, g)$ be an adelic $\mathbb K$-Cartier divisor on $X$.
Then the following are equivalent:
\begin{enumerate}[label=\rm(\arabic*)]
\item
$(D, g)$ is big.

\item
$D$ is big and $\widehat{\mu}_{\mathrm{\max}}^{\mathrm{asy}}(D, g) > 0$.
\end{enumerate}
\end{prop}

\begin{proof}
First of all,  
note that $\sup\limits_{x \in \Delta(D)^\circ} G_{(D,g)}(x) = \widehat{\mu}_{\mathrm{\max}}^{\mathrm{asy}}(D, g)$ by Lemma \ref{Lem: maximal value of G}.
Moreover, by Theorem~\ref{Thm: volume of arithmetic divisor},
\begin{equation}\label{eqn:prop:big:small:sec:01}
\widehat{\mathrm{vol}}(D, g) = (d+1)\int_{\Delta(D)^\circ}\max(G_{(D,g)}(x),0)\,\mathrm{d}x,
\end{equation}

\medskip
(1) $\Longrightarrow$ (2): By the above facts, one has
\begin{equation}\label{eqn:prop:big:small:sec:02}
\widehat{\mathrm{vol}}(D, g) \leqslant (d+1) \mathrm{vol}(D) \max (\widehat{\mu}_{\mathrm{\max}}^{\mathrm{asy}}(D, g), 0).
\end{equation}
Therefore, the assertion follows.

\medskip
(2) $\Longrightarrow$ (1): First of all, as $D$ is big, $\Delta(D)^\circ \not=
\emptyset$. Moreover, since $\widehat{\mu}_{\mathrm{\max}}^{\mathrm{asy}}(D, g) > 0$
and $G_{(D,g)}$ is continuous on $\Delta(D)^\circ$,
one can find a non-empty open set $U$ of $\Delta(D)^\circ$ such that
$G_{(D,g)} > 0$ on $U$, so that the assertion follows from \eqref{eqn:prop:big:small:sec:01}.
\end{proof}

\begin{defi}
An adelic $\mathbb K$-Cartier divisor $(D, g)$ is \emph{strongly big}\index{strongly big@strongly big} if
$D$ is big and $\nu_1^{\mathrm{asy}}(D, g) > 0$, that is, $D$ is big and there are positive integer $a$ and
$\phi \in H^0_{\mathbb K}(aD) \setminus \{ 0 \}$ such that $\widehat{\deg}_{\xi_{ag}}(\phi) > 0$.
Note that strong bigness implies bigness by Proposition~\ref{prop:Minkowski:nu:1:mu:max} and Proposition~\ref{prop:big:small:sec}.
Moreover if $S$ satisfies the Minkowski property of certain level, then
strong bigness is equivalent to bigness by Proposition~\ref{prop:Minkowski:nu:1:mu:max} and Proposition~\ref{prop:big:small:sec}.
\end{defi}

\begin{prop}\label{prop:big:R:Cartier:div:green:function:positive}
Let $(D, g)$ be an adelic $\mathbb K$-Cartier divisor on $X$ such that $D$ is big. Then there is an integrable function
$\varphi$ on $\Omega$ such that $(D, g + \varphi)$ is strongly big. 
\end{prop}

\begin{proof}
Since $D$ is big, there are a positive integer $a$ and $f \in K(X)^{\times}$ such that
$a D + (f)$ is effective. By Proposition~\ref{prop:adelic:R:Cartier:div:integrable},
a function given by $\omega \mapsto \ln \| f \|_{ag_{\omega}}$ is integrable.
Thus if we set
\[
\varphi(\omega) := \begin{cases}
{\displaystyle \frac{1}{a}(\ln \| f \|_{a g_{\omega}} + \ln 2 )} & \text{if $\omega \in \Omega_{\infty}$}, \\[1.5ex]
{\displaystyle \frac{1}{a}\ln \| f \|_{a g_{\omega}}} & \text{if $\omega \in \Omega \setminus \Omega_{\infty}$},
\end{cases}
\]
then $(\omega \in \Omega) \mapsto \varphi(\omega)$ is integrable.
Let $F_n$ be a vector subspace of $H^0_{\mathbb K}(naD)$ generated by $f^n$.
Then
\begin{align*}
\widehat{\deg}_{\xi_{na(g + \varphi)}}(F_n) & = - \int_{\Omega} \ln \| f^n \|_{na(g_{\omega} + \varphi(\omega))} \,\nu(\mathrm{d}\omega)= - n\int_{\Omega} \left( \ln \| f \|_{ag_{\omega}} - a \varphi(\omega) \right)\, \nu(\mathrm{d}\omega) \\
& = n \int_{\Omega_{\infty}} (\ln 2)\, \nu(\mathrm{d}\omega),
\end{align*}
so that 
\[
\nu_1
(H^0(naD), \xi_{na(g+\varphi)}) \geqslant 
\widehat{\deg}_{\xi_{na(g + \varphi)}}(F_n) = n \int_{\Omega_{\infty}} (\ln 2)\, \nu(\mathrm{d}\omega),
\]
which shows that $\nu_1^{\mathrm{asy}}
(D, g+\varphi) > 0$, so that $(D, g + \varphi)$
is strongly big.
\end{proof}

\begin{defi}
Let $(D, g)$ and $(D', g')$ be adelic $\mathbb K$-Cartier divisors on $X$.
We define $(D', g') \precsim (D, g)$ to be
\[
(D',g') \precsim (D, g) \quad\overset{\mathrm{def}}{\Longleftrightarrow}\quad
\text{either $(D', g') = (D, g)$ or $(D,g) - (D', g')$ is big.}
\]
\end{defi}

\begin{prop}\label{prop:partial:order:adelic:div}
\begin{enumerate}[label=\rm(\arabic*)]
\item
The relation $\precsim$ forms a partial order on the group of adelic $\mathbb K$-Cartier divisors on $X$.

\item
For adelic $\mathbb K$-Cartier divisors $(D,g), (D',g'), (E, h)$ and $(E', h')$ on $X$,
if $(D',g') \precsim (D,g)$ and $(E', h') \precsim (E, h)$, then $(D', g') + (E', h') \precsim (D,g) + (E, h)$ and
$a(D', g') \precsim a(D,g)$ for $a \in \mathbb K_{\geqslant 0}$.

\item
For adelic $\mathbb K$-Cartier divisors $(D,g)$ and $(D',g')$ on $X$,
if $(D',g') \precsim (D,g)$, then $\widehat{\mathrm{vol}}(D', g') \leqslant \widehat{\mathrm{vol}}(D,g)$.
\end{enumerate}
\end{prop}

\begin{proof}
(1) We assume that $(D',g') \precsim (D,g)$ and $(D,g) \precsim (D',g')$. If $(D',g') \not= (D,g)$,
then $(D,g) - (D', g')$ and $(D',g') - (D, g)$ are big, so that
\[
(0,0) = \left((D,g) - (D', g')\right) + \left((D',g') - (D, g)\right)
\]
is also big by Theorem~\ref{Thm: Brunn-Minkowski for adelic divisor}, which is a contradiction.
Next let us see that if $(D_1,g_1) \precsim (D_2,g_2)$ and $(D_2,g_2) \precsim (D_3,g_3)$, then
$(D_1,g_1) \precsim (D_3,g_3)$. Indeed, this is a consequence of Theorem~\ref{Thm: Brunn-Minkowski for adelic divisor} because
\[
(D_3, g_3) - (D_1, g_1) = \left( (D_3, g_3) - (D_2, g_2) \right) + \left( (D_2, g_2) - (D_1, g_1) \right).
\]

(2) is the consequences of Theorem~\ref{Thm: Brunn-Minkowski for adelic divisor} and Corollary~\ref{coro:homogeneous:volume:adelic:div} because
\[
\begin{cases}
\left( (D,g) + (E, h) \right) - \left( (D', g')  + (E', h') \right) = \left( (D,g) - (D', g') \right) + \left( (E, h) - (E', h') \right),\\
a (D,g)  - a (D', g') = a \left(  (D,g)  -  (D', g')\right).
\end{cases}
\]

(3) We may assume that $(D, g) - (D', g')$ is big. If $(D', g')$ is big, then the assertion follows from Theorem~\ref{Thm: Brunn-Minkowski for adelic divisor} because
$(D, g) = \left((D, g) - (D', g')\right) + (D', g')$. Otherwise, the assertion is obvious because $\widehat{\mathrm{vol}}(D', g') = 0$.
\end{proof}

\begin{prop}\label{prop:big:open:weak}
Let $(D, g)$ be a big adelic $\mathbb K$-Cartier divisor on $X$ and $(A, h)$ be an adelic $\mathbb R$-Cartier divisor on $X$.
Then there is a positive integer $n_0$ such that $n(D, g) + (A, h)$ is big for all $n \in \mathbb Z_{\geqslant n_0}$.
\end{prop}

\begin{proof}
It is sufficient to find a positive integer $n_0$ such that $n_0(D, g) + (A, h)$ is big because
$n(D, g) + (A, h) = n_0(D, g) + (A, h) + (n-n_0)(D, g)$.

As $D$ is big, one can find a positive integer $m$ such that $m D + A$ is big, so that, by Prioposition~\ref{prop:big:R:Cartier:div:green:function:positive},
$(m D +A, mg + h + \phi)$ is big for some non-negative integrable function $\phi$ on $\Omega$.
Let $a$ be a positive integer such that
\[
\widehat{\mathrm{vol}}(D, g) > \frac{(d+1)\mathrm{vol}(D)}{a} \int_{\Omega} \phi \nu(d\omega).
\]
Since
\[
\widehat{\mathrm{vol}}(D, g - \phi/a) \geqslant \widehat{\mathrm{vol}}(D, g) - \frac{(d+1)\mathrm{vol}(D)}{a} \int_{\Omega} \phi \nu(d\omega) > 0
\]
by using Proposition~\ref{prop:deg:plus:comp}, (2), 
one obtains $(D, g - \phi/a)$ is big.
Thus the assertion follows because
\[
(m+a) (D, g) + (A, h) = (mD + A, mg + h + \phi) + a(D, g - \phi/a).
\]
\end{proof}

\begin{theo}\label{thm:strong:cont:vol}
Let $(D,g), (D_1, g_1), \ldots, (D_n, g_n)$ be adelic $\mathbb R$-Cartier divisors on $X$.
Then
\[
\lim_{\varepsilon_1 \to 0, \ldots, \varepsilon_n \to 0} \widehat{\mathrm{vol}}((D, g) + \varepsilon_1 (D_1, g_1) + \cdots + \varepsilon_n (D_n, g_n)) = \widehat{\mathrm{vol}}(D, g).
\]
\end{theo}

\begin{proof}
Let us begin with the following Claim~\ref{claim:thm:strong:cont:vol:01}, Claim~\ref{claim:thm:strong:cont:vol:02}, Claim~\ref{claim:thm:strong:cont:vol:03} and Claim~\ref{claim:thm:strong:cont:vol:04}:

\begin{enonce}{Claim}\label{claim:thm:strong:cont:vol:01}
Let $(E, h)$ be an adelic $\mathbb R$-Cartier divisor on $X$.
Let $(0, f)$ be an adelic Cartier divisor on $X$.
Then $\lim_{\varepsilon \to 0} \widehat{\mathrm{vol}}(E, h + \varepsilon f) = \widehat{\mathrm{vol}}(E, h)$.
\end{enonce}

\begin{proof}
We set $\varphi_{1}(\omega) = \sup_{x \in X_{\omega}} \{ f_{\omega}(x) \}$ and $\varphi_{2}(\omega) = \sup_{x \in X_{\omega}} \{ -f_{\omega}(x) \}$.
Then, by Proposition~\ref{prop:adelic:R:Cartier:div:integrable}, 
$\varphi_{1}(\omega)$ and $\varphi_{2}(\omega)$ are integrable on $\Omega$, so that $\varphi(\omega) = \max \{ |\varphi_{1}(\omega)|, |\varphi_{2}(\omega)| \}$
is also integrable on $\Omega$ and $|f_{\omega}(x)| \leqslant \varphi(\omega)$ for all $x \in X_{\omega}$ and $\omega \in \Omega$.
Therefore, 
\[
h_{\omega} - |\varepsilon| \varphi(\omega) \leqslant h_{\omega} + \varepsilon f_{\omega} \leqslant h_{\omega} + |\varepsilon| \varphi(\omega),
\]
so that, by Proposition~\ref{prop:deg:plus:comp}, (1),
\[
\widehat{\deg}_+(H^0_{\mathbb K}(nE), \mathrm{e}^{n |\varepsilon| \varphi} \xi_{nh}) \leqslant
\widehat{\deg}_+(H^0_{\mathbb K}(nE), \xi_{n(h + \varepsilon f)}) \leqslant
\widehat{\deg}_+(H^0_{\mathbb K}(nE), \mathrm{e}^{-n |\varepsilon| \varphi} \xi_{nh}).
\]
Moreover, by Proposition~\ref{prop:deg:plus:comp}, (2),
\begin{multline*}
\widehat{\deg}_+(H^0_{\mathbb K}(nE), \xi_{nh}) \leqslant \widehat{\deg}_+(H^0_{\mathbb K}(nE), \mathrm{e}^{-n |\varepsilon| \varphi} \xi_{nh}) \\
\leqslant \widehat{\deg}_+(H^0_{\mathbb K}(nE), \xi_{nh}) + n |\varepsilon| \dim_K (H^0_{\mathbb K}(nE)) \int_{\Omega} \varphi \,\nu(\mathrm{d}\omega),
\end{multline*}
and
\begin{multline*}
\widehat{\deg}_+(H^0_{\mathbb K}(nE), \mathrm{e}^{n |\varepsilon| \varphi} \xi_{nh}) \leqslant \widehat{\deg}_+(H^0_{\mathbb K}(nE), \xi_{nh}) \\
\leqslant \widehat{\deg}_+(H^0_{\mathbb K}(nE), \mathrm{e}^{n |\varepsilon| \varphi} \xi_{nh}) + n |\varepsilon| \dim_K (H^0_{\mathbb K}(nE)) \int_{\Omega} \varphi\, \nu(\mathrm{d}\omega).
\end{multline*}
Therefore the assertion of the claim follows.
\end{proof}

\begin{enonce}{Claim}\label{claim:thm:strong:cont:vol:02}
Let $(B, f)$ be an adelic $\mathbb R$-Cartier divisor on $X$ such that
$(B, f) \pm (D_i, g_i)$ is big for every $i=1, \ldots, n$. Then
\begin{multline*}
\widehat{\mathrm{vol}}((D, g) - (|\varepsilon_1| + \cdots +  |\varepsilon_n|)(B, f)) \\ \leqslant
\widehat{\mathrm{vol}}((D, g) + \varepsilon_1 (D_1, g_1) + \cdots + \varepsilon_n (D_n, g_n)) \\ \leqslant
\widehat{\mathrm{vol}}((D, g) + (|\varepsilon_1| + \cdots +  |\varepsilon_n|)(B, f))
\end{multline*}
\end{enonce}

\begin{proof}
Since
\[
\begin{cases}
|\varepsilon_i|(B, f) - \varepsilon_i (D_i, g_i) = |\varepsilon_i|((B,f) \pm (D_i, g_i)), \\
\varepsilon_i (D_i, g_i) + |\varepsilon_i|(B, f) = |\varepsilon_i|((B, f) \pm (D_i, g_i)),
\end{cases}
\]
one has $-|\varepsilon_i|(B, f) \precsim \varepsilon_i (D_i, g_i) \precsim |\varepsilon_i|(B,f)$ by Proposition~\ref{prop:partial:order:adelic:div}, (2),
so that, by using Proposition~\ref{prop:partial:order:adelic:div}, (2) again,
\begin{multline*}
(D, g) - (|\varepsilon_1| + \cdots +  |\varepsilon_n|)(B, f) \\ \precsim
(D, g) + \varepsilon_1 (D_1, g_1) + \cdots + \varepsilon_n (D_n, g_n) \\ \precsim
(D, g) + (|\varepsilon_1| + \cdots +  |\varepsilon_n|)(B, f).
\end{multline*}
Therefore, by Proposition~\ref{prop:partial:order:adelic:div}, (3), one obtains the claim.
\end{proof}

\begin{enonce}{Claim}\label{claim:thm:strong:cont:vol:03}
Let $(H, g_H)$ be an adelic $\mathbb R$-Cartier divisor on $X$. 
Then there is an integrable function $\psi$ on $S$ such that
$(H, g_H - \psi)$ is not big.
\end{enonce}

\begin{proof}
Proposition~\ref{Pro: upper bound ess min 2} and Proposition~\ref{Pro: comparison ess minimum and asy max slope},
one obtains $\widehat{\mu}_{\max}^{\mathrm{asy}}(H,g_H) < \infty$, so that one can find
an integrable function $\psi$ on $S$ such that
\[
\widehat{\mu}_{\max}^{\mathrm{asy}}(H,g_H) < \int_{\Omega} \psi\, \nu(\mathrm{d}\omega).
\]
We choose a positive integer $n_0$ such that
\[
\widehat{\mu}_{\max}(H^0_{\mathbb R}(X, nH), \xi_{ng_H}) \leqslant \int_{\Omega} n \psi\, \nu(\mathrm{d}\omega)
\]
for all $n \geqslant n_0$. Thus, as $\xi_{n(g_H-\psi)} = \exp(n\psi)\xi_{ng_H}$, by Lemma~\ref{lem:mu:max:deg:plus}, (1),
\[
\widehat{\mu}_{\max}(H^0_{\mathbb R}(X, nH), \xi_{n(g_H-\psi)}) = \widehat{\mu}_{\max}(H^0_{\mathbb R}(X, nH), \xi_{ng_H}) - \int_{\Omega} n \psi \,\nu(\mathrm{d}\omega) \leqslant 0,
\]
so that the assertion follows from Lemma~\ref{lem:mu:max:deg:plus}, (2).

\end{proof}

\begin{enonce}{Claim}\label{claim:thm:strong:cont:vol:04}
Let $(H, g_H)$ be an adelic $\mathbb R$-Cartier divisor on $X$ and $\varphi$ be an integrable function on $\Omega$.
Then
\[
\widehat{\mathrm{vol}}(H, g_H + \varphi) \leqslant \widehat{\mathrm{vol}}(H, g_H) + (d+1)\mathrm{vol}(H) \int_{\Omega} | \varphi(\omega)|\, \nu(\mathrm{d}\omega).
\]
\end{enonce}

\begin{proof}
As $\xi_{n(g_H + \varphi)} = \exp(-n\varphi) \xi_{ng_H}$, by using Proposition~\ref{prop:deg:plus:comp}, (2),
\begin{multline*}
\widehat{\deg}_+(H^0_{\mathbb R}(X, nH), \xi_{n(g_H + \varphi)}) \leqslant \widehat{\deg}_+(H^0_{\mathbb R}(X, nH), \xi_{ng_H}) \\
+ n (\dim_K H^0_{\mathbb R}(nH)) \int_{\Omega} |\varphi(\omega)|\,\nu(\mathrm{d}\omega),
\end{multline*}
so that the assertion follows.
\end{proof}

First we assume that $D$ is big.
By Proposition~\ref{prop:big:R:Cartier:div:green:function:positive}, we can choose a $D$-Green functions family $g'$ such that $(D, g')$ is a big adelic
$\mathbb K$-Cartier divisor.
Then, by Proposition~\ref{prop:big:open:weak}, one can choose a positive integer $a$ such that
$a(D, g') \pm (D_i, g_i)$ is big for every $i=1, \ldots, n$. 
Then, by Claim~\ref{claim:thm:strong:cont:vol:02},
\begin{multline*}
\widehat{\mathrm{vol}}((D, g) - a (|\varepsilon_1| + \cdots +  |\varepsilon_n|)(D, g')) \\ \leqslant
\widehat{\mathrm{vol}}((D, g) + \varepsilon_1 (D_1, g_1) + \cdots + \varepsilon_n (D_n, g_n)) \\ \leqslant
\widehat{\mathrm{vol}}((D, g) + a (|\varepsilon_1| + \cdots +  |\varepsilon_n|)(D, g'))
\end{multline*}
If we set $f = g' - g$ and $\varepsilon = |\varepsilon_1| + \cdots +  |\varepsilon_n|$, then
\[
\begin{cases}
(D, g) - a (|\varepsilon_1| + \cdots +  |\varepsilon_n|)(D, g') = (1 - a \varepsilon) \left((D, g) + \left(0, \frac{a\varepsilon}{1 - a \varepsilon}f\right)\right),\\
(D, g) + a (|\varepsilon_1| + \cdots +  |\varepsilon_n|)(D, g') = (1 + a \varepsilon) \left((D, g) + \left(0, \frac{a\varepsilon}{1 + a \varepsilon}f\right)\right).
\end{cases}
\]
Therefore, by Claim~\ref{claim:thm:strong:cont:vol:01},
\begin{multline*}
\lim_{\varepsilon_1 \to 0, \ldots, \varepsilon_n \to 0} \widehat{\mathrm{vol}}((D, g) - a (|\varepsilon_1| + \cdots +  |\varepsilon_n|)(D, g')) \\
=  \lim_{\varepsilon \to 0}(1 - a \varepsilon)^{d+1} \widehat{\mathrm{vol}}\left((D, g) + \left(0, \frac{a\varepsilon}{1 - a \varepsilon} f\right)\right) = \widehat{\mathrm{vol}}(D, g).
\end{multline*}
In the same way,
\[
\lim_{\varepsilon_1 \to 0, \ldots, \varepsilon_n \to 0} \widehat{\mathrm{vol}}((D, g) + a (|\varepsilon_1| + \cdots +  |\varepsilon_n|)(D, g')) = \widehat{\mathrm{vol}}(D, g).
\]
One has the theorem in the case where $D$ is big.

\medskip

Next we assume that $D$ is not big. 
Let $(A, h)$ be a big adelic Cartier divisor on $X$ such that $D + A$ is big and $(A, h) \pm (D_i, g_i)$ are big for every $i=1, \ldots, n$.
Then, by Claim~\ref{claim:thm:strong:cont:vol:02}, if we set $\varepsilon = |\varepsilon_1| + \cdots + |\varepsilon_n|$, then
\[
0 \leqslant \widehat{\mathrm{vol}}((D, g) + \varepsilon_1 (D_1, g_1) + \cdots + \varepsilon_n (D_n, g_n)) \leqslant
\widehat{\mathrm{vol}}((D, g) + \varepsilon (A,h)),
\]
and hence one need to show that $\lim_{\varepsilon\downarrow 0} \widehat{\mathrm{vol}}((D, g) + \varepsilon (A,h)) = 0$.
By Claim~\ref{claim:thm:strong:cont:vol:03}, one can choose a non-negative integrable function $\varphi$ on $\Omega$
such that $(D,g) + (A, h) - (0, \varphi)$ is not big.
Then, as $(D, g) - (0, \varphi) + \varepsilon (A, h)  \precsim (D,g) + (A, h) - (0, \varphi) + \varepsilon (A, h)$, one has
\[
\widehat{\mathrm{vol}}((D, g) - (0, \varphi) + \varepsilon (A, h)) \leqslant  \widehat{\mathrm{vol}}((D,g) + (A, h) - (0, \varphi) + \varepsilon (A, h)).
\]
Since $D + A$ is big, by the previous case,
\[
\lim_{\varepsilon\downarrow 0} \widehat{\mathrm{vol}}((D,g) + (A, h) - (0, \varphi) + \varepsilon (A, h)) =
\widehat{\mathrm{vol}}((D,g) + (A, h) - (0, \varphi)) = 0,
\]
and hence
\begin{equation}\label{eqn:thm:strong:cont:vol:99}
\lim_{\varepsilon\downarrow 0} \widehat{\mathrm{vol}}((D, g) - (0, \varphi) + \varepsilon (A, h))  = 0.
\end{equation}
On the other hand, by Claim~\ref{claim:thm:strong:cont:vol:04},
\[
\widehat{\mathrm{vol}}((D, g) + \varepsilon (A, h)) \leqslant \widehat{\mathrm{vol}}((D, g) - (0, \varphi) + \varepsilon (A, h))
+ (d+1) \mathrm{vol}(D + \varepsilon A) \int_{\Omega} \varphi(\omega) \nu(d\omega).
\]
As $D$ is not big,  one obtains
\[
\lim_{\varepsilon\downarrow 0} \mathrm{vol}(D + \varepsilon A)  = \mathrm{vol}(D) = 0,
\]
and hence, by \eqref{eqn:thm:strong:cont:vol:99}, one has
$\lim_{\varepsilon\downarrow 0} \widehat{\mathrm{vol}}((D, g) + \varepsilon (A, h)) = 0$, as required.

\end{proof}

\begin{coro}\label{coro:big:open}
Let $H$ be a finite-dimensional vector space of $\widehat{\mathrm{Div}}_{\mathbb R}(X)$.
Then the set $\{ (D, g) \in H \mid \text{$(D, g)$ is big} \}$ is an open cone in $H$.
\end{coro}

\begin{proof}
The openness of it is a consequence of Theorem~\ref{thm:strong:cont:vol}.
One can check that it is a cone by Theorem~\ref{Thm: Brunn-Minkowski for adelic divisor} and Corollary~\ref{coro:homogeneous:volume:adelic:div}.
\end{proof}

\begin{coro}
The volume function $\widehat{\mathrm{vol}} : \widehat{\mathrm{Div}}_{\mathbb R}(X) \to \mathbb R$ factors through $\widehat{\mathrm{Div}}_{\mathbb R}(X)$ modulo
the vector subspace over $\mathbb R$ generated by principal Cartier divisors, that is, 
\[
\widehat{\mathrm{vol}}((D, g) + a_1 \widehat{(f_1)} + \cdots + a_r \widehat{(f_r)})
= \widehat{\mathrm{vol}}(D, g)
\]
for any $r \in \mathbb Z_{\geqslant 1}$, $(D, g) \in \widehat{\mathrm{Div}}_{\mathbb R}(X)$, $f_1, \ldots, f_r \in K(X)^{\times}$ and $a_1, \ldots, a_r \in \mathbb R$.
\end{coro}

\begin{proof}
If $a_1, \ldots, a_r \in \mathbb Z$, then the assertion is obvious. Next we assume that $a_1, \ldots, a_r \in \mathbb Q$.
We choose a positive integer $N$ such that $Na_1, \ldots, Na_r \in \mathbb Z$.
Then
\begin{align*}
N^{d+1} \widehat{\mathrm{vol}}(D, g) & = \widehat{\mathrm{vol}}(ND, Ng) = \widehat{\mathrm{vol}}((ND, Ng) + (Na_1) \widehat{(f_1)} + \cdots + (Na_r) \widehat{(f_r)})\\
& = N^{d+1} \widehat{\mathrm{vol}}((D, g) + a_1 \widehat{(f_1)} + \cdots + a_r \widehat{(f_r)}),
\end{align*}
as required. In general, take sequences $\{a_{1,n} \}_{n=1}^{\infty}, \ldots, \{a_{r,n} \}_{n=1}^{\infty}$ of rational numbers such that
$a_1 = \lim_{n\to\infty} a_{1, n}, \ldots, a_r = \lim_{n\to\infty} a_{r, n}$. Then, by Theorem~\ref{thm:strong:cont:vol},
\begin{align*}
\widehat{\mathrm{vol}}((D, g) + a_1 \widehat{(f_1)} + \cdots + a_r \widehat{(f_r)}) & = \lim_{n\to\infty}
\widehat{\mathrm{vol}}((D, g) + a_{1,n} \widehat{(f_1)} + \cdots + a_{r,n} \widehat{(f_r)}) \\
& = \widehat{\mathrm{vol}}(D, g),
\end{align*}
so that the assertion follows.
\end{proof}

%% file: ch7_2019_03_23.tex

\chapter{Nakai-Moishezon's criterion}

\IfChapVersion
\ChapVersion{Version of Chapter 7 : \\ \StrSubstitute{\DateChapSeven}{_}{\_}}
\fi

In this chapter, we fix an adelic curve $S=(K,(\Omega,\mathcal A,\nu),\phi)$. We assume that, either the $\sigma$-algebra $\mathcal A$ is discrete, or the field $K$ admits a countable subfield which is dense in each $K_\omega$, $\omega\in\Omega$.

\section{Graded algebra of adelic vector bundles}

Let $C$ be a non-negative real number. In this section, we assume that the adelic curve $S$ satisfies the \emph{tensorial minimal slope property}\index{tensorial minimal slope property} of level $\geqslant C$. Namely, for any pair $(\overline E,\overline F)$ of non-zero adelic vector bundles on $S$, one has
\begin{equation}\label{Equ: tensorial minimal slope superadditivity}\widehat{\mu}_{\min}(\overline E\otimes_{\varepsilon,\pi}\overline F)\geqslant\widehat{\mu}_{\min}(\overline E)+\widehat{\mu}_{\min}(\overline F)-C(\ln(\rang_K(E)\cdot\rang_K(F))).\end{equation}
Note that we have shown in Chapter \ref{Chap: Slsemistableopes of tensor product} that, if the field $K$ is perfect, then the adelic curve $S$ satisfies the tensorial minimal slope super-additivity of level $\geqslant \frac 23\nu(\Omega_\infty)$.

\begin{defi}\label{def:graded:algebra:adelic:vector:bundles}
Let $R_\sbullet=\bigoplus_{n\in\mathbb N}R_n$ be a graded $K$-algebra. We assume that,  for any $n\in\mathbb N$, $R_n$ is of finite rank over $K$. For any $n\in\mathbb N$, let $\xi_n=\{\|\ndot\|_{n,\omega}\}_{\omega\in\Omega}$ be a norm family on $R_n$. 
We say that $\overline R_{\sbullet}=\{(R_n,\xi_n)\}_{n\in\mathbb N}$ is a 
\emph{normed graded algebra on $S$}\index{normed graded algebra on S@normed graded algebra on $S$} if, for any $\omega \in \Omega$, 
$\overline R_{\sbullet,\omega} = \{ (R_{n,\omega}, \|\ndot\|_{n,\omega}) \}_{n \in \mathbb N}$ forms
a normed graded algebra over $K_{\omega}$, where $R_{n,\omega} := R_n \otimes_K K_{\omega}$
(cf. Subsection~\ref{Subsec:graded:algebra:norms}).
Moreover, if $(R_n, \xi_n)$ forms an adelic vector bundle on $S$ for all $n \in \mathbb N$, then
$\overline R_{\sbullet}$ is called a \emph{graded algebra of adelic vector bundles on $S$}\index{graded algebra of adelic vector bundles on S@graded algebra of adelic vector bundles on $S$}.
Furthermore, we say that $\overline R_\sbullet$ is \emph{of finite type}\index{of finite type@of finite type} if the underlying graded $K$-algebra $R_\sbullet$ is of finite type over $K$.
\end{defi}

\begin{prop}\label{Pro: sequence of normalised mu min} Let $C_0$ be a non-negative real number. Assume that the adelic curve $S$ satisfies the tensorial minimal slope property of level $\geqslant C_0$. 
Let $\overline R_\sbullet=\{(R_n,\xi_n)\}$ be a graded algebra of adelic vector bundles on $S$ such that $\overline R_0$ is the trivial adelic line bundle, namely $R_0=K$ and for any $\omega\in\Omega$, one has $\|1\|_\omega=1$. Suppose in addition that $R_\sbullet$ is generated as $R_0$-algebra by $R_1$. Then the sequence $\{\widehat{\mu}_{\min}(\overline R_n)/n\}_{n\in\mathbb N}$ converges to an element in $\mathbb R\cup\{+\infty\}$.
\end{prop}
\begin{proof}
Let $(n,m)$ be a couple of positive integers. Since $R_\sbullet$ is generated as $K$-algebra by $R_1$, the canonical $K$-linear map $f_{n,m}:R_n\otimes_KR_m\rightarrow R_{n+m}$ is surjective. Moreover, if we equip $R_n\otimes_KR_m$ with the $\varepsilon,\pi$-tensor product $\xi_n\otimes_{\varepsilon,\pi}\xi_m$, then, by the submultiplicativity condition, the homomorphism $f_{n,m}$ has height $\leqslant 1$. By Proposition \ref{Pro:inegalitdepente}, one has
\[\widehat{\mu}_{\min}(\overline R_n\otimes_{\varepsilon,\pi}\overline R_m)\leqslant\widehat{\mu}_{\min}(\overline R_{n+m}).\]
Moreover, by the assumption of tensorial minimal slope property, 
\[\widehat{\mu}_{\min}(\overline R_n\otimes_{\varepsilon,\pi}\overline R_m)\geqslant \widehat{\mu}_{\min}(\overline R_n)+\widehat{\mu}_{\min}(\overline R_m)-C_0(\ln(\rang_K(R_n))+\ln(\rang_K(R_m))).\]
Note that $R_\sbullet$ is a quotient $K$-algebra of $K[R_1]$. Hence $\rang_K(R_n)=O(n^{\rang_K(R_1)-1})$. By \cite[Proposition 1.3.5]{Chen10b}\footnote{In the statement of \cite[Proposition 1.3.5]{Chen10b}, we suppose given a \emph{positive} sequence $\{b_n\}_{n\in\mathbb N,\,n\geqslant 1}$ satisfying the \emph{weak subadditivity}\index{weak subadditivity@weak subadditivity} condition $b_{n+m}\leqslant b_n+b_m+f(n)+f(m)$, where $f:\mathbb N_{\geqslant 1}\rightarrow\mathbb R_+$ is a non-decreasing function such that $\sum_{\alpha\geqslant 0}f(2^\alpha)/2^{\alpha}<+\infty$. Then the sequence $\{b_n/n\}_{n\in\mathbb N,\,n\geqslant 1}$ converges in $\mathbb R_+$. However the same proof applies to a general (not necessarily positive) sequence satisfying the same weak subadditivity condition and leads to the convergence of the sequence $\{b_n/n\}_{n\in\mathbb N,\,n\geqslant 1}$ in $\mathbb R\cup\{-\infty\}$.}, the sequence $\{\widehat{\mu}_{\min}(\overline R_n)/n\}_{n\in\mathbb N}$ converges to an element in $\mathbb R\cup\{+\infty\}$. 
\end{proof}

\begin{defi}
Let $\overline{R}_\sbullet$ be a normed graded algebra  
on $S$. Let $M_\sbullet=\bigoplus_{n\in\mathbb Z}M_n$ be a $\mathbb Z$-graded $K$-linear space and $h$ be a positive integer. We say that $M_\sbullet$ is a \emph{$h$-graded $R_\sbullet$-module}\index{h-graded R-module@$h$-graded $R_\sbullet$-module} if $M_\sbullet$ is equipped with a structure of $R_\sbullet$-module such that
\[\forall\,(n,m)\in\mathbb N\times\mathbb Z,\quad\forall (a,x)\in R_n\times M_m,\quad ax\in M_{nh+m}.\]

Let $M_\sbullet$ be an $h$-graded $R_\sbullet$-module. Assume that each homogeneous component $M_n$ is of finite rank over $K$ and is equipped with a norm family $\xi_n'=\{\|\ndot\|_{n,\omega}'\}_{\omega\in\Omega}$. 
We say that $\overline{M}_\sbullet=\{\overline M_n\}_{n\in\mathbb Z}$ is a 
\emph{normed $h$-graded $\overline R_\sbullet$-module}\index{normed h-graded R-module@normed $h$-graded $\overline R_\sbullet$-module} if,
for any $\omega \in \Omega$, $\overline M_{\sbullet, \omega} = \{ (M_{n,\omega}, \|\ndot\|'_{n,\omega}) \}_{n \in \mathbb Z}$
forms a normed $h$-graded $\overline R_{\sbullet, \omega}$-module,
where $M_{n,\omega} := M_{n} \otimes_K K_{\omega}$ (cf. Subsetion~\ref{Subsec:graded:algebra:norms}).
We say that an $h$-graded $\overline R_\sbullet$-module $\overline M_\sbullet$ is \emph{of finite type}\index{of finite type@of finite type} if the underlying $h$-graded $R_\sbullet$-module $M_\sbullet$ is of finite type.
Moreover, if $(M_n, \xi'_n)$ forms an adelic vector bundle on $S$ for all $n \in \mathbb Z$, then
$\overline M_\sbullet$ is called a \emph{$h$-graded $\overline R_\sbullet$-module of adelic vector bundles on $S$}\index{h-graded R-module of adelic vector bundles on S@$h$-graded $\overline R_\sbullet$-module of adelic vector bundles on $S$}.
\end{defi}

\begin{prop}\label{Pro: estimate of asymptotique minimal slope} Let $C_0$ be a non-negative constant. We assume that the adelic curve $S$  satisfies the tensorial minimal slope property of level $\geqslant C_0$.
Let $\overline R_\sbullet$ be a graded algebra of adelic vector bundles which is of finite type, and $\overline{M}_\sbullet$ be an $h$-graded $\overline R_\sbullet$-module of adelic vector bundles on $S$ such that
$\overline{M}_\sbullet$ is of finite type, where $h$ is a positive integer. Then one has
\begin{equation}\label{Equ: asymptotic mu min of graded algebra}\liminf_{n\rightarrow+\infty}\frac{\widehat{\mu}_{\min}(\overline M_n)}{n}\geqslant \frac 1h\liminf_{n\rightarrow+\infty}\frac{\widehat{\mu}_{\min}(\overline R_n)}{n} >-\infty.\end{equation} 
\end{prop}
\begin{proof}
If we replace $\overline R_0$ by the trivial adelic line bundle, we obtain a new graded algebra of adelic vector bundles (denoted by $\overline{R}'_\sbullet$) and $\overline M_\sbullet$ is naturally equipped with a structure of module over this graded algebra of adelic vector bundles. Moreover, $R_\sbullet$ is a finite $R_\sbullet'$-algebra since $R_0$ is supposed to be of finite rank over $K$. In particular, $M_\sbullet$ is a module of finite type over $R_\sbullet'$. If $\{a_i\}_{i\in I}$ is a basis of $R_0$ over $K$ which contains $1\in R_0$ and if $\{b_j\}_{j\in J}$ is a finite family of homogeneous elements of positive degree in $R_\sbullet$, which generates $R_\sbullet$ as $R_0$-algebra, then $R_\sbullet'$ is generated as $K$-algebra by $\{a_ib_j\}_{(i,j)\in I\times J}$. This shows that $R_\sbullet'$ is a $K$-algebra of finite type. Therefore (by replacing $\overline R_\sbullet$ by $\overline R_{\sbullet}'$) we may assume without loss of generality that $\overline R_0$ is the trivial adelic line bundle.

We first prove the proposition in the particular case where $R_\sbullet$ is generated as $K$-algebra by $R_1$. Let $A$ be the infimum limit of the sequence $\{\widehat{\mu}_{\min}(\overline R_n)/n\}_{n\in\mathbb N,\,n\geqslant 1}$. By \cite[Lemma 2.1.6]{EGAII}, there exist integers $b_1$ and $m>0$ such that, for any integer $b$ with $b\geqslant b_1$ and any integer $\ell\geqslant 1$ the canonical $K$-linear map $R_{\ell m}\otimes_KM_{b}\rightarrow M_{b+\ell mh}$ is surjective. Hence by Proposition \ref{Pro:inegalitdepente}, one has \[\widehat{\mu}_{\min}(\overline M_{b+\ell mh})\geqslant\widehat{\mu}_{\min}(\overline R_{\ell m}\otimes_{\varepsilon,\pi}\overline M_b),\]
which leads to 
\[\widehat{\mu}_{\min}(\overline M_{b+\ell mh})\geqslant\widehat{\mu}_{\min}(\overline R_{\ell m})+\widehat{\mu}_{\min}(\overline M_b)-C_0\ln(\rang_K(R_{\ell m})\cdot\rang_K(M_b)).\]
Dividing the two sides of the inequality by $\ell mh$ and then letting $\ell$ tend to the infinity, we obtain
\[\liminf_{\ell\rightarrow+\infty}\frac{\widehat{\mu}_{\min}(\overline M_{b+\ell mh})}{\ell mh}\geqslant \frac{1}{h}A,\]
where we have used the fact that 
\[\lim_{\ell\rightarrow+\infty}\frac{\ln(\rang_K(R_{\ell m}))}{\ell}=0.\] 
Since $b\geqslant b_1$ is arbitrary, we obtain
\[\liminf_{n\rightarrow+\infty}\frac{\widehat{\mu}_{\min}(\overline M_n)}{n}\geqslant\frac 1hA.\]

We now consider the general case. By \cite[Lemma 2.1.6]{EGAII}, there exists a positive integer $u$ such that $R^{(u)}_\sbullet:=\bigoplus_{n\in\mathbb N}R_{un}$ is generated as $K$-algebra by $R^{(u)}_1=R_u$. Moreover, $R_\sbullet$ is a $u$-graded $R^{(u)}$-module of finite type and hence a finite $R^{(u)}_\sbullet$-algebra. Therefore $M_\sbullet$ is an $hu$-graded $R^{(u)}_\sbullet$-algebra of finite type. Let $B$ be the infimum limit of the sequence $\{\widehat{\mu}_{\min}(\overline R_{nu})/n\}_{n\in\mathbb N,\,n\geqslant 1}$. By applying the particular case of the proposition established above, we obtain
\[\liminf_{n\rightarrow+\infty}\frac{\widehat{\mu}_{\min}(\overline R_{n})}{n}\geqslant \frac{B}{u}\quad\text{and}\quad\liminf_{n\rightarrow+\infty}\frac{\widehat{\mu}_{\min}(\overline M_n)}{n}\geqslant\frac{B}{hu}.\]
Note that the first inequality actually implies that
\[\liminf_{n\rightarrow+\infty}\frac{\widehat{\mu}_{\min}(\overline R_n)}{n}=\frac{B}{u}\]
since $(\widehat{\mu}_{\min}(\overline R_{nu})/n)_{n\in\mathbb N,\,n\geqslant 1}$ converges to $B$. The inequality \eqref{Equ: asymptotic mu min of graded algebra} is thus proved. 

Finally, if the adelic curve $S$ satisfies the minimal slope property, by Proposition \ref{Pro: sequence of normalised mu min}  the sequence $(\widehat{\mu}_{\min}(\overline R_{nu})/n)_{n\in\mathbb N,\,n\geqslant 1}$ converges to some element $\mathbb R\cup\{+\infty\}$. Hence the last statement of the proposition is true. 
\end{proof}

\begin{rema}
Let $\overline R_\sbullet$ be a graded algebra of adelic vector bundles,  $I_\sbullet$ be a homogeneous ideal of $R_\sbullet$ and $R'_\sbullet$ be the quotient algebra $R_\sbullet/I_\sbullet$. If we equip $R'_n$ with the quotient norm family of that of $R_n$, then $\overline R_\sbullet'$ is a graded algebra of adelic vector bundles, denoted $\overline{R_\sbullet/I_\sbullet}$ (cf. Proposition~\ref{prop:norm:graded:module:quotient} \ref{Item: quotient algebra seminorms}).

More generally, let $\overline M_\sbullet$ be an $h$-graded $\overline R_\sbullet$-module and $Q_\sbullet$ is a graded quotient $R_\sbullet$-module of $M_\sbullet$. If we equip each $Q_n$ with the quotient norm family of that of $\overline M_n$, then $\overline{Q}_\sbullet$ becomes an $h$-graded $\overline R_\sbullet$-module (cf. Proposition~\ref{prop:norm:graded:module:quotient} \ref{Item: quotient module seminorm}).

Let $\overline{R}_\sbullet$ be a graded algebra of adelic vector bundles, $\overline M_\sbullet$ be an $h$-graded module, where $h\in\mathbb N$, $h\geqslant 1$. Assume that $I_\sbullet$ is a homogeneous ideal of $R_\sbullet$. Assume that $M_\sbullet$ is annihilated by $I_\sbullet$, then $\overline{M}_\sbullet$ is naturally equipped with a structure of $h$-graded $\overline{R_\sbullet/I_\sbullet}$-module (cf. Proposition~\ref{prop:norm:graded:module:quotient} \ref{Item: quotient module seminorm}).
\end{rema}

\begin{prop}\label{prop:lambda:graded:ring:module aa} We suppose that the adelic curve $S$ satisfies the tensorial minimal slope property. 
Let $\overline R_\sbullet$ be a graded algebra of adelic vector bundle, $I_\sbullet$,  $J_\sbullet$ and $M_\sbullet$ be homogeneous ideals of $R_\sbullet$ such that $J_\sbullet\subseteq M_\sbullet$ and $I_\sbullet\cdot M_\sbullet \subseteq J_\sbullet$. Let $R'_\sbullet=R_\sbullet/I_\sbullet$ and $Q_\sbullet=M_\sbullet/J_\sbullet$. For each $n\in\mathbb N$, we equip $R'_n$ and $Q_n$ with the quotient norm families of that of $\overline R_n$ and $\overline M_n$ respectively. Then one has
\[\liminf_{n\rightarrow+\infty}\frac{\widehat{\mu}_{\min}(\overline Q_n)}{n}\geqslant\liminf_{n\rightarrow+\infty}\frac{\widehat{\mu}_{\min}(\overline R'_n)}{n}.\]
\end{prop}
\begin{proof}
By the above remark, $\overline Q_\sbullet$ is equipped with a structure of graded $\overline{R}_\sbullet'$-module. Hence the statement follows from Proposition \ref{Pro: estimate of asymptotique minimal slope}.
\end{proof}

\section{Fundamental estimations}\label{sec:NM:Fundamental:estimation}

In this section, we prove some lower bound for the asymptotic minimal slope. We assume that the adelic curve $S$ satisfies the tensorial minimal slope property of level $\geqslant C_0$, where $C_0\geqslant 0$.
 
Let $\overline R_\sbullet=\{(R_n,\xi_n)\}_{}$ be a graded algebra of adelic vector bundles which is of finite type. We assume that $R_\sbullet$ is an integral ring. Let $X=\mathrm{Proj}(R_\sbullet)$ be the projective spectrum of $R_\sbullet$. If $Y$ is an integral closed subscheme of $X$ and $P_\sbullet\subseteq R_\sbullet$ is the defining homogeneous prime ideal of $Y$, we denote by $R_{Y,\sbullet}$ the quotient graded ring $R_\sbullet/P_\sbullet$. Note that each $R_{Y,n}$ is naturally equipped with the quotient norm family $\xi_{Y,n}$ of $\xi_n$ so that $\overline R_{Y,\sbullet}$ becomes a graded algebra of adelic vector bundles (cf. 
Proposition~\ref{prop:norm:graded:module:quotient}, Proposition~\ref{Pro:dominancealgebraic} and Proposition~\ref{Pro:mesurabilityofquotient}).

\begin{theo}\label{Thm:main estimates Nakai-Moiseshen}
Let $\mathfrak S_X$ be the set of all integral closed subschemes of $X$. To each $Y\in\mathfrak S_X$ we assigne a real number $\upsilon_Y$, a positive integer $n_Y$ and a non-zero element $s_Y$ in $R_{Y,n_Y}$ such that $\widehat{\deg}_{\xi_{Y,n_Y}}(s_Y)\geqslant n_Y\upsilon_Y$. Then there exists a finite subset $\mathfrak S$ of $\mathfrak S_X$ such that
\[\liminf_{n\rightarrow+\infty}\frac{\widehat{\mu}_{\min}(\overline R_n)}{n}\geqslant \min\{\upsilon(Y)\,:\,Y\in \mathfrak S\}.\]
\end{theo}
\begin{proof}

{\bf Step~1}: 
For any positive integer $h$, we set 
\[
R^{(h)}_n := R_{hn}\quad\text{and}\quad
R^{(h)} = \bigoplus_{n\in\mathbb N} R^{(h)}_n.
\]
If we assign $\upsilon_Y^h$, $hn_Y$ and $s_Y^h$ to each $Y \in \mathfrak{S}_X$, 
then $s_Y^h \in R_{Y, hn_Y} \setminus \{ 0 \}$ and
\[ \widehat{\deg}_{\xi_{Y,hn_Y}}(s_Y^h)\geqslant h\cdot\widehat{\deg}_{\xi_Y,n_Y}(s_Y)\geqslant  h\upsilon_Y,
\]
so that the above assignment satisfies the condition of the theorem for $R^{(h)}$.
Moreover, $R$ is a finitely generated $h$-graded $R^{(h)}$-module
(cf. \cite[Lemma~5.44]{Moriwaki2014}).
By using Proposition~\ref{prop:lambda:graded:ring:module aa}, 
we can see that if the theorem holds for $R^{(h)}$,
then it holds for $R$.
Therefore, by \cite[Chapitre~III, \S1, Proposition~3]{Bourbaki85}, 
we may assume that $R$ is generated by $R_1$ over $R_0$ and
$n_X=1$. Let $\mathscr O_{X}(1)$ be the tautological invertible sheaf of $X$ arising from $R_1$.

\medskip
We prove the theorem by induction on $d = \dim X$.

\medskip
{\bf Step~2}:
In the case where $d = 0$, $X = \Spec(F)$ for some finite extension field $F$ over $K$, so that
$R_n \subseteq H^0(X, \mathscr O_X(n)) \cong F$.
Therefore, $\dim_{K}(R_n)  \leqslant [F : K]$ for all $n\in\mathbb N$.
Let us consider the following sequence of homomorphisms:
\[
R_0 \overset{s_X\cdot}{\longrightarrow} R_1 \overset{s_X\cdot}{\longrightarrow} R_2 \overset{s_X\cdot}{\longrightarrow} R_3 
\overset{s_X\cdot}{\longrightarrow} \cdots
\overset{s_X\cdot}{\longrightarrow} R_{n-1} \overset{s_X\cdot}{\longrightarrow} R_n \overset{s_X\cdot}{\longrightarrow} \cdots,
\]
Note that each homomorphism is injective and $\dim_K(R_n)$ is bounded, so that
we can find a positive integer $N$ such that $R_{n} \overset{s_X\cdot}{\longrightarrow} R_{n+1}$ is an isomorphism for all $n \in\mathbb N_{N}$.
Therefore, by Proposition~\ref{Pro:inegalitdepente},
\[\widehat{\mu}_{\min}(\overline R_n)\geqslant
\widehat{\mu}_{\min}(\overline R_N)+(n-N)\,\widehat{\deg}_{\xi_1}(s_X)\geqslant \widehat{\mu}_{\min}(\overline R_N)+(n-N)\upsilon_X, 
\]
which leads to
\[\liminf_{n\rightarrow+\infty}\frac{\widehat{\mu}_{\min}(\overline R_n)}{n}\geqslant\upsilon_X.\]

\medskip
{\bf Step~3}:
We assume $d > 0$.
Let $I_\sbullet$ be the homogeneous ideal generated by $s_X$, that is, $I_\sbullet = R_\sbullet s_X$. By using the same ideas as in 
\cite[Chapter~I, Proposition~7.4]{Hart77},
we can find a sequence 
\[
I_\sbullet = I_{0,\sbullet} \subsetneq I_{1,\sbullet} \subsetneq \cdots \subsetneq I_{r,\sbullet} = R_\sbullet
\]
of homogeneous ideals of $R_\sbullet$ and non-zero homogeneous prime ideals $P_{1,\sbullet}, \ldots, P_{r,\sbullet}$ of $R_{\sbullet}$ such that
$P_{i,\sbullet} \cdot I_{i,\sbullet} \subseteq I_{i-1,\sbullet}$ for $i\in\{1, \ldots, r\}$.

\medskip
{\bf Step~4}:
Consider the following sequence:
\[
\begin{array}{ccccccc}
{R}_{0}  & \overset{\cdot s_X}{\longrightarrow} & {I}_{0, 1} & \hookrightarrow \cdots \hookrightarrow &
{I}_{i, 1} & \hookrightarrow \cdots \hookrightarrow & {I}_{r, 1} = {R}_1 \\
 & \vdots &  \vdots & \vdots &  \vdots & \vdots & \vdots \\
& \overset{\cdot s_X}{\longrightarrow} & {I}_{0, j} & \hookrightarrow \cdots \hookrightarrow &
{I}_{i, j} & \hookrightarrow \cdots \hookrightarrow & {I}_{r, j} = {R}_j \\
& \overset{\cdot s_X}{\longrightarrow} & {I}_{0, j+1} & \hookrightarrow \cdots \hookrightarrow &
{I}_{i, j+1} & \hookrightarrow \cdots \hookrightarrow & {I}_{r, j+1} = {R}_{j+1} \\
 & \vdots &  \vdots & \vdots &  \vdots & \vdots & \vdots \\
& \overset{\cdot s_X}{\longrightarrow} & {I}_{0, n} & \hookrightarrow \cdots \hookrightarrow &
{I}_{i, n} & \hookrightarrow \cdots \hookrightarrow & {I}_{r, n} = {R}_n
\end{array}
\]
By using Proposition~\ref{Pro:successive embeddings}, one has
\begin{equation}
\label{eqn:thm:lambda:estimate:01}
\widehat{\mu}_{\min}(\overline R_n)\geqslant\min\Big\{\min_{\begin{subarray}{c}i\in\{1,\ldots,r\}\\
j\in\{1,\ldots,n\}
\end{subarray}}\widehat{\mu}_{\min}(\overline{I_{i,j}/I_{{i-1},j}})+(n-j)\upsilon_X,\widehat{\mu}_{\min}(\overline R_0)+n\upsilon_X\Big\}.
\end{equation}
For any $i\in\{1,\ldots,r\}$ let $Y_i$ be the integral closed subscheme defined by $P_i$. By Proposition~\ref{prop:lambda:graded:ring:module aa}, 
one has \[\liminf_{m\rightarrow+\infty}\frac{\widehat{\mu}_{\min}(\overline{I_{i,m}/I_{i-1,m}})}{m}\geqslant\liminf_{m\rightarrow+\infty}\frac{\widehat{\mu}_{\min}(\overline R_{Y_i,m})}{m}
\]
Moreover, by the induction hypothesis,
there are a finite subset $\mathfrak S_i$ of $\mathfrak{S}_{Y_i}$ such that
\[\liminf_{m\rightarrow+\infty}\frac{\widehat{\mu}_{\min}(\overline R_{Y_i,m})}{m}\geqslant\min\{\upsilon_Z\,:\,Z\in\mathfrak S_i\}.
\] 
There the estimate \eqref{eqn:thm:lambda:estimate:01} leads to 
\[\liminf_{n\rightarrow+\infty}\frac{\widehat{\mu}_{\min}(\overline R_n)}{n}\geqslant \min\bigg\{\upsilon_Z\,:\,Z\in\{X\}\cup\bigcup_{i=1}^r\mathfrak S_i\bigg\}.\]
The theorem is thus proved.
\end{proof}

\begin{rema}\label{Rem: minoration of mu Rn}
Consider the following variante of the above theorem. Assume that $R_\sbullet$ is generated as $R_0$-algebra by $R_1$. By using Proposition \ref{Pro:suiteexactedeg}, we obtain that, for integers $n$ and $m$ such that $1\leqslant m\leqslant n$, one has 
\[
\begin{split}&\quad\;\widehat{\deg}(\overline R_n)\geqslant\sum_{j=1}^n\sum_{i=1}^r\widehat{\deg}(\overline{I_{i,j}/I_{i-1,j}})+\widehat{\deg}(\overline R_0)+\upsilon_X\sum_{k=0}^{n-1}\rang(R_k)\\
&\geqslant\sum_{j=1}^m\sum_{i=1}^r\widehat{\deg}(\overline{I_{i,j}/I_{i-1,j}})+\min_{i\in\{1,\ldots,r\}}\inf_{\ell\in\mathbb N_{\geqslant m}}\frac{\widehat{\mu}_{\min}(\overline{I_{i,\ell}/I_{i-1,\ell}})}{\ell}\sum_{j=m+1}^n j\rang(R_j/R_{j-1})
\\&
\qquad+\widehat{\deg}(\overline R_0)+\upsilon_X\sum_{k=0}^{n-1}\rang(R_k).
\end{split}\]
Dividing the two sides by $n\rang(R_n)$ and letting $n$ tend to the infinity, we obtain
\[\liminf_{n\rightarrow+\infty}\frac{\widehat{\mu}(R_n)}{n}\geqslant \frac{d}{d+1}\min_{i\in\{1,\ldots,r\}}\inf_{\ell\in\mathbb N_{\geqslant m}}\frac{\widehat{\mu}_{\min}(\overline{I_{i,\ell}/I_{i-1,\ell}})}{\ell}+\frac{1}{d+1}\upsilon_X,\]
where we have used the geometric Hilbert-Samuel theorem asserting that $\rang(R_n)=\deg(X)n^d+O(n^{d-1})$, with $d$ being the Krull dimension of the scheme $X$, which leads to
\[\lim_{n\rightarrow+\infty}\frac{1}{n\rang(R_n)}\sum_{j=0}^{n-1}\rang(R_j)=\frac{1}{d+1}.\]
Since $m$ is arbitrary, we obtain 
\begin{equation}\liminf_{n\rightarrow+\infty}\frac{\widehat{\mu}(R_n)}{n}\geqslant\frac{d}{d+1}\min_{i\in\{1,\ldots,r\}}\min_{Z\in\mathfrak S_i}\upsilon_Z+\frac{1}{d+1}\upsilon_X.\end{equation}
Note that in the general case where $R_\sbullet$ is not necessarily generated by $R_1$, the same argument leads to 
\[\limsup_{n\rightarrow+\infty}\frac{\widehat{\mu}(R_n)}{n}\geqslant\frac{d}{d+1}\min_{i\in\{1,\ldots,r\}}\min_{Z\in\mathfrak S_i}\upsilon_Z+\frac{1}{d+1}\upsilon_X.\]
\end{rema}

Under the strong tensorial minimal slope property (see Definition \ref{Def: strong tensorial minimal slope property}), Theorem \ref{Thm:main estimates Nakai-Moiseshen} admits the following analogue.

\begin{theo}\label{Thm: Nakai Moiseshen under strong tensor minimal slope}
We assume that the adelic curve $S$ satisfies the \emph{strong} tensorial minimal slope property of level $\geqslant C_1$, where $C_1\in\mathbb R_{\geqslant 0}$. Let $\mathfrak S_X$ be the set of all integral closed subschemes of $X$. Then one has
\[\liminf_{n\rightarrow+\infty}\frac{\widehat{\mu}_{\min}(\overline R_n)}{n}\geqslant\inf_{Y\in\mathfrak S_X}\limsup_{m\rightarrow+\infty}\frac{\widehat{\mu}_1(\overline{R}_{Y,m})}{m}.\]
\end{theo}
\begin{proof}
We reason by induction on the dimension $d$ of the scheme $X$. 

First we treat the case where $d=0$. Let $m$ be an integer, $m\geqslant 1$. Let $E$ be a vector subspace of $R_m$ such that $\widehat{\mu}_{\min}(\overline E)=\widehat{\mu}_1(\overline R_m)$. There exists any integer $N\in\mathbb N_{\geqslant 1}$ such that, for any $p\in\mathbb N_{\geqslant 1}$, the canonical $K$-linear map
\[R_{mN}\otimes E^{\otimes p}\longrightarrow R_{m(N+p)}\] 
is surjective. Therefore, by Proposition \ref{Pro:inegalitdepente} and the strong tensorial tensor minimal slope property (by an argument similar to the {\bf Step 2} of the proof of Theorem \ref{Thm:main estimates Nakai-Moiseshen}), one has
\[\begin{split}\widehat{\mu}_{\min}(\overline R_{m(N+p)})&\geqslant\widehat{\mu}_{\min}(\overline R_{mN})-C_1\ln(\rang(R_{mN}))+p\big(\widehat{\mu}_1(\overline R_m)-C_1\ln(\rang(E))\big)\\
&\geqslant\widehat{\mu}_{\min}(\overline R_{mN})-C_1\ln(\rang(R_{mN}))+p\big(\widehat{\mu}_1(\overline R_m)-C_1\ln(\rang(R_m))\big).
\end{split}\]
Dividing the two sides by $m(N+p)$ and letting $p$ tend to $+\infty$, by Proposition \ref{prop:lambda:graded:ring:module aa} we obtain 
\[\liminf_{n\rightarrow+\infty}\frac{\widehat{\mu}_{\min}(\overline{R}_n)}{n}\geqslant\frac{\widehat{\mu}_1(\overline R_m)}{m}-C_1\frac{\ln(\rang(R_m))}{m}.\]
Note that 
\[\lim_{m\rightarrow+\infty}\frac{\ln(\rang(R_m))}{m}=0.\]
Therefore, by taking the limsup when $m\rightarrow+\infty$, we obtain
\[\liminf_{n\rightarrow+\infty}\frac{\widehat{\mu}_{\min}(\overline{R}_n)}{n}\geqslant\limsup_{m\rightarrow+\infty}\frac{\widehat{\mu}_1(\overline R_m)}{m}.\]

We now assume that $d\geqslant 1$. Let $m$ be an integer such that  $m\geqslant 1$. Let $E$ be a vector subspace of $R_m$ such that $\widehat{\mu}_{\min}(\overline E)=\widehat{\mu}_1(\overline R_m)$. Let $I_\sbullet$ be the homogeneous ideal of $R^{(m)}_\sbullet=\bigoplus_{n\in\mathbb N}R_{mn}$ generated by $E$. That is, for any $n\in\mathbb N$, $I_n$ is the image of the canonical homomorphism
\[R_{(n-1)m}\otimes E\longrightarrow R_{nm}.\] As in the the proof of Theorem \ref{Thm:main estimates Nakai-Moiseshen}, we let 
\[I_\sbullet=I_{0,\sbullet}\subsetneq I_{1,\sbullet}\subsetneq\ldots\subsetneq I_{r,\sbullet}=R_\bullet^{(m)}\]
be a sequence of homogeneous ideals of $R^{(m)}_\sbullet$ and $P_{1,\sbullet},\ldots,P_{r,\sbullet}$ be non-zero homogeneous prime ideals of $R^{(m)}$ such that $P_{i,\sbullet}\cdot I_{i,\sbullet}\subset I_{i-1,\sbullet}$ for $i\in\{1,\ldots,r\}$. Let $p$ be an integer in $\mathbb N_{\geqslant 1}$. We denote by $F_p$ the image of the canonical $K$-linear map
\[R_0\otimes E^{\otimes p}\longrightarrow R_{mp}.\]
Consider the following sequence:
\[
\begin{array}{ccccccc}
F_p & = & {I}_{0, 1}E^{p-1} & \hookrightarrow \cdots \hookrightarrow &
{I}_{i, 1}E^{p-1} & \hookrightarrow \cdots \hookrightarrow & {I}_{r, 1}E^{p-1} \\
 & \vdots &  \vdots & \vdots &  \vdots & \vdots & \vdots \\
& = & {I}_{0, j}E^{p-j} & \hookrightarrow \cdots \hookrightarrow &
{I}_{i, j}E^{p-j} & \hookrightarrow \cdots \hookrightarrow & {I}_{r, j}E^{p-j} \\
& = & {I}_{0, j+1}E^{p-j-1} & \hookrightarrow \cdots \hookrightarrow &
{I}_{i, j+1}E^{p-j-1} & \hookrightarrow \cdots \hookrightarrow & {I}_{r, j+1}E^{p-j-1}\\
 & \vdots &  \vdots & \vdots &  \vdots & \vdots & \vdots \\
& = & {I}_{0, p} & \hookrightarrow \cdots \hookrightarrow &
{I}_{i, p} & \hookrightarrow \cdots \hookrightarrow & {I}_{r, p} = {R}_{mp}
\end{array}
\]
By Proposition \ref{Pro:successive embeddings} we obtain that 
\begin{equation}\label{Equ: minoration de mu min Rmp}\widehat{\mu}_{\min}(\overline{R}_{mp})\geqslant\min\Big\{\widehat{\mu}_{\min}(\overline F_p),\min_{\begin{subarray}{c}i\in\{1,\ldots,r\}\\
j\in\{1,\ldots,p\}\end{subarray}}\widehat{\mu}_{\min}(\overline{I_{i,j}E^{p-j}/I_{i-1,j}E^{p-j}})\Big\}.\end{equation}
By Proposition \ref{Pro:inegalitdepente} and the strong tensorial tensor minimal slope property, one has
\[\begin{split}\widehat{\mu}_{\min}(\overline F_p)&\geqslant\widehat{\mu}_{\min}(\overline R_0)-C_1\ln(\rang(R_0))+p\big(\widehat{\mu}_{\min}(\overline E)-C_1\ln(\rang(E))\big)\\
&\geqslant\widehat{\mu}_{\min}(\overline R_0)-C_1\ln(\rang(R_0))+p\big(\widehat{\mu}_{1}(\overline R_m)-C_1\ln(\rang(R_m))\big).
\end{split}\]
Similarly, for any $(i,j)\in\{1,\ldots,r\}\times\{1,\ldots,n\}$, one has
\[\begin{split}\widehat{\mu}_{\min}(\overline{I_{i,j}E^{p-j}/I_{i-1,j}E^{p-j}})\geqslant\widehat{\mu}_{\min}&(\overline{I_{i,j}/I_{i-1,j}})-C_1\ln(\rang(I_{i,j}/I_{i-1,j}))\\
&+(p-j)\big(\widehat{\mu}_{\min}(\overline E)-C_1\ln(\rang(E))\big).
\end{split}\]
For any $i\in\{1,\ldots,r\}$ let $Y_i$ be the integral closed subscheme defined by $P_i$. By Proposition~\ref{prop:lambda:graded:ring:module aa}, 
one has \[\liminf_{j\rightarrow+\infty}\frac{\widehat{\mu}_{\min}(\overline{I_{i,j}/I_{i-1,j}})}{j}\geqslant\liminf_{j\rightarrow+\infty}\frac{\widehat{\mu}_{\min}(\overline R_{Y_i,j})}{j}\geqslant\min_{Z\in\mathfrak S_{Y_i}}
\limsup_{k\rightarrow+\infty}\frac{\widehat{\mu}_1(\overline{R}_{Z,k})}{k},
\]
where the second inequality comes from the induction hypothesis.
Therefore, if we denote by $\upsilon$ the value
\[\inf_{Z\in\mathfrak S_X}\limsup_{k\rightarrow+\infty}\frac{\widehat{\mu}_{1}(\overline{R}_{Z,k})}{k},\]
then the inequality \eqref{Equ: minoration de mu min Rmp} leads to 
\[\liminf_{n\rightarrow+\infty}\frac{\widehat{\mu}_{\min}(\overline{R}_n)}{n}=\liminf_{p\rightarrow+\infty}\frac{\widehat{\mu}_{\min}(\overline{R}_{mp})}{mp}\geqslant\min\Big\{\frac{\widehat{\mu}_1(\overline{R}_m)}{m}-C_1\frac{\ln(\rang(R_m))}{m},\upsilon\Big\},\]
where the equality comes from Proposition~\ref{prop:lambda:graded:ring:module aa}. By taking the limsup when $m\rightarrow+\infty$, we obtain 
\[\liminf_{n\rightarrow+\infty}\frac{\widehat{\mu}_{\min}(\overline{R}_n)}{n}\geqslant\upsilon,\]
as desired.
\end{proof}


\section[A consequence of the extension property]{A consequence of the extension property of semipositive metrics}

Let $S=(K,(\Omega,\mathcal A,\nu),\phi)$ be an adelic curve which satisfies the tensorial minimal slope property.
We assume that, either (i) $\Omega_0$ is empty, or (ii) the field $K$ is countable, or (iii) $\Omega = \Omega_0$ and
$\#(\Omega_0) = 1$.
The purpose of this section is to prove the following theorem
as a consequence of the extension property of semipositive metrics
(cf. Theorem~\ref{thm:extension:property:over:Archimedean} and Theorem~\ref{Thm:extensionpropertynontrivial}).

\begin{theo}\label{thm:dominate:measure:lim}
Let $X$ be a geometrically reduced projective $K$-scheme, 
$L$ be a semiample invertible $\mathcal O_X$-module and
$\varphi = \{ \varphi_{\omega} \}_{\omega \in \Omega}$ be a metric family of $L$. 
Let $Y$ be a geometrically reduced closed subscheme of $X$,
$\xi_n := \{ \|\ndot\|_{n\varphi_{\omega}, \sup} \}_{\omega \in \Omega}$ and
$\rest{\xi_n}{Y} := \{ \|\ndot\|_{n\rest{\varphi_{\omega}}{Y_{\omega}}, \sup} \}_{\omega \in \Omega}$.
Let $R_{Y, n}$ be the image of $H^0(X, L^{\otimes n}) \to H^0(Y, \rest{L}{Y}^{\otimes n})$ and
$\xi_{Y,n}=\{\norm{\ndot}_{Y,n,\omega}\}_{\omega\in\Omega}$ be the quotient norm family on $R_{Y,n}$ obtained by
$H^0(X, L^{\otimes n}) \to R_{Y,n}$ and $\xi_n$. 
If $\varphi$ is dominated and measurable and $\varphi_{\omega}$ is semipositive for all $\omega \in \Omega$,
then we have the following:
\begin{enumerate}
\renewcommand{\labelenumi}{(\arabic{enumi})}
\item
$\xi_n$, $\rest{\xi_n}{Y}$ and $\xi_{Y,n}$ are dominated and measurable for all $n \geqslant 0$.
\item
For $s \in R_{Y,1} \setminus \{ 0 \}$, one has
\[
\lim_{n\to\infty} \frac{\widehat{\deg}_{\xi_{Y, n}}(s^{\otimes n})}{n} = \widehat{\deg}_{\rest{\xi_1}{Y}}(s).
\]
\end{enumerate}
\end{theo}

\begin{proof}
(1) In the case (iii), the assertion is obvious, so that we assume (i) and (ii).
First, by Theorem~\ref{Thm: Fubini-Study dominated} and Theorem~\ref{Thm: measurability of linear series},
$\xi_n$ is dominated and measurable for $n \geqslant 0$. 
Moreover, by Theorem~\ref{Thm: Fubini-Study dominated} and Theorem~\ref{Thm: measurability of linear series} together with Proposition~\ref{Pro: dominancy preserved by operators} and Proposition~\ref{Pro: measurability preserved by operators},  
$\rest{\xi_n}{Y}$ is dominated and measurable for $n \geqslant 0$.
Finally, by virtue of Proposition~\ref{Pro:dominancealgebraic} and Proposition~\ref{Pro:mesurabilityofquotient},
$\xi_{Y,n}$ is dominated and measurable for $n \geqslant 0$.

\medskip
Before starting the proof of (2), we need to prepare several facts.
Here we assume either (i) or (ii) or (iii).
We set $\xi_{Y,n} = \{ \|\ndot\|_{Y,n,\omega} \}_{\omega \in \Omega}$.
We claim the following:

\begin{enonce}{Claim}
\begin{enumerate}[label=\rm(\alph*)]

\item For all $\omega \in \Omega$, $n \geqslant 0$ and $s \in R_{Y, n, \omega}$,
\[ \| s \|_{n\rest{\varphi_{\omega}}{Y_{\omega}}, \sup} \leqslant \| s \|_{Y, n, \omega}.\]

\item For all $\omega \in \Omega$, $n \geqslant 1$ and $s \in R_{Y, 1, \omega} \setminus \{ 0\}$,
\[ \ln \| s \|_{\rest{\varphi_{\omega}}{Y_{\omega}},\sup} 
\leqslant \frac{\ln \| s^{\otimes n} \|_{Y,n,\omega}}{n} \leqslant \ln \| s \|_{Y, 1, \omega}.\]

\item For all $\omega \in \Omega$ and $s \in R_{Y, 1, \omega} \setminus \{ 0 \}$,
\[\lim_{n\to\infty} \frac{\ln \| s^{\otimes n} \|_{Y, n, \omega}}{n} = \ln \| s \|_{\rest{\varphi_{\omega}}{Y_{\omega}},\sup}.\]
\end{enumerate}
\end{enonce}

\begin{proof}
(a) Note that, for all $l \in H^0(X, L^{\otimes n})$ with $\rest{l}{Y} = s$,
one has $\|  s \|_{n\rest{\varphi_{\omega}}{Y_{\omega}}, \sup} \leqslant  \| l \|_{n\varphi_{\omega}, \sup}$, so that
the assertion follows.

(b) By Proposition~\ref{prop:norm:graded:module:quotient}, $\| s^{\otimes n} \|_{Y,n,\omega} \leqslant \left(\| s \|_{Y,1,\omega}\right)^n$. Moreover, by (a),
\[
\left(\| s \|_{\rest{\varphi_{\omega}}{Y_{\omega}},\sup}\right)^n = \| s^{\otimes n} \|_{\rest{n\varphi_{\omega}}{Y_{\omega}},\sup} \leqslant  \| s^{\otimes n} \|_{Y, n, \omega},
\]
so that one has (b).

(c) 
For a positive number $\epsilon$, by Theorem~\ref{thm:extension:property:over:Archimedean} and Theorem~\ref{Thm:extensionpropertynontrivial},
there is a positive integer $n_0$ such that, for all $n \geqslant n_0$,
we can find $l \in H^0(X_{\omega}, L_{\omega}^{\otimes n})$ such that $\rest{l}{Y_{\omega}} = s^{\otimes n}$ and
$\| l \|_{n\varphi_{\omega}, \sup} \leqslant \mathrm{e}^{n \epsilon} \left( \| s \|_{\rest{\varphi_{\omega}}{Y_{\omega}}, \sup}\right)^n$, and hence
\[
\ln \| s^{\otimes n} \|_{Y,n,\omega} \leqslant \ln \| l \|_{n\varphi_{\omega}, \sup}
\leqslant n \epsilon + n \ln \| s \|_{\rest{\varphi_{\omega}}{Y_{\omega}}, \sup}.
\]
Therefore, by (b),
\[
0 \leqslant \frac{\ln \| s^{\otimes n} \|_{Y,n,\omega}}{n} - \ln \| s \|_{\rest{\varphi_{\omega}}{Y_{\omega}},\sup} \leqslant \epsilon
\]
for all $n \geqslant n_0$, as required.
\end{proof}

(2)
By (1), $(\omega \in \Omega) \mapsto \Big|\ln \| s \|_{\rest{\varphi_{\omega}}{Y_{\omega}},\sup}\Big|$ and
$(\omega \in \Omega) \mapsto \Big|\ln \| s \|_{Y, 1, \omega}\Big|$ are integrable. Moreover, by (b), one has
\[
\left| \frac{\ln \| s^{\otimes n} \|_{Y,n,\omega}}{n}\right| \leqslant 
\max \left\{ \Big| \ln \| s \|_{\rest{\varphi_{\omega}}{Y_{\omega}},\sup}\Big|  , \Big| \ln \| s \|_{Y, 1, \omega} \Big| \right\},
\]
and hence, by Lebesgue's dominated convergence theorem together with (c),
\[
\lim_{n\to\infty} \frac{1}{n} \int_{\Omega} \ln \| s^{\otimes n} \|_{Y,n,\omega} \,\nu(d\omega) =
\int_{\Omega} \ln \| s \|_{\rest{\varphi_{\omega}}{Y_{\omega}},\sup} \,\nu(d\omega),
\]
which shows the assertion of the theorem.
\end{proof}

\section{Nakai-Moishezon's criterion in a general settings}
Let $S=(K,(\Omega,\mathcal A,\nu),\phi)$ be an adelic curve.
We assume that, either (i) $\Omega_0$ is empty, or (ii) the field $K$ is countable, or (iii) $\Omega = \Omega_0$ and
$\#(\Omega_0) = 1$.
In this section, let us consider the following Nakai-Moishezon's criterion in a general settings:

\begin{theo}\label{thm:Nakai:Moishezon:general:setting}
Let $X$ be a geometrically integral projective $K$-scheme, $L$ be an  
invertible $\mathcal O_X$-module and
$\varphi = \{ \varphi_{\omega} \}_{\omega \in \Omega}$ be a metric family of $L$.
Let $\xi_n := \{ \|\ndot\|_{n\varphi_{\omega}, \sup} \}_{\omega \in \Omega}$ and
$\rest{\xi_n}{Y} := \{ \|\ndot\|_{n\rest{\varphi_{\omega}}{Y_{\omega}}, \sup} \}_{\omega \in \Omega}$ for a subvariety $Y$ of $X$.
We assume the following:
\begin{enumerate}[label=\rm(\arabic*)]
\item {\rm (Dominancy and measurability)}
The metric family $\varphi$ is dominated and measurable.

\item {\rm (Semipositivity)}
$L$ is semiample and $\varphi_{\omega}$ is semipositive for all $\omega \in \Omega$.

\item {\rm (Bigness)}
For every subvariety $Y$ of $X$, $\rest{L}{Y}$ is big, and
there are a positive number $n_Y$ and $s_Y \in H^0(Y, \rest{L^{\otimes n_Y}}{Y}) \setminus \{ 0 \}$ such that
$\widehat{\deg}_{\rest{\xi_{n_Y}}{Y}}(s_Y) > 0$.
\end{enumerate}
Then one has 
\begin{equation}\label{eqn:thm:Nakai:Moishezon:general:setting:01}
\liminf_{n\to\infty} \frac{\hat{\mu}_{\min}\left(H^0(X, L^{\otimes n}), \xi_n \right)}{n} > 0.
\end{equation}
Moreover, if the adelic curve $S$ satisfies the strong Minkowski properties, then
\begin{equation}\label{eqn:thm:Nakai:Moishezon:general:setting:02}
\liminf_{n\to\infty} \frac{\nu_{\min}\left(H^0(X, L^{\otimes n}), \xi_n \right)}{n} > 0,
\end{equation}
so that, there are a positive integer $n$ and a basis $\{e_i\}_{i=1}^N$ of $H^0(X, L^{\otimes n})$ such that
$\widehat{\deg}_{\xi_n}(e_i) > 0$ for $i=1, \ldots, N$.
\end{theo}

\begin{proof}
First of all, for each subvariety $Y$ of $X$, as $\rest{L}{Y}$ is nef and big, one has
$(\rest{L}{Y}^{\dim Y}) > 0$, so that by the classical Nakai-Moishezon criterion, $L$ is ample.

For a subvariety $Y$ of $X$, we set
\[
\begin{cases}
R_{Y, n} :=  \text{the image of the natural homomorphism $H^0(X, L^{\otimes n}) \to H^0(Y, \rest{L}{Y}^{\otimes n})$},\\
R_{Y,n,\omega} := R_{Y,n} \otimes_K K_{\omega} \quad (\omega \in \Omega), \\
\|\ndot\|_{Y,n,\omega} := \text{the quotient norm of $\|\ndot\|_{n \varphi_{\omega}, \sup}$ on $R_{Y,n,\omega}$}\quad (\omega \in \Omega), \\
\xi_{Y, n} := \{ \|\ndot\|_{Y, n, \omega} \}_{\omega \in \Omega}.
\end{cases}
\]

\begin{enonce}{Claim}\label{claim:thm:Nakai:Moishezon:general:setting:01}
There are a positive number $n'_Y$ and $s'_Y \in R_{Y, n'_Y} \setminus \{ 0 \}$ such that
$\widehat{\deg}_{\xi_{Y, n'_Y}}(s'_Y) > 0$.
\end{enonce}

\begin{proof}
Fix a positive integer $n_0$ such that the natural homomorphism 
\[
H^0(X, L^{\otimes n}) \to H^0(Y, \rest{L}{Y}^{\otimes n})
\]
is surjective for all $n \geqslant n_0$, that is, $R_{Y,n} = H^0(Y, \rest{L}{Y}^{\otimes n})$ for all $n \geq n_0$, so
that $s_Y^{\otimes n_0} \in R_{Y, n_0n_Y} \setminus \{ 0 \}$.
By Theorem~\ref{thm:dominate:measure:lim}, (2),
one has
\[
\lim_{n\to\infty} \frac{\widehat{\deg}_{\xi_{Y, nn_0n_Y}}(s_Y^{\otimes nn_0})}{n} = \widehat{\deg}_{\rest{\xi_{n_0n_Y}}{Y}}(s_Y^{\otimes n_0}) = n_0 \widehat{\deg}_{\rest{\xi_{n_Y}}{Y}}(s_Y)> 0,
\]
so that there is a positive integer $n_1$ such that $\widehat{\deg}_{\xi_{Y, n_1n_0n_Y}}(s_Y^{\otimes n_1 n_0}) > 0$.
Therefore, of we set $n'_Y := n_1 n_0 n_Y$ and $s'_Y := s_Y^{\otimes n_1 n_0}$, one has the claim.
\end{proof}

The assertion follows from \eqref{eqn:thm:Nakai:Moishezon:general:setting:01} follows from the above claim together with Theorem~\ref{Thm:main estimates Nakai-Moiseshen}.
Further, if $S$ satisfies the strong Minkowski property, then there is a constant $C$ depending only on $S$ such that
\[
\nu_{\min}(H^0(X, L^{\otimes n}), \xi_n) + C \ln (\dim_K  H^0(X, L^{\otimes n})) \geqslant \widehat{\mu}_{\min}(H^0(X, L^{\otimes n}), \xi_n),
\]
and hence the assertion \eqref{eqn:thm:Nakai:Moishezon:general:setting:02} follows.
\end{proof}

\begin{rema}
In the case (iii) (i.e. $\Omega = \Omega_0$ and $\#(\Omega_0) = 1$),
$S$ satisfies the strong Minkowski property, that is,
if $E$ is a finite-dimensional vector space over $K$ and $\|\ndot\|$ is an ultrametric norm of $E$ over $(K,|\ndot|)$,
then $\nu_{\min}(E, \|\ndot\|) =
\widehat{\mu}_{\min}(E, \|\ndot\|)$, which can be checked as follows:

In general, one has $\nu_{\min}(E, \|\ndot\|) \leqslant
\widehat{\mu}_{\min}(E, \|\ndot\|)$ by Proposition~\ref{Pro: comparaison mu min et nu min}, so that it is sufficient to show that
$\nu_{\min}(E, \|\ndot\|) \geqslant
\widehat{\mu}_{\min}(E, \|\ndot\|)$. Let $(e_i)_{i=1}^r$ be an orthogonal basis of $E$ with respect to $\|\ndot\|$
(cf. Proposition~\ref{Pro:existenceepsorth}). 
Clearly we may assume that $\| e_r \| = \max \{ \|e_1\|, \ldots, \|e_r \| \}$.
Let $Q := E/(K e_1 + \cdots + Ke_{r-1})$ and $\|\ndot\|_Q$ be the quotient norm of $\|\ndot\|$ on $Q$.
Then $\|\pi(e_r)\|_Q = \| e_r \|$, where $\pi : E \to Q$ is the canonical homomorphism. 
Thus $-\log \| e_r \| \geqslant \widehat{\mu}_{\min}(E, \|\ndot\|)$, and hence
$\widehat{\deg}(e_i) \geqslant \widehat{\mu}_{\min}(E, \|\ndot\|)$ for all $i$. Therefore,
one has $\nu_{\min}(E, \|\ndot\|) \geqslant 
\widehat{\mu}_{\min}(E, \|\ndot\|)$.
\end{rema}

\begin{lemm}\label{Lem: asymptotic minima bounded from below by last minimum}
We assume that the adelic curve $S$ satisfies the strong Minkowski property. Let $X$ be a geometrically integral projective $K$-scheme, $L$ be an ample invertible $\mathcal O_X$-module, $\varphi=\{\varphi_\omega\}_{\omega\in\Omega}$ be a dominant and measurable metric family on $L$ such that $\varphi_\omega$ is semipositive for any $\omega\in\Omega$. If the height function $h_{(L,\varphi)}$ on the set $X(K^{\mathrm{ac}})$ of algebraic points of $X$ is bounded from below by a positive number, then one has
\[\liminf_{n\rightarrow+\infty}\frac{\widehat{\mu}_{\min}(H^0(X,L^{\otimes n}),\xi_n)}{n}> 0,\]
where for any $n\in\mathbb N$, $\xi_n=\{\norm{\ndot}_{n\varphi_\omega}\}_{\omega\in\Omega}$. 
\end{lemm}
\begin{proof}
We reason by induction on the dimension $d$ of the scheme $X$. The case where $d=0$ comes from Theorem \ref{thm:Nakai:Moishezon:general:setting}. In the following, we assume that the theorem is true for lower dimensional arithmetic varieties. Denote by 
\[\widehat{\mu}_{\min}^{\mathrm{asy}}(L,\varphi):=\liminf_{n\rightarrow+\infty}\frac{\widehat{\mu}_{\min}(H^0(X,L^{\otimes n}),\xi_n)}{n},\]
\[\widehat{\mu}_{\max}^{\mathrm{asy}}(L,\varphi):=\lim_{n\rightarrow+\infty}\frac{\widehat{\mu}_{\max}(H^0(X,L^{\otimes n}),\xi_n)}{n},\]
and
\[\widehat{\mu}^{\mathrm{asy}}(L,\varphi):=\liminf_{n\rightarrow+\infty}\frac{\widehat{\mu}(H^0(X,L^{\otimes n}),\xi_n)}{n}.\]
Assume by contradiction that $\widehat{\mu}_{\min}^{\mathrm{asy}}(L,\varphi)\leqslant 0$. If $\widehat{\mu}_{\max}^{\mathrm{asy}}(L,\varphi)>0$, then there exist a positive integer $n_X$ and $s_X\in H^0(X,L^{\otimes n_X})\setminus\{0\}$ such that $\widehat{\deg}_{\xi_{n_X}}(s_X)>0$. Moreover, by the induction hypothesis, for any subvariety $Y\subsetneq X$ one has
\[\widehat{\mu}_{\min}^{\mathrm{asy}}(L|_Y,\varphi|_Y)>0.\]
As a consequence, there are a positive number $n_Y$ and $s_Y \in H^0(Y, \rest{L^{\otimes n_Y}}{Y}) \setminus \{ 0 \}$ such that
$\widehat{\deg}_{\rest{\xi_{n_Y}}{Y}}(s_Y) > 0$. Therefore, by Theorem \ref{thm:Nakai:Moishezon:general:setting} one has $\widehat{\mu}_{\min}^{\mathrm{asy}}(L,\varphi)>0$, which leads to a contradiction. Therefore $\widehat{\mu}_{\max}^{\mathrm{asy}}(L,\varphi)\leqslant 0$. This observation actually leads to 
\[\widehat{\mu}_{\max}^{\mathrm{asy}}(L,\varphi)=\widehat{\mu}_{\min}^{\mathrm{asy}}(L,\varphi).\]
In fact, if $\widehat{\mu}_{\max}^{\mathrm{asy}}(L,\varphi)>\widehat{\mu}_{\min}^{\mathrm{asy}}(L,\varphi)$, then we can twist $(L,\varphi)$ by the pull-back of an adelic line bundle on $S$ to obtain a metric $\varphi'$ such that 
\[\widehat{\mu}_{\max}^{\mathrm{asy}}(L,\varphi')>0\geqslant \widehat{\mu}_{\min}^{\mathrm{asy}}(L,\varphi').\]
Note that $h_{(L,\varphi')}\geqslant h_{(L,\varphi)}$. By the above argument we still obtain a contradiction.

Finally, suppose that $\widehat{\mu}_{\max}^{\mathrm{asy}}(L,\varphi)=\widehat{\mu}_{\min}^{\mathrm{asy}}(L,\varphi)\leqslant 0$. By twisting $(L,\varphi)$  by the pull-back of an adelic line bundle on $S$, we may assume without loss of generality that $\widehat{\mu}_{\max}^{\mathrm{asy}}(L,\varphi)=\widehat{\mu}_{\min}^{\mathrm{asy}}(L,\varphi)= 0$. In this case one has $\widehat{\mu}^{\mathrm{asy}}(L,\varphi)=0$ since $\widehat{\mu}^{\mathrm{asy}}_{\min}(L,\varphi)\leqslant \widehat{\mu}^{\mathrm{asy}}(L,\varphi)\leqslant \widehat{\mu}_{\max}^{\mathrm{asy}}(L,\varphi)$.  By the induction hypothesis we obtain that, for any subvariety $Y\subsetneq X$, one has 
\[\widehat{\mu}_{\min}^{\mathrm{asy}}(L|_Y,\varphi|_Y)>0.\] 
However, by Remark \ref{Rem: minoration of mu Rn} we obtain that $\widehat{\mu}(L,\varphi)>0$, which leads to a contradiction. The theorem is thus proved.
\end{proof}

\begin{defi}
Let $X$ be a geometrically integral scheme over $\Spec K$, $L$ be an ample invertible $\mathcal O_X$-module, $\varphi=\{\varphi_\omega\}_{\omega\in\Omega}$ be a dominant and measurable metric family on $L$. We denote by $\nu_{\mathrm{abs}}(L,\varphi)$ the infimum of the height function $h_{D,g}$, called the \emph{absolute minimum}\index{absolute minimum@absolute minimum} of $(D,g)$.
\end{defi}

\begin{theo}
We assume that the adelic curve $S$ satisfies the strong Minkowski property. Let $X$ be a geometrically integral projective $K$-scheme, $L$ be an ample invertible $\mathcal O_X$-module, $\varphi=\{\varphi_\omega\}_{\omega\in\Omega}$ be a dominant and measurable metric family on $L$ such that $\varphi_\omega$ is semipositive for any $\omega\in\Omega$. Then the following inequality holds.
\begin{gather}\label{Equ: equality between asymptotic minimal slope and last minimum}\widehat{\mu}_{\min}^{\mathrm{asy}}(L,\varphi):=\liminf_{n\rightarrow+\infty}\frac{\widehat{\mu}_{\min}(H^0(X,L^{\otimes n}),\xi_n)}{n}=\nu_{\mathrm{abs}}(L,\varphi),
\end{gather}
where for any $n\in\mathbb N$, $\xi_n=\{\norm{\ndot}_{n\varphi_\omega}\}_{\omega\in\Omega}$, and $d$ is the Krull dimension of $X$.
\end{theo}
\begin{proof}For any $n\in\mathbb N$, let $E_n:=H^0(X,L^{\otimes n})$. Since the adelic curve $S$ satisfies the strong Minkowski property, one has
\[\widehat{\mu}_{\min}^{\mathrm{asy}}(L,\varphi)=\liminf_{n\rightarrow+\infty}\frac{\nu_{\min}(E_n,\xi_n)}{n}\leqslant \liminf_{n\rightarrow+\infty}\frac{\nu_{\min}^{\mathrm{a}}(E_n,\xi_n)}{n}\leqslant\liminf_{n\rightarrow+\infty}\frac{\widehat{\mu}_{\min}(E_{n,K^{\mathrm{ac}}},\xi_{n,K^{\mathrm{ac}}})}{n},\]
where the second inequality comes from Proposition \ref{Pro: comparison of nu under extension of fields} and the last inequality comes from Corollary \ref{Cor: comparison between mu i and nu i}.
Let $P$ be an algebraic point of $X$. For sufficiently positive integer $n$, the invertible $\mathcal O_X$-module $L^{\otimes n}$ is very ample hence defines a closed embedding $X\rightarrow\mathbb P(E_n)$. Let $\mathcal O_{E_n}(1)$ be the universal invertible sheaf on $\mathbb P(E_n)$. Then, viewed as a quotient vector space of rank $1$ of $E_{n,K^{\mathrm{ac}}}$, the Arakelov degree of $P^*(\mathcal O_{E_n}(1))$ (equipped with the quotient norm family) is bounded from above by $nh_{(L,\varphi)}(P)$ and bounded from below by $\widehat{\mu}_{\min}(E_{n,K^{\mathrm{ac}}},\xi_{n,K^{\mathrm{ac}}})$. Therefore we obtain $\widehat{\mu}_{\min}^{\mathrm{asy}}(L,\varphi)\leqslant h_{(L,\varphi)}(P)$. Since $P\in X(K^{\mathrm{ac}})$ is arbitrary, this leads to the inequality $\widehat{\mu}_{\min}^{\mathrm{asy}}(L,\varphi)\leqslant\nu_{\mathrm{abs}}(L,\varphi)$. Moreover, by Lemma \ref{Lem: asymptotic minima bounded from below by last minimum} the converse inequality also holds. Therefore the equality \eqref{Equ: equality between asymptotic minimal slope and last minimum} is proved.
  
\end{proof}

\if01
\section{Nakai-Moishezon's criterion over a trivially valued field}
We assume that the absolute value $|\ndot|$ is trivial.
Here we consider an adelic curve $S = (K, (\Omega, \mathcal{A}, \nu), \phi)$ as follows:
$\Omega$ has only one element $\omega$ (that is, $\Omega = \{ \omega \}$),
$\mathcal{A} = \{ \emptyset, \{ \omega \} \}$, $\nu(\Omega) = 1$ and $\phi(\omega) = |\ndot|$.
Let us begin with the following lemma:

\begin{lemm}\label{lem:strong:Minkowski:prop:trivial:value:case}
The adlic curve $S$ satisfies the strong Minkowski properties. More precisely,
if $E$ is a finite-dimensional vector space over $K$ and $\|\ndot\|$ is an ultrametric norm of $E$ over $(K,|\ndot|)$,
then $\nu_{\min}(E, \|\ndot\|) =
\widehat{\mu}_{\min}(E, \|\ndot\|)$.
\end{lemm}

\begin{proof}
In general, one has $\nu_{\min}(E, \|\ndot\|) \leqslant
\widehat{\mu}_{\min}(E, \|\ndot\|)$ by Proposition~\ref{Pro: comparaison mu min et nu min}, so that it is sufficient to show that
$\nu_{\min}(E, \|\ndot\|) \geqslant
\widehat{\mu}_{\min}(E, \|\ndot\|)$.
Let $(e_i)_{i=1}^r$ be an orthogonal basis of $E$ with respect to $\|\ndot\|$
(cf. Proposition~\ref{Pro:existenceepsorth}). 
Clearly we may assume that $\| e_r \| = \max \{ \|e_1\|, \ldots, \|e_r \| \}$.
Let $Q := E/(K e_1 + \cdots + Ke_{r-1})$ and $\|\ndot\|_Q$ be the quotient norm of $\|\ndot\|$ on $Q$.
Then $\|\pi(e_r)\|_Q = \| e_r \|$, where $\pi : E \to Q$ is the canonical homomorphism. 
Thus $-\log \| e_r \| \geqslant \widehat{\mu}_{\min}(E, \|\ndot\|)$, and hence
$\widehat{\deg}(e_i) \geqslant \widehat{\mu}_{\min}(E, \|\ndot\|)$ for all $i$. Therefore,
one has $\nu_{\min}(E, \|\ndot\|) \geqslant 
\widehat{\mu}_{\min}(E, \|\ndot\|)$.
\end{proof}

\begin{theo}
\label{them:Nakai:Moishezon:criterion:trivial:value:field}
Let $X$ be a projective variety over $K$, $L$ be a semiample invertible sheaf on $X$ and
$\varphi$ be a semipositive metric of $L$. We assume that, for every subvariety $Y$ of $X$,
$\rest{L}{Y}$ is big and there are a positive integer $m_Y$ and $\ell_Y \in H^0(Y, \rest{L}{Y}^{\otimes m_Y}) \setminus \{ 0 \}$
such that $\| \ell_Y \|_{m_Y \varphi_Y, \sup} < 1$. Then $L$ is ample and
there are a positive integer $n$ and a basis $s_1, \ldots, s_r$ of $H^0(X, L^{\otimes n})$ such that
$\| s_i \|_{n\varphi, \sup} < 1$ for all $i$.
\end{theo}

\begin{proof}
First of all, by using Nakai-Moishezon's criterion, one can see that $L$ is ample.
By our assumption, there are a positive integer $m_Y$ and $\ell_Y \in H^0(Y, \rest{L}{Y}^{\otimes m_Y}) \setminus \{ 0 \}$
such that $\| \ell_Y \|_{\rest{m_Y \varphi}{Y}, \sup} < 1$.
We choose $\epsilon > 0$ such that $e^{\epsilon} \| \ell_Y \|_{\rest{m_Y \varphi}{Y}, \sup} < 1$.
By Theorem~\ref{Thm:extensionpropertynontrivial},
there are a positive integer $n$ and $s \in H^0(X, L^{\otimes n m_Y}) \setminus \{ 0 \}$ such that
$\rest{s}{Y} = \ell_Y^{\otimes n}$ and $\|s \|_{nm_Y\varphi,\sup} \leqslant e^{n\epsilon} \|\ell_Y\|^n_{\rest{m_Y\varphi}{Y},\sup}$, so that $\|s \|_{nm_Y\varphi,\sup} < 1$, which means that $\| \ell_Y^{\otimes n} \|_{\xi_{Y, nm_Y}} < 1$.
Thus, by Theorem~\ref{thm:Nakai:Moishezon:criterion:adeic:curve} together with Lemma~\ref{lem:strong:Minkowski:prop:trivial:value:case},
the assertion follows.
\end{proof}
\fi

{
\section{Nakai-Moishezon's criterion over a number field}
Throughout this section, we fix a number field $K$ and the standard adelic curve $S =(K,(\Omega,\mathcal A,\nu),\phi)$ of $K$ as in Subsection~\ref{Subsec:Numberfields}.
Denote by $\Omega_{\mathrm{fin}}$ the set $\Omega\setminus\Omega_\infty$ of finite places of $K$, and by $\mathfrak o_K$ the ring of algebraic integers  in $K$. 
Note that} 
$S$ satisfies the strong Minkowski property (see \cite[Theorem 1.1]{chen17}).
Moreover, for $\omega \in \Omega_{\mathrm{fin}}$, the valuation ring of the completion $K_{\omega}$ of $K$ 
with respect to $\omega$ is denoted by $\mathfrak o_{\omega}$.

\if01
$\bullet$
For $\omega \in \Omega_{\mathrm{fin}}$, the valuation ring of the completion $K_{\omega}$ of $K$ 
with respect to $\omega$ is denoted by $\mathfrak o_{\omega}$.

$\bullet$
Let $\left(E, \xi = \{ \|\ndot\|_{\omega} \}_{\omega \in \Omega}\right)$ be an adelic vector bundle on $S$.
We always assume that $\xi$ is invariant by the complex conjugation structure (cf. Subsection~\ref{subsec:invariant:family:norms}), that is,
\[
\| x_1 \otimes_K^{\sigma} \lambda_1 + \cdots + x_l \otimes_K^{\sigma} \lambda_l \|_{\sigma} =
\| x_1 \otimes_K^{\overline{\sigma}} \overline{\lambda_1} + \cdots + x_l \otimes_K^{\overline{\sigma}} \overline{\lambda_l} \|_{\overline{\sigma}}
\]
for any $\sigma \in \Omega_{\infty,\mathrm c}$, $x_1, \ldots, x_l \in E$ and $\lambda_1, \ldots, \lambda_l \in K_{\sigma}$.
For convinience, if $\xi$ is invariant by the complex conjugation structure, then
$\xi$ is often said to be \emph{$F_{\infty}$-invariant}\index{F infty-invariant@$F_{\infty}$-invariant}.

$\bullet$
Let $X$ be a geometrically integral projective variety over $K$ and $L$ be an invertible sheaf on $X$.
For $\omega \in \Omega$, let $X_{\omega} := X \times_{\Spec K} K_{\omega}$ and
$L_{\omega} := L \otimes_{\mathcal O_X} \mathcal O_{X_{\omega}}$.
Let $\varphi_{\omega}$ be a continuous metric of $L_{\omega}$ on $X_{\omega}^{\mathrm{an}}$ for each $\omega \in \Omega$, and
$\varphi := \{ \varphi_{\omega} \}_{\omega \in \Omega}$. 
We always assume that the family of metrics $\{ \varphi_{\sigma} \}_{\sigma \in \Omega_{\infty}}$ on the Archimedean places is compatible with the complex conjugation, that is,
$F_{\sigma}^*(\varphi_{\bar{\sigma}}) = \varphi_{\sigma}$ for all $\sigma \in \Omega_{\infty, \mathrm c}$, 
where $F_{\sigma} : X_{\sigma} \to X_{\bar{\sigma}}$ is the complex conjugation map.
Note that the norm family $\{ \|\ndot\|_{n\varphi_{\omega}, \sup} \}_{\omega \in \Omega}$ of $H^0(X, L^{\otimes n})$
is $F_{\infty}$-invariant.
\fi

{
\subsection{Invariants $\bblambda$ and $\bbsigma$ for a graded algebra of adelic vector bundles}\label{subsection:inv:lambda:classical:settings:Nakai:Moishezon}
Let $\overline R_\sbullet=\{(R_n,\xi_n)\}_{n \in \mathbb Z_{\geqslant 0}}$ be a graded algebra of adelic vector bundles on $S$
such that $(R_n,\xi_n)$ is dominated and coherent 
for all $n\geqslant 0$. For the definition of the invariants $\bblambda$ and $\bbsigma$, see Subsection~\ref{subsection:inv:lambda:classical:settings}.

\begin{defi}
We say that $\overline R_\sbullet$ is \emph{asymptotically pure}\index{asymptotically pure@asymptotically pure} if 
\[ \limsup_{n\to\infty} \frac{\bbsigma(R_n, \xi_n)}{n} = 0. \]
\end{defi}

As a consequence of Proposition~\ref{prop:comp:lambda:nu:number:field}, we have the following:

\begin{prop}\label{prop:comp:lambda:nu:number:field:graded:ring}
One has the following inequalities:
\begin{multline*}
[K :\mathbb Q]\liminf_{n\to+\infty} \frac{\bblambda(R_n,\xi_n)}{n} \leqslant 
\liminf_{n\to+\infty}\frac{\nu_{\min}(R_n,\xi_n)}{n} \\
\leqslant 
[K :\mathbb Q]\liminf_{n\to+\infty}\frac{\bblambda(R_n,\xi_n)}{n} + 
\limsup_{n\to+\infty} \frac{\bbsigma(R_n,\xi_n)}{n}.
\end{multline*}
In particular, if $\overline R_\sbullet$ is asymptotically pure, then
\[
[K :\mathbb Q]\liminf_{n\to+\infty} \frac{\bblambda(R_n,\xi_n)}{n} = 
\liminf_{n\to+\infty}\frac{\nu_{\min}(R_n,\xi_n)}{n}.
\]
\end{prop}

\medskip
Let $\overline{M}_\sbullet=\{(M_n, \xi_{M_n}) \}_{n\in\mathbb Z}$ be a $h$-graded $\overline R_\sbullet$-module
such that $(M_n,\xi_n)$ is dominated and coherent 
for all $n\in \mathbb Z$.

\begin{prop}\label{prop:comp:lambda:R:M}
\begin{enumerate}[label=\rm(\arabic*)]
\item
If $R_\sbullet$ is generated by $R_1$ over $K$, then ${\displaystyle \lim_{n\to\infty} \frac{\bblambda(R_n, \xi_n)}{n}}$
exists in $\mathbb R \cup \{ \infty \}$.

\item
If $R_\sbullet = \bigoplus_{n=0}^{\infty} R_n$ is of finite type over $K$ and
$M_\sbullet = \bigoplus_{n \in \mathbb Z} M_n$ is finitely generated over $R_\sbullet$, then
\[
\frac{1}{h} \liminf_{n\to\infty} \frac{\bblambda(R_n,\xi_{n})}{n} \leqslant
\liminf_{n\to\infty} \frac{\bblambda(M_n,\xi_{M_n})}{n}.
\]
\end{enumerate}
\end{prop}

\begin{proof}
We set $\mathscr R_n :=  (R_n, \xi_{n})^{\mathrm{fin}}_{\leqslant 1}$ for $n \geqslant 0$, and
$\mathscr M_n := (M_n, \xi_{M_n})^{\mathrm{fin}}_{\leqslant 1}$ for $n \in \mathbb Z$.

\medskip
(1) For $\epsilon > 0$, we choose bases $e_1, \ldots, e_r$ and $e'_1, \ldots, e'_{r'}$ of
$R_n$ and $R_m$ over $K$, respectively, such that 
\[
\begin{cases}
e_1, \ldots, e_r \in \mathscr R_n, & \max \{ \| e_i \|_{\infty, n} \} \leqslant e^{-\bblambda(R_n, \xi_n) + \epsilon},\\
e'_1, \ldots, e'_{r'} \in \mathscr R_m, & \max \{ \| e'_j \|_{\infty, n} \} \leqslant e^{-\bblambda(R_m, \xi_m) + \epsilon}.
\end{cases}
\]
Then $e_i e'_j \in \mathscr R_{n+m}$ and
$\max \{ \| e_i e'_j \|_{\infty, n+m} \} \leqslant e^{-\bblambda(R_n, \xi_n) -\bblambda(R_n, \xi_n) + 2\epsilon}$.
Note that $\{ e_i e'_j \}$ forms generators of $R_{n+m}$ over $K$ because $R_n \otimes R_m \to R_{n+m}$ is surjective,
so that $e^{-\bblambda(R_{n+m}, \xi_{n+m})} \leqslant e^{-\bblambda(R_n, \xi_n) -\bblambda(R_n, \xi_n) + 2\epsilon}$.
Therefore, one has
\[
\bblambda(R_{n+m}, \xi_{n+m}) \geqslant \bblambda(R_n, \xi_n) + \bblambda(R_n, \xi_n)
\]
for all $n, m$. Thus the assertion follows from Fekete's lemma.

\medskip
(2) It can be proved in the similar way as in Proposition~\ref{Pro: estimate of asymptotique minimal slope}.
First we assume that $R_{\sbullet}$ is generated by $R_1$ over $K$. Then
there exist integers $b_1$ and $m>0$ such that, for any integer $b$ with $b\geqslant b_1$ and any integer $\ell\geqslant 1$ the canonical $K$-linear map $R_{\ell m}\otimes_KM_{b}\rightarrow M_{b+\ell mh}$ is surjective. 
For $\epsilon > 0$, we choose a basis $e_1, \ldots, e_r$ of $R_{\ell m}$ and a basis $m_1, \ldots, m_{r'}$ of $M_{b}$ such that
$e_1, \ldots, e_r \in \mathscr R_{\ell m}$, $m_1, \ldots, m_{r'} \in \mathscr M_{b}$,
$\max \{ \| e_i \|_{\infty, \ell m} \} \leqslant e^{-\bblambda(R_{\ell m}, \xi_{\ell m}) + \epsilon}$ and
$\max \{ \| m_j \|_{\infty, M_{b}} \} \leqslant e^{-\bblambda(M_{b}, \xi_{M_{b}}) + \epsilon}$.
Note that $e_i m_j \in \mathscr M_{b+\ell mh}$ and
\[
\| e_i m_j \|_{\infty, M_{b+\ell mh}} \leqslant e^{-\bblambda(R_{\ell m}, \xi_{\ell m}) -\bblambda(M_{b}, \xi_{M_{b}}) + 2\epsilon}.
\]
Moreover we can find a basis of $M_{b+\ell mh}$ among $\{ e_i m_j \}_{1 \leqslant i \leqslant r, 1 \leqslant j \leqslant r'}$, so that
\[
e^{-\bblambda(M_{b+\ell mh}, \xi_{M_{b+\ell mh}})} \leqslant e^{-\bblambda(R_{\ell m}, \xi_{\ell m}) -\bblambda(M_{b}, \xi_{M_{b}}) + 2\epsilon},
\]
and hence one has
\[
\bblambda(M_{b+\ell mh}, \xi_{M_{b+\ell mh}}) \geqslant \bblambda(R_{\ell m}, \xi_{\ell m}) + \bblambda(M_{b}, \xi_{M_{b}}).
\]
Therefore,
\[
\liminf_{l\to\infty} \frac{\bblambda(M_{b+\ell mh}, \xi_{M_{b+\ell mh}})}{\ell mh} \geqslant
\frac{1}{h} \liminf_{l\to\infty} \frac{\bblambda(R_{\ell m}, \xi_{\ell m})}{\ell m} \geqslant
\frac{1}{h} \liminf_{n\to\infty} \frac{\bblambda(R_{n}, \xi_{n})}{n},
\]
which implies 
\[
\liminf_{n\to\infty} \frac{\bblambda(M_{n}, \xi_{M_{n}})}{n} \geqslant
\frac{1}{h} \liminf_{n\to\infty} \frac{\bblambda(R_{n}, \xi_{n})}{n}
\]
because $b \geqslant b_1$ is arbitrary.

In general, we can find a positive integer $u$ such that
$R^{(u)}_{\sbullet} := \bigoplus_{n=0}^{\infty} R_{un}$ is generated by $R_1^{(u)} = R_u$ over $K$.
Note that $R_{\sbullet}$ is a finitely generated $R^{(u)}_{\sbullet}$-module. Therefore, by the previous observation, 
one has
\[
\begin{cases}
{\displaystyle \liminf_{n\to\infty} \frac{\bblambda(R_{n}, \xi_{n})}{n} \geqslant \frac{1}{u} \liminf_{n\to\infty} \frac{\bblambda(R_{un}, \xi_{un})}{n}}, \\[2ex]
{\displaystyle \liminf_{n\to\infty} \frac{\bblambda(M_{n}, \xi_{M_{n}})}{n} \geqslant \frac{1}{hu} \liminf_{n\to\infty} \frac{\bblambda(R_{un}, \xi_{un})}{n}}.
\end{cases}
\]
Moreover, as 
\[
\liminf_{n\to\infty} \frac{\bblambda(R_{un}, \xi_{un})}{un} \geqslant \liminf_{m\to\infty} \frac{\bblambda(R_{m}, \xi_{m})}{m},
\]
one obtains
\[
\liminf_{n\to\infty} \frac{\bblambda(R_{n}, \xi_{n})}{n} = \frac{1}{u} \liminf_{n\to\infty} \frac{\bblambda(R_{un}, \xi_{un})}{n}.
\]
Thus the assertion follows.
\end{proof}
}

{
\subsection{Dominancy and coherency of generically pure metric}
Let $X$ be a geometrically integral projective variety over $K$ and
$L$ be an invertible sheaf on $X$.
For $\omega \in \Omega$, let $X_{\omega} := X \times_{\Spec K} K_{\omega}$ and
$L_{\omega} := L \otimes_{\mathcal O_X} \mathcal O_{X_{\omega}}$.
Let $\varphi_{\omega}$ be a continuous metric of $L_{\omega}$ on $X_{\omega}^{\mathrm{an}}$ for each $\omega \in \Omega$, and
$\varphi := \{ \varphi_{\omega} \}_{\omega \in \Omega}$.

Let us begin with the definition of the generic purity of the metric family $\varphi$.

\begin{defi}\label{def:generic:purity}
We say that $\varphi$ is \emph{generically pure}\index{generically pure@generically pure}
if there are a non-empty Zariski open set $U$ of $\Spec(\mathfrak o_K)$,
a projective integral scheme $\mathscr X$ over $U$  
and an invertible $\mathcal O_{\mathscr X}$-module $\mathscr L$ such that $\mathscr X \times_{U} \Spec(K) = X$, $\rest{\mathscr L}{X} = L$ and,
for each $\omega \in U \cap \Omega_{\mathrm{fin}}$,
$\varphi_{\omega}$ coincides with the metric arising from $\mathscr X_{\omega}$ and $\mathscr L_{\omega}$,
where $\mathscr X_{\omega} = \mathscr X \times_{U} \Spec(\mathfrak o_{\omega})$ and
$\mathscr L_{\omega}$ is the pull-back of $\mathscr L$ to $\mathscr X_{\omega}$.
\end{defi}

{
\begin{prop}\label{prop:generic:purity:dominant:coherent}
\begin{enumerate}
\renewcommand{\labelenumi}{(\arabic{enumi})}
\item
If $L$ is generated by global sections
and $\varphi$ is generically pure, there exist a non-empty Zariski open set $U$ of
$\Spec(\mathfrak o_K)$ and a basis $\pmb{e} = (e_i)_{i=1}^r$ of $H^0(X, L)$ such that
$\varphi_{\omega} = \varphi_{\pmb{e}, \omega}$ for all $\omega \in U \cap \Omega_{\mathrm{fin}}$.

\item
If $L$ is semiample  
and $\varphi$ is generically pure, then $\varphi$ is dominated and
$\left(H^0(X, L^{\otimes n}), \{ \|\ndot\|_{n\varphi_{\omega}, \sup} \}_{\omega \in \Omega}\right)$
is coherent for all $n \geqslant 0$.
\end{enumerate}
\end{prop}

\begin{proof}
(1): We use the notation in Definition~\ref{def:generic:purity}.
Shrinking $U$ if necessarily, we may assume that
$H^0(\mathscr X, \mathscr L)$ is a free $\mathfrak o_U$-module and
$H^0(\mathscr X, \mathscr L) \otimes_{\mathfrak o_U} \mathcal O_{\mathscr X} \to
\mathscr L$ is surjective.
Let $\pmb{e} = (e_i)_{i=1}^r$ be a free basis of
$H^0(\mathscr X, \mathscr L)$ over $\mathfrak o_U$.
Then $(e_i)_{i=1}^r$ yields a free basis of $H^0(\mathscr X_{\omega}, \mathscr L_{\omega})$ over $\mathfrak o_{\omega}$
for any $\omega \in U \cap \Omega_{\mathrm{fin}}$.
Let $\|\ndot\|_{H^0(\mathscr X_{\omega}, \mathscr L_{\omega})}$ be the norm of $H^0(X_{\omega}, L_{\omega}^{\otimes r})$
arising from the lattice $H^0(\mathscr X_{\omega}, \mathscr L_{\omega})$. Then, by Proposition~\ref{Pro:orthogonallattice}, $\|\ndot\|_{H^0(\mathscr X_{\omega}, \mathscr L_{\omega})} = \|\ndot\|_{\pmb{e},\omega}$ for any $\omega \in U \cap \Omega_{\mathrm{fin}}$.
Moreover, $H^0(\mathscr X_{\omega}, \mathscr L_{\omega}) \otimes_{\mathfrak o_{\omega}} \mathcal O_{\mathscr X_{\omega}}\to \mathscr L_{\omega}$
is surjective, so that, by Proposition~\ref{Pro:modelFubiniStudy}, one has $\varphi_{\omega} = \varphi_{\pmb{e}, \omega}$, as required.

\medskip
(2) We choose a positive integer $m$ such that $L^{\otimes m}$ is generated by global sections and 
$\alpha_n : H^0(X, L^{\otimes m})^{\otimes n} \to H^0(X, L^{\otimes nm})$
is surjective for all $n \geqslant 1$.
Then, by (1), there are a non-empty Zariski open set $U$ of
$\Spec(\mathfrak o_K)$ and a basis $\pmb{e} = (e_i)_{i=1}^r$ of $H^0(X, L^{\otimes m})$ such that
$m \varphi_{\omega} = \varphi_{\pmb{e}, \omega}$ for all $\omega \in U \cap \Omega_{\mathrm{fin}}$. In particular, $m \varphi$ is dominated, so that
$\varphi$ is also dominated by Proposition~\ref{Pro: dominancy preserved by operators}.

For $\omega \in U \cap \Omega_{\mathrm{fin}}$, let $\|\ndot\|^{\otimes n}_{\pmb{e},\omega}$ be the $\varepsilon$-tensor products of $\|\ndot\|_{\pmb{e},\omega}$ on $H^0(X_{\omega}, L^{\otimes m}_{\omega})^{\otimes n}$. Note that,
by Proposition~\ref{Pro:alphatenso} together with \eqref{Equ: split tensor} in Remark~\ref{Rem:produittenrk1},
\[
\forall\, a_{i_1,\ldots,i_N} \in K_{\omega},\quad
\Bigg\| \sum_{\substack{(i_1, \ldots, i_r) \in \mathbb Z_{\geqslant 0}^r,\\ i_1 + \cdots + i_r = n}}  a_{i_1, \ldots, i_r} e_1^{\otimes i_1} \otimes \cdots \otimes e_{r}^{\otimes i_r}\Bigg\|^{\otimes n}_{\pmb{e},\omega} =
\max \{ |a_{i_1, \ldots, i_r}|_{\omega}\}.
\]
Moreover, by Remark~\ref{Rem: powers of quotient}, $n\varphi_{\pmb{e},\omega}$ coincides with
the quotient metric induced by the surjective homomorphism 
$H^0(X_{\omega}, L^{\otimes m}_{\omega})^{\otimes n} \otimes_{K_{\omega}} \mathcal O_{X_{\omega}}
\to L^{\otimes nm}_{\omega}$
and $\|\ndot\|^{\otimes n}_{\pmb{e},\omega}$.

Fix $s \in H^0(X, L^{\otimes n})$ ($n \geqslant 1$). Then $s^{\otimes m} \in H^0(X, L^{\otimes mn})$. As $\alpha_n$ is surjective, one can choose 
\[
f = \sum_{\substack{(i_1, \ldots, i_r) \in \mathbb Z_{\geqslant 0}^r,\\ i_1 + \cdots + i_r = n}}  f_{i_1, \ldots, i_r} e_1^{\otimes i_1} \otimes \cdots \otimes e_{r}^{\otimes i_r} \in H^0(X, L^{\otimes m})^{\otimes n}
\quad(f_{i_1, \ldots, i_r} \in K)
\]
such that $\alpha_n(f) = s^{\otimes m}$.
Then, by Proposition~\ref{Pro:positivityofquotientmetric},
\[
\big(\| s\|_{n\varphi_{\omega}, \sup}\big)^m = \| s^{\otimes m} \|_{nm\varphi_{\omega}, \sup} =
\| s^{\otimes m} \|_{n\varphi_{\pmb{e},\omega}, \sup} \leqslant \|f \|^{\otimes n}_{\pmb{e},\omega}
= \max \{ |f_{i_1, \ldots, i_r}|_{\omega}\},
\]
for all $\omega \in U \cap \Omega_{\mathrm{fin}}$, so that
so that $\| s\|_{n\varphi_{\omega}, \sup} \leqslant 1$ for all $\omega \in \Omega$
except finitely many $\omega$ because 
$\Omega \setminus (U \cap \Omega_{\mathrm{fin}})$ is finite and
$|f_{i_1, \ldots, i_r}|_{\omega} = 1$ for all
$i_1, \ldots, i_r$ and $\omega \in \Omega_{\mathrm{fin}}$ except finitely many $\omega$.
\end{proof}
}
\if01
As $X$ is geometrically integral, by \cite[Th\'{e}or\`{e}me~9.7.7]{EGAIV-3},
shrinking $U$ if necessarily, we may assume that,
for any $\omega \in U \cap \Omega_{\mathrm{fin}}$, the fiber of $\mathscr X \to U$ over $\omega$ is
geometrically integral over the residue field at $\omega$. Moreover, we may also assume that
$\mathscr L$ is ample over $U$. Note that the central fiber $\mathscr X_{\omega, o}$ of $\mathscr X_{\omega} \to \Spec(\mathfrak o_{\omega})$ coincides with the fiber of $\mathscr X \to U$ over $\omega$, so that,
by Proposition~\ref{Pro:latticeby model}, 
\[
\{ x \in R_n \mid \| x \|_{n\varphi_{\omega}, \sup} \leqslant 1 \} =
H^0(\mathscr X_{\omega}, \mathscr L_{\omega}^{\otimes n}) = H^0(\mathscr X, \mathscr L^{\otimes n}) \otimes_{\mathfrak o_U} \mathfrak o_{\omega}
\]
and $\|\ndot\|_{n\varphi_{\omega},\sup} = \|\ndot\|_{H^0(\mathscr X_{\omega}, \mathscr L_{\omega}^{\otimes n})}$ for all $n \geqslant 1$ and $\omega \in U \cap \Omega_{\mathrm{fin}}$.
In particular, by Proposition~\ref{prop:strongly:cohenent},
$\left(H^0(X, L^{\otimes n}), \{ \|\ndot\|_{n\varphi_{\omega}, \sup} \}_{\omega \in \Omega}\right)$
is dominated and coherent for all $n \geqslant 1$.
Note that $n\varphi_{\omega}$ coincides with
the metric arised from $\mathscr X_{\omega}$ and $\mathscr L_{\omega}^{\otimes n}$ (cf. Proposition~\ref{Pro:tensormodel}).
We choose a positive integer $n$ such that $\mathscr L^{\otimes n}$ is very ample over $U$.
For $\omega \in U \cap \Omega_{\mathrm{fin}}$, as $H^0(\mathscr X_{\omega}, \mathscr L_{\omega}^{\otimes n}) \otimes \mathcal O_{\mathscr X} \to \mathscr L_{\omega}^{\otimes n}$ is
surjective, by Proposition~\ref{Pro:modelFubiniStudy},
$n\varphi_{\omega}$ is equal to
the quotient metric induced by $\|\ndot\|_{H^0(\mathscr X_{\omega}, \mathscr L_{\omega}^{\otimes n})} = \|\ndot\|_{n\varphi_{\omega},\sup}$ and
$H^0(X_{\omega}, L_{\omega}^{\otimes n}) \otimes \mathcal O_{X_{\omega}} \to L_{\omega}^{\otimes n}$, which shows that $n\varphi$ is dominated.
Therefore, by Proposition~\ref{Pro: dominancy preserved by operators},
$\varphi$ is dominated.
\fi
\if01
\begin{enonce}{Claim}
For $\omega \in U \cap \Omega_{\mathrm{fin}}$,
$n\varphi_{\omega}$ coincides with
the metric induced by 
$\|\ndot\|_{H^0(\mathscr X_{\omega}, \mathscr L_{\omega}^{\otimes n})}$ and the surjection 
$H^0(X_{\omega}, L_{\omega}^{\otimes n}) \otimes \mathcal O_{X_{\omega}} \to L_{\omega}^{\otimes}$.
\end{enonce}

\begin{proof}
Let $\psi_{\omega}$ be the metric of $\mathcal O_{\mathbb P(H^0(X_{\omega},L_{\omega}^{\otimes n}))}(1)$
given by the model 
\[
\left(\mathbb P(H^0(\mathscr X_{\omega}, \mathscr L_{\omega}^{\otimes n})), \mathcal O_{\mathbb P(H^0(\mathscr X_{\omega}, \mathscr L_{\omega}^{\otimes n}))}(1)\right).
\]
As 
$\mathscr L_{\omega}^{\otimes n} = \rest{ \mathcal O_{\mathbb P(H^0(\mathscr X_{\omega}, \mathscr L_{\omega}^{\otimes n}))}(1)}{\mathscr X_{\omega}}$, one has 
\begin{equation}\label{eqn:prop:generic:purity:dominant:coherent:01}
\rest{\psi_{\omega}}{X_{\omega}^{\mathrm{an}}} = n\varphi_{\omega}.
\end{equation}
Let $\psi'_{\omega}$ be the quotient metric of $\mathcal O_{\mathbb P(H^0(X_{\omega},L_{\omega}^{\otimes n}))}(1)$ induced by
$\|\ndot\|_{H^0(\mathscr X_{\omega}, \mathscr L_{\omega}^{\otimes n})}$ and the surjection 
$H^0(X_{\omega}, L_{\omega}^{\otimes n}) \otimes O_{\mathbb P(H^0(X_{\omega}, L_{\omega}^{\otimes n}))} \to \mathcal O_{\mathbb P(H^0(X_{\omega}, L_{\omega}^{\otimes n}))}(1)$. Then, by
Proposition~\ref{Pro:modelFubiniStudy}, 
\begin{equation}\label{eqn:prop:generic:purity:dominant:coherent:02}
\psi'_{\omega} = \psi_{\omega}.
\end{equation} 
Obviously, $\rest{\psi'_{\omega}}{X_{\omega}^{\mathrm{an}}}$
is equal to the metric induced by 
$\|\ndot\|_{H^0(\mathscr X_{\omega}, \mathscr L_{\omega}^{\otimes n})}$ and the surjection 
$H^0(X_{\omega}, L_{\omega}^{\otimes n}) \otimes \mathcal O_{X_{\omega}} \to L_{\omega}^{\otimes}$.
Thus the assertion of the claim follows from \eqref{eqn:prop:generic:purity:dominant:coherent:01} and
\eqref{eqn:prop:generic:purity:dominant:coherent:02}.
\end{proof}

The above claim shows that $n\varphi$ is dominated.
Therefore, by Proposition~\ref{Pro: dominancy preserved by operators},
$\varphi$ is dominated
\fi

\subsection{Fine metric family}

Let $X$ be a geometrically integral projective variety over $K$ and
$L$ be an invertible sheaf on $X$.
For $\omega \in \Omega$, let $X_{\omega} := X \times_{\Spec K} K_{\omega}$ and
$L_{\omega} := L \otimes_{\mathcal O_X} \mathcal O_{X_{\omega}}$.
Let $\varphi_{\omega}$ be a continuous metric of $L_{\omega}$ on $X_{\omega}^{\mathrm{an}}$ for each $\omega \in \Omega$, and
$\varphi := \{ \varphi_{\omega} \}_{\omega \in \Omega}$.

\begin{defi}
We say that $\varphi$ is \emph{very fine}\index{very fine@very fine} if $\varphi$ is dominated and
there are a generically pure continuous metric family $\varphi' = \{ \varphi'_{\omega} \}_{\omega \in \Omega}$ of $L$ and a non-empty Zariski open set $U$ of
$\Spec(\mathfrak o_K)$ such that
$|\ndot|_{\varphi_{\omega}} \leqslant |\ndot|_{\varphi'_{\omega}}$ for all $\omega \in U \cap \Omega_{\mathrm{fin}}$.
Further, $\varphi$ is said to be \emph{fine}\index{fine@fine} if $r\varphi$ is very fine for some positive integer $r$.
\end{defi}

{
\begin{prop}\label{prop:fine:coherent}
Let $L$ and $M$ be invertible $\mathcal O_X$-module, and $\varphi$ and $\psi$ be continuous metric families of $L$ and $M$,
respectively.
\begin{enumerate}[label=\rm(\arabic*)]
\item
If $\varphi$ and $\psi$ are very fine, then $\varphi + \psi$ is very fine.

\item
If $\varphi$ and $\psi$ are fine, then $\varphi + \psi$ is fine.

\item
If $a\varphi$ is fine for some positive integer $a$, then $\varphi$ is fine.

\item
If $\varphi$ is fine, then $\left(H^0(X, L), \{ \|\ndot\|_{\varphi_{\omega}, \sup} \}_{\omega \in \Omega}\right)$ is
coherent.
\end{enumerate}
\end{prop}

\begin{proof}
(1) is obvious.

(2) We choose positive integers $r$ and $r'$ such that $r\varphi$ and $r'\psi$ are very fine.
Then, by (1), $rr'\varphi$ and $rr'\psi$ are very fine, so that
$rr'(\varphi+\psi)$ is very fine, as required.

(3) Since $a\varphi$ is fine, there is a positive integer $r$ such that $ra\varphi$ is very fine, so that
$\varphi$ is fine.

(4) Let $r$ be a positive integer such that $r\varphi$ is very fine.
Then there are a generically pure continuous metric family $\varphi' = \{ \varphi'_{\omega} \}_{\omega \in \Omega}$ of $L^{\otimes r}$ and a non-empty Zariski open set $U$ of
$\Spec(\mathfrak o_K)$ such that
$|\ndot|_{r\varphi_{\omega}} \leqslant |\ndot|_{\varphi'_{\omega}}$ for all $\omega \in U \cap \Omega_{\mathrm{fin}}$, so that, 
for $s \in H^0(X, L) \setminus \{ 0 \}$ and $\omega \in U \cap \Omega_{\mathrm{fin}}$,
$\| s^{\otimes r} \|_{r\varphi_{\omega}, \sup} \leqslant \| s^{\otimes r} \|_{\varphi'_{\omega},\sup}$.
By Proposition~\ref{prop:generic:purity:dominant:coherent}, $\| s^{\otimes r} \|_{\varphi'_{\omega},\sup} \leqslant 1$
expect finitely many $\omega \in \Omega$. Therefore, the same assertion holds for $\| s^{\otimes r} \|_{r\varphi_{\omega}, \sup}$. Note that $\| s^{\otimes r} \|_{r\varphi_{\omega}, \sup} = \| s \|_{\varphi_{\omega}, \sup}^r$, and hence
$\| s \|_{\varphi_{\omega},\sup} \leqslant 1$ expect finitely many $\omega \in \Omega$.
\end{proof}

For $n \geqslant 0$, we set $R_n := H^0(X, L^{\otimes n})$ and $\xi_n := \{ \|\ndot\|_{n\varphi_{\omega}, \sup} \}_{\omega \in \Omega}$. Note that $\overline R_{\sbullet} = \{ (R_n, \xi_n) \}_{n=0}^{\infty}$ forms
a graded algebra of adelic vector bundles over $S$ (cf. Definition~\ref{def:graded:algebra:adelic:vector:bundles}).
For $n \geqslant 0$ and $\omega \in \Omega$,
we denote $R_n \otimes_K K_{\omega} = H^0(X_{\omega}, L_{\omega}^{\otimes n})$ by $R_{n,\omega}$.
Moreover, for $n \geqslant 0$ and $\omega \in \Omega_{\mathrm{fin}}$, we set 
$\mathscr R_{n, \omega} := \{ x \in R_{n,\omega} \mid \| x \|_{n,\omega} \leqslant 1 \}$.
Note that $\mathscr R_{n, \omega}$ is a locally free $\mathfrak o_{\omega}$-module and
$\mathscr R_{n, \omega} \otimes_{\mathfrak o_{\omega}} K_{\omega} = R_{n,\omega}$
(cf. Proposition~\ref{prop:ultrametric:ball:lattice} and Proposition~\ref{Pro:normetreausauxdisc}). 
Further we set 
\[
\mathscr R_n := \{ x \in R_n \mid \text{$\| x \|_{n,\omega} \leqslant 1$ for all $\omega \in \Omega_{\mathrm{fin}}$} \}.
\]
If $(R_n,\xi_n)$ is dominated and coherent, then, 
by Proposition~\ref{prop:finiteness:cond:classical:setting} and Proposition~\ref{prop:finite:generation:dominated},
$\mathscr R_n$ is finitely generated over $\mathfrak o_K$, 
$\mathscr R_n \otimes_{\mathfrak o_K} K = R_n$,
$\mathscr R_n \otimes_{\mathfrak o_K} K_{\omega} = R_{n,\omega}$ and
$\mathscr R_n \otimes_{\mathfrak o_K} \mathfrak o_{\omega} = \mathscr R_{n, \omega}$ for all $\omega \in \Omega_{\mathrm{fin}}$.

\begin{prop}\label{prop:equiv:coherent}
We assume that $L$ is ample and $\varphi$ is dominated. Then the following are equivalent:
\begin{enumerate}[label=\rm(\arabic*)]
\item
The metric family $\varphi$ is fine.

\if01
{\color{green}
\item There exist a positive integer $r$, a basis $(e_i)_{i=1}^N$ of $H^0(X, L^{\otimes r})$ and a non-empty open
set $U$ of $\Spec(\mathfrak o_K)$ with the following properties:
\begin{enumerate}
\renewcommand{\labelenumii}{(\arabic{enumi}.\arabic{enumii})}
\item
$L^{\otimes r}$ is generated by global sections and $H^0(X, L^{\otimes r})^{\otimes n} \to H^0(X, L^{\otimes nr})$
is surjective for every $n \geqslant 1$.

\item
Let $\|\ndot\|'_{\omega}$ be the norm of $H^0(X_{\omega}, L_{\omega}^{\otimes r})$ given by
\[
\hskip4em \forall\, a_1, \ldots, a_N \in K_{\omega},\quad
\| a_1 e_1 + \cdots + a_N e_N \|'_{\omega} := \max \{ |a_1|_{\omega}, \ldots, |a_{N}|_{\omega} \}.
\]
If $\phi'_{\omega}$ is the metric of $L^{\otimes r}_{\omega}$ induced by 
$H^0(X_{\omega}, L_{\omega}^{\otimes r}) \otimes_{K_{\omega}}
\mathcal O_{X_{\omega}} \to L_{\omega}^{\otimes r}$
and the norm $\|\ndot\|'_{\omega}$,
then $r\varphi_{\omega} \leqslant \phi'_{\omega}$ for all $\omega \in U \cap \Omega_{\mathrm{fin}}$.
\end{enumerate}
}
\fi

\item
$( R_n, \xi_n )$ is
coherent for all $n \geqslant 0$.
\end{enumerate}
\end{prop}

\begin{proof}

(1) $\Longrightarrow$ (2): This is a consequence of Proposition~\ref{prop:fine:coherent}.

\if01
We choose a positive integer $r$ such that $r\varphi$ is very fine, $L^{\otimes r}$ is generated by global sections and
$H^0(X, L^{\otimes r})^{\otimes n} \to H^0(X, L^{\otimes nr})$
is surjective for every $n \geqslant 1$. As $r\varphi$ is very fine,
there are a non-empty Zariski open set $U$ of $\Spec(\mathfrak o_K)$,
a projective integral scheme $\mathscr X$ over $U$  
and an invertible $\mathcal O_{\mathscr X}$-module $\mathscr L$ such that $\mathscr X \times_{U} \Spec(K) = X$, $\rest{\mathscr L}{X} = L^{\otimes r}$ and,
for each $\omega \in U \cap \Omega_{\mathrm{fin}}$,
if $\varphi'_{\omega}$ is the metric of $L_{\omega}^{\otimes r}$ arising  from $\mathscr X_{\omega}$ and $\mathscr L_{\omega}$,
then $r\varphi_{\omega} \leqslant \varphi'_{\omega}$,
where $\mathscr X_{\omega} = \mathscr X \times_{U} \Spec(\mathfrak o_{\omega})$ and
$\mathscr L_{\omega}$ is the pull-back of $\mathscr L$ to $\mathscr X_{\omega}$.
Shrinking $U$ if necessarily, we may assume that
$H^0(\mathscr X, \mathscr L) \otimes_{\mathfrak o_U} \mathcal O_{\mathscr X} \to \mathscr L$ is surjective and
$H^0(\mathscr X, \mathscr L)$ is a free $\mathfrak o_U$-module.
Let $(e_i)_{i=1}^N$ be a free basis of $H^0(\mathscr X, \mathscr L)$ over $\mathfrak o_U$.
Then $(e_i)_{i=1}^N$ yields a free basis of $H^0(\mathscr X_{\omega}, \mathscr L_{\omega})$ over $\mathfrak o_{\omega}$
for any $\omega \in U \cap \Omega_{\mathrm{fin}}$.
Let $\|\ndot\|_{H^0(\mathscr X_{\omega}, \mathscr L_{\omega})}$ be the norm of $H^0(X_{\omega}, L_{\omega}^{\otimes r})$
arising from the lattice $H^0(\mathscr X_{\omega}, \mathscr L_{\omega})$. Then, by Proposition~\ref{Pro:orthogonallattice},
\[
\forall\, a_1, \ldots, a_N \in K_{\omega},\quad
\| a_1 e_1 + \cdots + a_N e_N \|_{H^0(\mathscr X_{\omega}, \mathscr L_{\omega})} := \max \{ |a_1|_{\omega}, \ldots, |a_{N}|_{\omega} \}.
\]
Let $\phi'_{\omega}$ be the quotient metric induced by
$\|\ndot\|_{H^0(\mathscr X_{\omega}, \mathscr L_{\omega})}$ and $H^0(X_{\omega}, L_{\omega}^{\otimes r}) \otimes \mathcal O_{X_{\omega}} \to L_{\omega}^{\otimes r}$.
Then as $H^0(\mathscr X_{\omega}, \mathscr L_{\omega}) \otimes_{\mathfrak o_{\omega}} \mathcal O_{\mathscr X_{\omega}}\to \mathscr L_{\omega}$
is surjective, by Proposition~\ref{Pro:modelFubiniStudy}, one has $\varphi'_{\omega} = \phi'_{\omega}$, as required.

\medskip
(2) $\Longrightarrow$ (3):
For $\omega \in U \cap \Omega_{\mathrm{fin}}$, let $\|\ndot\|'_{n,\varepsilon,\omega}$ be the $\varepsilon$-tensor products of $\|\ndot\|'_{\omega}$ on $H^0(X_{\omega r}, L^{\otimes r}_{\omega})^{\otimes n}$. Note that,
by Proposition~\ref{Pro:alphatenso} together with \eqref{Equ: split tensor} in Remark~\ref{Rem:produittenrk1},
\[
\forall\, a_{i_1,\ldots,i_N} \in K_{\omega},\quad
\Bigg\| \sum_{\substack{(i_1, \ldots, i_N) \in \mathbb Z_{\geqslant 0}^N,\\ i_1 + \cdots + i_N = n}}  a_{i_1, \ldots, i_N} e_1^{\otimes i_1} \otimes \cdots \otimes e_{N}^{i_N}\Bigg\|'_{n,\varepsilon,\omega} =
\max \{ |a_{i_1, \ldots, i_N}|_{\omega}\}.
\]
Moreover, by Remark~\ref{Rem: powers of quotient}, $n\phi'_{\omega}$ coincides with
the quotient metric induced by the surjective homomorphism 
$H^0(X_{\omega}, L^{\otimes r}_{\omega})^{\otimes n} \otimes_{K_{\omega}} \mathcal O_{X_{\omega}}
\to L^{\otimes rn}_{\omega}$
and $\|\ndot\|'_{n,\varepsilon,\omega}$.

Fix $s \in H^0(X, L^{\otimes n})$ ($n \geqslant 1$).
Then $s^{\otimes r} \in H^0(X, L^{\otimes nr})$. As $\|s^{\otimes r}\|_{nr\varphi_{\omega}, \sup} =
\|s\|_{n\varphi_{\omega}, \sup}^r$ and $\Omega \setminus (U \cap \Omega_{\mathrm{fin}})$ is finite,
it is sufficient to show that $\|s^{\otimes r}\|_{nr\varphi_{\omega}, \sup} \leqslant 1$ for all $\omega \in U \cap \Omega_{\mathrm{fin}}$ except finitely many $\omega$. We choose 
\[
f = \sum_{\substack{(i_1, \ldots, i_N) \in \mathbb Z_{\geqslant 0}^N,\\ i_1 + \cdots + i_N = n}}  a_{i_1, \ldots, i_N} e_1^{\otimes i_1} \otimes \cdots \otimes e_{N}^{i_N} \in H^0(X, L^{\otimes r})^{\otimes n}
\quad(a_{i_1, \ldots, i_N} \in K)
\]
such that the image of $f$ in $H^0(X, L^{\otimes rn})$ is $s^{\otimes r}$.
Then, by Proposition~\ref{Pro:positivityofquotientmetric},
\[
\| s^{\otimes r} \|_{nr\varphi_{\omega}, \sup} \leqslant 
\| s^{\otimes r} \|_{n\phi'_{\omega}, \sup} \leqslant \|f \|'_{n,\varepsilon,\omega}
= \max \{ |a_{i_1, \ldots, i_N}|_{\omega}\},
\]
so that the assertion follows because $|a_{i_1, \ldots, i_N}|_{\omega} \leqslant 1$ for all
$i_1, \ldots, i_N$ and $\omega \in U \cap \Omega_{\mathrm{fin}}$ except finitely many $\omega$.
}
\fi

\medskip
(2) $\Longrightarrow$ (1):
First note that $(R_n, \xi_n)$ is dominated for $n \geqslant 0$ by
Proposition~\ref{Pro: dominancy preserved by operators} and Theorem~\ref{Thm: Fubini-Study dominated}.
Moreover, by our assumption, $(R_n, \xi_n)$ is coherent for every $n \geqslant 0$.

Let $r$ be a positive integer such that $L^{\otimes r}$ is very ample.
Let $\mathscr X$ be the Zariski closure of $X$ in $\mathbb P(\mathscr R_r)$ and
$\mathscr L = \rest{\mathcal O_{\mathbb P(\mathscr R_r)}(1)}{\mathscr X}$.
Then $\rest{\mathscr L}{X} = L^{\otimes r}$. Moreover,
since $\mathscr R_r \otimes_{\mathfrak o_K} \mathcal O_{\mathbb P(\mathscr R_r)} \to
\mathcal O_{\mathbb P(\mathscr R_r)}(1)$ is surjective,
$\mathscr R_r \otimes_{\mathfrak o_K} \mathcal O_{\mathscr X} \to
\mathscr L$ is also surjective.
For each $\omega \in \Omega_{\mathrm{fin}}$,
let $\psi_{\omega}$ be the metric of $L_{\omega}^{\otimes r}$ 
arising from $\mathscr X_{\omega}$ and $\mathscr L_{\omega}$,
where $\mathscr X_{\omega} = \mathscr X \times_{\Spec(\mathfrak o_K)} \Spec(\mathfrak o_{\omega})$ and
$\mathscr L_{\omega}$ is the pull-back of $\mathscr L$ to $\mathscr X_{\omega}$.
Let $\varphi' = \{ \varphi'_{\omega} \}_{\omega \in \Omega}$ be the metric family of $L^{\otimes r}$ given by
\[
\varphi'_{\omega} := \begin{cases}
\psi_{\omega} & \text{if $\omega \in \Omega_{\mathrm{fin}}$}, \\
r\varphi_{\omega} & \text{otherwise}.
\end{cases}
\]
Here let us see 
\begin{equation}\label{eqn:prop:equiv:coherent:01}
\forall\, \omega \in \Omega,\, \forall\, x \in X_{\omega}^{\mathrm{an}}, \quad |\ndot|_{r\varphi_{\omega}}(x) \leqslant |\ndot|_{\varphi'_{\omega}}(x).
\end{equation}
Clearly we may assume that $\omega \in \Omega_{\mathrm{fin}}$. Note that
\[
\mathscr R_{r,\omega} = \{ s \in H^0(\mathscr X_{\omega}, \mathscr L_{\omega}) \mid \| s \|_{r\varphi_{\omega},\sup} \leqslant 1 \}
\]
and $\mathscr R_{r,\omega} \otimes \mathcal O_{\mathscr X_{\omega}} \to \mathscr L_{\omega}$ is surjective, so that,
by Proposition~\ref{Pro:modelFubiniStudy}, \eqref{eqn:prop:equiv:coherent:01} follows.
Therefore, $r\varphi$ is very fine, and hence $\varphi$ is fine.
\end{proof}

Finally we consider the following theorem:

\begin{theo}\label{thm:fine:asym:pure}
If $\varphi$ is very fine, then 
$\overline R_{\sbullet}$ is asymptotically pure.
\end{theo}
}

\begin{proof}
By our assumption, $\varphi$ is dominated and
there are a generically pure continuous metric family $\varphi' = \{ \varphi'_{\omega} \}_{\omega \in \Omega}$ of $L$ and a non-empty Zariski open set $U$ of
$\Spec(\mathfrak o_K)$ such that
$|\ndot|_{\varphi_{\omega}} \leqslant |\ndot|_{\varphi'_{\omega}}$ for all $\omega \in U \cap \Omega_{\mathrm{fin}}$.

First note that $(R_n, \xi_n)$ is dominated for $n \geqslant 0$ by
Proposition~\ref{Pro: dominancy preserved by operators} and Theorem~\ref{Thm: Fubini-Study dominated}.
Moreover, by Proposition~\ref{prop:fine:coherent} or Proposition~\ref{prop:equiv:coherent}, $(R_n, \xi_n)$ is coherent for every $n \geqslant 0$.

By the generic purity of $\varphi'$, there are a non-empty Zariski open set $U'$ of $\Spec(\mathfrak o_K)$,
a projective integral scheme $\mathscr X$ over $U'$ and an invertible $\mathcal O_{\mathscr X}$-module $\mathscr L$ such that $\mathscr X \times_{U'} \Spec(K) = X$, $\rest{\mathscr L}{X} = L$ and,
for each $\omega \in U' \cap \Omega_{\mathrm{fin}}$,
$\varphi'_{\omega}$ coincides with the metric arising from $\mathscr X_{\omega}$ and $\mathscr L_{\omega}$,
where $\mathscr X_{\omega} = \mathscr X \times_{U} \Spec(\mathfrak o_{\omega})$ and
$\mathscr L_{\omega}$ is the pull-back of $\mathscr L$ to $\mathscr X_{\omega}$.
Replacing $U$ and $U'$ by $U \cap U'$, we may assume that $U = U'$.
Moreover, as $X$ is geometrically integral over $K$, by virtue of \cite[Th\'{e}or\`{e}me~9.7.7]{EGAIV-3},
shrinking $U$ if necessarily, we may also assume that, for any $\omega \in U \cap \Omega_{\mathrm{fin}}$, 
the fiber of $\mathscr X \to U$ over $\omega$ is
geometrically integral over the residue field at $\omega$. 
Then, by Proposition~\ref{Pro:latticeby model}, 
\begin{equation}\label{eqn:prop:model:ample:dominant:01}
\begin{cases}
\{ x \in R_n \mid \| x \|_{n\varphi'_{\omega}, \sup} \leqslant 1 \} =
H^0(\mathscr X_{\omega}, \mathscr L_{\omega}^{\otimes n}) = H^0(\mathscr X, \mathscr L^{\otimes n}) \otimes_{\mathfrak o_U} \mathfrak o_{\omega}, \\
\|\ndot\|_{n\varphi'_{\omega},\sup} = \|\ndot\|_{H^0(\mathscr X_{\omega}, \mathscr L_{\omega}^{\otimes n})}
\end{cases}
\end{equation}
for all $n \geqslant 1$ and $\omega \in U \cap \Omega_{\mathrm{fin}}$.

\begin{enonce}{Claim}\label{claim:thm:dom:coh:asym:pure:02}
\begin{enumerate}[label=\rm(\alph*)]
\item
$|\ndot|_{n\varphi_{\omega}}(x) \leqslant |\ndot|_{n\varphi'_{\omega}}(x)$ for all $\omega \in U \cap \Omega_{\mathrm{fin}}$,
$x \in X_{\omega}^{\mathrm{an}}$ and $n \geqslant 1$.

\item
$\|\ndot\|_{n\varphi_{\omega}, \sup} \leqslant \|\ndot\|_{\mathscr R_{n,\omega}} \leqslant 
\min \left\{
|\varpi_{\omega}|_{\omega}^{-1} \|\ndot\|_{n\varphi_{\omega}, \sup},\ 
\|\ndot\|_{n\varphi'_{\omega}, \sup}\right\}$ for all $\omega \in U \cap \Omega_{\mathrm{fin}}$ and $n \geqslant 1$, where $\varpi_{\omega}$ is a uniformizing parameter of $\mathfrak o_{\omega}$.
\end{enumerate}
\end{enonce}

\begin{proof}
(a) is obvious.

(b) First of all, by Proposition~\ref{Pro:normetreausauxdisc},
\[
\|\ndot\|_{n\varphi_{\omega}, \sup} \leqslant \|\ndot\|_{\mathscr R_{n,\omega}} \leqslant 
|\varpi_{\omega}|_{\omega}^{-1} \|\ndot\|_{n\varphi_{\omega}, \sup}.
\]
By (a), one has $\|\ndot\|_{n\varphi_{\omega},\sup} \leqslant \|\ndot\|_{n\varphi'_{\omega},\sup}$,
so that, by \eqref{eqn:prop:model:ample:dominant:01}, one obtains
\[
\mathscr R_{n,\omega} \supseteq H^0(\mathscr X_{\omega}, \mathscr L_{\omega}^{\otimes n}).
\]
Therefore, by \eqref{eqn:prop:model:ample:dominant:01} again, (b) follows.
\end{proof}

\begin{enonce}{Claim}\label{claim:thm:dom:coh:asym:pure:03}
If we set $A_{\omega} = d_{\omega}(\varphi_{\omega}, \varphi'_{\omega})$ for $\omega \in \Omega$, then one has the following:
\begin{enumerate}[label=\rm(\alph*)]
\item
${\displaystyle \sup_{x \in R_n \setminus \{ 0 \}} \ln \frac{\| s \|_{\mathscr R_{n,\omega}}}{\| s \|_{n\varphi_{\omega}, \sup}} \leqslant A_{\omega} n}$ for all $\omega \in U \cap \Omega_{\mathrm{fin}}$ and
$n \geqslant 1$.
\item
${\displaystyle \int_{\Omega} A_{\omega} \nu(d\omega) =
\sum_{\omega \in \Omega} A_{\omega} \nu(\{ \omega \}) < \infty}$.
\end{enumerate}
\end{enonce}

\begin{proof}
(a) By using the inequality \eqref{Equ:comparaisondistancesup} in Subsection~\ref{subsec:distance:metric} together with
Claim~\ref{claim:thm:dom:coh:asym:pure:02}, 
\begin{align*}
\sup_{x \in R_n \setminus \{ 0 \}} \ln \frac{\| s \|_{\mathscr R_{n,\omega}}}{\| s \|_{n\varphi_{\omega}, \sup}} & \leqslant
\sup_{x \in R_n \setminus \{ 0 \}} \ln \frac{\| s \|_{n\varphi'_{\omega}, \sup}}{\| s \|_{n\varphi_{\omega}, \sup}} \\
& = d_{\omega}\left(\|\ndot\|_{n\varphi'_{\omega}, \sup}, \|\ndot\|_{n\varphi'_{\omega}, \sup}\right) \\
& \leqslant d_{\omega}(n\varphi_{\omega}, n\varphi'_{\omega}) = n A_{\omega}.
\end{align*}

(b) Note that $\varphi'$ is dominated by Proposition~\ref{prop:generic:purity:dominant:coherent}.
Moreover, $\varphi$ is dominated by our assumption, so that, by Proposition~\ref{Pro: dominancy preserved by operators},
the function $\omega \mapsto d_{\omega}(\varphi_{\omega}, \varphi'_{\omega}) = A_{\omega}$ is $\nu$-dominated. Thus
one obtains (b).
\end{proof}

Fix a positive number $\epsilon$. Then, by Claim~\ref{claim:thm:dom:coh:asym:pure:03}, 
there is a non-empty Zariski open set $U_{\epsilon}$ of $U$
such that 
\[
\sum_{\omega \in U_{\epsilon} \cap \Omega_{\mathrm{fin}}} A_{\omega} \nu(\{ \omega \}) \leqslant \epsilon.
\]
Thus, if we set $B = \sum_{\omega \in \Omega_{\mathrm{fin}} \setminus U_{\epsilon}}
-\ln |\varpi_{\omega}|_{\omega} \nu (\{\omega\})$, then, by Claim~\ref{claim:thm:dom:coh:asym:pure:02},
\begin{align*}
\bbsigma(R_n, \xi_n) & = \sum_{\omega \in \Omega_{\mathrm{fin}}}
\sup_{x \in R_{n,\omega} \setminus \{ 0 \}} \ln \left( \frac{\|x\|_{\mathscr R_{n,\omega}}}{\| x \|_{n\varphi_{\omega}, \sup}}\right) \nu (\{\omega\}) \\
& \leqslant \sum_{\omega \in U_{\epsilon} \cap \Omega_{\mathrm{fin}}}
\sup_{x \in R_{n,\omega} \setminus \{ 0 \}} \ln \left( \frac{\|x\|_{\mathscr R_{n,\omega}}}{\| x \|_{n\varphi_{\omega}, \sup}}\right) \nu (\{\omega\})
+ B \\
& \leqslant n \sum_{\omega \in U_{\epsilon} \cap \Omega_{\mathrm{fin}}}
A_{\omega}  \nu (\{\omega\})
+ B \leqslant n \epsilon + B
\end{align*}
for $n \geqslant 1$, and hence one has
\[
\limsup_{n\to\infty} \frac{\bbsigma(R_n, \xi_n)}{n} \leqslant \epsilon,
\]
so that the assertion of the theorem follows.
\end{proof}
}

{
\subsection{A generalization of Nakai-Moishezon's criterion}

Let $X$ be a geometrically integral projective variety over $K$ and $L$ be an invertible sheaf on $X$.
For $\omega \in \Omega$, let $X_{\omega} := X \times_{\Spec K} K_{\omega}$ and
$L_{\omega} := L \otimes_{\mathcal O_X} \mathcal O_{X_{\omega}}$.
Let $\varphi_{\omega}$ be a continuous metric of $L_{\omega}$ on $X_{\omega}^{\mathrm{an}}$ for each $\omega \in \Omega$, and
$\varphi := \{ \varphi_{\omega} \}_{\omega \in \Omega}$. 
For $n \geqslant 0$ and a subvariety $Y$ of $X$, let $\xi_n := \{ \|\ndot\|_{n\varphi_{\omega},\sup} \}_{\omega \in \Omega}$ and
$\rest{\xi_n}{Y} := \{ \|\ndot\|_{n\rest{\varphi_{\omega}}{Y_{\omega}},\sup} \}_{\omega \in \Omega}$.
In this subsection, let us consider the following Nakai-Moishezon's criterion over a number field, which gives a generalisation of Nakai-Moishezon's criterion due to Shouwu Zhang.

\begin{theo}
\label{theo:Nakai:Moishezon:criterion:number:field}
We assume the following:

\begin{enumerate}[label=\rm(\arabic*)]
\item {\rm (Fineness)} 
The metric family $\varphi$ is fine.

\item {\rm (Semipositivity)}
$L$ is semiample and $\varphi_{\omega}$ is semipositive for every $\omega \in \Omega$.

\item {\rm (Bigness)}
For every subvariety $Y$ of $X$, $\rest{L}{Y}$ is big, and 
there are a positive number $n_Y$ and $s_Y \in H^0(Y, \rest{L^{\otimes n_Y}}{Y}) \setminus \{ 0 \}$ such that
$\widehat{\deg}_{\rest{\xi_{n_Y}}{Y}}(s_Y) > 0$.
\end{enumerate} 
Then one has 
\[
\liminf_{n\to\infty} \frac{\bblambda\left(H^0(X, L^{\otimes n}), \xi_n\right)}{n} > 0.
\]
\end{theo}

\begin{proof}
We set $R_n := H^0(X, L^{\otimes n})$ for $n \geqslant 0$.
By Proposition~\ref{prop:comp:lambda:R:M}, for a positive number $h$,
\[
\frac{1}{h} \liminf_{n\to\infty} \frac{\bblambda\left(R_n, \xi_n \right)}{n}
\leqslant \liminf_{n\to\infty} \frac{\bblambda\left(R_n, \xi_n \right)}{n},
\]
so that, replacing $L$, $\varphi$, $n_Y$ and $s_Y$ by $L^{\otimes h}$, $h\varphi$, $hn_Y$ and $s_Y^{\otimes h}$
for a sufficiently large integer $h$, we may assume that $\varphi$ is very fine.
Moreover, by Remark~\ref{remark:number:field:measurable}, we can see that $\varphi$ is measurable.
Therefore, by Theorem~\ref{thm:Nakai:Moishezon:general:setting}, one has
\[
\liminf_{n\to\infty} \frac{\nu_{\min}\left(R_n, \xi_n \right)}{n} > 0.
\]
By Proposition~\ref{Pro: dominancy preserved by operators}, Theorem~\ref{Thm: Fubini-Study dominated} and Proposition~\ref{prop:fine:coherent}, $\left(R_n, \xi_n \right)$ is dominated and coherent
for all $n \geqslant 0$.
Therefore, by Proposition~\ref{prop:comp:lambda:nu:number:field:graded:ring} and
Theorem~\ref{thm:fine:asym:pure}, 
one can see the assertion of the theorem.
\end{proof}

\begin{rema}\label{rema:fineness:dominant:coherent:num:field}
In \cite{Chen_Moriwaki2017}, Nakai-Moishezon's criterion was proved under the following condition \eqref{eqn:rema:fineness:dominant:coherent:num:field}
instead of (Fineness) in Theroem~\ref{theo:Nakai:Moishezon:criterion:number:field}:
\begin{equation}\label{eqn:rema:fineness:dominant:coherent:num:field}
\begin{cases}
\text{The adelic vector bundle $\left(H^0(X, L^{\otimes n}), \xi_n \right)$ over $S$} \\
\text{is dominaited and coherent for every $n \geqslant 0$.}
\end{cases}
\end{equation}
\end{rema}
}

%% file: cha_2019_03_23.tex

\appendix
\chapter{Reminders on measure theory}

\IfChapVersion
\ChapVersion{Version of Appendix : \\ \StrSubstitute{\DateChapAppendix}{_}{\_}}
\fi

\section{Monotone class theorems}

We recall here a monotone class theorem in the functional form and several related results and we refer to \cite[\S I.2]{DM} and \cite[\S 2.2]{Yan} for reference. For convenience of readers, we include the proof here. We fix in this section a non-empty set $\Omega$. If $\mathcal H$ is a family of real-valued functions on $\Omega$, we denote by $\sigma(\mathcal H)$ the $\sigma$-algebra on $\Omega$ generated by $\mathcal H$. It is the smallest $\sigma$-algebra on $\Omega$ with respect to which all functions in $\mathcal H$ are measurable.

\begin{defi}\label{Def:lambda-family}
Let  $\mathcal H$ be a family of non-negative and bounded functions on $\Omega$. We say that $\mathcal H$ is a \emph{$\lambda$-family}\index{lambda-family@$\lambda$-family} if it verifies the following conditions:
\begin{enumerate}[label=\rm(\roman*)]
\item the constant function $1$ belongs to $\mathcal H$;
\item if $f$ and $g$ are two functions in $\mathcal H$, $a$ and $b$ are non-negative numbers, then $af+bg\in\mathcal H$;
\item if $f$ and $g$ are two functions in $\mathcal H$ such that $f\leqslant g$, then $g-f\in\mathcal H$;
\item if $\{f_n\}_{n\in\mathbb N}$ is an increasing and uniformly bounded sequence of functions in $\mathcal H$, then the limit of the sequence $\{f_n\}_{n\in\mathbb N}$ belongs to $\mathcal H$.
\end{enumerate}
\end{defi}

\begin{lemm}\label{Lem:monotoneclass}
Let $\mathcal H$ be a $\lambda$-family of non-negative and bounded functions on $\Omega$. If for any couple $(f,g)$ of functions in $\mathcal H$, one has $\min(f,g)\in\mathcal H$, then any non-negative, bounded and $\sigma(\mathcal H)$-measurable function on $\Omega$ belongs to $\mathcal H$. In particular, the $\sigma$-algebra $\sigma(\mathcal H)$ is equal to the set of all $A\subseteq\Omega$ such that $\indic_A\in\mathcal H$.
\end{lemm}
\begin{proof}
Let $\mathcal F$ be the set of all $A\subseteq\Omega$ such that $\indic_A\in\mathcal H$. Since $\mathcal H$ is a $\lambda$-family, we obtain that $\mathcal F$ is a $\lambda$-system\footnote{Namely $\mathcal F$ satisfies the following conditions: (i) $\Omega\in\mathcal F$; (ii) if $A\in\mathcal F$, $B\in\mathcal F$ and $A\subseteq B$, then $B\setminus A\in\mathcal F$, (iii) if $\{A_n\}_{n\in\mathbb N}$ is an increasing sequence of elements of $\mathcal F$, then the union $\bigcup_{n\in\mathbb N}A_n$ belongs to $\mathcal F$.} and at the same time a $\pi$-system (namely for all $A\in\mathcal F$ and $B\in\mathcal F$ one has $A\cap B\in\mathcal F$) since the family $\mathcal H$ is supposed to be stable by the operator $(f,g)\mapsto\min(f,g)$. Therefore $\mathcal F$ is actually a $\sigma$-algebra.

If $(f,g)$ is a couple of functions in $\mathcal H$, one has
\[\max(f,g)=f+g-\min(f,g)\in\mathcal H.\]
In particular, if $f\in\mathcal H$ and $a\in\mathbb R_+$, then 
\[\max(f-a,0)=\max(f,a)-a\in\mathcal H.\]
This property actually implies that, for any $f\in\mathcal H$ and any integer $n\geqslant 1$, one has $f^n\in\mathcal H$. In fact, the function $x\mapsto x^n$ is convex on $\mathbb R_+$, which can be written as the supremum of a countable family of functions of the form \[x\longmapsto \max(na^{n-1}x-(n-1)a^n,0)\] with $a\in\mathbb Q_+:=\mathbb Q\cap\mathbb R_+$. Therefore by the condition (iv) in Definition \ref{Def:lambda-family} one obtains
\[f^n=\sup_{a\in\mathbb Q_+}\max(na^{n-1}f-(n-1)a^n,0)\in\mathcal H.\]

If $f$ is an element of $\mathcal H$ and $t$ is a real number, $t>0$, one has $\min(t^{-1}f,1)\in\mathcal H$. Moreover, the sequence $\{1-\min(t^{-1}f,1)^n\}_{n\in\mathbb N,\,n\geqslant 1}$ is increasing and converges to $\indic_{\{f<t\}}$, which implies that $\indic_{\{f<t\}}\in\mathcal H$ and hence $\{f<t\}\in\mathcal F$. Therefore every function in $\mathcal H$ is $\mathcal F$-measurable, and thus $\sigma(\mathcal H)\subseteq\mathcal F$.

It remains to prove that any non-negative bounded $\mathcal F$-measurable function belongs to $\mathcal H$. Let $f$ be such a function. For any integer $n\geqslant 1$, let
\[f_n=\sum_{k=0}^{n2^n-1}\frac{k}{2^n}\indic_{\{k/2^n\leqslant f<(k+1)/2^n\}}+n\indic_{\{f\geqslant n\}}.\]
This is a function in $\mathcal H$. Moreover, the sequence $\{f_n\}_{n\in\mathbb N,\,n\geqslant 1}$ is increasing and converges to $f$. Therefore $f\in\mathcal H$.
\end{proof}

\begin{theo}\label{Thm:monotoneclass}
Let $\mathcal H$ be a $\lambda$-family of non-negative and bounded functions on $\Omega$ and $\mathcal C$ be a subset of $\mathcal H$. Assume that for any couple $(f,g)$ of functions in $\mathcal C$, the product function $fg$ belongs to $\mathcal C$. Then any non-negative and bounded  $\sigma(\mathcal C)$-measurable function belongs to $\mathcal H$. 
\end{theo}
\begin{proof}
By replacing $\mathcal H$ by the intersection of all $\lambda$-families containing $\mathcal C$ we may assume that $\mathcal H$ is the smallest $\lambda$-family which contains $\mathcal C$.

We first prove that $\mathcal H$ is stable by multiplication. Let $\mathcal H_1$ be the set of all non-negative and bounded functions $f$ on $\Omega$ such that $fg\in\mathcal H$ for any $g\in\mathcal C$. This is a $\lambda$-family containing $\mathcal C$. Hence one has $\mathcal H_1\supseteq\mathcal H$. Let $\mathcal H_2$ be the set of all non-negative and bounded functions $f$ on $\Omega$ such that $fg\in\mathcal H$ for any $g\in\mathcal H$. This is also a $\lambda$-family. Moreover, since $\mathcal H_1\supseteq\mathcal H$ one obtains $\mathcal H_2\supseteq\mathcal C$ and hence $\mathcal H_2\supseteq\mathcal H$, which implies that $\mathcal H$ is stable by multiplication.

Let $f$ and $g$ be two functions in $\mathcal H$. We will prove that $|f-g|\in\mathcal H$. By dilating the function $|f-g|$ by a positive constant, we may assume that $|f-g|$ is bounded from above by $1$. One has
\[(f-g)^2=f^2+g^2-2fg\in\mathcal H.\]
Let $\{f_n\}_{n\in\mathbb N}$ be the sequence of functions on $\Omega$ defined by the following recursive formula
\[f_0=0,\quad f_{n+1}=f_n+\frac{1}{2}((f-g)^2-f^2_n).\]
By induction on $n$, we can show that $f_n\in\mathcal H$ and $f_n\leqslant|f-g|$. In fact, these properties are trivially satisfied by $f_0$. If $f_n\in\mathcal H$ and $f_n\leqslant|f-g|$, then one has $f_{n+1}\in\mathcal H$. Moreover, by the relation $|f-g|\leqslant 1$ one obtains $f_{n+1}\leqslant |f-g|$ since the function $t\mapsto t-\frac12t^2$ is increasing on the interval $[0,1]$. The properties $f_n\in\mathcal H$ and $f_n\leqslant|f-g|$ show that the sequence $\{f_n\}_{n\in\mathbb N}$ is increasing and converges to $|f-g|$. Hence $|f-g|\in\mathcal H$, which implies that
\[\min(f,g)=\frac12(f+g-|f-g|)\in\mathcal H.\]
By Lemma \ref{Lem:monotoneclass}, any non-negative, bounded and $\sigma(\mathcal H)$-measurable function belongs to $\mathcal H$. The theorem is thus proved.
\end{proof}

{
\section{Measurable selection theorem}
In this section, we recall a measurable selection theorem due to  Kuratowski and Ryll-Nardzewski \cite{Kuratowski65}. See \cite[Chapter 5]{Parthasarathy72} for more details.

\begin{theo}\label{Thm: measurable selection}
Let $Y$ be a complete separable metric space and $\mathscr P(Y)$ be the set of subsets of $Y$. Let $(\Omega,\mathcal A)$ be a measurable space and $F:\Omega\rightarrow\mathscr P(Y)$ be a map. We assume that
\begin{enumerate}[label=\rm(\arabic*)]
\item for any $\omega\in\Omega$, the set $F(\omega)$ is a non-empty closed subset of $Y$,
\item for any open subset $U$ of $Y$, the set $\{\omega\in\Omega\,:\,F(\omega)\cap U\neq \varnothing\}$ belongs to $\mathcal A$. 
\end{enumerate}
Then there exist a measurable map $f:\Omega\rightarrow Y$ such that $f(\omega)\in F(\omega)$ for any $\omega\in\Omega$.
\end{theo}

{
\section{Vague convergence and weak convergence of measures}\label{Sec: vague convergence}

Let $X$ be a locally compact Hausdorff space. Recall that a \emph{Radon measure}\index{Radon measure@Radon measure} is by definition a Borel measure $\nu$ on $X$ which satisfies the following conditions:
\begin{enumerate}[label=\rm(\arabic*)]
\item $\nu$ is \emph{tight},\index{tight@tight} that is, for any Borel subset $B$ of $X$, $\nu(B)$ is equal to the supremum of $\nu(K)$, where $K$ runs over the set of compact subsets of $B$;
\item $\nu$ is \emph{outer regular}\index{outer regular@outer regular}, that is, for any Borel subset $B$ of $X$, $\nu(B)$ is equal to the infimum of $\nu(U)$, where $U$ runs over the set of open subsets of $X$ containing $B$; 
\item $\nu$ is \emph{locally finite}\index{locally finite@locally finite}, that is, for any $x\in X$ there exists a neighbourhood $U$ of $x$ such that $\nu(U)<+\infty$.
\end{enumerate} 
We denote by $\mathscr M(X)$ be the set of Radon measures on $X$. Let $C_c(X)$ be the vector space of continuous real-valued functions of compact support on $X$. We say that an $\mathbb R$-linear map $\varphi:C_c(X)\rightarrow\mathbb R$ is a \emph{positive functional}\index{positive functional@positive functional} if $\varphi(f)\geqslant 0$ for any non-negative function $f$ in $C_c(X)$. 
Recall the Riesz's representation theorem as follows. See \cite[\S56]{Halmos74} for a proof.
\begin{theo}
Let $X$ be a locally compact Hausdorff space. The map sending $\nu\in\mathscr M(X)$ to the positive functional
\[f\in C_c(X)\longrightarrow \int_Xf\,\mathrm{d}\nu\]
defines a bijection between the set $\mathscr M(X)$ and the set of all positive linear functionals on $C_c(X)$.
\end{theo}

The \emph{vague topology}\index{vague topology@vague topology} on $\mathscr M(X)$ is an example of weak-* topology if we identify $\mathscr M(X)$ with a subset of the dual space of $C_c(X)$. More precisely, we say that a sequence $\{\nu_n\}_{n\in\mathbb N}$ of Radon measures \emph{converges vaguely}\index{converges vaguely@converges vaguely} if for any function $f\in C_c(X)$, the sequence of integrals
$\{\int_{X}f\,\mathrm{d}\nu_n\}_{n\in\mathbb N}$
converges in $\mathbb R$. Note that the limit of the above sequence defines a positive linear functional on $C_c(X)$ when $f$ varies, which corresponds to a Radon measure, called the \emph{vague limit}\index{vague limit@vague limit} of $\{\nu_n\}_{n\in\mathbb N}$. 

If $\{\nu_n\}_{n\in\mathbb N} $ is a sequence of Radon probability measure which converges vaguely, the limite measure may have a total mass $<1$. In probability theory, the notion of weak convergence is also largely used. Let $\mathscr M_1(X)$ be the subset of $\mathscr M(X)$ of probability measures. Let $C_b(X)$ be the vector space of bounded continuous functions. We say that a sequence $\{\nu_n\}_{n\in\mathbb N}$ of  measures in $\mathscr M_1(X)$ (they are therefore \emph{probability} measures) \emph{converges weakly}\index{converges weakly@converges weakly} if for any bounded continuous function $f$ on $X$, the sequence of integrals
$\{\int_{X}f\,\mathrm{d}\nu_n\}_{n\in\mathbb N}$
converges in $\mathbb R$. Clearly, if the sequence $\{\nu_n\}_{n\in\mathbb N}$ converges weakly, then it also converges vaguely, and its vague limit is also called its \emph{weak limit}\index{weak limit@weak limit}. Note that in the weak convergence case the limit measure should be a probability measure. The following criterion provides a criterion of weak convergence for vaguely convergence sequence of Radon probability measures. We refer the readers to \cite[Theorem 13.16]{Klenke14} for the proof and for more details\footnote{In \cite[Theorem 13.16]{Klenke14}, it is assumed that the topological space is a locally compact Polish space.  This condition is satisfied notably when $X$ is a locally compact Hausdorff space with countable base, see \cite{Vaughan37}. However, it actually suffices that the topological space is locally compact and metrisable (see Lemma 13.10 of \cite{Klenke14} which is used in the proof of Theorem 13.16 of \emph{loc. cit.}).}.

\begin{theo}\label{Thm: criterion of weak convergence}
Let $X$ be a locally compact metrisable space and $\{\nu_n\}_{n\in\mathbb N}$ be a sequence of Radon probability measures on $X$, which converges vaguely to a limite measure $\nu$. Assume the limite measure $\nu$ is a probability measure. Then the sequence $\{\nu_n\}_{n\in\mathbb N}$ converge weakly to $\nu$.
\end{theo}

}

\section{Upper and lower integral}\label{Sec:ulint}

Let $(\Omega,\mathcal A,\nu)$ be  a measure space.
We denote by $\mathscr L^1(\Omega,{\mathcal A},\nu)$ the vector space of all real-valued $\nu$-integrable functions on $(\Omega,{\mathcal A})$. 
We say that a subset $A$ of $\Omega$ is $\nu$-\emph{negligible}\index{negligible@negligible} if there exists a set $B\in\mathcal A$ such that $\nu(B)=0$ and that $A\subseteq B$. We say that two functions $h_1$ and $h_2$ on $\Omega$ are $\nu$-\emph{indistinguishable}\index{indistinguishable@indistinguishable} if $\{h_1\neq h_2\}$ is a $\nu$-negligible set. Any function on $\Omega$ which is $\nu$-indistinguishable with the zero function is said to be $\nu$-\emph{negligible}\index{negligible@negligible}. In other words, a function $f$ on $\Omega$ is $\nu$-negligible if and only if $\{f\neq 0\}$ is a $\nu$-negligible set. If a formula depending on a variable $\omega\in\Omega$ is satisfied outside of a $\nu$-negligible set, we say that it holds $\nu$-almost everywhere (written in abbriviation as $\nu$-a.e.).

\begin{defi}\label{Def:upper and lower integral}
We construct two non-necessarily linear functional $\overline I_v(\ndot)$ and $\underline{I}_\nu(\ndot)$ as follows. For any function $h:\Omega\rightarrow \mathbb R$, let
\begin{align*}
\upint_{\Omega} h(\omega) \nu(d\omega):=\inf_{\begin{subarray}{c}
f\in\mathscr L^1(\Omega,\mathcal A,\nu)\\
f\geqslant h\,\nu\text{-a.e.}
\end{subarray}}\int_{\Omega} f(\omega) \nu(d\omega),\\
\lowint_{\Omega} h(\omega) \nu(d\omega):=\sup_{\begin{subarray}{c}
g\in\mathscr L^1(\Omega,\mathcal A,\nu)\\
g\leqslant h\,\nu\text{-a.e.}
\end{subarray}}\int_{\Omega} g(\omega) \nu(d\omega).
\end{align*}
If $h$ is not $\nu$-almost everywhere bounded from above by any integrable function, then $\upint_{\Omega} h(\omega) \nu(d\omega)$ is defined as $+\infty$ by convention. Similarly, if $h$ is not $\nu$-almost everywhere bounded from below by any integrable function, then $\lowint_{\Omega} h(\omega) \nu(d\omega)$ is defined as $-\infty$ by convention. The values $\upint_{\Omega} h(\omega) \nu(d\omega)$ and $\lowint_{\Omega} h(\omega) \nu(d\omega)$ are called \emph{upper integral}\index{upper integral@upper integral} and \emph{lower integral}\index{lower integral@lower integral} of the function $h$, respectively.
From now on, for simplicity,
\[
\upint_{\Omega} h(\omega) \nu(d\omega),\quad
\lowint_{\Omega} h(\omega) \nu(d\omega)
\quad\text{and}\quad
\int_{\Omega} f(\omega) \nu(d\omega)
\]
are denoted by $\overline I_\nu(h)$, $\underline I_\nu(h)$ and $I_\nu(f)$,
respectively, for any function $h$ on $\Omega$ and any integrable function $f$ on $\Omega$.
\end{defi}

The following properties are straightforward from the definition of upper and lower integrals.

\begin{prop}\phantomsection\label{Pro: preserve order int sup and inf}
\begin{enumerate}[label=\rm(\arabic*)]
\item For any function $h:\Omega\rightarrow\mathbb R$
\begin{equation}\underline I_\nu(h)\leqslant\overline I_\nu(h).\end{equation}
\item 
If $h_1$ and $h_2$ are two real-valued functions on $\Omega$ such that $h_1\leqslant h_2$, then
\begin{equation}\label{Equ:dominantordre}
\underline{I}_\nu(h_1)\leqslant\underline{I}_\nu(h_2),\quad
\overline{I}_\nu
(h_1)\leqslant\overline{I}_\nu(h_2).
\end{equation}
\end{enumerate}
\end{prop}

\begin{prop}\label{Pro:integrableupplower}
Let $h$ be a real-valued function on $\Omega$. Then $h$ is $\nu$-indistinguishable with a $\nu$-integrable function if and only if $\overline{I}_\nu(h)=\underline{I}_\nu(h)\in\mathbb R$.
\end{prop}
\begin{proof}
If $h$ is indistinguishable with a $\nu$-integrable function $\widetilde h$, then $\overline{I}_\nu(h)$ and $\underline{I}_\nu(h)$ are both equal to the integral of $\widetilde h$ with respect to the measure $\nu$, which is a real number.

Conversely, assume that $\overline{I}_\nu(h)=\underline{I}_\nu(h)$, then we can find two sequences $\{f_n\}_{n\in\mathbb N}$ and $\{g_n\}_{n\in\mathbb N}$ of functions in $\mathscr L^1(\Omega,\mathcal A,\nu)$ such that $g_n\leqslant h\leqslant f_n$ $\nu$-almost everywhere. and that 
\[\lim_{n\rightarrow+\infty}I_\nu(f_n)=\overline{I}_\nu(h)=\underline I_{\nu}(h)=\lim_{n\rightarrow+\infty}I_\nu(g_n).\] Without loss of generality, we may assume that the sequence $\{f_n\}_{n\in\mathbb N}$ is decreasing and $\{g_n\}_{n\in\mathbb N}$ is increasing (otherwise we replace $f_n$ by $\widetilde{f}_n=\min\{f_1,\ldots,f_n\}$ and $g_n$ by $\widetilde{g}_n=\max\{g_1,\ldots,g_n\}$). 
Let $f=\inf_{n\in\mathbb N}f_n$ and $g=\sup_{n\in\mathbb N}g_n$. By Lebesgue's dominant convergence theorem, we obtain that $f$ and $g$ are both $\nu$-integrable, and
\[I_\nu(f)=\overline{I}_\nu(h)=\underline{I}_\nu(h)=I_\nu(g).\]
Moreover, one has $g\leqslant h\leqslant f$ $\nu$-almost everywhere, which implies that $f=g=h$ $\nu$-almost everywhere.
\end{proof}

In general the operators $\overline{I}_\nu(\ndot)$ and $\underline{I}_\nu(\ndot)$ are not linear operators. However, they satisfies some convexity property.

\begin{prop}\label{Pro:subsuperadditive}
Let $h_1$ and $h_2$ be two real-valued functions on $\Omega$. 
\begin{enumerate}[label=\rm(\arabic*)]
\item Assume that $\{\overline I_\nu(h_1),\overline I_\nu(h_2)\}\neq\{+\infty,-\infty\}$. Then one has
\begin{equation}\label{Equ: subadditive}\overline{I}_\nu(h_1+h_2)\leqslant \overline{I}_\nu(h_1)+\overline{I}_\nu(h_2)\end{equation}
\item Assume that $\{\underline I_\nu(h_1),\underline I_\nu(h_2)\}\neq\{+\infty,-\infty\}$. Then one has \begin{equation}\label{Equ: superadditive}\underline{I}_\nu(h_1+h_2)\geqslant \underline{I}_\nu(h_1)+\underline{I}_\nu(h_2).\end{equation}
\end{enumerate}
\end{prop} 
\begin{proof}(1) We first treat the case where neither of $\overline I_\nu(h_1)$ and $\overline I_{\nu}(h_2)$ is $+\infty$.
If $f_1$ and $f_2$ be two $\nu$-integrable functions on $\Omega$ such that $f_1\geqslant h_1$ and $f_2\geqslant h_2$ $\nu$-almost everywhere. Then the sum $f_1+f_2$ is $\nu$-integrable, and $f_1+f_2\geqslant h_1+h_2$. Therefore $I_\nu(f_1)+I_\nu(f_2)=I_\nu(f_1+f_2)\geqslant \overline{I}_\nu(h_1+h_2)$. Since $f_1$ and $f_2$ are arbitrary, we obtain $\overline{I}_\nu(h_1)+\overline{I}_\nu(h_2)\geqslant\overline{I}_\nu(h_1+h_2)$.

If at least one of the upper integrals $\overline I_\nu(h_1)$ and $\overline I_\nu(h_2)$ is $+\infty$, then by the hypothesis $\{\overline I_\nu(h_1),\overline I_\nu(h_2)\}\neq\{+\infty,-\infty\}$ one has $\overline I_\nu(h_1)+\overline I_\nu(h_2)=+\infty$. Hence the inequality \eqref{Equ: subadditive} is trivial.

The proof of the statement (2) is very similar to that of (1). We omit the details.
\end{proof}

\begin{prop}\label{Pro:translationbyintegrable}
Let $h$ be a real-valued function on $\Omega$ and $\varphi$ be a $\nu$-integrable function. Then one has
\[\overline{I}_\nu(h+\varphi)=\overline{I}_\nu(h)+I_\nu(\varphi),\quad\underline{I}_\nu(h+\varphi)=\underline{I}_\nu(h)+I_\nu(\varphi).\] 
\end{prop}
\begin{proof} Since $\varphi$ is $\nu$-integrable, one has
\[\overline{I}_\nu(\varphi)=I_\nu(\varphi)=\underline{I}_\nu(\varphi)\in\mathbb R.\]
By Proposition \ref{Pro:subsuperadditive}, one has
\[\overline{I}_\nu(h+\varphi)\leqslant\overline{I}_\nu(h)+I_\nu(\varphi).\]
Moreover, if we apply this inequality to $h+\varphi$ and $-\varphi$, we obtain
\[\overline{I}_\nu(h)\leqslant\overline{I}_\nu(h+\varphi)-I_\nu(\varphi).\]
Therefore the first equality is true. The proof of the second equality is quite similar, we omit the details.
\end{proof}

\begin{prop}\label{Pro:multscal}
Let $h$ be a real-valued function on $\Omega$. If $a$ is a non-negative number, then one has 
\[\overline{I}_\nu(ah)=a\overline{I}_\nu(h),\quad\underline{I}_\nu(ah)=a\underline{I}_\nu(h).\]
\end{prop}
\begin{proof}
The assertions are trivial when $a=0$. In the following, we assume that $a>0$. If $f$ is a $\nu$-integrable function such that $h\leqslant f$ $\nu$-almost everywhere, then $af$ is a $\nu$-integrable function such that $ah\leqslant af$ $\nu$-almost everywhere. Therefore, we obtain that $\overline{I}_\nu(ah)\leqslant a\overline{I}(h)$. If we apply this inequality to $a^{-1}$ and $ah$ we get $\overline{I}_\nu(h)\leqslant a^{-1}\overline{I}_{\nu}(ah)$. Hence the first equality is true. The proof of the second equality is very similar, we omit the details.
\end{proof}

\begin{prop}\label{Pro:upperlownega}
Let $h$ be a real-valued function on $\Omega$. One has 
\begin{equation*}
\overline{I}_\nu(-h)=-\underline{I}_\nu(h),\quad\underline{I}_\nu(-h)=-\overline{I}_\nu(h).
\end{equation*}
\end{prop}
\begin{proof}
If $f$ is an $\nu$-integrable function such that $-h\leqslant f$ $\nu$-almost everywhere, then one has $-f\leqslant h$ $\nu$-almost everywhere. Since $f$ is arbitrary, we obtain $-\overline{I}_\nu(-h)\leqslant \underline{I}_\nu(h)$. Similarly, if $g$ is an $\nu$-integrable function such that $g\leqslant -h$ $\nu$-almost everywhere, then one has $h\leqslant -g$ $\nu$-almost everywhere Since $g$ is arbitrary, we obtain $-\underline{I}_\nu(-h)\geqslant \overline{I}_{\nu}(h)$. Finally, if we apply the obtained inequality to $-h$, we obtain $-\overline{I}_\nu(h)\leqslant \underline{I}_\nu(-h)$ and $-\underline{I}_\nu(h)\geqslant\overline{I}_\nu(-h)$. Therefore the equalities hold.
\end{proof}

\begin{prop}\label{Pro:differencedomi}
Let $h_1$ and $h_2$ be two real-valued functions on $\Omega$, and let $h=h_1+h_2$. Assume that $\{\overline{I}_\nu(h_1),\underline{I}_\nu(h_2)\}\neq\{+\infty,-\infty\}$. Then one has
\[\underline{I}_\nu(h)\leqslant\overline{I}_\nu(h_1)+\underline{I}_\nu(h_2)\leqslant\overline{I}_\nu(h).\]
\end{prop}
\begin{proof}
By the equality $h=h_1+h_2$ we obtain $h_1=(-h_2)+h$. Thus Proposition \ref{Pro:subsuperadditive} leads to 
\[\overline{I}_\nu(h_1)\leqslant\overline{I}_\nu(-h_2)+\overline I_{\nu}(h)=-\underline{I}_\nu(h_2)+\overline{I}_\nu(h),\]
where the equality comes from Proposition \ref{Pro:upperlownega}. Hence we obtain the inequality
\[\overline{I}_\nu(h)\geqslant\overline{I}_\nu(h_1)+\underline{I}_\nu(h_2).\]
We then apply this inequality to $-h$, $-h_2$ and $-h_1$ to get the other equality.
\end{proof}

\begin{defi}\label{Def:dominancefunc}
Let $h:\Omega\rightarrow\mathbb R$ be a real valued function on $\Omega$. We say that $h$ is \emph{$\nu$-dominated}\index{nu-dominated@$\nu$-dominated} if there exists a $\nu$-integrable function $f$ such that $\{\omega\in\Omega\,:\,|h(\omega)|\leqslant f(\omega)\}$ is a $\nu$-negligible set (in other words, $|h|\leqslant f$ $\nu$-almost everywhere). Note that this condition is equivalent to
\[\overline{I}_\nu(h)<+\infty\;\text{ and }\;\underline{I}_\nu(h)>-\infty.\]
We denote by $\mathscr D^1(\Omega,\mathcal A,\nu)$ the vector space of $\nu$-dominated functions on $\Omega$. Clearly one has $\mathscr D^1(\Omega,\mathcal A,\nu)\supseteq\mathscr L^1(\Omega,\mathcal A,\nu)$, and $\mathscr D^1(\Omega,\mathcal A,\nu)$ is invariant by the operator $f\mapsto |f|$ of taking the absolute value. Moreover, if $f$ and $g$ are real valued functions on $\Omega$ such that $|f|\leqslant |g|$ $\nu$-almost everywhere and that $g$ is $\nu$-dominated, then the function $f$ is also $\nu$-dominated.
\end{defi}

\begin{prop}\label{Pro:seminormed}
Let $\|\ndot\|_{\mathscr D^1_\nu}$ be the function on $\mathscr D^1(\Omega,\mathcal A,\nu)$ sending any $\nu$-dominated function $f$ to $\overline I_\nu(|f|)$. Then $\|\ndot\|_{\mathscr D^1_\nu}$ is a seminorm. Moreover, a function $f\in\mathscr D^1(\Omega,\mathcal A,\nu)$ satisfies $\|f\|_{\mathscr D_\nu^1}=0$ if and only if it is $\nu$-negligible. 
\end{prop}
\begin{proof}
Let $f$ be a $\nu$-dominated function and $a$ be a real number. One has $|af|=|a|\cdot|f|$. By Proposition \ref{Pro:multscal} we obtain that 
\[\|af\|_{\mathscr D_\nu^1}=\overline{I}_\nu(|af|)=\overline{I}_\nu(|a|\cdot|f|)=|a|\cdot\overline I_{\nu}(|f|)=|a|\cdot\|f\|_{\mathscr D_\nu^1}.\]
Moreover, if $f$ and $g$ are two $\nu$-dominated functions, then by Proposition \ref{Pro:subsuperadditive}, one has 
\[\|f+g\|_{\mathscr D_\nu^1}=\overline{I}_\nu(|f+g|)\leqslant\overline{I}_\nu(|f|+|g|)\leqslant\overline{I}_\nu(|f|)+\overline{I}_\nu(|g|)=\|f\|_{\mathscr D_\nu^1}+\|g\|_{\mathscr D_\nu^1},\]
where the first inequality comes from \eqref{Equ:dominantordre}, and the second inequality comes from \eqref{Pro:subsuperadditive}. Therefore $\|\ndot\|_{\mathscr D_\nu}$ is a seminorm on $\mathscr D_\nu^1(\Omega,\mathcal A,\nu)$.

Let $f$ be a $\nu$-negligible function. Then one has $|f|= 0$ $\nu$-almost everywhere. Hence one has $\|f\|_{\mathscr D_\nu^1}=\overline I_\nu(|f|)=0$. Conversely, if $f$ is a $\nu$-dominated function such that $\overline I_\nu(|f|)=0$, then one has $\underline I_\nu(|f|)=\overline I_{\nu}(|f|)=0$. By Proposition \ref{Pro:integrableupplower}, $|f|$ is $\nu$-indistinguishable with a $\nu$-integrable function $g$ of integral $0$. Moreover, since $|f|$ is non-negative, we obtain that the set $\{g<0\}$ is $\nu$-negligible. Therefore $g$ vanishes $\nu$-almost everywhere. Thus $f$ is $\nu$-negligible.
\end{proof}

\begin{prop}\label{Pro:Levi}
Let $\{f_n\}_{n\in\mathbb N}$ be an increasing sequence of non-negative functions on $\Omega$ and $f$ be the limit of $\{f_n\}_{n\in\mathbb N}$. Then one has
\[\lim_{n\rightarrow+\infty}\overline{I}_\nu(f_n)=\overline{I}_{\nu}(f).\]
\end{prop}
\begin{proof}Clearly one has $\overline{I}_\nu(f_n)\leqslant\overline{I}_\nu(f)$ for any $n\in\mathbb N$. Hence \[\lim_{n\rightarrow+\infty}\overline{I}_\nu(f_n)\leqslant\overline{I}_\nu(f).\]

If one of the functions $f_n$ is not dominated, then neither is $f$. Hence one has
\[\lim_{n\rightarrow+\infty}\overline{I}_\nu(f_n)=+\infty=\overline{I}_\nu(f).\] 
In the following, we assume that all the functions $f_n$ are dominated. Let $\epsilon>0$. For any $n\in\mathbb N$, let $g_n$ be an integrable function on $\Omega$ such that $f_n\leqslant g_n$ and that $\overline I_\nu(f_n)\geqslant I_\nu(g_n)-\epsilon$. Note that $\widetilde g_n:=\inf_{m\geqslant n}g_m$ is also an integrable function on $\Omega$ such that $f_n\leqslant \widetilde g_n$ and $\overline I_\nu(f_n)\geqslant I_\nu(\widetilde g_n)-\varepsilon$. Therefore, by replacing $g_n$ by $\widetilde{g}_n$, we may assume without loss of generality that the sequence $\{g_n\}_{n\in\mathbb N}$ is increasing. Let $g=\sup_{n\in\mathbb N}g_n$. By the monotone convergence theorem one has \[\overline{I}_\nu(g)=\lim_{n\rightarrow+\infty}I_\nu(g_n)\leqslant\lim_{n\rightarrow+\infty}\overline{I}_\nu(f_n)+\epsilon.\]
Moreover, since $g\geqslant f$, one has $\overline{I}_{\nu}(f)\leqslant\overline{I}_{\nu}(g)$. Therefore the proposition is proved.
\end{proof}

\begin{coro}\label{Cor:Levi}
Let $\{f_{n}\}_{n\in\mathbb N}$ be a sequence of non-negative functions on $\Omega$, and $f$ be the sum of the series $\sum_{n\in\mathbb N}f_n$. Then one has
\[\overline I_{\nu}(f)\leqslant\sum_{n\in\mathbb N}\overline{I}_\nu(f_n).\]
\end{coro}
\begin{proof}
For any $n\in\mathbb N$, let $g_n=\sum_{k=0}^nf_k$. The sequence $\{g_n\}_{n\in\mathbb N}$ is increasing, and converges to $f$. Therefore, by Proposition \ref{Pro:Levi}, one has 
\[\overline{I}_{\nu}(f)=\lim_{n\rightarrow+\infty}\overline{I}_\nu(g_n).\]
Moreover, by Proposition \ref{Pro:subsuperadditive}, for any $n\in\mathbb N$ one has
\[\overline{I}_\nu(g_n)\leqslant\sum_{k=0}^n\overline{I}_\nu(f_k).\]
Hence we obtain 
\[\overline I_{\nu}(f)\leqslant\sum_{n\in\mathbb N}\overline{I}_\nu(f_n).\]
\end{proof}

\begin{prop}\label{Pro:fatou}
Let $\{f_n\}_{n\in\mathbb N}$ be a sequence of non-negative functions on $\Omega$ and $f=\liminf_{n\rightarrow+\infty}f_n$. Then one has
\begin{equation}\label{Equ:Fatou}
\overline{I}_\nu(f)\leqslant\liminf_{n\rightarrow+\infty}\overline{I}_\nu(f_n).
\end{equation}
\end{prop}
\begin{proof}
For any $n\in\mathbb N$, let $g_n=\inf_{m\geqslant n}f_m$. Then the sequence $\{g_n\}_{n\in\mathbb N}$ is increasing and converges to $f$. By Proposition \ref{Pro:Levi}, one has
\[\overline{I}_\nu(f)=\lim_{n\rightarrow+\infty}\overline{I}_\nu(g_n)\leqslant\liminf_{n\rightarrow+\infty}\overline{I}_\nu(f_n),\]
where the inequality comes from the fact that $g_n\leqslant f_n$ for any $n\in\mathbb N$. The proposition is thus proved.
\end{proof}

\begin{prop}\label{Pro:completeness}
Let $D^1(\Omega,\mathcal A,\nu)$ be the quotient space of $\mathscr D^1(\Omega,\mathcal A,\nu)$ by the vector subspace of $\nu$-negligible functions. Then the seminorm $\|\ndot\|_{\mathscr D_\nu^1}$ on $\mathscr D^1(\Omega,\mathcal A,\nu)$ induced a norm $\|\ndot\|_{D_\nu^1}$ on $D^1(\Omega,\mathcal A,\nu)$ induced by the seminorm $\|\ndot\|_{\mathscr D_\nu^1}$, and the vector space $D^1(\Omega,\mathcal A,\nu)$ is complete with respect to this norm. 
\end{prop}
\begin{proof}
The first assertion is a direct consequence of Proposition \ref{Pro:seminormed}. In the following, we prove the second assertion.

Let $\{f_n\}_{n\in\mathbb N}$ be a Cauchy sequence in $\mathscr D^1(\Omega,\mathcal A,\nu)$. For any $\epsilon>0$ and any $m,n\in\mathbb N$, one has
$|f_n-f_m|\geqslant\epsilon\indic_{\{|f_n-f_m|>\epsilon\}}$, which implies that
\[\|f_n-f_m\|_{\mathscr D_\nu^1}=\overline I_\nu(|f_n-f_m|)\geqslant\epsilon\overline{I}_\nu(\indic_{\{|f_n-f_m|>\epsilon\}}).\]
Since $\{f_n\}_{n\in\mathbb N}$ is a Cauchy sequence, one has
\[\lim_{N\rightarrow+\infty}\sup_{\begin{subarray}{c}(n,m)\in\mathbb N^2\\ n\geqslant N,\,m\geqslant N\end{subarray}}\|f_n-f_m\|=0.\]
Therefore we can construct a subsequence $\{f_{n_k}\}_{k\geqslant 1}$ of $\{f_n\}_{n\in\mathbb N}$ such that
\[\forall\,k\in\mathbb N_{\geqslant 1},\quad \overline{I}_{\nu}(\indic_{\{|f_{n_k}-f_{n_{k+1}}|> 2^{-k}\}})<2^{-k}.\]
For any $m\in\mathbb N_{\geqslant 1}$, let $A_m=\bigcup_{k\geqslant m}\{|f_{n_k}-f_{n_{k+1}}|>2^{-k}\}$. Then the set
\[B:=\{\omega\in\Omega\,:\,\{f_{n_k}(\omega)\}_{k\geqslant 1}\text{ does not converge}\}\]
is contained in $\bigcap_{m\geqslant 1}A_m$. Moreover, for any $m\in\mathbb N_{\geqslant 1}$, by Corollary \ref{Cor:Levi} one has
\[\overline{I}_\nu(\indic_{A_m})\leqslant\sum_{k\geqslant m}\overline{I}_{\nu}(\indic_{\{|f_{n_k}-f_{n_{k+1}}|> 2^{-k}\}})\leqslant 2^{-m+1}.\]
Therefore one obtain $\overline I_\nu(\indic_B)=0$, which implies that $B$ is a $\nu$-negligible set. Thus we obtain that the sequence $\{f_{n_k}\}_{k\geqslant 1}$ converges $\nu$-almost everywhere to some function $f$ on $\Omega$. Note that by Proposition \ref{Pro:fatou} one has
\[\overline{I}_\nu(|f|)\leqslant\liminf_{k\rightarrow+\infty}\overline I_{\nu}(|f_{n_k}|).\]
Therefore $f$ is a dominated function. Finally, still by Proposition \ref{Pro:fatou}, for any $n\in\mathbb N$ one has
\[\overline I_\nu(|f_n-f|)\leqslant\liminf_{k\rightarrow+\infty}\overline I_{\nu}(|f_n-f_{n_k}|).\]
Hence one has 
\[\lim_{n\rightarrow+\infty}\overline I_\nu(|f_n-f|)=0.\]
The proposition is thus proved.
\end{proof}

\section{$L^1$ space}\phantomsection\label{Subsec: L1 space}
Let $\mathscr L^1(\Omega,\mathcal A,\nu)$ be the vector space of all real-valued $\nu$-integrable functions on the measurable space $(\Omega,\mathcal A)$. This vector space is equipped with the seminorm $\|\ndot\|_{\mathscr L^1_\nu}$ which sends a function $f\in\mathscr L^1(\Omega,\mathcal A,\nu)$ to
\[\|f\|_{\mathscr L^1_\nu}:=\int_{\Omega}|f(\omega)|\,\nu(\mathrm{d}\omega).\]
Note that the set of all functions $f\in\mathscr L^1(\Omega,\mathcal A,\nu)$ such that $\|f\|_{\mathscr L^1_\nu}=0$ forms a vector subspace of $\mathscr L^1(\Omega,\mathcal A,\nu)$. Such functions are said to be $\nu$-\emph{negligible}\index{negligible@negligible}. The quotient space of $\mathscr L^1(\Omega,\mathcal A,\nu)$ by the vector subspace of $\nu$-{negligible} functions is denoted by $L^1(\Omega,\mathcal A,\nu)$. The seminorm $\|\ndot\|_{\mathscr L^1_\nu}$ induces by quotient a norm on $L^1(\Omega,\mathcal A,\nu)$, which we denote by $\|\ndot\|_{L^1_\nu}$. Note that the vector space $L^1(\Omega,\mathcal A,\nu)$ is complete with respect to this norm, and the integration with respect to $\nu$ induces a continuous linear form on $L^1(\Omega,\mathcal A,\nu)$, which we denote by 
\[(\zeta\in L^1(\Omega,\mathcal A,\nu))\longmapsto\int_\Omega\zeta(\omega)\,\nu(\mathrm{d}\omega)\]
by abuse of notation.